\documentclass{report}[10pt, final]
\usepackage{latexsym, amscd, amsfonts, eucal, mathrsfs, amsmath, amssymb, amsthm, xypic, multind}

\def\Z{\mathbf{Z}}

\def\Q{\mathbf{Q}}

\DeclareMathOperator{\sNerve}{N}
\DeclareMathOperator{\Nerve}{N}
\DeclareMathOperator{\tNerve}{N}
\newcommand{\lNerve}{{\mathfrak F}}
\DeclareMathOperator{\sCoNerve}{\mathfrak{C}}

\DeclareMathOperator{\RFib}{RFib}
\renewcommand{\injlim}{\varinjlim}
\renewcommand{\hat}{\widehat}
\renewcommand{\projlim}{\varprojlim}
\renewcommand{\amalg}{\coprod}
\newcommand{\Rex}[1]{Rex(#1)}
\newcommand{\St}{St}
\DeclareMathOperator{\lex}{lex}
\DeclareMathOperator{\Mor}{Mor}

\newcommand{\Sh}{Sh}
\newcommand{\Un}{Un}
\DeclareMathOperator{\Gpd}{{\mathcal G}pd}
\DeclareMathOperator{\Post}{Post}
\newcommand{\h}[1]{\rm{h} \it{#1}}

\newcommand{\hn}[2]{\rm{h}_{#1} \! #2}

\newcommand{\Sum}{\sum}
 
\newcommand{\calKU}{{\mathcal K}{\mathcal U}}
\DeclareMathOperator{\Spf}{Spf}

\DeclareMathOperator{\rght}{R}
\DeclareMathOperator{\lft}{L}
\newcommand{\hyp}{\wedge}
\newcommand{\bfA}{{\mathbf A}}

\newcommand{\bfC}{{\mathbf C}}
\newcommand{\bfB}{{\mathbf B}}
\newcommand{\bfS}{{\mathbf S}}

\newcommand{\degree}{\circ}

\newcommand{\toposX}{\mathfrak X}
\newcommand{\Cech}{\v{C}ech\,}
\newcommand{\bigdot}{\bullet}

\newcommand{\join}{\star}
\DeclareMathOperator{\R}{\mathbf{R}}
\DeclareMathOperator{\calV}{\mathcal V}
\DeclareMathOperator{\rk}{rk}
\DeclareMathOperator{\mCech}{\check{C}}

\DeclareMathOperator{\Eff}{Eff}
\DeclareMathOperator{\Mod}{Mod}
\DeclareMathOperator{\Disc}{Disc}
\DeclareMathOperator{\Env}{Env}
\DeclareMathOperator{\Mult}{{\mathcal M}}
\DeclareMathOperator{\Multi}{{\mathcal M}}
\DeclareMathOperator{\Shv}{Shv}
\newcommand{\et}{\'{e}t}
\newcommand{\mathet}{\text{\et}}

\DeclareMathOperator{\LGeom}{\mathcal{L}{\mathcal T}op}
\DeclareMathOperator{\RGeom}{\mathcal{R}{\mathcal T}op}
\DeclareMathOperator{\Geo}{\mathcal{T}op}
\DeclareMathOperator{\GeoR}{R}

\DeclareMathOperator{\bHom}{{Map}}
\DeclareMathOperator{\Band}{Band}
\DeclareMathOperator{\Gerb}{Gerb}
\DeclareMathOperator{\Idem}{Idem}
\DeclareMathOperator{\Ret}{Ret}
\DeclareMathOperator{\calO}{\mathcal{O}}
\DeclareMathOperator{\calR}{\mathcal{R}}
\DeclareMathOperator{\Res}{\mathcal{R}es}
\DeclareMathOperator{\sSet}{\mathcal{S}et_{\Delta}}
\DeclareMathOperator{\mSet}{\mathcal{S}et_{\Delta}^{+}}
\DeclareMathOperator{\sCat}{\mathcal{C}at_{\Delta}}
\DeclareMathOperator{\SCat}{\mathcal{C}at_{\bfS}}

\DeclareMathOperator{\Inv}{Inv}

\DeclareMathOperator{\Grp}{\mathcal{G}pd}
\DeclareMathOperator{\Ab}{\mathcal{A}b}
\DeclareMathOperator{\EM}{\mathcal{E}\mathcal{M}}
\DeclareMathOperator{\Group}{\mathcal{G}rp}
\DeclareMathOperator{\tCat}{\mathcal{C}at_{\text{top}}}
\DeclareMathOperator{\Cat}{\mathcal{C}at}

\DeclareMathOperator{\Kan}{\mathcal{K}an}
\DeclareMathOperator{\calM}{\mathcal{M}}
\DeclareMathOperator{\calN}{\mathcal{N}}
\DeclareMathOperator{\calH}{\mathcal{H}}
\DeclareMathOperator{\Set}{\mathcal{S}et}

\DeclareMathOperator{\QCat}{\mathcal{Q}}
\DeclareMathOperator{\calQ}{\mathcal{Q}}
\DeclareMathOperator{\RPres}{\mathcal{P}r^{R}}
\DeclareMathOperator{\LPres}{\mathcal{P}r^{L}}
\newcommand{\LLPres}[1]{\mathcal{P}{\text r}^{\text L}_{#1}}
\newcommand{\RRPres}[1]{\mathcal{P}{\text r}^{\text R}_{#1}}
\DeclareMathOperator{\QC}{\QCat}
\DeclareMathOperator{\LFun}{Fun^{L}}
\DeclareMathOperator{\RFun}{Fun^{R}}
\DeclareMathOperator{\aHom}{Fun_{\mathcal{A}}}

\DeclareMathOperator{\CG}{\mathcal{C}\mathcal{G}}
\DeclareMathOperator{\Sing}{Sing}
\DeclareMathOperator{\Ob}{Ob}
\DeclareMathOperator{\Groth}{Groth}
\DeclareMathOperator{\cDelta}{{\bf \Delta}}
\DeclareMathOperator{\Sub}{Sub}

\DeclareMathOperator{\colim}{\varinjlim}

\DeclareMathOperator{\Top}{Top}

\DeclareMathOperator{\holim}{holim}
\DeclareMathOperator{\cosk}{cosk}
\DeclareMathOperator{\sk}{sk}
\DeclareMathOperator{\hocolim}{hocolim}
\DeclareMathOperator{\bd}{\partial}

\DeclareMathOperator{\calU}{\mathcal{U}}

\DeclareMathOperator{\calA}{\mathcal{A}}
\DeclareMathOperator{\calW}{\mathcal{W}}
\DeclareMathOperator{\calE}{\mathcal{E}}
\DeclareMathOperator{\calZ}{\mathcal{Z}}
\DeclareMathOperator{\calB}{\mathcal{B}}
\DeclareMathOperator{\calK}{\mathcal{K}}
 \DeclareMathOperator{\Pro}{Pro}

\DeclareMathOperator{\calF}{\mathcal{F}}
\DeclareMathOperator{\calG}{\mathcal{G}}
\DeclareMathOperator{\Hom}{Hom} 
\DeclareMathOperator{\HH}{H} 
\DeclareMathOperator{\id}{id} \DeclareMathOperator{\Fun}{Fun}
\DeclareMathOperator{\calC}{\mathcal{C}}
\DeclareMathOperator{\calI}{\mathcal{I}}
\DeclareMathOperator{\calJ}{\mathcal{J}}
\DeclareMathOperator{\SSet}{\mathcal{S}}

\DeclareMathOperator{\calX}{\mathcal{X}}

\DeclareMathOperator{\calY}{\mathcal{Y}}
\DeclareMathOperator{\op}{op}
\DeclareMathOperator{\calD}{\mathcal{D}}
\DeclareMathOperator{\Ind}{Ind} \DeclareMathOperator{\Acc}{Acc}
\DeclareMathOperator{\calP}{\mathcal{P}} 
\DeclareMathOperator{\PSL}{PSL}
\newcommand{\mset}[1]{(\Set^{+}_{\Delta})_{/#1}}

\newcommand{\push}[4]{\xymatrix{#1\ar[r]\ar[d] \ar@{}[rd]|(.7)*+{\lrcorner} & #2 \ar[d] \\ #3 \ar[r] & #4}}
\newcommand{\Push}[8]{\xymatrix{#1\ar[r]^{#5}\ar[d]_{#6} \ar@{}[rd]|(.7)*+{\lrcorner} & #2 \ar[d]^{#7} \\ #3 \ar[r]_{#8} & #4}}
\newcommand{\pull}[4]{\xymatrix{#1\ar[r]\ar[d] \ar@{}[rd]|(.3)*+{\ulcorner} & #2 \ar[d] \\ #3 \ar[r] & #4}}
\newcommand{\Pull}[8]{\xymatrix{#1\ar[r]^{#5}\ar[d]_{#6} \ar@{}[rd]|(.3)*+{\ulcorner} & #2 \ar[d]^{#7} \\ #3 \ar[r]_{#8} & #4}}
\newcommand{\adjoint}[2]{\xymatrix@1{#1\ar@<.4ex>[r] & #2 \ar@<.4ex>[l]}}
\newcommand{\Adjoint}[4]{\xymatrix@1{#2 \ar@<.4ex>[r]^{#1} & #3 \ar@<.4ex>[l]^{#4}}}

\newtheorem{theorem}{Theorem}[subsection]
\newtheorem{lemma}[theorem]{Lemma}

\newtheorem{proposition}[theorem]{Proposition}
\newtheorem{corollary}[theorem]{Corollary}
\newtheorem{fact}[theorem]{Fact}

\newtheorem*{proposition2}{Proposition}

\theoremstyle{definition}
\newtheorem{definition}[theorem]{Definition}
\newtheorem{variant}[theorem]{Variant}
\newtheorem{convention}[theorem]{Convention}
\newtheorem{example}[theorem]{Example}
\newtheorem{notation}[theorem]{Notation}
\newtheorem{counterexample}[theorem]{Counterexample}

\newtheorem{remark}[theorem]{Remark}
\newtheorem{warning}[theorem]{Warning}
\newtheorem*{remark2}{Remark}
\newtheorem*{warning2}{Warning}

\begin{document}

\title{Higher Topos Theory}
\author{Jacob Lurie}


\maketitle

\section*{Introduction}\label{intro}

\setcounter{theorem}{0}

Let $X$ be a nice topological space (for example, a CW complex). One goal of algebraic topology is to study the topology of $X$ by means of algebraic invariants, such as the singular cohomology groups
$\HH^{n}(X;G)$ of $X$ with coefficients in an abelian group $G$. These cohomology groups have proven to be an extremely useful tool, due largely to the fact that they enjoy excellent formal properties (which have been axiomatized by Eilenberg and Steenrod, see \cite{eilenbergsteenrod}), and the fact that they tend to be very computable. However, the usual definition of $\HH^{n}(X;G)$ in terms of singular $G$-valued cochains on $X$ is perhaps somewhat unenlightening. This raises the following question: can we understand the cohomology group $\HH^{n}(X;G)$ in more conceptual terms?

As a first step toward answering this question, we observe that $\HH^{n}(X;G)$ is a {\em representable} functor of $X$. That is, there exists an {\it Eilenberg-MacLane space} $K(G,n)$
and a universal cohomology class $\eta \in \HH^{n}( K(G,n); G)$ such that, for any nice topological space $X$, pullback of $\eta$ determines a bijection
$$ [X, K(G,n)] \rightarrow \HH^{n}(X;G).$$
Here $[X, K(G,n)]$ denotes the set of homotopy classes of maps from $X$ to $K(G,n)$.
The space $K(G,n)$ can be characterized up to homotopy equivalence by the above property, or by the the formula
$$ \pi_{k} K(G,n) \simeq \begin{cases} \ast & \text{if } k \neq n \\
G & \text{if } k = n. \end{cases}$$

In the case $n=1$, we can be more concrete. An Eilenberg MacLane space $K(G,1)$
is called a {\it classifying space} for $G$, and is typically denoted by $BG$. The universal cover
of $BG$ is a contractible space $EG$, which carries a free action of the group $G$ by covering transformations. We have a quotient map $\pi: EG \rightarrow BG$.
Each fiber of $\pi$ is a discrete topological space, on which the group $G$ acts simply transitively. We can summarize the situation by saying that $EG$ is a {\it $G$-torsor} over the classifying space $BG$.
For every continuous map $X \rightarrow BG$, the fiber product $\widetilde{X}: EG \times_{BG} X$ has the structure of a $G$-torsor on $X$: that is, it is a space endowed with a free action of $G$ and
a homeomorphism $\widetilde{X} / G \simeq X$. This construction determines a map
from $[X, BG]$ to the set of isomorphism classes of $G$-torsors on $X$. If $X$ is a well-behaved space (such as a CW complex), then this map is a bijection. We therefore have (at least) three different ways of thinking about a cohomology class $\eta \in \HH^{1}(X; G)$:
\begin{itemize}
\item[$(1)$] As a $G$-valued singular cocycle on $X$, which is well-defined up to coboundaries.
\item[$(2)$] As a continuous map $X \rightarrow BG$, which is well-defined up to homotopy.
\item[$(3)$] As a $G$-torsor on $X$, which is well-defined up to isomorphism.
\end{itemize}
These three points of view are equivalent if the space $X$ is sufficiently nice.
However, they are generally quite different from one another. The singular cohomology of a space $X$ is constructed using continuous maps from simplices $\Delta^k$ into $X$. If there are not many maps {\em into} $X$ (for example if every path in $X$ is constant), then we cannot expect singular cohomology to tell us very much about $X$. The second definition
uses maps from $X$ into the classifying space $BG$, which
(ultimately) relies on the existence of continuous real-valued
functions on $X$. If $X$ does not admit many real-valued
functions, then the set of homotopy classes $[X, BG]$ is also not a very useful invariant.
For such spaces, the third approach is the most powerful: there is a good theory of
$G$-torsors on an arbitrary topological space $X$. 

There is another reason for thinking about $\HH^{1}(X;G)$ in the language of $G$-torsors: it continues to make sense in situations where the traditional ideas of topology break down. If $\widetilde{X}$ is a $G$-torsor on a topological space $X$, then the projection map
$\widetilde{X} \rightarrow X$ is a local homeomorphism; we may therefore identify $\widetilde{X}$
with a sheaf of sets $\calF$ on $X$. The action of $G$ on $\widetilde{X}$ determines an action of
$G$ on $\calF$. The sheaf $\calF$ (with its $G$-action) and the space $\widetilde{X}$ (with its $G$-action) determine each other, up to canonical isomorphism. Consequently, we can formulate
the definition of a $G$-torsor in terms of the category $\Shv_{\Set}(X)$ of sheaves of sets
on $X$, without ever mentioning the topological space $X$ itself. The same definition makes
sense in any category which bears a sufficiently strong resemblance to the category of sheaves on a topological space: for example, in any {\em Grothendieck topos}. This observation allows us to construct a theory of torsors in a variety of nonstandard contexts, such as the \'{e}tale topology of algebraic varieties (see \cite{SGA}).

Describing the cohomology of $X$ in terms of the sheaf theory of $X$ has still another advantage, which comes into play even when the space $X$ is assumed to be a CW complex. For a general space
$X$, isomorphism classes of $G$-torsors on $X$ are classified not by the singular cohomology
$\HH^1_{\text{sing}}(X;G)$, but by the sheaf cohomology $\HH^{1}_{\text{sheaf}}(X; \calG)$ of
$X$ with coefficients in the constant sheaf $\calG$ associated to $G$. This sheaf cohomology is defined more generally for {\em any} sheaf of groups $\calG$ on $X$.
Moreover, we have a conceptual interpretation of $\HH^{1}_{\text{sheaf}}(X; \calG)$ in general: it classifies
$\calG$-torsors on $X$ (that is, sheaves $\calF$ on $X$ which carry an action of $\calG$ and locally admit a $\calG$-equivariant isomorphism $\calF \simeq \calG$) up to isomorphism. The general formalism of sheaf cohomology is extremely useful, even if we are interested only in the case where $X$ is a nice topological space: it includes, for example, the theory of cohomology with coefficients in a local system on $X$.

Let us now attempt to obtain a similar interpretation for cohomology classes $\eta \in \HH^{2}(X;G)$. 
What should play the role of a $G$-torsor in this case?
To answer this question, we return to the situation where $X$ is a CW complex, so that
$\eta$ can be identified with a continuous map $X \rightarrow K(G,2)$. 
We can think of
$K(G,2)$ as the classifying space of a group: not the discrete group $G$, but instead the classifying space $BG$ (which, if built in a sufficiently careful way, comes equipped with the structure of a topological abelian group). Namely, we can identify $K(G,2)$ with the quotient
$E/BG$, where $E$ is a contractible space with a free action of $BG$.
Any cohomology class $\eta \in \HH^{2}(X;G)$ determines a map $X \rightarrow K(G,2)$
(which is well-defined up to homotopy), and we can form the pullback $\widetilde{X} = E \times_{BG} X$. We now think of $\widetilde{X}$ as a torsor over $X$: not for the discrete group $G$, but instead for its classifying space $BG$.

To complete the analogy with our analysis in the case $n=1$, we would like to interpret
the fibration $\widetilde{X} \rightarrow X$ as defining some kind of sheaf $\calF$ on the space $X$.
This sheaf $\calF$ should have the property that for each $x \in X$, the stalk
$\calF_x$ can be identified with the fiber $\widetilde{X}_x \simeq BG$. Since the space $BG$ is not discrete (or homotopy equivalent to a discrete space), the situation cannot be adequately described in the usual language of set-valued sheaves. However, the classifying space $BG$ is 
{\em almost} discrete: since the homotopy groups $\pi_i BG$ vanish for $i > 1$, we
can recover $BG$ (up to homotopy equivalence) from its fundamental groupoid.
This suggests that we might try to think about $\calF$ as a ``groupoid-valued sheaf'' on $X$,
or a {\it stack}\index{gen}{stack} (in groupoids) on $X$.

\begin{remark2}\index{gen}{gerbe}
The condition that each stalk $\calF_x$ be equivalent to a classifying space $BG$
can be summarized by saying that $\calF$ is a {\it gerbe} on $X$: more precisely, it is a
gerbe banded by the constant sheaf $\calG$ associated to $G$.
We refer the reader to \cite{giraud} for an explanation of this terminology, and a proof that
such gerbes are indeed classified by the sheaf cohomology group $\HH^{2}_{\text{sheaf}}(X;\calG)$. 
\end{remark2}

For larger values of $n$, even the language of stacks is not sufficient to describe the nature of the sheaf $\calF$ associated to the fibration $\widetilde{X} \rightarrow X$. To address the situation,
Grothendieck proposed (in his infamous letter to Quillen; see \cite{pursuing}) that there
should be a theory of {\it $n$-stacks} on $X$, for every integer $n \geq 0$. Moreover, for every sheaf of abelian groups $\calG$ on $X$, the cohomology group
$\HH^{n+1}_{\text{sheaf}}(X;\calG)$ should have an interpreation as classifying a special type
of $n$-stack: namely, the class of $n$-gerbes banded by $\calG$ (for a discussion in the case $n=2$, we refer the reader to \cite{breen}; we will discuss the general case in \S \ref{chmdim}).
In the special case where the space $X$ is a point (and where we restrict our attention to
$n$-stacks in groupoids), the theory of $n$-stacks on $X$ should recover the classical homotopy theory of {\it $n$-types}: that is, CW complexes $Z$ such that the homotopy groups
$\pi_{i}(Z,z)$ vanish for $i > n$ (and every base point $z \in Z$). More generally, we should think of an $n$-stack (in groupoids) on a general space $X$ as a ``sheaf of $n$-types'' on $X$. 

When $n=0$, an $n$-stack on a topological space $X$ simply means a sheaf of sets on $X$.
The collection of all such sheaves can be organized into a category $\Shv_{\Set}(X)$, and this category is a prototypical example of a {\it Grothendieck topos}. The main goal of this book is to obtain an analogous understanding of the situation for $n > 0$. More precisely, we would like answers to the following questions:
\begin{itemize}
\item[$(Q1)$] Given a topological space $X$, what should we mean by a ``sheaf of $n$-types'' on $X$?

\item[$(Q2)$] Let $\Shv_{\leq n}(X)$ denote the collection of all sheaves of $n$-types on $X$.
What sort of a mathematical object is $\Shv_{\leq n}(X)$?

\item[$(Q3)$] What special features (if any) does $\Shv_{\leq n}(X)$ possess?
\end{itemize}

Our answers to questions $(Q2)$ and $(Q3)$ may be summarized as follows
(our answer to $(Q1)$ is more difficult to summarize, and we will avoid discussing it for the moment):

\begin{itemize}
\item[$(A2)$] The collection $\Shv_{\leq n}(X)$ has the structure of an
{\it $\infty$-category}.
\item[$(A3)$] The $\infty$-category $\Shv_{\leq n}(X)$ is an example of
an {\it $(n+1)$-topos}: that is, an $\infty$-category which satisfies
higher categorical analogues of Giraud's axioms for Grothendieck topoi
(see Theorem \ref{nchar}).\index{gen}{topos!$n$}\index{gen}{$n$-topos}
\end{itemize}

\begin{remark2}
Grothendieck's vision has been realized in various ways, thanks to the work of a number of mathematicians (most notably Brown, Joyal, and Jardine: see for example \cite{jardine}), and their work can also be used to provide answers to questions $(Q1)$ and $(Q2)$ (for more details, we refer the reader to \S \ref{hyperstacks}). Question $(Q3)$ has also been addressed (at least in limiting case $n = \infty$) by To\"{e}n and Vezzosi (see \cite{toen}) and in published work of Rezk. 
\end{remark2}

To provide more complete versions of the answers $(A2)$ and $(A3)$, we will need to develop the language of {\em higher category theory}. This is generally regarded as a technical and forbidding subject, but fortunately we will only need a small fragment of it. More precisely, we will need a theory of
{\it $(\infty,1)$-categories}: higher categories $\calC$ for which the $k$-morphisms of $\calC$ are required to be invertible for $k > 1$. In \S \ref{chap1}, we will present such a theory:
namely, one can define an {\it $\infty$-category}\index{gen}{$\infty$-category} to be a simplicial set
satisfying a weakened version of the Kan extension condition (see Definition \ref{qqcc}; simplicial sets satisfying this condition are also called {\it weak Kan complexes} or {\it quasicategories} in the literature). 
Our intention is that \S \ref{chap1} can be used as a short ``user's guide'' to $\infty$-categories:
it contains many of the basic definitions, and explains how many ideas from classical category theory can be extended to the $\infty$-categorical context. To simplify the exposition, we have deferred many proofs until later chapters, which contain a more thorough account of the theory.
The hope is that \S \ref{chap1} will be useful to readers who want to get the flavor of the subject, without becoming overwhelmed by technical details.

In \S \ref{chap2} we will shift our focus slightly: rather than study individual examples of $\infty$-categories, we consider {\em families} of $\infty$-categories $\{ \calC_{D} \}_{ D \in \calD}$ parametrized by the objects of another $\infty$-category $\calD$. We might expect such a family to be given by a map of $\infty$-categories $p: \calC \rightarrow \calD$: given such a map, we can then define each $\calC_{D}$ to be
the fiber product $\calC \times_{\calD} \{D\}$. This definition behaves poorly in general
(for example, the fibers $\calC_{D}$ need not be $\infty$-categories), but it behaves well if
we make suitable assumptions on the map $p$. Our goal in \S \ref{chap2} is to study some of these assumptions in detail, and show that they lead to a good {\em relative} version of higher category theory.

One motivation for the theory of $\infty$-categories is that it arises naturally in addressing questions like $(Q2)$ above. More precisely, given a collection of mathematical objects $\{ \calF_{\alpha} \}$ whose definition has a homotopy-theoretic flavor (like $n$-stacks on a topological space $X$), we can often organize the collection $\{ \calF_{\alpha} \}$ into an $\infty$-category
(in other words, we can find an $\infty$-category $\calC$ whose vertices correspond to the
objects $\calF_{\alpha}$). Another important example is provided by higher category theory itself: 
the collection of all $\infty$-categories can itself be organized into a (very large) $\infty$-category, which we will denote by $\Cat_{\infty}$. Our goal in \S \ref{chap4} is to study $\Cat_{\infty}$ and to show that it can be characterized by a universal property: namely, functors $\chi: \calD \rightarrow \Cat_{\infty}$ are classified, up to equivalence, by certain kinds of fibrations $\calC \rightarrow \calD$ (see Theorem \ref{straightthm} for a more precise statement). Roughly speaking, this correspondence assigns to a fibration
$\calC \rightarrow \calD$ the functor $\chi$ given by the formula $\chi(D) = \calC \times_{\calD} \{D\}$. 

Classically, category theory is a useful tool not so much because of the light it sheds on any particular
mathematical discipline, but instead because categories are so ubiquitous:
mathematical objects in many different settings (sets, groups, smooth manifolds, etcetera) can be organized into categories. Moreover, many elementary mathematical concepts can be described in purely categorical terms, and therefore make sense in each of these settings. For example, we can form products of sets, groups, and smooth manifolds: each of these notions can simply be described as a Cartesian product in the relevant category. Cartesian products are a special case of the more general notion of {\em limit}\index{gen}{limit}, which plays a central role in classical category theory. In \S \ref{chap3}, we will make a systematic study of limits (and the dual theory of colimits) in the $\infty$-categorical setting. We will also introduce the more general theory of {\em Kan extensions}, in both absolute and relative versions; this theory plays a key technical role throughout the later parts of the book.

In some sense, the material of \S \ref{chap1} through \S \ref{chap3} of this book should be regarded as purely formal. Our main results can be summarized as follows: there exists a reasonable theory of $\infty$-categories, and it behaves in more or less the same way as the theory of ordinary categories. Many of the ideas that we introduce are straightforward generalizations of their ordinary counterparts (though proofs in the $\infty$-categorical setting often require a bit of dexterity in manipulating simplicial sets), which will be familiar to mathematicians who are acquainted with ordinary category theory (as presented, for example, in \cite{maclane}). In \S \ref{chap5}, we introduce $\infty$-categorical analogues of more sophisticated concepts from classical category theory: presheaves, $\Pro$ and $\Ind$-categories, accessible and presentable categories, and localizations. The main theme is that most of the $\infty$-categories which appear ``in nature'' are large, but are determined by small subcategories. Taking careful advantage of this fact will allow us to deduce a number of pleasant results, such as an $\infty$-categorical version of the adjoint functor theorem (Corollary \ref{adjointfunctor}).

In \S \ref{chap6} we come to the heart of the book: the study of {\it $\infty$-topoi}, the $\infty$-categorical analogues of Grothendieck topoi. The theory of $\infty$-topoi is our answer
to the question $(Q3)$ in the limiting case $n = \infty$ (we will also study the analogous notion
for finite values of $n$). Our main result is an analogue of Giraud's theorem, which asserts the equivalence of ``extrinsic'' and ``intrinsic'' approaches to the subject (Theorem \ref{mainchar}). 
Roughly speaking, an $\infty$-topos is an $\infty$-category which ``looks like'' the $\infty$-category of all homotopy types. We will show that this intuition is justified in the sense that it is possible to reconstruct a large portion of classical homotopy theory inside an arbitrary $\infty$-topos. In other words, an $\infty$-topos
is a world in which one can ``do'' homotopy theory (much as an ordinary topos can be regarded as a world in which one can ``do'' other types of mathematics).




In \S \ref{chap7} we will discuss the relationship between our theory of $\infty$-topoi and ideas from classical topology. We will show that, if $X$ is a paracompact space, then the $\infty$-topos of ``sheaves of homotopy types'' on $X$ can be interpreted in terms of the classical homotopy theory of spaces {\em over} $X$. Another main theme is that various ideas from geometric topology (such as dimension theory and shape theory) can be described naturally in the language of $\infty$-topoi. We will also formulate and prove ``nonabelian'' generalizations of classical cohomological results, such as Grothendieck's vanishing theorem for the cohomology of Noetherian topological spaces, and the proper base change theorem.

\subsection*{Prerequisites and Suggested Reading}

We have made an effort to keep this book as self-contained as possible. The main prerequisite is familiarity with the classical homotopy theory of simplicial sets (good references include \cite{maysimp} and \cite{goerssjardine}; we have also provided a very brief review in \S \ref{simpset}). The remaining material that we need is either described in the appendix, or developed in the body of the text. However, our exposition of this background material is often somewhat terse, and the reader might benefit from consulting other treatments of the same ideas. Some suggestions for further reading are listed below.

\begin{warning2}
The list of references below is woefully incomplete. We have not attempted, either here or in the body of the text, to give a comprehensive survey of the literature on higher category theory. We have also not attempted to trace all of the ideas presented to their origins, or to present a detailed history of the subject. Many of the topics presented in this book have appeared elsewhere, or belong to the mathematical folklore; it should not be assumed that uncredited results are due to the author.
\end{warning2}

\begin{itemize}
\item {\bf Classical Category Theory:} Large portions of this book are devoted to providing $\infty$-categorical generalizations of the basic notions of category theory. A good reference for many of the concepts we use is MacLane's book \cite{maclane} (see also \cite{adamek} and \cite{makkai} for some of the more advanced material of \S \ref{chap5}).

\item {\bf Classical Topos Theory:} Our main goal in this book is to describe an $\infty$-categorical version of
topos theory. Familiarity with classical topos theory is not strictly necessary (we will define all of the relevant concepts as we need them), but will certainly be helpful. Good references include \cite{SGA} and \cite{where}.

\item {\bf Model Categories:} Quillen's theory of model categories provides a useful tool for studying
specific examples of $\infty$-categories, including the theory of $\infty$-categories itself.
We will summarize the theory of model categories in \S \ref{appmodelcat}; more complete
references include \cite{hovey}, \cite{hirschhorn}, and \cite{goerssjardine}.

\item {\bf Higher Category Theory:} There are many approaches to the theory of higher categories, some of which look quite different from the theory presented in this book. Several other possibilities are presented in the survey article \cite{leinster}. More detailed accounts can be found in \cite{leinster2}, 
\cite{simpson2}, and \cite{tamsamani}. 

In this book, we consider only {\it $(\infty,1)$-categories}: that is, higher categories in which all $k$-morphisms are assumed to be invertible for $k > 1$. There are a number of convenient ways to formalize this idea: via simplicial categories (see for example \cite{dwyerkan} and \cite{bergner}), via
Segal categories (\cite{simpson2}), via complete Segal spaces (\cite{completesegal}), or via
the theory of $\infty$-categories presented in this book (other references include \cite{joyalpub},
\cite{joyalnotpub}, \cite{nichols}, and \cite{quasicat}). The relationship between these various approaches is described in \cite{bergner2}, and an axiomatic framework which encompasses all of them is described in \cite{toenchar}.

\item {\bf Higher Topos Theory:} The idea of studying a topological space $X$ via the theory of
sheaves of $n$-types (or {\it $n$-stacks}) on $X$ goes back at least to Grothendieck
(\cite{pursuing}), and has been taken up a number of times in recent years. For small values of
$n$, we refer the reader to \cite{giraud}, \cite{street}, \cite{breen}, \cite{joyaltierney}, and \cite{polesello}. For the case $n=\infty$, we refer the reader to \cite{brown}, \cite{jardine}, \cite{hirschowitz}, and \cite{toen2}.
A very readable introduction to some of these ideas can be found in \cite{baezshul}.

Higher topos theory itself can be considered an abstraction of this idea: rather than studying
sheaves of $n$-types on a particular topological space $X$, we instead study general $n$-categories with the same formal properties. This idea has been implemented in the work of To\"{e}n and Vezzosi
(see \cite{toen} and \cite{toenvezz}), resulting in a theory which is essentially equivalent to
the one presented in this book. (A rather different variation on this idea in the case $n=2$ can
be also be found in \cite{ditopoi}.) The subject has also been greatly influenced by the unpublished ideas of Charles Rezk.
\end{itemize}

\section*{Terminology}

A few comments on some of the terminology which appears in this book:

\begin{itemize}
\item The word {\em topos} will always mean {\em Grothendieck} topos.\index{gen}{topos}

\item We let $\sSet$ denote the category of simplicial sets. If $J$ is a linearly ordered set, we
let $\Delta^{J}$ denote the simplicial set given by the nerve of $J$, so that
the collection of $n$-simplices of $\Delta^{J}$ can be identified with the collection
of all nondecreasing maps $\{0, \ldots, n \} \rightarrow J$. We will frequently apply this notation
when $J$ is a subset of $\{0, \ldots, n \}$; in this case, we can identify
$\Delta^{J}$ with a subsimplex of the standard $n$-simplex $\Delta^{n}$ (at
least if $J \neq \emptyset$; if $J = \emptyset$ then $\Delta^{J}$ is empty).

\item We will refer to a category $\calC$ as {\it accessible} or {\it presentable} if it is {\it locally accessible} or {\it locally presentable} in the terminology of \cite{makkai}.\index{gen}{accessible!category}\index{gen}{presentable!category}

\item Unless otherwise specified, the term {\it $\infty$-category} will be used to indicate a higher category in which all $n$-morphisms are invertible for $n > 1$.\index{gen}{$\infty$-category}

\item We will study higher category theory in Joyal's setting of {\it quasicategories}. However, we do not always follow Joyal's terminology.\index{gen}{quasicategory} In particular, we will use the term {\it $\infty$-category} to refer to what Joyal calls a {\it quasicategory} (which are, in turn, the same as the {\it weak Kan complex}\index{gen}{Kan complex!weak} of Boardman and Vogt); we will use the terms {\it inner fibration} and {\it inner anodyne map} where Joyal uses {\it mid-fibration} and {\it mid-anodyne map}.\index{gen}{weak Kan complex}

\item Let $n \geq 0$. We will say that a homotopy type $X$ (described by either a topological space or a Kan complex) is {\it $n$-truncated} if the homotopy groups $\pi_{i}(X,x)$ vanish for every point $x \in X$
and every $i > n$. By convention, we say that $X$ is $(-1)$-truncated if it is either empty or (weakly) contractible, and $(-2)$-truncated if $X$ is (weakly) contractible.

\item Let $n \geq 0$. We will say that a homotopy type $X$ (described either by a topological space or a Kan complex) is {\it $n$-connective} if $X$ is nonempty and the homotopy groups $\pi_{i}(X,x)$ vanish for every point $x \in X$ and every integer $i < n$. By convention, we will agree that every homotopy type $X$ is $(-1)$-connective. 

\item More generally, we will say that a map of homotopy types $f: X \rightarrow Y$ is
$n$-truncated ($n$-connective) if the homotopy fibers of $f$ are $n$-truncated ($n$-connective).
\end{itemize}

\begin{remark2}
For $n \geq 1$, a homotopy type $X$ is $n$-connective if and only if it is $(n-1)$-connected (in the usual terminology). In particular, $X$ is $1$-connective if and only if it is path connected.
\end{remark2}

\begin{warning2}
In this book, we will often be concerned with sheaves on a topological space $X$ (or some Grothendieck site) that take values in an $\infty$-category $\calC$. The most ``universal'' case is that in which $\calC$ is the $\infty$-category of $\SSet$ of spaces. Consequently, the term ``sheaf on $X$''
without any other qualifiers will typically refer to a sheaf of spaces on $X$, rather than a sheaf of sets on $X$. We will see that the collection of all $\SSet$-valued sheaves on $X$ can be organized into an
$\infty$-category, which we denote by $\Shv(X)$. In particular, $\Shv(X)$ will not denote the ordinary category of set-valued sheaves on $X$; if we need to consider this latter object, we will denote it
by $\Shv_{\Set}(X)$. 
\end{warning2}

\section*{Acknowledgements}

This book would never have come into existence without the advice and encouragement of many people. In particular, I would like to thank Vigleik Angeltveit, Vladimir Drinfeld, Matt Emerton, John Francis, Andre Henriques, James Parson, David Spivak, and James Wallbridge for many suggestions and corrections which have improved the readability of this book; Andre Joyal, who was kind enough to share with me a preliminary version of his work on the theory of quasi-categories; Charles Rezk, for explaining to me a very conceptual reformulation of the axioms for $\infty$-topoi (which we will describe in \S \ref{magnet}); Bertrand To\"{e}n and Gabriele Vezzosi, for many stimulating conversations about their work (which has considerable overlap with the material treated here); Mike Hopkins, for his advice and support throughout my time as a graduate student; Max Lieblich, for offering encouragement during early stages of this project; Josh Nichols-Barrer, for sharing with me some of his ideas about the foundations of higher category theory, and my editors Anna Pierrehumbert and Vickie Kearn at the Princeton Univesity Press, for helping to make this the best book that it can be. Finally, I would like to thank the American Institute of Mathematics for supporting me throughout the (seemingly endless) process of revising this book.

\tableofcontents

\chapter{An Overview of Higher Category Theory}\label{chap1}
\setcounter{theorem}{0}
\setcounter{subsection}{0}
This chapter is intended as a general introduction to higher category theory. We begin with what we feel is the most intuitive approach to the subject, based on {\it topological categories}. This approach is easy to understand, but difficult to work with when one wishes to perform even simple categorical constructions. As a remedy, we will introduce the more suitable formalism of {\it $\infty$-categories} (called
{\it weak Kan complexes} in \cite{quasicat} and {\it quasi-categories} in \cite{joyalpub}), which
provides a more convenient setting for adaptations of sophisticated category-theoretic ideas.
Our goal in \S \ref{highcat} is to introduce both approaches and to explain why they are equivalent to one another. The proof of this equivalence will rely on a crucial result (Theorem \ref{biggie}) which we will prove in \S \ref{valencequi}.

Our second objective in this chapter is to give the reader an idea of how to work with the formalism of $\infty$-categories. In \S \ref{langur} we will establish a vocabulary which includes $\infty$-categorical analogues (often direct generalizations) of most of the important concepts from ordinary category theory. To keep the exposition brisk, we will postpone the more difficult proofs until later chapters of this book. Our hope is that, after reading this chapter, a reader who does not wish to be burdened with the details will be able to understand (at least in outline) some of the more conceptual ideas described in \S \ref{chap5} and beyond.

\section{Foundations for Higher Category Theory}

\subsection{Goals and Obstacles}\label{highcat}

Recall that a {\it category}\index{gen}{category} $\calC$ consists of the following data:
\begin{itemize}
\item[$(1)$] A collection $\{ X, Y, Z, \ldots \}$ whose members are the {\it objects} of $\calC$. We typically
write $X \in \calC$ to indicate that $X$ is an object of $\calC$.
\item[$(2)$] For every pair of objects $X,Y \in \calC$, a set $\Hom_{\calC}(X,Y)$ of {\it morphisms} from
$X$ to $Y$. We will typically write $f: X \rightarrow Y$ to indicate that $f \in \Hom_{\calC}(X,Y)$, and say that $f$ {\it is a morphism from $X$ to $Y$}.
\item[$(3)$] For every object $X \in \calC$, an {\it identity morphism} $\id_{X} \in \Hom_{\calC}(X,X)$.

\item[$(4)$] For every triple of objects $X,Y, Z \in \calC$, a composition map
$$ \Hom_{\calC}(X,Y) \times \Hom_{\calC}(Y,Z) \rightarrow \Hom_{\calC}(X,Z).$$
Given morphisms $f: X \rightarrow Y$ and $g: Y \rightarrow Z$, we will usually denote the image of
the pair $(f,g)$ under the composition map by $gf$ or $g \circ f$.
\end{itemize}

These data are furthermore required to satisfy the following conditions, which guarantee that composition is unital and associative:

\begin{itemize}
\item[$(5)$] For every morphism $f: X \rightarrow Y$, we have
$\id_Y \circ f = f = f \circ \id_{X}$ in $\Hom_{\calC}(X,Y)$.
\item[$(6)$] For every triple of composable morphisms
$$ W \stackrel{f}{\rightarrow} X \stackrel{g}{\rightarrow} Y \stackrel{h}{\rightarrow} Z,$$
we have an equality $h \circ (g \circ f) = (h \circ g) \circ f$ in $\Hom_{\calC}(W,Z)$.
\end{itemize}

The theory of categories has proven to be a valuable organization tool in many areas of mathematics. Mathematical structures of virtually any type can be viewed as the objects of a suitable category $\calC$, where the morphisms in $\calC$ are given by structure-preserving maps. There is a veritable legion of examples of categories which fit this paradigm:
\begin{itemize}
\item The category $\Set$ whose objects are sets and whose morphisms are maps of sets.
\item The category $\Group$ whose objects are groups and whose morphisms are group homomorphisms.
\item The category $\Top$ whose objects are topological spaces and whose morphisms are continuous maps.
\item The category $\Cat$ whose objects are (small) categories and whose morphisms
are functors. (Recall that a functor $F$ from $\calC$ to $\calD$ is a map which assigns to each object
$C \in \calC$ another object $FC \in \calD$ and to each morphism $f: C \rightarrow C'$ in
$\calC$ a morphism $F(f): FC \rightarrow FC'$ in $\calD$, so that $F( \id_C) = \id_{FC}$ and
$F(g \circ f) = F(g) \circ F(f)$.)
\item \ldots
\end{itemize}

In general, the existence of a morphism $f: X \rightarrow Y$ in a category $\calC$ reflects some relationship that exists between the objects $X,Y \in \calC$. In some contexts, these relationships themselves become basic objects of study, and can themselves be fruitfully organized into categories:

\begin{example}\label{2cat1}
Let $\Group$ be the category whose objects are groups and whose morphisms are group homomorphisms. In the theory of groups, one is often concerned only with group homomorphisms
{\em up to conjugacy}. The relation of conjugacy can be encoded as follows: for every pair of
groups $G, H \in \Group$, there is a category $\bHom(G,H)$ whose objects are group homomorphisms
from $G$ to $H$ (that is, elements of $\Hom_{\Group}(G,H)$), where a morphism from $f: G \rightarrow H$ to $f': G \rightarrow H$ is an element $h \in H$ such that $h f(g) h^{-1} = f'(g)$ for all $g \in G$.
Note that two group homomorphisms $f,f': G \rightarrow H$ are conjugate if and only if they are isomorphic when viewed as objects of $\bHom(G,H)$.
\end{example}

\begin{example}\label{2cat2}
Let $X$ and $Y$ be topological spaces, and let $f_0, f_1: X \rightarrow Y$ be continuous maps.
Recall that a {\it homotopy} from $f_0$ to $f_1$ is a continuous map $f: X \times [0,1] \rightarrow Y$
such that $f | X \times \{0\}$ coincides with $f_0$ and $f | X \times \{1\}$ coincides with $f_1$.
In algebraic topology, one is often concerned not with the category $\Top$ of topological spaces,
but with its {\it homotopy category}: that is, the category obtained by identifying those pairs of
morphisms $f_0, f_1: X \rightarrow Y$ which are homotopic to one another. For many purposes, it is better to do something a little bit more sophisticated: namely, one can form a category
$\bHom(X,Y)$ whose objects are continuous maps $f: X \rightarrow Y$ and whose morphisms
are given by (homotopy classes of) homotopies.
\end{example}

\begin{example}\label{2cat3}
Given a pair of categories $\calC$ and $\calD$, the collection of all functors from
$\calC$ to $\calD$ is itself naturally organized into a category $\Fun( \calC, \calD)$, where
the morphisms are given by {\it natural transformations}. (Recall that, given a pair of functors
$F,G: \calC \rightarrow \calD$, a natural transformation $\alpha: F \rightarrow G$ is a collection
of morphisms $\{ \alpha_{C}: F(C) \rightarrow G(C) \}_{C \in \calC}$ which satisfy the following
condition: for every morphism $f: C \rightarrow C'$ in $\calC$, the diagram
$$ \xymatrix{ F(C) \ar[r]^{F(f)} \ar[d]^{\alpha_C} & F(C') \ar[d]^{\alpha_{C'} } \\
G(C) \ar[r]^{G(f)} & G(C') }$$
commutes in $\calD$.)
\end{example}

In each of these examples, the objects of interest can naturally be organized into what is
called a {\it $2$-category} (or {\it bicategory}): we have not only a collection of objects and a notion of morphisms between objects, but also a notion of morphisms between morphisms, which
are called {\it $2$-morphisms}. The vision of higher category theory is that there should exist a good notion of $n$-category for all $n \geq 0$, in which we have not only objects,
morphisms, and $2$-morphisms, but also $k$-morphisms for all $k
\leq n$. Finally, in some sort of limit we might hope to obtain a theory
of $\infty$-categories, where there are morphisms of all orders.\index{gen}{$2$-category}\index{gen}{bicategory}

\begin{example}\label{grape}
Let $X$ be a topological space, and $0 \leq n \leq \infty$. We
can extract an $n$-category $\pi_{\leq n} X$ (roughly) as follows.
The objects of $\pi_{\leq n} X$ are the points of $X$. If $x,y \in
X$, then the morphisms from $x$ to $y$ in $\pi_{\leq n} X$ are
given by continuous paths $[0,1] \rightarrow X$ starting at $x$
and ending at $y$. The $2$-morphisms are given by homotopies of
paths, the $3$-morphisms by homotopies between homotopies, and so
forth. Finally, if $n < \infty$, then two $n$-morphisms of
$\pi_{\leq n} X$ are considered to be the same if and only if they are homotopic to one another.\index{not}{pileqnX@$\pi_{\leq n}X$}\index{gen}{fundamental $n$-groupoid}

If $n = 0$, then $\pi_{\leq n} X$ can be identified with the set $\pi_0 X$ of path
components of $X$. If $n=1$, then our definition of $\pi_{\leq n}
X$ agrees with usual definition for the fundamental groupoid of
$X$. For this reason, $\pi_{\leq n} X$ is often called the {\it
fundamental $n$-groupoid of $X$}. It is an {\it $n$-groupoid} (rather than a mere $n$-category) 
because every $k$-morphism of $\pi_{\leq k} X$ has an inverse (at least ``up to
homotopy'').
\end{example}

There are many approaches to realizing the theory of higher categories.
We might begin by defining a $2$-category to be a
``category enriched over $\Cat$.'' In other words, we consider a
collection of objects together with a {\em category} of morphisms
$\Hom(A,B)$ for any two objects $A$ and $B$, and composition {\em
functors} $c_{ABC}: \Hom(A,B) \times \Hom(B,C) \rightarrow
\Hom(A,C)$ (to simplify the discussion, we will ignore identity
morphisms for a moment). These functors are required to satisfy an
associative law, which asserts that for any quadruple $(A,B,C,D)$
of objects, the diagram
$$ \xymatrix{ \Hom(A,B) \times \Hom(B,C) \times \Hom(C,D) \ar[d]
 \ar[r] & \Hom(A,C) \times \Hom(C,D) \ar[d] \\
 \Hom(A,B) \times \Hom(B,D) \ar[r] & \Hom(A,D) }$$
 commutes; in other words, one has an {\em equality} of functors
$$c_{ACD} \circ (c_{ABC} \times 1) = c_{ABD} \circ (1
\times c_{BCD})$$
from $\Hom(A,B) \times \Hom(B,C) \times
\Hom(C,D)$ to $\Hom(A,D)$. This leads to the definition of
a {\it strict $2$-category}.\index{gen}{$2$-category!strict}

At this point, we should object that the definition of a strict
$2$-category violates one of the basic philosophical principles of
category theory: one should never demand that two functors $F$ and $F'$ be
equal to one another. Instead one should postulate the existence of a natural
isomorphism between $F$ and $F'$. This means that the associative
law should not take the form of an equation, but of additional
structure: a collection of isomorphisms $\gamma_{ABCD}: c_{ACD} \circ
(c_{ABC} \times 1) \simeq c_{ABD} \circ (1 \times c_{BCD})$. We
should further demand that the isomorphisms $\gamma_{ABCD}$ be
functorial in the quadruple $(A,B,C,D)$ and satisfy
certain higher associativity conditions, which generalize the ``Pentagon axiom''
described in \S \ref{monoidaldef}. After formulating the
appropriate conditions, we arrive at the definition of a {\it weak
$2$-category}.\index{gen}{$2$-category!weak}

Let us contrast the notions of ``strict $2$-category'' and ``weak
$2$-category.'' The former is easier to define, since we do not
have to worry about the higher associativity conditions satisfied
by the transformations $\gamma_{ABCD}$. On the other hand, the
latter notion seems more natural if we take the philosophy of
category theory seriously. In this case, we happen to be lucky:
the notions of ``strict $2$-category'' and ``weak $2$-category''
turn out to be equivalent. More precisely, any weak $2$-category
is equivalent (in the relevant sense) to a strict $2$-category. The
choice of definition can therefore be regarded as a question of
aesthetics.

We now plunge onward to $3$-categories. Following the above
program, we might define a {\it strict $3$-category} to consist of a
collection of objects together with strict $2$-categories
$\Hom(A,B)$ for any pair of objects $A$ and $B$, together with a
strictly associative composition law. Alternatively, we could seek
a definition of {\it weak $3$-category} by allowing $\Hom(A,B)$ to
be only a weak $2$-category, requiring associativity only up to
natural $2$-isomorphisms, which satisfy higher associativity laws
up to natural $3$-isomorphisms, which in turn satisfy still higher
associativity laws of their own. Unfortunately, it turns out that
these notions are {\em not} equivalent.

Both of these approaches have serious drawbacks. The obvious
problem with weak $3$-categories is that an explicit definition is
extremely complicated (see \cite{tricat}, where a definition is given along these lines), to the point where it is
essentially unusable. On the other hand, strict $3$-categories
have the problem of not being the correct notion: most of the weak
$3$-categories which occur in nature are not equivalent to
strict $3$-categories. For example, the fundamental $3$-groupoid of
the $2$-sphere $S^2$ cannot be described using the language of
strict $3$-categories. The situation only gets worse (from either
point of view) as we pass to $4$-categories and beyond.

Fortunately, it turns out that major simplifications can be
introduced if we are willing to restrict our attention to
$\infty$-categories in which most of the higher morphisms are
invertible. Let us henceforth use the term {\it $(\infty,n)$-category}\index{gen}{$(\infty,n)$-category}
to refer to $\infty$-categories in which all $k$-morphisms are
invertible for $k > n$. The $\infty$-categories described in
Example \ref{grape} (when $n=\infty$) are all
$(\infty,0)$-categories. The converse, which asserts that every
$(\infty,0)$-category has the form $\pi_{\leq \infty} X$ for some
topological space $X$, is a generally accepted principle of higher
category theory. Moreover, the $\infty$-groupoid $\pi_{\leq \infty} X$ encodes the entire homotopy type of $X$. In other words, $(\infty,0)$-categories (that is,
$\infty$-categories in which {\em all} morphisms are invertible)
have been extensively studied from another point of view: they are
essentially the same thing as ``spaces'' in the sense of homotopy
theory, and there are many equivalent ways to describe them (for
example, we can use CW complexes or simplicial sets).

\begin{convention}\index{gen}{$\infty$-groupoid}\index{gen}{$\infty$-bicategory}\index{gen}{$\infty$-category}
We will often refer to $(\infty,0)$-categories as {\it $\infty$-groupoids} and $(\infty,2)$-categories as {\it $\infty$-bicategories}. Unless otherwise specified, the generic term {\it $\infty$-category} will mean $(\infty,1)$-category. 
\end{convention}

In this book, we will restrict our attention almost entirely to the theory of $\infty$-categories (in which we have only invertible $n$-morphisms for $n \geq 2$). Our reasons are threefold:

\begin{itemize}
\item[$(1)$] Allowing noninvertible $n$-morphisms for $n > 1$ introduces a number of additional complications to the theory, at both technical and conceptual levels. As we will see throughout this book, many ideas from category theory generalize to the $\infty$-categorical setting in a natural way. However, these generalizations are not so straightforward if we allow noninvertible $2$-morphisms. For example, one must distinguish between strict and lax fiber products, even in the setting of ``classical'' $2$-categories.

\item[$(2)$] For the applications studied in this book, we will not need to consider $(\infty,n)$-categories for $n > 2$. The case $n=2$ is of some relevance, because the collection of (small) $\infty$-categories can naturally be viewed as a (large) $\infty$-bicategory. However, we will generally be able to exploit this structure in an ad-hoc manner, without developing any general theory of $\infty$-bicategories.

\item[$(3)$] For $n > 1$, the theory of $(\infty,n)$-categories is most naturally viewed as a special case of {\em enriched} (higher) category theory. Roughly speaking, an $n$-category can be viewed as a category enriched over $(n-1)$-categories. As we explained above, this point of view is inadequate because it requires that composition satisfies an associative law up to equality, while in practice the associativity only holds up to isomorphism or some weaker notion of equivalence.
In other words, to obtain the correct definition we need to view the collection of $(n-1)$-categories
as an $n$-category, not as an ordinary category. Consequently, the naive approach is circular:  though it does lead to the correct theory of $n$-categories, we can only make sense of it if the theory of $n$-categories is already in place.

Thinking along similar lines, we can view an $(\infty,n)$-category as an $\infty$-category which is {\em enriched over $(\infty,n-1)$-categories}. The collection of $(\infty,n-1)$-categories is itself organized into an $(\infty,n)$-category $\Cat_{(\infty,n-1)}$, so at a first glance this definition suffers from the same problem of circularity. However, because the associativity properties of composition are required to hold up to {\em equivalence}, rather than up to arbitrary natural transformation, the noninvertible $k$-morphisms in $\Cat_{(\infty,n-1)}$ are irrelevant for $k > 1$. We may therefore view an $(\infty,n)$-category as a category enriched over $\Cat_{(\infty,n-1)}$, where the latter is regarded as an $\infty$-category by discarding noninvertible $k$-morphisms for $2 \leq k \leq n$.
In other words, the naive inductive definition of higher category theory is reasonable, {\em provided that we work in the $\infty$-categorical setting from the outset}. 
We refer the reader to \cite{tamsamani} for a definition of $n$-categories which follows this line of thought.

The theory of {\em enriched} $\infty$-categories is a useful and important one, but will not be treated in this book. Instead we refer the reader to \cite{DAG} for an introduction using the same language and formalism we employ here. 
\end{itemize}

Though we will not need a theory of $(\infty,n)$-categories for $n > 1$, the case $n=1$ is the main subject matter of this book. Fortunately, the above discussion suggests a definition. Namely, an $\infty$-category $\calC$ should be consist of a collection of objects, and an $\infty$-groupoid 
$\bHom_{\calC}(X,Y)$ for every pair of objects $X,Y \in \calC$. These $\infty$-groupoids can be identified with ``spaces'', and should be equipped with an associative composition law.
As before, we are faced with two choices as to how to make this
precise: do we require associativity on the nose, or only up to (coherent)
homotopy? Fortunately, the answer turns out to be irrelevant:
as in the theory of $2$-categories, any $\infty$-category
with a coherently associative multiplication can be replaced by an
equivalent $\infty$-category with a strictly associative
multiplication. We are led to the following:

\begin{definition}\label{ic}\index{gen}{topological category}\index{gen}{category!topological}
A {\it topological category} is a category which is enriched over
$\CG$, the category of compactly generated (and weakly Hausdorff) topological
spaces. The category of topological categories will be denoted by
$\tCat$.\index{not}{CG@$\CG$}\index{not}{Cattop@$\tCat$}
\end{definition}

More explicitly, a topological category $\calC$ consists of a
collection of objects, together with a (compactly generated)
topological space $\bHom_{\calC}(X,Y)$ for any pair of objects $X,Y
\in \calC$. These mapping spaces must be equipped with an
associative composition law, given by continuous maps
$$\bHom_{\calC}(X_0, X_1) \times \bHom_{\calC}(X_1, X_2) \times
\ldots \bHom_{\calC}(X_{n-1},X_n) \rightarrow
\bHom_{\calC}(X_0,X_n)$$ (defined for all $n \geq 0$). Here the
product is taken in the category of compactly generated
topological spaces.

\begin{remark}
The decision to work with compactly generated topological spaces,
rather than arbitrary spaces, is made in order to facilitate the comparison with more
combinatorial approaches to homotopy theory. This is a purely
technical point which the reader may safely ignore.
\end{remark}

It is possible to use Definition \ref{ic} as a foundation for higher category theory: that is, to {\em define} an $\infty$-category to be a topological category. However, this approach has a number of technical disadvantages. We will describe an alternative (though equivalent) formalism in the next section.

\subsection{$\infty$-Categories}\label{qqqc}

Of the numerous formalizations of higher category theory, Definition \ref{ic} is the quickest and most transparent. However, it is one of the most difficult to actually work with: many of the basic constructions of higher category theory give rise most naturally to $(\infty,1)$-categories for which the composition of
morphisms is associative only up to (coherent) homotopy (for several examples of this phenomenon, we refer the reader to \S \ref{langur}). In order to remain in the world of topological categories, it is necessary to combine these constructions with a ``straightening'' procedure which
produces a strictly associative composition law. Although it is always possible to do this
(see Theorem \ref{biggier}), it is much more technically convenient to work from the outset
within a more flexible theory of $(\infty,1)$-categories. Fortunately, there are many candidates
for such a theory, including the theory of Segal categories (\cite{simpson2}), the theory of complete Segal spaces (\cite{completesegal}), and the theory of model categories (\cite{hovey}, \cite{hirschhorn}).
To review all of these notions and their interrelationships would
involve too great a digression from the main purpose of this book.
However, the frequency with which we will encounter sophisticated
categorical constructions necessitates the use of {\em one} of
these more efficient approaches. We will employ the theory of {\it weak Kan complexes}, which goes back to Boardman-Vogt (\cite{quasicat}). These objects have subsequently been studied more extensively by Joyal (\cite{joyalpub} and \cite{joyalnotpub}), who calls them {\it quasicategories}. We will simply call them {\it $\infty$-categories}.\index{gen}{$\infty$-category}\index{gen}{quasicategory}

To get a feeling for what an $\infty$-category $\calC$ should be, it is useful to consider two extreme cases. If {\em every} morphism in $\calC$ is invertible, then $\calC$ is equivalent to the fundamental $\infty$-groupoid of a topological space $X$. In this case, higher category theory reduces to classical homotopy theory. On the other hand, if $\calC$ has no nontrivial $n$-morphisms for $n > 1$, then $\calC$ is equivalent to an ordinary category. A general formalism must capture the features of both of these examples. In other words, we need
a class of mathematical objects which can behave both like categories and like topological spaces. In \S \ref{highcat}, we achieved this by ``brute force'': namely, we directly amalgamated the theory of topological spaces and the theory of categories, by considering topological categories.
However, it is possible to approach the problem more directly using the theory of
{\em simplicial sets}. 
We will assume that the reader has some familiarity with the theory of simplicial sets; a brief review of this theory is included in \S \ref{simpset}, and a more extensive introduction can be found in \cite{goerssjardine}.

The theory of simplicial sets originated as a combinatorial approach to homotopy theory. Given any topological space $X$, one can associate a simplicial set $\Sing X$, whose $n$-simplices are precisely the continuous maps $| \Delta^n | \rightarrow X$, where $|\Delta^n| = \{ (x_0, \ldots, x_n) \in [0,1]^{n+1} | x_0 + \ldots + x_n =1 \}$ is the standard $n$-simplex. Moreover, the topological space $X$ is {\em determined}, up to weak homotopy equivalence, by $\Sing X$. More precisely, the singular complex functor $X \mapsto \Sing X$
admits a left adjoint, which carries every simplicial set $K$ to its {\it geometric realization} $|K|$. For every topological space $X$, the counit map \index{gen}{geometric realization!of simplicial sets}\index{not}{|K|@$|K|$}\index{not}{SingX@$\Sing X$} $| \Sing X | \rightarrow X$ is a weak homotopy equivalence. Consequently, if one is only interested in studying topological spaces up to weak homotopy equivalence, one might as well work simplicial sets instead.

If $X$ is a topological space, then the simplicial set $\Sing X$ has an important property, which is captured by the following definition:

\begin{definition}\label{strongkan}
Let $K$ be a simplicial set. We say that $K$ is a {\it Kan complex} if, for any $0 \leq i \leq n$ and any diagram of solid arrows\index{gen}{Kan complex}
$$ \xymatrix{ \Lambda^n_i \ar[r] \ar@{^{(}->}[d] & K \\
\Delta^n \ar@{-->}[ur],& }$$
there exists a dotted arrow as indicated rendering the diagram commutative. Here $\Lambda^n_i \subseteq \Delta^n$ denotes the $i$th horn, obtained from the simplex $\Delta^n$ by deleting the interior and the face opposite the $i$th vertex.
\end{definition}

The singular complex of any topological space $X$ is a Kan complex: this follows from the fact that the horn $| \Lambda^n_i |$ is a retract of the simplex $| \Delta^n |$ in the category of topological spaces. Conversely, any Kan complex $K$ ``behaves like'' a space: for example, there are simple combinatorial recipes for extracting homotopy groups from $K$ (which turn out be isomorphic to the homotopy groups of the topological space $|K|$). According to a theorem of Quillen (see \cite{goerssjardine} for a proof), the singular complex and geometric realization provide mutually inverse equivalences between the homotopy category of CW complexes and the homotopy category of Kan complexes.

The formalism of simplicial sets is also closely related to category theory.\index{gen}{nerve!of a category}\index{not}{NervecalC@$\Nerve(\calC)$} To any category $\calC$, we can associate a simplicial set $\Nerve(\calC)$, called the {\it nerve} of $\calC$. For each $n \geq 0$, we 
let $\Nerve(\calC)_{n} = \bHom_{\sSet}(\Delta^n, \Nerve(\calC))$ denote the set of all functors $[n]
\rightarrow \calC$. Here $[n]$ denotes the linearly ordered set $\{ 0, \ldots, n \}$, regarded as a category in the obvious way. More
concretely, $\Nerve(\calC)_n$ is the set of all composable
sequences of morphisms
$$ C_0 \stackrel{f_1}{\rightarrow} C_1 \ldots \stackrel{f_n}{\rightarrow} C_n$$ having length $n$.
In this description, the face map $d_i$ carries the above sequence
to
$$C_0 \stackrel{f_1}{\rightarrow} C_1 \ldots \stackrel{f_{i-1}}{\rightarrow}
C_{i-1} \stackrel{ f_{i+1} \circ f_i }{\rightarrow} C_{i+1}
\stackrel{f_{i+2}}{\rightarrow} \ldots
\stackrel{f_{n}}{\rightarrow} C_n$$ while the degeneracy $s_i$
carries it to $$C_0 \stackrel{f_1}{\rightarrow} C_1 \ldots
\stackrel{f_i}{\rightarrow} C_i \stackrel{\id_{C_i}}{\rightarrow}
C_i \stackrel{f_{i+1}}{\rightarrow} C_{i+1}
\stackrel{f_{i+2}}{\rightarrow} \ldots \stackrel{f_n}{\rightarrow}
C_n.$$

It is more or less clear from this description that the simplicial
set $\Nerve(\calC)$ is just a fancy way of encoding the structure
of $\calC$ as a category. More precisely, we note that the
category $\calC$ can be recovered (up to isomorphism) from its
nerve $\Nerve(\calC)$. The objects of $\calC$ are simply the {\it
vertices} of $\Nerve(\calC)$; that is, the elements of $\Nerve
(\calC)_0$. A morphism from $C_0$ to $C_1$ is given by an edge $\phi
\in \Nerve(\calC)_1$ with $d_1(\phi) = C_0$ and $d_0(\phi)=C_1$. The
identity morphism from an object $C$ to itself is given by the
degenerate simplex $s_0(C)$. Finally, given a diagram $C_0
\stackrel{\phi}{\rightarrow} C_1 \stackrel{\psi}{\rightarrow}
C_2$, the edge of $\Nerve(\calC)$ corresponding to $\psi \circ
\phi$ may be uniquely characterized by the fact that there exists
a $2$-simplex $\sigma \in \Nerve(\calC)_2$ with $d_2(\sigma)=\phi$,
$d_0(\sigma)=\psi$, and $d_1(\sigma)=\psi \circ \phi$.

It is not difficult to characterize those simplicial sets which arise as the nerve of a category:

\begin{proposition}\label{ruko}
Let $K$ be a simplicial set. Then the following conditions are equivalent:
\begin{itemize}
\item[$(1)$] There exists a small category $\calC$ and an isomorphism $K \simeq \Nerve(\calC)$.
\item[$(2)$] For each $0 < i < n$ and each diagram
$$ \xymatrix{ \Lambda^n_i \ar@{^{(}->}[d] \ar[r] & K \\
\Delta^n \ar@{-->}[ur], & \\}$$
there exists a {\em unique} dotted arrow rendering the diagram commutative.
\end{itemize}
\end{proposition}

\begin{proof}
We first show that $(1) \Rightarrow (2)$. Let $K$ be the nerve of a small category
$\calC$, and let $f_0: \Lambda^n_i \rightarrow K$ be a map of simplicial sets where
$0 < i < n$. We wish to show that $f_0$ can be extended uniquely to a map $f: \Delta^n \rightarrow K$.
For $0 \leq k \leq n$, let $X_i \in \calC$ be the image of the vertex $\{k\} \subseteq \Lambda^n_i$.
For $0 < k \leq n$, let $g_k: X_{k-1} \rightarrow X_{k}$ be the morphism in $\calC$ determined by the restriction $f_0 | \Delta^{ \{k-1,k\} }$. The composable chain of morphisms
$$ X_0 \stackrel{ g_1}{\rightarrow} X_1 \stackrel{g_2}{\rightarrow} \ldots \stackrel{g_n}{\rightarrow} X_n$$
determines an $n$-simplex $f: \Delta^n \rightarrow K$. We will show that $f$ is the desired solution to our extension problem (the uniqueness of this solution is evident: if $f': \Delta^n \rightarrow K$ is
any other map with $f' | \Lambda^n_i = f_0$, then $f'$ must correspond to the same chain of morphisms
in $\calC$, so that $f' = f$). It will suffice to prove the following, for every $0 \leq j \leq n$:
\begin{itemize}
\item[$(\ast_j)$] If $j \neq i$, then 
$$ f | \Delta^{ \{0, \ldots, j-1, j+1, \ldots, n \} } = f_0 | \Delta^{ \{0, \ldots, j-1, j+1, \ldots, n \} }.$$
\end{itemize}
To prove $(\ast_j)$, it will suffice to show that $f$ and $f_0$ have the same
restriction to $\Delta^{ \{k,k' \} }$, where $k$ and $k'$ adjacent
elements of the linearly ordered set $\{ 0, \ldots, j-1, j+1, \ldots, n \} \subseteq [n]$.
If $k$ and $k'$ are adjacent in $[n]$, then this follows by construction. In particular,
$(\ast)$ is automatically satisfied if $j=0$ or $j=n$.
Suppose instead that $k = j-1$ and $k' = j+1$, where $0 < j < n$. If $n = 2$, then $j=1=i$ and we obtain a contradiction. We may therefore assume that $n > 2$, so that either $j-1 > 0$ or
$j+1 < n$. Without loss of generality, $j-1 > 0$, so that $\Delta^{ \{ j-1, j+1 \} }
\subseteq \Delta^{ \{1, \ldots, n \} }$. The desired conclusion now follows from $(\ast_0)$.

We now prove the converse. Suppose that the simplicial set $K$ satisfies $(2)$; we claim that
$K$ is isomorphic to the nerve of a small category $\calC$. We construct the category
$\calC$ as follows:
\begin{itemize}
\item[$(i)$] The objects of $\calC$ are the vertices of $K$.

\item[$(ii)$] Given a pair of objects $x, y \in \calC$, we let
$\Hom_{\calC}(x,y)$ denote the collection of all edges
$e: \Delta^1 \rightarrow K$ such that $e|\{0\} =x$ and $e| \{1\} = y$.

\item[$(iii)$] Let $x$ be an object of $\calC$, Then the identity morphism
$\id_{x}$ is the edge of $K$ defined by the composition
$$ \Delta^1 \rightarrow \Delta^0 \stackrel{e}{\rightarrow} K.$$

\item[$(iv)$] Let $f: x \rightarrow y$ and $g: y \rightarrow z$ be morphisms
in $\calC$. Then $f$ and $g$ together determine a map
$\sigma_0: \Lambda^2_1 \rightarrow K$. In view of condition $(2)$, the map
$\sigma_0$ can be extended uniquely to a $2$-simplex $\sigma: \Delta^2 \rightarrow K$.
We define the composition $g \circ f$ to be the morphism from $x$ to $z$ in $\calC$ corresponding
to the edge given by the composition
$$ \Delta^1 \simeq \Delta^{ \{0,2\} } \subseteq \Delta^2 \stackrel{\sigma}{\rightarrow} K.$$
\end{itemize}

We first claim that $\calC$ is a category. To prove this, we must verify the following axioms:
\begin{itemize}
\item[$(a)$] For every object $y \in \calC$, the identity $\id_{y}$ is a unit with respect to composition.
In other words, for every morphism $f: x \rightarrow y$ in $\calC$ and every morphism
$g: y \rightarrow z$ in $\calC$, we have $\id_y \circ f = f$ and $g \circ \id_y = g$.
These equations are ``witnessed'' by the $2$-simplices $s_1(f), s_0(g) \in \Hom_{\sSet}( \Delta^2, K)$.

\item[$(b)$] Composition is associative. That is, for every sequence of composable morphisms
$$ w \stackrel{f}{\rightarrow} x \stackrel{g}{\rightarrow} y \stackrel{h}{\rightarrow} z,$$
we have $h \circ (g \circ f) = (h \circ g) \circ f$. To prove this, let us first choose
$2$-simplices $\sigma_{0 1 2}$ and $\sigma_{1 2 3}$ as indicated below:
$$ 
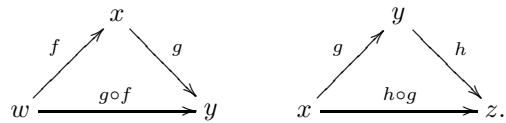
\xymatrix{ & x \ar[dr]^{g} & & & y \ar[dr]^{h} & \\
w \ar[ur]^{f} \ar[rr]^{g \circ f} & & y & x \ar[ur]^{g} \ar[rr]^{h \circ g} & & z. }$$
Now choose a $2$-simplex $\sigma_{0 2 3}$ corresponding to a diagram
$$ \xymatrix{ & y \ar[dr]^{h} & \\
w \ar[ur]^{ g \circ f} \ar[rr]^{ h \circ (g \circ f)} &  & z. }$$
These three $2$-simplices together define a map $\tau_0: \Lambda^3_2 \rightarrow K$.
Since $K$ satisfies condition $(2)$, we can extend $\tau_0$ to a $3$-simplex
$\tau: \Delta^3 \rightarrow K$. The composition
$$ \Delta^2 \simeq \Delta^{ \{0,1,3\}} \subseteq \Delta^3 \stackrel{\tau}{\rightarrow} K$$
corresponds to the diagram
$$ \xymatrix{ & x \ar[dr]^{h \circ g} & \\
w \ar[ur]^{f} \ar[rr]^{ h \circ (g \circ f) } & & z, }$$
which ``witnesses'' the associativity axiom $h \circ (g \circ f) = (h \circ g) \circ f$.
\end{itemize}

It follows that $\calC$ is a well-defined category. By construction, we have a canonical
map of simplicial sets $\phi: K \rightarrow \Nerve \calC$. To complete the proof, it will suffice
to show that $\phi$ is an isomorphism. We will prove, by induction on
$n \geq 0$, that $\phi$ induces a bijection $\Hom_{ \sSet}( \Delta^n, K) \rightarrow
\Hom_{\sSet}( \Delta^n, \Nerve \calC)$. For $n=0$ and $n=1$, this is obvious from
the construction. Assume therefore that $n \geq 2$, and choose an integer $i$ such
that $0 < i < n$. We have a commutative diagram
$$ \xymatrix{ \Hom_{\sSet}( \Delta^n, K ) \ar[r] \ar[d] & \Hom_{\sSet}( \Delta^n, \Nerve \calC) \ar[d] \\
\Hom_{ \sSet}( \Lambda^n_i, K) \ar[r] & \Hom_{\sSet}( \Lambda^n_i, \Nerve \calC ). }$$
Since $K$ and $\Nerve \calC$ both satisfy $(2)$ (for $\Nerve \calC$, this follows
from the first part of the proof), the vertical maps are bijective. It will therefore suffice to show that the lower horizontal map is bijective, which follows from the inductive hypothesis.
\end{proof}

We note that condition $(2)$ of Proposition \ref{ruko} is very similar to Definition \ref{strongkan}. However, it is different in two important respects. First, it requires the extension condition only for {\em inner} horns $\Lambda^n_i$ with $0 < i < n$. Second, the asserted condition is stronger in this case: not only does any map $\Lambda^n_i \rightarrow K$ extend to the simplex $\Delta^n$, but the extension is unique.\index{gen}{horn!inner}\index{gen}{inner horn}

\begin{remark}\label{nottes}
It is easy to see that it is not reasonable to expect condition $(2)$ of Proposition \ref{ruko} to hold for ``outer'' horns $\Lambda^n_i$, $i \in \{0,n\}$. Consider, for example, the case where $i=n=2$, and where $K$ is the nerve of a category $\calC$. Giving a map $\Lambda^2_2 \rightarrow K$
corresponds to supplying the solid arrows in the diagram
$$ \xymatrix{ & C_1 \ar[dr] & \\
C_0 \ar[rr] \ar@{-->}[ur] & & C_2,} $$
and the extension condition would amount to the assertion that one could always find a dotted arrow rendering the diagram commutative. This is true in general only when the category $\calC$ is a {\em groupoid}.
\end{remark}

We now see that the notion of a simplicial set is a flexible one: a simplicial set $K$ can be a good model for an $\infty$-groupoid (if $K$ is a Kan complex), or for an ordinary category (if it satisfies the hypotheses of Proposition \ref{ruko}). Based on these observations, we might expect that some more general class of simplicial sets could serve as models for $\infty$-categories in general.

Consider first an arbitrary simplicial set $K$. We can try to envision $K$ as a generalized category, whose objects are the vertices of $K$ (that is, the elements of $K_0$), and whose morphisms are the edges of $K$ (that is, the elements of $K_1$). A $2$-simplex
$\sigma: \Delta^2 \rightarrow K$ should be thought of as a diagram
$$ \xymatrix{ & Y \ar[dr]^{\psi} & \\
X \ar[ur]^{\phi} \ar[rr]^{\theta} & & Z }$$
together with an identification (or homotopy) between $\theta$ and $\psi \circ \phi$
which renders the diagram ``commutative''. (Note that, in higher category theory, 
this is not merely a condition: the homotopy $\theta \simeq \psi \circ \phi$ is an additional datum.)
Simplices of larger dimension may be thought of as verifying the
commutativity of certain higher-dimensional diagrams.

Unfortunately, for a general simplicial set $K$, the analogy
outlined above is not very strong. The essence of the problem is that, though we may refer
to the $1$-simplices of $K$ as ``morphisms'', there is in general no way to compose them.
Taking our cue from the example of $\Nerve(\calC)$,
we might say that a morphism $\theta:X \rightarrow Z$ is a composition of morphisms
$\phi: X \rightarrow Y$ and $\psi: Y \rightarrow Z$ if there exists a $2$-simplex
$\sigma: \Delta^2 \rightarrow K$ as in the diagram indicated above. 
We must now consider two potential difficulties: the $2$-simplex $\sigma$ may not exist, and
if it does it exist it may not be unique, so that we have more than one choice for the composition $\theta$.

The existence of $\sigma$ can be formulated as an extension condition on the simplicial set $K$.
We note that a composable pair of morphisms $(\psi, \phi)$ determines a map
of simplicial sets
$\Lambda^2_1 \rightarrow K$. Thus, the assertion that $\sigma$ can
always be found may be formulated as a extension property: any
map of simplicial sets $\Lambda^2_1 \rightarrow K$ can be extended to $\Delta^2$, as indicated in the following diagram:
$$ \xymatrix{ \Lambda^2_1\ar[r]\ar@{^{(}->}[d] & K \\
\Delta^2 \ar@{-->}[ur]}$$

The uniqueness of $\theta$ is another matter. It turns out to be
unnecessary (and unnatural) to require that $\theta$ be uniquely
determined. To understand this point, let us return to 
the example of the fundamental groupoid of a topological space
$X$. This is a category whose objects are the points $x \in X$.
The morphisms between a point $x \in X$ and a point $y \in X$ are given by
continuous paths $p: [0,1] \rightarrow X$ such that $p(0)=x$ and
$p(1)=y$. Two such paths are considered to be equivalent if there
is a homotopy between them. Composition in the fundamental
groupoid is given by concatenation of paths. Given paths $p,q:
[0,1] \rightarrow X$ with $p(0)=x$, $p(1)=q(0)=y$, and $q(1)=z$,
the composite of $p$ and $q$ should be a path joining $x$ to $z$.
There are many ways of obtaining such a path from $p$ and $q$. One
of the simplest is to define
$$r(t) = \begin{cases} p(2t) & \text{if } 0 \leq t \leq \frac{1}{2} \\
q(2t-1) & \text{if } \frac{1}{2} \leq t \leq 1. \end{cases}$$
However, we could just as well use the formula
$$r'(t) = \begin{cases} p(3t) & \text{if } 0 \leq t \leq \frac{1}{3} \\
q(\frac{3t-1}{2}) & \text{if } \frac{1}{3} \leq t \leq 1
\end{cases}$$
to define the composite path. Because the paths $r$ and $r'$ are homotopic to one another, it does not matter which one we choose.

The situation becomes more complicated if we try to think
$2$-categorically. We can capture more information about the space
$X$ by considering its {\it fundamental $2$-groupoid}. This is a
$2$-category whose objects are the points of $X$, whose morphisms
are paths between points, and whose $2$-morphisms are given by
homotopies between paths (which are themselves considered modulo
homotopy). In order to have composition of morphisms unambiguously defined, we would have to choose some formula once and for all.
Moreover, there is no particularly compelling choice; for example,
neither of the formulas written above leads to a strictly associative composition law.

The lesson to learn from this is that in higher-categorical
situations, we should not necessarily ask for a uniquely
determined composition of two morphisms. In the fundamental groupoid
example, there are many choices for a composite path but all of
them are homotopic to one another. Moreover, in keeping with the
philosophy of higher category theory, {\em any} path which is
homotopic to the composite should be just as good as the
composite itself. From this point of view, it is perhaps more natural to
view composition as a relation than as a function, and this is
very efficiently encoded in the formalism of simplicial sets: a
$2$-simplex $\sigma: \Delta^2 \rightarrow K$ should be viewed as
``evidence'' that $d_0(\sigma) \circ d_2(\sigma)$ is homotopic to $d_1(\sigma)$.

Exactly what conditions on a simplicial set $K$ will guarantee that it behaves like a higher category? Based on the above argument, it seems reasonable to require that $K$ satisfy an extension condition with respect to certain horn inclusions $\Lambda^n_i$, as in Definition \ref{strongkan}. However, as we observed in Remark \ref{nottes}, this is reasonable only for the inner horns where $0 < i < n$, which appear in the statement of Proposition \ref{ruko}. 

\begin{definition}\label{qqcc}\index{gen}{$\infty$-category}
An {\it $\infty$-category} is a simplicial set $K$ which has the
following property: for any $0 < i < n$, any map $f_0: \Lambda^n_i
\rightarrow K$ admits an extension $f: \Delta^n \rightarrow K$.
\end{definition}

Definition \ref{qqcc} was first formulated by Boardman and Vogt (\cite{quasicat}). They referred to $\infty$-catgories as {\it weak
Kan complexes},\index{gen}{Kan complex!weak} motivated by the obvious analogy with Definition \ref{strongkan}. Our terminology places more emphasis on the analogy with the characterization of ordinary categories given in Proposition \ref{ruko}: we require the same extension conditions, but drop the uniqueness assumption.

\begin{example}
Any Kan complex is an $\infty$-category. In particular, if $X$ is a topological space, then we may view its singular complex $\Sing X$ as an $\infty$-category: this is one way of defining the fundamental $\infty$-groupoid $\pi_{\leq \infty} X$ of $X$, introduced informally in Example \ref{grape}.
\end{example}

\begin{example}
The nerve of any category is an $\infty$-category. We will occasionally abuse terminology by
identifying a category $\calC$ with its nerve $\Nerve(\calC)$; by means of this identification, we may view ordinary category theory as a special case of the study of $\infty$-categories.
\end{example}

The weak Kan condition of Definition \ref{qqcc} leads to a very elegant and powerful version of higher category theory. This theory has been developed by Joyal in the papers \cite{joyalpub} and
\cite{joyalnotpub} (where simplicial sets satisfying the condition of Definition \ref{qqcc} are called
{\it quasi-categories}), and will be used throughout this book.

\begin{notation}
Depending on the context, we will use two different notations in
connection with simplicial sets. When emphasizing their
role as $\infty$-categories, we will often denote
them by
calligraphic letters such as $\calC$, $\calD$, and so forth. When
casting simplicial sets in their different (though related) role
of representing homotopy types, we will employ capital Roman
letters. To avoid confusion, we will also employ the latter notation
when we wish to contrast the theory of $\infty$-categories with some
other other approach to higher category theory, such as the theory
of topological categories.
\end{notation}

\subsection{Equivalences of Topological Categories}\label{stronghcat}

We have now introduced two approaches to higher category theory: one based on topological categories, and one based on simplicial sets. These two approaches turn out to be equivalent to one another. However, the equivalence itself needs to be understood in a higher-categorical sense. We take our cue from classical homotopy theory, in which we can take the basic objects to be either topological spaces or simplicial sets. It is not true that every Kan complex is isomorphic to the singular complex of a topological space, or that every CW complex is homeomorphic to the geometric realization of a simplicial set. However, both of these statements become true if we replace the words ``isomorphic to'' by ``homotopy equivalent to''. We would like to formulate a similar statement regarding our approaches to higher category theory. The first step is to find a concept which replaces ``homotopy equivalence''. If $F: \calC \rightarrow \calD$ is a functor between topological categories, under what circumstances should we regard $F$ as an ``equivalence'' (so that $\calC$ and $\calD$ really represent the same higher category)? 

The most naive answer is that $F$ should be regarded as an equivalence if it is an isomorphism of topological categories. This means that $F$ induces a bijection between the objects of $\calC$ and the objects of $\calD$, and a homeomorphism $\bHom_{\calC}(X,Y) \rightarrow \bHom_{\calD}(F(X),F(Y))$ for every pair of objects $X,Y \in \calC$. However, it is immediately obvious that this condition is far too strong; for example, in the case where $\calC$ and $\calD$ are ordinary categories (which we may view also as topological categories, where all morphism sets are endowed with the discrete topology), we recover the notion of an isomorphism between categories. This notion does not play an important role in category theory. One rarely asks whether or not two categories are isomorphic; instead, one asks whether or not they are equivalent. This suggests the following definition:

\begin{definition}\index{gen}{strong equivalence}
A functor $F: \calC \rightarrow \calD$ between topological categories is a  {\it strong equivalence} if it is an equivalence in the sense of enriched category theory. In other words, $F$ is a strong equivalence if it induces homeomorphisms $\bHom_{\calC}(X,Y) \rightarrow \bHom_{\calD}(F(X),F(Y))$ for every pair of objects $X,Y \in \calC$, and every object of $\calD$ is isomorphic (in $\calD$) to $F(X)$ for some $X \in \calC$.
\end{definition}

The notion of strong equivalence between topological categories has the virtue that, when restricted to ordinary categories, it reduces to the usual notion of equivalence. However, it is still not the right definition: for a pair of objects $X$ and $Y$ of a higher category $\calC$, the morphism space $\bHom_{\calC}(X,Y)$ should itself only be well-defined up to homotopy equivalence. 

\begin{definition}\label{vergen}\index{gen}{homotopy category!of a topological category}
Let $\calC$ be a topological category. The {\it homotopy category} $\h{\calC}$
is defined as follows:\index{not}{hcalC@$\h{\calC}$}
\begin{itemize}
\item The objects of $\h{\calC}$ are the objects of $\calC$.
\item If $X,Y \in \calC$, then we define $\Hom_{\h{\calC}}(X,Y)= \pi_0 \bHom_{\calC}(X,Y)$.
\item Composition of morphisms in $\h{\calC}$ is induced from the composition of morphisms
in $\calC$ by applying the functor $\pi_0$.
\end{itemize}
\end{definition}

\begin{example}\index{gen}{homotopy category!of spaces}
Let $\calC$ be the topological category whose objects are CW-complexes, where
$\bHom_{\calC}(X,Y)$ is the set of continuous maps from $X$ to $Y$, equipped with the
(compactly generated version of the) compact-open topology. We will denote the homotopy category of $\calC$ by $\calH$, and refer to $\calH$ as the {\it homotopy category of spaces}.\index{not}{Hcal@$\calH$}
\end{example}

There is a second construction of the homotopy category $\calH$, which will play an important role in what follows. First, we must recall a bit of terminology from classical homotopy theory.

\begin{definition}\index{gen}{weak homotopy equivalence!of topological spaces}
A map $f: X \rightarrow Y$ between topological spaces is said to be a {\it weak homotopy equivalence}
if it induces a bijection $\pi_0 X \rightarrow \pi_0 Y$, and if 
for every point $x \in X$ and every $i \geq 1$, the induced map of homotopy groups
$$ \pi_{i}(X,x) \rightarrow \pi_i(Y,f(x))$$ is an isomorphism.
\end{definition}

Given a space $X \in \CG$, classical homotopy theory ensures the existence of a CW-complex
$X'$ equipped with a weak homotopy equivalence $\phi: X' \rightarrow X$. Of course,
$X'$ is not uniquely determined; however, it is unique up to canonical homotopy equivalence,
so that the assignment
$$ X \mapsto [X] = X'$$
determines a functor $\theta: \CG \rightarrow \calH$. By construction, $\theta$ carries
weak homotopy equivalences in $\CG$ to isomorphisms in $\calH$. In fact, $\theta$
is universal with respect to this property. In other words, we may describe $\calH$
as the category obtained from $\CG$ by formally inverting all weak homotopy equivalences.
This is one version of Whitehead's theorem, which is usually stated as follows: every weak homotopy equivalence between CW complexes admits a homotopy inverse.\index{gen}{Whitehead's theorem}

We can now improve upon Definition \ref{vergen} slightly. We first observe that
the functor $\theta: \CG \rightarrow \calH$ preserves products. Consequently, we can apply the construction of Remark \ref{laxcon} to convert any topological category $\calC$ into a category enriched over $\calH$. We will denote this $\calH$-enriched category by $\h{\calC}$, and refer to it as the {\it homotopy category} of $\calC$. More concretely, the homotopy category
$\h{\calC}$ may be described as follows:\index{gen}{homotopy category!enriched over $\calH$}
\begin{itemize}
\item[$(1)$] The objects of $\h{\calC}$ are the objects of $\calC$.
\item[$(2)$] For $X,Y \in \calC$, we have
$$ \bHom_{ \h{\calC} }(X,Y) = [ \bHom_{\calC}(X,Y) ].$$
\item[$(3)$] The composition law on $\h{\calC}$ is obtained from the composition law on
$\calC$ by applying the functor $\theta: \CG \rightarrow \calH$.
\end{itemize}

\begin{remark}
If $\calC$ is a topological category, we have now defined $\h{\calC}$ in two different ways: first as an ordinary category, and then as a category enriched over $\calH$. These two definitions are compatible with one another, in the sense that $\h{\calC}$ (as an ordinary category) is the
underlying category of $\h{\calC}$ (as an $\calH$-enriched category). This follows immediately
from the observation that for every topological space $X$, there is a canonical bijection
$\pi_0 X \simeq \bHom_{\calH}( \ast, [X] ).$
\end{remark}

If $\calC$ is a topological category, we may imagine that $\h{\calC}$ is the object which is obtained by forgetting the topological morphism spaces of $\calC$ and remembering only their (weak) homotopy types. The following definition codifies the idea that these homotopy types should be ``all that really matter''.

\begin{definition}\label{defequiv}\index{gen}{equivalence!of topological categories}
Let $F: \calC \rightarrow \calD$ be a functor between topological categories. We will say that $F$ is a {\it weak equivalence}, or simply an {\it equivalence}, if the induced functor
$\h{ \calC} \rightarrow \h{ \calD}$ is an equivalence of $\calH$-enriched categories.
\end{definition}

More concretely, a functor $F$ is an equivalence if and only if:

\begin{itemize}
\item For every pair of objects $X,Y \in \calC$, the induced map
$$ \bHom_{\calC}(X,Y) \rightarrow \bHom_{\calD}(F(X),F(Y))$$ is a weak homotopy equivalence of topological spaces.

\item Every object of $\calD$ is isomorphic in $\h{ \calD}$ to $F(X)$, for some $X \in \calC$.
\end{itemize}

\begin{remark}\index{gen}{equivalence!in a topological category}
A morphism $f: X \rightarrow Y$ in $\calD$ is said to be an {\it equivalence} if the induced morphism in $\h{ \calD}$ is an isomorphism. In general, this is much weaker than the condition that $f$ be an isomorphism in $\calD$; see Proposition \ref{rooot}.
\end{remark}

It is Definition \ref{defequiv} which gives the correct notion of equivalence between topological categories (at least, when one is using them to describe higher category theory). We will agree that all relevant properties of topological categories are invariant under this notion of equivalence. We say that two topological categories are {\it equivalent} if there is an equivalence between them, or more generally if there is a chain of equivalences joining them. Equivalent topological categories should be regarded as ``the same'' for all relevant purposes.

\begin{remark}
According to Definition \ref{defequiv}, a functor $F: \calC \rightarrow \calD$ is an equivalence if and only if the induced functor $\h{\calC} \rightarrow \h{ \calD}$ is an equivalence. In other words, the homotopy category $\h{\calC}$ (regarded as a category which is enriched over $\calH$) is an invariant of $\calC$
which is sufficiently powerful to detect equivalences between $\infty$-categories.
This should be regarded as analogous to the more classical fact that the homotopy groups $\pi_i(X,x)$ of a CW complex $X$ are homotopy invariants which detect homotopy equivalences between CW complexes (by Whitehead's theorem). However, it is important to remember that $\h{ \calC}$ does not determine $\calC$ up to equivalence, just as the homotopy type of a CW complex is not determined by its homotopy groups.
\end{remark}

\subsection{Simplicial Categories}\label{compp1}

In the previous sections we introduced two very different approaches to the foundations of higher category theory: one based on topological categories, the other on simplicial sets. In order to prove that they are equivalent to one another, we will introduce a third approach,  which is closely related to the first but shares the combinatorial flavor of the second.

\begin{definition}\index{gen}{simplicial category}\index{gen}{category!simplicial}
A {\it simplicial category} is a category which is enriched over
the category $\sSet$ of simplicial sets. The category of simplicial
categories (where morphisms are given by simplicially enriched functors) will be
denoted by $\sCat$.\index{not}{Catsi@$\sCat$}
\end{definition}

\begin{remark}
Every simplicial category can be regarded as a simplicial object in the category $\Cat$. Conversely, a simplicial object of $\Cat$ arises from a simplicial category if and only if the underlying simplicial set of objects is constant.
\end{remark}

Like topological categories, simplicial categories can be used as models of higher category theory. If $\calC$ is a simplicial category, then we will generally think of the simplicial sets $\bHom_{\calC}(X,Y)$ as ``spaces'', or homotopy types. 

\begin{remark}
If $\calC$ is a simplicial category with the property that each of the simplicial sets
$\bHom_{\calC}(X,Y)$ is an $\infty$-category, then we may view $\calC$ itself as a kind of $\infty$-bicategory. We will not use this interpretation of simplicial categories in this book. Usually we will consider only {\em fibrant} simplicial categories: that is, simplicial categories for which the mapping objects $\bHom_{\calC}(X,Y)$ are Kan complexes.\index{gen}{fibrant!simplicial category}
\end{remark}

The relationship between simplicial categories and topological categories is easy to describe. Let $\sSet$ denote the
category of simplicial sets and $\CG$ the category of compactly
generated Hausdorff spaces. We recall that there exists a pair of
adjoint functors
$$ \Adjoint{||}{\sSet}{\CG}{\Sing}$$
which are called the {\it geometric realization} and {\it
singular complex} functors, respectively. Both of these functors commute
with finite products. Consequently, if $\calC$ is a simplicial
category, we may define a topological category $|\calC|$ in the
following way:\index{gen}{geometric realization!of simplicial categories}\index{not}{|calC|@$|\calC|$}
\index{not}{SingcalC@$\Sing \calC$}

\begin{itemize}
\item The objects of $|\calC|$ are the objects of $\calC$.

\item If $X,Y \in \calC$, then $\bHom_{|\calC|}(X,Y) = |
\bHom_{\calC}(X,Y)|$.

\item The composition law for morphisms in $|\calC|$ is obtained
from the composition law on $\calC$ by applying the geometric
realization functor.
\end{itemize}

Similarly, if $\calC$ is a topological category, we may obtain a
simplicial category $\Sing \calC$ by applying the singular complex
functor to each of the morphism spaces individually. The singular complex and geometric realization functors determine an adjunction between $\sCat$ and $\tCat$.
This adjunction should be understood as determining an ``equivalence'' between the
theory of simplicial categories and the theory of topological categories. This is essentially a formal consequence of the fact that the geometric realization and singular complex functors determine an equivalence between the homotopy theory of topological spaces and the homotopy theory of simplicial sets. More precisely, we recall that a map $f: S \rightarrow T$ of simplicial sets is said to be a {\it weak homotopy equivalence} if the induced map $|S| \rightarrow |T|$ of topological spaces is a weak homotopy equivalence. A theorem of Quillen (see \cite{goerssjardine} for a proof) asserts that the unit and counit morphisms
$$ S \rightarrow \Sing |S| $$
$$ | \Sing X | \rightarrow X$$
are weak homotopy equivalences, for every (compactly generated) topological space $X$ and every simplicial set $S$. It follows that the category obtained from $\CG$ by inverting weak homotopy equivalences (of spaces) is equivalent to the category obtained from $\sSet$ by inverting weak homotopy equivalences. We use the symbol $\calH$ to denote either of these (equivalent) categories.\index{not}{calH@$\calH$}

If $\calC$ is a simplicial category, we let $\h{\calC}$ denote the $\calH$-enriched category obtained by  applying the functor $\sSet \rightarrow \calH$ to each of the morphism spaces of $\calC$. We will refer to $\h{ \calC}$ as the {\it homotopy category of $\calC$}. We note that this is the same notation that was introduced in \S \ref{stronghcat} for the homotopy category of a topological category. However, there is little risk of confusion: the above remarks imply the existence of canonical isomorphisms\index{gen}{homotopy category!of a simplicial category}\index{not}{hcalC@$\h{\calC}$}
$$ \h{\calC} \simeq \h{ |\calC|}$$
$$ \h{\calD} \simeq \h{ \Sing \calD}$$
for every simplicial category $\calC$ and every topological category $\calD$.

\begin{definition}\index{gen}{equivalence!of simplicial categories}
A functor $\calC \rightarrow \calC'$ between simplicial categories is an {\it equivalence} if
the induced functor $\h{ \calC} \rightarrow \h{ \calC'}$ is an equivalence of $\calH$-enriched categories.
\end{definition}

In other words, a functor $\calC \rightarrow \calC'$ between simplicial categories is an equivalence if and only if the geometric realization $|\calC| \rightarrow |\calC'|$ is an equivalence of topological categories. In fact, one can say more. It follows easily from the preceding remarks that the unit and counit maps
$$ \calC \rightarrow \Sing | \calC |$$
$$ | \Sing \calD | \rightarrow \calD$$
induce {\em isomorphisms} between homotopy categories. Consequently, if we are working with topological or simplicial categories {\it up to equivalence}, we are always free to replace a simplicial category $\calC$ by $|\calC|$, or a topological category $\calD$ by $\Sing \calD$. In this sense, the notions of topological and simplicial category are equivalent and either can be used as a foundation for higher category theory.

\subsection{Comparing $\infty$-Categories with Simplicial Categories}\label{theequiv}

In \S \ref{compp1}, we introduced the theory of simplicial categories and explained why (for our purposes) it is equivalent to the theory of topological categories. In this section, we will show that the theory of simplicial categories is also closely related to the theory of $\infty$-categories.
Our discussion requires somewhat more elaborate constructions than were needed in the previous sections; a reader who does not wish to become bogged down in details is urged to skip ahead to \S \ref{working}.

We will relate simplicial categories with simplicial sets by means of the {\it simplicial nerve functor}
$$ \sNerve: \sCat \rightarrow \sSet,$$
originally introduced by Cordier (see \cite{coherentnerve}).
The nerve of an ordinary category $\calC$ is characterized by the formula
$$ \Hom_{ \sSet}( \Delta^n, \Nerve(\calC)) = \Hom_{\Cat}( [n], \calC );$$
here $[n]$ denotes the linearly ordered set $\{ 0, \ldots, n\}$, regarded
as a category. This
definition makes sense also when $\calC$ is a simplicial
category, but is clearly not very interesting: it makes no use of
the simplicial structure on $\calC$. In order to obtain a more interesting construction,
we need to replace the ordinary category $[n]$ by a suitable ``thickening'', a simplicial
category which we will denote by $\sCoNerve[\Delta^n]$\index{not}{CoNerve@$\sCoNerve[S]$}.

\begin{definition}\label{csimp1}
Let $J$ be a finite nonempty linearly ordered set. The simplicial category
$\sCoNerve[\Delta^J]$ is defined as follows:
\begin{itemize}
\item The objects of $\sCoNerve[\Delta^J]$ are the elements of
$J$.

\item If $i,j \in J$, then
$$\bHom_{\sCoNerve[\Delta^J]}(i,j) = \begin{cases} \emptyset & \text{if } j < i \\
\Nerve(P_{i,j}) & \text{if } i \leq j. \end{cases}$$
Here $P_{i,j}$ denotes the partially ordered set $\{ I \subseteq J: (i,j \in I) \wedge (
\forall k \in I) [i \leq k \leq j] ) \}$.

\item If $i_0 \leq i_1 \leq \ldots \leq i_n$, then the composition
$$ \bHom_{\sCoNerve[\Delta^J]}(i_0,i_1) \times \ldots \times
\bHom_{\sCoNerve[\Delta^J]}(i_{n-1},i_n) \rightarrow
\bHom_{\sCoNerve[\Delta^J]}(i_0,i_n)$$ is induced by the map of
partially ordered sets
$$ P_{i_0,i_1} \times \ldots  \times P_{i_{n-1},i_n} \rightarrow P_{i_0,i_n}$$
$$ ( I_1, \ldots, I_n ) \mapsto I_1 \cup \ldots \cup I_n.$$
\end{itemize}
\end{definition}

In order to help digest Definition \ref{csimp1}, let us analyze the structure of the
topological category $| \sCoNerve[ \Delta^n ] |$. The objects of this category
are elements of the set $[n] = \{ 0, \ldots, n\}$. For each
$0 \leq i \leq j \leq n$, the topological space $\bHom_{ |\sCoNerve[\Delta^n]| }(i,j)$ is homeomorphic to a cube; it may be identified with the set of all functions $p: \{ k \in [n]: i \leq k \leq j \} \rightarrow [0,1]$ which satisfy $p(i) = p(j) =
1$. The morphism space $\bHom_{ | \sCoNerve[\Delta^n] |}(i,j)$ is empty when
$j < i$, and composition of morphisms is given by concatenation of functions.

\begin{remark}\label{conervexp}
Let us try to understand better the simplicial category $\sCoNerve[\Delta^n]$ and its relationship to the ordinary category $[n]$. These categories have the same objects, namely the elements of $\{ 0, \ldots, n\}$.
In the category $[n]$, there is a unique morphism $q_{ij}: i \rightarrow j$ whenever $i \leq j$. In virtue of the uniqueness, these elements satisfy $q_{jk} \circ q_{ij} = q_{ik}$ for $i \leq j \leq k$.

In the simplicial category $\sCoNerve[\Delta^n]$, there is a vertex
$p_{ij} \in \bHom_{ \sCoNerve[\Delta^n] }(i,j)$, given by the element
$\{ i, j \} \in P_{ij}$. We note that $p_{jk} \circ p_{ij} \neq p_{ik}$ (unless we are in one of the degenerate cases where
$i =j$ or $j=k$). Instead, the collection of all compositions
$$ p_{i_n i_{n-1}} \circ p_{i_{n-1} i_{n-2}} \circ \ldots \circ p_{i_1 i_0},$$ where 
$i=i_0 < i_1 < \ldots <  i_{n-1} < i_n = j$ constitute all of the different vertices of the cube 
$\bHom_{\sCoNerve[\Delta^n]}(i,j)$. The weak contractibility of $\bHom_{\sCoNerve[\Delta^n]}(i,j)$
expresses the idea that although these compositions do not coincide, they are all canonically homotopic to one another. We observe that there is a (unique) functor
$\sCoNerve[\Delta^n] \rightarrow [n]$ which is the identity on objects, and that this functor
is an equivalence of simplicial categories. We can summarize the situation informally as follows: the simplicial category $\sCoNerve[\Delta^n]$ is a ``thickened version'' of $[n]$, where
we have dropped the strict associativity condition
$$ q_{jk} \circ q_{ij} = q_{ik}$$ and instead have imposed associativity only up to (coherent) homotopy. (We can formulate this idea more precisely by saying that $\sCoNerve[ \Delta^{\bigdot}]$ is a cofibrant replacement for $[\bigdot]$ with respect to a suitable model structure on the category of cosimplicial objects of $\sCat$.)
\end{remark}

The construction $J \mapsto \sCoNerve[\Delta^J]$ is functorial in $J$, as we now explain.

\begin{definition}\label{csimp2}
Let $f: J \rightarrow J'$ be a monotone map between linearly
ordered sets. The simplicial functor $\sCoNerve[f]:
\sCoNerve[\Delta^J] \rightarrow \sCoNerve[\Delta^{J'}]$ is defined
as follows:

\begin{itemize}
\item For each object $i \in \sCoNerve[\Delta^J]$,
$\sCoNerve[f](i) = f(i) \in \sCoNerve[\Delta^{J'}]$.

\item If $i \leq j$ in $J$, then the map
$ \bHom_{\sCoNerve[\Delta^J]}(i,j) \rightarrow
\bHom_{\sCoNerve[\Delta^{J'}]}(f(i),f(j))$ induced by $f$ is the nerve of the map $$P_{i,j}
\rightarrow P_{f(i),f(j)}$$
$$ I \mapsto f(I).$$
\end{itemize}
\end{definition}

\begin{remark}
Using the notation of Remark \ref{conervexp}, we note that Definition \ref{csimp2} has been rigged so that the functor $\sCoNerve[f]$ carries the vertex $p_{ij}
\in \bHom_{\sCoNerve[ \Delta^J ]}(i,j)$ to the vertex
$p_{f(i) f(j)} \in \bHom_{\sCoNerve[\Delta^{J'}]}(f(i), f(j))$.
\end{remark}

It is not difficult to check that the construction described in
Definition \ref{csimp2} is well-defined, and compatible with composition in $f$.
Consequently, we deduce that $\sCoNerve$ determines a functor
$$ \cDelta \rightarrow \sCat$$
$$ \Delta^n \mapsto \sCoNerve[ \Delta^n ],$$
which we may view as a cosimplicial object of $\sCat$.

\begin{definition}\index{gen}{simplicial nerve}\index{gen}{nerve!of a simplicial category}\index{not}{NervecalC@$\Nerve(\calC)$}\label{topnerve}
Let $\calC$ be a simplicial category. The {\it simplicial nerve}
$\Nerve(\calC)$ is the simplicial set described by the
formula
$$ \Hom_{ \sSet}( \Delta^n, \Nerve(\calC)) =
\Hom_{\sCat}( \sCoNerve[ \Delta^n ], \calC).$$ 

If $\calC$ is a topological category, we define the {\it
topological} {\it nerve} $\tNerve(\calC)$ of $\calC$ to be the
simplicial nerve of $\Sing \calC$.\index{gen}{topological nerve}\index{gen}{nerve!of a topological category}
\end{definition}

\begin{remark}
If $\calC$ is a simplicial (topological) category, we will often abuse terminology by referring to the 
simplicial (topological) nerve of $\calC$ simply as the {\em nerve} of $\calC$. 
\end{remark}

\begin{warning}
Let $\calC$ be a simplicial category. Then $\calC$ can be regarded as an ordinary category, by ignoring all simplices of positive dimension in the mapping spaces of $\calC$. The simplicial nerve of $\calC$ does {\em not} coincide with the nerve of this underlying ordinary category. Our notation is therefore potentially ambiguous. We will adopt the following convention: whenever $\calC$ is a simplicial category, $\Nerve(\calC)$ will denote the {\em simplicial} nerve of $\calC$, unless we specify otherwise. Similarly, if $\calC$ is a topological category, then the topological nerve
of $\calC$ does not generally coincide with the nerve of the underlying category; the notation
$\Nerve(\calC)$ will be used to indicate the topological nerve, unless otherwise specified.
\end{warning}

\begin{example}
Any ordinary category $\calC$ may be considered as a simplicial
category, by taking each of the simplicial sets
$\Hom_{\calC}(X,Y)$ to be {\em constant}. In this case, the set of
simplicial functors $\sCoNerve[\Delta^n] \rightarrow \calC$ may be
identified with the set of functors from $[n]$ into $\calC$.
Consequently, the simplicial nerve of $\calC$ agrees with the ordinary nerve of $\calC$, as defined in \S \ref{qqqc}. Similarly, the ordinary nerve of $\calC$ can be identified with the topological nerve of $\calC$, where $\calC$ is regarded as a topological category with discrete morphism spaces.
\end{example}

In order to get a feel for what the nerve of a topological
category $\calC$ looks like, let us explicitly describe its
low-dimensional simplices:

\begin{itemize}
\item The $0$-simplices of $\tNerve(\calC)$ may be identified with
the objects of $\calC$.

\item The $1$-simplices of $\tNerve(\calC)$ may be identified with
the morphisms of $\calC$.

\item To give a map from the boundary of a $2$-simplex into
$\tNerve(\calC)$ is to give a diagram (not necessarily commutative)
$$ \xymatrix{ & Y \ar[dr]^{f_{YZ}} & \\
X \ar[ur]^{f_{XY}} \ar[rr]^{f_{XZ}} & & Z. }$$
To give a $2$-simplex of $\tNerve(\calC)$ having this specified boundary is equivalent to
giving a path from $f_{XZ}$ to $f_{YZ} \circ f_{XY}$ in
$\bHom_{\calC}(X,Z)$.
\end{itemize}

The category $\sCat$ of simplicial categories admits (small)
colimits. Consequently, by formal nonsense, the functor
$\sCoNerve: \cDelta \rightarrow \sCat$ extends uniquely (up to unique isomorphism) to a
colimit-preserving functor $\sSet \rightarrow \sCat$, which we
will denote also by $\sCoNerve$. By construction, the functor
$\sCoNerve$ is left adjoint to the simplicial nerve functor $\sNerve$\index{not}{CoNerve@$\sCoNerve[S]$}. For each simplicial set $S$, we can view $\sCoNerve[S]$ as the simplicial category ``freely generated'' by $S$: every $n$-simplex $\sigma: \Delta^n \rightarrow S$ determines a functor $\sCoNerve[\Delta^n] \rightarrow \sCoNerve[S]$, which we can think of as a homotopy coherent diagram $[n] \rightarrow \sCoNerve[S]$. 

\begin{example}
Let $A$ be a partially ordered set. The simplicial category $\sCoNerve[\Nerve A]$ can be constructed using the following generalization of Definition \ref{csimp1}:
\begin{itemize}
\item The objects of $\sCoNerve[ \Nerve A]$ are the elements of $A$.
\item Given a pair of elements $a,b \in A$, the simplicial set $\bHom_{ \sCoNerve[ \Nerve A]}(a,b)$ can be identified with $\Nerve P_{a,b}$, where $P_{a,b}$ denotes the collection of linearly ordered subsets $S \subseteq A$ with least element $a$ and largest element $b$, partially ordered by inclusion.
\item Given a sequence of elements $a_0, \ldots, a_n \in A$, the composition map
$$ \bHom_{ \sCoNerve[ \Nerve A]}(a_0, a_1) \times \ldots \times \bHom_{ \sCoNerve[ \Nerve A]}( a_{n-1}, a_n) \rightarrow \bHom_{ \sCoNerve[ \Nerve A]}(a_0, a_n)$$
is induced by the map of partially ordered sets
$$ P_{a_0, a_1} \times \ldots \times P_{a_{n-1}, a_n} \rightarrow P_{ a_0, a_n}$$
$$ (S_1, \ldots, S_n) \mapsto S_1 \cup \ldots \cup S_n.$$
\end{itemize}
\end{example}

\begin{proposition}\label{toothy}
Let $\calC$ be a simplicial category having the property that, for every pair of objects
$X,Y \in \calC$, the simplicial set $\bHom_{\calC}(X,Y)$ is a Kan complex. Then the simplicial nerve $\sNerve(\calC)$ is an $\infty$-category.
\end{proposition}

\begin{proof}
We must show that if $0 < i < n$, then $\sNerve(\calC)$ has the right extension property with respect to the inclusion $\Lambda^n_i \subseteq \Delta^n$. Rephrasing this in the language of simplicial categories, we must show that $\calC$ has the right extension property with respect to the simplicial functor $\sCoNerve[\Lambda^n_i] \rightarrow \sCoNerve[\Delta^n].$
To prove this, we make use of the following observations concerning
$\sCoNerve[\Lambda^n_i]$, which we view as a simplicial subcategory 
of $\sCoNerve[\Delta^n]$:

\begin{itemize}
\item The objects of $\sCoNerve[\Lambda^n_i]$ are the objects of
$\sCoNerve[\Delta^n]$: that is, elements of the set $[n]$.

\item For $0 \leq j \leq k \leq n$, the simplicial set
$\bHom_{\sCoNerve[\Lambda^n_i]}(j,k)$ coincides with
$\bHom_{\sCoNerve[\Delta^n]}(j,k)$ unless $j=0$ and $k=n$
(note that this condition fails if $i=0$ or $i=n$).
\end{itemize}

Consequently, every extension problem
$$ \xymatrix{ \Lambda^n_i \ar@{^{(}->}[d] \ar[r]^{F} & \Nerve(\calC) \\
\Delta^n \ar@{-->}[ur] & }$$
is equivalent to
$$\xymatrix{ \bHom_{\sCoNerve[\Lambda^n_i]}(0,n) \ar[d] \ar[r] & \bHom_{\calC}(F(0), F(n)) \\
\bHom_{ \sCoNerve[\Delta^n]}(0,n) \ar@{-->}[ur]. & }$$
Since the simplicial set on the right is a Kan complex by assumption, it suffices to verify that the left vertical map is anodyne. This follows by inspection: the simplicial set $\bHom_{ \sCoNerve[ \Delta^n]}(0,n)$ can be identified with the cube $( \Delta^1 )^{ \{ 1, \ldots, n-1\} }$, and 
$\bHom_{\sCoNerve[ \Lambda^n_i]}(0,n)$ can be identified with the simplicial subset obtained by removing the interior of the cube together with one of its faces.
\end{proof}

\begin{remark}\label{goobrem}
The proof of Proposition \ref{toothy} yields a slightly stronger result: if $F: \calC \rightarrow \calD$ is a functor between simplicial categories which induces Kan fibrations
$\bHom_{\calC}(C,C') \rightarrow \bHom_{\calD}(F(C),F(C'))$ for every pair of objects $C,C' \in \calC$, then the associated map $\sNerve(\calC) \rightarrow \sNerve(\calD)$ is an inner fibration of simplicial sets (see Definition \ref{fibdeff}).
\end{remark}

\begin{corollary}\label{tooky}
Let $\calC$ be a topological category. Then the topological nerve $\tNerve(\calC)$ is an $\infty$-category.
\end{corollary}

\begin{proof}
This follows immediately from Proposition \ref{toothy}, since the singular complex of any topological space is a Kan complex.
\end{proof}

We now cite the
following theorem, which will be proven in \S \ref{compp2} and refined in \S \ref{compp3}:

\begin{theorem}\label{biggie}
Let $\calC$ be a topological category, and let $X, Y \in \calC$ be
objects. Then the counit map
$$|\bHom_{\sCoNerve[\tNerve(\calC)]}(X,Y)| \rightarrow \bHom_{\calC}(X,Y)$$
is a weak homotopy equivalence of topological spaces.
\end{theorem}

Assuming Theorem \ref{biggie}, we can now explain why the theory of $\infty$-categories is 
equivalent to the theory of topological categories (or, equivalently, simplicial categories).
The adjoint functors $\Nerve$ and $| \sCoNerve[ \bigdot ] |$ are not mutually inverse equivalences of categories. However, they {\em are} homotopy inverse to one another. To make this precise, we need to introduce a definition.

\begin{definition}\label{tulkas}\index{gen}{homotopy category!of a simplicial set}
Let $S$ be a simplicial set. The {\it homotopy category} $\h{S}$ is defined to be the homotopy category $\h{ \sCoNerve[S]}$ of the simplicial category $\sCoNerve[S]$.\index{not}{hS@$\h{S}$}
We will often view $\h{S}$ as a category enriched over the homotopy category
$\calH$ of spaces via the construction of \S \ref{compp1}: that is, for
every pair of vertices $x,y \in S$, we have $\bHom_{ \h{S}}(x,y) = [ \bHom_{ \sCoNerve[S]}(x,y) ]$.
A map $f: S \rightarrow T$ of simplicial sets is a {\it categorical equivalence} if 
the induced map $\h{S} \rightarrow \h{T}$ is an equivalence of $\calH$-enriched categories.\index{gen}{categorical equivalence}
\end{definition}

\begin{remark}
In \cite{joyalnotpub}, Joyal uses the term ``weak categorical equivalence'' for what we have called a ``categorical equivalence'', and reserves the term ``categorical equivalence'' for a stronger notion of equivalence.
\end{remark}

\begin{remark}
We have introduced the term ``categorical equivalence'', rather than simply ``equivalence'' or ``weak equivalence'', in order to avoid confusing the notion of categorical equivalence of simplicial sets with the (more classical) notion of weak homotopy equivalence of simplicial sets.
\end{remark}

\begin{remark}\label{gytyt}
It is immediate from the definition that $f: S \rightarrow T$ is a categorical equivalence if and only if
$\sCoNerve[S] \rightarrow \sCoNerve[T]$ is an equivalence (of simplicial categories), if and only if $|\sCoNerve[S]| \rightarrow |\sCoNerve[T]|$ is an equivalence (of topological categories).
\end{remark}

We now observe that the adjoint functors $(|\sCoNerve[ \bigdot ]|, \tNerve )$
determine an equivalence between the theory of simplicial sets (up to categorical equivalence)
and that of topological categories (up to equivalence). In other words, for any
topological category $\calC$ the counit map
$|\sCoNerve[ \tNerve(\calC)] | \rightarrow \calC$ is an equivalence of topological categories, and for any simplicial set $S$ the unit map
$S \rightarrow \tNerve |\sCoNerve[S]|$
is a categorical equivalence of simplicial sets. In view of Remark \ref{gytyt}, the second assertion is a formal consequence of the first. Moreover, the first assertion is merely a reformulation of Theorem \ref{biggie}.

\begin{remark}
The reader may at this point object that we have achieved a comparison between the theory of topological categories with the theory of simplicial sets, but that not every simplicial set is an $\infty$-category. However, every simplicial set is categorically equivalent to an $\infty$-category. In fact, Theorem \ref{biggie} implies that every simplicial set $S$ is categorically equivalent to the nerve of the topological category $| \sCoNerve[S] |$, which is an $\infty$-category (Corollary \ref{tooky}).
\end{remark}

\section{The Language of Higher Category Theory}\label{langur}

\setcounter{theorem}{0}

One of the main goals of this book is to demonstrate that many ideas from classical category theory can be adapted to the setting of higher categories. In this section, we will survey some of the simplest examples.

\subsection{The Opposite of an $\infty$-Category}\label{working}

If $\calC$ is an ordinary category, then the opposite category
$\calC^{op}$ is defined in the following way:

\begin{itemize}\index{gen}{opposite!of a category}
\item The objects of $\calC^{op}$ are the objects of $\calC$.
\item For $X,Y \in \calC$, we have $\Hom_{\calC^{op}}(X,Y) =
\Hom_{\calC}(Y,X)$. Identity morphisms and composition are defined
in the obvious way.
\end{itemize}

This definition generalizes without change to the setting of topological or simplicial categories. Adapting this definition to the setting of $\infty$-categories requires a few additional words.
We may define more generally the {\it opposite} of a simplicial set
$S$ as follows: For any finite, nonempty, linearly ordered set
$J$, we set $S^{op}(J) = S(J^{op})$, where $J^{op}$ denotes the
same set $J$ endowed with the opposite ordering. More concretely,
we have $S^{op}_n = S_n$, but the face and degeneracy maps on
$S^{op}$ are given by the formulas
$$ (d_i: S^{op}_n \rightarrow S^{op}_{n-1}) = (d_{n-i}: S_n
\rightarrow S_{n-1})$$
$$ (s_i: S^{op}_n \rightarrow S^{op}_{n+1}) = (s_{n-i}: S_n
\rightarrow S_{n+1}).$$\index{gen}{opposite!of a simplicial set}

The formation of opposite categories is fully compatible
with all of the constructions we have introduced for passing back
and forth between different models of higher category theory.

It is clear from the definition that a simplicial set $S$ is an $\infty$-category if and only if its opposite $S^{op}$ is an $\infty$-category: for $0 < i < n$, $S$ has the extension property
with respect to the horn inclusion $\Lambda^n_i \subseteq \Delta^n$ if
and only if $S^{op}$ has the extension property with respect to the horn inclusion
$\Lambda^n_{n-i} \subseteq \Delta^n$.

The construction $S \mapsto S^{op}$ determines an automorphism of the $\infty$-category
of $\infty$-categories. We will later see that this is (essentially) the {\em only} nontrivial automorphism
(see Theorem \ref{cabbi}).

\subsection{Mapping Spaces in Higher Category Theory}\label{prereq1}

If $X$ and $Y$ are objects of an ordinary category $\calC$, then one has a well-defined
set $\Hom_{\calC}(X,Y)$ of morphisms from $X$ to $Y$. In higher category theory, one has instead a morphism {\em space} $\bHom_{\calC}(X,Y)$. In the setting of topological or simplicial
categories, this morphism space (either a topological space or a simplicial set) is an inherent feature of the formalism. In the setting of $\infty$-categories, it is not so obvious
how $\bHom_{\calC}(X,Y)$ should be defined. However, it is at least clear what to do on the level of the homotopy category.

\begin{definition}\label{morspace}\index{not}{MapS@$\bHom_{S}(X,Y)$}
Let $S$ be a simplicial set containing vertices $x$ and $y$, and let
$\calH$ denote the homotopy category of spaces. We define
$\bHom_{S}(x,y) = \bHom_{\h{S}}(x,y) \in \calH$ to be the object of $\calH$ representing
the space of maps from $x$ to $y$ in $S$. Here $\h{S}$ denotes the homotopy category of $S$, regarded as a $\calH$-enriched category (Definition \ref{tulkas}). 
\end{definition}

\begin{warning}
Let $S$ be a simplicial set. The notation $\bHom_{S}( X, Y)$ has two {\em very} different meanings.
When $X$ and $Y$ are vertices of $S$, then our notation should be interpreted in the sense of Definition \ref{morspace}, so that $\bHom_{S}(X,Y)$ is an object of $\calH$. If $X$ and $Y$ are objects of $(\sSet)_{/S}$, then we instead let $\bHom_{S}(X,Y)$ denote the simplicial mapping object
$$ Y^{X} \times_{ S^{X} } \{ \phi \} \in \sSet,$$
where $\phi$ denotes the structural morphism $X \rightarrow S$. We trust that it will be clear in context which of these two definitions applies in a given situation.
\end{warning}

We now consider the following question: given a simplicial set $S$ containing a pair of vertices $x$ and $y$, how can we compute $\bHom_{S}(x,y)$? We have defined
$\bHom_{S}(x,y)$ as an object of the homotopy category $\calH$, but for many purposes it is important to choose a simplicial set $M$ which represents $\bHom_{S}(x,y)$. 
The most obvious candidate for $M$ is the simplicial set
$\bHom_{\sCoNerve[S]}(x,y)$. The advantages of this definition are that it works in all cases (that is, $S$ does not need to be an $\infty$-category), and comes equipped with an associative composition law. However, the construction of the simplicial set
$\bHom_{\sCoNerve[S]}(x,y)$ is quite complicated. Furthermore,
$\bHom_{\sCoNerve[S]}(x,y)$ is usually not a Kan complex, so it
can be difficult to extract algebraic invariants like homotopy
groups, even when a concrete description of its simplices is known. 

In order to address these shortcomings, we will introduce another simplicial set which
represents the homotopy type $\bHom_{S}(x,y) \in \calH$, at least when $S$ is an $\infty$-category. 
We define a new simplicial
set $\Hom^{\rght}_S(x,y)$, the space of {\it right morphisms} from
$x$ to $y$, by letting $\Hom_{\sSet}( \Delta^n, \Hom^{\rght}_S(x,y))$ denote the set of
all $z: \Delta^{n+1} \rightarrow S$ such that $z| \Delta^{ \{n+1 \}} = y$ and $z|
\Delta^{ \{0, \ldots, n\} }$ is a constant simplex at the vertex
$x$. The face and degeneracy operations on
$\Hom^{\rght}_S(x,y)_n$ are defined to coincide with corresponding operations on 
$S_{n+1}$.\index{not}{HomR@$\Hom^{\rght}_{S}(X,Y)$}

We first observe that when $S$ is an $\infty$-category, $\Hom^{\rght}_S(x,y)$ really is a ``space'':

\begin{proposition}\label{gura}
Let $\calC$ be an $\infty$-category containing a pair of objects $x$ and $y$. The simplicial set
$\Hom^{\rght}_{\calC}(x,y)$ is a Kan complex.
\end{proposition}

\begin{proof}
It is immediate from the definition that if $\calC$ is a
$\infty$-category, then $M=\Hom^{\rght}_{\calC}(x,y)$ satisfies the Kan
extension condition for every horn inclusion $\Lambda^n_i
\subseteq \Delta^n$ where $0 < i \leq n$. This implies that $M$ is
a Kan complex (Proposition \ref{greenwich}).
\end{proof}

\begin{remark}\label{needie}
If $S$ is a simplicial set and $x,y,z \in S_0$, then there is no
obvious composition law
$$\Hom^{\rght}_{S}(x,y) \times \Hom^{\rght}_{S}(y,z) \rightarrow \Hom^{\rght}_{S}(x,z).$$
We will later see that if $S$ is an $\infty$-category, then there is a
composition law which is well-defined up to a
contractible space of choices. The absence of a canonical choice for a composition
law is the main drawback of $\Hom^{\rght}_{S}(x,y)$, in comparison with
$\bHom_{\sCoNerve[S]}(x,y).$
The main goal of \S \ref{valencequi} is to show that, if $S$ is an $\infty$-category, then there is a 
(canonical) isomorphism between $\Hom^{\rght}_{S}(x,y)$ and $\bHom_{\sCoNerve[S]}(x,y)$
in the homotopy category $\calH$. In particular, we will conclude that
$\Hom^{\rght}_{S}(x,y)$ represents $\bHom_{S}(x,y)$, whenever $S$ is an $\infty$-category.
\end{remark}

\begin{remark}\index{not}{HomLS@$\Hom^{\lft}_{S}(X,Y)$}\label{swink}
The definition of $\Hom^{\rght}_{S}(x,y)$ is not self-dual: that
is, $\Hom^{\rght}_{S^{op}}(x,y) \neq \Hom^{\rght}_{S}(y,x)$ in
general. Instead we define $\Hom^{\lft}_S(x,y) =
\Hom^{\rght}_{S^{op}}(y,x)^{op}$, so that $\Hom^{\lft}_{S}(x,y)_n$ is
the set of all $z \in S_{n+1}$ such that $z|\Delta^{ \{0\} } = x$
and $z | \Delta^{ \{1, \ldots, n+1\} }$ is the constant simplex at
the vertex $y$. \end{remark}

Although the simplicial sets $\Hom^{\lft}_S(x,y)$ and
$\Hom^{\rght}_S(x,y)$ are generally not isomorphic to one another,
they are homotopy equivalent whenever $S$ is an $\infty$-category. To
prove this, it is convenient to define a third, self-dual, space
of morphisms: let $\Hom_{S}(x,y) = \{x\} \times_S S^{\Delta^1}
\times_S \{y\}$. In other words, to give an $n$-simplex of
$\Hom_{S}(x,y)$, one must give a map $f: \Delta^n \times \Delta^1
\rightarrow S$, such that $f| \Delta^n \times \{0\}$ is constant
at $x$ and $f| \Delta^n \times \{1\}$ is constant at $y$. We
observe that there exist natural inclusions
$$ \Hom^{\rght}_{S}(x,y) \hookrightarrow \Hom_S(x,y) \hookleftarrow
\Hom^{\lft}_{S}(x,y),$$ 
which are induced by retracting the cylinder $\Delta^n \times \Delta^1$ onto certain maximal dimensional simplices. We will later show (Corollary
\ref{homsetsagree}) that these inclusions are homotopy
equivalences, provided that $S$ is an $\infty$-category.\index{not}{HomS@$\Hom_{S}(X,Y)$}

\subsection{The Homotopy Category}\label{hcat}

For every ordinary category $\calC$, the nerve $\Nerve(\calC)$ is an $\infty$-category. 
Informally, we can describe the situation as follows: the nerve functor is a fully faithful inclusion from the bicategory of categories to the $\infty$-bicategory of $\infty$-categories.
Moreover, this inclusion has a left adjoint:

\begin{proposition}\label{leftadj}
The nerve functor $\Cat \rightarrow \sSet$ is right adjoint to the functor
$\h: \sSet \rightarrow \Cat$, which associates to every simplicial set $S$ its
homotopy category $\h{S}$ $($here we ignore the $\calH$-enrichment of $\h{S}${}$)$.
\end{proposition}

\begin{proof}
Let us temporarily distinguish between the nerve functor $\Nerve: \Cat \rightarrow \sSet$
and the simplicial nerve functor $\Nerve': \sCat \rightarrow \sSet$. These two functors are related by the fact that $\Nerve$ can be written as a composition
$$ \Cat \stackrel{i}{\subseteq} \sCat \stackrel{\Nerve'}{\rightarrow} \sSet.$$
The functor $\pi_0: \sSet \rightarrow \Set$ is a left adjoint to the inclusion functor
$\Set \rightarrow \sSet$, so the functor
$$ \sCat \rightarrow \Cat$$
$$ \calC \mapsto \h{\calC}$$
is left adjoint to $i$. It follows that $\Nerve = \Nerve' \circ i$ has a left adjoint, given by the composition
$$ \sSet \stackrel{ \sCoNerve[ \bigdot] }{\rightarrow} \sCat \stackrel{ \h{}}{\rightarrow} \Cat,$$
which coincides with the homotopy category functor $\h{}: \sSet \rightarrow \Cat$ by definition.
\end{proof}

\begin{remark}
The formation of the homotopy category is literally left adjoint to
the inclusion $\Cat \subseteq \sCat$. The analogous assertion is not
quite true in the setting of topological categories, since the
functor $\pi_0: \CG \rightarrow \Set$ is a left adjoint only when
restricted to locally path connected spaces.
\end{remark}

\begin{warning}
If $\calC$ is a simplicial category, then we do not necessarily
expect that $\h{\calC} \simeq \h{\sNerve(\calC)}$. However, this is always
the case when $\calC$ is {\it fibrant}, in the sense that every
simplicial set $\bHom_{\calC}(X,Y)$ is a Kan complex.
\end{warning}

\begin{remark}
If $S$ is a simplicial set, Joyal (\cite{joyalnotpub}) refers to the category $\h{S}$ as the {\it fundamental
category} of $S$. This is motivated by the observation that if $S$ is a Kan complex, then $\h{S}$ is
the fundamental groupoid of $S$ in the usual sense.
\end{remark}

Our objective, for the remainder of this section, is to obtain a more explicit understanding of
the homotopy category $\h{S}$ of a simplicial set $S$. Proposition \ref{leftadj} implies that $\h{S}$ 
admits the following presentation by generators and relations:

\begin{itemize}
\item The objects of $\h{S}$ are the vertices of $S$.

\item For every edge $\phi: \Delta^1 \rightarrow S$, there is a morphism $\overline{\phi}$ from
$\phi(0)$ to $\phi(1)$.

\item For each $\sigma: \Delta^2 \rightarrow S$, we have $\overline{d_0(\sigma)} \circ \overline{d_2(\sigma)} = \overline{d_1(\sigma)}$.

\item For each vertex $x$ of $S$, the morphism $\overline{s_0 x}$ is the identity $\id_x$. 
\end{itemize}

If $S$ is an $\infty$-category, there is a much more satisfying construction of the category $\h{S}$.
We will describe this construction in detail, since it nicely illustrates the utility of the weak Kan condition of Definition \ref{qqcc}.

Let $\calC$ be an $\infty$-category. We will construct a category $\pi(\calC)$\index{gen}{homotopy category!of an $\infty$-category}
(which we will eventually show to be equivalent to the homotopy category $\h{\calC}$). The objects of $\pi(\calC)$ are the vertices of $\calC$. Given an
edge $\phi: \Delta^1 \rightarrow \calC$, we shall say that $\phi$ has {\it source}
$C= \phi(0)$ and {\it target} $C' = \phi(1) $, and write $\phi:
C \rightarrow C'$. For each object $C$ of $\calC$, we
let $\id_{C}$ denote the degenerate edge $s_0(C): C \rightarrow C$.

Let $\phi: C \rightarrow C'$ and $\phi': C \rightarrow C'$ be a pair of edges of $\calC$ having
the same source and target. We will say that $\phi$ and $\phi'$ are {\it homotopic} if
there is a $2$-simplex $\sigma: \Delta^2 \rightarrow \calC$, which we depict as follows:
$$ \xymatrix{
& C' \ar[dr]^{\id_{C'}} & \\
C \ar[ur]^{\phi} \ar[rr]^{\phi'} & & C'.}$$
In this case, we say that $\sigma$ is a {\it homotopy} between $\phi$ and $\phi'$.\index{gen}{homotopy!between morphisms of $\calC$}

\begin{proposition}\label{extneeded}
Let $\calC$ be an $\infty$-category, and let $C$ and $C'$ be objects of
$\pi(\calC)$. Then the relation of homotopy is an equivalence relation
on the edges joining $C$ to $C'$.
\end{proposition}

\begin{proof}
Let $\phi: \Delta^1 \rightarrow \calC$ be an edge. Then $s_1(\phi)$ is a homotopy from $\phi$ to
itself. Thus homotopy is a reflexive relation.

Suppose next that $\phi, \phi', \phi'': C \rightarrow C'$ are edges with the same source and target. Let $\sigma$
be a homotopy from $\phi$ to $\phi'$ and $\sigma'$ a homotopy from
$\phi$ to $\phi''$. Let $\sigma'': \Delta^2 \rightarrow \calC$ denote the
constant map at the vertex $C'$. 
We have a commutative diagram
$$ \xymatrix{ \Lambda^3_1 \ar@{^{(}->}[d] \ar[rr]^{ (\sigma'', \bigdot, \sigma', \sigma) } & & \calC \\
\Delta^3 \ar@{-->}[urr]^{\tau}. & & }$$
Since $\calC$ is an $\infty$-category, there exists a $3$-simplex $\tau: \Delta^3 \rightarrow \calC$ as indicated by the dotted arrow in the diagram. It is easy to see that
$d_1(\tau)$ is a homotopy from $\phi'$ to $\phi''$.

As a special case, we may take $\phi=\phi''$; we then deduce that
the relation of homotopy is symmetric. It then follows immediately
from the above that the relation of homotopy is also transitive.
\end{proof}

\begin{remark}\label{cello}
The definition of homotopy that we have given is not evidently self-dual; in other words, 
it is not immediately obvious a homotopic pair of edges $\phi, \phi': C \rightarrow C'$ of an $\infty$-category $\calC$ remain homotopic when regarded as edges in the opposite $\infty$-category $\calC^{op}$. To prove this, let $\sigma$ be a homotopy from $\phi$ to $\phi'$, and consider
the commutative diagram
$$ \xymatrix{ \Lambda^3_2 \ar@{^{(}->}[d] \ar[rr]^{ ( \sigma, s_1 \phi, \bigdot, s_0 \phi)} & & \calC \\
\Delta^3 \ar@{-->}[urr]^{\tau}. & & }$$
The assumption that $\calC$ is an $\infty$-category guarantees a $3$-simplex
$\tau$ rendering the diagram commutative. The face $d_2 \tau$ may be regarded as a homotopy
from $\phi'$ to $\phi$ in $\calC^{op}$.
\end{remark}

We can now define the morphism sets of the category $\pi(\calC)$: given vertices $X$ and $Y$ of $\calC$, we let $\Hom_{\pi(\calC)}(X,Y)$ denote the set of homotopy
classes of edges $\phi: X \rightarrow Y$ in $\calC$. For each edge
$\phi: \Delta^1 \rightarrow \calC$, we let $[ \phi ]$ denote the corresponding morphism in
$\pi(\calC)$.

We define a composition law on $\pi(\calC)$ as follows. Suppose that $X$, $Y$, and $Z$ are vertices of $\calC$, and that we are given 
edges $\phi: X \rightarrow Y$, $ \psi : Y \rightarrow Z$.
The pair $(\phi, \psi)$ determines a map $\Lambda^2_1 \rightarrow \calC$. Since $\calC$ is an $\infty$-category, this map extends to a $2$-simplex $\sigma: \Delta^2 \rightarrow \calC$. We now define
$[ \psi] \circ [ \phi ] = [ d_1 \sigma ]$.

\begin{proposition}\label{trug}
Let $\calC$ be an $\infty$-category.
The composition law on $\pi(\calC)$ is well-defined. In other words, the homotopy class $[\psi] \circ [\phi]$ does not depend on the choice of $\psi$ representing $[\psi]$, the choice of $\phi$ representing $[ \phi]$, or the choice of the the $2$-simplex $\sigma$.
\end{proposition}

\begin{proof}
We begin by verifying the independence of the choice of $\sigma$. Suppose that we are given two $2$-simplices $\sigma, \sigma': \Delta^2 \rightarrow \calC$, satisfying
$$ d_0 \sigma = d_0 \sigma' = \psi $$
$$ d_2 \sigma = d_2 \sigma' = \phi.$$
Consider the diagram
$$ \xymatrix{ \Lambda^3_1 \ar@{^{(}->}[d] \ar[rr]^{( s_1 \psi, \bigdot, \sigma', \sigma)} & & \calC \\
\Delta^3 \ar@{-->}[urr]^{\tau}. & & }$$
Since $\calC$ is an $\infty$-category, there exists a $3$-simplex $\tau$ as indicated by the dotted arrow. It follows that $d_1 \tau$ is a homotopy from $d_1 \sigma$ to $d_1 \sigma'$.

We now show that $[\psi ] \circ [\phi]$ depends only on $\psi$ and $\phi$ only up to homotopy.
In view of Remark \ref{cello}, the assertion is symmetric with respect to $\psi$ and $\phi$; it will therefore suffice to show that $[ \psi ] \circ [ \phi ]$ does not change if we replace $\phi$ by a morphism $\phi'$ which is homotopic to $\phi$. Let $\sigma$ be a $2$-simplex with
$d_0 \sigma = \psi$, $d_2 \sigma = \phi$, and let $\sigma'$ be a homotopy from $\phi$ to $\phi'$.
Consider the diagram
$$ \xymatrix{ \Lambda^3_1 \ar@{^{(}->}[d] \ar[rr]^{ (s_0 \psi, \bigdot, \sigma, \sigma') } & & \calC \\
\Delta^3. \ar@{-->}[urr]^{\tau} & & }$$
Again, the hypothesis that $\calC$ is an $\infty$-category guarantees the existence of a $3$-simplex $\tau$ as indicated in the diagram. Let $\sigma'' = d_1 \tau$. Then 
$[ \psi ] \circ [ \phi' ] = [ d_1 \sigma' ]$. But $d_1 \sigma = d_1 \sigma'$ by construction, so that
$[ \psi ] \circ [ \phi ] = [ \psi ] \circ [ \phi' ]$ as desired.
\end{proof}

\begin{proposition}
If $\calC$ is an $\infty$-category, then $\pi(\calC)$ is a category.
\end{proposition}

\begin{proof}
Let $C$ be a vertex of $\calC$. We first verify that $[\id_{C}]$ is an identity with respect to the composition law on $\pi(\calC)$. For every edge $\phi: C' \rightarrow C$ in $\calC$, the
$2$-simplex $s_1(\phi)$ verifies the equation
$$ [ \id_{C} ] \circ [\phi] = [\phi].$$
This proves that $\id_{C}$ is a left
identity; the dual argument (Remark \ref{cello}) shows that $[\id_{C}]$ is a right
identity.

The only other thing we need to check is the associative law for
composition in $\pi(\calC)$. Suppose given a composable sequence of edges
$$ C \stackrel{\phi}{\rightarrow} C' \stackrel{\phi'}{\rightarrow} C'' \stackrel{\phi''}{\rightarrow} C'''.$$
Choose $2$-simplices $\sigma, \sigma', \sigma'': \Delta^2 \rightarrow \calC$, corresponding to diagrams
$$ \xymatrix{ & C' \ar[dr]^{\phi'} & \\
C \ar[ur]^{\phi} \ar[rr]^{\psi} & & C'' }$$
$$ \xymatrix{ & C'' \ar[dr]^{\phi''} & \\
C \ar[ur]^{\psi} \ar[rr]^{\theta} & & C''' }$$
$$ \xymatrix{ & C'' \ar[dr]^{\phi''} & \\
C' \ar[ur]^{\phi'} \ar[rr]^{\psi'} & & C''',}$$
respectively. Then $[\phi'] \circ [\phi] = [\psi]$, $[\phi''] \circ [\psi] = [\theta]$, and
$[\phi''] \circ [\phi'] = [\psi']$. Consider the diagram
$$ \xymatrix{
\Lambda^3_2 \ar[rr]^{(\sigma'', \sigma', \bigdot, \sigma)} \ar@{^{(}->}[d] & & \calC \\
\Delta^3. \ar@{-->}[urr]^{\tau} & & }$$
Since $\calC$ is an $\infty$-category, there exists a $3$-simplex $\tau$ rendering the diagram commutative. Then $d_2(\tau)$ verifies the equation
$[\psi'] \circ [\phi] = [\theta]$, so that
$$([\phi''] \circ [\phi']) \circ [\phi] = [\theta] = [\phi''] \circ
[\psi] = [\phi''] \circ ([\phi'] \circ [\phi])$$ as desired.
\end{proof}

We now show that if $\calC$ is an $\infty$-category, then $\pi(\calC)$ is
naturally equivalent (in fact isomorphic) to $\h{\calC}$.

\begin{proposition}
Let $\calC$ be an $\infty$-category. There exists a unique functor $F: \h{\calC} \rightarrow \pi(\calC)$ with the following properties:
\begin{itemize}
\item[$(1)$] On objects, $F$ is the identity map.
\item[$(2)$] For every edge $\phi$ of $\calC$, $F( \overline{\phi} ) = [\phi]$.
\end{itemize}
Moreover, $F$ is an isomorphism of categories.
\end{proposition}

\begin{proof}
The existence and uniqueness of $F$ follows immediately from our presentation
of $\h{\calC}$ by generators and relations. It is obvious that $F$ is bijective on objects and surjective on morphisms. To complete the proof, it will suffice to show that $F$ is faithful.

We first show that every morphism $f: x \rightarrow y$ in $\h{\calC}$ may be written as $\overline{\phi}$ for some $\phi \in \calC$. Since the morphisms in $\h{\calC}$ are generated by morphisms having the form $\overline{\phi}$ under composition, it suffices to show that the set of such morphisms contains all identity morphisms and is stable under composition. The first assertion is clear, since $\overline{ s_0 x}=\id_x$. For the second, we note that if $\phi: x \rightarrow y$ and $\phi': y \rightarrow z$ are composable edges, then there exists a 2-simplex $\sigma: \Delta^2 \rightarrow \calC$, which we may depict as follows:
$$ 
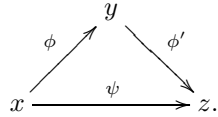
\xymatrix{ & y \ar[dr]^{\phi'} & \\
x \ar[ur]^{\phi} \ar[rr]^{\psi} & & z. }$$
Thus $\overline{ \phi' } \circ \overline{\phi} = \overline{ \psi }$.

Now suppose that $\phi,\phi': x \rightarrow y$ are such that $[\phi]=[\phi']$; we wish to show that $\overline{\phi}=\overline{\phi'}$. By definition, there exists a homotopy $\sigma: \Delta^2 \rightarrow \calC$ joining $\phi$ and $\phi'$. The existence of $\sigma$ entails the relation
$$ \id_y \circ \overline{\phi} = \overline{\phi'} $$ in the homotopy category $\h{S}$, so that
$\overline{\phi} = \overline{\phi'}$ as desired.
\end{proof}

\subsection{Objects, Morphisms and Equivalences}\label{obmor}

As in ordinary category theory, we may speak of {\it objects} and {\it morphisms} in a higher category $\calC$. If $\calC$ is a topological (or simplicial) category, these
should be understood literally as the objects and morphisms in the
underlying category of $\calC$. We may also apply this terminology
to $\infty$-categories (or even more general simplicial sets): if $S$ is a simplicial set, then the {\it
objects} of $S$ are the vertices $\Delta^0 \rightarrow S$, and the {\it morphisms}
of $S$ are edges $\Delta^1 \rightarrow S$. A morphism $\phi: \Delta^1 \rightarrow S$ is said to have {\it source} $X= \phi(0)$ and {\it target} $Y= \phi(1)$;
we will often denote this by writing $\phi: X \rightarrow Y$.\index{gen}{object!of an $\infty$-category}\index{gen}{morphism!in an $\infty$-category}
If $X: \Delta^0 \rightarrow S$ is an object of $S$, we will write
$\id_X = s_0(X): X \rightarrow X$ and refer to this as the {\it identity morphism} of $X$.

If $f,g: X \rightarrow Y$ are two morphisms in a higher category $\calC$, then $f$ and $g$
are {\it homotopic} if they determine the same morphism in the homotopy category $\h{\calC}$. In the setting of $\infty$-categories, this coincides with the notion of homotopy introduced in the previous section. In the setting of
topological categories, this simply means that $f$ and $g$ lie
in the same path component of $\bHom_{\calC}(X,Y)$. In either case,
we will sometimes indicate this relationship between $f$ and $g$ by writing $f
\simeq g$.

A morphism $f: X \rightarrow Y$ in an $\infty$-category $\calC$ is
said to be an {\it equivalence} if it determines an isomorphism in the homotopy category
$\h{\calC}$. We say that $X$ and $Y$ are {\it equivalent} if there is
an equivalence between them (in other words, if they are
isomorphic as objects of $\h{\calC}$).\index{gen}{equivalence!in an $\infty$-category}

If $\calC$ is a topological category, then the
requirement that a morphism $f: X \rightarrow Y$ be an equivalence
is quite a bit weaker than the requirement that $f$ be an
isomorphism. In fact, we have the following:

\begin{proposition}\label{rooot}\index{gen}{equivalence!in a topological category}
Let $f: X \rightarrow Y$ be a morphism in a topological category.
The following conditions are equivalent:

\begin{itemize}
\item[$(1)$] The morphism $f$ is an equivalence.

\item[$(2)$] The morphism $f$ has a {\it homotopy inverse} $g: Y
\rightarrow X$; that is, a morphism $g$ such that $f \circ g
\simeq \id_Y$ and $g \circ f \simeq \id_X$.

\item[$(3)$] For every object $Z \in \calC$, the induced map
$\bHom_{\calC}(Z,X) \rightarrow \bHom_{\calC}(Z,Y)$ is a homotopy
equivalence.

\item[$(4)$] For every object $Z \in \calC$, the induced map
$\bHom_{\calC}(Z,X) \rightarrow \bHom_{\calC}(Z,Y)$ is a weak
homotopy equivalence.

\item[$(5)$] For every object $Z \in \calC$, the induced map
$\bHom_{\calC}(Y,Z) \rightarrow \bHom_{\calC}(X,Z)$ is a homotopy
equivalence.

\item[$(6)$] For every object $Z \in \calC$, the induced map
$\bHom_{\calC}(Y,Z) \rightarrow \bHom_{\calC}(X,Z)$ is a weak
homotopy equivalence.
\end{itemize}
\end{proposition}

\begin{proof}
It is clear that $(2)$ is merely a reformulation of $(1)$. We will
show that $(2) \Rightarrow (3) \Rightarrow (4) \Rightarrow (1)$;
the implications $(2) \Rightarrow (5) \Rightarrow (6) \Rightarrow
(1)$ follow using the same argument.

To see that $(2)$ implies $(3)$, we note that if $g$ is a homotopy
inverse to $f$, then composition with $g$ gives a map
$\bHom_{\calC}(Z,Y) \rightarrow \bHom_{\calC}(Z,X)$ which is
homotopy inverse to composition with $f$. It is clear that $(3)$
implies $(4)$. Finally, if $(4)$ holds, then we note that $X$ and
$Y$ represent the same functor on $\h{\calC}$ so that $f$ induces an
isomorphism between $X$ and $Y$ in $\h{\calC}$.
\end{proof}

\begin{example}
Let $\calC$ be the category of CW-complexes, considered as a topological category by endowing
each of the sets $\Hom_{\calC}(X,Y)$ with the (compactly generated) compact open topology. A pair of objects $X,Y \in \calC$ are equivalent (in the sense defined above) if and only if they are homotopy equivalent (in the sense of classical topology).
\end{example}

If $\calC$ is an $\infty$-category (topological category, simplicial category), then we
shall write $X \in \calC$ to mean that $X$ is an object of
$\calC$. We will generally understand that all meaningful properties of
objects are invariant under equivalence. Similarly, all
meaningful properties of morphisms are invariant under
homotopy and under composition with equivalences.

In the setting of $\infty$-categories, there is a very useful characterization of equivalences which is due to Joyal.

\begin{proposition}[Joyal \cite{joyalnotpub}]\label{greenlem}\index{gen}{equivalence!in an $\infty$-category}
Let $\calC$ be an $\infty$-category, and $\phi: \Delta^1 \rightarrow \calC$ a morphism of $\calC$. Then $\phi$ is an equivalence if and only if, for every $n \geq 2$ and every map
$f_0: \Lambda^n_0 \rightarrow \calC$ such that $f_0 | \Delta^{\{0,1\}} = \phi$,
there exists an extension of $f_0$ to $\Delta^n$.
\end{proposition}

The proof requires some ideas which we have not yet introduced, and will be given in \S \ref{leftfib}.

\subsection{$\infty$-Groupoids and Classical Homotopy Theory}

Let $\calC$ be an $\infty$-category. We will say that $\calC$ is an {\it $\infty$-groupoid} if the homotopy category $\h{\calC}$ is a groupoid: in other words, if every morphism in $\calC$ is an equivalence. In \S \ref{highcat}, we asserted that the theory of $\infty$-groupoids is equivalent to classical homotopy theory. We can now formulate this idea in a very precise way:

\begin{proposition}[Joyal \cite{joyalpub}]\label{greenwich}\index{gen}{$\infty$-groupoid}
Let $\calC$ be a simplicial set. The following conditions are
equivalent:

\begin{itemize}
\item[$(1)$] The simplicial set $\calC$ is an $\infty$-category and its homotopy category $\h{\calC}$ is a groupoid.

\item[$(2)$] The simplicial set $\calC$ satisfies the extension condition
for all horn inclusions $\Lambda^n_i \subseteq \Delta^n$ for $0 \leq i < n$.

\item[$(3)$] The simplicial set $\calC$ satisfies the extension condition
for all horn inclusions $\Lambda^n_i \subseteq \Delta^n$ for $0 < i \leq n$.

\item[$(4)$] The simplicial set $\calC$ is a Kan complex; in other words, it
satisfies the extension condition for all horn inclusions
$\Lambda^n_i \subseteq \Delta^n$ for $0 \leq i \leq n$.
\end{itemize}
\end{proposition}

\begin{proof}
The equivalence $(1) \Leftrightarrow (2)$ follows immediately from
Proposition \ref{greenlem}.
Similarly, the equivalence $(1) \Leftrightarrow (3)$ follows by applying Proposition \ref{greenlem} to $\calC^{op}$. We conclude by observing that $(4) \Leftrightarrow (2) \wedge (3)$.
\end{proof}

\begin{remark}
The assertion that we can identify $\infty$-groupoids with spaces is less obvious in other formulations of higher category theory. For example, suppose that $\calC$ is a topological category whose homotopy category $\h{\calC}$ is a groupoid. For simplicity, we will assume furthermore that $\calC$ has a single object $X$. We may then identify $\calC$ with
the topological monoid $M=\Hom_{\calC}(X,X)$. The assumption that
$\h{ \calC}$ is a groupoid is equivalent to the assumption
that the discrete monoid $\pi_0 M$ is a group. In this case, one can show that the 
unit map $M \rightarrow \Omega BM$ is a weak homotopy
equivalence, where $BM$ denotes the classifying space of the
topological monoid $M$. In other words, up to equivalence,
specifying $\calC$ (together with the object $X$) is equivalent to
specifying the space $BM$ (together with its base point).
\end{remark}

Informally, we might say that the inclusion functor $i$ from Kan
complexes to $\infty$-categories exhibits the $\infty$-category of
(small) $\infty$-groupoids as a full subcategory of
the $\infty$-bicategory of (small) $\infty$-categories. Conversely, every
$\infty$-category $\calC$ has an ``underlying'' $\infty$-groupoid, which is obtained by discarding the noninvertible morphisms of $\calC$:\index{gen}{$\infty$-groupoid!underlying an $\infty$-category}

\begin{proposition}[\cite{joyalnotpub}]\label{lumba}
Let $\calC$ be an $\infty$-category. Let $\calC' \subseteq \calC$ be the largest simplicial
subset of $\calC$ having the property that every edge of $\calC'$ is an equivalence in $\calC$.
Then $\calC'$ is a Kan complex. It may be characterized by the
following universal property: for any Kan complex $K$, the induced
map $\Hom_{\sSet}(K,\calC') \rightarrow \Hom_{\sSet}(K,\calC)$ is a
bijection.
\end{proposition}

\begin{proof}
It is straightforward to check that $\calC'$ is an $\infty$-category. Moreover, if
$f$ is a morphism in $\calC'$, then $f$ has a homotopy inverse $g \in \calC$. Since
$g$ is itself an equivalence in $\calC$, we conclude that $g$ belongs to $\calC'$ and
is therefore a homotopy inverse to $f$ in $\calC'$. In other words, every morphism in $\calC'$ is an equivalence, so that $\calC'$ is a Kan complex by Proposition \ref{greenwich}. To prove the last assertion, we observe that if $K$ is an $\infty$-category, then any map of simplicial sets
$\phi: K \rightarrow \calC$ carries equivalences in $K$ to equivalences in $\calC$. In particular, if
$K$ is a Kan complex, then $\phi$ factors (uniquely) through $\calC'$.
\end{proof}

It follows from Proposition \ref{lumba} that the functor
$\calC \mapsto \calC'$ is right adjoint to
the inclusion functor from Kan complexes to $\infty$-categories. It is easy to see that this right adjoint is an invariant notion: that is, a categorical equivalence of $\infty$-categories $\calC \rightarrow \calD$ induces a homotopy equivalence
$\calC' \rightarrow \calD'$ of Kan complexes. 

\begin{remark}
It is easy to give analogous constructions in the case of topological or simplicial categories. For example, if $\calC$ is a topological category, then we can define $\calC'$ to be another topological category with the same objects as $\calC$, where $\bHom_{\calC'}(X,Y) \subseteq \bHom_{\calC}(X,Y)$ is the subspace consisting of equivalences in $\bHom_{\calC}(X,Y)$, equipped with the subspace topology.
\end{remark}

\begin{remark}
We will later introduce a relative version of the construction described in Proposition \ref{lumba}, which applies to certain families of $\infty$-categories (Corollary \ref{relativeKan}).
\end{remark}

Although the inclusion functor from Kan complexes to $\infty$-categories does not
literally have a left adjoint, it does have a left adjoint in a higher-categorical sense. This left adjoint is computed by any ``fibrant replacement'' functor (for the usual model structure) from $\sSet$ to itself, for
example the functor $S \mapsto \Sing |S|$.
The unit map $u: S \rightarrow
\Sing |S|$ is always a weak homotopy equivalence, but generally not a categorical equivalence. For example, if $S$ is an $\infty$-category, then $u$ is a categorical equivalence if and only if $S$ is a Kan complex. In general, $\Sing |S|$ may be regarded as the $\infty$-groupoid
obtained from $S$ by freely adjoining inverses to all the
morphisms in $S$.

\begin{remark}
The inclusion functor $i$ and its homotopy-theoretic left adjoint
may also be understood using the formalism of {\it
localizations of model categories}. In addition to its usual model
category structure, the category $\sSet$ of simplicial sets may be
endowed with the {\it Joyal model structure} which we will define in \S
\ref{compp3}. These model structures have the same cofibrations (in both cases, the
cofibrations are simply the monomorphisms of simplicial sets).
However, the Joyal model structure has fewer weak equivalences
(categorical equivalences, rather than weak homotopy equivalences)
and consequently more fibrant objects (all $\infty$-categories,
rather than only Kan complexes). It follows that the usual
homotopy theory of simplicial sets is a
localization of the homotopy theory of $\infty$-categories. The
identity functor from $\sSet$ to itself determines a Quillen
adjunction between these two homotopy theories, which plays the
role of $i$ and its left adjoint.
\end{remark}

\subsection{Homotopy Commutativity versus Homotopy Coherence}\label{comcoh}

Let $\calC$ be an $\infty$-category (topological category, simplicial category). 
To a first approximation,
working in $\calC$ is like working in its homotopy category $\h{
\calC}$: up to equivalence, $\calC$ and $\h{\calC}$ have the same
objects and morphisms. The main difference between $\h{\calC}$ and
$\calC$ is that in $\calC$, one must not ask whether or not
morphisms are {\em equal}; instead one should ask whether or not they are {\it homotopic}. If so, the homotopy itself is an additional datum which we will need to consider. Consequently, the notion of a commutative diagram in
$\h{\calC}$, which corresponds to a {\it homotopy commutative}
diagram in $\calC$, is quite unnatural and usually needs to be
replaced by the more refined notion of a {\it homotopy coherent}
diagram in $\calC$.\index{gen}{homotopy coherence}\index{gen}{diagram!homotopy commutative}\index{gen}{diagram!homotopy coherent}

To understand the problem, let us suppose that $F: \calI
\rightarrow \calH$ is a functor from an ordinary category $\calI$
into the homotopy category of spaces $\calH$. In other words, $F$
assigns to each object $X \in \calI$ a space (say, a CW complex)
$F(X)$, and to each morphism $\phi: X \rightarrow Y$ in $\calI$ a
continuous map of spaces $F(\phi): F(X) \rightarrow F(Y)$ (well-defined
up to homotopy), such that $F(\phi \circ \psi)$ is homotopic to
$F(\phi) \circ F(\psi)$ for any pair of composable morphisms $\phi,
\psi$ in $\calI$. In this situation, it may or may not be possible
to {\em lift} $F$ to an actual functor $\widetilde{F}$ from
$\calI$ to the ordinary category of topological spaces, such that
$\widetilde{F}$ induces a functor $\calI \rightarrow \calH$ which
is naturally isomorphic to $F$. In general there are obstructions
to both the existence and the uniqueness of the lifting
$\widetilde{F}$, even up to homotopy. To see this, let us suppose for a moment that
$\widetilde{F}$ exists, so that there exist homotopies
$k_{\phi}: \widetilde{F}(\phi) \simeq F(\phi)$. These homotopies determine {\em additional}
data on $F$: namely, one obtains a canonical homotopy $h_{\phi,\psi}$ from $F(\phi \circ \psi)$
to $F(\phi) \circ F(\psi)$ by composing
$$ F(\phi \circ \psi) \simeq \widetilde{F}(\phi \circ \psi) = \widetilde{F}(\phi) \circ \widetilde{F}(\psi)
\simeq F(\phi) \circ F(\psi).$$
The functor $F$ to the homotopy category
$\calH$ should be viewed as a first approximation to $\widetilde{F}$; we obtain a
second approximation when we take into account the homotopies
$h_{\phi, \psi}$. These homotopies are not arbitrary: the
associativity of composition gives a relationship between
$h_{\phi, \psi}, h_{\psi, \theta}, h_{\phi, \psi \circ \theta}$
and $h_{\phi \circ \psi, \theta}$, for a composable triple of
morphisms $(\phi, \psi, \theta)$ in $\calI$. This relationship may
be formulated in terms of the existence of a certain higher
homotopy, which is once again canonically determined by
$\widetilde{F}$ (and the homotopies $k_{\phi}$). To obtain the next approximation to
$\widetilde{F}$, we should take these higher homotopies into
account, and formulate the associativity properties that {\em
they} enjoy, and so on. Roughly speaking, a {\it homotopy coherent} diagram in
$\calC$ is a functor $F: \calI \rightarrow
\h{\calC}$, together with all of the extra data that would be
available if we were able to lift $F$ to a functor $\widetilde{F}: \calI \rightarrow \calC$.

The distinction between homotopy commutativity and homotopy coherence is arguably the {\em main} difficulty in working with higher categories. The idea of homotopy coherence is simple enough, and can be made precise in the setting of a general topological category. However, the amount of data required to specify a homotopy coherent diagram is considerable,
so the concept is quite difficult to employ in practical situations.

\begin{remark}
Let $\calI$ be an ordinary category and $\calC$ a topological category. Any functor
$F: \calI \rightarrow \calC$ determines a homotopy coherent diagram in $\calC$ (with all of the homotopies involved being constant). For many topological categories $\calC$, the converse fails: not every homotopy-coherent diagram in $\calC$ can be obtained in this way, even up to equivalence. In these cases, it is the notion of {\it homotopy coherent} diagram which is fundamental; a homotopy coherent diagram should be regarded as ``just as good'' as a strictly commutative diagram, for $\infty$-categorical purposes. As evidence for this, we remark that given an equivalence $\calC' \rightarrow \calC$, a strictly commutative diagram $F: \calI \rightarrow \calC$ cannot always be lifted to a strictly commutative diagram in $\calC'$; however it can always be lifted (up to equivalence) to a homotopy coherent diagram in $\calC'$.
\end{remark}

One of the advantages of working with $\infty$-categories is that the
definition of a homotopy coherent diagram is easy to formulate. We can simply define a homotopy coherent diagram in an $\infty$-category $\calC$ to
be a map of simplicial sets $f: \Nerve(\calI) \rightarrow \calC$. The
restriction of $f$ to simplices of low dimension encodes the
induced map on homotopy categories. Specifying $f$ on
higher-dimensional simplices gives precisely the ``coherence
data'' that the above discussion calls for.

\begin{remark}\label{remmt}
Another possible approach to the problem of homotopy coherence is
to restrict our attention to simplicial (or topological) categories $\calC$ in
which every homotopy coherent diagram is equivalent to a strictly commutative diagram. For example, this is always true when $\calC$ arises from a simplicial model category (Proposition \ref{gumby444}). Consequently, in the framework of model categories it is possible to ignore the theory of homotopy coherent diagrams, and work with strictly commutative diagrams instead. This approach is quite powerful, particularly when combined with the observation that every simplicial category $\calC$
admits a fully faithful embedding into a simplicial model category (for example, one can use
a simplicially enriched version of the Yoneda embedding). This idea can be used to show that
every homotopy coherent diagram in $\calC$ can be ``straightened'' to a commutative diagram, possibly after replacing $\calC$ by an equivalent simplicial category (for a more precise version
of this statement, we refer the reader to Corollary \ref{strictify}).
\end{remark}

\subsection{Functors between Higher Categories}\label{funcback}

The notion of a homotopy coherent diagram in an higher category
$\calC$ is a special case of the more general notion of a functor
$F: \calI \rightarrow \calC$ between higher categories
(specifically, it is the special case in which $\calI$ is assumed
to be an ordinary category). Just as the collection of all
ordinary categories forms a bicategory (with functors as
morphisms and natural transformations as $2$-morphisms), the
collection of all $\infty$-categories can
be organized into an $\infty$-bicategory. In particular, for any
$\infty$-categories $\calC$ and $\calC'$, we expect to be able to construct
an $\infty$-category $\Fun(\calC,\calC')$ of functors from $\calC$ to
$\calC'$.

In the setting of topological categories, the construction of an appropriate mapping object $\Fun(\calC, \calC')$ is quite difficult. The naive guess is that $\Fun(\calC, \calC')$ should be 
a category of topological functors from $\calC$ to $\calC'$: that is, functors which induce continuous maps between morphism spaces. However, we saw in \S \ref{comcoh} that this notion is generally too rigid, even in the special case where $\calC$ is an ordinary category.

\begin{remark}
Using the language of model categories, one might say that the
problem is that not every topological category is {\it cofibrant}.
If $\calC$ is a ``cofibrant'' topological category (for example, if $\calC =
|\sCoNerve[S]|$ where $S$ is a simplicial set), then the collection of topological
functors from $\calC$ to $\calC'$ is large enough to contain representatives for
every $\infty$-categorical functor from $\calC$ to $\calC'$. Most
ordinary categories are not cofibrant when viewed as topological categories.
More importantly, the property of being cofibrant is not stable under products, so that naive attempts to construct a mapping object $\Fun(\calC, \calC')$ need not give the correct answer even when $\calC$ itself is assumed cofibrant (if $\calC$ is cofibrant, then we are guaranteed to have ``enough'' topological functors $\calC \rightarrow \calC'$ to represent all functors between the underlying $\infty$-categories, but not necessarily enough natural transformations between them; note that the product $\calC \times [1]$ is usually not cofibrant, even in the simplest nontrivial case where $\calC = [1]$.) This is arguably the most important technical disadvantage of the theory of topological (or simplicial) categories as an approach to higher category theory.
\end{remark}

The construction of functor categories is much easier to describe in the framework of $\infty$-categories. If $\calC$ and $\calD$ are
$\infty$-categories, then we can simply define a {\it functor} from
$\calC$ to $\calD$ to be a map $p: \calC \rightarrow \calD$ of simplicial sets.\index{gen}{functor!between $\infty$-categories}

\begin{notation}\index{not}{Fun@$\Fun(\calC, \calC')$}
Let $\calC$ and $\calD$ be simplicial sets. We let $\Fun(\calC, \calD)$ denote the
simplicial set $\bHom_{\sSet}(\calC, \calD)$ parametrizing maps from $\calC$ to $\calD$.
We will use this notation only when $\calD$ is an $\infty$-category (the simplicial set $\calC$ will often, but not always, be an $\infty$-category as well). We will refer to $\Fun(\calC, \calD)$ as the
{\it $\infty$-category of functors from $\calC$ to $\calD$} (see Proposition \ref{tyty} below).
We will refer to morphisms in $\Fun(\calC, \calD)$ as {\it natural transformations} of functors, and
equivalences in $\Fun(\calC, \calD)$ as {\it natural equivalences}.\index{gen}{natural transformation}\index{gen}{natural equivalence}
\end{notation}

\begin{proposition}\label{tyty}
Let $K$ be an arbitrary simplicial set.
\begin{itemize}
\item[$(1)$] For every $\infty$-category $\calC$, the simplicial set $\Fun(K,\calC)$ is an $\infty$-category.

\item[$(2)$] Let $\calC \rightarrow \calD$ be a categorical equivalence of $\infty$-categories. Then the induced map $\Fun(K,\calC) \rightarrow \Fun(K,\calD)$ is a categorical equivalence.

\item[$(3)$] Let $\calC$ be an $\infty$-category, and $K \rightarrow K'$ a categorical equivalence of simplicial sets. Then the induced map $\Fun(K',\calC) \rightarrow \Fun(K,\calC)$ is a categorical equivalence.
\end{itemize}
\end{proposition}

The proof makes use of the Joyal model structure on $\sSet$, and will be given in \S \ref{compp3}.

\subsection{Joins of $\infty$-Categories}\label{join}

Let $\calC$ and $\calC'$ be ordinary categories. We will define a
new category $\calC \join \calC'$, called the {\it join} of $\calC$ and
$\calC'$. An object of $\calC \join \calC'$ is either an object of
$\calC$ or an object of $\calC'$. The morphism sets are given as follows:
 $$\Hom_{\calC \join \calC'}(X,Y) = \begin{cases} \Hom_{\calC}(X,Y) & \text{if } X,Y \in \calC \\
\Hom_{\calC'}(X,Y) & \text{if } X,Y \in \calC' \\
\emptyset & \text{if } X \in \calC', Y \in \calC \\
\ast & \text{if } X \in \calC, Y \in \calC'. \end{cases}$$\index{gen}{join!of categories}
Composition of morphisms in $\calC \join \calC'$ is defined in the
obvious way. 

The join construction described above is often useful when discussing diagram categories, limits, and colimits. In this section, we will introduce a generalization of this construction to the $\infty$-categorical setting.

\begin{definition}
If $S$ and $S'$ are
simplicial sets, then the simplicial set $S \star S'$\index{not}{Star@$S \star S'$} is defined as
follows: for each nonempty finite linearly ordered set $J$, we set
$$(S \star S')(J) = \coprod_{J = I \cup I'} S(I) \times
S'(I'),$$ where the union is taken over all decompositions of $J$ into disjoint subsets $I$ and $I'$, satisfying $i < i'$ for all $i \in I$, $i' \in I'$. Here we allow the
possibility that either $I$ or $I'$ is empty, in which case we agree to
the convention that $S(\emptyset) = S'(\emptyset) = \ast$.\index{gen}{join!of simplicial sets}
\end{definition}

More concretely, we have $$(S \star S')_{n} =
S_n \cup S'_n \cup \bigcup_{i+j = n-1} S_i \times S'_j.$$

The join operation endows $\sSet$ with the
structure of a monoidal category (see \S \ref{monoidaldef}).
The identity for the join operation is
the empty simplicial set $\emptyset = \Delta^{-1}$. More generally, we have
natural isomorphisms $\phi_{ij}: \Delta^{i-1} \star \Delta^{j-1} \simeq
\Delta^{(i+j)-1}$, for all $i, j \geq 0$.

\begin{remark}
The operation $\star$ is essentially determined by the isomorphisms
$\phi_{ij}$, together with its behavior under the formation of
colimits: for any fixed simplicial set $S$, the functors
$$ T \mapsto T \star S$$
$$ T \mapsto S \star T$$
commute with colimits, when regarded as functors from $\sSet$ to
the undercategory $(\sSet)_{S/}$ of simplicial sets {\em under} $S$.
\end{remark}

Passage to the nerve carries joins of
categories into joins of simplicial sets. More precisely, for every pair of
categories $\calC$ and $\calC'$, there is
a canonical isomorphism $$\Nerve( \calC \join \calC') \simeq
\Nerve(\calC) \join \Nerve(\calC').$$ (The existence of this
isomorphism persists when we allow $\calC$ and $\calC'$ to be a simplicial or
topological categories and apply the appropriate generalization of
the nerve functor.) This suggests that the join operation on
simplicial sets is the appropriate $\infty$-categorical analogue of
the join operation on categories.

We remark that the formation of joins does not commute with the
functor $\sCoNerve[\bigdot]$. However, the simplicial category $\sCoNerve[S \star S']$
contains $\sCoNerve[S]$ and $\sCoNerve[S']$ as full (topological)
subcategories, and contains no morphisms from objects of
$\sCoNerve[S']$ to objects of $\sCoNerve[S]$. Consequently, there is unique map $\phi: \sCoNerve[S \star S'] \rightarrow \sCoNerve[S] \star \sCoNerve[S']$ which reduces to the identity on
$\sCoNerve[S]$ and $\sCoNerve[S']$. We will later show that $\phi$ is an equivalence of simplicial categories (Corollary \ref{diamond3}).

We conclude by recording a pleasant property of the join
operation:

\begin{proposition}[Joyal \cite{joyalnotpub}]
If $S$ and $S'$ are $\infty$-categories, then $S \star S'$ is an $\infty$-category.
\end{proposition}

\begin{proof}
Let $p: \Lambda^n_i \rightarrow S \star S'$ be a map, with $0 < i
< n$. If $p$ carries $\Lambda^n_i$ entirely into $S \subseteq S
\star S'$ or into $S' \subseteq S \star S'$, then we deduce
the existence an extension of $p$ to $\Delta^n$ by invoking
the assumption that $S$ and $S'$ are $\infty$-categories. Otherwise,
we may suppose that $p$ carries the vertices $\{0, \ldots, j\}$
into $S$, and the vertices $\{ j+1, \ldots, n\}$ into $S'$. 
We may now restrict $p$ to obtain maps
$$ \Delta^{ \{ 0, \ldots, j \} } \rightarrow S$$
$$ \Delta^{ \{ j+1, \ldots, n \} } \rightarrow S',$$
which together determine a map $\Delta^n \rightarrow S \star S'$ extending $p$.
\end{proof}

\begin{notation}
Let $K$ be a simplicial set. The {\it left cone} $K^{\triangleleft}$ is defined to be
the join $\Delta^0 \join K$. Dually, the {\it right cone} $K^{\triangleright}$ is defined to be the join $K \join \Delta^0$. Either cone contains a distinguished vertex (belonging to $\Delta^0$), which we will refer to as the {\it cone point}.\index{gen}{cone point}\index{not}{Kleft@$K^{\triangleleft}$}\index{not}{Kright@$K^{\triangleright}$}
\end{notation}

\subsection{Overcategories and Undercategories}\label{slices}

Let $\calC$ be an ordinary category, and $X \in \calC$ an object.
The {\it overcategory} $\calC_{/X}$ is defined as follows:\index{gen}{overcategory}
the objects of $\calC_{/X}$ are
morphisms $Y \rightarrow X$ in $\calC$ having target $X$.
Morphisms are given by commutative triangles
$$\xymatrix{ Y \ar[dr] \ar[rr] & & Z \ar[dl] \\
& X }$$
and composition is defined in the obvious way.

One can rephrase the definition of the overcategory as follows.
Let $[0]$ denote the category with a single object, possessing
only an identity morphism. Then specifying an object $X \in \calC$
is equivalent to specifying a functor $x: [0] \rightarrow
\calC$. The overcategory $\calC_{/X}$ may then be described by
the following universal property: for any category $\calC'$, we
have a bijection
$$ \Hom(\calC', \calC_{/X}) \simeq \Hom_{x}(\calC' \join [0],
\calC),$$ where the subscript on the right hand side indicates
that we consider only those functors $\calC' \join [0]
\rightarrow \calC$ whose restriction to $[0]$ coincides with
$x$.

We would like to generalize the construction of overcategories to the $\infty$-categorical setting. 
Let us begin by working in the framework of topological categories. In this case, there is a natural candidate for the relevant overcategory. Namely, if $\calC$ is a topological category containing an object $X$, then the overcategory $\calC_{/X}$ (defined as above) has the structure of a topological category, where each morphism space $\bHom_{ \calC_{/X} }(Y,Z)$ is topologized as a subspace
of $\bHom_{\calC}(Y,Z)$ (here we are identifying an object of $\calC_{/X}$ with its image in $\calC$). This topological category is usually {\em
not} a model for the correct $\infty$-categorical slice
construction. The problem is that a morphism in $\calC_{/X}$ consists
of a commutative triangle
$$\xymatrix{ Y \ar[dr] \ar[rr] & & Z \ar[dl] \\
& X }$$
of objects over $X$. To obtain the correct notion, we should allow also
triangles which commute only up to homotopy. 

\begin{remark}
In some cases, the naive overcategory $\calC_{/X}$ is a good approximation
to the correct construction: see Lemma \ref{tulmand}.
\end{remark}

In the setting of $\infty$-categories, Joyal has given a much simpler description of the desired construction (see \cite{joyalpub}). This description will play a vitally important role throughout this book. We begin by noting the following:

\begin{proposition}[\cite{joyalpub}]\label{tyrii}
Let $S$ and $K$ be simplicial sets, and $p: K \rightarrow S$ an
arbitrary map. There exists a simplicial set $S_{/p}$ with the following universal property:
$$\Hom_{\sSet}(Y, S_{/p}) = \Hom_{p}( Y \join K, S),$$
where the subscript on the right hand side indicates that we
consider only those morphisms $f: Y \join K \rightarrow S$ such
that $f|K = p$.
\end{proposition}

\begin{proof}
One defines $(S_{/p})_n$ to be $\Hom_{p}(\Delta^n \join K, S)$.
The universal property holds by definition when $Y$ is a simplex.
It holds in general because both sides are compatible with the
formation of colimits in $Y$.
\end{proof}

Let $p: K \rightarrow S$ be as in Proposition \ref{tyrii}. If $S$ is an $\infty$-category, we will refer to $S_{/p}$ as an {\it overcategory} of $S$, or as the {\it $\infty$-category of objects of $S$ over $p$}.
The following result guarantees that the operation of passing to overcategories is well-behaved:\index{gen}{overcategory!of an $\infty$-category}\index{not}{calC/p@$\calC_{/p}$}

\begin{proposition}\label{gorban3}
Let $p: K \rightarrow \calC$ be a map of simplicial sets, and suppose
that $\calC$ is an $\infty$-category. Then $\calC_{/p}$ is an $\infty$-category.
Moreover, if $q: \calC \rightarrow \calC'$ is a categorical equivalence of
$\infty$-categories, then the induced map $\calC_{/p} \rightarrow \calC'_{/qp}$ is a categorical equivalence
as well.
\end{proposition}

The proof requires a number of ideas which we have not yet introduced, and will be postponed (see Proposition \ref{gorban4} for the first assertion and \S \ref{slim} for the second).

\begin{remark}
Let $\calC$ be an $\infty$-category.
In the particular case where $p: \Delta^n \rightarrow \calC$ classifies an
$n$-simplex $\sigma \in \calC_n$, we will often write $\calC_{/\sigma}$ in place of of $\calC_{/p}$.
In particular, if $X$ is an object of $\calC$, we let $\calC_{/X}$ denote the overcategory
$\calC_{/p}$, where $p: \Delta^0 \rightarrow \calC$ has image $X$.\index{not}{calC/X@$\calC_{/X}$}
\end{remark}

\begin{remark} Let $p: K \rightarrow \calC$ be a map of simplicial sets. The
preceding discussion can be dualized, replacing $Y \star K$
by $K \star Y$; in this case we denote the corresponding simplicial set
by $\calC_{p/}$ which (if $\calC$ is an $\infty$-category) we will refer to as an
{\it undercategory} of $\calC$. In the special case where $K = \Delta^n$ and $p$ classifies a simplex $\sigma \in \calC_{n}$, we will also write $\calC_{\sigma/}$ for $\calC_{p/}$; in particular, we will write $\calC_{X/}$ when $X$ is an object of $\calC$.\index{not}{calCp/@$\calC_{p/}$}\index{not}{calCX/@$\calC_{X/}$}\index{gen}{undercategory!of an $\infty$-category}
\end{remark}

\begin{remark}
If $\calC$ is an ordinary category and $X \in \calC$,
then there is a canonical isomorphism $\Nerve(\calC)_{/X} \simeq \Nerve (\calC_{/X})$. In other words, the overcategory construction for $\infty$-categories can be regarded as a {\em generalization} of the relevant construction from classical category theory.
\end{remark}

\subsection{Fully Faithful and Essentially Surjective Functors}

\begin{definition}\label{faulfa}\index{gen}{essentially surjective}
Let $F: \calC \rightarrow \calD$ be a functor between topological categories (simplicial categories, simplicial sets). We will say that $F$ is {\it essentially surjective} if the induced functor 
$\h{F}: \h{\calC} \rightarrow \h{\calD}$ is essentially surjective (that is, if every object of $\calD$ is
equivalent to $F(X)$ for some $X \in \calC$). 

We will say that $F$ is {\it fully faithful} if $\h{F}$ is a fully faithful functor of $\calH$-enriched categories. In other words, $F$ is fully faithful if and only if, for every pair of objects
$X, Y \in \calC$, the induced map
$\bHom_{\h{\calC}}(X,Y) \rightarrow \bHom_{\h{\calD}}(F(X),F(Y))$ is an isomorphism in the homotopy category $\calH$. \index{gen}{fully faithful}\index{gen}{functor!fully faithful}
\end{definition}

\begin{remark}
Because Definition \ref{faulfa} makes reference only to the homotopy categories of $\calC$ and $\calD$, it is invariant under equivalence and under operations which pass between
the various models for higher category theory that we have introduced.
\end{remark}

Just as in ordinary category theory, a functor $F$ is an equivalence if and only if it is fully faithful and essentially surjective.

\subsection{Subcategories of $\infty$-Categories}

Let $\calC$ be an $\infty$-category, and let $(\h{\calC})' \subseteq \h{\calC}$ be a subcategory of its homotopy category. We can then form a pullback diagram of simplicial sets
$$ \xymatrix{ \calC' \ar[r] \ar[d] & \calC \ar[d] \\
\Nerve (\h{\calC})' \ar[r] & \Nerve (\h{ \calC}). }$$
We will refer to $\calC'$ as the {\it subcategory of $\calC$ spanned by $(\h{\calC})'$}. In general, we will say that a simplicial subset $\calC' \subseteq \calC$ is a {\it subcategory} of $\calC$ if it arises via this construction.\index{gen}{subcategory!of an $\infty$-category}

\begin{remark}
We say ``subcategory'', rather than ``sub-$\infty$-category'', in
order to avoid awkward language. The terminology is not meant to
suggest that $\calC'$ is itself a category, or isomorphic to the nerve
of a category.
\end{remark}

In the case where $(\h{\calC})'$ is a full subcategory of $\h{\calC}$, we will say that
$\calC'$ is a {\it full subcategory} of $\calC$. In this case, $\calC'$ is determined by the set
$\calC'_0$ of those objects $X \in \calC$ which belong to $\calC'$. We will then say that
$\calC'$ is the {\it full subcategory of $\calC$ spanned by $\calC'_0$}.\index{gen}{subcategory!full}

It follows from Remark \ref{needie} that
the inclusion $\calC' \subseteq \calC$ is fully faithful. In general, any fully faithful functor $f: \calC'' \rightarrow \calC$ factors as a composition
$$ \calC'' \stackrel{f'}{\rightarrow} \calC' \stackrel{f''}{\rightarrow} \calC,$$
where $f'$ is an equivalence of $\infty$-categories and $f''$ is the inclusion of the full subcategory
$\calC' \subseteq \calC$ spanned by the set of objects $f( \calC''_0 ) \subseteq \calC_0$.

\subsection{Initial and Final Objects}

If $\calC$ is an ordinary category, then an object $X \in \calC$
is said to be {\it final} if $\Hom_{\calC}(Y,X)$ consists of a single element, for every $Y \in \calC$. Dually, an object $X \in \calC$ is {\it initial} if it is final when viewed as an object of $\calC^{op}$.
The goal of this section is to generalize these definitions to the $\infty$-categorical setting.\index{gen}{final object!of a category}\index{gen}{object!final}

If $\calC$ is a topological category, then a candidate definition immediately presents itself: we could ignore the topology on the morphism spaces, and consider those objects of $\calC$ which are final when $\calC$ is regarded as an ordinary category. This requirement is unnaturally strong. For example, the category $\CG$ of compactly generated Hausdorff spaces has a final object: the topological space $\ast$, consisting of a single point. However, there are objects of $\CG$ which are equivalent to $\ast$ (any contractible space) but
not isomorphic to $\ast$ (and therefore not final objects of $\CG$, at least in the classical sense). Since any reasonable $\infty$-categorical notion is stable under equivalence, we need to find a weaker condition.

\begin{definition}\label{inuy}\index{gen}{final object!of an $\infty$-category}\index{gen}{object!final}
Let $\calC$ be a topological category (simplicial category, simplicial set). An object $X \in \calC$ is {\it final} if it is final in the homotopy category
$\h{\calC}$, regarded as a category enriched over $\calH$. In other words, $X$ is final if and only if
for each $Y \in \calC$, the mapping space $\bHom_{ \h{\calC}}(Y,X)$ is weakly contractible
(that is, a final object of $\calH$).
\end{definition}

\begin{remark}
Since the Definition \ref{inuy} makes reference only to the homotopy category $\h{\calC}$, it is invariant under equivalence and under passing between the various models for higher category theory.
\end{remark}

In the setting of $\infty$-categories, it is convenient to employ a slightly more
sophisticated definition, which we borrow from \cite{joyalpub}.

\begin{definition}\label{strongfin}\index{gen}{strongly final}
Let $\calC$ be a simplicial set. A
vertex $X$ of $\calC$ is {\it strongly final} if the projection $\calC_{/X}
\rightarrow \calC$ is a trivial fibration of simplicial sets. 
\end{definition}

In other words, a vertex $X$ of $\calC$ is strongly final if and only if any
map $f_0: \bd \Delta^n \rightarrow \calC$ such that
$f_0(n) = X$ can be extended to a map $f: \Delta^n \rightarrow S$.

\begin{proposition}\label{harry}
Let $\calC$ be an $\infty$-category containing an object $Y$. The object
$Y$ is strongly final if and only if, for every object $X \in \calC$, the Kan
complex $\Hom^{\rght}_{\calC}(X,Y)$ is contractible.
\end{proposition}

\begin{proof}
The ``only if'' direction is clear: the space
$\Hom^{\rght}_{\calC}(X,Y)$ is the fiber of the projection $p: \calC_{/Y}
\rightarrow \calC$ over the vertex $X$. If $p$ is a trivial
fibration, then the fiber is a contractible Kan complex. Since $p$
is a right fibration (Proposition \ref{sharpen}), the converse holds as well (Lemma \ref{toothie}).
\end{proof}

\begin{corollary}
Let $\calC$ be a simplicial set. Every strongly final object of $\calC$ is a final object of $\calC$. The converse holds if $\calC$ is an $\infty$-category.
\end{corollary}

\begin{proof}
Let $[0]$ denote the category with a single object and a single morphism.
Suppose that $Y$ is a strongly final vertex of $\calC$. Then there exists a retraction of
$\calC^{\triangleright}$ onto $\calC$, carrying the cone point to $Y$. Consequently, we obtain a retraction of ($\calH$-enriched) homotopy categories from $(\h{\calC}) \star [0]$ to $\h{\calC}$, carrying the unique object of $[0]$ to $Y$. This implies that $Y$ is final in $\h{\calC}$, so that $Y$ is a final object of $\calC$.

To prove the converse, we note that if $\calC$ is an $\infty$-category then $\Hom_{\calC}^{\rght}(X,Y)$ represents the homotopy type $\bHom_{\calC}(X,Y) \in \calH$; by Proposition \ref{harry} this space is contractible for all $X$ if and only if $Y$ is strongly final.
\end{proof}

\begin{remark}
The above discussion dualizes in an evident way, so that we have a notion of {\em initial} objects of an $\infty$-category $\calC$.
\end{remark}

\begin{example}
Let $\calC$ be an ordinary category containing an object $X$. Then $X$ is a final (initial)
object of the $\infty$-category $\Nerve(\calC)$ if and only if it is a final (initial) object of $\calC$, in the usual sense.
\end{example}

\begin{remark}
Definition \ref{strongfin} is only natural in the case where $\calC$ is an $\infty$-category. For example, if $\calC$ is not an $\infty$-category, then the collection of strongly final vertices of $\calC$ need not be stable under equivalence.
\end{remark}

An ordinary category $\calC$ may have more than one final object,
but any two final objects are uniquely isomorphic to one another.
In the setting of $\infty$-categories, an analogous statement holds,
but is slightly more complicated because the word ``unique'' needs to be
interpreted in a homotopy theoretic sense:

\begin{proposition}[Joyal]\label{initunique}\index{gen}{final object!uniqueness}
Let $\calC$ be a $\infty$-category, and let $\calC'$ be the full subcategory
of $\calC$ spanned by the final vertices of $\calC$. Then $\calC'$ is either empty or
a contractible Kan complex.
\end{proposition}

\begin{proof}
We wish to prove that every map $p: \bd \Delta^n \rightarrow \calC'$
can be extended to an $n$-simplex of $\calC'$. If $n = 0$, this is
possible unless $\calC'$ is empty. For $n > 0$, the desired extension exists
because $p$ carries the $n$th vertex of $\bd \Delta^n$ to a final
object of $\calC$.
\end{proof}

\subsection{Limits and Colimits}\label{limitcolimit}

An important consequence of the distinction between homotopy
commutativity and homotopy coherence is that the appropriate
notions of limit and colimit in a higher category
$\calC$ do not coincide with the notion of limit and colimit in the homotopy category $\h{\calC}$ (where limits and colimits often do
not exist). Limits and colimits in
$\calC$ are often referred to as {\it homotopy limits} and
{\it homotopy colimits}, to avoid confusing them with ordinary limits
and colimits.

Homotopy limits and
colimits can be defined in a topological category, but the
definition is rather complicated. We will review a few special cases here, and discuss the general definition in the appendix (\S \ref{qlim7}).

\begin{example}\label{examprod}\index{gen}{homotopy product}\index{gen}{product!homotopy}
Let $\{ X_{\alpha} \}$ be a family of objects in a topological
category $\calC$. A {\it homotopy product} $X = \prod_{\alpha}
X_{\alpha}$ is an object of $\calC$ equipped with morphisms
$f_{\alpha}: X \rightarrow X_{\alpha}$ which induce a weak
homotopy equivalence
$$ \bHom_{\calC}(Y,X) \rightarrow \prod_{\alpha} \bHom_{\calC}(Y,
X_{\alpha})$$ for every object $Y \in \calC$.

Passing to path components and using the fact that $\pi_0$
commutes with products, we deduce that $$\Hom_{\h{\calC}}(Y,X) \simeq
\prod_{\alpha} \Hom_{\h{\calC}}(Y, X_{\alpha}),$$ so that any product in $\calC$ is
also a product in $\h{\calC}$. In particular, the object $X$ is
determined up to canonical isomorphism in $\h{\calC}$.

In the special case where the index set is empty, we recover the
notion of a final object of $\calC$: an object $X$ for which each
of the mapping spaces $\bHom_{\calC}(Y,X)$ is weakly contractible.
\end{example}

\begin{example}\label{exampull}\index{gen}{homotopy pullback}\index{gen}{pullback!homotopy}
Given two morphisms $\pi: X \rightarrow Z$ and $\psi: Y \rightarrow Z$
in a topological category $\calC$, let us define $\bHom_{\calC}(W,
X \times^h_Z Y)$ to be the space consisting of points $p \in
\bHom_{\calC}(W,X)$, $q \in \bHom_{\calC}(W,Y)$, together with a path $r: [0,1]
\rightarrow \bHom_{\calC}(W,Z)$ joining $\pi \circ p$ to $\psi \circ q$. We
endow $\bHom_{\calC}(W, X \times^h_Z Y)$ with the obvious topology,
so that $X \times^h_Z Y$ can be viewed presheaf of topological spaces
on $\calC$. A {\it homotopy fiber product for $X$ and $Y$ over
$Z$} is an object of $\calC$ which represents this presheaf, up to
weak homotopy equivalence. In other words, it is an object $P \in \calC$
equipped with a point $p \in \bHom_{\calC}(P, X \times^h_Z Y)$ which
induces weak homotopy equivalences $\bHom_{\calC}(W,P) \rightarrow
\bHom_{\calC}(W, X \times^h_Z Y)$ for every $W \in \calC$.

We note that, if there exists a fiber product (in the ordinary sense) $X \times_Z Y$ in the category
$\calC$, then this ordinary fiber product admits a (canonically determined) map to the homotopy fiber product (if the homotopy fiber product exists). This map need not be an equivalence, but it is an equivalence in many good cases. We also note that a homotopy fiber product $P$ comes equipped with a map
to the fiber product $X \times_Z Y$ taken in the category $\h{\calC}$ (if this fiber product exists); this map is usually not an isomorphism.
\end{example}

\begin{remark}
Homotopy limits and colimits in general may be described in
relation to homotopy limits of topological spaces. The homotopy
limit $X$ of a diagram of objects $\{X_{\alpha} \}$ in an
arbitrary topological category $\calC$ is determined, up to
equivalence, by the condition that there exist a natural weak homotopy
equivalence
$$\bHom_{\calC}(Y,X) \simeq \holim \{ \bHom_{\calC}(Y, X_{\alpha})
\}.$$ Similarly, the homotopy colimit of the diagram is characterized by
the existence of a natural weak homotopy equivalence
$$\bHom_{\calC}(X,Y) \simeq \holim \{ \bHom_{\calC}(X_{\alpha},Y)
\}.$$
For a more precise discussion, we refer the reader to Remark \ref{curble}.
\end{remark}

In the setting of $\infty$-categories, limits and colimits
are quite easy to define:

\begin{definition}[Joyal \cite{joyalpub}]\label{defcolim}\index{gen}{colimit}\index{gen}{limit}
Let $\calC$ be an $\infty$-category and let $p: K \rightarrow \calC$ be an
arbitrary map of simplicial sets. A {\it colimit} for $p$ is
an initial object of $\calC_{p/}$ and a {\it limit} for $p$ is a final
object of $\calC_{/p}$.
\end{definition}

\begin{remark}
According to Definition \ref{defcolim}, a colimit of a diagram $p: K \rightarrow \calC$
is an object of $\calC_{p/}$. We may identify this object with a map
$\overline{p}: K^{\triangleright} \rightarrow \calC$ extending $p$. In general, we will say that a map $\overline{p}: K^{\triangleright} \rightarrow \calC$ is a {\it colimit diagram} if it is a colimit
of $p = \overline{p} | K$. In this case, we will also abuse terminology by referring to
$\overline{p}(\infty) \in \calC$ as a {\it colimit of $p$}, where $\infty$ denotes the cone point of
$K^{\triangleright}$.\index{gen}{colimit!diagram}\index{gen}{limit!diagram}\index{gen}{diagram!(co)limit}
\end{remark}

If $p: K \rightarrow \calC$ is a diagram, we will sometimes
write $\varinjlim(p)$ to denote a colimit of $p$ (considered either as an object of
$\calC_{p/}$ or of $\calC$), and $\varprojlim(p)$ to denote a limit of $p$ (as either an object of
$\calC_{/p}$ or an object of $\calC$). This notation is slightly abusive, since $\varinjlim(p)$ is not uniquely determined by $p$. This phenomenon is familiar in classical category theory: the colimit of a diagram is not unique, but is determined up to canonical isomorphism. In the $\infty$-categorical setting, we have a similar uniqueness result: Proposition \ref{initunique} implies that the collection of candidates for $\varinjlim(p)$, if nonempty, is parametrized by a contractible Kan complex.

\begin{remark}
In \S \ref{quasilimit4}, we will show
that Definition \ref{defcolim} agrees with the classical theory of homotopy (co)limits, when
we specialize to the case where $\calC$ is the nerve of a topological category.
\end{remark}

\begin{remark}
Let $\calC$ be an $\infty$-category, $\calC' \subseteq \calC$ a full subcategory,
and $p: K \rightarrow \calC'$ a diagram. Then $\calC'_{p/} = \calC' \times_{\calC}
\calC_{p/}$. In particular, if $p$ has a colimit in $\calC$, and that
colimit belongs to $\calC'$, then the same object may be regarded as a
colimit for $p$ in $\calC'$.
\end{remark}

Let $f: \calC \rightarrow \calC'$ be a map between $\infty$-categories. Let $p: K
\rightarrow \calC$ be a diagram in $\calC$, having a colimit $x \in \calC_{p/}$.
The image $f(x) \in \calC'_{f p/}$ may or may not be a colimit for the composite map
$f \circ p$. If it is, we will say that $f$ {\it preserves} the colimit of the diagram $p$.\index{gen}{colimit!preservation of}\index{gen}{limits!preservation of}
Often we will apply this terminology not to a particular diagram $p$ but some class of diagrams: for example, we may speak of maps $f$ which
preserve coproducts, pushouts, or filtered colimits (see \S \ref{coexample} for a discussion of special classes of colimits). Similarly, we may ask whether or not a map $f$
preserves the limit of a particular diagram, or various families of diagrams.

We conclude this section by giving a simple example of a colimit-preserving functor.

\begin{proposition}\label{needed17}
Let $\calC$ be an $\infty$-category, $q: T \rightarrow \calC$ and $p: K \rightarrow \calC_{/q}$ two diagrams. Let $p_0$ denote the composition of $p$ with the projection 
$\calC_{/q} \rightarrow \calC$. Suppose that $p_0$ has a colimit in $\calC$. Then:
\begin{itemize}
\item[$(1)$] The diagram $p$ has a colimit in $\calC_{/q}$, and that colimit is preserved by the projection $\calC_{/q} \rightarrow \calC$.

\item[$(2)$] An extension $\widetilde{p}: K^{\triangleright} \rightarrow \calC_{/q}$ is a colimit
of $p$ if and only if the composition
$$K^{\triangleright} \rightarrow \calC_{/q} \rightarrow \calC$$
is a colimit of $p_0$.

\end{itemize}
\end{proposition}

\begin{proof}
We first prove the ``if'' direction of $(2)$. Let $\widetilde{p}: K^{\triangleright} \rightarrow \calC_{/q}$ be such that the composite map $\widetilde{p_0}: K^{\triangleright} \rightarrow \calC$ is a colimit of $p_0$. We wish to show that $\widetilde{p}$ is a colimit of $p$. We may identify $\widetilde{p}$ with a map $K \join \Delta^0 \join T \rightarrow \calC$. For this, it suffices to show that
for any inclusion $A \subseteq B$ of simplicial sets, it is possible to solve the lifting problem depicted in the following diagram:
$$ \xymatrix{ (K \join B \join T) \coprod_{ K \join A \join T} ( K \join \Delta^0 \join A \join T )
\ar@{^{(}->}[d] \ar[r] & \calC \\
K \join \Delta^0 \join B \join T. \ar@{-->}[ur] & }$$
Because $\widetilde{p_0}$ is a colimit of $p_0$, the projection
$$ \calC_{\widetilde{p_0}/} \rightarrow \calC_{p_0/}$$ is a trivial fibration of simplicial
sets and therefore has the right lifting property with respect to the inclusion
$A \join T \subseteq B \join T$.

We now prove $(1)$. Let $\widetilde{p_0}: K^{\triangleright} \rightarrow \calC$ be a colimit of $p_0$.
Since the projection $\calC_{\widetilde{p_0}/} \rightarrow \calC_{p_0/}$ is a trivial fibration, it has the right lifting property with respect to $T$: this guarantees the existence of an extension
$\widetilde{p}: K^{\triangleright} \rightarrow \calC$ lifting $\widetilde{p_0}$. The preceding analysis proves that $\widetilde{p}$ is a colimit of $p$.

Finally, the ``only if'' direction of $(2)$ follows from $(1)$, since any two colimits of $p$ are equivalent.
\end{proof}

\subsection{Presentations of $\infty$-Categories}

Like many other types of mathematical structures, $\infty$-categories can be described by generators and relations. In particular, it makes sense to speak of a {\it
finitely presented} $\infty$-category $\calC$. Roughly speaking, $\calC$ is finitely presented
if it has finitely many objects and its morphism spaces are determined
by specifying a finite number of generating morphisms, a finite
number of relations among these generating morphisms, a finite
number of relations among the relations, and so forth (a finite
number of relations in all).

\begin{example}\label{infinitemorphisms}
Let $\calC$ be the free higher category generated by a single
object $X$ and a single morphism $f: X \rightarrow X$. Then
$\calC$ is a finitely presented $\infty$-category with a single
object, and $\Hom_{\calC}(X,X) = \{ 1, f, f^2, \ldots \}$ is
infinite and discrete. In particular, we note that the finite
presentation of $\calC$ does not guarantee finiteness properties
of the morphism spaces.
\end{example}

\begin{example}
If we identify $\infty$-groupoids with spaces, then giving a
presentation for an $\infty$-groupoid corresponds to giving a cell
decomposition of the associated space. Consequently, the finitely
presented $\infty$-groupoids correspond precisely to the finite
cell complexes.
\end{example}

\begin{example}
Suppose that $\calC$ is a higher category with only two objects
$X$ and $Y$, and that $X$ and $Y$ have contractible endomorphism
spaces and that $\Hom_{\calC}(X,Y)$ is empty. Then $\calC$ is
completely determined by the morphism space $\Hom_{\calC}(Y,X)$,
which may be arbitrary. In this case, $\calC$ is finitely
presented if and only if $\Hom_{\calC}(Y,X)$ is a finite cell
complex (up to homotopy equivalence).
\end{example}

The idea of giving a presentation for an $\infty$-category is very
naturally encoded in the theory of simplicial sets; more
specifically, in Joyal's model structure on $\sSet$, which we will discuss in
\S \ref{compp2}. This model structure can be described as follows:

\begin{itemize}\index{gen}{model category!Joyal}
\item The fibrant objects of $\sSet$ are precisely the
$\infty$-categories.

\item The weak equivalences in $\sSet$ are precisely those maps
$p: S \rightarrow S'$ which induce equivalences $\sCoNerve[S] \rightarrow \sCoNerve[S']$
of simplicial categories.
\end{itemize}

If $S$ is an arbitrary simplicial set, we can
choose a ``fibrant replacement'' for $S$; that is, a categorical
equivalence $S \rightarrow \calC$ where $\calC$ is an $\infty$-category. 
For example, we can take $\calC$ to be the nerve of the topological
category $| \sCoNerve[S] |$. 
The $\infty$-category $\calC$ is
well-defined up to equivalence, and we may
regard it as an $\infty$-category which is ``generated by'' $S$. The simplicial set $S$ itself can be thought of as a ``blueprint'' for building $\calC$. We may view $S$ as generated from the empty (simplicial) set by adjoining nondegenerate simplices. Adjoining a $0$-simplex to $S$ has the effect of adding an object to the $\infty$-category $\calC$, and adjoining a $1$-simplex to $S$ has the effect of adjoining a morphism to $\calC$. Higher dimensional simplices can be thought of as encoding relations among the morphisms.

\subsection{Set-Theoretic Technicalities}

In ordinary category theory, one frequently encounters categories in which the collection of objects
is too large to form a set. Generally speaking, this does not create
any difficulties so long as we avoid doing anything which is obviously illegal
(such as considering the ``category of all categories'' as an object of itself).

The same issues arise in the setting of higher category theory, and are
in some sense even more of a nuisance. In ordinary category
theory, one generally allows a category $\calC$ to have a proper
class of objects, but still requires $\Hom_{\calC}(X,Y)$ to be a
{\em set} for fixed objects $X,Y \in \calC$. The formalism of $\infty$-categories treats
objects and morphisms on the same footing (they are both simplices of a simplicial set), and it is somewhat unnatural (though certainly possible) to directly impose the analogous condition; see \S \ref{locbrend} for a discussion.

There are several means of handling the technical difficulties
inherent in working with large objects (in either classical or higher category theory):

\begin{itemize}
\item[$(1)$] One can employ some set-theoretic device which enables one
to distinguish between ``large'' and ``small''. Examples include:
\begin{itemize}
\item Assuming the existence of a sufficient supply of
(Grothendieck) universes.

\item Working in an axiomatic framework which allows both sets and
{\it classes} (collections of sets which are possibly too large to
themselves be considered sets).

\item Working in a standard set-theoretic framework (such as
Zermelo-Frankel), but incorporating a theory of classes through
some ad-hoc device. For example, one can define a class to be a
collection of sets which is defined by some formula in the
language of set theory.
\end{itemize}

\item[$(2)$] One can work exclusively with ``small'' categories, and
mirror the distinction between ``large'' and ``small'' by keeping
careful track of relative sizes.

 \item[$(3)$] One can simply ignore the set-theoretic difficulties
 inherent in discussing ``large'' categories.

\end{itemize}

Needless to say, approach $(2)$ yields the most refined information. However, it has the disadvantage of burdening our exposition with an additional layer of technicalities. On the other hand, approach $(3)$ will sometimes be inadequate, since we will need to make arguments which play off the distinction between a ``large'' category and a ``small'' subcategory which determines it. Consequently, we shall officially adopt approach $(1)$ for the remainder of this book. More specifically, we assume that for every
cardinal $\kappa_0$, there exists a strongly inaccessible cardinal $\kappa \geq \kappa_0$.
We then let $\calU(\kappa)$ denote the collection of all sets having rank $< \kappa$, so that
$\calU(\kappa)$ is a {\it Grothendieck universe}: in other words, $\calU(\kappa)$ satisfies all of the usual axioms of set theory. We will refer to a mathematical object as {\it small} if it belongs to $\calU(\kappa)$ (or is isomorphic to such an object), and {\it essentially small} if it is equivalent (in whatever relevant sense) to a small object. Whenever it is convenient to do so, we will choose another strongly inaccessible cardinal $\kappa' > \kappa$, to obtain a larger Grothendieck universe
$\calU(\kappa')$ in which $\calU(\kappa)$ becomes small.\index{gen}{small}\index{gen}{essentially small}\index{gen}{Grothendieck universe}

For example, an $\infty$-category $\calC$ is essentially small if and only if it satisfies the following conditions:
\begin{itemize}
\item The set of isomorphism classes of objects in the homotopy
category $\h{\calC}$ has cardinality $< \kappa$.

\item For every morphism $f: X \rightarrow Y$ in $\calC$ and every $i \geq 0$, the homotopy
set $\pi_{i}( \Hom^{\rght}_{\calC}(X,Y), f)$ has cardinality $< \kappa$.
\end{itemize}

For a proof and further discussion, we refer the reader to \S \ref{locbrend}.

\begin{remark}
The existence of the strongly inaccessible cardinal $\kappa$ cannot be proven from the standard axioms of set theory, and the assumption that $\kappa$ exists cannot be proven consistent with the standard axioms for set theory. However, it should be clear that assuming the existence of $\kappa$ is merely the most convenient of the devices mentioned above; none of the results proven in this book will depend on this assumption in an essential way.
\end{remark}

\subsection{The $\infty$-Category of Spaces}\label{introducingspaces}

The category of sets plays a central role in classical category theory. The main reason
is that {\em every} category $\calC$ is enriched over sets: given a pair of objects
$X,Y \in \calC$, we may regard $\Hom_{\calC}(X,Y)$ as an object of $\Set$.
In the higher categorical setting, the proper analogue of $\Set$ is the
$\infty$-category $\SSet$ of {\it spaces}, which we will now
introduce.

\begin{definition}\label{defsset}\index{not}{Kan@$\Kan$}
Let $\Kan$ denote the full subcategory of $\sSet$ spanned by the collection of Kan complexes. 
We will regard $\Kan$ as a simplicial category. Let $\SSet = \Nerve(\Kan)$ denote the (simplicial) nerve of $\Kan$. We will refer to $\SSet$ as the {\it $\infty$-category of spaces}.\index{not}{SSet@$\SSet$}\index{gen}{$\infty$-category!of spaces}
\end{definition}

\begin{remark}
For every pair of objects $X,Y \in \Kan$, the simplicial set
$\bHom_{\Kan}(X,Y) = Y^X$ is a Kan complex. It follows that
$\SSet$ is an $\infty$-category (Proposition \ref{toothy}).
\end{remark}

\begin{remark}
There are many other ways to obtain a suitable ``$\infty$-category of
spaces''. For example, we could instead define $\SSet$ to be the (topological) nerve
of the category of CW-complexes and continuous maps.
All that really matters is that we have a $\infty$-category which is equivalent to $\SSet = \Nerve(\Kan)$. 
We have selected Definition \ref{defsset} for definiteness and to simplify our discussion of the
Yoneda embedding in \S \ref{presheaf1}.
\end{remark}

\begin{remark}
We will occasionally need to distinguish between ``large'' spaces and ``small'' spaces.
In such contexts, we will let $\SSet$ denote the $\infty$-category of small spaces (defined
by taking the simplicial nerve of the category of small Kan complexes), and $\widehat{\SSet}$ the $\infty$-category of large spaces (defined by taking the simplicial nerve of the
category of {\em all} Kan complexes). We observe that $\SSet$ is a large $\infty$-category, and that
$\widehat{\SSet}$ is even bigger.\index{not}{SSethat@$\widehat{\SSet}$}
\end{remark}




\chapter{Fibrations of Simplicial Sets}\label{chap2}

\setcounter{theorem}{0}
\setcounter{subsection}{0}


%

 Many classes of morphisms which play an important role in the homotopy theory of simplicial sets can be defined by their lifting properties (we refer the reader to \S \ref{liftingprobs} for a brief introduction and a summary of the terminology employed below). 

\begin{example} 
A morphism $p: X \rightarrow S$ of simplicial sets which has the right
lifting property with respect to every horn inclusion $\Lambda^n_i \subseteq \Delta^n$ is called a {\it Kan fibration}. A morphism $i: A \rightarrow B$ which has the left lifting property with respect to every Kan fibration is said to be {\it anodyne}.\index{gen}{Kan fibration}\index{gen}{fibration!Kan}\index{gen}{anodyne}
\end{example}

\begin{example}\index{gen}{trivial fibration}\index{gen}{fibration!trivial}
A morphism $p: X \rightarrow S$ of simplicial sets which has the right
lifting property with respect to every inclusion $\bd \Delta^n \subseteq \Delta^n$ is called a {\it trivial fibration}.
A morphism $i: A \rightarrow B$ has the {\em left} lifting property with respect to every trivial Kan fibration if and only if it is a {\it cofibration}; that is, if and only if $i$ is a monomorphism of simplicial sets.\index{gen}{cofibration!of simplicial sets}
\end{example}

By definition, a simplicial set $S$ is a $\infty$-category if it has
the extension property with respect to all horn inclusions
$\Lambda^n_i \subseteq \Delta^n$ with $0 < i < n$. As in classical homotopy theory, it is
convenient to introduce a {\em relative} version of this condition.

\begin{definition}[Joyal]\label{fibdeff}
A morphism $f: X \rightarrow S$ of simplicial sets is
\begin{itemize}
\item a {\it left fibration}\index{gen}{fibration!left} if $f$ has the right lifting property with respect to all horn inclusions $\Lambda^n_i \subseteq \Delta^n$, $0 \leq i < n$.\index{gen}{left fibration}
\item a {\it right fibration}\index{gen}{fibration!right} if $f$ has the right lifting property with respect to all horn inclusions $\Lambda^n_i \subseteq  \Delta^n$, $0 < i \leq n$.\index{gen}{right fibration}
\item an {\it inner fibration}\index{gen}{fibration!inner} if $f$ has the right lifting property with respect to all horn inclusions $\Lambda^n_i \subseteq \Delta^n$, $0 < i < n$.\index{gen}{inner fibration}
\end{itemize}

A morphism of simplicial sets $i: A \rightarrow B$ is
\begin{itemize}
\item {\it left anodyne} if $i$ has the left lifting property with respect to all left fibrations.\index{gen}{anodyne!left}\index{gen}{left anodyne}
\item {\it right anodyne} if $i$ has the left lifting property with respect to all right fibrations.\index{gen}{anodyne!right}\index{gen}{right anodyne}
\item {\it inner anodyne} if $i$ has the left lifting property with respect to all inner fibrations.\index{gen}{anodyne!inner}\index{gen}{inner anodyne}
\end{itemize}
\end{definition}

\begin{remark}
Joyal uses the terms {\it mid-fibration} and {\it mid-anodyne morphism} for what we have called
{\it inner fibrations} and {\it inner anodyne morphisms}.
\end{remark}

The purpose of this chapter is to study the notions of fibration defined above, which are basic tools in the theory of $\infty$-categories. In \S \ref{leftfibsec}, we study the theory of right (left) fibrations $p: X \rightarrow S$, which can be viewed as the $\infty$-categorical analogue of {\it categories
(co)fibered in groupoids} over $S$. We will apply these ideas in \S \ref{valencequi} to show
that the theory of $\infty$-categories is equivalent to the theory of simplicial categories.

There is also an analogue of the more general theory of (co)fibered categories, whose fibers
are not necessarily groupoids: this is the theory of {\it (co)Cartesian fibrations}, which we will introduce in \S \ref{cartfibsec}. Cartesian and coCartesian fibrations are both examples of inner fibrations, which we will study in \S \ref{midfibsec}.

\begin{remark}
To help orient the reader, we summarize the relationship between many of the classes of fibrations which we will study in this book.
If $f: X \rightarrow S$ is a map of simplicial sets, then we have the following implications:
$$ 
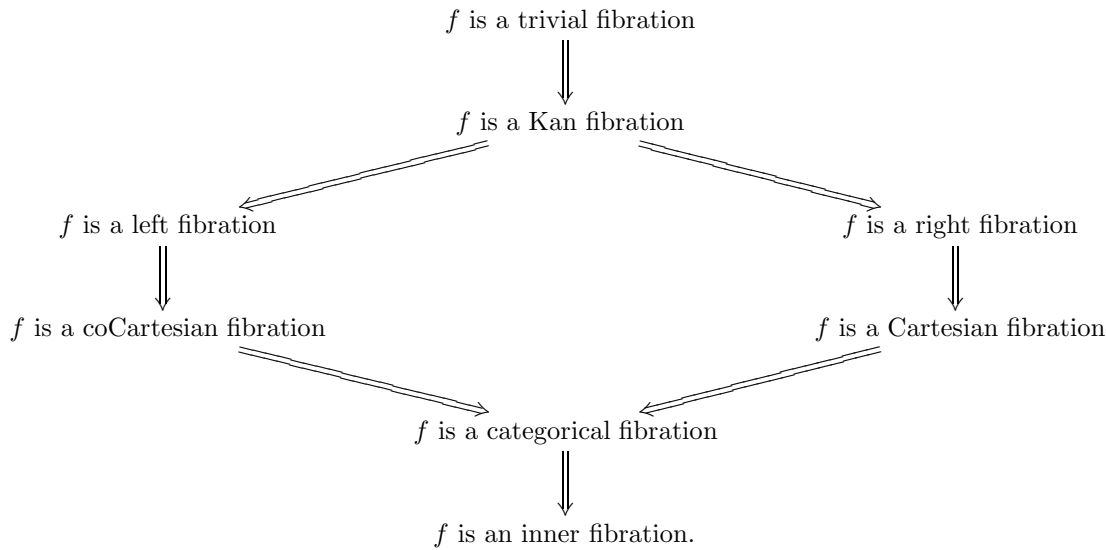
\xymatrix{ & \text{ $f$ is a trivial fibration} \ar@{=>}[d] & \\
& \text{ $f$ is a Kan fibration} \ar@{=>}[dr] \ar@{=>}[dl] & \\
\text{ $f$ is a left fibration} \ar@{=>}[d] & & \text{ $f$ is a right fibration} \ar@{=>}[d] \\
\text{ $f$ is a coCartesian fibration} \ar@{=>}[dr] & & \text{ $f$ is a Cartesian fibration} \ar@{=>}[dl] \\
& \text{ $f$ is a categorical fibration } \ar@{=>}[d] & \\
& \text{ $f$ is an inner fibration. } & }$$
In general, none of these implications is reversible.
\end{remark}

\begin{remark}
The small object argument (Proposition \ref{quillobj}) shows that every map $X \rightarrow Z$ of simplicial sets
admits a factorization $$X \stackrel{p}{\rightarrow} Y \stackrel{q}{\rightarrow} Z,$$ where
$p$ is anodyne (left anodyne, right anodyne, inner anodyne, a cofibration) and
$q$ is a Kan fibration (left fibration, right fibration, inner fibration, trivial fibration).
\end{remark}

\begin{remark}
The theory of left fibrations (left anodyne maps) is {\em dual} to the theory of right fibrations (right anodyne maps): a map $S \rightarrow T$ is a left fibration (left anodyne map) if and only if the induced map $S^{op} \rightarrow T^{op}$ is a right fibration (right anodyne map). Consequently, we will generally confine our remarks in \S \ref{leftfibsec} to the case of left fibrations; the analogous statements for right fibrations will follow by duality.
\end{remark}

\section{Left Fibrations}\label{leftfibsec}
 
\setcounter{theorem}{0}

In this section, we will study the class of {\em left fibrations} between simplicial sets. We begin
in \S \ref{scgp} with a review of some classical category theory: namely, the theory of categories cofibered in groupoids (over another category). We will see that the theory of left fibrations is a natural $\infty$-categorical generalization of this idea. 
In \S \ref{leftfib} we will show that the class of left fibrations is stable under various important constructions, such as the formation of slice $\infty$-categories.

It follows immediately from the definition that every Kan fibration of simplicial sets is a left fibration. The converse is false in general. However, it is possible to give a relatively simple criterion for
testing whether or not a left fibration $f: X \rightarrow S$ is a Kan fibration. We will establish this criterion in \S \ref{crit} and deduce some of its consequences.

The classical theory of Kan fibrations has a natural interpretation in the language of model categories: a map $p: X \rightarrow S$ is a Kan fibration if and only if $X$ is a fibrant object
of $(\sSet)_{/S}$, where the category $(\sSet)_{/S}$ is equipped with its usual model structure. There is a similar characterization of left fibrations: a map $p: X \rightarrow S$ is a left fibration if and only if $X$ is a fibrant object of $(\sSet)_{/S}$ with respect to certain model structure, which
we will refer to as the {\it covariant model structure}. We will define the covariant model structure in \S \ref{contrasec}, and give an overview of its basic properties. 

\subsection{Left Fibrations in Classical Category Theory}\label{scgp}

Before beginning our study of left fibrations, let us recall a bit of classical
category theory. Let $\calD$ be a small category, and suppose we are given a functor $$ \chi: \calD \rightarrow \Gpd,$$
where $\Gpd$ denotes the category of groupoids (where the morphisms are
given by functors). Using the functor $\chi$, we can extract a new category
$\calC_{\chi}$ via the classical {\it Grothendieck construction}: 

\begin{itemize}
\item The objects of $\calC_{\chi}$ are pairs $(D, \eta)$, where $D \in \calD$
and $\eta$ is an object of the groupoid $\chi(D)$.
\item Given a pair of objects $(D, \eta), (D', \eta') \in \calC_{\chi}$, a
morphism from $(D, \eta)$ to $(D', \eta')$ in $\calC_{\chi}$ is given by
a pair $(f, \alpha)$, where $f: D \rightarrow D'$ is a morphism in 
$\calD$, and $\alpha: \chi(f)(\eta) \simeq \eta'$ is an isomorphism in the groupoid $\chi(D')$.
\item Composition of morphisms is defined in the obvious way.
\end{itemize}

There is an evident forgetful functor $F: \calC_{\chi} \rightarrow \calD$, which
carries an object $(D, \eta) \in \calC_{\chi}$ to the underlying object $D \in \calD$.
Moreover, it is possible to reconstruct $\chi$ from the category $\calC_{\chi}$
(together with the forgetful functor $F$), at least up to equivalence; for example, if $D$ is an object of $\calD$, then the groupoid $\chi(D)$ is canonically equivalent to the fiber product $\calC_{\chi} \times_{\calD} \{D\}$. Consequently, the Grothendieck construction sets up a dictionary which relates 
functors $\chi: \calD \rightarrow \Gpd$ with categories $\calC_{\chi}$ admitting a functor $F: \calC_{\chi} \rightarrow \calD$. However, this dictionary is not perfect; not every functor $F: \calC \rightarrow \calD$ arises via the Grothendieck construction described above. To clarify the situation, we recall the following definition:

\begin{definition}\label{latelate}\index{gen}{category!cofibered in groupoids}
Let $F: \calC \rightarrow \calD$ be a functor between categories. We say that
{\it $\calC$ is cofibered in groupoids over $\calD$} if the following conditions are satisfied:
\begin{itemize}
\item[$(1)$] For every object $C \in \calC$ and every morphism
$\eta: F(C) \rightarrow D$ in $\calD$, there exists a morphism
$\widetilde{\eta}: C \rightarrow \widetilde{D}$ such that $F(\widetilde{\eta}) = \eta$.
\item[$(2)$] For every morphism $\eta: C \rightarrow C'$ in $\calC$ and every
object $C'' \in \calC$, the map
$$ \Hom_{\calC}(C',C'') \rightarrow \Hom_{\calC}(C,C'') \times_{ \Hom_{\calD}(F(C), F(C'')) } \Hom_{\calD}(F(C'), F(C''))$$ is bijective.
\end{itemize}
\end{definition}

\begin{example}\label{cik}
Let $\chi: \calD \rightarrow \Gpd$ be a functor from a category $\calD$
to the category of groupoids. Then the forgetful functor $\calC_{\chi} \rightarrow \calD$ exhibits $\calC_{\chi}$ as fibered in groupoids over $\calD$.
\end{example}

Example \ref{cik} admits a converse: suppose we begin with a category
$\calC$ fibered in groupoids over $\calD$. Then, for every 
every object $D \in \calD$, the fiber $\calC_{D} = \calC \times_{\calD} \{D\}$ is a
groupoid. Moreover, for every morphism $f: D \rightarrow D'$ in $\calD$, it
is possible to construct a functor $f_{!}: \calC_{D} \rightarrow \calC_{D'}$
as follows: for each $C \in \calC_{D}$, choose a morphism $\overline{f}: C \rightarrow C'$ covering the map $D \rightarrow D'$, and set $f_{!}(C) = C'$. The map $\overline{f}$ may not be uniquely determined, but it is unique up to isomorphism and depends functorially on $C$. Consequently, we obtain a
functor $f_{!}$, which is well-defined up to isomorphism. We can then
try to define a functor $\chi: \calD \rightarrow \Gpd$ by the formulas
$$ D \mapsto \calC_{D}$$
$$ f \mapsto f_{!}.$$
Unfortunately, this does not quite work: since the functor $f_{!}$ is determined only up to canonical isomorphism by $f$, the identity $(f \circ g)_{!} = f_{!} \circ g_{!}$
holds only up to canonical isomorphism, rather than up to equality. This is merely a technical inconvenience; it can be addressed in (at least) two ways:
\begin{itemize}
\item The groupoid $\chi(D) = \calC \times_{\calD} \{D\}$ can be described as the category of functors $G$ fitting into a commutative diagram
$$ \xymatrix{ & \calC \ar[d]^{F} \\
\{D\} \ar@{-->}[ur]^{G} \ar[r] & \calD. }$$
If we replace the one point category $\{D\}$ with the overcategory $\calD_{D/}$ in this definition, then we obtain a groupoid equivalent to $\chi(D)$ which
depends on $D$ in a strictly functorial fashion.
\item Without modifying the definition of $\chi(D)$, we can realize
$\chi$ as a functor from $\calD$ to an appropriate {\it bicategory} of groupoids.
\end{itemize}

We may summarize the above discussion informally by saying that the Grothendieck construction establishes an equivalence between functors
$\chi: \calD \rightarrow \Gpd$ and categories fibered in groupoids over $\calD$.

The theory of left fibrations should be regarded as an $\infty$-categorical generalization of Definition \ref{latelate}. As a preliminary piece of evidence for this assertion, we offer the following:

\begin{proposition}\label{stinkyer}
Let $F: \calC \rightarrow \calD$ be a functor between categories. Then $\calC$ is cofibered in groupoids over $\calD$ if and only if the induced map $\Nerve(F): \Nerve(\calC) \rightarrow \Nerve(\calD)$ is a left fibration of simplicial sets.
\end{proposition}

\begin{proof}
Proposition \ref{ruko} implies that $\Nerve(F)$ is an inner fibration. It follows that
$\Nerve(F)$ is a left fibration if and only if it has the right lifting property with respect to
$\Lambda^n_0 \subseteq \Delta^n$ for all $n > 0$. When $n = 1$, the relevant lifting property is equivalent to $(1)$ of Definition \ref{latelate}. When $n=2$ ($n=3$) the relevant lifting property is
equivalent to the surjectivity (injectivity) of the map described in $(2)$. For $n > 3$, the relevant lifting property is automatic (since a map $\Lambda^n_0 \rightarrow S$ extends {\em uniquely} to $\Delta^n$ when $S$ is isomorphic to the nerve of a category). 
\end{proof}

Let us now consider the structure of a general left fibration $p: X \rightarrow S$.
In the case where $S$ consists of a single vertex, Proposition \ref{greenwich} asserts that $p$ is a left fibration if and only if $X$ is a Kan complex. Since the class of left fibrations is stable under pullback, we deduce that for {\em any} left fibration $p: X \rightarrow S$ and any vertex $s$ of $S$, the fiber $X_{s} = X \times_S \{s\}$ is a Kan complex (which we can think of as the $\infty$-categorical analogue of a groupoid). Moreover, these Kan complexes are related to one another. More precisely, suppose that $f: s \rightarrow s'$ is an edge of the simplicial set $S$ and consider the inclusion
$i: X_s \simeq X_s \times \{0\} \subseteq X_s
\times \Delta^{1}$. In \S \ref{leftfib}, we will prove that $i$ is left anodyne (Corollary \ref{prodprod1}). 
It follows that we can solve the lifting problem
$$ \xymatrix{ \{0 \} \times X_{s} \ar@{^{(}->}[d] \ar@{^{(}->}[rr] & & X \ar[d]^{p} \\
\Delta^1 \times X_{s} \ar@{-->}[urr] \ar[r] & \Delta^1 \ar[r]^{f} & S.}$$
Restricting the dotted arrow to $\{1\} \times X_{s}$, we obtain a map
$f_{!}: X_{s} \rightarrow X_{s'}$. Of course, $f_{!}$ is not unique, but it is uniquely determined
up to homotopy.\index{not}{f!@$f_{!}$}

\begin{lemma}\label{functy}
Let $q: X \rightarrow S$ be a left fibration of simplicial sets. The assignment
$$s \in S_0 \mapsto X_{s} $$
$$ f \in S_1 \mapsto f_{!} $$
determines a (covariant) functor from the homotopy category $\h{S}$ into the homotopy category $\calH$ of spaces.
\end{lemma}

\begin{proof}
Let $f: s \rightarrow s'$ be an edge of $S$.
We note the following characterization of the morphism $f_{!}$ in $\calH$. Let $K$ be any simplicial set, and suppose given homotopy classes of maps $\eta \in \Hom_{\calH}(K,X_{s})$, $\eta' \in \Hom_{\calH}(K,X_{s'})$. Then
$\eta' = f_{!} \circ \eta$ if and only if there exists a map $p: K \times \Delta^1 \rightarrow X$ such that
$q \circ p$ is given by the composition $$K \times \Delta^1 \rightarrow \Delta^1 \stackrel{f}{\rightarrow} S,$$
$\eta$ is the homotopy class of $p | K \times \{0\}$, and $\eta'$ is the homotopy class of $p| K \times \{1\}$. 

Now consider any $2$-simplex $\sigma: \Delta^2 \rightarrow S$, which we will depict as
$$ \xymatrix{ & v \ar[dr]^{g} & \\
u \ar[ur]^{f} \ar[rr]^{h} & & w. }$$
We note that the inclusion
$X_{ u } \times \{0\} \subseteq X_{ u } \times \Delta^2$ is left-anodyne 
(Corollary \ref{prodprod1}). Consequently there exists a map $p: X_{ u } \times \Delta^2 \rightarrow X$ such that $p| X_{u} \times \{0\}$ is the inclusion $X_{ u } \subseteq X$ and $q \circ p$ is the composition
$ X_{ u } \times \Delta^2 \rightarrow \Delta^2 \stackrel{\sigma}{\rightarrow} S$. 
Then $f_{!} \simeq p| X_{ u } \times \{1\}$, $h_{!} = p| X_{u} \times \{2 \}$, and the map $p | X_{ u } \times \Delta^{ \{1,2\} }$ verifies the equation
$$ h_{!} = g_{!} \circ f_{!}$$ in $\Hom_{\calH}( X_{u}, X_{ w })$.
\end{proof}

We can summarize the situation informally as follows. Fix a simplicial set $S$. To give a left fibration $q: X \rightarrow S$, one must specify a Kan complex $X_{s}$ for each ``object'' of $S$, a map
$f_{!}: X_{s} \rightarrow X_{s'}$ for each ``morphism'' $f: s \rightarrow s'$ of $S$, and ``coherence data'' for these morphisms for each higher-dimensional simplex of $S$. In other words, giving a left fibration ought to be more or less the same thing as giving a functor from $S$ to the $\infty$-category $\SSet$ of spaces.
Lemma \ref{functy} can be regarded as a weak version of this assertion; we will prove something considerably more precise in \S \ref{contrasec} (see Theorem \ref{struns}). 

We close this section by establishing two simple properties of left fibrations, which will be neededΠin the proof of Proposition \ref{greenlem}:

\begin{proposition}\label{hamb1}
Let $p: \calC \rightarrow \calD$ be a left fibration of $\infty$-categories, and let
$f: X \rightarrow Y$ be a morphism in $\calC$ such that $p(f)$ is an equivalence in
$\calD$. Then $f$ is an equivalence in $\calC$.
\end{proposition}

\begin{proof}
Let $\overline{g}$ be a homotopy inverse to $p(f)$ in $\calD$, so that
there exists a $2$-simplex of $\calD$ depicted as follows:
$$ \xymatrix{ & p(Y) \ar[dr]^{\overline{g}} & \\
p(X) \ar[ur]^{p(f)} \ar[rr]^{\id_{p(X)}} & & p(X). }$$
Since $p$ is a left fibration, we can lift this to a diagram
$$ \xymatrix{ & Y \ar[dr]^{g} & \\
X \ar[ur]^{f} \ar[rr]^{\id_{X}} & & X }$$
in $\calC$. It follows that $g \circ f \simeq \id_{X}$, so that $f$ admits a left homotopy inverse.
Since $p(g) = \overline{g}$ is an equivalence in $\calD$, the same argument proves that
$g$ has a left homotopy inverse. This left homotopy inverse must coincide with $f$, since $f$ is a right homotopy inverse to $g$. Thus $f$ and $g$ are homotopy inverse in the $\infty$-category $\calC$, so that $f$ is an equivalence as desired.
\end{proof}

\begin{proposition}\label{hamb2}
Let $p: \calC \rightarrow \calD$ be a left fibration of $\infty$-categories, let
$Y$ be an object of $\calC$, and let
$\overline{f}: \overline{X} \rightarrow p(Y)$ be an equivalence in $\calD$. Then
there exists a morphism $f: X \rightarrow Y$ in $\calC$ such that $p(f) = \overline{f}$
(automatically an equivalence, in view of Proposition \ref{hamb1}). 
\end{proposition}

\begin{proof}
Let $\overline{g}: p(Y) \rightarrow \overline{X}$ be a homotopy inverse to $\overline{f}$
in $\calC$. Since $p$ is a left fibration, there exists a morphism
$g: Y \rightarrow X$ such that $\overline{g} = p(g)$. Since $\overline{f}$ and $\overline{g}$
are homotopy inverse to one another, there exists a $2$-simplex of $\calD$ which we can depict as follows:
$$ \xymatrix{ & p(X) \ar[dr]^{\overline{f}} & \\
p(Y) \ar[ur]^{ p(g) } \ar[rr]^{ \id_{p(Y)} } & & p(Y). }$$
Applying the assumption that $p$ is a left fibration once more, we can lift this to a diagram
$$ \xymatrix{ & X \ar[dr]^{f} & \\
Y \ar[ur]^{g} \ar[rr]^{\id_{Y}} & & Y, }$$
which proves the existence of $f$.
\end{proof}

\subsection{Stability Properties of Left Fibrations}\label{leftfib}

The purpose of this section is to show that left fibrations of simplicial sets exist in abundance. Our main results are Proposition \ref{sharpen} (which is our basic source of examples for left fibrations) and Corollary \ref{ichy} (which asserts that left fibrations are stable under the formation of functor categories). 

Let $\calC$ be an $\infty$-category, and let $\SSet$ denote the $\infty$-category of spaces.
One can think of a functor from $\calC$ to $\SSet$ as a ``cosheaf of spaces'' on $\calC$. By analogy with ordinary category theory, one might expect that the basic example of such a cosheaf would be the cosheaf corepresented by an object $C$ of $\calC$; roughly speaking this should be given by the functor
$$ D \mapsto \bHom_{\calC}(C,D).$$ 
As we saw in \S \ref{scgp}, it is natural to guess that such a functor can be encoded 
by a left fibration $\widetilde{\calC} \rightarrow \calC$. There is a natural candidate
for $\widetilde{\calC}$: the undercategory $\calC_{C/}$. 
Note that the fiber of the map
$$ f: \calC_{C/} \rightarrow \calC$$ over the object $D \in \calC$ is the Kan complex $\Hom_{\calC}^{\lft}(C,D)$. The assertion that $f$ is a left fibration is a consequence of the following more general result:

\begin{proposition}[Joyal]\label{sharpen}
Suppose given a diagram of simplicial sets
$$ K_0 \subseteq K \stackrel{p}{\rightarrow} X \stackrel{q}{\rightarrow} S$$
where $q$ is an inner fibration. Let $r = q \circ p: K \rightarrow S$,
$p_0 = p|K_0$, and $r_0 = r|K_0$. Then the induced map
$$ X_{p/} \rightarrow X_{p_0/} \times_{ S_{r_0/} } S_{r/}$$
is a left fibration. If the map $q$ is already a left fibration, then
the induced map
$$ X_{/p} \rightarrow X_{/ p_0} \times_{ S_{/r_0} } S_{/r}$$
is a left fibration as well.
\end{proposition}

Proposition \ref{sharpen} immediately implies the following half of
Proposition \ref{gorban3}, which we asserted earlier without proof:

\begin{corollary}[Joyal]\label{gorban4}\index{gen}{left fibration!and undercategories}
Let $\calC$ be an $\infty$-category and $p: K \rightarrow \calC$ an arbitrary diagram.
Then the projection $\calC_{p/} \rightarrow \calC$ is a left fibration. In particular, $\calC_{p/}$ is itself an $\infty$-category.
\end{corollary}

\begin{proof}
Apply Proposition \ref{sharpen} in the case where $X = \calC$, $S = \ast$, $A = \emptyset$, $B = K$.
\end{proof}

We can also use Proposition \ref{sharpen} to prove Proposition \ref{greenlem}, 
which was stated without proof in \S \ref{obmor}.

\begin{proposition2}
Let $\calC$ be an $\infty$-category, and $\phi: \Delta^1 \rightarrow \calC$ a morphism of $\calC$. Then $\phi$ is an equivalence if and only if, for every $n \geq 2$ and every map
$f_0: \Lambda^n_0 \rightarrow \calC$ such that $f_0 | \Delta^{\{0,1\}} = \phi$,
there exists an extension of $f_0$ to $\Delta^n$.
\end{proposition2}

\begin{proof}
Suppose first that $\phi$ is an equivalence, and let $f_0$ be as above. To find the desired extension of $f_0$, we must produce the dotted arrow in the associated diagram
$$ \xymatrix{ \{ 0 \} \ar@{^{(}->}[d] \ar[r] & \calC_{/ \Delta^{n-2} } \ar[d]^{q} \\
\Delta^1 \ar[r]^{\phi'} \ar@{-->}[ur] & \calC_{/ \bd \Delta^{n-2} }. }$$  
The projection map $p: \calC_{/ \bd \Delta^{n-2} } \rightarrow \calC$ is a right fibration (Proposition \ref{sharpen}). Since $\phi'$ is a preimage of $\phi$ under $p$, Proposition \ref{hamb1} implies that $\phi'$ is an equivalence. Because $q$ is a right fibration (Proposition \ref{sharpen} again), 
the existence of the dotted arrow follows from Proposition \ref{hamb2}.

We now prove the converse. Let $\phi: X \rightarrow Y$ be a morphism in $\calC$,
and consider the map $\Lambda^2_0 \rightarrow \calC$ indicated in the following diagram:
$$ \xymatrix{ & Y \ar@{-->}[dr]^{\psi} & \\
X \ar[ur]^{\phi} \ar[rr]^{\id_{X}} & & X. }$$
The assumed extension property ensures the existence of the dotted morphism
$\psi: Y \rightarrow X$ and a $2$-simplex $\sigma$ which verifies the identity
$\psi \circ \phi \simeq \id_{X}$. We now consider the map
$$ \tau_0: \Lambda^3_0 \stackrel{ ( \bigdot, s_0 \phi, s_1 \psi, \sigma) } \longrightarrow \calC.$$
Once again, our assumption allows us to extend $\tau_0$ to a $3$-simplex
$\tau: \Delta^3 \rightarrow \calC$, and the face $d_0 \tau$ verifies the identity
$\phi \circ \psi = \id_{Y}$. It follows that $\psi$ is a homotopy inverse to $\phi$, so that
$\phi$ is an equivalence in $\calC$.
\end{proof}

We now turn to the proof of Proposition \ref{sharpen}. It is an easy consequence of the following more basic observation:

\begin{lemma}[Joyal \cite{joyalnotpub}]\label{precough}
Let $f: A_0 \subseteq A$ and $g: B_0 \subseteq B$ be inclusions of simplicial sets.
Suppose either that $f$ is right anodyne, or that $g$ is left anodyne. Then the induced inclusion
$$ h: (A_0 \star B) \coprod_{ A_0 \star B_0 } (A \star B_0) \subseteq A \star B$$
is inner anodyne.
\end{lemma}

\begin{proof}
We will prove that $h$ is inner anodyne whenever $f$ is right anodyne; the other assertion follows by a dual argument. 

Consider the class of {\em all} morphisms $f$ for which the conclusion of the lemma holds (for any inclusion $g$). This class of morphisms is weakly saturated; to prove that it contains all right-anodyne morphisms, it suffices to show that it contains each of the inclusions
$f: \Lambda^n_{j} \subseteq \Delta^n$ for $0 < j \leq n$. We may therefore assume that $f$ is of this form.

Now consider the collection of all inclusions $g$ for which $h$ is inner anodyne (where $f$ is now fixed). This class of morphisms is also weakly saturated; to prove that it contains all inclusions, it suffices to show that the lemma holds when $g$ is of the form $\bd \Delta^m \subseteq \Delta^m$. In this case, $h$ can be identified with the inclusion $\Lambda^{n+m+1}_j \subseteq \Delta^{n+m+1}$, which is inner anodyne because $0 < j \leq n < n+m+1$.
\end{proof}

The following result can be proven by exactly the same argument:

\begin{lemma}[Joyal]\label{purcough}
Let $f: A_0 \rightarrow A$ and $g: B_0 \rightarrow B$ be inclusions of simplicial sets.
Suppose that $f$ is left anodyne. Then the induced inclusion
$$ (A_0 \star B) \coprod_{ A_0 \star B_0 } (A \star B_0) \subseteq A \star B$$
is left anodyne.
\end{lemma}

\begin{proof}[Proof of Proposition \ref{sharpen}]
After unwinding the definitions, the first assertion follows from Lemma \ref{precough} and the second from Lemma \ref{purcough}.
\end{proof}

For future reference, we record the following counterpart to Proposition \ref{sharpen}:

\begin{proposition}[Joyal]\label{sharpen2}
Let $\pi: S \rightarrow T$ be an inner fibration, $p: B \rightarrow S$ a map of simplicial sets,
$i: A \subseteq B$ an inclusion of simplicial sets, $p_0 = p | A$, $p' = \pi \circ p$, and $p'_0 = \pi \circ p_0 = p'|A$. Suppose that either $i$ is right anodyne, or $\pi$ is a left fibration. 
Then the induced map
$$ \phi: S_{p/} \rightarrow S_{p_0/} \times_{ T_{p'_0/}} T_{p'/}$$ is a trivial Kan fibration.
\end{proposition}

\begin{proof}
Consider the class of all cofibrations $i: A \rightarrow B$ for which $\phi$ is a trivial fibration for {\em every} inner fibration (right fibration) $p: S \rightarrow T$. It is not difficult to see that this is a weakly saturated class of morphisms; thus, it suffices to consider the case where $A = \Lambda^m_i$, $B = \Delta^m$, for $0 < i \leq m$ ($0 \leq i \leq m$).

Let $q: \bd \Delta^n \rightarrow S_{p/}$ be a map, and suppose given an extension of
$\phi \circ q$ to $\Delta^n$. We wish to find a compatible extension of $q$. Unwinding the definitions, we are given a map
$$ r: (\Delta^m \star \bd \Delta^n) \coprod_{\Lambda^m_i \star \bd \Delta^n} ( \Lambda^m_i \star \Delta^n) 
\rightarrow S$$ which we wish to extend to $\Delta^m \star \Delta^n$ in a manner that is compatible
with a given extension $\Delta^m \star \Delta^n \rightarrow T$ of the composite map $\pi \circ r$. 
The existence of such an extension follows immediately from the assumption that $p$ has the right lifting property with respect to the horn inclusion $\Lambda^{n+m+1}_i \subseteq \Delta^{n+m+1}$.
\end{proof}

The remainder of this section is devoted to the study of the behavior of left fibrations under exponentiation. Our goal is to prove an assertion of the following form: if $p: X \rightarrow S$ is a left fibration of simplicial sets, then so is the induced map $X^{K} \rightarrow S^{K}$, for every simplicial set $K$ (this is a special case of Corollary \ref{ichy} below). This is an easy consequence of the following characterization of left anodyne maps, which is due to Joyal:

\begin{proposition}[Joyal \cite{joyalnotpub}]\label{usejoyal}\index{gen}{left anodyne}\index{gen}{anodyne!left}
The following collections  of morphisms generate the same weakly saturated class of morphisms
of $\sSet$:
\begin{itemize}
\item[$(1)$] The collection $A_1$ of all horn inclusions $\Lambda^n_i
\subseteq \Delta^n$, $0 \leq i < n$.

\item[$(2)$] The collection $A_2$ of all inclusions $$(\Delta^m \times
\{0\}) \coprod_{ \bd \Delta^m \times \{0\} } (\bd \Delta^m \times
\Delta^1) \subseteq \Delta^m \times \Delta^1.$$

\item[$(3)$] The collection $A_3$ of all inclusions $$(S' \times \{0\})
\coprod_{S \times \{0\} } (S \times \Delta^1) \subseteq S' \times
\Delta^1,$$ where $S \subseteq S'$.

\end{itemize}
\end{proposition}

\begin{proof}
Let $S \subseteq S'$ be as in $(3)$. Working cell-by-cell on $S'$,
we deduce that every morphism in $A_3$ can be obtained as an
iterated pushout of morphisms belonging to $A_2$. Conversely,
$A_2$ is contained in $A_3$, which proves that they generate the
same weakly saturated collection of morphisms.

To proceed with the proof, we must first introduce a bit of
notation. The $(n+1)$-simplices of $\Delta^n \times \Delta^1$ are
indexed by order-preserving maps
$$ [n+1] \rightarrow [0, \ldots, n] \times
[0,1].$$ We let $\sigma_k$ denote the map
$$\sigma_k(m) =
\begin{cases} (m,0) & \text{if } m \leq k \\
(m-1,1) & \text{if } m > k. \end{cases}$$ We will also denote by $\sigma_k$
the corresponding $(n+1)$-simplex of $\Delta^n \times
\Delta^1$. We note that $\{ \sigma_k \}_{ 0 \leq k \leq n }$ are
precisely the nondegenerate $(n+1)$-simplices of $\Delta^n \times
\Delta^1$.

We define a collection $\{ X(k) \}_{ 0 \leq k \leq n+1 }$ of
simplicial subsets of $\Delta^n \times \Delta^1$ by descending
induction on $k$. We begin by setting
$$X(n+1) = (\Delta^n \times \{0\}) \coprod_{ \bd \Delta^n \times \{0\} } (\bd
\Delta^n \times \Delta^1).$$ Assuming that $X(k+1)$ has been
defined, we let $X(k) \subseteq \Delta^n \times \Delta^1$ be the
union of $X(k+1)$ and the simplex $\sigma_k$ (together with all the
faces of $\sigma_k$). We note that this description exhibits $X(k)$ as a
pushout $$ X(k+1) \coprod_{ \Lambda^{n+1}_k } \Delta^{n+1},$$ and
also that $X(0) = \Delta^{n} \times \Delta^1$. It follows that
each step in the chain of inclusions
$$ X(n+1) \subseteq X(n) \subseteq \ldots \subseteq X(1) \subseteq
X(0)$$ is contained in the class of morphisms generated by $A_1$,
so that the inclusion $X(n+1) \subseteq X(0)$ is generated by
$A_1$.

To complete the proof, we show that each inclusion in
$A_1$ is a retract of an inclusion in $A_3$. More specifically,
the inclusion $\Lambda^{n}_i \subseteq \Delta^{n}$ is a retract of
$$(\Delta^n \times \{0\}) \coprod_{ \Lambda^n_i \times \{0\} } (
\Lambda^n_i \times \Delta^1) \subseteq \Delta^n \times \Delta^1,$$
so long as $0 \leq i < n$. We will define the relevant maps
$$ \Delta^{n} \stackrel{j}{\rightarrow} \Delta^n \times \Delta^1
\stackrel{r}{\rightarrow} \Delta^{n}$$ and leave it to the reader
to verify that they are compatible with the relevant subobjects. The
map $j$ is simply the inclusion
$\Delta^n \simeq \Delta^n \times \{1\} \subseteq \Delta^n \times \Delta^1$.
The map $r$ is induced by a map of partially ordered sets, which we will also denote by $r$. It may be described by the formulae

$$r(m,0) =
\begin{cases} m & \text{if } m \neq i+1 \\
i & \text{if } m = i+1 \end{cases}$$
$$r(m,1) = m.$$
\end{proof}

\begin{corollary}\label{prodprod1}
Let $i: A \rightarrow A'$ be left-anodyne, and let $j: B \rightarrow B'$ be a cofibration. Then
the induced map $$(A \times B') \coprod_{A \times B} (A' \times B)
\rightarrow A' \times B'$$ is left-anodyne.
\end{corollary}

\begin{proof}
This follows immediately from Proposition \ref{usejoyal2}, which characterizes the class of left-anodyne maps as the class generated by $A_3$ (which is stable under smash products with any cofibration).
\end{proof}

\begin{remark}
A basic fact in the homotopy theory of simplicial sets is that the analogue of Corollary \ref{prodprod1} holds also for the class of {\em anodyne} maps of simplicial sets. Since the class of anodyne maps is generated (as a weakly saturated class of morphisms) by the class of left anodyne maps and the class of right anodyne maps, this classical fact follows from Corollary \ref{prodprod1} (together with the dual assertion concerning right anodyne maps).
\end{remark}

\begin{corollary}\label{ichy}\index{gen}{left fibration!and functor categories}
Let $p: X \rightarrow S$ be a left-fibration, and let $i: A \rightarrow B$ be any cofibration of simplicial sets. Then the induced map $q: X^{B} \rightarrow X^A \times_{ S^A } S^B$ is a left fibration. If $i$ is left anodyne, then $q$ is a trivial Kan fibration.
\end{corollary}

\begin{corollary}[Homotopy Extension Lifting Property]\label{helper}
Let $p: X \rightarrow S$ be a map of simplicial sets. Then $p$ is a left fibration if and only if the induced map
$$ X^{\Delta^1} \rightarrow X^{ \{0\} } \times_{ S^{ \{0\} } } S^{\Delta^1}$$
is a trivial Kan fibration of simplicial sets.
\end{corollary}

For future use, we record the following criterion for establishing that a morphism is left anodyne:

\begin{proposition}\label{trull11}
Let $p: X \rightarrow S$ be a map of simplicial sets, let
$s: S \rightarrow X$ be a section of $p$, and let
$h \in \Hom_{S}(X \times \Delta^1, X)$ be a (fiberwise)
simplicial homotopy from $s \circ p = h| X \times \{0\}$ to $\id_{X} = h| X \times \{1\}$. Then $s$
is left anodyne.
\end{proposition}

\begin{proof}
Consider a diagram
$$ \xymatrix{ S \ar[d]^{s} \ar[r]^{g} & Y \ar[d]^{q} \\
X \ar[r]^{g'} \ar@{-->}^{f}[ur] & Z }$$
where $q$ is a left fibration. We must show that it is possible to find a map
$f$ rendering the diagram commutative. Define
$F_0: (S \times \Delta^1) \coprod_{ S \times \{ 0 \} } (X \times \{0\})$
to be the composition of $g$ with the projection onto $S$.
Now consider the diagram
$$ \xymatrix{ (S \times \Delta^1) \coprod_{ S \times \{0\} }
(X \times \{0\}) \ar[d] \ar[rrr]^{F_0} & & & Y \ar[d]^{q} \\
X \times \Delta^1 \ar[rrr]^{g' \circ h} \ar@{-->}[urrr]^{F} & & & Z. }$$
Since $q$ is a left fibration and the left vertical map is left anodyne, it
is possible to supply the dotted arrow $F$ as indicated. Now we observe that
$f = F|X \times \{1\}$ has the desired properties.
\end{proof}

\subsection{A Characterization of Kan Fibrations}\label{crit}

Let $p: X \rightarrow S$ be a left fibration of simplicial sets. As we saw in \S \ref{scgp}, $p$ determines for each vertex $s$ of $S$ a Kan complex $X_{s}$, and for each edge $f: s \rightarrow s'$ a map of Kan complexes $f_{!}: X_{s} \rightarrow X_{s'}$ (which is well-defined up to homotopy). If $p$ is a Kan fibration, then the same argument allows us to construct a map $X_{s'} \rightarrow X_{s}$, which is a homotopy inverse to $f_{!}$. Our goal in this section is to prove the following converse:

\begin{proposition}\label{dent}\index{gen}{left fibration!and Kan fibrations}
Let $p: S \rightarrow T$ be a left fibration of simplicial sets. The following conditions are equivalent:

\begin{itemize}
\item[$(1)$] The map $p$ is a Kan fibration.
\item[$(2)$] For every edge $f: t \rightarrow t'$ in $T$, the map $f_{!}: S_{t} \rightarrow S_{t'}$ is an isomorphism
in the homotopy category $\calH$ of spaces.
\end{itemize}

\end{proposition}

\begin{lemma}\label{strike2}
Let $p: S \rightarrow T$ be a left fibration of simplicial sets. Suppose that $S$ and $T$ are Kan complexes, and that $p$ is a homotopy equivalence. Then $p$ induces a surjection from $S_0$ to $T_0$.
\end{lemma}

\begin{proof}
Fix a vertex $t \in T_0$. Since $p$ is a homotopy equivalence, there exists a vertex $s \in S_0$ and an edge $e$ joining $p(s)$ to $t$. Since $p$ is a left fibration, this edge lifts to an edge $e': s \rightarrow s'$ in $S$. Then $p(s')=t$.
\end{proof}

\begin{lemma}\label{toothie2}
Let $p: S \rightarrow T$ be a left fibration of simplicial sets. Suppose that $T$ is a Kan complex. Then $p$ is a Kan fibration.
\end{lemma}

\begin{proof} 
We note that the projection $S \rightarrow \ast$, being a composition of left fibrations
$S \rightarrow T$ and $T \rightarrow \ast$, is a left fibration, so that $S$ is also a Kan complex.
Let $A \subseteq B$ be an anodyne inclusion of simplicial sets. We must show that the map
$p: S^B \rightarrow S^A \times_{ T^A} T^B$ is surjective on vertices. Since $S$ and $T$ are Kan complexes, the maps $T^B \rightarrow T^A$ and $S^B \rightarrow S^A$ are trivial fibrations. It follows that $p$ is a homotopy equivalence and a left fibration. Now we simply apply Lemma \ref{strike2}.
\end{proof}

\begin{lemma}\label{toothie}
Let $p: S \rightarrow T$ be a left fibration of simplicial sets. Suppose that for every vertex $t \in T$, the fiber $S_{t}$ is contractible. Then $p$ is a trivial Kan fibration.
\end{lemma}

\begin{proof}
It will suffice to prove the analogous result for {\em right} fibrations (we do this in order to keep the notation we use below consistent with that employed in the proof of Proposition \ref{usejoyal}).

Since $p$ has nonempty fibers, it has the right lifting property with respect to the inclusion
$\emptyset = \bd \Delta^0 \subseteq \Delta^0$. Let $n > 0$, $f: \bd \Delta^n \rightarrow S$ any map, and $g: \Delta^n \rightarrow T$ an extension of $p \circ f$. We must show that there exists an extension $\widetilde{f}: \Delta^n \rightarrow S$ with $g = p \circ \widetilde{f}$.

Pulling back via the map $G$, we may suppose that $T = \Delta^n$ and $g$ is the identity map, so that $S$ is an $\infty$-category. Let $t$ denote the initial vertex of $T$. There is a unique map
$g': \Delta^n \times \Delta^1 \rightarrow T$ such that $g' | \Delta^n \times \{1\} = g$ and
$g' | \Delta^n \times \{0\}$ is constant at the vertex $t$.

Since the inclusion $\bd \Delta^n \times \{1\} \subseteq \bd \Delta^n \times \Delta^1$ is right anodyne, there exists an extension $f'$ of $f$ to $\bd \Delta^n \times \Delta^1$ which covers $g' | \bd \Delta^n \times \Delta^1$. To complete the proof, it suffices to show that we can extend $f'$ to a map $\widetilde{f}': \Delta^n \times \Delta^1 \rightarrow S$ (such an extension is automatically compatible with $g'$ in view of our assumptions that $T = \Delta^n$ and $n > 0$). Assuming this has been done, we simply define
$\widetilde{f} = \widetilde{f}' | \Delta^n \times \{1\}$.

Recall the notation of the proof of Proposition \ref{usejoyal}, and filter the simplicial set
$\Delta^n \times \Delta^1$ by the simplicial subsets
$$ X(n+1) \subseteq \ldots \subseteq X(0) = \Delta^n \times \Delta^1.$$
We extend the definition of $f'$ to $X(m)$ by a descending induction on $m$. When $m=n+1$, we note that $X(n+1)$ is obtained from $\bd \Delta^n \times \Delta^1$ by adjoining the interior of the simplex 
$ \bd \Delta^n \times \{0 \}$. Since the boundary of this simplex maps entirely into the contractible Kan complex $S_{t}$, it is possible to extend $f'$ to $X(n+1)$.

Now suppose the definition of $f'$ has been extended to $X(i+1)$. We note that $X(i)$ is obtained from $X(i+1)$ by pushout along a horn inclusion $\Lambda^{n+1}_i \subseteq \Delta^{n+1}$. If $i > 0$, then the assumption that $S$ is an $\infty$-category guarantees the existence of an extension of $f'$ to $X(i)$. When $i =0$, we note that $f'$ carries the initial edge of $\sigma_0$ into the fiber $S_{t}$. Since $S_{t}$ is a Kan complex, $f'$ carries the initial edge of $\sigma_0$ to an equivalence in $S$, and the desired extension of $f'$ exists by Proposition \ref{greenlem}.
\end{proof}

\begin{proof}[Proof of Proposition \ref{dent}]
Suppose first that $(1)$ is satisfied, and let $f: t \rightarrow t'$ be an edge in $T$. 
Since $p$ is a right fibration, the edge $f$ induces a map $f^{\ast}: S_{t'} \rightarrow S_{t}$, which is
well-defined up to homotopy. It is not difficult to check that the maps $f^{\ast}$ and
$f_{!}$ are homotopy inverse to one another; in particular, $f_{!}$ is a homotopy equivalence.
This proves that $(1) \Rightarrow (2)$.

Assume now that $(2)$ is satisfied. A map of simplicial sets is a Kan fibration if and only if it is both a right fibration and a left fibration; consequently, it will suffice to prove that $p$ is a right fibration. According to Corollary \ref{helper}, it will suffice to show that $$q: S^{\Delta^1} \rightarrow S^{ \{1\} } \times_{ T^{ \{1\} }} T^{ \Delta^1} $$ is a trivial Kan fibration. Corollary \ref{ichy} implies that $q$ is a left fibration. By Lemma \ref{toothie}, it suffices to show that the fibers of $q$ are contractible.

Fix an edge $f: t \rightarrow t'$ in $T$. Let $X$ denote the simplicial set of sections of the projection
$S \times_{ T} {\Delta^1} \rightarrow \Delta^1$, where $\Delta^1$ maps into $T$ via the edge $f$. Consider the fiber $q': X \rightarrow S_{ t'}$ of $q$ over the edge $f$. Since the $q$ and $q'$ have the same fibers (over points of $S^{ \{1\} }\times_{ T^{ \{1\} } } T^{ \Delta^1} $ whose second projection is the edge $f$), it will suffice to show that $q'$ is a trivial fibration for every choice of $f$.

Consider the projection $r: X \rightarrow S_{t}$. Since $p$ is a left fibration, $r$ is a trivial fibration.
Because $S_{t}$ is a Kan complex, so is $X$. Lemma \ref{toothie2} implies that $q'$ is a Kan fibration.
We note that $f_{!}$ is obtained by choosing a section of $r$ and then composing with $q'$. Consequently, assumption $(2)$ implies that $q'$ is a homotopy equivalence, and thus a trivial fibration, which completes the proof.
\end{proof}

\begin{remark}
Lemma \ref{toothie} is an immediate consequence of Proposition \ref{dent}, since any map between contractible Kan complexes is a homotopy equivalence. Lemma \ref{toothie2} also follows immediately, since if $T$ is a Kan complex, then its homotopy category is a groupoid, so that
{\em any} functor $\h{T} \rightarrow \calH$ carries edges of $T$ to invertible morphisms in $\calH$.
\end{remark}

\subsection{The Covariant Model Structure}\label{contrasec}

In \S \ref{leftfib}, we saw that a left fibration $p: X \rightarrow S$ determines a functor
$\chi$ from $\h{S}$ to the homotopy category $\calH$, carrying each 
vertex $s$ to the fiber $X_{s} = X \times_{S} \{s\}$. We would like to formulate
a more precise relationship between left fibrations over $S$ and functors
from $S$ into spaces. For this, it is convenient to employ Quillen's language of model categories. In this section, we will show that the category $(\sSet)_{/S}$ can be endowed with the structure of a simplicial model category, whose fibrant objects are precisely the left fibrations $X \rightarrow S$. In \S \ref{valencequi}, we will give an $\infty$-categorical version of the Grothendieck construction, provided by a suitable right Quillen functor
$$ (\sSet)^{\sCoNerve[S]} \rightarrow (\sSet)_{/S}$$
which we will prove to be a Quillen equivalence (Theorem \ref{struns}).

\begin{warning}
We will assume throughout this section that the reader is familiar with
the theory of model categories, as presented in \S \ref{appmodelcat}. We will also assume familiarity with the model structure on the category $\sCat$ of simplicial categories (see \S \ref{compp4}).
\end{warning}

\begin{definition}\index{gen}{cone!left}\index{gen}{left cone}\index{gen}{cone!right}\index{gen}{right cone}\index{gen}{mapping cone}\index{not}{Cl@$C^{\triangleleft}(f)$}\index{not}{Cr@$C^{\triangleleft}(f)$}
Let $f: X \rightarrow S$ be a map of simplicial sets. The
{\it left cone} of $f$ is the simplicial set $S \coprod_{X} X^{\triangleleft}$.
We will denote the left cone of $f$ by $C^{\triangleleft}(f)$.
Dually, we define the {\em right cone} of $f$ to be the simplicial
set $C^{\triangleleft}(f) = S \coprod_{X} X^{\triangleleft}$.
\end{definition}

\begin{remark}
Let $f: X \rightarrow S$ be a map of simplicial sets. There is a canonical monomorphism
of simplicial sets $S \rightarrow C^{\triangleleft}(f)$. We will generally identify $S$ with its image under this monomorphism, and thereby regard $S$ as a simplicial subset of
$C^{\triangleleft}(f)$. We note that there is a unique vertex of
$C^{\triangleleft}(f)$ which does not belong to $S$. We will refer to this vertex as
the {\it cone point} of $C^{\triangleleft}(f)$. 
\end{remark}

\begin{example}
Let $S$ be a simplicial set, and let $\id_{S}$ denote the identity map from $S$ to itself.
Then $C^{\triangleleft}(\id_S)$ and $C^{\triangleright}(\id_S)$ can be identified with
$S^{\triangleleft}$ and $S^{\triangleright}$, respectively.
\end{example}

\begin{definition}\index{gen}{covariant!cofibration}\index{gen}{covariant!fibration}\index{gen}{covariant!equivalence}\index{gen}{cofibration!covariant}\index{gen}{fibration!covariant}\index{gen}{equivalence!covariant}
Let $S$ be a simplicial set. We will say that a map $f: X \rightarrow Y$ in $(\sSet)_{/S}$ is
a:
\begin{itemize}
\item[$(C)$] {\it covariant cofibration} if it is a monomorphism of simplicial sets.
\item[$(W)$] {\it covariant equivalence} if the induced map
$$ X^{\triangleleft} \coprod_{X} S \rightarrow Y^{\triangleleft} \coprod_{Y} S $$
is a categorical equivalence.
\item[$(F)$] {\it covariant fibration} if it has the right lifting property with respect to every map
which is both a covariant cofibration and a covariant equivalence.
\end{itemize}
\end{definition}

\begin{lemma}\label{onehalff}
Let $S$ be a simplicial set. Then every left anodyne map in $(\sSet)_{/S}$ is a covariant
equivalence.
\end{lemma}

\begin{proof}
By general nonsense, it suffices to prove the result for a generating left anodyne inclusion of the form $\Lambda^n_i \subseteq \Delta^n$, where $0 \leq i < n$. In other words, we must show any map
$$ i: ( \Lambda^n_i)^{\triangleleft} \coprod_{ \Lambda^n_i} S \rightarrow (\Delta^n)^{\triangleleft} \coprod_{ \Delta^n} S$$
is a categorical equivalence. We now observe that $i$ is a pushout of the inner anodyne inclusion
$\Lambda^{n+1}_{i+1} \subseteq \Delta^{n+1}$.
\end{proof}

\begin{proposition}\label{covcech}\index{gen}{covariant!model structure}\index{gen}{model category!covariant}
Let $S$ be a simplicial set. The covariant cofibrations, covariant equivalences, and covariant fibrations determine a left proper, combinatorial model structure on $(\sSet)_{/S}$.
\end{proposition}

\begin{proof}
It suffices to show that conditions $(1)$, $(2)$, and $(3)$ of Proposition \ref{goot} are met.
We consider each in turn:
\begin{itemize}
\item[$(1)$] The class $(W)$ of weak equivalences is perfect. This follows from
Corollary \ref{perfpull}, since the functor $X \mapsto X^{\triangleleft} \coprod_{X} S$ commutes with filtered colimits.

\item[$(2)$] It is clear that the class $(C)$ of cofibrations is generated by a set. We must show that weak equivalences are stable under pushouts by cofibrations. In other words, suppose we are given a pushout diagram
$$ \xymatrix{ X \ar[r]^{j} \ar[d]^{i} & Y \ar[d] \\
X' \ar[r]^{j'} & Y' }$$
in $(\sSet)_{/S}$ where $i$ is a covariant cofibration and $j$ is a covariant equivalence. We must show that $j'$ is a covariant equivalence. We obtain a pushout diagram in $\sCat$
$$ \xymatrix{ \sCoNerve[X^{\triangleleft} \coprod_{X} S]  \ar[r] \ar[d] & \sCoNerve[Y^{\triangleleft} \coprod_{Y} S] \ar[d] \\
\sCoNerve[(X')^{\triangleleft} \coprod_{X'} S] \ar[r] & \sCoNerve[(Y')^{\triangleleft} \coprod_{Y'} S]}$$
which is homotopy coCartesian, since $\sCat$ is a left proper model category. Since the upper horizontal map is an equivalence, so it the bottom horizontal map; thus $j'$ is a covariant equivalence.

\item[$(3)$] We must show that a map $p: X \rightarrow Y$ in $\sSet$, which has the right lifting
property with respect to every map in $(C)$, belongs to $(W)$. We note in that case that $p$ is a trivial Kan fibration, and therefore admits a section $s: Y \rightarrow X$. We will show that
$p$ and $s$ induce mutually inverse isomorphisms between
$\sCoNerve[ X^{\triangleleft} \coprod_{X} S]$ and $\sCoNerve[ Y^{\triangleleft} \coprod_{Y} S]$
in the homotopy category $\h{\sCat}$; it will then follow that $p$ is a covariant equivalence.

Let $f: X \rightarrow X$ denote the composition $s \circ p$; we wish to show that the
map $\sCoNerve[ X^{\triangleleft} \coprod_{X} S]$ induced by $f$ is equivalent
to the identity in $\h{\sCat}$. We observe that $f$ is homotopic to the identity $\id_{X}$ via
a homotopy $h: \Delta^1 \times X \rightarrow X$. It will therefore suffice to show
that $h$ is a covariant equivalence. But $h$ admits a left inverse
$$X \simeq \{0\} \times X \subseteq \Delta^1 \times X$$
which is left anodyne (Corollary \ref{prodprod1}) and therefore
a covariant equivalence by Lemma \ref{onehalff}.
\end{itemize}
\end{proof}

\begin{proposition}\label{natsim}
The category $(\sSet)_{/S}$ is a simplicial model category $($with respect to the covariant model structure and the natural simplicial structure$)$.
\end{proposition}

\begin{proof}
We will deduce this from Proposition \ref{testsimpmodel}. The only nontrivial point is to verify that for any $X \in (\sSet)_{/S}$, the projection $X \times \Delta^n \rightarrow X$ is a covariant equivalence. But this map has a section $X \times \{0\} \rightarrow X \times \Delta^n$, which is left anodyne and therefore a covariant equivalence (Proposition \ref{onehalf}).
\end{proof}

We will refer to the model structure of Proposition \ref{covcech} as the {\it covariant model structure} on $(\sSet)_{/S}$. We will prove later that the covariantly fibrant objects of $(\sSet)_{/S}$ are precisely the left fibrations $X \rightarrow S$ (Corollary \ref{usewhere1}). For the time being, we will be content to make a much weaker observation:

\begin{proposition}\label{onehalf}
Let $S$ be a simplicial set.
\begin{itemize}
\item[$(1)$] Every left anodyne map in $(\sSet)_{/S}$ is a trivial cofibration with respect to the covariant model structure.
\item[$(2)$] Every covariant fibration in $(\sSet)_{/S}$ is a left fibration of simplicial sets.
\item[$(3)$] Every fibrant object of $(\sSet)_{/S}$ determines a left fibration $X \rightarrow S$.
\end{itemize}
\end{proposition}

\begin{proof}
Assertion $(1)$ follows from Lemma \ref{onehalff}, and 
the implications $(1) \Rightarrow (2) \Rightarrow (3)$ are obvious.
\end{proof}

Our next result expresses the idea that the covariant model structure on $(\sSet)_{/S}$ depends functorially on $S$:

\begin{proposition}\label{contrafunk}
Let $j: S \rightarrow S'$ be a map of simplicial sets. Let 
$j_{!}: (\sSet)_{/S} \rightarrow (\sSet)_{/S'}$ be the forgetful functor $($given by composition with
$j${}$)$, and let $j^{\ast}: (\sSet)_{/S'} \rightarrow (\sSet)_{/S}$ be its right adjoint, which is given by the formula
$$ j^{\ast} X' = X' \times_{S'} S.$$
Then we have a Quillen adjunction
$$ \Adjoint{ j_{!} }{ (\sSet)_{/S} }{ (\sSet)_{/S'}}{j^{\ast}}$$
(with the covariant model structures). 
\end{proposition}

\begin{proof}
It is clear that $j_{!}$ preserves cofibrations. For $X \in (\sSet)_{S}$, the pushout diagram
$$ \xymatrix{ S \ar[r] \ar[d] & S' \ar[d] \\
X^{\triangleleft} \coprod_{X} S \ar[r] & X^{\triangleleft} \coprod_{X} S'}$$
is a homotopy pushout (with respect to the Joyal model structure). Thus $j_{!}$ preserves covariant equivalences. It follows that $(j_{!}, j^{\ast})$ is a Quillen adjunction.
\end{proof}

\begin{remark}
Let $j: S \rightarrow S'$ be as in Proposition \ref{contrafunk}. If $j$ is a categorical
equivalence, then the Quillen adjunction $(j_{!}, j^{\ast})$ is a categorical
equivalence. This follows from Theorem \ref{struns} and Proposition \ref{lesstrick}.
\end{remark}

\begin{remark}\index{gen}{contravariant model structure}\index{gen}{model category!contravariant}
Let $S$ be a simplicial set. The covariant model structure on $(\sSet)_{/S}$
is usually not self-dual. Consequently, we may define a new model
structure on $(\sSet)_{/S}$ as follows:
\begin{itemize}
\item[$(C)$] A map $f$ in $(\sSet)_{/S}$ is a {\it contravariant cofibration} if
it is a monomorphism of simplicial sets.
\item[$(W)$] A map $f$ in $(\sSet)_{/S}$ is a {\it contravariant equivalence} if
$f^{op}$ is a covariant equivalence in $(\sSet)_{/S^{op}}$.\index{gen}{contravariant equivalence}\index{gen}{equivalence!contravariant} 
\item[$(F)$] A map $f$ in $(\sSet)_{/S}$ is a {\it contravariant fibration} if
$f^{op}$ is a covariant fibration in $(\sSet)_{/S^{op}}$.\index{gen}{contravariant fibration}\index{gen}{fibration!contravariant}
\end{itemize}
We will refer to this model structure on $(\sSet)_{/S}$ as the {\it contravariant} model structure. Propositions \ref{natsim}, \ref{onehalf} and \ref{contrafunk} have evident analogues in the contravariant setting.
\end{remark}

\section{Simplicial Categories and $\infty$-Categories}\label{valencequi}

\setcounter{theorem}{0}

For every topological category $\calC$ and every pair of objects $X,Y \in \calC$, Theorem \ref{biggie} asserts that the counit map
$$u: |\bHom_{\sCoNerve[\tNerve(\calC)] }(X,Y)| \rightarrow \bHom_{\calC}(X,Y)$$
is a weak homotopy equivalence of topological spaces. This result is the main ingredient needed to establish the equivalence between the theory of topological categories and the theory of $\infty$-categories. The goal of this section is to give a proof of Theorem \ref{biggie} and to develop some of its consequences. 

We first replace Theorem \ref{biggie} by a statement about {\em simplicial} categories. Consider the composition
$$ \bHom_{\sCoNerve[ \tNerve(\calC) ]}(X,Y) \stackrel{v}{\rightarrow}
\Sing \bHom_{| \sCoNerve[ \tNerve(\calC)] |}(X,Y) \stackrel{ \Sing(u)}{\rightarrow} \Sing \bHom_{\calC}(X,Y).$$
Classical homotopy theory ensures that $v$ is a weak homotopy equivalence. Moreover, $u$ is a weak homotopy equivalence of topological spaces if and only if $\Sing(u)$ is a weak homotopy equivalence of simplicial sets. Consequently, $u$ is a weak homotopy equivalence of topological spaces if and only if $\Sing(u) \circ v$ is a weak homotopy equivalence of simplicial sets.
It will therefore suffice to prove the following {\em simplicial} analogue of Theorem \ref{biggie}:

\begin{theorem}\label{biggiesimp}
Let $\calC$ be a fibrant simplicial category $($that is, a simplicial category in which each mapping
space $\bHom_{\calC}(x,y)$ is a Kan complex$)$, and let $x, y \in \calC$ be a pair of objects. The counit map
$$ u: \bHom_{ \sCoNerve[\sNerve(\calC)]}(x,y) \rightarrow \bHom_{\calC}(x,y)$$ is a weak
homotopy equivalence of simplicial sets. 
\end{theorem}

The proof will be given in \S \ref{compp2} (see Proposition \ref{wiretrack}).
Our strategy is as follows:

\begin{itemize}
\item[$(1)$] We will show that, for every simplicial set $S$, there is a close relationship between
{\it right fibrations} $S' \rightarrow S$ and {\it simplicial presheaves} $\calF: \sCoNerve[S]^{op} \rightarrow \sSet$. This relationship is controlled by the straightening and unstraightening
functors which we introduce in \S \ref{rightstraight}.
\item[$(2)$] Suppose that $S$ is an $\infty$-category. Then, for each object
$y \in S$, the projection $S_{/y} \rightarrow S$ is a right fibration, which corresponds to a simplicial presheaf $\calF: \sCoNerve[S]^{op} \rightarrow \sSet$. This simplicial presheaf
$\calF$ is related to $S_{/y}$ in two different ways:
\begin{itemize}
\item[$(i)$] As a simplicial presheaf, $\calF$ is weakly equivalent to the functor
$x \mapsto \bHom_{ \sCoNerve[S]}(x,y)$.
\item[$(ii)$] For each object $x$ of $S$, there is a canonical homotopy equivalence
$\calF(x) \rightarrow S_{/y} \times_{S} \{x\} \simeq \Hom_{S}^{\rght}(x,y)$.
Here the Kan complex $\Hom_{S}^{\rght}(x,y)$ is defined as in \S \ref{prereq1}.
\end{itemize}
\item[$(3)$] Combining observations $(i)$ and $(ii)$, we will conclude that
the mapping spaces $\Hom_{S}^{\rght}(x,y)$ are homotopy equivalent to the correpsonding
mapping spaces $\Hom_{\sCoNerve[S]}(x,y)$.
\item[$(4)$] In the special case where $S$ is the nerve of a fibrant simplicial category $\calC$,
there is a canonical map $\Hom_{\calC}(x,y) \rightarrow \Hom_{S}^{\rght}(x,y)$, which we will
show to be a homotopy equivalence in \S \ref{twistt}.
\item[$(5)$] Combining $(3)$ and $(4)$, we will obtain a canonical isomorphism
$\bHom_{\calC}(x,y) \simeq \bHom_{ \sCoNerve[ \sNerve(\calC) ] }(x,y)$ in the homotopy category of spaces. We will then show that this isomorphism is induced by the unit map appearing in the statement of Theorem \ref{biggiesimp}.
\end{itemize}

We will conclude this section with \S \ref{compp3}, where we apply Theorem \ref{biggiesimp} to 
construct the {\it Joyal model structure} on $\sSet$ and to establish a more refined version of the equivalence between $\infty$-categories and simplicial categories.

\subsection{The Straightening and Unstraightening Constructions (Unmarked Case)}\label{rightstraight}

In \S \ref{scgp}, we asserted that a left fibration $X \rightarrow S$ can be viewed as a functor
from $S$ into a suitable $\infty$-category of Kan complexes. Our goal in this section is to make this idea precise. For technical reasons, it will be somewhat more convenient to phrase our results
in terms of the dual theory of {\em right} fibrations $X \rightarrow S$.
Given any functor $\phi: \sCoNerve[S]^{op} \rightarrow \calC$ between simplicial categories,
we will define an {\it unstraightening functor} $\Un_{\phi}: \Set_{\Delta}^{\calC} \rightarrow
(\sSet)_{/S}$. If $\calF: \calC \rightarrow \sSet$ is a diagram taking values in Kan complexes, then
the associated map $\Un_{\phi} \calF \rightarrow S$ is a right fibration, whose fiber at
a point $s \in S$ is homotopy equivalent to the Kan complex $\calF( \phi(s) )$.

Fix a simplicial set $S$, a simplicial category $\calC$ and a functor $\phi: \sCoNerve[S] \rightarrow \calC^{op}$. Given an object $X \in (\sSet)_{/S}$, let $v$ denote the cone point of $X^{\triangleright}$. We can view the simplicial category $$\calM = \sCoNerve[X^{\triangleright}] \coprod_{ \sCoNerve[X]} \calC^{op}$$
as a correspondence from $\calC^{op}$ to $\{v\}$, which we can identify with a simplicial functor
$$ \St_{\phi} X: \calC \rightarrow \sSet.$$\index{gen}{straightening functor}\index{gen}{unstraightening functor}
\index{not}{Stphi@$\St_{\phi}$}\index{not}{Unphi@$\Un_{\phi}$}\index{not}{StS@$\St_{S}$}\index{not}{UnS@$\Un_{S}$}
This functor is described by the formula
$$ (\St_{\phi} X)(C) = \bHom_{\calM}(C,v).$$
We may regard $\St_{\phi}$ as a functor from $(\sSet)_{/S}$ to $(\sSet)^{\calC}$. We refer to
$\St_{\phi}$ as the {\it straightening functor} associated to $\phi$. In the special case where
$\calC = \sCoNerve[S]^{op}$ and $\phi$ is the identity map, we will write
$\St_{S}$ instead of $\St_{\phi}$.

By the adjoint functor theorem (or by direct construction), the straightening functor $\St_{\phi}$
associated to $\phi: \sCoNerve[S] \rightarrow \calC^{op}$ has a right adjoint, which we will denote by $\Un_{\phi}$ and refer to as the {\it unstraightening functor}. We now record the obvious functoriality properties of this construction.

\begin{proposition}\label{straightchange}
\begin{itemize}
\item[$(1)$] Let $p: S' \rightarrow S$ be a map of simplicial sets, $\calC$ a simplicial category, and
$\phi: \sCoNerve[S] \rightarrow \calC^{op}$ a simplicial functor, and let $\phi': \sCoNerve[S'] \rightarrow \calC^{op}$ denote the composition $\phi \circ \sCoNerve[p]$. 
Let $p_{!}: (\sSet)_{/S'} \rightarrow (\sSet)_{/S}$ denote the forgetful functor, given by composition with $p$. There is a natural isomorphism of functors
$$ \St_{\phi} \circ p_{!} \simeq \St_{\phi'}$$
from $(\sSet)_{/S'}$ to $\Set_{\Delta}^{\calC}$.

\item[$(2)$] Let $S$ be a simplicial set, $\pi: \calC \rightarrow \calC'$ a simplicial functor between simplicial categories, and $\phi: \sCoNerve[S] \rightarrow
\calC^{op}$ a simplicial functor. Then there is a natural isomorphism of functors $$\St_{\pi^{op} \circ \phi} \simeq \pi_{!} \circ \St_{\phi}$$
from $(\sSet)_{/S}$ to $\Set_{\Delta}^{\calC'}$. Here $\pi_{!}: \Set_{\Delta}^{\calC} \rightarrow \Set_{\Delta}^{\calC'}$
is the left adjoint to the functor $\pi^{\ast}: \Set_{\Delta}^{\calC'} \rightarrow \Set_{\Delta}^{\calC}$ given by composition with $\pi$.
\end{itemize}
\end{proposition}

Our main result is the following:

\begin{theorem}\label{struns}
Let $S$ be a simplicial set, $\calC$ a simplicial category, and 
$\phi: \sCoNerve[S] \rightarrow \calC^{op}$ a simplicial functor. The straightening and unstraightening functors determine a Quillen adjunction
$$ \Adjoint{\St_{\phi}}{ (\sSet)_{/S} }{ \Set_{\Delta}^{\calC} }{ \Un_{\phi} },$$
where $(\sSet)_{/S}$ is endowed with the contravariant model structure and
$\Set_{\Delta}^{\calC}$ with the projective model structure.
If $\phi$ is an equivalence of simplicial categories, then $(\St_{\phi}, \Un_{\phi})$ is a Quillen equivalence.
\end{theorem}

\begin{proof}
It is easy to see that $\St_{\phi}$ preserves cofibrations and weak equivalences, so that
the pair $( \St_{\phi}, \Un_{\phi})$ is a Quillen adjunction. The real content of
Theorem \ref{struns} is the final assertion. Suppose that $\phi$ is an equivalence of simplicial categories; then we wish to show that $( \St_{\phi}, \Un_{\phi})$ is a Quillen equivalence.
We will prove this result in \S \ref{fullun} as a consequence of Proposition \ref{fullmeal}.
\end{proof}

\subsection{Straightening Over a Point}\label{twistt}

In this section, we will study behavior of the straightening functor
$\St_{X}$ in the case where the simplicial set $X = \{x\}$ consists of a single vertex.
In this case, we can view $\St_{X}$ as a colimit-preserving functor
from the category of simplicial sets to itself. We begin with a few general remarks about such functors.

Let $\cDelta$ denote the category of combinatorial simplices and
$\sSet$ the category of simplicial sets, so that $\sSet$ may be
identified with the category of presheaves of sets on $\cDelta$.
If $\calC$ is {\em any} category which admits small colimits, then
any functor $f: \cDelta \rightarrow \calC$ extends to a
colimit-preserving functor $F: \sSet \rightarrow \calC$ (which is unique up to unique isomorphism). We may regard $f$ as a cosimplicial
object $C^{\bigdot}$ of $\calC$. In this case, we shall denote the\index{gen}{geometric realization}
functor $F$ by
$$ S \mapsto |S|_{C^{\bigdot}}.$$\index{not}{|S|_C@$|S|_{C^{\bigdot}}$}

\begin{remark}
Concretely, one constructs $|S|_{C^{\bigdot}}$ by taking the
disjoint union of $S_n \times C^n$ and making the appropriate
identifications along the ``boundaries''. In the language of category theory, the geometric realization is given by
the {\it coend} $$\int_{[n] \in \cDelta} S_n \times C^{n}.$$\index{gen}{coend}
\end{remark}

The functor $S \mapsto |S|_{C^{\bigdot}}$ has a right adjoint
which we shall denote by $\Sing_{C^{\bigdot}}$. It may be described
by the formula
$$ \Sing_{C^{\bigdot}}(X)_n = \Hom_{\calC}( C^n, X).$$\index{not}{Sing_C@$\Sing_{C^{\bigdot}}(X)$}

\begin{example}
Let $\calC$ be the category $\CG$ of compactly generated Hausdorff
spaces, and let $C^{\bigdot}$ be the cosimplicial space defined by
$$ C^n = \{ ( x_0, \ldots, x_n) \in [0,1]^{n+1} : x_0 + \ldots +
x_n = 1 \}.$$ Then $|S|_{C^{\bigdot}}$ is the usual {\it geometric
realization} $|S|$ of the simplicial set $S$ and
$\Sing_{C^{\bigdot}} = \Sing$ is the functor which assigns to each topological space $X$ its its singular complex.
\end{example}

\begin{example}\index{gen}{standard simplex}
Let $\calC$ be the category $\sSet$, and let $C^{\bigdot}$ be the {\it standard simplex} (the cosimplicial object of $\sSet$ given by the Yoneda embedding):
$$ C^n = \Delta^n.$$ Then $||_{C^{\bigdot}}$ and $\Sing_{C^{\bigdot}}$ are
both (isomorphic to) the identity functor on $\sSet$.
\end{example}

\begin{example}
Let $\calC = \Cat$, and let $f: \cDelta \rightarrow \Cat$ be the functor which associates
to each finite nonempty linearly ordered set $J$ the corresponding category.
Then $\Sing_{C^{\bigdot}}= \Nerve$ is the functor which associates to each category its nerve, and $||_{C^{\bigdot}}$ associates, to each simplicial set $S$, the homotopy category $\h{S}$ as defined in \S \ref{hcat}.
\end{example}

\begin{example}
Let $\calC = \sCat$, and let $C^{\bigdot}$ be the cosimplicial
object of $\calC$ given in Definitions \ref{csimp1} and
\ref{csimp2}. Then 
$\Sing_{C^{\bigdot}}$ is the simplicial nerve functor, and
$||_{C^{\bigdot}}$ is its left adjoint
$$ S \mapsto \sCoNerve[S]. $$
\end{example}

Let us now return to the case of the straightening functor
$\St_{X}$, where $X = \{x\}$ consists of a single vertex. The above remarks show that
we can identify $\St_{X}$ with the geometric realization functor
$| |_{Q^{\bigdot}}: \sSet \rightarrow \sSet$, for some cosimplicial object
$Q^{\bigdot}$ in $\sSet$. To describe $Q^{\bigdot}$ more explicitly, let
us first define a cosimplicial simplicial set $J^{\bigdot}$ by the formula
$$ J^{n} = (\Delta^n \star \{y\}) \coprod_{ \Delta^n } \{x\}.$$
The cosimplical simplicial set $Q^{\bigdot}$ can then be described by the formula
$Q^{n} = \bHom_{ \sCoNerve[J^{n}]}(x,y)$. 

In order to proceed with our analysis, we need to understand better the cosimplicial object $Q^{\bigdot}$ of $\sSet$. It admits the following description:\index{not}{Qdot@$Q^{\bigdot}$}

\begin{itemize}
\item For each $n \geq 0$, let $P_{[n]}$ denote the partially ordered set of
{\em nonempty} subsets of $[n]$, and $K_{[n]}$ the simplicial set $\Nerve(P)$ (which may be identified with a simplicial subset of the $(n+1)$-cube $(\Delta^1)^{n+1}$).
The simplicial set $Q^{n}$ is obtained by collapsing, for
each $0 \leq i \leq n$, the subset
$$ ( \Delta^1 )^{ \{ j: 0 \leq j < i \} } \times
\{1\} \times (\Delta^1)^{ \{ j: i < j \leq n \} } \subseteq K_{[n]}$$
to its quotient $( \Delta^1)^{ \{ j: i < j \leq n \} }$. 

\item A map $f: [n] \rightarrow [m]$ determines a map $P_f: P_{[n]}
\rightarrow P_{[m]}$, by setting $P_f(I)=f(I)$. The map $P_f$ in turn
induces a map of simplicial sets $K_{[n]} \rightarrow K_{[m]}$, which
determines a map of quotients $Q^n \rightarrow Q^m$
when $f$ is order-preserving.
\end{itemize}

\begin{remark}\index{not}{Qcaldot@$\calQ^{\bigdot}$}
Let $\calQ^{\bigdot} = |Q^{\bigdot}|$ denote the cosimplicial space obtained by
applying the (usual) geometric realization functor to
$Q^{\bigdot}$. The space $\calQ^n$ may be described as a
quotient of the cube of all functions $p: [n] \rightarrow [0,1]$ satisfying $p(0)=1$. This cube is to be divided by the following equivalence relation: $p \simeq p'$ if there
exists a nonnegative integer $i \leq n$ such that $p| \{i, \ldots
n\} = p'| \{i, \ldots, n\}$ and $p(i) = p'(i) = 1$.

Each $\calQ^n$ is homeomorphic to an $n$-simplex, and these
homeomorphisms may be chosen to be compatible with the face maps
of the cosimplicial space $\calQ^{\bigdot}$. However,
$\calQ^{\bigdot}$ is not isomorphic to the standard simplex
because it has very different degeneracies. For example, the
product of the degeneracy mappings $\calQ^n \rightarrow
(\calQ^1)^n$ is not injective for $n \geq 2$.
\end{remark}

Our goal for the remainder of this section is to study the functors
$\Sing_{Q^{\bigdot}}$ and $||_{Q^{\bigdot}}$ and to prove
that they are ``close'' to the identity functor. More precisely, there is a map $\pi: Q^{\bigdot} \rightarrow \Delta^{\bigdot}$ of
cosimplicial objects of $\sSet$. It is induced by a map $K_{[n]}
\rightarrow \Delta^n$, which the nerve of the map of partially ordered sets
$P_{[n]} \rightarrow [n]$ which carries each nonempty subset of
$[n]$ to its largest element.

\begin{proposition}\label{babyy}
Let $S$ be a simplicial set. Then the map $p_S: |S|_{Q^{\bigdot}}
\rightarrow S$ induced by $\pi$ is a weak homotopy equivalence.
\end{proposition}

\begin{proof}
Consider the collection $A$ of simplicial sets $S$ for which the assertion of Proposition \ref{babyy} holds. 
Since $A$ is stable under filtered colimits, it will suffice to prove that every simplicial set $S$ having only finitely many nondegenerate simplices belongs to $A$. We prove this by induction on the dimension $n$ of $S$, and the number of nondegenerate simplices of $S$ of dimension $n$. If $S = \emptyset$, there is nothing to prove; otherwise we may write
$$S \simeq S' \coprod_{ \bd \Delta^n } \Delta^n $$
$$ |S|_{Q^{\bigdot}} \simeq |S'|_{Q^{\bigdot}} \coprod_{ |\bd \Delta^n|_{Q^{\bigdot}} } |\Delta^n|_{Q^{\bigdot}}.$$
Since both of these pushouts are homotopy pushouts, it suffices to show that $p_{S'}$, $p_{\bd \Delta^n}$, and $p_{\Delta^n}$ are weak homotopy equivalences. For $p_{S'}$ and $p_{\bd \Delta^n}$, this follows from the inductive hypothesis; for $p_{\Delta^n}$, we need only observe that both $\Delta^n$ and $|\Delta^n|_{Q^{\bigdot}} = Q^n$ are weakly contractible.
\end{proof}

\begin{remark}
The strategy used to prove Proposition \ref{babyy} will reappear frequently throughout this book: it allows us to prove theorems about arbitrary simplicial sets by reducing to the case of simplices.
\end{remark}

\begin{proposition}\label{realremmy}
The adjoint functors
$\Adjoint{ | |_{Q^{\bigdot}}}{\sSet}{\sSet}{ \Sing_{Q^{\bigdot} } }$
determine a Quillen equivalence from the category $\sSet$
$($endowed with the Kan model structure$)$ to itself.
\end{proposition}

\begin{proof}
We first show that the functors $( | |_{Q^{\bigdot}}, \Sing_{Q^{\bigdot}} )$ determine a Quillen adjunction from
$\sSet$ to itself. For this, it suffices to prove that the functor
$S \mapsto |S|_{Q^{\bigdot}}$
preserves cofibrations and weak equivalences. The case of cofibrations is easy, and the second case follows from Proposition \ref{babyy}.
To complete the proof, it will suffice to show that the left derived functor
$L | |_{Q^{\bigdot}}$ determines an equivalence from the homotopy category
$\calH$ to itself. This follows immediately from Proposition \ref{babyy}, which implies
that $L | |_{Q^{\bigdot}}$ is equivalent to the identity functor.
\end{proof}

\begin{corollary}\label{remmy33}
Let $X$ be a Kan complex. Then the counit map
$$v: | \Sing_{Q^{\bigdot}} X|_{Q^{\bigdot}} \rightarrow X$$ is a weak
homotopy equivalence.
\end{corollary}

\begin{remark}\label{undef}
Let $S$ be a simplicial set containing a vertex $s$. Let $\calC$ be a simplicial category,
$\phi: \sCoNerve[S]^{op} \rightarrow \calC$ a simplicial functor, and $C = \phi(s) \in \calC$.
For every simplicial functor $\calF: \calC \rightarrow \sSet$, there is a canonical isomorphism
$$ (\Un_{\phi} \calF) \times_{S} \{s\} \simeq \Sing_{Q^{\bigdot}} \calF(C).$$
In particular, we have a canonical map from $\calF(C)$ to the fiber $(\Un_{\phi} \calF)_{s}$, which is a homotopy equivalence if $\calF(C)$ is fibrant.
\end{remark}

\begin{remark}\label{simfun}
Let $\calC$ and $\calC'$ be simplicial categories. Given a pair of simplicial functors
$\calF: \calC \rightarrow \sSet$, $\calF': \calC' \rightarrow \sSet$, we let $\calF \boxtimes \calF': \calC \times \calC' \rightarrow \sSet$ denote the functor described by the formula
$$(\calF \boxtimes \calF')(C,C') = \calF(C) \times \calF'(C').$$
Given a pair of simplicial functors $\phi: \sCoNerve[S]^{op} \rightarrow \calC$, $\phi': \sCoNerve[S']^{op} \rightarrow \calC'$, we let $\phi \boxtimes \phi'$ denote the induced map
$\sCoNerve[S \times S'] \rightarrow \calC \times \calC'$. We observe that there is a canonical isomorphism of functors
$$ \Un_{\phi \boxtimes \phi'}( \calF \boxtimes \calF') \simeq \Un_{\phi}(\calF) \times \Un_{\phi'}(\calF').$$
Restricting our attention to the case where $S' = \Delta^0$ and $\phi'$ is an isomorphism,
we obtain an isomorphism
$$ \Un_{\phi}( \calF \boxtimes K) \simeq \Un_{\phi}(\calF) \times \Sing_{Q^{\bigdot}} K,$$
for every simplicial set $K$. In particular, for every pair of functors
$\calF, \calG \in \Set_{\Delta}^{\calC}$, we have a chain of maps
\begin{eqnarray*}
\Hom_{\sSet}(K, \bHom_{ \Set^{\calC}_{\Delta}}( \calF, \calG) )
& \simeq & \Hom_{ \Set^{\calC}_{\Delta}}( \calF \boxtimes K, \calG ) \\
& \rightarrow & \Hom_{ (\sSet)_{/S} }( \Un_{\phi}( \calF \boxtimes K), \Un_{\phi} \calG ) \\
& \simeq & \Hom_{ (\sSet)_{/S} }( \Un_{\phi}(\calF) \times \Sing_{Q^{\bigdot}} K, \Un_{\phi} \calG) \\
& \rightarrow & \Hom_{ (\sSet)_{/S} }( \Un_{\phi}(\calF) \times K, \Un_{\phi} \calG) \\
& \simeq & \Hom_{\sSet}(K, \bHom_{ (\sSet)_{/S}}( \Un_{\phi}(\calF), \Un_{\phi}(\calG) ).
\end{eqnarray*}
This construction is natural in $K$, and therefore determines a map of simplicial sets
$$ \bHom_{ \Set^{\calC}_{\Delta}}(\calF, \calG) \rightarrow
\bHom_{ (\sSet)_{/S}}( \Un_{\phi} \calF, \Un_{\phi} \calG ).$$
Together, these maps endow the unstraightening functor
$\Un_{\phi}$ with the structure of a {\em simplicial} functor from
$\Set^{\calC}_{\Delta}$ to $(\sSet)_{/S}$. 
\end{remark}

The cosimplicial object $Q^{\bigdot}$ of $\sSet$ will play an important role in our proof
of Theorem \ref{biggie}. To explain this, let us suppose that $\calC$ is a simplicial category and
$S = \Nerve(\calC)$ is its simplicial nerve.
For every pair of vertices $\overline{x}, \overline{y} \in S$, we can consider the right mapping space
$\Hom^{\rght}_{S}(\overline{x},\overline{y})$. By definition, giving an $n$-simplex of $\Hom^{\rght}_{S}(\overline{x},\overline{y})$ is equivalent to giving a map of simplicial sets $J^{n} \rightarrow S$, which carries $x$ to $\overline{x}$ and $y$ to $\overline{y}$. Using the identification
$S \simeq \Nerve(\calC)$, we see that this is equivalent to giving a map
$\sCoNerve[ J^{n} ]$ into $\calC$, which again carries $x$ to $\overline{x}$ and
$y$ to $\overline{y}$. This is simply the data of a map of simplicial sets
$Q^{n} \rightarrow \bHom_{\calC}( \overline{x}, \overline{y} )$. Moreover, this identification
is natural with respect to $[n]$; we therefore have the following result:

\begin{proposition}\label{remmy22}
Let $\calC$ be a simplicial category, and let $X,Y \in \calC$ be
two objects. There is a natural isomorphism of simplicial sets
$\Hom^{\rght}_{\sNerve(\calC)}(X,Y) \simeq
\Sing_{Q^{\bigdot}} \bHom_{\calC}(X,Y)$.
\end{proposition}

\subsection{Straightening of Right Fibrations}\label{fullun}

Our goal in this section is to prove Theorem \ref{struns}, which asserts that
the Quillen adjunction $$\Adjoint{ \St_{\phi} }{ (\sSet)_{/S} }{\Set_{\Delta}^{\calC}}{ \Un_{\phi} }$$
is a Quillen equivalence when $\phi: \sCoNerve[S] \rightarrow \calC^{op}$ is an equivalence
of simplicial categories. We first treat the case where $S$ is a simplex.

\begin{lemma}\label{sticx}
Let $n$ be a nonnegative integer, let $[n]$ denote the linearly ordered set
$\{0, \ldots, n\}$, regarded as a $($discrete$)$ simplicial category, and let
$\phi: \sCoNerve[ \Delta^n ] \rightarrow [n]$ be the canonical functor. Then
the Quillen adjunction
$$\Adjoint{\St_{\phi}}{ (\sSet)_{/ \Delta^n } }{ \sSet^{[n]} }{\Un_{\phi} }$$
is a Quillen equivalence.
\end{lemma}

\begin{proof}
It follows from the definition of the contravariant model structure that the left derived
functor $L \St_{\phi}: \h{(\sSet)_{/ \Delta^n } } \rightarrow \h{ \Set_{\Delta}^{[n] }}$
is conservative. It will therefore suffice to show that the counit map
$L \St_{\phi} \circ R \Un_{\phi} \rightarrow \id$ is an isomorphism of functors
from $\h{ \Set_{\Delta}^{[n]} }$ to itself. For this, we must show that if
$\calF: [n] \rightarrow \sSet$ is projectively fibrant, then the counit map
$$ \St_{\phi} \Un_{\phi} \calF \rightarrow \calF$$
is an equivalence in $\Set_{\Delta}^{[n]}$. In other words, we may assume that
$\calF(i)$ is a Kan complex for $i \in [n]$, and we wish to prove that
each of the induced maps
$$ v_i: (\St_{\phi} \Un_{\phi} \calF)(i) \rightarrow \calF(i)$$
is a weak homotopy equivalence of simplicial sets. 

Let $\psi: [n] \rightarrow [1]$ be defined by the formula 
$$ \psi'(j) = \begin{cases} 0 & \text{if } 0 \leq j \leq i \\
1& \text{otherwise.} \end{cases}$$
Then, for every object $X \in (\sSet)_{/ \Delta^n }$, we have isomorphisms
$$ (\St_{\phi} X)(i) \simeq ( \St_{\psi \circ \phi} X)(0)
\simeq | X \times_{ \Delta^n } \Delta^{ \{n-i, \ldots, n \} } |_{Q^{\bigdot}},$$ where
the twisted geometric realization functor $| |_{Q^{\bigdot}}$ is as defined in \S \ref{twistt}.
Taking $X = \Un_{\phi} \calF$, we see that $v_i$ fits into a commutative diagram
$$ \xymatrix{ | X \times_{\Delta^n} \{n-i\} |_{Q^{\bigdot} } \ar[r]^{\sim} \ar@{^{(}->}[d]
& | \Sing_{Q^{\bigdot} } \calF(i) |_{Q^{\bigdot}} \ar[d] \\
| X \times_{\Delta^n} \Delta^{ \{n-i, \ldots, n\} } |_{Q^{\bigdot}} \ar[r]^{v_i} &
\calF(i). }$$
Here the upper horizontal map is an isomorphism supplied by Corollary \ref{undef}, 
and the right vertical map is a weak homotopy equivalence by Proposition \ref{remmy33}. Consequently, to prove that the map $v_i$ is a weak homotopy equivalence, it will suffice
to show that the left vertical map is a weak homotopy equivalence. In view of Proposition
\ref{babyy}, it will suffice to show that the inclusion
$$ X \times_{ \Delta^n } \{n-i\} \subseteq X \times_{ \Delta^n} \Delta^{ \{n-i, \ldots, n \} }$$
is a weak homotopy equivalence. In fact, $X \times_{ \Delta^n } \{n-i\}$ is a deformation
retract of $X \times_{ \Delta^n} \Delta^{ \{n-i, \ldots, n \} }$: this follows from
the observation that the projection $X \rightarrow \Delta^n$ is a right fibration
(Proposition \ref{onehalf}).
\end{proof}

It will be convenient to restate Lemma \ref{sticx} in a slightly modified form. First, we need to introduce a bit of terminology.

\begin{definition}\index{gen}{pointwise equivalence}\index{gen}{equivalence!pointwise}\label{scatterbrain}
Suppose given a commutative diagram of simplicial sets
$$ \xymatrix{ X \ar[rr]^{f} \ar[dr]^{p} & & Y \ar[dl]^{q} \\
& S, & }$$
where $p$ and $q$ are right fibrations. We will say that $f$ is a
{\it pointwise equivalence} if, for each vertex $s \in S$, the induced map
$X_{s} \rightarrow Y_{s}$ is a homotopy equivalence of Kan complexes.
\end{definition}

\begin{remark}
In the situation of Definition \ref{scatterbrain}, the following conditions are equivalent:
\begin{itemize}
\item[$(a)$] The map $f$ is a pointwise equivalence of right fibrations over $S$.
\item[$(b)$] The map $f$ is a contravariant equivalence in $(\sSet)_{/S}$. 
\item[$(c)$] The map $f$ is a categorical equivalence of simplicial sets.
\end{itemize}
The equivalence $(a) \Leftrightarrow (b)$ follows from Corollary \ref{prefibchar}
(see below), and the equivalence $(a) \Leftrightarrow (c)$ from Proposition \ref{apple1}.
\end{remark}

\begin{lemma}\label{gottapruve2}
Let $S' \subseteq S$ be simplicial sets. Let
$p: X \rightarrow S$ be any map, and let $q: Y \rightarrow S$ be a right fibration.
Let $X' = X \times_{S} S'$ and $Y' = Y \times_{S} S'$. The restriction map
$$ \phi: \bHom_{ (\sSet)_{/S}} (X, Y) \rightarrow \bHom_{ (\sSet)_{/S'}}( X', Y')$$ is a Kan fibration.
\end{lemma}

\begin{proof}
We first show that $\phi$ is a right fibration. It will suffice to show that
$\phi$ has the right lifting property with respect to every right anodyne inclusion
$A \subseteq B$. This follows from the fact that $q$ has the right lifting property with respect
to the induced inclusion
$$i: (B \times S') \coprod_{ A \times S'} ( A \times S) \subseteq B \times S,$$
since $i$ is again right anodyne (Corollary \ref{prodprod1}). 

Applying the preceding argument to the inclusion $\emptyset \subseteq S'$, we deduce
that the projection map $$\bHom_{ (\sSet)_{/S'}}(X',Y') \rightarrow \Delta^0$$ is a right fibration.
Proposition \ref{greenwich} implies that $\bHom_{ (\sSet)_{/S'}}(X',Y')$ is a Kan complex.
Lemma \ref{toothie2} now implies that $\phi$ is a Kan fibration as desired.
\end{proof}

\begin{lemma}\label{blem}
Let $\calU$ be a collection of simplicial sets. Suppose that:
\begin{itemize}
\item[$(i)$] The collection $\calU$ is stable under isomorphism. That is, if
$S \in \calU$ and $S' \simeq S$, then $S' \in \calU$.
\item[$(ii)$] The collection $\calU$ is stable under the formation of disjoint unions.
\item[$(iii)$] Every simplex $\Delta^n$ belongs to $\calU$.
\item[$(iv)$] Given a pushout diagram
$$ \xymatrix{ X \ar[r] \ar[d]^{f} & X' \ar[d] \\
Y \ar[r] & Y' }$$
in which $X$, $X'$, and $Y$ belong to $\calU$. If the map $f$ is a monomorphism, then
$Y'$ belongs to $\calU$.
\item[$(v)$] Suppose given a sequence of monomorphisms of simplicial sets
$$ X(0) \rightarrow X(1) \rightarrow \ldots $$
If each $X(i)$ belongs to $\calU$, then the colimit $\colim X(i)$ belongs to $\calU$.
\end{itemize}
Then every simplicial set belongs to $\calU$.
\end{lemma}

\begin{proof}
Let $S$ be a simplicial set; we wish to show that $S \in \calU$. In view of $(v)$, it will suffice
to show that each skeleton $\sk^{n} S$ belongs to $\calU$. We may therefore assume that
$S$ is finite dimensional. We now proceed by induction on the dimension $n$ of $S$.
Let $A$ denote the set of nondegenerate $n$-simplexes of $S$, so that we have a pushout
diagram
$$ \xymatrix{ \coprod_{\alpha \in A} \bd \Delta^n \ar[r] \ar[d] & \sk^{n-1} S \ar[d] \\
\coprod_{\alpha \in A} \Delta^n \ar[r] & S. }$$
Invoking assumption $(iv)$, we are reduced to proving that
$\sk^{n-1} S$, $\coprod_{ \alpha \in A} \bd \Delta^n$, and $\coprod_{\alpha \in A} \Delta^n$
belong to $\calU$. For the first two this follows from the inductive hypothesis, and for the
last it follows from assumptions $(ii)$ and $(iii)$.
\end{proof}

\begin{lemma}\label{postcuse}
Suppose given a commutative diagram of simplicial sets
$$ \xymatrix{ X \ar[rr]^{f} \ar[dr]^{p} & & Y \ar[dl]^{q} \\
& S, & }$$
where $p$ and $q$ are right fibrations. The following conditions are equivalent:
\begin{itemize}
\item[$(a)$] The map $f$ is a pointwise equivalence.
\item[$(b)$] The map $f$ is an equivalence in the simplicial category $( \sSet)_{/S}$
$($that is, $f$ admits a homotopy inverse{}$)$.
\item[$(c)$] For every object $A \in (\sSet)_{/S}$, composition with $f$ induces a homotopy
equivalence of Kan complexes $\bHom_{ (\sSet)_{/S}}( A, X) \rightarrow
\bHom_{ (\sSet)_{/S}}(A, Y)$.
\end{itemize}
\end{lemma}

\begin{proof}
The implication $(b) \Rightarrow (a)$ is clear (any homotopy inverse to $f$ determines
homotopy inverses for the maps $f_{s}: X_{s} \rightarrow Y_{s}$, for each vertex $s \in S$),
and the implication $(c) \Rightarrow (b)$ follows from Proposition \ref{rooot}.
We will prove that $(a) \Rightarrow (c)$. Let $\calU$ denote the collection of all simplicial
sets $A$ such that, for {\em every} map $A \rightarrow S$, composition with $f$ induces a homotopy equivalence of Kan complexes
$$\bHom_{ (\sSet)_{/S}}( A, X) \rightarrow \bHom_{ (\sSet)_{/S}}(A, Y).$$
We will show that $\calU$ satisfies the hypotheses of Lemma \ref{blem}, and therefore
contains {\em all} simplicial sets. Conditions $(i)$ and $(ii)$ are obvious, and
conditions $(iv)$ and $(v)$ follow from Lemma \ref{gottapruve2}. It will therefore suffice
to show that every simplex $\Delta^n$ belongs to $\calU$. For every
map $\Delta^n \rightarrow S$, we have a commutative diagram
$$ \xymatrix{ \bHom_{ (\sSet)_{/S}}( \Delta^n, X) \ar[r] \ar[d] & \bHom_{ (\sSet)_{/S}}( \Delta^n, Y) \ar[d] \\
\bHom_{ (\sSet)_{/S} }( \{n\}, X) \ar[r] & \bHom_{ ( \sSet)_{/S} }( \{n\}, Y). }$$
Since the inclusion $\{n\} \subseteq \Delta^n$ is right anodyne, the vertical maps are trivial
Kan fibrations. It will therefore suffice to show that the bottom horizontal map is
a homotopy equivalence, which follows immediately from $(a)$.
\end{proof}

\begin{lemma}\label{prestix}
Let $\phi: \sCoNerve[ \Delta^n] \rightarrow [n]$ be as in Lemma \ref{sticx}.
Suppose given a right fibration $X \rightarrow \Delta^n$, a
projectively fibrant diagram $\calF \in \Set_{\Delta}^{[n]}$, and a weak
equivalence of diagrams $\alpha: \St_{\phi} X \rightarrow \calF$. Then the adjoint map
$X \rightarrow \Un_{\phi} \calF$ is a pointwise equivalence of left fibrations
over $\Delta^n$.
\end{lemma}

\begin{proof}
For $0 \leq i \leq n$, let $X(i) = X \times_{\Delta^n} \Delta^{ \{n-i, \ldots, n \} } \subseteq X$.
We observe that $(\St_{\phi} X)(i)$ is canonically isomorphic to the twisted geometric realization
$| X(i) |_{Q^{\bigdot}}$, where $Q^{\bigdot}$ is defined as in \S \ref{twistt}. 
Since $X \rightarrow \Delta^n$ is a right fibration, the fiber
$X \times_{ \Delta^n } \{i\}$ is a deformation retract of $X(i)$. Using Proposition \ref{babyy},
we conclude that the induced inclusion
$| X \times_{ \Delta^n} \{n-i\} |_{Q^{\bigdot}} \rightarrow | X(i) |_{Q^{\bigdot}}$ is a
weak homotopy equivalence. Since $\alpha$ is a weak equivalence,
we get weak equivalences $| X \times_{ \Delta^n} \{n-i\} |_{Q^{\bigdot}} \rightarrow \calF(i)$
for each $0 \leq i \leq n$. Using Proposition \ref{realremmy}, we deduce that
the adjoint maps $X \times_{ \Delta^n } \{n-i\} \rightarrow \Sing_{Q^{\bigdot}} \calF(i)$
are again weak homotopy equivalences. The desired result now follows from
the observation that $\Sing_{Q^{\bigdot}} \calF(i) \simeq (\Un_{\phi} \calF) \times_{ \Delta^n} \{n-i\}$
(Remark \ref{undef}). 
\end{proof}

\begin{notation}\index{not}{RFibS@$\RFib(S)$}
For every simplicial set $S$, we let $\RFib(S)$ denote the full subcategory
of $( \sSet)_{/S}$ spanned by those maps $X \rightarrow S$ which are right fibrations.
\end{notation}

Proposition \ref{onehalf} implies that if $p: X \rightarrow S$ exhibits $X$ as a fibrant
object of the contravariant model category $(\sSet)_{/S}$, then $p$ is a right fibration. We will
prove the converse below (Corollary \ref{usewhere1}). For the moment, we will be content with the following weaker result:

\begin{lemma}\label{camine}
For every integer $n \geq 0$, the inclusion
$i: (\sSet)_{/\Delta^n}^{\degree} \subseteq \RFib( \Delta^n)$ is an equivalence of simplicial
categories.
\end{lemma}

\begin{proof}
It is clear that $i$ is fully faithful. To prove that $i$ is essentially surjective, consider
any left fibration $X \rightarrow \Delta^n$. Let $\phi: \sCoNerve[ \Delta^n ] \rightarrow [n]$
be defined as in Lemma \ref{sticx}, and choose a weak equivalence
$\St_{\phi} X \rightarrow \calF$, where $\calF \in \Set_{\Delta}^{[n]}$ is a projectively fibrant
diagram. Lemma \ref{prestix} implies that the adjoint map
$X \rightarrow \Un_{\phi} \calF$ is a pointwise equivalence of right fibrations in $\Delta^n$, and
therefore a homotopy equivalence in $\RFib( \Delta^n)$ (Lemma \ref{postcuse}). 
It now suffices to observe that $\Un_{\phi} \calF \in ( \sSet)_{/\Delta^n}^{\degree}$.
\end{proof}

\begin{lemma}\label{stucks}
For each integer $n \geq 0$, the unstraightening functor
$\Un_{\Delta^n}: (\Set_{\Delta}^{\sCoNerve[\Delta^n]})^{\degree} \rightarrow \RFib(\Delta^n)$ is
an equivalence of simplicial categories.
\end{lemma}

\begin{proof}
In view of Lemma \ref{camine} and Proposition \ref{weakcompatequiv}, it will suffice to show that
the Quillen adjunction $( \St_{\Delta^n}, \Un_{\Delta^n} )$ is a Quillen equivalence. This follows
immediately from Lemma \ref{sticx} and Proposition \ref{lesstrick}.
\end{proof}

\begin{proposition}\label{fullmeal}
For every simplicial set $S$, the unstraightening functor 
$\Un_{S}$ induces an equivalence of simplicial categories
$(\Set_{\Delta}^{\sCoNerve[S]^{op}})^{\degree} \rightarrow \RFib(S)$.
\end{proposition}

\begin{proof}
For each simplicial set $S$, let $( \Set^{\sCoNerve[S]^{op}}_{\Delta} )_{f}$ denote the
category of projectively fibrant objects of $\Set^{\sCoNerve[S]^{op}}_{\Delta}$, and let
$W_{S}$ be the class of weak equivalences in $( \Set^{\sCoNerve[S]^{op}}_{\Delta})_{f}$.
Let $W'_{S}$ be the collection of pointwise equivalences in $\RFib(S)$.
We have a commutative diagram of simplicial categories
$$ \xymatrix{ ( \Set^{\sCoNerve[S]^{op}}_{\Delta})^{\degree} \ar[r]^{ \Un_{S} } \ar[d] &
\RFib(S) \ar[d]^{\psi_{S}} \\
( \Set^{\sCoNerve[S]^{op}}_{\Delta})_{f}[W_{S}^{-1} ] \ar[r]^{\phi_S} & \RFib[ {W'}^{-1}_{S} ] }$$
(see Notation \ref{localdef}). Lemma \ref{kur} implies that the left vertical map is an equivalence.
Using Lemma \ref{postcuse} and Remark \ref{uppa}, we deduce that the right vertical
map is also an equivalence. It will therefore suffice to show that $\phi_{S}$ is an equivalence. 

Let $\calU$ denote the collection of simplicial sets $S$ for which $\phi_{S}$ is an equivalence.
We will show that $\calU$ satisfies the hypotheses of Lemma \ref{blem}, and therefore contains every simplicial set $S$. Conditions $(i)$ and $(ii)$ are obviously satisfied, and
condition $(iii)$ follows from Lemma \ref{stucks} and Proposition \ref{weakcompatequiv}.
We will verify condition $(iv)$; the proof of $(v)$ is similar.

Applying Corollary \ref{uspin}, we deduce:
\begin{itemize}
\item[$(\ast)$] The functor $S \mapsto ( \Set^{\sCoNerve[S]^{op}}_{\Delta})_{f} [W_S^{-1}]$ carries homotopy colimit diagrams indexed by a partially ordered set to homotopy limit diagrams in $\sCat$.
\end{itemize}

Suppose given a pushout diagram
$$ \xymatrix{ X \ar[r] \ar[d]^{f} & X' \ar[d] \\
Y \ar[r] & Y' }$$
in which $X, X', Y \in \calU$, where $f$ is a cofibration. We wish to prove that $Y' \in \calU$.
We have a commutative diagram
$$ \xymatrix{ ( \Set^{\sCoNerve[Y']^{op}}_{\Delta} )_{f}[W_{Y'}^{-1}] \ar[r]^{\phi_{Y'}} &
\RFib(Y')[ {W'}_{Y'}^{-1}] \ar[r]^{u} \ar[d]^{v} \ar[dr]^{w} & \RFib(Y)[ {W'}_{Y}^{-1}] \ar[d] \\
& \RFib(X')[ {W'}_{X'}^{-1}] \ar[r] & \RFib(X)[ {W'}_{X}^{-1} ]. }$$
Using $(\ast)$ and Corollary \ref{wspin}, we deduce that $\phi_{Y'}$ is an equivalence if and only if,
for every pair of objects $x,y \in \RFib(Y')[ {W'}_{Y'}^{-1}]$, the diagram
of simplicial sets
$$ \xymatrix{ \bHom_{ \RFib(Y')[ {W'}_{Y'}^{-1}] }( x, y) \ar[r] \ar[d] 
& \bHom_{ \RFib(Y)[ {W'}_{Y}^{-1}] }( u(x), u(y) ) \ar[d] \\
\bHom_{ \RFib(X')[ {W'}_{X'}^{-1}] }( v(x), v(y) ) \ar[r] &
\bHom_{ \RFib(X)[ {W'}_{X}^{-1}] }( w(x), w(y) ) }$$ 
is homotopy Cartesian. Since $\psi_{Y'}$ is a weak equivalence of simplicial categories, 
we may assume without loss of generality that $x = \psi_{Y'}( \overline{x} )$ and
$y = \psi_{Y'}( \overline{y} )$, for some $ \overline{x}, \overline{y} \in (\sSet)^{\degree}_{/Y'}$. 
It will therefore suffice to prove that the equivalent diagram
$$ \xymatrix{ \bHom_{ \RFib(Y') }( \overline{x}, \overline{y}) \ar[r] \ar[d] 
& \bHom_{ \RFib(Y) }( \overline{u}(\overline{x}), \overline{u}(\overline{y}) ) \ar[d] \\
\bHom_{ \RFib(X') }( \overline{v}(\overline{x}), \overline{v}(\overline{y}) ) \ar[r]^{g} &
\bHom_{ \RFib(X) }( \overline{w}(\overline{x}), \overline{w}(\overline{y}) ) }$$ 
is homotopy Cartesian. But this diagram is a pullback square, and the map $g$ is a Kan
fibration by Lemma \ref{gottapruve2}.
\end{proof}

We can now complete the proof of Theorem \ref{struns}.
Suppose that $\phi: \sCoNerve[S] \rightarrow \calC^{op}$
is an equivalence of simplicial categories; we wish to show that the adjoint functors
$( \St_{\phi}, \Un_{\phi})$ determine a Quillen equivalence between
$(\sSet)_{/S}$ and $\Set_{\Delta}^{\calC}$. Using Proposition \ref{lesstrick}, we can reduce to
the case where $\phi$ is an isomorphism. In view of Proposition \ref{weakcompatequiv}, it will
suffice to show that $\Un_{\phi}$ induces an equivalence of simplicial categories
$(\Set^{\sCoNerve[S]^{op}}_{\Delta})^{\degree} \rightarrow (\sSet)_{/S}^{\degree}$, which
follows immediately from Proposition \ref{fullmeal}.

\begin{corollary}\label{usewhere1}
Let $p: X \rightarrow S$ be a map of simplicial sets. The following conditions
are equivalent:
\begin{itemize}
\item[$(1)$] The map $p$ is a right fibration.
\item[$(2)$] The map $p$ exhibits $X$ as a fibrant object
of $(\sSet)_{/S}$ $($with respect to the contravariant model structure$)$.
\end{itemize}
\end{corollary}

\begin{proof}
The implication $(2) \Rightarrow (1)$ follows from Proposition \ref{onehalf}.
For the converse, let us suppose that $p$ is a right fibration. 
Proposition \ref{fullmeal} implies that the unstraightening functor
$\Un_{S}: ( \Set_{\Delta}^{\sCoNerve[S]^{op}})^{\degree} \rightarrow \RFun(S)$
is essentially surjective. Since $\Un_{S}$ factors through the inclusion
$i: ( \sSet)_{/S}^{\degree} \subseteq \RFun(S)$, we deduce that $i$ is essentially surjective.
Consequently, we can choose a simplicial homotopy equivalence $f: X \rightarrow Y$ in $(\sSet)_{/S}$, where $Y$ is fibrant. Let $g$ be a homotopy inverse to $X$, so that there exists a homotopy
$h: X \times \Delta^1 \rightarrow X$ from $\id_{X}$ to $g \circ f$.

To prove that $X$ is fibrant, we must show that every lifting problem
$$ \xymatrix{ A \ar@{^{(}->}[d]^{j} \ar[r]^{e_0} & X \ar[d]^{p} \\
B \ar@{-->}[ur]^{e} \ar[r] & S }$$
has a solution, provided that $j$ is a trivial cofibration in the contravariant model category
$(\sSet)_{/S}$. Since $Y$ is fibrant, the map $f \circ e_0$ can be extended to a map
$\overline{e}: B \rightarrow Y$ in $(\sSet)_{/S}$. Let $e' = g \circ \overline{e}$. The maps
$\overline{e}$ and $h \circ (e_0 \times \id_{\Delta^1})$ determine another lifting problem
$$ \xymatrix{ (A \times \Delta^1) \coprod_{ A \times \{1\} } ( B \times \{1\} ) \ar@{^{(}->}[d]^{j'} \ar[r] & X \ar[d]^{p} \\
B \times \Delta^1 \ar[r] \ar@{-->}[ur]^{E} & S. }$$
Proposition \ref{usejoyal} implies that $j'$ is right anodyne. Since $p$ is a right fibration,
there exists an extension $E$ as indicated in the diagram. The restriction
$e = E | B \times \{0\}$ is then a solution the original problem. 
\end{proof}

\begin{corollary}\label{prefibchar}\index{gen}{covariant!equivalence}
Suppose given a diagram of simplicial sets
$$ \xymatrix{ X \ar[dr]^{p} \ar[rr]^{f} & & Y \ar[dl]^{q} \\
& S & }$$ 
where $p$ and $q$ are right fibrations. Then $f$ is a contravariant equivalence
in $(\sSet)_{/S}$ if and only if $f$ is a pointwise equivalence.
\end{corollary}

\begin{proof}
Since $(\sSet)_{/S}$ is a simplicial model category, this follows immediately
from Corollary \ref{usewhere1} and Lemma \ref{postcuse}.
\end{proof}

Corollary \ref{usewhere1} admits the following generalization:

\begin{corollary}\label{usewhere2}
Suppose given a diagram of simplicial sets
$$ \xymatrix{ X \ar[rr]^{f} \ar[dr]^{p} & & Y \ar[dl]^{q} \\
& S, & }$$
where $p$ and $q$ are right fibrations. Then $f$ is a contravariant fibration in
$(\sSet)_{/S}$ if and only if $f$ is a right fibration.
\end{corollary}

\begin{proof}
The map $f$ admits a factorization
$$ X \stackrel{f'}{\rightarrow} X' \stackrel{f''}{\rightarrow} Y$$
where $f'$ is a contravariant equivalence and $f''$ is a contravariant fibration (in $(\sSet)_{/S}$). 
Proposition \ref{onehalf} implies that $f''$ is a right fibration, so the composition
$q \circ f''$ is a right fibration. Invoking Corollary \ref{prefibchar}, we conclude that
for every vertex $s \in S$, the map $f'$ induces a homotopy equivalence of fibers
$X_{s} \rightarrow X'_{s}$. Consider the diagram
$$ \xymatrix{ X_{s} \ar[rr]^{f'_s} \ar[dr] & & X'_{s} \ar[dl] \\
& Y_{s}. & }$$
The vertical maps in this diagram are right fibrations between Kan complexes, and therefore
Kan fibrations (Lemma \ref{toothie2}). Since $f_{s}$ is a homotopy equivalence, we conclude
that the induced map of fibers $f'_{y}: X_{y} \rightarrow X'_{y}$ is a homotopy equivalence
for each vertex $y \in Y$. Invoking Lemma \ref{postcuse}, we deduce that
$f'$ is an equivalence in the simplicial category $(\sSet)_{/Y}$.

We can now repeat the proof of Corollary \ref{usewhere1}. Let $g$ be a homotopy
inverse to $f'$ in the simplicial category $(\sSet)_{/Y}$, and let
$h: X \times \Delta^1 \rightarrow X$ be a homotopy from $\id_{X}$ to $g \circ f'$
(which projects to the identity on $Y$). 
To prove that $f$ is a covariant fibration, we must show that every lifting problem
$$ \xymatrix{ A \ar@{^{(}->}[d]^{j} \ar[r]^{e_0} & X \ar[d]^{f} \\
B \ar@{-->}[ur]^{e} \ar[r] & Y }$$
has a solution, provided that $j$ is a trivial cofibration in the contravariant model category
$(\sSet)_{/S}$. Since $f''$ is a contravariant fibration, the map $f' \circ e_0$ can be extended to a map
$\overline{e}: B \rightarrow X'$ in $(\sSet)_{/Y}$. Let $e' = g \circ \overline{e}$. The maps
$\overline{e}$ and $h \circ (e_0 \times \id_{\Delta^1})$ determine another lifting problem
$$ \xymatrix{ (A \times \Delta^1) \coprod_{ A \times \{1\} } ( B \times \{1\} ) \ar@{^{(}->}[d]^{j'} \ar[r] & X \ar[d]^{f} \\
B \times \Delta^1 \ar[r] \ar@{-->}[ur]^{E} & Y. }$$
Proposition \ref{usejoyal} implies that $j'$ is right anodyne. Since $f$ is a right fibration,
there exists an extension $E$ as indicated in the diagram. The restriction
$e = E | B \times \{0\}$ is then a solution the original problem. 
\end{proof}

We conclude this section with one more result which will be useful
in studying the Joyal model structure on $\sSet$. 
Suppose that $f: X \rightarrow S$ is any map of simplicial sets, and $\{s\}$ is a vertex of
$S$. Let $Q^{\bigdot}$ denote the cosimplicial object of $\sSet$ defined in
\S \ref{twistt}. Then we have a canonical map
$$| X_{s} |_{Q^{\bigdot}} \simeq (\St_{ \{s\} } X_{s})(s)
\rightarrow (\St_{S} X)(s).$$

\begin{proposition}\label{canuble}
Suppose that $f: X \rightarrow S$ is a right fibration of simplicial sets. Then for each vertex
$s$ of $S$, the canonical map
$\phi: | X_{s} |_{ Q^{\bigdot} } \rightarrow (\St_{S} X)(s)$ is a weak homotopy equivalence of simplicial sets.
\end{proposition}

\begin{proof}
Choose a weak equivalence $\St_{S} X \rightarrow \calF$, where
$\calF: \sCoNerve[S]^{op} \rightarrow \sSet$ is a projectively fibrant diagram.
Theorem \ref{struns} implies that the adjoint map
$X \rightarrow \Un_{S}( \calF)$ is a contravariant equivalence in
$(\sSet)_{/S}$. Applying the ``only if'' direction of Corollary \ref{usewhere1}, we conclude
that each of the induced maps
$$X_{s} \rightarrow (\Un_{S} \calF)_{s} \simeq \Sing_{ Q^{\bigdot} } \calF(s)$$ 
is a homotopy equivalence of Kan complexes. Using Proposition \ref{realremmy}, we deduce that the adjoint map $| X_{s} |_{Q^{\bigdot} } \rightarrow \calF(s)$ is a weak homotopy equivalence. 
It follows from the two-out-of-three property that $\phi$ is also a weak homotopy equivalence.
\end{proof}

\subsection{The Comparison Theorem}\label{compp2}

Let $S$ be an $\infty$-category containing a pair of objects $x$ and $y$, and let
$Q^{\bigdot}$ denote the cosimplicial object of $\sSet$ described in \S \ref{twistt}. 
We have a canonical map of simplicial sets
$$ f: | \Hom^{\rght}_{S}(x,y) |_{Q^{\bigdot}} \rightarrow \bHom_{\sCoNerve[S]}(x,y).$$
Moreover, in the special case where $S$ is the nerve of a fibrant simplicial category $\calC$, the composition
$$ | \Hom^{\rght}_{S}(x,y) |_{Q^{\bigdot}} \stackrel{f}{\rightarrow} \bHom_{\sCoNerve[S]}(x,y) \rightarrow \bHom_{\calC}(x,y)$$ can be identified with
the counit map 
$$ | \Sing_{Q^{\bigdot}} \bHom_{\calC}(x,y) |_{Q^{\bigdot}} \rightarrow \bHom_{\calC}(X,Y),$$
and is therefore a weak equivalence (Proposition \ref{remmy33}). 
Consequently, we may reformulate Theorem \ref{biggiesimp} in the following way:

\begin{proposition}\label{wiretrack}
Let $S$ be an $\infty$-category containing a pair of objects $x$ and $y$. Then the natural map
$$ f: | \Hom^{\rght}_{S}(x,y) |_{Q^{\bigdot}} \rightarrow \bHom_{\sCoNerve[S]}(x,y)$$
is a weak homotopy equivalence of simplicial sets. 
\end{proposition}

\begin{proof}
Let $C = S_{/y}^{\triangleright} \coprod_{ S_{/y} } S$, and let $v$ denote the image in $C$ of the cone point of $S_{/y}^{\triangleright}$. There is a canonical projection $\pi: C \rightarrow S$, which induces a map of simplicial sets 
$$f'': (\St_{S} S_{/y})(x) \rightarrow \bHom_{ \sCoNerve[S]}(x,y).$$
The map $f$ can be identified with the composition $f'' \circ f'$, where
$f'$ is the map
$$ | \Hom^{\rght}_{S}(x,y)|_{Q^{\bigdot}} \simeq 
(\St_{ \{x\} } S_{/y} \times_{S} \{x\})(x) \rightarrow (\St_{S} S_{/y})(y).$$
Since the projection $S_{/y} \rightarrow S$ is a right fibration, the map
$f'$ is a weak homotopy equivalence (Proposition \ref{canuble}). It will therefore suffice to show that
$f''$ is a weak homotopy equivalence. To see this, we consider the commutative diagram
$$ \xymatrix{ & (\St_{S} S_{/y})(x) \ar[dr]^{f''} & \\
(\St_{S} \{y\})(x) \ar[ur]^{g} \ar[rr]^{h} & & \bHom_{ \sCoNerve[S]}(x,y). }$$
The inclusion $i: \{ y\} \subseteq S_{/y}$ is a retract of the inclusion
$$ (S_{/y} \times \{1\} ) \coprod_{ \{ y \} \times \{1\} }
( \{ \id_{y} \} \times \Delta^1) \subseteq S_{/y} \times \Delta^1,$$
which is right anodyne by Corollary \ref{prodprod1}. It follows that $i$ is a contravariant equivalence
in $(\sSet)_{/S}$ (Proposition \ref{onehalf}), so the map $g$ is a weak homotopy equivalence of simplicial sets. Since the map $h$ is an isomorphism, the map $f''$ is also a weak homotopy equivalence by virtue of the two-out-of-three property.
\end{proof}

\subsection{The Joyal Model Structure}\label{compp3}

The category of simplicial sets can be endowed with a model structure for which the fibrant objects are precisely the $\infty$-categories. The original construction of this model structure is due to Joyal, who uses purely combinatorial arguments (\cite{joyalnotpub}). In this section, 
we will exploit the relationship between simplicial categories and
$\infty$-categories to give an alternative description of this model structure. 
Our discussion will make use of a model structure on the category $\sCat$ of simplicial categories, which we review in \S \ref{compp4}.

\begin{theorem}\label{biggier}\index{gen}{model category!Joyal}\index{gen}{simplicial set!Joyal model structure}
There exists a left proper, combinatorial model structure on
the category of simplicial sets with the following properties:
\begin{itemize}
\item[$(C)$] A map $p: S \rightarrow S'$ of simplicial sets is a {\it
cofibration} if and only if it is a monomorphism.

\item[$(W)$] A map $p: S \rightarrow S'$ is a {\it categorical
equivalence} if and only if the induced simplicial functor
$\sCoNerve[S] \rightarrow \sCoNerve[S']$ is an equivalence of simplicial categories.
\end{itemize}

Moreover, the adjoint functors $(\sCoNerve, \sNerve)$
determine a Quillen equivalence between $\sSet$ $($with the model structure defined above$)$ and $\sCat$.
\end{theorem}

Our proof will make use of the theory of {\em inner anodyne} maps of simplicial sets,
which we will study in detail in \S \ref{midfibsec}. We first establish a simple Lemma.

\begin{lemma}\label{rugg}
Every inner anodyne map $f: A \rightarrow B$ of simplicial sets is a categorical equivalence.
\end{lemma}

\begin{proof}
It will suffice to prove that if $f$ is inner anodyne, then the associated map
$\sCoNerve[f]$ is a trivial cofibration of simplicial categories. The collection of all
morphisms $f$ for which this statement holds is weakly saturated (Definition \ref{saturated}). Consequently, we may assume that $f$ is an inner horn inclusion $\Lambda^n_i \subseteq \Delta^n$, $0 < i < n$.
We now explicitly describe the map $\sCoNerve[f]$:

\begin{itemize}
\item The objects of $\sCoNerve[\bd \Lambda^n_i]$ are the objects of
$\sCoNerve[\Delta^n]$: namely, elements of the linearly ordered
set $[n] = \{0, \ldots, n\}$. 

\item For $0 \leq j \leq k \leq n$, the simplicial set
$\bHom_{\sCoNerve[\Lambda^n_i]}(j,k)$ is equal to
$\bHom_{\sCoNerve[\Delta^n]}(j,k)$ unless $j=0$ and $k=n$. In the
latter case, $$\bHom_{\sCoNerve[\Lambda^n_i]}(j,k) = K \subseteq
(\Delta^1)^{n-1} \simeq \bHom_{\sCoNerve[\Delta^n]}(j,k),$$
where $K$ is the simplicial subset of the cube $(\Delta^1)^{n-1}$ obtained
by removing the interior and a single face.
\end{itemize}

We observe that $\sCoNerve[f]$ is a pushout of the inclusion $\calE_{K} \subseteq
\calE_{ (\Delta^1)^{n-1} }$ (see \S \ref{compp4} for an explanation of this notation).
It now suffices to observe that the inclusion $K \subseteq (\Delta^1)^{n-1}$ is trivial fibration of simplicial sets (with respect to the usual model structure on $\sSet$).
\end{proof}

\begin{proof}[Proof of Theorem \ref{biggier}]
We first show that $\sCoNerve$ carries cofibrations of simplicial sets to cofibrations of simplicial categories. Since the class of all cofibrations of simplicial
sets is generated by the inclusions $\bd \Delta^n \subseteq
\Delta^n$, it suffices to show that each map $\sCoNerve[\bd
\Delta^n] \rightarrow \sCoNerve[\Delta^n]$ is a cofibration of
simplicial categories. If $n = 0$, then the inclusion $\sCoNerve[\bd \Delta^n] \subseteq
\sCoNerve[\Delta^n]$ is isomorphic to the inclusion $\emptyset
\subseteq \ast$ of simplicial categories, which is a
cofibration. In the case where $n
> 0$, we make use of the following explicit description of
$\sCoNerve[\bd \Delta^n]$ as a subcategory of
$\sCoNerve[\Delta^n]$:

\begin{itemize}
\item The objects of $\sCoNerve[\bd \Delta^n]$ are the objects of
$\sCoNerve[\Delta^n]$: namely, elements of the linearly ordered
set $[n] = \{0, \ldots, n\}$. 

\item For $0 \leq j \leq k \leq n$, the simplicial set
$\Hom_{\sCoNerve[\bd \Delta^n]}(j,k)$ is equal to
$\Hom_{\sCoNerve[\Delta^n]}(j,k)$ unless $j=0$ and $k=n$. In the
latter case, $\Hom_{\sCoNerve[\bd \Delta^n]}(j,k)$ consists of the
boundary of the cube 
$$(\Delta^1)^{n-1} \simeq
\Hom_{\sCoNerve[\Delta^n]}(j,k).$$
\end{itemize}

In particular, the inclusion $\sCoNerve[\bd \Delta^n] \subseteq
\sCoNerve[\Delta^n]$ is a pushout of the inclusion $\calE_{ \bd
(\Delta^1)^{n-1} } \subseteq \calE_{ (\Delta^1)^{n-1} }$, which is
a cofibration of simplicial categories (see \S \ref{compp4} for an explanation of our notation).

We now declare that a map $p: S \rightarrow S'$ of simplicial sets is a {\it categorical
fibration}\index{gen}{categorical fibration}\index{gen}{fibration!categorical} if it has the right lifting property with respect to
all maps which are cofibrations and categorical equivalences. We now claim that the cofibrations, categorical equivalences, and categorical fibrations determine a left proper, combinatorial model structure on $\sSet$. To prove this, it will suffice to show that the hypotheses of Proposition \ref{goot} are satisfied:

\begin{itemize}
\item[$(1)$] The class of categorical equivalences in $\sSet$ is perfect. This follows from Corollary \ref{perfpull}, since the functor $\sCoNerve$ preserves filtered colimits, and the class of equivalences between simplicial categories is perfect. 
\item[$(2)$] The class of categorical equivalences is stable under pushouts by cofibrations.
Since $\sCoNerve$ preserves cofibrations, this follows immediately from the left-properness of $\sCat$.
\item[$(3)$] A map of simplicial sets which has the right lifting property with respect to {\em all} cofibrations is a categorical equivalence. In other words, we must show that if $p: S \rightarrow S'$ is a trivial fibration of simplicial sets, then
the induced functor $\sCoNerve[p]: \sCoNerve[S] \rightarrow
\sCoNerve[S']$ is an equivalence of simplicial categories.

Since $p$ is a trivial fibration, it admits a section $s: S'
\rightarrow S$. It is clear that $\sCoNerve[p] \circ \sCoNerve[s]$
is the identity; it therefore suffices to show that $$\sCoNerve[s]
\circ \sCoNerve[p]: \sCoNerve[S] \rightarrow \sCoNerve[S]$$ is
homotopic to the identity.

Let $K$ denote the simplicial set $\bHom_{S'}(S,S)$. Then $K$ is a
contractible Kan complex, containing points $x$ and $y$ which
classify $s \circ p$ and $\id_S$. We note the existence of a
natural ``evaluation map'' $e: K \times S \rightarrow S$, such
that $s \circ p = e \circ (\{x\} \times \id_S)$, $\id_S = e \circ
( \{y\} \times \id_S )$. It therefore suffices to show that the
functor $\sCoNerve$ carries $\{x\} \times \id_{S}$ and $\{y\}
\times \id_{S}$ into homotopic morphisms. Since both of these maps
section the projection $K \times S \rightarrow S$, it suffices to
show that the projection $\sCoNerve[K \times S] \rightarrow
\sCoNerve[S]$ is an equivalence of simplicial categories. Replacing $S$ by $S \times K$
and $S'$ by $S$, we are reduced to the
special case where $S = S' \times K$ and $K$ is a contractible
Kan complex.

By the small object argument, we can find an inner anodyne map $S' \rightarrow V$, where $V$ is an $\infty$-category. The corresponding map $S' \times K \rightarrow V
\times K$ is also inner anodyne (Proposition \ref{usejoyal2}), so the maps
$\sCoNerve[S'] \rightarrow \sCoNerve[V]$ and $\sCoNerve[S' \times
K] \rightarrow \sCoNerve[V \times K]$ are both trivial
cofibrations (Lemma \ref{rugg}). It follows that we are free to replace $S'$ by $V$
and $S$ by $V \times K$. In other words, we may suppose that $S'$
is an $\infty$-category (and now we will have no further need of the
assumption that $S$ is isomorphic to the product $S' \times K$).

Since $p$ is surjective on vertices, it is clear that
$\sCoNerve[p]$ is essentially surjective. It therefore suffices to
show that for every pair of vertices $x,y \in S_0$, the induced map of simplicial
sets $\bHom_{\sCoNerve[S]}(x,y) \rightarrow
\bHom_{\sCoNerve[S']}(p(x),p(y))$ is a weak homotopy equivalence. Using Propositions
\ref{wiretrack} and \ref{babyy}, it suffices to show that the map
$\Hom^{\rght}_{S}(x,y) \rightarrow \Hom^{\rght}_{S'}(p(x),p(y))$ is a weak homotopy equivalence. 
This map is obviously a trivial fibration if $p$ is a trivial fibration. 
\end{itemize}

By construction, the functor $\sCoNerve$ preserves weak equivalences. We verified above that $\sCoNerve$ preserves cofibrations as well. It follows that the adjoint functors $(\sCoNerve, \Nerve)$ determine a Quillen adjunction
$$ \Adjoint{\sCoNerve}{\sSet}{\sCat}{\Nerve}.$$
To complete the proof, we  wish to show that this Quillen adjunction is a Quillen equivalence. According to Proposition \ref{quilleq}, we must show that for every simplicial set $S$ and every {\em fibrant} simplicial category $\calC$, a map
$$ u: S \rightarrow \Nerve(\calC)$$
is a categorical equivalence if and only if the adjoint map
$$ v: \sCoNerve[S] \rightarrow \calC$$
is an equivalence of simplicial categories. We observe that $v$ factors as a composition
$$ \sCoNerve[S] \stackrel{ \sCoNerve[u] }{\rightarrow} \sCoNerve[ \Nerve(\calC)]
\stackrel{w}{\rightarrow} \calC.$$
By definition, $u$ is a categorical equivalence if and only if $\sCoNerve[u]$ is an
equivalence of simplicial categories. We now conclude by observing that the counit
map $w$ is an equivalence of simplicial categories (Theorem \ref{biggiesimp}).
\end{proof}

We now establish a few pleasant properties enjoyed by the Joyal model structure on $\sSet$.
We first note that every object of $\sSet$ is cofibrant; in particular, the Joyal model structure is {\em left proper} (Proposition \ref{propob}).

\begin{remark}\label{rightprop}
The Joyal model structure is {\em not} right proper. To see this,
we note that the inclusion $\Lambda^2_1 \subseteq \Delta^2$ is a
categorical equivalence, but it does not remain so after pulling
back via the fibration $\Delta^{ \{0,2\} } \subseteq \Delta^2$.
\end{remark}

\begin{corollary}\label{equivstable}\index{gen}{categorical equivalence!and products}
Let $f: A \rightarrow B$ be a categorical equivalence of simplicial sets, and $K$ an arbitrary simplicial set. Then the induced map $A \times K \rightarrow B \times K$ is a categorical equivalence.
\end{corollary}

\begin{proof}
Choose an inner anodyne map $B \rightarrow Q$, where $Q$ is an $\infty$-category. Then $B \times K \rightarrow Q \times K$ is also inner anodyne, hence a categorical equivalence (Lemma \ref{rugg}). 
It therefore suffices to prove that $A \times K \rightarrow Q \times K$ is a categorical equivalence. In other words, we may suppose to begin with that $B$ is an $\infty$-category.

Now choose a factorization $A \stackrel{f'}{\rightarrow} R \stackrel{f''}{\rightarrow} B$ where $f'$ is an inner anodyne map and $f''$ is an inner fibration. Since $B$ is an $\infty$-category, $R$ is an $\infty$-category.
The map $A \times K \rightarrow R \times K$ is inner anodyne (since $f'$ is), and therefore a categorical equivalence; consequently, it suffices to show that $R \times K \rightarrow B \times K$ is a categorical equivalence. In other words, we may reduce to the case where $A$ is also an $\infty$-category.

Choose an inner anodyne map $K \rightarrow S$, where $S$ is an $\infty$-category. Then $A \times K \rightarrow A \times S$ and $B \times K \rightarrow B \times S$ are both inner anodyne, and therefore categorical equivalences. Thus, to prove that $A \times K \rightarrow B \times K$ is a categorical equivalence, it suffices to show that $A \times S \rightarrow B \times S$ is a categorical equivalence. In other words, we may suppose that $K$ is an $\infty$-category.

Since $A$ and $K$ are $\infty$-categories, $\h{(A \times K)} \simeq \h{A} \times \h{K}$; similarly
$\h{(B \times K)} \simeq \h{B} \times \h{K}$. It follows that $A \times K \rightarrow B \times K$ is essentially surjective, provided that $f$ is essentially surjective. Furthermore, for any pair of vertices $(a,k), (a',k') \in (A \times K)_0$, we have 
$$\Hom^{\rght}_{A \times K}( (a,k), (a',k') ) \simeq \Hom^{\rght}_{A}(a,a') \times \Hom^{\rght}_K(k,k')$$
$$\Hom^{\rght}_{B \times K}( (f(a),k), (f(a'),k')) \simeq \Hom^{\rght}_{B}(f(a),f(a')) \times \Hom^{\rght}_K(k,k').$$
It follows that $A \times K \rightarrow B \times K$ is fully faithful, provided that $f$ is fully faithful, which completes the proof.
\end{proof}

\begin{remark}\label{tokenn}
Since every inner anodyne map is a categorical equivalence, it follows that every categorical fibration $p: X \rightarrow S$ is a inner fibration (see Definition \ref{fibdeff}). The converse is false in general; however, it is true when $S$ is a point. In other words, the fibrant objects for the Joyal model structure on $\sSet$ are precisely the $\infty$-categories. The proof will be given in \S \ref{slin}, as Theorem \ref{joyalcharacterization}. We will assume this result for the remainder of the section. No circularity will result from this, since the proof of Theorem \ref{joyalcharacterization} will not use any of the results proven below.
\end{remark}

The functor $\sCoNerve[\bigdot]$ does not generally commute with products.
However, Corollary \ref{equivstable} implies that $\sCoNerve$ commutes with products in the following weak sense:

\begin{corollary}\label{prodcom}
Let $S$ and $S'$ be simplicial sets. The natural map
$$\sCoNerve[S \times S'] \rightarrow \sCoNerve[S] \times
\sCoNerve[S']$$ is an equivalence of simplicial categories.
\end{corollary}

\begin{proof}
Suppose first that there are fibrant simplicial categories
$\calC$, $\calC'$ with $S = \sNerve(\calC)$, $S' =
\sNerve(\calC')$. In this case, we have a diagram
$$ \sCoNerve[S \times S'] \stackrel{f}{\rightarrow} \sCoNerve[S] \times
\sCoNerve[S'] \stackrel{g}{\rightarrow} \calC \times \calC'.$$
By the two-out-of-three property, it suffices to show that $g$ and
$g \circ f$ are equivalences. Both of these assertions follow
immediately from the fact that the counit map
$\sCoNerve[\sNerve(\calD)] \rightarrow \calD$ is an equivalence for
{\em any} fibrant simplicial category $\calD$ (Theorem \ref{biggier}).

In the general case, we may choose categorical equivalences $S
\rightarrow T$, $S' \rightarrow T'$, where $T$ and $T'$ are nerves of fibrant simplicial categories. Since $S \times S'
\rightarrow T \times T'$ is a categorical equivalence, we reduce
to the case treated above.
\end{proof}

Let $K$ be a fixed simplicial set, and let $\calC$ be a simplicial set which is fibrant with respect to the Joyal model structure. 
Then $\calC$ has the extension property with respect to all inner anodyne maps, and is therefore a
$\infty$-category. It follows that $\Fun(K,\calC)$ is also an $\infty$-category.
We might call two morphisms $f,g: K
\rightarrow \calC$ {\it homotopic} if they are equivalent when viewed
as objects of $\Fun(K,\calC)$. On the other hand, the general
theory of model categories furnishes another notion of homotopy:
$f$ and $g$ are {\it left homotopic} if the map
$$ f \coprod g: K \coprod K \rightarrow \calC$$
can be extended over a mapping cylinder $I$ for $K$.

\begin{proposition}\label{fruity!}
Let $\calC$ be a $\infty$-category and $K$ an arbitrary simplicial set. A pair of morphisms $f,g: K \rightarrow \calC$ are homotopic if and only if they are left-homotopic.
\end{proposition}

\begin{proof}
Choose a contractible Kan complex $S$ containing a pair of distinct vertices, $x$ and $y$.
We note that the inclusion
$$ K \coprod K \simeq K \times \{x,y\} \subseteq K \times S$$
exhibits $K \times S$ as a mapping cylinder for $K$. It follows that $f$ and $g$ are left homotopic
if and only if the map $f \coprod g: K \coprod K \rightarrow \calC$ admits an extension to $K \times S$. 
In other words, $f$ and $g$ are left homotopic if and only if there exists a map $h: S \rightarrow \calC^K$ such that $h(x)=f$ and $h(y)=g$. We note that any such map factors through $Z$, where $Z \subseteq \Fun(K,\calC)$ is the largest Kan complex contained in $\calC^K$. Now, by classical homotopy theory, the map $h$ exists if and only if $f$ and $g$ belong to the same path component of $Z$. It is clear that this holds if and only if $f$ and $g$ are equivalent when viewed as objects of the $\infty$-category $\Fun(K,\calC)$. 
\end{proof}

We are now in a position to prove Proposition \ref{tyty}, which was asserted without proof in
\S \ref{funcback}. We first recall the statement.

\begin{proposition2}\index{gen}{$\infty$-category!of functors}
Let $K$ be an arbitrary simplicial set.
\begin{itemize}
\item[$(1)$] For every $\infty$-category $\calC$, the simplicial set $\Fun(K,\calC)$ is an $\infty$-category.

\item[$(2)$] Let $\calC \rightarrow \calD$ be a categorical equivalence of $\infty$-categories. Then the induced map $\Fun(K,\calC) \rightarrow \Fun(K,\calD)$ is a categorical equivalence.

\item[$(3)$] Let $\calC$ be an $\infty$-category, and $K \rightarrow K'$ a categorical equivalence of simplicial sets. Then the induced map $\Fun(K',\calC) \rightarrow \Fun(K,\calC)$ is a categorical equivalence.
\end{itemize}
\end{proposition2}

\begin{proof}
We first prove $(1)$. To show that $\Fun(K,\calC)$ is an $\infty$-category, it
suffices to show that it has the extension property with respect
to every inner anodyne inclusion $A \subseteq B$. This is equivalent
to the assertion that $\calC$ has the right lifting property with
respect to the inclusion $A \times K \subseteq B \times K$. But
$\calC$ is an $\infty$-category and $A \times K \subseteq B \times K$ is
inner anodyne (Corollary \ref{prodprod2}).

Let $\h{\sSet}$ denote the homotopy category of $\sSet$, taken with respect to the Joyal model structure. For each simplicial set $X$, we let $[X]$ denote the same simplicial set, considered
as an object of $\h{\sSet}$. For every pair of objects $X, Y \in \sSet$, $[X \times Y]$ is a product
for $[X]$ and $[Y]$ in $\h{\sSet}$. This is a general fact when $X$
and $Y$ are fibrant; in the general case, we choose fibrant
replacements $X \rightarrow X'$, $Y \rightarrow Y'$, and apply the fact that
the canonical map $X \times Y \rightarrow X' \times Y'$ is a categorical equivalence (Proposition \ref{fruity!}). 

If $\calC$ is an $\infty$-category, then $\calC$ is a fibrant object
of $\sSet$ (Theorem \ref{joyalcharacterization}). Proposition \ref{fruity!} allows us to identify
$ \Hom_{ \h{\sSet}}( [X], [\calC] )$ with the set of equivalence classes of objects in
the $\infty$-category $\Fun(X, \calC)$. In particular, we have a canonical bijections
$$ \Hom_{\h{\sSet}}( [X] \times [K], [\calC] ) 
\simeq \Hom_{ \h{\sSet}}( [X \times K] , [\calC] )
\simeq \Hom_{\h{ \sSet}}( [X], [\Fun(K,\calC)]).$$ 
It follows that $[ \Fun(K, \calC) ]$ is determined up to canonical isomorphism by
$[K]$ and $[\calC]$ (more precisely, it is an {\em exponential} $[\calC]^{[K]}$ in the homotopy
category $\h{\sSet}$), which proves $(2)$ and $(3)$.
\end{proof}

Our description of the Joyal model structure on $\sSet$ is different from the definition given in \cite{joyalnotpub}. Namely, Joyal defines a map $f: A \rightarrow B$ to be a
{\it weak categorical equivalence} if, for every $\infty$-category $\calC$,
the induced map
$$ \h{\Fun(B, \calC)} \rightarrow \h{\Fun(A, \calC)}$$
is an equivalence (of ordinary categories). To prove that our definition agrees with his, it will suffice to prove the following.\index{gen}{categorical equivalence!weak}

\begin{proposition}\label{joyaldef}
Let $f: A \rightarrow B$ be a map of simplicial sets. Then $f$ is a categorical equivalence if and only if it is a weak categorical equivalence.
\end{proposition}

\begin{proof}
Suppose first that $f$ is a categorical equivalence. If $\calC$ is an arbitrary $\infty$-category, Proposition \ref{tyty} implies that the induced map $\Fun(B,\calC) \rightarrow \Fun(A,\calC)$ is a categorical equivalence, so that $\h{\Fun(B,\calC)}\rightarrow \h{\Fun(A, \calC)}$ is an equivalence of categories. This proves that $f$ is a weak categorical equivalence.

Conversely, suppose that $f$ is a weak categorical equivalence. We wish to show that $f$ induces an isomorphism in the homotopy category of $\sSet$ with respect to the Joyal model structure. It will suffice to show that for any fibrant object $\calC$, $f$ induces a bijection $[B,\calC] \rightarrow [A,\calC]$, where $[X,\calC]$ denotes the set of homotopy classes of maps from $X$ to $\calC$. By Proposition \ref{fruity!}, $[X,\calC]$ may be identified with the set of isomorphism classes of objects in the category $\h{\Fun(X, \calC)}$. By assumption, $f$ induces an equivalence of categories 
$\h{\Fun(B, \calC)} \rightarrow \h{\Fun(A, \calC)}$, and therefore a bijection on isomorphism classes of objects.
\end{proof}

\begin{remark}
The proof of Proposition \ref{tyty} makes use of Theorem \ref{joyalcharacterization}, which asserts that the (categorically) fibrant objects of $\sSet$ are precisely the $\infty$-categories. Joyal proves the analogous assertion for his model structure in \cite{joyalnotpub}. We remark that one cannot formally deduce Theorem \ref{joyalcharacterization} from Joyal's result, since we {\em need} Theorem  \ref{joyalcharacterization} to prove that Joyal's model structure coincides with the one we have defined above. On the other hand, our approach {\em does} give a new proof of Joyal's theorem.
\end{remark}

\begin{remark}
Proposition \ref{joyaldef} permits us to define the Joyal model structure without reference to the theory of simplicial categories (this is Joyal's original point of view \cite{joyalnotpub}). Our approach 
is less elegant, but allows us to easily compare the theory of $\infty$-categories with other models of higher category theory, such as simplicial categories. There is another approach to obtaining comparison results, due to To\"{e}n. In \cite{toenchar}, he shows that if $\calC$ is a model category equipped with a cosimplicial object  $C^{\bigdot}$ satisfying certain conditions, then $\calC$ is (canonically) Quillen equivalent to Rezk's category of complete Segal spaces.
To\"{e}n's theorem applies in particular when $\calC$ is the category of simplicial sets, and
$C^{\bigdot}$ is the ``standard simplex'' $C^n = \Delta^n$.
In fact, $\sSet$ is in some sense universal with respect to this property, since it is generated by $C^{\bigdot}$ under colimits and the class of categorical equivalences is dictated by To\"{e}n's axioms. We refer the reader to \cite{toenchar} for details.
\end{remark}

\section{Inner Fibrations}\label{midfibsec}

\setcounter{theorem}{0}

In this section, we will study the theory of {\em inner fibrations} between simplicial sets. The meaning of this notion is somewhat difficult to motivate, because it has no counterpart in classical category theory: Proposition \ref{ruko} implies that {\em every} functor between ordinary categories $\calC \rightarrow \calD$ induces an inner fibration of nerves $\Nerve(\calC) \rightarrow \Nerve(\calD)$.

In the case where $S$ is a point, a map $p: X \rightarrow S$ is an inner fibration if and only if
$X$ is an $\infty$-category. Moreover, the class of inner fibrations is stable under base change: if 
$$\xymatrix{ X' \ar[d]^{p'} \ar[r] & X \ar[d]^{p} \\ 
S' \ar[r] & S }$$ 
is a pullback diagram of simplicial sets and $p$ is an inner fibration, then so is $p'$.
It follows that if $p: X \rightarrow S$ is an arbitrary inner fibration, then each fiber
$X_{s} = X \times_{S} \{s\}$ is an $\infty$-category. We may therefore think of $p$ as encoding a family of $\infty$-categories parametrized by $S$. However, the fibers $X_{s}$ depend functorially on $s$ only in a very weak sense.

\begin{example}
Let $F: \calC \rightarrow \calC'$ be a functor between ordinary categories. Then the map $\Nerve(\calC) \rightarrow \Nerve(\calC')$ is an inner fibration. Yet the fibers
$\Nerve(\calC)_{C} = \Nerve( \calC \times_{\calC'} \{C\} )$ and $\Nerve(\calC)_{D} = \Nerve( \calC \times_{ \calC'} \{D\} )$ over objects $C,D \in \calC'$ can have wildly different properties, even if $C$ and $D$ are isomorphic objects of $\calC'$.
\end{example}

In order to describe how the different fibers of an inner fibration are related to one another, we will introduce the notion of a {\it correspondence} between $\infty$-categories. We review the classical theory of correspondences in \S \ref{corresp}, and explain how to generalize this theory to the $\infty$-categorical setting. 
 
In \S \ref{joyalpr}, we will prove that the class of inner anodyne maps is stable under smash products with arbitrary cofibrations between simplicial sets. As a consequence, we will deduce that the class of inner fibrations (and hence the class of $\infty$-categories) is stable under the formation of mapping spaces.

In \S \ref{minin}, we will study the theory of {\em minimal} inner fibrations, a generalization of Quillen's theory of minimal Kan fibrations. In particular, we will define a class of minimal $\infty$-categories and show that every $\infty$-category $\calC$ is (categorically) equivalent to a minimal $\infty$-category $\calC'$, where $\calC'$ is well-defined up to (noncanonical) isomorphism. We will apply this theory in \S \ref{ncats} to develop a theory of $n$-categories for $n < \infty$.

\subsection{Correspondences}\label{corresp}

Let $\calC$ and $\calC'$ be categories. A {\it correspondence} from $\calC$ to $\calC'$ is a functor
$$ M: \calC^{op} \times \calC' \rightarrow \Set.$$\index{gen}{correspondence!between categories}
If $M$ is a correspondence from $\calC$ to $\calC'$, we can define a new category $\calC \star^{M} \calC'$ as follows. An object of $\calC \star^{M} \calC'$ is either an object of $\calC$ or an object of $\calC'$. For morphisms, we take
$$ \Hom_{\calC \star^{M} \calC'}(X,Y) = \begin{cases} \Hom_{\calC}(X,Y) & \text{if } X,Y \in \calC \\
\Hom_{\calC'}(X,Y) & \text{if } X,Y \in \calC' \\
M(X,Y) & \text{if } X \in \calC, Y \in \calC' \\
\emptyset & \text{if } X \in \calC', Y \in \calC. \end{cases}$$
Composition of morphisms is defined in the obvious way, using the composition laws in $\calC$ and $\calC'$, and the functoriality of $M(X,Y)$ in $X$ and $Y$.

\begin{remark}
In the special case where $F: \calC^{op} \times \calC' \rightarrow \Set$ is the constant functor taking the value $\ast$, the category $\calC \star^{F} \calC'$ coincides with the ordinary join 
$\calC \star \calC'$.
\end{remark}

For any correspondence $M: \calC \rightarrow \calC'$, there is an obvious functor $F: \calC \star^{M} \calC' \rightarrow [1]$ (here $[1]$ denotes the linearly ordered set
$\{0,1\}$, regarded as a category in the obvious way), uniquely determined by the condition that $F^{-1} \{0\} = \calC$ and $F^{-1} \{1\} = \calC'$. Conversely, given any category $\calM$ equipped with a functor $F: \calM \rightarrow [1]$, we can {\em define} $\calC = F^{-1} \{0\}$, $\calC' = F^{-1} \{1\}$, and a correspondence $M: \calC \rightarrow \calC'$
by the formula $M(X,Y) = \Hom_{\calM}(X,Y)$. We may summarize the situation as follows:

\begin{fact}\label{factus}
Giving a pair of categories $\calC$, $\calC'$ and a correspondence between them is equivalent to giving a category $\calM$ equipped with a functor $\calM \rightarrow [1]$.
\end{fact}

Given this reformulation, it is clear how to generalize the notion of a correspondence to the $\infty$-categorical setting.

\begin{definition}\label{quasicorresp}\index{gen}{correspondence!between $\infty$-categories}
Let $\calC$ and $\calC'$ be $\infty$-categories. A {\it correspondence} from $\calC$ to $\calC'$ is a $\infty$-category $\calM$ equipped with a map $F: \calM \rightarrow \Delta^1$ and identifications $\calC \simeq F^{-1} \{0\}$, $\calC' \simeq F^{-1} \{1\}$. 
\end{definition}

\begin{remark}
Let $\calC$ and $\calC'$ be $\infty$-categories. Fact \ref{factus} generalizes to the $\infty$-categorical setting in the following way: there is a canonical bijection between equivalence classes of correspondences from $\calC$ to $\calC'$ and equivalence classes of functors
$\calC^{op} \times \calC' \rightarrow \SSet$, where $\SSet$ denotes the $\infty$-category of spaces.
In fact, it is possible to prove a more precise result (a Quillen equivalence between certain model categories), but we will not need this.
\end{remark}

To understand the relevance of Definition \ref{quasicorresp}, we note the following:

\begin{proposition}
Let $\calC$ be an ordinary category, and let $p: X \rightarrow \Nerve(\calC)$ be a
map of simplicial sets. Then $p$ is an inner fibration if and only if
$X$ is an $\infty$-category.
\end{proposition}

\begin{proof}
This follows from the fact that any map $\Lambda^n_i \rightarrow
\Nerve(\calC)$, $0 < i < n$, admits a {\em unique} extension to
$\Delta^n$.
\end{proof}

It follows readily from the definition that an arbitrary map of simplicial sets $p: X \rightarrow S$
is an inner fibration if and only if the fiber of $p$ over any
simplex of $S$ is an $\infty$-category. In particular, an inner fibration $p$
associates to each vertex $s$ of $S$ an $\infty$-category $X_{s}$, and to each edge $f: s \rightarrow s'$ in $S$ a correspondence
between the $\infty$-categories $X_{s}$ and $X_{s'}$. Higher dimensional simplices give
rise to what may be thought of as compatible ``chains'' of
correspondences.

Roughly speaking, we might think of an inner fibration $p: X \rightarrow
S$ as a functor from $S$ into some kind of
$\infty$-category of $\infty$-categories, where the morphisms are
given by correspondences. However, this description is not quite
accurate, since the correspondences are required to ``compose''
only in a weak sense.
To understand the issue, let us return to the setting of {\em
ordinary} categories. If $\calC$ and $\calC'$ are two categories,
then the correspondences from $\calC$ to $\calC'$ themselves
constitute a category, which we may denote by $M(\calC, \calC')$.
There is a natural ``composition'' defined on correspondences. If
we view an object $F \in M(\calC, \calC')$ as a functor
$\calC^{op} \times \calC' \rightarrow \Set$, and $G \in M(\calC',
\calC'')$, then we can define $(G \circ F)(C,C'')$ to be the coend
$$ \int_{C' \in \calC'} F(C,C') \times G(C',C'').$$\index{gen}{coend}

If we view $F$ as determining a 
category $\calC \star^{F} \calC'$ and $G$ as determining a category
$\calC' \star^{G} \calC''$, then $\calC \star^{G \circ F} \calC''$ is obtained
by forming the pushout
$$ ( \calC \star^{F} \calC') \coprod_{ \calC' } ( \calC' \star^{G} \calC'')$$
and then discarding the objects of $\calC'$.

Now, giving a category equipped with a functor to $[2]$ is equivalent to giving a triple of categories $\calC$, $\calC'$,
$\calC''$, together with correspondences $F \in M(\calC,\calC')$,
$G \in M(\calC', \calC'')$, $H \in M(\calC, \calC'')$ and a map $\alpha: G \circ F \rightarrow H$. But the map $\alpha$ need not be an isomorphism. Consequently, the above data cannot
literally be interpreted as a functor from $[2]$ into a
category (or even a higher category) in which the morphisms are
given by correspondences.

If $\calC$ and $\calC'$ are categories, then a correspondence from $\calC$ to $\calC'$ can be regarded as a kind of generalized functor from $\calC$ to $\calC'$. More specifically, for any functor $f: \calC \rightarrow \calC'$, we can define a correspondence $M_f$ by the formula
$$ M_f(X,Y) = \Hom_{\calC'}(f(X),Y).$$\index{gen}{correspondence!associated to a functor}
This construction gives a fully faithful embedding $\bHom_{\Cat}(\calC, \calC') \rightarrow M(\calC,\calC')$. Similarly, any functor $g: \calC' \rightarrow \calC$ determines a correspondence $M_{g}$ given by the formula $M_{g}(X,Y) = \Hom_{\calC}(X,g(Y))$; we observe that $M_{f} \simeq M_{g}$ if and only if the functors $f$ and $g$ are adjoint to one another.

If an inner fibration $p: X \rightarrow S$ corresponds to a ``functor'' from $S$ to a higher category of $\infty$-categories with morphisms given by correspondences, then some special class of inner fibrations should correspond to functors from $S$ into an $\infty$-category of $\infty$-categories with morphisms given by actual functors. This is indeed the case, and the appropriate notion is that of a {\em (co)Cartesian fibration} which we will study in \S \ref{cartfibsec}.

\subsection{Stability Properties of Inner Fibrations}\label{joyalpr}

Let $\calC$ be an $\infty$-category and $K$ an arbitrary simplicial set. In \S \ref{funcback}, we asserted that $\Fun(K,\calC)$ is an $\infty$-category (Proposition \ref{tyty}). In the course of the proof, we invoked certain stability properties of the class of inner anodyne maps. The goal of this section is to establish the required properties, and deduce some of their consequences.
Our main result is the following analogue of Proposition \ref{usejoyal}:

\begin{proposition}[Joyal \cite{joyalnotpub}]\label{usejoyal2}\index{gen}{inner anodyne}
The following collections all generate the same class of morphisms
of $\sSet$:
\begin{itemize}
\item[$(1)$] The collection $A_1$ of all horn inclusions $\Lambda^n_i
\subseteq \Delta^n$, $0 < i < n$.

\item[$(2)$] The collection $A_2$ of all inclusions $$(\Delta^m \times
\Lambda^2_1) \coprod_{ \bd \Delta^m \times \Lambda^2_1 } (\bd
\Delta^m \times \Delta^2) \subseteq \Delta^m \times \Delta^2.$$

\item[$(3)$] The collection $A_3$ of all inclusions $$(S' \times
\Lambda^2_1) \coprod_{S \times \Lambda^2_1 } (S \times \Delta^2)
\subseteq S' \times \Delta^2,$$ where $S \subseteq S'$.

\end{itemize}
\end{proposition}

\begin{proof}
We will employ the strategy that we used to prove Proposition \ref{usejoyal}, though the details are slightly more complicated. Working cell-by-cell, we conclude that every morphism in $A_3$ belongs to the weakly saturated class of morphisms generated by $A_2$. We next show that every morphism in $A_1$ is a retract of a morphism belonging to $A_3$. More precisely, we will show that for
$0 < i < n$, the inclusion $\Lambda^n_i \subseteq \Delta^n$ is a retract of the inclusion
$$ (\Delta^n \times \Lambda^2_1) \coprod_{ \Lambda^n_i \times \Lambda^2_1 }
( \Lambda^n_i \times \Delta^2) \subseteq \Delta^n \times \Delta^2.$$
To prove this, we embed $\Delta^n$ into $\Delta^n \times \Delta^2$ via
the map of partially ordered sets
$ s: [n] \rightarrow [n] \times [2]$ given by
$$ s(j) = \begin{cases} (j,0) & \text{if } j < i \\
(j,1) & \text{if } j = i \\
(j,2) & \text{if } j > i. \end{cases}$$
and consider the retraction $\Delta^n \times \Delta^2 \rightarrow \Delta^n$ given
by the map
$$ r: [n] \times [2] \rightarrow [n]$$
$$ r(j,k) = \begin{cases} j & \text{if } j < i, k=0 \\
j & \text{if } j > i, k = 2 \\
i & \text{otherwise.} \end{cases}$$

We now show that every morphism in $A_2$ is inner anodyne (that is, it lies in the weakly saturated class of morphisms generated by $A_1$). Choose $m \geq 0$. For each $0 \leq i \leq j < m$, we let
$\sigma_{ij}$ denote the $(m+1)$-simplex of $\Delta^m \times \Delta^2$ corresponding to the map
$$ f_{ij}: [m+1] \rightarrow [m] \times [2]$$
$$f_{ij}(k) = \begin{cases} (k,0) & \text{if } 0 \leq k \leq i \\
(k-1, 1) & \text{if } i+1 \leq k \leq j+1 \\
(k-1, 2) & \text{if } j+2 \leq k \leq m+1. \end{cases}$$
For each $0 \leq i \leq j \leq m$, we let $\tau_{ij}$ denote the $(m+2)$-simplex of $\Delta^m \times \Delta^2$ corresponding to the map
$$ g_{ij}: [m+2] \rightarrow [m] \times [2] $$
$$g_{ij}(k) = \begin{cases} (k,0) & \text{if } 0 \leq k \leq i \\
(k-1, 1) & \text{if } i+1 \leq k \leq j+1 \\
(k-2, 2) & \text{if } j+2 \leq k \leq m+2. \end{cases}$$

Let $X(0) = (\Delta^m \times
\Lambda^2_1) \coprod_{ \bd \Delta^m \times \Lambda^2_1 } (\bd
\Delta^m \times \Delta^2)$. For $0 \leq j < m$, we let
$$ X(j+1) = X(j) \cup \sigma_{0j} \cup \ldots \cup \sigma_{jj}. $$
We have a chain of inclusions
$$ X(j) \subseteq X(j) \cup \sigma_{0j} \subseteq \ldots \subset X(j) \cup \sigma_{0j} \cup
\ldots \cup \sigma_{jj} = X(j+1),$$
each of which is a pushout of a morphism in $A_1$ and therefore inner anodyne. It follows
that each inclusion $X(j) \subseteq X(j+1)$ is inner-anodyne. Set $Y(0) = X(m)$, so that the inclusion $X(0) \subseteq Y(0)$ is inner anodyne. We now set $Y(j+1) = Y(j) \cup \tau_{0j} \cup \ldots \cup \tau_{jj}$
for $0 \leq j \leq m$. As before, we have a chain of inclusions
$$ Y(j) \subseteq Y(j) \cup \tau_{0j} \subseteq \ldots \subseteq Y_{j} \cup \tau_{0j} \cup 
\ldots \cup \tau_{jj} = Y(j+1)$$
each of which is a pushout of a morphism belonging to $A_1$. It follows that
each inclusion $Y(j) \subseteq Y(j+1)$ is inner anodyne. By transitivity, we conclude
that the inclusion $X(0) \subseteq Y(m+2)$ is inner anodyne. We conclude the proof by observing that $Y(m+2) = \Delta^m \times \Delta^2$.
\end{proof}

\begin{corollary}[Joyal \cite{joyalnotpub}]\label{berek}
A simplicial set $\calC$ is an $\infty$-category if and only if the restriction map $$
\Fun(\Delta^2,\calC) \rightarrow \Fun(\Lambda^2_1, \calC)$$ is a trivial fibration.
\end{corollary}

\begin{proof}
By Proposition \ref{usejoyal2}, $\calC \rightarrow \ast$ is an inner fibration if and only if $S$ has the extension property with respect to each of the inclusions in the class $A_2$.
\end{proof}

\begin{remark}
In \S \ref{qqqc}, we asserted that the main function of the weak Kan condition on a simplicial set $\calC$ is that it allows us to compose the edges of $\calC$. We can regard Corollary \ref{berek} as an affirmation of this philosophy: the class of $\infty$-categories $\calC$ is characterized by the requirement that one can compose morphisms in $\calC$, and the composition is well-defined up to a contractible space of choices.
\end{remark}

\begin{corollary}[Joyal \cite{joyalnotpub}]\label{prodprod2}
Let $i: A \rightarrow A'$ be an inner anodyne map of simplicial sets, and let $j: B \rightarrow B'$ be a cofibration. Then
the induced map $$(A \times B') \coprod_{A \times B} (A' \times B)
\rightarrow A' \times B'$$ is inner anodyne.
\end{corollary}

\begin{proof}
This follows immediately from Proposition \ref{usejoyal2}, which characterizes the class of inner anodyne maps as the class generated by $A_3$ (which is stable under smash products with any cofibration).
\end{proof}

\begin{corollary}[Joyal \cite{joyalnotpub}]\index{gen}{inner fibration!and functor categories}
Let $p: X \rightarrow S$ be an inner fibration, and let $i: A \rightarrow B$ be any cofibration of simplicial sets. Then the induced map $q: X^{B} \rightarrow X^A \times_{ S^A } S^B$ is an inner fibration. If $i$ is inner anodyne, then $q$ is a trivial fibration. In particular, if $X$ is a $\infty$-category, then so is $X^B$ for any simplicial set $B$.
\end{corollary}

\subsection{Minimal Fibrations}\label{minin}

One of the aims of homotopy theory is to understand the classification of spaces up to 
homotopy equivalence. In the setting of simplicial sets, this problem admits an attractive formulation in terms of Quillen's theory of {\em minimal} Kan complexes.
Every Kan complex $X$ is homotopy equivalent to a minimal Kan complex, and a map $X \rightarrow Y$ of minimal Kan complexes is a homotopy equivalence if and only if it is an isomorphism. Consequently, the classification of Kan complexes up to homotopy equivalence is equivalent to the classification of {\em minimal} Kan complexes up to isomorphism. Of course, in practical terms, this is not of much use for solving the classification problem. Nevertheless, the theory of minimal Kan complexes (and, more generally, minimal Kan fibrations) is a useful tool in the homotopy theory of simplicial sets. The purpose of this section is to describe a generalization of the theory of minimal models, in which Kan fibrations are replaced by inner fibrations. An exposition of this theory can also be found in \cite{joyalnotpub}.

We begin by introducing a bit of terminology. Suppose given a commutative diagram
$$ \xymatrix{ A \ar[d]^{i} \ar[r]^{u} & X \ar[d]^{p} \\
B \ar[r]^{v} \ar@{-->}[ur] & S }$$
of simplicial sets where $p$ is an inner fibration, and suppose also that we have a pair $f,f': B \rightarrow X$ of candidates for the dotted arrow which render the diagram commutative. We will say that $f$ and $f'$ are {\it homotopic relative to $A$ over $S$}
if they are equivalent when viewed as objects in the $\infty$-category given by
the fiber of the map
$$ X^B \rightarrow X^A \times_{ S^A} S^B.$$
Equivalently, $f$ and $f'$ are homotopic relative to $A$ over $S$ if
there exists a map $F: B \times \Delta^1 \rightarrow X$ such that
$F | B \times \{0\} = f$, $F | B \times \{1\} = f'$, $p \circ F = v \circ \pi_{B}$,
$F \circ (i \times \id_{\Delta^1}) = u \circ \pi_{A}$, and $F| \{b\} \times \Delta^1$
is an equivalence in the $\infty$-category $X_{v(b)}$ for every vertex $b$ of $B$.

\begin{definition}\label{trukea}\index{gen}{minimal!inner fibration}
Let $p: X \rightarrow S$ be an inner fibration of simplicial sets. We will say that $p$ is {\it minimal} if $f=f'$ for every pair of maps $f,f': \Delta^n \rightarrow X$ which are homotopic relative to $\bd \Delta^n$ over $S$.

We will say that an $\infty$-category $\calC$ is {\it minimal} if the associated inner fibration
$\calC \rightarrow \ast$ is minimal.\index{gen}{minimal!$\infty$-category}
\end{definition}

\begin{remark}
In the case where $p$ is a Kan fibration, Definition \ref{trukea} recovers the usual notion of a minimal Kan fibration. We refer the reader to \cite{goerssjardine} for a discussion of minimal fibrations in this more classical setting.
\end{remark}

\begin{remark}
Let $p: X \rightarrow \Delta^n$ be an inner fibration. Then $X$ is an $\infty$-category.
Moreover, $p$ is a minimal inner fibration if and only if $X$ is a minimal $\infty$-category.
This follows from the observation that for any pair of maps $f,f': \Delta^m \rightarrow X$, a homotopy between $f$ and $f'$ is automatically compatible with the projection to $\Delta^n$.
\end{remark}

\begin{remark}
If $p: X \rightarrow S$ is a minimal inner fibration and $T \rightarrow S$ is an arbitrary map of simplicial sets, then the induced map $X_{T} = X \times_{S} T \rightarrow T$ is a minimal inner fibration.
Conversely, if $p: X \rightarrow S$ is an inner fibration and if $X \times_{S} \Delta^n \rightarrow \Delta^n$ is minimal for {\em every} map $\sigma: \Delta^n \rightarrow S$, then $p$ is minimal. Consequently, for many purposes the study of minimal inner fibrations reduces to the study of minimal $\infty$-categories.
\end{remark}

\begin{lemma}\label{ststst}
Let $\calC$ be a minimal $\infty$-category, and let $f: \calC \rightarrow \calC$ be a functor
which is homotopic to the identity. Then $f$ is a monomorphism of simplicial sets.
\end{lemma}

\begin{proof}
Choose a homotopy $h: \Delta^1 \times \calC \rightarrow \calC$ from $\id_{\calC}$ to $f$.
We prove by induction on $n$ that the map $f$ induces an injection from the set of
$n$-simplices of $\calC$ to itself.
Let $\sigma, \sigma': \Delta^n \rightarrow \calC$ be such that
$f \circ \sigma = f \circ \sigma'$. By the inductive hypothesis, we deduce that
$\sigma | \bd \Delta^n = \sigma' | \bd \Delta^n = \sigma_0$. Consider the diagram
$$ \xymatrix{ 
(\Delta^2 \times \bd \Delta^n) \coprod_{ \Lambda^2_2 \times \bd \Delta^n }
( \Lambda^2_2 \times \Delta^n ) \ar[rrr]^-{G_0} \ar@{^{(}->}[d] & & & \calC \\
\Delta^2 \times \Delta^n \ar@{-->}[urrr]^{G} & & & }$$
where $G_0| \Lambda^2_2 \times \Delta^n$ is given by amalgamating
$h \circ (\id_{\Delta^1} \times \sigma)$ with $h \circ (\id_{ \Delta^1} \times \sigma')$, and
$G_0 | \Delta^2 \times \bd \Delta^n$ is given by the composition
$$ \Delta^2 \times \bd \Delta^n \rightarrow \Delta^1 \times \bd \Delta^n
\stackrel{\sigma_0}{\rightarrow} \Delta^1 \times \calC \stackrel{h}{\rightarrow} \calC.$$
Since $h| \Delta^1 \times \{X\}$ is an equivalence for every object
$X \in \calC$, Proposition \ref{goouse} implies the existence of the map $G$ indicated in the diagram. The restriction $G|\Delta^1 \times \Delta^n$ is a homotopy between
$\sigma$ and $\sigma'$ relative to $\bd \Delta^n$. Since $\calC$ is minimal, we deduce that $\sigma = \sigma'$.
\end{proof}

\begin{lemma}\label{stsst}
Let $\calC$ be a minimal $\infty$-category, and let $f: \calC \rightarrow \calC$ be a functor
which is homotopic to the identity. Then $f$ is an isomorphism of simplicial sets.
\end{lemma}

\begin{proof}
Choose a homotopy $h: \Delta^1 \times \calC \rightarrow \calC$ from $\id_{\calC}$ to $f$.
We prove by induction on $n$ that the map $f$ induces a {\em bijection} from the set of
$n$-simplices of $\calC$ to itself. The injectivity follows from Lemma \ref{ststst}, so it will suffice to prove the surjectivity. Choose an $n$-simplex $\sigma: \Delta^n \rightarrow \calC$.
By the inductive hypothesis, we may suppose that
$\sigma| \bd \Delta^n = f \circ \sigma'_0$, for some map $\sigma'_0: \bd \Delta^n \rightarrow \calC$.
Consider the diagram
$$ \xymatrix{ (\Delta^1 \times \bd \Delta^n) \coprod_{ \{1\} \times \bd \Delta^n } (\{1\} \times \Delta^n) \ar[rrr]^-{G_0} \ar@{^{(}->}[d] & & & \calC \\
\Delta^1 \times \Delta^n, \ar@{-->}[urrr]^{G} & &  & }$$
where $G_0| \Delta^1 \times \bd \Delta^n = h \circ (\id_{\Delta^1} \times \sigma'_0)$ and
$G_0 | \{1\} \times \Delta^n = \sigma$. If $n > 0$, then the existence of the map $G$ as indicated in the diagram follows from Proposition \ref{goouse}; if $n=0$ it is obvious. Now let
$\sigma' = G| \{0\} \times \Delta^n$. To complete the proof, it will suffice to show that
$f \circ \sigma' = \sigma$.

Consider now the diagram
$$ \xymatrix{ (\Lambda^2_0 \times \Delta^n) \coprod_{ \Lambda^2_0 \times \bd \Delta^n } ( \Delta^2 \times \bd \Delta^n) \ar[rrr]^-{H_0} \ar[d] & & & \calC \\
\Delta^2 \ar@{-->}[urrr]^{H} & & &  }$$
where $H_0 | \Delta^{ \{0,1\} } \times \Delta^n = h \circ (\id_{\Delta^1} \times \sigma')$, $H_0 | \Delta^{ \{1,2\}} \times \Delta^n = G$, and $H_0| (\Delta^2 \times \bd \Delta^n)$ given by the composition
$$ \Delta^2 \times \bd \Delta^n \rightarrow \Delta^1 \times \bd \Delta^n \stackrel{\sigma'_0}{\rightarrow} \Delta^1 \times \calC \stackrel{h}{\rightarrow} \calC.$$
The existence of the dotted arrow $H$ follows once again from Proposition \ref{goouse}.
The restriction $H| \Delta^{ \{1,2\} } \times \Delta^n$ is a homotopy from $f \circ \sigma'$
to $\sigma$ relative to $\bd \Delta^n$. Since $\calC$ is minimal, we conclude that
$f \circ \sigma' = \sigma$ as desired.
\end{proof}

\begin{proposition}
Let $f: \calC \rightarrow \calD$ be an equivalence of minimal $\infty$-categories. Then
$f$ is an isomorphism.
\end{proposition}

\begin{proof}
Since $f$ is a categorical equivalence, it admits a homotopy inverse $g: \calD \rightarrow \calC$.
Now apply Lemma \ref{stsst} to the compositions $f \circ g$ and $g \circ f$.
\end{proof}

The following result guarantees a good supply of minimal $\infty$-categories:

\begin{proposition}\label{minimod}
Let $p: X \rightarrow S$ be a inner fibration of simplicial sets. Then there exists
a retraction $r: X \rightarrow X$ onto a simplicial subset $X' \subseteq X$ with the following
properties:
\begin{itemize}
\item[$(1)$] The restriction $p|X': X' \rightarrow S$ is a minimal inner fibration.
\item[$(2)$] The retraction $r$ is compatible with the projection $p$, in the sense that
$p \circ r = p$.
\item[$(3)$] The map $r$ is homotopic over $S$ to $\id_{X}$ relative to $X'$.
\item[$(4)$] For every map of simplicial sets $T \rightarrow S$, the induced inclusion
$X' \times_{S} T \subseteq X \times_{S} T$ is a categorical equivalence.
\end{itemize}
\end{proposition}

\begin{proof}
For every $n \geq 0$, we define a relation on the set of $n$-simplices
of $X$: given two simplices $\sigma, \sigma': \Delta^n \rightarrow X$, we will write
$\sigma \sim \sigma'$ if $\sigma$ is homotopic to $\sigma'$ relative to $\bd \Delta^n$. We note that $\sigma \sim \sigma'$ if and only if $\sigma| \bd \Delta^n = \sigma'| \bd \Delta^n$ and
$\sigma$ is equivalent to $\sigma'$ where both are viewed as objects in the $\infty$-category given by a fiber of the map
$$ X^{\Delta^n} \rightarrow X^{ \bd \Delta^n} \times_{ S^{\bd \Delta^n} } S^{\Delta^n}.$$
Consequently, $\sim$ is an equivalence relation.

Suppose that $\sigma$ and $\sigma'$ are both degenerate, and $\sigma \sim \sigma'$.
From the equality $\sigma | \bd \Delta^n = \sigma' | \bd \Delta^n$ we deduce that $\sigma = \sigma'$. Consequently, there is at most one degenerate $n$-simplex of $X$ in each $\sim$-class. Let $Y(n) \subseteq X_n$ denote a set of representatives for the $\sim$-classes of $n$-simplices in $X$, which contains all degenerate simplices. We now define the simplicial subset
$X' \subseteq X$ recursively as follows: an $n$-simplex $\sigma: \Delta^n \rightarrow X$
belongs to $X'$ if $\sigma \in Y(n)$ and $\sigma| \bd \Delta^n$ factors through $X'$.

Let us now prove $(1)$. To show that $p|X'$ is an inner fibration, it suffices to prove that
every lifting problem of the form
$$ \xymatrix{ \Lambda^n_i \ar[r]^{s} \ar@{^{(}->}[d] & X' \ar[d] \\
\Delta^n \ar@{-->}[ur]^{\sigma} \ar[r] & S }$$
with $0 < i < n$ has a solution $f$ in $X'$. Since $p$ is an inner fibration, this lifting problem
has a solution $\sigma': \Delta^n \rightarrow X$ in the original simplicial set $X$. Let
$\sigma'_0 = d_i \sigma: \Delta^{n-1} \rightarrow X$ be the induced map.
Then $\sigma'_0 | \bd \Delta^{n-1}$ factors through $X'$. Consequently, $\sigma'_0$
is homotopic over $S$, relative to $\bd \Delta^{n-1}$ to some map
$\sigma_0: \Delta^{n-1} \rightarrow X'$. Let $g_0: \Delta^1 \times \Delta^{n-1} \rightarrow X$
be a homotopy from $\sigma'_0$ to $\sigma_0$, and let 
$g_1: \Delta^1 \times \bd \Delta^n \rightarrow X$ be the result of amalgamating
$g_0$ with the identity homotopy from $s$ to itself. Let $\sigma_1 = g_1| \{1\} \times \bd \Delta^n$.
Using Proposition \ref{goouse}, we deduce that $g_1$ extends to a homotopy from
$\sigma'$ to some other map $\sigma'': \Delta^{n} \rightarrow X$ with
$\sigma''| \bd \Delta^{n} = \sigma_1$. It follows that $\sigma''$ is homotopic
over $S$ relative to $\bd \Delta^n$ to a map $\sigma: \Delta^n \rightarrow X$ with the desired properties. This proves that $p|X'$ is an inner fibration. It is immediate from the construction that $p|X'$ is minimal.

We now verify $(2)$ and $(3)$ by constructing a map
$h: X \times \Delta^1 \rightarrow X$ such that $h|X \times \{0\}$ is the identity,
$h| X \times \{1\}$ is a retraction $r: X \rightarrow X$ with image $X'$, and $h$ is a homotopy
over $S$ and relative to $X'$. Choose an exhaustion
of $X$ by a transfinite sequence of simplicial subsets 
$$X' = X^0 \subseteq X^1 \subseteq \ldots $$
where each $X^{\alpha}$ is obtained from 
$$X^{< \alpha} = \bigcup_{\beta < \alpha} X^{\beta}$$ by adjoining a single nondegenerate
simplex, if such a simplex exists. We construct $h_{\alpha} = h| X^{\alpha} \times \Delta^1$
by induction on $\alpha$. By the inductive hypothesis, we may suppose that we have already
defined $h_{< \alpha} = h| X^{< \alpha} \times \Delta^1$. If $X = X^{< \alpha}$, then we are done.
Otherwise, we can write $X^{\alpha} = X^{< \alpha} \coprod_{ \bd \Delta^n } \Delta^n$
corresponding to some nondegenerate simplex $\tau: \Delta^n \rightarrow X$, and it
suffices to define $h_{\alpha} | \Delta^n \times \Delta^1$. If $\tau$ factors through $X'$, we define
$h_{\alpha} | \Delta^n \times \Delta^1$ to be the composition
$$ \Delta^n \times \Delta^1 \rightarrow \Delta^n \stackrel{\sigma}{\rightarrow} X.$$
Otherwise, we use Proposition \ref{goouse} to deduce the existence of the dotted arrow $h'$
in the diagram
$$ \xymatrix{ (\Delta^n \times \{0\}) \coprod_{ \bd \Delta^n \times \{0\} } (\bd \Delta^n \times \Delta^1)
\ar[rrrr]^{ (\tau, h_{< \alpha}) } \ar@{^{(}->}[d] & & & & X \ar[d]^{p} \\
\Delta^n \times \Delta^1 \ar[rrrr]^{p \circ \sigma} \ar@{-->}[urrrr]^{h_0} & & & & S. }$$
Let $\tau' = h'| \Delta^n \times \{1\}$. Then $\tau' | \bd \Delta^n$ factors through
$X'$. It follows that there is a homotopy $h'': \Delta^{n} \times \Delta^{ \{1,2\} } \rightarrow
X$ from $\tau'$ to $\tau''$, which is over $S$ and relative to $\bd \Delta^n$, and such that
$\tau''$ factors through $X'$. Now consider the diagram
$$ \xymatrix{ ( \Delta^n \times \Lambda^2_1 ) \coprod_{ \bd \Delta^n \times \Lambda^2_1 }
(\bd \Delta^n \times \Delta^2 ) \ar[rrrr]^{H_0} \ar@{^{(}->}[d] & & & & X \ar[d]^{p} \\
\Delta^n \times \Delta^2 \ar[rrrr] \ar@{-->}[urrrr]^{H} & & & &  S}$$
where $H_0|\Delta^n \times \Delta^{ \{0,1\} } = h'$, $H_0 | \Delta^n \times \Delta^{ \{1,2\} } =
h''$, and $H_0 | \bd \Delta^n \times \Delta^2$ is given by the composition
$$ \bd \Delta^n \times \Delta^2 \rightarrow \bd \Delta^n \times \Delta^1 \stackrel{h_{<\alpha}}{\rightarrow} X.$$
Using the fact that $p$ is an inner fibration, we deduce that there exists a dotted
arrow $H$ rendering the diagram commutative. We may now define
$h_{\alpha}|\Delta^n \times \Delta^1 = H|\Delta^n \times \Delta^{ \{0,2\} }$; it is easy to see that
this extension has all the desired properties.

We now prove $(4)$. Using Proposition \ref{tulky}, we can reduce to the case where $T = \Delta^n$. 
Without loss of generality, we can replace $S$ by $T = \Delta^n$, so that $X$ and $X'$ are $\infty$-categories. The above constructions show that $r: X \rightarrow X'$ is a homotopy inverse of the inclusion $i: X' \rightarrow X$, so that $i$ is an equivalence as desired.
\end{proof}

We conclude by recording a property of minimal $\infty$-categories which makes them very useful for certain applications. 

\begin{proposition}\label{minstrict}
Let $\calC$ be a minimal $\infty$-category, and let $\sigma: \Delta^n \rightarrow \calC$
be an $n$-simplex of $\calC$ such that $\sigma | \Delta^{ \{i, i+1\} } = \id_{C}: C \rightarrow C$
is a degenerate edge. Then $\sigma = s_i \sigma_0$ for some $\sigma_0 : \Delta^{n-1} \rightarrow \calC$.
\end{proposition}

\begin{proof}
We work by induction on $n$. Let $\sigma_0 = d_{i+1} \sigma$ and let $\sigma' = s_i \sigma_0$. 
We will prove that $\sigma = \sigma'$. Our first goal is to prove that $\sigma | \bd \Delta^n = \sigma' | \bd \Delta^n$; in other words, that $d_j \sigma = d_j \sigma'$ for $0 \leq j \leq n$.
If $j= i+1$ this is obvious; if
$j \notin \{ i, i+1\}$ then it follows from the inductive hypothesis. Let us consider the case
$i = j$, and set $\sigma_1 = d^i \sigma$. We need to prove that $\sigma_0 = \sigma_1$. The argument above establishes that $\sigma_0 | \bd \Delta^{n-1} = \sigma_1 | \bd \Delta^{n-1}$.
Since $\calC$ is minimal, it will suffice to show that $\sigma_0$ and $\sigma_1$ are homotopic relative to $\bd \Delta^{n-1}$. We now observe that 
$$( s_{n-1} \sigma_0, s_{n-2} \sigma_0, \ldots, s_{i+1} \sigma_0, \sigma, 
s_{i-1} \sigma_1, \ldots, s_0 \sigma_1)$$
provides the desired homotopy $\Delta^{n-1} \times \Delta^1 \rightarrow \calC$.

Since $\sigma$ and $\sigma'$ coincide on $\bd \Delta^n$, to prove that $\sigma = \sigma'$ it will suffice to prove that $\sigma$ and $\sigma'$ are homotopic relative to $\bd \Delta^{n}$.
We now observe that
$$ ( s_n \sigma', \ldots, s_{i+2} \sigma', s_{i} \sigma', s_{i} \sigma, s_{i-1} \sigma, \ldots, s_0 \sigma)$$
is a homotopy $\Delta^{n} \times \Delta^1 \rightarrow \calC$ with the desired properties.
\end{proof}

We can interpret Proposition \ref{minstrict} as asserting that in a minimal $\infty$-category
$\calC$, composition is ``strictly unital''. For example, in the special case where $n=2$ and $i=1$, Proposition \ref{minstrict} asserts that if $f: X \rightarrow Y$ is a morphism in an $\infty$-category $\calC$, then $f$ is the {\em unique} composition $\id_{Y} \circ f$.

\subsection{$n$-Categories}\label{ncats}

The theory of $\infty$-categories can be regarded as a generalization of classical category theory: 
if $\calC$ is an ordinary category, then its nerve $\Nerve(\calC)$ is an $\infty$-category which determines $\calC$ up to canonical isomorphism. Moreover, Proposition \ref{ruko} provides a precise characterization of those $\infty$-categories which can be obtained from ordinary categories. In this section, we will explain how to specialize the theory of $\infty$-categories to obtain a theory of $n$-categories, for every nonnegative integer $n$. (However, the ideas described here are appropriate for describing only those $n$-categories which have only invertible $k$-morphisms, for every $k \geq 2$.) 

Before we can give the appropriate definition, we need to introduce a bit of terminology.
Let $f,f': K \rightarrow \calC$ be two diagrams in an $\infty$-category $\calC$, and suppose
that $K' \subseteq K$ is a simplicial subset such that $f | K' = f'| K' = f_0$. We will say that
$f$ and $f'$ are {\it homotopic relative to $K'$} if they are equivalent when viewed as objects of the $\infty$-category $\Fun(K,\calC) \times_{ \Fun(K',\calC)} \{ f_0\}$. Equivalently, $f$ and $f'$
are homotopic relative to $K'$ if there exists a homotopy
$$ h: K \times \Delta^1 \rightarrow \calC$$\index{gen}{homotopy!relative to $K' \subseteq K$}
with the following properties:
\begin{itemize}
\item[$(i)$] The restriction $h | K' \times \Delta^1$ coincides with the composition
$$ K' \times \Delta^1 \rightarrow K' \stackrel{ f_0 }{\rightarrow} \calC. $$
\item[$(ii)$] The restriction $h | K \times \{0\}$ coincides with $f$.
\item[$(iii)$] The restriction $h | K \times \{1\}$ coincides with $f'$.
\item[$(iv)$] For every vertex $x$ of $K$, the restriction
$h | \{x\} \times \Delta^1$ is an equivalence in $\calC$.
\end{itemize}

We observe that if $K'$ contains every vertex of $K$, then condition $(iv)$ follows from condition $(i)$.

\begin{definition}\label{ncat}\index{gen}{$n$-category}
Let $\calC$ be a simplicial set and $n \geq -1$ an integer. We will say that
$\calC$ is an {\it $n$-category} if it is an $\infty$-category and the following additional conditions are satisfied:
\begin{itemize}
\item[$(1)$] Given a pair of maps $f, f': \Delta^n \rightarrow \calC$, if
$f$ and $f'$ are homotopic relative to $\bd \Delta^n$, then $f = f'$.

\item[$(2)$] Given $m > n$ and a pair of maps $f,f': \Delta^m \rightarrow \calC$, if
$f | \bd \Delta^m = f' | \bd \Delta^m$, then $f = f'$.

\end{itemize}
\end{definition}

It is sometimes convenient to extend Definition \ref{ncat} to the case where $n = -2$: we will say that
a simplicial set $\calC$ is a {\it $(-2)$-category} if it is a final object of $\sSet$: in other words, if it is isomorphic to $\Delta^0$.

\begin{example}\label{minuscat}
Let $\calC$ be a $(-1)$-category. Using condition $(2)$ of Definition \ref{ncat}, one shows by induction on $m$ that $\calC$ has at most one $m$-simplex. Consequently, we see that up to isomorphism there are precisely two $(-1)$-categories: $\Delta^{-1} \simeq \emptyset$ and $\Delta^0$.
\end{example}

\begin{example}\label{0catdef}
Let $\calC$ be a $0$-category, and let $X= \calC_0$ denote the set of objects of $\calC$.
Let us write $x \leq y$ if there is a morphism $\phi$ from $x$ to $y$ in $\calC$. Since $\calC$ is an $\infty$-category, this relation is reflexive and transitive. Moreover, condition $(2)$ of Definition \ref{ncat} guarantees that the morphism $\phi$ is unique if it exists. If $x \leq y$ and $y \leq x$, it follows that the morphisms relating $x$ and $y$ are mutually inverse equivalences. Condition $(1)$ then implies that $x = y$. 
We deduce that $(X, \leq)$ is a partially ordered set. It follows from Proposition \ref{huka} below that the map $\calC \rightarrow \Nerve(X)$ is an isomorphism. 

Conversely, it is easy to see that the nerve of any partially ordered set $(X, \leq)$ is a $0$-category in the sense of Definition \ref{ncat}. Consequently, the full subcategory of $\sSet$ spanned by the $0$-categories is equivalent to the category of partially ordered sets.
\end{example}

\begin{remark}\label{slurpper}
Let $\calC$ be an $n$-category, and let $m > n+1$. Then the restriction map 
$$\theta: \Hom_{\sSet}( \Delta^m, \calC) \rightarrow \Hom_{\sSet}( \bd \Delta^m, \calC)$$
is bijective.
If $n=-1$, this is clear from Example \ref{minuscat}; let us therefore suppose that $n \geq 0$, so that $m \geq 2$. The injectivity of $\theta$ follows immediately from part $(2)$ of Definition \ref{ncat}. 
To prove the surjectivity, we consider an arbitrary map $f_0: \bd \Delta^m \rightarrow \calC$.
Let $f: \Delta^m \rightarrow \calC$ be an extension of $f_0 | \Lambda^m_1$ (which exists
since $\calC$ is an $\infty$-category, and $0 < 1 < m$). Using condition $(2)$ again, we
deduce that $\theta(f) = f_0$.
\end{remark}

The following result shows that, in the case where $n=1$, Definition \ref{ncat} recovers the usual definition of a category:

\begin{proposition}\label{huka}
Let $S$ be a simplicial set. The following conditions are equivalent:
\begin{itemize}
\item[$(1)$] The unit map $u: S \rightarrow \Nerve(\h{S})$ is an isomorphism of simplicial sets.
\item[$(2)$] There exists a small category $\calC$ and an isomorphism $S \simeq \Nerve(\calC)$ of simplicial sets.
\item[$(3)$] The simplicial set $S$ is a $1$-category.
\end{itemize}
\end{proposition}

\begin{proof}
The implications $(1) \Rightarrow (2) \Rightarrow (3)$ are clear. Let us therefore assume that $(3)$ holds, and show that $f: S \rightarrow \Nerve(\h{S})$ is an isomorphism. We will prove, by induction on $n$, that the map $u$ is bijective on $n$-simplices. 

For $n=0$, this is clear. If $n=1$, the surjectivity of $u$ obvious. To prove the injectivity, we note that if $f(\phi) = f(\psi)$, then the edges $\phi$ and $\psi$ are homotopic in $S$. A simple application of condition $(2)$ of Definition \ref{ncat} then shows that $\phi = \psi$.

Now suppose $n > 1$. The injectivity of $u$ on $n$-simplices follows from condition $(3)$ of Definition \ref{ncat}, and the injectivity of $u$ on $(n-1)$-simplices. To prove the surjectivity, let us suppose given a map $s: \Delta^n \rightarrow \Nerve(\h{S})$. Choose $0 < i < n$. Since
$u$ is bijective on lower-dimensional simplices, the map $s| \Lambda^n_i$ factors uniquely
through $S$. Since $S$ is an $\infty$-category, this factorization extends to a map
$\widetilde{s}: \Delta^n \rightarrow S$. Since $\Nerve(\h{S})$ is the nerve of a category, a pair of maps from $\Delta^n$ into $\Nerve(\h{S})$ which agree on $\Lambda^n_i$ must be the same. We deduce that $u \circ \widetilde{s} = s$, and the proof is complete.
\end{proof}

\begin{remark}
The condition that an $\infty$-category $\calC$ be an $n$-category is not invariant under categorical equivalence. For example, if $\calD$ is a category with several objects, all of which are uniquely isomorphic to one another, then $\Nerve(\calD)$ is categorically equivalent to $\Delta^0$, but is not a $(-1)$-category. Consequently, there can be no intrinsic characterization of the class of $n$-categories itself. Nevertheless, there does exist a convenient description for the class of $\infty$-categories which are {\em equivalent} to $n$-categories; see Proposition \ref{tokerp}.
\end{remark}

Our next goal is to establish that the class of $n$-categories is stable under the formation of functor categories.
In order to do so, we need to reformulate Definition \ref{ncat} in a more invariant manner.
Recall that for any simplicial set $X$, the {\it $n$-skeleton} $\sk^n X$ is defined to be the simplicial subset of $X$ generated by all the simplices of $X$ having dimension $\leq n$.\index{not}{sknX@$\sk^n X$}\index{gen}{skeleton}

\begin{proposition}\label{ncatchar}
Let $\calC$ be an $\infty$-category and $n \geq -1$. The following are equivalent:
\begin{itemize}
\item[$(1)$] The $\infty$-category $\calC$ is an $n$-category.
\item[$(2)$] For every simplicial set $K$ and
every pair of maps $f,f': K \rightarrow \calC$ such that
$f | \sk^{n} K$ and $f'|\sk^{n} K$ are homotopic relative to
$\sk^{n-1} K$, we have $f = f'$.
\end{itemize}
\end{proposition}

\begin{proof}
The implication $(2) \Rightarrow (1)$ is obvious. Suppose that $(1)$ is satisfied and let
$f,f': K \rightarrow \calC$ be as in the statement of $(2)$. To prove that
$f = f'$ it suffices to show that $f$ and $f'$ agree on
every nondegenerate simplex of $K$. We may therefore reduce to the case where
$K = \Delta^m$. We now work by induction on $m$. If $m < n$, there is nothing to prove.
In the case $m=n$, the assumption that $\calC$ is an $n$-category immediately implies that $f=f'$. If $m > n$, the inductive hypothesis implies that $f| \bd \Delta^m = f'| \bd \Delta^m$, so that
$(1)$ implies that $f = f'$.
\end{proof}

\begin{corollary}\label{zooka}\index{gen}{$n$-categories!and functor categories}
Let $\calC$ be an $n$-category and $X$ a simplicial set. Then $\Fun(X,\calC)$ is an $n$-category.
\end{corollary}

\begin{proof}
Proposition \ref{tyty} asserts that $\Fun(X,\calC)$ is an $\infty$-category. 
We will show that $\Fun(X,\calC)$ satisfies condition $(2)$ of Proposition \ref{ncatchar}.
Suppose given a pair of maps $f,f': K \rightarrow \Fun(X,\calC)$ such that
$f| \sk^{n} K$ and $f'| \sk^{n} K$ are homotopic relative to $f| \sk^{n-1} K$. We wish to show
that $f = f'$. We may identify $f$ and $f'$ with maps $F,F': K \times X \rightarrow \calC$.
Since $\calC$ is an $n$-category, to prove that $F = F'$ it suffices to show that
$F | \sk^{n}(K \times X)$ and $F' | \sk^{n}(K \times X)$ are homotopic relative to
$\sk^{n-1}(K \times X)$. This follows at once, since
$\sk^{p}(K \times X) \subseteq (\sk^{p} K) \times X$ for every integer $p$.
\end{proof}

When $n=1$, Proposition \ref{ruko} asserts that the class of $n$-categories can be characterized by the uniqueness of certain horn fillers. We now prove a generalization of this result.

\begin{proposition}\label{ruko2}
Let $n \geq 1$, and let $\calC$ be an $\infty$-category. Then $\calC$ is an $n$-category if and only if it satisfies the following condition:
\begin{itemize}
\item For every $m > n$ and every diagram
$$ \xymatrix{ \Lambda^m_i \ar[r]^{f_0} \ar@{^{(}->}[d] & \calC \\
\Delta^m, \ar@{-->}[ur]^{f} & }$$
where $0 < i < m$, there exists a {\em unique} dotted arrow $f$ as indicated, which renders the diagram commutative.
\end{itemize}
\end{proposition}

\begin{proof}
Suppose first that $\calC$ is an $n$-category. Let $f,f': \Delta^m \rightarrow \calC$
be two maps with $f | \Lambda^m_i = f' | \Lambda^m_i$, where $0 < i < m$ and
$m > n$. We wish to prove that $f = f'$. Since $\Lambda^m_i$ contains the $(n-1)$-skeleton of
$\Delta^m$, it will suffice (by Proposition \ref{ncatchar}) to show that $f$ and $f'$
are homotopic relative to $\Lambda^m_i$. This follows immediately from the fact that the inclusion $\Lambda^m_i \subseteq \Delta^m$ is a categorical equivalence.

Now suppose that every map $f_0: \Lambda^m_i \rightarrow \calC$, where
$0 < i < m$ and $n < m$, extends uniquely to an $m$-simplex of $\calC$.
We will show that $\calC$ satisfies conditions $(1)$ and $(2)$ of Definition \ref{ncat}.
Condition $(2)$ is obvious: if $f,f': \Delta^m \rightarrow \calC$ are two maps which
coincide on $\bd \Delta^m$, then they coincide on $\Lambda^m_1$ and are
therefore equal to one another (here we use the fact that $m > 1$ because of our assumption that $n \geq 1$).
Condition $(1)$ is a bit more subtle. Suppose that $f,f': \Delta^n \rightarrow \calC$ are
homotopic via a homotopy $h: \Delta^n \times \Delta^1 \rightarrow \calC$
which is constant on $\bd \Delta^n \times \Delta^1$. For $0 \leq i \leq n$, let
$\sigma_i$ denote the $(n+1)$-simplex of $\calC$ obtained by composing
$h$ with the map
$$ [n+1] \rightarrow [n] \times [1]$$
$$ j \mapsto \begin{cases} (j,0) & \text{if } j \leq i \\
(j-1,1) & \text{if } j > i. \end{cases}$$
If $i < n$, then we observe that $\sigma_i | \Lambda^{n+1}_{i+1}$ is equivalent
to the restriction $( s_i d_i \sigma_i ) | \Lambda^{n+1}_{i+1}$. Applying our hypothesis, we conclude that $\sigma_i = s_i d_i \sigma_i$, so that $d_i \sigma_i = d_{i+1} \sigma_i$. A dual argument
establishes the same equality for $0 < i$. Since $n > 0$, we conclude that
$d_i \sigma_i = d_{i+1} \sigma_i$ for all $i$. Consequently, we have a chain of equalities
$$ f' = d_0 \sigma_0 = d_1 \sigma_0 = d_1 \sigma_1 = d_2 \sigma_1 = \ldots
= d_{n} \sigma_n = d_{n+1} \sigma_n = f$$
so that $f' = f$, as desired.
\end{proof}

\begin{corollary}\label{ncatsliceee}\index{gen}{$n$-categories!and overcategories}
Let $\calC$ be an $n$-category and let $p: K \rightarrow \calC$ be a diagram.
Then $\calC_{/p}$ is an $n$-category.
\end{corollary}

\begin{proof}
If $n \leq 0$, this follows easily from Examples \ref{minuscat} and \ref{0catdef}. We may therefore suppose that $n \geq 1$. Proposition \ref{gorban3} implies that $\calC_{/p}$ is an $\infty$-category.
According to Proposition \ref{ruko2}, it suffices to show that for every
$m > n$, $0 < i < m$, and every map $f_0: \Lambda^m_i \rightarrow \calC_{/p}$, there
exists a {\em unique} map $f: \Delta^m \rightarrow \calC_{/p}$ extending $f$.
Equivalently, we must show that there is a unique map $g$ rendering the diagram
$$ \xymatrix{ \Lambda^m_i \star K \ar@{^{(}->}[d] \ar[r]^{g_0} & \calC \\
\Delta^m \star K \ar@{-->}[ur]^{g} & }$$
commutative. The existence of $g$ follows from the fact that $\calC_{/p}$ is an $\infty$-category.
Suppose that $g': \Delta^m \star K \rightarrow \calC$ is another map which extends $g_0$.
Proposition \ref{ruko} implies that $g' | \Delta^m = g | \Delta^m$. We conclude that
$g$ and $g'$ coincide on the $n$-skeleton of $\Delta^m \star K$. Since $\calC$ is an $n$-category, we deduce that $g = g'$ as desired.
\end{proof}

We conclude this section by introducing a construction which allows us to pass from an arbitrary $\infty$-category $\calC$ to its ``underlying'' $n$-category, by discarding information about morphisms of order $> n$. In the case where $n=1$, we have already introduced the relevant construction: we simply replace $\calC$ by its homotopy category (or, more precisely, the nerve of its homotopy category).\index{gen}{$n$-category!underlying an $\infty$-category}

\begin{notation}
Let $\calC$ be an $\infty$-category and let $n \geq 1$. 
For every simplicial set $K$, let $[ K, \calC ]_{n} \subseteq \Fun( \sk^{n} K, \calC )$
be the subset consisting of those diagrams $\sk^n K \rightarrow \calC$ which extend
to the $(n+1)$-skeleton of $K$ (in other words, the image of the restriction map
$\Fun( \sk^{n+1} K, \calC) \rightarrow \Fun( \sk^n K, \calC)$). 
We define an equivalence relation $\sim$ on $[K, \calC]_{n}$ as follows: given two 
maps $f,g: \sk^n K \rightarrow \calC$, we write $f \sim g$ if $f$ and $g$ are homotopic
relative to $\sk^{n-1} K$.
\end{notation}

\begin{proposition}\label{undern}\index{not}{hncalC@$\hn{n}{\calC}$}
Let $\calC$ be an $\infty$-category and $n \geq 1$.
\begin{itemize}
\item[$(1)$] There exists a simplicial set $\hn{n}{\calC}$ with the following universal mapping property: $\Fun(K, \hn{n}{ \calC}) = [K, \calC]_{n} / \sim$.
\item[$(2)$] The simplicial set $\hn{n}{\calC}$ is an $n$-category.
\item[$(3)$] If $\calC$ is an $n$-category, then the natural map $\theta: \calC \rightarrow \hn{n}{\calC}$ is an isomorphism.
\item[$(4)$] For every $n$-category $\calD$, composition with $\theta$ induces an isomorphism of simplicial sets
$$ \psi: \Fun( \hn{n}{\calC}, \calD) \rightarrow \Fun( \calC, \calD). $$
\end{itemize}
\end{proposition}

\begin{proof}
To prove $(1)$, we begin by defining $\hn{n}{\calC}([m]) = [ \Delta^m, \calC ]_{n}/ \sim$, so that the desired universal property holds by definition whenever $K$ is a simplex. Unwinding the definitions, to check the universal property for a general simplicial set $K$ we must verify the following fact:
\begin{itemize}
\item[$(\ast)$] Given two maps $f,g: \bd \Delta^{n+1} \rightarrow \calC$ which are homotopic relative to
$\sk^{n-1} \Delta^{n+1}$, if $f$ extends to an $(n+1)$-simplex of $\calC$, then $g$ extends
to an $(n+1)$-simplex of $\calC$. 
\end{itemize} 
This follows easily from Proposition \ref{princex}.

We next show that $\hn{n}{\calC}$ is an $\infty$-category. Let
$\eta_0: \Lambda^m_i \rightarrow \hn{n}{\calC}$ be a morphism, where $0 < i < m$. We wish to show that $\eta_0$ extends to an $m$-simplex $\eta: \Delta^m \rightarrow \calC$. If $m \leq n+2$, then
$\Lambda^m_i = \sk^{n+1} \Lambda^m_i$, so that $\eta_0$ can be written as a composition
$$ \Lambda^m_i \rightarrow \calC \stackrel{\theta}{\rightarrow} \hn{n}{\calC}.$$
The existence of $\eta$ now follows from our assumption that $\calC$ is an $\infty$-category.
If $m > n+2$, then $\Hom_{\sSet}( \Lambda^m_i, \hn{n}{ \calC}) \simeq \Hom_{\sSet}( \Delta^m, \hn{n}{ \calC})$ by construction, so there is nothing to prove.

We next prove that $\hn{n}{ \calC}$ is an $n$-category. It is clear from the construction that for $m > n$, any two $m$-simplices of $\hn{n}{\calC}$ with the same boundary must coincide. Suppose next that we are given two maps $f,f': \Delta^n \rightarrow \hn{n}{ \calC}$ which are homotopic relative to
$\bd \Delta^n$. Let $F: \Delta^n \times \Delta^1 \rightarrow \hn{n}{\calC}$ be a homotopy from
$f$ to $f'$. Using $(\ast)$, we deduce that $F$ is the image under $\theta$ of a map $\widetilde{F}: \Delta^n \times \Delta^1 \rightarrow \hn{n}{\calC}$, where $\widetilde{F}| \bd \Delta^n \times \Delta^1$ factors through the projection $\bd \Delta^n \times \Delta^1 \rightarrow \bd \Delta^n$. Since $n > 0$, we conclude that $\widetilde{F}$ is a homotopy from $\widetilde{F} | \Delta^n \times \{0\}$ to $\widetilde{F} | \Delta^n \times \{1\}$, so that $f = f'$. This completes the proof of $(2)$.

To prove $(3)$, let us suppose that $\calC$ is an $n$-category; we prove by induction on $m$ that the map $\calC \rightarrow \hn{n}{\calC}$ is bijective on $m$-simplices. For $m < n$, this is clear. When $m = n$ it follows from part $(1)$ of Definition \ref{ncat}. When $m = n+1$, surjectivity follows from the construction of $h_n \calC$, and injectivity from part $(2)$ of Definition \ref{ncat}.
For $m > n+1$, we have a commutative diagram
$$ \xymatrix{ \Hom_{\sSet}( \Delta^m, \calC ) \ar[r] \ar[d] & \Hom_{\sSet}( \Delta^m, \hn{n}{ \calC}) \ar[d] \\
\Hom_{\sSet}( \bd \Delta^m, \calC ) \ar[r] & \Hom_{\sSet}( \bd \Delta^m, \hn{n}{\calC}) }$$
where the bottom horizontal map is an isomorphism by the inductive hypothesis, the
left vertical map is an isomorphism by construction, and the right vertical map is an isomorphism by Remark \ref{slurpper}; it follows that the upper horizontal map is an isomorphism as well.

To prove $(4)$, we observe that if $\calD$ is an $n$-category, then the composition
$$ \Fun( \calC, \calD) \rightarrow \Fun( \hn{n}{ \calC}, \hn{n}{ \calD})
\simeq \Fun( \hn{n}{ \calC}, \calD )$$
is an inverse to $\phi$, where the second isomorphism is given by $(3)$.
\end{proof}

\begin{remark}
The construction of Proposition \ref{undern} does not quite work if $n \leq 0$, since there
may exist equivalences in $h_n \calC$ which do not arise from equivalences in $\calC$. 
However, it is a simple matter to give an alternative construction in these cases which
satisfies conditions $(2)$, $(3)$, and $(4)$; we leave the details to the reader.
\end{remark}

\begin{remark}
In the case $n=1$, the $\infty$-category $\hn{1}{\calC}$ constructed in Proposition \ref{undern}
is isomorphic to the nerve of the homotopy category $\h{\calC}$.
\end{remark}

We now apply the theory of minimal $\infty$-categories (\S \ref{minin}) to obtain a characterization of the class of $\infty$-categories which are {\em equivalent} to $n$-categories. First, we need a definition from classical homotopy theory.

\begin{definition}\label{trunckan}\index{gen}{truncated!space}
Let $k \geq -1$ be an integer. A Kan complex $X$ is {\it $k$-truncated} if, for every $i > k$ and every point $x \in X$, we have $$ \pi_{i}(X,x) \simeq \ast.$$
By convention, we will also say that $X$ is {\it $(-2)$-truncated} if $X$ is contractible.
\end{definition}

\begin{remark}
If $X$ and $Y$ are homotopy equivalent Kan complexes, then $X$ is $k$-truncated if and only if $Y$ is $k$-truncated. In other words, we may view $k$-truncatedness as a condition on objects in the homotopy category $\calH$ of spaces.
\end{remark}

\begin{example}
A Kan complex $X$ is $(-1)$-truncated if it is either empty or contractible. It is $0$-truncated if the
natural map $X \rightarrow \pi_0 X$ is a homotopy equivalence (equivalently, $X$
is $0$-truncated if it is homotopy equivalent to a discrete space).
\end{example}

\begin{proposition}\label{tokerp}\index{gen}{$n$-category}
Let $\calC$ be an $\infty$-category and $n \geq -1$. The following conditions are equivalent:
\begin{itemize}
\item[$(1)$] There exists a minimal model
$\calC' \subseteq \calC$ such that $\calC'$ is an $n$-category.
\item[$(2)$] There exists a categorical equivalence $\calD \rightarrow \calC$, where
$\calD$ is an $n$-category.
\item[$(3)$] For every pair of objects $X,Y \in \calC$, the mapping space
$\bHom_{\calC}(X,Y) \in \calH$ is $(n-1)$-truncated.
\end{itemize}
\end{proposition}

\begin{proof}
It is clear that $(1)$ implies $(2)$. Suppose next that $(2)$ is satisfied; we will prove $(3)$. Without loss of generality, we may replace $\calC$ by $\calD$ and thereby assume that $\calC$ is an $n$-category. If $n=-1$, the desired result follows immediately from Example \ref{minuscat}.
Choose $m \geq n$ and an element $\eta \in \pi_m( \bHom_{\calC}(X,Y), f)$.
We can represent $\eta$ by a commutative diagram of simplicial sets
$$ \xymatrix{ \bd \Delta^m \ar@{^{(}->}[d] \ar[r] & \{f\} \ar[d] \\
\Delta^m \ar[r]^-{s} & \Hom^{\rght}_{\calC}(X,Y). }$$
We can identify $s$ with a map $\Delta^{m+1} \rightarrow \calC$ whose restriction to
$\bd \Delta^{m+1}$ is specified. Since $\calC$ is an $n$-category, the inequality
$m+1 > n$ shows that $s$ is uniquely determined. This proves that $\pi_m( \bHom_{\calC}(X,Y), f) \simeq \ast$, so that $(3)$ is satisfied.

To prove that $(3)$ implies $(1)$, it suffices to show that if $\calC$ is a {\em minimal} $\infty$-category which satisfies $(3)$, then $\calC$ is an $n$-category. We must show that the conditions of Definition \ref{ncat} are satisfied. The first of these conditions follows immediately from the assumption that $\calC$ is minimal. For the second, we must show that if $m > n$ and $f,f': \bd \Delta^m \rightarrow \calC$ are such that $f| \bd \Delta^m = f'| \bd \Delta^m$, then $f=f'$.
Since $\calC$ is minimal, it suffices to show that $f$ and $f'$ are homotopic relative
to $\bd \Delta^m$. We will prove that there is a map $g: \Delta^{m+1} \rightarrow \calC$
such that $d_{m+1} g = f$, $d_m g = f'$, and $d_i g = d_i s_m f = d_i s_m f'$ for $0 \leq i < m$.
Then the sequence $(s_0 f, s_1 f, \ldots, s_{m-1} f, g)$ determines a map
$\Delta^m \times \Delta^1 \rightarrow \calC$ which gives the desired homotopy between
$f$ and $f'$ (relative to $\bd \Delta^m$).

To produce the map $g$, it suffices to solve the lifting problem depicted in the diagram
$$ \xymatrix{ \bd \Delta^{m+1} \ar[r]^{ g } \ar@{^{(}->}[d] & \calC \\
\Delta^{m+1} \ar@{-->}[ur]. & }$$
Choose a fibrant simplicial category $\calD$ and an equivalence
of $\infty$-categories $\calC \rightarrow \Nerve(\calD)$.
According to Proposition \ref{princex}, it will suffice to prove that we can solve 
the associated lifting problem
$$ \xymatrix{ \sCoNerve[ \bd \Delta^{m+1}] \ar[r]^{G_0} \ar@{^{(}->}[d] & \calD \\
\sCoNerve[ \Delta^{m+1}] \ar@{-->}[ur]^{G}. }$$ 
Let $X$ denote the initial vertex of $\Delta^{m+1}$, considered as an object of
$\sCoNerve[\bd \Delta^{m+1}]$, and $Y$ the final vertex. Note that $G_0$
determines a map
$$ e_0: \bd (\Delta^1)^m \simeq \bHom_{ \sCoNerve[ \bd \Delta^{m+1}]}(X,Y) \rightarrow
\bHom_{\calD}(G_0(X), G_0(Y))$$
and that giving the desired extension $G$ is equivalent to extending $e_0$ to a map
$$ e: (\Delta^1)^m \simeq \bHom_{ \sCoNerve[ \Delta^{m+1}]}(X,Y) \rightarrow
\bHom_{ \calD}(G_0(X), G_0(Y)).$$
The obstruction to constructing $e$ lies in 
$\pi_{m-1}(\bHom_{\calD}(G_0(X), G_0(Y)), p)$ for an appropriately chosen base point
$p$. Since $(m-1) > (n-1)$, condition $(3)$ implies that this homotopy set is trivial, so that
the desired extension can be found.
\end{proof}

\begin{corollary}\index{gen}{truncated!space}
Let $X$ be a Kan complex. Then $X$ is $($categorically$)$ equivalent to an $n$-category if and only if it is $n$-truncated.
\end{corollary}

\begin{proof}
For $n=-2$ this is obvious. If $n \geq -1$, this follows from characterization $(3)$ of Proposition \ref{tokerp} and the following observation: a Kan complex $X$ is $n$-truncated if and only if, for every pair of vertices $x,y \in X_0$, the Kan complex
$$\{x\} \times_{X}  X^{\Delta^1} \times_{X} \{y\}$$
of paths from $x$ to $y$ is $(n-1)$-truncated.
\end{proof}

\begin{corollary}\label{zook}
Let $\calC$ be an $\infty$-category and $K$ a simplicial set. Suppose that, for every pair of objects
$C,D \in \calC$, the space $\bHom_{\calC}(C,D)$ is $n$-truncated. Then the
$\infty$-category $\Fun(K,\calC)$ has the same property.
\end{corollary}

\begin{proof}
This follows immediately from Proposition \ref{tokerp} and Corollary \ref{zooka}, since the functor
$$ \calC \mapsto \Fun(K,\calC)$$ preserves categorical equivalences between $\infty$-categories.
\end{proof}

\section{Cartesian Fibrations}\label{cartfibsec}

\setcounter{theorem}{0}

Let $p: X \rightarrow S$ be an inner fibration of simplicial sets. Each fiber of $p$ is an $\infty$-category, and each edge $f: s \rightarrow s'$ of $S$ determines a correspondence between
the fibers $X_{s}$ and $X_{s'}$.
In this section, we would like to study the case in which each of these correspondences is associated to a functor $f^{\ast}: X_{s'} \rightarrow X_{s}$. 
Roughly speaking, we can attempt to construct $f^{\ast}$ as follows: for each vertex $y \in X_{s'}$, we choose an edge $\widetilde{f}: x \rightarrow y$ lifting $f$, and set $f^{\ast} y = x$. However, this recipe does not uniquely determine $x$, even up to equivalence, since there might be many different choices for $\widetilde{f}$. To get a good theory, we need to make a good choice of $\widetilde{f}$. More precisely, we should require that
$\widetilde{f}$ be a {\it $p$-Cartesian} edge of $X$. In \S \ref{universalmorphisms}, we will introduce the definition of $p$-Cartesian edges and study their basic properties. In particular, we will see that a $p$-Cartesian edge $\widetilde{f}$ is determined up to equivalence by its target $y$ and its image in $S$. Consequently, if there is a sufficient supply of $p$-Cartesian edges of $X$, then we can use the above prescription to define the functor $f^{\ast}: X_{s'} \rightarrow X_{s}$.
This leads us to the notion of a {\it Cartesian fibration}, which we will study in \S \ref{funkymid}. 

In \S \ref{slib}, we will establish a few basic stability properties of the class of Cartesian fibrations (we will discuss other results of this type in \S \ref{chap4}, after we have developed the language of marked simplicial sets). In \S \ref{slik} we will show that if $p: \calC \rightarrow \calD$ is a Cartesian fibration of $\infty$-categories, then we can reduce many questions about $\calC$ to similar questions about the base $\calD$ and about the fibers of $p$. This technique has many applications, which we will discuss in \S \ref{slim} and \S \ref{slin}. Finally, in \S \ref{bifib}, we will study the theory of {\it bifibrations}, which is useful for constructing examples of Cartesian fibrations.

\subsection{Cartesian Morphisms}\label{universalmorphisms}

Let $\calC$ and $\calC'$ be ordinary categories, and let $M: \calC^{op} \times \calC' \rightarrow \Set$ be a correspondence between them. Suppose that we wish to know whether or not $M$ arises as the correspondence associated to some functor $g: \calC' \rightarrow \calC$. This is the case if and only if, for each object $C' \in \calC'$, we can find an object $C \in \calC$ and a point $\eta \in M(C,C')$
having the property that the ``composition with $\eta$'' map 
$$\psi: \Hom_{\calC}(D,C) \rightarrow M(D,C')$$ is bijective, for all $D \in \calC$. Note that
$\eta$ may be regarded as a morphism in the category
$ \calC \star^{M} \calC' $. We will say that $\eta$ is a {\it Cartesian} morphism in $\calC \star^{M} \calC'$ if $\psi$ is bijective for each $D \in \calC$. The purpose of this section is to generalize this notion to the $\infty$-categorical setting and to establish its basic properties.

\begin{definition}\label{univedge}\index{gen}{morphism!$p$-Cartesian}\index{gen}{Cartesian edge}
Let $p: X \rightarrow S$ be an inner fibration of simplicial sets.
Let $f: x \rightarrow y$ be an edge in $X$. We shall say that $f$
is {\it $p$-Cartesian} if the induced map
$$ X_{/f} \rightarrow X_{/y} \times_{ S_{/p(y)} } S_{/p(f)}$$ 
is a trivial Kan fibration.
\end{definition}

\begin{remark}
Let $\calM$ be an ordinary category, and let $p: \Nerve(\calM) \rightarrow \Delta^1$ be a map (automatically an inner fibration), and let $f: x \rightarrow y$ be a morphism in $\calM$ which projects isomorphically onto $\Delta^1$. Then $f$ is $p$-Cartesian in the sense of Definition \ref{univedge} if and only if it is Cartesian in the classical sense.
\end{remark}

We now summarize a few of the formal properties of Definition \ref{univedge}:

\begin{proposition}\label{stuch}
\begin{itemize}
\item[$(1)$] Let $p: X \rightarrow S$ be an isomorphism of simplicial sets. Then every edge of $X$ is $p$-Cartesian.

\item[$(2)$] Suppose given a pullback diagram
$$ \xymatrix{ X' \ar[d]^{p'} \ar[r]^{q} & X \ar[d]^{p} \\
S' \ar[r] & S }$$
of simplicial sets, where $p$ $($and therefore also $p'${}$)$ is an inner fibration. Let
$f$ be an edge of $X'$. If $q(f)$ is $p$-Cartesian, then $f$ is $p'$-Cartesian.

\item[$(3)$] Let $p: X \rightarrow Y$ and $q: Y \rightarrow Z$ be
inner fibrations, and let $f: x' \rightarrow x$ be an edge of $X$ such that $p(f)$ is $q$-Cartesian. Then
$f$ is $p$-Cartesian if and only if $f$ is $(q \circ p)$-Cartesian.
\end{itemize}
\end{proposition}

\begin{proof}
Assertions $(1)$ and $(2)$ follow immediately from the definition. To prove $(3)$, we consider the commutative diagram
$$ \xymatrix{ X_{/f} \ar[rr]^{\psi} \ar[dr]^{\psi'} & & X_{/x} \times_{ Z_{/ (q \circ p)(x)} } Z_{/ (q \circ p)(f)} \\
& X_{/x} \times_{ Y_{/p(x)} } Y_{/p(f)}. \ar[ur]^{\psi''} & }$$
The map $\psi''$ is a pullback of 
$$Y_{/p(f)} \rightarrow Y_{/p(x)} \times_{ Z_{/(q \circ p)(x)} } Z_{/(q \circ p)(f)}$$
and therefore a trivial fibration, in view of our assumption that $p(f)$ is $q$-Cartesian. If $\psi'$ is a trivial fibration, it follows that $\psi$ is a trivial fibration as well, which proves the ``only if'' direction of $(3)$.

For the converse, suppose that $\psi$ is a trivial fibration. Proposition \ref{sharpen} implies
that $\psi'$ is a right fibration. According to Lemma \ref{toothie}, it will suffice to prove that the fibers of $\psi'$ are contractible. Let $t$ be a vertex of $X_{/x} \times_{ Y_{/p(x)} } Y_{/p(f)}$, and let
$K = (\psi'')^{-1} \{ \psi''(t) \}$. Since $\psi''$ is a trivial fibration, $K$ is a contractible Kan complex.
Since $\psi$ is a trivial fibration, the simplicial set
$ (\psi')^{-1} K = \psi^{-1} \{ \psi''(t) \}$ is also a contractible Kan complex. It follows that the
fiber of $\psi'$ over the point $t$ is weakly contractible, as desired.
\end{proof}

\begin{remark}\label{univsay}
Let $p: X \rightarrow S$ be an inner fibration of simplicial sets. Unwinding the definition, we see that an edge $f: \Delta^1 \rightarrow X$ is $p$-Cartesian if and only if for every $n \geq 2$ and every commutative diagram
$$ \xymatrix{ \Delta^{ \{n-1, n\} } \ar[dr]^{f} \ar@{^{(}->}[d] & \\
\Lambda^n_n \ar[r] \ar@{^{(}->}[d] & X \ar[d]^{p} \\
\Delta^n \ar[r] \ar@{-->}[ur] & S, }$$
there exists a dotted arrow as indicated, rendering the diagram commutative.
\end{remark}

In particular, we note that Proposition \ref{greenlem} may be restated as follows:

\begin{itemize}
\item[$(\ast)$] Let $\calC$ be a $\infty$-category, and let $p: \calC \rightarrow \Delta^0$ denote the projection from $\calC$ to a point. A morphism $\phi$ of $\calC$ is $p$-Cartesian if and only if $\phi$ is an equivalence.
\end{itemize}

In fact, it is possible to strengthen assertion $(\ast)$ as follows:

\begin{proposition}\label{universalequiv}
Let $p: \calC \rightarrow \calD$ be an inner fibration between $\infty$-categories, and let $f: C \rightarrow C'$ be a morphism in $\calC$. The following conditions are equivalent:
\begin{itemize}
\item[$(1)$] The morphism $f$ is an equivalence in $\calC$.
\item[$(2)$] The morphism $f$ is $p$-Cartesian and $p(f)$ is an equivalence in $\calD$.
\end{itemize}
\end{proposition}

\begin{proof}
Let $q$ denote the projection from $\calD$ to a point. We note that both $(1)$ and $(2)$
imply that $p(f)$ is an equivalence in $\calD$, and therefore $q$-Cartesian by $(\ast)$.
The equivalence of $(1)$ and $(2)$ now follows from $(\ast)$ and the third part of Proposition \ref{stuch}.
\end{proof}

\begin{corollary}\label{corpal}
Let $p: \calC \rightarrow \calD$ be an inner fibration between $\infty$-categories. Every identity morphism of $\calC$ $($in other words, every degenerate edge of $\calC${}$)$ is $p$-Cartesian.
\end{corollary}

We now study the behavior of Cartesian edges under composition.

\begin{proposition}\label{protohermes}\index{gen}{Cartesian edge!and composition}
Let $p: \calC \rightarrow \calD$ be an inner fibration between simplicial sets, and let
$\sigma: \Delta^2 \rightarrow \calC$ be a $2$-simplex of $\calC$, which we will depict as
a diagram
$$ \xymatrix{ & C_1 \ar[dr]^{g} & \\
C_0 \ar[ur]^{f} \ar[rr]^{h} & & C_2. }$$
Suppose that $g$ is $p$-Cartesian. Then $f$ is $p$-Cartesian if and only if $h$ is $p$-Cartesian.
\end{proposition}

\begin{proof}
We wish to show that the map
$$ i_0: \calC_{/h} \rightarrow \calC_{/C_2} \times_{ \calD_{/p(C_2)}} \calD_{/p(h)}$$ is a trivial fibration
if and only if $$ i_1: \calC_{ /f } \rightarrow \calC_{/C_1} \times_{ \calD_{/p(C_1)}} \calD_{/p(f)}$$ is a trivial fibration. The dual of Proposition \ref{sharpen} implies that both maps are right fibrations. Consequently, by (the dual of) Lemma \ref{toothie}, it suffices to show that the fibers of $i_0$ are contractible if and only if the fibers of $i_1$ are contractible. 

For any simplicial subset $B \subseteq \Delta^2$, let $X_B = \calC_{/\sigma|B} \times_{ \calD_{\sigma|B}} \calD_{/\sigma}$. We note that $X_B$ is functorial in $B$, in the sense that an inclusion 
$A \subseteq B$ induces a map $j_{A,B}: X_B \rightarrow X_A$ (which is a right fibration, again by Proposition \ref{sharpen}). We note that $j_{ \Delta^{ \{2\} }, \Delta^{ \{0,2\} }}$
is the base change of $i_0$ by the map $\calD_{/p(\sigma)} \rightarrow \calD_{/p(h)}$, and that $j_{ \Delta^{ \{1\} }, \Delta^{ \{0,1\} }}$ is the base change of $i_1$ by the map
$\calD_{/\sigma} \rightarrow \calD_{/p(f)}$.  The maps $$ \calD_{ /p(f)} \leftarrow \calD_{/p(\sigma)} \rightarrow \calD_{/p(h)}$$ are both surjective on objects (in fact, both maps have sections).
Consequently, it suffices to prove that
$j_{ \Delta^{ \{1\}}, \Delta^{ \{0,1\} }}$ has contractible fibers if and only if $j_{ \Delta^{ \{2\}}, \Delta^{ \{0,2\} }}$ has contractible fibers. Now we observe that the compositions

$$ X_{\Delta^2} \rightarrow X_{ \Delta^{ \{0,2\} }} \rightarrow X_{ \Delta^{ \{2\} }}$$
$$ X_{\Delta^2} \rightarrow X_{ \Lambda^{2}_{1} } \rightarrow X_{\Delta^{ \{1,2\} }}  \rightarrow X_{ \Delta^{ \{2\} }}$$ 
coincide. By Proposition \ref{sharpen2}, $j_{A,B}$ is a trivial fibration whenever the inclusion $A \subseteq B$ is left anodyne. we deduce that $j_{ \Delta^{ \{2\} }, \Delta^{ \{0,2\} } }$ is a trivial fibration if and only if
$j_{\Delta^{ \{1,2\} }, \Lambda^2_1}$ is a trivial fibration. Consequently, it suffices to show that $j_{ \Delta^{ \{1,2\} }, \Lambda^2_1}$ is a trivial fibration if and only if $j_{ \Delta^{ \{1\} }, \Delta^{ \{0,1\} }}$ is a trivial fibration.

Since $j_{ \Delta^{ \{1,2\} }, \Lambda^2_1}$ is a pullback of $j_{ \Delta^{ \{1\} }, \Delta^{ \{0,1\} }}$, the ``if'' direction is obvious. For the converse, it suffices to show that the natural map
$$\calC_{/g} \times_{ \calD_{ /p(g)} } \calD_{/ p(\sigma)} \rightarrow \calC_{/ C_1} \times_{ \calD_{/ p(C_1)} } \calD_{/p(\sigma)}$$ is surjective on vertices. But this map is a trivial fibration, since the inclusion $\{1\} \subseteq \Delta^{ \{1,2\} }$ is left anodyne.
\end{proof}

Our next goal is to reformulate the notion of a Cartesian morphism in a form which will be useful later.
For convenience of notation, we will prove this result in a dual form. If $p: X \rightarrow S$ is an inner fibration and $f$ an edge of $X$, we will say that $f$ is {\it $p$-coCartesian} if is Cartesian with respect to the morphism $p^{op}: X^{op} \rightarrow S^{op}$.\index{gen}{coCartesian edge}\index{gen}{morphism!$p$-coCartesian}

\begin{proposition}\label{goouse}\index{gen}{Cartesian edge}
Let $p: Y \rightarrow S$ be an inner fibration of simplicial sets, and $e: \Delta^1 \rightarrow
Y$ an edge. Then $e$ is $p$-coCartesian if and only if for each $n \geq 1$ and each diagram
$$ \xymatrix{ \{0\} \times \Delta^1 \ar[drr]^{e} \ar@{^{(}->}[d] & & \\
(\Delta^n \times \{0\}) \coprod_{ \bd \Delta^n \times \{0\} } ( \bd \Delta^n \times \Delta^1) \ar[rr]^-{f} \ar@{^{(}->}[d] & & Y \ar[d]^{p} \\
\Delta^n \times \Delta^1 \ar@{-->}[urr]^{h} \ar[rr]^{g} & & S }$$
there exists a map $h$ as indicated, rendering the diagram commutative.
\end{proposition}

\begin{proof}
Let us first prove the ``only if'' direction. We recall a bit of the
notation used in the proof of Proposition \ref{usejoyal}; in particular, the filtration
$$X(n+1) \subseteq \ldots \subseteq X(0) = \Delta^n \times
\Delta^1$$ of $\Delta^n \times \Delta^1$. We construct $h|X(m)$
by descending induction on $m$. To begin, we set $h|X(n+1) = f$.
Now, for each $m$ the space $X(m)$ is obtained from $X(m+1)$ by
pushout along a horn inclusion $\Lambda^{n+1}_m \subseteq
\Delta^{m+1}$. If $m > 0$, the desired extension exists because
$p$ is an inner fibration. If $m = 0$, the desired extension exists
because of the hypothesis that $e$ is a $p$-coCartesian edge.

We now prove the ``if'' direction. Suppose that $e$ satisfies the
condition in the statement of the Proposition. We wish to show
that $e$ is $p$-coCartesian. In other words, we must show that for every $n \geq 2$ and
every diagram
$$ \xymatrix{ \Delta^{ \{0,1\} } \ar[dr]^{e} \ar@{^{(}->}[d] & \\
\Lambda^n_0 \ar[r] \ar@{^{(}->}[d] & X \ar[d]^{p} \\
\Delta^n \ar[r] \ar@{-->}[ur] & S }$$
there exists a dotted arrow as indicated, rendering the diagram commutative. 
Replacing $S$ by $\Delta^n$ and $Y$ by
$Y \times_{S} \Delta^n$, we may reduce to the case where $S$ is a $\infty$-category. We
again make use of the notation (and argument) employed in the proof
of Proposition \ref{usejoyal}. Namely, the inclusion $\Lambda^n_0
\subseteq \Delta^n$ is a retract of the inclusion
$$ (\Lambda^n_0 \times \Delta^1) \coprod_{ \Lambda^n_0 \times
\{0\} } (\Delta^n \times \{0\}) \subseteq \Delta^n \times
\Delta^1.$$ The retraction is implemented by maps
$$ \Delta^n \stackrel{j}{\rightarrow} \Delta^n \times \Delta^1
\stackrel{r}{\rightarrow} \Delta^n$$ which were defined in the
proof of Proposition \ref{usejoyal}. We now set $F = f \circ r$,
$G = g \circ r$.

Let $K = \Delta^{ \{1, 2, \ldots, n \} } \subseteq \Delta^n$.
Then $$F| (\bd K \times \Delta^1) \coprod_{ \bd K \times \{0\}} (K \times
\Delta^1)$$ carries $\{1\} \times \Delta^1$ into $e$. By assumption, there exists an
extension of $F$ to $K \times \Delta^1$ which is compatible with
$G$. In other words, there exists a compatible extension $F'$ of
$F$ to $$ \bd \Delta^n \times \Delta^1 \coprod_{ \bd \Delta^n
\times \{0\} } \Delta^n \times \{0 \}.$$ Moreover, $F'$ carries
$\{0\} \times \Delta^1$ to a degenerate edge; such an edge is
automatically coCartesian (by Corollary \ref{corpal}, since $S$ is an $\infty$-category), and therefore there exists an extension of $F'$ to all of $\Delta^n \times \Delta^1$ by the first part of
the proof.
\end{proof}

\begin{remark}\label{kermy}
Let $p: X \rightarrow S$ be an inner fibration of simplicial sets, $x$ a vertex of $X$, and 
$\overline{f}: \overline{x}' \rightarrow p(x)$ an edge of $S$ ending at $p(x)$. There may exist many $p$-Cartesian edges $f: x' \rightarrow x$ of $X$ with $p(f) = \overline{f}$. However, there is a sense in which any two such edges having the same target $x$ are equivalent to one another.
Namely, any $p$-Cartesian edge $f: x' \rightarrow x$ lifting $\overline{f}$ can be regarded as a final object of the $\infty$-category $X_{/x} \times_{ S_{/p(x)} } \{ \overline{f} \}$, and is therefore determined up to equivalence by $\overline{f}$ and $x$.
\end{remark}

We now spell out the meaning of Definition \ref{univedge} in the setting of simplicial categories.

\begin{proposition}\label{trainedg}\index{gen}{Cartesian edge!and simplicial categories}
Let $F: \calC \rightarrow \calD$ be a functor between simplicial categories.
Suppose that $\calC$ and $\calD$ are fibrant, and that for every pair of objects
$C,C' \in \calC$, the associated map
$$ \bHom_{\calC}(C,C') \rightarrow \bHom_{\calD}(F(C),F(C'))$$ is a Kan fibration.
Then:
\begin{itemize}
\item[$(1)$] The associated map $q: \sNerve(\calC) \rightarrow \sNerve(\calD)$ is an inner fibration between $\infty$-categories.

\item[$(2)$] A morphism $f: C' \rightarrow C''$ in $\calC$ is $q$-Cartesian if and only if,
for every object $C \in \calC$, the diagram of simplicial sets
$$ \xymatrix{ \bHom_{\calC}(C, C') \ar[r] \ar[d] & \bHom_{\calC}(C,C'') \ar[d] \\
\bHom_{\calD}(F(C), F(C')) \ar[r] & \bHom_{\calD}(F(C), F(C''))}$$
is homotopy Cartesian.
\end{itemize}

\end{proposition}

\begin{proof}
Assertion $(1)$ follows from Remark \ref{goobrem}. Let $f$ be a morphism in $\calC$. By definition, $f: C' \rightarrow C''$ is $q$-Cartesian if and only if 
$$\theta: \sNerve(\calC)_{/f} \rightarrow \sNerve(\calC)_{/C''} \times_{ \sNerve(\calD)_{/F(C'')}} \sNerve(\calD)_{/F(f)}$$ is a trivial fibration. Since $\theta$ is a right fibration between
right fibrations over $\calC$, $f$ is $q$-Cartesian if and only if for every object $C \in \calC$, 
the induced map
$$\theta_{C}: \{C\} \times_{\sNerve(\calC)} \sNerve(\calC)_{/f} \rightarrow \{C\} \times_{\sNerve(\calC)} \sNerve(\calC)_{/C''} \times_{ \sNerve(\calD)_{/F(C'')}} \sNerve(\calD)_{/F(f)} $$
is a homotopy equivalence of Kan complexes. This is equivalent to the assertion that the diagram
$$ \xymatrix{ \sNerve(\calC)_{/f} \times_{\calC} \{C\} \ar[r] \ar[d] & \sNerve(\calC)_{/C''} \times_{\sNerve(\calC)} \{C\} \ar[d] \\
\sNerve(\calD)_{/F(f)} \times_{\sNerve(\calD)} \{F(C)\} \ar[r] & \sNerve(\calD)_{/F(C'')} \times_{ \sNerve(\calD)} \{ F(C) \} }$$
is homotopy Cartesian. In view of Theorem \ref{biggie}, this diagram is equivalent to the diagram of simplicial sets
$$ \xymatrix{ \bHom_{\calC}(C, C') \ar[r] \ar[d] & \bHom_{\calC}(C,C'') \ar[d] \\
\bHom_{\calD}(F(C), F(C')) \ar[r] & \bHom_{\calD}(F(C), F(C'')).}$$
This proves $(2)$.
\end{proof}

In some contexts, it will be convenient to introduce a slightly larger class of edges:

\begin{definition}\index{gen}{Cartesian!locally}\index{gen}{locally Cartesian!edge}
Let $p: X \rightarrow S$ be an inner fibration, and let $e: \Delta^1 \rightarrow X$ be an edge.
We will say that $e$ is {\it locally $p$-Cartesian} if it is a $p'$-Cartesian edge of the fiber product
$X \times_{S} \Delta^1$, where $p': X \times_{S} \Delta^1 \rightarrow \Delta^1$ denotes the projection.
\end{definition}

\begin{remark}\label{intin}
Suppose given a pullback diagram
$$ \xymatrix{ X' \ar[r]^{f} \ar[d]^{p'} & X \ar[d]^{p} \\
S' \ar[r] & S }$$
of simplicial sets, where $p$ (and therefore also $p'$) is an inner fibration. An edge
$e$ of $X'$ is locally $p'$-Cartesian if and only if its image $f(e)$ is locally $p$-Cartesian.
\end{remark}

We conclude with a somewhat technical result which will be needed in \S \ref{bicat1}:

\begin{proposition}\label{sworkk}
Let $p: X \rightarrow S$ be an inner fibration of simplicial sets. Let
$f: x \rightarrow y$ be an edge of $X$ Suppose that there is a $3$-simplex
$\sigma: \Delta^3 \rightarrow X$ such that $d_1 \sigma = s_0 f$ and $d_2 \sigma = s_1 f$.
Suppose furthermore that there exists a $p$-Cartesian edge $\widetilde{f}:
\widetilde{x} \rightarrow y$ such that $p(\widetilde{f}) = p(f)$.
Then $f$ is $p$-Cartesian.
\end{proposition}

\begin{proof}
We have a diagram of simplicial sets
$$ \xymatrix{ \Lambda^2_2 \ar[rr]^{(\widetilde{f}, f, \bigdot)} \ar@{^{(}->}[d] & & X \ar[d]^{p} \\
\Delta^2 \ar[rr]^{s_0 p(f)} \ar@{-->}[urr]^{\tau} & & S.}$$
Because $\widetilde{f}$ is $p$-Cartesian, there exists a map $\tau$ rendering the diagram commutative. Let $g = d_2(\tau)$, which we regard as a morphism $x \rightarrow \widetilde{x}$ in the $\infty$-category $X_{p(x)} = X \times_{S} \{p(x)\}$. We will show that $g$ is an equivalence in $X_{p(x)}$. It will follow that $g$ is $p$-Cartesian and that $f$, being a composition of $p$-Cartesian edges, is $p$-Cartesian (Proposition \ref{protohermes}).

Now consider the diagram
$$ \xymatrix{ \Lambda^2_1 \ar[rr]^{(d_0 d_3 \sigma, \bigdot, g)} \ar@{^{(}->}[d] & & X \ar[d]^{p} \\
\Delta^2 \ar[rr]^{ d_3 p(\sigma) } \ar@{-->}[urr]^{\tau'} & & S.}$$
The map $\tau'$ exists since $p$ is an inner fibration. Let $g' = d_1 \tau'$. We will show that
$g': \widetilde{x} \rightarrow x$ is a homotopy inverse to $g$ in the $\infty$-category
$X_{p(x)}$. 

Using $\tau$ and $\tau'$, we construct a new diagram
$$ \xymatrix{ \Lambda^3_2 \ar[rr]^{ (\tau', d_3 \sigma, \bigdot, \tau) }\ar@{^{(}->}[d] & &  X \ar[d]^{p} \\
\Delta^3 \ar[rr]^{ s_0 d_3 p(\sigma) } \ar@{-->}[urr]^{\theta} & & S.}$$
Since $p$ is an inner fibration, we deduce the existence of $\theta: \Delta^3 \rightarrow X$ rendering the
diagram commutative. The simplex $d_2(\theta)$ exhibits $\id_{x}$ as a composition
$g' \circ g$ in the $\infty$-category $X_{p(s)}$. It follows that $g'$ is a left homotopy inverse to $g$.

We now have a diagram
$$ \xymatrix{ \Lambda^2_1 \ar[rr]^{(g, \bigdot, g')} \ar@{^{(}->}[d] & & X_{p(x)} \\
\Delta^2. \ar@{-->}[urr]^{\tau''} } $$
The indicated $2$-simplex $\tau''$ exists since $X_{p(x)}$ is an $\infty$-category, and
exhibits $d_1(\tau'')$ as a composition $g \circ g'$. To complete the proof, it will suffice to show that
$d_1(\tau'')$ is an equivalence in $X_{p(x)}$. 

Consider the diagrams
$$ \xymatrix{ \Lambda^3_1 \ar[rr]^{ (d_0 \sigma, \bigdot, s_1 \widetilde{f}, \tau')} \ar@{^{(}->}[d] & & X \ar[d]^{p} &  \Lambda^3_1 \ar[rr]^{ (\tau, \bigdot, d_1 \theta', \tau'') } \ar@{^{(}->}[d] & & X \ar[d]^{p} \\
\Delta^3 \ar@{-->}[urr]^{\theta'} \ar[rr]^{\sigma} & & S & \Delta^3 \ar[rr]^{ s_0 s_0 p(f)} \ar@{-->}[urr]^{\theta''} & & S. }$$
Since $p$ is an inner fibration, there exist $3$-simplices $\theta', \theta'': \Delta^3 \rightarrow X$ with the inducated properties. The $2$-simplex $d_1(\theta'')$ identifies $d_1(\tau'')$ as a map between two
$p$-Cartesian lifts of $p(f)$; it follows that $d_1(\tau'')$ is an equivalence, which completes the proof.
\end{proof}

\subsection{Cartesian Fibrations}\label{funkymid}

In this section, we will introduce the study of {\it Cartesian fibrations} between simplicial sets. The theory of Cartesian fibrations is a generalization of the theory of right fibrations studied in \S \ref{leftfibsec}. Recall that if $f: X \rightarrow S$
is a right fibration of simplicial sets, then the fibers $\{ X_{s} \}_{s \in S}$
are Kan complexes, which depend in a (contravariantly) functorial fashion on
the choice of vertex $s \in S$. The condition that $f$ be a Cartesian fibration has a similar flavor: we still require that $X_{s}$ depend functorially on $s$, but weaken the requirement that $X_{s}$ be a Kan complex; instead, we merely require that it is an $\infty$-category.

\begin{definition}\label{defcart}\index{gen}{Cartesian fibration}\index{gen}{fibration!Cartesian}
We will say that a map $p: X \rightarrow S$ of simplicial sets is a {\it Cartesian fibration} if the following conditions are satisfied:
\begin{itemize}
\item[$(1)$] The map $p$ is an inner fibration.
\item[$(2)$] For every edge $f: x \rightarrow y$ of $S$ and every vertex $\widetilde{y}$
of $X$ with $p(\widetilde{y}) = y$, there exists a $p$-Cartesian edge $\widetilde{f}: \widetilde{x} \rightarrow \widetilde{y}$ with $p(\widetilde{f}) = f$.
\end{itemize}

We say that $p$ is a {\it coCartesian fibration} if the opposite map $p^{op}: X^{op} \rightarrow S^{op}$ is a Cartesian fibration.\index{gen}{fibration!coCartesian}\index{gen}{coCartesian fibration}
\end{definition}

If a general inner fibration $p: X \rightarrow S$ associates to each vertex $s \in S$ an $\infty$-category $X_{s}$ and to each edge $s \rightarrow s'$ a correspondence from $X_{s}$ to $X_{s'}$, then $p$ is Cartesian if each of these correspondences arise from an (canonically determined) functor $X_{s'} \rightarrow X_{s}$. In other words, a Cartesian fibration with base $S$ ought to be roughly the same thing as a contravariant functor from $S$ into an $\infty$-category of $\infty$-categories, where the morphisms are given by {\em functors}. 
One of the main goals of \S \ref{chap4} is to give a precise formulation
(and proof) of this assertion.

\begin{remark}\label{gcart}
Let $F: \calC \rightarrow \calC'$ be a functor between (ordinary) categories. The induced map of simplicial sets $\Nerve(F): \Nerve(\calC) \rightarrow \Nerve(\calC')$ of simplicial sets is automatically an inner fibration; it is Cartesian if and only if $F$ is a {\it fibration} of categories in the sense of Grothendieck.
\end{remark}

The following formal properties follow immediately from the definition:

\begin{proposition}
\begin{itemize}
\item[$(1)$] Any isomorphism of simplicial sets is a Cartesian fibration.

\item[$(2)$] The class of Cartesian fibrations between simplicial sets is stable under base change.

\item[$(3)$] A composition of Cartesian fibrations is a Cartesian fibration.
\end{itemize}
\end{proposition}

Recall that an $\infty$-category $\calC$ is a Kan complex if and only if every morphism in $\calC$ is an equivalence. We now establish a relative version of this statement:

\begin{proposition}\label{goey}\index{gen}{Cartesian fibration!and right fibrations}
Let $p: X \rightarrow S$ be an inner fibration of simplicial sets. The following
conditions are equivalent:
\begin{itemize}
\item[$(1)$] The map $p$ is a Cartesian fibration and every edge in $X$ is
$p$-Cartesian.

\item[$(2)$] The map $p$ is a right fibration.

\item[$(3)$] The map $p$ is a Cartesian fibration and every fiber of $p$ is a Kan
complex.
\end{itemize}
\end{proposition}

\begin{proof}
In view of Remark \ref{univsay}, the assertion that every edge of $X$ is $p$-Cartesian is equivalent to the assertion that $p$ has the right lifting property with respect to $\Lambda^n_n \subseteq \Delta^n$ for all $n \geq 2$. The requirement that $p$ be a Cartesian fibration further imposes the right lifting property with respect to $\Lambda^1_1 \subseteq \Delta^1$. This proves that $(1) \Leftrightarrow (2)$.

Suppose that $(2)$ holds. Since we have established that $(2)$ implies $(1)$, we know that $p$ is Cartesian. Furthermore, we have already seen that the fibers of a right fibration are Kan complexes. Thus $(2)$ implies $(3)$.

We complete the proof by showing that $(3)$ implies that every edge $f: x \rightarrow y$
of $X$ is $p$-Cartesian. Since $p$ is a Cartesian fibration, there exists a $p$-Cartesian edge $f': x' \rightarrow y$
with $p(f') = p(f)$. Since $f'$ is $p$-Cartesian, 
there exists a $2$-simplex $\sigma: \Delta^2 \rightarrow X$ which we may depict as a diagram
$$ \xymatrix{ & x' \ar[dr]^{f'} & \\
x \ar[ur]^{g} \ar[rr]^{f} & & y, }$$
where $p(\sigma) = s_0 p(f)$. Then $g$ lies in the fiber $X_{p(x)}$, and is therefore an equivalence (since $X_{p(x)}$ is a Kan complex). It follows that $f$ is equivalent to $f'$ as objects of
$X_{/y} \times_{ S_{/p(y)} } \{p(f) \}$, so that $f$ is $p$-Cartesian as desired.
\end{proof}

\begin{corollary}\label{relativeKan}
Let $p: X \rightarrow S$ be a Cartesian fibration. Let $X'
\subseteq X$ consist of all those simplices $\sigma$ of $X$ such that
every edge of $\sigma$ is $p$-Cartesian. Then $p|X'$ is a right fibration.
\end{corollary}

\begin{proof}
We first show that $p|X'$ is an inner fibration. It suffices to show
that $p|X'$ has the right lifting property with respect to every
horn inclusion $\Lambda^n_i$, $0 < i < n$. If $n > 2$, then
this follows immediately from the fact the fact that $p$ has the
appropriate lifting property. If $n = 2$, then we must show that
if $f: \Delta^2 \rightarrow X$ is such that $f|\Lambda^2_1$
factors through $X'$, then $f$ factors through $X'$. This follows immediately from Proposition \ref{protohermes}.

We now wish to complete the proof by showing that $p$ is a right fibration. According to
Proposition \ref{goey}, it suffices to prove that every edge of $X'$ is $p|X'$-Cartesian. This follows immediately from the characterization given in Remark \ref{univsay}, since every edge of $X'$ is $p$-Cartesian when regarded as an edge of $X$.
\end{proof}

In order to verify that certain maps are Cartesian fibrations, it often convenient to work in a slightly more general setting. 

\begin{definition}\index{gen}{locally Cartesian!fibration}\index{gen}{fibration!locally Cartesian} A map $p: X \rightarrow S$ of simplicial sets is a {\it locally Cartesian fibration} if it is an inner fibration and, for every edge $\Delta^1 \rightarrow S$, the pullback $X \times_{S} \Delta^1 \rightarrow \Delta^1$ is a Cartesian fibration.
\end{definition}

In other words, an inner fibration $p: X \rightarrow S$ is a locally Cartesian fibration if and only if, for every vertex $x \in X$ and every edge $e: s \rightarrow p(x)$ in $S$, there exists a locally $p$-Cartesian edge $\overline{s} \rightarrow x$ which lifts $e$.

Let $p: X \rightarrow S$ be an inner fibration of simplicial sets. It is clear that every $p$-Cartesian morphism of $X$ is locally $p$-Cartesian. Moreover, Proposition \ref{protohermes} implies that the class of $p$-Cartesian edges of $X$ is stable under composition. Then following result can be regarded as a sort of converse:

\begin{lemma}\label{charloccart}
Let $p: X \rightarrow S$ be a locally Cartesian fibration of simplicial sets, and let
$f: x' \rightarrow x$ be an edge of $X$. The following conditions are equivalent:
\begin{itemize}
\item[$(1)$] The edge $e$ is $p$-Cartesian.
\item[$(2)$] For every $2$-simplex $\sigma$
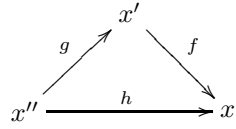
$$ \xymatrix{ & x' \ar[dr]^{f} & \\
x'' \ar[ur]^{g} \ar[rr]^{h} & & x }$$
in $X$, the edge $g$ is locally $p$-Cartesian if and only if the edge $h$ is locally $p$-Cartesian.
\item[$(3)$] For every $2$-simplex $\sigma$
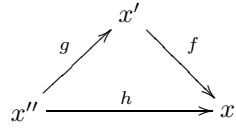
$$ \xymatrix{ & x' \ar[dr]^{f} & \\
x'' \ar[ur]^{g} \ar[rr]^{h} & & x }$$
in $X$, if $g$ is locally $p$-Cartesian, then $h$ is locally $p$-Cartesian.
\end{itemize}
\end{lemma}

\begin{proof}
We first show that $(1) \Rightarrow (2)$. Pulling back via the composition $p \circ \sigma: \Delta^2 \rightarrow S$, we can reduce to the case where $S = \Delta^2$. In this case, $g$ is locally $p$-Cartesian if and only if it is $p$-Cartesian, and likewise for $h$. We now conclude by applying Proposition \ref{protohermes}.

The implication $(2) \Rightarrow (3)$ is obvious. We conclude by showing that $(3) \Rightarrow (1)$. We must show that $\eta: X_{/f} \rightarrow
X_{/x} \times_{ S_{/p(x)} } S_{/p(f)}$ is a trivial fibration.
Since $\eta$ is a right fibration, it will suffice to
show that the fiber of $\eta$ over any vertex is contractible. Any such vertex determines a map
$\sigma: \Delta^2 \rightarrow S$ with $\sigma| \Delta^{ \{1,2\} } = p(f)$. Pulling back via
$\sigma$, we may suppose that $S= \Delta^2$. 

It will be convenient to introduce a bit of notation: for every map $q: K \rightarrow X$, let
$Y_{/q} \subseteq X_{/q}$ denote the full simplicial subset spanned by those vertices of
$X_{/q}$ which map to the initial vertex of $S$. 
We wish to show that the natural map
$Y_{/f} \rightarrow Y_{/x}$ is a trivial
fibration. By assumption, there exists a locally $p$-Cartesian morphism
$g: x'' \rightarrow x'$ in $X$ covering the edge $\Delta^{ \{0,1\} } \subseteq S$.
Since $X$ is an $\infty$-category, there exists a $2$-simplex $\tau: \Delta^2 \rightarrow X$ with
$d_2(\tau)=g$ and $d_0(\tau)=f$. Then $h = d_1(\tau)$ is a composite of $f$
and $g$, and assumption $(3)$ guarantees that $h$ is locally $p$-Cartesian. We have a commutative diagram
$$ 
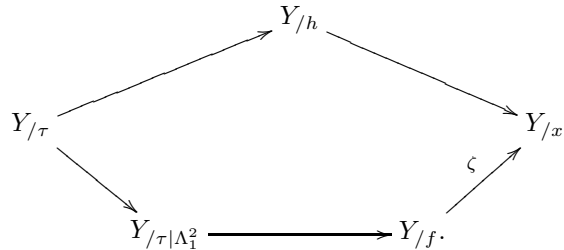
\xymatrix{ & & Y_{/h} \ar[drr]  & & \\
Y_{/\tau} \ar[urr] \ar[dr] & & & & Y_{/x} \\
& Y_{/\tau| \Lambda^2_1} \ar[rr] & & Y_{/f}. \ar[ur]^{\zeta} & }$$
Moreover, all of these maps in this diagram are trivial fibrations except possibly
$\zeta$, which is known to be a right fibration. It follows that $\zeta$ is a trivial fibration as well, which completes the proof.
\end{proof}

In fact, we have the following:

\begin{proposition}\label{gotta}
Let $p: X \rightarrow S$ be a locally Cartesian fibration. The following conditions are equivalent:
\begin{itemize}
\item[$(1)$] The map $p$ is a Cartesian fibration.
\item[$(2)$] Given a $2$-simplex
$$ \xymatrix{ x \ar[rr]^{f} \ar[dr]^{h} & & x' \ar[dl]^{g} \\
& z, & }$$
if $f$ and $g$ are locally $p$-Cartesian, then $h$ is locally $p$-Cartesian.
\item[$(3)$] Every locally $p$-Cartesian edge of $X$ is $p$-Cartesian.
\end{itemize}
\end{proposition} 
 
\begin{proof}
The equivalence $(2) \Leftrightarrow (3)$ follows from Lemma \ref{charloccart}, and the implication $(3) \Rightarrow (1)$ is obvious. To prove that $(1) \Rightarrow (3)$, let us suppose that
$e: x \rightarrow y$ is a locally $p$-Cartesian edge of $X$. Choose a $p$-Cartesian edge
$e': x' \rightarrow y$ lifting $p(e)$. The edges $e$ and $e'$ are both $p'$-Cartesian in
$X' = X \times_{S} \Delta^1$, where $p': X' \rightarrow \Delta^1$ denotes the projection. It follows that $e$ and $e'$ are equivalent in $X'$, and therefore also equivalent in $X$. Since $e'$ is $p$-Cartesian, we deduce that $e$ is $p$-Cartesian as well.
\end{proof}

\begin{remark}
If $p: X \rightarrow S$ is a locally Cartesian fibration, then we can associate to every edge
$s \rightarrow s'$ of $S$ a functor $X_{s'} \rightarrow X_{s}$, which is well-defined up to homotopy.
A $2$-simplex
$$ \xymatrix{ s \ar[rr] \ar[dr] & & s' \ar[dl] \\
& s'' & }$$
determines a triangle of $\infty$-categories
$$ \xymatrix{ X_{s} & & X_{s'} \ar[ll]^{F} \\
& X_{s''} \ar[ul]^{H} \ar[ur]^{G} & }$$
which commutes up to a (generally noninvertible) natural transformation $\alpha: F \circ G \rightarrow H$. Proposition \ref{gotta} implies that $p$ is a Cartesian fibration if and only if every such natural transformation is an equivalence of functors.
\end{remark}
 
\begin{corollary}
Let $p: X \rightarrow S$ be an inner fibration of simplicial sets. Then $p$ is Cartesian if and only if every pullback $X \times_{S} \Delta^n \rightarrow \Delta^n$ is
a Cartesian fibration, for $n \leq 2$.
\end{corollary}
 
One advantage the theory of locally Cartesian fibrations holds over the theory of Cartesian fibrations is the following ``fiberwise'' existence criterion:

\begin{proposition}\label{fibertest}
Suppose given a commutative diagram of simplicial sets.
$$ \xymatrix{ X \ar[dr]^{p} \ar[rr]^{r} & & Y \ar[dl]^{q} \\
& S & }$$
Suppose that:
\begin{itemize}
\item[$(1)$] The maps $p$ and $q$ are locally Cartesian fibrations, and $r$ is an inner fibration.
\item[$(2)$] The map $r$ carries locally $p$-Cartesian edges of $X$ to locally $q$-Cartesian edges of $Y$.
\item[$(3)$] For every vertex $s$ of $S$, the induced map $r_{s}: X_{s} \rightarrow Y_{s}$ is a
locally Cartesian fibration.
\end{itemize}
Then $r$ is a locally Cartesian fibration. Moreover, an edge $e$ of $X$ is locally $r$-Cartesian if and only if there exists a $2$-simplex $\sigma$
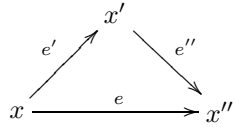
$$ \xymatrix{ & x' \ar[dr]^{e''} & \\
x \ar[ur]^{e'} \ar[rr]^{e} & & x'' }$$
with the following properties:
\begin{itemize}
\item[$(i)$] In the simplicial set $S$, we have $p( \sigma) = s^0( p(e))$.
\item[$(ii)$] The edge $e''$ is locally $p$-Cartesian.
\item[$(iii)$] The edge $e'$ is locally $r_{p(x)}$-Cartesian.
\end{itemize}
\end{proposition}

\begin{proof}
Suppose given a vertex $x'' \in X$ and an edge $e_0: y \rightarrow p(x'')$ in $Y$. It is clear that we can construct a $2$-simplex $\sigma$ in $X$ satisfying $(i)$ through $(iii)$, with
$p(e) = q(e_0)$. Moreover, $\sigma$ is uniquely determined up to equivalence. We will prove that $e$ is locally $r$-Cartesian. This will prove that $r$ is a locally Cartesian fibration, and the
``if'' direction of the final assertion. The converse will then follow from the uniqueness (up to equivalence) of locally $r$-Cartesian lifts of a given edge (with specified terminal vertex).
 
To prove that $e$ is locally $r$-Cartesian, we are free to pull back by the edge
$p(e): \Delta^1 \rightarrow S$, and thereby reduce to the case $S = \Delta^1$. Then
$p$ and $q$ are Cartesian fibrations. Since $e''$ is $p$-Cartesian and $r(e'')$ is $q$-Cartesian, Proposition \ref{stuch} implies that $e''$ is $r$-Cartesian. Remark \ref{intin} implies that
$e'$ is locally $p$-Cartesian. It follows from Lemma \ref{charloccart} that $e$ is locally $p$-Cartesian as well.
\end{proof}

\begin{remark}
The analogue of Proposition \ref{fibertest} for Cartesian fibrations is false.
\end{remark}

\subsection{Stability Properties of Cartesian Fibrations}\label{slib}

In this section, we will prove the class of Cartesian fibrations is stable under the formation of overcategories and undercategories. Since the definition of a Cartesian fibration is not self-dual, we must treat these results separately, using slightly different arguments (Propositions \ref{werylonger} and \ref{verylonger}). We begin with the following simple lemma.

\begin{lemma}\label{doweneed}
Let $A \subseteq B$ be an inclusion of simplicial sets. Then the inclusion
$$ (\{1\} \star B) \coprod_{ \{1\} \star A } (\Delta^1 \star A ) \subseteq \Delta^1 \star B$$ is inner anodyne.
\end{lemma}

\begin{proof}
Working by transfinite induction, we may reduce to the case where $B$ is obtained from $A$ by adjoining a single non-degenerate simplex, and therefore to the universal case
$B = \Delta^n$, $A = \bd \Delta^n$. Now the inclusion in question is isomorphic to $\Lambda^{n+2}_{1} \subseteq \Delta^{n+2}$.
\end{proof}

\begin{proposition}\label{werylonger}\index{gen}{Cartesian fibration!and overcategories}
Let $p: \calC \rightarrow \calD$ be a Cartesian fibration of simplicial sets, and let
$q: K \rightarrow \calC$ be a diagram. Then:
\begin{itemize}
\item[$(1)$] The induced map $p': \calC_{/q} \rightarrow \calD_{/pq}$ is a Cartesian fibration.
\item[$(2)$] An edge $f$ of $\calC_{/q}$ is $p'$-Cartesian if and only if the image
of $f$ in $\calC$ is $p$-Cartesian.
\end{itemize}
\end{proposition}

\begin{proof}
Proposition \ref{sharpen2} implies that $p'$ is an inner fibration. Let us call an edge
$f$ of $\calC_{q/}$ {\it special} if its image in $\calC$ is $p$Cartesian. To complete the proof, it will suffice to show that:
\begin{itemize}
\item[$(i)$] Given a vertex $\overline{q} \in \calC_{/q}$ and an edge
$\widetilde{f}: \overline{r}' \rightarrow p'( \overline{q} )$, there exists a special edge
$f: \overline{r} \rightarrow \overline{q}$ with $p'(f) = \widetilde{f}$. 
\item[$(ii)$] Every special edge of $\calC_{/q}$ is $p'$-Cartesian.
\end{itemize}

To prove $(i)$, let $\widetilde{f}'$ denote the image of $\widetilde{f}$ in $\calD$ and
$c$ the image of $\overline{q}$ in $\calC$. Using the assumption that $p$ is a coCartesian fibration, we can choose a $p$-coCartesian edge $f': c \rightarrow d$ lifting $\widetilde{f}'$. To extend this data to the desired edge $f$ of $\calC_{/q}$, it suffices to solve the lifting problem depicted in the diagram
$$ \xymatrix{ (\{1\} \star K) \coprod_{ \{1\} } \Delta^1 \ar[r] \ar@{^{(}->}[d]^{i} & \calC \ar[d]^{p} \\
\Delta^1 \star K  \ar[r] \ar@{-->}[ur] & \calD }$$
This lifting problem has a solution, since $p$ is an inner fibration and $i$ is inner anodyne (Lemma \ref{doweneed}).

To prove $(ii)$, it will suffice to show that if $n \geq 2$, then any lifting problem of the form
$$ \xymatrix{ \Lambda^n_n \star K \ar[r]^{g} \ar@{^{(}->}[d] & \calC \ar[d]^{p} \\
\Delta^n \star K \ar[r] \ar@{-->}[ur]^{G} & \calD }$$
has a solution, provided that $e=g( \Delta^{ \{n-1,n\} })$ is a $p$-Cartesian edge of $\calC$.
Consider the set $P$ of pairs $(K', G_{K'})$, where $K' \subseteq K$ and $G_{K'}$ fits in a commutative diagram 
$$ \xymatrix{ (\Lambda^n_n \star K) \coprod_{ \Lambda^n_n \star K'} (\Delta^n \star K') \ar[rrr]^-{G_{K'}} \ar@{^{(}->}[d] & & & \calC \ar[d]^{p} \\
\Delta^n \star K \ar[rrr] & & & \calD. }$$
Because $e$ is $p$-Cartesian, there exists an element
$(\emptyset, G_{\emptyset} ) \in P$. We regard $P$ as partially ordered, where
$(K', G_{K'} ) \leq (K'', G_{K''})$ if $K' \subseteq K''$ and $G_{K'}$ is a restriction of $G_{K''}$.
Invoking Zorn's lemma, we deduce the existence of a maximal element $(K', G_{K'})$ of $P$.
If $K' = K$, then the proof is complete. Otherwise, it is possible to enlarge $K'$ by adjoining a single nondegenerate $m$-simplex of $K$. Since $(K', G_{K''})$ is maximal, we conclude that the associated lifting problem
$$ \xymatrix{ (\Lambda^n_n \star \Delta^m) \coprod_{\Lambda^n_n \star \bd \Delta^m} (\Delta^n \star \bd \Delta^m) \ar[r] \ar@{^{(}->}[d] & \calC \ar[d]^{p} \\
\Delta^n \star \Delta^m \ar[r] \ar@{-->}[ur]^{\sigma} & \calD. }$$
has no solution.
The left vertical map is equivalent to the inclusion $\Lambda^{n+m+1}_{n+1} \subseteq \Delta^{n+m+1}$, which is inner anodyne. Since $p$ is an inner fibration by assumption, we obtain a contradiction.
\end{proof}

\begin{proposition}\label{verylonger}\index{gen}{coCartesian fibration!and overcategories}
Let $p: \calC \rightarrow \calD$ be a coCartesian fibration of simplicial sets, and let
$q: K \rightarrow \calC$ be a diagram. Then:
\begin{itemize}
\item[$(1)$] The induced map $p': \calC_{/q} \rightarrow \calD_{/pq}$ is a coCartesian fibration.
\item[$(2)$] An edge $f$ of $\calC_{/q}$ is $p'$-coCartesian if and only if the image
of $f$ in $\calC$ is $p$-coCartesian.
\end{itemize}
\end{proposition}

\begin{proof}
Proposition \ref{sharpen2} implies that $p'$ is an inner fibration. Let us call an edge
$f$ of $\calC_{/q}$ {\it special} if its image in $\calC$ is $p$-coCartesian. To complete the proof, it will suffice to show that:
\begin{itemize}
\item[$(i)$] Given a vertex $\overline{q} \in \calC_{/q}$ and an edge
$\widetilde{f}: p'( \overline{q} ) \rightarrow \overline{r}'$, there exists a special edge
$f: \overline{q} \rightarrow \overline{r}$ with $p'(f) = \widetilde{f}$. 
\item[$(ii)$] Every special edge of $\calC_{/q}$ is $p'$-coCartesian.
\end{itemize}

To prove $(i)$, we begin a commutative diagram
$$ \xymatrix{ \Delta^0 \star K \ar[r]^{\overline{q}} \ar@{^{(}->}[d] & \calC \ar[d] \\
\Delta^1 \star K \ar[r]^{\widetilde{f}} & \calD }. $$
Let $C \in \calC$ denote the image under $\overline{q}$ of the cone point of
$\Delta^0 \star K$, and choose a $p$-coCartesian morphism
$u: C \rightarrow C'$ lifting $\widetilde{f}| \Delta^1$. We now consider the
collection $P$ of all pairs $(L, f_L)$, where $L$ is a simplicial subset of $K$ and $f_L$ is a map fitting into a commutative diagram
$$ \xymatrix{ (\Delta^0 \star K) \coprod_{ \Delta^0 \star L} (\Delta^1 \star L) \ar[rrr]^-{f_L} \ar@{^{(}->}[d] & & & \calC \ar[d] \\
\Delta^1 \star K \ar[rrr]^{\widetilde{f}} & &  & \calD } $$
where $f_L | \Delta^1 = u$ and $f_L | \Delta^0 \star K = \overline{q}$. We partially order
the set $P$ as follows: $(L, f_L) \leq (L', f_{L'})$ if $L \subseteq L'$ and $f_{L}$ is equal
to the restriction of $f_{L'}$. The partially ordered set $P$ satisfies the hypotheses of Zorn's lemma, and therefore contains a maximal element $(L,f_L)$. If $L \neq K$, then we can choose a simplex
$\sigma: \Delta^n \rightarrow K$ of minimal dimension which does not belong to $L$. By maximality, we obtain a diagram
$$ \xymatrix{ \Lambda^{n+2}_{0} \ar[r] \ar@{^{(}->}[d] & \calC \ar[d] \\
\Delta^{n+2} \ar[r] \ar@{-->}[ur] & \calD }$$
in which the indicated dotted arrow cannot be supplied. This is a contradiction, since
the upper horizontal map carries the initial edge of $\Lambda^{n+2}_0$ to a $p$-coCartesian
edge of $\calC$. It follows that $L = K$, and we may take $f = f_L$. This completes the proof of $(i)$.

The proof of $(ii)$ is similar. Suppose given $n \geq 2$ and a diagram
$$ \xymatrix{ \Lambda^n_0 \star K \ar[r]^{f_0} \ar@{^{(}->}[d] & \calC \ar[d] \\
\Delta^n \star K \ar[r]^{g} \ar@{-->}[ur]^{f} & \calD }$$
be a commutative diagram, where $f_0 | K = q$ and $f_0 | \Delta^{ \{0,1\} }$ is a $p$-coCartesian edge of $\calC$. We wish to prove the existence of the dotted arrow $f$, indicated in the diagram.
As above, we consider the collection $P$ of all pairs $(L, f_L)$, where $L$ is a simplicial subset of $K$ and $f_L$ extends $f_0$ and fits into a commutative diagram
$$ \xymatrix{ (\Lambda^n_0 \star K) \coprod_{ \Lambda^n_0 \star L} (\Delta^n \star L) \ar[rrr]^-{f_L} \ar@{^{(}->}[d] & & &  \calC \ar[d] \\
\Delta^n \star K \ar[rrr]^{g} & & & \calD. }$$
We partially order $P$ as follows: $(L, f_L) \leq (L', f_{L'})$ if $L \subseteq L'$ and $f_L$ is a restriction of $f_{L'}$. Using Zorn's lemma, we conclude that $P$ contains a maximal element
$(L, f_L)$. If $L \neq K$, then we can choose a simplex $\sigma: \Delta^m \rightarrow K$ which does not belong to $L$, where $m$ is as small as possible. Invoking the maximality of $(L,f_L)$, we obtain a diagram
$$ \xymatrix{ \Lambda^{n+m+1}_{0} \ar[r]^{h} \ar@{^{(}->}[d] & \calC \ar[d] \\
\Delta^{n+m+1} \ar[r] \ar@{-->}[ur] & \calD }$$
where the indicated dotted arrow cannot be supplied. However, the map $h$ carries the initial edge of $\Delta^{n+m+1}$ to a $p$-coCartesian edge of $\calC$, so we obtain a contradiction. It follows that $L = K$, so that we can take $f = f_L$ to complete the proof.
\end{proof}

\subsection{Mapping Spaces and Cartesian Fibrations}\label{slik}

Let $p: \calC \rightarrow \calD$ be a functor between $\infty$-categories, and let
$X$ and $Y$ be objects of $\calC$. Then $p$ induces a map
$$ \phi: \bHom_{\calC}(X,Y) \rightarrow \bHom_{\calD}(p(X),p(Y)).$$
Our goal in this section is to understand the relationship between the fibers of $p$ and the {\em homotopy} fibers of $\phi$.

\begin{lemma}\label{sharpy}
Let $p: \calC \rightarrow \calD$ be an inner fibration of $\infty$-categories, and let $X,Y \in \calC$. The induced map $\phi: \Hom^{\rght}_{\calC}(X,Y) \rightarrow \Hom^{\rght}_{\calD}(p(X),p(Y))$ is a Kan fibration.
\end{lemma}

\begin{proof}
Since $p$ is an inner fibration, the induced map $\widetilde{\phi}: \calC_{/X} \rightarrow \calD_{/p(X)} \times_{\calD} \calC$ is a right fibration by Proposition \ref{sharpen}. We note that $\phi$ is obtained from $\widetilde{\phi}$ by restricting to the fiber over the vertex $Y$ of $\calC$. Thus $\phi$ is a right fibration; since the target of $\phi$ is a Kan complex, $\phi$ is a Kan fibration by Lemma \ref{toothie2}.
\end{proof}

Suppose the conditions of Lemma \ref{sharpy} are satisfied. Let us consider the problem of computing the fiber of $\phi$ over a vertex $\overline{e}: p(X) \rightarrow p(Y)$ of $\Hom^{\rght}_{\calD}(X,Y)$. 
Suppose that there is a $p$-Cartesian edge $e: X' \rightarrow Y$ lifting $\overline{e}$. By definition, we have
a trivial fibration
$$ \psi: \calC_{/e} \rightarrow \calC_{/Y} \times_{ \calD_{/p(Y)}} \calD_{/\overline{e}}.$$
Consider the $2$-simplex $\sigma = s_1(\overline{e})$, regarded as a vertex of $\calD_{/\overline{e}}$. Passing to the fiber, we obtain a trivial fibration
$$ F \rightarrow \phi^{-1}(e),$$ where $F$ denotes the fiber of $\calC_{/e} \rightarrow \calD_{/\overline{e}} \times_{\calD} \calC$ over the point $(\sigma,X)$.
On the other hand, we have a trivial fibration
$\calC_{/e} \rightarrow \calD_{/\overline{e}} \times_{ \calD_{/p(X)} } \calC_{/X'}$ by Proposition \ref{sharpen2}. Passing to the fiber again, we obtain a trivial fibration $F \rightarrow \Hom^{\rght}_{\calC_{p(X)}}(X,X')$. We may summarize the situation as follows:

\begin{proposition}\label{compspaces}
Let $p: \calC \rightarrow \calD$ be an inner fibration of $\infty$-categories. Let $X,Y \in \calC$, let
$\overline{e}: p(X) \rightarrow p(Y)$ be a morphism in $\calD$, and let $e: X' \rightarrow Y$ be a locally $p$-Cartesian morphism of $\calC$ lifting $\overline{e}$. Then in the homotopy category $\calH$ of spaces, there is
a fiber sequence
$$ \bHom_{\calC_{p(X)}}(X,X') \rightarrow \bHom_{\calC}(X,Y) \rightarrow \bHom_{\calD}(p(X),p(Y)).$$
Here the fiber is taken over the point classified by $\overline{e}: p(X) \rightarrow p(Y)$. 
\end{proposition}

\begin{proof}
The edge $\overline{e}$ defines a map $\Delta^1 \rightarrow \calD$. Note that the fiber
of the Kan fibration $\Hom^{\rght}_{\calC}(X,Y) \rightarrow \Hom^{\rght}_{\calD}(pX, pY)$ does not change if we replace $p$ by the induced projection $\calC \times_{\calD} \Delta^1 \rightarrow \Delta^1$.
We may therefore assume without loss of generality that $e$ is $p$-Cartesian, and the desired result follows from the above analysis.
\end{proof}

A similar assertion can be taken as a characterization of Cartesian morphisms:

\begin{proposition}\label{charCart}
Let $p: \calC \rightarrow \calD$ be an inner fibration of $\infty$-categories, and let
$f: Y \rightarrow Z$ be a morphism in $\calC$. The following are equivalent:
\begin{itemize}
\item[$(1)$] The morphism $f$ is $p$-Cartesian.
\item[$(2)$] For every object $X$ of $\calC$, composition with $f$ gives rise to a homotopy
Cartesian diagram
$$ \xymatrix{ \bHom_{\calC}(X,Y) \ar[r] \ar[d] & \bHom_{\calC}(X,Z) \ar[d] \\
\bHom_{\calD}( p(X), p(Y) ) \ar[r] & \bHom_{\calD}(p(X), p(Z)).} $$
\end{itemize}
\end{proposition}

\begin{proof}
Let $\phi: \calC_{/f} \rightarrow \calC_{/Z} \times_{\calD_{/p(Z)}} \calD_{/p(f)}$ be the canonical map; then $(1)$ is equivalent to the assertion that $\phi$ is a trivial fibration. According to Proposition \ref{sharpen}, $\phi$ is a right fibration. Thus, $\phi$ is a trivial fibration if and only if the fibers of $\phi$ are contractible Kan complexes. For each object $X \in \calC$, let 
$$\phi_X:  \calC_{/f} \times_{\calC} \{X \} \rightarrow \calC_{/Z} \times_{\calD_{/p(Z)}} \calD_{/p(f)}
\times_{\calC} \{X\}$$ be the induced map. Then $\phi_{X}$ is a right fibration between Kan complexes, and therefore a Kan fibration; it has contractible fibers if and only if it is a homotopy equivalence. Thus, $(1)$ is equivalent to the assertion that $\phi_{X}$ is a homotopy equivalence for every object $X$ of $\calC$.

We remark that $(2)$ is somewhat imprecise: although all the maps in the diagram are well defined in the homotopy category $\calH$ of spaces, we need to represent this by a commutative diagram in the category of simplicial sets before we can ask whether or not the diagram is homotopy Cartesian. We therefore rephrase $(2)$ more precisely: it asserts that the diagram of Kan complexes
$$ \xymatrix{ \calC_{/f} \times_{\calC} \{X\} \ar[r] \ar[d] & \calC_{/Z} \times_{\calC} \{X\} \ar[d] \\
\calD_{/p(f)} \times_{\calD} \{p(X)\} \ar[r] & \calD_{/p(Z)} \times_{\calD} \{p(X)\} }$$
is homotopy Cartesian. Lemma \ref{sharpy} implies that the right vertical map is a Kan fibration, so the homotopy limit in question is given by the fiber product $$\calC_{/Z} \times_{\calD_{/p(Z)}} \calD_{/p(f)} \times_{\calC} \{X\}.$$ 
Consequently, assertion $(2)$ is also equivalent to the condition that
$\phi_{X}$ be a homotopy equivalence for every object $X \in \calC$.
\end{proof}

\begin{corollary}\label{usefir}
Suppose given maps $\calC \stackrel{p}{\rightarrow} \calD \stackrel{q}{\rightarrow} \calE$ of $\infty$-categories, such that both $q$ and $q \circ p$ are locally Cartesian fibrations. Suppose that $p$ carries
locally $(q \circ p)$-Cartesian edges of $\calC$ to locally $q$-Cartesian edges of $\calD$, and that for every object $Z \in \calE$, the induced map $\calC_{Z} \rightarrow \calD_{Z}$ is a categorical equivalence. Then $p$ is a categorical equivalence.
\end{corollary}

\begin{proof}
Proposition \ref{compspaces} implies that $p$ is fully faithful. If $Y$ is any object of $\calD$, then $Y$ is equivalent in the fiber $\calD_{q(Y)}$ to the image under $p$ of some vertex of $\calC_{q(Y)}$. Thus $p$ is essentially surjective and the proof is complete.
\end{proof}

\begin{corollary}\label{usesec}
Let $p: \calC \rightarrow \calD$ be a Cartesian fibration of $\infty$-categories. Let $q: \calD' \rightarrow \calD$ be a categorical equivalence of $\infty$-categories. Then the induced map $q': \calC' = \calD' \times_{\calD} \calC \rightarrow \calC$ is a categorical equivalence.
\end{corollary}

\begin{proof}
Proposition \ref{compspaces} immediately implies that $q'$ is fully faithful. We claim that $q'$ is essentially surjective. Let $X$ be any object of $\calC$. Since $q$ is fully faithful, there exists an object $y$ of $T'$ and an equivalence $\overline{e}: q(Y) \rightarrow p(X)$. Since $p$ is Cartesian, we can choose a $p$-Cartesian edge $e: Y' \rightarrow X$ lifting $\overline{e}$. Since $e$ is $p$-Cartesian and $p(e)$ is an equivalence, $e$ is an equivalence. By construction, the object $Y'$ of $S$ lies in the image of $q'$.
\end{proof}

\begin{corollary}\label{heath}\index{gen}{Cartesian fibration!and trivial fibrations}
Let $p: \calC \rightarrow \calD$ be a Cartesian fibration of $\infty$-categories. Then $p$ is a categorical equivalence if and only if $p$ is a trivial fibration.
\end{corollary}

\begin{proof}
The ``if'' direction is clear. Suppose then that $p$ is a categorical equivalence. We first claim
that $p$ is surjective on objects. The essential surjectivity of $p$ implies that for each $Y \in \calD$ there is an equivalence $Y \rightarrow p(X)$, for some object $X$ of $\calC$. Since $p$ is Cartesian, this equivalence lifts to a $p$-Cartesian edge $\widetilde{Y} \rightarrow X$ of $S$, so that $p(\widetilde{Y}) = Y$.

Since $p$ is fully faithful, the map $\bHom_{\calC}(X,X') \rightarrow \bHom_{\calD}(p(X),p(X'))$ is a homotopy equivalence
for any pair of objects $X,X' \in \calC$. Suppose that $p(X) = p(X')$. Then, applying Proposition \ref{compspaces}, we deduce that $\bHom_{\calC_{p(X)}}(X,X')$ is contractible.
It follows that the $\infty$-category $\calC_{p(X)}$ is nonempty with contractible morphism spaces; it is therefore a contractible Kan complex. Proposition \ref{goey} now implies that $p$ is a right fibration. Since $p$ has contractible fibers, it is a trivial fibration by Lemma \ref{toothie}.
\end{proof}

We have already seen that if a $\infty$-category $S$ has an initial
object, then that initial object is essentially unique. We now establish a relative version of this
result. 


\begin{lemma}\label{sabreto}
Let $p: \calC \rightarrow \calD$ be a Cartesian fibration of $\infty$-categories, and let
$C$ be an object of $\calC$. Suppose that $D = p(C)$ is an initial object of $\calD$, and that
$C$ is an initial object of the $\infty$-category $\calC_{D} = \calC \times_{\calD} \{D\}$. 
Then $C$ is an initial object of $\calC$.
\end{lemma}

\begin{proof}
Let $C'$ be any object of $\calC$, and let $D' = p(C')$. Since $D$ is an initial object of
$\calD$, the space $\bHom_{\calD}(D,D')$ is contractible. In particular, there
exists a morphism $f: D \rightarrow D'$ in $\calD$. Let $\widetilde{f}: \widetilde{D} \rightarrow C'$
be a $p$-Cartesian lift of $f$. According to Proposition \ref{compspaces}, there exists a fiber sequence in the homotopy category $\calH$:
$$ \bHom_{\calC_{D}}(C, \widetilde{D}) \rightarrow \bHom_{\calC}(C,C') \rightarrow
\bHom_{\calD}(D,D').$$
Since the first and last space in the sequence are contractible, we deduce that $\bHom_{\calC}(C,C')$ is contractible as well, so that $C$ is an initial object of $\calC$. 
\end{proof}

\begin{lemma}\label{sabretooth}
Suppose given a diagram of simplicial sets
$$ \xymatrix{ \bd \Delta^n \ar[r]^{f_0} \ar@{^{(}->}[d] & X \ar[d]^{p} \\
\Delta^n \ar@{-->}[ur]^{f} \ar[r]^{g} & S }$$
where $p$ is a Cartesian fibration and $n > 0$. Suppose that
$f_0(0)$ is an initial object of the $\infty$-category $X_{g(0)} = X \times_{S} \{ g(0) \}$.
Then there exists a map $f: \Delta^n \rightarrow S$ as indicated by the dotted arrow in the diagram, which renders the diagram commutative.
\end{lemma}

\begin{proof}
Pulling back via $g$, we may replace $S$ by $\Delta^n$ and thereby reduce to the case where $S$ is an $\infty$-category and $g(0)$ is an initial object of $S$. It follows from Lemma \ref{sabreto} that 
$f_0(v)$ is an initial object of $S$, which implies the existence of the desired extension $f$.
\end{proof}

\begin{proposition}\label{topaz}
Let $p: X \rightarrow S$ be a Cartesian fibration of
simplicial sets. Assume that, for each vertex $s$ of $S$, the
$\infty$-category $X_{s} = X \times_{S} \{s\}$ has an initial object. 
\begin{itemize}
\item[$(1)$] Let $X' \subseteq X$
denote the full simplicial subset of $X$ spanned by those vertices $x$ which are initial objects of $X_{p(x)}$. Then $p|X'$ is a trivial fibration of simplicial sets.
\item[$(2)$] Let $\calC = \bHom_{S}(S, X)$ be the $\infty$-category of sections of $p$.
An arbitrary section $q: S \rightarrow X$ is an initial object of $\calC$ if and only if
$q$ factors through $X'$.
\end{itemize}
\end{proposition}

\begin{proof}
Since every fiber $X_{s}$ has an initial object, the map $p|X'$ has the right lifting
property with respect to the inclusion $\emptyset \subseteq \Delta^0$. If $n > 0$,
then Lemma \ref{sabretooth} shows that $p|X'$ has the right lifting property with
respect to $\bd \Delta^n \subseteq \Delta^n$. This proves $(1)$. In particular, we deduce
that there exists a map $q: S \rightarrow X'$ which is a section of $p$. In view of the uniqueness of initial objects, $(2)$ will follow if we can show that $q$ is an initial object of $\calC$.
Unwinding the definitions, we must show that for $n > 0$, any lifting problem
$$ \xymatrix{ S \times \bd \Delta^n \ar[r]^{f} \ar@{^{(}->}[d] & X \ar[d]^{q} \\
S \times \Delta^n \ar[r] \ar@{-->}[ur] & S}$$
can be solved, provided that $f | S \times \{0\} = q$. The desired extension can be constructed simplex-by-simplex, using Lemma \ref{sabretooth}. 
\end{proof}

\subsection{Application: Invariance of  Undercategories}\label{slim}

Our goal in this section is to complete the proof of Proposition \ref{gorban3} by proving the
following assertion:

\begin{itemize}\index{gen}{undercategory}
\item[$(\ast)$] Let $p: \calC \rightarrow \calD$ be an equivalence of $\infty$-categories, and let
$j: K \rightarrow \calC$ be a diagram. Then the induced map
$$ \calC_{j/} \rightarrow \calD_{p j/}$$
is a categorical equivalence.
\end{itemize}

We will need a lemma.

\begin{lemma}\label{blaha}
Let $p: \calC \rightarrow \calD$ be a fully faithful map of $\infty$-categories, and let $j: K \rightarrow \calC$
be any diagram in $\calC$. Then, for any object $x$ of $\calC$, the map of Kan complexes
$$ \calC_{j/} \times_{\calC} \{x\} \rightarrow \calD_{p j/ } \times_{\calD} \{p(x) \}$$ is a homotopy equivalence.
\end{lemma}

\begin{proof}
For any map $r: K' \rightarrow K$ of simplicial sets, let $C_{r} = \calC_{j r/} \times_{\calC} \{x\}$ and $D_{r} = \calD_{pjr/} \times_{\calD} \{p(x)\}$. 

Choose a transfinite sequence of simplicial subsets $K_{\alpha}$ of $K$, such that $K_{\alpha+1}$ is the result of adjoining a single nondegenerate simplex to $K_{\alpha}$, and $K_{\lambda} = \bigcup_{\alpha<\lambda} K_{\alpha}$ whenever $\lambda$ is a limit ordinal (we include the case where $\lambda = 0$, so that $K_0 = \emptyset$). Let $i_{\alpha}: K_{\alpha} \rightarrow K$ denote the inclusion. We claim the following:
\begin{itemize}
\item[$(1)$] For every ordinal $\alpha$, the map $\phi_{\alpha}: C_{i_{\alpha}} \rightarrow D_{i_{\alpha}}$ is a homotopy equivalence of simplicial sets.
\item[$(2)$] For every pair of ordinals $\beta \leq \alpha$, the maps $C_{i_{\alpha}} \rightarrow C_{i_{\beta}}$ and $D_{i_{\alpha}} \rightarrow D_{i_{\beta}}$ are Kan fibrations of simplicial sets.
\end{itemize}

We prove both of these claims by induction on $\alpha$. When $\alpha = 0$, $(2)$ is obvious and $(1)$ follows since both sides are isomorphic to $\Delta^0$. If $\alpha$ is a limit ordinal, $(2)$ is again obvious, while $(1)$ follows from the fact that both $C_{i_{\alpha}}$ and $D_{i_{\alpha}}$
are obtained as the inverse limit of a transfinite sequence of fibrations, and the map $\phi_{\alpha}$ is an inverse limit of maps which are individually homotopy equivalences.

Assume that $\alpha=\beta+1$ is a successor ordinal, so that $K_{\alpha} \simeq K_{\beta} \coprod_{ \bd \Delta^n } \Delta^n$. Let $f: \Delta^n \rightarrow K_{\alpha}$ be the induced map, so that
$$ C_{i_\alpha} = C_{i_{\beta}} \times_{ C_{ f | \bd \Delta^n } } C_{f} $$
$$ D_{i_\alpha} = D_{i_{\beta}} \times_{ D_{ f| \bd \Delta^n }} D_{f}.$$
We note that the projections $C_{f} \rightarrow D_{f | \bd \Delta^n }$ and $C_{f} \rightarrow D_{f| \bd \Delta^n}$ are left fibrations by Proposition \ref{sharpen}, and therefore Kan fibrations by Lemma \ref{toothie2}. This proves $(2)$, since the class of Kan fibrations is stable under pullback.
We also note that the pullback diagrams defining $X_{i_{\alpha}}$ and $Y_{i_{\alpha}}$ are also  homotopy pullback diagrams. Thus, to prove that $\phi_{\alpha}$ is a homotopy equivalence, it suffices to show that $\phi_{\beta}$ and the maps
$$ C_{f| \bd \Delta^n } \rightarrow D_{ f| \bd \Delta^n }$$
$$ C_{f} \rightarrow D_{f}$$
are homotopy equivalences. In other words, we may reduce to the case where $K$ is a {\em finite} complex. 

We now work by induction on the dimension of $K$. Suppose that the dimension of $K$ is $n$, and that the result is known for all simplicial sets having smaller dimension. Running through the above argument again, we can reduce to the case where $K = \Delta^n$. Let $v$ denote the final vertex of $\Delta^n$. By Proposition \ref{sharpen2}, the maps
$$ C_{j} \rightarrow C_{ j| \{v\} }$$
$$ D_{j} \rightarrow D_{ j| \{v\} }$$
are trivial fibrations. Thus, it suffices to consider the case where $K$ is a single point $\{v\}$. 
In this case, we have $C_{j} = \Hom^{\lft}_{\calC}( j(v), x)$ and $Y_{j} = \Hom^{\lft}_{\calD}( p(j(v)), p(x))$.
It follows that the map $\phi$ is a homotopy equivalence, since $p$ is assumed fully faithful.
\end{proof}

\begin{proof}[Proof of $(\ast)$]
Let $p: \calC \rightarrow \calD$ be a categorical equivalence of $\infty$-categories, and $j: K \rightarrow \calC$ any diagram.
We have a factorization
$$ \calC_{j/} \stackrel{f}{\rightarrow} \calD_{p j/} \times_{\calD} \calC \stackrel{g}{\rightarrow} \calD_{pj/}.$$
Lemma \ref{blaha} implies that $\calC_{j/}$ and $\calD_{p j/} \times_{\calD} \calC$ are fiberwise equivalent left-fibrations over $\calC$, so that $f$ is a categorical equivalence by Corollary \ref{usefir} (we note that the map $f$
automatically carries coCartesian edges to coCartesian edges, since {\em all} edges of the target $\calD_{p j/} \times_{\calD} \calC$ are coCartesian). The map $g$ is a categorical equivalence by Corollary \ref{usesec}. It follows that $g \circ f$ is a categorical equivalence, as desired.
\end{proof} 

\subsection{Application: Categorical Fibrations over a Point}\label{slin}

Our main goal in this section is to prove the following result:

\begin{theorem}\label{joyalcharacterization}\index{gen}{$\infty$-category!as a fibrant object of $\sSet$}
Let $\calC$ be a simplicial set. Then $\calC$ is fibrant for the Joyal model structure if and only if $\calC$ is an $\infty$-category.
\end{theorem}

The proof will require a few technical preliminaries.

\begin{lemma}\label{gorban2}
Let $p: \calC \rightarrow \calD$ be a categorical equivalence of $\infty$-categories and $m \geq 2$ an integer. Suppose given maps $f_0: \bd \Delta^{ \{1, \ldots, m\} } \rightarrow \calC$ and $h_0: \Lambda^m_0 \rightarrow \calD$ with
$h_0 | \bd \Delta^{ \{1, \ldots, m\} } = p \circ f_0$. Suppose further that the restriction of $h$ to $\Delta^{ \{0,1\} }$ is an equivalence in $\calD$. Then there exist maps $f: \Delta^{ \{1, \ldots, m\} } \rightarrow \calC$, $h: \Delta^{m} \rightarrow \calD$, with $h| \Delta^{ \{1, \ldots, n\} }= p \circ f$, $f_0 = f| \bd \Delta^{ \{1, \ldots, m\} }$, $h_0 = h| \Lambda^m_0$.
\end{lemma}

\begin{proof}
We may regard $h_0$ as a point of the simplicial set $\calD_{/ p \circ f_0}$. Since $p$ is a categorical equivalence, Proposition \ref{gorban3} implies that $p': \calC_{/f_0} \rightarrow \calD_{/p \circ f_0}$ is a categorical equivalence. It follows that $h_0$ lies in the essential image of $p'$. Consider the linearly ordered set
$\{ 0 < 0' < 1 < \ldots < n\}$ and the corresponding simplex $\Delta^{ \{0, 0', \ldots, n\} }$. By hypothesis, we can extend
$f_0$ to a map $f'_0: \Lambda^{ \{ 0', \ldots, m \} }_{0'} \rightarrow \calC$ and $h_0$ to a map $h'_0: \Delta^{ \{0,0'\} } \star \bd \Delta^{ \{1, \ldots, m \} } \rightarrow \calD$ such that $h'_0| \Delta^{ \{0,0'\} }$ is an equivalence
and $h'_0| \Lambda^{ \{0', \ldots, m\} }_0 = p \circ f'_0$. 

Since $h'_0| \Delta^{ \{0,0'\} }$ and $h'_0| \Delta^{ \{0,1\} }$ are both equivalences in $\calD$, we deduce that $h'_0| \Delta^{ \{0',1\} }$ is an equivalence in $\calD$. Since $p$ is a categorical equivalence, it follows that
$f'_0| \Delta^{ \{0',1\} }$ is an equivalence in $\calC$. Proposition \ref{greenlem} implies that $f'_0$ extends to a map $f': \Delta^{ \{0', \ldots, m\} } \rightarrow \calC$. The union of $p \circ f'$ and $h'_0$ determines a map
$\Lambda^{ \{0, 0', \ldots, m \} }_{0'} \rightarrow \calD$; since $\calD$ is an $\infty$-category, this extends to a map
$h': \Delta^{ \{0,0', \ldots, m\} } \rightarrow \calD$. We may now take $f = f' | \Delta^{ \{1, \ldots, m \} }$ and
$h = h' | \Delta^{m}$.
\end{proof}

\begin{lemma}\label{gorban}
Let $p: \calC \rightarrow \calD$ be a categorical equivalence of $\infty$-categories and $A \subseteq B$ any inclusion of simplicial sets. Let $f_0: A \rightarrow \calC$, $g: B \rightarrow \calD$ be any maps, and let $h_0: A \times \Delta^1 \rightarrow \calD$ be an equivalence from $g|A$  to $p \circ f_0$. Then there exists a map $f: B \rightarrow \calC$ and an equivalence $h: B \times \Delta^1 \rightarrow \calD$ from $g$ to $p \circ f$, such that
$f_0 = f|A$ and $h_0 = h| A \times \Delta^1$.
\end{lemma}

\begin{proof}
Working cell-by-cell with the inclusion $A \subseteq B$, we may reduce to the case where $B = \Delta^n$, $A = \bd \Delta^n$. If $n = 0$, the existence of the desired extensions is a reformulation of the assumption that $p$ is essentially surjective. Let us assume therefore that $n \geq 1$.

We consider the task of constructing $h: \Delta^n \times \Delta^1 \rightarrow \calD$. Consider the filtration
$$ X(n+1) \subseteq \ldots \subseteq X(0) = \Delta^n \times \Delta^1 $$
described in the proof of Proposition \ref{usejoyal}. We note that the value of $h$ on $X(n+1)$ is uniquely prescribed by $h_0$ and $g$. We extend the definition of $h$ to $X(i)$ by descending induction on $i$. We note that $X(i) \simeq X(i+1) \coprod_{ \Lambda^{n+1}_k} \Delta^{n+1}$. For $i > 0$, the existence of the required extension is guaranteed by the assumption that $\calD$ is an $\infty$-category. Since $n \geq 1$, Lemma \ref{gorban2}  allows us to extend $h$ over the simplex $\sigma_0$ and to define $f$ so that the desired conditions are satisfied.
\end{proof}

\begin{lemma}\label{gorbann}
Let $\calC \subseteq \calD$ be an inclusion of simplicial sets which is also a categorical equivalence. Suppose further that $\calC$ is an $\infty$-category. Then $\calC$ is a retract of $\calD$.
\end{lemma}

\begin{proof}
Enlarging $\calD$ by an inner anodyne extension if necessary, we may suppose that $\calD$ is an $\infty$-category. We now apply Lemma \ref{gorban} in the case
where $A = \calC$, $B = \calD$.
\end{proof}

\begin{proof}[Proof of Theorem \ref{joyalcharacterization}]
The ``only if'' direction has already been established (Remark \ref{tokenn}). 
For the converse, we must show that if $\calC$ is an $\infty$-category, then $\calC$ has the extension property with respect to every inclusion of simplicial sets
$A \subseteq B$ which is a categorical equivalence. Fix any map $A \rightarrow \calC$. Since the Joyal model structure is left-proper, the inclusion
$\calC \subseteq \calC \coprod_{A} B$ is a categorical equivalence. We now apply Lemma \ref{gorbann} to conclude that $\calC$ is a retract of $\calC \coprod_{A} B$.
\end{proof}

We can state Theorem \ref{joyalcharacterization} as follows: if $S$ is a point, then
$p: X \rightarrow S$ is a categorical fibration (in other words, a fibration with respect to the Joyal model structure on $\SSet$) if and only if it is an inner fibration. However, the class of inner fibrations does {\em not} coincide with the class of categorical fibrations in general. The following result describes the situation when $T$ is an $\infty$-category:

\begin{corollary}[Joyal]\index{gen}{categorical fibration!of $\infty$-categories}\label{gottaput}
Let $p: \calC \rightarrow \calD$ be a map of simpicial sets, where $\calD$ is an $\infty$-category.
Then $p$ is a categorical fibration if and only if the following conditions are satisfied:
\begin{itemize}
\item[$(1)$] The map $p$ is an inner fibration.
\item[$(2)$] For every equivalence $f: D \rightarrow D'$ in $\calD$, and every
object $C \in \calC$ with $p(C) = D$, there exists an equivalence
$\overline{f}: C \rightarrow C'$ in $\calC$ with $p( \overline{f}) = f$.
\end{itemize}
\end{corollary}

\begin{proof}
Suppose first that $p$ is a categorical fibration. Then $(1)$ follows immediately (since the inclusions $\Lambda^n_i \subseteq \Delta^n$ are categorical equivalences for $0 < i < n$).
To prove $(2)$, we let $\calD^{0}$ denote the largest Kan complex contained in $\calD$, so that
the edge $f$ belongs to $\calD$. There exists a contractible Kan complex $K$ containing
an edge $\widetilde{f}: \widetilde{D} \rightarrow \widetilde{D}'$ and a map $q: K \rightarrow \calD$ such that $q( \widetilde{f} ) = f$. Since the inclusion $\{ \widetilde{\calD} \} \subseteq K$ is a categorical equivalence, our assumption that $p$ is a categorical fibration allows us to lift
$q$ to a map $\widetilde{q}: K \rightarrow \calC$ such that $\widetilde{q}(\widetilde{D})=C$.
We can now take $\overline{f} = \widetilde{q}( \widetilde{f} )$; since $\widetilde{f}$ is an equivalence in $K$, $\overline{f}$ is an equivalence in $\calC$.

Now suppose that $(1)$ and $(2)$ are satisfied. We wish to show that $p$ is a categorical fibration. Consider a lifting problem
$$ \xymatrix{ A \ar@{^{(}->}[d]^{i} \ar[r]^{g_0} & \calC \ar[d]^{p} \\
B \ar[r]^{h} \ar@{-->}[ur]^{g} & \calD }$$
where $i$ is a cofibration and a categorical equivalence; we wish to show that 
there exists a morphism $g$ as indicated which renders the diagram commutative.
We first observe that condition $(1)$, together with our assumption that $\calD$ is an $\infty$-category, guarantee that $\calC$ is an $\infty$-category.
Applying Theorem \ref{joyalcharacterization}, we can extend $g_0$ to a map
$g': B \rightarrow \calC$, not necessarily satisfying $h = p \circ g'$.
Nevertheless, the maps $h$ and $p \circ g'$ have the same restriction to $A$.
Let 
$$H_0: (B \times \bd \Delta^1) \coprod_{ A \times \bd \Delta^1 } (A \times \Delta^1) \rightarrow \calD$$
be given by $(p \circ g',h)$ on $B \times \bd \Delta^1$, and by the composition
$$ A \times \Delta^1 \rightarrow A \subseteq B \stackrel{h}{\rightarrow} \calD$$
on $A \times \Delta^1$. Applying Theorem \ref{joyalcharacterization} once more, we deduce
that $H_0$ extends to a map $H: B \times \Delta^1 \rightarrow \calD$. The map $H$ carries
$\{a \} \times \Delta^1$ to an equivalence in $\calD$, for every vertex $a$ of $A$. Since
the inclusion $A \subseteq B$ is a categorical equivalence, we deduce that
$H$ carries $\{b \} \times \Delta^1$ to an equivalence, for every $b \in B$.

Let $$G_0: (B \times \{0\}) \coprod_{ A \times \{0\} } (A \times \Delta^1) \rightarrow \calC$$
be the composition of the projection to $B$ with the map $g'$. We have a commutative diagram
$$ \xymatrix{  (B \times \{0\}) \coprod_{ A \times \{0\} } (A \times \Delta^1) \ar[r]^-{G_0} \ar[d] & \calC \ar[d]^{p} \\
B \times \Delta^1 \ar[r]^{H} \ar@{-->}[ur]^{G} & \calD. }$$
To complete the proof, it will suffice to show that we can supply a map $G$ as indicated, rendering the diagram commutative; in this case, we can solve the original lifting problem by defining
$g = G | B \times \{1\}$.

We construct the desired extension $G$ working cell-by-cell on $B$. We start by applying assumption $(2)$ to construct the map $G| \{b\} \times \Delta^1$ for every vertex $b$ of $B$ (that does not already belong to $A$); moreover, we ensure that $G| \{b\} \times \Delta^1$ is an equivalence in $\calC$.

To extend $G_0$ to simplices of higher dimension, we encounter lifting problems of the type
$$ \xymatrix{ ( \Delta^n \times \{0\} ) \coprod_{ \bd \Delta^n \times \{0\} } ( \bd \Delta^n \times \Delta^1 ) \ar[r]^-{e} \ar@{^{(}->}[d] & \calC \ar[d]^{p} \\
\Delta^n \times \Delta^1 \ar[r] \ar@{-->}[ur] & \calD. }$$
According to Proposition \ref{goouse}, these lifting problems can be solved provided that
$e$ carries $\{ 0\} \times \Delta^1$ to a $p$-coCartesian edge of $\calC$. This follows immediately from Proposition \ref{universalequiv}.
\end{proof}

\subsection{Bifibrations}\label{bifib}

As we explained in \S \ref{leftfib}, left fibrations $p: X \rightarrow S$ can be thought of as {\em covariant} functors from $S$ into an $\infty$-category of spaces. Similarly, right fibrations $q: Y \rightarrow T$ can be thought of as {\em contravariant} functors from $T$ into an $\infty$-category of spaces. The purpose of this section is to introduce a convenient formalism which encodes covariant and contravariant functoriality simultaneously.

\begin{remark}
The theory of bifibrations will not play an important role in the remainder of the book. In fact, the only result from this section that we will actually use is Corollary \ref{tweezegork}, whose statement makes no mention of bifibrations. A reader who is willing to take Corollary \ref{tweezegork} on faith, or supply an alternative proof, may safely omit the material covered in this section.
\end{remark}

\begin{definition}\label{biv}\index{gen}{bifibration}
Let $S$, $T$, and $X$ be simplicial sets, and $p: X \rightarrow S
\times T$ a map. We shall say that $p$ is a {\it bifibration} if it is an inner fibration having the following
properties:

\begin{itemize}
\item For every $n\geq 1$ and every diagram of solid arrows
$$ \xymatrix{ \Lambda^n_0 \ar@{^{(}->}[d] \ar[r] & X \ar[d] \\
\Delta^n \ar@{-->}[ur] \ar[r]^{f} & S \times T}$$
such that $\pi_T \circ f$ maps $\Delta^{ \{0,1\} } \subseteq \Delta^n$ to a degenerate edge
of $T$, there exists a dotted arrow as indicated, rendering the diagram commutative. 
Here $\pi_T$ denotes the projection $S \times T \rightarrow T$.

\item For every $n\geq 1$ and every diagram of solid arrows
$$ \xymatrix{ \Lambda^n_n \ar@{^{(}->}[d] \ar[r] & X \ar[d] \\
\Delta^n \ar@{-->}[ur] \ar[r]^{f} & S \times T}$$
such that $\pi_S \circ f$ maps $\Delta^{ \{n-1,n\} } \subseteq \Delta^n$ to a degenerate edge
of $T$, there exists a dotted arrow as indicated, rendering the diagram commutative. 
Here $\pi_S$ denotes the projection $S \times T \rightarrow S$.

\end{itemize}
\end{definition}

\begin{remark}
The condition that $p$ be a bifibration is not a condition on $p$ alone, but refers also to a decomposition of the codomain of $p$ as a product $S \times T$. We note also that the definition is not symmetric in $S$ and $T$: instead, $p: X \rightarrow S \times T$ is a bifibration if and only if $p^{op}: X^{op} \rightarrow T^{op} \times S^{op}$ is a bifibration.
\end{remark}

\begin{remark}
Let $p: X \rightarrow S \times T$ be a map of simplicial sets. If $T = \ast$, then $p$ is a bifibration if and only if it is a {\em left} fibration. If $S = \ast$, then $p$ is a bifibration if and only if it is a {\em right} fibration.
\end{remark}

Roughly speaking, we can think of a bifibration $p: X \rightarrow S \times T$ as a bifunctor from $S \times T$ to an $\infty$-category of spaces; the functoriality is covariant in $S$ and contravariant in $T$.

\begin{lemma}\label{gork}
Let $p: X \rightarrow S \times T$ be a bifibration of simplicial sets. Suppose that
$S$ is an $\infty$-category. Then the composition $q=\pi_{T} \circ p$ is a Cartesian fibration of simplicial sets. Furthermore, an edge $e$ of $X$ is $q$-Cartesian if and only if
$\pi_{S}(p(e))$ is an equivalence.
\end{lemma}

\begin{proof}
The map $q$ is an inner fibration, since it is a composition of inner fibrations. Let us say that an edge $e: x \rightarrow y$ of $X$ is {\it quasi-Cartesian} if $\pi_S(p(e))$ is degenerate in $S$. 
Let $y \in X_0$ be any vertex of $X$, and $\overline{e}: \overline{x} \rightarrow q(y)$ an edge of $S$. The pair $(\overline{e}, s_0 q(y))$ is an edge of $S \times T$ whose projection to $T$ is degenerate; consequently, it lifts to a (quasi-Cartesian) edge $e: x \rightarrow y$ in $X$. It is immediate from Definition \ref{biv} that any quasi-Cartesian edge of $X$ is $q$-Cartesian. Thus, $q$ is a Cartesian fibration.

Now suppose that $e$ is a $q$-Cartesian edge of $X$. Then $e$ is equivalent to a quasi-Cartesian edge of $X$; it follows easily that $\pi_{S}(p(e))$ is an equivalence. Conversely, suppose that $e: x \rightarrow y$ is an edge of $X$ and that $\pi_{S}(p(e))$ is an equivalence. We wish to show that $e$ is $q$-Cartesian. Choose a quasi-Cartesian edge
$e': x' \rightarrow y$ with $q(e')=q(e)$. Since $e'$ is $q$-Cartesian, there exists a simplex $\sigma \in X_2$ with $d_0 \sigma = e'$, $d_1 \sigma = e$, and $q(\sigma) = s_0 q(e)$. Let $f = d_2(\sigma)$, so that $\pi_S (p(e')) \circ  \pi_{S} (p(f)) \simeq \pi_S p(e)$ in the $\infty$-category $S$. We note that $f$ lies in the fiber $X_{q(x)}$, which is left-fibered over $S$; since $f$ maps to an equivalence in $S$, it is an equivalence in $X_{q(x)}$. Consequently, $f$ is $q$-Cartesian, so that
$e = e' \circ f$ is $q$-Cartesian as well.
\end{proof}

\begin{proposition}\label{equivbifib}
Let $X \stackrel{p}{\rightarrow} Y \stackrel{q}{\rightarrow} S \times T$ be a diagram of simplicial sets.
Suppose that $q$ and $q \circ p$ are bifibrations, and that $p$ induces a homotopy equivalence $X_{(s,t)} \rightarrow Y_{(s,t)}$ of fibers over each vertex $(s,t)$ of $S \times T$. Then $p$ is a categorical equivalence.
\end{proposition}

\begin{proof}
By means of a standard argument (see the proof of Proposition \ref{babyy}) we may reduce to the case where $S$ and $T$ are simplices; in particular, we may suppose that $S$ and $T$ are $\infty$-categories.
Fix $t \in T_0$, and consider the map of fibers $p_{t}: X_{t} \rightarrow Y_{t}$. Both sides are left-fibered over $S \times \{t\}$, so that $p_{t}$ is a categorical equivalence by (the dual of) Corollary \ref{usefir}. We may then apply Corollary \ref{usefir} again (along with the characterization of Cartesian edges given in Lemma \ref{gork}) to deduce that $p$ is a categorical equivalence.
\end{proof}

\begin{proposition}\label{equivbifib2}
Let $p: X \rightarrow S \times T$ be a bifibration, let $f: S' \rightarrow S$,
$g: T' \rightarrow T$ be categorical equivalences {\em between $\infty$-categories}, and let $X' = X \times_{S \times T} (S' \times T')$.
Then the induced map $X' \rightarrow X$ is a categorical equivalence.
\end{proposition}

\begin{proof}
We will prove the result assuming that $f$ is an isomorphism. A dual argument will establish the result when $g$ is an isomorphism, and applying the result twice we will deduce the desired statement for arbitrary $f$ and $g$.

Given a map $i: A \rightarrow S$, let us say that $i$ is {\it good} if the induced map
$X \times_{S \times T} (A \times T') \rightarrow X \times_{S \times T} (A \times T')$
is a categorical equivalence. We wish to show that the identity map $S \rightarrow S$ is good; it will suffice to show that {\em all} maps $A \rightarrow S$ are good. Using the argument of Proposition \ref{babyy}, we can reduce to showing that every map $\Delta^n \rightarrow S$ is good. In other words, we may assume that $S = \Delta^n$, and in particular that $S$ is an $\infty$-category.
By Lemma \ref{gork}, the projection $X \rightarrow T$ is a Cartesian fibration. The desired result now follows from Corollary \ref{usesec}.
\end{proof}

We next prove an analogue of Lemma \ref{gorban}.

\begin{lemma}\label{gorbaniz}
Let $X \stackrel{p}{\rightarrow} Y \stackrel{q}{\rightarrow} S \times T$ satisfy the hypotheses of Proposition \ref{equivbifib}. Let $A \subseteq B$ be a cofibration of simplicial sets {\em over} $S \times T$. 
Let $f_0: A \rightarrow X$, $g: B \rightarrow Y$ be morphisms in $(\sSet)_{/ S \times T}$ 
and let $h_0: A \times \Delta^1 \rightarrow Y$ be a homotopy (again over $S \times T$) from $g|A$  to $p \circ f_0$.

Then there exists a map $f: B \rightarrow X$ (of simplicial sets over $S \times T$) and a homotopy $h: B \times \Delta^1 \rightarrow T$ (over $S \times T$) from $g$ to $p \circ f$, such that
$f_0 = f|A$ and $h_0 = h| A \times \Delta^1$.
\end{lemma}

\begin{proof}
Working cell-by-cell with the inclusion $A \subseteq B$, we may reduce to the case where $B = \Delta^n$, $A = \bd \Delta^n$. If $n = 0$, we may invoke the fact that $p$ induces a surjection
$\pi_0 X_{(s,t)} \rightarrow \pi_0 Y_{(s,t)}$ on each fiber. Let us assume therefore that $n \geq 1$.
Without loss of generality, we may pull back along the maps $B \rightarrow S$, $B \rightarrow T$, and reduce to the case where $S$ and $T$ are simplices.

We consider the task of constructing $h: \Delta^n \times \Delta^1 \rightarrow T$. We 
now employ the filtration
$$ X(n+1) \subseteq \ldots \subseteq X(0) $$
described in the proof of Proposition \ref{usejoyal}. We note that the value of $h$ on $X(n+1)$ is uniquely prescribed by $h_0$ and $g$. We extend the definition of $h$ to $X(i)$ by descending induction on $i$. We note that $X(i) \simeq X(i+1) \coprod_{ \Lambda^{n+1}_k} \Delta^{n+1}$. For $i > 0$, the existence of the required extension is guaranteed by the assumption that $Y$ is inner-fibered over $S \times T$.

We note that, in view of the assumption that $S$ and $T$ are simplices, any extension of of $h$ over the simplex $\sigma_0$ is automatically a map {\em over} $S \times T$. Since $S$ and $T$ are $\infty$-categories, Proposition \ref{equivbifib} implies that $p$ is a categorical equivalence of $\infty$-categories; the existence of the desired extension of $h$ (and the map $f$ now follows from Lemma \ref{gorban2}.
\end{proof}

\begin{proposition}\label{lemba}
Let $X \stackrel{p}{\rightarrow} Y \stackrel{q}{\rightarrow} S \times T$ satisfy the hypotheses of Proposition \ref{equivbifib}. Suppose that $p$ is a cofibration. Then there exists a retraction $r: Y \rightarrow X$ $($as a map of simplicial sets {\em over} $S \times T${}$)$ such that $r \circ p = \id_X$.
\end{proposition}

\begin{proof}
Apply Lemma \ref{gorbaniz} in the case $A = X$, $B = Y$.
\end{proof}

Let $q: \calM \rightarrow \Delta^1$ be an inner fibration, which we view as a correspondence from
$\calC = q^{-1} \{0\}$ to $\calD = q^{-1} \{1\}$. Evaluation at the endpoints of
$\Delta^1$ induces maps $\bHom_{\Delta^1}(\Delta^1, \calM) \rightarrow \calC$, $\bHom_{\Delta^1}(\Delta^1,\calM) \rightarrow \calD$.

\begin{proposition}\label{tweez}\index{gen}{bifibration!associated to a correspondence}
For every inner fibration $q: \calM \rightarrow \Delta^1$ as above, the map
$p: \bHom_{\Delta^1}(\Delta^1, \calM) \rightarrow \calC \times \calD$ is a bifibration.
\end{proposition}

\begin{proof}
We first show that $p$ is an inner fibration. It suffices to prove
that $q$ has the right-lifting property with respect to
$$ ( \Lambda^n_i \times \Delta^1 ) \coprod_{ \Lambda^n_i \times
\bd \Delta^1 } (\Delta^n \times \bd \Delta^1) \subseteq \Delta^n
\times \Delta^1, $$ for any $0 < i < n$. But this is a smash
product of $\bd \Delta^1 \subseteq \Delta^1$ with the inner anodyne
inclusion $\Lambda^n_i \subseteq \Delta^n$.

To complete the proof that $p$ is a bifibration, we verify that
every $n\geq 1$, $f_0: \Lambda^0_n \rightarrow X$ and
$g: \Delta^n \rightarrow S \times T$ with $g|\Lambda^n_0 = p \circ f_0$, if 
$(\pi_S \circ g) | \Delta^{ \{0,1\} }$ is degenerate, then there exists
$f: \Delta^n \rightarrow X$ with $g= p \circ f$ and $f_0 = f | \Lambda^n_0$. (The dual assertion,
regarding extensions of maps $\Lambda^n_n \rightarrow X$, is verified in the same way.)
The pair $(f_0,g)$ may be regarded as a map
$$h_{0}: (\Delta^n \times \{0,1\}) \coprod_{ \Lambda^n_0 \times \{0,1\} } (\Lambda^n_0 \times \Delta^1)
\rightarrow \calM$$
and our goal is to prove that $h_0$ extends to a map $h: \Delta^n \times \Delta^1 \rightarrow \calM$. 

Let $\{ \sigma_i \}_{0 \leq i \leq n}$ be the maximal-dimensional simplices of $\Delta^n \times \Delta^1$, as in the proof of Proposition \ref{usejoyal}. We set $$K(0) = (\Delta^n \times \{0,1\}) \coprod_{ \Lambda^n_0 \times \{0,1\} } (\Lambda^n_0 \times \Delta^1)$$ and, for $0 \leq i \leq n$, let $K(i+1) = K(i) \bigcup \sigma_i$. We construct maps $h_i: K_i \rightarrow \calM$, with $h_{i} = h_{i+1} | K_{i}$, by induction on $i$.  We note that for $i < n$, $K(i+1) \simeq K(i) \coprod_{ \Lambda^{n+1}_{i+1} } \Delta^{n+1}$, so that the desired extension exists in virtue of the assumption that $\calM$ is an $\infty$-category. If $i = n$, we have instead an isomorphism $\Delta^n \times \Delta^1 = K(n+1) \simeq K(n) \coprod_{ \Lambda^{n+1}_0 } \Delta^{n+1}$. The desired extension of $h_n$ can be found by Proposition \ref{greenlem}, since $h_0 | \Delta^{ \{0,1\} } \times \{0\}$ is an equivalence in $\calC \subseteq \calM$ by assumption.
\end{proof}

\begin{corollary}\label{tweeze}
Let $\calC$ be an $\infty$-category. Evaluation at the endpoints gives a bifibration
$\Fun(\Delta^1,\calC) \rightarrow \calC \times \calC$.
\end{corollary}

\begin{proof}
Apply Proposition \ref{tweez} to the correspondence $\calC \times \Delta^1$.
\end{proof}

\begin{corollary}\label{tweezegork}
Let $f: \calC \rightarrow \calD$ be a functor between $\infty$-categories. The projection
$$ \Fun(\Delta^1, \calD) \times_{\Fun( \{1\}, \calD) } \calC \rightarrow \Fun(\{0\}, \calD)$$
is a Cartesian fibration.
\end{corollary}

\begin{proof}
Combine Corollary \ref{tweeze} with Proposition \ref{gork}.
\end{proof}

\chapter{The $\infty$-Category of $\infty$-Categories}\label{chap4}

\setcounter{theorem}{0}
\setcounter{subsection}{0}

The power of category theory lies in its role as a unifying language for mathematics: nearly every class of mathematical structures (groups, manifolds, algebraic varieties, etcetera) can be organized into a category. This language is somewhat inadequate in situations where the 
structures need to be classified up to some notion of equivalence less rigid than isomorphism. For example, in algebraic topology one wishes to study topological spaces up to homotopy equivalence; in homological algebra one wishes to study chain complexes up to quasi-isomorphism. Both of these examples are most naturally described in terms of higher category theory (for example, the theory of $\infty$-categories developed in this book).

Another source of examples arises in category theory itself. In classical category theory, it is generally regarded as unnatural to ask whether two categories are isomorphic; instead, one asks whether or not they are equivalent. The same phenomenon arises in higher category theory. Throughout this book, we generally regard two $\infty$-categories $\calC$ and $\calD$ as ``the same'' if they are categorically equivalent, even if they are not isomorphic to one another as simplicial sets. In other words, we are not interested in the {\em ordinary} category of $\infty$-categories (a full subcategory of $\sSet$), but in an underlying $\infty$-category which we now define.

\begin{definition}\index{gen}{$\infty$-category!of $\infty$-categories}\index{not}{CatinftyD@$\Cat_{\infty}^{\Delta}$}\index{not}{Catinfty@$\Cat_{\infty}$}
The simplicial category $\Cat_{\infty}^{\Delta}$ is defined as follows:
\begin{itemize}
\item[$(1)$] The objects of $\Cat_{\infty}^{\Delta}$ are (small) $\infty$-categories.

\item[$(2)$] Given $\infty$-categories $\calC$ and $\calD$, we define $\bHom_{\Cat_{\infty}^{\Delta}}(\calC,\calD)$ to be the largest Kan complex contained in the $\infty$-category $\Fun(\calC, \calD)$.
\end{itemize}

We let $\Cat_{\infty}$ denote the simplicial nerve
$\Nerve(\Cat_{\infty}^{\Delta})$. We will refer to $\Cat_{\infty}$ as the {\it $\infty$-category
of $($small$)$ $\infty$-categories}.
\end{definition}

\begin{remark}
By construction, $\Cat_{\infty}$ arises as the nerve of a simplicial category
$\Cat_{\infty}^{\Delta}$, where composition is strictly associative. This is one advantage
of working with $\infty$-categories: the correct notion of functor is encoded by simply considering maps of simplicial sets (rather than homotopy coherent diagrams, say), so there is no difficulty in composing them.
\end{remark}

\begin{remark}
The mapping spaces in $\Cat^{\Delta}_{\infty}$ are Kan complexes, so that
$\Cat_{\infty}$ is an $\infty$-category (Proposition \ref{toothy}) as suggested by the terminology.
\end{remark}

\begin{remark}
By construction, the objects of $\Cat_{\infty}$ are $\infty$-categories, morphisms are given by functors, and $2$-morphisms are given by {\em homotopies} between functors. In other words, $\Cat_{\infty}$ discards all information about noninvertible natural transformations between functors. If necessary, we could retain this information by forming an {\it $\infty$-bicategory} of (small) $\infty$-categories. We do not wish to become involved in any systematic discussion of $\infty$-bicategories, so we will be content to consider only $\Cat_{\infty}$.
\end{remark}

Our goal in this chapter is to study the $\infty$-category $\Cat_{\infty}$. For example, we would like to show that $\Cat_{\infty}$ admits limits and colimits. There are two approaches to proving this assertion. We can attack the problem directly, by giving an explicit construction of the limits and colimits in question: see \S \ref{catlim} and \S \ref{catcolim}. Alternatively, we can try to
realize $\Cat_{\infty}$ as the $\infty$-category underlying a (simplicial) model category $\bfA$, and deduce the existence of limits and colimits in $\Cat_{\infty}$ from the existence of homotopy limits and homotopy colimits in $\bfA$ (Corollary \ref{limitsinmodel}). The objects of $\Cat_{\infty}$ can be identified with the fibrant-cofibrant objects of $\sSet$, with respect to the Joyal model structure. However, we cannot apply Corollary \ref{limitsinmodel} directly, because the Joyal model structure on $\sSet$ is not compatible with its (usual) simplicial structure. We will remedy this difficulty by introducing the category $\mSet$ of {\em marked} simplicial sets. We will explain how to endow $\mSet$ with the structure of a {\em simplicial} model category in such a way that there is an equivalence of simplicial categories $\Cat^{\Delta}_{\infty} \simeq (\mSet)^{\degree}$. This will allow us to identify $\Cat_{\infty}$ with the $\infty$-category underlying $\mSet$, so that Corollary \ref{limitsinmodel} can be invoked.

We will introduce the formalism of marked simplicial sets in \S \ref{twuf}. In particular, we will explain the construction of a model structure not only on $\mSet$ itself, but also for the category
$(\mSet)_{/S}$ of marked simplicial sets {\em over} a given simplicial set $S$. The fibrant objects of $(\mSet)_{/S}$ can be identified with Cartesian fibrations $X \rightarrow S$, which we can think of as contravariant functors from $S$ into $\Cat_{\infty}$. In \S \ref{strsec}, we will justify this intuition by introducing the {\it straightening} and {\it unstraightening} functors which will allow us to pass
back and forth between Cartesian fibrations over $S$ and functors from $S^{op}$ to $\Cat_{\infty}$.
This correspondence has applications to both the study of Cartesian fibrations and to the study of the $\infty$-category $\Cat_{\infty}$; we will survey some of these applications in \S \ref{hugr}.

\begin{remark}\index{not}{ChatCatinfty@$\widehat{\Cat}_{\infty}$}
In the later chapters of this book, it will be necessary to undertake a systematic study of $\infty$-categories which are not small. For this purpose, we introduce the following notational conventions:
$\Cat_{\infty}$ will denote the simplicial nerve of the category of {\em small} $\infty$-categories, while $\widehat{\Cat}_{\infty}$ denotes the the simplicial nerve of the category of $\infty$-categories which are not necessarily small. 
\end{remark}

\section{Marked Simplicial Sets}\label{twuf}

The Joyal model structure on $\sSet$ is a powerful tool in the study of  $\infty$-categories. However, in {\em relative} situations it is somewhat inconvenient. Roughly speaking, a categorical fibration $p: X \rightarrow S$ determines a family of $\infty$-categories $X_{s}$, parametrized by the vertices $s$ of $S$. However, we are generally more interested in those cases where $X_{s}$ can be regarded as a functor of $s$. As we explained in \S \ref{funkymid}, this naturally translates into the assumption that $p$ is a Cartesian (or coCartesian) fibration. According to Proposition \ref{funkyfibcatfib}, every Cartesian fibration is a categorical fibration, but the converse is false. Consequently, it is natural to to try to endow $(\sSet)_{/S}$ with some {\em other} model structure, in which the fibrant objects are precisely the Cartesian fibrations over $S$. 

Unfortunately, this turns out to be an unreasonable demand. In order to have a model category, we need to be able to form fibrant replacements: in other words, we need the ability to enlarge an arbitrary map $p: X \rightarrow S$ into a commutative diagram
$$ \xymatrix{ X \ar[dr]^{p} \ar[rr]^{\phi} & & Y \ar[dl]^{q} \\
& S & }$$
where $q$ is a Cartesian fibration {\em generated by $p$}. A question arises: for which edges $f$ of $X$ should $\phi(f)$ be $q$-Cartesian edge of $Y$? This sort of information is needed for the construction of $Y$; consequently, we need a formalism in which certain edges of $X$ have been distinguished, or {\it marked}.

\setcounter{theorem}{0}

\begin{definition}\index{gen}{marked!simplicial set}\index{gen}{simplicial set!marked}
A  {\it marked simplicial set} is a pair $(X,\calE)$ where $X$ is a simplicial set and
$\calE$ is a set of edges of $X$ which contains every degenerate edge. We will say that an edge of $X$ will be called {\it marked} if it belongs to $\calE$.\index{gen}{marked!edge}

A morphism $f: (X,\calE) \rightarrow (X', \calE')$ of marked simplicial sets is a map $f: X \rightarrow X'$ having the property that $f(\calE) \subseteq \calE'$. The category of marked simplicial sets will be denoted by $\mSet$.\index{not}{sSetm@$\mSet$}
\end{definition}

Every simplicial set $S$ may be regarded as a marked simplicial set, usually in many different ways. The two extreme cases deserve special mention: if $S$ is a simplicial set, we let $S^{\sharp} = (S,S_1)$ denote the marked simplicial set in which {\em every} edge of $S$ has been marked, and $S^{\flat} = (S, s_0(S_0))$ the marked simplicial set in which only the degenerate edges of $S$ have been marked.\index{not}{Xsharp@$X^{\sharp}$}\index{not}{Xflat@$X^{\flat}$}

\begin{notation}
Let $S$ be a simplicial set. We let $(\mSet)_{/S}$ denote the category of marked simplicial
sets equipped with a map to $S$ (which might otherwise be denoted as
$(\mSet)_{ / S^{\sharp} }$).\index{not}{sSetM/S@$(\mSet)_{/S}$}
\end{notation}

Our goal in this section is to study the theory of marked simplicial sets, and in particular to endow 
each $(\mSet)_{/S}$ with the structure of a model category. We will begin in \S \ref{bicat1} by introducing the notion of a {\it marked anodyne} morphism in $\mSet$. In \S \ref{bicat11}, we will establish a basic stability property of the class of marked anodyne maps, which implies the stability of Cartesian fibrations under exponentiation (Proposition \ref{doog}). In \S \ref{markmodel} we will introduce the {\it Cartesian model structure} on $(\mSet)_{/S}$, for every simplicial set $S$. In \S \ref{markprop}, we will study these model categories; in particular, we will see that each $(\mSet)_{/S}$ is a {\it simplicial} model category, whose fibrant objects are precisely the 
Cartesian fibrations $X \rightarrow S$ (with Cartesian edges of $X$ marked). Finally, we will conclude with \S \ref{compmodel}, where we compare the Cartesian model structure on $(\mSet)_{/S}$ with other model structures considered in this book (such as the Joyal and contravariant model structures).

\subsection{Marked Anodyne Morphisms}\label{bicat1}

In this section, we will introduce the class of {\em marked anodyne} morphisms in
$\mSet$. The definition is chosen so that the condition that a map
$\overline{X} \rightarrow \overline{S}$ have the right lifting property with
respect to all marked anodyne morphisms is closely related to the condition
that the underlying map of simplicial sets $X \rightarrow S$ be a Cartesian fibration (we refer the reader to Proposition \ref{dubudu} for a more precise statement). The theory of marked anodyne maps is a technical device which will prove useful when we discuss the Cartesian model structure in
\S \ref{markmodel}: every marked anodyne morphism is a trivial cofibration with respect to the Cartesian model structure, but not conversely. In this respect, the class of marked anodyne morphisms of $\mSet$ is analogous to the class of inner anodyne morphisms of $\sSet$. 

\begin{definition}\label{markanod}\index{gen}{anodyne!marked}\index{gen}{marked!anodyne}
The class of {\em marked anodyne} morphisms in $\mSet$ is the smallest weakly saturated (see \S \ref{liftingprobs}) class of morphisms such that:
\begin{itemize}
\item[$(1)$] For each $0 < i < n$, the inclusion $(\Lambda^n_i)^{\flat} \subseteq (\Delta^n)^{\flat}$ is marked anodyne.

\item[$(2)$] For every $n >0$, the inclusion
$$ ( \Lambda^n_n, \calE \cap (\Lambda^n_n)_{1} ) \subseteq ( \Delta^n, \calE )$$
is marked anodyne, where $\calE$ denotes the set of all degenerate edges of $\Delta^n$, together with the final edge $\Delta^{ \{n-1,n\} }$.

\item[$(3)$] The inclusion
$$ (\Lambda^2_1)^{\sharp} \coprod_{ (\Lambda^2_1)^{\flat} } (\Delta^2)^{\flat} \rightarrow (\Delta^2)^{\sharp}$$ is marked anodyne.

\item[$(4)$] For every Kan complex $K$, the map $K^{\flat} \rightarrow K^{\sharp}$ is marked anodyne.
\end{itemize}
\end{definition}

\begin{remark}
The definition of a marked simplicial set is self-dual. However, Definition
\ref{markanod} is not self-dual: if $A \rightarrow B$ is
marked anodyne, then the opposite morphism $A^{op} \rightarrow B^{op}$ need not be marked-anodyne. This reflects the fact that the theory of Cartesian fibrations is not self-dual.
\end{remark}

\begin{remark}
In part $(4)$ of Definition \ref{markanod}, it suffices to allow $K$ to range over a set of representatives for all isomorphism classes of Kan complexes with only countably many simplices. Consequently, we deduce that the class of marked anodyne morphisms in $\mSet$ is of small generation, so that the small object argument applies (see \S \ref{liftingprobs}). We will refine this observation further: see Corollary \ref{techycor}, below.
\end{remark}

\begin{remark}\label{sillin}
In Definition \ref{markanod}, we are free to replace $(1)$ by
\begin{itemize}
\item[$(1')$] For every inner anodyne map $A \rightarrow B$ of simplicial sets, the
induced map $A^{\flat} \rightarrow B^{\flat}$ is marked anodyne.
\end{itemize}
\end{remark}

\begin{proposition}\label{goldegg}
Consider the following classes of morphisms in $\mSet$:

\begin{itemize}
\item[$(2)$] All inclusions
$$ ( \Lambda^n_n, \calE \cap (\Lambda^n_n)_{1} ) \subseteq ( \Delta^n, \calE ),$$
where $n > 0$ and $\calE$ denotes the set of all degenerate edges of $\Delta^n$, together with the final edge $\Delta^{ \{n-1,n\} }$.

\item[$(2')$] All inclusions
$$ ( (\bd \Delta^n)^{\flat} \times (\Delta^1)^{\sharp})  \coprod_{ (\bd \Delta^n)^{\flat} \times
\{1\}^{\sharp} } ((\Delta^n)^{\flat} \times \{1\}^{\sharp}) \subseteq (\Delta^n)^{\flat} \times (\Delta^1)^{\sharp}.$$ 

\item[$(2'')$] All inclusions
$$ ( A^{\flat} \times (\Delta^1)^{\sharp})  \coprod_{ A^{\flat} \times
\{1\}^{\sharp} } (B^{\flat} \times \{1\}^{\sharp}) \subseteq B^{\flat} \times (\Delta^1)^{\sharp},$$
where $A \subseteq B$ is an inclusion of simplicial sets. 
\end{itemize}

The classes $(2')$ and $(2'')$ generate the same weakly saturated class of morphisms of $\mSet$, which contains the weakly saturated class generated by $(2)$. Conversely, the weakly saturated class of morphisms generated by $(1)$ and $(2)$ from Definition \ref{markanod} contains $(2')$ and $(2'')$.
\end{proposition}

\begin{proof}
To see that each of the morphisms specified in $(2'')$ is contained in the weakly saturated class generated by $(2')$, it suffices to work cell-by-cell with the inclusion $A \subseteq B$. The converse is obvious, since the class of morphisms of type $(2')$ is contained in the class of morphisms of type $(2'')$. To see that the weakly saturated class generated by $(2'')$ contains $(2)$, it suffices to show every morphism in $(2)$ is a retract of a morphism in $(2'')$. For this, we consider
maps
$$\Delta^{n} \stackrel{j}{\rightarrow} \Delta^n \times \Delta^1 \stackrel{r}{\rightarrow} \Delta^{n}.$$ 
Here $j$ is the composition of the identification $\Delta^n \simeq \Delta^n \times \{0\}$ with the inclusion $\Delta^n \times \{0\} \subseteq \Delta^n \times \Delta^1$, and $r$ may be identified with the map of partially ordered sets
$$r(m,i) =
\begin{cases} n & \text{if } m = n-1, i=1 \\
m & \text{otherwise.}  \end{cases}$$

Now we simply observe that $j$ and $r$ exhibit the inclusion
$$ ( \Lambda^n_n, \calE \cap (\Lambda^n_n)_{0} ) \subseteq ( \Delta^n, \calE ),$$
as a retract of
$$ ( (\Lambda^n_n)^{\flat} \times (\Delta^1)^{\sharp})  \coprod_{ (\Lambda_n^n)^{\flat} \times
\{1\}^{\sharp} } ((\Delta^n)^{\flat} \times \{1\}^{\sharp}) \subseteq (\Delta^n)^{\flat} \times (\Delta^1)^{\sharp}.$$ 

To complete the proof, we must show that each of the inclusions
$$ ( (\bd \Delta^n)^{\flat} \times (\Delta^1)^{\sharp})  \coprod_{ (\bd \Delta^n)^{\flat} \times
\{1\}^{\sharp} } ((\Delta^n)^{\flat} \times \{1\}^{\sharp}) \subseteq (\Delta^n)^{\flat} \times (\Delta^1)^{\sharp}$$ 
of type $(2')$ belongs to the weakly saturated class generated by $(1)$ and $(2)$. To see this, we
consider the filtration
$$ Y_{n+1} \subseteq \ldots \subseteq Y_0 = \Delta^n \times \Delta^1$$ which is the
{\em opposite} of the filtration defined in the proof of Proposition \ref{usejoyal}.
We let $\calE_i$ denote the class of all edges of $Y_i$ which are marked in
$(\Delta^n)^{\flat} \times (\Delta^1)^{\sharp}$. 
It will suffice to show that each inclusion $f_i: (Y_{i+1},\calE_{i+1}) \subseteq (Y_{i}, \calE_{i})$ lies in the weakly saturated class generated by $(1)$ and $(2)$. For $i \neq 0$, the map $f_i$ is a pushout of
$(\Lambda^{n+1}_{n+1-i})^{\flat} \subseteq (\Delta^{n+1})^{\flat}$. For $i=0$, $f_i$
is a pushout of 
$$ ( \Lambda^{n+1}_{n+1}, \calE \cap (\Lambda^{n+1}_{n+1})_{1} ) \subseteq ( \Delta^{n+1}, \calE ),$$
where and $\calE$ denotes the set of all degenerate edges of $\Delta^{n+1}$, together with $\Delta^{ \{n,n+1\}}$.
\end{proof}

We now characterize the class of marked-anodyne maps:

\begin{proposition}\label{dubudu}
A map $p: X \rightarrow S$ in $\mSet$ has the right lifting property with respect to all marked anodyne maps if and only if the following conditions are satisfied:
\begin{itemize}
\item[$(A)$] The map $p$ is an inner fibration of simplicial sets.
\item[$(B)$] An edge $e$ of $X$ is marked if and only if $p(e)$ is marked and $e$ is $p$-Cartesian.
\item[$(C)$] For every object $y$ of $X$ and every marked edge $\overline{e}: \overline{x} \rightarrow p(y)$ in $S$, there exists a marked edge $e: x \rightarrow y$ of $X$ with $p(e) = \overline{e}$. 
\end{itemize}
\end{proposition}

\begin{proof}
We first prove the ``only if'' direction. Suppose that $p$ has the right lifting property with respect to all marked anodyne maps. By considering maps of the form $(1)$ from Definition \ref{markanod}, we deduce that $(A)$ holds. Considering $(2)$ in the case $n=0$, we deduce that $(C)$ holds. Considering $(2)$ for $n > 0$, we deduce that every marked edge of $X$ is $p$-Cartesian.
For the converse, let us suppose that $e: x \rightarrow y$ is a $p$-Cartesian edge of $X$ and that $p(e)$ is marked in $S$. Invoking $(C)$, we deduce that there exists a marked edge $e': x' \rightarrow y$ with $p(e)=p(e')$. Since $e'$ is Cartesian,
we can find a $2$-simplex $\sigma$ of $X$ with $d_0(\sigma) = e'$, $d_1(\sigma)=e$, and 
$p(\sigma)= s_1 p(e)$. Then
$d_2(\sigma)$ an equivalence between $x$ and $x'$ in the $\infty$-category $X_{p(x)}$. Let $K$ denote the largest Kan complex contained in $X_{p(x)}$. Since $p$ has the right lifting property with respect to
$K^{\flat} \rightarrow K^{\sharp}$, we deduce that every edge of $K$ is marked; in particular, $d_2(\sigma)$ is marked. Since $p$ has the right lifting property with respect to the morphism described in $(3)$ of Definition \ref{markanod}, we deduce that $d_1(\sigma)=e$ is marked.

Now suppose that $p$ satisfied the hypotheses of the proposition. We must show that $p$ has the right lifting property with respect to the classes of morphisms $(1)$, $(2)$, $(3)$, and $(4)$ of Definition \ref{markanod}. For $(1)$, this follows from the assumption that $p$ is an inner fibration. 
For $(2)$, this follows from $(C)$ and from the assumption that every marked edge is $p$-Cartesian. For $(3)$, we are free to replace $S$ by $(\Delta^2)^{\sharp}$; then $p$ is a Cartesian fibration over an $\infty$-category $S$ and we may apply Proposition \ref{protohermes} to deduce that the class of $p$-Cartesian edges is stable under composition.

Finally, for $(4)$, we may replace $S$ by $K^{\sharp}$; 
then $S$ is a Kan complex and $p$ is a Cartesian fibration, so the $p$-Cartesian edges of $X$
are precisely the equivalences in $X$. Since $K$ is a Kan complex, any
map $K \rightarrow X$ carries the edges of $K$ to equivalences in $X$.
\end{proof}

By Quillen's small object argument, we 
deduce that a map $j: A \rightarrow B$ in $\mSet$ is marked anodyne if and only if it has the left lifting property with respect to all morphisms $p: X \rightarrow S$ satisfying the hypotheses of Proposition \ref{dubudu}. From this, we deduce:

\begin{corollary}\label{hermes}
The inclusion
$$ i: (\Lambda^2_2)^{\sharp} \coprod_{ (\Lambda^2_2)^{\flat} } (\Delta^2)^{\flat} \hookrightarrow (\Delta^2)^{\sharp}$$ is marked anodyne.
\end{corollary}

\begin{proof}
It will suffice to show that $i$ has the left lifting property with respect to any of the morphisms $p: X \rightarrow S$ described in Proposition \ref{dubudu}. Without loss of generality, we may replace $S$ by $(\Delta^2)^{\sharp}$; we now apply Proposition \ref{protohermes}.
\end{proof}

The following somewhat technical corollary will be needed in \S \ref{markmodel}:

\begin{corollary}\label{techycor}
In Definition \ref{markanod}, we can replace the class of morphisms $(4)$ by
\begin{itemize}
\item[$(4')$] the map $j: A^{\flat} \rightarrow (A, s_0 A_0 \bigcup \{f\} )$, where
$A$ is the quotient of $\Delta^3$ which co-represents the functor
$$ \Hom_{\sSet}(A,X) = \{ \sigma \in X_3, e \in X_1: d_1 \sigma = s_0 e, 
d_2 \sigma = s_1 e \}$$ and $f \in A_1$ is the image of
$\Delta^{ \{0,1\}} \subseteq \Delta^3$ in $A$.
\end{itemize}
\end{corollary}

\begin{proof}
We first show that for every Kan complex $K$, the map $i: K^{\flat} \rightarrow K^{\sharp}$
lies in the weakly saturated class of morphisms generated by $(4')$. We note that $i$ can be obtained 
as an iterated pushout of morphisms having the form
$K^{\flat} \rightarrow (K, s_0 K_0 \bigcup \{e\})$, where $e$ is an edge of $K$. It therefore suffices to show that there exists a map $p: A \rightarrow K$ such that $p(f) = e$. In other words, we must prove that there exists a $3$-simplex $\sigma: \Delta^3 \rightarrow K$ with $d_1 \sigma = s_0 e$ and $d_2 \sigma = s_1 e$. This follows immediately from the Kan extension condition.

To complete the proof, it will suffice to show that the map $j$ is marked anodyne.
To do so, it suffices to prove that for any diagram
$$ \xymatrix{ A^{\flat} \ar@{^{(}->}[d] \ar[r] & X \ar[d]^{p} \\
(A, s_0 A_0 \cup \{f\}) \ar[r] \ar@{-->}[ur] & S }$$
for which $p$ satisfies the conditions of Proposition \ref{dubudu}, there
exists a dotted arrow as indicated, rendering the diagram commutative. This is
simply a reformulation of Proposition \ref{sworkk}.
\end{proof}

\begin{definition}\label{conf1}\index{not}{Xnatural@$X^{\natural}$}
Let $p: X \rightarrow S$ be a Cartesian fibration of simplicial sets. We let
$X^{\natural}$ denote the marked simplicial set $(X, \calE)$, where $\calE$ is the set of
$p$-Cartesian edges of $X$.
\end{definition}

\begin{remark}\label{abus}
Our notation is slightly abusive, since $X^{\natural}$ depends not only on $X$ but also
on the map $X \rightarrow S$. 
\end{remark}

\begin{remark}\label{abuss}
According to Proposition \ref{dubudu}, a map
$(Y, \calE) \rightarrow S^{\sharp}$ has the right lifting property with respect to all marked anodyne maps if and only if the underlying map $Y \rightarrow S$ is a Cartesian fibration and
$(Y, \calE) = Y^{\natural}$. 
\end{remark}

We conclude this section with the following easy result, which will be needed later:

\begin{proposition}\label{eggwhite}
Let $p: X \rightarrow S$ be an inner fibration of simplicial sets, and let $f: A \rightarrow B$
be a marked anodyne morphism in $\mSet$, let
$q: B \rightarrow X$ be map of simplicial sets which carries each marked edge of $B$
to a $p$-Cartesian edge of $X$, and $q_0 = q \circ f$. 
Then the induced map
$$ X_{/q} \rightarrow X_{/q_0} \times_{S_{/pq_0}} S_{/pq}$$
is a trivial fibration of simplicial sets. 
\end{proposition}

\begin{proof}
It is easy to see that the class of all morphisms $f$ of $\mSet$ which satisfy the desired conclusion is weakly saturated. It therefore suffices to prove that this class contains collection of generators for the weakly saturated class of marked anodyne morphisms.
If $f$ induces a left anodyne map on the underlying simplicial sets, then the desired result is automatic. It therefore suffices to consider the case where $f$ is the inclusion
$$ ( \Lambda^n_n, \calE \cap (\Lambda^n_n)_{1} ) \subseteq ( \Delta^n, \calE )$$
as described in $(2)$ of Definition \ref{markanod}. In this case, a lifting problem
$$ \xymatrix{ \bd \Delta^m \ar[r] \ar@{^{(}->}[d] & X_{/q} \ar[d] \\
\Delta^m \ar[r] \ar@{-->}[ur] & X_{/q_0} \times_{S_{/pq_0}} S_{/pq} }$$
can be reformulated as an equivalent lifting problem
$$ \xymatrix{ \Lambda^{n+m+1}_{n+m+1} \ar[r]^{\sigma_0} \ar@{^{(}->}[d] & X \ar[d]^{p} \\
\Delta^{n+m+1} \ar[r] \ar@{-->}[ur] & S. }$$
This lifting problem admits a solution, since the hypothsis on $q$ guarantees
that $\sigma_0$ carries $\Delta^{ \{n+m, n+m+1\} }$ to a $p$-Cartesian edge of $X$.
\end{proof}

\subsection{Stability Properties of Marked Anodyne Morphisms}\label{bicat11}

Our main goal in this section is to prove the following stability result:

\begin{proposition}\label{doog}\index{gen}{Cartesian fibration!and functor categories}
Let $p: X \rightarrow S$ be a Cartesian fibration of
simplicial sets, and let $K$ be an arbitrary simplicial set. Then:
\begin{itemize}
\item[$(1)$] The induced map $p^K: X^{K} \rightarrow S^{K}$ is a
Cartesian fibration.

\item[$(2)$] An edge $\Delta^1 \rightarrow X^K$ is $p^K$-Cartesian if
and only if, for every vertex $k$ of $K$, the induced edge $\Delta^1
\rightarrow X$ is $p$-Cartesian.
\end{itemize}
\end{proposition}

We could easily have given an ad-hoc proof of this result in \S \ref{slib}. However, we have opted instead to give a proof using the language of marked simplicial sets. 

\begin{definition}
A morphism $(X, \calE) \rightarrow (X', \calE')$ in $\mSet$ is a {\it cofibration} if the underlying map $X \rightarrow X'$ of simplicial sets is a cofibration.
\end{definition}

The main ingredient we will need to prove Proposition \ref{doog} is the following:

\begin{proposition}\label{markanodprod}
The class of marked anodyne maps in $\mSet$ is stable under smash products with arbitrary cofibrations. In other words, if $f: X \rightarrow X'$ is marked anodyne, and 
$g: Y \rightarrow Y'$ is a cofibration, then the induced map
$$ (X \times Y') \coprod_{ X \times Y} (X' \times Y) \rightarrow X' \times Y'$$ is marked anodyne.
\end{proposition}

\begin{proof}
The argument is tedious, but straightforward. 
Without loss of generality, we may suppose that $f$ belongs either to the class $(2')$ of
Proposition \ref{goldegg}, or one of the classes specified in $(1)$, $(3)$, or $(4)$ of Definition \ref{markanod}.
The class of cofibrations is generated by the inclusions $(\bd \Delta^n)^{\flat} \subseteq (\Delta^n)^{\flat}$ and $(\Delta^1)^{\flat} \subseteq (\Delta^1)^{\sharp}$; thus we may suppose that $g: Y \rightarrow Y'$ is one of these maps. There are eight cases to consider:

\begin{itemize}
\item[(A1)] Let $f$ be the inclusion $(\Lambda^n_i)^{\flat} \subseteq (\Delta^n)^{\flat}$ and $g$ the inclusion $(\bd \Delta^n)^{\flat} \rightarrow (\Delta^n)^{\flat}$, where $0 < i < n$. Since the class of inner anodyne maps between simplicial sets is stable under smash products with inclusions, the smash product of $f$ and $g$ is marked-anodyne (see Remark \ref{sillin}). 

\item[(A2)] Let $f$ denote the inclusion $(\Lambda^n_i)^{\flat} \rightarrow (\Delta^n)^{\flat}$, and $g$ the map $(\Delta^1)^{\flat} \rightarrow (\Delta^1)^{\sharp}$, where $0 < i < n$. Then the smash product of $f$ and $g$ is an isomorphism (since $\Lambda^n_i$ contains all vertices of $\Delta^n$).

\item[(B1)] Let $f$ be the inclusion $$(\{ 1\}^{\sharp} \times (\Delta^n)^{\flat} ) \coprod_{ \{1\}^{\sharp} \times (\bd \Delta^n)^{\flat} } ((\Delta^1)^{\sharp} \times (\bd \Delta^n)^{\flat} ) \subseteq (\Delta^1)^{\sharp} \times (\Delta^n)^{\flat},$$ and let $g$ be the inclusion $(\bd \Delta^n)^{\flat} \rightarrow (\Delta^n)^{\flat}$. Then the smash product of $f$ and $g$ belongs to the class $(2'')$ of Proposition \ref{goldegg}.

\item[(B2)] Let $f$ be the inclusion $$(\{ 1\}^{\sharp} \times (\Delta^n)^{\flat} ) \coprod_{ \{1\}^{\sharp} \times (\bd \Delta^n)^{\flat} } ((\Delta^1)^{\sharp} \times (\bd \Delta^n)^{\flat} ) \subseteq (\Delta^1)^{\sharp} \times (\Delta^n)^{\flat},$$ and let $g$ denote the map
$(\Delta^1)^{\flat} \rightarrow (\Delta^1)^{\sharp}$. If $n > 0$, then the smash product of $f$ and $g$ is an isomorphism. If $n = 0$, then the smash product may be identified with
the map $ (\Delta^1 \times \Delta^1, \calE) \rightarrow (\Delta^1 \times \Delta^1)^{\sharp}$, where $\calE$ consists of all degenerate edges together with $\{0\} \times \Delta^1$, $\{1\} \times \Delta^1$, and $\Delta^1 \times \{1\}$. This map may be obtained as a composition of
two marked anodyne maps: the first is of type $(3)$ in Definition \ref{markanod} (adjoining the ``diagonal'' edge to $\calE$) and the second is the map described in Corollary \ref{hermes}
(adjoining the edge $\Delta^1 \times \{0\}$ to $\calE$).

\item[(C1)] Let $f$ be the inclusion $$(\Lambda^2_1)^{\sharp} \coprod_{ (\Lambda^2_1)^{\flat} } (\Delta^2)^{\flat} \rightarrow (\Delta^2)^{\sharp},$$ and let $g$ the inclusion
$(\bd \Delta^n)^{\flat} \subseteq (\Delta^n)^{\flat}$. Then the smash product of $f$ and $g$ is an isomorphism for $n > 0$, and isomorphic to $f$ for $n=0$.

\item[(C2)] Let $f$ be the inclusion $$(\Lambda^2_1)^{\sharp} \coprod_{ (\Lambda^2_1)^{\flat} } (\Delta^2)^{\flat} \rightarrow (\Delta^2)^{\sharp},$$ and let $g$ be the canonical map $(\Delta^1)^{\flat} \rightarrow (\Delta^1)^{\sharp}.$ Then the smash product of $f$ and $g$ is a pushout of the map $f$.

\item[(D1)] Let $f$ be the map $K^{\flat} \rightarrow K^{\sharp}$, where $K$ is a Kan complex, and let $g$ the inclusion $(\bd \Delta^n)^{\flat} \subseteq (\Delta^n)^{\flat}$. Then the smash product of $f$ and $g$ is an isomorphism for $n > 0$, and isomorphic to $f$ for $n=0$.

\item[(D2)] Let $f$ be the map $K^{\flat} \rightarrow K^{\sharp}$, where $K$ is a Kan complex, and let $g$ be the map $(\Delta^1)^{\flat} \rightarrow (\Delta^1)^{\sharp}$. The smash product of $f$ and $g$ can be identified with the inclusion
$$(K \times \Delta^1, \calE) \subseteq (K \times \Delta^1)^{\sharp},$$ where $\calE$
denotes the class of all edges $e=(e',e'')$ of $K \times \Delta^1$ for which either $e': \Delta^1 \rightarrow K$ or $e'': \Delta^1 \rightarrow \Delta^1$ is degenerate. This inclusion can be obtained as a transfinite composition of pushouts of the map $$(\Lambda^2_1)^{\sharp} \coprod_{ (\Lambda^2_1)^{\flat} } (\Delta^2)^{\flat} \rightarrow (\Delta^2)^{\sharp}.$$
\end{itemize}
\end{proof}

We now return to our main objective:

\begin{proof}[Proof of Proposition \ref{doog}]
Since $p$ is a Cartesian fibration, it induces a map
$X^{\natural} \rightarrow S^{\sharp}$ which has the right lifting property with respect to all marked anodyne maps. By Proposition \ref{markanodprod}, the induced map
$$ (X^{\natural})^{K^{\flat}} \rightarrow (S^{\sharp})^{K^{\flat}}=(S^K)^{\sharp}$$
has the right lifting property with respect to all marked anodyne morphisms. The desired result now follows from Remark \ref{abuss}.
\end{proof}

\subsection{The Cartesian Model Structure}\label{markmodel}

Let $S$ be a simplicial set. Our goal in this section is to introduce
the {\it Cartesian model structure} on the category $(\mSet)_{/S}$ of marked simplicial sets over $S$. We will eventually show that the fibrant objects of
$(\mSet)_{/S}$ correspond precisely to Cartesian fibrations $X \rightarrow S$, and that they encode (contravariant) functors from $S$ into the $\infty$-category
$\Cat_{\infty}$.

The category $\mSet$ is {\it Cartesian-closed}; that is, for any two objects
$X,Y \in \mSet$, there exists an internal mapping object $Y^X$ equipped with an
``evaluation map'' $Y^X \times X \rightarrow Y$ which induces bijections
$$\Hom_{\mSet}(Z,Y^X) \rightarrow \Hom_{\mSet}(Z \times X, Y)$$
for every $Z \in \mSet$. We let $\bHom^{\flat}(X,Y)$ denote the underlying simplicial
set of $Y^X$, and $\bHom^{\sharp}(X,Y) \subseteq \bHom^{\flat}(X,Y)$ the simplicial subset
consisting of all simplices $\sigma \subseteq \bHom^{\flat}(X,Y)$ such that every edge of
$\sigma$ is a marked edge of $Y^X$. Equivalently, we may describe these simplicial sets by the mapping properties
$$ \Hom_{\sSet}(K, \bHom^{\flat}(X,Y)) \simeq \Hom_{\mSet}( K^{\flat} \times X,Y) $$
$$ \Hom_{\sSet}(K, \bHom^{\sharp}(X,Y)) \simeq \Hom_{\mSet}(K^{\sharp} \times X,Y).$$\index{not}{Homflat@$\bHom^{\flat}(X,Y)$}\index{not}{Homsharp@$\bHom^{\sharp}(X,Y)$}

If $X$ and $Y$ are objects of $(\mSet)_{/S}$, then we let
$\bHom_{S}^{\sharp}(X,Y)$ and $\bHom_{S}^{\flat}(X,Y)$ denote the simplicial subsets
of $\bHom^{\sharp}(X,Y)$ and $\bHom^{\flat}(X,Y)$ classifying those maps which are compatible with the projections to $S$.\index{not}{MapSflat@$\bHom_{S}^{\flat}(X,Y)$}\index{not}{MapSSharp@$\bHom_{S}^{\sharp}(X,Y)$}

\begin{remark}
If $X \in (\mSet)_{/S}$ and $p: Y \rightarrow S$ is a Cartesian fibration, then
$\bHom^{\flat}_{S}(X,Y^{\natural})$ is an $\infty$-category, and $\bHom^{\sharp}_{S}(X,Y^{\natural})$ is the largest Kan complex contained in $\bHom^{\flat}_S(X,Y^{\natural})$.
\end{remark}

\begin{lemma}\label{hoopsp}
Let $f: \calC \rightarrow \calD$ be a functor between $\infty$-categories.
The following are equivalent:
\begin{itemize}
\item[$(1)$] The functor $f$ is a categorical equivalence.
\item[$(2)$] For every simplicial set $K$, the induced map
$\Fun(K,\calC) \rightarrow \Fun(K,\calD)$ is a categorical equivalence.
\item[$(3)$] For every simplicial set $K$, the functor $\Fun(K,\calC) \rightarrow \Fun(K,\calD)$ induces
a homotopy equivalence between the largest Kan complex contained in $\Fun(K,\calC)$ and the largest Kan complex contained in $\Fun(K,\calD)$.
\end{itemize}
\end{lemma}

\begin{proof}
The implications $(1) \Rightarrow (2) \Rightarrow (3)$ are obvious. Suppose that $(3)$ is satisfied.
Let $K = \calD$. According to $(3)$, there exists an object $x$ of $\Fun(K,\calC)$ whose image
in $\Fun(K,\calD)$ is equivalent to the identity map $K \rightarrow \calD$. We may identify
$x$ with a functor $g: \calD \rightarrow \calC$ having the property that $f \circ g$ is homotopic
to the identity $\id_{\calD}$. It follows that $g$ also has the property asserted by $(3)$, so the same argument shows that there is a functor $f': \calC \rightarrow \calD$ such that
$g \circ f'$ is homotopic to $\id_{\calC}$. It follows that $f \circ g \circ f'$ is homotopic to both
$f$ and $f'$, so that $f$ is homotopic to $f'$. Thus $g$ is a homotopy inverse to $f$, which proves that $f$ is an equivalence.
\end{proof}

\begin{proposition}\label{markdefeq}
Let $S$ be a simplicial set, and let $p: X \rightarrow Y$ be a morphism in $(\mSet)_{/S}$. The following are equivalent:
\begin{itemize}
\item[$(1)$] For every Cartesian fibration $Z \rightarrow S$, the induced map 
$$ \bHom^{\flat}_{S}(Y, Z^{\natural}) \rightarrow \bHom^{\flat}_{S}(X, Z^{\natural})$$ is an
equivalence of $\infty$-categories.
\item[$(2)$] For every Cartesian fibration $Z \rightarrow S$, the induced map
$$ \bHom^{\sharp}_{S}(Y, Z^{\natural}) \rightarrow \bHom^{\sharp}_{S}(X, Z^{\natural})$$ is a homotopy equivalence of Kan complexes.
\end{itemize}
\end{proposition}

\begin{proof}
Since $\bHom^{\sharp}_{S}(M,Z^{\natural})$ is the largest Kan complex contained in 
$\bHom^{\flat}_{S}(M, Z^{\natural})$, it is clear that $(1)$ implies $(2)$. Suppose that $(2)$ is satisfied, and let $Z \rightarrow S$ be a Cartesian fibration. We wish to show that
$$ \bHom^{\flat}_{S}(Y, Z^{\natural}) \rightarrow \bHom^{\flat}_{S}(X, Z^{\natural})$$
is an equivalence of $\infty$-categories. According to Lemma \ref{hoopsp}, it suffices to show that
$$ \bHom^{\flat}_{S}(Y, Z^{\natural})^K \rightarrow \bHom^{\flat}_{S}(X, Z^{\natural})^K$$ 
induces a homotopy equivalence on the maximal Kan complexes contained in each side.
Let $Z(K) = Z^K \times_{S^K} S$. Proposition \ref{doog} implies that $Z(K) \rightarrow S$
is a Cartesian fibration, and that there is a natural identification
$$ \bHom^{\flat}_{S}(M, Z(K)^{\natural}) \simeq \bHom^{\flat}_{S}(M, Z(K)^{\natural}).$$
The largest Kan complex contained in the right hand side is $\bHom^{\sharp}_{S}(M,Z(K)^{\natural})$. On the other hand, the natural map
$$ \bHom^{\sharp}_{S}(Y,Z(K)^{\natural}) \rightarrow \bHom^{\sharp}_{S}(X, Z(K)^{\natural})$$ is homotopy equivalence by assumption $(2)$.
\end{proof}

We will say that a map $X \rightarrow Y$ in $(\mSet)_{/S}$ is a {\it Cartesian equivalence} if it satisfies the equivalent conditions of Proposition \ref{markdefeq}.\index{gen}{equivalence!Cartesian}\index{gen}{Cartesian!equivalence}

\begin{remark}
Let $f: X \rightarrow Y$ be a morphism in $(\mSet)_{/S}$ which is {\em marked anodyne} when regarded as a map of marked simplicial sets. Since the smash product of $f$ with any inclusion $A^{\flat} \subseteq B^{\flat}$ is also marked-anodyne, we deduce that the map
$$ \phi: \bHom^{\flat}_{S}(Y, Z^{\natural}) \rightarrow \bHom^{\flat}_{S}(X, Z^{\natural})$$ is a trivial fibration for {\em every} Cartesian fibration $Z \rightarrow S$. Consequently, $f$ is a Cartesian equivalence.
\end{remark}

Let $S$ be a simplicial set and let $X,Y \in (\mSet)_{/S}$. We will say a pair of morphisms
$f,g: X \rightarrow Y$ are {\it strongly homotopic} if there exists a contractible Kan complex
$K$ and a map $K \rightarrow \bHom^{\flat}_{S}(X,Y)$, whose image contains both of the vertices $f$ and $g$. If $Y = Z^{\natural}$, where $Z \rightarrow S$ is a Cartesian fibration, then this simply means that $f$ and $g$ are equivalent when viewed as objects of the $\infty$-category $\bHom^{\flat}_{ S}(X,Y)$.

\begin{proposition}\label{crispy}
Let $X \stackrel{p}{\rightarrow} Y \stackrel{q}{\rightarrow} S$ be a diagram of simplicial sets, where both $q$ and $q \circ p$ are Cartesian fibrations. The following assertions are equivalent:

\begin{itemize}
\item[$(1)$] The map $p$ induces a Cartesian equivalence $X^{\natural} \rightarrow Y^{\natural}$ in
$(\mSet)_{/S}$.

\item[$(2)$] There exists a map $r: Y \rightarrow X$ which is a strong homotopy inverse to $p$, in the sense that $p \circ r$ and $r \circ p$ are both strongly homotopic to the identity.

\item[$(3)$] The map $p$ induces a categorical equivalence $X_{s} \rightarrow Y_{s}$, for each vertex $s$ of $S$.
\end{itemize}
\end{proposition}

\begin{proof}
The equivalence between $(1)$ and $(2)$ is easy, as is the assertion that $(2)$ implies $(3)$. It therefore suffices to show that $(3)$ implies $(2)$. We will construct $r$ and a homotopy from $r \circ p$ to the identity. It then follows that the map $r$ satisfies $(3)$, so the same argument will show that $r$ has a right homotopy inverse; by general nonsense this right homotopy inverse is automatically homotopic to $p$ and the proof will be complete.

Choose a transfinite sequence of simplicial subsets $S(\alpha) \subseteq S$, where each
$S(\alpha)$ is obtained from $\bigcup_{ \beta < \alpha} S(\beta)$ by adjoining a single
nondegenerate simplex (if such a simplex exists).
We construct
$r_{\alpha}: Y \times_{S} S(\alpha) \rightarrow X$ and an equivalence $h_{\alpha}: (X \times_{S} S(\alpha) ) \times \Delta^1\rightarrow X \times_{S} S(\alpha)$ from $r_{\alpha} \circ p$ to the identity, by induction on $\alpha$. By this device we may reduce to the case where $S = \Delta^n$, and the maps
$$r^0: Y' \rightarrow X$$
$$h^0: X' \times \Delta^1 \rightarrow X$$ 
are already specified, where $Y' = Y \times_{\Delta^n} \bd \Delta^n \subseteq Y$ and
$X' = X \times_{\Delta^n} \bd \Delta^n \subseteq X$. We may regard $r'$ and $h'$ together as defining a map
$\psi_0: Z'  \rightarrow X$, where
$$Z' = Y' \coprod_{ X' \times \{0\} } (X' \times \Delta^1) \coprod_{ X' \times \{1\} } X.$$
Let $Z = Y \coprod_{ X \times \{0\} } X \times \Delta^1$; then our goal is to solve the lifting problem depicted in the following diagram:
$$ \xymatrix{ Z' \ar[r]^{\psi_0} \ar@{^{(}->}[d] & X \ar[d] \\
Z \ar[r] \ar@{-->}[ur]^{\psi} & \Delta^n }$$
in such a way that $\psi$ carries $\{x\} \times \Delta^1$ to an equivalence in $X$, for
every object $x$ of $X$. We note that this last condition is vacuous for $n > 0$.

If $n=0$, the problem amounts to constructing a map $Y \rightarrow X$ which is homotopy inverse to $p$: this is possible in view of the assumption that $p$ is a categorical equivalence.
For $n > 0$, we note that any map $\phi: Z \rightarrow X$ extending $\phi_0$ is automatically compatible with the projection to $S$ (since $S$ is a simplex and $Z'$ contains all vertices of $Z$). 
Since the inclusion $Z' \subseteq Z$ is a cofibration between cofibrant objects in the model category $\sSet$ (with the Joyal model structure), and $X$ is a $\infty$-category (since $q$ is an inner fibration and $\Delta^n$ is a $\infty$-category), Proposition \ref{princex} asserts that it is sufficient to show that the extension $\phi$ exists up to homotopy. Since Corollary \ref{usefir} implies that $p$ is an equivalence, we are free to replace the inclusion $Z' \subseteq Z$ with the weakly equivalent inclusion
$$ (X \times \{1\}) \coprod_{ X \times_{ \Delta^n} \bd\Delta^n \times \Delta^1 }
(X \times_{\Delta^n} \bd \Delta^n \times \{1\}) \subseteq X \times \Delta^1.$$
Since $\phi_0$ carries $\{x\} \times \Delta^1$ to a $(q \circ p)$-Cartesian edge of $X$, for
every vertex $x$ of $X$, the existence of $\phi$ follows from Proposition \ref{goldegg}.
\end{proof}

\begin{lemma}\label{insdod}
Let $S$ be a simplicial set, let $i: X \rightarrow Y$ be a cofibration in $(\mSet)_{/S}$, and let $Z \rightarrow S$ be a Cartesian fibration. Then the associated map $p: \bHom_{S}^{\sharp}(Y,Z^{\natural}) \rightarrow \bHom_{S}^{\sharp}(X, Z^{\natural})$ is a Kan fibration.
\end{lemma}

\begin{proof}
Let $A \subseteq B$ be an anodyne inclusion of simplicial sets. We must show that $p$ has the right lifting property with respect to $p$. Equivalently, we must show that $Z^{\natural} \rightarrow S$ has the right lifting property with respect to the inclusion
$$ (B^{\sharp} \times X) \coprod_{A^{\sharp} \times X} (A^{\sharp} \times Y) \subseteq
B^{\sharp} \times Y.$$
This follows from Proposition \ref{markanodprod}, since the inclusion $A^{\sharp} \subseteq B^{\sharp}$ is marked anodyne.
\end{proof}

\begin{proposition}\label{markmodell}\index{gen}{model category!Cartesian}\index{gen}{Cartesian model structure}
Let $S$ be a simplicial set. There exists a left proper, combinatorial model structure on $(\mSet)_{/S}$, which may be described as follows:

\begin{itemize}
\item[$(C)$] The cofibrations in $(\mSet)_{/S}$ are those morphisms $p: X \rightarrow Y$ in $(\mSet)_{/S}$ which are cofibrations when regarded as morphisms of simplicial sets.

\item[$(W)$] The weak equivalences in $(\mSet)_{/S}$ are the Cartesian equivalences.

\item[$(F)$] The fibrations in $(\mSet)_{/S}$ are those maps which have the right lifting property with respect to every map which is simultaneously a cofibration and a Cartesian equivalence.
\end{itemize}
\end{proposition}

\begin{proof}
It suffices to show that the hypotheses of Proposition \ref{goot} are satisfied by the class $(C)$ of cofibrations and the class $(W)$.
\begin{itemize}
\item[(1)] The class $(W)$ of Cartesian equivalences is perfect, in the sense of Definition \ref{perfequiv}. To prove this, we first observe that the class of marked anodyne maps
is generated by the classes $(1)$, $(2)$, $(3)$ of Definition \ref{markanod} and
$(4')$ of Corollary \ref{techycor}. By Proposition \ref{quillobj}, there exists a functor
$T$ from $(\mSet)_{/S}$ to itself and a (functorial) factorization
$$X \stackrel{i_X}{\rightarrow} T(X) \stackrel{j_X}{\rightarrow} S^{\sharp}$$
where $i_X$ is marked anodyne (and therefore a Cartesian equivalence) and
$j_{X}$ has the right lifting property with respect to all marked anodyne maps, and therefore
corresponds to a Cartesian fibration over $S$. Moreover, the functor $T$ commutes
with filtered colimits. According to Proposition \ref{crispy}, a map $X \rightarrow Y$
in $(\mSet)_{/S}$ is a Cartesian equivalence if and only if, for each vertex $s \in S$, 
the induced map $T(X)_{s} \rightarrow T(Y)_{s}$ is a categorical equivalence.
It follows from Corollary \ref{perfpull} that $(W)$ is a perfect class of morphisms.

\item[(2)] The class of weak equivalences is stable under pushouts by cofibrations.
Suppose given a pushout diagram
$$ \xymatrix{ X \ar[r]^{p} \ar[d]^{i} & Y \ar[d] \\
X' \ar[r]^{p'} & Y' }$$ where $i$ is a cofibration and $p$ is a Cartesian equivalence. We wish to show that $p'$ is also a Cartesian equivalence. In other words, we must show that for any Cartesian fibration $Z \rightarrow S$, the associated map $\bHom^{\sharp}_{S}(Y', Z^{\natural}) \rightarrow
\bHom^{\sharp}_{S}(X', Z^{\natural})$ is a homotopy equivalence. Consider the pullback diagram
$$ \xymatrix{ \bHom^{\sharp}_{S}(Y', Z^{\natural}) \ar[r] \ar[d] & \bHom^{\sharp}_{S}(X',Z^{\natural}) \ar[d] \\
\bHom^{\sharp}_{S}(Y, Z^{\natural}) \ar[r] & \bHom^{\sharp}_{S}(X,Z^{\natural}). }$$
Since $p$ is a Cartesian equivalence, the bottom horizontal arrow is a homotopy equivalence.
According to Lemma \ref{insdod}, the right vertical arrow is a Kan fibration; it follows that the
diagram is homotopy Cartesian and so the top horizontal arrow is an equivalence as well.

\item[(3)] A map $p: X \rightarrow Y$ in $(\mSet)_{/S}$ which has the right lifting property with
respect to every map in $(C)$ belongs to $(W)$. Unwinding the definition, we see that
$p$ is a trivial fibration of simplicial sets, and that an edge $e$ of $X$ is marked if and only if
$p(e)$ is a marked edge of $Y$. It follows that $p$ has a section $s$, with $s \circ p$ fiberwise homotopic to $\id_{X}$. From this, we deduce easily that $p$ is a Cartesian equivalence.

\end{itemize}

\end{proof}

\begin{warning}
Let $S$ be a simplicial set. We must be careful to distinguish between
{\it Cartesian fibrations} of simplicial sets (in the sense of Definition \ref{defcart}) and
fibrations with respect to the Cartesian model structure on $(\mSet)_{/S}$ (in the sense of
Proposition \ref{markmodell}). Though distinct, these notions are closely related: for example, the fibrant objects of $(\mSet)_{/S}$ are precisely those objects of the form $X^{\natural}$, where
$X \rightarrow S$ is a Cartesian fibration (Proposition \ref{markedfibrant}).
\end{warning}

\begin{remark}\label{twuff}
The definition of the Cartesian model structure on $(\mSet)_{/S}$ is not self-opposite. Consequently, we can define another model structure on $(\mSet)_{/S}$ as follows:
\begin{itemize}
\item[$(C)$] The cofibrations in $(\mSet)_{/S}$ are precisely the monomorphisms.
\item[$(W)$] The weak equivalences in $(\mSet)_{/S}$ are precisely the
{\it coCartesian equivalences}\index{gen}{equivalence!coCartesian}\index{gen}{coCartesian equivalence}:
that is, those morphisms $f: \overline{X} \rightarrow \overline{Y}$ such that the induced map
$f^{op}: \overline{X}^{op} \rightarrow \overline{Y}^{op}$ is a Cartesian equivalence in
$(\mSet)_{/S^{op}}$. 
\item[$(F)$] The fibrations in $(\mSet)_{/S}$ are those morphisms which have the right lifting property with respect to every morphism satisfying both $(C)$ and $(W)$.
\end{itemize}
We will refer to this model structure on $(\mSet)_{/S}$ as the {\it coCartesian model structure}.\index{gen}{coCartesian model structure}\index{gen}{model category!coCartesian}
\end{remark}

\subsection{Properties of the Cartesian Model Structure}\label{markprop}

In this section, we will establish some of the basic properties of Cartesian model structures
on $(\mSet)_{/S}$ which was introduced in \S \ref{markmodel}. In particular, we will show that each $(\mSet)_{/S}$ is a {\it simplicial} model category, and characterize its fibrant objects.

\begin{proposition}\label{markedfibrant}
An object $X \in (\mSet)_{/S}$ is fibrant $($with respect to the Cartesian model structure$)$ if and only if $X \simeq Y^{\natural}$, where $Y \rightarrow S$ is a Cartesian fibration.
\end{proposition}

\begin{proof}
Suppose first that $X$ is fibrant. The small object argument implies that there exists a marked anodyne map $j: X \rightarrow Z^{\natural}$ for some Cartesian fibration $Z \rightarrow S$. Since $j$ is marked anodyne, it is a Cartesian equivalence. Since $X$ is fibrant, it has the extension property with respect to the trivial cofibration $j$; thus $X$ is a retract of $Z^{\natural}$. It follows that $X$ is isomorphic to $Y^{\natural}$, where $Y$ is a retract of $Z$.

Now suppose that $Y \rightarrow S$ is a Cartesian fibration; we claim that $Y^{\natural}$ has
the right lifting property with respect to any trivial cofibration $j: A \rightarrow B$ in $(\mSet)_{/S}$.
Since $j$ is a Cartesian equivalence, the map $\eta: \bHom_{S}^{\sharp}(B,Y^{\natural}) \rightarrow
\bHom_{S}^{\sharp}(A, Y^{\natural})$ is a homotopy equivalence of Kan complexes.
Hence, for any map $f: A \rightarrow Z^{\natural}$, there is a map $g: B \rightarrow Z^{\natural}$
such that $g|A$ and $f$ are joined by an edge $e$ of $\bHom_{S}^{\sharp}(A,Z^{\natural})$.
Let $M = (A \times (\Delta^1)^{\sharp}) \coprod_{ A \times \{1\}^{\sharp} } (B \times \{1\}^{\sharp})
\subseteq B \times (\Delta^1)^{\sharp}$. We observe that $e$ and $g$ together determine a map $M \rightarrow Z^{\natural}$. Consider the diagram
$$ \xymatrix{ M \ar[r] \ar[d] & Z^{\natural} \ar[d] \\
B \times (\Delta^1)^{\sharp} \ar[r] \ar@{-->}[ur]^{F} & S^{\sharp}. }$$
The left vertical arrow is marked anodyne, by Proposition \ref{markanodprod}. Consequently, there exists a dotted arrow $F$ as indicated. We note that $F|B \times \{0\}$ is an extension
of $f$ to $B$, as desired.
\end{proof}

We now study the behavior of the Cartesian model structures with respect to products.

\begin{proposition}\label{urlt}
Let $S$ and $T$ simplicial sets, and let $Z$ be an object of $(\mSet)_{/T}$. Then the functor $$
(\mSet)_{/S} \rightarrow (\mSet)_{/S \times T}$$
$$ X \mapsto X \times Z$$
preserves Cartesian equivalences.
\end{proposition}

\begin{proof}
Let $f: X \rightarrow Y$ be a Cartesian equivalence in $(\mSet)_{/S}$. We wish to show that $f \times \id_{Z}$ is a Cartesian equivalence in $(\mSet)_{/S \times T}$. Let $X \rightarrow X'$ be a marked anodyne map where $X' \in (\mSet)_{/S}$ is fibrant. Now choose a marked-anodyne map $X' \coprod_{X} Y \rightarrow Y'$, where $Y' \in (\mSet)_{/S}$ is fibrant. Since the product maps $X \times Z \rightarrow X' \times Z$ and $Y \times Z \rightarrow Y' \times Z$ are also marked anodyne
(by Proposition \ref{markanodprod}), it suffices to show that $X' \times Z \rightarrow Y' \times Z$ is a Cartesian equivalence. In other words, we may reduce to the situation where $X$ and $Y$ are fibrant. By Proposition \ref{crispy}, $f$ has a homotopy inverse $g$; then
$g \times \id_Y$ is a homotopy inverse to $f \times \id_Y$.
\end{proof}

\begin{corollary}\label{compatprod}
Let $f: A \rightarrow B$ be a cofibration in $(\mSet)_{/S}$ and $f': A' \rightarrow B'$ a cofibration in $(\mSet)_{/T}$. Then the smash product map
$$(A \times B' ) \coprod_{ A \times B } (A' \times B) \rightarrow A' \times B'$$ is a cofibration in
$(\mSet)_{/S \times T}$, which is trivial if either $f$ or $g$ is trivial.
\end{corollary}

\begin{corollary}
Let $S$ be a simplicial set, and regard $(\mSet)_{/S}$ as a simplicial category with mapping objects given by $\bHom^{\sharp}_{S}(X,Y)$. Then $(\mSet)_{/S}$ is a {\it simplicial} model category.
\end{corollary}

\begin{proof}
Unwinding the definitions, we are reduced to proving the following: given a cofibration
$i: X \rightarrow X'$ in $(\mSet)_{/S}$ and a cofibration $j: Y \rightarrow Y'$ in $\sSet$, the induced cofibration
$$ (X' \times Y^{\sharp}) \coprod_{ X \times Y^{\sharp} } (X \times {Y'}^{\sharp})
\subseteq X' \times {Y'}^{\sharp}$$
in $(\mSet)_{/S}$ is trivial if either $i$ is a Cartesian equivalence of $j$ is a weak homotopy equivalence. If $i$ is trivial, this follows immediately
from Corollary \ref{compatprod}. If $j$ is trivial, the same argument applies, provided that
we can verify that $Y^{\sharp} \rightarrow {Y'}^{\sharp}$ is a Cartesian equivalence in $\mSet$.
Unwinding the definitions, we must show that for every $\infty$-category $Z$, the restriction map
$$ \theta: \bHom^{\sharp}( {Y'}^{\sharp}, Z^{\natural} ) \rightarrow \bHom^{\sharp}( Y^{\sharp}, Z^{\natural}) $$
is a homotopy equivalence of Kan complexes. Let $K$ be the largest Kan complex contained in $Z$, so that $\theta$ can be identified with the restriction map
$$ \bHom_{\sSet}(Y', K) \rightarrow \bHom_{\sSet}(Y,K).$$
Since $j$ is a weak homotopy equivalence, this map is a trivial fibration.
\end{proof}

\begin{remark}
There is a second simplicial structure on $(\mSet)_{/S}$, where the simplicial mapping spaces are given by $\bHom_{S}^{\flat}(X,Y)$. This simplicial structure is {\em not} compatible with
the Cartesian model structure: for fixed $X \in (\mSet)_{/S}$, the functor
$$A \mapsto A^{\flat} \times X$$ does not carry weak homotopy equivalences (in the $A$-variable) to Cartesian equivalences. It does, however, carry {\em categorical} equivalences (in $A$) to Cartesian equivalences, and consequently $(\mSet)_{/S}$ is endowed with the structure of a $\sSet$-enriched model category, where we regard $\sSet$ as equipped with the Joyal model structure. This second simplicial structure reflects the fact that $(\mSet)_{/S}$ is really a model for an $\infty$-bicategory.
\end{remark}

\begin{remark}\label{duality}
Suppose $S$ is a Kan complex. A map $p: X \rightarrow S$ is a Cartesian fibration if and only if it is a coCartesian fibration (this follows in general from Proposition \ref{groob}; if $S = \Delta^0$, the main case of interest for us, it is obvious). Moreover, the class $p$-coCartesian edges of $X$ coincides with the class of $p$-Cartesian edges of $X$: both may be described as the class of equivalences in $X$. Consequently, if $A \in (\mSet)_{/S}$, then
$$ \bHom^{\flat}_{S}(A, X^{\natural}) \simeq \bHom^{\flat}_{S^{op}}(A^{op}, (X^{op})^{\natural})^{op},$$
where $A^{op}$ is regarded as a marked simplicial set in the obvious way. It follows that a map
$A \rightarrow B$ is a Cartesian equivalence in $(\mSet)_{/S}$ if and only if $A^{op} \rightarrow B^{op}$ is a Cartesian equivalence in $(\mSet)_{/S^{op}}$. In other words, the Cartesian model structure on $(\mSet)_{/S}$ is self-dual {\em when $S$ is a Kan complex}. In particular, if $S = \Delta^0$, we deduce that the functor $$ A \mapsto A^{op}$$ determines an autoequivalence of the model category $\mSet \simeq (\mSet)_{/\Delta^0}$.
\end{remark}

\subsection{Comparison of Model Categories}\label{compmodel}

Let $S$ be a simplicial set. We now have a plethora of model structures on categories of simplicial sets over $S$:

\begin{itemize}
\item[(0)] Let $\calC_0$ denote the category $(\sSet)_{/S}$ of simplicial sets over $S$ endowed with the {\em Joyal} model structure defined in \S \ref{compp3}: the cofibrations are monomorphisms of simplicial sets, and the weak equivalences are categorical equivalences.
\item[(1)] Let $\calC_1$ denote the category $(\mSet)_{/S}$ of marked simplicial sets over $S$,
endowed with the {\em marked} model structure of Proposition \ref{markmodell}: the cofibrations are
maps $(X,\calE_X) \rightarrow (Y, \calE_Y)$ which induce monomorphisms $X \rightarrow Y$, and the weak equivalences are the Cartesian equivalences.
\item[(2)] Let $\calC_2$ denote the category $(\mSet)_{/S}$ of marked simplicial sets over $S$,
endowed with the following {\em localization} of the Cartesian model structure: a map
$f: (X, \calE_X) \rightarrow (Y, \calE_Y)$ is a cofibration if the underlying map $X \rightarrow Y$ is a monomorphism, and a weak equivalence if $f: X^{\sharp} \rightarrow Y^{\sharp}$ is a marked
equivalence in $(\mSet)_{/S}$.
\item[(3)] Let $\calC_3$ denote the category $(\sSet)_{/S}$ of simplicial sets over $S$, which
is endowed with the {\em contravariant} model structure described in \S \ref{contrasec}: the cofibrations are the monomorphisms, and the weak equivalences are the contravariant equivalences.
\item[(4)] Let $\calC_4$ denote the category $(\sSet)_{/S}$ of simplicial sets over $S$, endowed with the usual homotopy-theoretic model structure: the cofibrations are the monomorphisms of simplicial sets, and the weak equivalences are the weak homotopy equivalences of simplicial sets.
\end{itemize}

The goal of this section is to study the relationship between these five model categories.
We may summarize the situation as follows:

\begin{theorem}\label{bigdiag}
There exists a sequence of Quillen adjunctions
$$ \calC_0 \stackrel{F_0}{\rightarrow} \calC_1 \stackrel{F_1}{\rightarrow} \calC_2
\stackrel{F_2}{\rightarrow} \calC_3 \stackrel{F_3}{\rightarrow} \calC_4$$
$$ \calC_0 \stackrel{G_0}{\leftarrow} \calC_1 \stackrel{G_1}{\leftarrow} \calC_2
\stackrel{G_2}{\leftarrow} \calC_3 \stackrel{G_3}{\leftarrow} \calC_4$$
which may be described as follows:
\begin{itemize}
\item[$(A0)$] The functor $G_0$ is the forgetful functor from $(\mSet)_{/S}$ to $(\sSet)_{/S}$, which
ignores the collection of marked edges. The functor $F_0$ is the left adjoint to $G_0$, which is given by
$X \mapsto X^{\flat}$. The Quillen adjunction $(F_0,G_0)$ is a Quillen equivalence if $S$ is a Kan complex.

\item[$(A1)$] The functors $F_1$ and $G_1$ are the identity functors on $(\mSet)_{/S}$.

\item[$(A2)$] The functor $F_2$ is the forgetful functor from $(\mSet)_{/S}$ to $(\sSet)_{/S}$, which ignores
the collection of marked edges. The functor $G_2$ is the right adjoint to $F_2$, which is given by
$X \mapsto X^{\sharp}$. The Quillen adjunction $(F_2,G_2)$ is a Quillen equivalence for
{\em every} simplicial set $S$.

\item[$(A3)$] The functors $F_3$ and $G_3$ are the identity functors on $(\mSet)_{/S}$. The Quillen
adjunction $(F_3, G_3)$ is a Quillen equivalence whenever $S$ is a Kan complex.
\end{itemize}
\end{theorem}

The rest of this section is devoted to giving a proof of Theorem \ref{bigdiag}. We will organize our efforts as follows. First, we verify that the model category $\calC_{2}$ is well-defined (the analogous results for the other model structures have already been established). We then consider each of the adjunctions $(F_i,G_i)$ in turn, and show that it has the desired properties.

\begin{proposition}\label{lmark}
Let $S$ be a simplicial set. There exists a left proper, combinatorial model structure on the category
$(\mSet)_{/S}$ which may be described as follows:
\begin{itemize}
\item[$(C)$] A map $f: (X, \calE_{X}) \rightarrow (Y, \calE_Y)$ is a cofibration if and only if the underlying map $X \rightarrow Y$ is a monomorphism of simplicial sets.
\item[$(W)$] A map $f: (X, \calE_{X}) \rightarrow (Y, \calE_{Y})$ is a weak equivalence
if and only if the induced map $X^{\sharp} \rightarrow Y^{\sharp}$ is a Cartesian equivalence
in $(\mSet)_{/S}$.
\item[$(F)$] A map $f: (X, \calE_{X}) \rightarrow (Y, \calE_Y)$ is a fibration if and only if
it has the right lifting property with respect to all trivial cofibrations.
\end{itemize}
\end{proposition}

\begin{proof}
It suffices to show that the conditions of Proposition \ref{goot} are satisfied. We check them in turn:

\begin{itemize}
\item[(1)]  The class $(W)$ of Cartesian equivalences is perfect, in the sense of Definition \ref{perfequiv}. This follows from Corollary \ref{perfpull}, since the class of Cartesian equivalences is perfect, and the functor $(X, \calE_X) \rightarrow X^{\sharp}$ commutes with filtered colimits.

\item[(2)] The class of weak equivalences is stable under pushouts by cofibrations. This
follows from the analogous property of the Cartesian model structure, since the functor
$(X, \calE_X) \mapsto X^{\sharp}$ preserves pushouts.

\item[(3)] A map $p: (X,\calE_X) \rightarrow (Y,\calE_Y)$ which has the right lifting property with respect to
every cofibration is a weak equivalence. In this case, the underlying map of simplicial
sets is a trivial fibration, so the induced map $X^{\sharp} \rightarrow Y^{\sharp}$ has the right lifting property with respect to all trivial cofibrations, and is a Cartesian equivalence as observed in the proof of Proposition \ref{markmodell}.

\end{itemize}
\end{proof} 

\begin{proposition}\label{markedjoyal}
Let $S$ be simplicial set. Consider the adjoint functors
$$ \Adjoint{F_0}{(\sSet)_{/S}}{(\mSet)_{/S}}{G_0}$$
described by the formulas
$$ F_0(X) = X^{\flat}$$
$$ G_0(X,\calE) = X.$$
The adjoint functors $(F_0,G_0)$ determine a Quillen adjunction between $(\sSet)_{/S}$ (with the Joyal model structure) and $(\mSet)_{/S}$ (with the Cartesian model structure). If $S$ is a Kan complex, then $(F_0,G_0)$ is a Quillen equivalence.
\end{proposition}

\begin{proof}
To prove that $(F_0,G_0)$ is a Quillen adjunction, it will suffice to show that $F_1$ preserves cofibrations and trivial cofibrations. The first claim is obvious. For the second, we must show that if
$X \subseteq Y$ is a categorical equivalence of simplicial sets over $S$, then the induced map
$X^{\flat} \rightarrow Y^{\flat}$ is a Cartesian equivalence in $(\mSet)_{/S}$. For this, it suffices to
show that for any Cartesian fibration $p: Z \rightarrow S$, the restriction map
$$ \bHom_{S}^{\flat}(Y^{\flat}, Z^{\natural}) \rightarrow \bHom_{S}^{\flat}(X^{\flat},Z^{\natural})$$ is
a trivial fibration of simplicial sets. In other words, we must show that for every inclusion
$A \subseteq B$ of simplicial sets, it is possible to solve any lifting problem of the form
$$ \xymatrix{ A \ar[r] \ar@{^{(}->}[d] & \bHom_{S}^{\flat}(Y^{\flat}, Z^{\natural}) \ar[d] \\
B \ar[r] \ar@{-->}[ur] & \bHom_{S}^{\flat}(X^{\flat}, Z^{\natural}). }$$
Replacing $Y$ by $Y \times B$ and $X$ by $(X \times B) \coprod_{ X \times A} (Y \times A)$, 
we may suppose that $A = \emptyset$ and $B = \ast$. Moreover, we may rephrase the lifting problem as the problem of constructing the dotted arrow indicated in the following diagram:
$$ \xymatrix{ X \ar@{^{(}->}[d] \ar[r] & Z \ar[d]^{p} \\
Y \ar[r] \ar@{-->}[ur] \ar[r] & S }$$
By Proposition \ref{funkyfibcatfib}, $p$ is a categorical fibration, and the lifting problem has a solution in virtue of the assumption that $X \subseteq Y$ is a categorical equivalence.

Now suppose that $S$ is a Kan complex. We want to prove that $(F_0,G_0)$ is a Quillen equivalence. In other words, we must show that for any fibrant object of $(\mSet)_{/S}$ corresponding to a Cartesian fibration $Z \rightarrow S$, a map $X \rightarrow Z$ in $(\sSet)_{/S}$ is a categorical equivalence if and only if the associated map $X^{\flat} \rightarrow Z^{\natural}$ is a Cartesian equivalence.

Suppose first that $X \rightarrow Z$ is a categorical equivalence. Then the induced map
$X^{\flat} \rightarrow Z^{\flat}$ is a Cartesian equivalence, by the argument given above. It therefore suffices to show that $Z^{\flat} \rightarrow Z^{\natural}$ is a Cartesian equivalence.
Since $S$ is a Kan complex, $Z$ is an $\infty$-category; let $K$ denote the largest Kan complex contained in $Z$. The marked edges of $Z^{\natural}$ are precisely the edges which belong to $K$, so we have a pushout diagram
$$ \xymatrix{ K^{\flat} \ar[r] \ar[d] & K^{\sharp} \ar[d] \\
Z^{\flat} \ar[r] & Z^{\natural}. }$$
It follows that $Z^{\flat} \rightarrow Z^{\natural}$ is marked anodyne, and therefore a Cartesian equivalence.

Now suppose that $X^{\flat} \rightarrow Z^{\natural}$ is a Cartesian equivalence. Choose a factorization $X \stackrel{f}{\rightarrow} Y \stackrel{g}{\rightarrow} Z$, where $f$ is a categorical equivalence and $g$ is a categorical fibration. We wish to show that $g$ is a categorical equivalence. 
Proposition \ref{groob} implies that $Z \rightarrow S$ is a categorical fibration, so that $X' \rightarrow S$ is a categorical fibration. Applying Proposition \ref{groob} again, we deduce that
$Y \rightarrow S$ is a Cartesian fibration. Thus we have a factorization
$$ X^{\flat} \rightarrow Y^{\flat} \rightarrow Y^{\natural} \rightarrow Z^{\natural}$$
where the first two maps are Cartesian equivalences by the arguments given above, and the 
composite map is a Cartesian equivalence. Thus $Y^{\natural} \rightarrow Z^{\natural}$ is an equivalence between fibrant objects of $(\mSet)_{/S}$, and therefore admits a homotopy inverse. The existence of this homotopy inverse proves that $g$ is a categorical equivalence, as desired.
\end{proof}

\begin{proposition}\label{marklocal}
Let $S$ be a simplicial set, and let $F_1$ and $G_1$ denote the identity functor from
$(\mSet)_{/S}$ to itself. Then $(F_1,G_1)$ determines a Quillen adjunction between
$\calC_1$ and $\calC_2$.
\end{proposition}

\begin{proof}
We must show that $F_1$ preserves cofibrations and trivial cofibrations. The first claim is obvious. For the second,
let $B: (\mSet)_{/S} \rightarrow (\mSet)_{/S}$ be the functor defined by
$$B(M, \calE_M) = M^{\sharp}.$$
We wish to show that if $X \rightarrow Y$ is a Cartesian equivalence in $(\mSet)_{/S}$, then
$B(X) \rightarrow B(Y)$ is a Cartesian equivalence.

We first observe that if $X \rightarrow Y$ is marked anodyne, then the induced map $B(X) \rightarrow B(Y)$ is also marked anodyne: by general nonsense, it suffices to check this for the generators described in Definition \ref{markanod}, for which it is obvious. Now return to the case of a general Cartesian equivalence $p: X \rightarrow Y$, and choose a diagram
$$ \xymatrix{ X \ar[r]^{i} \ar[d]^{p} & X' \ar[d] \ar[dr]^{q}\\
Y \ar[r] & X' \coprod_{X} Y \ar[r]^{j} & Y' }$$
in which $X'$ and $Y'$ are (marked) fibrant and $i$ and $j$ are marked anodyne. It follows that
$B(i)$ and $B(j)$ are marked anodyne, and therefore Cartesian equivalences. Thus, to prove that
$B(p)$ is a Cartesian equivalence, it suffices to show that $B(q)$ is a Cartesian equivalence.
But $q$ is a Cartesian equivalence between fibrant objects of $(\mSet)_{/S}$, and therefore has a homotopy inverse. It follows that $B(q)$ also has a homotopy inverse, and is therefore a Cartesian equivalence as desired.
\end{proof}

\begin{remark}
In the language of model categories, we may summarize Proposition \ref{marklocal} by saying that the model structure of Proposition \ref{lmark} is a {\em localization} of the Cartesian model structure on $(\mSet)_{/S}$.
\end{remark}

\begin{proposition}\label{romb}
Let $S$ be a simplicial set, and consider the adjunction
$$ \Adjoint{F_2}{ (\mSet)_{/S}}{ (\sSet)_{/S}}{G_2}$$
determined by the formulas
$$ F_2( X, \calE) = X$$
$$G_2(X) = X^{\sharp}.$$
The adjoint functors $(F_2,G_2)$ determines a Quillen equivalence between $\calC_2$ and
$\calC_3$.
\end{proposition}

\begin{proof}
We first claim that $F_2$ is conservative: that is, a map $f: (X, \calE_X) \rightarrow
(Y, \calE_Y)$ is a weak equivalence in $\calC_2$ if and only if the induced map
$X \rightarrow Y$ is a weak equivalence in $\calC_3$. Unwinding the definition, $f$ is a weak equivalence if and only if $X^{\sharp} \rightarrow Y^{\sharp}$ is a Cartesian equivalence. This holds if and only if, for every Cartesian fibration $Z \rightarrow S$, the induced map
$$\phi: \bHom^{\sharp}_S(Y^{\sharp}, Z^{\natural}) \rightarrow \bHom^{\sharp}_S(X^{\sharp}, Z^{\natural})$$ is a homotopy equivalence. Let $Z^0 \rightarrow S$ be the right fibration
associated to $Z \rightarrow S$ (see Corollary \ref{relativeKan}). There are natural identifications $\bHom^{\sharp}_S(Y^{\sharp}, Z^{\natural}) \simeq \bHom_{S}(Y, Z^0)$, $\bHom^{\sharp}_S(X^{\sharp}, Z^{\natural}) \simeq \bHom_{S}(X,Z^0)$. Consequently, $f$ is a weak equivalence if and only if, for every right fibration $Z^0 \rightarrow S$, the associated map
$$ \bHom_{S}(Y, Z^0) \rightarrow \bHom_{S}(X, Z^0)$$ is a homotopy equivalence. Since $\calC_3$ is a simplicial model category for which the fibrant objects are precisely the right fibrations $Z^0 \rightarrow S$ (Corollary \ref{usewhere1}), this is equivalent to the assertion that $X \rightarrow Y$ is a weak equivalence in $\calC_3$. 

To prove that $(F_2, G_2)$ is a Quillen adjunction, it suffices to show that $F_2$ preserves cofibrations and trivial cofibrations. The first claim is obvious, and the second follows because $F_2$ preserves all weak equivalences, by the above argument.

To show that $(F_2, G_2)$ is a Quillen equivalence, we must show that the unit and counit
$$ LF_2 \circ RG_2 \rightarrow \id$$
$$ \id \rightarrow RG_2 \circ LF_2$$
are weak equivalences. In view of the fact that $F_2 = LF_2$ is conservative, the second assertion follows from the first. As to the first, it suffices to show that if $X$ is a fibrant object of
$\calC_3$, then the counit map $(F_2 \circ G_2)(X) \rightarrow X$ is a weak equivalence.
But this map is an isomorphism.
\end{proof}

\begin{proposition}\label{strstr}
Let $S$ be a simplicial set, and let $F_3$ and $G_3$ denote the identity functor
from $(\sSet)_{/S}$ to itself. Then $(F_3, G_3)$ gives a Quillen adjunction between
$\calC_3$ and $\calC_4$. If $S$ is a Kan complex, then $(F_3, G_3)$ is a Quillen equivalence 
$($in other words, the model structures on $\calC_3$ and $\calC_4$ coincide$)$.
\end{proposition}

\begin{proof}
To prove that $(F_3,G_3)$ is a Quillen adjunction, it suffices to prove that $F_3$ preserves cofibrations and weak equivalences. The first claim is obvious (the cofibrations in $\calC_3$ and $\calC_4$ are the same). For the second, we note that both $\calC_3$ and $\calC_4$ are simplicial model categories in which every object is cofibrant. Consequently, a map $f: X \rightarrow Y$
is a weak equivalence if and only if, for every fibrant object $Z$, the associated map
$\bHom(Y,Z) \rightarrow \bHom(X,Z)$ is a homotopy equivalence of Kan complexes.
Thus, to show that $F_3$ preserves weak equivalences, it suffices to show that
$G_3$ preserves fibrant objects. A map $p: Z \rightarrow S$ is fibrant as an object of 
$\calC_4$ if and only if $p$ is a Kan fibration, and fibrant as an object of $\calC_3$ if and only if $p$ is a right fibration (Corollary \ref{usewhere1}). Since every Kan fibration is a right fibration, it follows that $F_3$ preserves weak equivalences. If $S$ is a Kan complex, then the converse holds: according to Lemma \ref{toothie}, every right fibration $p: Z \rightarrow S$ is a Kan fibration. It follows that $G_3$ preserves weak equivalences as well, so that the two model structures under consideration coincide.
\end{proof}

\section{Straightening and Unstraightening}\label{strsec}

\setcounter{theorem}{0}

Let $\calC$ be a category, and let $\chi: \calC^{op} \rightarrow \Cat$
be a functor from $\calC$ to the category $\Cat$ of small categories.
To this data, we can associate (by means of the {\it Grothendieck construction} discussed in \S \ref{scgp}) a new category $\widetilde{\calC}$ which may be described as follows:
\begin{itemize}
\item The objects of $\widetilde{\calC}$ are pairs $(C, \eta)$ where
$C \in \calC$ and $\eta \in \chi(C)$.
\item Given a pair of objects $(C, \eta), (C', \eta') \in \widetilde{\calC})$, a morphism
from $(C, \eta)$ to $(C', \eta')$ in $\widetilde{\calC}$ is a pair $(f, \alpha)$, where
$f: C \rightarrow C'$ is a morphism in the category $\calC$ and $\alpha: \eta \rightarrow \chi(f)(\eta')$ is a morphism in the category $\chi(C)$.
\item Composition is defined in the obvious way.
\end{itemize}
This construction establishes an equivalence between
$\Cat$-valued functors on $\calC^{op}$ and categories which are {\it fibered over
$\calC$}. (To formulate the equivalence precisely, it is best to view $\Cat$ as
a {\it bicategory}, but we will not dwell on this technical point here.)

The goal of this section is to establish an $\infty$-categorical version of the equivalence described above. We will replace the category $\calC$ by a simplicial set $S$, the category $\Cat$ by the $\infty$-category
$\Cat_{\infty}$, and the notion of ``fibered category'' with the notion of
``Cartesian fibration''. In this setting, we will obtain an equivalence of $\infty$-categories, which arises from a Quillen equivalence of simplicial model categories. On one side, we have the category
$(\mSet)_{/S}$, equipped with the Cartesian model structure (a simplicial model category whose fibrant objects are precisely the Cartesian fibrations $X \rightarrow S$; see \S \ref{markprop}). On the other, we have the category of simplicial functors
$$ \sCoNerve[S]^{op} \rightarrow \mSet,$$
equipped with the projective model structure (see \S \ref{quasilimit3}), whose underlying $\infty$-category is equivalent to $\Fun(S^{op}, \Cat_{\infty})$ (Proposition \ref{gumby444}).
The situation may be summarized as follows:

\begin{theorem}\label{straightthm}
Let $S$ be a simplicial set, $\calC$ a simplicial category, and $\phi: \sCoNerve[S] \rightarrow \calC^{op}$ a functor between simplicial categories. Then there exists a pair of adjoint functors
$$ \Adjoint{ \St^{+}_{\phi}}{(\mSet)_{/S}}{(\mSet)^{\calC}}{\Un^{+}_{\phi}} $$
with the following properties:
\begin{itemize}
\item[$(1)$] The functors $(\St^{+}_{\phi}, \Un^{+}_{\phi})$ determine a Quillen adjunction between
$(\mSet)_{/S}$ $($with the Cartesian model structure$)$ and $(\mSet)^{\calC}$ $($with the projective model structure$)$. 

\item[$(2)$] If $\phi$ is an equivalence of simplicial categories, then $(\St^{+}_{\phi}, \Un^{+}_{\phi})$ is a
Quillen equivalence.
\end{itemize}
\end{theorem}

We will refer to $\St^{+}_{\phi}$ and $\Un^{+}_{\phi}$ as the {\it straightening} and {\it unstraightening} functors, respectively. We will give a construct these functors in \S \ref{markmodel2}, and establish part $(1)$ of Theorem \ref{straightthm}. Part $(2)$ is more difficult and requires some preliminary work; we will begin in \S \ref{funkystructure} by analyzing the structure of Cartesian fibrations $X \rightarrow \Delta^n$.
We will apply these analyses in \S \ref{markmodel24} to complete the proof of Theorem \ref{straightthm} in the case where $S$ is a simplex. In \S \ref{markmodel25}, we will deduce the general result, using formal arguments to reduce to the special case of a simplex. 

In the case where $\calC$ is an ordinary category, the straightening and unstraightening procedures
of \S \ref{markmodel2} can be substantially simplified. We will discuss the situation in
\S \ref{altstr}, where we provide an analogue of Theorem \ref{straightthm} (see Propositions \ref{kudd} and \ref{sulken}).

\subsection{The Straightening Functor}\label{markmodel2}

Let $S$ be a simplicial set, and let $\phi: \sCoNerve[S] \rightarrow \calC^{op}$ be a functor between simplicial categories, which we regard as fixed throughout this section. Our objective is
to define the {\it straightening functor} $\St_{\phi}^{+}: (\mSet)_{/S} \rightarrow (\mSet)^{\calC}$ and its right adjoint $\Un_{\phi}^{+}$. 
The intuition is that an object $X$ of $(\mSet)_{/S}$ associates $\infty$-categories to vertices of $S$ in a homotopy coherent fashion, and the functor $\St^{+}_{\phi}$ ``straightens'' this diagram to obtain an $\infty$-category valued functor on $\calC$. The right adjoint $\Un^{+}_{\phi}$ should be viewed as a forgetful functor, which
takes a strictly commutative diagram and retains the underlying homotopy coherent diagram.

The functors $\St^{+}_{\phi}$ and $\Un^{+}_{\phi}$ are more elaborate versions of the straightening and unstraightening functors introduced in \S \ref{rightstraight}. We begin by recalling the unmarked version of the construction. For each object $X \in (\sSet)_{/S}$, form a pushout diagram of simplicial categories
$$ \xymatrix{ \sCoNerve[X] \ar[r] \ar[d]^{\phi} & \sCoNerve[X^{\triangleright}] \ar[d] \\
\calC^{op} \ar[r] & \calC^{op}_X }$$
where the left vertical map is given by composing $\phi$ with the map
$\sCoNerve[X] \rightarrow \sCoNerve[S]$. The functor
$\St_{\phi} X: \calC \rightarrow \sSet$ is defined by the formula
$$ (\St_{\phi}X)(C) = \bHom_{ \calC^{op}_X }(C, \ast)$$
where $\ast$ denotes the cone point of $X^{\triangleright}$.

We will define $\St^{+}_{\phi}$ by designating certain marked edges on the simplicial sets
$(\St_{\phi}X)(C)$, which depend in a natural way on the marked edges  of $X$.
In order to describe this dependence, we need to introduce a bit of notation.

\begin{notation}
Let $X$ be an object of $(\sSet)_{/S}$. Given an $n$-simplex $\sigma$ of the simplicial set
$\bHom_{\calC^{op}}(C,D)$, we let $\sigma^{\ast}: (\St_{\phi}X)(D)_n \rightarrow
(\St_{\phi}X)(C)_n$ denote the associated map on $n$-simplices.

Let $c$ be a vertex of $X$, and $C = \phi(c) \in \calC$.
We may identify $c$ with a map $c: \Delta^0 \rightarrow X$.
Then $c \star \id_{\Delta^0}: \Delta^1 \rightarrow X^{\triangleright}$ is an edge of
$X^{\triangleright}$, which determines a morphism $C \rightarrow \ast$ in
$\calC^{op}_X$, which we may identify with a vertex $\widetilde{c} \in (\St_{\phi}X)(C)$.

Similarly, suppose that $f: c \rightarrow d$ is an edge of $X$, corresponding to a morphism
$$ C \stackrel{F}{\rightarrow} D$$ in the simplicial category $\calC^{op}$.
We may identify $f$ with a map $f: \Delta^1\rightarrow X$. Then $f \star \id_{\Delta^1}: \Delta^2 \rightarrow X^{\triangleright}$ determines a map $\sCoNerve[\Delta^2] \rightarrow \calC_X$, which we may identify with a diagram (not strictly commutative)
$$ \xymatrix{ C \ar[rr]^{F} \ar[dr]^{\widetilde{c}} & & D \ar[dl]^{\widetilde{d}} \\
& \ast & }$$
together with an edge $$\widetilde{f}: \widetilde{c} \rightarrow \widetilde{d} \circ F = F^{\ast} 
\widetilde{d}$$
in the simplicial set $\bHom_{\calC^{op}_X}(C, \ast) = (\St_{\phi}X)(C)$.
\end{notation}

\begin{definition}\index{gen}{straightening functor}\index{not}{Stphi+@$\St_{\phi}^{+}$}
Let $S$ be a simplicial set, $\calC$ a simplicial category, and $\phi: \sCoNerve[S] \rightarrow \calC^{op}$ a simplicial functor. Let $(X, \calE)$ be an object of $(\mSet)_{/S}$. Then
$$ \St^{+}_{\phi}(X, \calE): \calC \rightarrow \mSet$$
is defined by the formula
$$ \St^{+}_{\phi}(X, \calE)(C) = ((\St_{\phi}X)(C), \calE_{\phi}(C))$$
where $\calE_{\phi}(C)$ is the set of all edges of $(\St_{\phi}X)(C)$ having the form
$$ G^{\ast} \widetilde{f},$$ where $f: d \rightarrow e$ is a marked edge of $X$, giving rise to an
edge $\widetilde{f}: \widetilde{d} \rightarrow F^{\ast} \widetilde{e}$ in $(\St_{\phi}X)(D)$, and $G$ belongs to $\bHom_{\calC^{op}}(C,D)_1$.
\end{definition}

\begin{remark}
The construction
$$ (X, \calE) \mapsto \St^{+}_{\phi}(X, \calE) = (\St_{\phi}X, \calE_{\phi})$$
is obviously functorial in $X$. Note that we may characterize the subsets
$\{ \calE_{\phi}(C) \subseteq (\St_{\phi} X)(C)_1 \}$ as the smallest collection of sets which contain
$\widetilde{f}$, for every $f \in \calE$, and depend functorially on $C$.
\end{remark}

The following formal properties of the straightening functor follow immediately from the definition:

\begin{proposition}\label{formall}

\begin{itemize}
\item[$(1)$] Let $S$ be a simplicial set, $\calC$ a simplicial category, and $\phi: \sCoNerve[S] \rightarrow
\calC^{op}$ a simplicial functor; then the associated straightening functor
$$\St^{+}_{\phi}: (\mSet)_{/S} \rightarrow (\mSet)^{\calC}$$ preserves colimits.

\item[$(2)$] Let $p: S' \rightarrow S$ be a map of simplicial sets, $\calC$ a simplicial category, and
$\phi: \sCoNerve[S] \rightarrow \calC^{op}$ a simplicial functor, and let $\phi': \sCoNerve[S'] \rightarrow \calC^{op}$ denote the composition $\phi \circ \sCoNerve[p]$. 
Let $p_{!}: (\mSet)_{/S'} \rightarrow (\mSet)_{/S}$ denote the forgetful functor, given by composition with $p$. There is a natural isomorphism of functors
$$ \St^{+}_{\phi} \circ p_{!} \simeq \St^{+}_{\phi'}$$
from $(\mSet)_{/S'}$ to $(\mSet)^{\calC}$.

\item[$(3)$] Let $S$ be a simplicial set, $\pi: \calC \rightarrow \calC'$ a simplicial functor between simplicial categories, and $\phi: \sCoNerve[S] \rightarrow
\calC^{op}$ a simplicial functor. Then there is a natural isomorphism of functors $$\St^{+}_{\pi \circ \phi} \simeq \pi_{!} \circ \St^{+}_{\phi}$$
from $(\mSet)_{/S}$ to $(\mSet)^{\calC'}$. Here $\pi_{!}: (\mSet)^{\calC} \rightarrow (\mSet)^{\calC'}$
is the left adjoint to the functor $\pi^{\ast}: (\mSet)^{\calC'} \rightarrow (\mSet)^{\calC}$ given by composition with $\pi$: see \S \ref{quasilimit3}.
\end{itemize}
\end{proposition}

\begin{corollary}\label{giraf}\index{gen}{unstraightening functor}\index{not}{Unphi+@$\Un^{+}_{\phi}$}
Let $S$ be a simplicial set, $\calC$ a simplicial category, and $\phi: \sCoNerve[S] \rightarrow \calC^{op}$ any simplicial functor. The straightening functor $\St^{+}_{\phi}$ has a right adjoint 
$$\Un^{+}_{\phi}: (\mSet)^{\calC} \rightarrow (\mSet)_{/S}.$$
\end{corollary}

\begin{proof}
This follows from part $(1)$ of Proposition \ref{formall} and the adjoint functor theorem. (Alternatively, one can construct $\Un^{+}_{\phi}$ directly; we leave details to the reader.)
\end{proof}

\begin{notation}\index{not}{StS+@$\St_{S}^{+}$}\index{not}{UnS+@$\Un_{S}^{+}$}
Let $S$ be a simplicial set, let $\calC = \sCoNerve[S]^{op}$, and let $\phi: \sCoNerve[S] \rightarrow \calC^{op}$ be the identity map. In this case, we will denote $\St^{+}_{\phi}$ by $\St^{+}_{S}$ and $\Un^{+}_{\phi}$ by $\Un^{+}_{S}$. 
\end{notation}

Our next goal is to show that the straightening and unstraightening functors $(\St^{+}_{\phi}, \Un^{+}_{\phi})$ give a {\em Quillen} adjunction between the model categories $(\mSet)_{/S}$ and
$(\mSet)^{\calC}$. The first step is to show that $\St^{+}_{\phi}$ preserves cofibrations.

\begin{proposition}\label{cougherup}
Let $S$ be a simplicial set, $\calC$ a simplicial category, and $\phi: \sCoNerve[S] \rightarrow \calC^{op}$ a simplicial functor. The functor $\St^{+}_{\phi}$ carries cofibrations $($with respect to the Cartesian model structure on $(\mSet)_{/S}${}$)$ to cofibrations $($with respect to the projective model structure on $(\mSet)^{\calC})${}$)$.
\end{proposition}

\begin{proof}
Let $j: A \rightarrow B$ be a cofibration in $(\mSet)_{/S}$; we wish to show that
$\St^{+}_{\phi}(j)$ is a cofibration. By general nonsense, we may suppose that $j$ is a
generating cofibration, either having the form $(\bd \Delta^n)^{\flat} \subseteq (\Delta^n)^{\flat}$
or $(\Delta^1)^{\flat} \rightarrow (\Delta^1)^{\sharp}$. Using Proposition \ref{formall}, we may reduce to the case where $S=B$, $\calC = \sCoNerve[S]$, and $\phi$ is the identity map. The result now follows from a straightforward computation.
\end{proof}

To complete the proof that $(\St^{+}_{\phi}, \Un^{+}_{\phi})$ is a Quillen adjunction, it suffices to show that
$\St^{+}_{\phi}$ preserves trivial cofibrations. Since every object of $(\mSet)_{/S}$ is cofibrant, this is equivalent to the apparently stronger claim that if $f: X \rightarrow Y$ is a Cartesian equivalence
in $(\mSet)_{/S}$, then $\St^{+}_{\phi}(f)$ is a weak equivalence in $(\mSet)^{\calC}$. The main step
is to establish this in the case where $f$ is marked anodyne. First, we need a few lemmas.

\begin{lemma}\label{blurgh}
Let $\calE$ be the set of all degenerate edges of $\Delta^n \times \Delta^1$, together with the edge
$\{n\} \times \Delta^1$. Let $B \subseteq \Delta^n \times \Delta^1$ be the coproduct
$$ ( \Delta^n \times \{1\} ) \coprod_{ \bd \Delta^n \times \{1\} } (\bd \Delta^n \times \Delta^1).$$
Then the map
$$i: ( B, \calE \cap B_1 ) \subseteq ( \Delta^n \times \Delta^1, \calE)$$ is marked anodyne.
\end{lemma}

\begin{proof}
We must show that $i$ has the left lifting property with respect to every map $p: X \rightarrow S$ satisfying the hypotheses of Proposition \ref{dubudu}. This is simply a reformulation of Proposition \ref{goouse}.
\end{proof}

\begin{lemma}\label{blughel}
Let $K$ be a simplicial set, $K' \subseteq K$ a simplicial subset, and $A$ a set of vertices of $K$.
Let $\calE$ denote the set of all degenerate edges of $K \times \Delta^1$, together with the edges
$\{a\} \times \Delta^1$ where $a \in A$. Let $B = (K' \times \Delta^1) \coprod_{ K' \times \{1\} } (K \times \{1\}) \subseteq K \times \Delta^1$. Suppose that, for every nondegenerate simplex $\sigma$ of $K$, either $\sigma$ belongs to $K'$, or the final vertex of $\sigma$ belongs to $A$. Then the inclusion
$$ (B, \calE \cap B_1) \subseteq (K \times \Delta^1, \calE)$$ is marked anodyne.
\end{lemma}

\begin{proof}
Working cell-by-cell, we reduce to Lemma \ref{blurgh}.
\end{proof}

\begin{lemma}\label{brend}
Let $X$ be a simplicial set, and let $\calE \subseteq \calE'$ be sets of edges of $X$ containing all degenerate edges. The following conditions are equivalent:
\begin{itemize}
\item[$(1)$] The inclusion $(X, \calE) \rightarrow (X, \calE')$ is trivial cofibration in $\mSet$ (with respect to the Cartesian model structure). 
\item[$(2)$] For every $\infty$-category $\calC$ and every map $f: X \rightarrow \calC$ which carries each edge of
$\calE$ to an equivalence in $\calC$, $f$ also carries each edge of $\calE'$ to an equivalence in $\calC$.
\end{itemize}
\end{lemma}

\begin{proof}
By definition, $(1)$ holds if and only if for every $\infty$-category $\calC$, the inclusion
$$j: \bHom^{\flat}( (X, \calE'), \calC^{\natural}) \rightarrow \bHom^{\flat}( (X, \calE), \calC^{\natural} )$$ is a categorical equivalence. Condition $(2)$ is the assertion that $j$ is an isomorphism. Thus $(2)$ implies $(1)$. Suppose that $(1)$ is satisfied, and let $f: X \rightarrow \calC$ be a vertex of
$\bHom^{\flat}((X,\calE), \calC^{\natural})$. By hypothesis, there exists an equivalence $f \simeq f'$, where $f'$ belongs to the image of $j$. Let $e \in \calE'$; then $f'(e)$ is an equivalence in $\calC$. Since $f$ and $f'$ are equivalent, $f(e)$ is also an equivalence in $\calC$. Consequently, $f$ also belongs to the image of $j$, and the proof is complete.
\end{proof}

\begin{proposition}\label{spec2}
Let $S$ be a simplicial set, $\calC$ a simplicial category, and $\phi: \sCoNerve[S] \rightarrow \calC^{op}$ a simplicial functor. The functor $\St^{+}_{\phi}$ carries marked anodyne maps in $(\mSet)_{/S}$ (with respect to the Cartesian model structure) to trivial cofibrations in 
$(\mSet)^{\calC}$ (with respect to the projective model structure).
\end{proposition}

\begin{proof}
Let $f: A \rightarrow B$ be a marked anodyne map in $(\mSet)_{/S}$. We wish to prove that
$\St^{+}_{\phi}(f)$ is a trivial cofibration. It will suffice to prove this under the assumption that $f$ is one of the generators for the class of marked anodyne maps, as given in Definition \ref{markanod}.
Using Proposition \ref{formall}, we may reduce to the case where $S$ is the underyling simplicial set of $B$, $\calC = \sCoNerve[S]^{op}$, and $\phi$ is the identity. There are four cases to consider:

\begin{itemize}
\item[(1)]  Suppose first that $f$ is among the morphisms listed in $(1)$ of Definition \ref{markanod}; that is, $f$ is an inclusion $(\Lambda^n_i)^{\flat} \subseteq (\Delta^n)^{\flat}$, where $0 < i < n$. Let $v_k$ denote the $k$th vertex of $\Delta^n$, which we may also think of as an object of the simplicial category $\calC$. 
We note that $\St^{+}_{\phi}(f)$ is an isomorphism when evaluated at $v_k$ for $k \neq 0$. Let $K$ denote the cube $(\Delta^1)^{ \{ j: 0 < j \leq n, j \neq i \} }$, let $K' = \bd K$, let $A$ denote the set of all vertices of $K$ corresponding to subsets of $\{ j: 0 < j \leq n, j \neq i\}$ which contain an element $> i$, and let
$\calE$ denote the set of all degenerate edges of $K \times \Delta^1$ together with all edges of the form
$\{a\} \times \Delta^1$, where $a \in A$. Finally, let $B = (K \times \{1\} ) \coprod_{ K' \times \{1\} } (K' \times \Delta^1)$. The morphism
$\St^{+}_{\phi}(f)(v_n)$ is a pushout of $g: (B, \calE \cap B_1) \subseteq (K \times \Delta^1, \calE)$. 
Since $i > 0$, we may apply Lemma \ref{blughel} to deduce that $g$ is marked-anodyne, and therefore a trivial cofibration in $\mSet$.

\item[(2)]
Suppose that $f$ is among the morphisms of part $(2)$ in Definition \ref{markanod}; that is, $f$
is an inclusion $$( \Lambda^n_n, \calE \cap (\Lambda^n_n)_{1} ) \subseteq ( \Delta^n, \calF ),$$
where $\calF$ denotes the set of all degenerate edges of $\Delta^n$, together with the final edge $\Delta^{ \{n-1,n\} }$. If $n > 1$, then one can repeat the argument given above in case $(1)$, except that the set of vertices $A$ needs to be replaced by the set of all vertices of $K$ which correspond to subsets of $\{j: 0 < j < n\}$ which contain $n-1$. If $n=1$, then we observe that $\St^{+}_{\phi}(f)(v_n)$ is
isomorphic to the inclusion $\{1\}^{\sharp} \subseteq (\Delta^1)^{\sharp}$, which is again a marked anodyne map and therefore a trivial cofibration in $\mSet$.

\item[(3)]
Suppose next that $f$ is the morphism $$ (\Lambda^2_1)^{\sharp} \coprod_{ (\Lambda^2_1)^{\flat} } (\Delta^2)^{\flat} \rightarrow (\Delta^2)^{\sharp}$$
specified in $(3)$ of Definition \ref{markanod}. Simple computation shows that
$\St^{+}_{\phi}(f)(v_n)$ is an isomorphism for $n \neq 0$, and $\St^{+}_{\phi}(f)(v_0)$ is may be identified with the inclusion $$(\Delta^1 \times \Delta^1, \calE) \subseteq (\Delta^1 \times \Delta^1)^{\sharp},$$ where $\calE$ denotes the set of all degenerate edges of $\Delta^1 \times \Delta^1$ together with
$\Delta^1 \times \{0\}$, $\Delta^1 \times \{1\}$, and $\{1\} \times \Delta^1$. This inclusion may be obtained as a pushout of  $$(\Lambda^2_1)^{\sharp} \coprod_{ (\Lambda^2_1)^{\flat} } (\Delta^2)^{\flat} \rightarrow (\Delta^2)^{\sharp}$$ followed by a pushout of
 $$ (\Lambda^2_2)^{\sharp} \coprod_{ (\Lambda^2_2)^{\flat} } (\Delta^2)^{\flat} \rightarrow (\Delta^2)^{\sharp}.$$ The first of these maps is marked-anodyne by definition; the second is marked anodyne by Corollary \ref{hermes}.

\item[(4)] 
Suppose that $f$ is the morphism $K^{\flat} \rightarrow K^{\sharp}$, where $K$ is a Kan complex, as in $(4)$ of Definition \ref{markanod}. For each vertex $v$ of $K$, let $\St^{+}_{\phi}(K^{\flat})(v)=(X_v, \calE_v)$, so that $\St^{+}_{\phi}(K^{\sharp})= X_v^{\sharp}$. Given a morphism $g \in \bHom_{ \sCoNerve[K]}(v,v')_n$, we let $g^{\ast}: X_v \times \Delta^n \rightarrow X_{v'}$ denote the induced map. We wish to show that the natural map $(X_v, \calE_v) \rightarrow X_{v}^{\sharp}$ is an equivalence in $\mSet$. 
By Lemma \ref{brend}, it suffices to show that for every $\infty$-category $Z$, if $h: X_v \rightarrow Z$
carries each edge belonging to $\calE_{v}$ into an equivalence, then $h$ carries {\em every} edge of $X_v$ to an equivalence. 

We first show that $h$ carries $\widetilde{e}$ to an equivalence, for every edge
$e: v \rightarrow v'$ in $K$. Let $m_{e}: \Delta^1 \rightarrow \bHom_{\calC^{op}}(v,v')$ denote the degenerate edge at the vertex corresponding to $e$.
Since $K$ is a Kan complex, the edge $e: \Delta^1 \rightarrow K$
extends to a $2$-simplex $\sigma: \Delta^2 \rightarrow K$ depicted as follows
$$ \xymatrix{ & v' \ar[dr]^{e'} & \\
v \ar[ur]^{e} \ar[rr]^{\id_{v}} & & v. }$$ 
Let $m_{e'}: \Delta^1 \rightarrow \bHom_{\calC}(v',v)$ denote the degenerate edge corresponding to $e'$. The map $\sigma$ gives rise to a diagram a diagram
$$ \xymatrix{ \widetilde{v} \ar[r]^{ \widetilde{e} } \ar[d]^{\id_{\widetilde{v}}} & e^{\ast} \widetilde{v}' \ar[d]^{ m_e^{\ast} \widetilde{e}' } \\
\widetilde{v} \ar[r] & e^{\ast} (e')^{\ast} \widetilde{v} }$$
in the simplicial set $X_{v}$. Since $h$ carries the left vertical arrow and the bottom horizontal arrow into equivalences, it follows that $h$ carries the composition
$(m_e^{\ast} \widetilde{e'}) \circ \widetilde{e}$ to an equivalence in $Z$; thus
$h(\widetilde{e})$ has a left homotopy inverse. A similar argument shows that $h(\widetilde{e})$ has a right homotopy inverse, so that $h(\widetilde{e})$ is an equivalence.

We observe that every edge of $X_{v}$ has the form
$g^{\ast} \widetilde{e}$, where $g$ is an edge of $\bHom_{\calC^{op}}(v,v')$ and
$e: v' \rightarrow v''$ is an edge of $K$. We wish to show that $h( g^{\ast} \widetilde{e} )$
is an equivalence in $Z$. Above, we have shown that this is true if $v=v'$ and $g$ is the identity.
We now consider the more general case where $g$ is not necessarily the identity, but is a degenerate edge corresponding to some map $v' \rightarrow v$ in $\calC$. Let $h'$ denote the composition
$$ X_{v'} \rightarrow X_{v} \stackrel{h}{\rightarrow} Z.$$
Then $h( g^{\ast} \widetilde{e}) = h'( \widetilde{e} )$ is an equivalence in $Z$ by the argument given above.

Now consider the case where $g: \Delta^1 \rightarrow \bHom_{\calC^{op}}(v,v')$ is nondegenerate.
In this case, there is a simplicial homotopy $G: \Delta^1 \times \Delta^1 \rightarrow \bHom_{\calC}(v,v')$ with $g = G| \Delta^1 \times \{0\}$ and $g' = G | \Delta^1 \times \{1\}$ a degenerate edge of $\bHom_{\calC^{op}}(v,v')$ (for example, we can arrange that $g'$ is
the constant edge at an endpoint of $g$). The map $G$ induces a simplicial homotopy 
$G(e)$ from $g^{\ast} \widetilde{e}$ to $(g')^{\ast} \widetilde{e}$. Moreover, the edges $G(e)| \{0\} \times \Delta^1$ and $G(e)| \{1\} \times \Delta^1$ belong to $\calE_{v}$, and are therefore carried by $h$ into equivalences in $Z$. Since $h$ carries $(g')^{\ast} \widetilde{e}$ into an equivalence of
$Z$, it carries $g^{\ast} \widetilde{e}$ into an equivalence of $Z$, as desired.
\end{itemize}
\end{proof}

We now study the behavior of straightening functors with respect to products.

\begin{notation}
Given two simplicial functors $\calF: \calC \rightarrow \mSet$, $\calF': \calC' \rightarrow \mSet$, we let $\calF \boxtimes \calF': \calC \times \calC' \rightarrow \mSet$ denote the functor described by the formula
$$(\calF \boxtimes \calF')(C,C') = \calF(C) \times \calF'(C').$$
\end{notation}

\begin{proposition}\label{spek3}
Let $S$ and $S'$ be simplicial sets, $\calC$ and $\calC'$ simplicial categories, and $\phi: \sCoNerve[S] \rightarrow \calC^{op}$, $\phi': \sCoNerve[S'] \rightarrow (\calC')^{op}$ simplicial functors; let $\phi \boxtimes \phi'$ denote the induced functor $\sCoNerve[S \times S'] \rightarrow (\calC \times \calC')^{op}$. For every $M \in (\mSet)_{/S}$, $M' \in (\mSet)_{/S'}$, the natural map
$$ s_{M,M'}: \St^{+}_{\phi \boxtimes \phi'}(M \times M') \rightarrow \St^{+}_{\phi}(M) \boxtimes \St^{+}_{\phi'}(M')$$
is a weak equivalence of functors $\calC \times \calC' \rightarrow \mSet$. 
\end{proposition}

\begin{proof}
Since both sides are compatible with the formations of filtered colimits in $M$, we may suppose that $M$ has only finitely many nondegenerate simplices. 
We work by induction on the dimension $n$ of $M$ and the number of $n$-dimensional simplices of $M$. If $M = \emptyset$ there is nothing to prove. If $n \neq 1$, we may choose a nondegenerate simplex of $M$ having maximal dimension and thereby write 
$M = N \coprod_{ (\bd \Delta^n)^{\flat} } (\Delta^n)^{\flat}$. By the inductive hypothesis we may suppose that the result is known
for $N$ and $(\bd \Delta^n)^{\flat}$. The map $s_{M,M'}$ is a pushout of the maps
$s_{N,M'}$ and $s_{ (\Delta^n)^{\flat}, M'}$ over $s_{ (\bd \Delta^n)^{\flat}, M' }$.
Since $\mSet$ is left-proper, this pushout is a homotopy pushout; it therefore suffices to prove the result after replacing $M$ by $N$, $(\bd \Delta^n)^{\flat}$, or $(\Delta^n)^{\flat}$. In the first two cases, the inductive hypothesis implies that $s_{M,M'}$ is an equivalence; we are therefore reduced to the case $M = (\Delta^n)^{\flat}$. If $n=0$, the result is obvious. If $n>2$, we set $$K = \Delta^{ \{0,1\} } \coprod_{ \{1\} } \Delta^{ \{1,2\} } \coprod_{ \{2\} } \ldots \coprod_{ \{n-1\} } \Delta^{ \{n-1, n\} } \subseteq \Delta^n.$$
The inclusion $K \subseteq \Delta^n$ is inner anodyne, so that $K^{\flat} \subseteq M$ is marked-anodyne. By Proposition \ref{spec2}, we deduce that $s_{M,M'}$ is an equivalence if and only if $s_{K^{\flat}, M'}$ is an equivalence, which follows from the inductive hypothesis since $K$ is $1$-dimensional.

We may therefore suppose that $n=1$. Using the above argument, we may reduce to the case where $M$ consists of a single edge, either marked or unmarked. Repeating the above argument with the roles of $M$ and $M'$ interchanged, we may suppose that $M'$ also consists of a single edge. Applying Proposition \ref{formall}, we may reduce to the case where $S = M$, $S' = M'$, $\calC = \sCoNerve[S]^{op}$, and $\calC' = \sCoNerve[S']^{op}$.

Let us denote the vertices of $M$ by $x$ and $y$, and the unique edge joining them by
$e: x \rightarrow y$. Similarly, we let $x'$ and $y'$ denote the vertices of $M'$, and $e': x' \rightarrow y'$ the edge which joins them. We note that the map
$s_{M,M'}$ induces an isomorphism when evaluated on any object of $\calC \times \calC'$
{\em except} $(x,x')$. Moreover, the map
$$s_{M,M'}(x,x'): \St^{+}_{\phi \boxtimes \phi'}(M \times M')(x,x') \rightarrow \St^{+}_{\phi}(M)(x) \times \St^{+}_{\phi'}(M')(x')$$
obtained from $s_{ (\Delta^1)^{\flat}, (\Delta^1)^{\flat}}$ by successive pushouts
along cofibrations of the form $(\Delta^1)^{\flat} \subseteq (\Delta^1)^{\sharp}$. Since $\mSet$ is left proper, we may reduce to the case where $M = M' = (\Delta^1)^{\flat}$. The result now follows from a simple explicit computation.
\end{proof}

We now study the situation in which $S = \Delta^0$, $\calC = \sCoNerve[S]$, and $\phi$ is the identity map. In this case, $\St^{+}_{\phi}$ may be regarded as a functor $T: \mSet \rightarrow \mSet$. The underlying functor of simplicial sets
is familiar: we have
$$ T( X, \calE) = ( | X |_{Q^{\bigdot} } , \calE' ),$$ where $Q$ denotes the cosimplicial object of $\sSet$ considered in \S \ref{twistt}. In that section, we exhibited a natural map $|X|_{Q^{\bigdot}} \rightarrow X$ which we proved to be a weak homotopy equivalence. We now prove a stronger version of that result:

\begin{proposition}\label{spek4}
For any marked simplicial set $M = (X, \calE)$, the natural map
$|X|_{Q^{\bigdot}} \rightarrow X$ induces a Cartesian equivalence
$$ T(M) \rightarrow M.$$
\end{proposition} 

\begin{proof}
As in the proof of Proposition \ref{spek3}, we may reduce to the case where $M$ consists of a simplex of dimension $\leq 1$ (either marked or unmarked). In these cases, the map $T(M) \rightarrow M$ is an isomorphism in $\mSet$.
\end{proof}

\begin{corollary}\label{spek5}
Let $S$ be a simplicial set, $\calC$ a simplicial category, $\phi: \sCoNerve[S] \rightarrow \calC^{op}$ a simplicial
functor, and $X \in (\mSet)_{/S}$ an object. For every $K \in \mSet$, there is a natural equivalence 
$$ \St^{+}_{\phi}(M \times K) \rightarrow \St^{+}_{\phi}(M) \boxtimes K$$
of functors from $\calC$ to $\mSet$.
\end{corollary}

\begin{proof}
Combine the equivalences of Proposition \ref{spek4} (in the case where
$S' = \Delta^0$, $\calC' = \sCoNerve[S']^{op}$, and $\phi'$ is the identity ) and Proposition \ref{spek5}.
\end{proof}

We can now complete the proof that $(\St^{+}_{\phi}, \Un^{+}_{\phi})$ is a Quillen adjunction:

\begin{corollary}\label{spek6}
Let $S$ be a simplicial set, $\calC$ a simplicial category, and $\phi: \sCoNerve[S]^{op} \rightarrow \calC$ a simplicial functor. The straightening functor $\St^{+}_{\phi}$ carries Cartesian equivalences in $(\mSet)_{/S}$ to $($objectwise$)$ Cartesian equivalences in $(\mSet)^{\calC}$.
\end{corollary}

\begin{proof}
Let $f: M \rightarrow N$ be a Cartesian equivalence in $(\mSet)_{/S}$. Choose
a marked anodyne map $M \rightarrow M'$, where $M'$ is fibrant; then choose a marked anodyne map $M' \coprod_M N \rightarrow N'$, with $N'$ fibrant. Since $\St^{+}_{\phi}$ carries marked anodyne maps to equivalences by Proposition \ref{spec2}, it suffices to prove that the induced map $\St^{+}_{\phi}(M') \rightarrow \St^{+}_{\phi}(N')$ is an equivalence. In other words, we may replace $M$ by $M'$ and $N$ by $N'$, thereby reducing to the case where $M$ and $N$ are fibrant.

Since $f$ is an Cartesian equivalence of fibrant objects, it has a homotopy inverse
$g$. We claim that $\St^{+}_{\phi}(g)$ is an inverse to $\St^{+}_{\phi}(f)$ in the homotopy category of $( \mSet )^{\calC}$. We will show that $\St^{+}_{\phi}(f) \circ \St^{+}_{\phi}(g)$ is homotopic to the identity; applying the same argument with the roles of $f$ and $g$ reversed will then establish the desired result.

Since $f \circ g$ is homotopic to the identity, there is a map $h: N \times K^{\sharp} \rightarrow N$, where $K$ is a contractible Kan complex containing vertices $x$ and $y$, such that
$f \circ g = h| N \times \{x\}$ and $\id_N = h|N \times \{y\}$. The map $\St^{+}_{\phi}(h)$ factors as
$$ \St^{+}_{\phi}(N \times K^{\sharp}) \rightarrow \St^{+}_{\phi}(N) \boxtimes K^{\sharp} \rightarrow \St^{+}_{\phi}(N)$$
where the left map is an equivalence by Corollary \ref{spek5} and the right map because $K$ is contractible. Since $\St^{+}_{\phi}(f \circ g)$ and $\St^{+}_{\phi}( \id_N)$ are both sections of $\St^{+}_{\phi}(h)$, they represent the same morphism in the homotopy category of $(\mSet)^{\calC}$.
\end{proof}

\subsection{Cartesian Fibrations over a Simplex}\label{funkystructure}

A map of simplicial sets $p: X \rightarrow S$ is a Cartesian fibration if and only if the pullback map $X \times_S \Delta^n \rightarrow \Delta^n$ is a Cartesian fibration, for each simplex of $S$.
Consequently, we might imagine that Cartesian fibrations $X \rightarrow \Delta^n$ are the ``primitive building blocks'' out of which other Cartesian fibrations are built. The goal of this section is to prove a structure theorem for these building blocks. This result has a number of consequences, and will play a vital role in the proof of Theorem \ref{straightthm}.

Note that $\Delta^n$ is the nerve of the category associated to the linearly ordered set
$$[n] = \{ 0 < 1 < \ldots < n\} .$$
Since a Cartesian fibration $p: X \rightarrow S$ can be thought of as giving a (contravariant) functor from $S$ to $\infty$-categories, it is natural to expect a close relationship between Cartesian fibrations $X \rightarrow \Delta^n$ and composable sequences of maps between $\infty$-categories $$ A^0 \leftarrow A^1 \leftarrow \ldots \leftarrow A^n.$$
In order to establish this relationship, we need to introduce a few definitions.

Suppose given a composable sequence of maps
$$ \phi: A^0 \leftarrow A^1 \leftarrow \ldots \leftarrow A^n$$
of simplicial sets. The {\it mapping simplex} $M(\phi)$ of $\phi$
is defined as follows. If $J$ is a nonempty finite linearly
ordered set with greatest element $j$, then to specify a map
$\Delta^J \rightarrow M(\phi)$ one must specify an
order-preserving map $f: J \rightarrow [n]$ together
with a map $\sigma: \Delta^J \rightarrow A^{f(j)}$. Given an order-preserving map
$p: J \rightarrow J'$ of partially ordered sets containing largest elements $j$ and $j'$,
there is natural map $M(\phi)(\Delta^{J'}) \rightarrow M(\phi)(\Delta^J)$ which carries
$(f,\sigma)$ to $(f \circ p, e \circ \sigma)$, where $e: A^{f(j')} \rightarrow A^{f(p(j))}$ is obtained
from $\phi$ in the obvious way.\index{gen}{mapping simplex}

\begin{remark}\label{megveg}\index{not}{Mphi@$M(\phi)$}
The mapping simplex $M(\phi)$ is equipped with a natural map $p:
M(\phi) \rightarrow \Delta^n$; the fiber of $p$ over the vertex
$j$ is isomorphic to the simplicial set $A^j$.
\end{remark}

\begin{remark}\label{megvegg}
More generally, let $f: [m] \rightarrow [n]$ be an
order-preserving map, inducing a map $\Delta^m \rightarrow
\Delta^n$. Then $M(\phi) \times_{\Delta^n} \Delta^m$ is naturally
isomorphic to $M(\phi')$, where the sequence $\phi'$ is given by
$$ A^{f(0)} \leftarrow \ldots \leftarrow A^{f(m)}.$$
\end{remark}

\begin{notation}\label{conf2}
Let $\phi: A^0 \leftarrow \ldots \leftarrow A^n$ be a composable sequence of maps of simplicial sets. To give an edge $e$ of $M(\phi)$, one must give a pair of integers
$0 \leq i \leq j \leq n$ and an edge $\overline{e} \in A^{j}$. We will say that
$e$ is {\em marked} if $\overline{e}$ is degenerate; let $\calE$ denote the set of all marked edges of $M(\phi)$. Then the pair $(M(\phi), \calE)$ is a marked simplicial set which we will denote by
$M^{\natural}(\phi)$.\index{gen}{mapping simplex!marked}\index{not}{Mnatphi@$M^{\natural}(\phi)$}
\end{notation}

\begin{remark}
There is a potential ambiguity between the terminology of Definition \ref{conf1} and that of
Notation \ref{conf2}. Suppose that 
$ \phi: A^0 \leftarrow \ldots \leftarrow A^n$ is a composable sequence of maps and that
$p: M(\phi) \rightarrow \Delta^n$ is a Cartesian fibration. Then $M(\phi)^{\natural}$ (Definition \ref{conf1})
and $M^{\natural}(\phi)$ (Notation \ref{conf2}) do not generally coincide as marked simplicial sets. We feel that there is little danger of confusion, since it is very rare that $p$ is a Cartesian fibration.
\end{remark}

\begin{remark}\label{funkytok}
The construction of the mapping simplex is functorial, in the sense that a commutative ladder
$$ \xymatrix{ \phi:A^0 \ar[d]^{f_0} & \ldots \ar[l] \ar[d] & A^n \ar[l] \ar[d]^{f_n} \\
\psi:B^0 & \ldots \ar[l] & B^n \ar[l] }$$
induces a map $M(f): M(\phi) \rightarrow M(\psi)$. Moreover, if each $f_i$ is a categorical equivalence, then $f$ is a categorical equivalence (this follows by induction on $n$, using the fact that the Joyal model structure is left proper).
\end{remark}

\begin{definition}
Let $p: X \rightarrow \Delta^n$ be a Cartesian fibration, and let
$$\phi: A^0 \leftarrow \ldots \leftarrow A^n$$ be a composable sequence of maps.\index{gen}{quasi-equivalence}
A map $q: M(\phi) \rightarrow X$ is a {\it quasi-equivalence} if it has the following properties:
\begin{itemize}
\item[$(1)$] The diagram 
$$ \xymatrix{ M(\phi) \ar[rr]^{q} \ar[dr] & & X \ar[dl]^{p} \\
& \Delta^n & }$$ is commutative.

\item[$(2)$] The map $q$ carries marked edges of $M(\phi)$ to $p$-Cartesian edges of $S$; in
other words, $q$ induces a map $M^{\natural}(\phi) \rightarrow X^{\natural}$ of marked simplicial sets.

\item[$(3)$] For $0 \leq i \leq n$, the induced map
$A^i \rightarrow p^{-1} \{i\}$ is a categorical equivalence.
\end{itemize}
\end{definition}

The goal of this section is to prove the following:

\begin{proposition}\label{simplexplay}
Let $p: X \rightarrow \Delta^n$ be a Cartesian fibration. 
\begin{itemize}
\item[$(1)$] There exists a composable sequence of maps
$$ \phi: A^0 \leftarrow A^1 \leftarrow \ldots \leftarrow A^n $$
and a quasi-equivalence $q: M(\phi) \rightarrow X$.

\item[$(2)$] Let 
$$ \phi: A^0 \leftarrow A^1 \leftarrow \ldots \leftarrow A^n $$
be a composable sequence of maps and $q: M(\phi) \rightarrow X$ a quasi-equivalence.
For any map $T \rightarrow \Delta^n$, the induced map
$$ M(\phi) \times_{\Delta^n} T \rightarrow X \times_{\Delta^n} T$$
is a categorical equivalence.
\end{itemize}
\end{proposition}

We first show that, to establish $(2)$ of Proposition \ref{simplexplay}, it suffices to consider
the case where $T$ is a simplex:

\begin{proposition}\label{tulky}
Suppose given a diagram $$ X \rightarrow Y \rightarrow Z$$ of
simplicial sets. For any map $T \rightarrow Z$, we let $X_T$
denote $X \times_Z T$ and $Y_T$ denote $Y \times_Z T$. The
following statements are equivalent:
\begin{itemize}
\item[$(1)$] For any map $T \rightarrow Z$, the induced map $X_T
\rightarrow Y_T$ is a categorical equivalence.

\item[$(2)$] For any $n \geq 0$ and any map $\Delta^n \rightarrow Z$, the
induced map $X_{\Delta^n} \rightarrow Y_{\Delta^n}$ is a
categorical equivalence.
\end{itemize}

\end{proposition}

\begin{proof}
It is clear that $(1)$ implies $(2)$. Let us prove the converse.
Since the class of categorical equivalences is stable under
filtered colimits, it suffices to consider the case where $T$ has
only finitely many nondegenerate simplices. We now work by
induction on the dimension of $T$, and the number of nondegenerate
simplices contained in $T$. If $T$ is empty, there is nothing to
prove. Otherwise, we may write $T = T' \coprod_{\bd \Delta^n}
\Delta^n$. By the inductive hypothesis, the maps
$$ X_{T'} \rightarrow Y_{T'}$$
$$ X_{\bd \Delta^n} \rightarrow Y_{\bd \Delta^n}$$
are categorical equivalences, and by assumption $X_{\Delta^n}
\rightarrow Y_{\Delta^n}$ is a categorical equivalence as well. We
note that $$X_{T} = X_{T'} \coprod_{ X_{\bd \Delta^n} }
X_{\Delta^n}$$ $$Y_{T} = Y_{T'} \coprod_{ Y_{\bd \Delta^n} }
Y_{\Delta^n}.$$ Since the Joyal model structure is left-proper,
these pushouts are homotopy pushouts, and therefore categorically
equivalent to one another.
\end{proof}

Suppose $p: X \rightarrow \Delta^n$ is a Cartesian fibration,
and $q: M(\phi) \rightarrow X$ is a quasi-equivalence. 
Let $f: \Delta^m \rightarrow \Delta^n$ be any
map. We note (see Remark \ref{funkytok}) that $M(\phi)
\times_{\Delta^n} \Delta^m$ may be identified with a mapping
simplex $M(\phi')$, and that the induced map
$$ M(\phi') \rightarrow X \times_{\Delta^n} \Delta^m$$ is again a quasi-equivalence.
Consequently, to establish $(2)$ of Proposition \ref{simplexplay}, it suffices to prove
that every quasi-equivalence is a categorical equivalence. First, we need the following lemma.

\begin{lemma}\label{coraveg}
Let $$ \phi: A^0 \leftarrow \ldots \leftarrow A^n$$
be a composable sequence of maps between simplicial sets, where $n > 0$. Let $y$ be
a vertex of $A^n$, and let the edge $e: y' \rightarrow y$ be the image of
$\Delta^{ \{n-1, n\} } \times \{y\}$ under the map $\Delta^n \times A^n \rightarrow M(\phi)$.
Let $x$ be any vertex of $M(\phi)$ which does not belong to the fiber $A^n$. Then
composition with $e$ induces a weak homotopy equivalence of simplicial sets
$$ \bHom_{ \sCoNerve[M(\phi)] }(x,y') \rightarrow \bHom_{ \sCoNerve[M(\phi)] }(x,y).$$
\end{lemma}

\begin{proof}
Replacing $\phi$ by an equivalent diagram if necessary (using Remark \ref{funkytok}), we may suppose that the map $A^n \rightarrow A^{n-1}$ is a cofibration. Let $\phi'$ denote the composable subsequence
$$ A^0 \leftarrow \ldots \leftarrow A^{n-1}.$$
Let $\calC = \sCoNerve[M(\phi)]$ and let
$\calC_{-} = \sCoNerve[M(\phi')] \subseteq \calC$. There is a pushout diagram in $\sCat$
$$ \xymatrix{ \sCoNerve[A^n \times \Delta^{n-1}] \ar[r] \ar[d] & \sCoNerve[ A^n \times \Delta^n] \ar[d] \\
\calC_{-} \ar[r] & \calC. }$$
This diagram is actually a homotopy pushout, since $\sCat$ is a left proper model category
and the top horizontal map is a cofibration. Form now the pushout
$$ \xymatrix{ \sCoNerve[ A^n \times \Delta^{n-1}] \ar[d] \ar[r] & \sCoNerve[ A^n \times ( \Delta^{n-1} 
\coprod_{ \{n-1\} } \Delta^{ \{n-1, n\} } )] \ar[d] \\
\calC_{-} \ar[r] & \calC_0. }$$
This diagram is also a homotopy pushout. Since the diagram of simplicial sets
$$ \xymatrix{  \{n-1\} \ar[r] \ar[d] & \Delta^{ \{n-1,n\} } \ar[d] \\
\Delta^{n-1} \ar[r] & \Delta^n } $$
is homotopy coCartesian (with respect to the Joyal model structure), we deduce that the natural
map $\calC_0 \rightarrow \calC$ is an equivalence of simplicial categories. It therefore suffices
to prove that composition with $e$ induces a weak homotopy equivalence
$$ \bHom_{\calC_0}(x,y') \rightarrow \bHom_{\calC}(x,y).$$

Form a pushout square
$$ \xymatrix{ \sCoNerve[ A^{n} \times \{n-1,n\} ] \ar[r] \ar[d] & \sCoNerve[ A^n] \times \sCoNerve[ \Delta^{ \{n-1,n\} } ] \ar[d] \\
\calC_0 \ar[r]^{F} & \calC'. }$$
The left vertical map is a cofibration (since $A^n \rightarrow A^{n-1}$ is a cofibration of simplicial sets), and the upper horizontal map is an equivalence of simplicial categories (Corollary \ref{prodcom}). Invoking the left-properness of $\sCat$, we conclude that $F$ is an equivalence of simplicial categories. Consequently, it will suffice to prove
that $\bHom_{\calC'}( F(x), F(y')) \rightarrow \bHom_{\calC'}(F(x),F(y))$ is a weak homotopy equivalence. We now observe that this map is an isomorphism of simplicial sets.
\end{proof}

\begin{proposition}\label{qequiv}
Let $p: X \rightarrow \Delta^n$ be a Cartesian fibration, let
$$ \phi: A^0 \leftarrow \ldots \leftarrow A^n$$ be a composable
sequence of maps of simplicial sets, and let $q: M(\phi)
\rightarrow X$ be a quasi-equivalence. Then $q$ is a categorical
equivalence.
\end{proposition}

\begin{proof}
We proceed by induction on $n$. The result is obvious if $n = 0$,
so let us assume that $n > 0$. Let $\phi'$ denote the composable
sequence of maps
$$ A^0 \leftarrow A^1 \leftarrow \ldots \leftarrow A^{n-1}$$
which is obtained from $\phi$ by omitting $A^n$. Let $v$ denote
the final vertex of $\Delta^n$, and let $T = \Delta^{ \{0,
\ldots, n-1\} }$ denote the face of $\Delta^n$ which is opposite
$v$. Let $X_{v} = X \times_{ \Delta^n} \{v\}$ and $X_T = X
\times_{\Delta^n} T$.

We note that $M(\phi) = M(\phi')\coprod_{ A^n \times
T } (A^n \times \Delta^n) $. We wish to show that the simplicial functor
$$F: \calC \simeq \sCoNerve[M(\phi)] \simeq  \sCoNerve[M(\phi')] \coprod_{
\sCoNerve[A^n \times T]} \sCoNerve[A^n \times \Delta^n] \rightarrow
\sCoNerve[X]$$ is an equivalence of simplicial categories. We note that $\calC$ decomposes
naturally into full subcategories $\calC_{+} = \sCoNerve[ A^n
\times \{v\} ]$ and $\calC_{-} = \sCoNerve[M(\phi')]$, having
the property that $\bHom_{\calC}(X,Y) = \emptyset$ if $x \in
\calC_{+}$, $y \in \calC_{-}$.

Similarly, $\calD = \sCoNerve[X]$ decomposes into full
subcategories $\calD_{+}= \sCoNerve[X_v]$ and $\calD_{-} =
\sCoNerve[X_T]$, satisfying $\bHom_{\calD}(x,y) = \emptyset$ if $x
\in \calD_{+}$ and $y \in \calD_{-}$. We observe that $F$
restricts to give an equivalence between $\calC_{-}$ and
$\calD_{-}$ by assumption, and gives an equivalence between
$\calC_{+}$ and $\calD_{+}$ by the inductive hypothesis. To
complete the proof, it will suffice to show that if $x \in
\calC_{-}$ and $y \in \calC_{+}$, then $F$ induces a homotopy
equivalence
$$ \bHom_{\calC}(x,y) \rightarrow \bHom_{\calD}(F(x),F(y)).$$

We may identify the object $y \in \calC_{+}$ with a vertex of $A^n$.
Let $e$ denote the edge of $M(\phi)$ which is the image of $\{y\} \times
\Delta^{ \{n-1,n\} }$ under the map $A^n \times \Delta^n \rightarrow
M(\phi)$. We let $[e]: y' \rightarrow y$ denote the corresponding
morphism in $\calC$. We have a commutative diagram
$$ \xymatrix{ \bHom_{\calC_{-}}(x,y') \ar[r] \ar[d] & \bHom_{\calC}(x,y) \ar[d] \\
\bHom_{\calD_{-}}(F(x),F(y')) \ar[r] & \bHom_{\calD}(F(x),F(y)). }$$
Here the left vertical arrow is a weak homotopy equivalence by the inductive hypothesis,
and the bottom horizontal arrow (which is given by composition with $[e]$) is a weak homotopy equivalence because $q(e)$ is $p$-Cartesian. Consequently, to complete the proof, it suffices to show that the top horizontal arrow (given by composition with $e$) is a weak homotopy equivalence. This follows immediately from Lemma \ref{coraveg}.
\end{proof}

To complete the proof of Proposition \ref{simplexplay}, it now suffices to show that for any Cartesian fibration $p: X \rightarrow \Delta^n$, there exists a quasi-equivalence $M(\phi) \rightarrow X$.
In fact, we will prove something slightly stronger (in order to make
our induction work):

\begin{proposition}\label{sharpsimplex}
Let $p: X \rightarrow \Delta^n$ be a Cartesian fibration of simplicial sets and $A$ another simplicial set.
Suppose given a commutative diagram of marked simplicial sets
$$ \xymatrix{ A^{\flat} \times (\Delta^n)^{\sharp} \ar[dr] \ar[rr]^{s} & & X^{\natural} \ar[dl] \\
& (\Delta^n)^{\sharp}. & }$$

Then there exists a sequence of composable
morphisms
$$ \phi: A^0 \leftarrow \ldots \leftarrow A^n,$$
a map $A \rightarrow A^n$, and an extension
$$ \xymatrix{ A^{\flat} \times (\Delta^n)^{\sharp} \ar[dr] \ar[r] & M^{\natural}(\phi) \ar[r]^{f} \ar[d] & X^{\natural} \ar[dl] \\
& (\Delta^n)^{\sharp}. & }$$
of the previous diagram, such that $f$ is a quasi-equivalence.
\end{proposition}

\begin{proof}
The proof goes by induction on $n$. We begin by considering the
fiber $s$ over the final vertex $v$ of $\Delta^n$. The map
$s_v: A \rightarrow X_v = X \times_{\Delta^n} \{v\}$ admits a
factorization
$$ A \stackrel{g}{\rightarrow} A^n \stackrel{h}{\rightarrow} S_v$$
where $g$ is a cofibration and $h$ is a trivial Kan fibration. The smash product inclusion
$$ (\{v\}^{\sharp} \times (A^n)^{\flat}) \coprod_{ \{v\}^{\sharp} \times A^{\flat}} ((\Delta^n)^{\sharp} \times A^{\flat}) \subseteq (\Delta^n)^{\sharp} \times (A^n)^{\flat}$$
is marked anodyne (Proposition \ref{markanodprod}). 
Consequently, we deduce the existence of a dotted arrow $f_0$ as indicated in the diagram
$$ \xymatrix{ A^{\flat} \times (\Delta^n)^{\sharp} \ar@{^{(}->}[d] \ar[r] & X^{\natural} \ar[d] \\
(A^n)^{\flat} \times (\Delta^n)^{\sharp} \ar@{-->}[ur]^{f_0} \ar[r] & (\Delta^n)^{\sharp} }$$
of marked simplicial sets, where $f_0 |(A^n \times \{n\}) = h$.

If $n=0$, we are now done. If $n > 0$, then we apply the inductive hypothesis to the diagram
$$ \xymatrix{ (A^n)^{\flat} \times (\Delta^{n-1})^{\sharp} \ar[dr] \ar[rr]^{f_0 | A^n \times \Delta^{n-1}} & & (X \times_{\Delta^n} \Delta^{n-1})^{\natural} \ar[dl] \\
& (\Delta^{n-1})^{\sharp} & }$$
to deduce the existence of a 
composable sequence of maps
$$ \phi': A^0 \leftarrow \ldots \leftarrow A^{n-1},$$ a map
$A^n \rightarrow A^{n-1}$, and a commutative diagram
$$ \xymatrix{ (A^n)^{\flat} \times (\Delta^{n-1})^{\sharp} \ar[dr] \ar[r] & M^{\natural}(\phi') \ar[r]^{f'} & (X \times_{\Delta^n} \Delta^{n-1})^{\natural} \ar[dl] \\
& (\Delta^{n-1})^{\sharp} & }$$
where $f'$ is a quasi-equivalence. We now define
$\phi$ to be the result of appending the map $A^n
\rightarrow A^{n-1}$ to the beginning of $\phi'$, and let $f: M(\phi) \rightarrow X$ be the
map obtained by amalgamating $f_0$ and $f'$.
\end{proof}

\begin{corollary}\label{presalad}
Let $p: X \rightarrow S$ be a Cartesian fibration of simplicial sets, and let
$q: Y \rightarrow Z$ be a coCartesian fibration. Define new simplicial sets
$Y'$ and $Z'$ equipped with maps $Y' \rightarrow S$, $Z' \rightarrow S$ via the formulas
$$ \Hom_{S}(K, Y') \simeq \Hom( X \times_{S} K, Y)$$
$$ \Hom_{S}(K,Z') \simeq \Hom( X \times_{S} K, Z).$$
Then:
\begin{itemize}
\item[$(1)$] Composition with $q$ determines a coCartesian fibration
$q': Y' \rightarrow Z'$.
\item[$(2)$] An edge $\Delta^1 \rightarrow Y'$ is $q'$-coCartesian if and only if 
the induced map $\Delta^1 \times_{S} X \rightarrow Y$
carries $p$-Cartesian edges to $q$-coCartesian edges.
\end{itemize}
\end{corollary}

\begin{proof}
Let us say that an edge of $Y'$ is {\it special} if it satisfies the hypothesis of $(2)$. Our first goal is to show that there is a sufficient supply of special edges in $Y'$. More precisely, we claim that given any edge $e: z \rightarrow z'$ in $Z'$ and any vertex $\widetilde{z} \in Y'$ covering $z$, there exists a special edge $\widetilde{e}: \widetilde{z} \rightarrow \widetilde{z}'$ of $Y'$ which covers $e$. 

Suppose that the edge $e$ covers an edge $e_0: s \rightarrow s'$ in $S$. We can identify
$\widetilde{z}$ with a map from $X_{s}$ to $Y$.
Using Proposition \ref{simplexplay}, we can choose a morphism $\phi: X'_{s} \leftarrow X'_{s'}$ and a quasi-equivalence $M(\phi) \rightarrow X \times_{S} \Delta^1$. Composing with $\widetilde{z}$, we obtain a map $X'_{s} \rightarrow Y$. Using Propositions \ref{funkyfibcatfib} and \ref{princex}, we may reduce to the problem of providing a dotted arrow in the diagram
$$ \xymatrix{ X'_{s} \ar@{^{(}->}[d] \ar[r] & Y \ar[d]_{q} \\
M(\phi) \ar@{-->}[ur] \ar[r] & Z }$$
which carries the marked edges of $M^{\natural}(\phi)$ to $q$-coCartesian edges of $Y$.
This follows from the the fact that $q^{X_{s}}: Y^{X_{s}} \rightarrow Z^{X_{s}}$ is a coCartesian fibration, and the description of the $q^{X_{s}}$-coCartesian edges (Proposition \ref{doog}).

To complete the proofs of $(1)$ and $(2)$, it will suffice to show that $q'$ is an inner fibration and that every special edge of $Y'$ is $q'$-coCartesian. For this, we must show that every lifting problem
$$ \xymatrix{ \Lambda^n_i \ar[r]^{\sigma_0} \ar@{^{(}->}[d] & Y' \ar[d]^{q'} \\
\Delta^n \ar[r] \ar@{-->}[ur] & Z' }$$
has a solution, provided that either $0 < i < n$, or $i =0$, $n \geq 2$, and
$\sigma_0 | \Delta^{ \{0,1\} }$ is special. We can reformulate this lifting problem
using the diagram
$$ \xymatrix{ X \times_{S} \Lambda^n_i \ar[r] \ar@{^{(}->}[d] & Y \ar[d]^{q} \\
X \times_{S} \Delta^n \ar[r] \ar@{-->}[ur] & Z. }$$
Using Proposition \ref{simplexplay}, we can choose a composable sequence of morphisms
$$ \psi: X'_0 \leftarrow \ldots \leftarrow X'_{n} $$
and a quasi-equivalence $M(\psi) \rightarrow X \times_{S} \Delta^n$. Invoking
Propositions \ref{funkyfibcatfib} and \ref{princex}, we may reduce to the associated mapping problem
$$ \xymatrix{ M(\psi) \times_{ \Delta^n} \Lambda^n_i \ar[r] \ar[d] & Y \ar[d]^{q} \\
M(\psi) \ar[r] \ar@{-->}[ur] & Z.}$$
Since $i < n$, this is equivalent to the mapping problem
$$ \xymatrix{ X'_{n} \times \Lambda^n_i \ar[r] \ar@{^{(}->}[d] & Y \ar[d]^{q} \\
X'_{n} \times \Delta^n \ar[r] & Z, }$$
which admits a solution in virtue of Proposition \ref{doog}.
\end{proof}

\begin{corollary}\label{skinnysalad}
Let $p: X \rightarrow S$ be a Cartesian fibration of simplicial sets, and let
$q: Y \rightarrow S$ be a coCartesian fibration. Define a new simplicial set $T$
equipped with a map $T \rightarrow S$ by the formula
$$ \Hom_{S}(K, T) \simeq \Hom_{S}( X \times_{S} K, Y).$$
Then:
\begin{itemize}
\item[$(1)$] The projection $r: T \rightarrow S$ is a coCartesian fibration.
\item[$(2)$] An edge $\Delta^1 \rightarrow Z$ is $r$-coCartesian if and only if 
the induced map $\Delta^1 \times_{S} X \rightarrow \Delta^1 \times_{S} Y$
carries $p$-Cartesian edges to $q$-coCartesian edges.
\end{itemize}
\end{corollary}

\begin{proof}
Apply Corollary \ref{presalad} in the case where $Z = S$.
\end{proof}

We conclude by noting the following additional property of quasi-equivalences, using the terminology of \S \ref{markmodel}:

\begin{proposition}\label{halfy}
Let $S = \Delta^n$, let $p: X \rightarrow S$ be a Cartesian fibration, let
$$\phi: A^0 \leftarrow \ldots \leftarrow A^n$$ be a composable sequence of maps, and let
$q: M(\phi) \rightarrow X$ be a quasi-equivalence. The induced map $M^{\natural}(\phi) \rightarrow X^{\natural}$ is a Cartesian equivalence in $(\mSet)_{/S}$.
\end{proposition}

\begin{proof}
We must show that for any Cartesian fibration $Y \rightarrow S$, the induced map
of $\infty$-categories $$\bHom^{\flat}_{S}(X^{\natural}, Y^{\natural}) \rightarrow \bHom^{\flat}_S( M^{\natural}(\phi), Y^{\natural} )$$
is a categorical equivalence. Because $S$ is a simplex, the left side may be identified with a full subcategory of $Y^X$ and the right side with a full subcategory of $Y^{M(\phi)}$. Since $q$ is a categorical equivalence, the natural map
$Y^X \rightarrow Y^{M(\phi)}$ is a categorical equivalence; thus, to complete the proof, it suffices to observe that a map of simplicial sets $f: X \rightarrow Y$ is compatible with the projection to $S$ and preserves marked edges if and only if $q \circ f$ has the same properties.
\end{proof}

\subsection{Straightening over a Simplex}\label{markmodel24}

Let $S$ be a simplicial set, $\calC$ a simplicial category, and $\phi: \sCoNerve[S]^{op} \rightarrow \calC$ a simplicial functor. In \S \ref{markmodel2}, we introduced the straightening and unstraightening functors
$$ \Adjoint{ \St^{+}_{\phi} }{ (\mSet)_{/S}}{(\mSet)^{\calC}}{ \Un^{+}_{\phi}}.$$
In this section, we will prove that $( \St^{+}_{\phi}, \Un^{+}_{\phi} )$ is a Quillen equivalence provided that $\phi$ is a categorical equivalence and $S$ is a simplex (the case of a general simplicial set $S$ will be treated in \S \ref{markmodel25}). 

Our first step is to prove the result in the case where $S$ is a point and $\phi$ is an isomorphism
of simplicial categories. We can identify the functor $\St^{+}_{\Delta^0}$ with the functor
$T: \mSet \rightarrow \mSet$ studied in \S \ref{markmodel2}. Consequently, Theorem \ref{straightthm} is an immediate consequence of Proposition \ref{spek4}:

\begin{lemma}\label{utest}
The functor $T: \mSet \rightarrow \mSet$ has a right adjoint $U$, and the pair
$(T,U)$ is a Quillen equivalence from $\mSet$ to itself.
\end{lemma}

\begin{proof}
We have already established the existence of the unstraightening functor $U$ in
\S \ref{markmodel2}, and proved that $(T,U)$ is a Quillen adjunction. To complete the proof,
it suffices to show that the left derived functor of $T$ (which we may identify with $T$, since
every object of $\mSet$ is cofibrant) is an equivalence from the homotopy category
of $\mSet$ to itself. But Proposition \ref{spek4} asserts that $T$ is isomorphic to the identity functor on the homotopy category of $\mSet$.
\end{proof}

Let us now return to the case of a general equivalence $\phi: \sCoNerve[S] \rightarrow \calC^{op}$.
Since we know that $(\St^{+}_{\phi}, \Un^{+}_{\phi})$ give a Quillen adjunction between $(\mSet)_{/S}$ and
$(\mSet)^{\calC}$, it will suffice to prove that the unit and counit
$$ u: \id \rightarrow R \Un^{+}_{\phi} \circ L \St^{+}_{\phi}$$
$$ v: L \St^{+}_{\phi} \circ R \Un^{+}_{\phi} \rightarrow \id$$
are weak equivalences. Our first step is to show that $R \Un^{+}_{\phi}$ detects weak equivalences: this reduces the problem of proving that $v$ is an equivalence to the problem of proving that $u$ is an equivalence. 

\begin{lemma}\label{garbz}
Let $S$ be a simplicial set, $\calC$ a simplicial category, and $\phi: \sCoNerve[S] \rightarrow \calC^{op}$ an essentially surjective functor. Let $p: \calF \rightarrow \calG$ be a map between (weakly) fibrant objects
of $(\mSet)^{\calC}$. Suppose that $\Un^{+}_{\phi}(p): \Un^{+}_{\phi} \calF \rightarrow \Un^{+}_{\phi} \calG$ is a Cartesian equivalence. Then $p$ is an equivalence.
\end{lemma}

\begin{proof}
Since $\phi$ is essentially surjective, it suffices to prove that $\calF(C) \rightarrow \calF(D)$
is a Cartesian equivalence for every object $C \in \calC$ which lies in the image of $\phi$. 
Let $s$ be a vertex of $S$ with $\psi(s) = C$. Let $i: \{s\} \rightarrow S$ denote the inclusion, and
$i^{\ast}: (\mSet)_{/S} \rightarrow \mSet$ denote the functor of passing to the fiber over $s$:
$$ i^{\ast} X = X_{s} = X \times_{ S^{\sharp} } \{s\}^{\sharp}.$$
Let $i_{!}$ denote the left adjoint to $i^{\ast}$. Let $\{C\}$ denote the trivial
category with one object (and only the identity morphism), and let $j: \{C\} \rightarrow \calC$ be the simplicial functor corresponding to the inclusion of $C$ as an object of $\calC$. According to Proposition \ref{formall}, we have a natural identification of functors
$$ \St^{+}_{\phi} \circ i_{!} \simeq j_{!} \circ T.$$
Passing to adjoints, we get another identification
$$ i^{\ast} \circ \Un^{+}_{\phi} \simeq U \circ j^{\ast}$$
from $(\mSet)^{\calC}$ to $\mSet$. Here $U$ denotes the right adjoint of $T$.

According to Lemma \ref{utest}, the functor $U$ detects equivalences between fibrant objects
of $\mSet$. Thus, it suffices to prove that $U( j^{\ast} \calF) \rightarrow U( j^{\ast} \calG)$ is a Cartesian equivalence. Using the identification above, we are reduced to proving that
$$ \Un^{+}_{\phi}(\calF)_{s} \rightarrow \Un^{+}_{\phi}(\calG)_{s}$$ is a Cartesian equivalence.
But $\Un^{+}_{\phi}(\calF)$ and $\Un^{+}_{\phi}(\calG)$ are fibrant objects of $(\mSet)_{/S}$, and therefore correspond to Cartesian fibrations over $S$: the desired result now follows from
Proposition \ref{crispy}.
\end{proof}

We have now reduced the proof of Theorem \ref{straightthm} to the problem of showing that
if $\phi: \sCoNerve[S] \rightarrow \calC^{op}$ is an equivalence of simplicial categories, then
the unit transformation
$$ u: \id \rightarrow R \Un^{+}_{\phi} \circ \St^{+}_{\phi}$$
is an isomorphism of functors from the homotopy category $\h{(\mSet)_{/S}}$ to itself. 

Our first step is to analyze the effect of the straightening functor $\St^{+}_{\phi}$ on a mapping simplex.
We will need a bit of notation. For any $X \in (\mSet)_{/S}$ and any vertex $s$ of $S$,
we let $X_{s}$ denote the fiber $X \times_{ S^{\sharp} } \{s\}^{\sharp}$, and let $i^{s}$
denote the composite functor
$$ \{s\} \hookrightarrow \sCoNerve[S] \stackrel{\phi}{\rightarrow} \calC^{op}$$
of simplicial categories. According to Proposition \ref{formall}, there is a natural identification
$$ \St^{+}_{\phi}(X_s) \simeq i^{s}_{!} T(X_s),$$ and consequently an induced map
$$ \psi^X_{s}: T(X_s) \rightarrow \St^{+}_{\phi}(X)(s).$$

\begin{lemma}\label{piecemeal}
Let $$\theta: A^0 \leftarrow \ldots \leftarrow A^n$$ 
be a composable sequence of maps of simplicial sets, and let
$M^{\natural}(\theta) \in (\mSet)_{\Delta^n}$ be its mapping simplex.
For each $0 \leq i \leq n$, the map
$$ \psi^{M^{\natural}(\theta)}_{i}: T(A^i)^{\flat} \rightarrow \St^{+}_{\Delta^n}(M^{\natural}(\theta))(i) $$  
is a Cartesian equivalence in $\mSet$.
\end{lemma}

\begin{proof}
The proof goes by induction on $n$. We first observe that $\psi^{M^{\natural}(\theta)}_{n}$
is an isomorphism; we may therefore restrict our attention to $i < n$. 
Let $\theta'$ be the composable sequence
$$ A^0 \leftarrow \ldots \leftarrow A^{n-1},$$
and $M^{\natural}(\theta')$ its mapping simplex, which we may regard either as an object
of $(\mSet)_{/\Delta^n}$ or $(\mSet)_{/\Delta^{n-1}}$.

For $i < n$, we have a commutative diagram
$$ \xymatrix{ & \St^{+}_{\Delta^n}(M^{\natural}(\theta'))(i) \ar[dr]^{f_{i}} & \\
T((A^i)^{\flat}) \ar[ur]^{ \psi^{M^{\natural}(\theta')}_i} \ar[rr]
 & & \St^{+}_{\Delta^n}(M^{\natural}(\theta))(i).}$$
By Proposition \ref{formall}, $\St^{+}_{\Delta^n} M^{\natural}(\theta') \simeq j_! \St^{+}_{\Delta^{n-1}} M^{\natural}(\theta')$, where $j: \sCoNerve[\Delta^{n-1}] \rightarrow \sCoNerve[\Delta^n]$ denotes the inclusion. Consequently, the inductive hypothesis implies that the maps
$$ T(A^i)^{\flat} \rightarrow \St^{+}_{\Delta^{n-1}}(M^{\natural}(\theta'))(i) $$ are Cartesian equivalences for $i < n$. It now suffices to prove that $f_i$ is a Cartesian equivalence, for $i < n$.

We observe that there is a (homotopy) pushout diagram
$$ \xymatrix{ (A^n)^{\flat} \times (\Delta^{n-1})^{\sharp} \ar[r] \ar[d] & (A^n)^{\flat} \times (\Delta^n)^{\sharp} \ar[d] \\
M^\natural(\theta') \ar[r] & M^\natural(\theta) }.$$
Since $\St^{+}_{\Delta^n}$ is a left Quillen functor, it induces a homotopy pushout diagram
$$ \xymatrix{ \St^{+}_{\Delta^n} ( (A^n)^{\flat} \times (\Delta^{n-1})^\sharp) \ar[r]^{g} \ar[d] &
\St^{+}_{\Delta^n} ( (A^n)^{\flat} \times (\Delta^n)^{\sharp} ) \ar[d] \\
\St^{+}_{\Delta^n} M^{\natural}(\theta') \ar[r] & \St^{+}_{\Delta^n} M^\natural(\theta).}$$
in $(\mSet)^{\calC}$. We are therefore reduced to proving that $g$ induces a Cartesian equivalence 
after evaluation at any $i < n$.

According to Proposition \ref{spek3}, the vertical maps of the diagram
$$ \xymatrix{ 
\St^{+}_{\Delta^n} ((A^n)^{\flat} \times (\Delta^{n-1})^{\sharp}) \ar[r] \ar[d] & \St^{+}_{\Delta^n} ((A^n)^{\flat} \times (\Delta^n)^{\sharp}) \ar[d] \\ T(A^n)^{\flat} \boxtimes
\St^{+}_{\Delta^n} (\Delta^{n-1})^{\sharp} \ar[r] & T(A^n)^{\flat} \boxtimes \St^{+}_{\Delta^n} (\Delta^n)^{\sharp}}$$
are Cartesian equivalences. To complete the proof we must show that
$$\St^{+}_{\Delta^n}( \Delta^{n-1})^{\sharp} \rightarrow \St^{+}_{\Delta^n}( \Delta^n )^{\sharp}$$ induces a Cartesian equivalence when evaluated at any $i < n$. Consider the diagram
$$ \xymatrix{ \{ n-1 \}^{\sharp} \ar[r] \ar[d] & ( \Delta^{n-1} )^{\sharp} \ar[d] \\
 ( \Delta^{ \{n-1,n\}} )^{\sharp} \ar[r] & (\Delta^n)^{\sharp}.}$$
The horizontal arrows are marked anodyne. It therefore suffices to show that
$$\St^{+}_{\Delta^n} \{n-1\}^{\sharp} \rightarrow \St^{+}_{\Delta^n} (\Delta^{ \{n-1,n\} })^{\sharp}$$
induces Cartesian equivalences when evaluated at any $i < n$. This follows from an easy computation.
\end{proof}

\begin{proposition}\label{speccase}
Let $n \geq 0$. Then the Quillen adjunction 
$$\Adjoint{\St^{+}_{\Delta^n}}{ (\mSet)_{/ \Delta^n} }{ (\mSet)^{\sCoNerve[\Delta^n]}}{ \Un^{+}_{\Delta^n}}$$ is a Quillen equivalence.
\end{proposition}

\begin{proof}
As we have argued above, it suffices to show that the unit
$$ \id \rightarrow R \Un^{+}_{\phi} \circ \St^{+}_{\Delta^n}$$
is an isomorphism of functors from $\h{ (\mSet)_{\Delta^n}}$ to itself. In other words, we must show that given an object
$X \in (\mSet)_{/ \Delta^n}$ and a weak equivalence
$$ \St^{+}_{\Delta^n} X \rightarrow \calF,$$
where $\calF \in (\mSet)^{\sCoNerve[\Delta^n]}$ is fibrant, the adjoint map
$$ j: X \rightarrow \Un^{+}_{\Delta^n} \calF$$ is a Cartesian equivalence in $(\mSet)_{/ \Delta^n}$.

Choose a fibrant replacement for $X$: that is, a Cartesian equivalence  $X \rightarrow Y^{\natural}$ where $Y \rightarrow \Delta^n$ is a Cartesian fibration. According to Proposition \ref{simplexplay}, there exists a composable sequence of maps
$$ \theta: A^0 \leftarrow \ldots \leftarrow A^n$$
and a quasi-equivalence $M^{\natural}(\theta) \rightarrow Y^{\natural}$. Proposition \ref{halfy}
implies that $M^{\natural}(\theta) \rightarrow Y^{\natural}$ is a Cartesian equivalence. Thus,
$X$ is equivalent to $M^{\natural}(\theta)$ in the homotopy category of $(\mSet)_{/\Delta^n}$ and we are free to replace $X$ by $M^{\natural}(\theta)$, thereby reducing to the case where $X$ is a mapping simplex.

We wish to prove that $j$ is a Cartesian equivalence. Since $\Un^{+}_{\Delta^n} \calF$ is fibrant,
Proposition \ref{halfy} implies that it suffices to show that $j$ is a quasi-equivalence: in other words, 
we need to show that the induced map of fibers 
$j_s: X_{s} \rightarrow (\Un^{+}_{\Delta^n} \calF)_s$ is a Cartesian equivalence, for each vertex $s$ of $\Delta^n$. 
As in the proof of Lemma \ref{garbz}, we may identify $(\Un^{+}_{\Delta^n} \calF)_s$ with
$U(\calF(s))$, where $U$ is the right adjoint to $T$. By Lemma \ref{utest}, $X_{s} \rightarrow U (\calF(s))$ is a Cartesian equivalence if and only if the adjoint map
$T(X_{s}) \rightarrow \calF(s)$ is a Cartesian equivalence. This map factors as a composition
$$ T(X_{s}) \rightarrow \St^{+}_{\Delta^n}(X)(s) \rightarrow \calF(s).$$
The map on the left is a Cartesian equivalence by Lemma \ref{piecemeal}, and the map on
the right in virtue of the assumption that $\St^{+}_{\Delta^n} X \rightarrow \calF$ is a weak equivalence.
\end{proof}

\subsection{Straightening in the General Case}\label{markmodel25}

Let $S$ be a simplicial set and $\phi: \sCoNerve[S] \rightarrow \calC^{op}$ an equivalence of simplicial categories. Our goal in this section is to complete the proof of Theorem \ref{straightthm} by showing that $( \St_{\phi}^{+}, \Un_{\phi}^{+})$ is a Quillen equivalence between
$(\mSet)_{/S}$ and $(\mSet)^{\calC}$. In \S \ref{markmodel24}, we handled the case where $S$ was a simplex (and $\phi$ an isomorphism), by verifying that the unit map
$\id \rightarrow R \Un_{\phi}^{+} \circ \St_{\phi}^{+}$ is an isomorphism of functors
from $\h{(\mSet)_{/S}}$ to itself. 

Here is the idea of the proof.  Without loss of generality, we may suppose that $\phi$ is an isomorphism (since the pair $(\phi_{!}, \phi^{\ast})$ is a Quillen equivalence between $(\mSet)^{ \sCoNerve[S]^{op}}$ and $(\mSet)^{ \calC }$, by Proposition \ref{lesstrick}). We wish to show that $\Un^{+}_{\phi}$
induces an equivalence from the homotopy category of $(\mSet)^{\calC}$ to
the homotopy category of $(\mSet)_{/S}$. According to Proposition \ref{speccase}, this is true
whenever $S$ is a simplex. In the general case, we would like to regard $(\mSet)^{\calC}$
and $(\mSet)_{/S}$ as somehow built out of pieces which are associated to simplices, and deduce that $\Un^{+}_{\phi}$ is an equivalence because it is an equivalence on each piece. In order to make this argument work, it is necessary to work not just with the homotopy categories of
$(\mSet)^{\calC}$ and $(\mSet)_{/S}$, but with the simplicial categories which give rise to them.

We recall that both $(\mSet)^{\calC}$ and $(\mSet)_{/S}$ are {\em simplicial} model categories with respect to the simplicial mapping spaces defined by
$$ \Hom_{\sSet}(K, \bHom_{ (\mSet)^{\calC}}(\calF, \calG)) = \Hom_{(\mSet)^{\calC}}(
\calF \boxtimes K^{\sharp}, \calG)$$
$$ \Hom_{\sSet}(K, \bHom_{(\mSet)_S}(X,Y)) = \Hom_{\sSet}(K, \bHom_{S}^{\sharp}(X,Y))
= \Hom_{(\mSet)_{/S}}(X \times K^{\sharp}, Y).$$
The functor $\St^{+}_{\phi}$ is not a simplicial functor. However, it is {\em weakly} compatible with the simplicial structure in the sense that there is a natural map
$$ \St^{+}_{\phi} (X \boxtimes K^{\sharp}) \rightarrow (\St^{+}_{\phi} X) \boxtimes K^{\sharp}$$
for any $X \in (\mSet)_{/S}$, $K \in \sSet$ (according to Corollary \ref{spek5}, this map is a weak equivalence in $(\mSet)^{\calC}$). Passing to adjoints, we get natural maps
$$ \bHom_{ (\mSet)^{\calC} }( \calF, \calG) \rightarrow \bHom_{S}^{\sharp}( \Un^{+}_{\phi} \calF, \Un^{+}_{\phi} \calG ).$$ In other words, $\Un^{+}_{\phi}$ {\em does} have the structure of a simplicial functor.
We now invoke Proposition \ref{weakcompatequiv} to deduce the following:

\begin{lemma}\label{gottaprove0}
Let $S$ be a simplicial set, $\calC$ a simplicial category, and $\phi: \sCoNerve[S] \rightarrow \calC^{op}$ a simplicial functor. The following are equivalent:
\begin{itemize}
\item[$(1)$] The Quillen adjunction $( \St^{+}_{\phi}, \Un^{+}_{\phi} )$ is a Quillen equivalence.

\item[$(2)$] The functor $\Un^{+}_{\phi}$ induces an equivalence of simplicial categories
$$ (\Un^{+}_{\phi})^{\degree} : ((\mSet)^{\calC})^{\degree} \rightarrow ((\mSet)_{/S})^{\degree},$$
where $((\mSet)^{\calC})^{\degree}$ denotes the full (simplicial) subcategory of
$( ( \mSet)^{\calC})$ consisting of fibrant-cofibrant objects, and $((\mSet)_{/S})^{\degree}$ denotes the full (simplicial) subcategory of $(\mSet)_{/S}$ consisting of fibrant-cofibrant objects.
\end{itemize}
\end{lemma}

Consequently, to complete the proof of Theorem \ref{straightthm}, it will suffice to show that
if $\phi$ is an equivalence of simplicial categories, then $(\Un^{+}_{\phi})^{\degree}$ is an equivalence of simplicial categories. The first step is to prove that $(\Un^{+}_{\phi})^{\degree}$ is fully faithful.

\begin{lemma}\label{gottaprove2}
Let $S' \subseteq S$ be simplicial sets, and let
$p: X \rightarrow S$, $q: Y \rightarrow S$ be Cartesian fibrations.
Let $X' = X \times_{S} S'$ and $Y' = Y \times_{S} S'$. The restriction map
$$ \bHom_{S}^{\sharp}(X^{\natural}, Y^{\natural}) \rightarrow \bHom_{S'}^{\sharp}( {X'}^{\natural}, {Y'}^{\natural})$$ is a Kan fibration.
\end{lemma}

\begin{proof}
It suffices to show that the map $Y^{\natural} \rightarrow S$ has the right lifting property with
respect to the inclusion
$$ ( {X'}^{\natural} \times B^{\sharp} ) \coprod_{ {X'}^{\natural} \times A^{\sharp} } (X^{\natural} \times A^{\sharp}) \subseteq X^{\natural} \times B^{\sharp},$$
for any anodyne inclusion of simplicial sets $A \subseteq B$.

But this is a smash product of a marked cofibration ${X'}^{\natural} \rightarrow X^{\natural}$ 
(in $(\mSet)_{/S}$) and a trivial marked cofibration $A^{\sharp} \rightarrow B^{\sharp}$
( in $\mSet$), and is therefore a trivial marked cofibration. We conclude by observing that $Y^{\natural}$ is a fibrant object of $(\mSet)_{/S}$ (Proposition \ref{markedfibrant}).
\end{proof}

\begin{proof}[Proof of Theorem \ref{straightthm}]
For each simplicial set $S$, let $(\mSet)^{\sCoNerve[S]^{op}}_{f}$ denote the
category of projectively fibrant objects of $(\mSet)^{\sCoNerve[S]^{op}}$, and let
$W_{S}$ be the class of weak equivalences in $( \mSet)^{\sCoNerve[S]^{op}}_{f}$.
Let $W'_{S}$ be the collection of pointwise equivalences in $(\mSet)_{/S}^{\degree}$.
We have a commutative diagram of simplicial categories
$$ \xymatrix{ ((\mSet)^{\sCoNerve[S]^{op}})^{\degree} \ar[r]^{ \Un^{+}_{S} } \ar[d] &
(\mSet)_{/S}^{\degree} \ar[d]^{\psi_{S}} \\
( \mSet)^{\sCoNerve[S]^{op}}_{f}[W_{S}^{-1} ] \ar[r]^{\phi_S} & ( \mSet)_{/S}^{\degree}[ {W'}^{-1}_{S} ] }$$ (see Notation \ref{localdef}). In view of Lemma \ref{gottaprove0}, it will suffice to show
that the upper horizontal map is an equivalence of simplicial categories.
Lemma \ref{kur} implies that the left vertical map is an equivalence.
Using Lemma \ref{postcuse} and Remark \ref{uppa}, we deduce that the right vertical
map is also an equivalence. It will therefore suffice to show that $\phi_{S}$ is an equivalence. 

Let $\calU$ denote the collection of simplicial sets $S$ for which $\phi_{S}$ is an equivalence.
We will show that $\calU$ satisfies the hypotheses of Lemma \ref{blem}, and therefore contains every simplicial set $S$. Conditions $(i)$ and $(ii)$ are obviously satisfied, and
condition $(iii)$ follows from Lemma \ref{gottaprove0} and Proposition \ref{speccase}.
We will verify condition $(iv)$; the proof of $(v)$ is similar.

Applying Corollary \ref{uspin}, we deduce:
\begin{itemize}
\item[$(\ast)$] The functor $S \mapsto ( \mSet)^{\sCoNerve[S]^{op}}_{f} [W_S^{-1}]$ carries homotopy colimit diagrams indexed by a partially ordered set to homotopy limit diagrams in $\sCat$.
\end{itemize}

Suppose given a pushout diagram
$$ \xymatrix{ X \ar[r] \ar[d]^{f} & X' \ar[d] \\
Y \ar[r] & Y' }$$
in which $X, X', Y \in \calU$, where $f$ is a cofibration. We wish to prove that $Y' \in \calU$.
We have a commutative diagram
$$ \xymatrix{ ( \mSet)^{\sCoNerve[Y']^{op}}_{f}[W_{Y'}^{-1}] \ar[dr]^{\phi_{Y'}} & & \\
& (\mSet)_{/Y'}^{\degree}[ {W'}_{Y'}^{-1}] \ar[r]^{u} \ar[d]^{v} \ar[dr]^{w} & (\mSet)_{/Y}^{\degree}[ {W'}_{Y}^{-1}] \ar[d] \\
& (\mSet)_{/X'}^{\degree}[ {W'}_{X'}^{-1}] \ar[r] & (\mSet)_{/X}^{\degree}[ {W'}_{X}^{-1} ]. }$$
Using $(\ast)$ and Corollary \ref{wspin}, we deduce that $\phi_{Y'}$ is an equivalence if and only if,
for every pair of objects $x,y \in (\mSet)_{/Y'}^{\degree}[ {W'}_{Y'}^{-1}]$, the diagram
of simplicial sets
$$ \xymatrix{ \bHom_{ (\mSet)_{/Y'}^{\degree}[ {W'}_{Y'}^{-1}] }( x, y) \ar[r] \ar[d] 
& \bHom_{ (\mSet)_{/Y}^{\degree}[ {W'}_{Y}^{-1}] }( u(x), u(y) ) \ar[d] \\
\bHom_{ (\mSet)_{/X'}^{\degree}[ {W'}_{X'}^{-1}] }( v(x), v(y) ) \ar[r] &
\bHom_{ (\mSet)_{/X}^{\degree}[ {W'}_{X}^{-1}] }( w(x), w(y) ) }$$ 
is homotopy Cartesian. Since $\psi_{Y'}$ is a weak equivalence of simplicial categories, 
we may assume without loss of generality that $x = \psi_{Y'}( \overline{x} )$ and
$y = \psi_{Y'}( \overline{y} )$, for some $ \overline{x}, \overline{y} \in (\mSet)^{\degree}_{/Y'}$. 
It will therefore suffice to prove that the equivalent diagram
$$ \xymatrix{ \bHom^{\sharp}_{ Y' }( \overline{x}, \overline{y}) \ar[r] \ar[d] 
& \bHom^{\sharp}_{Y }( \overline{u}(\overline{x}), \overline{u}(\overline{y}) ) \ar[d] \\
\bHom^{\sharp}_{ X' }( \overline{v}(\overline{x}), \overline{v}(\overline{y}) ) \ar[r]^{g} &
\bHom^{\sharp}_{X }( \overline{w}(\overline{x}), \overline{w}(\overline{y}) ) }$$ 
is homotopy Cartesian. But this diagram is a pullback square, and the map $g$ is a Kan
fibration by Lemma \ref{gottaprove2}.
\end{proof}

\subsection{The Relative Nerve}\label{altstr}

In \S \ref{markmodel}, we defined {\em straightening} and {\em unstraightening} functors, which
give rise to a Quillen equivalence of model categories
$$ \Adjoint{\St^{+}_{\phi}}{(\mSet)_{/S}}{(\mSet)^{\calC}}{\Un^{+}_{\phi}}$$
whenever $\phi: \sCoNerve[S] \rightarrow \calC^{op}$ is a weak equivalence of simplicial categories.
For many purposes, these constructions are unnecessarily complicated. For example, suppose that
$\calF: \calC \rightarrow \mSet$ is a (weakly) fibrant diagram, so that $\Un^{+}_{\phi}(\calF)$ is a
fibrant object of $(\mSet)_{/S}$ corresponding to a Cartesian fibration of simplicial sets $X \rightarrow S$. For every vertex $s \in S$, the fiber $X_{s}$ is an $\infty$-category which is equivalent
to $\calF( \phi(s) )$, but usually not isomorphic to $\calF( \phi(s) )$. In the special case where
$\calC$ is an ordinary category and $\phi: \sCoNerve[ \Nerve(\calC)^{op} ] \rightarrow \calC^{op}$ is the counit map, there is another version of unstraightening construction $\Un^{+}_{\phi}$ which does not share this defect. Our goal in this section is to introduce this simpler construction, which we call
the {\it marked relative nerve} $\calF \mapsto \Nerve^{+}_{\calF}(\calC)$, and to study its basic properties.

\begin{remark}
To simplify the exposition which follows, the relative nerve functor introduced below will actually be
an alternative to the {\em opposite} of the unstraightening functor
$$ \calF \mapsto ( \Un^{+}_{\phi} \calF^{op} )^{op},$$
which is a right Quillen functor from the projective model structure on
$(\mSet)^{\calC}$ to the coCartesian model structure on $\mset{ \Nerve(\calC)}$. 
\end{remark}


\begin{definition}\label{sulke}\index{gen}{relative nerve}\index{gen}{nerve!relative}\index{not}{NervefI@$\Nerve^{R}_{\calF}(\calC)$}
Let $\calC$ be a small category, and let $f: \calC \rightarrow \sSet$ be a functor. We define a simplicial set $\Nerve_{f}(\calC)$, the {\it nerve of $\calC$ relative to $f$}, as follows. For every nonempty finite linearly ordered set $J$, a map $\Delta^J \rightarrow \Nerve_{f}(\calC)$ consists of the following data:
\begin{itemize}
\item[$(1)$] A functor $\sigma$ from $J$ to $\calC$. 
\item[$(2)$] For every nonempty subset $J' \subseteq J$ having a maximal element $j'$, a
map $\tau(J'): \Delta^{J'} \rightarrow \calF( \sigma(j'))$. 
\item[$(3)$] For nonempty subsets $J'' \subseteq J' \subseteq J$, with maximal elements $j'' \in J''$, $j' \in J'$, the diagram $$ \xymatrix{ \Delta^{J''} \ar[r]^{\tau(J'')} \ar@{^{(}->}[d] & f( \sigma(j'') ) \ar[d] \\
\Delta^{J'} \ar[r]^{\tau(J')} & f( \sigma(j') ) }$$
is required to commute.
\end{itemize}
\end{definition}

\begin{remark}
Let $\calI$ be denote the linearly ordered set $[n]$, regarded as a category, and let
$f: \calI \rightarrow \sSet$ correspond to a composable sequence of morphisms
$\phi: X_0 \rightarrow \ldots \rightarrow X_n$.
Then $\Nerve_{f}(\calI)$ is closely related to the mapping simplex $M^{op}(\phi)$ introduced in \S \ref{funkystructure}. More precisely, there is a canonical map $\Nerve_{f}(\calI) \rightarrow M^{op}(\phi)$ compatible with the projection to $\Delta^{n}$, which induces an isomorphism on each fiber.
\end{remark}

\begin{remark}\label{staplur}
The simplicial set $\Nerve_{f}(\calC)$ of Definition \ref{sulke} depends functorially on $f$. When
$f$ takes the constant value $\Delta^0$, there is a canonical isomorphism
$\Nerve_{f}(\calC) \simeq \Nerve(\calC)$. In particular, for {\em any} functor $f$, there is a canonical map $\Nerve_{f}(\calC) \rightarrow \Nerve(\calC)$; the fiber of this map over an object $C \in \calC$ can be identified with the simpicial set $f(C)$. 
\end{remark}

\begin{remark}\label{pear1}\index{not}{FF@$\lNerve_{X}(\calC)$}
Let $\calC$ be a small $\infty$-category. The construction $f \mapsto \Nerve_{f}(\calC)$ determines a functor from $(\sSet)^{\calC}$ to $(\sSet)_{/ \Nerve(\calC)}$. This functor admits a left adjoint, which we will denote by $X \mapsto \lNerve_X(\calC)$ (the existence of this functor follows from the adjoint functor theorem). If $X \rightarrow \Nerve(\calC)$ is a left fibration, then $\calF_{X}(\calC)$ is a functor
$\calC \rightarrow \sSet$ which assigns to each $C \in \calC$ a simplicial set which is weakly equivalent to the fiber $X_{C} = X \times_{ \Nerve(\calC)} \{C\}$; this follows from Proposition \ref{kudd} below. 
\end{remark}

\begin{example}\label{spacek}
Let $\calC$ be a small category, and regard $\Nerve(\calC)$ as an object of
$(\sSet)_{/ \Nerve(\calC)}$ via the identity map. Then $\lNerve_{\Nerve(\calC)}(\calC) \in (\sSet)^{\calC}$ can be identified with the functor $C \mapsto \Nerve(\calC_{/C})$. 
\end{example}

\begin{remark}\label{pear2}
Let $g: \calC \rightarrow \calD$ be a functor between small categories and let $f: \calD \rightarrow \sSet$ be a diagram. There is a canonical isomorphism of simplicial sets
$\Nerve_{f \circ g}(\calC) \simeq \Nerve_{f}(\calD) \times_{ \Nerve(\calD) } \Nerve(\calC).$
In other words, the diagram of categories
$$ \xymatrix{ (\sSet)^{\calD} \ar[r]^{g^{\ast}} \ar[d]^{\Nerve_{\bigdot}(\calD)} & (\sSet)^{\calC} \ar[d]^{ \Nerve_{\bigdot}(\calC) } \\
(\sSet)_{/ \Nerve(\calD)} \ar[r]^{ \Nerve(g)^{\ast}} & (\sSet)_{/ \Nerve(\calC)} }$$
commutes up to canonical isomorphism. Here $g^{\ast}$ denotes the functor given by composition with
$g$, and $\Nerve(g)^{\ast}$ the functor given by pullback along the map of simplicial sets
$\Nerve(g): \Nerve(\calC) \rightarrow \Nerve(\calD)$.
\end{remark}

\begin{remark}\label{saltine}
Combining Remarks \ref{pear1} and \ref{pear2}, we deduce that for any functor
$g: \calC \rightarrow \calD$ between small categories, the diagram of left adjoints
$$ \xymatrix{ (\sSet)^{\calD} & (\sSet)^{\calC} \ar[l]_{g_{!} } \\
(\sSet)_{/ \Nerve(\calD)} \ar[u]_{ \lNerve_{\bigdot}(\calD)} & (\sSet)_{/\Nerve(\calC)} \ar[u]_{\lNerve_{\bigdot}(\calC)} \ar[l] }$$
commutes up to canonical isomorphism; here $g_{!}$ denotes the functor of left Kan extension along $g$, and the bottom arrow is the forgetful functor given by composition with
$\Nerve(g): \Nerve(\calC) \rightarrow \Nerve(\calD)$.
\end{remark}

\begin{notation}
Let $\calC$ be a small category, and let $f: \calC \rightarrow \sSet$ be a functor. We let
$f^{op}$ denote the functor $\calC \rightarrow \sSet$ described by the formula
$f^{op}(C) = f(C)^{op}$. We will use a similar notation in the case where
$f$ is a functor from $\calC$ to the category $\mSet$ of marked simplicial sets.
\end{notation}

\begin{remark}\label{scuz}
Let $\calC$ be a small category, let $S = \Nerve(\calC)^{op}$, and let $\phi: \sCoNerve[S] \rightarrow \calC^{op}$ be the counit map. For each $X \in (\sSet)_{/ \Nerve(\calC)}$, there is a canonical map
$$\alpha_{\calC}(X): \St_{\phi} X^{op} \rightarrow \lNerve_{X}(\calC)^{op}.$$
The collection of maps $\{ \alpha_{\calC}(X) \}$ is uniquely determined by the following requirements:

\begin{itemize}
\item[$(1)$] The morphism $\alpha_{\calC}(X)$ depends functorially on $X$. More precisely, suppose
given a commutative diagram of simplicial sets
$$ \xymatrix{ X \ar[rr]^{f} \ar[dr] & & Y \ar[dl] \\
& \Nerve(\calC). & }$$
Then the diagram
$$ \xymatrix{ \St_{\phi} X^{op} \ar[r]^{\alpha_{\calC}(X)} \ar[d]^{ \St_{\phi} f^{op}} & \lNerve_{X}(\calC)^{op} \ar[d]^{ \lNerve_{f}(\calC)^{op}} \\
\St_{\phi} Y^{op} \ar[r]^{ \alpha_{\calC}(Y) } & \lNerve_{Y}(\calC)^{op} }$$
commutes. 

\item[$(2)$] The transformation $\alpha_{\calC}$ depends functorially on $\calC$ in the following sense: for every functor $g: \calC \rightarrow \calD$,
if $\phi': \sCoNerve[(\Nerve \calD)^{op}] \rightarrow \calD^{op}$ denotes the counit map and
$X \in (\sSet)_{/ \Nerve(\calC)}$, then the diagram
$$ \xymatrix{ \St_{\phi} X^{op} \ar[d] \ar[r]^{g_! \alpha_{\calC}} & g_! \lNerve_{X}(\calC)^{op} \ar[d] \\
\St_{\phi'} X^{op} \ar[r]^{ \alpha_{\calD}} & \lNerve_{X}(\calD)^{op} }$$
commutes, where the vertical arrows are the isomorphisms provided by Remark \ref{saltine} and
Proposition \ref{straightchange}.

\item[$(3)$] Let $\calC$ be the category associated to a partially ordered set $P$, and let
$X = \Nerve(\calC)$, regarded as an object of
$(\sSet)_{/ \Nerve(\calC)}$ via the identity map. Then $(\St_{\phi} X^{op}) \in (\sSet)^{\calC}$
can be identified with the functor $p \mapsto \Nerve X_{p}$, where for each $p \in P$ we let
$X_{p}$ denote the collection of nonempty finite chains in $P$ having largest element $p$. Similarly, Example \ref{spacek} allows us to identify $\lNerve_{X}(\calC) \in (\sSet)^{\calC}$ with the functor $p \mapsto \Nerve \{ q \in P: q \leq p \}$. The map
$\alpha_{\calC}(X): (\St_{\phi} X^{op}) \rightarrow \lNerve_{X}(\calC)^{op}$ is induced by the map
of partially ordered sets $X_p \rightarrow \{ q \in P: q \leq p \}$ which carries every chain to its
smallest element.
\end{itemize}

To see that the collection of maps $\{ \alpha_{\calC}(X) \}_{ X \in (\sSet)_{/ \Nerve(\calC)}}$ is determined by these properties, we first note that because the functors $\St_{\phi}$ and $\lNerve_{\bigdot}(\calC)$ commute with colimits, any natural transformation $\beta_{\calC}: \St_{\phi}(\bigdot^{op}) \rightarrow \lNerve_{\bigdot}(\calC)^{op}$ is determined by its values $\beta_{\calC}(X): \St_{\phi}(X^{op}) \rightarrow \lNerve_{X}(\calC)^{op}$ in the case where $X = \Delta^n$ is a simplex. In this case, any map $X \rightarrow \Nerve \calC$ factors through the isomorphism $X \simeq \Nerve [n]$, so we can use property $(2)$ to reduce to the case where the category $\calC$ is a partially ordered set and the map $X \rightarrow \Nerve(\calC)$ is an isomorphism. The behavior of the natural transformation $\alpha_{\calC}$ is then dictated by property $(3)$. This proves the uniqueness of the natural transformations $\alpha_{\calC}$; the existence follows by a similar argument.
\end{remark}

The following result summarizes some of the basic properties of the relative nerve functor:

\begin{lemma}\label{sulken2}
Let $\calI$ be a category and
let $\alpha: f \rightarrow f'$ be a natural transformation of functors
$f,f': \calC \rightarrow \sSet$.
\begin{itemize}
\item[$(1)$] Suppose that, for each $I \in \calC$, the map $\alpha(I): f(I) \rightarrow f'(I)$ is
an inner fibration of simplicial sets. Then the induced map $\Nerve_{f}(\calC) \rightarrow
\Nerve_{f'}(\calC)$ is an inner fibration.
\item[$(2)$] Suppose that, for each $I \in \calI$, the simplicial set $f(I)$ is an $\infty$-category. Then $\Nerve_{f}(\calC)$ is an $\infty$-category.
\item[$(3)$] Suppose that, for each $I \in \calC$, the map $\alpha(I): f(I) \rightarrow f'(I)$ is
a categorical fibration of $\infty$-categories. Then the induced map
$\Nerve_{f}(\calC) \rightarrow \Nerve_{f'}(\calC)$ is a categorical fibration of $\infty$-categories.
\end{itemize}
\end{lemma}

\begin{proof}
Consider a commutative diagram
$$ \xymatrix{ \Lambda^n_i \ar[r] \ar@{^{(}->}[d] & \Nerve_{f}(\calI) \ar[d]^{p} \\
\Delta^n \ar@{-->}[ur] \ar[r] & \Nerve_{f'}(\calC), }$$
and let $I$ be the image of $\{n\} \subseteq \Delta^n$ under the
bottom map. If $0 \leq i < n$, then the lifting problem depicted in the diagram above is equivalent to the existence of a dotted arrow in an associated diagram
$$ \xymatrix{ \Lambda^n_i \ar[r]^{g} \ar@{^{(}->}[d] & f(I) \ar[d]^{\alpha(I)} \\
\Delta^n \ar@{-->}[ur] \ar[r] & f'(I).}$$
If $\alpha(I)$ is an inner fibration and $0 < i < n$, then we conclude that this lifting problem admits a solution. This proves $(1)$. To prove $(2)$, we apply $(1)$ in the special case where $f'$ is the constant functor taking the value $\Delta^0$. It follows that $\Nerve_{f}(\calC) \rightarrow \Nerve(\calC)$ is an inner fibration, so that $\Nerve_{f}(\calC)$ is an $\infty$-category.

We now prove $(3)$. According to Corollary \ref{gottaput}, an inner fibration
$\calD \rightarrow \calE$ of $\infty$-categories is a categorical fibration if and only if the following condition is satisfied:
\begin{itemize}
\item[$(\ast)$] For every equivalence $e: E \rightarrow E'$ in $\calE$, and every
object $D \in \calD$ lifting $E$, there exists an equivalence $\overline{e}: D \rightarrow D'$
in $\calD$ lifting $e$.
\end{itemize}

We can identify equivalences in $\Nerve_{f'}(\calC)$ with triples
$(g: I \rightarrow I', X, e: X' \rightarrow Y)$ where $g$ is an isomorphism in $\calC$, $X$ is an object
of $f'(I)$, $X'$ is the image of $X$ in $f'(I')$, and $e: X' \rightarrow Y$ is an equivalence in
$f'(I')$. Given a lifting $\overline{X}$ of $X$ to $f(I)$, we can apply the assumption that
$\alpha(I')$ is a categorical fibration (and Corollary \ref{gottaput}) to lift $e$ to an equivalence $\overline{e}: \overline{X}' \rightarrow \overline{Y}$ in $f(I')$. This produces the desired equivalence
$(g: I \rightarrow I', \overline{X}, \overline{e}: \overline{X}' \rightarrow \overline{Y})$ in
$\Nerve_{f}(\calC)$.
\end{proof}

We now introduce a slightly more elaborate version of the relative nerve construction.

\begin{definition}\index{not}{Nervefplus@$\Nerve^{+}_{\calF}(\calC)$}
Let $\calC$ be a small category and $\calF: \calC \rightarrow \mSet$ a functor. We let
$\Nerve^{+}_{\calF}(\calC)$ denote the marked simplicial set
$( \Nerve_{f}(\calC), M)$, where $f$ denotes the composition
$\calC \stackrel{\calF}{\rightarrow} \mSet \rightarrow \sSet$ and
$M$ denotes the collection of all edges $\overline{e}$ of $\Nerve_{f}(\calC)$ with the following property:
if $e: C \rightarrow C'$ is the image of $\overline{e}$ in $\Nerve(\calC)$, and $\sigma$ denotes the
edge of $f(C')$ determined by $\overline{e}$, then $\sigma$ is a marked edge of $\calF(C')$.
We will refer to $\Nerve^{+}_{\calF}(\calC)$ as the {\it marked relative nerve} functor.\index{gen}{relative nerve!marked}\index{gen}{nerve!marked relative}\index{gen}{marked relative nerve}
\end{definition}

\begin{remark}
Let $\calC$ be a small category. We will regard the construction $\calF \mapsto \Nerve^{+}_{\calF}(\calC)$ as determining a functor from $(\mSet)^{\calC}$ to $\mset{\Nerve(\calC) }$ (see Remark \ref{staplur}). This functor admits a left adjoint, which we will denote by $\overline{X} \mapsto \lNerve^{+}_{\overline{X}}(\calC)$.\index{not}{FFplus@$\lNerve^{+}_{\overline{X}}(\calC)$}
\end{remark}

\begin{remark}\label{saltine2}
Remark \ref{saltine} has an evident analogue for the functors $\lNerve^{+}$: 
for any functor $g: \calC \rightarrow \calD$ between small categories, the diagram of left adjoints
$$ \xymatrix{ (\mSet)^{\calD} & (\mSet)^{\calC} \ar[l]_{g_{!} } \\
\mset{\Nerve(\calD)} \ar[u]_{ \lNerve^{+}_{\bigdot}(\calD)} & \mset{\calC} \ar[u]_{\lNerve^{+}_{\bigdot}(\calC)} \ar[l] }$$
commutes up to canonical isomorphism.
\end{remark}

\begin{lemma}\label{coughup}
Let $\calC$ be a small category. Then:
\begin{itemize}
\item[$(1)$] The functor $X \mapsto \lNerve_{X}(\calC)$ carries
cofibrations in $(\sSet)_{/ \Nerve(\calC)}$
to cofibrations in $(\sSet)^{\calC}$ $($with respect to the projective model structure$)$.
\item[$(2)$] The functor $\overline{X} \mapsto \lNerve^{+}_{\overline{X}}(\calC)$ carries
cofibrations $($with respect to the coCartesian monoidal structure on  $\mset{\Nerve(\calC)}${}$)$ to cofibrations in $(\mSet)^{\calC}$ $($with respect to the projective model structure$)$.
\end{itemize}
\end{lemma}

\begin{proof}
We will give the proof of $(2)$; the proof of $(1)$ is similar. It will suffice to show that the right adjoint functor $\Nerve^{+}_{\bigdot}(\calC): (\mSet)^{\calC} \rightarrow \mSet{ \Nerve(\calC) }$ preserves trivial fibrations. Let $\calF \rightarrow \calF'$ be a trivial fibration in $(\mSet)^{\calC}$ with respect to the projective model structure, so that for each $C \in \calC$ the induced map $\calF(C) \rightarrow \calF'(C)$ is a trivial fibration of marked simplicial sets. We wish to prove that the induced map
$\Nerve^{+}_{\calF}(\calC) \rightarrow \Nerve^{+}_{\calF'}(\calC)$ is also a trivial fibration of marked simplicial sets. Let $f$ denote the composition $\calC \stackrel{\calF}{\rightarrow} \mSet \rightarrow \sSet$, and let $f'$ be defined likewise. We must verify two things:
\begin{itemize}
\item[$(1)$] Every lifting problem of the form
$$ \xymatrix{ \bd \Delta^n \ar[r] \ar@{^{(}->}[d] & \Nerve_{f}(\calC) \ar[d] \\
\Delta^n \ar[r]^{u} & \Nerve_{f'}(\calC) }$$
admits a solution. Let $C \in \calC$ denote the image of the final vertex of $\Delta^n$ under the map $u$.
Then it suffices to solve a lifting problem of the form
$$ \xymatrix{ \bd \Delta^n \ar[r] \ar@{^{(}->}[d] & f(C) \ar[d] \\
\Delta^n \ar[r] & f'(C), }$$
which is possible since the right vertical map is a trivial fibration of simplicial sets.

\item[$(2)$] If $\overline{e}$ is an edge of $\Nerve^{+}_{\calF}(\calC)$ whose image $\overline{e}'$ in
$\Nerve^{+}_{\calF'}(\calC)$ is marked, then $\overline{e}$ is itself marked. Let
$e: C \rightarrow C'$ be the image of $\overline{e}$ in $\Nerve(\calC)$, and let
$\sigma$ denote the edge of $\calF(C')$ determined by $\overline{e}$. Since
$\overline{e}'$ is a marked edge of $\Nerve^{+}_{\calF'}(\calC)$, the image of
$\sigma$ in $\calF'(C')$ is marked. Since the map $\calF(C') \rightarrow \calF'(C')$ is a trivial fibration of marked simplicial sets, we deduce that $\sigma$ is a marked edge of $\calF(C')$, so that
$\overline{e}$ is a marked edge of $\Nerve^{+}_{\calF}(\calC)$ as desired.
\end{itemize}
\end{proof}

\begin{remark}\label{scuz2}
Let $\calC$ be a small category, let $S = \Nerve(\calC)^{op}$, and let $\phi: \sCoNerve[S] \rightarrow \calC^{op}$ be the counit map. For every $\overline{X} = (X,M) \in \mset{\Nerve(\calC)}$, the 
morphism $\alpha_{\calC}(X): \St_{\phi}(X^{op}) \rightarrow \lNerve_{X}(\calC)^{op}$ of Remark \ref{scuz} induces a natural transformation $\St^{+}_{\phi} \overline{X}^{op} \rightarrow \lNerve^{+}_{ \overline{X}}(\calC)^{op}$, which we will denote by $\alpha^{+}_{\calC}( \overline{X})$. We will regard the collection
of morphisms $\{ \alpha^{+}_{\calC}( \overline{X}) \}_{ \overline{X} \in \mset{\Nerve(\calC)}}$ as determining a natural transformation of functors 
$$\alpha_{\calC}: \St^{+}_{\phi}( \bigdot^{op}) \rightarrow \lNerve^{+}_{\bigdot}(\calC)^{op}.$$
\end{remark}

\begin{lemma}\label{standrum}
Let $\calC$ be small category, let $S = \Nerve(\calC)^{op}$, and let $\phi: \sCoNerve[S] \rightarrow \calC^{op}$ be the counit map, and let $C \in \calC$ be an object. Then:
\begin{itemize}
\item[$(1)$] For every $X \in (\sSet)_{/ \Nerve(\calC)}$, the map 
$\alpha_{\calC}(X): \St_{\phi}(X^{op}) \rightarrow \lNerve_{X}(\calC)^{op}$ of Remark \ref{scuz} induces
a weak homotopy equivalence of simplicial sets $\St_{\phi}(X^{op})(C) \rightarrow \lNerve_{X}(\calC)(C)^{op}$.
\item[$(2)$] For every $\overline{X} \in \mset{\Nerve(\calC)}$, the map
$\alpha^{+}_{\calC}( \overline{X}): \St^{+}_{\phi}( \overline{X}^{op}) \rightarrow \lNerve^{+}_{\overline{X}}(\calC)^{op}$ of Remark \ref{scuz2} induces a Cartesian equivalence $\St^{+}_{\phi}( \overline{X}^{op})(C) \rightarrow \lNerve^{+}_{\overline{X}}(\calC)(C)^{op}$. 
\end{itemize}
\end{lemma}

\begin{proof}
We will give the proof of $(2)$; the proof of $(1)$ is similar but easier. Let us say that an object
$\overline{X} \in \mset{\Nerve(\calC)}$ is {\it good} if the map $\alpha^{+}_{\calC}( \overline{X})$ is a weak equivalence. We wish to prove that every object $\overline{X} = (X,M) \in \mset{ \Nerve(\calC)}$ is good. The proof proceeds in several steps.
\begin{itemize}
\item[$(A)$] Since the functors $\St^{+}_{\phi}$ and $\lNerve^{+}_{\bigdot}(\calC)$ both commute with filtered colimits, the collection of good objects of $\mset{ \Nerve(\calC)}$ is stable under filtered colimits.
We may therefore reduce to the case where the simplicial set $X$ has only finitely many nondegenerate simplices.
\item[$(B)$] Suppose given a pushout diagram
$$ \xymatrix{ \overline{X} \ar[r]^{f} \ar[d]^{g} & \overline{X}' \ar[d] \\
\overline{Y} \ar[r] & \overline{Y}' }$$
in the category $\mset{ \Nerve(\calC) }$. Suppose that either $f$ or $g$ is a cofibration, and that
the objects $\overline{X}, \overline{X}'$, and $\overline{Y}$ are good. Then $\overline{Y}'$ is good.
This follows from the fact that the functors $\St^{+}_{\phi}$ and $\lNerve^{+}_{\bigdot}(\calC)$ preserve cofibrations (Proposition \ref{cougherup} and Lemma \ref{coughup}), together with the observation that the projective model structure on $(\mSet)^{\calC}$ is left proper. 
\item[$(C)$] Suppose that $X \simeq \Delta^n$ for $n \leq 1$. In this case,
the map $\alpha^{+}_{\calC}( \overline{X})$ is an isomorphism (by direct calculation), so that $\overline{X}$ is good.
\item[$(D)$] We now work by induction on the number of nondegenerate marked edges of $\overline{X}$.
If this number is nonzero, then there exists a pushout diagram
$$ \xymatrix{ (\Delta^1)^{\flat} \ar[r] \ar[d] & (\Delta^1)^{\sharp} \ar[d] \\
\overline{Y} \ar[r] & \overline{X} }$$
where $\overline{Y}$ has fewer nondegenerate marked edges than $\overline{X}$, so that
$\overline{Y}$ is good by the inductive hypothesis. The marked simplicial sets
$(\Delta^1)^{\flat}$ and $(\Delta^1)^{\sharp}$ are good by virtue of $(C)$, so that $(B)$ implies that
$\overline{X}$ is good. We may therefore reduce to the case where $\overline{X}$ contains no nondegenerate marked edges, so that $\overline{X} \simeq X^{\flat}$.
\item[$(E)$] We now argue by induction on the dimension $n$ of $X$ and the number of nondegenerate $n$-simplices of $X$. If $X$ is empty, there is nothing to prove; otherwise, we have a pushout diagram
$$ \xymatrix{ \bd \Delta^n \ar[r] \ar[d] & \Delta^n \ar[d] \\
Y \ar[r] & X. }$$
The inductive hypothesis implies that $( \bd \Delta^n)^{\flat}$ and $Y^{\flat}$ are good.
Invoking step $(B)$, we can reduce to the case where $X$ is an $n$-simplex. In view of
$(C)$, we may assume that $n \geq 2$.

Let $Z = \Delta^{ \{0,1\} } \coprod_{ \{ 1\} } \Delta^{ \{1,2\} } \coprod_{ \{1\} } \ldots
\coprod_{ \{n-1\}} \Delta^{ \{n-1, n \}}$ so that $Z \subseteq X$ is an inner anodyne inclusion.
We have a commutative diagram
$$ \xymatrix{ \St^{+}_{\phi} (Z^{op})^{\flat} \ar[r]^{u} \ar[d]^{v} & \St^{+}_{\phi}(X^{op})^{\flat} \ar[d] \\
\lNerve^{+}_{ Z^{\flat}}(\calC)^{op} \ar[r]^{w} & \lNerve^{+}_{ X^{\flat}}(\calC)^{op}. }$$
The inductive hypothesis implies that $v$ is a weak equivalence, and Proposition \ref{spec2} implies that $u$ is a weak equivalence. To complete the proof, it will suffice to show that
$w$ is a weak equivalence. 

\item[$(F)$] The map $X \rightarrow \Nerve(\calC)$ factors as a composition
$$ \Delta^n \simeq \Nerve( [n] ) \stackrel{g}{\rightarrow} \Nerve(\calC).$$
Using Remark \ref{saltine2} (together with the fact that the left
Kan extension functor $g_!$ preserves weak equivalences between projectively cofibrant objects),
we can reduce to the case where $\calC = [n]$ and the map $X \rightarrow \Nerve(\calC)$ is an isomorphism.

\item[$(G)$] Fix an object $i \in [n]$. A direct computation shows that the map
$\lNerve^{+}_{Z^{\flat}}(\calC)(i) \rightarrow \lNerve^{+}_{X^{\flat}}(\calC)(i)$ can be identified with the inclusion 
$$ (  \Delta^{ \{0,1\} } \coprod_{ \{ 1\} } \Delta^{ \{1,2\} } \coprod_{ \{1\} } \ldots
\coprod_{ \{i-1\}} \Delta^{ \{i-1, i \}})^{op, \flat} \subseteq (\Delta^{i})^{op, \flat}.$$
This inclusion is marked anodyne, and therefore an equivalence of marked simplicial sets as desired.
\end{itemize}
\end{proof}

\begin{proposition}\label{kudd}
Let $\calC$ be a small category. Then:
\begin{itemize}
\item[$(1)$] The functors $\lNerve_{\bigdot}(\calC)$ and $\Nerve_{\bigdot}(\calC)$ determine
a Quillen equivalence between $( \sSet)_{/ \Nerve(\calC)}$ $($endowed with the covariant model structure$)$
and $(\sSet)^{\calC}$ $($endowed with the projective model structure$)$.
\item[$(2)$] The functors $\lNerve^{+}_{\bigdot}(\calC)$ and $\Nerve^{+}_{\bigdot}(\calC)$ determine
a Quillen equivalence between $\mset{ \Nerve(\calC)}$ $($endowed with the coCartesian model structure) and $(\mSet)^{\calC}$ $($endowed with the projective model structure$)$.
\end{itemize}
\end{proposition}

\begin{proof}
We will give the proof of $(2)$; the proof of $(1)$ is similar but easier. We first show that
the adjoint pair $( \lNerve^{+}_{\bigdot}(\calC), \Nerve^{+}_{\bigdot}(\calC) )$ is a Quillen adjunction.
It will suffice to show that the functor $\lNerve^{+}_{\bigdot}(\calC)$ preserves cofibrations and weak equivalences. The case of cofibrations follows from Lemma \ref{coughup}, and the case of weak
equivalences from Lemma \ref{standrum} and Corollary \ref{spek6}. To prove that
$( \lNerve^{+}_{\bigdot}(\calC), \Nerve^{+}_{\bigdot}(\calC) )$ is a Quillen equivalence, it will suffice to show that the left derived functor $L \lNerve^{+}_{\bigdot}(\calC)$ induces an equivalence from the homotopy category $\h{ \mset{\Nerve(\calC)}}$ to the homotopy category $\h{ (\mSet)^{\calC}}$. In view of Lemma \ref{standrum}, it will suffice to prove an analogous result for the straightening functor
$\St^{+}_{\phi}$, where $\phi$ denotes the counit map
$\sCoNerve[ \Nerve(\calC)^{op}] \rightarrow \calC^{op}$. We now invoke Theorem \ref{straightthm}.
\end{proof}

\begin{corollary}
Let $\calC$ be a small category, and let $\alpha: f \rightarrow f'$ be a natural transformation of functors
$f,f': \calC \rightarrow \sSet$. Suppose that, for each $C \in \calC$, the induced map
$f(C) \rightarrow f'(C)$ is a Kan fibration. Then the induced map
$\Nerve_{f}(\calC) \rightarrow \Nerve_{f'}(\calC)$ is a covariant fibration in
$(\sSet)_{/ \Nerve(\calC)}$. In particular, if each $f(C)$ is Kan complex, then the map
$\Nerve_{f}(\calC) \rightarrow \Nerve(\calC)$ is a left fibration of simplicial sets.
\end{corollary}

\begin{corollary}\label{sandcor}
Let $\calC$ be a small category and $\calF: \calC \rightarrow \mSet$ a fibrant object of
$(\mSet)^{\calC}$. Let $S = \Nerve(\calC)$, and let $\phi: \sCoNerve[S^{op}] \rightarrow \calC^{op}$ denote the counit map. Then the natural transformation $\alpha_{\calC}^{+}$ of Remark \ref{scuz2}
induces a weak equivalence
$\Nerve_{\calF}(\calC)^{op} \rightarrow (\Un^{+}_{\phi} \calF^{op})$
$($with respect to the Cartesian model structure on $\mset{ S^{op}}${}$)$.
\end{corollary}

\begin{proof}
It suffices to show that $\alpha_{\calC}^{+}$ induces an isomorphism of right derived functors
$R \Nerve_{\bigdot}(\calC)^{op} \rightarrow R (\Un^{+}_{\phi} \bigdot^{op})$, which follows immediately from Lemma \ref{standrum}. 
\end{proof}

\begin{proposition}\label{sulken}
Let $\calC$ be a category, and let $f: \calC \rightarrow \sSet$ be a functor such that
$f(C)$ is an $\infty$-category for each $C \in \calC$. Then:
\begin{itemize}
\item[$(1)$] The projection $p: \Nerve_{f}(\calC) \rightarrow \Nerve(\calC)$ is a coCartesian fibration of simplicial sets. 

\item[$(2)$] Let $e$ be an edge of $\Nerve_{f}(\calC)$, covering a morphism $C \rightarrow C'$
in $\calC$. Then $e$ is $p$-coCartesian if and only if the corresponding edge of $f(C')$ is an equivalence.

\item[$(3)$] The coCartesian fibration $p$ is associated to the functor $\Nerve(f): \Nerve(\calC) \rightarrow \Cat_{\infty}$ $($see \S \ref{universalfib}$)$. 
\end{itemize}
\end{proposition}

\begin{proof}
Let $\calF: \calC \rightarrow \mSet$ be the functor described by the formula
$\calF(C) = f(C)^{\natural}$. Then $\calF$ is a projectively fibrant object of $(\mSet)^{\calC}$.
Invoking Proposition \ref{kudd}, we deduce that $\Nerve^{+}_{\calF}(\calC)$ is a fibrant object
of $\mset{ \Nerve(\calC)}$. Invoking Proposition \ref{markedfibrant}, we deduce that the underlying
map $p: \Nerve_{f}(\calC) \rightarrow \Nerve(\calC)$ is a coCartesian fibration of simplicial sets, and that
the $p$-coCartesian morphisms of $\Nerve_{f}(\calC)$ are precisely the marked wedges of
$\Nerve^{+}_{\calF}(\calC)$. This proves $(1)$ and $(2)$. To prove $(3)$, we let
$S = \Nerve(\calC)$ and $\phi: \sCoNerve[S]^{op} \rightarrow \calC^{op}$ be the counit map.
By definition, a coCartesian fibration $X \rightarrow \Nerve(\calC)$ is associated to $f$ if and only if
it is equivalent to $(\Un_{\phi} f^{op})^{op}$; the desired equivalence is furnished by Corollary \ref{sandcor}.

\end{proof}

\section{Applications}\label{hugr}
\setcounter{theorem}{0}

The purpose of this section is to survey some applications of technology developed in
\S \ref{twuf} and \S \ref{strsec}. 
In \S \ref{apps}, we give some applications to the theory of Cartesian fibrations. In \S \ref{universalfib}, we will introduce the language of {\it classifying maps} which will allow us to exploit the Quillen equivalence provided by Theorem \ref{straightthm}. Finally, in \S \ref{catlim} and \S \ref{catcolim}, we will use Theorem \ref{straightthm} to give explicit constructions of limits and colimits in the $\infty$-category $\Cat_{\infty}$ (and also in the $\infty$-category $\SSet$ of spaces).

\subsection{Structure Theory for Cartesian Fibrations}\label{apps}

The purpose of this section is to prove that Cartesian fibrations between simplicial sets enjoy several pleasant properties. For example, every Cartesian fibration is a categorical fibration (Proposition \ref{funkyfibcatfib}), and categorical equivalences are stable under pullbacks by Cartesian fibrations (Proposition \ref{basechangefunky}). These results are fairly easy to prove for Cartesian fibrations $X \rightarrow S$ in the case where $S$ is an $\infty$-category. Theorem \ref{straightthm} provides a method for reducing to this special case:

\begin{proposition}\label{pretrokee}
Let $p: S \rightarrow T$ be a categorical equivalence of simplicial sets. Then the forgetful functor $$p_{!}: (\mSet)_{/S} \rightarrow (\mSet)_{/T}$$ and its right adjoint $p^{\ast}$ induce a Quillen equivalence between $(\mSet)_{/S}$ and $(\mSet)_{/T}$.
\end{proposition}

\begin{proof}
Let $\calC = \sCoNerve[S]^{op}$ and $\calD = \sCoNerve[T]^{op}$.
Consider the diagram of model categories and left Quillen functors:
$$\xymatrix{ (\mSet)_{/S} \ar[r]^{p_{!}} \ar[d]^{ \St^{+}_{S} } & (\mSet)_{/T} \ar[d]^{\St^{+}_{T}} \\
\calC \ar[r]^{ \sCoNerve[p]_{!} } & \calD }.$$ According to Proposition \ref{formall}, this diagram commutes (up to natural isomorphism). Theorem \ref{straightthm} implies that the vertical arrows are Quillen equivalences. Since $p$ is a categorical equivalence, $\sCoNerve[p]$ is an equivalence of simplicial categories, so that $\sCoNerve[p]_{!}$ is a Quillen equivalence (Proposition \ref{lesstrick}). It follows that $(p_!, p^{\ast})$ is a Quillen equivalence as well. 
\end{proof}

\begin{corollary}\label{tttroke}
Let $p: X \rightarrow S$ be a Cartesian fibration of simplicial sets, and let
$S \rightarrow T$ be a categorical equivalence. Then there exists a Cartesian fibration $Y \rightarrow T$, and an equivalence of $X$ with $S \times_T Y$ $($as Cartesian fibrations over $X${}$)$.
\end{corollary}

\begin{proof}
Proposition \ref{pretrokee} implies that the right derived functor $R p^{\ast}$ is essentially surjective.
\end{proof}

As we explained in Remark \ref{rightprop}, the Joyal model structure on $\sSet$ is {\em not} right proper. In other words, it is possible to have a categorical fibration $X \rightarrow S$ and 
a categorical equivalence $T \rightarrow S$ such that the induced map $X \times_{S} T \rightarrow X$ is not a categorical equivalence. This poor behavior of categorical fibrations is one of the reason that they do not play a prominent role in the theory of $\infty$-categories. Working with a stronger notion of fibration corrects the problem:

\begin{proposition}\label{basechangefunky}\index{gen}{Cartesian fibration!and pullbacks}
Let $p: X \rightarrow S$ be a Cartesian fibration, and let $T \rightarrow S$ be a categorical equivalence. Then the induced map $X \times_{S} T \rightarrow X$ is a categorical equivalence.
\end{proposition}

\begin{proof}
We first suppose that the map $T \rightarrow S$ is inner anodyne. By means of a simple argument, we may reduce to the case where $T \rightarrow S$ is a middle horn inclusion $\Lambda^n_i \subseteq \Delta^n$, where $0 < i < n$. 
According Proposition \ref{simplexplay}, there exists a sequence of maps
$$ \phi: A^0 \leftarrow \ldots \leftarrow A^n $$ and a map
$M(\phi) \rightarrow X$ which is a categorical equivalence, such that
$M(\phi) \times_{S} T \rightarrow X \times_{S} T$ is also a categorical equivalence.
Consequently, it suffices to show that the inclusion $M(\phi) \times_{S} T \subseteq M(\phi)$ is a categorical equivalence. But this map is a pushout of the inclusion
$A^n \times \Lambda^n_i \subseteq A^n \times \Delta^n$, which is inner anodyne.

We now treat the general case. Choose an inner anodyne map $T \rightarrow T'$ where $T'$ is an $\infty$-category. Then choose an inner anodyne map $T' \coprod_{T} S \rightarrow S'$, where $S'$ is also an $\infty$-category. The map $S \rightarrow S'$ is inner anodyne; in particular it is a categorical equivalence, so by
Corollary \ref{tttroke} there is a Cartesian fibration $X' \rightarrow S'$ and an equivalence
$X \rightarrow X' \times_{S'} S$ of Cartesian fibrations over $S$. We have a commutative diagram
$$ \xymatrix{ & X' \times_{S'} T \ar[r]^{u'} & X' \times_{S'} T' \ar[dr]^{u''} & \\
X \times_{S} T \ar[ur]^{u} \ar[dr]^{v} & & & X'. \\
& X \ar[r]^{v'} & X' \times_{S'} S \ar[ur]^{v''} & \\}$$
Consequently, to prove that $v$ is a categorical equivalence, it suffices to show that every other arrow in the diagram is a categorical equivalence. The maps $u$ and $v'$ are equivalences of Cartesian fibrations, and therefore categorical equivalences. The other three maps are special cases of the assertion we are trying to prove., For the map $u''$, we are in the special case of the map $S' \rightarrow T'$, which is an equivalence of $\infty$-categories: in this case we simply apply Corollary \ref{usesec}. For the maps $u'$ and $v''$, we need to know that the assertion of the proposition is valid in the special case of the maps $S \rightarrow S'$ and $T \rightarrow T'$. Since these maps are inner anodyne, the proof is complete.
\end{proof}

\begin{corollary}\label{basety}
Let 
$$ \xymatrix{ X \ar[r] \ar[d] & X' \ar[d]^{p'} \\
S \ar[r] & S' }$$
be a pullback diagram of simplicial sets, where $p'$ is a Cartesian fibration. Then the diagram is homotopy Cartesian $($with respect to the Joyal model structure$)$.
\end{corollary}

\begin{proof}
Choose a categorical equivalence $S' \rightarrow S''$, where $S''$ is an $\infty$-category.
Using Proposition \ref{pretrokee}, we may assume without loss of generality that
$X' \simeq X'' \times_{S''} S'$, where $X'' \rightarrow S''$ is a Cartesian fibration. Now choose
a factorization
$$ S \stackrel{\theta'}{\rightarrow} T \stackrel{\theta''}{\rightarrow} S''$$
where $\theta'$ is a categorical equivalence and $\theta''$ is a categorical fibration. 
The diagram
$$ T \rightarrow S'' \leftarrow X''$$ is fibrant. Consequently, the desired conclusion is equivalent to the assertion that the map $X \rightarrow T \times_{S''} X''$ is a categorical equivalence, which follows immediately from Proposition \ref{basechangefunky}.
\end{proof}

We now prove a stronger version of Corollary \ref{usefir}, which does not require that the base $S$ is a $\infty$-category.

\begin{proposition}\label{apple1}
Suppose given a diagram of simplicial sets
$$ \xymatrix{ X \ar[dr]^{p} \ar[rr]^{f} & & Y \ar[dl]^{q} \\
& S & }$$ 
where $p$ and $q$ are Cartesian fibrations, and $f$ carries $p$-Cartesian edges
to $q$-Cartesian edges.
The following conditions are equivalent:
\begin{itemize}
\item[$(1)$] The map $f$ is a categorical equivalence.
\item[$(2)$] For each vertex $s$ of $S$, $f$ induces a categorical equivalence $X_{s} \rightarrow Y_{s}$.
\item[$(3)$] The map $X^{\natural} \rightarrow Y^{\natural}$ is a Cartesian equivalence in
$(\mSet)_{/S}$. 
\end{itemize}

\end{proposition}

\begin{proof}
The equivalence of $(2)$ and $(3)$ follows from Proposition \ref{crispy}. We next show that $(2)$ implies $(1)$. In virtue of Proposition \ref{tulky}, we may reduce to the case where $S$ is a simplex. Then $S$ is an $\infty$-category and the desired result follows from Corollary \ref{usefir}. (Alternatively, we could observe that $(2)$ implies that $f$ has a homotopy inverse.)

To prove that $(1)$ implies $(3)$, we choose an inner anodyne map $j: S \rightarrow S'$, where $S'$ is an $\infty$-category. Let $X^{\natural}$ denote the object of $(\mSet)_{/S}$ associated to the Cartesian fibration $p: X \rightarrow S$, and let $j_{!} X^{\natural}$ denote the same marked simplicial set, regarded as an object of $(\mSet)_{/T}$. Choose a marked anodyne map
$j_{!} X^{\natural} \rightarrow {X'}^{\natural}$, where $X' \rightarrow S'$ is a Cartesian fibration.
By Proposition \ref{pretrokee}, the map $X^{\natural} \rightarrow j^{\ast} {X'}^{\natural}$ is a Cartesian equivalence, so that $X \rightarrow X' \times_{S'} S$ is a categorical equivalence. According to Proposition \ref{basechangefunky}, the map $X' \times_{S'} S \rightarrow X'$ is a categorical equivalence; thus the composite map $X \rightarrow X'$ is a categorical equivalence.

Similarly, we may choose a marked anodyne map $${X'}^{\natural} \coprod_{ j_{!} X^{\natural}} j_{!} Y^{\natural} \rightarrow {Y'}^{\natural}$$ for some Cartesian fibration $Y' \rightarrow S'$. Since
the Cartesian model structure is left-proper, the map $j_{!} Y^{\natural} \rightarrow {Y'}^{\natural}$
is a Cartesian equivalence, so we may argue as above to deduce that $Y \rightarrow Y'$ is a categorical equivalence. Now consider the diagram
$$ \xymatrix{ X \ar[r]^{f} \ar[d] & Y \ar[d] \\
X' \ar[r]^{f'} \ar[r] & Y'. }$$
We have argued that the vertical maps are categorical equivalences. The map $f$ is a categorical equivalence by assumption. It follows that $f'$ is a categorical equivalence. Since $S'$ is an $\infty$-category, we may apply Corollary \ref{usefir} to deduce that $X'_{s} \rightarrow Y'_{s}$
is a categorical equivalence for each object $s$ of $S'$. It follows that ${X'}^{\natural}
\rightarrow {Y'}^{\natural}$ is a Cartesian equivalence in $(\mSet)_{/S}$, so that we have a commutative diagram
$$ \xymatrix{ X^{\natural} \ar[r] \ar[d] & Y^{\natural} \ar[d] \\
j^{\ast} {X'}^{\natural} \ar[r] & j^{\ast} {Y'}^{\natural} }$$
where the vertical and bottom horizontal arrows are Cartesian equivalences in $(\mSet)_{/S}$. It follows that the top horizontal arrow is a Cartesian equivalence as well, so that $(3)$ is satisfied.
\end{proof}

\begin{corollary}\label{ruy}
Let
$$ \xymatrix{ W \ar[r] \ar[d] & X \ar[d] &  \\
Y \ar[r] & Z \ar[r] & S }$$
be a diagram of simplicial sets. Suppose that every morphism in this diagram is a right fibration, and that the square is a pullback. Then the diagram is homotopy Cartesian with respect to the contravariant model structure on $(\sSet)_{/S}$.
\end{corollary}

\begin{proof}
Choose a fibrant replacement
$$ X' \rightarrow Y' \leftarrow Z'$$
for the diagram
$$ X \rightarrow Y \leftarrow Z$$
in $(\sSet)_{/S}$, and let $W' = X' \times_{Z'} Y'$. We wish to show that the induced map
$i: W \rightarrow W'$ is a covariant equivalence in $(\sSet)_{/S}$. According to Corollary \ref{prefibchar}, it will suffice to show that for each vertex $s$ of $S$, the map of fibers
$W_{s} \rightarrow W'_{s}$ is a homotopy equivalence of Kan complexes. 
To prove this, we observe that we have a natural transformation of diagrams from
$$ \xymatrix{ W_{s} \ar[r] \ar[d] & X_s \ar[d] \\
Y_{s} \ar[r] & Z_{s} }$$
to
$$ \xymatrix{ W'_{s} \ar[r] \ar[d] & X'_s \ar[d] \\
Y'_{s} \ar[r] & Z'_{s} }$$
which induces homotopy equivalences 
$$X_{s} \rightarrow X'_{s} \quad \quad Y_{s} \rightarrow Y'_{s} \quad \quad Z_{s} \rightarrow Z'_{s}$$ 
(Corollary \ref{prefibchar}), where both diagrams are homotopy Cartesian (Proposition \ref{dent}).
\end{proof}

\begin{proposition}\label{funkyfibcatfib}\index{gen}{Cartesian fibration!and categorical fibrations}
Let $p: X \rightarrow S$ be a Cartesian fibration of simplicial sets. Then $p$ is a categorical fibration.
\end{proposition}

\begin{proof}
Consider a diagram
$$ \xymatrix{ A \ar[r] \ar@{^{(}->}[d]^{i} & X \ar[d]^{p} \\
B \ar[r] \ar@{-->}[ur]^{f} \ar[r]^{g}  & S} $$
of simplicial sets where $i$ is an inclusion and a categorical equivalence. We must demonstrate the existence of the indicated dotted arrow. Choose a categorical equivalence
$j: S \rightarrow T$, where $T$ is an $\infty$-category. By Corollary \ref{tttroke}, there exists
a Cartesian fibration $q: Y \rightarrow T$
such that $Y \times_{T} S$ is equivalent to $X$. Thus, there exist maps $$ u: X \rightarrow Y \times_{T} S $$
$$ v: Y \times_{T} S \rightarrow X$$
such that $u \circ v$ and $v \circ u$ are homotopic to the identity (over $S$).

Consider the induced diagram
$$ \xymatrix{ A \ar[r] \ar@{^{(}->}[d]^{i} & Y \\
B. \ar@{-->}[ur]_{f'} &} $$
Since $Y$ is an $\infty$-category, there exists a dotted arrow $f'$ making the diagram commutative. Let $g' = q \circ f': B \rightarrow T$. We note that $g'|A = (j \circ g)|A$. Since $T$ is an $\infty$-category and $i$ is a categorical equivalence, there exists a homotopy
$B \times \Delta^1 \rightarrow T$ from $g'$ to $j \circ g$ which is fixed on $A$. Since
$q$ is a Cartesian fibration, this homotopy lifts to a homotopy from $f'$ to some map
$f'': B \rightarrow Y$, so that we have a commutative diagram 
$$ \xymatrix{ A \ar[r] \ar@{^{(}->}[d]^{i} & Y \ar[d]^{q} \\
B \ar[r] \ar@{-->}[ur]_{f''} \ar[r] & T.} $$

Consider the composite map
$$ f''': B \stackrel{(f'',g)}{\rightarrow} Y \times_{T} S \stackrel{v}{\rightarrow} X.$$
Since $f'$ is homotopic to $f''$, and $v \circ u$ is homotopic to the identity, we conclude that
$f'''|A$ is homotopic to $f_0$ (via a homotopy which is fixed over $S$). Since $p$ is a Cartesian fibration, we can extend $h$ to a homotopy from $f'''$ to the desired map $f$.
\end{proof}

In general, the converse to Proposition \ref{funkyfibcatfib} fails: a categorical fibration of simplicial sets $X \rightarrow S$ need not be a Cartesian fibration. This is clear, since the property of being a categorical fibration is self-dual while the condition of being a Cartesian fibration is not. However, in the case where $S$ is a Kan complex, the theory of Cartesian fibrations {\em is} self-dual, and we have the following result:

\begin{proposition}\label{groob}
Let $p: X \rightarrow S$ be a map of simplicial sets, where $S$ is a Kan complex.
The following assertions are equivalent:
\begin{itemize}
\item[$(1)$] The map $p$ is a Cartesian fibration.
\item[$(2)$] The map $p$ is a coCartesian fibration.
\item[$(3)$] The map $p$ is a categorical fibration.
\end{itemize}
\end{proposition}

\begin{proof}
We will prove that $(1)$ is equivalent to $(3)$; the equivalence of $(2)$ and $(3)$ follows from a dual argument. Proposition \ref{funkyfibcatfib} shows that $(1)$ implies $(3)$ (for this implication, the assumption that $S$ is a Kan complex is not needed).

Now suppose that $(3)$ holds. Then $X$ is an $\infty$-category. Since every edge of $S$ is an equivalence, the $p$-Cartesian edges of $X$ are precisely the equivalences in $X$. It therefore suffices to show that if if $y$ is a vertex of $X$ and $\overline{e}: \overline{x} \rightarrow p(y)$ is an edge of $S$, then $\overline{e}$ lifts to an equivalence $e: x \rightarrow y$ in $S$. Since
$S$ is a Kan complex, we can find a contractible Kan complex $K$ and a map
$\overline{q}: K \rightarrow S$ such that $\overline{e}$ is the image of an edge $e': x' \rightarrow y'$ in $K$. 
The inclusion $\{y'\} \subseteq K$ is a categorical equivalence; since $p$ is a categorical fibration, we can lift $\overline{q}$ to a map $q: K \rightarrow X$ with $q(y')=y$. Then $e=q(e')$ has the desired properties.
\end{proof}

\subsection{Universal Fibrations}\label{universalfib}

In this section, we will apply Theorem \ref{straightthm} to construct a {\em universal} Cartesian fibration. Recall that $\Cat_{\infty}$ is defined to be the nerve of the simplicial category
$\Cat_{\infty}^{\Delta} = ( \mSet)^{\degree}$ of $\infty$-categories. In particular, we may regard the inclusion $\Cat_{\infty}^{\Delta} \hookrightarrow \mSet$ as a (projectively) fibrant object 
$\calF \in (\mSet)^{ \Cat_{\infty}^{\Delta} }$. Applying the unstraightening functor
$\Un^{+}_{\Cat_{\infty}^{op}}$, we obtain a fibrant object of $(\mSet)_{/\Cat_{\infty}^{op}}$, which we may identify with Cartesian fibration $q: \calZ \rightarrow \Cat_{\infty}^{op}$. We will refer to $q$
as the {\it universal Cartesian fibration}.\index{gen}{Cartesian fibration!universal} We observe that the objects of $\Cat_{\infty}$ can be identified with $\infty$-categories, and that the fiber of $q$ over an $\infty$-category $\calC$ can be identified with $U(\calC)$, where $U$ is the functor described in Lemma \ref{utest}. In particular, there is a canonical equivalence of $\infty$-categories $$\calC \rightarrow U(\calC) =
\calZ \times_{ \Cat_{\infty}^{op} } \{ \calC \}.$$ Thus we may think of $q$ as a Cartesian fibration which associates to each object of $\Cat_{\infty}$ the associated $\infty$-category.\index{gen}{universal!Cartesian fibration}

\begin{remark}
The $\infty$-categories $\Cat_{\infty}$ and $\calZ$ are {\em large}. However, the universal Cartesian fibration $q$ is small in the sense that for any small simplicial set $S$ and any map $f: S \rightarrow \Cat_{\infty}^{op}$, the fiber product $S \times_{ \QC^{op}} \calZ$ is small. This is because the fiber product can be identified with $\Un^{+}_{\phi}(\calF| \sCoNerve[S])$, where $\phi: \sCoNerve[S] \rightarrow \mSet$ is the composition of $\sCoNerve[f]$ with the inclusion.
\end{remark}

\begin{definition}\label{classer}\index{gen}{Cartesian fibration!classified by $f: S \rightarrow \Cat_{\infty}^{op}$}\index{gen}{classifying map!for a (co)Cartesian fibration}
Let $p: X \rightarrow S$ be a Cartesian fibration of simplicial sets. We will say that
a functor $f: S \rightarrow \Cat_{\infty}^{op}$ {\it classifies $p$} if there is an equivalence of Cartesian fibrations $X \rightarrow \calZ \times_{ \Cat^{op}_{\infty}} S
\simeq \Un^{+}_{S} f$. 

Dually, if $p: X \rightarrow S$ is a coCartesian fibration, then we will say that a functor
$f: S \rightarrow \Cat_{\infty}$ {\it classifies $p$} if $f^{op}$ classifies the
Cartesian fibration $p^{op}: X^{op} \rightarrow S^{op}$.\index{gen}{coCartesian fibration!classified by $f: S \rightarrow \Cat_{\infty}$}
\end{definition}

\begin{remark}\label{uniright}
Every Cartesian fibration $X \rightarrow S$ between {\em small} simplicial sets admits a classifying map $\phi: S \rightarrow \Cat_{\infty}^{op}$, which is uniquely determined up to equivalence. 
This is one expression of the idea that $\calZ \rightarrow \Cat_{\infty}^{op}$ is a {\it universal} Cartesian fibration. However, it is not immediately obvious that this property characterizes $\Cat_{\infty}$ up to equivalence, because $\Cat_{\infty}$ is not itself small. To remedy the situation, let us consider an arbitrary uncountable regular cardinal $\kappa$, and let $\Cat_{\infty}(\kappa)$ denote the full subcategory of $\Cat_{\infty}$ spanned by the $\kappa$-small $\infty$-categories. We then deduce the following:
\begin{itemize}
\item[$(\ast)$] Let $p: X \rightarrow S$ be a Cartesian fibration between small simplicial sets.
Then $p$ is classified by a functor $\chi: S \rightarrow \Cat_{\infty}(\kappa)^{op}$ if and only if,
for every vertex $s \in S$, the fiber $X_{s}$ is essentially $\kappa$-small. In this case, $\chi$ is determined uniquely up to homotopy.
\end{itemize}
Enlarging the universe and applying $(\ast)$ in the case where $\kappa$ is the supremum of all small cardinals, we deduce the following property:
\begin{itemize}
\item[$(\ast')$] Let $p: X \rightarrow S$ be a Cartesian fibration between simplicial sets which are not necessarily small. Then $p$ is classified by a functor $\chi: S \rightarrow \Cat_{\infty}^{op}$ if and only if, for every vertex $s \in S$, the fiber $X_{s}$ is essentially small. In this case, $\chi$ is determined uniquely up to homotopy.
\end{itemize}
This property evidently determines the $\infty$-category $\Cat_{\infty}$ (and the
Cartesian fibration $q: \calZ \rightarrow \Cat_{\infty}^{op}$) up to equivalence.
\end{remark}

\begin{warning}
The terminology of Definition \ref{classer} has the potential to cause confusion in the case where
$p: X \rightarrow S$ is both a Cartesian fibration and a coCartesian fibration. In this case,
$p$ is classified both by a functor $S \rightarrow \Cat_{\infty}^{op}$ (as a Cartesian fibration)
and by a functor $S \rightarrow \Cat_{\infty}$ (as a coCartesian fibration).
\end{warning}

The category $\Kan$ of Kan complexes can be identified with a full (simplicial) subcategory of $\Cat_{\infty}^{\Delta}$. Consequently we may identify the $\infty$-category
$\SSet$ of spaces with the full simplicial subset of $\Cat_{\infty}$, spanned by the vertices which represent $\infty$-groupoids. We let $\calZ^0 = \calZ \times_{ \Cat_{\infty}^{op}} \SSet^{op}$ be the restriction of the universal Cartesian fibration. The fibers of $q^0: \calZ^0 \rightarrow \SSet^{op}$
are Kan complexes (since they are equivalent to the $\infty$-categories represented by the vertices of $\SSet$). It follows from Proposition \ref{goey} that $q^0$ is a right fibration. We will refer to
$q^0$ as the {\it universal right fibration}.\index{gen}{universal!right fibration}\index{gen}{right fibration!universal}

Proposition \ref{goey} translates immediately into the following characterization of right fibrations:

\begin{proposition}
Let $p: X \rightarrow S$ be a Cartesian fibration of simplicial sets. The following conditions are equivalent:
\begin{itemize}
\item[$(1)$] The map $p$ is a right fibration.
\item[$(2)$] Every functor $f: S \rightarrow \Cat_{\infty}^{op}$ which classifies $p$ factors
through $\SSet^{op} \subseteq \Cat_{\infty}^{op}$.
\item[$(3)$] There exists a functor $f: S \rightarrow \SSet^{op}$ which classifies $p$.
\end{itemize}
\end{proposition}

Consequently, we may speak of right fibrations $X \rightarrow S$ being classified by functors
$S \rightarrow \SSet^{op}$, and left fibrations being classified by functors $S \rightarrow \SSet$.\index{gen}{right fibration!classifed by $S \rightarrow \SSet^{op}$}\index{gen}{left fibration!classified by $S \rightarrow \SSet$}\index{gen}{classifying map!for a right fibration}\index{gen}{classifying map!for a left fibration}

The $\infty$-category $\Delta^0$ corresponds to a vertex of $\Cat_{\infty}$ which we will denote by $\ast$. The fiber of $q$ over this point may be identified with $U \Delta^0 \simeq \Delta^0$; consequently, there is a unique vertex $\ast_{\calZ}$ of $\calZ$ lying over $\ast$. 
We note that $\ast$ and $\ast_{\calZ}$ belong to the subcategories $\SSet$ and $\calZ^0$. Moreover, we have the following:

\begin{proposition}\label{unifinal}
Let $q^0: \calZ^0 \rightarrow \SSet^{op}$ be the universal right fibration. The vertex $\ast_{\calZ}$ is a final object of the $\infty$-category $\calZ^0$.
\end{proposition}

\begin{proof}
Let $n > 0$, and let $f_0: \bd \Delta^n \rightarrow \calZ^0$ have the property that $f_0$ carries the final vertex of $\Delta^n$ to $\ast_{\calZ}$. We wish to show that there exists an extension
$$ \xymatrix{ \bd \Delta^n \ar[r]^{f_0} \ar@{^{(}->}[d] & \calZ \\
\Delta^n \ar@{-->}[ur]^{f} }$$
(in which case the map $f$ automatically factors through $\calZ^0$).

Let $\calD$ denote the simplicial category containing $\SSet^{op}_{\Delta}$ as a full subcategory, together with one additional object $X$, with the morphisms given by
$$ \bHom_{\calD}( K, X) = K $$
$$ \bHom_{\calD}(X,X) =  \ast $$
$$ \bHom_{\calD}(X,K) = \emptyset$$
for all $K \in \SSet^{op}_{\Delta}$. Let $\calC = \sCoNerve[ \Delta^n \star \Delta^0 ]$, and let
$\calC_0$ denote the full subcategory $\calC_0 = \sCoNerve[ \bd \Delta^n \star \Delta^0 ]$. 
We will denote the objects of $\calC$ by $\{ v_0, \ldots, v_{n+1} \}$.
Giving the map $f_0$ is tantamount to giving a simplicial functor $F_0: \calC_0 \rightarrow \calD$
with $F_0(v_{n+1})=X$, and constructing $f$ amounts to giving a simplicial functor $F: \calC \rightarrow \calD$ which extends $F_0$.

We note that the inclusion $\bHom_{\calC_0}( v_i, v_j) \rightarrow \bHom_{\calC}(v_i,v_j)$ is an isomorphism, unless $i=0$ and $j \in \{n, n+1\}$. Consequently, to define $F$, it suffices to find extensions
$$ \xymatrix{ \bHom_{\calC_0}(v_0, v_n) \ar[r] \ar@{^{(}->}[d] & \bHom_{\calD}(F_0(v_0), F_0(v_n)) \\ \bHom_{\calC}(v_0,v_n) \ar@{-->}[ur]^{j} }$$
$$ \xymatrix{ \bHom_{\calC_0}(v_0, v_{n+1}) \ar[r] \ar@{^{(}->}[d] & \bHom_{\calD}(F_0(v_0), F_0(v_{n+1})) \\ \bHom_{\calC}(v_0,v_{n+1}) \ar@{-->}[ur]^{j'} }$$
such that the following diagram commutes:
$$ \xymatrix{ \bHom_{\calC}(v_0,v_n) \times \bHom_{\calC}(v_n,v_{n+1}) \ar[r] \ar[d] &
\bHom_{\calD}(F_0(v_0), F_0(v_n)) \times \bHom_{\calD}(F_0(v_n), F_0(v_{n+1})) \ar[d] \\
\bHom_{\calC}(v_0, v_{n+1}) \ar[r] & \bHom_{\calD}(F_0(v_0), F_0(v_{n+1})). }$$

We note that $\bHom_{\calC}(v_{n},v_{n+1})$ is a point. In view of the assumption that
$f_0$ carries the final vertex of $\Delta^n$ to $\ast_{\calZ}$, we see that
$\bHom_{\calD}( F(v_n), F(v_{n+1}))$ is a point. It follows that, for any fixed choice of $j'$, there is a unique choice of $j$ for which the above diagram commutes. It therefore suffices to show that $j'$ exists. Since $\bHom_{ \calD }( F_0(v_0), X )$ is a Kan complex, it will suffice to show that the inclusion $\bHom_{\calC_0}( v_0, v_{n+1}) \rightarrow \bHom_{\calC}(v_0,v_{n+1})$ is an anodyne map of simplicial sets. In fact, it is isomorphic to the inclusion
$$ (\{1\} \times (\Delta^1)^{n-1}) \coprod_{ \{1\} \times \bd (\Delta^1)^{n-1} }
(\Delta^1 \times \bd (\Delta^1)^{n-1}) \subseteq \Delta^1 \times \Delta^{n-1},$$
which is the smash product of the cofibration $\bd (\Delta^1)^{n-1} \subseteq (\Delta^1)^{n-1}$ with the anodyne inclusion $\{1\} \subseteq \Delta^1$.
\end{proof}

\begin{corollary}\label{grt}
The universal right fibration $q^0: \calZ^0 \rightarrow \SSet^{op}$ is representable
by the final object of $\SSet$.
\end{corollary}

\begin{proof}
Combine Propositions \ref{unifinal} and \ref{reppfunc}.
\end{proof}

\begin{corollary}\label{unipull}
Let $p: X \rightarrow S$ be a left fibration between small simplicial sets. Then there exists a map
$S \rightarrow \SSet$ and an equivalence of left fibrations $X \simeq S \times_{\SSet} \SSet_{\ast/}$. 
\end{corollary}

\begin{proof}
Combine Corollary \ref{grt} with Remark \ref{uniright}.
\end{proof}

\subsection{Limits of $\infty$-Categories}\label{catlim}

The $\infty$-category $\Cat_{\infty}$ can be identified with the simplicial nerve of
$(\mSet)^{\degree}$. It follows from Corollary \ref{limitsinmodel} that $\Cat_{\infty}$ admits (small) limits and colimits, which can be computed in terms of homotopy (co)limits in the model category $\mSet$. For many applications, it is convenient to be able to construct limits and colimits while working entirely in the setting of $\infty$-categories. We will describe the construction of limits in this section; the case of colimits will be discussed in \S \ref{catcolim}.

Let $p: S^{op} \rightarrow \Cat_{\infty}$ be a diagram in $\Cat_{\infty}$. Then $p$ classifies
a Cartesian fibration $q: X \rightarrow S$. We will show (Corollary \ref{blurt} below) that
the limit $\projlim(p) \in \Cat_{\infty}$ can be identified with the $\infty$-category of {\em Cartesian} sections of $q$. We begin by proving a more precise assertion:

\begin{proposition}\label{charcatlimit}\index{gen}{limit!of $\infty$-categories}
Let $K$ be a simplicial set, $\overline{p}: K^{\triangleright} \rightarrow \Cat^{op}_{\infty}$ a diagram
in the $\infty$-category of spaces, $\overline{X} \rightarrow K^{\triangleright}$ a Cartesian fibration classified by $\overline{p}$, and $X = \overline{X} \times_{K^{\triangleright}} K$. 
The following conditions are equivalent:
\begin{itemize}
\item[$(1)$] The diagram $\overline{p}$ is a colimit of $p = \overline{p} | K$.

\item[$(2)$] The restriction map
$$ \theta: \bHom^{\flat}_{K^{\triangleright}}( (K^{\triangleright})^{\sharp}, \overline{X}^{\natural}) \rightarrow \bHom^{\flat}_{K}(K^{\sharp}, X^{\natural})$$
is an equivalence of $\infty$-categories.
\end{itemize}
\end{proposition}

\begin{proof}
According to Proposition \ref{cofinalcategories}, there exists a small category $\calC$ and a cofinal map $f: \Nerve(\calC) \rightarrow K$; let $\overline{\calC}= \calC \star [0]$ be the category obtained from $\calC$ by adjoining a new final object, and let
$\overline{f}: \Nerve(\overline{\calC}) \rightarrow K^{\triangleright}$ be the induced map (which is also cofinal). The maps $f$ and $\overline{f}$ are contravariant equivalences in
$(\sSet)_{/ K^{\triangleright} }$, and therefore induce Cartesian equivalences
$$ \Nerve(\calC)^{\sharp} \rightarrow K^{\sharp}$$
$$ \Nerve(\overline{\calC})^{\sharp} \rightarrow (K^{\triangleright})^{\sharp}.$$
We have a commutative diagram
$$ \xymatrix{ \bHom^{\flat}_{K^{\triangleright}}( (K^{\triangleright})^{\sharp}, 
\overline{X}^{\natural} ) \ar[r]^{\theta} \ar[d] & 
\bHom^{\flat}_{K^{\triangleright}}( K^{\sharp}, 
\overline{X}^{\natural} ) \ar[d] \\
\bHom^{\flat}_{K^{\triangleright}}( \Nerve(\overline{\calC})^{\sharp}, 
\overline{X}^{\natural} ) \ar[r]^{\theta'} & 
\bHom^{\flat}_{K^{\triangleright}}( \Nerve(\calC)^{\sharp}, 
\overline{X}^{\natural} ). }$$
The vertical arrows are categorical equivalences. Consequently, condition $(2)$ holds for 
$\overline{p}: K^{\triangleright} \rightarrow \Cat_{\infty}^{op}$ if and only if condition $(2)$
holds for the composition $\Nerve(\overline{\calC}) \rightarrow K^{\triangleright} \rightarrow
\Cat_{\infty}^{op}$. We may therefore assume without loss of generality that
$K = \Nerve(\calC)$. 

Using Corollary \ref{strictify}, we may further suppose that $\overline{p}$ is obtained as the simplicial nerve of a functor $\overline{\calF}: \overline{\calC}^{op} \rightarrow (\mSet)^{\degree}$. 
Changing $\overline{\calF}$ if necessary, we may suppose that it is a {\em strongly} fibrant
diagram in $\mSet$. Let $\calF = \overline{\calF}|\calC^{op}$. 
Let $\overline{\phi}: \sCoNerve[K^{\triangleright}]^{op} \rightarrow \overline{\calC}^{op}$
be the counit map, and $\phi: \sCoNerve[K]^{op} \rightarrow \calC^{op}$ the restriction of $\overline{\phi}$. We may assume without loss of generality that $\overline{X} = \St^{+}_{\phi} \overline{\calF}$. We have a (not strictly commutative) diagram of categories and functors
$$ \xymatrix{ \mSet \ar[d]^{\St^{+}_{\ast}} \ar[r]^{\times K^{\sharp}} & (\mSet)_{/K} \ar[d]^{\St^{+}_{\phi}} \\
\mSet \ar[r]^{\delta} & (\mSet)^{\calC^{op}}, }$$
where $\delta$ denotes the diagonal functor. This diagram commutes up to a natural transformation
$$ \St^{+}_{\phi}( K^{\sharp} \times Z) \rightarrow
\St^{+}_{\phi}(K^{\sharp}) \boxtimes \St^{+}_{\ast}(Z) \rightarrow \delta( \St^{+}_{\ast} Z ).$$
Here the first map is a weak equivalence by Proposition \ref{spek3}, and the second map is a weak  equivalence because $L \St^{+}_{\phi}$ is an equivalence of categories (Theorem \ref{straightthm}) and therefore carries the final object $K^{\sharp} \in \h{(\mSet)_{/K}}$ to a final object of
$\h{(\mSet)^{\calC^{op}}}$. We therefore obtain a diagram of {\em right} derived functors
$$ \xymatrix{ \h{\mSet} & \h{(\mSet)_{/K}} \ar[l]^{\Gamma} \\
\h{\mSet} \ar[u]^{R \Un^{+}_{\ast}} & \h{(\mSet)^{\calC^{op}}} \ar[l] \ar[u]^{R \Un^{+}_{\phi}}, }$$
which commutes up to natural isomorphism, where we regard $(\mSet)^{\calC^{op}}$ as equipped with the {\em injective} model structure described in \S \ref{quasilimit3}. Similarly, we have a commutative diagram
$$ \xymatrix{ \h{\mSet}& \h{(\mSet)_{/K^{\triangleright}}} \ar[l]^{\Gamma'} \\
\h{\mSet} \ar[u]^{R \Un^{+}_{\ast}} & \h{(\mSet)^{\overline{\calC}^{op}}} \ar[l] \ar[u]^{R \Un^{+}_{\overline{\phi}}}. }$$
Condition $(2)$ is equivalent to the assertion that the restriction map
$\Gamma'(\overline{X}^{\natural}) \rightarrow \Gamma(X^{\natural})$ is an isomorphism in
$\h{\mSet}$. Since the vertical functors in both diagrams are equivalences of categories (Theorem \ref{straightthm}), this is equivalent to the assertion that the map
$$ \varprojlim \overline{\calF} \rightarrow \varprojlim \calF$$
is a weak equivalence in $\mSet$. Since $\overline{\calC}$ has an initial object $v$, $(2)$ is equivalent to the assertion that $\overline{\calF}$ exhibits $\overline{\calF}(v)$ as a homotopy limit
of $\calF$ in $(\mSet)^{\degree}$. Using Theorem \ref{colimcomparee}, we conclude that
$(1) \Leftrightarrow (2)$ as desired.
\end{proof}

It follows from Proposition \ref{charcatlimit} that limits in $\Cat_{\infty}$ are computed by forming $\infty$-categories of Cartesian sections:

\begin{corollary}\label{blurt}
Let $p: K \rightarrow \Cat^{op}_{\infty}$ be a diagram in the $\infty$-category $\Cat_{\infty}$ of spaces and let
$X \rightarrow K$ be a Cartesian fibration classified by $p$. There
is a natural isomorphism 
$$ \projlim(p) \simeq \bHom^{\flat}_{K}( K^{\sharp}, X^{\natural} ) $$
in the homotopy category $\h{\Cat_{\infty}}$.
\end{corollary}

\begin{proof}
Let $\overline{p}: (K^{\triangleright})^{op} \rightarrow \Cat_{\infty}^{op}$ be a limit of $p$, and let
$X' \rightarrow K^{\triangleright}$ be a Cartesian fibration classified by
$\overline{p}$. Without loss of generality we may suppose $X \simeq X' \times_{K^{\triangleright}} K$. We have maps
$$ \bHom^{\flat}_{K}( K^{\sharp}, X^{\natural}) \leftarrow \bHom^{\flat}_{K^{\triangleright}}( (K^{\triangleright})^{\sharp}, {X'}^{\natural})
\rightarrow \bHom^{\flat}_{K^{\triangleright} } ( \{v\}^{\sharp}, {X'}^{\natural}),$$
where $v$ denotes the cone point of $K^{\triangleright}$. Proposition \ref{charcatlimit} implies
that the left map is an equivalence of $\infty$-categories. Since the inclusion
$\{v\}^{\sharp} \subseteq (K^{\triangleright})^{\sharp}$ is marked anodyne, the map on the right is a trivial fibration. We now conclude by observing that the space
$ \bHom^{\flat}_{K^{\triangleright} }( \{v\}^{\sharp}, {X'}^{\natural}) \simeq
X' \times_{ K^{\triangleright} } \{v\}$ can be identified with $\overline{p}(v) = \varprojlim(p)$.
\end{proof}

Using Proposition \ref{charcatlimit}, we can easily deduce an analogous characterization of
limits in the $\infty$-category of spaces.

\begin{corollary}\label{charspacelimit}\index{gen}{limit!of spaces}
Let $K$ be a simplicial set, $\overline{p}: K^{\triangleleft} \rightarrow \SSet$ a diagram
in the $\infty$-category of spaces, and $X \rightarrow K^{\triangleleft}$ a left fibration classified by $\overline{p}$. The following conditions are equivalent:
\begin{itemize}
\item[$(1)$] The diagram $\overline{p}$ is a limit of $p = \overline{p} | K$.

\item[$(2)$] The restriction map 
$$ \bHom_{K^{\triangleleft}}( K^{\triangleleft}, X) \rightarrow
\bHom_{K^{\triangleleft}}(K, X)$$
is a homotopy equivalence of Kan complexes.
\end{itemize}
\end{corollary}

\begin{proof}
The usual model structure on $\sSet$ is a localization of the Joyal model structure. It follows that the inclusion $\Kan \subseteq \Cat^{\Delta}_{\infty}$ preserves homotopy limits (of diagrams indexed by categories). Using Theorem \ref{colimcomparee}, Proposition \ref{cofinalcategories}, and Corollary \ref{strictify}, we conclude that the inclusion $\SSet \subseteq \Cat_{\infty}$ preserves (small) limits.
The desired equivalence now follows immediately from Proposition \ref{charcatlimit}.
\end{proof}

\begin{corollary}\label{needta}
Let $p: K \rightarrow \SSet$ be a diagram in the $\infty$-category $\SSet$ of spaces, and let
$X \rightarrow K$ be a left fibration classified by $p$. There
is a natural isomorphism
$$ \projlim(p) \simeq \bHom_{K}( K, X ) $$
in the homotopy category $\calH$ of spaces.
\end{corollary}

\begin{proof}
Apply Corollary \ref{blurt}.
\end{proof}

\begin{remark}
It is also possible to adapt the proof of Proposition \ref{charcatlimit} to give a direct proof of
Corollary \ref{charspacelimit}. We leave the details to the reader.
\end{remark}

\subsection{Colimits of $\infty$-Categories}\label{catcolim}

In this section, we will address the problem of constructing {\em colimits} in
the $\infty$-category $\Cat_{\infty}$. Let $p: S^{op} \rightarrow \Cat_{\infty}$ be
diagram, classifying a Cartesian fibration $f: X \rightarrow S$. In \S \ref{catlim}, we saw that
$\projlim(p)$ can be identified with the $\infty$-category of Cartesian sections of $f$. To construct the colimit $\injlim(p)$, we need to find an $\infty$-category which admits a map {\em from} each fiber $X_{s}$. The natural candidate, of course, is $X$ itself. However, because $X$ is generally not an $\infty$-category, we must take some care to formulate a correct statement.

\begin{lemma}\label{wilkins}
Let 
$$\xymatrix{
X' \ar[r] \ar[d] & X \ar[d]^{p} \\
S' \ar[r]^{q} & S }$$
be a pullback diagram of simplicial sets, where $p$ is a Cartesian fibration and
$q^{op}$ is cofinal. The induced map ${X'}^{\natural} \rightarrow X^{\natural}$ is
a Cartesian equivalence $($in $\mSet${}$)$.
\end{lemma}

\begin{proof}
Choose a cofibration $S' \rightarrow K$, where $K$ is a contractible Kan complex.
The map $q$ factors as a composition
$$ S' \stackrel{q'}{\rightarrow} S \times K \stackrel{q''}{\rightarrow} S.$$
It is obvious that the projection $X^{\natural} \times K^{\sharp} \rightarrow X^{\natural}$ is a Cartesian equivalence. We may therefore replace $S$ by $S \times K$ and $q$ by $q'$, thereby reducing to the case where $q$ is a cofibration. Proposition \ref{cofbasic} now implies that
$q$ is left-anodyne. It is easy to see that the collection of cofibrations $q: S' \rightarrow S$ for which the desired conclusion holds is weakly saturated. We may therefore reduce to the case where
$q$ is a horn inclusion $\Lambda^n_i \subseteq \Delta^n$, where $0 \leq i < n$.

We now apply Proposition \ref{simplexplay} to choose a sequence of composable maps
$$ \phi: A^0 \leftarrow \ldots \leftarrow A^n $$ and
a quasi-equivalence $M(\phi) \rightarrow X$. We have a commutative diagram of marked simplicial sets
$$ \xymatrix{ M^{\natural}(\phi) \times_{ (\Delta^n)^{\sharp} } (\Lambda^n_i)^{\sharp}
\ar@{^{(}->}[d]^{i} \ar[r] & {X'}^{\natural} \ar@{^{(}->}[d] \\
M^{\natural}(\phi) \ar[r] & X, }$$
Using Proposition \ref{halfy}, we deduce that the horizontal maps are Cartesian equivalences. To complete the proof, it will suffice to show that $i$ is a Cartesian equivalence. We now observe that $i$ is a pushout of the inclusion $i'': (\Lambda^n_i)^{\sharp} \times (A^n)^{\flat}
\subseteq (\Delta^n)^{\sharp} \times (A^n)^{\flat}$. It will therefore suffice to prove that $i''$
is a Cartesian equivalence. Using Proposition \ref{urlt}, we are reduced to proving that the inclusion
$(\Lambda^n_i)^{\sharp} \subseteq (\Delta^n)^{\sharp}$ is a Cartesian equivalence. According to Proposition \ref{strstr}, this is equivalent to the assertion that the horn inclusion
$\Lambda^n_i \subseteq \Delta^n$ is a weak homotopy equivalence, which is obvious.
\end{proof}

\begin{proposition}\label{charcatcolimit}\index{gen}{colimit!of $\infty$-categories}
Let $K$ be a simplicial set, $\overline{p}: K^{\triangleleft} \rightarrow \Cat_{\infty}^{op}$ be a diagram in the $\infty$-category $\Cat_{\infty}$, 
$\overline{X} \rightarrow K^{\triangleleft}$ a Cartesian fibration classified
by $\overline{p}$, and $X = \overline{X} \times_{ K^{\triangleleft} } K$.
The following conditions are equivalent:
\begin{itemize}
\item[$(1)$] The diagram $\overline{p}$ is a limit of $p = \overline{p} | K$.
\item[$(2)$] The inclusion $X^{\natural} \subseteq \overline{X}^{\natural}$ is a Cartesian equivalence in $(\mSet)_{/K^{\triangleleft}}$.
\item[$(3)$] The inclusion $X^{\natural} \subseteq \overline{X}^{\natural}$
is a Cartesian equivalence in $\mSet$.
\end{itemize}
\end{proposition}

\begin{proof}
Using the small object argument, we can construct a factorization
$$ X \stackrel{i}{\rightarrow} Y \stackrel{j}{\rightarrow} K^{\triangleleft} $$
where $j$ is a Cartesian fibration, $i$ induces a marked anodyne map
$X^{\natural} \rightarrow Y^{\natural}$, and 
$X \simeq Y \times_{ K^{\triangleleft} } K$. 
Since $i$ is marked anodyne, we can solve the lifting problem
$$ \xymatrix{ X^{\natural} \ar@{^{(}->}[d]^{i} \ar[r] & \overline{X}^{\natural} \ar[d] \\
Y^{\natural} \ar[r] \ar@{-->}[ur]^{q} & (K^{\triangleleft})^{\sharp}. }$$
Since $i$ is a Cartesian equivalence in $(\mSet)_{/K^{\triangleleft}}$, condition $(2)$ is equivalent to the assertion that $q$ is an equivalence of Cartesian fibrations over $K^{\triangleleft}$. Since $q$
induces an isomorphism over each vertex of $K$, this is equivalent to:
\begin{itemize}
\item[$(2')$] The map $q_{v}: Y_{v} \rightarrow \overline{X}_{v}$ is an equivalence of $\infty$-categories, where $v$ denotes the cone point of $K^{\triangleleft}$.
\end{itemize}
We have a commutative diagram
$$ \xymatrix{ Y^{\natural}_{v} \ar[r]^{q_v} \ar@{^{(}->}[d] & \overline{X}^{\natural}_v \ar@{^{(}->}[d] \\
Y^{\natural} \ar[r]^{q} & \overline{X}^{\natural}. }$$
Lemma \ref{wilkins} implies that the vertical maps are Cartesian equivalences. It follows that $(2') \Leftrightarrow (3)$, so that $(2) \Leftrightarrow (3)$.

To complete the proof, we will show that $(1) \Leftrightarrow (2)$.
According to Proposition \ref{cofinalcategories}, there exists a small category $\calC$ and a map $p: \Nerve(\calC) \rightarrow K$ such that $p^{op}$ is cofinal. Let $\overline{\calC} = [0] \star \calC$ be the category obtained by adjoining an initial object to $\calC$. 
Consider the diagram
$$ \xymatrix{ (X \times_{K} \Nerve(\calC))^{\natural} \ar@{^{(}->}[r] \ar[d] & (\overline{X} \times_{K^{\triangleleft}} \Nerve(\overline{\calC}))^{\natural} \ar[d] \\
X^{\natural} \ar@{^{(}->}[r] & \overline{X}^{\natural}. }$$
Lemma \ref{wilkins} implies that the vertical maps are Cartesian equivalences (in $\mSet$).
It follows that the upper horizontal inclusion is a Cartesian equivalence if and only if the lower horizontal inclusion is a Cartesian equivalence. Consequently, it will suffice to prove the
equivalence $(1) \Leftrightarrow (2)$ after replacing $K$ by $\Nerve(\calC)$.

Using Corollary \ref{strictify}, we may further suppose that $\overline{p}$ is the nerve
of a functor $\calF: \overline{\calC} \rightarrow (\mSet)^{\degree}$. Let $\overline{\phi}: \sCoNerve[K^{\triangleleft}] \rightarrow \overline{\calC}$ be the counit map, and let $\phi: \sCoNerve[K] \rightarrow \calC$ be the restriction of $\overline{\phi}$. Without loss of generality, we may suppose that $\overline{X} = \Un_{\overline{\phi}} \calF$. We have a commutative diagram of homotopy categories and right derived functors
$$ \xymatrix{ \h{(\mSet)^{\overline{\calC}}} \ar[r]^{G} \ar[d]^{R \Un^{+}_{\overline{\phi}}} &  \h{(\mSet)^{\calC}}
\ar[d]^{ R \Un^{+}_{\phi} } \\
\h{(\mSet)_{/(K^{\triangleleft})}} \ar[r]^{G'} & \h{(\mSet)_{/K}} }$$
where $G$ and $G'$ are restriction functors. Let $F$ and $F'$ be the left adjoints
to $G$ and $G'$, respectively. According to Theorem \ref{colimcomparee}, assumption $(1)$ is equivalent to the assertion that $\calF$ lies in the essential image of $F$. Since each of the vertical functors is equivalence of categories (Theorem \ref{straightthm}), this is equivalent to the assertion that $\overline{X}$ lies in the essential image of $F'$. 
Since $F'$ is fully faithful, this is equivalent to the assertion that the counit map
$$ F' G' \overline{X} \rightarrow \overline{X}$$
is an isomorphism in $\h{(\mSet)_{/K^{\triangleleft})}}$, which is clearly a reformulation of $(2)$.
\end{proof}

\begin{corollary}\label{tolbot}
Let $p: K^{op} \rightarrow \Cat_{\infty}$ be a diagram, classifying a Cartesian fibration
$X \rightarrow K$. Then there is a natural isomorphism
$\varinjlim(p) \simeq X^{\natural}$ in the homotopy category
in $\h{\Cat_{\infty}}$.
\end{corollary}

\begin{proof}
Let $\overline{p}: (K^{op})^{\triangleright} \rightarrow \Cat_{\infty}$ be a colimit of $p$, classifying a Cartesian fibration $\overline{X} \rightarrow K^{\triangleleft}$. Let $v$ denote the cone point of
$K^{\triangleleft}$, so that $\varinjlim(p) \simeq \overline{X}_{v}$. We now observe that the inclusions
$$\overline{X}^{\natural}_{v} \hookrightarrow \overline{X}^{\natural} \hookleftarrow X^{\natural}$$
are both Cartesian equivalences (Lemma \ref{wilkins} and Proposition \ref{charcatcolimit}). 
\end{proof}

\begin{warning}
In the situation of Corollary \ref{tolbot}, the marked simplicial set $X^{\natural}$ is usually not a fibrant object of $\mSet$, even when $K$ is an $\infty$-category.
\end{warning}

Using exactly the same argument, we can establish a version of Proposition \ref{charcatcolimit} which describes colimits in the $\infty$-category of spaces:

\begin{proposition}\label{charspacecolimit}\index{gen}{colimit!of spaces}
Let $K$ be a simplicial set, $\overline{p}: K^{\triangleright} \rightarrow \SSet$ be a diagram in the $\infty$-category of spaces, $\overline{X} \rightarrow K^{\triangleright}$ a left fibration classified
by $\overline{p}$, and $X = \overline{X} \times_{ K^{\triangleright} } K$.
The following conditions are equivalent:
\begin{itemize}
\item[$(1)$] The diagram $\overline{p}$ is a colimit of $p = \overline{p} | K$.
\item[$(2)$] The inclusion $X \subseteq \overline{X}$ is a covariant equivalence in
$(\sSet)_{/K^{\triangleright}}$.
\item[$(3)$] The inclusion $X \subseteq \overline{X}$
is a weak homotopy equivalence of simplicial sets.
\end{itemize}
\end{proposition}

\begin{proof}
Using the small object argument, we can construct a factorization
$$ X \stackrel{i}{\hookrightarrow} Y \stackrel{j}{\rightarrow} K^{\triangleright} $$
where $i$ is left anodyne, $j$ is a left fibration, and the inclusion
$X \subseteq Y \times_{K^{\triangleright}} K$ is an isomorphism. Choose a dotted arrow
$q$ as indicated in the diagram
$$ \xymatrix{ X \ar@{^{(}->}[d]^{i} \ar[r] & \overline{X} \ar[d] \\
Y \ar[r] \ar@{-->}[ur]^{q} & K^{\triangleright}. }$$
Since $i$ is a covariant equivalence in $(\sSet)_{/K^{\triangleright}}$, condition $(2)$ is equivalent to the assertion that $q$ is an equivalence of left fibrations over $K^{\triangleright}$. Since $q$
induces an isomorphism over each vertex of $K$, this is equivalent to the assertion that
$q_{v}: Y_{v} \rightarrow \overline{X}_{v}$ is an equivalence, where $v$ denotes the cone point
of $K^{\triangleright}$. We have a commutative diagram
$$ \xymatrix{ Y_{v} \ar[r]^{q_v} \ar[d] & \overline{X}_v \ar[d] \\
Y \ar[r]^{q} & \overline{X}. }$$
Proposition \ref{strokhop} implies that the vertical maps are right anodyne, and therefore weak homotopy equivalences. Consequently, $q_{v}$ is a weak homotopy equivalence if and only if $q$ is a weak homotopy equivalence. Since the inclusion $X \subseteq Y$ is a weak homotopy equivalence, this proves that $(2) \Leftrightarrow (3)$.

To complete the proof, we will show that $(1) \Leftrightarrow (2)$.
According to Proposition \ref{cofinalcategories}, there exists a small category $\calC$ and a cofinal map $\Nerve(\calC) \rightarrow K$. Let $\overline{\calC} = \calC \star [0]$ be the category obtained from $\calC$ by adjoining a new final object.
Consider the diagram
$$ \xymatrix{ X \times_{K} \Nerve(\calC) \ar@{^{(}->}[r] \ar[d] & \overline{X} \times_{K^{\triangleright}} \Nerve(\overline{\calC}) \ar[d] \\
X \ar@{^{(}->}[r] & \overline{X}. }$$
Proposition \ref{strokhop} implies that $\overline{X} \rightarrow K^{\triangleright}$ is smooth, so that the vertical arrows in the above diagram are cofinal. In particular, the vertical arrows are weak homotopy equivalences, so that the upper horizontal inclusion is a weak homotopy equivalence
if and only if the lower horizontal inclusion is a weak homotopy equivalence. Consequently, 
it will suffice to prove the equivalence $(1) \Leftrightarrow (2)$ after replacing $K$ by $\Nerve (\calC)$.

Using Corollary \ref{strictify}, we may further suppose that $\overline{p}$ is obtained as the nerve of a functor $\calF: \overline{\calC} \rightarrow \Kan$. Let $\overline{\phi}: \sCoNerve[K^{\triangleright}] \rightarrow \overline{\calC}$ be the counit map, and let $\phi: \sCoNerve[K] \rightarrow \calC$ be the restriction of $\overline{\phi}$. Without loss of generality, we may suppose that $\overline{X}^{op} = \Un_{\overline{\phi}} \calF$. We have a commutative diagram of homotopy categories and right derived functors
$$ \xymatrix{ \h{(\sSet)^{\overline{\calC}}} \ar[r]^{G} \ar[d]^{R \Un_{\overline{\phi}}} & \h{(\sSet)^{\calC}}
\ar[d]^{ R \Un_{\phi} } \ar[d] \\
\h{ (\sSet)_{/ (K^{\triangleright})^{op}}} \ar[r]^{G'} & \h{ (\sSet)_{/K}} }$$
where $G$ and $G'$ are restriction functors. Let $F$ and $F'$ be the left adjoints
to $G$ and $G'$, respectively. According to Theorem \ref{colimcomparee}, assumption $(1)$ is equivalent to the assertion that $\calF$ lies in the essential image of $F$. Since each of the vertical functors is equivalence of categories (Theorem \ref{struns}), this is equivalent to the assertion that
$\overline{X}^{op}$ lies in the essential image of $F'$. Since $F'$ is fully faithful, this is equivalent to the assertion that the counit map
$$ F' G' \overline{X}^{op} \rightarrow \overline{X}^{op}$$
is an isomorphism in $\h{ (\sSet)_{/ (K^{\triangleright})^{op}}}$, which is clearly equivalent to $(2)$.
This shows that $(1) \Leftrightarrow (2)$ and completes the proof.
\end{proof}

\begin{corollary}\label{needka}
Let $p: K \rightarrow \SSet$ be a diagram which classifies a left fibration
$\widetilde{K} \rightarrow K$, and let $X \in \SSet$ be a colimit of $p$. Then
there is a natural isomorphism
$$ \widetilde{K} \simeq X$$
in the homotopy category $\calH$.
\end{corollary}

\begin{proof}
Let $\overline{p}: K^{\triangleright} \rightarrow \SSet$ be a colimit diagram which extends $p$, and
$\widetilde{K}' \rightarrow K^{\triangleright}$ a left fibration classified by $\overline{p}$. Without loss of generality, we may suppose that $\widetilde{K} = \widetilde{K}' \times_{K^{\triangleright}} K$ and
$X = \widetilde{K}' \times_{K^{\triangleright}} \{v\}$, where $v$ denotes the cone point of $K^{\triangleright}$. Since the inclusion $\{v\} \subseteq K^{\triangleright}$ is right anodyne
and the map $\widetilde{K}' \rightarrow K^{\triangleright}$ is a left fibration, Proposition \ref{strokhop} implies that the inclusion
$X \subseteq \widetilde{K}'$ is right anodyne, and therefore a weak homotopy equivalence.
On the other hand, Proposition \ref{charspacecolimit} implies that the inclusion
$\widetilde{K} \subseteq \widetilde{K}'$ is a weak homotopy equivalence. The
composition
$$ X \simeq \widetilde{K}' \simeq \widetilde{K}$$
is the desired isomorphism in $\calH$.
\end{proof}

\chapter{Limits and Colimits}\label{chap3}

\setcounter{theorem}{0}
\setcounter{subsection}{0}

This chapter is devoted to the study of limits and colimits in the $\infty$-categorical setting. Our goal is to provide tools for proving the existence of limits and colimits, for analyzing them, and for comparing them to the (perhaps more familiar) notion of homotopy limits and colimits in simplicial categories. We will generally confine our remarks to colimits; analogous results for limits can be obtained by passing to the opposite $\infty$-categories.

We begin in \S \ref{chap3cofinal} by introducing the notion of a {\it cofinal} map between simplicial sets. If $f: A \rightarrow B$ is a cofinal map of simplicial sets, then we can identify colimits of
a diagram $p: B \rightarrow \calC$ with colimits of the induced diagram $p \circ f: A \rightarrow \calC$. This is a basic maneuver which will appear repeatedly in the later chapters of this book. 
Consequently, it is important to have a good supply of cofinal maps. This is guaranteed by Theorem \ref{hollowtt}, which can be regarded as an $\infty$-categorical generalization of Quillen's ``Theorem A''.

In \S \ref{c4s2}, we introduce a battery of additional techniques for analyzing colimits. We will explain how to analyze colimits of complicated diagrams in terms of colimits of simpler diagrams. Using these ideas, we can often reduce questions about the behavior of arbitrary colimits to questions about a few basic constructions, which we will analyze explicitly in \S \ref{coexample}. We will also explain the relationship between the $\infty$-categorical theory of colimits and the more classical theory of homotopy colimits, which can be studied very effectively using the language of model categories.

The other major topic of this chapter is the theory of {\em Kan extensions}, which can be viewed as relative versions of limits and colimits. We will study the properties of Kan extensions in \S \ref{relacoim}, and prove some fundamental existence theorems which we will need throughout the later chapters of this book.

 \section{Cofinality}\label{chap3cofinal}

\setcounter{theorem}{0}

Let $\calC$ be an $\infty$-category, and let $p: K \rightarrow \calC$ be a diagram in $\calC$ indexed by a simplicial set $K$. In \S \ref{limitcolimit} we introduced the definition of a {\em colimit} $\injlim(p)$ for the diagram $p$. In practice, it is often possible to replace $p$ by a simpler diagram without changing the colimit $\injlim(p)$. In this section, we will introduce a general formalism which will allow us to make replacements of this sort: the theory of {\em cofinal} maps between simplicial sets. We begin in \S \ref{cofinal} with a definition of the class of cofinal maps, and show (Proposition \ref{gute}) that if a map $q: K' \rightarrow K$ is cofinal, then there is an equivalence
$\injlim(p) \simeq \injlim(p \circ q)$ (provided that either colimit exists). In \S \ref{smoothness} we will reformulate the definition of cofinality, using the formalism of contravariant model categories (\S \ref{contrasec}). We conclude in \S \ref{quillA} by establishing an important recognition criterion for cofinal maps, in the special case where $K$ is an $\infty$-category. This result can be regarded as a refinement of Quillen's ``Theorem A''. 
 
\subsection{Cofinal Maps}\label{cofinal}

The goal of this section is to introduce the definition of a cofinal map $p: S \rightarrow T$ of simplicial sets, and study the basic properties of this notion. Our main result is Proposition \ref{gute}, which characterizes cofinality in terms of the behavior of $T$-indexed colimits.

\begin{definition}[Joyal \cite{joyalnotpub}]\index{gen}{cofinal}
Let $p: S \rightarrow T$ be a map of simplicial sets. We shall
say that $p$ is {\it cofinal} if, for any right fibration $X
\rightarrow T$, the induced map of of simplicial sets
$$ \bHom_{T}(T,X) \rightarrow \bHom_{T}(S, X)$$
is a homotopy equivalence.
\end{definition}

\begin{remark}
The simplicial set $\bHom_{T}(S,X)$ parametrizes sections of the right fibration $X \rightarrow T$.
It may be described as the fiber of the induced map $X^{S} \rightarrow T^{S}$ over the vertex of $T^S$ corresponding to the map $p$. Since $X^{S} \rightarrow T^{S}$ is a right fibration, the fiber $\bHom_{T}(S,X)$ is a Kan complex. Similarly, $\bHom_{T}(T,X)$ is a Kan complex.
\end{remark}

We begin by recording a few simple observations about the class of cofinal maps:

\begin{proposition}\label{cofbasic}
\begin{itemize}
\item[$(1)$] Any isomorphism of simplicial sets is cofinal.

\item[$(2)$] Let $f: K \rightarrow K'$ and $g: K' \rightarrow K''$ be
maps of simplicial sets. Suppose that $f$ is cofinal. Then $g$ is cofinal if and only if $g \circ f$ is cofinal.

\item[$(3)$] If $f: K \rightarrow K'$ is a cofinal map between simplicial
sets, then $f$ is a weak homotopy equivalence.

\item[$(4)$] An inclusion $i: K \subseteq K'$ of simplicial sets is
cofinal if and only if it is right anodyne.
\end{itemize}
\end{proposition}

\begin{proof}
Assertions $(1)$ and $(2)$ are obvious. We prove $(3)$. Let $S$ be a Kan complex.
Since $f$ is cofinal, the composition
$$ \bHom_{\sSet}(K',S) = \bHom_{K}(K', S \times K) \rightarrow \bHom_{K}(K,S \times
K) = \bHom_{\sSet}(K,S)$$ is a homotopy equivalence. Passing to connected components, we deduce that $K$ and $K'$ co-represent the same functor in the homotopy category $\calH$ of spaces. It follows that $f$ is a weak homotopy equivalence, as desired.

We now prove $(4)$. Suppose first that $i$ is right-anodyne. Let
$X \rightarrow K'$ be a right fibration. Then the induced map $\Hom_{K'}(K',X)
\rightarrow \Hom_{K'}(K,X)$ is a trivial fibration, and in
particular a homotopy equivalence.

Conversely, suppose that $i$ is a cofinal inclusion of simplicial sets. 
We wish to show that $i$ has the left
lifting property with respect to any right fibration. In other
words, we must show that given any diagram of solid arrows
$$ \xymatrix{ K \ar@{^{(}->}[d] \ar[r]^{s} & X \ar[d] \\
K' \ar@{=}[r] \ar@{-->}[ur] & K', }$$
for which the right-vertical map is a right fibration, there exists a dotted arrow as indicated, rendering the diagram commutative. Since $i$ is cofinal, the map $s$ is homotopic
to a map which extends over $K'$. In other words, there exists a map
$$ s': (K \times \Delta^1) \coprod_{ K \times \{1\} } (K'
\times \{1\}) \rightarrow X,$$
compatible with the projection to $K'$, such that $s'| K \times \{0\}$ coincides with $s$.
Since the inclusion $$ (K \times \Delta^1) \coprod_{K \times \{1\} } (K' \times \{1\}) \subseteq
K' \times \Delta^1$$ is right-anodyne, there exists a map $s'': K' \times \Delta^1 \rightarrow
X$ which extends $s'$, and is compatible with the projection to $K'$. The map
$s'' | K \times \{0\}$ has the desired properties.
\end{proof}

\begin{warning}
The class of cofinal maps does {\em not} satisfy the
``two-out-of-three'' property. If $f: K \rightarrow K'$ and $g: K'
\rightarrow K''$ are such that $g \circ f$ and $g$ are cofinal,
then $f$ need not be cofinal.
\end{warning}

Our next goal is to establish an alternative characterization of cofinality, in terms of the behavior of colimits (Proposition \ref{gute}). First, we need a lemma.

\begin{lemma}\label{cogh}
Let $\calC$ be an $\infty$-category, and let $p: K \rightarrow \calC$, $q: K'
\rightarrow \calC$ be diagrams. Define simplicial sets $M$ and $N$ by the
formulas
$$ \Hom(X,M) = \{ f: (X \times K) \star K' \rightarrow \calC :
f|(X \times K) = p \circ \pi_{K}, f|K' = q \}$$
$$ \Hom(X,N) = \{ g: K \star (X \times K') \rightarrow \calC :
f|K = p, f|(X \times K') = q \circ \pi_{K'} \}. $$ Here $\pi_K$
and $\pi_{K'}$ denote the projection from a product to the factor
indicated by the subscript.

Then $M$ and $N$ are Kan complexes, which are (naturally) homotopy
equivalent to one another.
\end{lemma}

\begin{proof}
We define a simplicial set $\calD$ as follows. For every finite, nonempty, linearly ordered set $J$, to give a map $\Delta^{J} \rightarrow \calD$ is to supply the following data:

\begin{itemize}
\item A map $\Delta^{J} \rightarrow \Delta^1$, corresponding to a decomposition of
$J$ as a disjoint union $J_{-} \coprod J_{+}$, where $J_{-} \subseteq J$ is closed downwards and
$J_{+} \subseteq J$ is closed upwards.

\item A map $e: (K \times \Delta^{J_{-}}) \star (K' \times
\Delta^{J_{+}}) \rightarrow \calC$ such that $e| K \times
\Delta^{J_{-}} = p \circ \pi_K$ and $e| K' \times \Delta^{J_{+}} =
q \circ \pi_{K'}$.
\end{itemize}

We first claim that $\calD$ is an $\infty$-category. Fix a finite linearly ordered set $J$ as above, and let $j \in J$ be neither the largest nor the smallest element of $J$. Let 
$f_0: \Lambda^J_{j} \rightarrow \calD$ be any map; we wish to show that there exists a map
$f: \Delta^J \rightarrow \calD$ which extends $f_0$. We first observe that the induced projection
$\Lambda^J_{j} \rightarrow \Delta^1$ extends {\em uniquely} to $\Delta^J$ (since $\Delta^1$ is isomorphic to the nerve of a category). Let $J = J_{-} \coprod J_{+}$ be the induced decomposition of $J$. Without loss of generality, we may suppose that $j \in J_{-}$. In this case, we may identify
$f_0$ with a map
$$ (( K \times \Lambda^{J_{-}}_j ) \star (K' \times \Delta^{J_+} ))
\coprod_{ (K \times \Lambda^{J_-}_j ) \star (K' \times \bd \Delta^{J_+}) }
(( K \times \Delta^{J_{-}} ) \star (K' \times \bd \Delta^{J_+})) \rightarrow 
\calC$$
and our goal is to find an extension 
$$f: ( K \times \Delta^{J_{-}} ) \star (K' \times \Delta^{J_+} ) \rightarrow \calC.$$
Since $\calC$ is an $\infty$-category, it will suffice to show that the inclusion
$$ (( K \times \Lambda^{J_{-}}_j ) \star (K' \times \Delta^{J_+} ))
\coprod_{ (K \times \Lambda^{J_-}_j ) \star (K' \times \bd \Delta^{J_+}) }
(( K \times \Delta^{J_{-}} ) \star (K' \times \bd \Delta^{J_+})) \subseteq 
 ( K \times \Delta^{J_{-}} ) \star (K' \times \Delta^{J_+} )$$ is inner anodyne.
According to Lemma \ref{precough}, it suffices to check that the inclusion $K \times \Lambda^{J_{-}}_j \subseteq K \times \Delta^{J_{-}}$ is 
right anodyne. This follows from Corollary \ref{prodprod1}, since $\Lambda^{J_{-}}_j \subseteq \Delta^{J_{-}}$ is right anodyne.

The $\infty$-category $\calD$ has just two objects, which we will denote by
$x$ and $y$. We observe that $M = \Hom^{\rght}_{\calD}(x,y)$ and $N = \Hom^{\lft}_{\calD}(x,y)$.
Proposition \ref{gura} implies that $M$ and $N$ are Kan complexes.
Propositions \ref{babyy} and \ref{wiretrack} imply each these Kan complexes is weakly homotopy equivalent to $\bHom_{ \sCoNerve[\calD]}(x,y)$, so that $M$ and $N$ are homotopy equivalent to one another as desired.
\end{proof}

\begin{remark}\label{coughi}
In the situation of Lemma \ref{cogh}, the homotopy equivalence
between $M$ and $N$ is furnished by the composition of a chain of weak homotopy
equivalences
$$ M \leftarrow |M|_{Q^{\bigdot}} \rightarrow
\Hom_{\sCoNerve[\calD]}(x,y) \leftarrow |N|_{Q^{\bigdot}} \rightarrow
N,$$ which is functorial in the triple $(\calC,p: K \rightarrow \calC,q: K' \rightarrow \calC)$.
\end{remark}

\begin{proposition}\label{coughing}
Let $v: K' \rightarrow K$ be a cofinal map and $p: K \rightarrow \calC$
a diagram in an $\infty$-category $\calC$. Then the map $\phi: \calC_{p/}
\rightarrow \calC_{pv/}$ is an equivalence of left fibrations
over $\calC$: in other words, it induces a homotopy equivalence of Kan
complexes after passing to the fiber over every object $x$ of $\calC$.
\end{proposition}

\begin{proof}
We wish to prove that the map
$$ \calC_{p/} \times_{\calC} \{x\} \rightarrow \calC_{pv/} \times_{\calC} \{x\}$$
is a homotopy equivalence of Kan complexes. Lemma \ref{cogh} implies that the left
hand side is homotopy equivalent $\bHom_{\calC}(K, \calC_{/x})$. Similarly, the right hand
side can be identified with $\bHom_{\calC}(K', \calC_{/x})$. Using the functoriality implicit in the proof of Lemma \ref{cogh} (see Remark \ref{coughi}), it suffices to show that the restriction map
$$ \bHom_{\calC}(K, \calC_{/x}) \rightarrow \bHom_{\calC}(K', \calC_{/x})$$ is a homotopy equivalence. Since $v$ is cofinal, this follows immediately from the fact that the projection
$\calC_{/x} \rightarrow \calC$ is a right fibration.
\end{proof}

\begin{proposition}\label{gute}
Let $v: K' \rightarrow K$ be a map of (small) simplicial sets. The following conditions are equivalent:
\begin{itemize}
\item[$(1)$] The map $v$ is cofinal.
\item[$(2)$] Given any $\infty$-category $\calC$ and any diagram $p: K \rightarrow \calC$, the induced map $\calC_{p/} \rightarrow \calC_{p'/}$ is an equivalence of $\infty$-categories, where $p' = p \circ v$.
\item[$(3)$] For every $\infty$-category $\calC$ and every diagram $\overline{p}: K^{\triangleright} \rightarrow \calC$ which is a colimit of $p = \overline{p}|K$, the induced map $\overline{p}': {K'}^{\triangleright} \rightarrow \calC$
is a colimit of $p' = \overline{p}'|K'$.
\end{itemize}
\end{proposition}

\begin{proof}
Suppose first that $(1)$ is satisfied. Let $p: K \rightarrow \calC$ be as in $(2)$. Proposition \ref{coughing} implies that the induced map $\calC_{p/} \rightarrow \calC_{p'/}$ induces a homotopy equivalence of Kan complexes, after passing to the fiber over any object of $\calC$. Since
both $\calC_{p/}$ and $\calC_{p'/}$ are left-fibered over $\calC$, Corollary \ref{usefir} implies that
$\calC_{p/} \rightarrow \calC_{p'/}$ is a categorical equivalence. This proves that $(1) \Rightarrow (2)$.

Now suppose that $(2)$ is satisfied, and let $\overline{p}: K^{\triangleright} \rightarrow \calC$
be as in $(3)$. Then we may identify $\overline{p}$ with an initial object of the $\infty$-category
$\calC_{p/}$. The induced map $\calC_{p/} \rightarrow \calC_{p'/}$ is an equivalence, and therefore carries the initial object $\overline{p}$ to an initial object $\overline{p}'$ of $\calC_{p'/}$; thus
$\overline{p}'$ is a colimit of $p'$. This proves that $(2) \Rightarrow (3)$.

It remains to prove that $(3) \Rightarrow (1)$. For this, we make use of the theory of 
classifying right fibrations (\S \ref{universalfib}). Let $X \rightarrow K$ be a right fibration. We wish to show that composition with $v$ induces a homotopy equivalence $\bHom_{K}(K,X) \rightarrow \bHom_{K}(K',X)$. It will suffice to prove this result after replacing $X$ by any equivalent right fibration. Let $\SSet$ denote the $\infty$-category of spaces. According to Corollary \ref{unipull}, there is a classifying map $p: K \rightarrow \SSet^{op}$ and an equivalence of right fibrations between $X$ and $(\SSet_{\ast/})^{op} \times_{\SSet^{op}} K$, where $\ast$ denotes a final object
of $\SSet$.

The $\infty$-category $\SSet$ admits small limits (Corollary \ref{limitsinmodel}). It follows that there exists a map
$\overline{p}: K^{\triangleright} \rightarrow \SSet^{op}$ which is a colimit of $p = \overline{p}|K$. Let $x$
denote the image in $\SSet$ of the cone point of $K^{\triangleright}$. Let $\overline{p}': {K'}^{\triangleright} \rightarrow \SSet^{op}$ be the induced map. Then, by hypothesis, $\overline{p}'$ is a colimit of 
$p' = \overline{p}'|K'$. According to Lemma \ref{cogh}, there is a (natural) chain of weak homotopy equivalences relating $\bHom_{K}(K,X)$ with $(\SSet^{op})_{p/} \times_{ \SSet^{op}} \{y\}$. 
Similarly, there is a chain of weak homotopy equivalences connecting $\bHom_{K}(K',X)$ with
$(\SSet^{op})_{p'/} \times _{\SSet^{op}} \{y\}$. Consequently, we are reduced to proving that the left vertical map in the diagram
$$ \xymatrix{ (\SSet^{op})_{p/} \times_{\SSet^{op}} \{y\} \ar[d] & (\SSet^{op})_{\overline{p}/} \times_{\SSet^{op}} \{y\} \ar[l] \ar[r] \ar[d] & (\SSet^{op})_{x/} \times_{ \SSet^{op} } \{y\} \ar[d] \\
(\SSet^{op})_{p'/} \times_{\SSet^{op}} \{y\} & (\SSet^{op})_{\overline{p}'/} \times_{\SSet^{op}} \{y\} \ar[l] \ar[r]  & (\SSet^{op})_{x/} \times_{ \SSet^{op} } \{y\} } $$
is a homotopy equivalence. Since $\overline{p}$ and $\overline{q}$ are colimits of $p$ and $q$, the left horizontal maps are trivial fibrations. Since the inclusions of the cone points into $K^{\triangleright}$ and
${K'}^{\triangleright}$ are right anodyne, the right horizontal maps are also trivial fibrations.
It therefore suffices to prove that the right vertical map is a homotopy equivalence. But this map is an isomorphism of simplicial sets.
\end{proof}

\begin{corollary}\label{stoog}
Let $p: K \rightarrow K'$ be a map of simplicial sets, and
$q: K' \rightarrow K''$ a categorical equivalence. Then $p$
is cofinal if and only if $q \circ p$ is cofinal. In particular, $($taking $p = \id_{S'}${}$)$
$q$ itself is cofinal.
\end{corollary}

\begin{proof}
Let $\calC$ be an $\infty$-category, $r'': K'' \rightarrow \calC$ a diagram, and
set $r' = r'' \circ q$, $r = r' \circ p$. Since $q$ is a categorical equivalence, $\calC_{r''/}
\rightarrow \calC_{r'/}$ is a categorical equivalence. It follows that $\calC_{r/} \rightarrow \calC_{r''/}$ is a categorical equivalence if and only if $\calC_{r/} \rightarrow \calC_{r'/}$ is a categorical equivalence. We now apply the characterization $(2)$ of Proposition \ref{gute}.
\end{proof}

\begin{corollary}\label{cofinv}
The property of cofinality is homotopy invariant. In other words,
if two maps $f,g: K \rightarrow K'$ have the same image in the
homotopy category of $\sSet$ obtained by inverting all categorical
equivalences, then $f$ is cofinal if and only if $g$ is cofinal.
\end{corollary}

\begin{proof}
Choose a categorical equivalence $K' \rightarrow \calC$, where $\calC$ is an $\infty$-category.
In view of Corollary \ref{stoog}, we may replace $K'$ by $\calC$ and thereby assume that
$K'$ is itself an $\infty$-category. Since $f$ and $g$ are homotopic, there exists a cylinder
object $S$ equipped with a trivial fibration $p: S
\rightarrow K$, a map $q: S \rightarrow \calC$, and two sections
$s,s': K \rightarrow S$ of $p$, such that $f = q \circ s$, $g = q
\circ s'$. Since $p$ is a categorical equivalence, so is every
section of $p$. Consequently, $s$ and $s'$ are cofinal. We now
apply Proposition \ref{cofbasic} to deduce that $f$ is cofinal if
and only if $q$ is cofinal. Similarly, $g$ is cofinal if and only
if $q$ is cofinal.
\end{proof}

\begin{corollary}\label{twork}
Let $p: X \rightarrow S$ be a map of simplicial sets. The following are equivalent:
\begin{itemize}
\item[$(1)$] The map $p$ is a cofinal right fibration.
\item[$(2)$] The map $p$ is a trivial fibration.
\end{itemize}
\end{corollary}

\begin{proof}
Clearly any trivial fibration is a right fibration. Furthermore, any trivial fibration is
a categorical equivalence, hence cofinal by Corollary \ref{stoog}. Thus $(2)$ implies $(1)$. Conversely, suppose that $p$ is a cofinal right fibration. Since $p$ is cofinal, the natural map
$\bHom_{S}(S,X) \rightarrow \bHom_{S}(X,X)$ is a homotopy equivalence of Kan complexes.
In particular, there exists a section $f: S \rightarrow X$ of $p$ such that
$f \circ p$ is (fiberwise) homotopic to the identity map of $X$. Consequently, for each
vertex $s$ of $S$, the fiber $X_{s} = X \times_{S} \{s\}$ is a contractible Kan complex
(since the identity map $X_{s} \rightarrow X_{s}$ is homotopic to the constant map with value $f(s)$). The dual of Lemma \ref{toothie} now shows that $p$ is a trivial fibration.
\end{proof}

\begin{corollary}\label{stufe}
A map $X \rightarrow Z$ of simplicial sets is cofinal if and only
if it admits a factorization $$X \stackrel{f}{\rightarrow} Y
\stackrel{g}{\rightarrow} Z,$$ where $X \rightarrow Y$ is
right-anodyne and $Y \rightarrow Z$ is a trivial fibration.
\end{corollary}

\begin{proof}
The ``if'' direction is clear: if such a factorization exists,
then $f$ is cofinal (since it is right anodyne), $g$ is cofinal
(since it is a categorical equivalence), and consequently $g \circ
f$ is cofinal (since it is a composition of cofinal maps).

For the ``only if'' direction, let us suppose that $X \rightarrow
Z$ is a cofinal map. By the small object argument (Proposition \ref{quillobj}), there is a
factorization $$X \stackrel{f}{\rightarrow} Y
\stackrel{g}{\rightarrow} Z$$ where $f$ is right-anodyne and $g$
is a right fibration. The map $g$ is cofinal by Proposition \ref{cofbasic}, and therefore a trivial fibration by Corollary \ref{twork}.
\end{proof}

\begin{corollary}\label{prodcofinal}
Let $p: S \rightarrow S'$ be a cofinal map, and $K$ any simplicial
set. Then the induced map $K \times S \rightarrow K \times S'$ is
cofinal.
\end{corollary}

\begin{proof}
Using Corollary \ref{stufe}, we may suppose that $p$ is either
right anodyne or a trivial fibration. Then the induced map $K \times S \rightarrow K \times S'$ has the same property.
\end{proof}

\subsection{Smoothness and Right Anodyne Maps}\label{smoothness}

In this section, we explain how to characterize the classes of right anodyne and cofinal morphisms in terms of the contravariant model structures studied in \S \ref{contrasec}. We also introduce a third class of maps between simplicial sets, which we call {\it smooth}.

We begin with the following characterization of right anodyne maps:

\begin{proposition}\label{hunef}\index{gen}{right anodyne}
Let $i: A \rightarrow B$ be a map of simplicial sets. The following conditions are equivalent:
\begin{itemize}
\item[$(1)$] The map $i$ is right anodyne.
\item[$(2)$] For any map of simplicial sets $j: B \rightarrow C$, the map $i$ is a trivial cofibration with respect
to the contravariant model structure on $(\sSet)_{/C}$.
\item[$(3)$] The map $i$ is a trivial cofibration with respect to the contravariant model structure on
$(\sSet)_{/B}$.
\end{itemize}
\end{proposition}

\begin{proof}
The implication $(1) \Rightarrow (2)$ follows immediately from Proposition \ref{onehalf}, and the implication $(2) \Rightarrow (3)$ is obvious. Suppose that $(3)$ holds. To prove $(1)$, it suffices to show that given any diagram
$$ \xymatrix{ A \ar@{^{(}->}[d]^{i} \ar[r] & X \ar[d]^p \\
B \ar[r] \ar@{-->}[ur]^{f} & Y }$$
such that $p$ is a right fibration, one can supply the dotted arrow $f$ as indicated. Replacing
$p: X \rightarrow Y$ by the pullback $X \times_{Y} B \rightarrow B$, we may reduce to the case where
$Y = B$. Corollary \ref{usewhere1} implies that $X$ is a fibrant object of $(\sSet)_{/B}$ (with respect to contravariant model structure) so that the desired map $f$ can be found.
\end{proof}

\begin{corollary}\label{nonobcomp}
Suppose given maps $A \stackrel{i}{\rightarrow} B \stackrel{j}{\rightarrow} C$ of simplicial
sets. If $i$ and $j \circ i$ are right anodyne, and $j$ is a cofibration, then $j$ is right-anodyne.
\end{corollary}

\begin{proof}
By Proposition \ref{hunef}, $i$ and $j \circ i$ are contravariant equivalences in $(\sSet)_{/C}$. It follows that $j$ is a trivial cofibration in $(\sSet)_{/C}$, so that $j$ is right anodyne (by Proposition \ref{hunef} again).
\end{proof}

\begin{corollary}\label{anothernonob}
Let $$ \xymatrix{ A' \ar[d]^{f'} & A \ar[r] \ar[l]^{u} \ar[d]^{f} & A'' \ar[d]^{f''} \\
B' & B \ar[r] \ar[l]^{v} & B'' }$$
be a diagram of simplicial sets. Suppose that $u$ and $v$ are monomorphisms,
and that $f, f'$, and $f''$ are right anodyne. Then the induced map
$$ A' \coprod_{A} A'' \rightarrow B' \coprod_{B} B''$$
is right anodyne.
\end{corollary}

\begin{proof}
According to Proposition \ref{hunef}, each of the maps $f$, $f'$, and $f''$ is a contravariant
equivalence in $(\sSet)_{/B' \coprod_{B} B''}$. The assumption on $u$ and $v$ guarantees that
$f' \coprod_{f} f''$ is also a contravariant equivalence in $(\sSet)_{/B' \coprod_{B} B''}$, so that
$f' \coprod_{f} f''$ is right anodyne by Proposition \ref{hunef} again.
\end{proof}

\begin{corollary}\label{filtanodyne}
The collection of right anodyne maps of simplicial sets is stable under filtered colimits.
\end{corollary}

\begin{proof}
Let $f: A \rightarrow B$ be a filtered colimit of right anodyne morphisms $f_{\alpha}: A_{\alpha} \rightarrow B_{\alpha}$. According to Proposition \ref{hunef}, each $f_{\alpha}$ is a contravariant equivalence in $(\sSet)_{/B}$. Since contravariant equivalences are stable under filtered colimits, we conclude that $f$ is a contravariant equivalence in $(\sSet)_{/B}$ so that $f$ is right anodyne by Proposition \ref{hunef}.
\end{proof}

Proposition \ref{hunef} has an analogue for cofinal maps:

\begin{proposition}\label{huneff}
Let $i: A \rightarrow B$ be a map of simplicial sets. The following conditions are equivalent:
\begin{itemize}
\item[$(1)$] The map $i$ cofinal.
\item[$(2)$] For any map $j: B \rightarrow C$, the inclusion $i$ is a contravariant
equivalence in $(\sSet)_{/C}$.
\item[$(3)$] The map $i$ is a contravariant equivalence in
$(\sSet)_{/B}$.
\end{itemize}
\end{proposition}

\begin{proof}
Suppose $(1)$ is satisfied. By Corollary \ref{stufe}, $i$ admits a factorization as a right anodyne map followed by a trivial fibration. Invoking Proposition \ref{hunef}, we conclude that $(2)$ holds. 
The implication $(2) \Rightarrow (3)$ is obvious. If $(3)$ holds, then we can choose a factorization
$$ A \stackrel{i'}{\rightarrow} A' \stackrel{i''}{\rightarrow} B$$ of $i$, where $i'$ is right anodyne and $i''$ is a right fibration. Then $i''$ is a contravariant fibration (in $\sSet_{/B}$) and a contravariant weak equivalence, and is therefore a trivial fibration of simplicial sets. 
We now apply Corollary \ref{stufe} to conclude that $i$ is cofinal.
\end{proof}

\begin{corollary}\label{weakcont}
Let $p: X \rightarrow S$ be a map of simplicial sets, where $S$ is a Kan complex. Then
$p$ is cofinal if and only if it is a weak homotopy equivalence.
\end{corollary}

\begin{proof}
By Proposition \ref{huneff}, $p$ is cofinal if and only if it is a contravariant equivalence
in $(\sSet)_{/S}$. If $S$ is a Kan complex, then Proposition \ref{strstr} asserts that
the contravariant equivalences are precisely the weak homotopy equivalences.
\end{proof}

Let $p: X \rightarrow Y$ be an arbitrary map of simplicial sets. In \S \ref{contrasec} we showed that $p$ induces a Quillen adjunction $(p_{!}, p^{\ast})$ between the contravariant model categories
$(\sSet)_{/X}$ and $(\sSet)_{/Y}$. The functor $p^{\ast}$ itself has a right adjoint, which we will denote by $p_{\ast}$; it is given by
$$ p_{\ast}(M) =  \bHom_{Y}(X,M).$$
The adjoint functors $(p^{\ast}, p_{\ast})$ are not Quillen adjoints in general. Instead we have:

\begin{proposition}\label{smoothdef}
Let $p: X \rightarrow Y$ be a map of simplicial sets. The following conditions are equivalent:

\begin{itemize}
\item[$(1)$] For any right-anodyne map $i: A \rightarrow B$ in $(\sSet)_{/Y}$, the induced map
$A \times_Y X \rightarrow B \times_{Y} X$ is right-anodyne.

\item[$(2)$] For every Cartesian diagram
$$\xymatrix{ X' \ar[r] \ar[d]^{p'} & X \ar[d]^{p} \\
Y' \ar[r] & Y, },$$ the functor ${p'}^{\ast}: (\sSet)_{/Y'} \rightarrow (\sSet)_{/X'}$ preserves contravariant equivalences. 

\item[$(3)$] For every Cartesian diagram
$$\xymatrix{ X' \ar[r] \ar[d]^{p'} & X \ar[d]^{p} \\
Y' \ar[r] & Y, },$$ the adjoint functors $( {p'}^{\ast}, p'_{\ast})$ give rise to a Quillen adjunction between the contravariant model categories $(\sSet)_{/Y'}$ and $(\sSet)_{/X'}$. 
\end{itemize}

\end{proposition}

\begin{proof}
Suppose that $(1)$ is satisfied; let us prove $(2)$. Since property $(1)$ is clearly stable under base change, we may suppose that $p' = p$. Let $u: M \rightarrow N$ be a contravariant equivalence in
$(\sSet)_{/Y}$. If $M$ and $N$ are fibrant, then $u$ is a homotopy equivalence, so that $p^{\ast}(u): p^{\ast} M \rightarrow p^{\ast} N$ is also a homotopy equivalence. In the general case, we may select a diagram
$$ \xymatrix{ M \ar[r]^i \ar[d]^{u} & M' \ar[d] \ar[dr]^{v} \\
N \ar[r]^-{i'} & N \coprod_{M} M' \ar[r]^-{j} &  N' } $$
where $M'$ and $N'$ are fibrant, and the maps $i$ and $j$ are right anodyne (and therefore $i'$ is also right anodyne). Then $p^{\ast}(v)$ is a contravariant equivalence, while the maps
$p^{\ast}(i)$, $p^{\ast}(j)$, and $p^{\ast}(i')$ are all right anodyne; by Proposition \ref{hunef} they are contravariant equivalences as well. It follows that $p^{\ast}(u)$ is a contravariant equivalence.

To prove $(3)$, it suffices to show that ${p'}^{\ast}$ preserves cofibrations and trivial cofibrations. The first statement is obvious, and the second follows immediately from $(2)$. Conversely the existence of a Quillen adjunction $({p'}^{\ast}, p_{\ast})$ implies that ${p'}^{\ast}$ preserves contravariant equivalences between cofibrant objects. Since every object of $(\sSet)_{/Y'}$ is cofibrant, we deduce that $(3)$ implies $(2)$. 

Now suppose that $(2)$ is satisfied, and let $i: A \rightarrow B$ be a right-anodyne map in $(\sSet)_{/Y}$ as in $(1)$. Then $i$ is a contravariant equivalence in $(\sSet)_{/B}$. Let $p': X \times_{Y} B \rightarrow B$ be base change of $p$; then $(2)$ implies that the induced map
$i': {p'}^{\ast} A \rightarrow {p'}^{\ast} B$ is a contravariant equivalence in $(\sSet)_{/B \times_{Y} X}$. By Proposition \ref{hunef}, the map $i'$ is right anodyne. Now we simply note that $i'$ may be identified with the map $A \times_Y X \rightarrow B \times_{Y} X$ in the statement of $(1)$.
\end{proof}

\begin{definition}\index{gen}{smooth}
We will say that a map $p: X \rightarrow Y$ of simplicial sets is {\em smooth} if it satisfies the (equivalent) conditions of Proposition \ref{smoothdef}.
\end{definition}

\begin{remark}\label{gonau}
Let 
$$ \xymatrix{ X' \ar[d] \ar[r]^{f'} & X \ar[d]^{p} \\
S' \ar[r]^{f} & S }$$
be a pullback diagram of simplicial sets. Suppose that $p$ is smooth and that $f$ is cofinal. Then $f'$ is cofinal: this follows immediately from characterization $(2)$ of Proposition \ref{smoothdef} and characterization $(3)$ of Proposition \ref{huneff}.
\end{remark}

We next give an alternative characterization of smoothness. Let
$$ \xymatrix{ X' \ar[d]^{p'} \ar[r]^{q'} & X \ar[d]^{p} \\
Y' \ar[r]^{q} & Y }$$
be a Cartesian diagram of simplicial sets. Then we obtain an isomorphism
$R {p'}^{\ast} R q^{\ast} \simeq R {q'}^{\ast} R p^{\ast}$ of right-derived functors, which induces
a natural transformation
$$ \psi_{p,q}: L q'_{!} R {p'}^{\ast} \rightarrow R p^{\ast} L q_{!}.$$

\begin{proposition}\label{smoothbase}
Let $p: X \rightarrow Y$ be a map of simplicial sets. The following conditions are equivalent:
\begin{itemize}
\item[$(1)$] The map $p$ is smooth.
\item[$(2)$] For every Cartesian rectangle
$$ \xymatrix{ X'' \ar[d]^{p''} \ar[r]^{q'} & X' \ar[d]^{p'} \ar[r] & X \ar[d]^{p} \\
Y'' \ar[r]^{q} & Y' \ar[r] & Y, }$$ the natural transformation
$\psi_{p',q}$ is an isomorphism of functors from the homotopy category of
$(\sSet)_{/Y''}$ to the homotopy category of $(\sSet)_{/X'}$ (here all categories are
endowed with the contravariant model structure).
\end{itemize}
\end{proposition}

\begin{proof}
Suppose that $(1)$ is satisfied, and consider any Cartesian rectangle as in $(2)$. Since $p$ is smooth, $p'$ and $p''$ are also smooth. It follows that ${p'}^{\ast}$ and ${p''}^{\ast}$ preserve weak equivalences, so they may be identified with their right derived functors. Similarly, $q_{!}$ and $q'_{!}$ preserve weak equivalences, so they may be identified with their left derived functors. Consequently, the natural transformation $\psi_{p',q}$ is simply obtained by passage to the homotopy category from the natural transformation
$$ q'_{!} {p''}^{\ast} \rightarrow {p'}^{\ast} q_{!}.$$
But this is an isomorphism of functors before passage to the homotopy categories.

Now suppose that $(2)$ is satisfied. Let $q: Y'' \rightarrow Y'$ be a right-anodyne map
in $(\sSet)_{/Y}$, and form the Cartesian square as in $(2)$. Let us compute the value of the functors $L q'_{!} R {p''}^{\ast}$ and $R {p'}^{\ast} L q_{!}$ on the object $Y''$ of $(\sSet)_{/Y''}$. The composite $L q'_{!} R {p''}^{\ast}$ is easy: because $Y''$ is fibrant and
$X'' = {p''}^{\ast} Y''$ is cofibrant, the result is $X''$, regarded as an object of $(\sSet)_{/X'}$. The other composition is slightly trickier: $Y''$ is cofibrant, but $q_{!} Y''$ is not fibrant when viewed as an object of $(\sSet)_{/Y'}$. However, in view of the assumption that $q$ is right anodyne, Proposition \ref{hunef} ensures that $Y'$ is a fibrant replacement for $q_{!} Y'$; thus we may identify  $R {p'}^{\ast} L q_{!}$ with the object ${p'}^{\ast} Y' = X'$ of $(\sSet)_{/ X'}$. Condition $(2)$ now implies that the natural map $X'' \rightarrow X'$ is a contravariant equivalence in $(\sSet)_{/X'}$. Invoking Proposition \ref{hunef}, we deduce that $q'$ is right anodyne, as desired.
\end{proof}

\begin{remark}
The terminology ``smooth'' is suggested by the analogy of Proposition \ref{smoothbase} with the {\em smooth base change theorem} in the theory of \'{e}tale cohomology (see, for example, \cite{freitag}). 
\end{remark}

\begin{proposition}\label{usenonob}
Suppose given a commutative diagram
$$ \xymatrix{ X \ar[r]^{i} \ar[dr]^{p} \ar[d] & X' \ar[d]^{p'} \\
X'' \ar[r]^{p''} & S }$$
of simplicial sets. Assume that $i$ is a cofibration, and that $p,p'$, and $p''$ are smooth. Then
the induced map $X' \coprod_{X} X'' \rightarrow S$ is smooth.
\end{proposition}

\begin{proof}
This follows immediately from Corollary \ref{anothernonob} and characterization $(1)$ of Proposition \ref{smoothdef}.
\end{proof}

\begin{proposition}\label{usefiltanodyne}
The collection of smooth maps $p: X \rightarrow S$ is stable under filtered colimits in
$(\sSet)_{/S}$. 
\end{proposition}

\begin{proof}
Combine Corollary \ref{filtanodyne} with characterization $(1)$ of Proposition \ref{smoothdef}.
\end{proof}

\begin{proposition}\label{strokhop}\index{gen}{coCartesian fibration!and smoothness}
Let $p: X \rightarrow S$ be a coCartesian fibration of simplicial sets. Then $p$ is smooth.
\end{proposition}

\begin{proof}
Let $i: B' \rightarrow B$ be a right anodyne map in $(\sSet)_{/S}$; we wish to show that the induced map $B' \times_{S} X \rightarrow B \times_{S} X$ is right anodyne. By general nonsense, we may reduce ourselves to the case where $i$ is an inclusion $\Lambda^n_i \subseteq \Delta^n$ where
$0 <  i \leq n$. Making a base change, we may suppose that $S = B$. By Proposition \ref{simplexplay}, there exists a composable sequence of maps
$$ \phi: A^0 \rightarrow \ldots \rightarrow A^n $$ and a quasi-equivalence
$M^{op}(\phi) \rightarrow X$. Consider the diagram
$$ \xymatrix{ M^{op}(\phi) \times_{\Delta^n} \Lambda^n_i \ar@{^{(}->}[d] \ar[r] \ar[dr]^{f} & 
X \times_{\Delta^n} \Lambda^n_i \ar@{^{(}->}[d]^{h} \\
M^{op}(\phi) \ar[r]^{g} & X }$$
The left vertical map is right-anodyne, since it is a pushout of the inclusion
 $A^0 \times \Lambda^n_i \subseteq A^0 \times \Delta^n$. It follows that $f$ is cofinal, being a composition of a right-anodyne map and a categorical equivalence. Since $g$ is cofinal (being a categorical equivalence) we deduce from Proposition \ref{cofbasic} that $h$ is cofinal. Since
$h$ is a monomorphism of simplicial sets, it is right-anodyne by Proposition \ref{cofbasic}.
\end{proof}

\begin{proposition}\label{longwait5}\index{gen}{bifibration!and smoothness}
Let $p: X \rightarrow S \times T$ be a bifibration. Then the composite map
$\pi_{S} \circ p: X \rightarrow S$ is smooth.
\end{proposition}

\begin{proof}
For every map $T' \rightarrow T$, let $X_{T'} = X \times_{T} T'$. We note that
$X$ is a filtered colimit of $X_{T'}$, as $T'$ ranges over the finite simplicial subsets
of $T$. Using Proposition \ref{usefiltanodyne}, we can reduce to the case where $T$ is finite.
Working by induction on the dimension and the number of nondegenerate simplices of
$T$, we may suppose that $T = T' \coprod_{ \bd \Delta^n } \Delta^n$, where the result is known
for $T'$ and for $\bd \Delta^n$. Applying Proposition \ref{usenonob}, we can reduce to the case $T = \Delta^n$. We now apply Lemma \ref{gork} to deduce that $p$ is a coCartesian fibration, and therefore smooth by Proposition \ref{strokhop}.
\end{proof}

\begin{lemma}\label{covg}
Let $\calC$ be an $\infty$-category containing an object $C$, and let
$f: X \rightarrow Y$ be a covariant equivalence in $(\sSet)_{/\calC}$. The induced map
$$ X \times_{\calC} \calC^{/C} \rightarrow Y \times_{\calC} \calC^{/C}$$ is also a covariant equivalence in $\calC^{/C}$.
\end{lemma}

\begin{proof}
It will suffice to prove that for every object $Z \rightarrow \calC$ of $(\sSet)_{/\calC}$, the fiber
product $Z \times_{\calC} \calC^{/C}$ is a homotopy product of $Z$ with
$\calC^{/C}$ in $(\sSet)_{/\calC}$ (with respect to the covariant model structure). Choose a factorization
$$ Z \stackrel{i}{\rightarrow} Z' \stackrel{j}{\rightarrow} \calC,$$
where $i$ is left anodyne and $j$ is a left fibration. According to Corollary \ref{usewhere1}, we may regard $Z'$ as a fibrant replacement for $Z$ in $(\sSet)_{/\calC}$. It therefore suffices to prove
that the map $i': Z \times_{\calC} \calC^{/C} \rightarrow Z' \times_{\calC} \calC^{/C}$ is a covariant equivalence. According to Proposition \ref{huneff}, it will suffice to prove that $i'$ is left anodyne.
The map $i'$ is a base change of $i$ by the projection $p: \calC^{/C} \rightarrow \calC$; it therefore suffices to prove that $p^{op}$ is smooth. This follows from Proposition \ref{strokhop}, since
$p$ is a right fibration of simplicial sets.
\end{proof}

\begin{proposition}\label{longwait44}
Let $\calC$ be an $\infty$-category, and
$$\xymatrix{ X \ar[rr]^{f} \ar[dr]^{p} & & Y \ar[dl]_{q} \\
& \calC & }$$ be a commutative diagram of simplicial sets. Suppose that
$p$ and $q$ are smooth. The following conditions are equivalent:
\begin{itemize}
\item[$(1)$] The map $f$ is a covariant equivalence in $(\sSet)_{/\calC}$.
\item[$(2)$] For each object $C \in \calC$, the induced map of fibers
$X_{C} \rightarrow Y_{C}$ is a weak homotopy equivalence.
\end{itemize} 
\end{proposition}

\begin{proof}
Suppose that $(1)$ is satisfied, and let $C$ be an object of $\calC$.
We have a commutative diagram of simplicial sets
$$ \xymatrix{ X_{C} \ar[r] \ar[d] & Y_{C} \ar[d] \\
X \times_{\calC} \calC^{/C} \ar[r] & Y \times_{\calC} \calC^{/C}.} $$
Lemma \ref{covg} implies that the bottom horizontal map is a covariant equivalence. The vertical maps are both pullbacks of the right anodyne inclusion 
$ \{ C\} \subseteq \calC^{/C}$ along smooth maps, and are therefore right anodyne. In particular, the vertical arrows and the bottom horizontal arrow are all weak homotopy equivalences; it follows that the map $X_{C} \rightarrow Y_{C}$ is a weak homotopy equivalence as well.

Now suppose that $(2)$ is satisfied. Choose a commutative diagram
$$ \xymatrix{ X \ar[rr]^{f} \ar[d] & & Y \ar[d] \\
X' \ar[rr]^{f'} \ar[dr]^{p'}  & & Y' \ar[dl]^{q'} \\
& \calC & }$$
in $(\sSet)_{/\calC}$, where the vertical arrows are left anodyne and the maps
$p'$ and $q'$ are left fibrations. Using Proposition \ref{strokhop}, we conclude that
$p'$ and $q'$ are smooth. Applying $(1)$, we deduce that for each object $C \in \calC$,
the maps $X_{C} \rightarrow X'_{C}$ and $Y_{C} \rightarrow Y'_{C}$ are weak homotopy equivalences. It follows that each fiber $f'_{C}: X'_{C} \rightarrow Y'_{C}$ is a homotopy equivalence of Kan complexes, so that $f'$ is an equivalence of left fibrations and therefore
a covariant equivalence. Inspecting the above diagram, we deduce that $f$ is also a covariant equivalence, as desired.
\end{proof}

\subsection{Quillen's Theorem A for $\infty$-Categories}\label{quillA}

Suppose that $f: \calC \rightarrow \calD$ is a functor between $\infty$-categories, and that we wish to determine whether or not $f$ is cofinal. According to Proposition \ref{gute}, the cofinality of $f$ is equivalent to the assertion that for any diagram $p: \calD \rightarrow \calE$, $f$ induces an equivalence
$$ \injlim(p) \simeq \injlim(p \circ f).$$
One can always define a morphism
$$ \phi: \injlim(p \circ f) \rightarrow \injlim(p)$$
(provided that both sides are defined); the question is whether or not we can define an inverse $\psi = \phi^{-1}$. Roughly speaking, this involves defining a compatible family of maps
$\psi_{D}: p(D) \rightarrow \injlim(p \circ f)$, indexed by $D \in \calD$. The only reasonable
candidate for $\psi_{D}$ is a composition
$$ p(D) \rightarrow (p \circ f)(C) \rightarrow \injlim(p \circ f),$$
where the first map arises from a morphism $D \rightarrow f(C)$ in $\calC$. Of course, the existence of $C$ is not automatic. Moreover, even if $C$ exists, it may is usually not unique. The collection of candidates for $C$ is parametrized by the $\infty$-category $\calC^{D/} = \calC \times_{\calD} \calD^{D/}$. In order to make the above construction work, we need the $\infty$-category
$\calC^{D/}$ to be weakly contractible. More precisely, we will prove the following result:

\begin{theorem}[Joyal \cite{joyalnotpub}]\label{hollowtt}\index{gen}{Quillen's theorem A!for $\infty$-categories}
Let $f: \calC \rightarrow \calD$ be a map of simplicial sets, where $\calD$ is an $\infty$-category. The following
conditions are equivalent:
\begin{itemize}
\item[$(1)$] The functor $f$ is cofinal.
\item[$(2)$] For every object $D \in \calD$, the simplicial set
$\calC \times_{ \calD } \calD_{D/}$ is weakly contractible.
\end{itemize}
\end{theorem}

We first need to establish the following lemma:

\begin{lemma}\label{trull6prime}
Let $p: U \rightarrow S$ be a Cartesian fibration of simplicial sets. Suppose
that for every vertex $s$ of $S$, the fiber $X_{s} = p^{-1} \{s\}$ is weakly contractible. Then $p$ is cofinal.
\end{lemma}

\begin{proof}
Let $q: N \rightarrow S$ be a right fibration. For every map of simplicial sets $T \rightarrow S$,
let $X_{T} = \bHom_{S}(T,N)$ and $Y_{T} = \bHom_{S}(T \times_{S} U, N)$. Our goal is to prove that the natural map $X_{S} \rightarrow Y_{S}$ is a homotopy equivalence of Kan complexes.
We will prove, more generally, that for any map $T \rightarrow S$, the map
$\phi_{T}: Y_{T} \rightarrow Z_{T}$ is a homotopy equivalence. The proof goes by induction on the
(possibly infinite) dimension of $T$. Choose a transfinite sequence of simplicial
subsets $T(\alpha) \subseteq T$, where each $T(\alpha)$ is obtained from
$T(< \alpha) = \bigcup_{\beta < \alpha} T(\beta)$ by adjoining a single nondegenerate simplex
of $T$ (if such a simplex exists). We prove that $\phi_{T(\alpha)}$ is a homotopy equivalence
by induction on $\alpha$. Assuming that $\phi_{T(\beta)}$ is a homotopy equivalence for every
$\beta < \alpha$, we deduce that $\phi_{T(< \alpha)}$ is the homotopy inverse limit of a tower of equivalences, and therefore a homotopy equivalence. If $T(\alpha) = T(< \alpha)$, we are done. Otherwise, we may write $T(\alpha) = T(< \alpha) \coprod_{ \bd \Delta^n} \Delta^n$. Then
$\phi_{T(\alpha)}$ can be written as a homotopy pullback of $\phi_{T(< \alpha)}$ with
$\phi_{\Delta^n}$ over $\phi_{ \bd \Delta^n}$. The third map is a homotopy equivalence
by the inductive hypothesis. Thus, it suffices to prove that $\phi_{\Delta^n}$ is an equivalence.
In other words, we may reduce to the case $T = \Delta^n$.

By Proposition \ref{simplexplay}, there exists a composable sequence of maps
$$ \theta: A^0 \leftarrow \ldots \leftarrow A^n$$
and a quasi-equivalence $f: M(\theta) \rightarrow X \times_{S} T$, where
$M(\theta)$ denotes the mapping simplex of the sequence $\theta$. 
Given a map $T' \rightarrow T$, we let $Z_{T'} = \bHom_{S}(M(\theta) \times_{T} T', N)$. 
Proposition \ref{funkyfibcatfib} implies that $q$ is a categorical fibration. It follows that, for
any map $T' \rightarrow T$, the categorical equivalence 
$M(\theta) \times_{T} T' \rightarrow U \times_{S} T'$ induces another categorical equivalence
$\psi_{T'} = Y_{T'} \rightarrow Z_{T'}$. Since $Y_{T'}$ and $Z_{T'}$ are Kan complexes, the map $\psi_{T'}$ is a homotopy equivalence. Consequently, to prove that $\phi_{T}$ is an equivalence, it suffices to show that the composite map
$$ X_{T} \rightarrow Y_{T} \rightarrow Z_{T}$$ is an equivalence.

Consider the composition
$$ u: X_{ \Delta^{n-1} } \stackrel{u'}{\rightarrow} Z_{ \Delta^{n-1} } \stackrel{u''}{\rightarrow} \bHom_{S}( \Delta^{n-1} \times A^n, N) \stackrel{u'''}{\rightarrow} \bHom_{S}( \{n-1\} \times A^n, N)$$
Using the fact that $q$ is a right fibration and that $A^n$ is weakly contractible, we deduce that $u$ and $u'''$ are homotopy equivalences. The inductive hypothesis implies that $u'$ is a homotopy equivalence. Consequently, $u''$ is also a homotopy equivalence. 
The space $Z_{T}$ fits into a homotopy Cartesian diagram
$$ \xymatrix{ Z_{T} \ar[r] \ar[d]^{v''} & Z_{\Delta^{n-1}} \ar[d]^{u''} \\
\bHom_{S}( \Delta^n \times A^n,N ) \ar[r] & \bHom_{S}(\Delta^{n-1} \times A^n, N).}$$
It follows that $v''$ is a homotopy equivalence. Now consider the composition
$$ v: X_{\Delta^n} \stackrel{v'}{\rightarrow} Z_{\Delta^n} \stackrel{v''}{\rightarrow}
\bHom_{S}( \Delta^n \times A^n, N) \stackrel{v'''}{\rightarrow} \bHom_{S}( \{n\} \times A^n, N).$$
Again, because $q$ is a right fibration and $A^n$ is weakly contractible, the maps
$v$ and $v'''$ are homotopy equivalences. Since $v''$ is a homotopy equivalence, we deduce
that $v'$ is a homotopy equivalence, as desired.
\end{proof}

\begin{proof}[Proof of Theorem \ref{hollowtt}]
Using the small object argument, we can factor $f$ as a composition
$$ \calC \stackrel{f'}{\rightarrow} \calC' \stackrel{f''}{\rightarrow} \calD$$
where $f'$ is a categorical equivalence and $f''$ is an inner fibration. Then $f''$ is cofinal if and only if $f$ is cofinal (Corollary \ref{cofinv}). For every $D \in \calD$, the map
$\calD_{D/} \rightarrow \calD$ is a left fibration, so the induced map
$\calC_{D/} \rightarrow {\calC'}_{D/}$ is a categorical equivalence (Proposition \ref{basechangefunky}). Consequently, it will suffice to prove that $(1) \Leftrightarrow (2)$ for the morphism $f'': \calC' \rightarrow \calD$. In other words, we may assume that the simplicial set $\calC$ is an $\infty$-category.

Suppose first that $(1)$ is satisfied, and choose $D \in \calD$. The projection
$\calD_{D/} \rightarrow \calD$ is a left fibration, and therefore smooth (Proposition \ref{strokhop}). Applying Remark \ref{gonau}, we deduce that the projection
$\calC \times_{\calD} \calD_{D/} \rightarrow \calD_{D/}$ is cofinal, and therefore a weak homotopy equivalence (Proposition \ref{cofbasic}). Since $\calD_{D/}$ has an initial object, it is weakly contractible. Therefore $\calC \times_{\calD} \calD_{D/}$ is weakly contractible, as desired.

We now prove that $(2) \Rightarrow (1)$.
Let $\calM = \Fun(\Delta^1,\calD) \times_{ \Fun( \{1\}, \calD )} \calC$. Then the map $f$ factors as a composition
$$ \calC \stackrel{f'}{\rightarrow} \calM \stackrel{f''}{\rightarrow} \calD$$
where $f'$ is the obvious map and $f''$ is given by evaluation at the vertex $\{0\} \subseteq
\Delta^1$. Note that there is a natural projection map $\pi: \calM \rightarrow \calC$, that
$f'$ is a section of $\pi$, and that there is a simplicial homotopy
$h: \Delta^1 \times \calM \rightarrow \calM$ from $\id_{\calM}$ to $f' \circ \pi$
which is compatible with the projection to $\calC$. It follows from
Proposition \ref{trull11} that $f'$ is right anodyne.

Corollary \ref{tweezegork} implies that $f''$ is a Cartesian fibration. 
The fiber of $f''$ over an object $D \in \calD$ is isomorphic to $\calC \times_{\calD} \calD^{D/}$, 
which is equivalent to $\calC \times_{\calD} \calD_{D/}$ and therefore weakly contractible (Proposition \ref{certs}). By assumption, the fibers of $f''$ are weakly contractible.
Lemma \ref{trull6prime} asserts that $f''$ is cofinal. It follows that $f$, as a composition of cofinal maps, is also cofinal.
\end{proof}

Using Theorem \ref{hollowtt} we can easily deduce the following classical result of Quillen:

\begin{corollary}[Quillen's Theorem A]\index{gen}{Quillen's theorem A}
Let $f: \calC \rightarrow \calD$ be a functor between ordinary categories. Suppose that, for
every object $D \in \calD$, the fiber product category $\calC \times_{\calD} \calD_{D/}$ has
weakly contractible nerve. Then $f$ induces a weak homotopy equivalence of simplicial sets
$\Nerve(\calC) \rightarrow \Nerve(\calD)$.
\end{corollary}

\begin{proof}
The assumption implies that $\Nerve(f): \Nerve(\calC) \rightarrow \Nerve(\calD)$ satisfies the hypotheses of Theorem \ref{hollowtt}. It follows that $\Nerve(f)$ is a cofinal map of simplicial sets, and therefore a weak homotopy equivalence (Proposition \ref{cofbasic}).
\end{proof}

\section{Techniques for Computing Colimits}\label{c4s2}

\setcounter{theorem}{0}

In this section, we will introduce various techniques for computing, analyzing, and manipulating colimits. We begin in \S \ref{quasilimit2} by introducing a variant on the join construction of \S \ref{langur}. The new join construction is (categorically) equivalent to the version we are already familiar with, but has better formal behavior in some contexts. For example, they permit us to define a {\em parametrized} version of overcategories and undercategories, which we will analyze in \S \ref{consweet}.

In \S \ref{quasilimit1}, we address the following question: given a diagram $p: K \rightarrow \calC$ and a decomposition of $K$ into ``pieces'', how is the colimit
$\injlim(p)$ related to the colimits of those pieces? For example, if $K = A \cup B$, then it seems reasonable expect an equation of the form
$$\colim(p) = ( \colim p|A ) \coprod_{ \colim(p| A \cap B)} (\colim p|B).$$ 
Of course there are many variations on this theme; we will lay out a general framework in \S \ref{quasilimit1}, and apply it to specific situations in \S \ref{coexample}.

Although the $\infty$-categorical theory of colimits is elegant and powerful, it can be be difficult to work with because the colimit $\colim(p)$ of a diagram $p$ is only well-defined up to equivalence.
This problem can sometimes be remedied by working in the more rigid theory of model categories, where the notion of $\infty$-categorical colimit should be replaced by the notion of {\it homotopy colimit} (see \S \ref{quasilimit3}). In order to pass smoothly between these two settings, we need to know that the $\infty$-categorical theory of colimits agrees with the more classical theory of homotopy colimits. A precise statement of this result (Theorem \ref{colimcomparee}) will be formulated and proven in \S \ref{quasilimit4}. 

\subsection{Alternative Join and Slice Constructions}\label{quasilimit2}

In \S \ref{join}, we introduced the {\em join} functor $\star$ on simplicial sets. In \cite{joyalnotpub}, Joyal introduces a closely related operation $\diamond$ on simplicial sets. This operation is equivalent to $\star$ (Proposition \ref{rub3}) but is more technically convenient in certain contexts. In this section we will review the definition of the operation $\diamond$ and to establish some of its basic properties (see also \cite{joyalnotpub} for a discussion).

\begin{definition}[\cite{joyalnotpub}]\index{not}{diamond@$\diamond$}
Let $X$ and $Y$ be simplicial sets. The simplicial set $X \diamond Y$ is defined to be pushout
$$ X \coprod_{ X \times Y \times \{0\} } (X \times Y \times \Delta^1) \coprod_{X \times Y \times \{1\} } Y.$$
\end{definition}

We note that since $X \times Y \times (\bd \Delta^1) \rightarrow X \times Y \times \Delta^1$ is a monomorphism, the pushout diagram defining $X \diamond Y$ is a homotopy pushout in $\sSet$ (with respect to the Joyal model structure). Consequently, we deduce that categorical equivalences $X \rightarrow X'$, $Y \rightarrow Y'$ induce a categorical equivalence $X \diamond Y \rightarrow X' \diamond Y'$.

The simplicial set $X \diamond Y$ admits a map $p: X \diamond Y \rightarrow \Delta^1$, with
$X \simeq p^{-1} \{0\}$ and $Y \simeq p^{-1} \{1\}$. Consequently, there is a unique map
$X \diamond Y \rightarrow X \star Y$ which is compatible with the projection to $\Delta^1$ and induces the identity maps on $X$ and $Y$.

\begin{proposition}\label{rub3}
For any simplicial sets $X$ and $Y$, the natural map $\phi: X
\diamond Y \rightarrow X \star Y$ is a categorical equivalence.
\end{proposition}

\begin{proof}
Since both sides are compatible with the formation of filtered
colimits in $X$, we may suppose that $X$ contains only finitely
many nondegenerate simplices. If $X$ is empty, then $\phi$ is an isomorphism and the result is
obvious. Working by induction on the dimension of $X$ and the
number of nondegenerate simplices in $X$, we may write
$$ X = X' \coprod_{\bd \Delta^{n}} \Delta^n,$$ and we may assume
that the statement is known for the pairs $(X',Y)$ and $(\bd
\Delta^n, Y)$. Since the Joyal model structure on $\sSet$ is
left-proper, we have a map of homotopy pushouts
$$ (X' \diamond Y) \coprod_{ \bd \Delta^n \diamond Y } (\Delta^n
\diamond Y) \rightarrow (X' \star Y) \coprod_{ \bd \Delta^n
\star Y} ( \Delta^n \star Y),$$ and we are therefore reduced
to proving the assertion in the case where $X = \Delta^n$.
The inclusion $$ \Delta^{ \{0,1\} } \coprod_{ \{1\} } \ldots \coprod_{ \{n-1\} }
\Delta^{ \{n-1,n\} } \subseteq \Delta^n$$ is inner anodyne. Thus, if
$n > 1$, we can conclude by induction. Thus we may suppose that
$X = \Delta^0$ or $X = \Delta^1$. By a similar argument, we may reduce to the case where
$Y = \Delta^0$ or $Y = \Delta^1$. The desired result now follows from an explicit calculation.
\end{proof}

\begin{corollary}\label{gyyyt}
Let $S \rightarrow T$ and $S' \rightarrow T'$ be categorical
equivalences of simplicial sets. Then the induced map
$$ S \star S' \rightarrow T \star T'$$ is a categorical
equivalence.
\end{corollary}

\begin{proof}
This follows immediately from Proposition \ref{rub3}, since the operation $\diamond$
has the desired property.
\end{proof}

\begin{corollary}\label{diamond3}
Let $X$ and $Y$ be simplicial sets. Then the natural map
$$ \sCoNerve[X \star Y] \rightarrow \sCoNerve[X] \star
\sCoNerve[Y]$$ is an equivalence of simplicial categories.
\end{corollary}

\begin{proof}
Using Corollary \ref{gyyyt}, we may reduce to the case where $X$ and $Y$ are $\infty$-categories.
We note that $\sCoNerve[X \star Y]$ is a correspondence from $\sCoNerve[X]$ to
$\sCoNerve[Y]$. To complete the proof, it suffices to show that 
$\bHom_{\sCoNerve[X \star Y]}(x,y)$ is weakly contractible, for any
pair of objects $x \in X$, $y \in Y$. Since $X \star Y$ is an $\infty$-category, we can apply Theorem \ref{biggie} to deduce that the mapping space $\bHom_{\sCoNerve[X \star Y]}(x,y)$ is weakly homotopy equivalent to $\Hom^{\rght}_{X \star Y}(x,y)$, which consists of a single point.
\end{proof}

For fixed $X$, the functor
$$ Y \mapsto X \diamond Y$$
$$ \sSet \rightarrow (\sSet)_{X/}$$
preserves all colimits. By the adjoint functor theorem (or by direct construction), this functor has a right adjoint $$ ( p: X \rightarrow \calC) \mapsto \calC^{p/}. $$
Since the functor $Y \mapsto X \diamond Y$ preserves cofibrations and categorical equivalences, we deduce that
the passage from $\calC$ to $\calC^{p/}$ preserves categorical fibrations and categorical equivalences between $\infty$-categories. Moreover, Proposition \ref{rub3} has the following consequence:\index{not}{calC^p/@$\calC^{p/}$}

\begin{proposition}\label{certs}
Let $\calC$ be an $\infty$-category, and let $p: X \rightarrow \calC$ be a diagram. Then the natural map $$ \calC_{p/} \rightarrow \calC^{p/}$$ is an equivalence of $\infty$-categories.
\end{proposition}

According to Definition \ref{defcolim}, a colimit for a diagram $p: X \rightarrow \calC$ is an initial object of the $\infty$-category $\calC_{p/}$. In view of the above remarks, an object of $\calC_{p/}$ is a colimit for $p$ if and only if its image in $\calC^{p/}$ is an initial object; in other words, we can replace $\calC_{p/}$ by $\calC^{p/}$ (and $\star$ by $\diamond$) in Definition \ref{defcolim}.

By Proposition \ref{sharpen}, for any $\infty$-category $\calC$ and any
map $p: X \rightarrow \calC$, the induced map $\calC_{p/} \rightarrow \calC$
is a left fibration. We now show that $\calC^{p/}$ has the same
property:

\begin{proposition}\label{sharpenn}
Suppose given a diagram of simplicial sets
$$ K_0 \subseteq K \stackrel{p}{\rightarrow} X \stackrel{q}{\rightarrow} S$$
where $q$ is a categorical fibration. Let $r = q \circ p: K \rightarrow S$,
$p_0 = p|K_0$, and $r_0 = r|K_0$. Then the induced map
$$ \phi: X^{p/} \rightarrow X^{p_0/} \times_{ S^{r_0/} } S^{r/}$$
is a left fibration.
\end{proposition}

\begin{proof}
We must show that $q$ has the right lifting property with respect to every left-anodyne inclusion
$A_0 \subseteq A$. Unwinding the definition, this amounts to proving that $q$ has the right
lifting property with respect to the inclusion
$$ i: (A_0 \diamond K) \coprod_{ A_0 \diamond K_0} (A \diamond K_0) \subseteq
A \diamond K.$$
Since $q$ is a categorical fibration, it suffices to show that $i$ is a categorical equivalence.
The above pushout is a homotopy pushout, so it will suffice to prove the analogous statement for the weakly equivalent inclusion
$$ (A_0 \star K) \coprod_{ A_0 \star K_0 } (A \star K_0) \subseteq A \star K.$$
But this map is inner anodyne (Lemma \ref{precough}).
\end{proof}

\begin{corollary}\label{measles}
Let $\calC$ be an $\infty$-category, and let $p: K \rightarrow \calC$ be any
diagram. For every vertex $v$ of $\calC$, the map $\calC_{p/} \times_{\calC}
\{v\} \rightarrow \calC^{p/} \times_{\calC} \{v\}$ is a homotopy
equivalence of Kan complexes.
\end{corollary}

\begin{proof}
The map $\calC_{p/} \rightarrow \calC^{p/}$ is a categorical equivalence
of left fibrations over $\calC$; now apply Proposition
\ref{apple1}.
\end{proof}

\begin{corollary}\label{homsetsagree}
Let $\calC$ be an $\infty$-category containing vertices $x$ and $y$. The maps
$$ \Hom^{\rght}_{\calC}(x,y) \rightarrow \Hom_{\calC}(x,y) \leftarrow \Hom^{\lft}_{\calC}(x,y)$$
are homotopy equivalences of Kan complexes (see \S \ref{prereq1} for an explanation of this notation).
\end{corollary}

\begin{proof}
Apply Corollary \ref{measles} (the dual of Corollary \ref{measles}) to the case where $p$ is the
inclusion $\{x \} \subseteq \calC$ (the inclusion $\{y\} \subseteq \calC$).
\end{proof}

\begin{remark}
The above ideas dualize in an evident way; given a map of simplicial sets
$p: K \rightarrow X$, we can define a simplicial set $X^{/p}$ with the universal
mapping property
$$ \Hom_{\sSet}(K', X^{/p} ) = \Hom_{(\sSet)_{K/}}( K' \diamond K, X).$$\index{not}{calC^/p@$\calC^{/p}$}
\end{remark}

\subsection{Parametrized Colimits}\label{consweet}

Let $p: K \rightarrow \calC$ be a diagram in an $\infty$-category $\calC$. The goal of this section is to make precise the idea that the colimit $\injlim(p)$ depends functorially on $p$ (provided that $\injlim(p)$ exists). We will prove this in a very general context, in which not only the diagram $p$ but also the simplicial set $K$ is allowed to vary. We begin by introducing a {\em relative} version of the $\diamond$-operation.

\begin{definition}
Let $S$ be a simplicial set, and let $X,Y \in (\sSet)_{/S}$. We define
$$X \diamond_S Y = X \coprod_{X \times_{S} Y \times \{0\} } (X \times_{S} Y \times \Delta^1) \coprod_{ X \times_{S} Y \times \{1\}} Y \in (\sSet)_{/S}.$$\index{not}{diamondS@$\diamond_{S}$}
\end{definition}

We observe that the operation $\diamond_S$ is compatible with base change in the following sense: for any map $T \rightarrow S$ of simplicial sets and any objects $X,Y \in (\sSet)_{/S}$, there is a natural isomorphism
$$ (X_T \diamond_T Y_T) \simeq (X \diamond_S Y)_T,$$
where we let $Z_T$ denote the fiber product $Z \times_{S} T$. 
We also note that in the case where $S$ is a point, the operation $\diamond_S$ coincides
with the operation $\diamond$ introduced in \S \ref{quasilimit2}. 

Fix $K \in (\sSet)_{/S}$. We note that functor $(\sSet)_{/S} \rightarrow ( (\sSet)_{/S})_{K/}$ defined by
$$ X \mapsto K \diamond_S S$$
has a right adjoint; this right adjoint associates to a diagram
$$ \xymatrix{ K \ar[dr] \ar[rr]^{p_S} & & Y \ar[dl] \\
& S }$$
the simplicial set $Y^{p_S/}$, defined by the property that
$ \Hom_{S}(X, Y^{p_S/})$ classifies commutative diagrams
$$ \xymatrix{ K  \ar[r]^{p_S} \ar@{^{(}->}[d] & Y \ar[d] \\  
K \diamond_S X \ar[ur] \ar[r] & S.}$$

The base-change properties of the operation $\diamond_S$ imply similar base-change properties for the relative slice construction: given a map $p_S: K \rightarrow Y$ in $(\sSet)_{/S}$ and any map $T \rightarrow S$, we have a natural isomorphism
$$ Y^{p_S/} \times_{S} T \simeq (Y \times_{S} T)^{p_T/}$$
where $p_T$ denotes the induced map $K_{T} \rightarrow Y_{T}$. 
In particular, the fiber of $Y^{p_S/}$ over a vertex $s$ of $S$ can be identified with the absolute slice construction $Y_s^{p_s/}$ studied in \S \ref{quasilimit2}.\index{not}{X^p_S/@$X^{p_{S}/}$}

\begin{remark}
Our notation is somewhat abusive: the simplicial set $Y^{p_S/}$ depends not only on the map
$p_S: K \rightarrow Y$, but also on the simplicial set $S$. We will attempt to avoid confusion by always indicating the simplicial set $S$ by a subscript in the notation for the map in question; we will only omit this subscript in the case $S = \Delta^0$, in which case the functor described above
coincides with the definition given in \S \ref{quasilimit2}.
\end{remark}

\begin{lemma}\label{stumble}
Let $n > 0$, and let 
$$B = ( \Lambda^n_n \times \Delta^1 ) \coprod_{ \Lambda^n_n \times \bd \Delta^1 }
( \Delta^n \times \bd \Delta^1 ) \subseteq \Delta^n \times \Delta^1.$$
Suppose given a diagram of simplicial sets
$$ \xymatrix{ A \times B \ar[r]^{f_0} \ar@{^{(}->}[d] & Y \ar[d]^{q} \\
A \times \Delta^n \times \Delta^1 \ar[r] \ar@{-->}[ur]^{f} & S}$$
in which $q$ is a Cartesian fibration, and that $f_0$ carries $\{ a\} \times \Delta^{ \{n-1,n\} } \times \{1\}$ to a
$q$-Cartesian edge of $Y$, for each vertex $a$ of $A$. Then there exists a morphism $f$
rendering the diagram commutative. 
\end{lemma}

\begin{proof}
Invoking Proposition \ref{doog}, we may replace $q: Y \rightarrow S$ by the induced map
$Y^A \rightarrow S^A$, and thereby reduce to the case where $A = \Delta^0$. 
We now recall the notation introduced in the proof of Proposition \ref{usejoyal}: more specifically, the family $\{ \sigma_i\}_{ 0 \leq i \leq n}$ of nondegenerate simplices of $\Delta^n \times \Delta^1$.
Let $B(0) = B$, and more generally set $B(n) = B \cup \sigma_n \cup \ldots \cup \sigma_{n+1-i}$
so that we we have a filtration
$$ B(0) \subseteq \ldots \subseteq B(n+1) = \Delta^n \times \Delta^1.$$
A map $f_0: B(0) \rightarrow Y$ has been prescribed for us already; we construct extensions $f_i: B(i) \rightarrow Y$ by induction on $i$. For $i < n$, there is a pushout diagram
$$ \xymatrix{ \Lambda^{n+1}_{n-i} \ar[r] \ar@{^{(}->}[d] & B(i) \ar@{^{(}->}[d] \\
\Delta^{n+1} \ar[r] & B(i+1) }.$$
Thus, the extension $f_{i+1}$ can be found in virtue of the assumption that $q$ is an inner fibration. For $i = n$, we obtain instead a pushout diagram
$$ \xymatrix{ \Lambda^{n+1}_{n+1} \ar[r] \ar@{^{(}->}[d] & B(n) \ar@{^{(}->}[d] \\
\Delta^{n+1} \ar[r] & B(n+1) },$$
and the desired extension can be found in virtue of the assumption that $f_0$ carries the edge $\Delta^{ \{n-1,n\}} \times \{1\}$ to a $q$-Cartesian edge of $Y$.
\end{proof}

\begin{proposition}\label{colimfam}
Suppose given a diagram of simplicial sets
$$ \xymatrix{ K \ar[drr]^{t} \ar[r]^{p_S} & X \ar[r]^{q} \ar[dr] & Y \ar[d] \\
& & S. }$$
Let $p'_{S} = q \circ p_{S}$. Suppose further that:
\begin{itemize}
\item[$(1)$] The map $q$ is a Cartesian fibration.
\item[$(2)$] The map $t$ is a coCartesian fibration.
\end{itemize}
Then the induced map $r: X^{p_S/} \rightarrow Y^{p'_S/}$ is a Cartesian fibration; moreover an edge of $X^{p_S/}$ is $r$-Cartesian if and only if its image in $X$ is $q$-Cartesian.
\end{proposition}

\begin{proof}
We first show that $r$ is an inner fibration. Suppose given $0 < i < n$ and a diagram
$$ \xymatrix{ \Lambda^n_i \ar[r] \ar@{^{(}->}[d] & X^{p_S/} \ar[d] \\
\Delta^n \ar[r] \ar@{-->}[ur] & Y^{p'_{S}/},}$$
we must show that it is possible to provide the dotted arrow. Unwinding the definitions, we see that it suffices to produce the indicated arrow in the diagram
$$ \xymatrix{ K \diamond_S \Lambda^n_i \ar[r] \ar@{^{(}->}[d] & X \ar[d]^{q} \\
K \diamond_S \Delta^n \ar[r] \ar@{-->}[ur] & Y.}$$
Since $q$ is a Cartesian fibration, it is a categorical fibration by Proposition \ref{funkyfibcatfib}.
Consequently, it suffices to show that the inclusion
$$ K \diamond_S \Lambda^n_i \subseteq K \diamond_S \Delta^n$$ is a categorical equivalence. In view of the definition of $K \diamond_S M$ as a (homotopy) pushout
$$ K \coprod_{ K \times_S M \times \{0\}} (K \times_S M \times \Delta^1) \coprod_{ K \times_S M \times \{1\} } M,$$ it suffices to verify that the inclusions
$$ \Lambda^n_i \subseteq \Delta^n$$
$$ K \times_S \Lambda^n_i \subseteq K \times_S \Delta^n$$ 
are categorical equivalences. The first statement is obvious; the second follows from (the dual of) Proposition \ref{basechangefunky}.

Let us say that an edge of $X^{p_S/}$ is {\it special} if its image in $X$ is $q$-Cartesian. To complete the proof, it will suffice to show that every special edge of $X^{p_S/}$ is $r$-Cartesian, and that there are sufficiently many special edges of $X^{p_S/}$. More precisely, consider any $n \geq 1$ and any diagram
$$ \xymatrix{ \Lambda^n_n \ar[r]^{h} \ar@{^{(}->}[d] & X^{p_S/} \ar[d] \\
\Delta^n \ar[r] \ar@{-->}[ur] & Y^{p'_{S}/}.}$$
We must show that:
\begin{itemize}
\item If $n = 1$, then there exists a dotted arrow rendering the diagram commutative, classifying a special edge of $X^{p_S/}$. 
\item If $n > 1$ and $h| \Delta^{ \{n-1,n\} }$ classifies a special edge of $X^{p_S/}$, then there exists a dotted arrow rendering the diagram commutative.
\end{itemize}
Unwinding the definitions,  we have a diagram
$$ \xymatrix{ K \diamond_S \Lambda^n_n \ar[r]^{f_0} \ar@{^{(}->}[d] & X \ar[d]^{q} \\
K \diamond_S \Delta^n \ar[r] \ar@{-->}[ur]^{f} & Y}$$
and we wish to prove the existence of the indicated arrow $f$. As a first step, we consider the restricted diagram
$$ \xymatrix{ \Lambda^n_n \ar[r]^{f_0| \Lambda^n_n} \ar@{^{(}->}[d] & X \ar[d]^{q} \\
\Delta^n \ar[r] \ar@{-->}[ur]^{f_1} & Y}.$$
By assumption, $f_0 | \Lambda^n_n$ carries $\Delta^{ \{n-1,n\} }$ to a $q$-Cartesian edge of $X$ (if $n > 1$), so there exists a map $f_1$ rendering the diagram commutative (and classifying a $q$-Cartesian edge of $X$ if $n = 1$).
It now suffices to produce the dotted arrow in the diagram
$$ \xymatrix{ (K \diamond_S \Lambda^n_n) \coprod_{ \Lambda^n_n } \Delta^n  \ar[r] \ar@{^{(}->}[d]^{i} & X \ar[d]^{q} \\
K \diamond_S \Delta^n \ar[r] \ar@{-->}[ur]^{f} & Y,}$$
where the top horizontal arrow is the result of amalgamating $f_0$ and $f_1$.

Without loss of generality, we may replace $S$ by $\Delta^n$. By (the dual of) Proposition \ref{simplexplay}, there exists a composable sequence of maps
$$ \phi: A^0 \rightarrow \ldots \rightarrow A^n$$ and a quasi-equivalence
$M^{op}(\phi) \rightarrow K$. We have a commutative diagram
$$ \xymatrix{ 
(M^{op}(\phi) \diamond_S \Lambda^n_n) \coprod_{ \Lambda^n_n } \Delta^n \ar@{^{(}->}[d]^{i'} \ar[r] & (K \diamond_S \Lambda^n_n) \coprod_{ \Lambda^n_n } \Delta^n \ar[d]^{i} \\
 M^{op}(\phi) \diamond_S \Delta^n \ar[r] & K \diamond_S \Delta^n}.$$
Since $q$ is a categorical fibration, Proposition \ref{princex} shows that it suffices to produce a dotted arrow $f'$ in the induced diagram
$$ \xymatrix{ (M^{op}(\phi) \diamond_S \Lambda^n_n) \coprod_{ \Lambda^n_n } \Delta^n  \ar[r] \ar@{^{(}->}[d]^{i} & X \ar[d]^{q} \\
M^{op}(\phi) \diamond_S \Delta^n \ar[r] \ar@{-->}[ur]^{f'} & Y}.$$
Let $B$ be as the statement of Lemma \ref{stumble}; then we have a pushout diagram
$$ \xymatrix{ 
A^0 \times B \ar[r] \ar@{^{(}->}[d]^{i''} & (M^{op}(\phi) \diamond_S \Lambda^n_n) \coprod_{ \Lambda^n_n} \Delta^n \ar[d] \\
A^0 \times \Delta^n \times \Delta^1 \ar[r] & M^{op}(\phi) \diamond_S \Delta^n.}$$
Consequently, it suffices to prove the existence of the map $f''$ in the diagram
$$ \xymatrix{ A^0 \times B \ar[r]^{g} \ar@{^{(}->}[d]^{i''} & X \ar[d]^{q} \\
A^0 \times \Delta^n \times \Delta^1 \ar[r] \ar@{-->}[ur]^{f''} & Y}.$$
Here the map $g$ carries $\{a\} \times \Delta^{ \{n-1,n\} } \times \{1\}$ to a $q$-Cartesian edge of $Y$, for each vertex $a$ of $A^0$. The existence of $f''$ now follows from Lemma \ref{stumble}.
\end{proof}

\begin{remark}\label{superfam}
In most applications of Proposition \ref{colimfam}, we will have $Y = S$. In that case,
$Y^{p'_{S}/}$ can be identified with $S$, and the conclusion is that the projection
$X^{p_S/} \rightarrow S$ is a Cartesian fibration.
\end{remark}

\begin{remark}\label{notnec}
The hypothesis on $s$ in Proposition \ref{colimfam} can be weakened: all we need in the proof is existence of maps $M^{op}(\phi) \rightarrow K \times_{S} \Delta^n$ which are universal categorical equivalences (that is, induce categorical equivalences $M^{op}(\phi) \times_{ \Delta^n } T \rightarrow K \times_{S} T$ for any $T \rightarrow \Delta^n$). Consequently, Proposition \ref{colimfam} remains valid when $K \simeq S \times K^0$, for {\em any} simplicial set $K^0$ (not necessarily an $\infty$-category). It seems likely that Proposition \ref{colimfam} remains valid whenever $s$ is a smooth map of simplicial sets, but we have not been able to prove this.
\end{remark}

We can now express the idea that the colimit a diagram should depend functorially on the diagram (at least for ``smoothly parametrized'' families of diagrams):

\begin{proposition}\label{familycolimit}\index{gen}{colimit!in families}
Let $q: Y \rightarrow S$ be a Cartesian fibration, let
$p_S: K \rightarrow Y$ be a diagram. Suppose that:
\begin{itemize}
\item[$(1)$] For each vertex $s$ of $S$, the restricted diagram $p_s: K_s \rightarrow Y_s$
has a colimit in the $\infty$-category $Y_s$.
\item[$(2)$] The composition $q \circ p_S$ is a coCartesian fibration.
\end{itemize}

There exists a map $p'_S$ rendering the diagram
$$ \xymatrix{ K \ar@{^{(}->}[d] \ar[r]^{p_S} & Y \ar[d]^{q} \\
K \diamond_S S \ar[ur]^{p'_S} \ar[r] & S }$$
commutative, and having the property that for each vertex $s$ of $S$, the restriction $p'_s: K_s \diamond \{s\} \rightarrow Y_s$ is a colimit of $p_s$.
Moreover, the collection of all such maps is parametrized by a contractible Kan complex.
\end{proposition}

\begin{proof}
Apply Proposition \ref{topaz} to the Cartesian fibration $Y^{p_S/}$, and observe that the collection of sections of a trivial fibration constitutes a contractible Kan complex.
\end{proof}

\subsection{Decomposition of Diagrams}\label{quasilimit1}
 
Let $\calC$ be an $\infty$-category, and $p: K \rightarrow \calC$ a diagram indexed by a simplicial set $K$. In this section, we will try to analyze the colimit $\colim(p)$ (if it exists) in terms of the colimits $\{ \colim(p | K_{I}) \}$, where $\{ K_{I} \}$ is some family of simplicial subsets of $K$. In fact, it will be useful to work in slightly more generality: we will allow each $K_{I}$ to be an arbitrary simplicial set mapping to $K$ (not necessarily via a monomorphism).

Throughout this section, we will fix a simplicial set $K$, an ordinary category $\calI$, and a functor
$F: \calI \rightarrow (\sSet)_{/K}$. It may be helpful to imagine that $\calI$
is a partially ordered set and that $F$ is an order-preserving map
from $\calI$ to the collection of simplicial subsets of $K$; this
will suffice for many but not all of our applications. We will
denote $F(I)$ by $K_I$, and the tautological map $K_I \rightarrow
K$ by $\pi_{I}$.

Our goal is to show that, under appropriate hypotheses, we can recover the colimit
of a diagram $p: K \rightarrow \calC$ in terms of the colimits of diagrams
$p \circ \pi_I: K_I \rightarrow \calC$. Our first goal is to show that the construction of these colimits is suitably functorial in $I$. For this, we need an auxiliary construction.

\begin{notation}\label{nixx}\index{not}{KF@$K_{F}$}
We define a simplicial set $K_{F}$ as follows. A map $\Delta^n \rightarrow K_{F}$
is determined by the following data:

\begin{itemize}
\item[$(i)$] A map $\Delta^{n} \rightarrow \Delta^{1}$, corresponding to a decomposition
$[n] = \{ 0, \ldots, i \} \cup \{ i+1, \ldots, n\}$.

\item[$(ii)$] A map $e_{-}: \Delta^{ \{ 0, \ldots, i\} } \rightarrow K$.

\item[$(iii)$] A map $e_{+}: \Delta^{\{ i+1, \ldots, n\} } \rightarrow \Nerve(\calI)$, which we may view as a chain of composable morphisms
$$ I(i+1) \rightarrow \ldots \rightarrow I(n)$$
in the category $\calI$.

\item[$(iv)$] For each $j \in \{i+1, \ldots, n\}$, a map
$e_{j}$ which fits into a commutative diagram
$$ \xymatrix{ & K_{I(j)} \ar[d]^{\pi_{I(j)}} \\
\Delta^{ \{0, \ldots, i\} } \ar[r]^{e_{-}} \ar[ur]^{e_j} & K. }$$
Moreover, for $j \leq k$ we require that
$e_k$ is given by the composition
$$ \Delta^{ \{0, \ldots, i\} } \stackrel{ e_j }{\rightarrow} 
K_{I(j)} \rightarrow K_{I(k)}.$$
\end{itemize}
\end{notation}

\begin{remark}
In the case where $i < n$, the maps $e_{-}$ and $\{ e_j \}_{j > i}$ are completely
determined by $e_{i+1}$, which can be arbitrary.
\end{remark}

The simplicial set $K_{F}$ is equipped with a map $K_{F} \rightarrow \Delta^1$. Under this map, the preimage of the vertex $\{0\}$ is $K \subseteq K_{F}$, and the preimage of the vertex $\{1\}$
is $\Nerve(\calI) \subseteq K_{F}$. For $I \in \calI$, we will denote the
corresponding vertex of $\Nerve(\calI) \subseteq K_F$ by $X_I$. We
note that, for each $I \in \calI$, there is a commutative diagram
$$ \xymatrix{ K_I \ar[r]^{\pi_I} \ar@{^{(}->}[d] & K \ar@{^{(}->}[d] \\
K_I^{\triangleright} \ar[r]^{\pi'_I} & K_F }$$
where $\pi'_I$ carries the cone point of $K_I^{\triangleright}$ to
the vertex $X_I$ of $K_F$.

Let us now suppose that $p: K \rightarrow \calC$ is a diagram in an $\infty$-category $\calC$. Our next goal is to prove Proposition \ref{extet}, which will allow us to extend $p$ to a larger diagram $K_{F} \rightarrow \calC$ which carries each vertex $X_I$ to a colimit of $p \circ \pi_I: K_I \rightarrow \calC$. First, we need a lemma.

\begin{lemma}\label{chort}
Let $\calC$ be an $\infty$-category, and let
$\sigma: \Delta^n \rightarrow \calC$ be a simplex having the property that
$\sigma(0)$ is an initial object of $\calC$. Let $\bd \sigma = \sigma | \bd \Delta^n$. 
The natural map $\calC_{\sigma/} \rightarrow \calC_{\bd \sigma/}$ is a trivial fibration.
\end{lemma}

\begin{proof}
Unwinding the definition, we are reduced to solving the extension problem depicted in the diagram
$$ \xymatrix{ (\bd \Delta^n \star \Delta^m) \coprod_{
\bd \Delta^n \star \bd \Delta^m } (\Delta^n \star \bd \Delta^m) \ar[r]^-{f_0} 
\ar@{^{(}->}[d] & \calC \\
\Delta^n \star \Delta^m. \ar@{-->}[ur]^-{f} & }$$
We can identify the domain of $f_0$ with $\bd \Delta^{n+m+1}$. Our hypothesis
guarantees that $f_0(0)$ is an initial object of $\calC$, which in turn guarantees the existence of $f$.
\end{proof}

\begin{proposition}\label{extet}
Let $p: K \rightarrow \calC$ be a diagram in an $\infty$-category $\calC$, let $\calI$ be an ordinary category, and let $F: \calI \rightarrow (\sSet)_{/K}$ be a functor. Suppose that, for each $I \in \calI$, the induced diagram $p_I = p \circ \pi_I: K_I \rightarrow \calC$ has a colimit
$q_I: K_I^{\triangleright} \rightarrow \calC$.

There exists a map $q: K_F \rightarrow \calC$ such that $q \circ \pi'_I = q_{I}$ and $q|K = p$. 
Furthermore, for any such $q$, the induced map $\calC_{q/} \rightarrow \calC_{p/}$ is a trivial fibration.
\end{proposition}

\begin{proof}
For each $X \subseteq \Nerve(\calI)$, we let $K_{X}$ denote the
simplicial subset of $K_F$ consisting of all simplices $\sigma \in
K_{F}$ such that $\sigma \cap \Nerve(\calI) \subseteq X$. We note that
$K_{\emptyset} = K$ and that $K_{ \Nerve(\calI) } = K_{F}$.

Define a transfinite sequence $Y_{\alpha}$ of simplicial subsets
of $\Nerve(\calI)$ as follows. Let $Y_{0} = \emptyset$, and let
$Y_{\lambda} = \bigcup_{ \gamma < \lambda } Y_{\gamma}$ when
$\lambda$ is a limit ordinal. Finally, let $Y_{\alpha+1}$ be
obtained from $Y_{\alpha}$ by adjoining a single nondegenerate
simplex, provided that such a simplex exists. We note that for
$\alpha$ sufficiently large, such a simplex will not exist and we
set $Y_{\beta} = Y_{\alpha}$ for all $\beta > \alpha$.

We define a sequence of maps $q_{\beta}: K_{Y_{\beta}} \rightarrow
\calC$ so that the following conditions are satisfied:

\begin{itemize}
\item[$(1)$] We have $q_{0} = p: K_{\emptyset} = K \rightarrow \calC$.

\item[$(2)$] If $\alpha < \beta$, then $q_{\alpha} = q_{\beta} |
K_{Y_{\alpha}}$.

\item[$(3)$] If $\{X_I\} \subseteq Y_{\alpha}$, then $q_{\alpha} \circ \pi'_I = q_I: K_I^{\triangleright} \rightarrow \calC$.

\end{itemize}
Provided that such a sequence can be constructed, we may conclude
the proof by setting $q = q_{\alpha}$ for $\alpha$ sufficiently
large.

The construction of $q_{\alpha}$ goes by induction on $\alpha$. If
$\alpha = 0$, then $q_{\alpha}$ is determined by condition $(1)$;
if $\alpha$ is a (nonzero) limit ordinal, then $q_{\alpha}$ is
determined by condition $(2)$. Suppose that $q_{\alpha}$ has been
constructed; we give a construction of $q_{\alpha+1}$.

There are two cases to consider. Suppose first that $Y_{\alpha+1}$
is obtained from $Y_{\alpha}$ by adjoining a vertex $X_I$. In this
case, $q_{\alpha+1}$ is uniquely determined by conditions $(2)$
and $(3)$.

Now suppose that $X_{\alpha+1}$ is obtained from $X_{\alpha}$ by
adjoining a nondegenerate simplex $\sigma$ of positive dimension, corresponding
to a sequence of composable maps
$$ I_0 \rightarrow \ldots \rightarrow I_n$$
in the category $\calI$. We
note that the inclusion $K_{Y_{\alpha}} \subseteq
K_{Y_{\alpha+1}}$ is a pushout of the inclusion
$$ K_{I_0} \star \bd \sigma \subseteq K_{I_0} \star \sigma.$$
Consequently, constructing the map $q_{\alpha+1}$ is tantamount to
finding an extension of a certain map $s_0: \bd \sigma \rightarrow \calC_{p_I/}$ to the whole of the simplex $\sigma$. By assumption, $s_0$ carries the initial vertex of $\sigma$ to an initial
object of $\calC_{p_I/}$, so that the desired extension $s$ can be found. For use below, we record
a further property of our construction: the projection $\calC_{q_{\alpha+1}/} \rightarrow \calC_{q_{\alpha}/}$ is a pullback of the map $( \calC_{ p_I/} )_{s/} \rightarrow (\calC_{p_I/})_{s_0/}$, which is a trivial fibration.

We now wish to prove that for any extension $q$ with the above properties, the induced map
$\calC_{q/} \rightarrow \calC_{p/}$ is a trivial fibration. We first observe that the map $q$ can be obtained by the inductive construction given above: namely, we take $q_{\alpha}$ to be
the restriction of $q$ to $K_{Y_{\alpha}}$. It will therefore suffice to show that, for every pair of ordinals $\alpha \leq \beta$, the induced map $\calC_{q_{\beta}/} \rightarrow \calC_{q_{\alpha}/}$ is a trivial fibration. The proof of this goes by induction on $\beta$: the case $\beta = 0$ is clear, and if $\beta$ is a limit ordinal we observe that the inverse limit of transfinite tower of trivial fibrations is itself a trivial fibration. We may therefore suppose that $\beta = \gamma + 1$ is a successor ordinal.
Using the factorization
$$ \calC_{q_{\beta}/} \rightarrow \calC_{q_{\gamma}/} \rightarrow \calC_{q_{\alpha}/}$$
and the inductive hypothesis, we are reduced to proving this in the case where $\beta$ is the successor of $\alpha$, which was treated above.
\end{proof}

Let us now suppose that we are given diagrams $p: K \rightarrow \calC$, $F: \calI \rightarrow (\sSet)_{/K}$ as in the statement of Proposition \ref{extet}, and let $q: K_{F} \rightarrow \calC$ be a map which satisfies the conclusions of the Proposition. Since $\calC_{q/} \rightarrow \calC_{p/}$ is a trivial fibration, we may identify colimits of the diagram $q$ with colimits of the diagram $p$ (up to equivalence). Of course, this is not useful in itself,
since the diagram $q$ is {\em more} complicated than $p$. Our objective
now is to show that, under the appropriate hypotheses, we may
identify the colimits of $q$ with the colimits of $q|\Nerve(\calI)$.
First, we need a few lemmas.

\begin{lemma}[Joyal \cite{joyalnotpub}]\label{chotle}
Let $f: A_0 \subseteq A$ and $g: B_0 \subseteq B$ be inclusions of simplicial sets, and suppose
that $g$ is a weak homotopy equivalence. Then the induced map
$$h: (A_0 \star B) \coprod_{ A_0 \star B_0} (A \star B_0) \subseteq A \star B$$ is right anodyne.
\end{lemma}

\begin{proof}
Our proof follows the pattern of Lemma \ref{precough}. The collection of all maps $f$ which
satisfy the conclusion (for {\em any} choice of $g$) forms a weakly saturated class of morphisms. It will therefore suffice to prove that the $h$ is right anodyne when $f$ is the inclusion $\bd \Delta^n \subseteq \bd \Delta^n$. Similarly, the collection of all maps $g$ which satisfy the conclusion (for fixed $f$) forms a weakly saturated class. We may therefore reduce to the case where $g$ is a horn inclusion $\Lambda^m_i \subseteq \Delta^m$. In this case, we may identify $h$ with the horn inclusion $\Lambda^{m+n+1}_{i+n+1} \subseteq \Delta^{m+n+1}$, which is clearly right-anodyne.
\end{proof}

\begin{lemma}\label{chotle2}
Let $A_0 \subseteq A$ be an inclusion of simplicial sets, and let
$B$ be weakly contractible. Then the inclusion $A_0 \star B
\subseteq A \star B$ is right anodyne.
\end{lemma}

\begin{proof}
As above, we may suppose that the inclusion $A_0 \subseteq A$ is identified with
$ \bd \Delta^n \subseteq \Delta^n$. If $K$ is a point, then the inclusion $A_0 \times B \subseteq A \times B$ is 
isomorphic to $\Lambda^{n+1}_{n+1}
\subseteq \Delta^{n+1}$, which is clearly right-anodyne.

In the general case, $B$ is nonempty, so we may choose a vertex
$b$ of $B$. Since $B$ is weakly contractible, the inclusion $\{b\}
\subseteq B$ is a weak homotopy equivalence. We have already shown
that $A_0 \star \{b\} \subseteq A \star \{b\}$ is right anodyne. It
follows that the pushout inclusion
$$A_0 \star B \subseteq (A \star \{b\}) \coprod_{ A_0 \star \{b\} } (A_0
\star B)$$ is right anodyne. To complete the proof, we apply Lemma
\ref{chotle} to deduce that the inclusion
$$  (A \star \{b\}) \coprod_{ A_0 \star \{b\} } (A_0
\star B) \subseteq A \star B$$ is right anodyne.
\end{proof}

\begin{notation}
Let $\sigma \in K_n$ be a simplex of $K$. We define a category
$\calI_{\sigma}$ as follows. The objects of $\calI_{z}$ are pairs $(I,
\sigma')$, where $I \in \calI$, $\sigma' \in (K_I)_n$, and $\pi_I(\sigma') = \sigma$. A
morphism from $(I',\sigma')$ to $(I'',\sigma'')$ in $\calI_{\sigma}$ consists of
a morphism $\alpha: I' \rightarrow I''$ in $\calI$ with the
property that $F(\alpha)(\sigma') = \sigma''$. We let $\calI'_{\sigma} \subseteq
\calI_{\sigma}$ denote the full subcategory consisting of pairs
$(I,\sigma')$ where $\sigma'$ is a degenerate simplex in $K_{I}$. Note that
if $\sigma$ is nondegenerate, $\calI'_{\sigma}$ is empty.
\end{notation}

\begin{proposition}\label{utl}
Let $K$ be a simplicial set, $\calI$ an ordinary category, and $F: \calI \rightarrow (\sSet)_{/K}$ a functor.
Suppose further that:
\begin{itemize}
\item[$(1)$] For each nondegenerate simplex $\sigma$ of $K$, the category
$\calI_{\sigma}$ is acyclic (that is, the simplicial set $\Nerve(\calI_{\sigma})$
is weakly contractible).

\item[$(2)$] For each degenerate simplex $\sigma$ of $K$, the inclusion $\Nerve
(\calI'_{\sigma}) \subseteq \Nerve(\calI_{\sigma})$ is a weak homotopy equivalence.
\end{itemize}

Then the inclusion $\Nerve(\calI) \subseteq K_{F}$ is right anodyne.
\end{proposition}

\begin{proof}
Consider any family of subsets $\{ L_{n} \subseteq K_{n}\}$ which is
stable under the ``face maps'' $d_i$ on $K$ (but not necessarily
the degeneracy maps $s_i$, so that the family $\{ L_{n} \}$ does
not necessarily have the structure of a simplicial set). We define
a simplicial subset $L_{F} \subseteq K_F$ as follows: a {\em nondegenerate} simplex
$\Delta^n \rightarrow K_{F}$ belongs to $L_{F}$ if and only if the
corresponding (possibly degenerate) simplex $\Delta^{ \{0, \ldots, i\} } \rightarrow K$ belongs to
$L_i \subseteq K_i$ (see Notation \ref{nixx}). 

We note that if $L = \emptyset$, then $L_{F} = \Nerve(\calI)$. If $L = K$, then $L_{F}= K_{F}$ (so that our notation is unambiguous).
Consequently, it will suffice to prove that for any $L \subseteq
L'$, the inclusion $L_{F} \subseteq L'_{F}$ is right-anodyne. By
general nonsense, we may reduce to the case where $L'$ is obtained
from $L$ by adding a single simplex $\sigma \in K_n$.

We now have two cases to consider. Suppose first that the simplex
$\sigma$ is nondegenerate. In this case, it is not difficult to see
that the inclusion $L_{F} \subseteq L'_{F}$ is a pushout of $ \bd
\sigma \star \Nerve(\calI_{\sigma}) \subseteq \sigma \star \Nerve(\calI_{\sigma})$. By hypothesis,
$N \calI_{z}$ is weakly contractible, so that the inclusion $L_{F}
\subseteq L'_{F}$ is right anodyne by Lemma \ref{chotle2}.

In the case where $\sigma$ is degenerate, we observe that $L_{F}
\subseteq L'_{F}$ is a pushout of the inclusion
$$ (\bd \sigma \star \Nerve( \calI_{\sigma} )) \coprod_{ \bd \sigma \star \Nerve(\calI'_{\sigma})} (\sigma \star \Nerve(\calI'_{\sigma})) \subseteq \sigma
\star \Nerve( \calI_{\sigma} ),$$ which is right anodyne by Lemma
\ref{chotle}.
\end{proof}

\begin{remark}
Suppose that $\calI$ is a partially ordered set, and that $F$ is
an order-preserving map from $\calI$ to the collection of
simplicial subsets of $K$. In this case, we observe that
$\calI'_{\sigma} = \calI_{\sigma}$ whenever $\sigma$ is a degenerate simplex of $K$, and that
$\calI_{\sigma} = \{ I \in \calI: \sigma \in K_I \}$ for any $\sigma$. Consequently, the
conditions of Proposition \ref{utl} hold if and only if each of
the partially ordered subsets $\calI_{\sigma} \subseteq \calI$ has a
contractible nerve. This holds automatically if $\calI$ is
directed and $K = \bigcup_{I \in \calI} K_I$.
\end{remark}

\begin{corollary}\label{util}
Let $K$ be a simplicial set, $\calI$ a category, and $F: \calI \rightarrow (\sSet)_{/K}$ a functor which satisfies the hypotheses of Proposition \ref{utl}. Let $\calC$ be an $\infty$-category, $p: K \rightarrow \calC$ any diagram, and let $q: K_{F} \rightarrow \calC$ be an extension of $p$ which satisfied the conclusions of Proposition \ref{extet}. The natural maps
$$ \calC_{p/} \leftarrow \calC_{q/} \rightarrow \calC_{q | \Nerve(\calI)/}$$ are trivial fibrations.
In particular, we may identify colimits of $p$ with colimits of $q| \Nerve(\calI)$.
\end{corollary}

\begin{proof}
This follows immediately from Proposition \ref{utl}, since the right anodyne inclusion
$\Nerve \calI \subseteq K_{F}$ is cofinal and therefore induces a trivial fibration $\calC_{q/} \rightarrow \calC_{q|\Nerve(\calI)/}$ by Proposition \ref{gute}.
\end{proof}

We now illustrate the usefulness of Corollary \ref{util} by giving a sample application. First, a bit of terminology. If $\kappa$ and $\tau$ are regular cardinals, we will write $\tau \ll \kappa$ if, for any cardinals $\tau_0 < \tau$, $\kappa_0 < \kappa$, we have $\kappa_0^{\tau_0} < \kappa$ (we refer the reader to Definition \ref{ineq} and the surrounding discussion for more details concerning this condition).\index{not}{kappalltau@$\kappa \ll \tau$}

\begin{corollary}\label{uterrr}
Let $\calC$ be an $\infty$-category and $\tau \ll \kappa$ regular cardinals. Then
$\calC$ admits $\kappa$-small colimits if and only if $\calC$ admits $\tau$-small
colimits and colimits indexed by (the nerves of) $\kappa$-small, $\tau$-filtered partially ordered sets.
\end{corollary}

\begin{proof}
The ``only if'' direction is obvious. Conversely, let $p: K
\rightarrow \calC$ be any $\kappa$-small diagram. Let $\calI$ denote the partially
ordered set of $\tau$-small simplicial subsets of $K$. Then
$\calI$ is directed and $\bigcup_{I \in \calI} K_I = K$, so that the hypotheses of
Proposition \ref{utl} are satisfied. Since each $p_I = p \circ \pi_I$ has a colimit
in $\calC$, there exists a map $q: K_{F} \rightarrow \calC$ satisfying the Proposition \ref{extet}.
Since $\calC_{q/} \rightarrow \calC_{p/}$ is an equivalence of $\infty$-categories, $p$ has a colimit if and only if $q$ has a colimit. By Corollary \ref{util}, $q$ has a colimit if and only if
$q| \Nerve(\calI)$ has a colimit. It is clear that $\calI$ is a $\tau$-filtered partially ordered set.
Furthermore, it is $\kappa$-small provided that $\tau \ll \kappa$.
\end{proof}

Similarly, we have:

\begin{corollary}
Let $f: \calC \rightarrow \calC'$ be a functor between $\infty$-categories, and let
$\tau \ll \kappa$ be regular cardinals. Suppose that
$\calC$ admits $\kappa$-small colimits. Then $f$ preserves $\kappa$-small colimits if and only if it preserves $\tau$-small colimits, and all colimits indexed by (the nerves of) $\kappa$-small, $\tau$-filtered partially ordered sets.
\end{corollary}

We will conclude this section with another application of Proposition
\ref{utl}, in which $\calI$ is not a partially ordered
set, and the maps $\pi_I: K_I \rightarrow K$ are not (necessarily)
injective. Instead, we take $\calI$ to be the {\it category of\index{gen}{category!of simplices}
simplices of $K$}. In other words, an object of $I \in \calI$ consists
of a map $\sigma_I: \Delta^n \rightarrow K$, and a morphism
from $I$ to $I'$ is given by a commutative diagram
$$ \xymatrix{ \Delta^n \ar[dr]^{\sigma_I} \ar[rr] & & \Delta^{n'} \ar[dl]_{\sigma'_{I'}} \\
& K. & }$$
For each $I \in \calI$, we let $K_I$ denote the domain $\Delta^n$ of $\sigma_I$, and we
let $\pi_{I} = \sigma_I: K_I \rightarrow K$.

\begin{lemma}\label{snick}
Let $K$ be a simplicial set, and let $\calI$ denote the
category of simplices of $K$ (as defined above). Then there is a
retraction $r: K_{F} \rightarrow K$ which fixes $K \subseteq K_F$.
\end{lemma}

\begin{proof}
Given a map $e: \Delta^n \rightarrow K_{F}$, we will describe the
composite map $r \circ e: \Delta^n \rightarrow K$. The map $e$ classifies the following data:
\begin{itemize}
\item[$(i)$] A decomposition $[n] = \{0, \ldots, i\} \cup \{i+1, \ldots, n \}$. 
\item[$(ii)$] A map $e_{-}: \Delta^{i} \rightarrow K$.
\item[$(iii)$] A string of morphisms
$$ \Delta^{ m_{i+1} } \rightarrow \ldots \rightarrow \Delta^{ m_n } \rightarrow K.$$
\item[$(iv)$] A compatible family of maps
$\{ e_j: \Delta^{i} \rightarrow \Delta^{m_j} \}_{ j > i }$, having the property that each composition $\Delta^i \stackrel{e_j}{\rightarrow} \Delta^{m_j} \rightarrow K$
coincide with $e_{-}$.
\end{itemize}

If $i = n$, we set $r \circ e = e_{-}$. Otherwise, we let $r \circ e$ denote the composition
$$ \Delta^n \stackrel{f}{\rightarrow} \Delta^{m_n} \rightarrow K$$
where $f: \Delta^n \rightarrow \Delta^{m_n}$ is defined as follows:

\begin{itemize}
\item The restriction $f | \Delta^{i}$ coincides with $e_n$.

\item For $i < j \leq n$, we let $f(j)$ denote the image in
$\Delta^{m_n}$ of the final vertex of $\Delta^{m_j}$.
\end{itemize}
\end{proof}

\begin{proposition}\label{cofinalcategories}
For every simplicial set $K$, there exists a category
$\calI$ and a cofinal map $f: \Nerve(\calI) \rightarrow K$.
\end{proposition}

\begin{proof}
We take $\calI$ to be the category of simplices of $K$, as defined
above, and $f$ to the composition of the inclusion $\Nerve(\calI)
\subseteq K_{F}$ with the retraction $r$ of Lemma \ref{snick}. 
To prove that $f$ is cofinal, it suffices to show that the inclusion
$\Nerve(\calI) \subseteq K_{F}$ is right anodyne, and that the retraction $r$ is cofinal.

To show that $\Nerve( \calI) \subseteq K_{F}$ is right anodyne, it suffices to show that the hypotheses of Proposition \ref{utl} are satisfied.  Let $\sigma: \Delta^J \rightarrow K$ be a simplex of $K$. We observe
that the category $\calI_{\sigma}$ may be described as follows: its
objects consist of pairs of maps $(s: \Delta^J \rightarrow
\Delta^{M}, t: \Delta^{M} \rightarrow K)$ with $t \circ s = \sigma$. A
morphism from $(s,t)$ to $(s',t')$ consists of a map
$$ \alpha: \Delta^M \rightarrow \Delta^{M'}$$
with $s' = \alpha \circ s$, $t = t' \circ \alpha$. In particular,
we note that $\calI_{\sigma}$ has an initial object $(\id_{\Delta^J},
\sigma)$. It also has a final object: namely, a pair $(s,t)$ such
that $s$ is surjective and $t: \Delta^M \rightarrow K$ is
nondegenerate. It follows that $\Nerve(\calI_{\sigma})$ is weakly
contractible for {\em any} simplex $\sigma$ of $K$. Moreover, if $z$ is
degenerate, then any final object of $\calI_{\sigma}$ belongs to
$\calI'_{\sigma}$ (and is therefore a final object of $\calI'_{\sigma}$). We conclude that $\Nerve(\calI'_{\sigma})$ is weakly contractible when $\sigma$ is
degenerate, so that the inclusion $\Nerve(\calI'_{\sigma}) \subseteq
\Nerve(\calI_{\sigma})$ is a weak homotopy equivalence. This completes the verification of the hypotheses of Proposition \ref{utl}.

We now show that $r$ is cofinal. According to Proposition \ref{gute}, it suffices to show that
for any $\infty$-category $\calC$ and any map $p: K \rightarrow \calC$, the induced map
$\calC_{q/} \rightarrow \calC_{p/}$ is a categorical equivalence, where $q = p \circ r$.
This follows from Proposition \ref{extet}.
\end{proof}

\begin{variant}\label{baryvar}\index{gen}{barycentric subdivision}
Let $K$ be a simplicial set, and let $\calI$ be the category of simplices of $K$ as above.
Let $\calI'$ be the full subcategory of $\calI$ spanned by the nondegenerate simplices of $K$. 
The inclusion $\calI' \subseteq \calI$ has a left adjoint $L$. It follows immediately from Theorem \ref{hollowtt} that the inclusion $\Nerve(\calI') \subseteq \Nerve(\calI)$ is cofinal. Consequently,
we obtain also a cofinal map $f: \Nerve(\calI') \rightarrow K$. The simplicial set
$\Nerve(\calI')$ can be identified with the {\it barycentric subdivision} of $K$. The assertion that $f$ is cofinal can be regarded as a generalization of the classical fact that barycentric subdivision does not change the weak homotopy type of a simplicial set. 

Note the category of nondegenerate simplices of $\Nerve(\calI')$ can be identified with a partially ordered set. The nerve of this partially ordered set can be identified with the {\it second barycentric subdivision $K^{(2)}$ of $K$}. Applying the above argument twice, we conclude that
there is a cofinal map $K^{(2)} \rightarrow K$. Consequently, we obtain the following refinement of 
Proposition \ref{cofinalcategories}: for every simplicial set $K$, there exists a partially ordered set $A$ and a cofinal map $\Nerve(A) \rightarrow K$.
\end{variant}

\subsection{Homotopy Colimits}\label{quasilimit4}

Our goal in this section is to compare the $\infty$-categorical theory of colimits with
the more classical theory of homotopy colimits in simplicial categories (see
Remark \ref{curble}). Our main result is the following:

\begin{theorem}\label{colimcomparee}\index{gen}{colimit!homotopy}\index{gen}{homotopy colimit}
Let $\calC$ and $\calI$ be fibrant simplicial categories and
$F: \calI \rightarrow \calC$ a simplicial functor. Suppose given an object $C \in \calC$ and a compatible
family of maps $\{ \eta_{I}: F(I) \rightarrow C \}_{I \in \calI}$. The following conditions are
equivalent:
\begin{itemize}
\item[$(1)$] The maps $\eta_I$ exhibit $C$ as a homotopy colimit of the diagram $F$.
\item[$(2)$] Let $f: \Nerve(\calI) \rightarrow \Nerve(\calC)$ be the simplicial nerve of $F$,
and $\overline{f}: \Nerve(\calI)^{\triangleright} \rightarrow \Nerve(\calC)$ the
extension of $f$ determined by the maps $\{ \eta_I \}$. Then
$\overline{f}$ is a colimit diagram in $\Nerve(\calC)$.
\end{itemize}
\end{theorem}

\begin{remark}
For an analogous result (in a slightly different setting), we refer
the reader to \cite{hirschowitz}.
\end{remark}

The proof of Theorem \ref{colimcomparee} will occupy the remainder of this section. 
We begin with a convenient criterion for detecting colimits in $\infty$-categories:

\begin{lemma}\label{kamma}
Let $\calC$ be an $\infty$-category, $K$ a simplicial set, and
$\overline{p}: K^{\triangleright} \rightarrow \calC$ a diagram. The following conditions are equivalent:
\begin{itemize}
\item[$(i)$] The diagram $\overline{p}$ is a colimit of $p = \overline{p} | K$.
\item[$(ii)$] Let $X \in \calC$ denote the image under $\overline{p}$ of the cone point
of $K^{\triangleright}$, let $\delta: \calC \rightarrow \Fun(K,\calC)$ denote the diagonal
embedding, and let $\alpha: p \rightarrow \delta(X)$ denote the natural transformation determined by $\overline{p}$. Then, for every object $Y \in \calC$, composition with $\alpha$ induces
a homotopy equivalence
$$ \phi_{Y}: \bHom_{\calC}(X, Y) \rightarrow \bHom_{ \Fun(K,\calC) }( p, \delta(Y) ).$$
\end{itemize}
\end{lemma}

\begin{proof}
Using Corollary \ref{homsetsagree}, we can identify the mapping space
$\bHom_{\Fun(K,\calC) }( p, \delta(Y) )$ with the fiber
$\calC^{p/} \times_{\calC} \{Y\}$, for each object $Y \in \calC$.
Under this identification, the map $\phi_{Y}$ can be identified with
the fiber over $Y$ of the composition
$$ \calC^{X/} \stackrel{\phi'}{\rightarrow} \calC^{ \overline{p}/ } \stackrel{\phi''}{\rightarrow} \calC^{p/},$$
where $\phi'$ is a section to the trivial fibration $\calC^{\overline{p}/} \rightarrow \calC^{X/}$. 
The map $\phi''$ is a left fibration (Proposition \ref{sharpenn}). Condition $(i)$ is equivalent to the requirement that $\phi''$ be a trivial Kan fibration, and condition $(ii)$ is equivalent to the
requirement that each of the maps
$$ \phi''_{Y}: \calC^{ \overline{p}/} \times_{\calC} \{Y \} \rightarrow \calC^{p/} \times_{\calC} \{Y\}.$$
is a homotopy equivalence of Kan compexes (which, in view of Lemma \ref{toothie2}, is equivalent to the requirement that $\phi''_{Y}$ be a trivial Kan fibration). The equivalence of these two conditions now follows from Lemma \ref{toothie}.
\end{proof}

The key to Theorem \ref{colimcomparee} is the following result, which compares
the construction of diagram categories in the $\infty$-categorical and simplicial settings:

\begin{proposition}\label{gumby444}\index{gen}{straightening of diagrams}
Let $S$ be a small simplicial set, $\calC$ a small simplicial category, and 
$u: \sCoNerve[S] \rightarrow \calC$ an equivalence. Suppose that
$\bfA$ is a combinatorial simplicial model category, and let
$\calU$ be a $\calC$-chunk of $\bfA$ (see Definition \ref{cattusi}). 
Then the induced map
$$ \sNerve ( \calU^{\calC} )^{\degree}  \rightarrow \Fun(S,
\sNerve(\calU^{\degree}) )$$ is a categorical equivalence of
simplicial sets.
\end{proposition}

\begin{remark}
In the statement of Proposition \ref{gumby444}, it makes no difference whether we regard
$\bfA^{\calC}$ as endowed with the projective or injective model structure.
\end{remark}

\begin{remark}
An analogous result was proved by Hirschowitz and Simpson; see \cite{hirschowitz}.
\end{remark}

\begin{proof}
Choose a regular cardinal $\kappa$ such that $S$ and $\calC$ are $\kappa$-small.
Using Lemma \ref{exchunk}, we can write $\calU$ as a $\kappa$-filtered colimit
of small $\calC$-chunks $\calU'$ contained in $\calU$. Since the collection of
categorical equivalences is stable under filtered colimits, it will suffice to prove the
result after replacing $\calU$ by each $\calU'$; in other words, we may suppose that
$\calU$ is small. 

According to Theorem \ref{biggier}, we may identify
the homotopy category of $\sSet$ (with respect to the Joyal model
structure) with the homotopy category of $\sCat$. We now observe that, because
$\sNerve( \calU^{\degree})$ is an $\infty$-category, the simplicial set
$\Fun(S,\sNerve( \calU^{\degree}))$ can be identified with an exponential
$[ \Nerve(\calU^{\degree}) ]^{ [S] }$ in the homotopy category $\h{\sSet}$. We now conclude by applying Corollary \ref{sniffle}.
\end{proof}

One consequence of Proposition \ref{gumby444} is that every homotopy coherent diagram
in a suitable model category $\bfA$ can be ``straightened'', as we indicated in Remark \ref{remmt}.

\begin{corollary}\label{strictify}
Let $\calI$ be a fibrant simplicial category, $S$ a simplicial
set, and $p: \sNerve(\calI) \rightarrow S$ a map. Then it is possible to find the following:
\begin{itemize}
\item[$(1)$] A fibrant simplicial category $\calC$.
\item[$(2)$] A simplicial functor $P: \calI \rightarrow \calC$.
\item[$(3)$] A categorical equivalence of simplicial sets
$j: S \rightarrow \sNerve(\calC)$.
\item[$(4)$] An equivalence between $j \circ p$ and $\sNerve(P)$, as objects of the
$\infty$-category $\Fun( \Nerve(\calI), \sNerve(\calC))$.
\end{itemize}
\end{corollary}

\begin{proof}
Choose an equivalence $i: \sCoNerve[S] \rightarrow \calC_0$, where
$\calC_0$ is fibrant; let $\bfA$ denote the model category of
simplicial presheaves on $\calC_0$ (endowed with the {\em injective} model structure). Composing $i$ with the Yoneda
embedding of $\calC_0$, we obtain a fully faithful simplicial
functor $\sCoNerve[S] \rightarrow \bfA^{\degree}$, which we may
alternatively view as a morphism $j_0: S \rightarrow \sNerve
(\bfA^{\degree})$.

We now apply Proposition \ref{gumby444} to the case where $u$ is the
counit map $\sCoNerve[\sNerve(\calI)] \rightarrow
\calI$. We deduce that the natural map
$$ \sNerve (\bfA^{\calI})^{\degree} \rightarrow \Fun( \Nerve(\calI), \sNerve(\bfA^{\degree}
) )$$ is an equivalence. From the essential
surjectivity, we deduce that $j_0 \circ p$ is equivalent to
$\sNerve(P_0)$, where $P_0: \calI \rightarrow \bfA^{\degree}$ is a
simplicial functor.

We now take $\calC$ to be the essential image of $\sCoNerve[S]$ in
$\bfA^{\degree}$, and note that $j_0$ and $P_0$ factor uniquely
through maps $j: S \rightarrow \sNerve(\calC)$, $P: \calI
\rightarrow \calC$ which possess the desired properties.
\end{proof}

We now return to our main result.

\begin{proof}[Proof of Theorem \ref{colimcomparee}:]
Let $\bfA$ denote the category $\Set_{\Delta}^{\calC}$, endowed with the projective
model structure. Let $j: \calC^{op} \rightarrow \bfA$ denote the Yoneda embedding, and let $\calU$ denote the full subcategory of $\bfA$ spanned by those objects which are weakly equivalent
to $j(C)$ for some $C \in \calC$, so that $j$ induces an equivalence of simplicial categories
$\calC^{op} \rightarrow \calU^{\degree}$. Choose a trivial injective cofibration
$j \circ F \rightarrow F'$, where $F'$ is a injectively fibrant object of $\bfA^{\calI^{op}}$. 
Let $f': \Nerve(\calI)^{op} \rightarrow \Nerve( \calU^{\degree} )$
be the nerve of $F'$, and let $C' = j(C)$, so that the maps 
$\{ \eta_{I}: F(I) \rightarrow C \}_{I \in \calI}$ induce a natural transformation $\alpha: \delta(C') \rightarrow f'$, where $\delta: \Nerve( \calU^{\degree}) \rightarrow \Fun( \Nerve(\calI)^{op}, \Nerve(\calU^{\degree}) )$ denotes the diagonal embedding. In view of Lemma \ref{kamma}, condition
$(1)$ admits the following reformulation:

\begin{itemize}
\item[$(1')$] For every object $A \in \calU^{\degree}$, composition with
$\alpha$ induces a homotopy equivalence 
$$ \bHom_{ \Nerve(\calU^{\degree}) }( A,C') \rightarrow \bHom_{ \Fun( \Nerve(\calI)^{op}, \Nerve(\calU^{\degree})) }( \delta(A),f' ).$$
\end{itemize}

Using Proposition \ref{gumby444}, we can reformulate this condition again:

\begin{itemize}
\item[$(1'')$] For every object $A \in \calU^{\degree}$, the canonical map
$$ \bHom_{ \bfA }( A,C') \rightarrow \bHom_{ \bfA^{ \calI^{op}}}( \delta'(A), F')$$
is a homotopy equivalence, where $\delta': \bfA \rightarrow \bfA^{\calI^{op}}$ denotes
the diagonal embedding. 
\end{itemize}

Let $B \in \bfA$ be a limit of the diagram $F'$, so we have a canonical map
$\beta: C' \rightarrow B$ between fibrant objects of $\bfA$. Condition $(2)$ is
equivalent to the assertion that $\beta$ is a weak equivalence in $\bfA$, while
condition $(1'')$ is equivalent to the assertion that composition with
$\beta$ induces a homotopy equivalence
$$ \bHom_{\bfA}(A, C') \rightarrow \bHom_{\bfA}(A, B)$$
for each $A \in \calU^{\degree}$. The implication $(2) \Rightarrow (1'')$ is clear.
Conversely, suppose that $(1'')$ is satisfied. For each $X \in \calC$, the
object $j(X)$ belongs to $\calU^{\degree}$, so that $\beta$ induces a
homotopy equivalence
$$ C'(X) \simeq \bHom_{\bfA}( j(X), C') \rightarrow \bHom_{\bfA}( j(X), B)
\simeq B(X).$$
It follows that $\beta$ is a weak equivalence in $\bfA$ as desired.
\end{proof}

\begin{corollary}\label{limitsinmodel}
Let $\bfA$ be a combinatorial simplicial model category. The associated $\infty$-category
$S = \sNerve( \bfA^{\degree})$ admits $($small$)$ limits and colimits.
\end{corollary}

\begin{proof}
We give the argument for colimits; the case of limits follows by a dual argument. Let $p: K \rightarrow S$ be a (small) diagram in $S$. 
By Proposition \ref{cofinalcategories}, there exists a (small) category
$\calI$ and a cofinal map $q: \Nerve(\calI) \rightarrow K$. Since $q$ is cofinal, $p$ has a colimit in $S$ if and only if $p \circ q$ has a colimit in $S$; thus we may reduce to the case where $K = \Nerve(\calI)$. 

Using Proposition \ref{gumby444}, we may suppose that $p$ is the nerve of a injectively fibrant diagram $p': \calI \rightarrow \bfA^{\degree}$. Let $\overline{p'}: \calI \star \{x\} \rightarrow \bfA^{\calI}$ be a limit of $p'$, so that $\overline{p}'$ is a homotopy limit diagram in $\bfA$. Now choose a trivial fibration $\overline{p}'' \rightarrow \overline{p}'$ in $\bfA^{\calI}$, where $\overline{p}''$ is cofibrant. The simplicial nerve of $\overline{p}''$ determines a colimit diagram
$\overline{f}: \Nerve(\calI)^{\triangleright} \rightarrow S$, by Theorem \ref{colimcomparee}. We now
observe that $f = \overline{f} | \Nerve(\calI)$ is equivalent to $p$, so that $p$ also admits a colimit in $S$.
\end{proof}

\section{Kan Extensions}\label{relacoim}
\setcounter{theorem}{0}

Let $\calC$ and $\calI$ be ordinary categories. There is an obvious ``diagonal'' functor $\delta: \calC \rightarrow \calC^{\calI}$, which carries an object $C \in \calC$ to the constant diagram $\calI \rightarrow \calC$ taking the value $C$. If $\calC$ admits small colimits, then the functor $\delta$ has a left adjoint $\calC^{\calI} \rightarrow \calC$. This left adjoint admits an explicit description: it carries an arbitrary diagram $f: \calI \rightarrow \calC$ to the colimit $\colim(f)$. Consequently, we can think of the theory of colimits as the study of left adjoints to diagonal functors.\index{gen}{diagonal functor}

More generally, if one is given a functor $i: \calI \rightarrow \calI'$ between diagram categories, then composition with $i$ induces a functor $i^{\ast}: \calC^{\calI'} \rightarrow \calC^{\calI}$. Assuming that $\calC$ has a sufficient supply of colimits, one can construct a left adjoint to $i^{\ast}$.
We then refer to this left adjoint as {\it left Kan extension along $i$}.\index{gen}{Kan extension}

In this section, we will study the $\infty$-categorical analogue of the theory of left Kan extensions. 
In the extreme case where $\calI'$ is the one-object category $\ast$, this theory simply reduces to the theory of colimits introduced in \S \ref{limitcolimit}. Our primary interest will be at the opposite extreme, when $i$ is a fully faithful embedding; this is the subject of \S \ref{kanex}. We will treat the general case in \S \ref{bigkanext}.

With a view toward later applications, we will treat not only the theory of {\em absolute} left Kan extensions, but also a relative notion which works over a base simplicial set $S$. The most basic example is the case of a {\it relative colimit} which we study in \S \ref{relcol}.

\subsection{Relative Colimits}\label{relcol}

In \S \ref{limitcolimit}, we introduced the notions of limit and colimit for a diagram
$p: K \rightarrow \calC$ in an $\infty$-category $\calC$. For many applications, it is convenient to have a {\em relative} version of these notions, which makes reference not to an $\infty$-category $\calC$ but to an arbitrary inner fibration of simplicial sets.  

\begin{definition}\label{relcoldef}\index{gen}{colimit!relative}\index{gen}{$f$-colimit}
Let $f: \calC \rightarrow \calD$ be an inner fibration of simplicial sets, let
$\overline{p}: K^{\triangleright} \rightarrow \calC$ be diagram, and let $p = \overline{p}|K$.
We will say that $\overline{p}$ is an {\it $f$-colimit} of $p$ if the map
$$ \calC_{\overline{p}/} \rightarrow \calC_{p/} \times_{ \calD_{f p/} } \calD_{ f\overline{p}/} $$
is a trivial fibration of simplicial sets. In this case, we will also say that
$\overline{p}$ is an {\it $f$-colimit diagram}.
\end{definition}

\begin{remark}\label{suppwolf}
Let $f: \calC \rightarrow \calD$ and $\overline{p}: K^{\triangleright} \rightarrow \calC$ be
as in Definition \ref{relcoldef}. Then $\overline{p}$ is an $f$-colimit of $p = \overline{p}|K$ if and only if the map 
$$ \phi: \calC_{\overline{p}/} \rightarrow \calC_{p/} \times_{ \calD_{f p/} } \calD_{ f \overline{p}/} $$
is a categorical equivalence. The ``only if'' direction is clear. The converse follows from
Proposition \ref{sharpen} (which implies that $\phi$ is a left fibration), Proposition \ref{funkyfibcatfib} (which implies that $\phi$ is a categorical fibration), and the fact that a categorical fibration which is a categorical equivalence is a trivial Kan fibration.

Observe that Proposition \ref{sharpen} also implies that the map
$$ \calD_{f  \overline{p}/} \rightarrow \calD_{f p/}$$
is a left fibration. Using Propositions \ref{basechangefunky} and \ref{funkyfibcatfib}, we conclude that $ \calC_{p/} \times_{ \calD_{f p/} } \calD_{f \overline{p}/}$ is a homotopy fiber product of
$\calC_{p/}$ with $\calD_{f \overline{p}/}$ over $\calD_{f p/}$ (with respect to the Joyal model structure on $\sSet$). 
Consequently, we deduce that $\overline{p}$ is an $f$-colimit diagram if and only if the
diagram of simplicial sets
$$ \xymatrix{ \calC_{ \overline{p}/} \ar[r] \ar[d] & \calD_{ f \overline{p} / } \ar[d] \\
\calC_{ p/ } \ar[r] & \calD_{f p/ } }$$
is homotopy Cartesian.
\end{remark}

\begin{example}
Let $\calC$ be an $\infty$-category and $f: \calC \rightarrow \ast$ the projection of $\calC$ to a point. Then a diagram $\overline{p}: K^{\triangleright} \rightarrow \calC$ is an $f$-colimit if and only if it is a colimit in the sense of Definition \ref{defcolim}.
\end{example}

\begin{example}\label{exex1}
Let $f: \calC \rightarrow \calD$ be an inner fibration of simplicial sets, and let
$e: \Delta^1 = (\Delta^0)^{\triangleright} \rightarrow \calC$ be an edge of $\calC$. Then $e$ is an $f$-colimit if and only if it is $f$-coCartesian.
\end{example}

The following basic stability properties follow immediately from the definition:

\begin{proposition}\label{basrel}
\begin{itemize}

\item[$(1)$] Let $f: \calC \rightarrow \calD$ be a trivial fibration of simplicial sets. Then
every diagram $\overline{p}: K^{\triangleright} \rightarrow \calC$ is an $f$-colimit.

\item[$(2)$] Let $f: \calC \rightarrow \calD$ and $g: \calD \rightarrow \calE$ be inner fibrations of simplicial sets, and let $\overline{p}: K^{\triangleright} \rightarrow \calC$ be a diagram.
Suppose that $f \circ \overline{p}$ is a $g$-colimit. Then $\overline{p}$ is an $f$-colimit if and only if $\overline{p}$ is a $g \circ f$-colimit.

\item[$(3)$] Let $f: \calC \rightarrow \calD$ be an inner fibration of $\infty$-categories, and let
$\overline{p}, \overline{q}: K^{\triangleright} \rightarrow \calC$ be diagrams which are equivalent when viewed as objects of the $\infty$-category $\Fun(K^{\triangleright}, \calC)$. Then $\overline{p}$ is 
an $f$-colimit if and only if $\overline{q}$ is an $f$-colimit.

\item[$(4)$] Suppose given a Cartesian diagram
$$ \xymatrix{ \calC' \ar[d]^{f'} \ar[r]^{g} & \calC \ar[d]^{f} \\
\calD' \ar[r] & \calD }$$
of simplicial sets, where $f$ (and therefore also $f'$) is an inner fibration. Let $\overline{p}: K^{\triangleright} \rightarrow \calC'$ be a diagram. If $g \circ \overline{p}$ is an $f$-colimit, then
$\overline{p}$ is an $f'$-colimit.
\end{itemize}
\end{proposition}

\begin{proposition}\label{summertoy}
Suppose give a commutative diagram of $\infty$-categories
$$ \xymatrix{ \calC \ar[r]^{f} \ar[d]^{p} & \calC' \ar[d]^{p'} \\
\calD \ar[r] & \calD' }$$
where the horizontal arrows are categorical equivalences and the vertical arrows are inner fibrations.
Let $\overline{q}: K^{\triangleright} \rightarrow \calC$ be a diagram and let $q = \overline{q}|K$
Then $\overline{q}$ is a $p$-colimit of $q$ if and only if
$f \circ \overline{q}$ is a $p'$-colimit of $f \circ q$.
\end{proposition}

\begin{proof}
Consider the diagram
$$ \xymatrix{ \calC_{ \overline{q}/ } \ar[r] \ar[d] & \calC'_{ f  \overline{q}/} \ar[d] \\
\calC_{q/} \times_{\calD_{p  q/}} \calD_{p  \overline{q}/} \ar[r] &
\calC'_{f  q/} \times_{ \calD'_{p'  f  q/}} \calD'_{p'  f  \overline{q}/} }.$$
According to Remark \ref{suppwolf}, it will suffice to show that the left vertical map is a categorical equivalence if and only if the right vertical map is a categorical equivalence. For this, it suffices to show that both of the horizontal maps are categorical equivalences. Proposition \ref{gorban3}
implies that the maps $ \calC_{ \overline{q}/} \rightarrow \calC'_{ f  \overline{q}/ }$,
$ \calC_{q/ } \rightarrow \calC'_{f  q/}$, $\calD_{p  \overline{q}/} \rightarrow \calD'_{p'  f  \overline{q}/}$, and $ \calD_{p  q/} \rightarrow \calD'_{p'  f  q/}$
are categorical equivalences. It will therefore suffice to show that the diagrams
$$ \xymatrix{ \calC_{q/} \times_{ \calD_{p  q/}} \calD_{p  \overline{q}/} \ar[r] \ar[d] &
\calC_{q/} \ar[d] &
\calC'_{ f  q/} \times_{ \calD'_{p'  f  q/} } \calD'_{ p'  f  \overline{q}/}
\ar[r] \ar[d] & \calC'_{ f  q /} \ar[d] \\
\calD_{p  \overline{q}/} \ar[r]^{\psi} & \calD_{p  q/} & \calD'_{p'  f  \overline{q}/}
\ar[r]^{\psi'} & \calD'_{p'  f  q/} }$$
are homotopy Cartesian (with respect to the Joyal model structure). This follows from
Proposition \ref{basechangefunky}, since $\psi$ and $\psi'$ are coCartesian fibrations.
\end{proof}

The next pair of results can be regarded as a generalization of Proposition \ref{gute}. They assert that, when computing relative colimits, we are free to replace any diagram by a cofinal subdiagram.

\begin{proposition}\label{relexists}
Let $p: \calC \rightarrow \calD$ be an inner fibration of $\infty$-categories, let $i: A \rightarrow B$ be a cofinal map, and let $\overline{q}: B^{\triangleright} \rightarrow \calC$ be a diagram.
Then $\overline{q}$ is a $p$-colimit if and only if $\overline{q} \circ i^{\triangleright}$
is a $p$-colimit.
\end{proposition}

\begin{proof}
Recall (Remark \ref{suppwolf}) that $\overline{q}$ is a relative colimit diagram if and only if the diagram
$$ \xymatrix{ \calC_{\overline{q}/} \ar[r] \ar[d] & \calC_{q/} \ar[d] \\
\calD_{ \overline{q}_0/} \ar[r] & \calD_{q_0/} }$$
is homotopy Cartesian with respect to the Joyal model structure. Since $i$ and $i^{\triangleright}$ are both cofinal, this is equivalent to the assertion that the diagram
$$ \xymatrix{ \calC_{\overline{q}  i^{\triangleright} /} \ar[r] \ar[d] & \calC_{q  i/} \ar[d] \\
\calD_{\overline{q}_0  i^{\triangleright} /} \ar[r] & \calD_{q_0  i/} }$$
is homotopy Cartesian, which (by Remark \ref{suppwolf}) is equivalent to the assertion that
$\overline{q} \circ i^{\triangleright}$ is a relative colimit diagram.
\end{proof}

\begin{proposition}\label{relexist}
Let $p: \calC \rightarrow \calD$ be a coCartesian fibration of $\infty$-categories, let $i: A \rightarrow B$ be a cofinal map, and let 
$$ \xymatrix{ B \ar[r]^{q} \ar[d] & \calC \ar[d]^{p} \\
B^{\triangleright} \ar[r]^{\overline{q}_0} & \calD }$$
be a diagram. Suppose that $q \circ i$
has a relative colimit lifting $\overline{q}_0 \circ i^{\triangleright}$. Then $q$ has a relative colimit lifting $\overline{q}_0$.
\end{proposition}

\begin{proof}
Let $q_0 = \overline{q}_0 | B$. We have a commutative diagram
$$ \xymatrix{ \calC_{q/} \ar[r]^-{f} \ar[d] & \calC_{q  i/} \times_{ \calD_{p  q  i/} } \calD_{p  q/} \ar[r] \ar[d] & \calC_{q  i/} \ar[d] \\
\calD_{q_0/} \ar[r] & \calD_{q_0/} \ar[r] & \calD_{q_0  i/} } $$
where the horizontal maps are categorical equivalences (since $i$ is cofinal, and by
Proposition \ref{basechangefunky}). Proposition \ref{werylonger} implies that the vertical maps are coCartesian fibrations, and that $f$ preserves coCartesian edges. Applying Proposition \ref{apple1} to $f$, we deduce that the map
$\phi: \calC_{q/} \times_{ \calD_{q_0/} } \{ \overline{q}_0 \} \rightarrow
\calC_{q  i/ } \times_{ \calD_{q_0  i/} } \{ \overline{q}_0  i^{\triangleright} \}$
is a categorical equivalence. Since $\phi$ is essentially surjective, we conclude that
there exists an extension $\overline{q}: B^{\triangleright} \rightarrow \calC$ of $q$
which covers $\overline{q}_0$, such that $\overline{q} \circ i^{\triangleright}$ is a $p$-colimit diagram. We now apply Proposition \ref{relexists} to conclude that $\overline{q}$ is itself a $p$-colimit diagram.
\end{proof}

Let $p: X \rightarrow S$ be a coCartesian fibration. 
The following results will allow us to reduce the theory of $p$-colimits to the theory of ordinary colimits in the fibers of $p$.

\begin{proposition}\label{chocolatelast}
Let $p: X \rightarrow S$ be an inner fibration of $\infty$-categories, $K$ a simplicial set, and
$\overline{h}: \Delta^1 \times K^{\triangleright} \rightarrow X$ a natural transformation from
$\overline{h}_0 = \overline{h} | \{0\} \times K^{\triangleright}$ to $\overline{h}_1 = \overline{h} | \{1\} \times K^{\triangleright}$.
Suppose that:
\begin{itemize}
\item[$(1)$] For every vertex $x$ of $K^{\triangleright}$, the restriction
$\overline{h} | \Delta^1 \times \{x\}$ is a $p$-coCartesian edge of $X$.
\item[$(2)$] The composition
$$ \Delta^1 \times \{\infty\} \subseteq \Delta^1 \times K^{\triangleright}
\stackrel{\overline{h}}{\rightarrow} X \stackrel{p}{\rightarrow} S$$
is a degenerate edge of $S$, where $\infty$ denotes the cone point of
$K^{\triangleright}$.
\end{itemize}

Then $h_0$ is a $p$-colimit if and only if $h_1$ is a $p$-colimit.
\end{proposition}

\begin{proof}
Let $h = \overline{h} | \Delta^1 \times K$, $h_0 = h | \{0\} \times K$, and $h_1 = h | \{1\} \times K$.
Consider the diagram
$$ \xymatrix{ X_{\overline{h}_0/} \ar[d] & X_{\overline{h}/} \ar[l]_{\phi} \ar[r] \ar[d] & X_{\overline{h}_1/} \ar[d] \\
X_{h_0/} \times_{S_{p  h_0/} } S_{p  \overline{h}_0} &
X_{h/} \times_{S_{p  h/}} S_{p  \overline{h}/} \ar[r] \ar[l]_{\psi} &
X_{h_1/} \times_{S_{p  h_1/}} S_{p  \overline{h}_1/} }$$
According to Remark \ref{suppwolf}, it will suffice to show that the left vertical map is a categorical equivalence if and only if the right vertical map is a categorical equivalence. For this, it will suffice to show that each of the horizontal arrows is a categorical equivalence. Because the inclusions
$\{1\} \times K \subseteq \Delta^1 \times K$ and $\{1\} \times K^{\triangleright} \subseteq
\Delta^1 \times K^{\triangleright}$ are right anodyne, the horizontal maps on the right are trivial fibrations. We are therefore reduced to proving that $\phi$ and $\psi$ are categorical equivalences.

Let $f: x \rightarrow y$ denote the edge of $X$ obtained by restricting $\overline{h}$ to the cone point of $K^{\triangleright}$. The map $\phi$ fits into a commutative diagram
$$ \xymatrix{ X_{\overline{h}/} \ar[r]^{\phi} \ar[d] & X_{h_0/} \ar[d] \\
X_{f/} \ar[r] & X_{x/}. }$$
Since the inclusion of the cone point into $K^{\triangleright}$ is right anodyne, the vertical arrows are trivial fibrations. Moreover, hypotheses $(1)$ and $(2)$ guarantee that $f$ is an equivalence in $X$, so that the map $X_{f/} \rightarrow X_{x/}$ is a trivial fibration. This proves that $\phi$ is a categorical equivalence.

The map $\psi$ admits a factorization
$$ X_{h/} \times_{S_{p  h/}} S_{p  \overline{h}/}
\stackrel{\psi'}{\rightarrow} 
X_{h_0/} \times_{ S_{p  h_0/}} S_{p  \overline{h}/}
\stackrel{\psi''}{\rightarrow}
X_{h_0} \times_{ S_{p  h_0/}} S_{ p  \overline{h}_0/}.$$
To complete the proof, it will suffice to show that $\psi'$ and $\psi''$ are trivial fibrations of simplicial sets. We first observe that $\psi'$ is a pullback of the map
$$X_{h/} \rightarrow X_{h_0/} \times_{S_{p  h_0/} } S_{p  h/},$$
which is a trivial fibration (Proposition \ref{eggwhite}). The map $\psi''$
is a pullback of the left fibration $\psi''_0: S_{p  \overline{h}/} \rightarrow S_{p  \overline{h}_0/}$. It therefore suffices to show that $\psi''_0$ is a categorical equivalence.
To prove this, we consider the diagram
$$ \xymatrix{ S_{p  \overline{h}/} \ar[r]^{\psi''_0} \ar[d] & S_{p  \overline{h}_0/} \ar[d] \\
S_{p(f)/} \ar[r]^{\psi''_1} & S_{p(x)/} }$$
As above, we observe that the vertical arrows are trivial fibrations, and $\psi''_1$ is a trivial fibration because the morphism $p(f)$ is an equivalence in $S$. It follows that $\psi''_0$ is a categorical equivalence, as desired.
\end{proof}

\begin{proposition}\label{relcolfibtest}
Let $q: X \rightarrow S$ be a locally coCartesian fibration of $\infty$-categories, let $s$ be an object of $S$, and let $\overline{p}: K^{\triangleright} \rightarrow X_{s}$ be a diagram. The following conditions are equivalent:
\begin{itemize}
\item[$(1)$] The map $\overline{p}$ is a $q$-colimit diagram.
\item[$(2)$] For every morphism $e: s \rightarrow s'$ in $S$, the associated functor
$e_{!}: X_{s} \rightarrow X_{s'}$ has the property that $e_{!} \circ \overline{p}$ is a colimit
diagram in the $\infty$-category $X_{s'}$.
\end{itemize}
\end{proposition}

\begin{proof}
Assertion $(1)$ is equivalent to the statement that the map
$$ \theta: X_{\overline{p}/} \rightarrow X_{p/} \times_{ S_{qp/} } S_{q \overline{p}/}$$
is a trivial fibration of simplicial sets. Since $\theta$ is a left fibration, it will suffice to show that the fibers of $\theta$ are contractible. Consider an arbitrary vertex of $S_{q \overline{p}/}$, corresponding to a morphism $t: K \star \Delta^1 \rightarrow S$. Since $K \star \Delta^1$
is categorically equivalent to $( K \star \{0\} ) \coprod_{ \{0\} } \Delta^1$ and
$t | K \star \{0\}$ is constant, we may assume without loss of generality that
$t$ factors as a composition
$$ K \star \Delta^1 \rightarrow \Delta^1 \stackrel{e}{\rightarrow} S.$$
Here $e: s \rightarrow s'$ is an edge of $S$. Pulling back by the map $e$, we can reduce to the problem of proving the following analogue of $(1)$ in the case where $S = \Delta^1$:
\begin{itemize}
\item[$(1')$] The projection $h_0: X_{ \overline{p}/} \times_{S} \{s' \} \rightarrow
X_{p/} \times_{S} \{ s' \}$ is a trivial fibration of simplicial sets.
\end{itemize}

Choose a coCartesian transformation $\overline{\alpha}: K^{\triangleright} \times \Delta^1 \rightarrow X$ from $\overline{p}$ to $\overline{p}'$, which covers the projection
$$ K^{\triangleright} \times \Delta^1 \rightarrow \Delta^1 \simeq S.$$
Consider the diagram
$$ \xymatrix{ X_{ \overline{p}/ } \times_{S} \{s' \} 
\ar[d]^{h_0} & X_{ \overline{\alpha} / } \times_{S} \{ s' \}
\ar[l] \ar[r] \ar[d]^{h} & X_{ \overline{p}' /} \times_{ S} \{ s' \}
\ar[d]^{h_1} \\
X_{p/} \times_{S } \{s'\} &
X_{\alpha/} \times_{ S } \{s'\} \ar[r] \ar[l] &
X_{p'/} \times_{ S } \{s'\}. }$$
Note that the vertical maps are left fibrations (Proposition \ref{sharpen}). Since the inclusion
$K^{\triangleright} \times \{1\} \subseteq K^{\triangleright} \times \Delta^1$ is right anodyne,
the upper right horizontal map is a trivial fibration. Similarly, the lower right horizontal map is a trivial fibration. Since $\overline{\alpha}$ is a coCartesian transformation, we deduce that the left  horizontal maps are also trivial fibrations (Proposition \ref{eggwhite}).
Condition $(2)$ is equivalent to the assertion that
$h_1$ is a trivial fibration (for each edge $e: s \rightarrow s'$ of the original simplicial set $S$). Since $h_1$ is a left fibration, and therefore a categorical fibration (Proposition \ref{funkyfibcatfib}), this is equivalent to the assertion that $h_1$ is a categorical equivalence. Chasing through the diagram, we deduce that
$(2)$ is equivalent to the assertion that $h_0$ is a categorical equivalence, which (by the same argument) is equivalent to the assertion that $h_0$ is a trivial fibration.
\end{proof}

\begin{corollary}\label{constrel}
Let $p: X \rightarrow S$ be a coCartesian fibration of $\infty$-categories, and let $K$ be a simplicial set.
Suppose that:
\begin{itemize}
\item[$(1)$] For each vertex $s$ of $S$, the fiber $X_{s} = X \times_{S} \{s\}$ admits colimits
for all diagrams indexed by $K$.
\item[$(2)$] For each edge $f: s \rightarrow s'$, the associated functor
$X_{s} \rightarrow X_{s'}$ preserves colimits of $K$-indexed diagrams.
\end{itemize}
Then for every diagram
$$ \xymatrix{ K \ar[r]^{q} \ar@{^{(}->}[d] & X \ar[d]^{p} \\
K^{\triangleright} \ar[r]^{f} \ar@{-->}[ur]^{\overline{q}} & S }$$
there exists a map $\overline{q}$ as indicated, which is a $p$-colimit.
\end{corollary}

\begin{proof}
Consider the map $K \times \Delta^1 \rightarrow K^{\triangleright}$ which is the identity on $K \times \{0\}$ and carries $K \times \{1\}$ to the cone point of $K^{\triangleright}$. 
Let $F$ denote the composition
$$ K \times \Delta^1 \rightarrow K^{\triangleright} \stackrel{f}{\rightarrow} S,$$
and let $Q: K \times \Delta^1 \rightarrow X$ be a coCartesian lifting of $F$
to $X$, so that $Q$ is a natural transformation from $q$ to a map
$q': K \rightarrow X_{s}$, where $s$ is the image under $f$ of the cone point of
$K^{\triangleright}$. In view of assumption $(1)$, there exists a map
$\overline{q}': K^{\triangleright} \rightarrow X_{s}$ which is a colimit of $q'$.
Assumption $(2)$ and Proposition \ref{relcolfibtest} guarantee that
$\overline{q}'$ is also a $p$-colimit diagram, when regarded as a map
from $K^{\triangleright}$ to $X$. 

We have a commutative diagram
$$ \xymatrix{ (K \times \Delta^1) \coprod_{ K \times \{1\} } (K^{\triangleright}
\times \{1\} ) \ar[rr]^-{(Q,\overline{q}')} \ar@{^{(}->}[d] & & X \ar[d]^{p} \\
(K \times \Delta^1)^{\triangleright} \ar@{-->}[urr]^{r} \ar[rr] & & S. }$$
The left vertical map is an inner fibration, so there exists a morphism
$r$ as indicated, rendering the diagram commutative. We now consider the map
$ K^{\triangleright} \times \Delta^1 \rightarrow (K \times \Delta^1)^{\triangleright}$
which is the identity on $K \times \Delta^1$ and carries the other vertices of
$K^{\triangleright} \times \Delta^1$ to the cone point of $(K \times \Delta^1)^{\triangleright}$.
Let $\overline{Q}$ denote the composition
$$ K^{\triangleright} \times \Delta^1 \rightarrow (K \times \Delta^1)^{\triangleright}
\stackrel{r}{\rightarrow} X, $$
and let $\overline{q} = \overline{Q} | K^{\triangleright} \times \{0\}$.
Then $\overline{Q}$ can be regarded as a natural transformation
$\overline{q} \rightarrow \overline{q}'$ of diagrams $K^{\triangleright} \rightarrow X$.
Since $\overline{q}'$ is a $p$-colimit diagram, Proposition \ref{chocolatelast}
implies that $\overline{q}$ is a $p$-colimit diagram as well.
\end{proof}

\begin{proposition}\label{timal}
Let $p: X \rightarrow S$ be a coCartesian fibration of $\infty$-categories, and let
$\overline{q}: K^{\triangleright} \rightarrow X$ be a diagram. Assume that:
\begin{itemize}
\item[$(1)$] The map $\overline{q}$ carries each edge of $K$ to a $p$-coCartesian
edge of $K$.
\item[$(2)$] The simplicial set $K$ is weakly contractible.
\end{itemize}

Then $\overline{q}$ is a $p$-colimit diagram if and only if 
it carries every edge of $K^{\triangleright}$ to a $p$-coCartesian edge of $X$.
\end{proposition}

\begin{proof}
Let $s$ denote the image under $p \circ \overline{q}$ of the cone point of $K^{\triangleright}$.
Consider the map $K^{\triangleright} \times \Delta^1 \rightarrow K^{\triangleright}$
which is the identity on $K^{\triangleright} \times \{0\}$ and collapses $K^{\triangleright} \times \{1\}$ to the cone point of $K^{\triangleright}$. Let $h$ denote the composition
$$ K^{\triangleright} \times \Delta^1 \rightarrow K^{\triangleright}
\stackrel{\overline{q}}{\rightarrow} X \stackrel{p}{\rightarrow} S,$$
which we regard as a natural transformation from $p \circ \overline{q}$ to the constant
map with value $s$. Let $H: \overline{q} \rightarrow \overline{q}'$ be a coCartesian transformation
from $\overline{q}$ to a diagram $\overline{q}': K^{\triangleright} \rightarrow X_{s}$.
Using Proposition \ref{protohermes}, we conclude that $\overline{q}'$ carries each
edge of $K$ to a $p$-coCartesian edge of $X$, which is therefore an equivalence
in $X_{s}$. 

Let us now suppose that $\overline{q}$ carries {\em every} edge of
$K^{\triangleright}$ to a $p$-coCartesian edge of $X$. Arguing as above, we conclude
that $\overline{q}'$ carries each edge of $K^{\triangleright}$ to an equivalence in $X_{s}$.
Let $e: s \rightarrow s'$ be an edge of $S$ and $e_{!}: X_{s} \rightarrow X_{s'}$ an associated functor. The composition
$$ K^{\triangleright} \stackrel{\overline{q}'}{\rightarrow} X_{s} \stackrel{e_{!}}{\rightarrow}
X_{s'}$$
carries each edge of $K^{\triangleright}$ to an equivalence in $X_{s}$, and 
is therefore a colimit diagram in $X_{s'}$ (Corollary \ref{silt}). Proposition \ref{relcolfibtest} implies
that $\overline{q}'$ is a $p$-colimit diagram, so that Proposition \ref{chocolatelast} implies that $\overline{q}$ is a $p$-colimit diagram as well.

For the converse, let us suppose that $\overline{q}$ is a $p$-colimit diagram.
Applying Proposition \ref{chocolatelast}, we conclude that $\overline{q}'$ is a $p$-colimit diagram. In particular, $\overline{q}'$ is a colimit diagram in the $\infty$-category
$X_{s}$. Applying Corollary \ref{silt}, we conclude that $\overline{q}'$ carries each
edge of $K^{\triangleright}$ to an equivalence in $X_{s}$. Now consider an arbitrary
edge $f: x \rightarrow y$ of $K^{\triangleright}$. If $f$ belongs to $K$, then
$\overline{q}(f)$ is $p$-coCartesian by assumption. Otherwise, we may suppose that
$y$ is the cone point of $K$. The map $H$ gives rise to a diagram
$$ \xymatrix{ \overline{q}(x) \ar[r]^{\overline{q}(f)} \ar[d]^{\phi} & \overline{q}(y) \ar[d]^{\phi'} \\
\overline{q}'(x) \ar[r]^{\overline{q}'(f)} & \overline{q}'(y) }$$
in the $\infty$-category $X \times_{S} \Delta^1$. Here 
$\overline{q}'(f)$ and $\phi'$ are equivalences in $X_{s}$, so that
$\overline{q}(f)$ and $\phi$ are equivalent as morphisms
$\Delta^1 \rightarrow X \times_{S} \Delta^1$. Since $\phi$ is $p$-coCartesian, we conclude
that $\overline{q}(f)$ is $p$-coCartesian, as desired.
\end{proof}

\begin{lemma}\label{gooodbar}
Let $p: \calC \rightarrow \calD$ be an inner fibration of $\infty$-categories, let
$C \in \calC$ be an object, and let $D = p(C)$.  Then $C$ is a $p$-initial object of $\calC$ if and only if $(C, \id_{D})$ is an initial object of $\calC \times_{ \calD} \calD_{D/}$.
\end{lemma}

\begin{proof}
We have a commutative diagram
$$ \xymatrix{ \calC_{C/} \times_{ \calD_{D/} } \calD_{ \id_{D}/} \ar[r]^{\psi} \ar[d]^{\phi} & \calC_{C/} \ar[d]^{\phi'} \\
\calC \times_{ \calD} \calD_{D/} \ar@{=}[r] & \calC \times_{\calD} \calD_{D/} }$$
where the vertical arrows are left fibrations, and therefore categorical fibrations (Proposition \ref{funkyfibcatfib}). We wish to show that $\phi$ is a trivial fibration if and only if $\phi'$ is a trivial fibration. This is equivalent to proving that $\phi$ is a categorical equivalence if and only if $\phi'$ is a categorical equivalence. For this, it will suffice to show that $\psi$ is a categorical equivalence.
But $\psi$ is a pullback of the trivial fibration $\calD_{\id_{D}/} \rightarrow \calD_{D/}$, and therefore itself a trivial fibration.
\end{proof}

\begin{proposition}\label{panna}
Suppose given a diagram of $\infty$-categories
$$ \xymatrix{ \calC \ar[dr]^{q} \ar[rr]^{p} & & \calD \ar[dl]_{r} \\
& \calE & }$$
where $p$ and $r$ are inner fibrations, $q$ is a Cartesian fibration, and
$p$ carries $q$-Cartesian morphisms to $r$-Cartesian morphisms.

Let $C \in \calC$ be an object, $D = p(C)$, and $E = q(C)$.
Let $\calC_{E} = \calC \times_{\calE} \{E\}$, $\calD_{E} = \calD \times_{\calE} \{E\}$, and
$p_{E}: \calC_{E} \rightarrow \calD_{E}$ the induced map. Suppose that
$C$ is a $p_{E}$-initial object of $\calC_{E}$. Then $C$ is a $p$-initial object of $\calC$.
\end{proposition}

\begin{proof}
Our hypothesis, together with Lemma \ref{gooodbar}, implies that
$(C, \id_{D})$ is an initial object of
$$\calC_{E} \times_{\calD_{E}} (\calD_{E})_{D/}
\simeq (\calC \times_{\calD} \calD_{D/} ) \times_{\calE_{E/}} \{ \id_{E} \}.$$
We will prove that the map
$\phi: \calC \times_{\calD} \calD_{D/} \rightarrow \calE_{E/}$ is a Cartesian fibration.
Since $\id_{E}$ is an initial object of $\calE_{E/}$, Lemma \ref{sabreto} will allow us to conclude that
$(C, \id_{D})$ is an initial object of $\calC \times_{\calD} \calD_{D/}$. We can then conclude the proof by applying Lemma \ref{gooodbar} once more.

It remains to prove that $\phi$ is a Cartesian fibration. Let us say that a morphism of
$\calC \times_{\calD} \calD_{D/}$ is {\it special} if its image in $\calC$ is $q$-Cartesian.
Since $\phi$ is obviously an inner fibration, it will suffice to prove the following assertions:

\begin{itemize}
\item[$(1)$] Given an object $X$ of $\calC \times_{\calD} \calD_{D/}$ and
a morphism $\overline{f}: \overline{Y} \rightarrow \phi(X)$ in $\calE_{E/}$, we can write
$\overline{f} = \phi(f)$ where $f$ is a special morphism of $\calC \times_{\calD} \calD_{D/}$.

\item[$(2)$] Every special morphism in $\calC \times_{\calD} \calD_{D/}$ is $\phi$-Cartesian.
\end{itemize}

To prove $(1)$, we first identify $X$ with a pair consisting of an object $C'' \in \calC$ and
a morphism $D \rightarrow p(C'')$ in $\calD$, and $\overline{f}$ with a $2$-simplex
$\overline{\sigma}: \Delta^2 \rightarrow \calE$ which we depict as a diagram:
$$ \xymatrix{ & E' \ar[dr]^{\overline{g}} & \\
E \ar[ur] \ar[rr] & & q(C''). }$$
Since $q$ is a Cartesian fibration, the morphism $\overline{g}$ can be written
as $q(g)$ for some morphism $g: C' \rightarrow C''$ in $\calC$. We now have
a diagram
$$ \xymatrix{ & p(C') \ar[dr]^{p(g)} & \\
D \ar[rr] & & p(C'')  }$$
in $\calD$. Since $p$ carries $q$-Cartesian morphisms to $r$-Cartesian morphisms,
we conclude that $p(g)$ is $r$-Cartesian, so that the above diagram can be completed to a $2$-simplex $\sigma: \Delta^2 \rightarrow \calD$ such that $r(\sigma) = \overline{\sigma}$.

We now prove $(2)$. Suppose $n \geq 2$, and we have a commutative diagram
$$ \xymatrix{ \Lambda^n_n \ar[r]^-{\sigma_0} \ar@{^{(}->}[d] & \calC \times_{\calD} \calD_{D/} \ar[d] \\
\Delta^n \ar[r] \ar@{-->}[ur]^{\sigma} & \calE_{E/} }$$
where $\sigma_0$ carries the final edge of $\Lambda^n_n$ to a special morphism of
$\calC \times_{\calD} \calD_{D/}$. We wish to prove the existence of the morphism $\sigma$
indicated in the diagram. We first let $\tau_0$ denote the composite map
$$ \Lambda^n_n \stackrel{\sigma_0}{\rightarrow} \calC \times_{\calD} \calD_{D/} \rightarrow \calC.$$
Consider the diagram
$$ \xymatrix{ \Lambda^n_n \ar[r]^{\tau_0} \ar@{^{(}->}[d] & \calC \ar[d]^{q} \\
\Delta^n \ar[r] \ar@{-->}[ur]^{\tau} & \calE. }$$
Since $\tau_0( \Delta^{ \{n-1, n\} })$ is $q$-Cartesian, there exists an extension
$\tau$ as indicated in the diagram. The morphisms $\tau$ and $\sigma_0$ together determine a map $\theta_0$ which fits into a diagram
$$ \xymatrix{ \Lambda^{n+1}_{n+1} \ar[r]^{\theta_0} \ar@{^{(}->}[d] & \calD \ar[d]^{r} \\
\Delta^{n+1} \ar[r] \ar@{-->}[ur]^{\theta} & \calE. }$$
To complete the proof, it suffices to prove the existence of the indicated arrow $\theta$.
This follows from the fact that $\theta_0( \Delta^{ \{n,n+1 \} }) = (p \circ \tau_0)( \Delta^{ \{n-1,n\} })$
is an $r$-Cartesian morphism of $\calD$.
\end{proof}

Proposition \ref{panna} immediately implies the following slightly stronger statement:

\begin{corollary}\label{pannaheave}
Suppose given a diagram of $\infty$-categories
$$ \xymatrix{ \calC \ar[dr]^{q} \ar[rr]^{p} & & \calD \ar[dl]_{r} \\
& \calE & }$$
where $q$ and $r$ are Cartesian fibrations, $p$ is an inner fibration, and 
$p$ carries $q$-Cartesian morphisms to $r$-Cartesian morphisms.

Suppose given another $\infty$-category $\calE_0$ equipped with a functor
$s: \calE_0 \rightarrow \calE$. Set $\calC_0 = \calC \times_{\calE} \calE_0$,
$\calD_0 = \calD \times_{\calE} \calE_0$, and let $p_0: \calC_0 \rightarrow \calD_0$ be the functor induced by $p$.
Let $\overline{f}_0: K^{\triangleright} \rightarrow \calC_0$ be a diagram and let
$\overline{f}$ denote the composition $K^{\triangleright} \stackrel{\overline{f}_0}{\rightarrow}
\calC_0 \rightarrow \calC$. Then $\overline{f}_0$ is a $p_0$-colimit diagram if and only if
$\overline{f}$ is a $p$-colimit diagram.
\end{corollary}

\begin{proof}
Let $f_0 = \overline{f}_0 | K$ and $f = \overline{f} | K$. Replacing our diagram by
$$ \xymatrix{ \calC_{f/} \ar[rr] \ar[dr] & & \calD_{pf/} \ar[dl] \\
& \calE_{qf/}, & }$$
we can reduce to the case where $K = \emptyset$. Then $\overline{f}_0$ determines
an object $C \in \calC_0$. Let $E$ denote the image of $C$ in $\calE_0$. We have a commutative diagram
$$ \xymatrix{ \{E\} \ar[rr]^{s'} \ar[dr]^{s''} & & \calE_0 \ar[dl]^{s} \\
& \calE. & }$$
Consequently, to prove Corollary \ref{pannaheave} for the map $s$, it will suffice to prove the
analogous assertions for $s'$ and $s''$; these follow from Proposition \ref{panna}.
\end{proof}

\begin{corollary}\label{superduck}
Let $p: \calC \rightarrow \calE$ be a Cartesian fibration of $\infty$-categories, 
$E \in \calE$ an object, and $\overline{f}: K^{\triangleright} \rightarrow \calC_{E}$ a diagram.
Then $\overline{f}$ is a colimit diagram in $\calC_{E}$ if and only if
it is a $p$-colimit diagram in $\calC$.
\end{corollary}

\begin{proof}
Apply Corollary \ref{pannaheave} in the case where $\calD = \calE$.
\end{proof}

\subsection{Kan Extensions along Inclusions}\label{kanex}

In this section, we introduce the theory of {\em left Kan extensions}. Let $F: \calC \rightarrow \calD$
be a functor between $\infty$-categories, and let $\calC^{0}$ be a full subcategory of $\calC$.
Roughly speaking, the functor $F$ is a left Kan extension of its restriction $F_0 = F | \calC^{0}$ if
the values of $F$ are as ``small'' as possible, given the values of $F_0$. In order to make this precise, we need to introduce a bit of terminology.

\begin{notation}
Let $\calC$ be an $\infty$-category, and let $\calC^{0}$ be a full subcategory. If $p:
K \rightarrow \calC$ is a diagram, we let $\calC^{0}_{/p}$ denote the fiber product
$\calC_{/p} \times_{ \calC} \calC^{0}$. In particular, if $C$ is an object of $\calC$,
then $\calC^{0}_{/C}$ denotes the full subcategory of $\calC_{/C}$ spanned by the morphisms $C' \rightarrow C$ where $C' \in \calC^{0}$.
\end{notation}

\begin{definition}\label{defKan}
Suppose given a commutative diagram of $\infty$-categories
$$ \xymatrix{ \calC^{0} \ar@{^{(}->}[d] \ar[r]^{F_0} & \calD \ar[d]^{p} \\
\calC \ar[r] \ar[ur]^{F} & \calD', }$$
where $p$ is an inner fibration and
the left vertical map is the inclusion of a full subcategory $\calC^{0} \subseteq \calC$.

We will say that $F$ is a {\it $p$-left Kan extension of $F_0$ at $C \in \calC$} if
the induced diagram
$$ \xymatrix{ (\calC^{0}_{/C}) \ar@{^{(}->}[d] \ar[r]^{F_C} & \calD \ar[d]^{p} \\
(\calC_{/C}^{0})^{\triangleright} \ar[ur] \ar[r] & \calD' }$$
exhibits $F(C)$ as a $p$-colimit of $F_{C}$.

We will say that $F$ is a {\it $p$-left Kan extension of $F_0$} if it is a $p$-left Kan extension
of $F_0$ at $C$, for every object $C \in \calC$.

In the case where $\calD' = \Delta^0$, we will omit mention of $p$ simply say that $F$ is a {\it left Kan extension of $F_0$} if the above condition is satisfied.\index{gen}{Kan extension}
\end{definition}

\begin{remark}\label{mozartwatch}
Consider a diagram
$$ \xymatrix{ \calC^{0} \ar@{^{(}->}[d] \ar[r]^{F_0} & \calD \ar[d]^{p} \\
\calC \ar[r] \ar[ur]^{F} & \calD' }$$
as in Definition \ref{defKan}. If
$C$ is an object of $\calC^{0}$, then the functor $F_{C}: (\calC_{/C}^{0})^{\triangleright} \rightarrow \calD$ is automatically a $p$-colimit. To see this, we observe that $\id_{C}: C \rightarrow C$ is a final object of $\calC^{0}_{/C}$. Consequently, the inclusion $\{ \id_C \} \rightarrow (\calC_{/C}^{0})$
is cofinal and we are reduced to proving that $F(\id_{C}): \Delta^1 \rightarrow \calD$ is a colimit
of its restriction to $\{0\}$, which is obvious.
\end{remark}

\begin{example}
Consider a diagram
$$ \xymatrix{ \calC \ar@{^{(}->}[d] \ar[r]^{q} & \calD \ar[d]^{p} \\
\calC^{\triangleright} \ar[r] \ar[ur]^{\overline{q}} & \calD'. }$$
The map $\overline{q}$ is a $p$-left Kan extension of $q$ if and only if it is a $p$-colimit of $q$.
The ``only if'' direction is clear from the definition, and the converse follows immediately from 
Remark \ref{mozartwatch}.
\end{example}

We first note a few basic stability properties for the class of left Kan extensions.

\begin{lemma}\label{switcher}
Consider a commutative diagram of $\infty$-categories
$$ \xymatrix{ \calC^{0} \ar@{^{(}->}[d] \ar[r]^{F_0} & \calD \ar[d]^{p} \\
\calC \ar[r] \ar[ur]^{F} & \calD' }$$
as in Definition \ref{defKan}. Let $C$ and $C'$ equivalent objects of $\calC$.
Then $F$ is a $p$-left Kan extension of $F_0$ at $C$ if and only if 
$F$ is a $p$-left Kan extension of $F_0$ at $C'$.
\end{lemma}

\begin{proof}
Let $f: C \rightarrow C'$ be an equivalence, so that the restriction maps
$$ \calC_{/C} \leftarrow \calC_{/f} \rightarrow \calC_{/C'}$$
are trivial fibrations of simplicial sets. Let $\calC^{0}_{/f} = \calC^{0} \times_{\calC} \calC_{/f}$, so that we have trivial fibrations
$$ \calC^{0}_{/C} \stackrel{g}{\leftarrow} \calC^{0}_{/f} \stackrel{g'}{\rightarrow} \calC^{0}_{/C'}.$$
Consider the associated diagram
$$ \xymatrix{ & (\calC^{0}_{/C})^{\triangleright} \ar[dr]^{F_{C}} & \\
(\calC^{0}_{/f})^{\triangleright} \ar[ur]^{G} \ar[dr]^{G'} & & \calD \\
& (\calC^{0}_{/C'})^{\triangleright} \ar[ur]^{F_{C'}} & .}$$
This diagram does not commute, but the functors
$F_{C} \circ G$ and $F_{C'} \circ G'$ are equivalent in the $\infty$-category
$\calD^{ (\calC^{0}_{/f})^{\triangleright}}$. Consequently,
$F_{C} \circ G$ is a $p$-colimit diagram if and only if $F_{C'} \circ G'$ is a $p$-colimit diagram
(Proposition \ref{basrel}).
Since $g$ and $g'$ are cofinal, we conclude that $F_{C}$ is a $p$-colimit diagram if and only if $F_{C'}$ is a $p$-colimit diagram (Proposition \ref{relexists}).
\end{proof}

\begin{lemma}\label{basekann}
\begin{itemize}
\item[$(1)$] Let $\calC$ be an $\infty$-category, $p: \calD \rightarrow \calD'$ an inner fibration of $\infty$-categories, and $F,F': \calC \rightarrow \calD$ be two functors which are equivalent
in $\calD^{\calC}$. Let $\calC^{0}$ be a full subcategory of $\calC$. Then $F$ is
a $p$-left Kan extension of $F| \calC^{0}$ if and only if $F'$ is a $p$-left Kan extension of $F'| \calC^{0}$.

\item[$(2)$] Suppose given a commutative diagram of $\infty$-categories
$$ \xymatrix{ \calC^{0} \ar[d]^{G_0} \ar[r] & \calC \ar[r]^{F} \ar[d]^{G} & \calD \ar[d] \ar[r]^{p} & \calE \ar[d] \\
{\calC'}^{0} \ar[r] & \calC' \ar[r]^{F'} & \calD' \ar[r]^{p'} & \calE' }$$
be a commutative diagram of $\infty$-categories, where the left horizontal maps are inclusions of full subcategories, the right horizontal maps are inner fibrations,
and the vertical maps are categorical equivalences. Then
$F$ is a $p$-left Kan extension of $F| \calC^{0}$ if and only if $F'$ is a $p'$-left Kan extension of $F' | {\calC'}^{0}$.
\end{itemize}
\end{lemma}

\begin{proof}
Assertion $(1)$ follows immediately from Proposition \ref{basrel}.
Let us prove $(2)$. Choose an object $C \in \calC$, and consider the diagram
$$ \xymatrix{ (\calC^{0}_{/C})^{\triangleright} \ar[r] \ar[d] & \calD \ar[d] \ar[r]^{p} & \calE \ar[d] \\
({\calC'}^{0}_{/G(C)})^{\triangleright} \ar[r] & \calD' \ar[r]^{p'} & \calE' }$$
We claim that the upper left horizontal map is a $p$-colimit diagram if and only if the bottom left horizontal map is a $p'$-colimit diagram.
In view of Proposition \ref{summertoy}, it will suffice to show that each of the vertical maps is an equivalence
of $\infty$-categories. For the middle and right vertical maps, this holds by assumption.
To prove that the left vertical map is a categorical equivalence, we consider the diagram
$$ \xymatrix{ \calC^{0}_{/C} \ar[r] \ar[d] & {\calC'}^{0}_{/G(C)} \ar[d] \\
\calC_{/C} \ar[r] & \calC'_{/G(C)}. }$$
The bottom horizontal map is a categorical equivalence by Proposition \ref{gorban3}, and the vertical maps are inclusions of full subcategories. It follows that the top horizontal map is fully faithful, and its essential image consists of those morphisms $C' \rightarrow G(C)$ where
$C'$ is equivalent (in $\calC'$) to the image of an object of $\calC^{0}$. Since $G_0$ is essentially surjective, this is the whole of ${\calC'}^{0}_{/G(C)}$. 

It follows that if $F'$ is a $p'$-left Kan extension of $F' | { \calC'}^{0}$, then $F$ is a $p$-left Kan extension of $F | \calC^{0}$. Conversely, if $F$ is a $p$-left Kan extension of $F| \calC^{0}$, then $F'$ is a $p'$-left Kan extension of $F' | { \calC'}^{0} $ at $G(C)$, for every object $C \in \calC$. Since $G$ is essentially surjective, Lemma \ref{switcher} implies that $F'$ is a $p'$-left Kan extension of $F' | {\calC'}^{0}$ at every object of $\calC'$. This completes the proof of $(2)$.
\end{proof}

\begin{lemma}\label{kan0}
Suppose given a diagram of $\infty$-categories
$$ \xymatrix{ \calC^{0} \ar@{^{(}->}[d] \ar[r]^{F_0} & \calD \ar[d]^{p} \\
\calC \ar[r] \ar[ur]^{F} & \calD' }$$
as in Definition \ref{defKan}, where $F$ is a left Kan extension of $F_0$ relative
to $p$. Then the induced map
$$ \calD_{F/} \rightarrow \calD'_{p F/ } \times_{ \calD'_{p F_0/} } \calD_{F_0/}$$
is a trivial fibration of simplicial sets. In particular, we may identify $p$-colimits of $F$ with $p$-colimits of $F_0$.
\end{lemma}

\begin{proof}
Using Lemma \ref{basekann}, we may reduce to the case where $\calC$ is minimal.
Let us call a simplicial subset $\calE \subseteq \calC$ {\it complete} if it has the following property:
for any simplex $\sigma: \Delta^n \rightarrow \calC$, if $\sigma| \Delta^{ \{0, \ldots, i\} }$ factors through $\calC^{0}$ and $\sigma| \Delta^{ \{i+1, \ldots, n \} }$ factors through $\calE$, then $\sigma$
factors through $\calE$. Note that if $\calE$ is complete, then $\calC^{0} \subseteq \calE$.
We next define a transfinite sequence of {\em complete} simplicial subsets of $\calC$
$$ \calC^0 \subseteq \calC^1 \subseteq \ldots $$
as follows: if $\lambda$ is a limit ordinal, we let $\calC^{\lambda} = \bigcup_{ \alpha < \lambda} \calC^{\alpha}$. If $\calC^{\alpha} = \calC$, then we set $\calC^{\alpha+1} = \calC$. Otherwise, choose some simplex $\sigma: \Delta^n \rightarrow \calC$ which does not belong to $\calC^{\alpha}$, where the dimension $n$ of $\sigma$ is chosen as small as possible, and let
$\calC^{\alpha+1}$ be the smallest complete simplicial subset of $\calC$ containing
$\calC^{\alpha}$ and the simplex $\sigma$.

Let $F_{\alpha} = F | \calC^{\alpha}$. We will prove that for every
$\beta \leq \alpha$ the projection
$$ \phi_{\alpha,\beta}: \calD_{F_{\alpha}/} \rightarrow \calD'_{p F_{\alpha}/}
\times_{ \calD'_{ p F_{\beta} / }} \calD_{F_{\beta} / }$$
is a trivial fibration of simplicial sets. Taking $\alpha \gg \beta = 0$, we have $\calC^{\alpha} = \calC$ and the proof will be complete.

Our proof proceeds by induction on $\alpha$. If $\alpha = \beta$, then $\phi_{\alpha,\beta}$ is an isomorphism and there is nothing to prove. If $\alpha > \beta$ is a limit ordinal, then the inductive hypothesis implies that $\phi_{\alpha,\beta}$ is the inverse limit of a transfinite tower of trivial fibrations, and therefore a trivial fibration. 
It therefore suffices to prove that if $\phi_{\alpha,\beta}$ is a trivial fibration, then $\phi_{\alpha+1, \beta}$ is a trivial fibration. We observe that $\phi_{\alpha+1,\beta} = \phi'_{\alpha,\beta} \circ \phi_{\alpha+1,\alpha}$, where $\phi'_{\alpha,\beta}$ is a pullback of $\phi_{\alpha,\beta}$ and
therefore a trivial fibration by the inductive hypothesis. Consequently, it will suffice to prove
that $\phi_{\alpha+1,\alpha}$ is a trivial fibration. The result is obvious if
$\calC^{\alpha+1} = \calC^{\alpha}$, so we may assume without loss of generality that 
$\calC^{\alpha+1}$ is the smallest complete simplicial subset of $\calC$ containing
$\calC^{\alpha}$ together with a simplex $\sigma: \Delta^n \rightarrow \calC$,
where $\sigma$ does not belong to $\calC^{\alpha}$. Since $n$ is chosen to be minimal, we may suppose that $\sigma$ is nondegenerate, and that the boundary of $\sigma$ already belongs to
$\calC^{\alpha}$.

Form
a pushout diagram
$$ \xymatrix{ \calC^0_{/ \sigma} \star \bd \Delta^n \ar[r] \ar[d] & \calC^{\alpha} \ar[d] \\
\calC^{0}_{/\sigma} \star \Delta^n \ar[r] & \calC'. }$$
By construction there is an induced map $\calC' \rightarrow \calC$, which is easily shown to be a monomorphism of simplicial sets; we may therefore identify $\calE'$ with its image in $\calC$. Since $\calC$ is minimal, we can apply Proposition \ref{minstrict} to deduce that $\calC'$ is complete, so that $\calC' = \calC^{\alpha+1}$. Let $G$ denote the composition
$$ \calC^0_{/\sigma} \star \Delta^n \rightarrow \calC \stackrel{F}{\rightarrow} \calD $$
and $G_{\bd} = G | \calC^{0}_{/\sigma} \star \bd \Delta^n$.
It follows that $\phi_{\alpha+1,\alpha}$ is a pullback of the induced map
$$ \psi: \calD_{G/} \rightarrow \calD'_{p G/} \times_{ \calD'_{p G_{\bd}/} }
\calD_{ G_{\bd}/}.$$
To complete the proof, it will suffice to show that $\psi$ is a trivial fibration of simplicial sets.

Let $G_0 = G | \calC^0_{/\sigma}$. Let $\calE = \calD_{G_0/}$, 
$\calE' = \calD'_{p \circ G_0/}$, and let $q: \calE \rightarrow \calE'$ be the induced map.
We can identify $G$ with a map $\sigma': \Delta^n \rightarrow \calE$. Let
$\sigma'_0 = \sigma' | \bd \Delta^n$. Then we wish to prove that
the map
$$ \psi': \calE_{ \sigma'/ } \rightarrow \calE'_{ q \sigma' / }
\times_{ \calE'_{ q \sigma'_0 / } } \calE_{ q \sigma'_0 /}$$
is a trivial fibration. Let $C = \sigma(0)$. 

The projection $\calC^0_{/\sigma} \rightarrow \calC^{0}_{/C}$ is a trivial fibration of simplicial sets, and therefore cofinal. Since $F$ is a $p$-left Kan extension of $F_0$ at $C$, we conclude that $\sigma'(0)$ is a $q$-initial object of $\calE$.

To prove that $\psi$ is a trivial fibration, it will suffice to prove that $\psi$ has the right lifting property with respect to the inclusion $\bd \Delta^m \subseteq \Delta^m$, for each $m \geq 0$. Unwinding the definitions, this amounts to the existence of a dotted arrow as indicated in the diagram
$$ \xymatrix{ \bd \Delta^{n+m+1} \ar[r]^{s} \ar@{^{(}->}[d] & \calE \ar[d]^{q} \\
\Delta^{n+m+1}. \ar@{-->}[ur] \ar[r] & \calE' }$$
However, the map $s$ carries the initial vertex of $\Delta^{n+m+1}$ to a vertex of $\calE$ which
is $q$-initial, so that the desired extension can be found.
\end{proof}

\begin{proposition}\label{acekan}
Let $F: \calC \rightarrow \calD$ be a functor between $\infty$-categories, 
$p: \calD \rightarrow \calD'$ an inner fibration of $\infty$-categories, and
$\calC^{0} \subseteq \calC^{1} \subseteq \calC$ full subcategories. Suppose that
$F| \calC^{1}$ is a $p$-left Kan extension of $F| \calC^{0}$. Then $F$ is a $p$-left Kan extension of $F|\calC^{1}$ if and only if $F$ is a $p$-left Kan extension of $F| \calC^{0}$.
\end{proposition}

\begin{proof}
Let $C$ be an object of $\calC$; we will show that $F$ is a $p$-left Kan extension of $F|\calC^{0}$ at $C$ if and only if $F$ is a $p$-left Kan extension of $F| \calC^{1}$ at $C$. Consider the composition
$$ F^0_{C}: (\calC^{0}_{/C})^{\triangleright} \subseteq
(\calC^{1}_{/C})^{\triangleright} \stackrel{F^1_{C}}{\rightarrow} \calD.$$
We wish to show that $F^0_{C}$ is a $p$-colimit diagram
if and only if $F^1_{C}$ is a $p$-colimit
diagram. According to Lemma \ref{kan0}, it will suffice to show that $F^1_{C} | \calC^{1}_{/C}$ is a left Kan extension of $F^0_{C}$. Let $f: C' \rightarrow C$ be an object of
$\calC^{1}_{/C}$. We wish to show that the composite map
$$ (\calC^{0}_{/f})^{\triangleright} \rightarrow (\calC^{0}_{/C'})^{\triangleright} \stackrel{F^{0}_{C'}}{\rightarrow} \calD$$ is a $p$-colimit diagram. Since the projection
$\calC^{0}_{/f} \rightarrow \calC^{0}_{/C'}$ is cofinal (in fact, a trivial fibration), it will suffice to show that $F^0_{C'}$ is a $p$-colimit diagram (Proposition \ref{relexists}). This follows from our hypothesis that $F| \calC^{1}$ is a $p$-left Kan extension of $F| \calC^{0}$.
\end{proof}

\begin{proposition}\label{stormus}
Let $F: \calC \times \calC' \rightarrow \calD$ be a functor between $\infty$-categories, $p:
\calD \rightarrow \calD'$ an inner fibration of $\infty$-categories, 
and $\calC^{0} \subseteq \calC$ a full subcategory. The following conditions are
equivalent:
\begin{itemize}
\item[$(1)$] The functor $F$ is a $p$-left Kan extension of $F| \calC^{0} \times \calC'$.
\item[$(2)$] For each object $C' \in \calC'$, the induced functor
$F_{C'}: \calC \times \{C'\} \rightarrow \calD$ is a $p$-left Kan extension of
$F_{C'}| \calC^{0} \times \{C'\}$.
\end{itemize}
\end{proposition}

\begin{proof}
It suffices to show that $F$ is a $p$-left Kan extension of $F|\calC^{0} \times \calC'$
at an object $(C,C') \in \calC \times \calC'$ if and only if $F_{C'}$ is a $p$-left Kan extension
of $F_{D} | \calC^{0} \times \{D\}$ at $C$. This follows from the observation that the inclusion
$\calC^{0}_{/C} \times \{ \id_{C'} \} \subseteq \calC^{0}_{/C} \times \calC'_{/C'}$ is cofinal
(because $\id_{C'}$ is a final object of $\calC'_{/C'}$). 
\end{proof}

\begin{lemma}\label{kanhalf}
Let $m \geq 0$, $n \geq 1$ be integers, and let
$$ \xymatrix{ (\bd \Delta^m \times \Delta^n ) \coprod_{ \bd \Delta^m \times \bd \Delta^n}
( \Delta^m \times \bd \Delta^n) \ar@{^{(}->}[d] \ar[rr]^-{f_0} & & X \ar[d]^{p} \\
\Delta^m \times \Delta^n \ar[rr] \ar@{-->}[urr]^{f} & & S }$$
be a diagram of simplicial sets, where $p$ is an inner fibration and
$f_0(0,0)$ is a $p$-initial vertex of $X$.
Then there exists a morphism $f: \Delta^m \times \Delta^n \rightarrow X$ rendering the diagram commutative. 
\end{lemma}

\begin{proof}
Choose a sequence of simplicial sets
$$( \bd \Delta^m \times \Delta^n ) \coprod_{ \bd \Delta^m \times \bd \Delta^n }
( \Delta^m \times \bd \Delta^n) = Y(0) \subseteq \ldots \subseteq Y(k) = \Delta^m \times \Delta^n,$$
where each $Y(i+1)$ is obtained from $Y(i)$ by adjoining a single nondegenerate simplex whose boundary already lies in $Y(i)$. We prove by induction on $i$ that $f_0$ can be
extended to a map $f_i$ such that the diagram 
$$ \xymatrix{ Y(i) \ar@{^{(}->}[d] \ar[r]^{f_i} & X \ar[d]^{p} \\
\Delta^m \times \Delta^n \ar[r] & S }$$
is commutative. Having done so, we can then complete the proof by choosing $i = k$.

If $i = 0$, there is nothing to prove. Let us therefore suppose that $f_i$ has been constructed, and consider the problem of constructing $f_{i+1}$ which extends $f_i$. This is equivalent to the lifting problem
$$ \xymatrix{ \bd \Delta^{r} \ar@{^{(}->}[d] \ar[r]^{\sigma_0} & X \ar[d]^{p} \\
\Delta^r \ar@{-->}[ur]^{\sigma} \ar[r] & S. }$$
It now suffices to observe that
where $r > 0$ and $\sigma_0(0) = f_0(0,0)$ is a $p$-initial vertex of $X$
(since every simplex of $\Delta^m \times \Delta^n$ which violates one of these conditions
already belongs to $Y(0)$ ).
\end{proof}

\begin{lemma}\label{sillytech}
Suppose given a diagram of simplicial sets
$$ \xymatrix{ X \ar[rr]^{p} \ar[dr] & & Y \ar[dl] \\
& S, & }$$
where $p$ is an inner fibration. Let $K$ be a simplicial set, let
$q_{S} \in \bHom_{S}(K \times S, X)$, and let $q'_S = p \circ q_{S}$.
Then the induced map
$$ X^{q_{S}/} \rightarrow Y^{q_{S}/}$$
is an inner fibration $($where the above simplicial sets are defined as in \S \ref{consweet}{}$)$.
\end{lemma}

\begin{proof}
Unwinding the definitions, we see that every lifting problem
$$ \xymatrix{ A \ar@{^{(}->}[d]^{i} \ar[r] & X^{q_{S}/} \ar[d] \\
B \ar[r] \ar@{-->}[ur] & Y^{q_{S}/} }$$
is equivalent to a lifting problem
$$ \xymatrix{ (A \times (K \diamond \Delta^0) ) \coprod_{ A \times K} (B \times K) \ar[r] \ar@{^{(}->}[d]^{i'} & X \ar[d]^{p} \\
B \times (K \diamond \Delta^0) \ar@{-->}[ur] \ar[r] & Y. }$$
We wish to show that this lifting problem has a solution, provided that $i$ is inner anodyne.
Since $p$ is an inner fibration, it will suffice to prove that $i'$ is inner anodyne, which follows from 
Corollary \ref{prodprod2}.
\end{proof}

\begin{lemma}\label{kan1}
Consider a diagram of $\infty$-categories
$$ \calC \rightarrow \calD' \stackrel{p}{\leftarrow} \calD,$$
where $p$ is an inner fibration. Let $\calC^{0} \subseteq \calC$ be a full subcategory.
Suppose given $n > 0$ and a commutative diagram
$$ \xymatrix{ \bd \Delta^n \ar@{^{(}->}[d] \ar[r]^-{f_0} & \bHom_{\calD'}(\calC, \calD) \ar[d] \\
\Delta^n \ar[r]^-{g} \ar@{-->}[ur]^-{f} & \bHom_{\calD'}(\calC^{0},\calD) }$$
with the property that the functor $F: \calC \rightarrow \calD$, determined by
evaluating $f_0$ at the vertex $\{0\} \subseteq \bd \Delta^n$, is a $p$-left Kan extension of
$F| \calC^{0}$. Then there exists a dotted arrow $f$ rendering the diagram commutative.
\end{lemma}

\begin{proof}
The proof uses the same strategy as that of Lemma \ref{kan0}. Using Lemma \ref{basekann} and Proposition \ref{princex}, we may replace $\calC$ by a minimal model and thereby assume that $\calC$ is minimal. As in the proof of Lemma \ref{kan0}, let us call a simplicial subset $\calE \subseteq \calC$ {\it complete} if it has the following property:
for any simplex $\sigma: \Delta^n \rightarrow \calC$, if $\sigma| \Delta^{ \{0, \ldots, i\} }$ factors through $\calC^{0}$ and $\sigma| \Delta^{ \{i+1, \ldots, n \} }$ factors through $\calE$, then $\sigma$
factors through $\calE$. Let $P$ denote the partially ordered set of pairs $(\calE, f_{\calE})$, where
$\calE \subseteq \calC$ is complete and $f_{\calE}$ is a map rendering commutative the diagram
$$ \xymatrix{ \bd \Delta^n \ar@{^{(}->}[d] \ar[r]^-{f_0} & \bHom_{\calD'}(\calC,\calD) \ar[d] \\
\Delta^n \ar@{=}[d] \ar[r]^-{f_{\calE}} & \bHom_{\calD'}(\calE, \calD) \ar[d] \\
\Delta^n \ar[r]^-{g} & \bHom_{\calD'}(\calC^{0},\calD). }$$
We partially order $P$ as follows: $(\calE, f_{\calE}) \leq (\calE', f_{\calE'} )$ if
$\calE \subseteq \calE'$ and $f_{\calE} = f_{\calE'} | \calE$. Using Zorn's lemma, we deduce that $P$ has a maximal element $(\calE, f_{\calE})$. If $\calE = \calC$, we may take $f = f_{\calE}$ and the proof is complete. Otherwise, choose a simplex $\sigma: \Delta^m \rightarrow \calC$ which does not belong to $\calE$, where $m$ is as small as possible. It follows that $\sigma$ is nondegenerate, and that the boundary of $\sigma$ belongs to $\calE$. Form a pushout diagram
$$ \xymatrix{ \calC^{0}_{/\sigma} \star \bd \Delta^m \ar[r] \ar@{^{(}->}[d] & \calE \ar[d] \\
\calC^{0}_{/\sigma} \star \Delta^m \ar[r] & \calE'. }$$
As in the proof of Lemma \ref{kan0}, we may identify $\calE'$ with a complete simplicial subset of $\calC$, which strictly contains $\calE$. Since $(\calE, f_{\calE})$ is maximal, we conclude that
$f_{\calE}$ does not extend to $\calE'$. Consequently, we deduce that there does not exist a dotted arrow rendering the diagram
$$ \xymatrix{ \calC^{0}_{/\sigma} \star \bd \Delta^m \ar@{^{(}->}[d] \ar[r] & \Fun(\Delta^n,\calD) \ar[d] \\
\calC^{0}_{/\sigma} \star \Delta^m \ar[r] \ar@{-->}[ur] & \Fun(\Delta^n, \calD')
\times_{ \Fun(\bd \Delta^n, \calD') } \Fun( \bd \Delta^n, \calD) }$$
commutative. Let 
$q: \calC^{0}_{/\sigma} \rightarrow \Fun( \Delta^n, \calD)$ be the restriction of the upper horizontal map,
and let $q': \calC^{0}_{/\sigma} \rightarrow \Fun(\Delta^n, \calD')$, 
$q_{\bd}: \calC^{0}_{/\sigma} \rightarrow \Fun(\bd \Delta^n, \calD)$,
$q'_{\bd}: \calC^{0}_{/\sigma} \rightarrow \Fun(\bd \Delta^n, \calD')$
be defined by composition with $q$. It follows that there exists no solution
to the associated lifting problem
$$ \xymatrix{ \bd \Delta^m \ar[r] \ar@{^{(}->}[d] & \Fun(\Delta^n,\calD)_{q/} \ar[d] \\
\Delta^m \ar[r] \ar@{-->}[ur] & \Fun(\Delta^n, \calD')_{q'/} \times_{
\Fun( \bd \Delta^n, \calD')_{q'_{\bd}/} } \Fun( \bd \Delta^n,\calD)_{q_{\bd}/}. }$$
Applying Proposition \ref{princex}, we deduce also the insolubility
of the equivalent lifting problem
$$ \xymatrix{ \bd \Delta^m \ar[r] \ar[d] & \Fun(\Delta^n, \calD)^{q/} \ar[d] \\
\Delta^m \ar[r] \ar@{-->}[ur] & \Fun(\Delta^n,\calD')^{q'/} \times_{
\Fun(\bd \Delta^n, \calD')^{q'_{\bd}/} } \Fun(\bd \Delta^n,\calD)^{q_{\bd}/}. }$$

Let $q_{\Delta^n}$ denote the map $\calC^0_{/\sigma} \times \Delta^n \rightarrow
\calD \times \Delta^n$ determined by $q$, and let
and let $\calX = (\calD \times \Delta^n)^{q_{\Delta^n}/}$ be the simplicial set constructed in \S \ref{consweet}. Let $q'_{\Delta^n}: \calC^{0}_{/\sigma} \times \Delta^n \rightarrow
\calD' \times \Delta^n$ and $\calX' = ( \calD' \times \Delta^n)^{q'_{\Delta^n}/}$
be defined similarly. We have natural isomorphisms
$$ \Fun(\Delta^n,\calD)^{q/} \simeq \bHom_{\Delta^n}(\Delta^n, \calX)$$
$$ \Fun(\bd \Delta^n, \calD)^{q_{\bd}/} \simeq \bHom_{\Delta^n}( \bd \Delta^n, \calX).$$
$$ \Fun(\Delta^n,\calD')^{q'/} \simeq \bHom_{\Delta^n}(\Delta^n, \calX')$$
$$ \Fun(\bd \Delta^n, \calD')^{q'_{\bd}/} \simeq \bHom_{\Delta^n}( \bd \Delta^n, \calX').$$
These identifications allow us reformulate our insoluble lifting problem once more:
$$ \xymatrix{ ( \bd \Delta^m \times \Delta^n ) \coprod_{ \bd \Delta^m \times \bd \Delta^n }
( \Delta^m \times \bd \Delta^n) \ar[rr]^-{g_0} \ar@{^{(}->}[d] & & \calX \ar[d]^{\psi} \\
\Delta^m \times \Delta^n \ar@{-->}[urr]^{g} \ar[rr] & & \calX'. }$$
We have a commutative diagram
$$ \xymatrix{ \calX \ar[rr]^{\psi} \ar[dr]^{r} & & \calX' \ar[dl]_{r'} \\
& \Delta^n. & }$$
Proposition \ref{colimfam} implies that $r$ and $r'$ are 
Cartesian fibrations, and that $\psi$ carries $r$-Cartesian edges
to $r'$-Cartesian edges. Lemma \ref{sillytech} implies that $\psi$ is an inner fibration.
Let $\psi_0: \calX_{ \{0\} } \rightarrow \calX'_{ \{0\} }$ be the diagram induced by taking
the fibers over the vertex $\{0\} \subseteq \Delta^n$. We have a commutative diagram
$$ \xymatrix{ \calD_{ \calC^{0}_{/\sigma(0)}/ } \ar[d]^{\theta} & \calD_{ \calC^{0}_{/\sigma} }\ar[l] \ar[d] \ar[r] &
\calX_{ \{0\} } \ar[d]^{\psi_0} \\
\calD'_{ \calC^{0}_{/ \sigma(0)}/} & \calD'_{ \calC^{0}_{/\sigma}/} \ar[r] \ar[l] & \calX'_{ \{0\} }}$$
in which the horizontal arrows are categorical equivalences. The vertex
$g_0(0,0) \in \calX'_{ \{0\} }$ lifts to a vertex of 
$\calD_{ \calC^{0}_{/\sigma}/}$ whose image in $\calD_{ \calC^{0}_{/\sigma(0)}/ }$
is $\theta$-initial (in virtue of our assumption that $F$ is a $p$-left Kan extension
of $F| \calC^{0}$). It follows that $g_0(0,0)$ is $\psi_0$-initial when regarded as a vertex of $\calX_{ \{0\} }$. 
Applying Proposition \ref{panna}, we deduce
that $g_0(0,0)$ is $\psi$-initial when regarded as a vertex of $\calX$. 
Lemma \ref{kanhalf} now guarantees the existence of the dotted arrow $g$, contradicting the maximality of $(\calE, f_{\calE})$.
\end{proof}

The following result addresses the existence problem for left Kan extensions:

\begin{lemma}\label{kan2}\label{Kan extension!existence of}
Suppose given a diagram of $\infty$-categories
$$ \xymatrix{ \calC^{0} \ar@{^{(}->}[d] \ar[r]^{F_0} & \calD \ar[d]^{p} \\
\calC \ar[r] \ar@{-->}[ur]^{F} & \calD' }$$
where $p$ is an inner fibration, and the left vertical arrow is the inclusion
of a full subcategory. The following conditions are equivalent:

\begin{itemize}
\item[$(1)$] There exists a functor $F: \calC \rightarrow \calD$ rendering the diagram
commutative, such that $F$ is a $p$-left Kan extension of $F_0$. 

\item[$(2)$] For every object $C \in \calC$, the diagram given by the composition
$$ \calC^{0}_{/C} \rightarrow \calC^{0} \stackrel{F_0}{\rightarrow} \calD$$ admits a $p$-colimit.
\end{itemize}

\end{lemma}

\begin{proof}
It is clear that $(1)$ implies $(2)$. Let us therefore suppose that $(2)$ is satisfied; we wish to prove that $F_0$ admits a left Kan extension. We will follow the basic strategy used in the proofs
of Lemmas \ref{kan0} and \ref{kan1}. 
Using Proposition \ref{princex} and Lemma \ref{basekann}, we can replace the inclusion $\calC^{0} \subseteq \calC$ by any categorically equivalent inclusion
${\calC'}^{0} \subseteq \calC'$. Using Proposition \ref{minimod}, we can choose
$\calC'$ to be a minimal model for $\calC$; we thereby reduce to the case where
$\calC$ is itself a minimal $\infty$-category.

We will say that a simplicial subset $\calE \subseteq \calC$ is {\it complete} if it has the following property:
for any simplex $\sigma: \Delta^n \rightarrow \calC$, if $\sigma| \Delta^{ \{0, \ldots, i\} }$ factors through $\calC^{0}$ and $\sigma| \Delta^{ \{i+1, \ldots, n \} }$ factors through $\calE$, then $\sigma$
factors through $\calE$. Note that if $\calE$ is complete, then $\calC^{0} \subseteq \calE$.
Let $P$ be the set of all pairs $(\calE, f_{\calE})$ where $\calE \subseteq \calC$ is complete, $f_{\calE}$ is a map of simplicial
sets which fits into a commutative diagram
$$ \xymatrix{ \calC^{0} \ar@{^{(}->}[d] \ar[r]^{F_0} & \calD \ar@{=}[d] \\
\calE \ar[r] \ar@{^{(}->}[d] \ar[r]^{f_{\calE}} & \calD \ar[d]^{p} \\
\calC \ar[r] & \calD', }$$
and every object $C \in \calE$, the composite map
$$ (\calC^{0}_{/C})^{\triangleright} \subseteq (\calE_{/C})^{\triangleright} \rightarrow \calE \stackrel{f_{\calE}}{\rightarrow} \calD$$
is a $p$-colimit diagram. We view $P$ as a partially ordered set, with
$(\calE, f_{\calE}) \leq (\calE', f_{\calE'})$ if $\calE \subseteq \calE'$ and
$f_{\calE'}|\calE = f_{\calE}$. This partially ordered set satisfies the hypotheses of Zorn's lemma,
and therefore has a maximal element which we will denote by $(\calE, f_{\calE})$. If
$\calE = \calC$, then $f_{\calE}$ is a $p$-left Kan extension of $F_0$ and the proof is complete.

Suppose that $\calE \neq \calC$. Then there is a simplex $\sigma: \Delta^n \rightarrow \calC$
which does not factor through $\calE$; choose such a simplex where $n$ is as small as possible.
The minimality of $n$ guarantees that $\sigma$ is nondegenerate, that $\sigma | \bd \Delta^n$ factors through $\calE$, and (if $n > 0$) that $\sigma(0) \notin \calC^{0}$. Form
a pushout diagram
$$ \xymatrix{ \calC^0_{/ \sigma} \star \bd \Delta^n \ar[r] \ar@{^{(}->}[d] & \calE \ar[d] \\
\calC^{0}_{/\sigma} \star \Delta^n \ar[r] & \calE'. }$$
This diagram induces a map $\calE' \rightarrow \calC$, which is easily shown to be a monomorphism of simplicial sets; we may therefore identify $\calE'$ with its image in $\calC$.
Since $\calC$ is minimal, we can apply Proposition \ref{minstrict} to deduce that $\calE' \subseteq \calC$ is complete. Since $(\calE, F_{\calE}) \in P$ is maximal, it follows that we cannot extend
$F_{\calE}$ to a functor $F_{\calE'}: \calE' \rightarrow \calD$ such that $(\calE', F_{\calE'}) \in P$.

Let $q$ denote the composition
$$ \calC^{0}_{/\sigma} \rightarrow \calC^{0} \stackrel{F_0}{\rightarrow} \calD.$$
The map $f_{\calE}$ determines a commutative diagram
$$ \xymatrix{ \bd \Delta^n \ar@{^{(}->}[d] \ar[r]^{g_0} & \calD_{q/} \ar[d]^{p'} \\
\Delta^n \ar[r] \ar@{-->}[ur]^{g} & \calD'_{p  q/}. }$$
Extending $f_{\calE}$ to a map $f_{\calE'}$ such that $(\calE', f_{\calE'}) \in P$ is
equivalent to producing a morphism $g: \Delta^n \rightarrow \calD_{q/}$ rendering the above
diagram commutative which, if $n=0$, is a $p$-colimit of $q$. In the case
$n=0$, the existence of such an extension follows from assumption $(2)$.
If $n > 0$, let $C = \sigma(0)$; then the projection $\calC^{0}_{/\sigma} \rightarrow
\calC^{0}_{/C}$ is a trivial fibration $\infty$-categories and $q$ factors as a composition
$$ \calC^{0}_{/\sigma} \rightarrow \calC^{0}_{/C} \stackrel{q'}{\rightarrow} \calD.$$ 
We obtain therefore a commutative diagram
$$ \xymatrix{ \calD_{q/} \ar[r]^{r} \ar[d]^{p'} & \calD_{q'/} \ar[d]^{p''} \\
\calD'_{p  q/} \ar[r] & \calD'_{p  q'/} }$$
where the horizontal arrows are categorical equivalences.
Since $(\calE, f_{\calE}) \in P$, $(r \circ g_0)(0)$ is a $p''$-initial vertex of
$\calD_{q'/}$. Applying Proposition \ref{summertoy}, we conclude that $g_0(0)$ is a $p'$-initial
vertex of $\calD_{q/}$, which guarantees the existence of the desired extension $g$.
This contradicts the maximality of $(\calE, f_{\calE})$ and completes the proof.
\end{proof}

\begin{corollary}\label{kanexistleft}
Let $p: \calD \rightarrow \calE$ be a coCartesian fibration of $\infty$-categories. Suppose
that each fiber of $p$ admits small colimits, and that for every morphism
$E \rightarrow E'$ in $\calE$, the associated functor $\calD_{E} \rightarrow \calD_{E'}$
preserves small colimits. Let $\calC$ be a small $\infty$-category, and 
$\calC^{0} \subseteq \calC$ a full subcategory. Then every functor
$F_0: \calC^{0} \rightarrow \calD$ admits a left Kan extension relative
to $p$.
\end{corollary}

\begin{proof}
This follows immediately from Lemma \ref{kan2} and Corollary \ref{constrel}.
\end{proof}

Combining Lemmas \ref{kan1} and \ref{kan2}, we deduce:

\begin{proposition}\label{lklk}
Suppose given a diagram of $\infty$-categories
$$ \calC \rightarrow \calD' \stackrel{p}{\leftarrow} \calD,$$
where $p$ is an inner fibration. Let $\calC^{0}$ be a full subcategory of $\calC$.
Let $\calK \subseteq \bHom_{\calD'}(\calC, \calD)$ be the full subcategory spanned by those functors $F: \calC \rightarrow \calD$ which are $p$-left Kan extensions of $F | \calC^{0}$.
Let $\calK' \subseteq \bHom_{\calD'}(\calC^{0}, \calD)$ be the full subcategory spanned
by those functors
$F_0: \calC^{0} \rightarrow \calD$ with the property that, for each object
$C \in \calC$, the induced diagram
$ \calC^{0}_{/C} \rightarrow \calD$ has a $p$-colimit. Then the restriction functor
$\calK \rightarrow \calK'$ is a trivial fibration of simplicial sets.
\end{proposition}

\begin{corollary}\label{leftkanextdef}
Suppose given a diagram of $\infty$-categories
$$ \calC \rightarrow \calD' \stackrel{p}{\leftarrow} \calD,$$
where $p$ is an inner fibration. Let $\calC^{0}$ be a full subcategory of $\calC$.
Suppose further that, for every functor $F_0 \in \bHom_{\calD'}(\calC^{0}, \calD)$,
there exists a functor $F \in \bHom_{\calD'}(\calC,\calD)$ which is a $p$-left Kan extension
of $F_0$.
Then the restriction map $i^{\ast}: \bHom_{\calD'}(\calC,\calD) \rightarrow 
\bHom_{\calD'}(\calC^{0},\calD)$ admits a section
$i_{!}$, whose essential image consists of precisely of those functors $F$
which are $p$-left Kan extensions of $F | \calC^{0}$.
\end{corollary}

In the situation of Corollary \ref{leftkanextdef}, we will refer to $i_{!}$ as a {\it left Kan extension functor}. We note that Proposition \ref{lklk} proves not only the existence of $i_{!}$, but also its uniqueness up to homotopy (the collection of all such functors is parametrized by a contractible Kan complex). The following characterization of $i_{!}$ gives a second explanation for its uniqueness:

\begin{proposition}\label{leftkanadj}
Suppose given a diagram of $\infty$-categories
$$ \calC \rightarrow \calD' \stackrel{p}{\leftarrow} \calD,$$
where $p$ is an inner fibration. Let $i: \calC^{0} \subseteq \calC$ be the inclusion of a full subcategory,
and suppose that every functor $F_0 \in \bHom_{\calD'}(\calC^{0},\calD)$ admits a $p$-left Kan extension. Then the left Kan extension functor
$i_{!}: \bHom_{ \calD'}( \calC^{0}, \calD) \rightarrow \bHom_{\calD'}(\calC, \calD)$
is a left adjoint to the restriction functor $i^{\ast}: \bHom_{\calD'}(\calC, \calD) \rightarrow
\bHom_{\calD'}(\calC^{0}, \calD)$.
\end{proposition}

\begin{proof}
Since $i^{\ast} \circ i_{!}$ is the identity functor on $\bHom_{\calD'}(\calC^{0}, \calD)$, there is an obvious candidate for the unit
$$ u: \id \rightarrow i^{\ast} \circ i_{!}$$
of the adjunction: namely, the identity. According to Proposition \ref{storut}, it will suffice to prove
that for every $F \in \bHom_{\calD'}(\calC^{0}, \calD)$,
$G \in \bHom_{\calD'}(\calC, \calD)$, composition with $u$ induces a homotopy equivalence
$$ \bHom_{\bHom_{\calD'}(\calC,\calD)
}(i_{!} F, G) \rightarrow \bHom_{\bHom_{\calD'}(\calC^{0},\calD)}( i^{\ast} i_{!} F,
i^{\ast} G) \stackrel{u}{\rightarrow}\bHom_{\bHom_{\calD'}(\calC^{0},\calD)}( F,
i^{\ast} G)$$
in the homotopy category $\calH$. This morphism in $\calH$ is represented by the restriction map 
$$ \Hom^{\rght}_{ \bHom_{\calD'}(\calC,\calD)}( i_{!}F,G) \rightarrow \Hom^{\rght}_{\bHom_{\calD'}(\calC^{0},\calD)}(F, i^{\ast} G)$$
which is a trivial fibration by Lemma \ref{kan1}.
\end{proof}

\begin{remark}
Throughout this section we have focused our attention on the theory of (relative) {\em left} Kan extensions. There is an entirely dual theory of {\em right} Kan extensions in the $\infty$-categorical setting, which can be obtained from the theory of left Kan extensions by passing to opposite $\infty$-categories.
\end{remark}

\subsection{Kan Extensions along General Functors}\label{bigkanext}

Our goal in this section is to generalize the theory of Kan extensions to the case where the
change of diagram category is not necessarily given by a fully faithful inclusion 
$\calC^{0} \subseteq \calC$. As in \S \ref{kanex}, we will discuss only the theory of {\em left} Kan extensions; a dual theory of right Kan extensions can be obtained by passing to opposite $\infty$-categories.

The ideas introduced in this section are relatively elementary extensions of the ideas of \S \ref{kanex}. However, we will encounter a new complication.
Let $\delta: \calC \rightarrow \calC'$ be a change of diagram $\infty$-category,
$f: \calC \rightarrow \calD$ a functor, and $\delta_{!}(f): \calC' \rightarrow \calD$
its left Kan extension along $\delta$ (to be defined below). Then one does not generally
expect that $\delta^{\ast} \delta_{!}(f)$ to be equivalent to the original functor $f$. Instead, one
has only a unit transformation $f \rightarrow \delta^{\ast} \delta_{!}(f)$. To set up the theory, this unit transformation must be taken as part of the data. Consequently, the theory of Kan extensions in general requires more elaborate notation and terminology than the special case treated in \S \ref{kanex}. We will compensate for this by considering only the case of {\em absolute} left Kan extensions. It is straightforward to set up a relative theory as in \S \ref{kanex}, but we will not need such a theory in this book.

\begin{definition}\index{gen}{left extension}
Let $\delta: K \rightarrow K'$ be a map of simplicial sets, let $\calD$ be an $\infty$-category, and let
$f: K \rightarrow \calD$ be a diagram. A {\it left extension of $f$ along $\delta$} consists of a map
$f': K' \rightarrow \calD$ and a morphism $f \rightarrow f' \circ \delta$ in the $\infty$-category
$\Fun(K,\calD)$. 
\end{definition}

Equivalently, we may view a left extension of $f: K \rightarrow \calD$ along $\delta: K \rightarrow K'$
as a map $F: M^{op}(\delta) \rightarrow \calD$ such that $F|K=f$, where $M^{op}(\delta) = M(\delta^{op})^{op} = (K \times \Delta^1) \coprod_{ K \times \{1\} } K'$ denotes the mapping cylinder of $\delta$.

\begin{definition}\label{genkan}\index{gen}{Kan extension}
Let $\delta: K \rightarrow K'$ be a map of simplicial sets, and let
$F: M^{op}(\delta) \rightarrow \calD$ be a diagram in an $\infty$-category $\calD$
(which we view as a left extension of $f = F|K$ along $\delta$). We will say that $F$ is a {\it left
Kan extension of $f$ along $\delta$} if there exists a commutative diagram
$$ \xymatrix{ M^{op}(\delta) \ar[r]^{F''} \ar[dr] & \calK \ar[r]^{F'} \ar[d]^{p} & \calD \\
& \Delta^1 & }$$
where $F''$ is a categorical equivalence, $\calK$ is an $\infty$-category, 
$F= F' \circ F''$, and $F'$ is a left Kan extension of $F' | \calK \times_{ \Delta^1} \{0\}$. 
\end{definition}

\begin{remark}\label{genka}
In the situation of Definition \ref{genkan}, the map $p: \calK \rightarrow \Delta^1$ is {\em automatically} a coCartesian fibration. To prove this, choose a factorization
$$ M(\delta^{op})^{\natural} \stackrel{i}{\rightarrow} (\calK')^{\sharp} \rightarrow (\Delta^1)^{\sharp}$$
where $i$ is marked anodyne, and $\calK' \rightarrow \Delta^1$ is a Cartesian fibration. Then $i$ is a quasi-equivalence, so that Proposition \ref{simplexplay} implies that $M(\delta^{op}) \rightarrow \calK'$ is a categorical equivalence. It follows that $\calK$ is equivalent to $(\calK')^{op}$ (via an equivalence which respects the projection to $\Delta^1$), so that the projection
$p$ is a coCartesian fibration.
\end{remark}

The following result asserts that the condition of Definition \ref{genkan} is essentially independent of the choice of $\calK$.

\begin{proposition}\label{oave}
Let $\delta: K \rightarrow K'$ be a map of simplicial sets, and let
$F: M^{op}(\delta) \rightarrow \calD$ be a diagram in an $\infty$-category $\calD$ which
is a left Kan extension along $\delta$. Let
$$ \xymatrix{ M^{op} \ar[r]^{F''} \ar[dr] & \calK \ar[d]^{p} \\
& \Delta^1 }$$
be a diagram where $F''$ is both a cofibration and a categorical equivalence of simplicial sets. Then $F= F' \circ F''$, for some map $F': \calK \rightarrow \calD$
which is a left Kan extension of $F' | \calK \times_{\Delta^1} \{0\}$.
\end{proposition}

\begin{proof}
By hypothesis, there exists a commutative diagram
$$ \xymatrix{ M^{op}(\delta) \ar[r]^{G''} \ar[d]^{F''} & \calK' \ar[r]^{G'} \ar[d]^{q} & \calD \\
\calK \ar[r]^{p} \ar@{-->}[ur]^{r} & \Delta^1 & }$$
where $\calK'$ is an $\infty$-category, $F = G' \circ G''$, and $G''$ is a categorical equivalence, and
$G'$ is a left Kan extension of $G' | \calK' \times_{\Delta^1} \{0\}$. Since $\calK'$ is an
$\infty$-category, there exists a map $r$ as indicated in the diagram such that
$G'' = r \circ F''$. We note that $r$ is a categorical equivalence so that the commutativity
of the lower triangle $p = q \circ r$ follows automatically. We now define
$F' = G' \circ r$, and note that part $(2)$ of Lemma \ref{basekann} implies that $F'$ 
is a left Kan extension of $F'| \calK \times_{\Delta^1} \{0\}$.
\end{proof}

We have now introduced two different definitions of left Kan extensions: Definition \ref{defKan}, which applies in the situation of an inclusion $\calC^{0} \subseteq \calC$ of a full subcategory into an $\infty$-category $\calC$, and Definition \ref{genkan} which applies in the case of a general map $\delta: K \rightarrow K'$ of simplicial sets. These two definitions are essentially the same. More precisely, we have the following assertion:

\begin{proposition}\label{compkan}
Let $\calC$ and $\calD$ be $\infty$-categories, and let $\delta: \calC^{0} \rightarrow \calC$ denote the inclusion of a full subcategory.

\begin{itemize}
\item[$(1)$] Let $f: \calC \rightarrow \calD$ be a functor, $f_0$ its restriction to $\calC^{0}$, so that $( f, \id_{f_0})$ can be viewed as a left extension of $f_0$ along $\delta$.
Then $(f, \id_{f_0})$ is a left Kan extension of $f_0$ along $\delta$ if and only if $f$ is a left Kan extension of $f_0$.  

\item[$(2)$] A functor $f_0: \calC^{0} \rightarrow \calD$ has a left Kan extension
if and only if it has a left Kan extension along $\delta$.
\end{itemize}
\end{proposition}

\begin{proof}
Let $\calK$ denote the full subcategory of $\calC \times \Delta^{1}$ spanned by the objects
$(C, \{i\})$ where either $C \in \calC^{0}$ or $i=1$, so that we have inclusions
$$ M^{op}(\delta) \subseteq \calK \subseteq \calC \times \Delta^{1}.$$
To prove $(1)$, suppose that $f: \calC \rightarrow \calD$ is a left Kan extension of 
$f_0 = f | \calC^{0}$ and let $F$ denote the composite map
$$ \calK \subseteq \calC \times \Delta^{1} \rightarrow \calC \stackrel{f}{\rightarrow} \calD.$$
It follows immediately that $F$ is a left Kan extension of $F | \calC^{0} \times \{0\}$, so that
$F| M^{op}(\delta)$ is a left Kan extension of $f_0$ along $\delta$.

To prove $(2)$, we observe that the ``only if'' follows from $(1)$; the converse follows from the existence criterion of Lemma \ref{kan2}.
\end{proof}

Suppose that $\delta: K^{0} \rightarrow K^{1}$ is a map of simplicial sets, $\calD$ an $\infty$-category, and that every diagram $K^{0} \rightarrow \calD$ admits a left Kan extension along $\delta$. Choose a diagram
$$\xymatrix{ M^{op}(\delta) \ar[rr]^{j} \ar[dr] & & \calK \ar[dl] \\
& \Delta^1 & }$$
where $j$ is inner anodyne and $\calK$ is an $\infty$-category, which we regard
as a correspondence from $\calK^{0} = \calK \times_{\Delta^1} \{0\}$ to 
$\calK^{1} = \calK \times_{ \Delta^1} \{1\}$. Let $\calC$ denote the full subcategory
of $\Fun(\calK, \calD)$ spanned by those functors $F: \calK \rightarrow \calD$ such that
$F$ is a left Kan extension of $F_0 = F|\calK^{0}$. The restriction map
$p: \calC \rightarrow \Fun(K^0, \calD)$ can be written as a composition of
$\calC \rightarrow \calD^{\calK^0}$ (a trivial fibration by Proposition \ref{lklk}) and
$\Fun(\calK^{0},\calD) \rightarrow \Fun(K^0,\calD)$ (a trivial fibration since $K^0 \rightarrow \calK^0$ is inner anodyne), and is therefore a trivial fibration. Let $\overline{\delta}_{!}$ be the composition of
a section of $p$ with the restriction map $\calC \subseteq \Fun(\calK, \calD) \rightarrow \Fun( M^{op}(\delta), \calD)$, and let $\delta_{!}$ denote the composition of 
$\overline{\delta}_{!}$ with the restriction map $\Fun( M^{op}(\delta),\calD) \rightarrow 
\Fun(K^1,\calD)$. Then $\overline{\delta}_{!}$ and $\delta_{!}$ are well-defined up to equivalence, at least once $\calK$ has been fixed (independence of the choice of $\calK$ will follow from the characterization given in Proposition \ref{charleftkangen}). We will abuse terminology by referring to {\em both} $\overline{\delta}_{!}$ and $\delta_{!}$ as {\it left Kan extension along $\delta$} (it should be clear from context which of these functors is meant in a given situation). We observe that $\overline{\delta}_{!}$ assigns to each object $f_0: K^0 \rightarrow \calD$ a 
left Kan extension of $f_0$ along $\delta$. 

\begin{example}
Let $\calC$ and $\calD$ be $\infty$-categories, and let
$i: \calC^{0} \rightarrow \calC$ be the inclusion of a full subcategory. Suppose that
$i_{!}: \Fun(\calC^{0}, \calD) \rightarrow \Fun(\calC, \calD)$ is a section of $i^{\ast}$, which satisfies the conclusion of Corollary \ref{leftkanextdef}. Then $i_{!}$ is a left Kan extension along $i$ in the sense defined above; this follows easily from Proposition \ref{compkan}. 
\end{example}

Left Kan extension functors admit the following characterization:

\begin{proposition}\label{charleftkangen}
Let $\delta: K^0 \rightarrow K^1$ be a map of simplicial sets, let $\calD$ be an $\infty$-category, let $\delta^{\ast}: \Fun(K^1,\calD) \rightarrow \Fun(K^0,\calD)$ be the restriction functor, and let
$\delta_{!}: \Fun(K^0,\calD) \rightarrow \Fun(K^1,\calD)$ be a functor of left Kan extension along $\delta$. Then $\delta_{!}$ is a left adjoint of $\delta^{\ast}$.
\end{proposition}

\begin{proof}
The map $\delta$ can be factored as a composition
$$ K^0 \stackrel{i}{\rightarrow} M^{op}(\delta) \stackrel{r}{\rightarrow} K^1$$ where
$r$ denotes the natural retraction of $M^{op}(\delta)$ onto $K^1$. Consequently, $\delta^{\ast} = i^{\ast} \circ r^{\ast}$. Proposition \ref{leftkanadj} implies that the left Kan extension functor $\overline{\delta}_{!}$ is
a left adjoint to $i^{\ast}$. By Proposition \ref{compadjoint}, it will suffice to prove that $r^{\ast}$ is a right adjoint to the restriction functor $j^{\ast}: \Fun(M^{op}(\delta),\calD) \rightarrow \Fun(K^1,\calD)$. 
Using Corollary \ref{tweezegork}, we deduce that $j^{\ast}$ is a coCartesian fibration. Moreover, there is a simplicial homotopy 
$\Fun(M^{op}(\delta), \calD) \times \Delta^1 \rightarrow \Fun(M^{op}(\delta), \calD)$ from the identity to $r^{\ast} \circ j^{\ast}$, which is a fiberwise homotopy over $\Fun(K^1, \calD)$. 
It follows that for every object $F$ of $\Fun(K^1,\calD)$, $r^{\ast} F$ is a final object of the $\infty$-category $\Fun(M^{op}(\delta), \calD) \times_{ \Fun(K^1,\calD)} \{F\}$. Applying Proposition \ref{quuquu}, we deduce that $r^{\ast}$ is right adjoint to $j^{\ast}$ as desired.
\end{proof}

Let $\delta: K^{0} \rightarrow K^{1}$ be a map of simplicial sets and $\calD$ an $\infty$-category which which that left Kan extension $\delta_{!}: \Fun(K^0,\calD) \rightarrow \delta_{!} \Fun(K^1,\calD)$ is defined. In general, the terminology ``Kan extension'' is perhaps somewhat unfortunate: if
$F: K^0 \rightarrow \calD$ is a diagram, then $\delta^{\ast} \delta_{!} F$ need not coincide with $F$, even up to equivalence. If $\delta$ is fully faithful, then the unit map
$F \rightarrow \delta^{\ast} \delta_{!} F$ is an equivalence: this follows from
Proposition \ref{compkan}. We will later need the following more precise assertion:

\begin{proposition}\label{timeless}
Let $\delta: \calC^{0} \rightarrow \calC^{1}$ and $f_0: \calC^{0} \rightarrow \calD$ be functors between $\infty$-categories, and let $f_{1}: \calC^{1} \rightarrow \calD$, $\alpha: f_{0} \rightarrow \delta^{\ast} f_1 = f_1 \circ \delta$ be a left Kan extension of $f_0$ along $\delta$.
Let $C$ be an object of $\calC^{0}$ such that, for every $C' \in \calC^{0}$, the functor
$\delta$ induces an isomorphism
$$ \bHom_{\calC^{0}}(C',C) \rightarrow \bHom_{\calC^{1}}( \delta C', \delta C)$$
in the homotopy category $\calH$. Then the morphism
$\alpha(C): f_0(C) \rightarrow f_1(\delta C)$ is an equivalence in $\calD$.
\end{proposition}

\begin{proof}
Choose a diagram
$$ \xymatrix{ M^{op}(\delta) \ar[r]^{G} \ar[dr] & \calM \ar[d] \ar[r]^{F} & \calD \\
& \Delta^1 & }$$
where $\calM$ is a correspondence from $\calC^{0}$ to $\calC^{1}$ associated to $\delta$, $F$ is a left Kan extension of $f_0 = F| \calC^{0}$, and $F \circ G$ is the map
$M^{op}(\delta) \rightarrow \calD$ determined by $f_0$, $f_1$, and $\alpha$.
Let $u: C \rightarrow \delta C$ be the morphism in $\calM$ given by the image
of $\{C\} \times \Delta^1 \subseteq M^{op}(\delta)$ under $G$. Then
$\alpha(C) = F(u)$, so it will suffice to prove that $F(u)$ is an equivalence. Since
$F$ is a left Kan extension of $f_0$ at $\delta C$, the composition
$$ (\calC^{0}_{/\delta C})^{\triangleright} \rightarrow
\calM \stackrel{F}{\rightarrow} \calD$$
is a colimit diagram. Consequently, it will suffice to prove that
$u: C \rightarrow \delta C$ is a final object of
$\calC^{0}_{/ \delta C}$. Consider the diagram
$$ \calC^{0}_{/C} \leftarrow \calC^{0}_{/u} \stackrel{q}{\rightarrow} \calC^{0}_{/\delta C}.$$
The $\infty$-category on the left has a final object $\id_{C}$, and the map on the left
is a trivial fibration of simplicial sets. We deduce that $s^0 u$ is a final object of\
$\calC^{0}_{/u}$. Since $q( s^0 u) = u \in \calC^{0}_{/\delta C}$, it will suffice to show that
$q$ is an equivalence of $\infty$-categories. We observe that $q$ is a map of right fibrations over $\calC^{0}$. According to Proposition \ref{apple1}, it will suffice to show that for each object $C'$ in $\calC^{0}$, the map $q$ induces a homotopy equivalence of Kan complexes
$$ \calC^{0}_{/u} \times_{ \calC^{0}} \{C'\} \rightarrow \calC^{0}_{/ \delta C} \times_{\calC^{0}} \{C'\}.$$
This map can be identified with the map
$$ \bHom_{\calC^0}( C',C) \rightarrow \bHom_{\calM}(C', \delta C)
\simeq \bHom_{\calC^1}(\delta C', \delta C),$$
in the homotopy category $\calH$, and is therefore a homotopy equivalence by assumption.
\end{proof}

We conclude this section by proving that the construction of left Kan extensions behaves well in families.

\begin{lemma}\label{longerwait}
Suppose given a commutative diagram
$$ \xymatrix{ \calC^{0} \ar[r]^{q} \ar[dr]^{i} & \calC \ar[d]^{p} \ar[r]^{F} & \calD \\
& \calE & } $$
of $\infty$-categories, where $p$ and $q$ are coCartesian fibrations, 
$i$ is the inclusion of a full subcategory, and $i$ carries $q$-coCartesian morphisms
of $\calC^{0}$ to $p$-coCartesian morphisms of $\calC$. The following conditions are equivalent:
\begin{itemize}
\item[$(1)$] The functor $F$ is a left Kan extension of $F| \calC^{0}$. 
\item[$(2)$] For each object $E \in \calE$, the induced functor
$F_{E}: \calC_{E} \rightarrow \calD$ is a left Kan extension of $F_E | \calC^{0}_{E}$. 
\end{itemize}
\end{lemma}

\begin{proof}
Let $C$ be an object of $\calC$ and let $E = p(C)$. Consider the composition
$$ (\calC^0_{E})^{\triangleright}_{/C} \stackrel{G^{\triangleright}}{\rightarrow} (\calC^{0}_{/C})^{\triangleright} \stackrel{F_{C}}{\rightarrow} \calD.$$
We will show that $F_C$ is a colimit diagram if and only if $F_C \circ G^{\triangleright}$ is a colimit diagram. For this, it suffices to show that the inclusion $G: (\calC^{0}_{E})_{/C} \subseteq
\calC^{0}_{/C}$ is cofinal. According to Proposition \ref{verylonger}, the projection $p': \calC_{/C} \rightarrow \calE_{/E}$
is a coCartesian fibration, and a morphism
$$ \xymatrix{ C' \ar[rr]^{f} \ar[dr] & & C'' \ar[dl] \\
& C & }$$
in $\calC_{/C}$ is $p'$-coCartesian if and only if $f$ is $p$-coCartesian. It follows that
$p'$ restricts to a coCartesian fibration $\calC'_{/C} \rightarrow \calE_{/E}$. 
We have a pullback diagram of simplicial sets
$$ \xymatrix{ (\calC^{0}_{E})_{/C} \ar[r]^{G} \ar[d] & \calC^0_{/C} \ar[d] \\
\{ \id_{E} \} \ar[r]^{G_0} & \calE_{/E}. }$$
The right vertical map is smooth (Proposition \ref{strokhop}) and $G_0$ is right anodyne, so that $G$ is right anodyne as desired.
\end{proof}

\begin{proposition}\label{longwait2}\label{Kan extension!in families}
Let
$$ \xymatrix{ X \ar[dr]^{p} \ar[rr]^{\delta} & & Y \ar[dl]^{q} \\
& S & }$$
be a commutative diagram of simplicial sets, where $p$ and $q$ are coCartesian fibrations,
and $\delta$ carries $p$-coCartesian edges to $q$-coCartesian edges.
Let $f_0: X \rightarrow \calC$ be a diagram in an $\infty$-category $\calC$, and let
$f_1: Y \rightarrow \calC$, $\alpha: f_0 \rightarrow f_1 \circ \delta$ be a left extension of 
$f_0$. The following conditions are equivalent:
\begin{itemize}
\item[$(1)$] The transformation $\alpha$ exhibits $f_1$ as a left Kan extension of $f_0$ along
$\delta$.
\item[$(2)$] For each vertex $s \in S$, the restriction $\alpha_{s}: f_0|X_{s} \rightarrow
(f_1 \circ \delta)|X_s$ exhibits $f_1|Y_{s}$ as a left Kan extension of $f_0|X_{s}$ along
$\delta_{s}: X_{s} \rightarrow Y_{s}$.
\end{itemize}
\end{proposition}

\begin{proof}
Choose an equivalence of simplicial categories $\sCoNerve(S) \rightarrow \calE$, where $\calE$ is fibrant, and let $[1]$ denote the linearly ordered set $\{0,1\}$, regarded as a category. Let
$\phi'$ denote the induced map $\sCoNerve( S \times \Delta^1) \rightarrow \calE \times [1]$.
Let $M$ denote the marked simplicial set
$$((X^{op})^{\natural} \times (\Delta^1)^{\sharp}) \coprod_{ (X^{op})^{\natural} \times \{0\} }
( Y^{op})^{\natural}.$$
Let $\St_{\phi}^{+}: (\mSet)_{(S \times \Delta^1)^{op}} \rightarrow (\mSet)^{\calE \times [1]}$ denote the straightening functor defined in \S \ref{markmodel2}, and choose a fibrant replacement
$$ \St_{\phi}^{+} M \rightarrow Z$$
in $(\mSet)^{\calE \times [1]}$. 
Let $S' = \sNerve(\calE)$, so that $S' \times \Delta^1 \simeq
\sNerve(\calE \times [1])$, and let $\psi: \sCoNerve(S' \times \Delta^1) \rightarrow \calE \times [1]$
be the counit map. Then $$\Un^{+}_{\psi}(Z)$$
is a fibrant object of $(\mSet)_{/ (S' \times \Delta^1)^{op}}$, which we may identify with a coCartesian fibration of simplicial sets $\calM \rightarrow S' \times \Delta^1$. 

We may regard $\calM$ as a correspondence from $\calM^{0} = \calM \times_{\Delta^1} \{0\}$ to
$\calM^{1} = \calM \times_{\Delta^1} \{1\}$. By construction, we have a unit map
$$ u: M^{op}(\delta) \rightarrow \calM \times_{S'} S.$$
Theorem \ref{straightthm} implies that the induced maps
$u_0: X \rightarrow \calM^{0} \times_{S'} S$, $u_1: Y \rightarrow \calM^{1} \times_{S'} S$ are equivalences of coCartesian fibrations. Proposition \ref{basechangefunky} implies that the maps
$\calM^{0} \times_{S'} S \rightarrow \calM^{0}$, $\calM^{1} \times_{S'} S \rightarrow \calM^{1}$ are categorical equivalences. 

Let $u'$ denote the composition
$$ M^{op}(\delta) \stackrel{u}{\rightarrow} \calM \times_{S'} S \rightarrow \calM,$$
and let $u'_0: X \rightarrow \calM^{0}$, $u'_1: Y \rightarrow \calM^{1}$ be defined similarly. The above argument shows that $u'_0$ and $u'_1$ are categorical equivalences. 
Consequently, the map $u'$ is a quasi-equivalence of coCartesian fibrations over $\Delta^1$, and therefore a categorical equivalence (Proposition \ref{simplexplay}). Replacing $\calM$ by the product
$\calM \times K$ if necessary, where $K$ is a contractible Kan complex, we may suppose that $u'$ is a cofibration of simplicial sets. Since $\calD$ is an $\infty$-category, there exists a functor $F: \calM \rightarrow \calD$ as indicated in the diagram below:
$$ \xymatrix{ M^{op}(\delta) \ar[rr]^{ (f_0, f_1, \alpha)} \ar[d] & & \calD \\
\calM. \ar@{-->}[urr]^{F} & & }$$
Consequently, we may reformulate condition $(1)$ as follows:
\begin{itemize}
\item[$(1')$] The functor $F$ is a left Kan extension of
$F | \calM^{0}$. 
\end{itemize}

Proposition \ref{apple1} now implies that, for each
vertex $s$ of $S$, the map $X_{s} \rightarrow \calM^{0}_{s}$ is a categorical equivalence.
Similarly, for each vertex $s$ of $S$, the inclusion $Y_{s} \rightarrow \calM^{1}_{s}$ is a categorical equivalence. It follows that the inclusion $M^{op}(\delta)_{s} \rightarrow \calM_{s}$ is a
quasi-equivalence, and therefore a categorical equivalence (Proposition \ref{simplexplay}). 
Consequently, we may reformulate condition $(2)$ as follows:

\begin{itemize}
\item[$(2')$] For each vertex $s \in S$, the functor $F| \calM_{s}$ is a left Kan extension of
$F| \calM^{0}_{s}$. 
\end{itemize}

Using Lemma \ref{basekann}, it is easy to see that the collection of objects $s \in S'$ such that
$F| \calM_{s}$ is a left Kan extension of $F| \calM^{0}_{s}$ is stable under equivalence. Since
the inclusion $S \subseteq S'$ is a categorical equivalence, we conclude that $(2')$ is equivalent to the following apparently stronger condition:

\begin{itemize}
\item[$(2'')$] For every object $s \in S'$, the functor $F| \calM_{s}$ is a left Kan extension of
$F| \calM^{0}_{s}$. 
\end{itemize}

The equivalence of $(1')$ and $(2'')$ follows from Lemma \ref{longerwait}. 
\end{proof}

\section{Examples of Colimits}\label{coexample}

\setcounter{theorem}{0}

In this section, we will analyze in detail the colimits of some very simple diagrams.
Our first three examples are familiar from classical category theory: coproducts (\S \ref{quasilimit5}), 
pushouts (\S \ref{quasilimit6}), and coequalizers (\S \ref{quasilimit8}). 

Our fourth example is slightly more unfamiliar. Let $\calC$ be an ordinary category which admits coproducts. Then $\calC$ is naturally {\em tensored} over the category of sets. Namely, for each
$C \in \calC$ and $S \in \Set$, we can define $C \otimes S$ to be the coproduct of a collection
of copies of $C$, indexed by the set $S$. The object $C \otimes S$ is characterized by the
following universal mapping property:
$$ \Hom_{\calC}(C \otimes S, D) \simeq \Hom_{\Set}( S, \Hom_{\calC}(C,D) ).$$
In the $\infty$-categorical setting, it is natural to try to generalize this definition by allowing
$S$ to an object of $\SSet$. In this case, $C \otimes S$ can again be viewed as a kind of colimit, but cannot be written as a {\em coproduct} unless $S$ is discrete. We will study the situation in \S \ref{quasilimit7}.

Our final objective in this section is to study the theory of {\em retracts} in an $\infty$-category $\calC$. In \S \ref{retrus}, we will see that there is a close relationship between retracts in $\calC$ and idempotent endomorphisms, just as in classical homotopy theory. Namely, any retract of an object $C \in \calC$ determines an idempotent endomorphism of $C$; conversely, if $\calC$ is {\it idempotent complete}, then every idempotent endomorphism of $C$ determines a retract of $C$.
We will return to this idea in \S \ref{surot}.

\subsection{Coproducts}\label{quasilimit5}

In this section, we discuss the simplest type of colimit: namely, coproducts.
Let $A$ be a set; we may regard $A$ as a category with
$$ \Hom_{A}(I,J) = \begin{cases} \ast & \text{if } I=J \\ \emptyset & \text{if } I \neq J. \end{cases}$$
We will also identify $A$ with the (constant) simplicial set which is the nerve of this category.
We note a functor $G: A \rightarrow \sSet$ is injectively fibrant if and only if it takes values in the category $\Kan$ of Kan complexes. If this condition is satisfied, then the product
$\prod_{\alpha \in A} G(\alpha)$ is a homotopy limit for $G$.

Let $F: A \rightarrow \calC$ be a functor from $A$ to a fibrant simplicial category; in other words, $F$ specifies a collection $\{ X_{\alpha} \}_{\alpha \in A}$ of objects in $\calC$. A homotopy colimit for $F$ will be referred to as a {\it homotopy coproduct} of the objects $\{ X_{\alpha} \}_{\alpha \in A}$. Unwinding the definition, we see that a homotopy coproduct is an object $X \in \calC$\index{gen}{coproduct!homotopy}\index{gen}{homotopy!coproduct} equipped with morphisms
$\phi_{\alpha}: X_{\alpha} \rightarrow X$ such that the induced map
$$ \bHom_{\calC}(X,Y) \rightarrow \prod_{\alpha \in A} \bHom_{\calC}(X_{\alpha},Y)$$ 
is a homotopy equivalence for every object $Y \in \calC$. 
Consequently, we recover the description given in Example \ref{examprod}. As we noted earlier, this characterization can be stated entirely in terms of the enriched homotopy category $\h{\calC}$: the maps $\{ \phi_{\alpha} \} $ exhibit $X$ as a homotopy coproduct of the family $\{ X_{\alpha} \}_{\alpha \in A}$
if and only if the induced map
$$ \bHom_{\calC}(X,Y) \rightarrow \prod_{\alpha \in A} \bHom_{\calC}(X_{\alpha},Y)$$ is an isomorphism in the homotopy category $\calH$ of spaces, for each $Y \in \calC$.

Now suppose that $\calC$ is an $\infty$-category, and let $p: A \rightarrow \calC$ be a map. 
As above, we may identify this with a collection of objects $\{ X_{\alpha} \}_{\alpha \in A}$ of $\calC$. 
To specify an object of $\calC_{p/}$ is to give an object $X \in \calC$ together with morphisms
$\phi_{\alpha}: X_{\alpha} \rightarrow X$ for each $\alpha \in A$. Using Theorem \ref{colimcomparee}, we deduce that
$X$ is a colimit of the diagram $p$ if and only if the induced map
$$ \bHom_{\calC}(X,Y) \rightarrow \prod_{\alpha \in A} \bHom_{\calC}(X_{\alpha},Y)$$
is an isomorphism in $\calH$, for each object $Y \in \calC$. In this case, we say that $X$ is a {\it coproduct} of the family $\{ X_{\alpha} \}_{\alpha \in A}$.

In either setting, we will denote the (homotopy) coproduct of a family of objects $\{ X_\alpha \}_{\alpha \in A}$ by
$$ \coprod_{\alpha \in A} X_I.$$
It is well-defined up to (essentially unique) equivalence.

Using Corollary \ref{util}, we deduce the following:

\begin{proposition}\index{gen}{coproduct}\index{gen}{product}\label{makerus}
Let $\calC$ be an $\infty$-category, and let $\{ p_{\alpha}: K_{\alpha} \rightarrow \calC\}_{\alpha \in A}$ be a family of diagrams in $\calC$ indexed by a set $A$. Suppose that each $p_{\alpha}$ has a colimit $X_{\alpha}$ in $\calC$. Let $K = \coprod K_{\alpha}$, and let $p: K \rightarrow \calC$ be the result of amalgamating the maps $p_{\alpha}$. Then $p$ has a colimit in $\calC$ if and only if the family $\{ X_{\alpha} \}_{\alpha \in A}$ has a coproduct in $\calC$; in this case, one may identify colimits of $p$ with coproducts $ \coprod_{\alpha \in A} X_{\alpha}.$
\end{proposition}

\subsection{Pushouts}\label{quasilimit6}

Let $\calC$ be an $\infty$-category. A {\it square} in $\calC$ is a map
$\Delta^1 \times \Delta^1 \rightarrow \calC$. We will typically denote
squares in $\calC$ by diagrams
$$ \xymatrix{ X' \ar[r]^{p'} \ar[d]^{q'} & X \ar[d]^{q} \\
Y' \ar[r]^{p} & Y, }$$
with the ``diagonal'' morphism $r: X' \rightarrow Y$ and homotopies
$r \simeq q \circ p'$, $r \simeq p \circ q'$ being implicit.\index{gen}{square}

We have isomorphisms of simplicial sets
$$ (\Lambda^2_0)^{\triangleright} \simeq \Delta^1 \times \Delta^1 \simeq (\Lambda^2_2)^{\triangleleft}.$$
Consequently, given a square $\sigma: \Delta^1 \times \Delta^1 \rightarrow \calC$, it makes sense to ask whether or not $\sigma$ is a limit or colimit diagram. If $\sigma$ is a limit diagram, we will also say that $\sigma$ is a {\it pullback square} or a {\it Cartesian square}, and we will informally write $X' = X \times_{Y} Y'$. Dually, if $\sigma$ is a colimit diagram, we will say that $\sigma$ is a {\it pushout square} or a {\it coCartesian square}, and write $Y = X \coprod_{X'} Y'$.\index{gen}{square!pullback}\index{gen}{square!pushout}\index{gen}{pullback}\index{gen}{pushout}\index{gen}{square!Cartesian}\index{gen}{square!coCartesian}

Now suppose that $\calC$ is a (fibrant) simplicial category. By definition, a commutative diagram
$$ \xymatrix{ X' \ar[r]^{p'} \ar[d]^{q'} & X \ar[d]^{q} \\
Y' \ar[r]^{p} & Y }$$
is a homotopy pushout square if, for every object $Z \in \calC$, the diagram
$$ \xymatrix{ \bHom_{\calC}(Y,Z) \ar[r] \ar[d] & \bHom_{\calC}(Y',Z) \ar[d] \\
\bHom_{\calC}(X,Z) \ar[r] & \bHom_{\calC}(X',Z) }$$
is a homotopy pullback square in $\Kan$. Using Theorem \ref{colimcomparee}, we can reduce questions about pushout diagrams in an arbitrary $\infty$-category to questions about homotopy pullback squares in $\Kan$.\index{gen}{pushout!homotopy}

The following basic transitivity property for pushout squares will be used repeatedly throughout this book: 

\begin{lemma}\label{transplantt}
Let $\calC$ be an $\infty$-category, and suppose given a map 
$\sigma: \Delta^2 \times \Delta^1 \rightarrow \calC$ which we will depict as a diagram
$$ \xymatrix{ X \ar[r] \ar[d] & Y \ar[d] \ar[r] & Z \ar[d] \\
X' \ar[r] & Y' \ar[r] & Z'. }$$
Suppose that the left square is a pushout in $\calC$. Then the right square is a pushout
if and only if the outer square is a pushout.
\end{lemma}

\begin{proof}
For every subset $A$ of $\{ x,y,z, x', y', z'\}$, let $\calD(A)$ denote the corresponding
full subcategory of $\Delta^2 \times \Delta^1$, and let $\sigma(A)$ denote the restriction of $\sigma$ to $\calD(A)$. We may regard $\sigma$ as determining an object
$\widetilde{\sigma} \in \calC_{\sigma(\{ y,z,x',y',z' \})/}$. Consider the maps
$$ \calC_{ \sigma(\{ z,x',z' \} ) /} \stackrel{\phi}{\leftarrow}
\calC_{ \sigma(\{ y,z, x',y',z'\})/} \stackrel{\psi}{\rightarrow} \calC_{\sigma(\{ y,x',y' \})/} $$
The map on $\phi$ is the composition of the trivial fibration
$$ \calC_{\sigma(\{ z,x',y',z' \} )/} \rightarrow \calC_{\sigma(\{ z,x',z' \}) /}$$ with a pullback of 
$$\calC_{\sigma( \{y,z,y',z'\} ) /} \rightarrow \calC_{ \sigma( \{z,y',z' \}) /},$$ also a trivial fibration in virtue of our assumption that the square
$$ \xymatrix{ Y \ar[d] \ar[r] & Z \ar[d] \\
Y' \ar[r] & Z' }$$
is a pullback in $\calC$. The map $\psi$ is a trivial fibration because the inclusion
$\calD( \{ y,x',y' \}) \subseteq \calD( \{y,z,x',y',z' \})$ is left anodyne. It follows that
$\phi(\widetilde{\sigma})$ is an initial object of $ \calC_{ \sigma(\{ z,x',z' \} ) /}$
if and only if $\psi(\widetilde{\sigma})$ is an initial object of 
$\calC_{\sigma(\{ y,x',y' \}) /}$, as desired.
\end{proof}

Our next objective is to apply Proposition \ref{utl} to show that
in many cases, complicated colimits may be decomposed as pushouts
of simpler colimits. Suppose given a pushout diagram of simplicial sets
$$ \xymatrix{ L' \ar[r]^{i} \ar[d] & L \ar[d] \\
K' \ar[r] & K }$$
and a diagram $p: K \rightarrow \calC$, where $\calC$ is an $\infty$-category.
Suppose furthermore that $p|K'$, $p|L'$, and $p|L$ admit colimits in $\calC$, which we will denote by $X$, $Y$, and $Z$, respectively. If we suppose further that the map $i$ is a cofibration of simplicial sets, then the hypotheses of Proposition \ref{extet} are satisfied. 
Consequently, we deduce:

\begin{proposition}\label{train}
Let $\calC$ be an $\infty$-category, and let $p: K \rightarrow \calC$ be a map of simplicial
sets. Suppose given a decomposition $K = K' \coprod_{L'} L$, where $L' \rightarrow L$ is a monomorphism of simplicial sets. Suppose further that
$p|K'$ has a colimit $X \in \calC$, $p|L'$ has a colimit $Y \in \calC$, and $p|L$ has
a colimit $Z \in \calC$. Then one may identify colimits for $p$ with
pushouts $X \coprod_{Y} Z$.
\end{proposition}

\begin{remark}
The statement of Proposition \ref{train} is slightly vague.
Implicit in the discussion is that identifications of $X$ with the
colimit of $p|K'$ and $Y$ with the colimit of $p|L'$ induce a
morphism $Y \rightarrow X$ in $\calC$ (and similarly for $Y$ and $Z$).
This morphism is not uniquely determined, but it is determined up
to a contractible space of choices: see the proof of Proposition
\ref{extet}.
\end{remark}

It follows from Proposition \ref{train} that any finite colimit
can be built using initial objects and pushout squares. For example,
we have the following:

\begin{corollary}\label{allfin}\index{gen}{colimit!finite}
Let $\calC$ be an $\infty$-category. Then $\calC$ admits all finite colimits if
and only if $\calC$ admits pushouts and has an initial object.
\end{corollary}

\begin{proof}
The ``only if'' direction is clear. For the converse, let us
suppose that $\calC$ has pushouts and an initial object. Let $p: K
\rightarrow \calC$ be any diagram, where $K$ is a finite simplicial
set: that is, $K$ has only finitely many nondegenerate simplices.
We will prove that $p$ has a colimit. The proof goes by induction:
first on the dimension of $K$, then on the number of simplices of
$K$ having the maximal dimension.

If $K$ is empty, then an initial object of $\calC$ is a colimit for
$p$. Otherwise, we may fix a nondegenerate simplex of $K$ having
the maximal dimension, and thereby decompose $K \simeq K_0
\coprod_{ \bd \Delta^n } \Delta^n$. By the inductive hypothesis,
$p|K_0$ has a colimit $X$ and $p| \bd \Delta^n$ has a colimit $Y$.
The $\infty$-category $\Delta^n$ has a final object, so $p| \Delta^n$
has a colimit $Z$ (which we may take to be $p(v)$, where $v$ is
the final vertex of $\Delta^n$). Now we simply apply Proposition
\ref{train} to deduce that $X \coprod_Y Z$ is a colimit for $p$.
\end{proof}

Using the same argument, one can show:

\begin{corollary}\label{allfinn}
Let $f: \calC \rightarrow \calC'$ be a functor between $\infty$-categories. Assume
that $\calC$ has all finite colimits. Then $f$ preserves all finite
colimits if and only if $f$ preserves initial objects and pushouts.
\end{corollary}

We conclude by showing how {\em all} colimits can be constructed
out of simple ones.

\begin{proposition}\label{alllimits}
Let $\calC$ be an $\infty$-category. Suppose that $\calC$ admits pushouts and
$\kappa$-small coproducts. Then $\calC$ admits colimits for all $\kappa$-small diagrams.
\end{proposition}

\begin{proof}
If $\kappa = \omega$, we have already shown this as Corollary
\ref{allfin}. Let us therefore suppose that $\kappa > \omega$, and
that $\calC$ has pushouts and $\kappa$-small sums.

Let $p: K \rightarrow \calC$ be a diagram, where $K$ is
$\kappa$-small. We first suppose that the dimension $n$ of $K$ is
finite: that is, $K$ has no nondegenerate simplices of dimension
$> n$. We prove that $p$ has a colimit, working by induction on
$n$.

If $n=0$, then $K$ consists of a finite disjoint union of fewer
than $\kappa$ vertices. The colimit of $p$ exists by the
assumption that $\calC$ has $\kappa$-small sums.

Now suppose that every diagram indexed by a $\kappa$-small
simplicial set of dimension $n$ has a colimit. Let $p: K
\rightarrow \calC$ be a diagram, with the dimension of $K$ equal to
$n+1$. Let $K^n$ denote the $n$-skeleton of $K$, and $K'_{n+1} \subseteq K_{n+1}$ the
set of all nondegenerate $(n+1)$-simplices of $K$, so that there is a pushout
diagram of simplicial sets
$$ K^n \coprod_{ K'_{n+1} \times \bd \Delta^{n+1}} (K'_{n+1} \times
\Delta^{n+1}) \simeq K.$$ By Proposition \ref{train}, we can
construct a colimit of $p$ as a pushout, using colimits for
$p|K^n$, $p|(K'_{n+1} \times \bd \Delta^{n+1})$, and $p|(K'_{n+1}
\times \Delta^{n+1})$. The first two exist by the inductive
hypothesis; the last, because it is a sum of diagrams which
possess colimits.

Now let us suppose that $K$ is not necessarily finite dimensional.
In this case, we can filter $K$ by its skeleta $\{ K^n \}$. This
is a family of simplicial subsets of $K$ indexed by the set
$\Z_{\geq 0}$ of nonnegative integers. By what we have shown
above, each $p|K^n$ has a colimit $x_n$ in $\calC$. Since this family
is directed and covers $K$, Corollary \ref{util} shows that we may identify colimits of
$p$ with colimits of a diagram $\Nerve(\Z_{\geq 0}) \rightarrow \calC$ which we may write informally as
$$ x_0 \rightarrow x_1 \rightarrow \ldots $$

Let $L$ be the simplicial subset of
$\Nerve(\Z_{\geq 0})$ which consists of all vertices, together with the
edges which join consecutive integers. A simple computation shows
that the inclusion $L \subseteq \Nerve(\Z_{\geq 0})$ is a categorical
equivalence, and therefore cofinal. Consequently, it suffices to
construct the colimit of a diagram $L \rightarrow \calC$. But $L$ is
$1$-dimensional, and is $\kappa$-small since $\kappa > \omega$.
\end{proof}

The same argument proves also the following:

\begin{proposition}\label{allimits}
Let $\kappa$ be a regular cardinal, and let $f: \calC \rightarrow \calD$ be a functor between $\infty$-categories, where $\calC$ admits $\kappa$-small colimits. Then
$f$ preserves $\kappa$-small colimits if and only if $f$ preserves pushout squares
and $\kappa$-small coproducts.
\end{proposition}

Let $\calD$ be an $\infty$-category containing an object $X$, and suppose that $\calD$
admits pushouts. Then $\calD_{X/}$ admits pushouts, and these pushouts map be computed in $\calD$. In other words, the projection $f: \calD_{X/} \rightarrow \calD$ preserves pushouts. In fact, this is a special case of a very general result; it requires only that $f$ is a left fibration and the
simplicial set $\Lambda^2_0$ is weakly contractible.

\begin{lemma}\label{pregoes}
Let $f: \calC \rightarrow \calD$ be a left fibration of $\infty$-categories, and let
$K$ be a weakly contractible simplicial set. Then any map
$\overline{p}: K^{\triangleright} \rightarrow \calC$ is an $f$-colimit diagram.
\end{lemma}

\begin{proof}
Let $p = \overline{p} | K$. We must show that the map
$$\phi:  \calC_{\overline{p} /} \rightarrow \calC_{p/} \times_{\calD_{f \circ p/}} \calD_{f \circ \overline{p}/}$$ is a trivial fibration of simplicial sets. In other words, we must show that
we can solve any lifting problem of the form
$$ \xymatrix{ (K \star A) \coprod_{ K \star A_0} (K^{\triangleright} \star A) \ar[r] \ar@{^{(}->}[d] & \calC \ar[d]^{f} \\
K^{\triangleright} \star A \ar[r] \ar@{-->}[ur] & \calD. }$$
Since $f$ is a left fibration, it will suffice to prove that the left vertical map is left anodyne, which follows immediately from Lemma \ref{chotle}. 
\end{proof}

\begin{proposition}\label{goeselse}
Let $f: \calC \rightarrow \calD$ be a left fibration of $\infty$-categories, and let
$p: K \rightarrow \calC$ be a diagram. Suppose that
$K$ is weakly contractible. Then:
\begin{itemize}
\item[$(1)$] Let $\overline{p}: K^{\triangleright} \rightarrow \calC$ be an extension of $p$.
Then $\overline{p}$ is a colimit of $p$ if and only if $f \circ \overline{p}$ is a colimit of $f \circ p$.

\item[$(2)$] Let $\overline{q}: K^{\triangleright} \rightarrow \calD$ be a colimit of $f \circ p$.
Then $\overline{q} = f \circ \overline{p}$, where $\overline{p}$ is an extension $($ automatically a colimit, in virtue of $(1)$ $)$ of $p$.
\end{itemize}
\end{proposition}

\begin{proof}
To prove $(1)$, fix an extension $\overline{p}: K^{\triangleright} \rightarrow \calC$. 
We have a commutative diagram
$$ \xymatrix{ \calC_{\overline{p}/} \ar[r]^-{\phi} & \calC_{p/} \times_{\calD_{fp/}} \calD_{f\overline{p}/} \ar[r]^-{\psi'} \ar[d] & \calC_{p/} \ar[d]^{\theta} \\
& \calD_{f \overline{p}/} \ar[r]^-{\psi} & \calD_{fp/}. }$$
Lemma \ref{pregoes} implies that $\phi$ is a trivial Kan fibration.
If $f \circ \overline{p}$ is a colimit diagram, the map $\psi$ is a trivial fibration. Since $\psi'$ is a pullback of $\psi$, we conclude that $\psi'$ is a trivial fibration. It follows that $\psi' \circ \phi$ is a trivial fibration, so that $\overline{p}$ is a colimit diagram. This proves the ``if'' direction of $(1)$.

To prove the converse, let us suppose that $\overline{p}$ is a colimit diagram. The maps
$\phi$ and $\psi' \circ \phi$ are both trivial fibrations. It follows that the fibers of $\psi'$ are contractible. Using Lemma \ref{chotle2}, we conclude that the map $\theta$ is a trivial fibration, and therefore surjective on vertices. It follows that the fibers of $\psi$ are contractible. Since $\psi$ is a left fibration with contractible fibers, it is a trivial fibration (Lemma \ref{toothie}). Thus $f \circ \overline{p}$ is a colimit diagram and the proof is complete.

To prove $(2)$, it suffices to show that $f$ has the right lifting property with respect to the inclusion
$i: K \subseteq K^{\triangleright}$. Since $f$ is a left fibration, it will suffice to show that $i$ is left anodyne, which follows immediately from Lemma \ref{chotle2}.
\end{proof}

\subsection{Coequalizers}\label{quasilimit8}

Let $\calI$ denote the category depicted by the diagram
$$\xymatrix{ X \ar@<.4ex>[r]^{F} \ar@<-.4ex>[r]_{G} & Y}.$$
In other words, $\calI$ has two objects, $X$ and $Y$, with
$ \Hom_{\calI}(X,X) = \Hom_{\calI}(Y,Y) = \ast$, $ \Hom_{\calI}(Y,X) = \emptyset$, and
$ \Hom_{\calI}(X,Y) = \{ F,G \}.$

To give a diagram $p: \Nerve(\calI) \rightarrow \calC$ in an $\infty$-category $\calC$, one must give a pair of morphisms $f = p(F)$, $g = p(G)$ in $\calC$, having the same domain $x=p(X)$ and the same codomain $y = p(Y)$. A colimit for the diagram $p$ is said to be a {\em coequalizer} of $f$ and $g$.\index{gen}{coequalizer}\index{gen}{equalizer}

Applying Corollary \ref{util}, we deduce the following:

\begin{proposition}\label{uty}
Let $K$ and $A$ be a simplicial sets, and let $i_0, i_1: A \rightarrow K$ be embeddings having disjoint images in $K$. Let $K'$ denote the coequalizer of $i_0$ and $i_1$; in other words, the simplicial set obtained from $K$ by identifying the image of $i_0$ with the image of $i_1$.
Let $p: K' \rightarrow \calC$ be a diagram in an $\infty$-category $S$, and let $q: K \rightarrow \calC$
be the composition
$$ K \rightarrow K' \stackrel{p}{\rightarrow} S.$$
Suppose that the diagrams $q \circ i_0= q \circ i_1$ and $q$ possess colimits $x$ and $y$ in $S$. 
Then $i_0$ and $i_1$ induce maps $j_0, j_1: x \rightarrow y$ $($well-defined up to homotopy$)$; 
colimits for $p$ may be identified with coequalizers of $j_0$ and $j_1$.
\end{proposition}

Like pushouts, coequalizers are a basic construction out of which other colimits can be built.
More specifically, we have the following:

\begin{proposition}\label{appendixdiagram}
Let $\calC$ be an $\infty$-category and $\kappa$ a regular cardinal. Then $\calC$ has all $\kappa$-small colimits if and only if $\calC$ has coequalizers and $\kappa$-small coproducts.
\end{proposition}

\begin{proof}
The ``only if'' direction is obvious. For the converse, suppose that $\calC$ has coequalizers and $\kappa$-small coproducts. In view of Proposition \ref{alllimits}, it suffices to show that $\calC$ has pushouts.
Let $p: \Lambda^2_0$ be a pushout diagram in $\calC$. We note that $\Lambda^2_0$ is the quotient of $\Delta^{ \{0,1\} } \coprod \Delta^{ \{0,2\} }$ obtained by identifying the initial vertex of $\Delta^{ \{0,1\} }$ with the initial vertex of $\Delta^{ \{0,2\} }$. In view of Proposition \ref{uty}, it suffices to show that $p| ( \Delta^{ \{0,1\} } \coprod \Delta^{ \{0,2\} } )$ and $p| \{0\}$ possess colimits in $\calC$. The second assertion is obvious. Since $\calC$ has finite sums, to prove that there exists a colimit for
$p| ( \Delta^{ \{0,1\} } \coprod \Delta^{ \{0,2\} } )$, it suffices to prove that $p| \Delta^{ \{0,1\}}$ and
$p|\Delta^{ \{0,2\} }$ possess colimits in $\calC$. This is immediate, since $\Delta^{ \{0,1\} }$ and $\Delta^{ \{0,2\} }$ both have final objects.
\end{proof}

Using the same argument, we deduce:

\begin{proposition}\label{appendicites}
Let $\kappa$ be a regular cardinal and $\calC$ be an $\infty$-category which admits $\kappa$-small colimits. A full subcategory $\calD \subseteq \calC$ is stable under $\kappa$-small colimits in $\calC$ if and only if $\calD$ is stable under coequalizers and under $\kappa$-small sums.
\end{proposition}

\subsection{Tensoring with Spaces}\label{quasilimit7}

Every ordinary category $\calC$ can be regarded as a category enriched over $\Set$.
Moreover, if $\calC$ admits coproducts, then $\calC$ can be regarded as {\em tensored} over $\Set$, in an essentially unique way. In the $\infty$-categorical setting, one has a similar situation: if $\calC$ is an $\infty$-category which admits all small limits, then $\calC$ may be regarded as tensored over the $\infty$-category $\SSet$ of spaces. To make this idea precise, we would need a good theory of {\em enriched} $\infty$-categories, which lies outside the scope of this book. We will instead settle for a slightly ad-hoc point of view which nevertheless allows us to construct the relevant tensor products. We begin with a few remarks concerning representable functors in the $\infty$-categorical setting.

\begin{definition}\label{repfunc}
Let $\calD$ be a closed monoidal category, and let $\calC$ be a category enriched over $\calD$.
We will say that a $\calD$-enriched functor $G: \calC^{op} \rightarrow \calD$ is {\it representable} if there exists an object $C \in \calC$ and a map $\eta: 1_{\calD} \rightarrow G(C)$ such that the induced map
$$ \bHom_{\calC}(X,C) \simeq \bHom_{\calC}(X,C) \otimes 1_{\calD}
\rightarrow \bHom_{\calC}(X,C) \otimes G(C) \rightarrow G(X)$$
is an isomorphism, for every object $X \in \calC$. In this case, we will say that 
$(C,\eta)$ {\it represents} the functor $F$. 
\end{definition}

\begin{remark}
In the situation of Definition \ref{repfunc}, we will sometimes abuse terminology and simply say that the functor $F$ is {\it represented} by the object $C$.\index{gen}{representable!functor}
\end{remark}

\begin{remark}
The dual notion of a {\it corepresentable functor} is may be defined in an obvious way.
\end{remark}

\begin{definition}\index{gen}{functor!representable by an object}
Let $\calC$ be an $\infty$-category, and let $\SSet$ denote the $\infty$-category of spaces.
We will say that a functor $F: \calC^{op} \rightarrow \SSet$ is {\it representable} if
the underlying functor
$$ \h{F}: \h{\calC}^{op} \rightarrow \h{\SSet} \simeq \calH$$
of $($ $\calH$-enriched $)$ homotopy categories is representable. We will say that a pair
$C \in \calC$, $\eta \in \pi_0 F(C)$ {\it represents} $F$ if the pair $(C,\eta)$ represents
$h F$. 
\end{definition}

\begin{proposition}\label{reppfunc}
Let $f: \widetilde{\calC} \rightarrow \calC$ be a right fibration of $\infty$-categories, let
$\widetilde{C}$ be an object of $\widetilde{\calC}$, $C = f(\widetilde{C}) \in \calC$, and let
$F: \calC^{op} \rightarrow \SSet$ be a functor which classifies $f$ (\S \ref{universalfib}).
The following conditions are equivalent:
\begin{itemize}
\item[$(1)$] Let $\eta \in \pi_0 F(C) \simeq \pi_0 
( \widetilde{\calC} \times_{\calC} \{C\} )$ be the connected component containing
$\widetilde{C}$. Then the pair $(C, \eta)$ represents the functor $F$. 

\item[$(2)$] The object $\widetilde{C} \in \widetilde{\calC}$ is final.

\item[$(3)$] The inclusion $\{ \widetilde{\calC} \} \subseteq \widetilde{\calC}$ is a contravariant equivalence in $(\sSet)_{/\calC}$. 
\end{itemize}
\end{proposition}

\begin{proof}
We have a commutative diagram of right fibrations
$$ \xymatrix{ \widetilde{\calC}_{/ \widetilde{C}} \ar[r]^{\phi} \ar[d] & \widetilde{\calC} \ar[d] \\
\calC_{/C} \ar[r] & \calC. }$$
Observe that the left vertical map is actually a trivial fibration.
Fix an object $D \in \calC$. The fiber of the upper horizontal map
$$ \phi_{D}: \widetilde{\calC}_{/ \widetilde{C}} \times_{\calC} \{D\} \rightarrow \widetilde{\calC} \times_{\calC} \{D\}$$ can be identified, in the homotopy category $\calH$, with the map
$\bHom_{\calC}(D,C) \rightarrow F(C).$ The map $\phi_{D}$ is a right fibration of Kan complexes, and therefore a Kan fibration. If $(1)$ is satisfied, then $\phi_{D}$ is a homotopy equivalence, and therefore a trivial fibration. It follows that the fibers of $\phi$ are contractible. Since $\phi$ is a right fibration, it is a trivial fibration (Lemma \ref{toothie}). This proves that $\widetilde{C}$ is a final object of $\widetilde{\calC}$. Conversely, if $(2)$ is satisfied, then $\phi_{D}$ is a trivial Kan fibration and therefore a weak homotopy equivalence. Thus $(1) \Leftrightarrow (2)$. 

If $(2)$ is satisfied, then the inclusion $\{ \widetilde{C} \} \subseteq \widetilde{\calC}$ is right anodyne, and therefore a contravariant equivalence by Proposition \ref{hunef}. Thus $(2) \Rightarrow (3)$. Conversely, suppose that $(3)$ is satisfied. The inclusion
$\{\id_{C} \} \subseteq \calC_{/C}$ is right anodyne, and therefore a contravariant equivalence. It follows that the lifting problem
$$ \xymatrix{ \{\id_{C} \} \ar[r]^{\widetilde{C}} \ar@{^{(}->}[d] & \widetilde{\calC} \ar[d]^{f} \\
\calC_{/C} \ar[r] \ar@{-->}[ur]^{e} & \calC }$$
has a solution. We observe that $e$ is a contravariant equivalence of right fibrations over $\calC$, and therefore a categorical equivalence. By construction, $e$ carries a final object
of $\calC_{/C}$ to $\widetilde{C}$, so that $\widetilde{C}$ is a final object of $\widetilde{\calC}$. 
\end{proof}

We will say that a right fibration $\widetilde{\calC} \rightarrow \calC$ is {\it representable} if
$\widetilde{\calC}$ has a final object.\index{gen}{representable!right fibration}

\begin{remark}
Let $\calC$ be an $\infty$-category, and let $p: K \rightarrow \calC$ be a diagram. Then
the right fibration $\calC_{/p} \rightarrow \calC$ is representable if and only if $p$ has a limit in $\calC$.
\end{remark}

\begin{remark}\index{gen}{corepresentable!functor}\index{gen}{corepresentable!left fibration}\index{gen}{functor!corepresentable}
All of the above ideas dualize in an evident way, so that we may speak of {\it corepresentable functors} and {\it corepresentable left fibrations} in the setting of $\infty$-categories.
\end{remark}

\begin{notation}
For each diagram $p: K \rightarrow \calC$ in an $\infty$-category $\calC$, we let
$\calF_{p}: \h{\calC} \rightarrow \calH$ denote the $\calH$-enriched functor
corresponding to the left fibration $\calC^{p/} \rightarrow \calC$.

If $p: \ast \rightarrow \calC$ is the inclusion of an object $X$ of $\calC$, then we write $\calF_{X}$ for $\calF_{p}$. We note that $\calF_{X}$ is the functor corepresented by $X$:
$$ \calF_X(Y) = \bHom_{\calC}(X,Y).$$
\end{notation}

Now suppose that $X$ is an object in an $\infty$-category $\calC$, and let $p: K \rightarrow \calC$
be a constant map taking the value $X$. For every object $Y$ of $\calC$, we have an isomorphism of simplicial sets
$(\calC^{p/}) \times_{\calC} \{Y\} \simeq ( \calC^{X/} \times_{\calC} \{Y\} )^K$. This identification is functorial up to homotopy, so we actually obtain an equivalence
$$\calF_{p}(Y) \simeq \bHom_{\calC}(X,Y)^{[K]}$$
in the homotopy category $\calH$ of spaces, where $[K]$ denotes the simplicial set $K$ regarded as an object of $\calH$. Applying Proposition \ref{reppfunc}, we deduce the following:

\begin{corollary}\label{charext}
Let $\calC$ be an $\infty$-category, $X$ and object of $\calC$, and $K$ a simplicial set.
Let $p: K \rightarrow \calC$ be the constant map taking the value $X$. The objects of the fiber
$\calC^{p/} \times_{\calC} \{Y\}$ are classified, up to equivalence, by maps
$\psi: [K] \rightarrow \bHom_{\calC}(X,Y)$
in the homotopy category $\calH$. Such a map $\psi$ classifies a colimit for $p$ if and only if
it induces isomorphisms
$$ \bHom_{\calC}(Y,Z) \simeq \bHom_{\calC}(X,Z)^{[K]} $$ 
in the homotopy category $\calH$, for every object $Z$ of $\calC$.
\end{corollary}

In the situation of Corollary \ref{charext}, we will denote a colimit for $p$ by $X \otimes K$, if such a colimit exists. We note that $X \otimes K$ is well defined up to (essentially unique) equivalence, and that it depends (up to equivalence) only on the weak homotopy type of the simplicial set $K$.\index{not}{XotimesK@$X \otimes K$}\index{gen}{tensor product with spaces}

\begin{corollary}\label{silt}
Let $\calC$ be an $\infty$-category, let $K$ be a weakly contractible simplicial set, and let
$p: K \rightarrow \calC$ be a diagram which carries each edge of $K$ to an equivalence in $\calC$.
Then:
\begin{itemize}
\item[$(1)$] The diagram $p$ has a colimit in $\calC$.
\item[$(2)$] An arbitrary extension $\overline{p}: K^{\triangleright} \rightarrow \calC$
is a colimit for $\calC$ if and only if $\overline{p}$ carries each edge of
$K^{\triangleright} \rightarrow \calC$ to an equivalence in $\calC$.
\end{itemize}
\end{corollary}

\begin{proof}
Let $\calC' \subseteq \calC$ be the largest Kan complex contained in $\calC$. By assumption, $p$ factors through $\calC'$. Since $K$ is weakly contractible, we conclude that $p: K \rightarrow \calC'$ is homotopic to a constant map $p': K \rightarrow \calC'$. Replacing $p$ by $p'$ if necessary, we may reduce to the case where $p$ is constant, taking value equal to some
fixed object $C \in \calC$.

Let $\overline{p}: K^{\triangleright} \rightarrow \calC$ be the constant map with value $C$. Using the characterization of colimits in Corollary \ref{charext}, we deduce that $\overline{p}$ is a colimit diagram in $\calC$. This proves $(1)$, and (in view of the uniqueness of colimits up to equivalence) the ``only if'' direction of $(2)$. To prove the converse, we suppose that $\overline{p}'$ is an arbitrary extension of $p$ which carries each edge of $K^{\triangleright}$ to an equivalence in $\calC$. Then $\overline{p}'$ factors through $\calC'$. Since $K^{\triangleright}$ is weakly contractible, we conclude as above that $\overline{p}'$ is homotopic to a constant map, and is therefore a colimit diagram.
\end{proof}

\subsection{Retracts and Idempotents}\label{retrus}

Let $\calC$ be a category. An object $Y \in \calC$ is said to be a {\it retract} of an object\index{gen}{retract} $X \in \calC$ if there is a commutative diagram
$$ \xymatrix{ & X \ar[dr]^{r} & \\
Y \ar[rr]^{\id_{Y}} \ar[ur]^{i} & & Y}$$
in $\calC$. In this case we can identify $Y$ with a subobject of $X$ via the monomorphism $i$, and think of $r$ as a retraction from $X$ onto $Y \subseteq X$. We observe also that the map
$i \circ r: X \rightarrow X$ is idempotent. Moreover, this idempotent determines $Y$ up to canonical isomorphism: we can recover $Y$ as the equalizer of the pair of maps $(\id_{X}, i \circ r): X \rightarrow X$ (or, dually, as the coequalizer of the same pair of maps). Consequently, we obtain an injective map from the collection of isomorphism classes of retracts of $X$ to the set of idempotent maps $f: X \rightarrow X$. We will say that $\calC$ is {\it idempotent complete} if this correspondence is bijective for every $X \in \calC$: that is, if every idempotent map $f: X \rightarrow X$ comes from a (uniquely determined) retract of $X$. If $\calC$ admits equalizers (or coequalizers), then $\calC$ is 
idempotent complete.\index{gen}{idempotent!completeness}\index{gen}{idempotent!in an ordinary category}

These ideas can be adapted to the $\infty$-categorical setting in a straightforward way. If $X$ and $Y$ are objects of an $\infty$-category $\calC$, then we say that $Y$ is a {\it retract} of $X$ if it is a retract of $X$ in the homotopy category $\h{\calC}$. Equivalently, $Y$ is a retract of $X$ if there exists a $2$-simplex $\Delta^2 \rightarrow \calC$ corresponding to a diagram
$$ \xymatrix{ & X \ar[dr]^{r} & \\
Y \ar[rr]^{\id_{Y}} \ar[ur]^{i} & & Y.}$$
As in the classical case, there is a correspondence between retracts $Y$ of $X$ and idempotent maps $f: X \rightarrow X$. However, there are two important differences: first, the notion of an idempotent map needs to be interpreted in an $\infty$-categorical sense. It is not enough to require that $f = f \circ f$ in the homotopy category $\h{\calC}$. This would correspond to the condition that there is a path $p$ joining $f$ to $f \circ f$ in the endomorphism space of $X$, which would give rise to {\em two} paths from $f$ to $f \circ f \circ f$. In order to have a hope of recovering $Y$, we need these paths to be homotopic. This condition does not even make sense unless $p$ is specified; thus we must take $p$ as part of the data of an idempotent map. In other words, in the $\infty$-categorical setting, idempotence is not merely a condition, but involves additional data (see Definition \ref{diags}).

The second important difference between the classical and $\infty$-categorical theory of retracts is that in the $\infty$-categorical case, one cannot recover a retract $Y$ of $X$ as the limit (or colimit) of a {\em finite} diagram involving $X$.

\begin{example}
Let $R$ be a commutative ring, and let $C_{\bigdot}(R)$ be the category of complexes of finite free $R$-modules, so that an object of $C_{\bigdot}(R)$ is a chain complex
$$ \ldots \rightarrow M_{1} \rightarrow M_0 \rightarrow M_{-1} \rightarrow \ldots $$
such that each $M_{i}$ is a finite free $R$-module, and $M_{i} = 0$ for $|i| \gg 0$; morphisms in $C_{\bigdot}(R)$ are given by morphisms of chain complexes. There is a natural simplicial structure on the category $C_{\bigdot}(R)$, for which the mapping spaces are Kan complexes; let $\calC = \sNerve( C_{\bigdot}(R) )$ be the associated $\infty$-category. Then $\calC$ admits all finite limits and colimits ($\calC$ is an example of a {\em stable} $\infty$-category; see \cite{DAG}). However, $\calC$ is idempotent complete if and only if every finitely generated projective $R$-module is stably free.
\end{example}

The purpose of this section is to define the notion of an {\it idempotent} in an $\infty$-category $\calC$, and to obtain a correspondence between idempotents and retracts in $\calC$.

\begin{definition}\index{not}{Idem+@$\Idem^{+}$}\index{not}{Idem@$\Idem$}\index{not}{Ret@$\Ret$}
The simplicial set $\Idem^{+}$ is defined as follows: 
for every nonempty, finite, linearly ordered set $J$, $\Hom_{\sSet}(\Delta^{J}, \Idem^{+})$ can be identified with the set of pairs $(J_0, \sim)$, where $J_0 \subseteq J$
and $\sim$ is an equivalence relation on $J_0$ which satisfies the following condition:
\begin{itemize}
\item[$(\ast)$] Let $i \leq j \leq k$ be elements of $J$ such that $i,k \in J_0$, and
$i \sim k$. Then $j \in J_0$, and $i \sim j \sim k$.
\end{itemize}

Let $\Idem$ denote the simplicial subset of $\Idem^{+}$, corresponding to those pairs
$(J_0, \sim)$ such that $J=J_0$. Let $\Ret \subseteq \Idem^{+}$ denote the simplicial subset corresponding to those pairs $(J_0, \sim)$ such that the quotient $J_0 / \sim$ has at most one element.
\end{definition}

\begin{remark}
The simplicial set $\Idem$ has exactly one nondegenerate simplex in each dimension $n$ (corresponding to the equivalence relation $\sim$ on $\{0, 1, \ldots, n\}$ given by
$( i \sim j) \Leftrightarrow (i = j)$ ), and the set of nondegenerate simplices of $\Idem$ is stable under passage to faces. In fact, $\Idem$ is characterized up to unique isomorphism by these two properties.
\end{remark}

\begin{definition}\label{diags}\index{gen}{idempotent!in an $\infty$-category}\index{gen}{retraction diagram!weak}\index{gen}{retraction diagram!small}
Let $\calC$ be an $\infty$-category.
\begin{itemize}
\item[$(1)$] An {\it idempotent} in $\calC$ is a map of simplicial sets
$\Idem \rightarrow \calC$. We will refer to $\Fun(\Idem,\calC)$ as the {\it $\infty$-category of
idempotents in $\calC$}.

\item[$(2)$] A {\it weak retraction diagram} in $\calC$ is a map of simplicial sets
$\Ret \rightarrow \calC$. We will refer to $\Fun(\Ret,\calC)$ as the {\it $\infty$-category of weak retraction diagrams in $\calC$}.

\item[$(3)$] A {\it strong retraction diagram} in $\calC$ is a map of simplicial sets
$\Idem^{+} \rightarrow \calC$. We will refer to $\Fun(\Idem^{+},\calC)$ as the {\it $\infty$-category of strong retraction diagrams in $\calC$}. 
\end{itemize}
\end{definition}

We now spell out Definition \ref{diags} in more concrete terms. We first observe that
$\Idem^{+}$ has precisely two vertices. Once of these vertices, which we will denote by $x$, belongs to $\Idem$, and the other, which we will denote by $y$, does not. The simplicial
set $\Ret$ can be identified with the quotient of $\Delta^{2}$ obtained by collapsing
$\Delta^{ \{0,2\} }$ to the vertex $y$. A weak retraction diagram $F: \Ret \rightarrow \calC$
in an $\infty$-category $\calC$ can therefore be identified with a $2$-simplex
$$ \xymatrix{ & X \ar[dr] & \\
Y \ar[ur] \ar[rr]^{\id_{Y}} & & Y }$$
where $X = F(x)$ and $Y = F(y)$. In other words, it is precisely the datum that we need in order to exhibit $Y$ as a retract of $X$ in the homotopy category $\h{\calC}$.

To give an idempotent $F: \Idem \rightarrow \calC$ in $\calC$, it suffices to specify the image under $F$ of each nondegenerate simplex of $\Idem$ in each dimension $n \geq 0$. Taking $n=0$, we obtain an object $X = F(x) \in \calC$.
Taking $n=1$, we get a morphism $f: X \rightarrow X$. Taking $n=2$, we get a $2$-simplex of $\calC$ corresponding to a diagram
$$ \xymatrix{ & X \ar[dr]^{f} & \\
X \ar[ur]^{f} \ar[rr]^{f} & & X }$$
which verifies the equation $f = f \circ f$ in the homotopy category $\h{\calC}$. Taking $n > 2$, we
get higher-dimensional diagrams which express the idea that $f$ is not only idempotent ``up to homotopy'', but ``up to coherent homotopy''.

The simplicial set $\Idem^{+}$ can be thought of as ``interweaving'' its simplicial subsets
$\Idem$ and $\Ret$, so that giving a strong retraction diagram $F: \Idem^{+} \rightarrow \calC$
is equivalent to giving a weak retraction diagram
$$ \xymatrix{ & X \ar[dr]^{r} & \\
Y \ar[ur]^{i} \ar[rr]^{\id_{Y}} & & Y }$$
together with a coherently idempotent map $f = i \circ r: X \rightarrow X$. Our next result makes precise the sense in which $f$ really is ``determined'' by $Y$.

\begin{lemma}\label{sturbie}
Let $J \subseteq \{ 0, \ldots, n\}$, and let $K \subseteq \Delta^n$ be the simplicial subset spanned by the nondegenerate simplices of $\Delta^n$ which do not contain $\Delta^{J}$. Suppose
that there exist $0 \leq i < j < k \leq n$ such that $i,k \in J$, $j \notin J$. Then the inclusion
$K \subseteq \Delta^n$ is inner anodyne.
\end{lemma}

\begin{proof}
Let $P$ denote the collection of all subset $J' \subseteq \{0, \ldots, n\}$ which contain $J \cup \{j\}$. Choose a linear ordering
$$ \{ J(1) \leq \ldots \leq J(m) \}$$
of $P$, with the property that if $J(i) \subseteq J(j)$, then $i \leq j$. Let
$$K(k) = K \cup \bigcup_{ 1\leq i \leq k} \Delta^{ J(i) }.$$
Note that there are pushout diagrams
$$ \xymatrix{ \Lambda^{J(i)}_{j} \ar[r] \ar[d] & \Delta^{J(i)} \ar[d] \\
K(i-1) \ar[r] & K(i). }$$
It follows that the inclusions $K(i-1) \subseteq K(i)$ are inner anodyne. Therefore
the composite inclusion $K = K(0) \subseteq K(m) = \Delta^n$ is also inner anodyne.
\end{proof}

\begin{proposition}
The inclusion $\Ret \subseteq \Idem^{+}$ is an inner anodyne map of simplicial sets.
\end{proposition}

\begin{proof}
Let $\Ret_{m} \subseteq \Idem^{+}$ be the simplicial subset defined so that
$(J_0, \sim): \Delta^{J} \rightarrow \Idem^{+}$ factors through $\Ret_m$ if and only if
the quotient $J_0 / \sim$ has cardinality $\leq m$. We observe that there is a pushout diagram
$$ \xymatrix{ K \ar[r] \ar[d] & \Delta^{2m} \ar[d] \\
\Ret_{m-1} \ar[r] & \Ret_{m} }$$
where $K \subseteq \Delta^{2m}$ denote the simplicial subset spanned by those faces
which do not contain $\Delta^{ \{1, 3, \ldots, 2m-1\} }$. If $m \geq 2$, Lemma \ref{sturbie} implies that the upper horizontal arrow is inner anodyne, so that the inclusion
$\Ret_{m-1} \subseteq \Ret_{m}$ is inner anodyne. The inclusion
$\Ret \subseteq \Idem^{+}$ can be identified with an infinite composition
$$ \Ret = \Ret_{1} \subseteq \Ret_{2} \subseteq \ldots $$
of inner anodyne maps, and is therefore inner anodyne.
\end{proof}

\begin{corollary}\label{gurgh}
Let $\calC$ be an $\infty$-category. Then the restriction map
$$ \Fun(\Idem^{+}, \calC) \rightarrow \Fun(\Ret, \calC)$$ from strong retraction diagrams
to weak retraction diagrams is a trivial fibration of simplicial sets. In particular,
every weak retraction diagram in $\calC$ can be extended to a strong retraction diagram.
\end{corollary}

We now study the relationship between strong retraction diagrams and idempotents in an $\infty$-category $\calC$. We will need the following lemma, whose proof is somewhat tedious.

\begin{lemma}
The simplicial set $\Idem^{+}$ is an $\infty$-category.
\end{lemma}

\begin{proof}
Suppose given $0 < i < n$ and a map $\Lambda^n_i \rightarrow \Idem^{+}$, corresponding
to a compatible family of pairs $\{ ( J_{k}, \sim_{k}) \}_{k \neq i}$, where $J_k \subseteq \{ 0, \ldots, k-1, k+1, \ldots, n \}$ and $\sim_k$ is an equivalence relation $J_{k}$ defining an element of
$\Hom_{\sSet}( \Delta^{ \{0, \ldots, k-1, k+1, \ldots, n\} }, \Idem^{+})$. Let $J = \bigcup J_{k}$, and
define a relation $\sim$ on $J$ as follows: if $a,b \in J$, then $a \sim b$ if and only if either
$$( \exists k \neq i) [ (a, b \in J_{k}) \wedge (a \sim_{k} b) ]$$
or
$$ ( a \neq b \neq i \neq a) \wedge (\exists c \in J_{a} \cap J_{b}) [ (a \sim_{b} c) \wedge (b \sim_{a} c)].$$
We must prove two things: that $(J, \sim) \in \Hom_{\sSet}(\Delta^n, \Idem^{+})$, and that the restriction of $(J, \sim)$ to $\{ 0, \ldots, k-1, k+1, \ldots, n\}$ coincides with $(J_k, \sim_k)$ for $k \neq i$.

We first check that $\sim$ is an equivalence relation. It is obvious that $\sim$ is reflexive and symmetric. Suppose that $a \sim b$ and that $b \sim c$; we wish to prove that $a \sim c$. There are several cases to consider:

\begin{itemize}
\item Suppose that there exists $j \neq i$, $k \neq i$ such that $a,b \in J_{j}$, 
$b,c \in J_{k}$, and $a \sim_{j} b \sim_{k} c$. If $a \neq k$, then also $a \in J_{k}$
and $a \sim_{k} b$, and we may conclude that $a \sim c$ by invoking the transitivity
of $\sim_{k}$. Therefore we may suppose that $a = k$. By the same argument, we may suppose
that $b = j$; we therefore conclude that $a \sim c$.

\item Suppose that there exists $k \neq i$ with $a,b \in J_k$, that $b \neq c \neq i \neq b$
and there exists $d \in J_{b} \cap J_{c}$ with $a \sim_{k} b \sim_{c} d \sim_{b} c$. If $a = b$ or $a=c$ there is nothing to prove; assume therefore that $a \neq b$ and $a \neq c$. Then $a \in J_{c}$
and $a \sim_{c} b$, so by transitivity $a \sim_{c} d$. Similarly, $a \in J_{b}$ and $a \sim_{b} d$
so that $a \sim_{b} c$ by transitivity. 

\item Suppose that $a \neq b \neq i \neq a$, $b \neq c \neq i \neq b$, and that there
exist $d \in J_{a} \cap J_{b}$ and $e \in J_{b} \cap J_{c}$ such that
$a \sim_{b} d \sim_{a} b \sim_{c} e \sim_{b} c$. It will suffice to prove that $a \sim_{b} c$. If $c = d$, this is clear; let us therefore assume that $c \neq d$.
By transitivity, it suffices to show that $d \sim_{b} e$. Since $c \neq d$, we have
$d \in J_{c}$ and $d \sim_{c} b$, so that $d \sim_{c} e$ by transitivity, and therefore
$d \sim_{b} e$.
\end{itemize}

To complete the proof that $(J, \sim)$ belongs to $\Hom_{\sSet}(\Delta^n, \Idem^{+})$, we must
show that if $a < b < c$, $a \in J$, $c \in J$, and $a \sim c$, then also $b \in J$ and
$a \sim b \sim c$. There are two cases to consider. Suppose first that there exists $k \neq j$ such that $a,c \in J_{k}$ and $a \sim_{k} c$. These relations hold for any $k \notin \{i,a,c\}$. If it is possible to choose $k \neq b$, then we conclude that $b \in J_{k}$ and $a \sim_{k} b \sim_{k} c$ as desired. 
Otherwise, we may suppose that the choices $k=0$ and $k=n$ are impossible, so that
$a = 0$ and $c = n$. Then $a < i < c$, so that $i \in J_{b}$ and $a \sim_{b} i \sim_{b} c$.
Without loss of generality we may suppose $b < i$. Then $a \sim_{c} i$, so that $b \in J_{c}$
and $a \sim_{c} b \sim_{c} i$ as desired. 

We now claim that $(J,\sim): \Delta^n \rightarrow \Idem^{+}$ is an extension of the original
map $\Lambda^n_i \rightarrow \Idem^{+}$. In other words, we claim that for $k \neq i$, 
$J_{k} = J \cap \{0, \ldots, k-1, k+1, \ldots, n\}$ and $\sim_{k}$ is the restriction of $\sim$ to
$J_{k}$. The first claim is obvious. For the second, let us suppose that $a,b \in J_{k}$ and
$a \sim b$. We wish to prove that $a \sim_{k} b$. It will suffice to prove that $a \sim_{j} b$ for any
$j \notin \{i, a, b\}$. Since $a \sim b$, either such a $j$ exists, or $a \neq b \neq i \neq a$
and there exists $c \in J_{a} \cap J_{b}$ such that $a \sim_{b} c \sim_{a} b$. If there exists
$j \notin \{a,b,c,i\}$, then we conclude that $a \sim_{j} c \sim_{j} b$ and hence $a \sim_{j} b$
by transitivity. Otherwise, we conclude that $c = k \neq i$ and that $0,n \in \{a,b,c\}$. 
Without loss of generality, $i < c$; thus $0 \in \{a,b\}$ and we may suppose without loss of generality that $a < i$. Since $a \sim_{b} c$, we conclude that $i \in J_{b}$ 
and $a \sim_{b} i \sim_{b} c$. Consequently,
$i \in J_{a}$ and $i \sim_{a} c \sim_{a} b$, so that $i \sim_{a} b$ by transitivity
and therefore $i \sim_{c} b$. We now have $a \sim_{c} i \sim_{c} b$ so that $a \sim_{c} b$
as desired.
\end{proof}

\begin{remark}
It is clear that $\Idem \subseteq \Idem^{+}$ is the full simplicial subset spanned by the vertex $x$, and therefore an $\infty$-category as well.
\end{remark}

According to Corollary \ref{gurgh}, every weak retraction diagram
$$ 
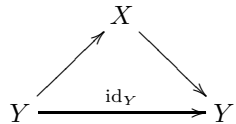
\xymatrix{ & X \ar[dr] & \\
Y \ar[ur] \ar[rr]^{\id_{Y}} & & Y }$$
in an $\infty$-category $\calC$ can be extended to a strong retraction diagram $F: \Idem^{+} \rightarrow \calC$, which restricts to give an idempotent in $\calC$. Our next goal is to show that
$F$ is canonically determined by the restriction $F| \Idem$.

Our next result expresses the idea that if an idempotent in $\calC$ arises in this manner, then $F$ is essentially unique.

\begin{lemma}\label{streaka}
The $\infty$-category $\Idem$ is weakly contractible.
\end{lemma}

\begin{proof}
An explicit computation shows that the topological space $|\Idem|$ is connected, simply connected, and has vanishing homology in degrees greater than zero. (Alternatively, we can deduce this from Proposition \ref{slanger} below.) 
\end{proof}

\begin{lemma}\label{linkdink}
The inclusion $\Idem \subseteq \Idem^{+}$ is a cofinal map of simplicial sets.
\end{lemma}

\begin{proof}
According to Theorem \ref{hollowtt}, it will suffice to prove that the simplicial sets
$\Idem_{x/}$ and $\Idem_{y/}$ are weakly contractible. The simplicial set $\Idem_{x/}$ is an $\infty$-category with an initial object, and therefore weakly contractible. The projection
$\Idem_{y/} \rightarrow \Idem$ is an isomorphism, and $\Idem$ is weakly contractible by Lemma \ref{streaka}.
\end{proof}

\begin{proposition}\label{autokan}
Let $\calC$ be an $\infty$-category, and let $F: \Idem^{+} \rightarrow \calC$ be a strong retraction diagram. Then $F$ is a left Kan extension of $F| \Idem$.
\end{proposition}

\begin{remark}
Passing to opposite $\infty$-categories, it follows that a strong retraction diagram $F: \Idem^{+} \rightarrow \calC$ is also a {\em right} Kan extension of $F|\Idem$.
\end{remark}

\begin{proof}
We must show that the induced map
$$ (\Idem_{/y})^{\triangleright} \rightarrow (\Idem^{+}_{/y})^{\triangleright} \stackrel{G}{\rightarrow} \Idem^{+}
\stackrel{F}{\rightarrow} \calC$$
is a colimit diagram. Consider the commutative diagram
$$ \xymatrix{ \Idem_{/y} \ar[r] \ar[d] & \Idem^{+}_{/y} \ar[d] \\
\Idem \ar[r] & \Idem^{+}. }$$
The lower horizontal map is cofinal by Lemma \ref{linkdink}, and the vertical maps are isomorphisms: therefore the upper horizontal map is also cofinal. Consequently, it will suffice to prove that $F \circ G$ is a colimit diagram, which is obvious.
\end{proof}

We will say that an idempotent $F: \Idem \rightarrow \calC$ in an $\infty$-category
$\calC$ is {\it effective} if it extends to a map $\Idem^{+} \rightarrow \calC$. According to Lemma \ref{kan2}, $F$ is effective if and only if it has a colimit in $\calC$. We will say that $\calC$ is {\it idempotent complete} if every idempotent in $\calC$ is effective.\index{gen}{idempotent complete}\index{gen}{idempotent!effective}

\begin{corollary}
Let $\calC$ be an $\infty$-category, and let $\calD \subseteq \Fun(\Idem,\calC)$ be the full subcategory spanned by the effective idempotents in $\calC$. The restriction map
$\Fun(\Idem^{+},\calC) \rightarrow \calD$ is a trivial fibration. In particular, if
$\calC$ is idempotent complete, then we have a diagram
$$ \Fun(\Ret,\calC) \leftarrow \Fun(\Idem^{+}, \calC) \rightarrow \Fun(\Idem,\calC)$$
of trivial fibrations.
\end{corollary}

\begin{proof}
Combine Proposition \ref{autokan} with Proposition \ref{lklk}.
\end{proof}

By definition, an $\infty$-category $\calC$ is idempotent complete if and only if every idempotent
$\Idem \rightarrow \calC$ has a colimit. In particular, if $\calC$ admits all small colimits, then it is idempotent complete. As we noted above, this is not necessarily true if $\calC$ admits only finite colimits. However, it turns out that filtered colimits do suffice: this assertion is not entirely obvious, since the $\infty$-category $\Idem$ itself is not filtered.

\begin{proposition}\label{slanger}
Let $A$ be a linearly ordered set with no largest element. Then there exists a cofinal map
$p: \Nerve(A) \rightarrow \Idem$.
\end{proposition}

\begin{proof}
Let $p: \Nerve(A) \rightarrow \Idem$ be the unique map which carries nondegenerate simplices to nondegenerate simplices. Explicitly, this map carries a simplex $\Delta^{J} \rightarrow \Nerve(A)$ corresponding to a map $s: J \rightarrow A$ of linearly ordered sets to the equivalence relation
$( i \sim j) \Leftrightarrow ( s(i) = s(j) )$. We claim that $p$ is cofinal. According to Theorem \ref{hollowtt}, it will suffice to show that the fiber product $\Nerve(A) \times_{\Idem} \Idem_{x/}$ is weakly contractible. We observe that $\Nerve(A) \times_{ \Idem} \Idem_{x} \simeq \Nerve(A')$, where
$A'$ denote the set $A \times \{0,1\}$ equipped with the partial ordering
$$ (\alpha, i) < (\alpha', j) \Leftrightarrow ( j = 1) \wedge ( \alpha < \alpha' ).$$

For each $\alpha \in A$, let $A_{< \alpha} = \{ \alpha' \in A: \alpha' < \alpha \}$ and let
$$A'_{\alpha} = \{ (\alpha',i) \in A' : (\alpha' < \alpha) \vee ( (\alpha',i) = (\alpha,1) ) \}.$$
By hypothesis, we can write
$A$ as a filtered union $\bigcup_{\alpha \in A} A_{< \alpha}$. It therefore suffices to prove
that for each $\alpha \in A$, the map
$$f: \Nerve(A_{< \alpha}) \times_{\Idem} \Idem_{x/} \rightarrow \Nerve(A) \times_{\Idem} \Idem_{x/}$$
has a nullhomotopic geometric realization $|f|$. But this map factors through
$\Nerve(A'_{\alpha})$, and $|\Nerve(A'_{\alpha})|$ is contractible because $A'_{\alpha}$ has a largest element.
\end{proof}

\begin{corollary}\label{swwe}
Let $\kappa$ be a regular cardinal, and suppose that $\calC$ is an $\infty$-category which admits $\kappa$-filtered colimits. Then $\calC$ is idempotent complete.
\end{corollary}

\begin{proof}
Apply Proposition \ref{slanger} to the linearly ordered set consisting of all ordinals 
less than $\kappa$ (and observe that this linearly ordered set is $\kappa$-filtered).
\end{proof}

\chapter{Presentable and Accessible $\infty$-Categories}\label{chap5}

\setcounter{theorem}{0}
\setcounter{subsection}{0}

Many categories which arise naturally, such as the category $\calA$ of abelian groups, are
{\em large}: they have a proper class of objects, even when the objects are considered only up to isomorphism. However, though $\calA$ itself is large, it is in some sense determined by the much smaller category $\calA_0$ of finitely generated abelian groups: $\calA$ is naturally equivalent to the category of $\Ind$-objects of $\calA_0$. This remark carries more than simply philosophical significance. When properly exploited, it can be used to prove statements such as the following:

\begin{proposition}\label{usepresent}
Let $F: \calA \rightarrow \Set$ be a contravariant functor from $\calA$ to the category of sets. Then $F$ is representable by an object of $\calA$ if and only if it carries colimits in $\calA$ to limits in $\Set$.
\end{proposition}

Proposition \ref{usepresent} is valid not only for the category $\calA$ of abelian groups, but for any {\it presentable} category: that is, any category which possess all (small) colimits and satisfies mild set-theoretic assumptions (such categories are referred to as {\it locally presentable} in \cite {adamek}). Our goal in this chapter is to develop an $\infty$-categorical generalization of the theory of presentable categories, and to obtain higher-categorical analogues of Proposition
\ref{usepresent} and related results (such as the adjoint functor theorem).

The most basic example of a presentable $\infty$-category is the $\infty$-category $\SSet$
of spaces. More generally, we can define an $\infty$-category $\calP(\calC)$ of {\em presheaves}
(of spaces) on an arbitrary small $\infty$-category $\calC$. We will study the properties of $\calP(\calC)$ in \S \ref{c5s1}; in particular, we will see that there exists a Yoneda embedding $j: \calC \rightarrow \calP(\calC)$ which is fully faithful, just as in ordinary category theory. Moreover, we 
give a characterization of $\calP(\calC)$ in terms of $\calC$: it is, in some sense, freely generated by the essential image of $j$ under (small) colimits.

The presheaf $\infty$-categories $\calP(\calC)$ are all presentable. Conversely, any presentable $\infty$-category can be obtained as a {\em localization} of some presheaf $\infty$-category $\calP(\calC)$ (Proposition \ref{pretop}). To make sense of this statement, we need a theory of localizations of $\infty$-categories. We will develop such a theory in \S \ref{c5s2}, as part of a more general theory of adjoint functors between $\infty$-categories.

In \S \ref{c5s3} we will introduce, for every small $\infty$-category $\calC$, an $\infty$-category $\Ind(\calC)$ of {\it $\Ind$-objects of $\calC$}. Roughly speaking, this is an $\infty$-category which is obtained from $\calC$ by freely adjoining colimits for all {\em filtered} diagrams. It is characterized up to equivalence by the fact that $\Ind(\calC)$ contains a full subcategory equivalent to $\calC$, which generates $\Ind(\calC)$ under filtered colimits and consists of {\em compact} objects. 

The construction of $\Ind$-categories will be applied in \S \ref{c5s5} to the study of {\em accessible} $\infty$-categories. Roughly speaking, an $\infty$-category $\calC$ is {\it accessible} if it is generated under (sufficiently) filtered colimits by a small subcategory $\calC^{0} \subseteq \calC$. We will prove that the class of accessible $\infty$-categories is stable under
a various categorical constructions. Results of this type will play an important technical role later this book: they generally allow us to dispense with the set-theoretic aspects
of an argument (such as cardinality estimation), and to focus instead on the more conceptual 
aspects.

We will say that an $\infty$-category $\calC$ is {\it presentable} if $\calC$ is accessible and admits (small) colimits.\index{gen}{presentable!$\infty$-category}\index{gen}{$\infty$-category!presentable} In \S \ref{c5s6}, we will describe the theory of presentable $\infty$-categories in detail. In particular, we will generalize Proposition \ref{usepresent} to the $\infty$-categorical setting, and prove an analogue of the adjoint functor theorem. 
We will also study localizations of presentable $\infty$-categories, following ideas of Bousfield.
The theory of presentable $\infty$-categories will play a vital role in the study of $\infty$-topoi, which is the subject of the next chapter.

\section{$\infty$-Categories of Presheaves}\label{c5s1}

\setcounter{theorem}{0}

The category of sets plays a central role in classical category theory. The primary reason for this is Yoneda's lemma, which asserts that for any category $\calC$, the ``Yoneda embedding''
$$j: \calC \rightarrow \Set^{\calC^{op}}$$
$$ C \mapsto \Hom_{\calC}(\bigdot, C)$$\index{gen}{Yoneda embedding!classical}
is fully faithful. Consequently, objects in $\calC$ can be thought of as a kind of ``generalized sets'', and various questions about the category $\calC$ can be reduced to questions about the category of sets.

If $\calC$ is an $\infty$-category, then the mapping {\em sets} of the above discussion should be replaced by mapping {\em spaces}. Consequently, one should expect the Yoneda embedding to take values in presheaves of {\em spaces}, rather than presheaves of sets. To formalize this, 
we introduce the following notation:

\begin{definition}\index{not}{PcalC@$\calP(\calC)$}
Let $S$ be a simplicial set. We let $\calP(S)$ denote the simplicial set
$\Fun(S^{op}, \SSet)$; here $\SSet$ denotes the $\infty$-category of spaces defined in \S \ref{introducingspaces}. We will refer
to $\calP(S)$ as the {\it $\infty$-category of presheaves on $S$}.\index{gen}{presheaf}
\end{definition}

\begin{remark}
More generally, for any $\infty$-category $\calC$, we might refer to
$\Fun( S^{op}, \calC)$ as the {\it $\infty$-category of $\calC$-valued presheaves on $S$}. 
Unless otherwise specified, the word ``presheaf'' will always refer to a $\SSet$-valued presheaf.
This is somewhat nonstandard terminology: one usually understands
the term ``presheaf'' to refer to a presheaf of sets, rather than
a presheaf of spaces. The shift in terminology is justified by the
fact that the important role of $\Set$ in ordinary category theory
is taken on by $\SSet$ in the $\infty$-categorical setting.
\end{remark}

Our goal in this section is to establish the basic properties of $\calP(S)$. We begin in \S \ref{presheaf3} by reviewing two other possible definitions of $\calP(S)$: one via the theory of right fibrations over $S$, another via simplicial presheaves on the category $\sCoNerve[S]$. Using the ``straightening'' results of \S \ref{fullun} and \S \ref{quasilimit4}, we will show that all three of these definitions are equivalent.

The presheaf $\infty$-categories $\calP(S)$ are examples of {\em presentable} $\infty$-categories (see \S \ref{c5s6}). In particular, each $\calP(S)$ admits small limits and colimits. We will give a proof of this assertion in \S \ref{presheaf2}, by reducing to the case where $S$ is a point.

The main question regarding the $\infty$-category $\calP(S)$ is how it relates to the original simplicial set $S$. In \S \ref{presheaf1} we will construct a map $j: S \rightarrow \calP(S)$, which is an $\infty$-categorical analogue of the usual Yoneda embedding. Just as in classical category theory, the Yoneda embedding is fully faithful. In particular, we note that any $\infty$-category $\calC$ can be embedded in a larger $\infty$-category which admits limits and colimits; this observation allows us to construct an {\it idempotent completion} of $\calC$, which we will study in \S \ref{surot}.

In \S \ref{presheaf4}, we will characterize the $\infty$-category $\calP(S)$ in terms of the Yoneda embedding $j: S \rightarrow \calP(S)$. Roughly speaking, we will show that $\calP(S)$ is freely generated by $S$ under colimits (Theorem \ref{charpresheaf}). In particular, if $\calC$ is a category which admits colimits, then any diagram $f: S \rightarrow \calC$ extends uniquely (up to homotopy) to a functor $F: \calP(S) \rightarrow \calC$. In \S \ref{completecomp}, we will give a criterion for determining whether or not $F$ is an equivalence.

\subsection{Other Models for $\calP(S)$}\label{presheaf3}

Let $S$ be a simplicial set. We have defined the $\infty$-category $\calP(S)$ of presheaves on $S$ to be the mapping space $\Fun(S^{op}, \SSet)$. However, there are several equivalent models which would serve equally well; we discuss two of them in this section.

Let $\calP'_{\Delta}(S)$ denote the full subcategory of
$(\sSet)_{/S}$ spanned by the right fibrations $X \rightarrow S$. We define $\calP'(S)$ to be the simplicial nerve $\sNerve(\calP'_{\Delta}(S))$. Because $\calP'_{\Delta}(S)$ is a fibrant simplicial category, $\calP'(S)$ is an $\infty$-category. We will see in a moment that $\calP'(S)$ is (naturally) equivalent to $\calP(S)$. In order to do this, we need to introduce a third model.

Let $\phi: \sCoNerve[S]^{op} \rightarrow \calC$ be an equivalence of simplicial categories.
Let $\Set_{\Delta}^{\calC}$ denote the category of simplicial functors $\calC \rightarrow \sSet$ (which we may view as simplicial presheaves on $\calC^{op}$). We regard $\Set_{\Delta}^{\calC}$ as endowed with the {\em projective} model structure defined in \S \ref{quasilimit3}. With respect to this structure,
$\Set_{\Delta}^{\calC}$ is a simplicial model category; we let $\calP''_{\Delta}(\phi) =  (\Set_{\Delta}^{\calC})^{\degree}$ denote the full simplicial subcategory consisting of fibrant-cofibrant objects, and we define $\calP''(\phi)$ to be the simplicial nerve $\sNerve(\calP''_{\Delta}(\phi))$.

We are now ready to describe the relationship between these different models:

\begin{proposition}\label{othermod}
Let $S$ be a simplicial set, and let $\phi: \sCoNerve[S]^{op} \rightarrow \calC$ be
an equivalence of simplicial categories. Then there are $($canonical$)$ equivalences
of $\infty$-categories
$$ \calP(S) \stackrel{f}{\leftarrow} \calP''(\phi) \stackrel{g}{\rightarrow} \calP'(S).$$
\end{proposition}

\begin{proof}
The map $f$ was constructed in Proposition \ref{gumby444}; it therefore suffices to give a construction of $g$. 

Recall that the category $(\sSet)_{/S}$ of simplicial sets over $S$ is endowed model structure, the {\it contravariant} model structure defined in \S \ref{contrasec}. Moreover, this model structure is simplicial (Proposition \ref{natsim}) and the fibrant objects are precisely the right fibrations over $S$ (Corollary \ref{usewhere1}). Thus, we may identify
$\calP'_{\Delta}(S)$ with the simplicial category $(\sSet)_{/S}^{\degree}$ of fibrant-cofibrant objects of $(\sSet)_{/S}$. 

According to Theorem \ref{struns}, the straightening and unstraightening functors
$(\St_{\phi}, \Un_{\phi})$ determine a Quillen equivalence between
$(\sSet)^{\calC}$ and $(\sSet)_{/S}$. Moreover, for any $X \in (\sSet)_{/S}$ and any simplicial set $K$, there is a natural chain of equivalences
$$ \St_{\phi} (X \times K) \rightarrow (\St_{\phi} X) \otimes |K|_{Q^{\bigdot}} \rightarrow
(\St_{\phi} X) \otimes K.$$
(The fact that the first map is an equivalence follows easily from Proposition \ref{spek3}.)
It follows from Proposition \ref{weakcompatequiv} that $\Un_{\phi}$ is endowed with the structure of a simplicial functor, and induces an equivalence of simplicial categories
$$ (\Set_{\Delta}^{\calC})^{\degree} \rightarrow  (\sSet)_{/S}^{\degree}.$$
We obtain the desired equivalence $g$ by passing to the simplicial nerve.
\end{proof}

\subsection{Colimits in $\infty$-Categories of Functors}\label{presheaf2}

Let $S$ be an arbitrary simplicial set. Our goal in this section is to prove that the $\infty$-category $\calP(S)$ of presheaves on $S$ admits small limits and colimits. There are (at least) three approaches to proving this:

\begin{itemize}
\item[$(1)$] According to Proposition \ref{othermod}, we may identify $\calP(S)$ with the
$\infty$-category underlying the simplicial model category $\Set_{\Delta}^{\sCoNerve[S]^{op}}$. We can then deduce the existence of limits and colimits in $\calP(S)$ by invoking Corollary \ref{limitsinmodel}.

\item[$(2)$] Since the $\infty$-category $\SSet$ classifies left fibrations, the $\infty$-category
$\calP(S)$ classifies left fibrations over $S^{op}$: in other words, homotopy classes of maps
$K \rightarrow \calP(S)$ can be identified with equivalence classes of left fibrations
$X \rightarrow K \times S^{op}$. It is possible to generalize Proposition \ref{charspacecolimit} and
Corollary \ref{charspacelimit} to describe limits and colimits in $\calP(S)$ entirely in the language of left fibrations. The existence problem can then be solved by exhibiting explicit constructions of left fibrations.

\item[$(3)$] Applying either $(1)$ or $(2)$ in the case where $S$ is a point, we can deduce that
the $\infty$-category $\SSet \simeq \calP( \ast)$ admits limits and colimits. We can then attempt to deduce the same result for $\calP(S) = \Fun( S^{op}, \SSet)$ using a general result about (co)limits in functor categories (Proposition \ref{limiteval}). 
\end{itemize}

Although approach $(1)$ is probably the quickest, we will adopt approach $(3)$ because it gives additional information: our proof will show that the formation of limits and colimits in $\calP(S)$ are computed pointwise. The same proof will also apply to $\infty$-categories of $\calC$-valued presheaves in the case where $\calC$ is not necessarily the $\infty$-category $\SSet$ of spaces.

\begin{lemma}\label{topaz2}
Let $q: Y \rightarrow S$ be a Cartesian fibration of simplicial sets, and let
$\calC = \bHom_{S}(S,Y)$ denote the $\infty$-category of sections of $q$. Let
$p: S \rightarrow Y$ be an object of $\calC$ having the property that $p(s)$ is an initial object of the 
fiber $Y_{s}$ for each vertex $s$ of $S$. Then $p$ is an initial object of $\calC$.
\end{lemma}

\begin{proof}
By Proposition \ref{colimfam}, the map $Y^{p_S/} \rightarrow S$ is a Cartesian fibration. By hypothesis, for each vertex $s$ of $S$, the map $Y^{p_S/} \times_{S} \{s\}  \rightarrow  Y_s$
is a trivial fibration. It follows that the projection $Y^{p_S/} \rightarrow Y$ is an equivalence of Cartesian fibrations over $S$, and therefore a categorical equivalence; taking sections over $S$ we obtain another categorical equivalence
$$ \bHom_{S}( S, Y^{p_S/} ) \rightarrow \bHom_{S}(Y,S).$$
But this map is just the left fibration $j: \calC^{p/} \rightarrow \calC$; it follows that $j$ is a categorical equivalence. Applying Propostion \ref{apple1} to the diagram
$$ \xymatrix{ \calC^{p/} \ar[dr]^{j} \ar[rr]^{j} & & \calC \ar[dl]^{\id_{\calC}} \\
& \calC,}$$
we deduce that $j$ induces categorical equivalences $\calC_{p/} \times_{\calC} \{t\} \rightarrow \{t\}$ for each vertex $t$ of $Q$. Thus the fibers of $j$ are contractible Kan complexes, so that $j$ is a trivial fibration (by Lemma \ref{toothie}) and $p$ is an initial object of $\calC$, as desired.
\end{proof}

\begin{proposition}\label{limiteval}\index{gen}{colimit!in a functor category}\index{gen}{limit!in a functor category}
Let $K$ be a simplicial set, $q: X \rightarrow S$ a Cartesian fibration, and
$p: K \rightarrow \bHom_{S}(S,X)$ a diagram.
For each vertex $s$ of $S$, 
we let $p_s: K \rightarrow X_{s}$ be the induced map. Suppose,
furthermore, that each $p_s$ has a colimit in the $\infty$-category $X_{s}$. Then:

\begin{itemize}
\item[$(1)$] There exists a map $\overline{p}: K \diamond \Delta^0 \rightarrow \bHom_{S}(S,X)$
which extends $p$ and induces a colimit diagram $\overline{p}: K \diamond
\Delta^0 \rightarrow X_{s}$, for each vertex $s \in S$.
                                                                                                                                                                                                                                                                                                                                                                                                                                                                                                       \item[$(2)$] An arbitrary extension $\overline{p}: K \diamond \Delta^0 \rightarrow \bHom_{S}(S,X)$
of $p$ is a colimit for $p$ if and only if each $\overline{p}_s: K \diamond
\Delta^0 \rightarrow X_{s}$ is a colimit for $p_s$.
\end{itemize}
\end{proposition}

\begin{proof}
Choose a factorization $K \rightarrow K' \rightarrow \bHom_{S}(S,X)$ of $p$,
where $K \rightarrow K'$ is inner anodyne (and therefore a
categorical equivalence) and $K' \rightarrow \calC^S$ is an inner fibration (so that $K'$ is an $\infty$-category). The map $K \rightarrow K'$ is a categorical equivalence, and therefore cofinal. We are free to replace $K$ by $K'$, and may thereby assume 
that $K$ is an $\infty$-category.

We apply Proposition \ref{familycolimit} to the Cartesian
fibration $X \rightarrow S$ and the diagram $p_S: K
\times S \rightarrow X$ determined by the map $p$. We
deduce that there exists a map $$\overline{p}_S: (K \times S) \diamond_S S = (K
\diamond \Delta^0) \times S \rightarrow X$$ having the property that
its restriction to the fiber over each $s \in S$ is a colimit of
$p_s$; this proves $(1)$.

The ``if'' direction of $(2)$ follows immediately from Lemma \ref{topaz2}. The ``only if'' follows
from $(1)$ and the fact that colimits, when they exist, are unique up to equivalence.
\end{proof}

\begin{corollary}
Let $K$ and $S$ be simplicial sets, and let $\calC$ be an $\infty$-category which admits $K$-indexed colimits. Then:
\begin{itemize}
\item[$(1)$] The $\infty$-category $\Fun(S, \calC)$ admits $K$-indexed colimits.
\item[$(2)$] A map $K^{\triangleright} \rightarrow \Fun(S,\calC)$ is a colimit diagram if and only if,
for each vertex $s \in S$, the induced map $K^{\triangleright} \rightarrow \calC$ is a colimit diagram.
\end{itemize}
\end{corollary}

\begin{proof}
Apply Proposition \ref{limiteval} to the projection $\calC \times S \rightarrow S$.
\end{proof}

\begin{corollary}\label{storum}
Let $S$ be a simplicial set. The $\infty$-category $\calP(S)$ of presheaves on $S$ admits all small limits and colimits.
\end{corollary}

\subsection{Yoneda's Lemma}\label{presheaf1}

In this section, we will construct the $\infty$-categorical analogue of the Yoneda embedding, and prove that it is fully faithful. We begin with a somewhat naive approach, based on the formalism of simplicial categories. We note that an analogoue of Yoneda's Lemma is valid in enriched category
theory (with the usual proof). Namely, suppose that $\calC$ is a category enriched over another category $\calE$. Then there is an ``enriched Yoneda embedding"
$$ i: \calC \rightarrow \calE^{\calC^{op}}$$
$$ X \mapsto \bHom_{\calC}( \bigdot, X).$$\index{gen}{Yoneda embedding!simplicial}

Consequently, for any simplicial
category $\calC$, one obtains a fully faithful embedding $i$ of
$\calC$ into the simplicial category $\bHom_{\sCat}(\calC^{op}, \sSet)$ of
simplicial functors from $\calC^{op}$ into $\sSet$. In fact, $i$
is fully faithful in the strong sense that it induces {\em
isomorphisms} of simplicial sets$$ \bHom_{\calC}(X,Y) \rightarrow
\bHom_{\Set_{\Delta}^{\calC^{op}}}(i(X), i(Y)),$$ rather than merely
weak homotopy equivalences. Unfortunately, this assertion does not
necessarily have any $\infty$-categorical content, because the 
simplicial category $\Set_{\Delta}^{\calC^{op}}$ does not generally
represent the correct $\infty$-category of functors from
$\calC^{op}$ to $\sSet$.

Let us describe an analogous construction in the setting of $\infty$-categories. Let $K$ be a simplicial set, and let
$\calC = \sCoNerve[K]$. Then $\calC$ is a simplicial category, so
$$ (X,Y) \mapsto \Sing|\Hom_{\calC}(X,Y)|$$
determines a simplicial functor from
$ \calC^{op} \times \calC$ to the category $\Kan$.
The functor $\sCoNerve$ does not commute with products, but there
exists a natural map $\sCoNerve[K^{op} \times K] \rightarrow
\calC^{op} \times \calC$. Composing with this map, we obtain a map
of simplicial sets
$$ \sCoNerve[K^{op} \times K] \rightarrow \Kan.$$
Passing to the adjoint, we obtain a map of simplicial sets
$K^{op} \times K \rightarrow \SSet,$ which we can identify
with $$j: K \rightarrow
\Fun(K^{op}, \SSet) = \calP(K).$$
We shall refer to $j$ (or, more generally, any map equivalent to
$j$) as the {\it Yoneda embedding}.\index{gen}{Yoneda embedding}

\begin{proposition}[$\infty$-Categorical Yoneda Lemma]\label{fulfaith}\index{gen}{Yoneda's Lemma}
Let $K$ be a simplicial set. Then the Yoneda embedding $j: K
\rightarrow \calP(K)$ is fully faithful.
\end{proposition}

\begin{proof}
Let $\calC' = \Sing | \sCoNerve[K^{op}] |$ be the ``fibrant
replacement'' for $\calC=\sCoNerve[K^{op}]$. We endow
$\Set_{\Delta}^{\calC'}$ with the {\em projective} model structure described in \S \ref{quasilimit3}.

We note that the Yoneda embedding factors as a composition
$$ K \stackrel{j'}{\rightarrow} \sNerve( (\Set_{\Delta}^{\calC'})^{\degree} )
\stackrel{j''}{\rightarrow} \Fun( K^{op}, \SSet),$$ where $j''$ is the map of Proposition \ref{gumby444} and consequently a categorical equivalence. It therefore suffices to
prove that $j'$ is fully faithful. For this, we need only show
that the adjoint map
$$ J: \sCoNerve[K] \rightarrow \Set_{\Delta}^{\calC'}.$$
is a fully-faithful functor between simplicial categories. We now observe that
$J$ is the composition of an equivalence $\sCoNerve[K] \rightarrow (\calC')^{op}$
with the (simplicial enriched) Yoneda embedding 
$(\calC')^{op} \rightarrow \Set_{\Delta}^{\calC'}$, which is fully faithful
in virtue of the classical (simplicially enriched) version of Yoneda's Lemma.
\end{proof}

We conclude by establishing another pleasant property of the Yoneda
embedding:

\begin{proposition}\label{yonedaprop}\index{gen}{Yoneda embedding!and limits}
Let $\calC$ be a small $\infty$-category, and $j: \calC \rightarrow \calP(\calC)$ the Yoneda embedding. Then $j$ preserves all limits which exist in $\calC$.
\end{proposition}

\begin{proof}
Let $p: K \rightarrow \calC$ be a diagram having a limit in $\calC$. We
wish to show that $j$ carries any limit for $p$ to a limit of
$j \circ p$. Choose a category $\calI$ and a cofinal map
$N(\calI^{op}) \rightarrow K^{op}$ (the existence of which is guaranteed by Proposition \ref{cofinalcategories}) Replacing $K$ by $\Nerve(\calI)$, we may suppose that $K$ is
the nerve of a category. Let $\overline{p}: \Nerve(\calI)^{\triangleleft} \rightarrow
\calC$ be a limit for $p$.

We recall the definition of the Yoneda embedding. It involves the
choice of an equivalence $\sCoNerve[\calC] \rightarrow
\calD$, where $\calD$ is a fibrant simplicial category. For
definiteness, we took $\calD$ to be $\Sing |\sCoNerve[\calC]|$. However, we could just as well
choose some other fibrant simplicial category $\calD'$ equivalent to $\sCoNerve[\calC]$ and
obtain a ``modified Yoneda embedding'' $j': \calC \rightarrow \calP(\calC)$; it is easy to see that
$j'$ and $j$ are equivalent functors, so it suffices to show that $j'$ preserves the limit of $p$.
Using Corollary \ref{strictify}, we may suppose that $\overline{p}$ is obtained
from a functor between simplicial categories
$\overline{q}: \{x\} \star \calI \rightarrow \calD$ by passing to the nerve. According to Theorem \ref{colimcomparee}, $\overline{q}$ is a homotopy limit of $q = \overline{q} | \calI$.
Consequently, for each object $Z \in \calD$, the induced functor
$$ \overline{q}_Z: I \mapsto \Hom_{\calD}(Z, \overline{q}(I))$$
is a homotopy limit of $q_Z = \overline{q}_Z|\calI$. Taking $Z$ to be the image of an object
$C$ of $\calC$, we deduce that
$$ \Nerve(\calI)^{\triangleleft} \rightarrow \calC \stackrel{j'}{\rightarrow} \calP(\calC) \rightarrow \SSet$$
is a limit for its restriction to $\Nerve(\calI)$, where the map on the right is given by ``evaluation at $C$''. Proposition \ref{limiteval} now implies that $j' \circ \overline{p}$ is a limit for $j' \circ p$, as desired.
\end{proof}

\subsection{Idempotent Completions}\label{surot}

Recall that an $\infty$-category $\calC$ is said to be {\it idempotent complete} if every
functor $\Idem \rightarrow \calC$ admits a colimit in $\calC$ (see \S \ref{retrus}).
If an $\infty$-category $\calC$ is not idempotent complete, then we can attempt to correct the situation by passing to a larger $\infty$-category.

\begin{definition}\index{gen}{idempotent completion}
Let $f: \calC \rightarrow \calD$ be a functor between $\infty$-categories. We will say that $f$ {\it exhibits $\calD$ as an idempotent completion of $\calC$} if $\calD$ is idempotent complete, $f$ is fully faithful, and every object of $\calD$ is a retract of $f(C)$, for some object $C \in \calC$.
\end{definition}

Our goal in this section is to show that $\infty$-category $\calC$ has an idempotent completion $\calD$, which is unique up to equivalence. The uniqueness is a consequence of Proposition \ref{charidemcomp}, proven below. The existence question is much easier to address.

\begin{proposition}\label{idmcoo}
Let $\calC$ be an $\infty$-category. Then $\calC$ admits an idempotent completion.
\end{proposition}

\begin{proof}
Enlarging the universe if necessary, we may suppose that $\calC$ is small. Let
$\calC'$ denote the full subcategory of $\calP(\calC)$ spanned by those objects which are retracts of objects which belong to the image of the Yoneda embedding $j: \calC \rightarrow \calP(\calC)$.
Then $\calC'$ is stable under retracts in $\calP(\calC)$. Since $\calP(\calC)$ admits all small colimits, Corollary \ref{swwe} implies that
$\calP(\calC)$ is idempotent complete. It follows that $\calC'$ is idempotent complete. Proposition \ref{fulfaith} implies that the Yoneda embedding $j: \calC \rightarrow \calC'$ is fully faithful, and therefore exhibits $\calC'$ as an idempotent completion of $\calC$.
\end{proof}

We now address the question of uniqueness for idempotent completions. First, we need a few preliminary results.

\begin{lemma}\label{sweeble}
Let $\calC$ be an $\infty$-category which is idempotent complete, and let $p: K \rightarrow \calC$ be a diagram. Then $\calC_{/p}$ and $\calC_{p/}$ are also idempotent complete.
\end{lemma}

\begin{proof}
By symmetry, it will suffice to prove that $\calC_{/p}$ is idempotent complete. Let
$q: \calC_{/p} \rightarrow \calC$ be the associated right fibration, and let
$F: \Idem \rightarrow \calC_{/p}$ be an idempotent. We will show that $F$ has a limit.
Since $\calC$ is idempotent complete,
$q \circ F$ has a limit $\overline{q \circ F}: \Idem^{\triangleleft} \rightarrow \calC$. Consider the lifting problem
$$ \xymatrix{ \Idem \ar@{^{(}->}[d] \ar[r]^{F} & \calC_{/p} \ar[d]^{q} \\
\Idem^{\triangleleft} \ar[r]^{ \overline{q \circ F} } \ar@{-->}[ur]^{\overline{F}} & \calC. }$$
The right vertical map is a right fibration, and the left vertical map is right anodyne (Lemma \ref{chotle2}), so that there exists a dotted arrow $\overline{F}$ as indicated. Using Proposition \ref{goeselse}, we deduce that $\overline{F}$ is a limit of $F$. 
\end{proof}

\begin{lemma}\label{sweerum}
Let $f: \calC \rightarrow \calD$ be a functor between $\infty$-categories which exhibits $\calD$ as an idempotent completion of $\calC$, and let $p: K \rightarrow \calD$ be a diagram. Then the induced map $f_{/p}: \calC \times_{\calD} \calD_{/p} \rightarrow \calD_{/p}$ exhibits $\calD_{/p}$ as an idempotent completion of $\calC \times_{\calD} \calD_{/p}$.
\end{lemma}

\begin{proof}
Lemma \ref{sweeble} asserts that $\calD_{/p}$ is idempotent complete. We must show that every object $\overline{D} \in \calD_{/p}$ is a retract of $f_{/p}( \overline{C})$, for some
$\overline{C} \in \calC \times_{\calD} \calD_{/p}$. Let $q: \calD_{/p} \rightarrow \calD$ be the projection, and let $D = q( \overline{D})$. Since $f$ exhibits $\calD$ as an idempotent completion of $\calC$, there is a diagram
$$ \xymatrix{ & f(C) \ar[dr] & \\
D' \ar[rr]^{g} \ar[ur] & & D }$$
in $\calD$, where $g$ is an equivalence. Since $q$ is a right fibration, we can lift this to a diagram
$$ \xymatrix{ & \overline{f(C)} \ar[dr] & \\
\overline{D}' \ar[rr]^{\overline{g}} \ar[ur] & & \overline{D} }$$
in $\calD_{/q}$. Since $\overline{g}$ is $q$-Cartesian and $g$ is an equivalence,
$\overline{g}$ is ann equivalence. It follows that $\overline{D}$ is a retract of
$\overline{f(C)}$. By construction, $\overline{f(C)} = f_{/p}(\overline{C})$ for 
an appropriately chosen object $\overline{C} \in \calC \times_{\calD} \calD_{/p}$.  
\end{proof}

\begin{lemma}\label{wequivvv}
Let $f: \calC \rightarrow \calD$ be a functor between $\infty$-categories which exhibits $\calD$ as an idempotent completion of $\calC$. Suppose that $\calD$ has an initial object $\emptyset$. Then
$\calC$ is weakly contractible as a simplicial set.
\end{lemma}

\begin{proof}
Without loss of generality, we may suppose that $\calC$ is a full subcategory of $\calD$ and that $f$ is the inclusion. Since $f$ exhibits $\calD$ as an idempotent completion of $\calC$, the initial object $\emptyset$ of $\calD$ admits a map $f: C \rightarrow \emptyset$, where $C \in \calC$.
The $\infty$-category $\calC_{C/}$ has an initial object, and
is therefore weakly contractible. Since composition
$$ \calC_{f/} \rightarrow \calC_{C/} \rightarrow \calC$$
is both a weak homotopy equivalence (in fact, a trivial fibration) and weakly nullhomotopic, we
conclude that $\calC$ is weakly contractible.
\end{proof}

\begin{lemma}\label{sweetrum}
Let $f: \calC \rightarrow \calD$ be a functor between $\infty$-categories which exhibits $\calD$ as an idempotent completion of $\calC$. Then $f$ is cofinal.
\end{lemma}

\begin{proof}
According to Theorem \ref{hollowtt}, it suffices to prove that for every object $D \in \calD$, 
simplicial set $\calC \times_{\calD} \calD_{D/}$ is weakly contractible. Lemma \ref{sweerum} asserts that $f_{D/}$ is also an idempotent completion, and Lemma \ref{wequivvv} completes the proof.
\end{proof}

\begin{lemma}\label{honeybeen}
Let $F: \calC \rightarrow \calD$ be a functor between $\infty$-categories, and let
$\calC^{0} \subseteq \calC$ be a full subcategory such that the inclusion exhibits
$\calC$ as an idempotent completion of $\calC^{0}$. Then $F$ is a left Kan extension
of $F|\calC^{0}$. 
\end{lemma}

\begin{proof}
We must show that for every object $C \in \calC$, the composite map
$$ (\calC^0_{/C})^{\triangleright} \rightarrow (\calC_{/C})^{\triangleright}
\stackrel{G}{\rightarrow} \calC \stackrel{F}{\rightarrow} \calD$$
is a colimit diagram in $\calD$. Lemma \ref{sweerum} guarantees that 
$\calC^0_{/C} \subseteq \calC_{/C}$ is an idempotent completion, and therefore cofinal by Lemma \ref{sweetrum}. Consequently, it suffices prove that $F \circ G$ is a colimit diagram, which is obvious.
\end{proof}

\begin{lemma}\label{beenhoney}
Let $\calC$ and $\calD$ be $\infty$-categories which are idempotent complete, and let
$\calC^{0} \subseteq \calC$ be a full subcategory such that the inclusion exhibits
$\calC$ as an idempotent completion of $\calC^{0}$. Then any functor
$F_0: \calC^{0} \rightarrow \calD$ has an extension $F: \calC \rightarrow \calD$.
\end{lemma}

\begin{proof}
We will suppose that the $\infty$-categories $\calC$ and $\calD$ are small. Let $\calP(\calD)$ be the $\infty$-category of presheaves on $\calD$ (see \S \ref{c5s1}), $j: \calD \rightarrow \calP(\calD)$ the Yoneda embedding, and $\calD'$ the essential image of $j$. According to Proposition \ref{princex}, it will suffice to prove that $j \circ F_0$ can be extended to a functor $F': \calC \rightarrow \calD'$.
Since $\calP(\calD)$ admits small colimits, we can choose $F': \calC \rightarrow \calP(\calD)$
to be a left Kan extension of $j \circ F_0$. Every object of $\calC$ is a retract of an object of $\calC^{0}$, so that every object in the essential image of $F'$ is a retract of the Yoneda image of an object of $\calD$. Since $\calD$ is idempotent complete, it follows that the $F'$ factors through $\calD'$.
\end{proof}

\begin{proposition}\label{charidemcomp}\index{gen}{idempotent completion!universal property of}
Let $f: \calC \rightarrow \calD$ be a functor which exhibits $\calD$ as the idempotent completion of $\calC$, and let $\calE$ be an $\infty$-category which is idempotent complete. Then composition with $f$ induces an equivalence of $\infty$-categories
 $f^{\ast}: \Fun(\calD, \calE) \rightarrow \Fun(\calC, \calE)$.
\end{proposition}

\begin{proof}
Without loss of generality, we may suppose that $f$ is the inclusion of a full subcategory. In this case, we combine Lemmas \ref{honeybeen}, \ref{beenhoney}, and Proposition \ref{lklk} to deduce that $f^{\ast}$ is a trivial fibration.
\end{proof}

\begin{remark}
Let $\calC$ be a small $\infty$-category, and let $f: \calC \rightarrow \calC'$ be an idempotent completion of $\calC$. The proof of Proposition \ref{idmcoo} shows that $\calC'$ is equivalent to a full subcategory of $\calP(\calC)$, and therefore locally small (see \S \ref{locbrend}). Moreover, every object of $\h{\calC'}$ 
is the image of some retraction map in $\h{\calC}$; it follows that the set of equivalence classes of objects in $\calC'$ is bounded in size. It follows that $\calC'$ is essentially small.
\end{remark}

\subsection{The Universal Property of $\calP(S)$}\label{presheaf4}

Let $S$ be a (small) simplicial set. We have defined $\calP(S)$ to be the $\infty$-category of maps from $S^{op}$ into the $\infty$-category $\SSet$ of spaces. Informally, we may view $\calP(S)$ as the limit of a diagram in the $\infty$-bicategory of (large) $\infty$-categories: namely, the constant diagram carrying $S^{op}$ to $\SSet$. In more concrete terms, our definition of $\calP(S)$ leads immediately to a characaterization of $\calP(S)$ by a universal mapping property: for every $\infty$-category $\calC$, there is an equivalence of $\infty$-categories (in fact an isomorphism of simplicial sets)
$$ \Fun(\calC, \calP(S)) \simeq \Fun(\calC \times S^{op}, \SSet).$$
The goal of this section is to give a dual characterization of $\calP(S)$: it may also be viewed
as a {\em colimit} of copies of $\SSet$, indexed by $S$. However, this colimit needs to be understood in an appropriate $\infty$-bicategory of $\infty$-categories where the morphisms are given by {\em colimit preserving} functors. In other words, we will show that $\calP(S)$ is in some sense ``freely generated'' by $S$ under small colimits (Theorem \ref{charpresheaf}). First, we need to introduce a bit of notation.

\begin{notation}
Let $\calC$ be an $\infty$-category and $S$ a simplicial set. We will let
$\LFun( \calP(S), \calC)$ denote the full subcategory of
$\Fun( \calP(S), \calC)$ spanned by those functors $\calP(S) \rightarrow \calC$
which preserve small colimits. 

The motivation for this notation is as follows: in \S \ref{afunc5}, we will use the notation
$\LFun( \calD, \calC)$ to denote the full subcategory of $\Fun( \calD, \calC)$ spanned by those functors which are {\em left adjoints}. In \S \ref{aftt}, we will see that when $\calD = \calP(S)$ (or, more generally, when $\calD$ is presentable), then a functor $\calD \rightarrow \calC$ is a left adjoint if and only if it preserves small colimits (see Corollary \ref{adjointfunctor} and Remark \ref{afi}). 
\end{notation}

We wish to prove that if $\calC$ is an $\infty$-category which admits small colimits, then any map $S \rightarrow \calC$ extends in an essentially unique fashion to a colimit-preserving functor $\calP(S) \rightarrow \calC$. To prove this, we need a second characterization of the colimit-preserving functors $f: \calP(S) \rightarrow \calC$: they are precisely those functors which are left Kan extensions of their restriction to the essential image of the Yoneda embedding.

\begin{lemma}\label{repco}\index{gen}{corepresentable!functor}
Let $S$ be a small simplicial set, let $s$ be a vertex of $S$, let
$e: \calP(S) \rightarrow \SSet$ be the map given by evaluation at $s$, and let
$f: \calC \rightarrow \calP(S)$ be the associated left fibration (see \S \ref{universalfib}). Then
$f$ is corepresentable by the object $j(s) \in \calP(S)$, where $j: S \rightarrow \calP(S)$ denotes the Yoneda embedding.
\end{lemma}

\begin{proof}
Without loss of generality, we may suppose that $S$ is an $\infty$-category.
We make use of the equivalent model $\calP'(S)$ of \S \ref{presheaf3}. Observe that the functor
$f: \calP(S) \rightarrow \SSet$ is equivalent to $f': \calP'(S) \rightarrow \SSet$, where
$f'$ is the nerve of the simplicial functor $\calP'_{\Delta}(S) \rightarrow \Kan$ which
associates to each left fibration $Y \rightarrow S$ the fiber $Y_{s} = Y \times_{S} \{s\}$. 
Furthermore, under the equivalence of $\calP(S)$ with $\calP'(S)$, the object $j(s)$
corresponds to a left fibration $X(s) \rightarrow S$ which is corepresented by $s$. Then
$X(s)$ contains an initial object $x$ lying over $s$. The choice of $x$ determines a point 
$\eta \in \pi_0 f'(X(s))$. According to Proposition \ref{reppfunc}, to show that $X(s)$ corepresents $f'$, it suffices to show that for every left fibration $X \rightarrow S$, the map
$$ \bHom_{S}( X(s), Y) \rightarrow Y_{s}, $$
given by evaluation at $x$, is a homotopy equivalence of Kan complexes. 
We may rewrite the space on the right hand side as $\bHom_{S}( \{x\}, Y)$. According
to Proposition \ref{natsim}, the covariant model structure on $(\sSet)_{/S}$ is compatible with the simplicial structure. It therefore suffices to prove that the inclusion $i: \{ x\} \subseteq X(s)$ is a
covariant equivalence. But this is clear, since $i$ is the inclusion of an initial object and therefore left anodyne.
\end{proof}

\begin{lemma}\label{longwait0}
Let $S$ be a small simplicial set, and let $j: S \rightarrow \calP(S)$ denote the Yoneda embedding. Then $\id_{\calP(S)}$ is a left Kan extension of $j$ along itself.
\end{lemma}

\begin{proof}
Let $\calC \subseteq \calP(S)$ denote the essential image of $j$. According to Proposition \ref{fulfaith}, $j$ induces an equivalence $S \rightarrow \calC$. It therefore suffices to prove that
$\id_{\calP(S)}$ is a left Kan extension of its restriction to $\calC$. Let $X$ be an object of $\calP(S)$; we must show that the natural map
$$ \phi: \calC_{/X}^{\triangleright} \subseteq \calP(S)_{/X}^{\triangleright} \rightarrow \calP(S)$$
is a colimit diagram.

According to Proposition \ref{limiteval}, it will suffice to prove that for each vertex $s$ of $S$, the map
$$ \phi_{s}: \calC_{/X}^{\triangleright} \rightarrow \SSet$$ given by composing $\phi$ with the evaluation map is a colimit diagram in $\SSet$. Let 
$\calD \rightarrow \calC_{/X}^{\triangleright}$ be the pullback of the universal left fibration along $\phi_{s}$, and let $\calD^0 \subseteq \calD$ be the preimage in $\calD$ of
$\calC_{/X} \subseteq \calC_{/X}^{\triangleright}$. According to Proposition \ref{charspacecolimit}, it will suffice to prove that the inclusion $\calD^{0} \subseteq \calD$ is a weak homotopy equivalence of simplicial sets.

Let $C = j(s)$. Let $\calE = \calC_{/X}^{\triangleright} \times_{\calP(S)} \calP(S)_{C/}$,
let $\calE^{0} = \calC_{/X} \times_{\calC} \calC_{C/} \subseteq \calE$, and let
$\calE^{1} = \calC_{/X} \times_{\calC} \{ \id_{C} \} \subseteq \calE^{0}$.
Lemma \ref{repco} implies that the left fibrations
$$ \calD \rightarrow \calC_{/X}^{\triangleright} \leftarrow \calE$$ are
equivalent. It therefore suffices to show that the inclusion $\calE^{0} \subseteq \calE$ is a weak homotopy equivalence. To prove this, we observe that both $\calE$ and $\calE^{0}$ contain
$\calE^{1}$ as a deformation retract (that is, there is a retraction $r: \calE \rightarrow \calE^{1}$ and a homotopy $\calE \times \Delta^1 \rightarrow \calE$ from $r$ to $\id_{\calE}$, so that
the inclusion $\calE^{1} \subseteq \calE$ is a homotopy equivalence; the situation for
$\calE^0$ is similar).
\end{proof}

\begin{lemma}\label{natee}
Let $$\xymatrix{ A \ar[rr]^{f} \ar[dr] & & B \ar[dl]^{g} \\
& S & }$$
be a diagram of simplicial sets. The following conditions are equivalent:

\begin{itemize}
\item[$(1)$] The map $f$ is a covariant equivalence in $(\sSet)_{/S}$.

\item[$(2)$] For every diagram $p: S \rightarrow \calC$ taking values in an $\infty$-category $\calC$, and every limit $\overline{p \circ g}: B^{\triangleleft} \rightarrow \calC$ of
$p \circ g$, the composition $\overline{p \circ g} \circ f^{\triangleleft}: A^{\triangleleft} \rightarrow \calC$ is a limit diagram.

\item[$(3)$] For every diagram $p: S \rightarrow \SSet$ taking values in the $\infty$-category
$\SSet$ of spaces, and every limit $\overline{p \circ g}: B^{\triangleleft} \rightarrow \SSet$ of
$p \circ g$, the composition $\overline{p \circ g} \circ f^{\triangleleft}: A^{\triangleleft} \rightarrow \SSet$ is a limit diagram.
\end{itemize}
\end{lemma}

\begin{proof}
The equivalence of $(1)$ and $(3)$ follows from Corollary \ref{needta} (and the definition of a contravariant equivalence). The implication $(2) \Rightarrow (3)$ is obvious. We show that $(3) \Rightarrow (2)$. Let $p: S \rightarrow \calC$ and $\overline{p \circ g}$ be as in $(2)$. 
Passing to a larger universe if necessary, we may suppose that $\calC$ is small.
For each object $C \in \calC$, let $j_C: \calC \rightarrow \SSet$ denote the composition of the Yoneda embedding $j: \calC \rightarrow \calP(\calC)$ with the map $\calP(\calC) \rightarrow \SSet$ given by evaluation at $C$. Combining Proposition \ref{yonedaprop} with Proposition \ref{limiteval}, we deduce that each $j_{C} \circ \overline{p \circ g}$ is a limit diagram. Applying $(3)$, we conclude that each $j_{C} \circ \overline{p \circ g} \circ f^{\triangleleft}$ is a limit diagram. We now apply Propositions \ref{yonedaprop} and \ref{limiteval} to conclude that $\overline{p \circ g} \circ f^{\triangleleft}$ is a limit diagram, as desired.
\end{proof}

\begin{lemma}\label{longwait1}\index{gen}{Yoneda embedding!and left Kan extensions}
Let $S$ be a small simplicial set, $j: S \rightarrow \calP(S)$ the Yoneda embedding,
let $\calC$ denote the full subcategory of $\calP(S)$ spanned by the objects $j(s)$, where
$s$ is a vertex of $S$, and let $\calD$ be an arbitrary $\infty$-category.

\begin{itemize}
\item[$(1)$] Let $f: \calP(S) \rightarrow \calD$ be a functor. Then $f$ is a left Kan extension of $f|\calC$ if and only if $f$ preserves small colimits.
\item[$(2)$] Suppose that $\calD$ admits small colimits, and let $f_0: \calC \rightarrow \calD$
be an arbitrary functor. There exists an extension $f: \calP(S) \rightarrow \calD$ which is a left Kan extension of $f_0 = f|\calC$.
\end{itemize}
\end{lemma}

\begin{proof}
Assertion $(2)$ follows immediately from Lemma \ref{kan2}, since the $\infty$-category
$\calC_{/X}$ is small for each object $X \in \calP(S)$. We will prove $(1)$. Suppose first that
$f$ preserves small colimits. We must show that for each $X \in \calP(S)$, the composition
$$ \calC_{/X}^{\triangleright} \stackrel{\delta}{\rightarrow} \calP(S) \stackrel{f}{\rightarrow} \calD$$
is a colimit diagram. Lemma \ref{longwait0} implies that $\delta$ is a colimit diagram; if $f$ preserves small colimits, then $f \circ \delta$ is also a colimit diagram.

Now suppose that $f$ is a left Kan extension of $f_0 = f | \calC$. We wish to prove that $f$
preserves small colimits. Let $K$ be a small simplicial set, and let
$\overline{p}: K^{\triangleright} \rightarrow \calP(S)$ be a colimit diagram. We must show that
$f \circ \overline{p}$ is also a colimit diagram.

Let $$\overline{\calE} = \calC \times_{ \Fun( \{0\}, \calP(S) } \Fun(\Delta^1,\calP(S) ) \times_{\Fun( \{1\}, \calP(S))} K^{\triangleright},$$
and let $\calE = \overline{\calE} \times_{ K^{\triangleright} } K \subseteq \overline{\calE}$. 
We have a commutative diagram
$$ \xymatrix{ \calE \ar[d] \ar[r] & \overline{\calE} \ar[d] \\
K \ar[r] & K^{\triangleright}. }$$
where the vertical arrows are coCartesian fibrations (Corollary \ref{tweezegork}). 
Let $\overline{\eta}: \overline{\calE} \diamond_{ K^{\triangleright} } K^{\triangleright} 
\rightarrow \calP(S)$ be the natural map, and $\eta = \overline{\eta} | \calE \diamond_{K} K$. 
Proposition \ref{longwait2} implies that $f \circ \eta$ exhibits $f \circ p$ as a left Kan extension of 
$f \circ (\eta | \calE)$ along $q|\calE$. Similarly, $f \circ \overline{\eta}$ exhibits $f \circ \overline{p}$ as a left Kan extension of $f \circ (\overline{\eta} | \overline{\calE})$. It will therefore suffice to prove
that every colimit of $f \circ (\overline{\eta} | \overline{\calE})$ is also a colimit of
$f \circ (\eta | \calE)$. According to Lemma \ref{natee}, it suffices to show that the inclusion
$\calE \subseteq \overline{\calE}$ is a contravariant equivalence in $(\sSet)_{/\calC}$.

Since the map $\overline{\calE} \rightarrow K^{\triangleright} \times \calC$ is a bivariant fibration, 
we can apply Proposition \ref{longwait5} to deduce that the map
$\overline{\calE}^{op} \rightarrow \calC^{op}$ is smooth. Similarly, $\calE^{op} \rightarrow \calC^{op}$ is smooth. According to Proposition \ref{longwait44}, the inclusion
$\calE \subseteq \overline{\calE}$ is a contravariant equivalence if and only if, for every
object $C \in \calC$, the inclusion of fibers $\calE_{C} \subseteq \overline{\calE}_{C}$ 
is a weak homotopy equivalence. Lemma \ref{repco} implies that
$\overline{\calE}_{C} \rightarrow K^{\triangleright}$
is equivalent to the left fibration given by the pullback of the universal left fibration
along the map
$$ K^{\triangleright} \stackrel{\overline{p}}{\rightarrow} \calP(S) \stackrel{s}{\rightarrow} \SSet.$$
We now conclude by applying Proposition \ref{charspacecolimit}, noting that
$\overline{p}$ is a colimit diagram by assumption and that $s$ preserves colimits by
Proposition \ref{limiteval}.
\end{proof}

\begin{theorem}\label{charpresheaf}\index{gen}{presheaf!universal property of $\calP(S)$}
Let $S$ be a small simplicial set, and let $\calC$ be an $\infty$-category which admits small colimits. Composition with the Yoneda embedding $j: S \rightarrow \calP(S)$ induces
an equivalence of $\infty$-categories
$$ \LFun( \calP(S), \calC) \rightarrow \Fun(S,\calC).$$
\end{theorem}

\begin{proof}
Combine Corollary \ref{leftkanextdef} with Lemma \ref{longwait1}.
\end{proof}

\begin{definition}
Let $\calC$ be an $\infty$-category. A full subcategory $\calC' \subseteq \calC$ is {\it stable under colimits} if, for any small diagram $p: K \rightarrow \calC'$ which has a colimit
$\overline{p}: K^{\triangleright} \rightarrow \calC$ in $\calC$, the map $\overline{p}$ factors through $\calC'$.

Let $\calC$ be an $\infty$-category which admits all small colimits. Let $A$ be a collection of objects of $\calC$. We will say that $A$ {\it generates $\calC$ under colimits} if the following condition is satisfied: for any full subcategory $\calC' \subseteq \calC$ containing every element of $A$, if $\calC'$ is stable under colimits, then $\calC = \calC'$. 

We say that a map $f: S \rightarrow \calC$ {\it generates $\calC$ under colimits} if the image
$f(S_0)$ generates $\calC$ under colimits.\index{gen}{generation under colimits}
\end{definition}

\begin{corollary}\label{gencolcot}
Let $S$ be a small simplicial set. Then the Yoneda embedding $j: S \rightarrow \calP(S)$ generates $\calP(S)$ under small colimits.
\end{corollary}

\begin{proof}
Let $\calC$ be the smallest full subcategory of $\calP(S)$ which contains the essential image of
$j$ and is stable under small colimits. Applying Theorem \ref{charpresheaf}, we deduce that
the diagram $j: S \rightarrow \calC$ is equivalent to $F \circ j$, for some colimit-preserving
functor $F: \calP(S) \rightarrow \calC$. We may regard $F$ as a colimit preserving functor
from $\calP(S)$ to itself. Applying Theorem \ref{charpresheaf} again, we deduce that $F$ is
equivalent to the identity functor from $\calP(S)$ to itself. It follows that every object
of $\calP(S)$ is equivalent to an object which lies in $\calC$, so that $\calC = \calP(S)$ as desired.
\end{proof}

\subsection{Complete Compactness}\label{completecomp}

Let $S$ be a small simplicial set, and $f: S \rightarrow \calC$ a diagram in an $\infty$-category $\calC$. Our goal in this section is to analyze the following question: when is the diagram $f: S \rightarrow \calC$ equivalent to the Yoneda embedding $j: S \rightarrow \calP(S)$?
An obvious necessary condition is that $\calC$ admit small colimits (Corollary \ref{storum}). 
Conversely, if $\calC$ admits small colimits, then Theorem \ref{charpresheaf} implies that $f$ is equivalent to $F \circ j$, where $F: \calP(S) \rightarrow \calC$ is a colimit-preserving functor.
We are now reduced to the question of deciding whether or not the functor $F$ is an equivalence.
There are two obvious necessary conditions for this to be so: $f$ must be fully faithful (Proposition \ref{fulfaith}), and $f$ must generate $\calC$ under colimits (Corollary \ref{gencolcot}). We will show that the converse holds, provided that the essential image of $f$ consists of {\em completely compact} objects of $\calC$ (see Definition \ref{complcompdef} below).

We begin by considering an arbitrary simplicial set $S$ and a vertex $s$ of $S$.
Composing the Yoneda embedding $j: S \rightarrow \calP(S)$ with the ``evaluation map''
$$\calP(S)  = \Fun(S^{op}, \SSet) \rightarrow
\Fun( \{s\}, \SSet) \simeq \SSet,$$
we obtain a map $j_{s}: S \rightarrow \SSet$. We will refer to $j_{s}$
as the {\it functor corepresented by $s$}.\index{gen}{corepresentable!functor}

\begin{remark}
The above definition makes sense even when the simplicial set $S$ is not small. However, in
this case we need to replace $\SSet$ (the simplicial nerve of the category of {\em small} Kan complexes) by the (very large) $\infty$-category $\widehat{\SSet}$, where $\hat{\SSet}$ is the simplicial nerve of the category of {\em all} Kan complexes (not necessarily small).
\end{remark}

\begin{definition}\label{complcompdef}\index{gen}{compact object!completely}\index{gen}{completely compact}
Let $\calC$ be an $\infty$-category which admits small colimits. We will say that an object
$C \in \calC$ is {\it completely compact} if the functor $j_{C}: \calC \rightarrow \widehat{\SSet}$
corepresented by $C$ preserves small colimits.
\end{definition}

The requirement that an object $C$ of an $\infty$-category $\calC$ be completely compact
is {\em very} restrictive (see Example \ref{tryu} below). We introduce this notion not because it is a generally useful one, but because it is relevant for the purpose of characterizing $\infty$-categories of presheaves.

Our first goal is to establish that the class of completely compact objects of $\calC$ is stable under retracts.

\begin{lemma}\label{compcompcomp}
Let $\calC$ be an $\infty$-category, $K$ a simplicial set, and $\overline{p}, \overline{q}: K^{\triangleright} \rightarrow \calC$ a pair of diagrams. Suppose $\overline{q}$ is a colimit diagram, and $\overline{p}$ is a retract of $\overline{q}$ in the $\infty$-category $\Fun(K^{\triangleright}, \calC)$. Then $\overline{p}$ is a colimit diagram.
\end{lemma}

\begin{proof}
Choose a map $\sigma: \Delta^2 \times K^{\triangleright} \rightarrow \calC$ such that
$\sigma | \{1\} \times K^{\triangleright} = \overline{q}$ and
$\sigma | \Delta^{ \{0,2\} } \times K^{\triangleright} = \overline{p} \circ \pi_{K^{\triangleright}}$. 
We have a commutative diagram of simplicial sets:
$$ \xymatrix{ \calC_{\sigma/} \ar[r] \ar[d] & \calC_{\sigma| \Delta^2 \times K/} \ar[d] \\
\calC_{\sigma | \Delta^{ \{1,2\} /} \times K^{\triangleright}} \ar[r]^{f} \ar[d] & \calC_{\sigma| \Delta^{ \{1,2\} } \times K/} \ar[d] \\
\calC_{\sigma| \{2\} \times K^{\triangleright}/} \ar[r]^{f'} & \calC_{\sigma| \{2\} \times K/}.}$$

We first claim that both vertical compositions are categorical equivalences. We give the argument for the right vertical composition; the other case is similar. We have a factorization
$$ \calC_{\sigma | \Delta^2 \times K/} \stackrel{g'}{\rightarrow}
\calC_{\sigma| \Delta^{ \{0,2\}} \times K/} \stackrel{g''}{\rightarrow} \calC_{\sigma | \{2\} \times K/}$$
where the $g'$ is a trivial fibration, and $g''$ admits a section $s$, where $s$ is also a section
of the trivial fibration $\calC_{/ \sigma| \Delta^{ \{0,2\} \times K}} \rightarrow \calC_{/ \sigma| \{0\} \times K}$. Consequently, $s$ and therefore also $g''$ are categorical equivalences.
It follows that the map $f'$ is a retract of $f$ in the homotopy category of $\sSet$ (taken with respect to the Joyal model structure). 

The map $f$ sits in a commutative diagram
$$ \xymatrix{ \calC_{\sigma | \Delta^{ \{1,2\}/ } \times K^{\triangleright}} \ar[r]^{f} \ar[d] & \calC_{ \sigma| \Delta^{ \{1,2\}/ } \times K} \ar[d] \\
\calC_{\overline{q}/} \ar[r] & \calC_{q/} }$$
where the vertical maps and the lower horizontal map are trivial fibrations. It follows that
$f$ is a categorical equivalence. Since $f'$ is a retract of $f$, $f'$ is also a categorical equivalence. Since $f'$ is a left fibration, we deduce that $f'$ is a trivial fibration (Corollary \ref{heath}), so that $\overline{p}$ is a colimit diagram as desired.
\end{proof}

\begin{lemma}\label{retcompact}
Let $\calC$ be an $\infty$-category which admits small colimits. Let $C$ and
$D$ be objects of $\calC$. Suppose that $C$ is completely compact, and that $D$
is a retract of $C$ (that is, there exist maps $f: D \rightarrow C$ and $r: C \rightarrow D$
with $r \circ f \simeq \id_{D}$. Then $D$ is completely compact. In particular, if $C$ and
$D$ are equivalent, then $D$ is completely compact.
\end{lemma}

\begin{proof}
Let $j: \calC^{op} \rightarrow \SSet^{\calC}$ denote the Yoneda embedding (for $\calC^{op})$. Since $D$ is a retract of $C$, $j(D)$ is a retract of $j(C)$. Let $\overline{p}: K^{\triangleright} \rightarrow \calC$ be a diagram. Then $j(D) \circ \overline{p}: K^{\triangleright} \rightarrow \SSet$
is a retract of $j(C) \circ \overline{p}: K^{\triangleright} \rightarrow \SSet$ in the $\infty$-category
$\Fun(K^{\triangleright}, \SSet)$. If $\overline{p}$ is a colimit diagram, then $j(C) \circ \overline{p}$ is a colimit diagram (since $C$ is completely compact). Lemma \ref{compcompcomp} now implies
that $j(D) \circ \overline{p}$ is a colimit diagram as well.
\end{proof}

In order to study the condition of complete compactness in more detail, it is convenient to introduce a slightly more general notion.

\begin{definition}\index{gen}{completely compact!left fibration}
Let $\calC$ be an $\infty$-category which admits small colimits, and let
$\phi: \widetilde{\calC} \rightarrow \calC$ be a left fibration. We will say that $\phi$ is
{\it completely compact} if it is classified by a functor $\calC \rightarrow \widehat{\SSet}$ that preserves small colimits.
\end{definition}

\begin{lemma}\label{bstick}
Let $\calC$ be an $\infty$-category which admits small colimits, $f: X' \rightarrow X$
a map of Kan complexes, and
$$ \xymatrix{ \calF' \ar[r] \ar[d] & \calF \ar[d] \\
X' \times \calC \ar[r]^{f \times \id_{\calC}} & X \times \calC }$$
be a diagram of left fibrations over $\calC$, which is a homotopy pullback square
$($ with respect to the covariant model structure on $(\sSet)_{/\calC}$ $)$. 
If $\calF \rightarrow \calC$ is completely compact, then $\calF' \rightarrow \calC$ is completely compact.
\end{lemma}

\begin{proof}
Replacing the diagram by an equivalent one if necessary, we may suppose that 
it is Cartesian and that $f$ is a Kan fibration. Let $\overline{p}: K^{\triangleright} \rightarrow \calC$ be a colimit diagram, and let $F: \calC \rightarrow \hat{\SSet}$ be a functor which classifies
the left fibration $\calF'$. We wish to show that $F \circ \overline{p}$ is a colimit diagram in
$\hat{\SSet}$.

We have a pullback diagram
$$ \xymatrix{ K \times_{\calC} \calF' \ar[r] \ar[d]^{\psi'} & K \times_{\calC} \calF \ar[d]^{\psi} \\
K^{\triangleright} \times_{\calC} \calF' \ar[r] & K^{\triangleright} \times_{\calC} \calF }$$
of simplicial sets, which is homotopy Cartesian (with respect to the usual model structure on $\sSet$) since the horizontal maps are pullbacks of $f$.
Since $\calF$ is completely compact, Proposition \ref{charspacecolimit} implies that the inclusion $\psi$ is a weak homotopy equivalence. It follows that $\psi'$ is also a weak homotopy equivalence. Applying Proposition \ref{charspacecolimit} again, we deduce that $F \circ \overline{p}$ is a colimit diagram as desired.
\end{proof}

\begin{lemma}\label{compactslice}
Let $\calC$ be a presentable $\infty$-category, $p: K \rightarrow \calC$
be a small diagram, and let $X \in \calC_{/p}$ be an object whose image in $\calC$ is completely compact. Then $X$ is completely compact.
\end{lemma}

\begin{proof}
Let $\overline{p}: K^{\triangleleft} \rightarrow \calC$ be a limit of $p$, carrying the cone point to an object $Z \in \calC$. Then we have trivial fibrations
$$ \calC_{/Z} \leftarrow \calC_{/\overline{p} } \rightarrow \calC_{/p}.$$
Consequently, we may replace the diagram $p: K \rightarrow \calC$ with the inclusion
$\{Z\} \rightarrow \calC$.

We may identify the object $X \in \calC_{/Z}$ with a morphism $f: Y \rightarrow Z$ in $\calC$.
We have a commutative diagram of simplicial sets
$$ \xymatrix{ (\calC_{/Z})_{f/} \ar[dr]^{\psi} \ar[r]^{\theta} & (\calC_{/Y})_{f/} \ar[d] \ar[r]^{\theta'} & (\calC_{/Y})^{f/} \ar[d]^{\psi'} \\
& \calC_{/Z} \ar[r]^{\theta'_0} & \calC^{/Z} }$$
where $\theta$ is an isomorphism, the maps $\theta'$ and $\theta'_0$ are categorical equivalences (see \S \ref{quasilimit2}), and the vertical maps are left fibrations. We wish to prove
that $\psi$ is a completely compact left fibration. It will therefore suffice to prove that
$\psi'$ is completely compact. We have a (homotopy) pullback diagram
$$ \xymatrix{ \calC_{Y/}^{/f} \ar[r] \ar[d] & \calC_{Y/}^{\Delta^1} \times_{\calC^{\{1\}}} \{Z\} \ar[d] \\
\calC^{/Z} \ar[r] & (\calC_{Y/} \times_{\calC} \{Z\}) \times \calC^{/Z} }$$
of left fibrations over $\calC^{/Z}$. We observe that the left fibrations in the lower part of the diagram are constant. According to Lemma \ref{bstick}, to prove that $\psi'$ is completely compact, it will suffice to prove that the left fibration $\calC_{Y/}^{\Delta^1} \times_{ \calC^{ \{1\} } } \{Z\} \stackrel{\psi''}{\rightarrow} \calC^{/Z}$ is completely compact. We observe that $\psi''$ admits a factorization
$$ \calC_{Y/}^{\Delta^1} \times_{ \calC^{ \{1\} }} \{Z\} \stackrel{\phi}{\rightarrow}
\calC_{Y/} \times_{ \calC^{ \{0\} } } \calC^{/Z} \stackrel{\phi'}{\rightarrow} \calC^{/Z}$$
where $\phi$ is a trivial fibration, and $\phi'$ is a pullback of the left fibration
left fibration $\phi'': \calC_{Y/} \rightarrow \calC$. Since $Y$ is completely compact, $\phi''$ is completely compact. The projection $\calC^{/Z} \rightarrow \calC$ is equivalent to
$\calC_{/Z} \rightarrow \calC$, and therefore commutes with colimits by Proposition \ref{needed17}.
It follows that $\phi'$ is completely compact, which completes the proof.
\end{proof}

\begin{proposition}\label{dda}
Let $S$ be a small simplicial set, and let $j: S \rightarrow \calP(S)$ denote the Yoneda embedding.
Let $C$ be an object of $\calP(S)$. The following conditions are equivalent:
\begin{itemize}
\item[$(1)$] The object $C \in \calP(S)$ is completely compact.
\item[$(2)$] There exists a vertex $s$ of $S$ such that $C$ is a retract of $j(s)$.
\end{itemize}
\end{proposition}

\begin{proof}
Suppose first that $(1)$ is satisfied.
Let $S_{/C} = S \times_{\calP(S)} \calP(S)_{/C}$. According to Lemma \ref{longwait0}, 
the natural map
$$ S^{\triangleright}_{/C} \stackrel{j'}{\rightarrow} \calP(S)_{/C}^{\triangleright} \rightarrow \calP(S)$$
is a colimit diagram. Let $f: \calP(S) \rightarrow \SSet$ be the functor corepresented by $C$.
Since $C$ is completely compact. $f(C)$ can be identified with a colimit of the diagram $f | S_{/C}$. 
The space $f(C)$ is homotopy equivalent to $\bHom_{\calP(S)}(C,C)$, and therefore contains a point corresponding to $\id_{C}$. It follows that $\id_{C}$ lies in the image of
$\bHom_{\calP(S)}(C, j'(\widetilde{s}) ) \rightarrow \bHom_{\calP(S)}(C,C)$, for some
vertex $\widetilde{s}$ of $S_{/C}$. The vertex $\widetilde{s}$ classifies a vertex
$s \in S$ equipped with a morphism $\alpha: j(s) \rightarrow C$. It follows that there
is a commutative triangle
$$ \xymatrix{ & j(s) \ar[dr]^{\alpha} & \\
C \ar[ur] \ar[rr]^{\id_{C}} & & C }$$
in the $\infty$-category $\calP(S)$, so that $C$ is a retract of $j(s)$.

Now suppose that $(2)$ is satisfied. According to Lemma \ref{retcompact}, it suffices to prove that $j(s)$ is completely compact. Using Lemma \ref{repco}, we may identify the functor
$\calP(S) \rightarrow \SSet$ co-represented by $j(s)$ with the functor given by evaluation at $s$.
Proposition \ref{limiteval} implies that this functor preserves all limits and colimits that exist in
$\calP(S)$.
\end{proof}

\begin{example}\label{tryu}
Let $\calC$ be the $\infty$-category $\SSet$ of spaces. Then an object
$C \in \SSet$ is completely compact if and only if it is equivalent to $\ast$, the final
object of $\SSet$.
\end{example}

We now use the theory of completely compact objects to give a characterization of
presheaf $\infty$-categories.

\begin{proposition}\label{trumptow}\index{gen}{completely compact!and presheaf $\infty$-categories}
Let $S$ be a small simplicial set and $\calC$ an $\infty$-category which admits small colimits. Let $F: \calP(S) \rightarrow \calC$ be functor which preserves small colimits, and
$f = F \circ j$ its composition with the Yoneda embedding $j: S \rightarrow \calP(S)$. 
Suppose further that:
\begin{itemize}
\item[$(1)$] The functor $f$ is fully faithful.

\item[$(2)$] For every vertex $s$ of $S$, the object $f(s) \in \calC$ is
completely compact.
\end{itemize}

Then $F$ is fully faithful.
\end{proposition}

\begin{proof}
Let $C$ and $D$ be objects of $\calP(S)$. We wish to prove that the natural map
$$ \eta_{C,D}: \bHom_{\calP(S)}(C,D) \rightarrow \bHom_{\calC}(F(C), F(D))$$ is an isomorphism in the homotopy category $\calH$.
Suppose first that $C$ belongs to the essential image
of $j$. Let $G: \calP(S) \rightarrow \SSet$ be a functor co-represented by $C$, and let
$G': \calC \rightarrow \SSet$ be a functor co-represented by $F(C)$. Then we have a natural transformation of functors $G \rightarrow G' \circ F$. Assumption $(2)$ implies that $G'$ preserves small colimits, so that $G' \circ F$ preserves small colimits. Proposition \ref{dda} implies that
$G$ preserves small colimits. It follows that the collection of objects $D \in \calP(S)$ such
that $\eta_{C,D}$ is an equivalence is stable under small colimits. If $D$ belongs to the essential image of $j$, then assumption $(1)$ implies that $\eta_{C,D}$ is an equivalence. It follows from Lemma \ref{longwait0} that the essential image of $j$ generates $\calP(S)$ under small colimits; thus $\eta_{C,D}$ is an isomorphism in $\calH$ for every object $D \in \calP(S)$.

We now prove the result in general. Fix $D \in \calP(S)$. Let $H: \calP(S)^{op} \rightarrow \SSet$
be a functor represented by $D$, and let $H': \calC^{op} \rightarrow \SSet$ be a functor represented by $FD$. Then we have a natural transformation of functors $H \rightarrow H' \circ F^{op}$, which we wish to prove is an equivalence. By assumption, $F^{op}$ preserves small limits. Proposition \ref{yonedaprop} implies that $H$ and $H'$ preserve small limits. It follows that the collection $P$ of objects $C \in \calP(S)$ such that $\eta_{C,D}$ is an equivalence is stable under small colimits.
The special case above established that $P$ contains the essential image of the Yoneda embedding. We once again invoke Lemma \ref{longwait0} to deduce every object of $\calP(S)$ belongs to $P$, as desired.
\end{proof}

\begin{corollary}\label{charpr}\index{gen}{equivalence!of presheaf $\infty$-categories}
Let $\calC$ be an $\infty$-category which admits small colimits. Let $S$ be a small
simplicial set and $F: \calP(S) \rightarrow \calC$ a colimit preserving functor. Then
$F$ is an equivalence if and only if the following conditions are satisfied:
\begin{itemize}
\item[$(1)$] The composition $f =F \circ j: S \rightarrow \calC$ is fully faithful.
\item[$(2)$] For every vertex $s \in S$, the object $f(s) \in \calC$ is completely compact.
\item[$(3)$] The set of objects $\{ f(s): s \in S_0\}$ generates $\calC$ under colimits.
\end{itemize}
\end{corollary}

\begin{proof}
If $(1)$, $(2)$, and $(3)$ are satisfied, then $F$ is fully faithful (Proposition \ref{trumptow}). 
Since $\calP(S)$ is admits small colimits, and $F$ preserves small colimits, the essential image of $F$ is stable under small colimits. Using $(3)$, we conclude that $F$ is essentially surjective and therefore an equivalence of $\infty$-categories. For the converse, it suffices to check that
$\id_{\calP(S)}: \calP(S) \rightarrow \calP(S)$ satisfies $(1)$, $(2)$, and $(3)$. For this, we invoke Propsition \ref{fulfaith}, Proposition \ref{dda}, and Lemma \ref{longwait0}, respectively.
\end{proof}

\begin{corollary}\label{swapKK}\index{gen}{presheaf!and overcategories}
Let $\calC$ be a small $\infty$-category, and let $p: K \rightarrow \calC$ be a diagram, and
let $p': K \rightarrow \calP(\calC)$ be the composition of $p$ with the Yoneda embedding
$j: \calC \rightarrow \calP(\calC)$, and let $f: \calC_{/p} \rightarrow \calP(\calC)_{/p'}$ be the induced map. Let $F: \calP(\calC_{/p}) \rightarrow \calP(\calC)_{/p'}$ be a colimit-preserving functor such that $F \circ j'$ is equivalent to $f$, where $j': \calC_{/p} \rightarrow \calP(\calC_{/p})$ denotes the Yoneda embedding for $\calC_{/p}$ $($according to Theorem \ref{charpresheaf}, $F$ exists and is unique up to equivalence$)$. Then $F$ is an equivalence of $\infty$-categories.
\end{corollary}

\begin{proof}
We will show that the $f$ satisfies conditions $(1)$ through $(3)$ of Corollary \ref{charpr}.
The assertion that $f$ is fully faithful follows immediately from the assertion that $j$ is fully faithful (Proposition \ref{fulfaith}). To prove that the essential image of $f$ consists of completely compact objects, we use Lemma \ref{compactslice} to reduce to proving that the essential image of $j$ consists of completely compact objects of $\calP(\calC)$, which follows from Proposition \ref{dda}.
It remains to prove that $\calP(\calC)_{/p'}$ is generated under colimits by $f$. Let
$\overline{X}$ be an object of $\calP(\calC)_{/p'}$ and $X$ its image in $\calP(\calC)$. 
Let $\calD \subseteq \calP(\calC)$ be the essential image of $j$, and $\overline{\calD}$ the inverse image of $\calD$ in $\calP(\calC)_{/p'}$, so that $\overline{\calD}$ is the essential image of $f$.
Using Lemma \ref{longwait0}, we can choose a colimit diagram $\overline{q}: L^{\triangleright} \rightarrow \calP(\calC)$ which carries the cone point to $X$ such that $q = \overline{q}|L$ factors through $\calD$. Since the inclusion of the cone point into $L^{\triangleright}$ is right anodyne,
there exists a map $\overline{q}': L^{\triangleright} \rightarrow \calP(\calC)_{/p'}$ lifting
$\overline{q}$, which carries the cone point of $L^{\triangleright}$ to $\overline{X}$. 
Proposition \ref{needed17} implies that $\overline{q}'$ is a colimit diagram, so that 
$\overline{X}$ can be written as the colimit of a diagram $L \rightarrow \overline{\calD}$.
\end{proof}

\section{Adjoint Functors}\label{c5s2}

\setcounter{theorem}{0}

Let $\calC$ and $\calD$ be (ordinary) categories. Two functors
$$ \Adjoint{F}{ \calC }{\calD}{G}$$
are said to be {\it adjoint} to one another if there is a functorial bijection
$$ \Hom_{\calD}(F(C), D) \simeq \Hom_{\calC}( C, G(D))$$
defined for $C \in \calC$, $D \in \calD$. Our goal in this section is
to extend the theory of adjoint functors to the $\infty$-categorical setting.\index{gen}{adjoint functor!between ordinary categories}

By definition, a pair of functors $F$ and $G$ (as above) are adjoint if and only if they
determine the same correspondence $$ \calC^{op} \times \calD \rightarrow \Set.$$ 
In \S \ref{corresp}, we introduced an $\infty$-categorical generalization of the notion of a correspondence. In certain cases, a correspondence $\calM$ from an $\infty$-category
$\calC$ to an $\infty$-category $\calD$ determines a functor $F: \calC \rightarrow \calD$, which we say is a functor {\it associated} to $\calM$. We will study these associated functors in \S \ref{afunc1}. The notion of a correspondence is self-dual, so it is possible that the correspondence $\calM$ also determines an associated functor $G: \calD \rightarrow \calC$. In this case, we will say that
$F$ and $G$ are adjoint. We will study the basic properties of adjoint functors in \S \ref{afunc2}.

One of the most important features of adjoint functors is their behavior with respect to limits and colimits: left adjoints preserve colimits, while right adjoints preserve limits. We will prove an $\infty$-categorical analogue of this statement in \S \ref{afunc3}. In certain situations, the {\it adjoint functor theorem} provides a converse to this statement: see \S \ref{aftt}.

The theory of model categories provides a host of examples of adjoint functors between $\infty$-categories. In \S \ref{afunc4}, we will show that a simplicial Quillen adjunction
between a pair of model categories $(\bfA,\bfA')$ determines an adjunction between the associated $\infty$-categories $(\sNerve(\bfA^{\degree}), \sNerve({\bfA'}^{\degree}))$. We will also consider some other examples of situations which give rise to adjoint functors.

In \S \ref{afunc4half}, we study the behavior of adjoint functors when restricted to overcategories.
Our main result (Proposition \ref{curpse}) can summarized as follows: suppose that
$F: \calC \rightarrow \calD$ is a functor between $\infty$-categories which admits a right adjoint
$G$. Assume further that the $\infty$-category $\calC$ admits pullbacks. Then for every object
$C$, the induced functor $\calC_{/C} \rightarrow \calD_{/FC}$ admits a right adjoint, given by
the formula 
$$(D \rightarrow FC) \mapsto (GD \times_{GFC} C \rightarrow C).$$

If a functor $F: \calC \rightarrow \calD$ has a right adjoint $G$, then $G$ is uniquely determined up to equivalence. In \S \ref{afunc5}, we will prove a strong version of this statement, phrased as an (anti)equivalence of functor categories.

In \S \ref{locfunc}, we will restrict the theory of adjoint functors to the special case in which one of the functors is the inclusion of a full subcategory. In this case, we obtain the theory of {\em localizations} of $\infty$-categories. This theory will play a central role in our study of presentable $\infty$-categories (\S \ref{c5s6}), and later in the study of $\infty$-topoi (\S \ref{chap6}). It is also useful in the study of {\it factorization systems} on $\infty$-categories, which we will discuss in
\S \ref{factgen1}.

Finally, in \S \ref{cataut}, we will apply some of the ideas of this and earlier sections to analyze the $\infty$-category of equivalences from $\Cat_{\infty}$ to itself. The main result (Theorem \ref{cabbi}) is that there is essentially only one nontrivial self-equivalence of $\Cat_{\infty}$: namely, the operation which carries each $\infty$-category $\calC$ to its opposite $\calC^{op}$.

\subsection{Correspondences and Associated Functors}\label{afunc1}

Let $p: X \rightarrow S$ be a Cartesian fibration of simplicial sets. In \S \ref{universalfib}, we
saw that $p$ is classified by a functor $S^{op} \rightarrow \Cat_{\infty}$. In particular, if
$S = \Delta^1$, then $p$ determines a diagram
$$ G: \calD \rightarrow \calC$$
in the $\infty$-category $\Cat_{\infty}$, which is well-defined up to equivalence.
We can obtain this diagram by applying the straightening
functor $\St_{S}^{+}$ to the marked simplicial set $X^{\natural}$, and then taking a fibrant replacement. In general, this construction is rather complicated. However, in the special case where
$S = \Delta^1$, it is possible to give a direct construction of $G$; that is our goal in this section.

\begin{definition}\label{fibas}\index{gen}{functor!associated to a correspondence}\index{gen}{correspondence!associated functor}
Let $p: \calM \rightarrow \Delta^1$ be a Cartesian fibration, and suppose given
equivalences of $\infty$-categories $h_0: \calC \rightarrow p^{-1} \{0\}$ and $h_1: \calD \rightarrow p^{-1} \{1\}$. We will say that a functor $g: \calD \rightarrow \calC$ is {\it associated to $\calM$} if
there is a commutative diagram
$$ \xymatrix{ \calD \times \Delta^1 \ar[dr] \ar[rr]^{s} & & \calM \ar[dl] \\
& \Delta^1 & }$$
such that $s| \calD \times \{1\} = h_1$, $s| \calD \times \{0\} = h_0 \circ g$, and
$s| \{x\} \times \Delta^1$ is a $p$-Cartesian edge of $\calM$ for every object $x$ of $\calD$.
\end{definition}

\begin{remark}
The terminology of Definition \ref{fibas} is slightly abusive: it would be more accurate to say that $g$ is associated to the triple $(p: \calM \rightarrow \Delta^1, h_0: \calC \rightarrow p^{-1} \{0\},
h_1: \calD \rightarrow p^{-1} \{1\} )$. 
\end{remark}

\begin{proposition}\label{candi}
Let $\calC$ and $\calD$ be $\infty$-categories, and let $g: \calD \rightarrow \calC$
be a functor. 
\begin{itemize}
\item[$(1)$] There exists a diagram
$$ \xymatrix{ \calC \ar[r] \ar[d] & \calM \ar[d]^{p} & \calD \ar[l] \ar[d] \\
\{0\} \ar[r] & \Delta^1 & \{1\} \ar[l] }$$
where $p$ is a Cartesian fibration, the associated maps $\calC \rightarrow p^{-1} \{0\}$ and $\calD \rightarrow p^{-1} \{1\}$ are isomorphisms, and $g$ is associated to $\calM$.

\item[$(2)$] Suppose given a commutative diagram
$$ \xymatrix{ \calC \ar[r] \ar[dd] & \calM' \ar[d]^{s} & \calD \ar[l] \ar[dd] \\
& \calM \ar[d]^{p} & \\
\{0\} \ar[r] & \Delta^1 & \{1\} \ar[l] }$$ 
where $s$ is a categorical equivalence, $p$ and $p' = p \circ s$ are Cartesian fibrations,
and the maps $\calC \rightarrow p^{-1} \{0\}$, $\calD \rightarrow p^{-1} \{1\}$ are categorical equivalences. The functor $g$ is associated to $\calM$ if and only if it is associated to $\calM'$.

\item[$(3)$] Suppose given diagrams 
$$ \xymatrix{ \calC \ar[r] \ar[d] & \calM' \ar[d]^{p'} & \calD \ar[l] \ar[d] \\
\{0\} \ar[r] & \Delta^1 & \{1\} \ar[l] }$$
$$ \xymatrix{ \calC \ar[r] \ar[d] & \calM'' \ar[d]^{p''} & \calD \ar[l] \ar[d] \\
\{0\} \ar[r] & \Delta^1 & \{1\} \ar[l] }$$
as above, such that $g$ is associated to both $\calM'$ and $\calM''$. Then there
exists a third such diagram
$$ \xymatrix{ \calC \ar[r] \ar[d] & \calM \ar[d]^{p} & \calD \ar[l] \ar[d] \\
\{0\} \ar[r] & \Delta^1 & \{1\} \ar[l] }$$
and a diagram 
$$ \calM' \leftarrow \calM \rightarrow \calM''$$
of categorical equivalences in $(\sSet)_{\calC \coprod \calD/ \,/\Delta^1}$.
\end{itemize}
\end{proposition}

\begin{proof}
We begin with $(1)$. Let $\calC^{\natural}$ and $\calD^{\natural}$ denote
the simplicial sets $\calC$ and $\calD$ considered as marked simplicial sets, where the marked edges are precisely the equivalences. We set
$$N = ( \calD^{\natural} \times (\Delta^{1})^{\sharp}) \coprod_{ \calD^{\natural} \times \{0\}^{\sharp} } \calC^{\natural}.$$ 
The small object argument implies the existence of a factorization
$$ N \rightarrow N(\infty) \rightarrow (\Delta^1)^{\sharp},$$
where the left map is marked anodyne and the right map has the right lifting property with respect to all marked anodyne morphisms. We remark
that we can obtain $N(\infty)$ as the colimit of a transfinite sequence of simplicial sets $N(\alpha)$, where $N(0) = N$, $N(\alpha)$ is the colimit of the sequence $\{ N(\beta) \}_{\beta < \alpha}$ when
$\alpha$ is a limit ordinal, and each $N(\alpha+1)$ fits into a pushout diagram
$$ \xymatrix{ A \ar@{^{(}->}[d] \ar[r] & N(\alpha) \ar[d] \ar[dr] &  \\
B \ar[r] \ar@{-->}[ur] & N(\alpha+1) \ar[r] & (\Delta^1)^{\sharp} }$$
where the left vertical map is one of the generators for the class of marked anodyne maps
given in Definition \ref{markanod}. We may furthermore assume that there does {\em not} exist a dotted arrow as indicated in the diagram. It follows by induction on $\alpha$ that
$N(\alpha) \times_{\Delta^1} \{0\} \simeq \calC^{\natural}$ and $N(\alpha) \times_{ \Delta^1} \{1\} \simeq \calD^{\natural}$. According to Proposition \ref{dubudu}, $N(\infty) \simeq \calM^{\natural}$ for
some Cartesian fibration $\calM \rightarrow \Delta^1$. It follows immediately that
$\calC \simeq \calM_{\{0\}}$, $\calD \simeq \calM_{\{1\}}$, and that $g$ is associated to $\calM$.

We now prove $(2)$. The ``if'' direction is immediate from the definition. Conversely, suppose that $g$ is associated to $\calM$. To show that $g$ is associated to $\calM'$, we need to produce the dotted arrow indicated in the diagram
$$ \xymatrix{ \calD \times \bd \Delta^1 \ar[r] \ar@{^{(}->}[d] & \calM' \\
\calD \times \Delta^1. \ar@{-->}[ur] & }$$
According to Proposition \ref{princex}, we may replace $\calM'$ by the equivalent $\infty$-category
$\calM$; the desired result then follows form the assumption that $g$ is associated to $\calM$.

To prove $(3)$, we take $\calM$ to be the correspondence constructed in the course of proving $(1)$. It will suffice to construct an appropriate categorical equivalence $\calM \rightarrow \calM'$; the same argument will construct the desired map $\calM \rightarrow \calM''$. Consider the diagram
$$ \xymatrix{ N \ar@{^{(}->}[d]^{s'} \ar[r]^{s} & \calM' \ar[d] \\
\calM \ar[r] \ar@{-->}[ur]^{s''} & \Delta^1. }$$
(Here we identify $N$ with its underlying simplicial set by forgetting the class of marked edges, and
the top horizontal map exhibits $g$ as associated to $\calM'$.) In the terminology of \S \ref{funkystructure}, the maps $s$ and $s'$ are both quasi-equivalences. By Proposition
\ref{qequiv}, they are categorical equivalences. The projection $\calM' \rightarrow \Delta^1$ is
a categorical fibration and $s'$ is a trivial cofibration, which ensures the existence of the arrow $s''$. 
The factorization $s = s'' \circ s'$ shows that $s''$ is a categorical equivalence, and completes the proof.
\end{proof}

Proposition \ref{candi} may be informally summarized by saying that every functor
$g: \calD \rightarrow \calC$ is associated to some Cartesian fibration $p: \calM \rightarrow \Delta^1$, and that $\calM$ is determined up to equivalence. Conversely, the Cartesian fibration also determines $g$:

\begin{proposition}\label{funcas}
Let $p: \calM \rightarrow \Delta^1$ be a Cartesian fibration, and let
$h_0: \calC \rightarrow p^{-1} \{0\}$ and $h_1: \calD \rightarrow p^{-1} \{1\}$ be
categorical equivalences. There exists a functor $g: \calD \rightarrow \calC$ associated to $\calM$.
Any other functor $g': \calC \rightarrow \calD$ is associated to $p$ if and only if $g$ is equivalent to $g'$ as objects of the $\infty$-category $\calC^{\calD}$.
\end{proposition}

\begin{proof}
Consider the diagram
$$ \xymatrix{ \calD^{\flat} \times \{1\} \ar[d] \ar[r] & \calM^{\natural} \ar[d] \\
\calD^{\flat} \times (\Delta^1)^{\sharp} \ar[r] \ar@{-->}^{s}[ur] & (\Delta^1)^{\sharp}.}$$
By Proposition \ref{markanodprod}, the left vertical map is marked anodyne, so the
dotted arrow exists. Consider the map $s_0: s| \calD \times \{0\}: \calD \rightarrow p^{-1} \{0\}$. 
Since $h_0$ is a categorical equivalence, there exists a map $g: \calD \rightarrow \calC$
such that the functions $h_0 \circ g$ and $s_0$ are equivalent. Let $e: \calD \times \Delta^1 \rightarrow \calM$ be an equivalence from $h_0 \circ g$ to $s_0$. Let $e': \calD \times \Lambda^2_1 \rightarrow \calM$ be the result of amalgamating $e$ with $s$. Then we have a commutative diagram of marked simplicial sets
$$ \xymatrix{ \calD^{\flat} \times (\Lambda^2_1)^{\sharp} \ar@{^{(}->}[d] \ar[r]^{e'} & \calM^{\natural} \ar[d] \\
\calD^{\flat} \times (\Delta^2)^{\sharp} \ar[r] \ar@{-->}[ur]^{e''} & (\Delta^1)^{\sharp}. }$$
Because left vertical map is marked anodyne there exists a morphism $e''$ as indicated, rendering the diagram commutative. The restriction
$e''| \calD \times \Delta^{ \{0,2\} }$ exhibits $g$ as associated to $\calM$.

Now suppose that $g'$ is another functor associated to $p$. Then there exists a commutative diagram of marked simplicial sets
$$ \xymatrix{ \calD^{\flat} \times \{1\} \ar@{^{(}->}[d] \ar[r] & \calM^{\natural} \ar[d] \\
\calD^{\flat} \times (\Delta^1)^{\sharp} \ar[r] \ar@{-->}[ur]^{s'} & (\Delta^1)^{\sharp},}$$
with $g' = s'| \calD \times \{0\}$. Let $s''$ be the map obtained by amalgamating
$s$ and $s'$. Consider the diagram
$$ \xymatrix{ \calD^{\flat} \times (\Lambda^2_2)^{\sharp} \ar@{^{(}->}[d] \ar[r]^{s''} & \calM^{\natural} \ar[d] \\
\calD^{\flat} \times (\Delta^2)^{\sharp} \ar[r] \ar@{-->}[ur]^{s'''} & (\Delta^1)^{\sharp}. }$$
Since the left vertical map is marked anodyne, the indicated dotted arrow $s''$ exists.
The restriction $s''| \calD \times \Delta^{ \{0,1\} }$ is an equivalence between $h_0 \circ g$ and $h_0 \circ g'$. Since $h_0$ is a categorical equivalence, $g$ and $g'$ are themselves homotopic.

Conversely, suppose that $f: \calD \times \Delta^1 \rightarrow \calC$ is an equivalence from $g'$ to $g$. The maps $s$ and $h_0 \circ f$ amalgamate to give a map $f': \calD \times \Lambda^2_1 \rightarrow \calC$ which fits into a commutative diagram of marked simplicial sets:
$$ \xymatrix{ \calD^{\flat} \times (\Lambda^2_1)^{\sharp} \ar@{^{(}->}[d] \ar[r]^{f'} & \calM^{\natural} \ar[d] \\
\calD^{\flat} \times (\Delta^2)^{\sharp} \ar[r] \ar@{-->}[ur]^{f''} & (\Delta^1)^{\sharp}. }$$
The left vertical map is marked anodyne, so there exists a dotted arrow $f''$ as indicated; then
the map $f'' | \calD \times \Delta^{ \{0,2\} }$ exhibits that $g'$ is associated to $p$.
\end{proof}

\begin{proposition}\label{compass}
Let $p: \calM \rightarrow \Delta^2$ be a Cartesian fibration, and suppose given equivalences
of $\infty$-categories $\calC \rightarrow p^{-1} \{0\}$, $\calD \rightarrow p^{-1} \{1\}$, and
$\calE \rightarrow p^{-1} \{2\}$. Suppose that $\calM \times_{\Delta^2} \Delta^{ \{0,1\} }$
is associated to a functor $f: \calD \rightarrow \calC$, and that $\calM \times_{ \Delta^2} \Delta^{ \{1,2\} }$ is associated to a functor $g: \calE \rightarrow \calD$. Then $\calM \times_{\Delta^2} \Delta^{ \{0,2\} }$ is associated to the composite functor $f \circ g$.
\end{proposition}

\begin{proof}
Let $X$ be the mapping simplex of the sequence of functors
$$ \calE \stackrel{g}{\rightarrow} \calD \stackrel{f}{\rightarrow} \calC.$$
Since $f$ and $g$ are associated to restrictions of $\calM$, we obtain a commutative diagram
$$ \xymatrix{ X \times_{\Delta^2} \Lambda^2_1 \ar@{^{(}->}[d] \ar[r] & \calM \ar[d] \\
X \ar[r] \ar@{-->}[ur]^{s} & \Delta^2.}$$
The left vertical inclusion is a pushout of $\calE \times \Lambda^2_1 \subseteq \calE \times \Delta^2$, which is inner anodyne. Since $p$ is inner anodyne, there exists a dotted arrow $s$ as indicated in the diagram. The restriction $s| X \times_{\Delta^2} \Delta^{ \{0,2\}}$ exhibits
that the functor $f \circ g$ is associated to the correspondence $\calM \times_{\Delta^2} \Delta^{ \{0,2\}}$.
\end{proof}

\begin{remark}
Taken together, Propositions \ref{candi} and \ref{funcas} assert that there is a bijective correspondence between equivalence classes of functors $\calD \rightarrow \calC$
and equivalence classes of Cartesian fibrations $p: \calM \rightarrow \Delta^1$ equipped with
equivalences $\calC \rightarrow p^{-1} \{0\}$, $\calD \rightarrow p^{-1} \{1\}$.
\end{remark}

\subsection{Adjunctions}\label{afunc2}

In \S \ref{afunc1}, we established a dictionary that allows us to pass back and forth between functors
$g: \calD \rightarrow \calC$ and Cartesian fibrations $p: \calM \rightarrow \Delta^1$. The dual argument shows if $p$ is a coCartesian fibration it also determines a functor
$f: \calC \rightarrow \calD$. In this case, we will say that $f$ and $g$ are {\em adjoint} functors.\index{gen}{adjoint functor!between $\infty$-categories}

\begin{definition}\index{gen}{adjunction}\index{gen}{correspondence!adjunction}
Let $\calC$ and $\calD$ be $\infty$-categories. An {\it adjunction} between $\calC$ and $\calD$ is a map $q: \calM \rightarrow \Delta^1$ which is both a Cartesian fibration and a coCartesian fibration, together with equivalences $\calC \rightarrow \calM_{\{0\}}$ and $\calD \rightarrow \calM_{\{1\}}$.

Let $\calM$ be an adjunction between $\calC$ and $\calD$, and let $f: \calC \rightarrow \calD$ and $g: \calD \rightarrow \calC$ be functors associated to $\calM$. In this case, we will say that $f$ is {\it left adjoint} to $g$ and $g$ is {\it right adjoint} to $f$.\index{gen}{right adjoint}\index{gen}{left adjoint}
\end{definition}

\begin{remark}
Propositions \ref{candi} and \ref{funcas} imply that if a functor $f: \calC \rightarrow \calD$ has a right adjoint $g: \calD \rightarrow \calC$, then $g$ is uniquely determined up to homotopy. In fact, we will later see that $g$ is determined up to a contractible ambiguity.
\end{remark}

We now verify a few basic properties of adjunctions:

\begin{lemma}\label{gruft}
Let $p: X \rightarrow S$ be a locally Cartesian fibration of simplicial sets. Let $e: s \rightarrow s'$ be an edge of $S$ with the following property:
\begin{itemize}
\item[$(\ast)$] For every $2$-simplex
$$ \xymatrix{ & x' \ar[dr]^{\overline{e}'} & \\
x \ar[rr]^{\overline{e}''} \ar[ur]^{\overline{e}} & & x''}$$
in $X$ such that $p( \overline{e}) = e$, if $\overline{e}$ and $\overline{e}'$ are locally $p$-Cartesian, then $\overline{e}''$ is locally $p$-Cartesian.
\end{itemize}
Let $\overline{e}: x \rightarrow y$ be a locally $p$-coCartesian edge such that
$p( \overline{e} ) = e$. Then $\overline{e}$ is $p$-coCartesian.
\end{lemma}

\begin{proof}
We must show that for any $n \geq 2$ and any diagram
$$ \xymatrix{ \Lambda^n_0 \ar[r]^{f} \ar@{^{(}->}[d] & X \ar[d] \\
\Delta^n \ar[r] \ar@{-->}[ur] & S }$$
such that $f|\Delta^{ \{0,1\} }=\overline{e}$, there exists a dotted arrow as indicated. Pulling back along the bottom horizontal map, we may reduce to the case $S = \Delta^n$; in particular, $X$ and
$S$ are both $\infty$-categories.

According to (the dual of) Proposition \ref{charCart}, it suffices to show that composition with $\overline{e}$ gives a homotopy Cartesian diagram
$$ \xymatrix{ \bHom_{X}(y,z) \ar[r] \ar[d] & \bHom_{X}(x,z) \ar[d] \\
\bHom_{S}(p(y),p(z)) \ar[r] & \bHom_{S}(p(x),p(z))}.$$
There are two cases to consider: if $\bHom_{S}(p(y), p(z)) = \emptyset$, there is nothing to prove. Otherwise, we must show that composition with $f$ induces a homotopy equivalence $\bHom_{X}(y,z) \rightarrow \bHom_{X}(x,z)$.

In view of the assumption that $S = \Delta^n$, there is a unique morphism $g_0: p(y) \rightarrow p(z)$. Let $g: y' \rightarrow z$ be a locally $p$-Cartesian edge lifting $g_0$. We have a commutative diagram
$$ \xymatrix{ \bHom_{X}(y,y') \ar[r] \ar[d] & \bHom_{X}(x,y') \ar[d] \\
\bHom_{X}(y,z) \ar[r] & \bHom_{X}(x,z). }$$
Since $g$ is locally $p$-Cartesian, the left vertical arrow is a homotopy equivalence. Since
$e$ is locally $p$-coCartesian, the top horizontal arrow is a homotopy equivalence. It will therefore suffice to show that the map $\bHom_{X}(x,y') \rightarrow \bHom_{X}(x,z)$ is a homotopy equivalence.

Choose a locally $p$-Cartesian edge $\overline{e}': x' \rightarrow y'$ in $X$ with $p( \overline{e}')=e$, so that we have another commutative diagram
$$ \xymatrix{ & \bHom_{X}(x,x') \ar[dl] \ar[dr] & \\
\bHom_{X}(x,y') \ar[rr] & & \bHom_{X}(x,z). }$$
Using the two-out-of-three property, we are reduced to proving that both of the diagonal
arrows are homotopy equivalences. For the diagonal arrow on the left, this follows from our assumption that $\overline{e}'$ is locally $p$-Cartesian. For the arrow on the right, it suffices to show that the composition $g \circ \overline{e}'$ is locally $p$-coCartesian, which follows
from assumption $(\ast)$.
\end{proof}

\begin{corollary}\label{grutt1}
Let $p: X \rightarrow S$ be a Cartesian fibration of simplicial sets. An edge $e: x \rightarrow y$ of $X$ is $p$-coCartesian if and only if it is locally $p$-coCartesian (see the discussion preceding Proposition \ref{gotta}).
\end{corollary}

\begin{corollary}\label{getcocart}
Let $p: X \rightarrow S$ be a Cartesian fibration of simplicial sets. The following conditions are equivalent:
\begin{itemize}
\item[$(1)$] The map $p$ is a coCartesian fibration. 
\item[$(2)$] For every edge $f: s \rightarrow s'$ of $S$, the induced functor
$f^{\ast}: X_{s'} \rightarrow X_{s}$ has a left adjoint.
\end{itemize}
\end{corollary}

\begin{proof}
By definition, the functor corresponding to an edge $f: \Delta^1 \rightarrow S$ has a left adjoint if and only if the pullback $X \times_{S} \Delta^1 \rightarrow \Delta^1$ is a coCartesian fibration.
In other words, condition $(2)$ is equivalent to the assertion that for every
edge $f: s \rightarrow s'$ and every vertex $\widetilde{s}$ of $X$ lifting $s$, there
exists a {\em locally} $p$-coCartesian edge $\widetilde{f}: \widetilde{s} \rightarrow \widetilde{s}'$ lifting $f$. Using Corollary \ref{grutt1}, we conclude that $\widetilde{f}$ is automatically
$p$-coCartesian, so that $(2)$ is equivalent to $(1)$.
\end{proof}

\begin{proposition}\label{compadjoint}\index{gen}{adjoint functor!and composition}
Let $f: \calC \rightarrow \calD$ and $f': \calD \rightarrow \calE$ be functors between $\infty$-categories. Suppose that $f$ has a right adjoint $g$ and that $f'$ has a right adjoint $g'$. Then
$g \circ g'$ is right adjoint to $f' \circ f$.
\end{proposition}

\begin{proof}
Let $\phi$ denote the composable sequence of morphisms
$$ \calC \stackrel{g}{\leftarrow} \calD \stackrel{g'}{\leftarrow} \calE.$$
Let $M(\phi)$ denote the mapping simplex, and choose a factorization
$$ M(\phi) \stackrel{s}{\rightarrow} X \stackrel{q}{\rightarrow} \Delta^2$$
where $s$ is a quasi-equivalence and $X \rightarrow \Delta^2$ is a Cartesian fibration (using
Proposition \ref{sharpsimplex}). We first show that $q$ is a coCartesian fibration. In other words, we must show that for every object $\overline{x} \in \calC$ and every morphism $e: q(\overline{x}) \rightarrow y$, there is a $q$-Cartesian edge $\overline{e}: \overline{x} \rightarrow \overline{y}$ lifting $e$. This is clear if $e$ is degenerate. If $e = \Delta^{ \{0,1\} } \subseteq \Delta^2$, then
the existence of a left adjoint to $g$ implies that $e$ has a locally $q$-coCartesian lift $\overline{e}$. Lemma \ref{gruft} implies that $\overline{e}$ is $q$-coCartesian. Similarly,
if $e = \Delta^{ \{1,2\} }$, then we can find a $q$-coCartesian lift of $e$. Finally, if $e$ is the long edge $\Delta^{ \{0,2\} }$, then we may write $e$ as a composite $e' \circ e''$; the existence of a $q$-coCartesian lift of $e$ follows from the existence of $q$-coCartesian lifts of $e'$ and $e''$.
We now apply Proposition \ref{compass} and deduce that the adjunction $X \times_{\Delta^2} \Delta^{ \{0,2\} }$ is associated to both $g \circ g'$ and $f' \circ f$.
\end{proof}

In classical category theory, one can spell out the relationship between a pair of adjoint
functors $f: \calC \rightarrow \calD$ and $g: \calD \rightarrow \calC$ by specifying a {\it unit} transformation $\id_{\calC} \rightarrow g \circ f$ (or, dually, a {\it counit} $f \circ g \rightarrow \id_{\calD}$). This concept generalizes to the $\infty$-categorical setting as follows:

\begin{definition}
Suppose given a pair of functors
$$ \Adjoint{f}{\calC}{\calD}{g}$$
between $\infty$-categories. A {\it unit transformation} for $(f,g)$ is a morphism
$u: \id_{\calC} \rightarrow g \circ f$ in $\Fun(\calC,\calC)$ with the following property:
for every pair of objects $C \in \calC$, $D \in \calD$, the composition
$$ \bHom_{\calD}(f(C), D) \rightarrow \bHom_{\calC}(g(f(C)), g(D))
\stackrel{u(C)}{\rightarrow} \bHom_{\calC}(C, g(D))$$ is an isomorphism in the homotopy category $\calH$.\index{gen}{transformation!unit}\index{gen}{transformation!counit}\index{gen}{unit transformation}\index{gen}{counit transformation}
\end{definition}

\begin{proposition}\label{storut}\index{gen}{adjoint functor!and unit transformations}
Let $f: \calC \rightarrow \calD$ and $g: \calD \rightarrow \calC$ be a pair of functors between $\infty$-categories $\calC$ and $\calD$. The following conditions are equivalent:
\begin{itemize}
\item[$(1)$] The functor $f$ is a left adjoint to $g$.
\item[$(2)$] There exists a unit transformation $u: \id_{\calC} \rightarrow g \circ f$.
\end{itemize}
\end{proposition}

\begin{proof}
Suppose first that $(1)$ is satisfied.
Choose an adjunction $p: M \rightarrow \Delta^1$ which is associated to $f$ and $g$; according to $(1)$ of Proposition \ref{candi} we may identify $M_{ \{0\} }$ with $\calC$ and $M_{ \{1\} }$ with $\calD$. Since $f$ is associated to $M$, there is a map $F: \calC \times \Delta^1 \rightarrow M$
such that $F | \calC \times \{0\} = \id_{\calC}$ and $F| \calC \times \{1\} = f$, with each edge
$F| \{c\} \times \Delta^1$ $p$-coCartesian. Similarly, there is a map $G: \calD \times \Delta^1 \rightarrow M$ with $G | \calD \times \{1\} = \id_{\calD}$, $G| \calD \times \{0\} = g$, and 
$G| \{d\} \times \Delta^1$ is $p$-Cartesian for each object $d \in \calD$. Let 
$F' : \Lambda^2_2 \times \calC \rightarrow M$ be such that $F' | \Delta^{ \{0,2\} } \times \calC = F$ and $F' | \Delta^{ \{1,2\} } \times \calC = G \circ ( f \times \id_{\Delta^1} )$. Consider the diagram
$$ \xymatrix{ \Lambda^2_2 \times \calC \ar@{^{(}->}[d] \ar[r]^{F'} & M \ar[d] \\
\Delta^2 \times \calC \ar[r] \ar@{-->}[ur]^{F''} & \Delta^1. }$$ Using the fact $F' | \{c\} \times \Delta^{ \{1,2\} }$ is $p$-Cartesian for every object $c \in C$, we deduce the existence of the indicated dotted arrow $F''$. We now define $u = F' | \calC \times \Delta^{ \{0,1\} }$. We may regard $u$ as a natural transformation $\id_{\calC} \rightarrow g \circ f$. We claim that $u$ is a unit transformation.
 In other words, we must show that for any objects $C \in \calC$, $D \in \calD$, the composite map
$$ \bHom_{\calD}(fC, D) \rightarrow \bHom_{\calC}(gfC, gD)
\stackrel{u}{\rightarrow} \bHom_{\calC}(C,gD)$$
is an isomorphism in the homotopy category $\calH$ of spaces. This composite map
fits into a commutative diagram
$$ \xymatrix{ \bHom_{\calD}( f(C), D) \ar[r] \ar[d] &  \bHom_{\calD}(g(f(C)), g(D)) \ar[r] &
\bHom_{\calD}(C,g(D)) \ar[d] \\
\bHom_{M}(C,D) \ar[rr] & & \bHom_{M}(C,D). }$$
The left and right vertical arrows in this diagram are given by composition with
a $p$-coCartesian and a $p$-Cartesian morphism in $M$, respectively. Proposition \ref{compspaces} implies that these maps are homotopy equivalences.

We now prove that $(2) \Rightarrow (1)$. Choose a correspondence $p:M \rightarrow \Delta^1$ from
$\calC$ to $\calD$ which is associated to the functor $g$, via a map
$G: \calD \times \Delta^1 \rightarrow M$ as above. We have natural transformations
$$ \id_{\calC} \stackrel{u}{\rightarrow} g \circ f \stackrel{G \circ (f \times \id_{\Delta^1})}{\longrightarrow}
f.$$
Let $F: \calC \times \Delta^1 \rightarrow M$ be a composition of these transformations. We will complete the proof by showing that $F$ exhibits $M$ as a correspondence associated to the functor $f$. It will suffice to show that for each object $C \in \calC$, $F(C): C \rightarrow fC$ is 
$p$-coCartesian. According to Proposition \ref{charCart}, it will suffice to show that for each object $D \in \calD$, composition with $F(C)$ induces a homotopy equivalence
$ \bHom_{\calD}(f(C),D) \rightarrow \bHom_{M}(C,D)$. As above, this map fits into a commutative diagram 
$$ \xymatrix{ \bHom_{\calD}( f(C), D) \ar[r] \ar[d] &  \bHom_{\calD}(g(f(C)), g(D)) \ar[r] &
\bHom_{\calD}(C,g(D)) \ar[d] \\
\bHom_{M}(C,D) \ar[rr] & & \bHom_{M}(C,D) }$$
where the upper horizontal composition is an equivalence (since $u$ is a unit transformation)
and the right vertical arrow is an equivalence (since it is given by composition with a $p$-Cartesian morphism). It follows that the left vertical arrow is also a homotopy equivalence, as desired.
\end{proof}

\begin{proposition}\label{adjhom}
Let $\calC$ and $\calD$ be $\infty$-categories, and let $f: \calC \rightarrow \calD$ and
$g: \calD \rightarrow \calC$ be adjoint functors. Then $f$ and $g$ induce adjoint functors
$\h{f}: \h{\calC} \rightarrow \h{\calD}$ and $\h{g}:\h{\calD} \rightarrow \h{\calC}$ between $(${}$\calH$-enriched$)$ homotopy categories.
\end{proposition}

\begin{proof}
This follows immediately from Proposition \ref{storut}, since a unit transformation
$\id_{\calC} \rightarrow g \circ f$ induces a unit transformation
$\id_{\h{\calC}} \rightarrow (\h{g}) \circ (\h{f})$.
\end{proof}

The converse to Proposition \ref{adjhom} is false. 
If $f: \calC \rightarrow \calD$ and $g: \calD \rightarrow \calC$ are functors such that
$\h{f}$ and $\h{g}$ are adjoint to one another, then $f$ and $g$ are not necessarily adjoint.
Nevertheless, the existence of adjoints can be tested at the level of (enriched) homotopy categories.

\begin{lemma}\label{storkk}
Let $p: \calM \rightarrow \Delta^1$ be an inner fibration of simplicial sets, giving a correspondence
between the $\infty$-categories $\calC = \calM_{ \{0\} }$ and $\calD = \calM_{ \{1\} }$. Let $c$ be an object of $\calC$, $d$ an object of $\calD$, and $f: c \rightarrow d$ a morphism. The following are equivalent:
\begin{itemize}
\item[$(1)$] The morphism $f$ is $p$-Cartesian.
\item[$(2)$] The morphism $f$ gives rise to a Cartesian morphism in the enriched homotopy category $\h{\calM}$; in other words, composition with $p$ induces homotopy equivalences
$$ \bHom_{\calC }(c',c) \rightarrow \bHom_{\calM}(c',d)$$ for every
object $c' \in \calC$.
\end{itemize}
\end{lemma}

\begin{proof}
This follows immediately from Proposition \ref{charCart}.
\end{proof}

\begin{lemma}\label{storkkkk}
Let $p: \calM \rightarrow \Delta^1$ be an inner fibration, so that $\calM$ can be identified with a correspondence from $\calC = p^{-1} \{0\}$ to $\calD = p^{-1} \{1\}$. The following conditions are equivalent:
\begin{itemize}
\item[$(1)$] The map $p$ is a Cartesian fibration.
\item[$(2)$] There exists a $\calH$-enriched functor functor $g: \h{ \calD} \rightarrow \h{\calC}$ and a functorial identification $$\bHom_{\calM}(c,d) \simeq \bHom_{\calC}(c,g(d)).$$
\end{itemize}
\end{lemma}

\begin{proof}
If $p$ is a Cartesian fibration, then there is a functor $\calD \rightarrow \calC$ associated to $\calM$; we can then take $g$ to be the associated functor on enriched homotopy categories. Conversely, suppose that there exists a functor $g$ as above. We wish to show that $p$ is a Cartesian fibration. In other words, we must show that for every object $d \in \calD$, there is an object $c \in \calC$
and a $p$-Cartesian morphism $f: c \rightarrow d$. We take $c=g(d)$; in view of the identification
$\bHom_{\calM}(c,d) \simeq \bHom_{\calC}(c,c)$, there exists a morphism $f: c \rightarrow d$ corresponding to the identity $\id_{c}$. Lemma \ref{storkk} implies that $f$ is $p$-Cartesian, as desired.
\end{proof}

\begin{proposition}\index{gen}{adjoint functor!existence of}\label{sumpytump}
Let $f: \calC \rightarrow \calD$ be a functor between $\infty$-categories. Suppose that
the induced functor of $\calH$-enriched categories $\h{f}: \h{\calC} \rightarrow \h{\calD}$ admits a right adjoint. Then $f$ admits a right adjoint.
\end{proposition}

\begin{proof}
According to $(1)$ of Proposition \ref{candi}, there is a coCartesian fibration
$p: \calM \rightarrow \Delta^1$ associated to $f$. Let $\h{g}$ be the right adjoint of
$\h{f}$. Applying Lemma \ref{storkkkk}, we deduce that $p$ is a Cartesian fibration. Thus
$p$ is an adjunction, so that $f$ has a right adjoint as desired.
\end{proof}

\subsection{Preservation of Limits and Colimits}\label{afunc3}

Let $\calC$ and $\calD$ be ordinary categories, and let $F: \calC \rightarrow \calD$ be a functor.
If $F$ has a right adjoint $G$, then $F$ preserves colimits; we have a chain of natural isomorphisms
\begin{eqnarray*}
\Hom_{\calD}( F (\varinjlim C_{\alpha}), D) & \simeq & \Hom_{\calC}( \varinjlim C_{\alpha}, G(D)) \\
& \simeq & \varprojlim \Hom_{\calC}( C_{\alpha}, G(D)) \\
& \simeq & \varprojlim \Hom_{\calD}(F(C_{\alpha}), D) \\
& \simeq & \Hom_{\calD}( \varinjlim F(C_{\alpha}), D).
\end{eqnarray*}
In fact, this is in some sense the {\em defining} feature of left adjoints: under suitable set-theoretic assumptions, the {\em adjoint functor theorem} asserts that any colimit preserving functor admits a right adjoint. We will prove an $\infty$-categorical version of the adjoint functor theorem in \S \ref{aftt}. Our goal in this section is to lay the groundwork, by showing that left adjoints preserve colimits in the $\infty$-categorical setting. We will first need to establish several lemmas.

\begin{lemma}\label{lotusss}
Suppose given a diagram
$$ \xymatrix{ K \times \Delta^1 \ar[rr]^P \ar[dr] & & \calM \ar[dl]^q \\
& \Delta^1 & }$$
of simplicial sets, where $\calM$ is an $\infty$-category and $P| \{k\} \times \Delta^1$
is $q$-coCartesian for every vertex $k$ of $K$. Let $p = P| K \times \{0\}$. Then the induced map
$$ \psi: \calM_{P/} \rightarrow \calM_{p/}$$ induces a trivial fibration
$$ \psi_1: \calM_{P/} \times_{\Delta^1} \{1\} \rightarrow \calM_{p/} \times_{\Delta^1} \{1\}.$$
\end{lemma}

\begin{proof}
If $K$ is a point, then the assertion of the Lemma reduces immediately to the definition of a coCartesian edge. In the general case, we note that
$\psi$ and $\psi_1$ are both left fibrations between $\infty$-categories. Consequently, it suffices to show that $\psi_1$ is a categorical equivalence. In doing so, we are free to replace
$\psi$ by the equivalent map $\psi': \calM^{P/} \rightarrow \calM^{p/}$. To prove that $\psi'_1:
\calM^{P/} \times_{ \Delta^1} \{1\} \rightarrow \calM^{p/} \times_{\Delta^1} \{1\}$ is a trivial fibration, we must show that for every inclusion $A \subseteq B$ of simplicial sets and any map
$$ k_0: ((K \times \Delta^1) \diamond A) \coprod_{ (K \times \{0\}) \diamond A }
((K \times \{0\}) \diamond B) \rightarrow \calM$$
with $k_0| K \times \Delta^1 = P$ and $k_0(B) \subseteq q^{-1} \{1\}$, there exists an extension
of $k_0$ to a map $k: (K \times \Delta^1) \diamond B \rightarrow \calM$. Let
$$X = (K \times \Delta^1) \coprod_{K \times \Delta^1 \times B \times \{0\}} (K \times \Delta^1 \times B \times \Delta^1)$$ and let $h: X \rightarrow K \diamond B$ be the natural map. 
Let 
$$ X' = h^{-1} ((K \times \Delta^1) \diamond A) \coprod_{ (K \times \{0\}) \diamond A }
((K \times \{0\}) \diamond B) \subseteq X,$$
and let $\widetilde{k}_0: X' \rightarrow \calM$ be the composition $k_0 \circ h$. It suffices to
prove that there exists an extension of $\widetilde{k}_0$ to a map $\widetilde{k}: X \rightarrow \calM$.
Replacing $\calM$ by $\bHom_{\Delta^1}(K,\calM)$, we may reduce to the case where $K$ is a point, which we already treated above.
\end{proof}

\begin{lemma}\label{storka}
Let $q: \calM \rightarrow \Delta^1$ be a correspondence between $\infty$-categories
$\calC = q^{-1} \{0\}$ and $\calD = q^{-1} \{1\}$, and let $p: K \rightarrow \calC$ be a diagram in $\calC$. Let $f: c \rightarrow d$ be a $q$-Cartesian morphism in $\calM$ from $c \in \calC$
to $d \in \calD$. Let $r: \calM_{p/} \rightarrow \calM$ be the projection, and let
$\overline{d}$ be an object of $\calM_{p/}$ with $r(\overline{d})=d$. Then:
\begin{itemize}
\item[$(1)$] There exists a morphism $\overline{f}: \overline{c} \rightarrow \overline{d}$
in $\calM_{p/}$ satisfying $f = r( \overline{f} )$.
\item[$(2)$] Any morphism $\overline{f}: \overline{c} \rightarrow \overline{d}$ which
satisfies $r( \overline{f}) = f$ is $r$-Cartesian.
\end{itemize}
\end{lemma}

\begin{proof}
We may identify $\overline{d}$ with a map $\overline{d}: K \rightarrow \calM_{/d}$. Consider the set of pairs $( L, s )$ where $L \subseteq K$ and $s: L \rightarrow \calM_{/f}$ sits in a commutative diagram
$$ \xymatrix{ L \ar[r] \ar@{^{(}->}[d] & \calM_{/f} \ar[d] \\
K \ar[r] & \calM_{/d}. }$$
We order these pairs by setting $(L,s) \leq (L',s')$ if $L \subseteq L'$ and $s = s' | L$. By Zorn's lemma, there exists a pair $(L,s)$ which is maximal with respect to this ordering. To prove $(1)$, it suffices to show that $L = K$. Otherwise, we may obtain a larger simplicial subset $L' = L \coprod_{ \bd \Delta^n } \Delta^n \subseteq K$ by adjoining a single nondegenerate simplex. By maximality, there is no solution to the associated lifting problem
$$ \xymatrix{ \bd \Delta^n \ar[r] \ar@{^{(}->}[d] & \calM_{/f} \ar[d] \\
\Delta^n \ar[r] \ar@{-->}[ur] & \calM_{/d},}$$
nor to the associated lifting problem
$$ \xymatrix{ \Lambda^{n+2}_{n+2} \ar@{^{(}->}[d] \ar[r]^{s} & \calM \ar[d]^{q} \\
\Delta^{n+2} \ar[r] \ar@{-->}[ur] & \Delta^1, }$$
which contradicts the fact that $s$ carries $\Delta^{ \{n+1, n+2\} }$ to the $q$-Cartesian
morphism $f$ in $\calM$.

Now suppose that $\overline{f}$ is a lift of $f$. To prove that $\overline{f}$ is $r$-Cartesian, it suffices to show that for every $m \geq 2$ and every diagram
$$ \xymatrix{ \Lambda^{m}_m \ar[r]^{g_0} \ar@{^{(}->}[d] & \calM_{p/} \ar[d] \\
\Delta^m \ar[r] \ar@{-->}[ur]^{g} \ar[r] & \calM }$$
such that $g_0 | \Delta^{ \{m-1, m\} } = \widetilde{f}$, there exists a dotted arrow
$g$ as indicated, rendering the diagram commutative. We can identify the diagram with a map
$$t_0: ( K \star \Lambda^m_m) \coprod_{ \Lambda^m_m} \Delta^m \rightarrow \calM.$$
Consider the set of all pairs $(L, t)$, where $L \subseteq K$ and 
$$t: (K \star \Lambda^m_m) \coprod_{ L \star \Lambda^m_m } (L \star \Delta^m) \rightarrow \calM$$
is an extension of $t_0$. As above, we order the set of such pairs by declaring $(L,t) \leq (L',t')$ if $L \subseteq L'$
and $t = t' | L$. Zorn's lemma guarantees the existence of a maximal pair $(L,t)$. If $L = K$, we are done; otherwise let $L'$ be obtained from $L$ by adjoining a single nondegenerate $n$-simplex of $K$. By maximality, the map $t$ does not extend to $L'$; consequently the associated mapping problem
$$ \xymatrix{ (\Delta^n \star \Lambda^m_m) \coprod_{ \bd \Delta^n \star \Lambda^m_m }
(\bd \Delta^n \star \Delta^m) \ar[r] \ar@{^{(}->}[d] & \calM \ar[d] \\
\Delta^n \star \Delta^m \ar[r] \ar@{-->}[ur] & \Delta^1 }$$
has no solution. But this contradicts the assumption that $r(\widetilde{f})=f$ is a $q$-Cartesian edge of $\calM$.
\end{proof}

\begin{lemma}\label{lotusk}
Let $q: \calM \rightarrow \Delta^1$ be a correspondence between the $\infty$-categories
$\calC = q^{-1} \{0\}$ and $\calD = q^{-1} \{1 \}$. Let $f: c \rightarrow d$ be a morphism
in $\calM$ between objects $c \in \calC$, $d \in \calD$. Let $p: K \rightarrow \calC$
be a diagram, and consider an associated map
$$ k: \calM_{p/} \times_{\calM} \{c\} \rightarrow \calM_{p/} \times_{\calM} \{d\}$$
$($ the map $k$ is well-defined up to homotopy, according to Lemma \ref{functy} $)$.
If $f$ is $q$-Cartesian, then $k$ is a homotopy equivalence.
\end{lemma}

\begin{proof}
Let $X = (\calM_{p/})^{\Delta^1} \times_{\calM^{\Delta^1} } \{f\}$, and consider the diagram
$$ \xymatrix{ & X \ar[dr]^{u} \ar[dl]^{v} & \\
\calM_{p/} \times_{\calM} \{c\} & & \calM_{p/} \times_{\calM} \{d\}.}$$
The map $u$ is a homotopy equivalence, and $k$ is defined as the composition
of $v$ with a homotopy inverse to $u$. Consequently, it will suffice to show that
$v$ is a trivial fibration. To prove this, we must show that $v$ has the right lifting
property with respect to $\bd \Delta^n \subseteq \Delta^n$, which is equivalent to solving a lifting problem
$$ \xymatrix{ (\bd \Delta^n \times \Delta^1) \coprod_{ \bd \Delta^n \times \{1\} } (\Delta^n \times \{1\}) \ar[r] \ar@{^{(}->}[d] & \calM_{p/} \ar[d]^{r} \\
\Delta^n \times \Delta^1 \ar[r] \ar@{-->}[ur] & \calM}.$$
If $n=0$, we invoke $(1)$ of Lemma \ref{storka}. If $n > 0$, then Proposition \ref{goouse} implies
that it suffices to show that the upper horizontal map carries $\{n \} \times \Delta^1$
to an $r$-Cartesian edge of $\calM_{p/}$, which also follows from assertion $(2)$ of Lemma \ref{storka}.
\end{proof}

\begin{lemma}\label{lotuss}
Let $q: \calM \rightarrow \Delta^1$ be a Cartesian fibration, and let $\calC = q^{-1} \{0\}$.
The inclusion $\calC \subseteq \calM$ preserves all colimits which exist in $\calC$.
\end{lemma}

\begin{proof}
Let $\overline{p}: K^{\triangleright} \rightarrow \calC$ be a colimit of $p = \overline{p}|K$.
We wish to show that $\calM_{\overline{p}/} \rightarrow \calM_{p/}$ is a trivial fibration. Since
we have a diagram
$$ \calM_{\overline{p}/} \rightarrow \calM_{p/} \rightarrow \calM$$
of left fibrations, it will suffice to show that the induced map
$$ \calM_{\overline{p}/} \times_{\calM} \{d\} \rightarrow \calM_{p/} \times_{\calM} \{d\}$$ is a homotopy
equivalence of Kan complexes, for each object $d$ of $\calM$. If $d$ belongs to
$\calC$, this is obvious. In general, we may choose a $q$-Cartesian morphism
$f: c \rightarrow d$ in $\calM$. Composition with $f$ gives a commutative
diagram
$$ \xymatrix{ [\calM_{\overline{p}/} \times_{\calM} \{c\}] \ar[r] \ar[d] & [\calM_{p/} \times_{\calM} \{c\}] \ar[d] \\
[\calM_{\overline{p}/} \times_{\calM} \{d\}] \ar[r] & [\calM_{p/} \times_{\calM} \{d\}] }$$
in the homotopy category $\calH$ of spaces. The upper horizontal map is a homotopy equivalence since $\overline{p}$ is a colimit of $p$ in $\calC$. The vertical maps are homotopy equivalences by Lemma \ref{lotusk}. Consequently, the bottom horizontal map is also a homotopy equivalence, as desired.
\end{proof}

\begin{proposition}\label{adjointcol}\index{gen}{adjoint functor!and (co)limits}
Let $f: \calC \rightarrow \calD$ be a functor between $\infty$-categories which has a right
adjoint $g: \calD \rightarrow \calC$. Then $f$ preserves all colimits which exist in $\calC$, and $g$ preserves all limits which exist in $\calD$.
\end{proposition}

\begin{proof}
We will show that $f$ preserves colimits; the analogous statement for $g$ follows by a dual argument. Let $\overline{p}: K^{\triangleright} \rightarrow \calC$ be a colimit for
$p = \overline{p}|K$. We must show that $f \circ \overline{p}$ is a colimit of
$f \circ p$.

Let $q: \calM \rightarrow \Delta^1$ be an adjunction between $\calC = \calM_{\{0\}}$ and $\calD = \calM_{ \{1\} }$ which is associated to $f$ and $g$.  We wish to show that
$$\phi_1: \calD_{f \overline{p}/} \rightarrow \calD_{f p/}$$
is a trivial fibration. Since $\phi_1$ is a left fibration, it suffices to show that
$\phi_1$ is a categorical equivalence.

Since $\calM$ is associated to $f$, there is a map
$F: \calC \times \Delta^1 \rightarrow \calM$ with $F| \calC \times \{0\} = \id_{\calC}$, 
$F| \calC \times \{1\} = f$, and $F| \{c\} \times \Delta^1$ a $q$-coCartesian morphism of $\calM$
for every object $c \in \calC$. 
Let $\overline{P} = F \circ (\overline{p} \times \id_{\Delta^1})$
be the induced map $K^{\triangleright} \times \Delta^1 \rightarrow \calM$, and let
$P = \overline{P}|K \times \Delta^1$.

Consider the diagram
$$\xymatrix{ \calM_{\overline{p}/} \ar[r]^{\phi'} & \calM_{p/} \\
\calM_{\overline{P}/} \ar[r] \ar[u]^{\overline{v}} \ar[d]^{\overline{u}} & \calM_{P/} \ar[u]^{v} \ar[d]^{u}\\
\calM_{f \overline{p}/} \ar[r]^{\phi} & \calM_{f p/}. }$$
We note that every object in this diagram is an $\infty$-category with a map to $\Delta^1$; moreover, the map $\phi_{1}$ is obtained from $\phi$ by passage to the fiber over $\{1\} \subseteq \Delta^1$. Consequently, to prove that $\phi_1$ is a categorical equivalence, it suffices to verify three things:

\begin{itemize}
\item[$(1)$] The bottom vertical maps $u$ and $\overline{u}$ are trivial fibrations.
This follows from the fact that $K \times \{1\} \subseteq K \times \Delta^1$ and
$K^{\triangleright} \times \{ 1\} \subseteq K^{\triangleright} \times \Delta^1$ are right anodyne inclusions (Proposition \ref{sharpen2}).

\item[$(2)$] The upper vertical maps $v$ and $\overline{v}$ are trivial fibrations when restricted to
$\calD \subseteq \calM$. This follows from Lemma \ref{lotusss}, since $F$ carries
each $\{c\} \times \Delta^1$ to a $q$-coCartesian edge of $\calM$.

\item[$(3)$] The map $\phi'$ is a trivial fibration, since $\overline{p}$ is a colimit of
$p$ in $\calM$ according to Lemma \ref{lotuss}.

\end{itemize}

\end{proof}

\begin{remark}
Under appropriate set-theoretic hypotheses, one can prove a converse to Proposition \ref{adjointcol}. See Corollary \ref{adjointfunctor}. 
\end{remark}

\subsection{Examples of Adjoint Functors}\label{afunc4}

In this section, we describe a few simple criteria for establishing the existence of adjoint functors.

\begin{lemma}\label{adjfunclemma}
Let $q: \calM \rightarrow \Delta^1$ be a coCartesian fibration associated to a functor
$f: \calC \rightarrow \calD$, where $\calC = q^{-1} \{0\}$ and $\calD = q^{-1} \{1\}$. Let $D$ be an object of $\calD$. The following are equivalent:
\begin{itemize}
\item[$(1)$] There exists a $q$-Cartesian morphism $g: C \rightarrow D$ in $\calM$, where
$C \in \calC$.
\item[$(2)$] The right fibration $\calC \times_{\calD} \calD^{/D} \rightarrow \calC$ is
representable.
\end{itemize}
\end{lemma}

\begin{proof}
Let $F: \calC \times \Delta^1 \rightarrow \calM$ be
a $p$-coCartesian natural transformation from $\id_{\calC}$ to $f$. Define a simplicial
set $X$ so that for every simplicial set $K$, $\Hom_{\sSet}(K,X)$ parametrizes
maps $H: K \times \Delta^2 \rightarrow \calM$ such that
$h=H | K \times \{0\}$ factors through $\calC$, $H | K \times \Delta^{ \{0,1\} } = F \circ 
(h|(K \times \{0\}) \times \id_{\Delta^1})$, and $H | K \times \{2\}$ is the constant map at the vertex $D$. We have restriction maps
$$ \xymatrix{ & X \ar[dr] \ar[dl] & \\
\calC \times_{\calM} \calM^{/D} & & \calC \times_{\calD} \calD^{/D}.}$$
which are both trivial fibrations (the map on the right because $\calM$ is an $\infty$-category, the map on the left because $F$ is a $p$-coCartesian transformation). Consequently, $(2)$ is equivalent to 
the assertion that the $\infty$-category $\calC \times_{\calM} \calM^{/D}$ has a final object.
It now suffices to observe that a final object of $\calC \times_{\calM} \calM^{/D}$ is {\em precisely}
a $q$-Cartesian morphism $C \rightarrow D$, where $C \in \calC$.
\end{proof}

\begin{proposition}\label{adjfuncbaby}
Let $F: \calC \rightarrow \calD$ be a functor between $\infty$-categories.
The following are equivalent:
\begin{itemize}
\item[$(1)$] The functor $F$ has a right adjoint.
\item[$(2)$] For every pullback diagram
$$ \xymatrix{ \overline{\calC} \ar[r] \ar[d]^{p'} & \overline{\calD} \ar[d]^{p} \\
\calC \ar[r]^{F} & \calD, }$$
if $p$ is a representable right fibration, then $p'$ is also a representable right fibration.
\end{itemize}
\end{proposition}

\begin{proof}
Let $\calM$ be a correspondence from $\calC$ to $\calD$ associated to $F$, and apply
Lemma \ref{adjfunclemma} to each object of $D$.
\end{proof}

\begin{proposition}\label{quuquu}
Let $p: \calC \rightarrow \calD$ be a Cartesian fibration of $\infty$-categories, and let
$s: \calD \rightarrow \calC$ be a section of $p$ such that $s(D)$ is an initial object
of $\calC_{D} = \calC \times_{\calD} \{D\}$ for every object $D \in \calD$. Then
$s$ is a left adjoint of $p$.
\end{proposition}

\begin{proof}
Let $\calC^{0} \subseteq \calC$ denote the full subcategory of $\calC$ spanned by those
objects $C \in \calC$ such that $C$ is initial in the $\infty$-category $\calC_{p(C)}$. According to Proposition \ref{topaz}, the restriction $p| \calC^0$ is a trivial fibration from
$\calC^0$ to $\calD$. Consequently, it will suffice to show that the inclusion
$\calC^{0} \subseteq \calC$ is left adjoint to the composition $s \circ p: \calC \rightarrow \calC^{0}$.
Let $\calM \subseteq \calC \times \Delta^1$ be the full subcategory spanned by the vertices
$(C, \{i\})$ where $i = 1$ or $C \in \calC^0$. Let $q: \calM \rightarrow \Delta^1$ be the projection. It is clear that $q$ is a coCartesian fibration which is associated to the inclusion $\calC^0 \subseteq \calC$. To complete the proof, it will suffice to show that $q$ is also a Cartesian fibration
which is associated to $s \circ p$.

We first show that $q$ is a Cartesian fibration. It will suffice to show that for any object
$C \in \calC$, there is a $q$-Cartesian edge $(C',0) \rightarrow (C,1)$ in $\calM$. By assumption,
$C' = (s \circ p)(C)$ is an initial object of $\calC_{p(C)}$. Consequently, there exists
a morphism $f: C' \rightarrow C$ in $\calC_{p(C)}$; we will show that $f \times \id_{\Delta^1}$
is a $q$-Cartesian edge of $\calM$. To prove this, it suffices to show that for every $n \geq 2$ and every diagram
$$ \xymatrix{ \Lambda^n_n \ar@{^{(}->}[d] \ar[r]^{G_0} & \calM \ar[d] \\ 
\Delta^n \ar[r] \ar@{-->}[ur]^{G} & \Delta^1 }$$
such that $F_0 | \Delta^{ \{n-1, n\} } = f \times \id_{\Delta^1}$, there exists a dotted arrow $F: \Delta^n \rightarrow \calM$
as indicated, rendering the diagram commutative. We may identify $G_0$ with a map
$g_0: \Lambda^n_n \rightarrow \calC$. The composite map $p \circ g_0$ carries
$\Delta^{ \{n-1,n\} }$ to a degenerate edge of $\calD$, and therefore admits an extension
$\overline{g}: \Delta^n \rightarrow \calD$. Consider the diagram
$$ \xymatrix{ \Lambda^n_n \ar[r]^{g_0} \ar@{^{(}->}[d] & \calC \ar[d]^{p} \\
\Delta^n \ar[r]^{ \overline{g} } \ar@{-->}[ur]^{g} & \calD. }$$
Since $g_0$ carries the initial vertex $v$ of $\Delta^n$ to an initial object of the fiber
$\calC_{\overline{g}(v)}$, Lemma \ref{sabretooth} implies the existence of the indicated
map $g$ rendering the diagram commutative. This gives rise to a map $G: \Delta^n \rightarrow \calM$
with the desired properties, and completes the proof that $q$ is a Cartesian fibration.

We now wish to show that $s \circ p$ is associated to $q$. To prove this, it suffices to prove the existence of a map $H: \calC \times \Delta^1 \rightarrow \calC$ such that
$p \circ H = p \circ \pi_{\calC}$, $H | \calC \times \{1\} = \id_{\calC}$, and $H| C \times \{0\} = s \circ p$. We construct the map $H$ inductively, working cell-by-cell on $\calC$. Suppose that we have a nondegenerate simplex $\sigma: \Delta^n \rightarrow \calC$ and that $H$ has already been defined on $\sk^{n-1} \calC \times \Delta^1$. To define $H \circ (\sigma \times \id_{\Delta^1})$, we must solve a lifting problem that may be depicted as follows:
$$ \xymatrix{ (\bd \Delta^n \times \Delta^1) \coprod_{ \bd \Delta^n \times \bd \Delta^1} (\Delta^n \times \bd \Delta^1) \ar[rrrr]^{h_0} \ar@{^{(}->}[d] & & & & \calC \ar[d]^{p} \\
\Delta^n \times \Delta^1 \ar[rrrr] \ar@{-->}[urrrr]^{h} & & & & \calD.}$$
We now consider the filtration
$$ X(n+1) \subseteq X(n) \subseteq \ldots \subseteq X(0) = \Delta^n \times \Delta^1$$
defined in the proof of Proposition \ref{usejoyal}. Let $Y(i) = X(i) \coprod_{ \bd \Delta^n \times \{0\} } (\Delta^n \times \{1\})$. For $i > 0$, the inclusion $Y(i+1) \subseteq Y(i)$ is a pushout of the inclusion $X(i+1) \subseteq X(i)$, and therefore inner anodyne. Consequently, we may use the assumption that $p$ is an inner fibration to extend $h_0$ to a map defined on $Y(1)$. The
inclusion $Y(1) \subseteq \Delta^n \times \Delta^1$ is a pushout of $\bd \Delta^{n+1} \subseteq \Delta^{n+1}$; we then obtain the desired extension $h$ by applying Lemma \ref{sabretooth}.
\end{proof}

\begin{proposition}\label{simpex}
Let $\calM$ be a fibrant simplicial category equipped with a functor $p: \calM \rightarrow \Delta^1$ $($here we identify $\Delta^1$ with the two-object category whose nerve is $\Delta^1${}$)$, so that we may view $\calM$ as a correspondence between the simplicial categories $\calC = p^{-1} \{0\}$ and
$\calD = p^{-1} \{1\}$. The following are equivalent
\begin{itemize}
\item[$(1)$] The map $p$ is a Cartesian fibration.
\item[$(2)$] For every object $D \in \calD$, there exists a morphism $f: C \rightarrow D$ in $\calM$
which induces homotopy equivalences
$$ \bHom_{\calC}(C',C) \rightarrow \bHom_{\calM}(C',D)$$
for every $C' \in \calC$.
\end{itemize}
\end{proposition}

\begin{proof}
This follows immediately from Proposition \ref{trainedg}, since nonempty morphism spaces
in $\Delta^1$ are contractible.
\end{proof}

\begin{corollary}\label{sturk}
Let $\calC$ and $\calD$ be fibrant simplicial categories, and let 
$$\Adjoint{F}{\calC}{\calD}{G}$$
be a pair of adjoint functors
$F: \calC \rightarrow \calD$ $($ in the sense of enriched category theory, so that there is a natural isomorphism of simplicial sets $\bHom_{\calC}(F(C),D) \simeq \bHom_{\calD}(C,G(D))$ for $C \in \calC$, $D \in \calD$ $)$. Then the induced functors
$$ \Adjoint{f}{\Nerve(\calC)}{\Nerve(\calD)}{g}$$
are also adjoint to one another.
\end{corollary}

\begin{proof}
Let $\calM$ be the correspondence associated to the adjunction $(F,G)$. In other words, $\calM$ is a simplicial category containing $\calC$ and $\calD$ as full (simplicial) subcategories, with
$$\bHom_{\calM}(C,D) = \bHom_{\calC}(C,G(D)) = \bHom_{\calD}(F(C),D)$$
$$ \bHom_{\calM}(D,C) = \emptyset$$
for every pair of objects $C \in \calC$, $D \in \calD$. Let $M = \sNerve(\calM)$. Then
$M$ is a correspondence between $\sNerve(\calC)$ and $\sNerve(\calD)$. By Proposition
\ref{simpex}, it is an adjunction. It is easy to see that this adjunction is associated to both $f$ and $g$.
\end{proof}

The following variant on the situation of Corollary \ref{sturk} arises very often in practice:

\begin{proposition}\label{quiladj}\index{gen}{adjoint functor!Quillen}\index{gen}{Quillen adjunction}
Let $\bfA$ and $\bfA'$ be simplicial model categories, and let
$$ \Adjoint{F}{\bfA}{\bfA'}{G}$$
be a $($simplicial$)$ Quillen adjunction. Let $\calM$ be the simplicial category defined as in the proof of Corollary $\ref{sturk}$, and let $\calM^{\degree}$ be the full subcategory of $\calM$ consisting of those objects which are fibrant-cofibrant $($either as objects of $\bfA$ or as objects of $\bfA'${}$)$. Then $\sNerve(\calM^{\degree})$ determines an adjunction between $\sNerve( \bfA^{\degree})$ and $\sNerve( {\bfA'}^{\degree})$.
\end{proposition}

\begin{proof}
We need to show that $\sNerve(\calM^{\degree})\rightarrow \Delta^1$ is both a Cartesian fibration and a coCartesian fibration. We will argue the first point; the second follows from a dual argument.
According to Proposition \ref{princex}, it suffices to show that for every fibrant-cofibrant object
$D$ of $\bfA'$, there is a fibrant-cofibrant object $C$ of $\bfA$ and a morphism $f: C \rightarrow D$
in $\calM^{\degree}$ which induces weak homotopy equivalences
$$ \bHom_{ \bfA}(C',C) \rightarrow \bHom_{\calM}(C',D)$$ for every
fibrant-cofibrant object $C' \in \bfA$. We define $C$ to be a cofibrant replacement for $GD$: in other words, we choose a cofibrant object $C$ with a trivial fibration $C \rightarrow G(D)$ in the model category $\bfA$. Then $\bHom_{\bfA}(C',C) \rightarrow \bHom_{\calM}(C',D) = \bHom_{\bfA}(C',G(D))$ is a trivial fibration of simplicial sets, whenever $C'$ is a cofibrant object of $\bfA$.
\end{proof}

\begin{remark}
Suppose that $F: \bfA \rightarrow \bfA'$ and $G: \bfA' \rightarrow \bfA$ are as in Proposition \ref{quiladj}. We may associate to the adjunction $\sNerve(M^{\degree})$ a pair of adjoint functors $f: \sNerve( \bfA^{\degree}) \rightarrow \sNerve( {\bfA'}^{\degree})$ and $g: \sNerve( { \bfA'}^{\degree}) \rightarrow \sNerve( \bfA^{\degree} )$. In this situation, $f$
is often called a (nonabelian) {\it left derived functor} of $F$, and $g$ a (nonabelian) {\it right derived functor} of $G$. On the level of homotopy categories, $f$ and $g$ reduce to the usual derived functors associated to the Quillen adjunction (see \S \ref{quilladj}).\index{gen}{functor!derived}\index{gen}{derived functor}\index{gen}{left derived functor}\index{gen}{right derived functor}
\end{remark}

\subsection{Adjoint Functors and Overcategories}\label{afunc4half}

Our goal in this section is to prove the following result:

\begin{proposition}\label{curpse}
Suppose given an adjunction of $\infty$-categories
$$ \Adjoint{F}{\calC}{\calD}{G}.$$
Assume that the $\infty$-category $\calC$ admits pullbacks, and let $C$ be an object of $\calC$.
Then:
\begin{itemize}
\item[$(1)$] The induced functor $f: \calC^{/C} \rightarrow \calD^{/FC}$ admits a right adjoint $g$.
\item[$(2)$] The functor $g$ is equivalent to the composition
$$ \calD^{/FC} \stackrel{g'}{\rightarrow} \calC^{/GFC} \stackrel{g''}{\rightarrow} \calC^{/C},$$ 
where $g'$ is induced by $G$ and $g''$ is induced by pullback along the unit map
$C \rightarrow GFC$.
\end{itemize}
\end{proposition}

Proposition \ref{curpse} is an immediate consequence of the following more general result, which we will prove at the end of this section:

\begin{lemma}\label{starfi}
Suppose given an adjunction between $\infty$-categories
$$ \Adjoint{F}{\calC}{\calD}{G}.$$
Let $K$ be a simplicial set, and suppose given a pair of diagrams
$p_0: K \rightarrow \calC$, $p_1: K \rightarrow \calD$, and a natural
transformation $h: F \circ p_0 \rightarrow p_1$. 
Assume that $\calC$ admits pullbacks and $K$-indexed limits. Then:
\begin{itemize}
\item[$(1)$] 
Let
$f: \calC^{/p_0} \rightarrow \calD^{/p_1}$ denote the composition
$$ \calC^{/p_0} \rightarrow \calD^{/Fp_0} \stackrel{\circ 
\alpha}{\rightarrow}
\calD^{/p_1}.$$
Then $f$ admits a right adjoint $g$.
\item[$(2)$] The functor $g$ is equivalent the composition
$$ \calD^{/p_1} \stackrel{g'}{\rightarrow}
\calC^{/Gp_1} \stackrel{g''}{\rightarrow} \calC^{/p_0}.$$
Here $g''$ is induced by pullback along the natural transformation
$p_0 \rightarrow Gp_1$ adjoint to $h$ $($see below$)$.
\end{itemize}
\end{lemma}

We begin by recalling a bit of notation which will be needed in the proof. Suppose that $q: X \rightarrow S$ is an inner fibration of simplicial sets and
$p_S: K \rightarrow X$ is an arbitrary map, then we have defined a map of simplicial sets
$X^{/p_S} \rightarrow S$, which is characterized by the following universal property:
for every simplicial set $Y$ equipped with a map to $S$, there is a pullback diagram
$$ \xymatrix{ \Hom_{S}( Y, X^{/p_S}) \ar[r] \ar[d] & \Hom_{S}( Y \diamond_{S} S, X) \ar[d] \\
\{ p \} \ar[r] & \Hom_{S}(S, X). }$$ We refer the reader to \S \ref{consweet} for a more detailed discussion.



\begin{lemma}\label{curtwell}
Let $q: \calM \rightarrow \Delta^1$ be a coCartesian fibration of simplicial sets, classifying a functor
$F$ from $\calC = \calM \times_{ \Delta^1 } \{0\}$ to $\calD = \calM \times_{ \Delta^1 } \{1\}$.
Let $K$ be a simplicial set, and suppose given a commutative diagram
$$ \xymatrix{ K \times \Delta^1 \ar[dr] \ar[rr]^{g_{\Delta^1} } & & 
\ar[dl] \calM \\
& \Delta^1, & }$$
which restricts to give a pair of diagrams
$$ \calC \stackrel{g_0}{\leftarrow} K \stackrel{g_1}{\rightarrow} \calD.$$
 
Then:
\begin{itemize}
\item[$(1)$] The projection $q': \calM^{/g_{\Delta^1}} \rightarrow 
\Delta^1$ is a coCartesian fibration of simplicial sets, classifying a 
functor $F': \calC^{/g_0} \rightarrow \calD^{/g_1}$. Moreover, an
edge of $\calM^{/ g_{\Delta^1} }$ is $q'$-coCartesian if and only if its image in 
$\calM$ is $q$-coCartesian.

\item[$(2)$] Suppose that for every vertex $k$ in $K$, the map 
$g_{\Delta}$ carries $\{ k\} \times \Delta^1$
to a $q$-coCartesian morphism in $\calM$, so that $g_{\Delta^1}$ 
determines an equivalence $g_1 \simeq F \circ g_0$. Then $F'$ is homotopic 
to the composite functor
$$ \calC^{/g_0} \rightarrow \calD^{/ Fg_0} \simeq \calD^{/ g_1}.$$

\item[$(3)$] Suppose that $\calM = \calD \times \Delta^1$, and that $q$ is 
the projection onto the second factor, so that we can identify $F$ with 
the identity functor from $\calD$ to itself. Let $\overline{g}: K \times 
\Delta^1 \rightarrow \calD$ 
denote the composition $g_{\Delta^1}$ with the projection map
$\calM \rightarrow \calD$, so that we can regard
$\overline{g}$ as a morphism from $g_0$ to $g_1$ in $\Fun(K, \calD)$. Then 
the functor $F': \calD^{/g_0} \rightarrow \calD^{/g_1}$ is induced by 
composition with $\overline{g}$.

\end{itemize}
\end{lemma}

\begin{proof}
Assertion $(1)$ follows immediately from Proposition \ref{colimfam}.

We now prove $(2)$. Since $F$ is associated to the correspondence $\calM$, 
there exists a natural transformation $\alpha: \calC \times \Delta^1 
\rightarrow \calM$ from $\id_{\calC}$ to $F$, such that for each $C' \in 
\calC$, the induced map $\alpha_{C}: C' \rightarrow FC'$ is 
$q$-coCartesian. Without loss of generality, we may assume that 
$g_{\Delta^1}$ is given by the composition
$$ K \times \Delta^1 \stackrel{ g_0}{\rightarrow} \calC \times \Delta^1
\stackrel{\alpha}{\rightarrow} \calM.$$
In this case, $\alpha$ induces a map $\alpha': \calC^{/g_0} 
\times \Delta^1 \rightarrow \calM^{/ g_{\Delta^1}}$, which we may identify 
with a natural transformation from $\id_{\calC^{/g_0}}$ to the functor 
$\calC^{/g_0} \rightarrow \calD^{/Fg_0}$ determined by $F$. To show that 
this functor coincides with $F'$, it will suffice to show that $\alpha'$ 
carries each object of $\calC^{/g_0}$ to a $q'$-coCartesian morphism in 
$\calM^{ /g_{\Delta^1}}$. This follows immediately from the description of 
the $q'$-coCartesian edges given in assertion $(1)$.

We next prove $(3)$. Consider the diagram
$$ \calD^{ / g_0} \stackrel{p}{\leftarrow} \calD^{/\overline{g}} 
\stackrel{p'}{\rightarrow} \calD^{/g_1}.$$
By definition, ``composition with $\overline{g}$'' refers to a functor 
from $\calD_{/g_0}$ 
to $\calD^{/g_1}$ obtained by composing $p'$ with a section to the trivial 
fibration $p$.
To prove that this functor is homotopic to $F'$, it will suffice to show that
$F' \circ p$ is homotopic to $p'$. For this, we must produce a map
$\beta: \calD^{/\overline{g}} \times \Delta^1 \rightarrow \calM^{ / 
g_{\Delta^1}}$ from
$p$ to $p'$, such that $\beta$ carries each object of 
$\calD^{/\overline{g}}$ to a 
$q'$-coCartesian edge of $\calM^{/g_{\Delta^1}}$. We observe that
$\calD^{/\overline{g}} \times \Delta^1$ can be identified with 
$\calM^{/h_{\Delta^1}}$, where
$h: \Delta^1 \times \Delta^1 \rightarrow \calM \simeq \calD \times \Delta^1$ is the product
of $\overline{g}$ with the identity map. We now take $\beta$ to be the 
restriction map
$\calM^{ / h_{\Delta^1} } \rightarrow \calM^{ / g_{\Delta^1} }$ induced by the diagonal
inclusion $\Delta^1 \subseteq \Delta^1 \times \Delta^1$. Using $(1)$, we readily deduce that
$\beta$ has the desired properties.

\end{proof}

We will also need the following counterpart to Proposition \ref{colimfam}:

\begin{lemma}\label{limfam}
Suppose given a commutative diagram of simplicial sets
$$ \xymatrix{ K \times S \ar[r]^{p_S} \ar[dr] & X \ar[r]^{q} \ar[d] & Y \ar[dl] \\
& S & }$$
where the left diagonal arrow is projection onto the second factor, and
$q$ is an Cartesian fibration. Assume further that:
\begin{itemize}
\item[$(\ast)$] For every vertex $k \in K$, the map $p_S$ carries each
edge of $\{k \} \times S$ to a $q$-Cartsian edge in $X$.
\end{itemize}
Let $p'_S = q \circ p_S$. Then the map $q': X^{/p_S} \rightarrow Y^{/p'_S}$ is a Cartesian fibration. Moreover, an edge of $X^{/p_S}$ is $q'$-Cartesian if and only if its image in $X$ is $q$-Cartesian.
\end{lemma}

\begin{proof} To give the proof, it is convenient to use the language of 
marked simplicial sets (see \S \ref{twuf}). Let $X^{\natural}$ denote 
the marked simplicial set whose underlying simplicial set is $X$, where we 
consider an edge of $X^{\natural}$ to be marked if it is $q$-Cartesian. 
Let $\overline{X}^{\natural}$ denote the marked simplicial set whose 
underlying simplicial set is $X^{/p_S}$, where we consider an edge to be 
marked if and only if its image in $X$ is marked. According to Proposition 
\ref{dubudu}, it will suffice to show that the map 
$\overline{X}^{\natural} \rightarrow ( Y^{/p'_S} )^{\sharp}$ has the right 
lifting property with respect to every marked anodyne map $i: A 
\rightarrow B$. Let $\overline{A}$ and $\overline{B}$ denote the 
simplicial sets underlying $A$ and $B$, respectively. Suppose given a 
diagram of marked simplicial sets $$ \xymatrix{ A \ar@{^{(}->}[d]^{i} \ar[r] & 
\overline{X}^{\natural} \ar[d] \\ B \ar[r] \ar@{-->}[ur] & ( 
Y^{/p'_S})^{\sharp}. }$$ We wish to show that there exists a dotted arrow, 
rendering the diagram commutative. We begin by choosing a solution to the 
associated lifting problem $$ \xymatrix{ A \ar@{^{(}->}[d] \ar[r] & X^{\natural} 
\ar[d] \\ B \ar[r] \ar@{-->}[ur] & Y^{\sharp}, }$$ which is possible in 
view of our assumption that $q$ is a Cartesian fibration. To extend this 
to a solution to the original problem, it suffices to solve another 
lifting problem $$ \xymatrix{ (\overline{A} \times K \times \Delta^1) 
\coprod_{ (\overline{A} \times K \times \bd \Delta^1) } ( \overline{B} 
\times K \times \bd \Delta^1) \ar[r]^-{f} \ar@{^{(}->}[d]^{j} & X 
\ar[d]^{q} \\ 
\overline{B} \times K \times \Delta^1 \ar[r] \ar@{-->}[ur] & Y. }$$ By 
construction, the map $f$ induces a map of marked simplicial sets from $B 
\times K^{\flat} \times \{0\}$ to $X^{\natural}$. Using assumption 
$(\ast)$, we conclude that $f$ also induces a map of marked simplicial 
sets from $B \times K^{\flat} \times \{1\}$ to $X^{\natural}$. Using 
Proposition \ref{dubudu} again (and our assumption that $q$ is a 
Cartesian fibration), we are reduced to proving that the map $j$ induces a 
marked anodyne map $$ (A \times (K \times \Delta^1)^{\flat}) \coprod_{ A 
\times (K \times \bd \Delta^1)^{\flat} } ( B \times (K \times \bd 
\Delta^1)^{\flat}) \rightarrow B \times (K \times \Delta^1)^{\flat}.$$ 
Since $i$ is marked anodyne by assumption, this follows immediately from 
Proposition \ref{markanodprod}. \end{proof}

\begin{lemma}\label{spitfure} Let $q: \calM \rightarrow \Delta^1$ be a 
Cartesian fibration of simplicial sets, associated to a functor $G$ from 
$\calD = \calM \times_{ \Delta^1} \{1\}$ to $\calC = \calM 
\times_{\Delta^1} \{0\}$. Suppose given a simplicial set $K$ 
and a commutative diagram
$$ \xymatrix{ K \times \Delta^1 \ar[rr]^{g_{\Delta^1} } \ar[dr] & & \calM 
\ar[dl] \\
& \Delta^1, & }$$
so that $g_{\Delta^1}$ restricts to a pair of functors
$$ \calC \stackrel{g_0}{\leftarrow} K \stackrel{g_1}{\rightarrow} \calD.$$
Suppose furthermore that, for every vertex $k$ of $K$, the corresponding 
morphism $g_0(k) \rightarrow g_1(k)$ is $q$-Cartesian. Then:
\begin{itemize}
\item[$(1)$] The induced map $q': \calM^{ / f_{\Delta^1} } \rightarrow \Delta^1$ is a Cartesian fibration.
Moreover, an edge of $\calM^{/ f_{\Delta^1}}$ is $q'$-Cartesian if and only if its image in
$\calM$ is $q$-Cartesian.
\item[$(2)$] The associated functor $\calD^{/g_1} \rightarrow 
\calC^{/g_0}$ is 
homotopic to the composition of the functor $G': \calD^{/g_1} \rightarrow 
\calC^{/Gg_1}$ induced by $G$ and the equivalence
$\calC^{/Gg_1} \simeq \calC^{/g_0}$ determined by the map $g_{\Delta^1}$.
\end{itemize}
\end{lemma}

\begin{proof} Assertion $(1)$ follows immediately from Lemma \ref{limfam}. 
We will prove $(2)$. Since the functor $G$ is associated to $q$, there 
exists a map $\alpha: \calD \times \Delta^1 \rightarrow \calM$ which is a 
natural transformation from $G$ to $\id_{\calD}$, such that for every 
object $D \in \calD$ the induced map $\alpha_{D}: \{D \} \times \Delta^1 
\rightarrow \calM$ is a $q$-Cartesian edge of $\calM$. Without loss of 
generality, we may assume that $g$
coincides with the composition
$$ K \times \Delta^1 \stackrel{g_1}{\rightarrow} \calD \times \Delta^1
\stackrel{\alpha}{\rightarrow} \calM.$$
In this case, $\alpha$ 
induces a map $\alpha': \calD^{/g_1} \times \Delta^1 \rightarrow 
\calM^{/f_{\Delta^1}}$, which is a natural transformation from $G'$ to the 
identity. Using $(1)$, we deduce that $\alpha'$ carries each object of 
$\calD^{/g_1}$ to a $q'$-Cartesian edge of $\calM^{/ f_{\Delta^1}}$. It 
follows that $\alpha'$ exhibits $G'$ as the functor associated to the 
Cartesian fibration $q'$, as desired. \end{proof}

\begin{proof}[Proof of Lemma \ref{starfi}]
Let $q: \calM \rightarrow \Delta^1$ be a correspondence from $\calC = 
\calM \times_{\Delta^1} \{0\}$ to
$\calD = \calM \times_{\Delta^1} \{1\}$, which is associated to the pair of adjoint functors $F$ and $G$.
The natural transformation $h$ determines a map
map $\alpha$ determines a map
$\alpha: K \times \Delta^1 \rightarrow \calM$, which is a natural 
transformation from $p_0$ to $p_1$. Using the fact that $q$ is both a 
Cartesian and a coCartesian fibration, we can form a commutative
square $\sigma$:
$$ \xymatrix{ & Gp_1 \ar[dr]^{\phi} & \\
p_0 \ar[ur] \ar[rr]^{\alpha} \ar[dr]^{\psi} & & p_1 \\
& Fp_0 \ar[ur] & }$$
in the $\infty$-category $\Fun(K, \calM)$, where the morphism $\phi$ is
$q$-Cartesian and the morphism $\psi$ is $q$-coCartesian.

Let $\calN = \calM \times \Delta^1$. We can identify $\sigma$ with a map
$\sigma_{\Delta^1 \times \Delta^1}: K \times \Delta^1 \times \Delta^1 
\rightarrow \calM \times \Delta^1$. Let $\calN' = \calN^{/ 
\sigma_{\Delta^1 \times \Delta^1}}$. Proposition \ref{colimfam} 
implies that the projection $\calN' \rightarrow \Delta^1 \times \Delta^1$ 
is a coCartesian fibration, associated to some diagram of 
$\infty$-categories
$$ \xymatrix{ & \calC^{/Gp_1} \ar[dr]^{f'} & \\
\calC^{/p_0} \ar[rr] \ar[dr] \ar[ur]^{f''} & & \calD^{/p_1} \\
& \calD^{/Fp_0} \ar[ur]. & }$$
Lemma \ref{curtwell} allows us to identify the functors in the lower 
triangle, so we see that the horizontal composition is homotopic to the 
functor $f$. To complete the proof of $(1)$, it will suffice to show that
the functors $f'$ and $f''$ admit right adjoints. To prove $(2)$, it 
suffices to show that those right adjoints are given by $g'$ and $g''$, 
respectively. The adjointness of $f'$ and $g'$ follows from Lemma 
\ref{spitfure}.

It follows from Lemma \ref{curtwell} that the functor $f'': \calC^{/p_0}
\rightarrow \calC^{/Gp_1}$ is given by composition with the transformation
$h': p_0 \rightarrow Gp_1$ which is adjoint to $h$.
The pullback functor 
$g''$ 
is right adjoint to $f''$ by definition; the only nontrivial point is to 
establish the existence of $g''$. Here we must use our hypotheses on the 
$\infty$-category $\calC$. Let $\overline{p}_0: K^{\triangleleft} 
\rightarrow \calC$ be a limit of $p_0$, let $\overline{Gp}_1: 
K^{\triangleleft} \rightarrow \calC$ be a limit of $Gp_1$. Let us identify 
$h'$ with a map $K \times \Delta^1 \rightarrow \calC$, and choose an 
extension $\overline{h}': K^{\triangleleft} \times \Delta^1 \rightarrow 
\calC$ which is a natural transformation from $\overline{p}_0$ to
$\overline{Gp}_1$. Let $C \in \calC$ denote the image under
$\overline{p}_0$ of the cone point of $K^{\triangleleft}$, let
$C' \in \calC$ denote the image under $\overline{Gp}_1$ of the cone point
of $K^{\triangleleft}$, and let $j: C \rightarrow C'$ be the morphism
induced by $\overline{h}'$. We have a commutative diagram of 
$\infty$-categories:
$$ \xymatrix{ \calC^{/ p_0} & \calC^{ / h' } \ar[r]^{f''_1} \ar[l]^{f''_0} 
& \calC^{/ 
Gp_1} \\
\calC^{/ \overline{p}_0 } \ar[u] \ar[d] & \calC^{/ \overline{h}'} \ar[l] 
\ar[r] \ar[u] \ar[d] & \calC^{/ \overline{Gp}_1} \ar[u] \ar[d] \\
\calC^{/C} & \calC^{/ j} \ar[r] \ar[l] & \calC^{ /C '}. } $$
In this diagram, the left horizontal arrows are trivial Kan fibrations, as 
are all of the vertical arrows. The functor $f''$ is obtained by composing
$f''_0$ with a section to the trivial Kan fibration $f''_1$. Utilizing
the vertical equivalences, we can identify $f''$ with the functor
$\calC^{/C} \rightarrow \calC^{/C'}$ given by composition with $j$.
But this functor admits a right adjoint, in view of our assumption that
$\calC$ admits pullbacks.
\end{proof}

\subsection{Uniqueness of Adjoint Functors}\label{afunc5}

We have seen that if $f: \calC \rightarrow \calD$ is a functor which admits a right adjoint
$g: \calD \rightarrow \calC$, then $g$ is uniquely determined up to homotopy. Our next result is a slight refinement of this assertion. 

\begin{definition}
Let $\calC$ and $\calD$ be $\infty$-categories. We let $\LFun(\calC, \calD)
\subseteq \Fun(\calC, \calD)$ denote the full subcategory of $\Fun(\calC,\calD)$ spanned by those functors $F: \calC \rightarrow \calD$ which are left adjoints. Similarly, we define
$\RFun(\calC,\calD)$ to be the full subcategory of $\Fun(\calC, \calD)$ spanned by those functors which are right adjoints.\index{not}{FunL@$\LFun(\calC, \calD)$}\index{not}{FunR@$\RFun(\calC, \calD)$}
\end{definition}

\begin{proposition}\label{switcheroo}
Let $\calC$ and $\calD$ be $\infty$-categories. Then the $\infty$-categories
$\LFun(\calC, \calD)$ and $\RFun(\calD, \calC)^{op}$ are $($canonically$)$ equivalent to one another.
\end{proposition}

\begin{proof}
Enlarging the universe if necessary, we may assume without loss of generality that $\calC$ and $\calD$ are small. Let $j: \calD \rightarrow \calP(\calD)$ be the Yoneda embedding.
Composition with $j$ induces a fully faithful embedding
$$ i: \Fun(\calC,\calD) {\rightarrow} \Fun(\calC, \calP(\calD)) \simeq
\Fun( \calC \times \calD^{op}, \SSet).$$
The essential image of $i$ consists of those functors $G: \calC \times \calD^{op} \rightarrow \SSet$
with the property that, for each $C \in \calC$, the induced functor
$G_{C}: \calD^{op} \rightarrow \SSet$ is representable by an object $D \in \calD$. 
The functor $i$ induces a fully faithful embedding
$$ i_0: \RFun(\calC, \calD) \rightarrow \Fun(\calC \times \calD^{op}, \SSet)$$
whose essential image consists of those functors $G$ which belong to the essential image
of $i$, and furthermore satisfy the additional condition that for each $D \in \calD$, the induced
functor $G_{D}: \calC \rightarrow \SSet$ is corepresentable by an object $C \in \calC$ (this
follows from Proposition \ref{adjfuncbaby}). Let $\calE \subseteq \Fun(\calC \times \calD^{op}, \SSet)$
be the full subcategory spanned by those functors which satisfy these two conditions, so that
the Yoneda embedding induces an equivalence
$$ \RFun(\calC,\calD) \rightarrow \calE.$$
We note that the above conditions are self-dual, so that the same reasoning gives an equivalence of $\infty$-categories
$$ \RFun(\calD^{op}, \calC^{op}) \rightarrow \calE.$$ 
We now conclude by observing that there is a natural equivalence of $\infty$-categories
$\RFun(\calD^{op}, \calC^{op}) \simeq \LFun( \calD, \calC)^{op}$.
\end{proof}

We will later need a slight refinement of Proposition \ref{switcheroo}, which exhibits some functoriality in $\calC$. We begin with a few preliminary remarks concerning the construction of presheaf $\infty$-categories.

Let $f: \calC \rightarrow \calC'$ be a functor between small $\infty$-categories. Then composition with $f$ induces a restriction functor $G: \calP(\calC') \rightarrow \calP(\calC)$. However, there
is another slightly less evident functoriality of the construction $\calC \mapsto \calP(\calC)$. Namely, according to Theorem \ref{charpresheaf}, there is a colimit-preserving functor $\calP(f): \calP(\calC) \rightarrow \calP(\calC')$, uniquely determined up to equivalence, such that the diagram
$$ \xymatrix{ \calC \ar[d] \ar[r]^{f} & \calC' \ar[d] \\
\calP(\calC) \ar[r]^{ \calP(f) } & \calP(\calC') } $$\index{not}{Pcalf@$\calP(f)$}
commutes up to homotopy (here the vertical arrows are given by the Yoneda embeddings).

The functor $\calP(f)$ has an alternative characterization in the language of adjoint functors:

\begin{proposition}\label{adjobs}
Let $f: \calC \rightarrow \calC'$ be a functor between small $\infty$-categories, let
$G: \calP(\calC') \rightarrow \calP(\calC)$ be the functor given by composition with $f$.
Then $G$ is right adjoint to $\calP(f)$.
\end{proposition}

\begin{proof}
We first prove that $G$ admits a left adjoint. Let $e: \calP(\calC) \rightarrow \Fun( \calP(\calC)^{op}, \widehat{\SSet})$ denote the Yoneda embedding.
According to Proposition \ref{adjfuncbaby}, it will suffice to show that for each
$M \in \calP(\calC)$, the composite functor $e(M) \circ G$ is corepresentable.
Let $\calD$ denote the full subcategory of $\calP(\calC)$ spanned by those objects $M$ such that $G \circ e_M$ is corepresentable. Since $\calP(\calC)$ admits small colimits, Proposition \ref{yonedaprop}
implies that the collection of corepresentable functors on $\calP(\calC)$ is stable under small colimits. According to Propositions \ref{yonedaprop} and \ref{limiteval}, the functor
$M \mapsto e(M) \circ G$ preserves small colimits. It follows that $\calD$ is stable under small colimits in $\calP(\calC)$. Since $\calP(\calC)$ is generated under small colimits by the Yoneda embedding $j_{\calC'}: \calC \rightarrow \calP(\calC)$ (Corollary \ref{gencolcot}), it will suffice to show that $j_{\calC}(C) \in \calD$ for each $C \in \calC$. According to 
Lemma \ref{repco}, $e(j_{\calC}(C))$ is equivalent to the functor $\calP(\calC) \rightarrow \widehat{\SSet}$ given by evaluation at $C$. Then $e( j_{\calC}(C) ) \circ G$ is equivalent to the functor given by evaluation at $f(C) \in \calC'$, which is corepresentable (Lemma \ref{repco} again).
We conclude that $G$ has a left adjoint $F$.

To complete the proof, we must show that $F$ is equivalent to $\calP(f)$. To prove this, it will suffice to show that $F$ preserves small colimits and that the diagram
$$ \xymatrix{ \calC \ar[r]^{f} \ar[d] & \calC' \ar[d] \\
\calP(\calC) \ar[r]^{F} & \calP(\calC') }$$
commutes up to homotopy. The first point is obvious: since $F$ is a left adjoint, it preserves all colimits which exist in $\calP(\calC)$ (Proposition \ref{adjointcol}). For the second, choose a counit map $v: F \circ G \rightarrow \id_{\calP(\calC')}$. By construction, the functor $f$ induces
a natural transformation $u: j_{\calC} \rightarrow G \circ j_{\calC'} \circ f$. To complete
the proof, it will suffice to show that the composition
$$ \theta: F \circ j_{\calC} \stackrel{u}{\rightarrow} F \circ G \circ j_{\calC'} \circ f
\stackrel{v}{\rightarrow} j_{\calC'} \circ f$$
is an equivalence of functors from $\calC$ to $\calP(\calC')$. Fix objects
$C \in \calC$, $M \in \calP(\calC')$. We have a commutative diagram
$$ \xymatrix{ \bHom_{\calP(\calC')}( j_{\calC'}(f(C)), M) \ar[r] \ar@{=}[d]  
& \bHom_{\calP(\calC)}( G(j_{\calC'}(f(C))), G(M)) \ar[r] \ar[d] & 
\bHom_{\calP(\calC)}( j_{\calC}(C), G(M) ) \ar[d] \\
\bHom_{\calP(\calC')}( j_{\calC'}(f(C)), M) \ar[r] & 
\bHom_{\calP(\calC')}( F(G( j_{\calC'}( f(C) ))), M) \ar[r] & \bHom_{\calP(\calC')}(F(j_{\calC}(C)), M) }$$
in the homotopy category $\calH$ of spaces, where the vertical arrows are isomorphisms. 
Consequently, to prove that the lower horizontal composition is an isomorphism, it suffices to prove that the upper horizontal composition is an isomorphism. Using Lemma \ref{repco}, we reduce to the assertion that $M( f(C)) \rightarrow (G(M))(C)$ is an isomorphism in $\calH$, which follows immediately from the definition of $G$.
\end{proof}

\begin{remark}\label{switcheroo2}
Suppose given a functor $f: \calD \rightarrow \calD'$
which admits a right adjoint $g$. 
Let $\calE \subseteq \Fun(\calC \times \calD^{op}, \SSet)$
and $\calE' \subseteq \Fun(\calC \times (\calD')^{op}, \SSet)$
be defined as in the proof of Proposition \ref{switcheroo}, and consider the diagram
$$ \xymatrix{ \RFun(\calC, \calD) \ar[d]^{\circ g} \ar[r] & \calE \ar[d] & \LFun(\calD, \calC)^{op} \ar[d]^{\circ f} \ar[l] \\
\RFun(\calC, \calD') \ar[r] & \calE' & \LFun(\calD', \calC)^{op}. \ar[l] }$$
Here the middle vertical map is given by composition with $\id_{\calC} \times f$. The square on the right is manifestly commutative, but the square on the left commutes only up to homotopy. To verify the second point, we observe that the square in question is given by applying the functor
$\bHom(\calC, \bigdot)$ to the diagram
$$ \xymatrix{ \calD \ar[r] \ar[d]^{g} & \calP(\calD) \ar[d]^{G} \\
\calD' \ar[r] & \calP(\calD') }$$
where $G$ is given by composition with $f$ and the horizontal arrows are given by the Yoneda embedding. Let $\calP^0(\calD) \subseteq \calP(\calD)$ and $\calP^0(\calD')$ denote the essential images of the Yoneda embeddings. Proposition \ref{adjfuncbaby} asserts that $G$ carries $\calP^{0}(\calD')$ into $\calP^{0}(\calD)$, so that it will suffice to verify that the
diagram
$$ \xymatrix{ \calD \ar[r] \ar[d]^{g} & \calP^0(\calD) \ar[d]^{G^0} \\
\calD' \ar[r] & \calP^0(\calD') }$$
is homotopy commutative. In view of Proposition \ref{compadjoint}, it will suffice to show that
$G^0$ admits a left adjoint $F^0$ and that the diagram
$$ \xymatrix{ \calD \ar[r]  & \calP^0(\calD) \\
\calD' \ar[r] \ar[u] & \calP^0(\calD') \ar[u]^{F_0} }$$
is homotopy commutative. 
According to Proposition \ref{adjobs}, the functor $G$ has a left adjoint $\calP(f)$ which fits into a commutative diagram $$ \xymatrix{ \calD \ar[r] & \calP(\calD) \\
\calD' \ar[r] \ar[u]^{f} & \calP(\calD') \ar[u]^{\calP(f)}. }$$ 
In particular, $\calP(f)$ carries $\calP^0(\calD)$ into $\calP^{0}(\calD')$ and therefore
restricts to give a left adjoint $F^0: \calP^{0}(\calD) \rightarrow \calP^{0}(\calD')$ which verifies the desired commutativity.
\end{remark}

We conclude this section by establishing the following consequence of
Proposition \ref{adjobs}: 

\begin{corollary}\label{coughspaz}
Let $\calC$ be a small $\infty$-category, $\calD$ a locally small $\infty$-category which admits small colimits. Let $F: \calP(\calC) \rightarrow \calD$ be a colimit-preserving functor, let
$f: \calC \rightarrow \calD$ denotes the composition of $F$ with the Yoneda embedding of
$\calC$, and let $G: \calD \rightarrow \calP(\calC)$ be the functor given by the composition
$$ \calD \stackrel{j'}{\rightarrow} \Fun( \calD^{op}, \SSet) \stackrel{\circ f}{\rightarrow} \calP(\calC).$$
Then $G$ is a right adjoint to $F$. Moreover, the map 
$$f = F \circ j \rightarrow (F \circ (G \circ F)) \circ j = (F \circ G) \circ f$$
exhibits $F \circ G$ as a left Kan extension of $f$ along itself.
\end{corollary}

The proof requires a few preliminaries:

\begin{lemma}\label{stimp}
Suppose given a pair of adjoint functors
$$ \Adjoint{f}{\calC}{\calD}{g}$$
between $\infty$-categories. Let $T: \calC \rightarrow \calX$ be any functor. Then
$T \circ g: \calD \rightarrow \calX$ is a left Kan extension of $T$ along $f$.
\end{lemma}

\begin{proof}
Let $p: \calM \rightarrow \Delta^1$ be a correspondence associated to the pair of adjoint functors
$f$ and $g$. Choose a $p$-Cartesian homotopy $h$ from $r$ to $\id_{M}$, where
$r$ is a functor from $\calM$ to $\calC$; thus $r | \calD$ is homotopic to $g$. 
It will therefore suffice to show that the composition
$$ \overline{T}: \calM \stackrel{r}{\rightarrow} \calC \stackrel{T}{\rightarrow} \calX$$
is a left Kan extension of $\overline{T} | \calC \simeq T$. For this, we must show that
for each $D \in \calD$, the functor $\overline{T}$ exhibits $\overline{T}(D)$
as a colimit of the diagram
$$ (\calC \times_{\calM} \calM_{/D}) \rightarrow \calM \stackrel{\overline{T}}{\rightarrow} \calX.$$
We observe that $\calC \times_{\calM} \calM_{/D}$ has a final object, given by any
$p$-Cartesian morphism $e: C \rightarrow D$. It therefore suffices to show that
$\overline{T}(e)$ is an equivalence in $\calX$, which follows immediately from the construction of $\overline{T}$.
\end{proof}

\begin{lemma}\label{kanspaz}
Let $f: \calC \rightarrow \calC'$ be a functor between small $\infty$-categories and
$\calX$ an $\infty$-category which admits small colimits.
Let $H: \calP(\calC) \rightarrow \calX$ be a functor which preserves small colimits, and 
$h: \calC \rightarrow \calX$ the composition of $F$ with the Yoneda embedding $j_{\calC}: \calC \rightarrow \calP(\calC)$. Then the composition 
$$\calC' \stackrel{j_{\calC'}}{\rightarrow} \calP(\calC') \stackrel{ \circ f}{\rightarrow} \calP(\calC)
\stackrel{H}{\rightarrow} \calX$$
is a left Kan extension of $h$ along $f$.
\end{lemma}

\begin{proof}
Let $G: \calP(\calC') \rightarrow \calP(\calC)$ be the functor given by composition with $f$.
In view of Proposition \ref{acekan}, it will suffice to show that $H \circ G$ is a 
left Kan extension of $h$ along $j_{\calC} \circ f$.
 
Theorem \ref{charpresheaf} implies the existence of a functor
$F: \calP(\calC) \rightarrow \calP(\calC')$ which preserves small colimits, such that
$F \circ j_{\calC} \simeq j_{\calC'} \circ f$. Moreover, Lemma \ref{longwait1} ensures that
$F$ is a left Kan extension of $f$ along the fully faithful Yoneda embedding $j_{\calC}$. 
Using Proposition \ref{acekan} again, we are reduced to proving that
$H \circ G$ is a left Kan extension of $H$ along $F$. This follows immediately from Proposition \ref{adjobs} and Lemma \ref{stimp}.
\end{proof}

\begin{proof}[Proof of Corollary \ref{coughspaz}]
The first claim follows from Proposition \ref{adjobs}. To prove the second, we may assume without loss of generality that $\calD$ is minimal, so that $\calD$ is union of small full subcategories
$\{ \calD_{\alpha} \}$. It will suffice to show that, for each index $\alpha$ such that $f$ factors
through $\calD_{\alpha}$, the restricted transformation $f \rightarrow ((F \circ G)| \calD_{\alpha}) \circ f$ exhibits $(F \circ G)| \calD_{\alpha}$ as a left Kan extension of
$f$ along the induced map $\calC \rightarrow \calD_{\alpha}$, which follows from Lemma \ref{kanspaz}.
\end{proof}

\subsection{Localization Functors}\label{locfunc}

Suppose we are given a $\infty$-category $\calC$ and a collection $S$ of morphisms
of $\calC$ which we would like to invert. In other words, we wish to find an $\infty$-category $S^{-1} \calC$ equipped with a functor $\eta: \calC \rightarrow S^{-1} \calC$ which carries each morphism in $S$ to an equivalence, and is in some sense universal with respect to these properties.
One can give a general construction of $S^{-1} \calC$ using the formalism of \S \ref{bicat1}. Without loss of generality, we may suppose that $S$ contains all the identity morphisms in $\calC$. Consequently, the pair $(\calC, S)$ may be regarded as a marked simplicial set, and we can choose a marked anodyne map $( \calC, S) \rightarrow ( S^{-1} \calC, S')$, where $S^{-1} \calC$ is an $\infty$-category and $S'$ is the collection of all equivalences in $S^{-1} \calC$. However, this construction is generally very difficult to analyze, and the properties of $S^{-1} \calC$ are very difficult to control. For example, it might be the case that $\calC$ is locally small and $S^{-1} \calC$ is not.

Under suitable hypotheses on $S$ (see \S \ref{invloc}), there is a drastically simpler approach: we can find the desired $\infty$-category $S^{-1} \calC$ {\em inside} of $\calC$, as the full subcategory of {\em $S$-local} objects of $\calC$.\index{not}{SinvC@$S^{-1} \calC$}

\begin{example}\label{excom}
Let $\calC$ be the (ordinary) category of abelian groups, $p$ a
prime number, and let $S$ denote the collection of morphisms $f$
whose kernel and cokernel consist entirely of $p$-power torsion elements. A
morphism $f$ lies in $S$ if and only if it induces an isomorphism
after inverting the prime number $p$. In this case, we may
identify $S^{-1} \calC$ with the full subcategory of $\calC$ consisting of those abelian groups which
are {\em uniquely $p$-divisible}. The functor $\calC \rightarrow S^{-1} \calC$ is given
by $$ M \mapsto M \otimes_{\Z} \Z[ \frac{1}{p}].$$
\end{example}

In Example \ref{excom}, the functor $\calC \rightarrow
S^{-1} \calC$ is actually left adjoint to an inclusion functor. We will take this as our starting point.

\begin{definition}\label{swagga}
A functor $f: \calC \rightarrow \calD$ between $\infty$-categories is a {\it localization}
if $f$ has a fully faithful right adjoint.\index{gen}{localization}\index{gen}{functor!localization}
\end{definition}

\begin{warning}
Let $f: \calC \rightarrow \calD$ be a localization functor, and let $S$ denote the collection of all morphisms $\alpha$ in $\calC$ such that $f(\alpha)$ is an equivalence. Then, for any
$\infty$-category $\calE$, composition with $f$ induces a fully faithful functor
$$ \Fun( \calD, \calE) \stackrel{ \circ f}{\rightarrow} \Fun(\calC, \calE)$$
whose essential image consists of those functors $p: \calC \rightarrow \calE$ which
carry each $\alpha \in S$ to an equivalence in $\calE$ (Proposition \ref{unlap}). We may describe the situation more
informally by saying that $\calD$ is obtained from $\calC$ by inverting the morphisms of $S$.

Some authors use the term ``localization'' in a more general sense, to describe any
functor $f: \calC \rightarrow \calD$ in which $\calD$ is obtained by inverting some collection $S$ of morphisms in $\calC$. Such a morphism $f$ need not be a localization in the sense of Definition \ref{swagga}; however, it is in many cases (see Proposition \ref{local}).
\end{warning}

If $f: \calC \rightarrow \calD$ is a localization of $\infty$-categories, then we will also say that $\calD$ is a {\it localization} of $\calC$. In this case, a right adjoint $g: \calD \rightarrow \calC$
of $f$ gives an equivalence between $\calD$ and a full subcategory of $\calC$ (the essential image of $g$). We let $L: \calC \rightarrow \calC$ denote the composition $g \circ f$. We will abuse terminology by referring to $L$ as a {\it localization functor} if it arises in this way.

The following result will allow us to recognize localization functors:

\begin{proposition}\label{recloc}
Let $\calC$ be a $\infty$-category, and let $L: \calC \rightarrow \calC$ be a functor with essential
image $L\calC \subseteq \calC$. The following conditions are equivalent:
\begin{itemize}
\item[$(1)$] There exists a functor $f: \calC \rightarrow \calD$ with a fully faithful 
right adjoint $g: \calD \rightarrow \calC$, and an equivalence between $g \circ f$ and $L$.

\item[$(2)$] When regarded as a functor from $\calC$ to $L\calC$, $L$ is a left adjoint of the
inclusion $L\calC \subseteq \calC$.

\item[$(3)$] There exists a natural transformation $\alpha: \calC \times \Delta^1 \rightarrow \calC$ from $\id_{\calC}$ to $L$ such that, for every object $C$ of $\calC$, the morphisms
$L(\alpha(C)), \alpha(LC): LC \rightarrow LLC$ of $\calC$ are equivalences.
\end{itemize}
\end{proposition}

\begin{proof}
It is obvious that $(2)$ implies $(1)$ (take $\calD= L\calC$, $f = L$, and $g$ to be the inclusion).
The converse follows from the observation that, since $g$ is fully faithful, we are free to replace $\calD$ by the essential image of $g$ (which is equal to the essential image of $L$).

We next prove that $(2)$ implies $(3)$. Let $\alpha: \id_{\calC} \rightarrow L$ be a unit
for the adjunction. Then, for each pair of objects $C \in \calC$, $D \in L \calC$, composition with $\alpha(C)$ induces a homotopy equivalence
$$ \bHom_{\calC}( LC, D) \rightarrow \bHom_{\calC}( C, D),$$
and in particular a bijection $\Hom_{\h{\calC}}(LC, D) \rightarrow \Hom_{\h{\calC}}(C,D)$.
If $C$ belongs to $L \calC$, then Yoneda's lemma implies that
$\alpha(C)$ is an isomorphism in $\h{\calC}$. This proves that $\alpha(LC)$ is an equivalence for every $C \in \calC$. Since $\alpha$ is a natural transformation, we have obtain a diagram
$$ \xymatrix{ C \ar[r]^{\alpha(C)} \ar[d]^{\alpha(C)} & LC \ar[d]^{L \alpha(C)}  \\
LC \ar[r]^{\alpha(LC)} & LLC.}$$
Since composition with $\alpha(C)$ gives an injective map
from $\Hom_{\h{\calC}}(LC,LLC)$ to $\Hom_{\h{\calC}}(C,LLC)$, we conclude that
$\alpha(LC)$ is homotopic to $L \alpha(C)$; in particular, $\alpha(LC)$ is also an equivalence. This proves $(3)$.

Now suppose that $(3)$ is satisfied; we will prove that $\alpha$ is the unit of an adjunction between $\calC$ and $L \calC$. In other words, we must show that for each $C \in \calC$ and
$D \in \calC$, composition with $\alpha(C)$ induces a homotopy equivalence 
$$ \phi: \bHom_{\calC}( LC, LD) \rightarrow \bHom_{\calC}( C, LD).$$
By Yoneda's lemma, it will suffice to show that for every Kan complex $K$, the induced map
$$ \Hom_{\calH}( K, \bHom_{\calC}(LC, LD)) \rightarrow \Hom_{\calH}(K, \bHom_{\calC}(C,LD))$$
is a bijection of sets, where $\calH$ denotes the homotopy category of spaces. Replacing
$\calC$ by $\Fun(K, \calC)$, we are reduced to proving the following:

\begin{itemize}
\item[$(\ast)$] Suppose that $\alpha: \id_{\calC} \rightarrow L$ satisfies $(3)$. Then, for every
pair of objects $C, D \in \calC$, composition with $\alpha(C)$ induces a bijection of sets
$$ \phi: \Hom_{ \h{\calC}}( LC, LD) \rightarrow \Hom_{ \h{\calC}}(C, LD).$$
\end{itemize}

We first show that $\phi$ is surjective. Let $f$ be a morphism from $C$ to $LD$. We then have a commutative diagram
$$ \xymatrix{ C \ar[r]^{f} \ar[d]^{\alpha(C)} & LD \ar[d]^{\alpha(LD)} \\
LC \ar[r]^{Lf} & LLD, }$$
so that $f$ is homotopic to the composition $(\alpha(LD)^{-1} \circ Lf) \circ \alpha(C)$; this proves that the homotopy class of $f$ lies in the image of $\phi$.

We now show that $\phi$ is injective. Let $g: LC \rightarrow LD$ be an arbitrary morphism. 
We have a commutative diagram
$$ \xymatrix{ LC \ar[r]^{g} \ar[d]^{\alpha(LC)} & LD \ar[d]^{ \alpha(LD) } \\
LLC \ar[r]^{ Lg } & LLD, }$$
so that $g$ is homotopic to the composition
\begin{eqnarray*}
\alpha(LD)^{-1} \circ Lg \circ \alpha(LC) & 
\simeq & \alpha(LD)^{-1} \circ Lg \circ L \alpha(C) \circ (L\alpha(C))^{-1} \circ \alpha(LC) \\
& \simeq & \alpha(LD)^{-1} \circ L( g \circ \alpha(C) ) \circ (L \alpha(C))^{-1} \circ \alpha(LC)
\end{eqnarray*}
In particular, $g$ is determined by $g \circ \alpha(C)$ up to homotopy.

\end{proof}

\begin{remark}\label{localcolim}\index{gen}{localization!and colimits}
Let $L: \calC \rightarrow \calD$ be a localization functor and $K$ a simplicial set. Suppose
that every diagram $p: K \rightarrow \calC$ admits a colimit in $\calC$. Then the $\infty$-category $\calD$ has the same property. Moreover, we can give an explicit prescription for computing colimits in $\calD$. Let $q: K \rightarrow \calD$ be a diagram, and let
$p: K \rightarrow \calC$ be the composition of $q$ with a right adjoint to $L$. 
Choose a colimit $\overline{p}: K^{\triangleright} \rightarrow \calC$. Since $L$ is a left
adjoint, $L \circ \overline{p}$ is a colimit diagram in $\calD$, and $L \circ p$ is equivalent
to the diagram $q$.
\end{remark}

We conclude this section by introducing a few ideas which will allow us to recognize localization functors, when they exist.

\begin{definition}\label{locaobj}\index{gen}{localization!of an object}
Let $\calC$ be an $\infty$-category and $\calC^{0} \subseteq \calC$ a full subcategory. We will say that a morphism $f: C \rightarrow D$ in $\calC$ {\it exhibits $D$ as a $\calC^0$-localization of $C$} if $D \in \calC^{0}$, and composition with $f$ induces an isomorphism
$$ \bHom_{\calC^{0}}(D,E) \rightarrow \bHom_{\calC}(C,E)$$ in the homotopy category $\calH$, for each object $E \in \calC^{0}$.
\end{definition}

\begin{remark}\label{initrem}
In the situation of Definition \ref{locaobj}, a morphism $f: C \rightarrow D$ exhibits
$D$ as a localization of $C$ if and only if $f$ is an initial object of the $\infty$-category
$\calC^{0}_{C/} = \calC_{C/} \times_{\calC} \calC^{0}$. In particular, $f$ is uniquely determined up to equivalence.
\end{remark}

\begin{proposition}\label{testreflect}
Let $\calC$ be an $\infty$-category and $\calC^0 \subseteq \calC$ a full subcategory.
The following conditions are equivalent:
\begin{itemize}
\item[$(1)$] For every object $C \in \calC$, there exists a localization
$f: C \rightarrow D$ relative to $\calC^{0}$.
\item[$(2)$] The inclusion $\calC^{0} \subseteq \calC$ admits a left adjoint.
\end{itemize}
\end{proposition}

\begin{proof}
Let $\calD$ be the full subcategory of $\calC \times \Delta^1$ spanned by objects
of the form $(C,i)$, where $C \in \calC^0$ if $i=1$. Then the projection $p: \calD \rightarrow \Delta^1$ is a correspondence from $\calC$ to $\calC^0$ which is associated to the inclusion functor $i: \calC^0 \subseteq \calC$. It follows that $i$ admits a left adjoint if and only if
$p$ is a coCartesian fibration. It now suffices to observe that if $C$ is an object of $\calC$, then we may identify $p$-coCartesian edges $f: (C,0) \rightarrow (D,1)$ of $\calD$ with
localizations $C \rightarrow D$ relative to $\calC^{0}$.
\end{proof}

\begin{remark}\label{reflective}\index{gen}{subcategory!reflective}\index{gen}{reflective subcategory}
By analogy with classical category theory, we will say that a full subcategory $\calC^0$ of an $\infty$-category $\calC$ is a {\it reflective subcategory} if the hypotheses of Proposition \ref{testreflect} are satisfied by the inclusion $\calC^0 \subseteq \calC$. 
\end{remark}

\begin{example}\label{swang}
Let $\calC$ be an $\infty$-category which has a final object, and let $\calC^{0}$ be the full subcategory of $\calC$ spanned by the final objects. Then the inclusion $\calC^{0} \subseteq \calC$ admits a left adjoint.
\end{example}

\begin{corollary}\label{sweng}
Let $p: \calC \rightarrow \calD$ be a coCartesian fibrations between $\infty$-categories, let
$\calD^{0} \subseteq \calD$ be a full subcategory, and let $\calC^{0} = \calC \times_{\calD} \calD^{0}$.
If the inclusion $\calD^{0} \subseteq \calD$ admits a left adjoint, then the inclusion
$\calC^{0} \subseteq \calC$ admits a left adjoint. 
\end{corollary}

\begin{proof}
In view of Proposition \ref{testreflect}, it will suffice to show that for every object $C \in \calC$, there
a morphism $f: C \rightarrow C_0$ which is a localization of $C$ relative to $\calC^{0}$. Let
$D = p(C)$, let $\overline{f}: D \rightarrow D_0$ be a localization of $D$ relative to $\calD_0$, and let
$f: C \rightarrow C_0$ be a $p$-coCartesian morphism in $\calC$ lifting $\overline{f}$. We claim that
$f$ has the desired property. Choose any object $C' \in \calC^{0}$, and let $D' = p(C') \in \calD^{0}$.
We obtain a diagram of spaces
$$ \xymatrix{ \bHom_{\calC}( C_0, C') \ar[r]^{\phi} \ar[d] &  \bHom_{\calC}(C, C') \ar[d] \\
\bHom_{\calD}( D_0, D') \ar[r]^{\psi} & \bHom_{\calD}( D, D') }$$
which commutes up to preferred homotopy. By assumption, the map $\psi$ is a homotopy equivalence. Since $f$ is $p$-coCartesian, the map $\phi$ induces a homotopy equivalence after passing to the homotopy fibers over any pair of points $\eta \in \bHom_{\calD}(D_0, D')$, $\psi(\eta)
\in \bHom_{\calD}(D,D')$. Using the long exact sequence of homotopy groups associated to the vertical fibrations, we conclude that $\phi$ is a homotopy equivalence, as desired.
\end{proof}

\begin{proposition}\label{unlap}
Let $\calC$ be an $\infty$-category, let $L: \calC \rightarrow \calC$ be a localization functor with essential image $L\calC$. Let $S$ denote the collection of all morphisms $f$ in $\calC$ such that $Lf$ is an equivalence. Then, for every $\infty$-category $\calD$, composition with
$f$ induces a fully faithful functor
$$ \psi: \Fun( L \calC , \calD) \rightarrow \Fun( \calC, \calD).$$
Moreover, the essential image of $\psi$ consists of those functors $F: \calC \rightarrow \calD$
such that $F(f)$ is an equivalence in $\calD$, for each $f \in S$.
\end{proposition}

\begin{proof}
Let $S_0$ be the collection of all morphisms $C \rightarrow D$ in $\calC$ which exhibit
$D$ as an $L\calC$-localization of $C$. We first claim that, for any functor
$F: \calC \rightarrow \calD$, the following conditions are equivalent:
\begin{itemize}
\item[$(a)$] The functor $F$ is a right Kan extension of $F | L\calC$.
\item[$(b)$] The functor $F$ carries each morphism in $S_0$ to an equivalence in $\calD$.
\item[$(c)$] The functor $F$ carries each morphism in $S$ to an equivalence in $\calD$.
\end{itemize}

The equivalence of $(a)$ and $(b)$ follows immediately from the definitions
(since a morphism $f: C \rightarrow D$ exhibits $D$ as an $L \calC$-localization
of $C$ if and only if $f$ is an initial object of $(L \calC) \times_{\calC} \calC_{C/}$ ), and the
implication $(c) \Rightarrow (b)$ is obvious. To prove that $(b) \Rightarrow (c)$, let us
consider any map $f: C \rightarrow D$ which belongs to $S$. We have a commutative diagram
$$ \xymatrix{ C \ar[r]^{f} \ar[d] & LC \ar[d]^{f} \\
D \ar[r] & LD. }$$
Since $f \in S$, the map $Lf$ is an equivalence in $\calC$. If $F$ satisfies $(b)$, then $F$ carries
each of the horizontal maps to an equivalence in $\calD$. It follows from the two-out-of-three property that $Ff$ is an equivalence in $\calD$ as well, so that $F$ satisfies $(c)$.

Let $\Fun^{0}(\calC, \calD)$ denote the full subcategory of $\Fun(\calC, \calD)$ spanned by those functors which satisfy $(a)$, $(b)$, and $(c)$. Using Proposition \ref{lklk}, we deduce
that the restriction functor $\phi: \Fun^{0}(\calC, \calD) \rightarrow \Fun( L \calC, \calD)$ is fully faithful.
We now observe that $\psi$ is a right homotopy inverse to $\phi$. It follows that
$\phi$ is essentially surjective, and therefore an equivalence. Being right homotopy inverse to an equivalence, the functor $\psi$ must itself be an equivalence.
\end{proof}

\subsection{Factorization Systems}\label{factgen1}

Let $f: X \rightarrow Z$ be a map of sets. Then $f$ can be written as a composition
$$ X \stackrel{f'}{\rightarrow} Y \stackrel{f''}{\rightarrow} Z$$
where $f'$ is surjective and $f''$ is injective. This factorization is uniquely determined up to (unique) isomorphism: the set $Y$ can be characterized either as the image of the map $f$, or as the quotient
of $X$ by the equivalence relation $R = \{ (x,y) \in X^2: f(x) = f(y) \}$. We can describe the situation formally by saying that the collections of surjective and injective maps form a {\em factorization system} on the category $\Set$ of sets (see Definition \ref{spanhun}). In this section, we will describe a theory of factorization systems in the $\infty$-categorical setting. These ideas are due to Joyal and we refer the reader to \cite{joyalnotpub} for further details.

\begin{definition}\label{defperp}\index{gen}{orthogonal}\index{gen}{orthogonal!left}\index{gen}{orthogonal!right}\index{gen}{left orthogonal}\index{gen}{right orthogonal}
Let $f: A \rightarrow B$ and $g: X \rightarrow Y$ be morphisms in an $\infty$-category $\calC$.
We will say that $f$ is {\it left orthogonal} to $g$, or that $g$ is {\it right orthogonal to $f$},
if the following condition is satisfied:
\begin{itemize}
\item[$(\ast)$] For every commutative diagram
$$ \xymatrix{ A \ar[r] \ar[d]^{f} & X \ar[d]^{g} \\
B \ar[r] & Y }$$
in $\calC$, the mapping space $\bHom_{\calC_{A/ \, /Y}}(B,X )$
is contractible. (Here we abuse notation by identifying $B$ and $X$ with the corresponding objects
of $\calC_{A/\,/Y}$.)
\end{itemize}
In this case, we will write $f \perp g$.\index{not}{fperpg@$f \perp g$}
\end{definition}

\begin{remark}
More informally: a morphism $f: A \rightarrow B$ in an $\infty$-category $\calC$ is left orthogonal to another morphism $g: X \rightarrow Y$ if, for every commutative diagram
$$ \xymatrix{ A \ar[d]^{f} \ar[r] & X \ar[d]^{g} \\
B \ar[r] \ar@{-->}[ur] & Y, }$$
the space of dotted arrows rendering the diagram commutative is contractible.
\end{remark}

\begin{remark}\label{spack}
Let $f: A \rightarrow B$ and $g: X \rightarrow Y$ be morphisms in an $\infty$-category $\calC$.
Fix a morphism $A \rightarrow Y$, which we can identify with an object
$\overline{Y} \in \calC_{A/}$. Lifting $g: X \rightarrow Y$ to an object of
$\widetilde{X} \in \calC_{A/ \, /Y}$ is equivalent to lifting $g$ to a morphism
$\overline{g}: \overline{X} \rightarrow \overline{Y}$ in $\calC_{A/}$.
The map $f: A \rightarrow B$ determines an object $\overline{B} \in \calC_{A/}$, 
and lifting $f$ to an object $\widetilde{B} \in \calC_{A/ \, /Y}$ is equivalent to giving a map 
$h: \overline{B} \rightarrow \overline{Y}$ in $\calC_{A/}$. We therefore have a fiber
sequence of spaces
$$ \bHom_{ \calC_{ A/ \, /Y} }( \widetilde{B}, \widetilde{X} )
\rightarrow \bHom_{\calC_{A/}}( \overline{B}, \overline{X} )
\rightarrow \bHom_{\calC_{A/}}( \overline{B}, \overline{Y} ),$$
where the fiber is taken over the point $h$. Consequently, condition
$(\ast)$ of Definition \ref{defperp} can be reformulated as follows:
for every morphism $\overline{g}: \overline{X} \rightarrow \overline{Y}$ in
$\calC_{A/}$ lifting $g$, composition with $\overline{g}$ induces a homotopy equivalence
$$ \bHom_{\calC_{A/}}( \overline{B}, \overline{X} ) \rightarrow \bHom_{\calC_{A/}}( \overline{B}, \overline{Y}).$$
\end{remark}

\begin{notation}\index{not}{Sperp@$S^{\perp}$}\index{not}{perpS@$^{\perp}S$}
Let $\calC$ be an $\infty$-category, and let $S$ be a collection of morphisms in $\calC$. We let
$S^{\perp}$ denote the collection of all morphisms in $\calC$ which are right orthogonal to $S$, and
$^{\perp}S$ the collection of all morphisms in $\calC$ which are left orthogonal to $S$.
\end{notation}

\begin{remark}
Let $\calC$ be an ordinary category containing a pair of morphisms $f$ and $g$. If $f \perp g$, then
$f$ has the left lifting property with respect to $g$, and $g$ has the right lifting property with respect to $f$. It follows that for any collection $S$ of morphisms in $\calC$, we have inclusions
$S^{\perp} \subseteq S_{\perp}$ and $^{\perp} S \subseteq _{\perp}S$, where the latter classes of morphisms are defined in \S \ref{liftingprobs}.
\end{remark}

Applying Remark \ref{spack} to an $\infty$-category $\calC$ and its opposite, we obtain the following result:

\begin{proposition}\label{swimmm}
Let $\calC$ be an $\infty$-category and $S$ a collection of morphisms in $\calC$.
\begin{itemize}
\item[$(1)$] The sets of morphisms $S^{\perp}$ and $^{\perp}S$ contain every equivalence in $\calC$.
\item[$(2)$] The sets of morphisms $S^{\perp}$ and $^{\perp}S$ are closed under the formation of retracts.
\item[$(3)$] Suppose given a commutative diagram
$$ \xymatrix{ & Y \ar[dr]^{g} & \\
X \ar[ur]^{f} \ar[rr]^{h} & & Z }$$
in $\calC$, where $g \in S^{\perp}$. Then $f \in S^{\perp}$ if and only if $h \in S^{\perp}$. In particular, $S^{\perp}$ is closed under composition.

\item[$(4)$] Suppose given a commutative diagram
$$ \xymatrix{ & Y \ar[dr]^{g} & \\
X \ar[ur]^{f} \ar[rr]^{h} & & Z }$$
in $\calC$, where $f \in {}^{\perp}S$. Then $g \in {}^{\perp}S$ if and only if $h \in {}^{\perp}S$. In particular, $^{\perp}S$ is closed under composition.

\item[$(5)$] The set of morphisms $S^{\perp}$ is stable under pullbacks: that is, given a pullback diagram
$$ \xymatrix{ X' \ar[d]^{g'} \ar[r] & X \ar[d]^{g} \\
Y' \ar[r] & Y }$$
in $\calC$, if $g$ belongs to $S^{\perp}$, then $g'$ belongs to $S^{\perp}$.
\item[$(6)$] The set of morphisms $^{\perp}S$ is stable under pushouts: that is, given a pushout diagram
$$ \xymatrix{ A \ar[d]^{f} \ar[r] & A' \ar[d]^{f'} \\
B \ar[r] & B', }$$
if $f$ belongs to $^{\perp}S$, then so does $f'$.
\item[$(7)$] Let $K$ be a simplicial set such that $\calC$ admits $K$-indexed colimits.
Then the full subcategory of $\Fun( \Delta^1, \calC)$ spanned by the elements of
$^{\perp}S$ is closed under $K$-indexed colimits.
\item[$(8)$] Let $K$ be a simplicial set such that $\calC$ admits $K$-indexed limits. Then the full subcategory of $\Fun( \Delta^1, \calC)$ spanned by the elements of $S^{\perp}$ is closed under $K$-indexed limits.
\end{itemize}
\end{proposition}

\begin{remark}\label{smule}
Suppose given a pair of adjoint functors $\Adjoint{F}{\calC}{\calD.}{G}$
Let $f$ be a morphism in $\calC$ and $g$ a morphism in $\calD$. Then
$f \perp G(g)$ if and only if $F(f) \perp g$.
\end{remark}

\begin{definition}[Joyal]\label{spanhun}\index{gen}{factorization system}
Let $\calC$ be an $\infty$-category. A {\it factorization system} on $\calC$ is a pair
$(S_L, S_R)$, where $S_L$ and $S_R$ are collections of morphisms of $\calC$ which satisfy the following axioms:
\begin{itemize}
\item[$(1)$] The collections $S_L$ and $S_R$ are stable under the formation of retracts.
\item[$(2)$] Every morphism in $S_L$ is left orthogonal to every morphism in $S_R$.
\item[$(3)$] For every morphism $h: X \rightarrow Z$ in $\calC$, there exists a commutative
triangle
$$ \xymatrix{ & Y \ar[dr]^{g} & \\
X \ar[ur]^{f} \ar[rr]^{h} & & Z }$$
where $f \in S_L$ and $g \in S_R$. 
\end{itemize}
We will call $S_L$ the {\it left set} of the factorization system, and $S_R$ the
{\it right set} of the factorization system.
\end{definition}

\begin{example}\label{scumm}
Let $\calC$ be an $\infty$-category. Then $\calC$ admits a factorization system
$(S_L, S_R)$, where $S_L$ is the collection of all equivalences in $\calC$, and $S_R$ consists of all morphisms of $\calC$.
\end{example}

\begin{remark}\label{spill}
Let $(S_L, S_R)$ be a factorization system on an $\infty$-category $\calC$. Then
$(S_R, S_L)$ is a factorization system on the opposite $\infty$-category $\calC^{op}$. 
\end{remark}

\begin{proposition}\label{swin}
Let $\calC$ be an $\infty$-category, and let $(S_L, S_R)$ be a factorization system on $\calC$. 
Then $S_L = {^{\perp}S_R}$ and $S_R = S_L^{\perp}$.
\end{proposition}

\begin{proof}
By symmetry, it will suffice to prove the first assertion. The inclusion $S_L \subseteq {^{\perp}S_R}$ follows immediately from the definition. To prove the reverse inclusion, let us suppose that $h: X \rightarrow Z$ is a morphism in $\calC$ which is left orthogonal to every morphism in $S_R$. Choose a commutative triangle
$$ \xymatrix{ & Y \ar[dr]^{g} & \\
X \ar[ur]^{f} \ar[rr]^{h} & & Z }$$
where $f \in S_L$ and $g \in S_R$, and consider the associated diagram
$$ \xymatrix{ X \ar[r]^{f} \ar[d]^{h} & Y \ar[d]^{g} \\
Z \ar[r]^{\id} \ar@{-->}[ur] & Z. }$$
Since $h \perp g$, we can complete this diagram to a $3$-simplex of $\calC$ as indicated. This $3$-simplex exhibits $h$ as a retract of $f$, so that $h \in S_L$ as desired.
\end{proof}

\begin{remark}
It follows from Proposition \ref{swin} that a factorization system $(S_L, S_R)$ on an $\infty$-category $\calC$ is completely determined by {\em either} the left set $S_L$ or the right set $S_R$.
\end{remark}

\begin{corollary}\label{spen}
Let $(S_L, S_R)$ be a factorization system on an $\infty$-category $\calC$. Then
the collections of morphisms $S_L$ and $S_R$ contain all equivalences and are stable under composition.
\end{corollary}

\begin{proof}
Combine Propositions \ref{swin} and \ref{swimmm}.
\end{proof}

\begin{remark}
It follows from Corollary \ref{spen} that a factorization system $(S_L, S_R)$ on $\calC$ determines a pair of subcategories $\calC^{L}, \calC^{R} \subseteq \calC$, each containing all the objects of $\calC$: the morphisms of $\calC^{L}$ are the elements of $S_L$, and the morphisms of $\calC^{R}$ are the elements of $S_R$.
\end{remark}

\begin{example}
Let $p: \calC \rightarrow \calD$ be a coCartesian fibration of $\infty$-categories. Then there is an associated factorization system $(S_L, S_R)$ on $\calC$, where $S_L$ is the class of $p$-coCartesian morphisms of $\calC$, and $S_R$ is the class of morphisms $g$ of $\calC$ such that
$p(g)$ is an equivalence in $\calD$. If $\calD \simeq \Delta^0$, this recovers the factorization system of Example \ref{scumm}; if $p$ is an isomorphism, this recovers the opposite of the factorization system of Example \ref{scumm}.
\end{example}

\begin{example}\label{spink}
Let $\calX$ be an $\infty$-topos and let $n \geq -2$ be an integer. Then there exists a factorization system $(S_L, S_R)$ on $\calX$, where $S_L$ denotes the collection of $(n+1)$-connective morphisms
of $\calX$ and $S_R$ the collection of $n$-truncated morphisms of $\calC$. See \S 
\ref{homotopysheaves}.
\end{example}

Let $(S_L, S_R)$ be a factorization system on an $\infty$-category $\calC$, so that any morphism
$h: X \rightarrow Z$ factors as a composition
$$ \xymatrix{ & Y \ar[dr]^{g} & \\
X \ar[ur]^{f} \ar[rr]^{h} & & Z }$$
where $f \in S_L$ and $g \in S_R$. For many purposes, it is important to know that this factorization is {\em canonical}. More precisely, we have the following result:

\begin{proposition}\label{canfact}
Let $\calC$ be an $\infty$-category and let $S_L$ and $S_R$ be collections of morphisms
in $\calC$. Suppose that $S_L$ and $S_R$ are stable under equivalence in $\Fun( \Delta^1, \calC)$, and contain every equivalence in $\calC$.
The following conditions are equivalent:
\begin{itemize}
\item[$(1)$] The pair $(S_L, S_R)$ is a factorization system on $\calC$.
\item[$(2)$] The restriction map $p: \Fun'( \Delta^2, \calC) \rightarrow \Fun( \Delta^{ \{0,2\} }, \calC)$
is a trivial Kan fibration. Here $\Fun'( \Delta^2, \calC)$ denotes the full subcategory of
$\Fun( \Delta^2, \calC)$ spanned by those diagrams
$$ \xymatrix{ & Y \ar[dr]^{g} & \\
X \ar[ur]^{f} \ar[rr]^{h} & & Z }$$
such that $f \in S_L$ and $g \in S_R$. 
\end{itemize}
\end{proposition}

\begin{corollary}\label{funcsys}
Let $\calC$ be an $\infty$-category equipped with a factorization system $(S_L, S_R)$, and let $K$ be an arbitrary simplicial set. Then the $\infty$-category $\Fun(K, \calC)$ admits a factorization system
$(S_L^{K}, S_R^{K})$, where $S_{L}^{K}$ denotes the collection of all morphisms $f$ in 
$\Fun(K, \calC)$ such that $f(v) \in S_L$ for each vertex $v$ of $K$, and
$S_{R}^{K}$ is defined likewise.
\end{corollary}

The remainder of this section is devoted to the proof of Proposition \ref{canfact}. We begin with a few preliminary results.

\begin{lemma}\label{prefukt}
Let $\calC$ be an $\infty$-category and let $(S_L, S_R)$ be a factorization system on $\calC$.
Let $\calD$ be the full subcategory of $\Fun( \Delta^1, \calC)$ spanned by the elements of $S_R$.
Then:
\begin{itemize}
\item[$(1)$] The $\infty$-category $\calD$ is a localization of $\Fun( \Delta^1, \calC)$; in other words, the inclusion $\calD \subseteq \Fun( \Delta^1, \calC)$ admits a left adjoint.

\item[$(2)$] A morphism $\alpha: h \rightarrow g$ in $\Fun( \Delta^1, \calC)$, corresponding to a commutative diagram
$$ \xymatrix{ X \ar[d]^{h} \ar[r]^{f} & Y \ar[d]^{g} \\
Z' \ar[r]^{e} & Z }$$
exhibits $g$ as a $\calD$-localization of $h$ $($see Definition \ref{locaobj}$)$ if and only if $g \in S_R$, $f \in S_L$, and 
$e$ is an equivalence.
\end{itemize}
\end{lemma}

\begin{proof}
We will prove the ``if'' direction of assertion $(2)$.
It follows from the definition of a factorization system that for every object
$h \in \Fun( \Delta^1, \calC)$, there exists a morphism $\alpha: h \rightarrow g$ satisfying
the condition stated in $(2)$, which therefore exhibits $g$ as a $\calD$-localization of $h$.
Invoking Proposition \ref{testreflect}, we will deduce $(1)$. Because a $\calD$-localization of
$h$ is uniquely determined up to equivalence, we will also deduce the ``only if'' direction of assertion $(2)$.

Suppose given a commutative diagram
$$ \xymatrix{ X \ar[d]^{h} \ar[r]^{f} & Y \ar[d]^{g} \\
Z' \ar[r]^{e} & Z }$$
where $f \in S_L$, $g \in S_R$, and $e$ is an equivalence, and let $\overline{g}: \overline{Y} \rightarrow \overline{Z}$ be another element of $S_R$. We have a diagram of spaces
$$ \xymatrix{ \bHom_{\Fun(\Delta^1, \calC)}( g, \overline{g}) \ar[r]^{\psi} \ar[d] & \bHom_{ \Fun( \Delta^1, \calC)}( h, \overline{g} ) \ar[d] \\
\bHom_{\calC}( Z, \overline{Z}) \ar[r]^{\psi_0} & \bHom_{\calC}( Z', \overline{Z} )}$$
which commutes up to canonical homotopy. We wish to prove that $\psi$ is a homotopy equivalence.

Since $e$ is an equivalence in $\calC$, the map $\psi_0$ is a homotopy equivalence.
It will therefore suffice to show that $\psi$ induces a homotopy equivalence after passing to the homotopy fibers over any point of $\bHom_{\calC}(Z, \overline{Z}) \simeq \bHom_{\calC}(Z', \overline{Z})$. These homotopy fibers can be identified with the homotopy fibers of the vertical arrows in the diagram
$$ \xymatrix{ \bHom_{ \calC}( Y, \overline{Y} ) \ar[r] \ar[d] & \bHom_{\calC}( X, \overline{Y}) \ar[d] \\
\bHom_{ \calC}( Y, \overline{Z} ) \ar[r] & \bHom_{ \calC}( X, \overline{Z} ). }$$
It will therefore suffice to show that this diagram (which commutes up to specified homotopy) is a homotopy pullback. Unwinding the definition, this is equivalent to the assertion that $f$
is left orthogonal to $\overline{g}$, which is part of the definition of a factorization system.
\end{proof}

\begin{lemma}\label{hulfer}
Let $K$, $A$, and $B$ be simplicial sets. Then the diagram
$$ \xymatrix{ K \times B \ar[r] \ar[d] & K \times (A \star B) \ar[d] \\
B \ar[r] & ( K \times A) \star B }$$
is a homotopy pushout square of simplicial sets $($with respect to the Joyal model structure$)$.
\end{lemma}

\begin{proof}
We consider the larger diagram
$$ \xymatrix{ K \times B \ar[r] \ar[d] & K \times (A \diamond B) \ar[r] \ar[d] & K \times (A \star B) \ar[d] \\
B \ar[r] & (K \times A) \diamond B \ar[r] & (K \times A) \star B. }$$
The square on the left is a pushout square in which the horizontal maps are monomorphisms
of simplicial sets, and therefore a homotopy pushout square (since the Joyal model structure is left proper). The square on the right is a homotopy pushout square, since the horizontal arrows are both categorical equivalences (Proposition \ref{rub3}). It follows that the outer square is also a homotopy pushout, as desired.
\end{proof}

\begin{notation}
In the arguments which follow, we let $Q$ denote the simplicial subset of
$\Delta^3$ spanned by all simplices which do not contain $\Delta^{ \{1,2 \} }$. Note that
$Q$ is isomorphic to the product $\Delta^1 \times \Delta^1$ as a simplicial set.
\end{notation}

\begin{lemma}\label{sidewise}
Let $\calC$ be an $\infty$-category, and let $\sigma: Q \rightarrow \calC$ be a diagram, which
we depict as
$$ \xymatrix{ A \ar[r] \ar[d] & X \ar[d] \\
B \ar[r] & Y. }$$
Then there is a canonical categorical equivalence
$$\theta: \Fun( \Delta^3, \calC) \times_{ \Fun( Q, \calC) } \{ \sigma \} \rightarrow \bHom_{ \calC_{ A/ \, /Y }}( B, X)$$
In particular, $\Fun( \Delta^3, \calC) \times_{ \Fun(Q, \calC)} \{ \sigma \}$ is a Kan complex.
\end{lemma}

\begin{proof}
We will identify $\bHom_{ \calC_{A/ \, /Y}}(B,X)$ with the simplicial set $Z$ defined by the following universal property: for every simplicial set $K$, we have a pullback diagram of sets
$$ \xymatrix{ \Hom_{\sSet}(K, Z) \ar[r] \ar[d] & \Hom_{\sSet}( \Delta^0 \star (K \times \Delta^1) \star \Delta^0, \calC) \ar[d] \\
\Delta^0 \ar[r] & \Hom_{\sSet}( \Delta^0 \star (K \times \bd \Delta^1) \star \Delta^0, \calC). }$$
The map $\theta$ is then induced by the natural transformation
$$ K \times \Delta^3 \simeq K \times ( \Delta^0 \star \Delta^1 \star \Delta^0)
\rightarrow \Delta^0 \star (K \times \Delta^1) \star \Delta^0.$$

We wish to prove that $\theta$ is a categorical equivalence. Since $\calC$ is an
$\infty$-category, it will suffice to show that for every simplicial set $K$, the rightmost square of the diagram
$$ \xymatrix{ K \times ( \Delta^{ \{0\} } \coprod \Delta^{ \{3 \} }) \ar[r] \ar[d] & K \times C \ar[r] \ar[d] & K \times \Delta^3 \ar[d] \\
\Delta^{ \{0\} } \coprod \Delta^{ \{3\} } \ar[r] & \Delta^0 \star (K \times \bd \Delta^1) \star \Delta^0 \ar[r] & \Delta^0
\star (K \times \Delta^1) \star \Delta^0 }$$
is a homotopy pushout square (with respect to the Joyal model structure). For this we need only verify that the left and outer squares are homotopy pushout diagrams, which follows from repeated application of Lemma \ref{hulfer}.
\end{proof}

\begin{proof}[Proof of Proposition \ref{canfact}]
We first show that $(1) \Rightarrow (2)$. Assume that $(S_L, S_R)$ is a factorization system on $\calC$.
The map $p: \Fun'( \Delta^2, \calC) \rightarrow \Fun( \Delta^{ \{0,2\} }, \calC)$ is obviously a categorical fibration. It will therefore suffice to show that $p$ is a categorical equivalence.

Let $\calD$ be the full subcategory of $\Fun( \Delta^1 \times \Delta^1, \calC)$
spanned by those diagrams of the form
$$\xymatrix{ X \ar[d]^{h} \ar[r]^{f} & Y \ar[d]^{g} \\
Z' \ar[r]^{e} & Z }$$
where $f \in S_L$, $g \in S_R$, and $e$ is an equivalence in $\calC$. The map $p$ factors as a composition
$$ \Fun'(\Delta^2, \calC) \stackrel{p'}{\rightarrow} \calD \stackrel{p''}{\rightarrow} \Fun( \Delta^1, \calC)$$
where $p'$ carries a diagram
$$ \xymatrix{ & Y \ar[dr]^{g} & \\
X \ar[rr]^{h} \ar[ur]^{f} & & Z }$$
to the partially degenerate square
$$\xymatrix{ X \ar[d]^{h} \ar[dr]^{h} \ar[r]^{f} & Y \ar[d]^{g} \\
Z \ar[r]^{\id} & Z, }$$
and $p''$ is given by restriction to the left vertical edge of the diagram. To complete the proof, it will suffice to show that $p'$ and $p''$ are categorical equivalences.

We first show that $p'$ is a categorical equivalence.
The map $p'$ admits a left inverse $q$, given by composition with an inclusion
$\Delta^2 \subseteq \Delta^1 \times \Delta^1$. We note that $q$ is a pullback of the restriction map
$q': \Fun''( \Delta^2, \calC) \rightarrow \Fun( \Delta^{ \{0,2 \} }, \calC),$
where $\Fun''( \Delta^2, \calC)$ is the full subcategory spanned by diagrams of the form
$$ \xymatrix{ X \ar[dr] \ar[d] & \\
Z' \ar[r]^{e} & Z }$$
where $e$ is an equivalence. Since $q'$ is a trivial Kan fibration (Proposition \ref{lklk}), $q$ is a trivial Kan fibration, so that $p'$ is a categorical equivalence as desired.

We now complete the proof by showing that $p''$ is a trivial Kan fibration. Let
$\calE$ denote the full subcategory of $\Fun( \Delta^1, \calC) \times \Delta^1$ spanned by those pairs
$(g,i)$ where either $i=0$ or $g \in S_R$. The projection map $r: \calE \rightarrow \Delta^1$ is a 
Cartesian fibration associated to the inclusion $\Fun'( \Delta^1, \calC) \subseteq \Fun( \Delta^1, \calC)$, where $\Fun'(\Delta^1, \calC)$ is the full subcategory spanned by the elements of $S_R$. 
Using Lemma \ref{prefukt}, we conclude that $r$ is also a coCartesian fibration. Moreover,
we can identify
$$ \calD \subseteq \Fun( \Delta^1 \times \Delta^1, \calC) \simeq
\bHom_{ \Delta^1}( \Delta^1, \calE)$$
with the full subcategory spanned by the coCartesian sections of $r$. In terms of this identification,
$p''$ is given by evaluation at the initial vertex $\{0\} \subseteq \Delta^1$, and is therefore
a trivial Kan fibration as desired. This completes the proof that $(1) \Rightarrow (2)$.

Now suppose that $(2)$ is satisfied, and choose a section $s$ of the trivial Kan fibration $p$.
Let $s$ carry each morphism $f: X \rightarrow Z$ to a commutative diagram
$$ \xymatrix{ & Y \ar[dr]^{s_R(f)} & \\
X \ar[ur]^{ s_L(f)} \ar[rr]^{f} & & Z. }$$
If $s_{R}(f)$ is an equivalence, then $f$ is equivalent to $s_L(f)$ and therefore belongs to
$S_L$. Conversely, if $f$ belongs to $S_L$, then the diagram
$$ \xymatrix{ & Z \ar[dr]^{\id} & \\
X \ar[ur]^{f} \ar[rr]^{f} & & Z }$$
is a preimage of $f$ under $p$, and therefore equivalent to $s(f)$; this implies that
$s_L(f)$ is an equivalence. We have proved the following:
\begin{itemize}
\item[$(\ast)$] A morphism $f$ of $\calC$ belongs to $S_L$ if and only if $s_L(f)$ is an equivalence in $\calC$.
\end{itemize}

It follows immediately from $(\ast)$ that $S_L$ is stable under the formation of retracts; similarly, $S_R$ is stable under the formation of retracts. To complete the proof, it will suffice to show that $f \perp g$ whenever $f \in S_L$ and $g \in S_R$. Fix a commutative diagram $\sigma:$
$$ \xymatrix{ A \ar[d]^{f} \ar[r] & X \ar[d]^{g} \\
B \ar[r] & Y }$$
in $\calC$. In view of Lemma \ref{sidewise}, it will suffice to show that the Kan complex
$\Fun( \Delta^3, \calC) \times_{ \Fun( Q, \calC ) } \{ \sigma \}$ is contractible.

Let $\calD$ denote the full subcategory of $\Fun( \Delta^2 \times \Delta^1, \calC)$ spanned by those diagrams
$$ \xymatrix{ C \ar[r] \ar[d]^{u'} & Z \ar[d]^{v'} \\
C' \ar[r] \ar[d]^{u''} & Z' \ar[d]^{v''} \\
C'' \ar[r] & Z'' }$$
for which $u' \in S_L$, $v'' \in S_R$, and the maps $v'$ and $u''$ are equivalences.
Let us identify $\Delta^3$ with the full subcategory of $\Delta^2 \times \Delta^1$ spanned by
all those vertices except for $(2,0)$ and $(0,1)$. Applying Proposition \ref{lklk} twice,
we deduce that the restriction functor $\Fun( \Delta^2 \times \Delta^1, \calC) \rightarrow \Fun( \Delta^3, \calC)$ induces a trivial Kan fibration from $\calD$ to the full subcategory
$\calD' \subseteq \Fun(\Delta^3, \calC)$ spanned by those diagrams
$$ \xymatrix{ C \ar[d]^{u'} \ar[r] & Z' \ar[d]^{v''} \\
C' \ar[r] \ar[ur] & Z'' }$$
such that $u' \in S_L$ and $v'' \in S_R$. It will therefore suffice to show that the fiber
$\calD \times_{ \Fun(Q, \calC)} \{ \sigma \}$ is contractible.

By construction, the restriction functor $\calD \rightarrow \Fun(Q, \calC)$ is equivalent to the composition
$$ q: \calD \subseteq \Fun( \Delta^2 \times \Delta^1, \calC) \rightarrow
\Fun( \Delta^{ \{0,2\} } \times \Delta^1, \calC).$$
It will therefore suffice to show that $q^{-1} \{ \sigma \}$ is a contractible Kan complex.
Invoking assumption $(2)$ and $(\ast)$, we deduce that $q$ induces an equivalence from
$\calD$ to the full subcategory of $\Fun( \Delta^{ \{0,2\} } \times \Delta^1, \calC)$ spanned by
those diagrams
$$ \xymatrix{ C \ar[d]^{u} \ar[r] & Z \ar[d]^{v} \\
C'' \ar[r] & Z'' }$$
such that $u \in S_L$ and $v \in S_R$. The desired result now follows from our assumption that
$f \in S_L$ and $g \in S_R$.
\end{proof}

\subsection{Application: Automorphisms of $\Cat_{\infty}$}\label{cataut}

In \S \ref{working}, we saw that for every $\infty$-category $\calC$, the opposite simplicial set
$\calC^{op}$ is again an $\infty$-category. We will see below that the construction
$\calC \mapsto \calC^{op}$ determines an equivalence from the $\infty$-category
$\Cat_{\infty}$ to itself. In fact, this is essentially the {\em only} nontrivial self-equivalence of
$\Cat_{\infty}$. More precisely, we have the following result due to To\"{e}n (see \cite{toenchar}):

\begin{theorem}\label{cabbi}\index{gen}{$\infty$-category!opposite}\index{gen}{opposite!$\infty$-category}
Let $\calE$ denote the full subcategory of $\Fun( \Cat_{\infty}, \Cat_{\infty})$ spanned by the equivalences. Then $\calE$ is equivalent to the $($nerve of the$)$ discrete category
$\{ \id, r \}$, where $r: \Cat_{\infty} \rightarrow \Cat_{\infty}$ is a functor which associates to every
$\infty$-category its opposite. 
\end{theorem}

Our goal in this section is to give a proof of Theorem \ref{cabbi}. We first outline the basic strategy. Fix an equivalence $f$ from $\Cat_{\infty}$ to itself. The first step is to argue that $f$ is determined by its restriction to a reasonably small subcategory of $\Cat_{\infty}$. To prove this, we will introduce the notion of a subcategory $\calC^{0} \subseteq \calC$ which {\it strongly generates} $\calC$
(Definitions \ref{coughball} and \ref{ballcough}). We will then show that $\Cat_{\infty}$ is strongly generated by the subcategory consisting of nerves of partially ordered sets (in fact, it is generated
by an even smaller subcategory: see Theorem \ref{upsquare}). This will allow us to reduce to the problem of understanding the category of self-equivalences of the category of partially ordered sets, which is easy to tackle directly: see Proposition \ref{cape}. 

We begin by introducing some definitions.

\begin{definition}\label{coughball}\index{gen}{strongly generates}\index{gen}{generates!strongly}
Let $f: \calC \rightarrow \calD$ be a functor between $\infty$-categories.
We will say that $f$ is {\em strongly generates} the $\infty$-category $\calD$ if the identity transformation
$\id: f \rightarrow f$ exhibits the identity functor $\id_{\calD}$ as a left Kan extension of $f$
along $f$.
\end{definition}

\begin{remark}
In other words, a functor $f: \calC \rightarrow \calD$ strongly generates the $\infty$-category
$\calD$ if and only if, for every object $D \in \calD$, the evident diagram
$( \calC \times_{\calD} \calD_{/D} )^{\triangleright} \rightarrow \calD_{/D}^{\triangleright} \rightarrow \calD$
exhibits $D$ as a colimit of the diagram
$(\calC \times_{\calD} \calD_{/D}) \rightarrow \calC \stackrel{f}{\rightarrow} \calD.$
In particular, this implies that every object of $\calD$ can be obtained as the colimit of a diagram
which factors through $\calC$. Moreover, if $\calC$ is small and $\calD$ is locally small, then the diagram can be assumed small.
\end{remark}

\begin{remark}\label{copse}
Let $f: \calC \rightarrow \calD$ be a functor between $\infty$-categories, where $\calC$ is small, $\calD$ is locally small, and $\calD$ admits small colimits. In view of Theorem \ref{charpresheaf}, we may assume without loss of generality that $f$ factors as a composition
$$ \calC \stackrel{j}{\rightarrow} \calP(\calC) \stackrel{F}{\rightarrow} \calD,$$
where $j$ denotes the Yoneda embedding and $F$ preserves small colimits.
Corollary \ref{coughspaz} implies that $F$ has a right adjoint $G$, given by the composition
$$ \calD \stackrel{j'}{\rightarrow} \Fun( \calD^{op}, \SSet) \stackrel{\circ f}{\rightarrow} \calP(\calC),$$
where $j'$ denotes the Yoneda embedding for $\calD$; moreover, the transformation
$$ f = F \circ j \rightarrow (F \circ (G \circ F)) \circ j \simeq (F \circ G) \circ f$$
exhibits $(F \circ G)$ as a left Kan extension of $f$ along itself. It follows that $f$
strongly generates $\calD$ if and only if the counit map $F \circ G \rightarrow \id_{\calD}$ is an equivalence of functors. This is equivalent to the requirement that the functor $G$ is fully faithful.

In other words, the functor $f: \calC \rightarrow \calD$ strongly generates $\calD$ if and only if the
induced functor $F: \calP(\calC) \rightarrow \calD$ exhibits $\calD$ as a localization of
$\calP(\calC)$. In particular, the Yoneda embedding $\calC \rightarrow \calP(\calC)$
strongly generates $\calP(\calC)$, for any small $\infty$-category $\calC$.
\end{remark}

\begin{remark}\label{copo}
Let $f: \calC \rightarrow \calD$ be as in Remark \ref{copse}, and let $\calE$ be an $\infty$-category which admits small colimits. Let $\Fun^{0}(\calD, \calE)$ denote the full subcategory of
$\Fun(\calD, \calE)$ spanned by those functors which preserve small colimits. Then composition
with $f$ induces a fully faithful functor $\Fun^{0}(\calD, \calE) \rightarrow \Fun(\calC, \calE)$. 
This follows from Theorem \ref{charpresheaf}, Proposition \ref{unichar}, and Remark \ref{copse}.
\end{remark}

\begin{definition}\index{gen}{generates!strongly}\index{gen}{strongly generates}\label{ballcough}
Let $\calC$ be an $\infty$-category. We will say that a full subcategory $\calC^{0} \subseteq \calC$ 
{\it strongly generates} $\calC$ if the inclusion functor $\calC^{0} \rightarrow \calC$ strongly generates $\calC$, in the sense of Definition \ref{coughball}.
\end{definition}

\begin{remark}\label{sobre}
In other words, $\calC^0$ strongly generates $\calC$ if and only if the identity functor
$\id_{\calC}$ is a left Kan extension of $\id_{\calC} | \calC^{0}$. It follows from Proposition \ref{acekan} that if $\calC^{0} \subseteq \calC^{1} \subseteq \calC$ are full subcategories and
$\calC^{0}$ strongly generates $\calC$, then $\calC^{1}$ also strongly generates $\calC$.
\end{remark}

\begin{example}
The $\infty$-category $\SSet$ of spaces is strongly generated by its final object;
this follows immediately from Remark \ref{copse}.
\end{example}

The main ingredient in the proof of Theorem \ref{cabbi} is the following result, which is of interest in its own right:

\begin{theorem}\label{upsquare}
The $\infty$-category $\Cat_{\infty}$ is strongly generated by the full subcategory consisting of the objects $\{ \Delta^n \}_{n \geq 0}$. 
\end{theorem}

\begin{remark}
It follows from Theorem \ref{upsquare} that the theory of $\infty$-categories can be obtained as a suitable localization of the model category $\Fun( \cDelta^{op}, \sSet)$ of bisimplicial sets. It is possible to describe this localization precisely: this leads to Rezk's theory of {\em complete Segal spaces}, which is another model for the theory of higher categories where every $k$-morphism is assumed to be invertible for $k > 1$. For more details, we refer the reader to \cite{completesegal}.
\end{remark}

\begin{proof}[Proof of Theorem \ref{upsquare}]
We can identify the full subcategory in question with $\Nerve( \cDelta)$, the (nerve of) the category of simplices. The fully faithful embedding $f: \Nerve(\cDelta) \rightarrow \Cat_{\infty}$ can be extended (up to equivalence) to a colimit-preserving functor $F: \calP( \Nerve(\cDelta) ) \rightarrow \Cat_{\infty}$, which admits a right adjoint $G$ (this follows from Corollary \ref{adjointfunctor}, but an explicit description will be given below). In view of Remark \ref{copse}, it will suffice to show that the unit transformation $\id \rightarrow F \circ G$ is an equivalence.

We now reformulate the desired conclusion in the language of model categories. We can identify
$\Cat_{\infty}$ with the underlying $\infty$-category $\bfA^{\degree}$ of the simplicial model category
$\bfA = \mSet$ of marked simplicial sets, with the Cartesian model structure described in \S \ref{twuf}. The diagram $f$ is then obtained from a diagram $\overline{f}: \cDelta \rightarrow \bfA$, given by the cosimplicial object $[n] \mapsto (\Delta^n)^{\flat}.$
Since $\bfA$ is a simplicial model category (with respect to the simplicial structure
given by the formula $X \otimes K = X \times K^{\sharp}$), we can extend $\overline{f}$ to a colimit preserving functor
$$\overline{F}: \Fun( \cDelta^{op}, \sSet) \rightarrow \bfA.$$
Here $\Fun( \cDelta^{op}, \sSet)$ can be identified with the category of bisimplicial sets.
Since the cosimplicial object $\overline{f} \in \Fun( \cDelta, \bfA)$ is Reedy cofibrant (see \S \ref{coreed}), the functor $\overline{F}$ is a left Quillen functor if we endow $\Fun( \cDelta^{op}, \sSet)$ with the injective model structure (Example \ref{cabletome}). The functor $\overline{F}$ has a right adjoint
$\overline{G}$, given by the formula
$$ \overline{G}(X)^{m,n} = \Hom_{ \bfA}( (\Delta^m)^{\flat} \times (\Delta^n)^{\sharp}, X).$$
This right adjoint induces a functor from $\bfA^{\degree}$ to $\Fun(\cDelta^{op}, \sSet)^{\degree}$, which
(after passing to the simplicial nerve) is equivalent to the functor $G: \Cat_{\infty}
\rightarrow \calP( \Nerve(\cDelta))$ considered above. Consequently, it will suffice to show that
the counit map $L \overline{F} \circ R \overline{G} \rightarrow \id_{ \h{\bfA}}$ is an equivalence of functors, where $L \overline{F}$ and $R \overline{G}$ denote the left and right derived functors of
$\overline{F}$ and $\overline{G}$, respectively. Since every object of $\Fun( \cDelta^{op}, \sSet)$ is cofibrant, we can identify $\overline{F}$ with its left derived functor. We are therefore reduced to proving the following:
\begin{itemize}
\item[$(\ast)$] Let $\overline{X} = (X,M)$ be a fibrant object of the category $\bfA$ of marked simplicial sets. Then the counit map $\eta_{\overline{X}}: \overline{F} \overline{G} \overline{X} \rightarrow \overline{X}$ is a weak equivalence in $\bfA$.
\end{itemize}
Since $\overline{X}$ is fibrant, the simplicial set $X$ is an $\infty$-category and $M$ is the collection of all equivalences in $X$. Unwinding the definitions, we can identify
$\overline{F} \overline{G} X$ with the marked simplicial set $(Y, N)$ described as follows:
\begin{itemize}
\item[$(a)$] An $n$-simplex of $Y$ is a map of simplicial sets
$\Delta^n \times \Delta^n \rightarrow X$, which carries every morphism of
$\{i\} \times \Delta^n$ to an equivalence in $\calC$, for $0 \leq i \leq n$.
\item[$(b)$] An edge $\Delta^1 \rightarrow Y$ belongs to $N$ if and only if the corresponding
map $\Delta^1 \times \Delta^1 \rightarrow X$ factors through the projection onto the second
factor.
\end{itemize}

In terms of this identification, the map $\eta_{\overline{X}}: (Y,N) \rightarrow (X,M)$ is defined on $n$-simplices by composing with the diagonal map $\Delta^n \rightarrow \Delta^n \times \Delta^n$.

Let $N'$ denote the collection of all edges of $Y$ which correspond to maps from
$(\Delta^1 \times \Delta^1)^{\sharp}$ into $\overline{X}$. The map $\eta_{\overline{X}}$ factors as a composition
$$ (Y, N) \stackrel{i}{\rightarrow} (Y,N') \stackrel{ \eta'_{\overline{X}}}{\rightarrow} (X,M).$$
We claim that the map $i$ is a weak equivalence of marked simplicial sets. To prove this, it
will suffice to show that for every edge $\alpha$ which belongs to $N'$, there exists a 
$2$-simplex  $\sigma:$
$$ \xymatrix{ & y' \ar[dr]^{\alpha''} & \\
y \ar[ur]^{\alpha'} \ar[rr]^{\alpha} & & y'' }$$
in $Y$, where $\alpha'$ and $\alpha''$ belong to $N$. To see this, let us suppose that
$\alpha$ classifies a commutative diagram
$$ \xymatrix{ A \ar[d]^{p} \ar[r]^{q} \ar[dr]^{r} & A' \ar[d]^{p'} \\
B \ar[r]^{q'} & B' }$$
in the $\infty$-category $X$. We wish to construct an appropriate $2$-simplex
$\sigma$ in $Y$, corresponding to a map $\widetilde{\sigma}: \Delta^2 \times \Delta^2 \rightarrow X^0$, where $X^0$ denotes the largest Kan complex contained in $X$. Let $T$ denote the full subcategory
of $\Delta^2 \times \Delta^2$ spanned by all vertices except for $(0,2)$, and let
$\widetilde{\sigma}_0: T \rightarrow X^0$ be the map described by the diagram
$$ \xymatrix{ A \ar[r]^{\id} \ar[d]^{q} & A \ar[d]^{p} \ar[r]^{q} & A' \ar[d]^{\id} \\
A' \ar[r]^{\id} & A' \ar[r]^{\id} \ar[d]^{p'} & A' \ar[d]^{p'} \\
& B' \ar[r]^{\id} & B'. }$$
To prove that $\widetilde{\sigma}_0$ can be extended to a map $\widetilde{\sigma}$ with the desired properties, it suffices to solve an extension problem of the form
$$ \xymatrix{ T \coprod_{ \Delta^{ \{0,2\} } } \Delta^2 \ar[r] \ar[d] & X^0 \\
\Delta^2 \times \Delta^2. \ar@{-->}[ur] & }$$
This is possible because $X^0$ is a Kan complex and the left vertical map is a weak homotopy equivalence. This completes the proof that $i$ is a weak equivalence. By the two-out-of-three property, it will now suffice to show that $\eta'_{\overline{X}}: (Y,N') \rightarrow (X,M)$ is an equivalence of marked simplicial sets.

We now define maps $R_{\leq}, R_{\geq}: \Delta^1 \times Y \rightarrow Y$ as follows.
Consider a map $g: \Delta^n \rightarrow \Delta^1 \times Y$, corresponding to a partition
$[n] = [n]_{-} \cup [n]_{+}$ and a map $\widetilde{g}: \Delta^n \times \Delta^n \rightarrow X$.
We then define $R_{\leq} \circ g$ to be the $n$-simplex of $Y$ corresponding to the map
$\widetilde{g} \circ \tau: \Delta^n \times \Delta^n \rightarrow X$, where
$\tau: \Delta^n \times \Delta^n \rightarrow \Delta^n \times \Delta^n$ is defined on vertices by the formula
$$ \tau(i,j) = \begin{cases} (i,j) & \text{if } i \leq j \\
(i,j) & \text{if } j \in [n]_{-} \\
(i,i) & \text{otherwise.} \end{cases}$$
Similarly, we let $R_{\geq} \circ g$ correspond to the map $\widetilde{g} \circ \tau'$, where
$\tau'$ is given on vertices by the formula
$$ \tau_{i,j} = \begin{cases} (i,j) & \text{if } i \geq j \\
(i,j) & \text{if } j \in [n]_{+} \\
(i,I) & \text{otherwise.} \end{cases}$$
The map $R_{\leq}$ defines a homotopy from $\id_{Y}$ to an idempotent map
$r_{\leq}: Y \rightarrow Y$. Similarly, $R_{\geq}$ defines a homotopy from an idempotent map
$r_{\geq}: Y \rightarrow Y$ to the identity map $\id_{Y}$. Let $Y_{\leq}, Y_{\geq} \subseteq Y$
denote the images of the maps $r_{\leq}$ and $r_{\geq}$, respectively. Let
$N'_{\leq}$ denote the collection of all edges of $Y$ which belong to $N'$, and define
$N'_{\geq}$ similarly. The map $R_{\leq}$ determines a map
$(Y,N') \times (\Delta^1)^{\sharp} \rightarrow (Y,N')$, which exhibits
$(Y_{\leq}, N'_{\leq})$ as a deformation retract of $(Y, N')$ in the category of marked simplicial sets.
Similarly, the map $R_{\geq}$ exhibits $(Y_{\leq} \cap Y_{\geq}, N'_{\leq} \cap N'_{\geq})$ as a
deformation retract of $(Y_{\leq}, N'_{\leq})$. It will therefore suffice to show that the composite map
$$ (Y_{\leq} \cap Y_{\geq}, N'_{\leq} \cap N'_{\geq}) \subseteq (Y, N') \rightarrow (X,M)$$
is a weak equivalence of marked simplicial sets. We now complete the proof by observing that this composite map is an isomorphism.
\end{proof}



It follows from Remark \ref{sobre} and Theorem \ref{upsquare} that $\Cat_{\infty}$ is strongly generated by the full subcategory spanned by those $\infty$-categories of the form $\Nerve P$, where
$P$ is a partially ordered set. Our next step is to describe this subcategory in more intrinsic terms.

\begin{proposition}
Let $\calC$ be an $\infty$-category. The following conditions are equivalent:
\begin{itemize}
\item[$(1)$] The $\infty$-category $\calC$ is equivalent to the nerve of a partially ordered set $P$.
\item[$(2)$] For every $\infty$-category $\calD$ and every pair of functors
$F, F': \calD \rightarrow \calC$ such that $F(x) \simeq F'(x)$ for each object $x \in \calD$, 
the functors $F$ and $F'$ are equivalent as objects of $\Fun(\calD, \calC)$.
\item[$(3)$] For every $\infty$-category $\calD$, the map of sets
$$ \pi_0 \bHom_{\Cat_{\infty}}( \calD, \calC)
\rightarrow \Hom_{ \Set}( \pi_0 \bHom_{ \Cat_{\infty}}( \Delta^0, \calD),
\pi_0 \bHom_{ \Cat_{\infty} }( \Delta^0, \calC) )$$
is injective.
\end{itemize}
\end{proposition}

\begin{proof}
The implication $(1) \Rightarrow (2)$ is obvious, and $(3)$ is just a restatement of $(2)$.
Assume $(2)$; we will show that $(1)$ is satisfied. Let $P$ denote the collection of equivalence
classes of objects of $\calC$, where $x \leq y$ if the space $\bHom_{\calC}(x,y)$ is nonempty. 
There is a canonical functor $\calC \rightarrow \Nerve P$. To prove that this functor is an equivalence, it will suffice to show the following:
\begin{itemize}
\item[$(\ast)$] For every pair of objects $x,y \in \calC$, the space $\bHom_{\calC}(x,y)$ is either empty or contractible.
\end{itemize}
To prove $(\ast)$, we may assume without loss of generality that $\calC$ is the nerve of a fibrant
simplicial category $\overline{\calC}$. Let $x$ and $y$ be objects of $\overline{\calC}$ such that
the Kan complex $K = \bHom_{ \overline{\calC} }(x,y)$ is nonempty.
We define a new (fibrant) simplicial category $\overline{\calD}$ so that $\overline{\calD}$ consists of a pair of objects $\{x', y'\}$, with
$$ \bHom_{\overline{\calD}}(x',x') \simeq \bHom_{\overline{\calD}}(y',y') \simeq \Delta^0$$
$$ \bHom_{\overline{\calD}}(x',y') \simeq K \quad \bHom_{\overline{\calD}}(y',x') \simeq \emptyset.$$
We let $\overline{F}, \overline{F}': \overline{\calD} \rightarrow \overline{\calC}$ be simplicial functors
such that $\overline{F}(x') = \overline{F}'(x') = x$, $\overline{F}(y') = \overline{F}'(y')=y$, where
$\overline{F}$ induces the identity map from $\bHom_{\overline{\calD}}(x',y') = K = \bHom_{\overline{\calC}}(x,y)$ to itself, while $\overline{F}'$ induces a constant map from $K$ to itself.
Then $\overline{F}$ and $\overline{F}'$ induce functors $F$ and $F'$ from $\sNerve( \overline{\calD} )$ to $\calC$. It follows from assumption $(2)$ that the functors $F$ and $F'$ are equivalent, which implies that the identity map from $K$ to itself is homotopic to a constant map; this proves that $K$ is contractible. 
\end{proof}

\begin{corollary}\label{bigegg}
Let $\sigma: \Cat_{\infty} \rightarrow \Cat_{\infty}$ be an equivalence of $\infty$-categories, and let
$\calC$ be an $\infty$-category $($which we regard as an object of $\Cat_{\infty}${}$)$. Then $\calC$ is equivalent to the nerve of a partially ordered set if and only if $\sigma(\calC)$ is equivalent to the nerve of a partially ordered set.
\end{corollary}

\begin{lemma}\label{calcul}
Let $\sigma, \sigma' \in \{ \id_{\cDelta}, r \} \subseteq \Fun( \cDelta, \cDelta)$, where
$r$ denotes the reversal functor from $\cDelta$ to itself. Then
$$ \Hom_{ \Fun(\cDelta, \cDelta)}(\sigma, \sigma') = \begin{cases} \emptyset & \text{if } \sigma= \sigma' \\
\{ \id \} & \text{if } \sigma = \sigma'. \end{cases}$$
\end{lemma}

\begin{proof}
Note that $\sigma$ and $\sigma'$ are both the identity at the level of objects. 
Let $\alpha: \sigma \rightarrow \sigma'$ be a natural transformation. Then, for each $n \geq 0$,
$\alpha_{[n]}$ is a map from $[n]$ to itself. We claim that $\alpha_{[n]}$ is given by the formula
$$ \alpha_{[n]}(i) = \begin{cases} i & \text{if } \sigma = \sigma' \\
n-i & \text{if } \sigma \neq \sigma'.\end{cases}$$
To prove this, we observe that a choice of $i \in [n]$ determines a map $[0] \rightarrow [n]$, which allows us to reduce to the case $n=0$ (where the result is obvious) by functoriality.

It follows from the above argument that the natural transformation $\alpha$ is uniquely determined, if it exists. Moreover, $\alpha$ is a well-defined natural transformation if and only if each $\alpha_{[n]}$ is an order-preserving map from $[n]$ to itself; this is true if and only if $\sigma = \sigma'$.
\end{proof}

\begin{proposition}\label{cape}
Let $\calP$ denote the category of partially ordered sets, and let
$\sigma: \calP \rightarrow \calP$ be an equivalence of categories. Then
$\sigma$ is isomorphic either to the identity functor $\id_{\calP}$ or the functor
$r$ which carries every partially ordered set $X$ to the same set with the opposite ordering.
\end{proposition}

\begin{proof}
Since $\sigma$ is an equivalence of categories, it carries the final object $[0] \in \calP$ to itself (up to canonical isomorphism). It follows that for every
partially ordered set $X$, we have a canonical bijection of sets
$$ \eta_{X}: X \simeq \Hom_{ \calP}( [0], X) \simeq
\Hom_{\calP}( \sigma([0]), \sigma(X) ) \simeq \Hom_{\calP}( [0], \sigma(X)) \simeq \sigma(X).$$

We next claim that $\sigma([1])$ is isomorphic to $[1]$ as a partially ordered set.
Since $\eta_{[1]}$ is bijective, the partially ordered set $\sigma([1])$ has precisely two elements.
Thus $\sigma([1])$ is isomorphic either to $[1]$ or to a partially ordered set $\{x,y\}$ with two elements, neither larger than the other. In the second case, the set $\Hom_{\calP}( \sigma([1]), \sigma([1]))$
has four elements. This is impossible, since $\sigma$ is an equivalence of categories and
$\Hom_{\calP}( [1], [1])$ has only three elements. Let $\alpha: \sigma([1]) \rightarrow [1]$
be an isomorphism (automatically unique, since the ordered set $[1]$ has no automorphisms in $\calP$). 

The map $\alpha \circ \eta_{[1]}$ is a bijection from the set $[1]$ to itself. We will assume that
this map is the identity, and prove that $\sigma$ is isomorphic to the identity functor $\id_{\calP}$.
The same argument, applied to $\sigma \circ r$, will show that if $\alpha \circ \eta_{[1]}$ is not the identity, then $\sigma$ is isomorphic to $r$.

To prove that $\sigma$ is equivalent to the identity functor, it will suffice to show that for every
partially ordered set $X$, the map $\eta_{X}$ is an isomorphism of partially ordered sets.
In other words, we must show that both $\eta_{X}$ and $\eta_{X}^{-1}$ are maps of partially ordered sets. We will prove that $\eta_{X}$ is a map of partially ordered sets; the same argument, applied to 
an inverse to the equivalence $\sigma$, will show that $\eta^{-1}_{X}$ is a map of partially ordered sets.
Let $x,y \in X$ satisfy $x \leq y$; we wish to prove that $\eta_{X}(x) \leq \eta_{X}(y)$ in $\sigma(X)$.
The pair $(x,y)$ defines a map of partially ordered sets $[1] \rightarrow X$. By functoriality, we
may replace $X$ by $[1]$, and thereby reduce to the problem of proving that $\eta_{[1]}$ is a map of partially ordered sets. This follows from our assumption that $\alpha \circ \eta_{[1]}$ is the identity map.
\end{proof}

\begin{proof}[Proof of Theorem \ref{cabbi}]
Let $\calC$ be the full subcategory of $\Cat_{\infty}$ spanned by those $\infty$-categories
which are equivalent to the nerves of partially ordered sets, and let $\calC^{0}$ denote the full subcategory of $\calC$ spanned by the objects $\{ \Delta^n \}_{n \geq 0}$. 
Corollary \ref{bigegg} implies that every object $\sigma \in \calE$ restricts to an equivalence from $\calC$ to itself. According to Proposition \ref{cape}, $\sigma| \calC$ is equivalent either to the identity functor, or to the restriction $r|\calC$. In either case, we conclude that $\sigma$ also induces an equivalence from $\calC^{0}$ to itself.

Using Theorem \ref{upsquare} and Remark \ref{copo}, we deduce that the restriction functor
$\calE \rightarrow \Fun(\calC^{0}, \calC^{0})$ is fully faithful. In particular, any object $\sigma \in \calE$ is determined by the restriction $\sigma | \calC$, so that $\sigma$ is equivalent to either $\id$ or
$r$ by virtue of Proposition \ref{cape}. Since $\calC^{0}$ is equivalent to the nerve of the category
$\cDelta$, Lemma \ref{calcul} implies the existence of a fully faithful embedding from
$\calE$ to the nerve of the discrete category $\{ \id, r \}$. To complete the proof, it will suffice to show that this functor is essentially surjective. In other words, we must show that there exists a functor
$R: \Cat_{\infty} \rightarrow \Cat_{\infty}$ whose restriction to $\calC$ is equivalent to $r$. 

To carry out the details, it is convenient to replace $\Cat_{\infty}$ by an equivalent
$\infty$-category with a slightly more elaborate definition. Recall that $\Cat_{\infty}$ is defined to be the simplicial nerve of a simplicial category $\Cat_{\infty}^{\Delta}$, whose objects are $\infty$-categories, where $\bHom_{\Cat_{\infty}^{\Delta}}(X,Y)$ is the largest Kan complex contained in
$\Fun(X,Y)$. We would like to define $R$ to be induced by the functor $X \mapsto X^{op}$, but this
is not a simplicial functor from $\Cat_{\infty}^{\Delta}$ to itself; instead we have a canonical isomorphism $\bHom_{\Cat_{\infty}^{\Delta}}( X^{op}, Y^{op} )
\simeq \bHom_{\Cat_{\infty}^{\Delta}}(X,Y)^{op}.$
However, if we let $\Cat_{\infty}^{\top}$ denote the topological category obtained by geometrically realizing the morphism spaces in 
$\Cat_{\infty}^{\Delta}$, then $i$ induces an autoequivalence of $\Cat_{\infty}^{\top}$ as a topological category (via the natural homeomorphisms $| K | \simeq |K^{op}|$, which is defined for every simplicial set $K$). We now define $\Cat'_{\infty}$ to be the topological nerve of $\Cat_{\infty}^{\top}$ (see 
Definition \ref{topnerve}). Then $\Cat'_{\infty}$ is an $\infty$-category equipped with a canonical equivalence $\Cat_{\infty} \rightarrow \Cat'_{\infty}$, and the involution $i$ induces
an involution $I$ on $\Cat'_{\infty}$, which carries each object $\calD \in \Cat'_{\infty}$ to 
the opposite $\infty$-category $\calD^{op}$. We now define $R$ to be the composition
$$ \Cat_{\infty} \rightarrow \Cat'_{\infty} \stackrel{I}{\rightarrow} \Cat'_{\infty} \rightarrow \Cat_{\infty},$$
where the last map is a homotopy inverse to the equivalence $\Cat_{\infty} \rightarrow \Cat'_{\infty}$.
It is easy to see that $R$ has the desired properties (moreover, we note that for every
object $\calD \in \Cat_{\infty}$, the image $R \calD$ is canonically equivalent with the opposite $\infty$-category $\calD^{op}$).
\end{proof}

 \section{$\infty$-Categories of Inductive Limits}\label{c5s3}

\setcounter{theorem}{0}

Let $\calC$ be a category. An {\it $\Ind$-object} of $\calC$ is a diagram\index{gen}{Ind-object}
$f: \calI \rightarrow \calC$ where $\calI$ is a small filtered category. We will informally denote the $\Ind$-object $f$ by $$ [\colim X_i] $$
where $X_i = f(i)$. The collection of all $\Ind$-objects of $\calC$ forms a category, where
the morphisms are given by the formula
$$ \Hom_{\Ind(\calC)}( [\colim X_i], [\colim Y_j]) =
\varprojlim \varinjlim \Hom_{\calC}(X_i, Y_j).$$
We note that $\calC$ may be identified with a full subcategory of $\Ind(\calC)$, corresponding
to diagrams indexed by the one-point category $\calI = \ast$. 
The idea is that $\Ind(\calC)$ is obtained from $\calC$ by formally adjoining colimits of filtered diagrams. More precisely, $\Ind(\calC)$ may be described by the following universal property:
for any category $\calD$ which admits filtered colimits, and any functor $F: \calC \rightarrow \calD$, there exists a functor $\widetilde{F}: \Ind(\calC) \rightarrow \calD$, whose restriction to $\calC$ is isomorphic to $F$, and which commutes with filtered colimits. Moreover, $\widetilde{F}$ is determined up to (unique) isomorphism.\index{not}{IndC@$\Ind(\calC)$}

\begin{example}
Let $\calC$ denote the category of finitely presented groups. Then $\Ind(\calC)$ is equivalent to the category of groups. (More generally, one could replace ``group'' by any type of mathematical structure described by algebraic operations which are required to satisfy equational axioms.)
\end{example}

Our objective in this section is to generalize the definition of $\Ind(\calC)$ to the case where $\calC$ is an $\infty$-category. If we were to work in the setting of simplicial (or topological) categories, we could apply the definition given above directly. However, this leads to a number of problems:

\begin{itemize}
\item[$(1)$] The construction of $\Ind$-categories does not preserve equivalences between simplicial categories.

\item[$(2)$] The obvious generalization of the right hand side in equation above is given by
$$ \varprojlim \colim \bHom_{\calC}(X_i,Y_j).$$
While the relevant limits and colimits certainly exist in the category of simplicial sets, they
are not necessarily the correct objects: really one should replace the limit by a homotopy limit.

\item[$(3)$] In the higher-categorical setting, we should really allow the indexing diagram $\calI$
to be a higher category as well. While this does not result in any additional generality (Corollary \ref{rot}), the restriction to the diagrams indexed by ordinary categories is a technical inconvenience.
\end{itemize}

Although these difficulties are not insurmountable, it is far more convenient to proceed differently, using the theory of $\infty$-categories. In \S \ref{c5s1}, we showed that if $\calC$ is a $\infty$-category, then $\calP(\calC)$ can be interpreted as an $\infty$-category which is freely generated by $\calC$ under colimits. We might therefore hope to find $\Ind(\calC)$ {\it inside} of $\calP(\calC)$, as a full subcategory. The problem, then, is to characterize this subcategory, and to prove that it has the appropriate universal mapping property.

We will begin in \S \ref{smallfilt}, by introducing the definition of a {\em filtered $\infty$-category}.
Let $\calC$ be a small $\infty$-category. In \S \ref{indlim}, we will define $\Ind(\calC)$ to be the smallest full subcategory of $\calP(\calC)$ which contains all representable presheaves on $\calC$ and is stable under filtered colimits. There is also a more direct characterization of which presheaves $F: \calC \rightarrow \SSet^{op}$ belong to $\Ind(\calC)$: they are precisely the {\em right exact} functors, which we will study in \S \ref{rexex}.

In \S \ref{indlim}, we will define the $\Ind$-categories $\Ind(\calC)$ and study their properties. In particular, we will show that morphism spaces in $\Ind(\calC)$ {\em are} computed by the naive formula
$$ \Hom_{\Ind(\calC)}( [\colim X_i], [\colim Y_j]) =
\varprojlim \varinjlim \Hom_{\calC}(X_i, Y_j).$$
Unwinding the definitions, this amounts to two conditions:
\begin{itemize}
\item[$(1)$] The (Yoneda) embedding of $j: \calC \rightarrow \Ind(\calC)$ is fully faithful (Proposition \ref{fulfaith}).
\item[$(2)$] For each object $C \in \calC$, the corepresentable functor
$$ \Hom_{\Ind(\calC)}( j(C), \bigdot) $$ commutes with filtered colimits.
\end{itemize}
It is useful to translate condition $(2)$ into a definition: an object $D$ of an $\infty$-category
$\calD$ is said to be {\it compact} if the functor $\calD \rightarrow \SSet$ corepresented by $D$ commutes with filtered colimits. We will study this compactness condition in \S \ref{compobj}.

One of our main results asserts that the $\infty$-category $\Ind(\calC)$ is obtained from $\calC$ by freely adjoining colimits of filtered diagrams (Proposition \ref{intprop}). In \S \ref{agileco}, we will describe a similar construction in the case where the class of filtered diagrams has been replaced by {\em any} class of diagrams. We will revisit this idea in \S \ref{stable11}, where we will study the $\infty$-category obtained from $\calC$ by freely adjoining colimits of {\em sifted} diagrams.

\subsection{Filtered $\infty$-Categories}\label{smallfilt}

Recall that a partially ordered\index{gen}{filtered!partially ordered set}
set $A$ is {\it filtered} if every finite subset of $A$ has an
upper bound in $A$. Diagrams indexed by directed partially ordered sets are extremely common in mathematics. For example, if $A$ is the set $$\Z_{\geq 0} = \{0, 1, \ldots \}$$ of natural
numbers, then a diagram indexed by $A$ is a sequence
$$ X_0 \rightarrow X_1 \rightarrow \ldots .$$
The formation of direct limits for such sequences is one of the most basic constructions in mathematics.

In classical category theory, it is convenient to consider not
only diagrams indexed by filtered partially ordered sets, but also
more general diagrams indexed by filtered categories. A category
$\calC$ is said to be {\it filtered} if it satisfies following
conditions:\index{gen}{filtered!category}

\begin{itemize}
\item[$(1)$] For every finite collection $\{ X_i \}$ of objects of $\calC$, there exists an object
$X \in \calC$ equipped with morphisms $\phi_i: X_i \rightarrow X$.

\item[$(2)$]  Given any two morphisms $f,g: X \rightarrow Y$ in $\calC$,
there exists a morphism $h: Y \rightarrow Z$ such that $h \circ f
= h \circ g$.
\end{itemize}

Condition $(1)$ is analogous to the requirement that any
finite part of $\calC$ admits an ``upper bound'', while condition $(2)$ guarantees that the upper bound is unique in some asymptotic sense.

If we wish to extend the above definition to the $\infty$-categorical setting, it is natural to strengthen the second condition.

\begin{definition}\label{topfilt}\index{gen}{filtered!topological category}
Let $\calC$ be a topological category. We will say that $\calC$ is {\it filtered} if it satisfies
the following conditions:
\begin{itemize}
\item[$(1')$] For every finite set $\{X_i\}$ of objects of $\calC$, there exists an object $X \in \calC$
and morphisms $\phi_i: X_i \rightarrow X$.
\item[$(2')$] For every pair $X, Y \in \calC$ of objects of $\calC$, every nonnegative integer $n \geq 0$, and every continuous map $S^n \rightarrow \bHom_{\calC}(X,Y)$, there exists a morphism $Y \rightarrow Z$ such that the induced map $S^n \rightarrow \bHom_{\calC}(X,Z)$ is nullhomotopic.
\end{itemize}
\end{definition}

\begin{remark}
It is easy to see that an ordinary category $\calC$ is filtered in the usual sense if and only if it is filtered when regarded as a topological category with discrete mapping spaces. Conversely, 
if $\calC$ is a filtered topological category, then its homotopy category $\h{\calC}$ is filtered (when viewed as an ordinary category). 
\end{remark}

\begin{remark}
Condition $(2')$ of Definition \ref{topfilt} is a reasonable analogue of condition $(2)$ in the definition of a filtered category. In the special case $n=0$, condition $(2')$ asserts that any pair of morphisms
$f,g: X \rightarrow Y$ become {\em homotopic} after composition with some map $Y \rightarrow Z$.
\end{remark}

\begin{remark}
Topological spheres $S^n$ need not play any distinguished role in the definition of a filtered topological category. Condition $(2')$ is equivalent to the following apparently stronger condition:
\begin{itemize}
\item[$(2'')$] For every pair $X, Y \in \calC$ of objects of $\calC$, every finite cell complex $K$, and every continuous map $K \rightarrow \bHom_{\calC}(X,Y)$, there exists a morphism $Y \rightarrow Z$ such that the induced map $K \rightarrow \bHom_{\calC}(X,Z)$ is nullhomotopic.
\end{itemize}
\end{remark}

\begin{remark}
The condition that a topological category $\calC$ be filtered depends only on the homotopy category $\h{\calC}$, viewed as a $\calH$-enriched category. Consequently if $F: \calC \rightarrow \calC'$ is an equivalence of topological categories, then $\calC$ is filtered if and only if $\calC'$ is filtered.
\end{remark}

\begin{remark}
Definition \ref{topfilt} has an obvious analogue for (fibrant) simplicial categories: one simply replaces the topological $n$-sphere $S^n$ by the simplicial $n$-sphere $\bd \Delta^n$. It is easy to see that a topological category $\calC$ is filtered if and only if the simplicial category $\Sing \calC$ is filtered. Similarly, a (fibrant) simplicial category $\calD$ is filtered if and only if the topological category $|\calD|$ is filtered.
\end{remark}

We now wish to study the analogue of Definition \ref{topfilt} in the setting of $\infty$-categories. It will be convenient to introduce a slightly more general notion:

\begin{definition}\label{filtquas}\index{gen}{filtered!$\infty$-category}\index{gen}{$\kappa$-filtered}
Let $\kappa$ be a regular cardinal, and let $\calC$ be a $\infty$-category. We will say that
$\calC$ is {\it $\kappa$-filtered} if, for every $\kappa$-small simplicial set $K$ and every 
map $f: K \rightarrow \calC$, there exists a map $\overline{f}: K^{\triangleright} \rightarrow \calC$ extending $f$. (In other words, 
$\calC$ is $\kappa$-filtered if it has the extension property with respect to the inclusion
$K \subseteq K^{\triangleright}$, for every $\kappa$-small simplicial set $K$.)

We will say that $\calC$ is {\it filtered} if it is $\omega$-filtered.
\end{definition}

\begin{example}
Let $\calC$ be the nerve of a partially ordered set $A$. Then $\calC$ is $\kappa$-filtered if and only if every $\kappa$-small subset of $A$ has an upper bound in $A$.
\end{example}

\begin{remark}\label{falg}
One may rephrase Definition \ref{filtquas} as follows: an $\infty$-category $\calC$ is $\kappa$-filtered if and only if,
for every diagram $p: K \rightarrow \calC$, where $K$ is $\kappa$-small, the slice $\infty$-category
$\calC_{p/}$ is nonempty. 
Let $q: \calC \rightarrow \calC'$ be a categorical equivalence of $\infty$-categories. Proposition \ref{gorban3} asserts that the induced map $\calC_{p/} \rightarrow \calC'_{q \circ p/}$ is a categorical equivalence. Consequently $\calC_{p/}$ is nonempty if and only if $\calC'_{q \circ p/}$ is nonempty. It follows that $\calC$ is $\kappa$-filtered if and only if $\calC'$ is $\kappa$-filtered.
\end{remark}

\begin{remark}\label{tweeny}
An $\infty$-category $\calC$ is $\kappa$-filtered if and only if, for every
$\kappa$-small partially ordered set $A$, $\calC$ has the right lifting property with respect
to the inclusion $\Nerve(A) \subseteq \Nerve(A)^{\triangleright}  \simeq \Nerve(A \cup \{ \infty\} )$. The ``only if'' direction is obvious. For the converse, we observe that for every $\kappa$-small diagram $p: K \rightarrow \calC$, the $\infty$-category $\calC_{p/}$ is equivalent to
$\calC_{q/}$, where $q$ denotes the composition
$K'' \stackrel{p'}{\rightarrow} K \stackrel{p}{\rightarrow} \calC$. Here $K''$ is the second barycentric subdivision of $K$ and $p'$ is the map described in Variant \ref{baryvar}. We now observe that $K''$ is equivalent to the nerve of a $\kappa$-small partially ordered set. 
\end{remark}

\begin{remark}\index{gen}{filtered!simplicial set}
We will say that an arbitrary simplicial set $S$ is {\it $\kappa$-filtered} if there exists a categorical equivalence $j: S \rightarrow \calC$, where $\calC$ is a $\kappa$-filtered $\infty$-category. In view of the Remark \ref{falg}, this condition is independent of the choice of $j$.
\end{remark}

Our next major goal is to prove Proposition \ref{stook}, which asserts that an $\infty$-category
$\calC$ is filtered if and only if the associated topological category $| \sCoNerve[\calC] |$ is filtered.
First, we need a lemma.

\begin{lemma}\label{goony}
Let $\calC$ be an $\infty$-category. Then $\calC$ is filtered if and only if it has the right extension
property with respect to every inclusion $\bd \Delta^n \subseteq \Lambda^{n+1}_{n+1}$, $n \geq 0$.
\end{lemma}

\begin{proof}
The ``only if'' direction is clear: we simply take $K = \bd \Delta^n$ in Definition \ref{filtquas}.
For the converse, let us suppose that the assumption of Definition \ref{filtquas} is satisfied whenever $K$ is the boundary of a simplex; we must then show that it remains satisfied for {\em any} $K$ which has only finitely many nondegenerate simplices.

We work by induction on the dimension of $K$, and then by the number of nondegenerate simplices of $K$. If $K$ is empty, there is nothing to prove (since it is the boundary of a $0$-simplex). Otherwise, we may write $K = K' \coprod_{ \bd \Delta^n } \Delta^n$, where $n$ is the dimension of $K$. 

Choose a map $p: K \rightarrow \calC$; we wish to show that $p$ may be extended to a map
$\widetilde{p}: K \star \{y\} \rightarrow \calC$. We first consider the restriction $p|K'$; by the inductive hypothesis it admits an extension $q: K' \star \{x\} \rightarrow \calC$. The
restriction $q | \bd \Delta^n \star \{x\}$ together with $p|\Delta^n$ assemble to give a map
$$ r: \bd \Delta^{n+1} \simeq ( \bd \Delta^n \star \{x\} ) \coprod_{ \bd \Delta^n } \Delta^n \rightarrow
\calC.$$
By assumption, the map $r$ admits an extension
$$\widetilde{r}: \bd \Delta^{n+1} \star \{y\} \rightarrow \calC.$$

Let $$s: ( K' \star \{x\} ) \coprod_{ \bd \Delta^{n+1}  } ( \bd \Delta^{n+1} \star \{y\})$$ denote the result of amalgamating $r$ with $\widetilde{p}$. We note that the inclusion
$$ (K' \star \{x\} ) \coprod_{ \bd \Delta^n \star \{x\} } ( \bd \Delta^{n+1} \star \{y\} )
\subseteq (K' \star \{x\} \star \{y\}) \coprod_{ \bd \Delta^n \star \{x\} \star \{y\} } ( \Delta^n \star \{y\} )$$
is a pushout of
$$ (K' \star \{x\}) \coprod_{ \bd \Delta^n \star \{x\} } (\bd \Delta^n \star \{x\} \star \{y\} )
\subseteq K' \star \{x\} \star \{y\},$$ and therefore a categorical equivalence by Lemma \ref{doweneed}. It follows that $s$ admits an extension
$$ \widetilde{s}: (K' \star \{x\} \star \{y\}) \coprod_{ \bd \Delta^n \star \{x\} \star \{y\} } ( \Delta^n \star \{y\} ) \rightarrow \calC, $$ and we may now define $\widetilde{p} = \widetilde{s} | K \star \{y\}$. 
\end{proof}

\begin{proposition}\label{stook}
Let $\calC$ be a topological category. Then $\calC$ is filtered if and only if the $\infty$-category
$\tNerve(\calC)$ is filtered.
\end{proposition}

\begin{proof}
Suppose first that $\tNerve(\calC)$ is filtered. We verify conditions $(1')$ and $(2')$ of Definition \ref{topfilt}:

\begin{itemize}
\item[$(1')$] Let $\{X_i\}_{i \in I}$ be a finite collection of
objects of $\calC$, corresponding to a map $p: I \rightarrow \tNerve(\calC)$, where $I$ is regarded as a discrete simplicial set. If $\tNerve(\calC)$ is filtered, then $p$ extends to a map
$\widetilde{p}: I \star \{x\} \rightarrow \tNerve(\calC)$, corresponding to an object $X = p(x)$
equipped with maps $X_i \rightarrow X$ in $\calC$.

\item[$(2')$] Let $X, Y \in \calC$ be objects, $n \geq 0$, and $S^n \rightarrow \bHom_{\calC}(X,Y)$ a map. We note that this data may be identified with a topological functor $F: | \sCoNerve[K] | \rightarrow \calC$, where $K$ is the simplicial set obtained from $\bd \Delta^{n+2}$ by collapsing the initial face
$\Delta^{n+1}$ to a point. If $\tNerve(\calC)$ is filtered, then $F$ extends to a functor $\widetilde{F}$ defined on $| \sCoNerve[ K \star \{z\}] |$; this gives an object $Z = \widetilde{F}(z)$ and a morphism
$Y \rightarrow Z$ such that the induced map $S^n \rightarrow \bHom_{\calC}(X,Z)$ is nullhomotopic.
\end{itemize}

For the converse, let us suppose that $\calC$ is filtered. We wish to show that $\tNerve(\calC)$ is filtered. By Lemma \ref{goony}, it will suffice to prove that $\tNerve(\calC)$ has the extension property with respect to the inclusion $\bd \Delta^n \subseteq \Lambda^{n+1}_{n+1}$, for each $n \geq 0$. Equivalently, it suffices to show that $\calC$ has the right extension property with
respect to the inclusion $| \sCoNerve[ \bd \Delta^n] | \subseteq | \sCoNerve[ \Lambda^{n+1}_{n+1} ] |$.
If $n = 0$, this is simply the assertion that $\calC$ is nonempty; if $n = 1$, this is the assertion that for any pair of objects $X,Y \in \calC$ there exists an object $Z$ equipped with morphisms
$X \rightarrow Z$, $Y \rightarrow Z$. Both of these conditions follow from part $(1)$ of Definition \ref{topfilt}; we may therefore assume that $n > 1$.

Let $\calA_0 = | \sCoNerve[\bd \Delta^n] |$, $\calA_1 = | \sCoNerve[ \bd \Delta^n \coprod_{ \Lambda^n_n}  \Lambda^n_n \star \{ n+1 \}]|$, $\calA_2 = | \sCoNerve[ \Lambda^{n+1}_{n+1}] |$, and
$\calA_3 = | \sCoNerve[ \Delta^{n+1}] |$, so that we have inclusions of topological categories
$$ \calA_0 \subseteq \calA_1 \subseteq \calA_2 \subseteq \calA_3.$$

We will make use of the description of $\calA_3$ given in Remark \ref{conervexp}: its objects are integers $i$ satisfying $0 \leq i \leq n+1$, with $\bHom_{\calA_3}(i,j)$ given by the cube of all functions $p: \{i, \ldots, j\} \rightarrow [0,1]$ satisfying $p(i)=p(j)=1$ for $i \leq j$, and$\Hom_{\calA_3}(i,j) = \emptyset$ for $j < i$. Composition is given by amalgamation of functions.

We note that $\calA_1$ and $\calA_2$ are subcategories of $\calA_3$ having the same objects, where:
\begin{itemize}
\item $\bHom_{\calA_1}(i,j) = \bHom_{\calA_2}(i,j) = \bHom_{\calA_3}(i,j)$ unless $i=0$ and
$j \in \{n,n+1\}$. 

\item $\bHom_{\calA_1}(0,n) = \bHom_{\calA_2}(0,n)$ is the boundary of the cube $\bHom_{\calA_3}(0,n) = [0,1]^{n-1}$.

\item $\bHom_{\calA_1}(0,n+1)$ consists of all functions $p: [n+1] \rightarrow [0,1]$
satisfying $p(0) = p(n+1) =1$ and $(\exists i) [ (1 \leq i \leq n-1) \wedge p(i) \in \{0,1\} ]$.

\item $\bHom_{\calA_2}(0,n+1)$ is the union of $\bHom_{\calA_1}(0,n+1)$ with the collection
of functions $p: \{0, \ldots, n+1\} \rightarrow [0,1]$ satisfying $p(0) = p(n) = p(n+1) = 1$.
\end{itemize}

Finally, we note that $\calA_0$ is the full subcategory of $\calA_1$ (or $\calA_2$) whose set of objects is $\{0, \ldots, n\}$.

We wish to show that any topological functor $F: \calA_0 \rightarrow \calC$ can be extended to a 
functor $\widetilde{F}: \calA_2 \rightarrow \calC$. Let $X = F(0)$, $Y = F(n)$. Then $F$ induces a map $S^{n-1} \simeq \bHom_{\calA_0}(0,n) \rightarrow \bHom_{\calC}(X,Y)$. Since $\calC$ is filtered, there exists a map $\phi: Y \rightarrow Z$ such that the induced map
$f: S^{n-1} \rightarrow \bHom_{\calC}(X,Z)$ is nullhomotopic.

Now set $\widetilde{F}(n+1) = Z$; for $p \in \bHom_{\calA_1}(i, n+1)$, we set
$\widetilde{F}(p) = \phi \circ F(q)$, where $q \in \bHom_{\calA_1}(i,n)$ is such that
$q| \{i, \ldots, n-1 \} = p | \{i, \ldots, n-1\}$. Finally, we note that the assumption that
$f$ is nullhomotopic allows us to extend $\widetilde{F}$ from $\bHom_{\calA_1}(0,n+1)$ to the
whole of $\bHom_{\calA_2}(0,n+1)$.
\end{proof}

\begin{remark}\label{elisa}
Suppose that $\calC$ is a $\kappa$-filtered $\infty$-category, and let $K$ be a simplicial set which is categorically equivalent to a $\kappa$-small simplicial set. Then $\calC$ has the extension property with respect to the inclusion $K \subseteq K^{\triangleright}$. This follows from Proposition \ref{princex}: to test whether or not a map $K \rightarrow S$ extends over $K^{\triangleright}$, it suffices to check in the homotopy category of $\sSet$ (with respect to the Joyal model structure), where we may replace $K$ by an equivalent $\kappa$-small simplicial set. 
\end{remark}

\begin{proposition}\label{smallity}
Let $\calC$ be a $\infty$-category with a final object. Then $\calC$ is
$\kappa$-filtered for every regular cardinal $\kappa$. Conversely, if $\calC$ is $\kappa$-filtered and there exists a categorical
equivalence $K \rightarrow \calC$, where $K$ is a $\kappa$-small simplicial set, then $\calC$ has a final object.
\end{proposition}

\begin{proof}
We remark that $\calC$ has a final object if and only if there exists
a retraction $r$ of $\calC^{\triangleright}$ onto $\calC$. If $\calC$ is $\kappa$-filtered and categorically equivalent to a $\kappa$-small simplicial set, then the existence of such a retraction follows from Remark \ref{elisa}. On the other hand, if
the retraction $r$ exists, then any map $p: K \rightarrow \calC$
admits an extension $K^{\triangleright} \rightarrow \calC$: one merely
considers the composition $ K^{\triangleright} \rightarrow \calC^{\triangleright} \stackrel{r}{\rightarrow} \calC.$
\end{proof}

A useful observation from classical category theory is that, if we are only interested in using filtered categories to index colimit diagrams, then in fact we do not need the notion of a filtered category at all: we can work instead with diagrams indexed by filtered partially ordered sets. We now prove an $\infty$-categorical analogue of this statement.

\begin{proposition}\label{rot}
Suppose that $\calC$ is a $\kappa$-filtered $\infty$-category. Then there
exists a $\kappa$-filtered partially ordered set $A$ and a cofinal map
$\Nerve(A) \rightarrow \calC$.
\end{proposition}

\begin{proof}
The proof uses the ideas introduced in \S \ref{quasilimit1}, and in particular Proposition \ref{utl}.
Let $X$ be a set of size $\geq \kappa$, and regard $X$ as a
category with a unique isomorphism between any pair of objects. We
note that $\Nerve(X)$ is a contractible Kan complex; consequently the
projection $\calC \times \Nerve(X) \rightarrow \calC$ is cofinal. Hence, it
suffices to produce a cofinal map $\Nerve(A) \rightarrow \calC \times \Nerve(X)$
with the desired properties.

Let $\{ K_{\alpha} \}_{\alpha \in A}$ be the collection of all
simplicial subsets of $K = \calC \times \Nerve(X)$ which are $\kappa$-small and
possess a final vertex. Regard $A$ as a partially ordered by inclusion.
We first claim that $A$ is $\kappa$-filtered and that
$\bigcup_{\alpha \in A} K_{\alpha} = K$. To prove both
of these assertions, it suffices to show that any $\kappa$-small
simplicial subset $L \subseteq K$ is contained in a
$\kappa$-small simplicial subset $L'$ which has a final
vertex.

Since $\calC$ is $\kappa$-filtered, the composition
$$L \rightarrow \calC
\times \Nerve(X) \rightarrow \calC$$ extends to a map $p: L^{\triangleright} \rightarrow \calC$.
Since $X$ has cardinality $\geq \kappa$, there exists an element
$x \in X$ which is not in the image of $L_0 \rightarrow N(X)_0 = X$.
Lift $p$ to a map $\widetilde{p}: L^{\triangleright} \rightarrow K$ which extends the inclusion $L \subseteq K \times
\Nerve(X)$ and carries the cone point to the element $x \in X = N(X)_0$. It
is easy to see that $\widetilde{p}$ is injective, so that we may
regard $L^{\triangleright}$ as a simplicial subset of $K \times
\Nerve(X)$. Moreover, it is clearly $\kappa$-small and has a
final vertex, as desired.

Now regard $A$ as a category, and let $F: A \rightarrow (\sSet)_{/K}$ be the functor
which carries each $\alpha \in A$ to the simplicial set $K_{\alpha}$. For each
$\alpha \in A$, choose a final vertex $x_\alpha$ of $K_{\alpha}$.
Let $K_F$ be defined as in \S \ref{quasilimit1}. We claim next that there 
exists a retraction $r: K_F \rightarrow K$ with
the property that $r(X_{\alpha}) = x_{\alpha}$ for each $I \in \calI$.

The construction of $r$ proceeds as in the proof of Proposition
\ref{extet}. Namely, we well-order the finite linearly ordered subsets $B
\subseteq A$, and define $r|K'_B$ by induction on $B$.
Moreover, we will select $r$ so that it has the property that if
$B$ is nonempty with largest element $\beta$, then $r(K'_B)
\subseteq K_{\beta}$.

If $B$ is empty, then $r|K'_B=r|K$ is the identity map. Otherwise,
$B$ has a least element $\alpha$ and a largest element $\beta$. We are
required to construct a map $K_{\alpha} \star \Delta^B \rightarrow
K_{\beta}$, or a map $r_B: \Delta^B \rightarrow K_{\id|K_{\alpha}/}$, where
the values of this map on $\bd \Delta^B$ have already been
determined. If $B$ is a singleton, we define this map to carry the
vertex $\Delta^B$ to a final object of $K_{\id|K_{\alpha}/}$ lying
over $x_{\beta}$. Otherwise, we are guaranteed that {\em some}
extension exists by the fact that $r_B|\bd \Delta^B$ carries the
final vertex of $\Delta^B$ to a final object of
$K_{\id|K_{\alpha}/}$.

Now let $j: \Nerve(A) \rightarrow K$ denote the restriction
of the retraction of $r$ to $\Nerve(A)$. Using Propositions \ref{extet}
and \ref{utl}, we deduce that $j$ is a cofinal map as desired.
\end{proof}

A similar technique can be used to prove the following characterization of
$\kappa$-filtered $\infty$-categories:

\begin{proposition}\label{charfiltt}
Let $S$ be a simplicial set. The following conditions are equivalent:
\begin{itemize}
\item[$(1)$] The simplicial set $S$ is $\kappa$-filtered.
\item[$(2)$] There exists a diagram of simplicial sets $\{ Y_{\alpha} \}_{ \alpha \in \calI}$ having
colimit $Y$ and a categorical equivalence $S \rightarrow Y$, 
where each $Y_{\alpha}$ is $\kappa$-filtered and the indexing category $\calI$ is $\kappa$-filtered.
\item[$(3)$] There exists a categorical equivalence $S \rightarrow \calC$ where $\calC$ is a $\kappa$-filtered union of simplicial subsets $\calC_{\alpha} \subseteq \calC$ such that each $\calC_{\alpha}$ is an $\infty$-category with a final object.
\end{itemize}
\end{proposition}

\begin{proof}
Let $T: \sSet \rightarrow \sSet$ be the ``fibrant replacement'' functor given by
$$ T(X) = \tNerve( | \sCoNerve[X] |).$$
There is a natural transformation $j_{X}: X \rightarrow T(X)$ which is a categorical equivalence every simplicial set $X$ the map $j_{X}$ is a categorical equivalence. Moreover, each $T(X)$ is an $\infty$-category. Furthermore, the functor $T$ preserves inclusions and commutes with filtered colimits.

It is clear that $(3)$ implies $(2)$. Suppose that $(2)$ is satisfied. Replacing the diagram
$\{ Y_{\alpha} \}_{ \alpha \in \calI}$ by $\{ T(Y_{\alpha}) \}_{\alpha \in \calI}$ if necessary, we may suppose that each $Y_{\alpha}$ is an $\infty$-category. It follows that $Y$ is an $\infty$-category.
If $p: K \rightarrow Y$ is a diagram indexed by a $\kappa$-small simplicial set, then $p$ factors
through a map $p_{\alpha}: K \rightarrow Y_{\alpha}$ for some $\alpha \in \calI$, in virtue of the assumption that $\calI$ is $\kappa$-filtered. Since $Y_{\alpha}$ is a $\kappa$-filtered $\infty$-category, we can find an extension $K^{\triangleright} \rightarrow Y_{\alpha}$ of $p_{\alpha}$, hence an extension $K^{\triangleright} \rightarrow Y$ of $p$.

Now suppose that $(1)$ is satisfied. Replacing $S$ by $T(S)$ if necessary, we may suppose that $S$ is an $\infty$-category. Choose a set $X$ of cardinality $\geq \kappa$, and let $\Nerve(X)$ be defined as in the proof of Proposition \ref{rot}. The proof of Proposition \ref{rot} shows that we may write $S \times \Nerve(X)$ as a $\kappa$-filtered union of simplicial subsets $\{ Y_{\alpha} \}$, where
each $Y_{\alpha}$ has a final vertex. We now take $\calC = T(S \times \Nerve(X) )$, and let
$\calC_{\alpha} = T(Y_{\alpha})$: these choices satisfy $(3)$, which completes the proof.
\end{proof}

By definition, a $\infty$-category $\calC$ is $\kappa$-filtered if any map $p: K \rightarrow \calC$, whose source $K$ is $\kappa$-small, extends over the cone $K^{\triangleright}$. We now consider the possibility of constructing this extension uniformly in $p$. First, we need a few lemmas.

\begin{lemma}\label{stull2}
Let $\calC$ be a filtered $\infty$-category. Then $\calC$ is weakly contractible.
\end{lemma}

\begin{proof}
Since $\calC$ is filtered, it is nonempty. Fix an object $C \in \calC$. Let $| \calC | $ denote the geometric realization of $\calC$ as a simplicial set. We identify $C$ with a point of the topological space $| \calC |$. By Whitehead's theorem, to show that $\calC$ is weakly contractible, it suffices to show that for every $i \geq 0$, the homotopy set $\pi_{i}( |\calC|, C)$ consists of a single point.
If not, we can find a finite simplicial subset $K \subseteq \calC$ containing $C$ such that the map $f: \pi_{i}( |K|, C) \rightarrow \pi_{i}( | \calC|,C)$ has nontrivial image. But $\calC$ is filtered, so the inclusion $K \subseteq \calC$ factors through a map $K^{\triangleright} \rightarrow \calC$.
It follows that $f$ factors through $\pi_{i}( |K^{\triangleright}|, C)$. But this homotopy set is trivial, since $K^{\triangleright}$ is weakly contractible.
\end{proof}

\begin{lemma}\label{forfilt}
Let $\calC$ be a $\kappa$-filtered $\infty$-category, and let $p: K \rightarrow \calC$
be a diagram indexed by a $\kappa$-small simplicial set $K$. Then $\calC_{p/}$ is $\kappa$-filtered.
\end{lemma}

\begin{proof}
Let $K'$ be a $\kappa$-small simplicial set, and $p': K' \rightarrow \calC_{p/}$ a $\kappa$-small diagram. Then we may identify $p'$ with a map $q: K \star K' \rightarrow \calC$, and we get an isomorphism $( \calC_{p/} )_{p'/} \simeq \calC_{q/}$. Since $K \star K'$ is $\kappa$-small,
the $\infty$-category $\calC_{q/}$ is nonempty.
\end{proof}

\begin{proposition}\label{undertruck}
Let $\calC$ be an $\infty$-category and $\kappa$ a regular cardinal. Then $\calC$ is $\kappa$-filtered if and only if, for each $\kappa$-small simplicial set $K$, the diagonal map
$d: \calC \rightarrow \Fun(K,\calC)$ is cofinal.
\end{proposition}

\begin{proof}
Suppose first that the diagonal map $d: \calC \rightarrow \Fun(K,\calC)$ is cofinal, for any $\kappa$-small simplicial set $K$. Choose any map $j: K \rightarrow \calC$; we wish to show that $j$ can be extended to $K^{\triangleright}$. By Proposition \ref{princex}, it suffices to show that $j$ can be extended to the equivalent simplicial set $K \diamond \Delta^0$. In other words, we must produce
an object $C \in \calC$ and a morphism $j \rightarrow d(C)$ in $\Fun(K,\calC)$. It will suffice to prove that
the $\infty$-category $\calD = \calC \times_{ \Fun(K,\calC) } \Fun(K,\calC)_{j/}$ is nonempty. We now invoke Theorem \ref{hollowtt} to deduce that $\calD$ is weakly contractible.

Now suppose that $S$ is $\kappa$-filtered, and that $K$ is a $\kappa$-small simplicial set.
We wish to show that the diagonal map $d: \calC \rightarrow \Fun(K,\calC)$ is cofinal. By Theorem \ref{hollowtt}, it suffices to prove that for every object $X \in \Fun(K,\calC)$, the $\infty$-category
$\Fun(K,\calC)^{X/} \times_{\Fun(K,\calC)} \calC$ is weakly contractible. But if we identify $X$ with a map
$x: K \rightarrow \calC$, then we get a natural identification
$$ \Fun(K,\calC)^{X/} \times_{ \Fun(K,\calC) } \calC \simeq \calC^{x/},$$ which is $\kappa$-filtered
by Lemma \ref{forfilt} and therefore weakly contractible by Lemma \ref{stull2}.
\end{proof}

\subsection{Right Exactness}\label{rexex}

Let $\calA$ and $\calB$ be abelian categories. In classical homological algebra, a functor $F: \calA \rightarrow \calB$ is said to be
{\it right exact} if it is additive, and whenever
$$A' \rightarrow A \rightarrow A'' \rightarrow 0$$
is an exact sequence in $\calA$, the induced sequence
$$F(A') \rightarrow F(A) \rightarrow F(A'') \rightarrow 0$$
is exact in $\calB$.\index{gen}{right exact functor}\index{gen}{functor!right exact}

The notion of right exactness generalizes in a natural way to functors between categories which are not assumed to be abelian. Let $F: \calA \rightarrow \calB$ be a functor between abelian categories, as above. Then $F$ is additive if and only if $F$ preserves finite coproducts. Furthermore, an additive functor $F$ is right exact if and only if it preserves coequalizer diagrams. Since every finite colimit can be built out of finite coproducts and coequalizers, right exactness is equivalent to the requirement that $F$ preserves all finite colimits. This condition makes sense whenever the category $\calA$ admits finite colimits.

It is possible to generalize even further, to the case of a functor between arbitrary categories. To simplify the discussion, let us suppose that $\calB = \Set^{op}$. Then we may regard a functor $F: \calA \rightarrow \calB$ as a presheaf of sets on the category $\calA$. Using this presheaf we can define a new category $\calA_{F}$, whose objects are pairs $(A, \eta)$ where $A \in \calA$ and
$\eta \in F(A)$, and morphisms from $(A, \eta)$ to $(A', \eta')$ are maps $f: A \rightarrow A'$
such that $f^{\ast}(\eta') = \eta$, where $f^{\ast}$ denotes the induced map
$F(A') \rightarrow F(A)$. If $\calA$ admits finite colimits, then the functor $F$ preserves finite colimits if and only if the category $\calA_{F}$ is filtered.

Our goal in this section is to adapt the notion of right-exact functors to the $\infty$-categorical context. We begin with the following:

\begin{definition}\index{gen}{$\kappa$-right exact}\index{gen}{functor!$\kappa$-right exact}\label{spuss}
Let $F: \calA \rightarrow \calB$ be a functor between $\infty$-categories and $\kappa$ a regular cardinal. We will say that
$F$ is {\it $\kappa$-right exact} if, for any right fibration $\calB' \rightarrow \calB$
where $\calB'$ is $\kappa$-filtered, the $\infty$-category $\calA' = \calA \times_{\calB} \calB'$ is also $\kappa$-filtered. We will say that $F$ is {\it right exact} if it is $\omega$-right exact.
\end{definition}

\begin{remark}
We also have an evident dual notion of {\em left exact} functor.\index{gen}{left exact!functor}
\index{gen}{functor!left exact}
\end{remark}

\begin{remark}
If $\calA$ admits finite colimits, then a functor $F: \calA \rightarrow \calB$ is right exact
if and only if $F$ preserves finite colimits (see Proposition \ref{swarmmy} below).
\end{remark}

We note the following basic stability properties of $\kappa$-right exact maps.

\begin{proposition}
Let $\kappa$ be a regular cardinal.
\begin{itemize}
\item[$(1)$] If $F: \calA \rightarrow \calB$ and $G: \calB \rightarrow \calC$ are $\kappa$-right
exact functors between $\infty$-categories, then $G \circ F: \calA \rightarrow \calC$ is $\kappa$-right exact.
\item[$(2)$] Any equivalence of $\infty$-categories is $\kappa$-right exact.
\item[$(3)$] Let $F: \calA \rightarrow \calB$ be a $\kappa$-right exact functor, and let
$F': \calA \rightarrow \calB$ be homotopic to $F$. Then $F'$ is $\kappa$-right exact.
\end{itemize}
\end{proposition}

\begin{proof}
Property $(1)$ is immediate from the definition. We will establish $(2)$ and $(3)$ as a consequence of the following more general assertion: if $F: \calA \rightarrow \calB$ and $G: \calB \rightarrow \calC$ are functors such that $F$ is a categorical equivalence, then $G$ is $\kappa$-right exact if and only if $G \circ F$ is $\kappa$-right exact. To prove this, let $\calC' \rightarrow \calC$ be a right fibration. Proposition \ref{basechangefunky} implies that the induced map
$$ \calA' = \calA \times_{\calC} \calC' \rightarrow \calB \times_{\calC} \calC' = \calB'$$
is a categorical equivalence. Thus $\calA'$ is $\kappa$-filtered if and only if $\calB'$ is $\kappa$-filtered.

We now deduce $(2)$ by specializing to the case where $G$ is the identity map. To prove $(3)$,
we choose a contractible Kan complex $K$ containing a pair of vertices $\{x,y\}$ and a map $g: K \rightarrow \calB^{\calA}$ with $g(x) = F$, $g(y) = F'$. Applying the above argument to the composition
$$ \calA \simeq \calA \times \{x\} \subseteq \calA \times K \stackrel{G}{\rightarrow} \calB,$$
we deduce that $G$ is $\kappa$-right exact. Applying the converse to the diagram
$$ \calA \simeq \calA \times \{y\} \subseteq \calA \times K \stackrel{G}{\rightarrow} \calB$$
we deduce that $F'$ is $\kappa$-right exact.
\end{proof}

The next result shows that the $\kappa$-right exactness of a functor $F: \calA \rightarrow \calB$ can be tested on a very small collection of right fibrations $\calB' \rightarrow \calB$.

\begin{proposition}\label{swarmy}
Let $F: \calA \rightarrow \calB$ be a functor between $\infty$-categories and $\kappa$ a regular cardinal. The following are equivalent:
\begin{itemize}
\item[$(1)$] The functor $F$ is $\kappa$-right exact.
\item[$(2)$] For every object $B$ of $\calB$, the $\infty$-category
$ \calA \times_{\calB} \calB_{/B}$ is $\kappa$-filtered.
\end{itemize}
\end{proposition}

\begin{proof}
We observe that for every object $B \in \calB$, the $\infty$-category $\calB_{/B}$ is right-fibered
over $\calB$ and is $\kappa$-filtered (since it has a final object). Consequently, $(1)$ implies $(2)$. Now suppose that $(2)$ is satisfied. Let $T: (\sSet)_{/\calB} \rightarrow (\sSet)_{/\calB}$
denote the composite functor
$$ (\sSet)_{/\calB} \stackrel{\St_{\calB}}{\rightarrow} (\sSet)^{\sCoNerve[\calB^{op}]}
\stackrel{\Sing |\bigdot|}{\rightarrow} (\sSet)^{\sCoNerve[\calB^{op}]} \stackrel{\Un_{\calB}}{\rightarrow} (\sSet)_{/\calB}.$$
We will use the following properties of $T$:
\begin{itemize}
\item[$(i)$] There is a natural transformation $j_{X}: X \rightarrow T(X)$, where $j_{X}$ is
a contravariant equivalence in $(\sSet)_{/\calB}$ for every $X \in (\sSet)_{/\calB}$.
\item[$(ii)$] For every $X \in (\sSet)_{/\calB}$, the associated map $T(X) \rightarrow \calB$ is
a right fibration.
\item[$(iii)$] The functor $T$ commutes with filtered colimits.
\end{itemize}
We will say that an object $X \in (\sSet)_{/\calB}$ is {\it good} if the $\infty$-category
$T(X) \times_{\calB} \calA$ is $\kappa$-filtered. We now make the following observations:
\begin{itemize}
\item[$(A)$] If $X \rightarrow Y$ is a contravariant equivalence in $(\sSet)_{/\calB}$, then
$X$ is good if and only if if $Y$ is good. This follows from the fact that $T(X) \rightarrow T(Y)$
is an equivalence of right fibrations, so that the induced map $T(X) \times_{\calB} \calA
\rightarrow T(Y) \times_{\calB} \calA$ is an equivalence of right fibrations and consequently a categorical equivalence of $\infty$-categories.
\item[$(B)$] If $X \rightarrow Y$ is a categorical equivalence in $(\sSet)_{/\calB}$, then
$X$ is good if and only if $Y$ is good. This follows $(A)$, since every
categorical equivalence is a contravariant equivalence.
\item[$(C)$] The collection of good objects of $(\sSet)_{\calB}$ is stable under $\kappa$-filtered colimits. This follows from the fact that the functor $X \mapsto T(X) \times_{\calB} \calA$ commutes with $\kappa$-filtered colimits (in fact, with all filtered colimits) and Proposition \ref{charfiltt}.
\item[$(D)$] If $X \in (\sSet)_{/\calB}$ corresponds to a right fibration $X \rightarrow \calB$, then
$X$ is good if and only if $X \times_{ \calB} \calA$ is $\kappa$-filtered.
\item[$(E)$] For every object $B \in \calB$, the overcategory $\calB_{/B}$ is a good
object of $(\sSet)_{/\calB}$. In view of $(D)$, this is equivalent to the assumption $(2)$.
\item[$(F)$] If $X$ consists of a single vertex $x$, then $X$ is good. To see this, let
$B \in \calB$ denote the image of $X$. The natural map $X \rightarrow \calB_{/B}$
can be identified with the inclusion of a final vertex; this map is right anodyne and therefore
a contravariant equivalence. We now conclude by applying $(A)$ and $(E)$.
\item[$(G)$] If $X \in (\sSet)_{/\calB}$ is an $\infty$-category with a final object $x$, then
$X$ is good. To prove this, we note that $\{x\}$ is good by $(F)$ and the inclusion
$\{x\} \subseteq X$ is right anodyne, hence a contravariant equivalence. We conclude by applying $(A)$.
\item[$(H)$] If $X \in (\sSet)_{/\calB}$ is $\kappa$-filtered, then $X$ is good. To prove this, we apply
Proposition \ref{charfiltt} to deduce the existence of a categorical equivalence $i:X \rightarrow \calC$, where $\calC$ is a $\kappa$-filtered union of $\infty$-categories with final objects. Replacing $\calC$ by $\calC \times K$ if necessary, where $K$ is a contractible Kan complex, we may suppose that $i$ is a cofibration. Since $\calB$ is an $\infty$-category, the lifting problem
$$ \xymatrix{ S \ar[r] \ar[d]^{i} & \calB \\
\calC \ar@{-->}[ur] }$$
has a solution. Thus we may regard $\calC$ as an object of $(\sSet)_{/\calB}$.
According to $(B)$, it suffices to show that $\calC$ is good. But $\calC$ is a $\kappa$-filtered
colimit of good objects of $(\sSet)_{\calB}$ (by $(G)$), and is therefore itself good (by $(C)$).
\end{itemize}

Now let $\calB' \rightarrow \calB$ be a right fibration, where $\calB'$ is $\kappa$-filtered.
By $(H)$, $\calB'$ is a good object of $(\sSet)_{/\calB}$. Applying $(D)$, we deduce that
$\calA' = \calB' \times_{\calB} \calA$ is $\kappa$-filtered. This proves $(1)$.
\end{proof}

Our next goal is to prove Proposition \ref{swarmmy}, which gives a very concrete characterization of right exactness under the assumption that there is a sufficient supply of colimits. We first need a few preliminary results.

\begin{lemma}\label{devass}
Let $\calB' \rightarrow \calB$ be a Cartesian fibration. Suppose that $\calB$ has an initial object
$B$ and that $\calB'$ is filtered. Then the fiber $\calB'_{B} = \calB' \times_{\calB} \{B\}$ is a contractible Kan complex.
\end{lemma}

\begin{proof}
Since $B$ is an initial object of $\calB$, the inclusion $\{ B\}^{op} \subseteq \calB^{op}$ is cofinal. Proposition \ref{strokhop} implies that the inclusion $(\calB'_{B})^{op} \subseteq (\calB')^{op}$ is also cofinal, and therefore a weak homotopy equivalence. It now suffices to prove that $\calB'$ is weakly contractible, which follows from Lemma \ref{stull2}.
\end{proof}

\begin{lemma}\label{druv}
Let $f: \calA \rightarrow \calB$ be a right exact functor between $\infty$-categories, and let
$A \in \calA$ be an initial object. Then $f(A)$ is an initial object of $\calB$.
\end{lemma}

\begin{proof}
Let $B$ be an object of $\calB$. Proposition \ref{swarmy} implies that
$\calA' = \calB_{/B} \times_{\calB} \calA$ is filtered. We may identify
$\bHom_{\calB}(f(A), B)$ with the fiber of the right fibration $\calA' \rightarrow \calA$ over the object $A$. We now apply Lemma \ref{devass} to deduce that $\bHom_{\calB}(f(A),B)$ is contractible.
\end{proof}

\begin{lemma}\label{devic}
Let $\kappa$ be a regular cardinal, $f: \calA \rightarrow \calB$ a $\kappa$-right exact functor between $\infty$-categories, and $p: K \rightarrow \calA$ be a diagram indexed by a $\kappa$-small simplicial set $K$.
The induced map $\calA_{p/} \rightarrow \calB_{f p/}$ is $\kappa$-right exact.
\end{lemma}

\begin{proof}
According to Proposition \ref{swarmy}, it suffices to prove that for each object
$\overline{B} \in \calB_{f \circ p/}$, the $\infty$-category
$\calA' = \calA_{p/} \times_{\calB_{f p/}} ( \calB_{f p/} )_{/\overline{B}}$ is $\kappa$-filtered.
Let $B$ denote the image of $\overline{B}$ in $\calB$, and let
$q: K' \rightarrow \calA'$ be a diagram indexed by a $\kappa$-small simplicial set $K'$;
we wish to show that $q$ admits an extension to ${K'}^{\triangleright}$. We may regard $p$ and $q$ together as defining a diagram
$K \star K' \rightarrow \calA \times_{\calB} \calB_{/B}$. Since $f$ is $\kappa$-filtered,
we can extend this to a map
$(K \star K')^{\triangleright} \rightarrow \calA \times_{\calB} \calB_{/B}$, which can be identified with an extension $\overline{q}: {K'}^{\triangleright} \rightarrow \calA'$ of $q$.
\end{proof}

\begin{proposition}\label{swarmmy}\index{gen}{right exact!and colimits}
Let $f: \calA \rightarrow \calB$ be a functor between $\infty$-categories and let $\kappa$ be a regular cardinal.

\begin{itemize}
\item[$(1)$] If $f$ is $\kappa$-right exact, then $f$ preserves all $\kappa$-small colimits
which exist in $\calA$.
\item[$(2)$] Conversely, if $\calA$ admits $\kappa$-small colimits and $f$ preserves $\kappa$-small colimits, then $f$ is right exact.
\end{itemize}
\end{proposition}

\begin{proof}
Suppose first that $f$ is $\kappa$-right exact. Let $K$ be a $\kappa$-small simplicial set, and
let $\overline{p}: K^{\triangleright} \rightarrow \calA$ be a colimit of $p = \overline{p}|K$. 
We wish to show that $f \circ \overline{p}$ is a colimit diagram. Using Lemma \ref{devic}, we may replace $\calA$ by $\calA_{p/}$ and $\calB$ by $\calB_{f p/}$, and thereby reduce to the case $K = \emptyset$. We are then reduced to proving that $f$ preserves initial objects, which follows from Lemma \ref{druv}.

Now suppose that $\calA$ admits $\kappa$-small colimits, and that $f$ preserves $\kappa$-small colimits. We wish to prove that $f$ is $\kappa$-right exact. Let $B$ be an object of $\calB$ and set
$\calA' = \calA \times_{\calB} \calB_{/B}$. We wish to prove that $\calA'$ is $\kappa$-filtered.
Let $p': K \rightarrow \calA'$ be a diagram indexed by a $\kappa$-small simplicial set $K$; we wish to prove that $p'$ extends to a map $\overline{p}': K^{\triangleright} \rightarrow \calA'$. Let $p: K \rightarrow \calA$ be the composition of $p'$ with the projection $\calA' \rightarrow \calA$, and
let $\overline{p}: K^{\triangleright} \rightarrow \calA$ be a colimit of $p$. We may identify
$f \circ \overline{p}$ and $p'$ with objects of $\calB_{f p/}$. Since $f$ preserves $\kappa$-small colimits, $f \circ \overline{p}$ is an initial object of $\calB_{f p/}$, so that there exists a morphism
$\alpha: f \circ \overline{p} \rightarrow p'$ in $\calB_{f \circ p/}$. The morphism $\alpha$
can be identified with the desired extension $\overline{p}': K^{\triangleright} \rightarrow \calA'$.
\end{proof}

\begin{remark}
The results of this section all dualize in an evident way: a functor $G: \calA \rightarrow \calB$
is said to be {\it $\kappa$-left exact} if the induced functor $G^{op}: \calA^{op} \rightarrow \calB^{op}$ is $\kappa$-right exact. In the case where $\calA$ admits $\kappa$-small limits, this is equivalent to the requirement that $G$ preserves $\kappa$-small limits.
\end{remark}

\begin{remark}
Let $\calC$ be an $\infty$-category, and let $F: \calC \rightarrow \SSet^{op}$ be a functor, and let $\widetilde{\calC} \rightarrow \calC$ be the associated right fibration (the pullback of the universal right fibration $\calQ^0 \rightarrow \SSet^{op}$). If $F$ is $\kappa$-right exact, then
$\widetilde{\calC}$ is $\kappa$-filtered (since $\calQ^0$ has a final object). If
$\calC$ admits $\kappa$-small colimits, then the converse holds: if
$\widetilde{\calC}$ is $\kappa$-filtered, then $F$ preserves $\kappa$-small colimits by
Proposition \ref{geort}, and is therefore $\kappa$-right exact by Proposition \ref{swarmy}.
The converse does not hold in general: it is possible to give an example of right fibration
$\widetilde{\calC} \rightarrow \calC$ such that $\widetilde{\calC}$ is filtered, yet the classifying functor $F: \calC \rightarrow \SSet^{op}$ is not right exact.
\end{remark}

\subsection{Filtered Colimits}\label{fcolm}

Filtered categories tend not to be very interesting in themselves. Instead, they are primarily useful for indexing diagrams in other categories. 
This is because the colimits of filtered diagrams enjoy certain exactness properties which are not shared by colimits in general. In this section, we will formulate and prove these exactness properties in the $\infty$-categorical setting. First, we need a few definitions.

\begin{definition}\label{bicard}\index{gen}{$\kappa$-closed}
Let $\kappa$ be a regular cardinal. We will say that an
$\infty$-category $\calC$ is {\it $\kappa$-closed} if every diagram
$p: K \rightarrow \calC$ indexed by a $\kappa$-small simplicial set $K$
admits a colimit $\overline{p}: K^{\triangleright} \rightarrow \calC$.
\end{definition}

In a $\kappa$-closed $\infty$-category, it is possible to construct $\kappa$-small colimits
functorially. More precisely, suppose that $\calC$ is an $\infty$-category and that $K$ is a simplicial set with the property that every diagram $p: K \rightarrow \calC$ has a colimit in $\calC$. Let $\calD$ denote the full subcategory of $\Fun(K^{\triangleright}, \calC)$ spanned by the colimit diagrams. Proposition \ref{lklk} implies that the restriction functor
$\calD \rightarrow \Fun(K,\calC)$ is a trivial fibration. It therefore admits a section $s$ (which is unique up to a contractible ambiguity). Let $e: \Fun(K^{\triangleright}, \calC) \rightarrow \calC$ be the functor given by evaluation at the cone point of $K^{\triangleright}$. We will refer to the composition
$$ \Fun(K,\calC) \stackrel{s}{\rightarrow} \calD \subseteq \Fun(K^{\triangleright}, \calC) \stackrel{e}{\rightarrow} \calC$$
as a {\it colimit} functor; it associates to each diagram $p: K \rightarrow \calC$ a colimit of $p$
in $\calC$. We will generally denote colimit functors by $\colim_{K}: \Fun(K,\calC) \rightarrow \calC$.\index{gen}{colimit!functor}\index{gen}{functor!colimit}\index{not}{colimK@$\colim_{K}$}

\begin{lemma}\label{tractab}
Let $F \in \Fun(K,\SSet)$ be a {\em corepresentable} functor (that is, $F$ lies in the essential image of the Yoneda embedding $K^{op} \rightarrow \Fun(K,\SSet)$), and let $X \in \SSet$ be a colimit of $F$. Then $X$ is contractible.
\end{lemma}

\begin{proof}
Without loss of generality, we may suppose that $K$ is an $\infty$-category. Let $\widetilde{K} \rightarrow K$ be a left fibration classified by $F$. Since $F$ is corepresentable, $\widetilde{K}$ has an initial object and is therefore weakly contractible. Corollary \ref{needka} implies that
there is an isomorphism $\widetilde{K} \simeq X$ in the homotopy category $\calH$, so that $X$ is also contractible.
\end{proof}

\begin{proposition}\label{frent}\index{gen}{filtered colimit!left exactness of}
Let $\kappa$ be a regular cardinal and let $\calI$ be an $\infty$-category. The following conditions are equivalent:
\begin{itemize}
\item[$(1)$] The $\infty$-category $\calI$ is $\kappa$-filtered.
\item[$(2)$] The colimit functor $\colim_{\calI}: \Fun(\calI, \SSet) \rightarrow \SSet$
preserves $\kappa$-small limits.
\end{itemize}
\end{proposition}

\begin{proof}
Suppose that $(1)$ is satisfied. According to Proposition \ref{rot}, there exists a $\kappa$-filtered partially ordered set $A$ and a cofinal map $i: \Nerve(A) \rightarrow \SSet$. Since $i$ is cofinal, the colimit functor for $\calI$ admits a factorization
$$ \Fun(\calI,\SSet) \stackrel{i^{\ast}}{\rightarrow} \Fun(\Nerve(A), \SSet) {\rightarrow} \SSet.$$
Proposition \ref{limiteval} implies that $i^{\ast}$ preserves limits. We may therefore replace
$\calI$ by $\Nerve(A)$ and thereby reduce to the case where $\calI$ is itself the nerve of a $\kappa$-filtered partially ordered set $A$.

We note that the functor $\colim_{\calI}: \Fun(\calI, \SSet) \rightarrow \SSet$ can be characterized as the
left adjoint to the diagonal functor $\delta: \SSet \rightarrow \Fun(\calI,\SSet)$. Let $\bfA$ denote the category
of all functors from $A$ to $\sSet$; we regard $\bfA$ as a simplicial model category with respect to the {\em projective} model structure described in \S \ref{quasilimit3}. Let $\phi^{\ast}: \sSet \rightarrow \bfA$ denote the diagonal functor which associates to each simplicial set $K$ the constant functor $A \rightarrow \sSet$ with value $K$, and let $\phi_{!}$ be a left adjoint of
$\phi^{\ast}$, so that the pair $(\phi^{\ast}, \phi_{!})$ gives a Quillen adjunction between
$\bfA$ and $\sSet$. Proposition \ref{gumby444} implies that there is an equivalence of $\infty$-categories $\sNerve(\bfA^{\degree}) \rightarrow
\Fun(\calI,\SSet)$, and $\delta$ may be identified with the right derived functor of $\phi^{\ast}$. Consequently, the functor $\colim_{\calI}$ may be identified with the left derived functor of $\phi_{!}$. To prove that $\colim_{\calI}$ preserves $\kappa$-small limits, it suffices to prove that $\colim_{\calI}$ preserves fiber products and $\kappa$-small products. According to Theorem \ref{colimcomparee}, it suffices to prove that
$\phi_{!}$ preserves homotopy fiber products and $\kappa$-small homotopy products. For fiber products, this reduces to the classical assertion that if we are given a family of homotopy
Cartesian squares
$$ \xymatrix{ W_{\alpha} \ar[r] \ar[d] & X_{\alpha} \ar[d] \\
Y_{\alpha} \ar[r] & Z_{\alpha} }$$
in the category of Kan complexes, indexed by a filtered partially ordered set $A$, then the colimit square
$$ \xymatrix{ W \ar[r] \ar[d] & X \ar[d] \\
Y \ar[r] & Z }$$
is also homotopy Cartesian. The assertion regarding homotopy products is handled similarly.

Now suppose that $(2)$ is satisfied. Let $K$ be a $\kappa$-small simplicial set and
$p: K \rightarrow \calI^{op}$ a diagram; we wish to prove that $\calI^{op}_{/p}$ is nonempty.
Suppose otherwise. Let $j: \calI^{op} \rightarrow \Fun(\calI,\SSet)$ be the Yoneda embedding, let
$q = j \circ p$, and let $\overline{q}: K^{\triangleleft} \rightarrow \Fun(\calI,\SSet)$ be a limit of
$q$, and let $X \in \Fun(\calI,\SSet)$ be the image of the cone point of $K^{\triangleleft}$
under $\overline{q}$. Since $j$ is fully faithful and $\calI^{op}_{/p}$ is empty, we have
$\bHom_{ \SSet^{\calI} }(j(I), X) = \emptyset$ for each $I \in \calI$. Using Lemma \ref{repco}, we may identify $\bHom_{ \SSet^{\calI} }(j(I),X)$ with $X(I)$ in the homotopy category $\calH$ of spaces. We therefore conclude that $X$ is an initial object of $\Fun(\calI,\SSet)$. Since the functor $\colim_{\calI}: \Fun(\calI,\SSet) \rightarrow \SSet$ is a left adjoint, it preserves initial objects. We conclude that $\colim_{\calI} X$ is an initial object of $\SSet$.
On the other hand, if $\colim_{\calI}$ preserves $\kappa$-small limits, then
$\colim_{\calI} \circ \overline{q}$ exhibits $\colim_{\calI} X$ as the limit of the diagram
$\colim_{\calI} \circ q: K \rightarrow \SSet$. For each vertex $k$ in $K$, Lemmas \ref{repco} and \ref{tractab} imply that $\colim_{\calI} q(k)$ is contractible, and therefore a final object of $\SSet$. It follows that
$\colim_{\calI} X$ is also a final object of $\SSet$. This is a contradiction, since the initial object
of $\SSet$ is not final.
\end{proof}

\subsection{Compact Objects}\label{compobj}

Let $\calC$ be a category which admits filtered colimits. An object $C \in \calC$ is said to be {\it compact} if the corepresentable functor $$ \Hom_{\calC}(C, \bigdot)$$ commutes with filtered colimits.\index{gen}{compact object!of a category}

\begin{example}
Let $\calC = \Set$ be the category of sets. An object $C \in \calC$ is compact if and only if
is finite.
\end{example}

\begin{example}\label{compactgroup}
Let $\calC$ be the category of groups. An object $G$ of $\calC$ is compact if and only if 
it is finitely presented (as a group).
\end{example}

\begin{example}\label{compactspacee}
Let $X$ be a topological space, and let $\calC$ be the category of open sets of $X$ (with morphisms given by inclusions). Then an object $U \in \calC$ is compact if and only if $U$ is compact when viewed as a topological space: that is, every open cover of $U$ admits a finite subcover.
\end{example}

\begin{remark}
Because of Example \ref{compactgroup}, many authors call an object $C$ of a category $\calC$ {\it finitely presented} if $\Hom_{\calC}(C, \bigdot)$ preserves filtered colimits. Our terminology is motivated instead by Example \ref{compactspacee}.
\end{remark}

\begin{definition}\label{kcontdef}\index{gen}{functor!$\kappa$-continuous}\index{gen}{functor!continuous}\index{gen}{continuous functor}\index{gen}{$\kappa$-continuous functor}
Let $\calC$ be an $\infty$-category which admits small, $\kappa$-filtered colimits. We will say a functor $f: \calC \rightarrow \calD$ is {\it $\kappa$-continuous} if it preserves $\kappa$-filtered colimits.

Let $\calC$ be an $\infty$-category containing an object $C$, and let
$j_{C}: \calC \rightarrow \hat{\SSet}$ denote the functor corepresented by $C$.
If $\calC$ admits $\kappa$-filtered colimits, then we will say that $C$ is {\it $\kappa$-compact}
if $j_{C}$ is $\kappa$-continuous. We will say that $C$ is {\it compact} if
it is $\omega$-compact (and $\calC$ admits filtered colimits).\index{gen}{compact object!of an $\infty$-category}\index{gen}{$\kappa$-compact!object}

Let $\kappa$ be a regular cardinal, and let $\calC$ be an $\infty$-category which admits small, $\kappa$-filtered colimits. We will say that a left fibration $\widetilde{\calC} \rightarrow \calC$
is {\it $\kappa$-compact} if it is classified by a $\kappa$-continuous functor
$\calC \rightarrow \hat{\SSet}$.\index{gen}{$\kappa$-compact!left fibration}
\end{definition}

\begin{notation}
Let $\calC$ be an $\infty$-category and $\kappa$ a regular cardinal. We will generally let
$\calC^{\kappa}$ denote the full subcategory spanned by the $\kappa$-compact objects of $\calC$. \index{not}{calCkappa@$\calC^{\kappa}$}
\end{notation}

\begin{lemma}\label{misst}
Let $\calC$ be an $\infty$-category which admits small $\kappa$-filtered colimits,
and let $\calD \subseteq \Fun(\calC, \widehat{\SSet})$ be the full subcategory spanned by the $\kappa$-continuous functors $f: \calC \rightarrow \hat{\SSet}$. Then $\calD$ is stable under $\kappa$-small limits in $\hat{\SSet}^{\calC}$.
\end{lemma}

\begin{proof}
Let $K$ be a $\kappa$-small simplicial set, and let $p: K \rightarrow \Fun(\calC, \widehat{\SSet})$ be a diagram, which we may identify with a map $p': \calC \rightarrow \Fun(K, \widehat{\SSet})$. Using Proposition \ref{limiteval}, we may obtain a limit of the diagram $p$ by composing $p'$ with a limit functor
$$\varprojlim: \Fun(K,\widehat{\SSet}) \rightarrow \widehat{\SSet}$$
(that is, a right adjoint to the diagonal functor $\widehat{\SSet} \rightarrow \Fun(K,\widehat{\SSet})$; see \S \ref{fcolm}). It therefore suffices to show that the functor $\varprojlim$ is $\kappa$-continuous. This is simply a reformulation of Proposition \ref{frent}.
\end{proof}

The basic properties of $\kappa$-compact left fibrations are summarized in the following Lemma::

\begin{lemma}\label{hardstuff0}
Let $\kappa$ be a regular cardinal.

\begin{itemize}
\item[$(1)$] Let $\calC$ be an $\infty$-category which admits small, $\kappa$-filtered colimits, and let $C \in \calC$ be an object. Then $C$ is $\kappa$-compact if and only if the left fibration $\calC_{C/} \rightarrow \calC$ is $\kappa$-compact.

\item[$(2)$] Let $f: \calC \rightarrow \calD$ be a $\kappa$-continuous functor between
$\infty$-categories which admit small, $\kappa$-filtered colimits, and let
$\widetilde{\calD} \rightarrow \calD$ be a $\kappa$-compact left fibration. Then
the associated left fibration $\calC \times_{\calD} \widetilde{\calD} \rightarrow \calC$
is also $\kappa$-compact.

\item[$(3)$] Let $\calC$ be an $\infty$-category which admits small, $\kappa$-filtered colimits, and let $\bfA \subseteq (\sSet)_{/\calC}$ denote the full subcategory spanned by the
$\kappa$-compact left fibrations over $\calC$. Then $\bfA$ is stable under
$\kappa$-small homotopy limits (with respect to the covariant model structure on $(\sSet)_{/\calC}$. 
In particular, $\bfA$ is stable under the formation homotopy pullbacks, $\kappa$-small products, and $($ if $\kappa$ is uncountable $)$ the formation of homotopy inverse limits of towers.
\end{itemize}
\end{lemma}

\begin{proof}
Assertions $(1)$ and $(2)$ are obvious. To prove $(3)$, let us suppose that
$\widetilde{\calC}$ is a $\kappa$-small homotopy limit of $\kappa$-compact
left fibrations $\widetilde{\calC}_{\alpha} \rightarrow \calC$. Let
$\calJ$ be a small, $\kappa$-filtered $\infty$-category, and 
$\overline{p}: \calJ^{\triangleright} \rightarrow \calC$ a colimit diagram.
We wish to prove that the composition of $\overline{p}$ with the functor
$\calC \rightarrow \hat{\SSet}$ classifying $\widetilde{\calC}$ is a colimit diagram.
Applying Proposition \ref{rot}, we may reduce to the case where $\calJ$ is the nerve
of a $\kappa$-filtered partially ordered set $A$. According to Theorem \ref{struns},
it will suffice to show that the collection of homotopy colimit diagrams
$$ A \cup \{\infty \} \rightarrow \Kan$$
is stable under $\kappa$-small homotopy limits in the diagram category $(\sSet)^{A \cup \{\infty\} }$, which follows easily from our assumption that $A$ is $\kappa$-filtered.
\end{proof}

Our next goal is to prove a very useful stability result for $\kappa$-compact objects (Proposition \ref{placeabovee}). We first need to establish a few technical lemmas.

\begin{lemma}\label{hardstuff1}
Let $\kappa$ be a regular cardinal, let $\calC$ be an $\infty$-category which admits small, $\kappa$-filtered colimits, and let $f: C \rightarrow D$ be a morphism in $\calC$. Suppose that $C$ and $D$ are $\kappa$-compact objects of $\calC$. Then $f$ is a $\kappa$-compact object of
$\Fun(\Delta^1, \calC)$.
\end{lemma}

\begin{proof}
Let $X = \Fun(\Delta^1, \calC) \times_{ \Fun(\{1\}, \calC) } \calC_{f/}$,
$Y = \Fun(\Delta^1, \calC_{C/})$, and $Z = \Fun(\Delta^1, \calC) \times_{ \Fun(\{1\}, \calC) } \calC_{C/},$
so that we have a (homotopy) pullback diagram
$$ \xymatrix{ \Fun(\Delta^1, \calC)_{f/} \ar[r] \ar[d] & X \ar[d] \\
Y \ar[r] & Z }$$
of left fibrations over $\Fun(\Delta^1, \calC)$.
According to Lemma \ref{hardstuff0}, it will suffice to show that $X$, $Y$, and $Z$
are $\kappa$-compact left fibrations. To show that $X$ is a $\kappa$-compact left fibration, it suffices to show that $\calC_{f/} \rightarrow \calC$ is a $\kappa$-compact left fibration, which
follows since we have a trivial fibration $\calC_{f/} \rightarrow \calC_{D/}$, where $D$ is $\kappa$-compact by assumption. Similarly, we have a trivial fibration
$Y \rightarrow \Fun(\Delta^1, \calC) \times_{ \calC^{(0)}} \calC_{C/}$, so that the 
$\kappa$-compactness of $C$ implies that $Y$ is a $\kappa$-compact left fibration. Lemma \ref{hardstuff0} and the compactness of $C$ immediately imply that $Z$ is a $\kappa$-compact left fibration, which completes the proof.
\end{proof}

\begin{lemma}\label{hardstuff2}
Let $\kappa$ be a regular cardinal, and let $\{ \calC_{\alpha} \}$ be a $\kappa$-small family of $\infty$-categories having product $\calC$. Suppose that each $\calC$ admits small, $\kappa$-filtered colimits. Then:
\begin{itemize}
\item[$(1)$] The $\infty$-category $\calC$ admits $\kappa$-filtered colimits.

\item[$(2)$] If $C \in \calC$ is an object whose image in each $\calC_{\alpha}$ is $\kappa$-compact, then $C$ is $\kappa$-compact as an object of $\calC$.
\end{itemize}
\end{lemma}

\begin{proof}
The first assertion is obvious, since colimits in a product can be computed pointwise. For the second, choose an object $C \in \calC$ whose images $\{ C_{\alpha} \in \calC_{\alpha} \}$ are $\kappa$-compact. 

The left fibration $\calC_{C/} \rightarrow \calC$ can obtained as a $\kappa$-small product of the left fibrations
$\calC \times_{\calC_{\alpha}} (\calC_{\alpha})_{C_{\alpha}/} \rightarrow \calC$. Lemma \ref{hardstuff0} implies that each factor is $\kappa$-compact, so that the product is also $\kappa$-compact.
\end{proof}

\begin{lemma}\label{showtop}
Let $S$ be a simplicial set, and suppose given a tower
$$ \ldots X(1) \stackrel{f_1}{\rightarrow} X(0) \stackrel{f_0}{\rightarrow} S$$
where each $f_i$ is a left fibration. Then the inverse limit
$X(\infty)$ is a homotopy inverse limit of the tower $\{ X(i) \}$ with respect to the covariant model structure on $(\sSet)_{/S}$.
\end{lemma}

\begin{proof}
Construct a ladder
$$ \xymatrix{ \ldots \ar[r] & X(1) \ar[r]^{f_1} \ar[d] & X(0) \ar[r]^{f_0} \ar[d] & S \ar[d] \\
\ldots \ar[r] & X'(1) \ar[r]^{f'_1} & X'(0) \ar[r]^{f'_0} & S }$$
where the vertical maps are covariant equivalences and the tower
$\{ X'(i) \}$ is fibrant, in the sense that each of the maps $f'_i$ is a covariant fibration.
We wish to show that the induced map on inverse limits $X(\infty) \rightarrow X'(\infty)$
is a covariant equivalence. Since both $X(\infty)$ and $X'(\infty)$ are left-fibered over $S$,
this can be tested by passing to the fibers over each vertex $s$ of $S$. We may therefore reduce to the case where $S$ is a point, in which case the tower $\{ X(i) \}$ is already fibrant (since
a left fibration over a Kan complex is a Kan fibration; see Lemma \ref{toothie2}).
\end{proof}

\begin{lemma}\label{hardstuff3}
Let $\kappa$ be an uncountable regular cardinal, and let
$$ \ldots \rightarrow \calC^2 \stackrel{f_2}{\rightarrow} \calC^1 \stackrel{f_1}{\rightarrow} \calC^0$$ be a tower of $\infty$-categories. Suppose that each $\calC^i$ admits small $\kappa$-filtered colimits, and that each of the functors $f_i$ is a categorical fibration which preserves $\kappa$-filtered colimits.
Let $\calC$ denote the inverse limit of the tower. Then:
\begin{itemize}
\item[$(1)$] The $\infty$-category $\calC$ admits small $\kappa$-filtered colimits, and the projections
$p_n: \calC \rightarrow \calC^n$ are $\kappa$-continuous.

\item[$(2)$] If $C \in \calC$ has $\kappa$-compact image in $\calC^i$ for each $i \geq 0$, then
$C$ is a $\kappa$-compact object of $\calC$.
\end{itemize}
\end{lemma}

\begin{proof}
Let $\overline{q}: K^{\triangleright} \rightarrow \calC$ be a diagram indexed by an arbitrary simplicial set, let $q = \overline{q}|K$, and set $\overline{q}_n = p_n \circ \overline{q}$, 
$q_n = p_n \circ q$. Suppose that each $\overline{q}_n$ is a colimit diagram in $\calC^n$.
Then the map $\calC_{\overline{q}/} \rightarrow \calC_{q/}$ is the inverse limit of a tower of trivial fibrations $\calC^n_{\overline{q}_n/} \rightarrow \calC^n_{q_n/}$, and therefore a a trivial fibration.

To complete the proof of $(1)$, it will suffice to show that if $K$ is a $\kappa$-filtered $\infty$-category, then any diagram $q: K \rightarrow \calC$ can be extended to a map
$\overline{q}: K^{\triangleright} \rightarrow \calC$ with the property described above.
To construct $\overline{q}$, it suffices to construct a compatible family $\overline{q}_n:
K^{\triangleright} \rightarrow \calC^n$. We begin by selecting arbitrary colimit diagrams
$\overline{q}'_n: K^{\triangleright} \rightarrow \calC^n$ which extend $q_n$. 
We now explain how to adjust these choices to make them compatible with one another, using induction on $n$. Set $\overline{q}_0 = \overline{q}'_0$. Suppose next that $n > 0$.
Since $f_n$ preserves $\kappa$-filtered colimits, we may identify
$\overline{q}_{n-1}$ and $f_n \circ \overline{q}'_{n}$ with initial objects of
$\calC^{n-1}_{q_{n-1}/}$. It follows that there exists an equivalence
$e: \overline{q}_{n-1} \rightarrow f_n \circ \overline{q}'_{n}$ in $\calC^{n-1}_{q_{n-1}/}$. The map
$f_n$ induces a categorical fibration $\calC^n_{q_n/} \rightarrow \calC^{n-1}_{q_{n-1}/}$, so that
$e$ lifts to an equivalence $\overline{e}: \overline{q}_n \rightarrow \overline{q}'_n$ in
$\calC^n_{q_n/}$. The existence of the equivalence $\overline{e}$ proves that
$\overline{q}_n$ is a colimit diagram in $\calC^n$, and we have
$\overline{q}_{n-1} = f_n \circ \overline{q}_n$ by construction. This proves $(1)$.

Now suppose that $C \in \calC$ is as in $(2)$, and let $C^n = p_n(C) \in \calC^n$. 
The left fibration $\calC_{/C}$ is the inverse limit of a tower of left fibrations
$$ \ldots \rightarrow \calC^1_{C^1/} \times_{\calC^1} \calC \rightarrow \calC^0_{C^0/} \times_{\calC^0} \calC.$$
Using Lemma \ref{hardstuff0}, we deduce that each term in this tower is a $\kappa$-compact left fibration over $\calC$. Proposition \ref{sharpen} implies that each map in the tower is a left fibration, so that $\calC_{C/}$ is a homotopy inverse limit of a tower of $\kappa$-compact left fibrations, by Lemma \ref{showtop}. We now apply Lemma \ref{hardstuff0} again to deduce that
$\calC_{C/}$ is a $\kappa$-compact left fibration, so that $C \in \calC$ is $\kappa$-compact as desired.
\end{proof}

\begin{proposition}\label{placeabovee}\index{gen}{compact object!of a functor $\infty$-category}
Let $\kappa$ be a regular cardinal, let $\calC$ be an $\infty$-category which admits
small $\kappa$-filtered colimits, and let $f: K \rightarrow \calC$ be a diagram indexed by a $\kappa$-small simplicial set $K$. Suppose that for each vertex $x$ of $K$, $f(x) \in \calC$ is $\kappa$-compact. Then $f$ is a $\kappa$-compact object of $\Fun(K,\calC)$.
\end{proposition}

\begin{proof}
Let us say that a simplicial set $K$ is {\em good} if it satisfies the conclusions of the lemma.
We wish to prove that all $\kappa$-small simplicial sets are good. The proof proceeds in several steps:

\begin{itemize}
\item[$(1)$] Given a pushout square
$$ \xymatrix{ K' \ar[r] \ar[d]^{i} & K \ar[d] \\
L' \ar[r] & L }$$
where $i$ is a cofibration and the simplicial sets $K'$, $K$, and $L'$ are good, the simplicial
set $L$ is also good. To prove this, we observe that the associated diagram of $\infty$-categories
$$ \xymatrix{ \Fun(L,\calC) \ar[r] \ar[d] & \Fun(L',\calC) \ar[d] \\
\Fun(K,\calC) \ar[r] & \Fun(K',\calC) }$$ is homotopy Cartesian, and every arrow in the diagram preserves
$\kappa$-filtered colimits (by Proposition \ref{limiteval}). Now apply Lemma \ref{yoris}.

\item[$(2)$] If $K \rightarrow K'$ is a categorical equivalence and $K$ is good, then $K'$ is good:
the forgetful functor $\Fun(K',\calC) \rightarrow \Fun(K,\calC)$ is an equivalence of $\infty$-categories, and therefore detects $\kappa$-compact objects.

\item[$(3)$] Every simplex $\Delta^n$ is good. To prove this, we observe that the inclusion
$$ \Delta^{ \{0,1\} } \coprod_{ \{1\} } \ldots \coprod_{ \{n-1\} } \Delta^{ \{n-1,n\}} \subseteq \Delta^n$$
is a categorical equivalence. Applying $(1)$ and $(2)$, we can reduce to the case $n \leq 1$.
If $n=0$ there is nothing to prove, and if $n = 1$ we apply Lemma \ref{hardstuff1}.

\item[$(4)$] If $\{ K_{\alpha} \}$ is a $\kappa$-small collection of good simplicial sets having coproduct $K$, then $K$ is also good. To prove this, we observe that
$\Fun(\calC) \simeq \prod_{\alpha} \Fun(K_{\alpha}, \calC)$ and apply Lemma \ref{hardstuff2}.

\item[$(5)$] If $K$ is a $\kappa$-small simplicial set of dimension $\leq n$, then $K$ is good. The proof is by induction on $n$. Let
$K^{(n-1)} \subseteq K$ denote the $(n-1)$-skeleton of $K$, so that we have a pushout diagram
$$ \xymatrix{ \coprod_{ \sigma \in K_n} \bd \Delta^n \ar[r] \ar[d] & K^{(n-1)} \ar[d] \\
\coprod_{\sigma \in K_n} \Delta^n \ar[r] & K.}$$
The inductive hypothesis implies that $\coprod_{\sigma \in K_n} \bd \Delta^n$ and $K^{(n-1)}$
are good. Applying $(3)$ and $(4)$, we deduce that $\coprod_{\sigma \in K_n} \Delta^n$ is good. We now apply $(1)$ to deduce that $K$ is good.

\item[$(6)$] Every $\kappa$-small simplicial set $K$ is good. If $\kappa = \omega$, then this
follows immediately from $(5)$, since every $\kappa$-small simplicial set is finite dimensional.
If $\kappa$ is uncountable, then we have an increasing filtration
$$ K^{(0)} \subseteq K^{(1)} \subseteq \ldots $$
which gives rise to a tower of $\infty$-categories
$$ \ldots \Fun(K^{(1)}, \calC) \rightarrow \Fun(K^{(0)}, \calC) $$
having (homotopy) inverse limit $\Fun(K,\calC)$. Using Proposition \ref{limiteval}, we deduce that the hypotheses of Lemma \ref{hardstuff3} are satisfied, so that $K$ is good.
\end{itemize}
\end{proof}

\begin{corollary}\label{jurman}
Let $\kappa$ be a regular cardinal, and let $\calC$ be an $\infty$-category which admits
small, $\kappa$-filtered colimits. Suppose that $p: K \rightarrow \calC$ is a $\kappa$-small diagram with the property that for every vertex $x$ of $K$, $p(x)$ is a $\kappa$-compact object of $\calC$.
Then the left fibration $\calC_{p/} \rightarrow \calC$ is $\kappa$-compact.
\end{corollary}

\begin{proof}
It will suffice to show that the equivalent left fibration $\calC^{p/ } \rightarrow \calC$
is $\kappa$-compact. Let $P$ be the object of $\Fun(K,\calC)$ corresponding to $p$. Then
we have an isomorphism of simplicial sets $$\calC^{p/} \simeq \calC \times_{ \Fun(K,\calC) } \Fun(K,\calC)^{P/}.$$ Proposition \ref{placeabovee} asserts that $P$ is a $\kappa$-compact object
of $\Fun(K,\calC)$, so that the left fibration $$\Fun(K,\calC)^{P/} \rightarrow \Fun(K,\calC)$$ is $\kappa$-compact.
Proposition \ref{limiteval} implies that the diagonal map $\calC \rightarrow \Fun(K,\calC)$ preserves $\kappa$-filtered colimits, so we can apply part $(2)$ of Lemma \ref{hardstuff0} to deduce that
$\calC^{p/} \rightarrow \calC$ is $\kappa$-compact, as well.
\end{proof}

\begin{corollary}\label{tyrmyrr}
Let $\calC$ be an $\infty$-category which admits small, $\kappa$-filtered colimits, and let
$\calC^{\kappa}$ denote the full subcategory of $\calC$ spanned by the $\kappa$-compact objects. Then $\calC^{\kappa}$ is stable under the formation of all $\kappa$-small colimits which exist in $\calC$.
\end{corollary}

\begin{proof}
Let $K$ be a $\kappa$-small simplicial set, and let $\overline{p}: K^{\triangleright} \rightarrow \calC$ be a colimit diagram. Suppose that, for each vertex $x$ of $K$, the object $\overline{p}(x) \in \calC$ is $\kappa$-compact. We wish to show that $C= \overline{p}(\infty) \in \calC$ is $\kappa$-compact, where $\infty$ denotes the cone point of $K^{\triangleright}$. Let $p = \overline{p}|K$, and consider the maps
$$ \calC_{p/} \leftarrow \calC_{ \overline{p}/ } \rightarrow \calC_{ C/}. $$
Both are trivial fibrations (the first because $\overline{p}$ is a colimit diagram, and the second because the inclusion $\{ \infty\} \subseteq K^{\triangleright}$ is right anodyne). Corollary \ref{jurman} asserts that the left fibration $\calC_{p/} \rightarrow \calC$ is $\kappa$-compact. It follows that the equivalent left fibration $\calC_{C/}$ is $\kappa$-compact, so that $C$ is a $\kappa$-compact object of $\calC$ as desired.
\end{proof}

\begin{remark}\label{trmyr}
Let $\kappa$ be a regular cardinal, and let $\calC$ be an $\infty$-category which admits $\kappa$-filtered colimits. Then the full subcategory $\calC^{\kappa} \subseteq \calC$ of $\kappa$-compact objects is stable under retracts. If $\kappa > \omega$, this follows from 
Proposition \ref{slanger} and Corollary \ref{tyrmyrr} (since every retract can be obtained as a $\kappa$-small colimit). We give an alternative argument that works also in the most important case $\kappa = \omega$. Let $C$ be $\kappa$-compact, and let $D$ be a retract of $C$. Let
$j: \calC^{op} \rightarrow \Fun(\calC, \widehat{\SSet})$ be the Yoneda embedding. Then
$j(D) \in \Fun(\calC, \widehat{\SSet})$ is a retract of $j(C)$. Since $j(C)$ preserves $\kappa$-filtered colimits, then Lemma \ref{compcompcomp} implies that $j(D)$ preserves $\kappa$-filtered colimits, so that $D$ is $\kappa$-compact.
\end{remark}

The following result gives a convenient description of the compact objects of an $\infty$-category of presheaves:

\begin{proposition}\label{charsmallpre}
Let $\calC$ be a small $\infty$-category, $\kappa$ a regular cardinal, and 
$C \in \calP(\calC)$ an object. The following are equivalent:
\begin{itemize}
\item[$(1)$] There exists a diagram $p: K \rightarrow \calC$ indexed by a $\kappa$-small
simplicial set, such that $j \circ p$ has a colimit $D$ in $\calP(\calC)$, and $C$ is a retract of $D$.
\item[$(2)$] The object $C$ is $\kappa$-compact.
\end{itemize}
\end{proposition}

\begin{proof}
Proposition \ref{dda} asserts that for every object $A \in \calC$, $j(A)$ is completely compact, and in particular $\kappa$-compact. According to Corollary \ref{tyrmyrr} and Remark \ref{trmyr}, the collection of $\kappa$-compact objects of $\calP(\calC)$ is stable under $\kappa$-small colimits and retracts. Consequently, $(1) \Rightarrow (2)$.

Now suppose that $(2)$ is satisfied. Let $\calC_{/C} = \calC \times_{\calP(\calC)} \calP(\calC)_{/C}$.
Lemma \ref{longwait0} implies that the composition
$$ \overline{p}: \calC_{/C}^{\triangleright} \rightarrow \calP(S)_{/C}^{\triangleright} \rightarrow \calP(S)$$ is a colimit diagram. As in the proof of Corollary \ref{uterrr}, we can write $C$ as the colimit
of a $\kappa$-filtered diagram $q: \calI \rightarrow \calP(\calC)$, where each object
$q(I)$ is the colimit of $\overline{p} | \calC^{0}$, where $\calC^{0}$ is a $\kappa$-small
simplicial subset of $\calC_{/C}$. Since $C$ is $\kappa$-compact, we may argue as in the proof of Proposition \ref{dda} to deduce that $C$ is a retract of $q(I)$, for some object $I \in \calI$.
This proves $(1)$.
\end{proof}

We close with a result which we will need in \S \ref{c5s6}. First, a bit of notation: if $\calC$ is a small $\infty$-category and $\kappa$ a regular cardinal, we let $\calP^{\kappa}(\calC)$ denote the full subcategory consisting of $\kappa$-compact objects of $\calP(\calC)$.

\begin{proposition}\label{kcolim}
Let $\calC$ be a small, idempotent complete $\infty$-category and $\kappa$ a regular cardinal. The following conditions are equivalent:
\begin{itemize}
\item[$(1)$] The $\infty$-category $\calC$ admits $\kappa$-small colimits.
\item[$(2)$] The Yoneda embedding $j: \calC \rightarrow \calP^{\kappa}(\calC)$ has a left adjoint.
\end{itemize}
\end{proposition}

\begin{proof}
Suppose that $(1)$ is satisfied. For each object
$M \in \calP(\calC)$, let $F_M: \calP(\calC) \rightarrow \hat{\SSet}$ denote the associated corepresentable functor. Let $\calD \subseteq \calP(\calC)$
denote the full subcategory of $\calP(\calC)$ spanned by those objects $M$
such that $F_M \circ j: \calC \rightarrow \hat{\SSet}$ is corepresentable.
According to Proposition \ref{limiteval}, composition with $j$ induces a {\em limit-preserving} functor
$$ \Fun(\calP(\calC), \widehat{\SSet}) \rightarrow \Fun(\calC, \widehat{\SSet}).$$
Applying Proposition \ref{yonedaprop} to $\calC^{op}$, we conclude that the collection
of corepresentable functors on $\calC$ is stable under retracts and $\kappa$-small limits.
A second application of Proposition \ref{yonedaprop} (this time to 
$\calP(\calC)^{op}$) now shows that $\calD$ is stable under retracts
and $\kappa$-small colimits in $\calP(\calC)$. Since $j$ is fully faithful, $\calD$ contains the essential image of $j$. It follows from Proposition \ref{charsmallpre} that $\calD$ contains
$\calP^{\kappa}(\calC)$. We now apply Proposition \ref{adjfuncbaby} to deduce that
$j: \calC \rightarrow \calP^{\kappa}(\calC)$ admits a left adjoint.

Conversely, suppose that $(2)$ is satisfied. Let $L$ denote a left adjoint to the Yoneda embedding, let
$p: K \rightarrow \calC$ be a $\kappa$-small diagram, and let
$q = j \circ p$. Using Corollary \ref{tyrmyrr}, we deduce that $q$ has a colimit
$\overline{q}: K^{\triangleright} \rightarrow \calP^{\kappa}(\calC)$. Since $L$ is a left adjoint,
$L \circ \overline{q}$ is a colimit of $L \circ q$. Since $j$ is fully faithful, the diagram
$p$ is equivalent to $L \circ q$, so that $p$ has a colimit as well.
\end{proof}

\subsection{Ind-Objects}\label{indlim}

Let $S$ be a simplicial set. In \S \ref{presheaf4}, we proved that the $\infty$-category
$\calP(S)$ is {\em freely} generated under small colimits by the image of the Yoneda embedding $j: S \rightarrow \calP(S)$ (Theorem \ref{charpresheaf}). Our goal in this section is to study
the an analogous construction, where we allow only {\em filtered} colimits.

\begin{definition}\label{indsmall}\index{not}{IndC@$\Ind(\calC)$}\index{not}{IndkC@$\Ind_{\kappa}(\calC)$}\index{gen}{$\infty$-category!of $\Ind$-objects}
Let $\calC$ be a small $\infty$-category and let $\kappa$ be a regular cardinal. 
We let $\Ind_{\kappa}(\calC)$ denote the full subcategory of $\calP(\calC)$ spanned by those functors $f: \calC^{op} \rightarrow \SSet$ which classify right fibrations
$\widetilde{\calC} \rightarrow \calC$ where the $\infty$-category
$\widetilde{\calC}$ is $\kappa$-filtered. In the case where $\kappa = \omega$, we will simply write $\Ind(\calC)$ for $\Ind_{\kappa}(\calC)$. We will refer to $\Ind(\calC)$ as the $\infty$-category of {\it $\Ind$-objects} of $\calC$.\index{gen}{Ind-object!of an $\infty$-category}
\end{definition}

\begin{remark}\label{repareind}
Let $\calC$ be a small $\infty$-category and $\kappa$ a regular cardinal. Then the Yoneda embedding $j: \calC \rightarrow \calP(\calC)$ factors through $\Ind_{\kappa}(\calC)$. This follows immediately from Lemma \ref{repco}, since $j(C)$ classifies the right fibration $\calC_{/C} \rightarrow \calC$. The $\infty$-category $\calC_{/C}$ has a final object and is therefore $\kappa$-filtered (Proposition \ref{smallity}).
\end{remark}

\begin{proposition}\label{geort}
Let $\calC$ be a small $\infty$-category, and let $\kappa$ be a regular cardinal.
The full subcategory $\Ind_{\kappa}(\calC) \subseteq \calP(\calC)$ is stable under $\kappa$-filtered colimits.
\end{proposition}

\begin{proof}
Let $\calP'_{\Delta}(\calC)$ denote the full subcategory of $(\sSet)_{/\calC}$ spanned
by the right fibrations $\widetilde{\calC} \rightarrow \calC$. According to Proposition \ref{othermod},
the $\infty$-category $\calP(\calC)$ is equivalent to the simplicial nerve $\sNerve (\calP'_{\Delta}(\calC))$. Let $\Ind'_{\kappa}(\calC)$ denote the full subcategory of
$\calP'_{\Delta}(\calC)$ spanned by right fibrations $\widetilde{\calC} \rightarrow \calC$ where
$\widetilde{\calC}$ is $\kappa$-filtered. It will suffice to prove that for any diagram
$p: \calI \rightarrow \sNerve(\Ind'_{\Delta}(\calC))$ indexed by a small $\kappa$-filtered $\infty$-category $\calI$, the colimit of $p$ in $\sNerve(\calP'_{\Delta}(\calC))$ also belongs to
$\Ind'_{\kappa}(\calC)$. Using Proposition \ref{rot}, we may reduce to the case where $\calI$ is the nerve of a $\kappa$-filtered partially ordered set $A$. Using Proposition \ref{gumby444}, we may further reduce to the case where $p$ is the simplicial nerve of a diagram taking values in the ordinary category $\Ind'_{\kappa}(\calC)$. In virtue of Theorem \ref{colimcomparee}, it will suffice to prove that $\Ind'_{\kappa}(\calC) \subseteq \calP'_{\Delta}(\calC)$ is stable under $\kappa$-filtered homotopy colimits. We may identify $\calP'_{\Delta}$ with the collection of fibrant objects of
$(\sSet)_{/\calC}$ with respect to the contravariant model structure. Since the class of contravariant equivalences is stable under filtered colimits, any $\kappa$-filtered colimit in
$(\sSet)_{/\calC}$ is also a homotopy colimit. Consequently, it will suffice to prove that
$\Ind'_{\kappa}(\calC) \subseteq \calP'_{\Delta}(\calC)$ is stable under $\kappa$-filtered colimits.
This follows immediately from the definition of a $\kappa$-filtered $\infty$-category.
\end{proof}

\begin{corollary}\label{indpr}
Let $\calC$ be a small $\infty$-category, let $\kappa$ be a regular cardinal, and let
$F: \calC^{op} \rightarrow \SSet$ be an object of $\calP(\calC)$. The following conditions are equivalent:

\begin{itemize}
\item[$(1)$] There exists a $($small$)$ $\kappa$-filtered $\infty$-category $\calI$, a diagram
$p: \calI \rightarrow \calC$ such that $F$ is a colimit of the composition
$j \circ p: \calI \rightarrow \calP(\calC)$.

\item[$(2)$] The functor $F$ belongs to $\Ind_{\kappa}(\calC)$.

\end{itemize}

If $\calC$ admits $\kappa$-small colimits, then $(1)$ and $(2)$ are equivalent to
\begin{itemize}

\item[$(3)$] The functor $F$ preserves $\kappa$-small limits.
\end{itemize}
\end{corollary}

\begin{proof}
Lemma \ref{longwait0} implies that $F$ is a colimit of the diagram
$$ \calC_{/F} \rightarrow \calC \stackrel{j}{\rightarrow} \calP(\calC),$$
and Lemma \ref{repco} allows us to identify $\calC_{/F} = \calC \times_{\calP(\calC)} \calP(\calC)_{/F}$ with the right fibration associated to $F$. Thus $(2) \Rightarrow (1)$. The converse follows from Proposition \ref{geort}, since every representable functor belongs to $\Ind_{\kappa}(\calC)$ (Remark \ref{repareind}).

Now suppose that $\calC$ admits $\kappa$-small colimits. If $(3)$ is satisfied, then
$F^{op}: \calC \rightarrow \SSet^{op}$ is $\kappa$-right exact by Proposition \ref{frent}. The right fibration associated to $F$ is the pullback of the universal right fibration by $F^{op}$. Using Corollary \ref{grt}, the universal right fibration over $\SSet^{op}$ is representable by the final object of $\SSet$. Since $F$ is $\kappa$-right exact, the fiber product
$(\SSet^{op})_{/\ast} \times_{\SSet^{op}} \calC$ is $\kappa$-filtered. Thus $(3) \Rightarrow (2)$.

We now complete the proof by showing that $(1) \Rightarrow (3)$. First suppose that
$F$ lies in the essential image of the Yoneda embedding $j: \calC \rightarrow \calP(\calC)$. According to Lemma \ref{repco}, $j(C)$ is equivalent to the composition of the opposite Yoneda embedding $j': \calC^{op} \rightarrow \Fun(\calC, \SSet)$ with the evaluation functor
$e: \Fun(\calC, \SSet) \rightarrow \SSet$ associated to the object $C \in \calC$. Propositions \ref{yonedaprop} and \ref{limiteval} imply that $j'$ and $e$ preserve $\kappa$-small limits, so that $j(C)$ preserves $\kappa$-small limits. To conclude the proof, it will suffice to show that the collection of functors $F: \calC^{op} \rightarrow \SSet$ which satisfy $(3)$ is stable under $\kappa$-filtered colimits: this follows easily from Proposition \ref{frent}.
\end{proof}

\begin{proposition}\label{justcut}
Let $\calC$ be a small $\infty$-category, let $\kappa$ be a regular cardinal, and let
$j: \calC \rightarrow \Ind_{\kappa}(\calC)$ be the Yoneda embedding. For each object $C \in \calC$, $j(C)$ is a $\kappa$-compact object of $\Ind_{\kappa}(\calC)$.
\end{proposition}

\begin{proof}
The functor $\Ind_{\kappa}(\calC) \rightarrow \SSet$ co-represented by $j(C)$ is equivalent to the composition 
$$ \Ind_{\kappa}(\calC) \subseteq \calP(\calC) \rightarrow \SSet$$
where the first map is the canonical inclusion and the second is given by evaluation at
$C$. The second map preserves all colimits (Proposition \ref{limiteval}), and the
first preserves $\kappa$-filtered colimits since $\Ind_{\kappa}(\calC)$ is stable under
$\kappa$-filtered colimits in $\calP(\calC)$ (Proposition \ref{geort}).
\end{proof}

\begin{remark}
Let $\calC$ be a small $\infty$-category and $\kappa$ a regular cardinal. Suppose
that $\calC$ is equivalent to an $n$-category, so that the Yoneda embedding
$j: \calC \rightarrow \calP(\calC)$ factors through $\calP_{\leq n-1}(\calC) = \Fun(\calC^{op}, \tau_{\leq n-1} \SSet)$, where
$\tau_{\leq n-1} \SSet$ denotes the full subcategory of $\SSet$ spanned by the
$(n-1)$-truncated spaces: that is, spaces whose homotopy groups vanish in dimensions $n$ and above. The class of $(n-1)$-truncated spaces is stable under filtered colimits, so that
$\calP_{\leq n-1}(\calC)$ is stable under filtered colimits in $\calP(\calC)$. Corollary \ref{indpr} implies that $\Ind(\calC) \subseteq \calP_{\leq n-1}(\calC)$. In particular, $\Ind(\calC)$ is itself equivalent to an $n$-category. In particular, if $\calC$ is the nerve of an ordinary category $\calI$, then $\Ind(\calC)$ is equivalent to the nerve of an ordinary category $\calJ$, which is uniquely determined up to equivalence. Moreover, $\calJ$ admits filtered colimits, and there is a fully faithful embedding $\calI \rightarrow \calJ$ which generates $\calJ$ under filtered colimits, whose essential image consists of compact objects of $\calJ$. It follows
that $\calJ$ is equivalent to the category of $\Ind$-objects of $\calI$, in the sense of ordinary category theory.
\end{remark}

According to Corollary \ref{indpr}, we may characterize $\Ind_{\kappa}(\calC)$ as the smallest full subcategory of $\calP(\calC)$ which contains the image of the Yoneda embedding $j: \calC \rightarrow \calP(\calC)$ and is stable under $\kappa$-filtered colimits. Our goal is to obtain a more precise characterization of $\Ind_{\kappa}(\calC)$: namely, we will show that it is {\em freely} generated by $\calC$ under $\kappa$-filtered colimits. 

\begin{lemma}\label{diverti}
Let $\calD$ be an $\infty$-category $($not necessarily small$)$. There exists a fully faithful functor
$i: \calD \rightarrow \calD'$ with the following properties:
\begin{itemize}
\item[$(1)$] The $\infty$-category $\calD'$ admits small colimits.
\item[$(2)$] A small diagram $K^{\triangleright} \rightarrow \calD$ is a colimit if and only if the composite map $K^{\triangleright} \rightarrow \calD'$ is a colimit.
\end{itemize}
\end{lemma}

\begin{proof}
Let $\calD' = \Fun( \calD, \widehat{\SSet})^{op}$, and let $i$
be the opposite of the Yoneda embedding. Then $(1)$ follows from Proposition \ref{limiteval} and $(2)$ from Proposition \ref{yonedaprop}.
\end{proof}

We will need the following analogue of Lemma \ref{longwait1}:

\begin{lemma}\label{waitlong1}
Let $\calC$ be a small $\infty$-category, $\kappa$ a regular cardinal,  
$j: \calC \rightarrow \Ind_{\kappa}(\calC)$ the Yoneda embedding, and $\calC' \subseteq \calC$ the essential image of $j$. Let $\calD$ be an $\infty$-category which admits small $\kappa$-filtered colimits. Then:
\begin{itemize}
\item[$(1)$] Every functor $f_0: \calC' \rightarrow \calD$ admits a left Kan extension
$f: \Ind_{\kappa}(\calC) \rightarrow \calD$.

\item[$(2)$] An arbitrary functor $f: \Ind_{\kappa}(\calC) \rightarrow \calD$ is a left Kan extension of $f| \calC'$ if and only if $f$ is $\kappa$-continuous.

\end{itemize}
\end{lemma}

\begin{proof}
Fix an arbitrary functor $f_0: \calC' \rightarrow \calD$. Without loss of generality, we may assume that $\calD$ is a full subcategory of a larger $\infty$-category $\calD'$, satisfying the conclusions of Lemma \ref{diverti}; in particular, $\calD$ is stable under small $\kappa$-filtered colimits in $\calD'$. We may further assume that $\calD$ coincides with its essential image in $\calD'$. Lemma \ref{longwait1} guarantees the existence of a functor $F: \calP(\calC) \rightarrow \calD'$ which is a left Kan extension of $f_0 = F | \calC'$, and such that $F$ preserves small colimits. 
Since $\Ind_{\kappa}(\calC)$ is generated by $\calC'$ under $\kappa$-filtered colimits (Corollary \ref{indpr}), the restriction $f = F | \Ind_{\kappa}(\calC)$ factors through $\calD$. It is then clear that $f: \Ind_{\kappa}(\calC) \rightarrow \calD$ is a left Kan extension of $f_0$, and that $f$ is $\kappa$-continuous. This proves $(1)$ and the ``only if'' direction of $(2)$ (since left Kan extensions of $f_0$ are unique up to equivalence).

We now prove the ``if'' direction of $(2)$. Let $f: \Ind_{\kappa}(\calC) \rightarrow \calD$ be the functor constructed above, and let $f': \Ind_{\kappa}(\calC) \rightarrow \calD$ be an arbitrary
$\kappa$-continuous functor such that $f | \calC'  = f' | \calC'$. We wish to prove that $f'$ is a left Kan extension of $f' | \calC'$. Since $f$ is a left Kan extension of
$f| \calC'$, there exists a natural transformation $\alpha: f \rightarrow f'$ which is an equivalence when restricted to $\calC'$. Let $\calE \subseteq \Ind_{\kappa}(\calC)$ be the full subcategory spanned by those objects $C$ for which the morphism $\alpha_{C}: f(C) \rightarrow f'(C)$ is an equivalence in $\calD$. By hypothesis, $\calC' \subseteq \calE$. Since both $f$ and $f'$ are $\kappa$-continuous, $\calE$ is stable under $\kappa$-filtered colimits in $\Ind_{\kappa}(\calC)$. We now apply Corollary \ref{indpr} to conclude that $\calE = \Ind_{\kappa}(\calC)$. It follows that $f'$ and $f$ are equivalent, so that $f'$ is a left Kan extension of $f' | \calC'$ as desired.
\end{proof}

\begin{remark}\label{poweryoga}
The proof of Lemma \ref{waitlong1} is very robust, and can be used to establish a number of analogous results. Roughly speaking, given any class $S$ of colimits, one can consider the smallest full subcategory $\calC''$ of $\calP(\calC)$ which contains the essential image $\calC'$ of the Yoneda embedding and is stable under colimits of type $S$. Given any functor $f_0: \calC' \rightarrow \calD$, where $\calD$ is an $\infty$-category which admits colimits of type $S$, one can show that there exists a functor $f: \calC'' \rightarrow \calD$ which is a left Kan extension of $f_0 = f| \calC'$, and that $f$ is characterized by the fact that it preserves colimits of type $S$. Taking $S$ to be the class of all small colimits, we recover Lemma \ref{longwait1}. Taking $S$ to be the class of all small $\kappa$-filtered colimits, we recover Lemma \ref{waitlong1}. Other variations are possible as well: we will exploit this idea further in \S \ref{agileco}.
\end{remark}

\begin{proposition}\label{intprop}
Let $\calC$ and $\calD$ be $\infty$-categories, and let $\kappa$ be a regular cardinal.
Suppose that $\calC$ is small and that $\calD$ admits small $\kappa$-filtered colimits.
Then composition with the Yoneda embedding induces an equivalence of $\infty$-categories
$$ \bHom_{\kappa}( \Ind_{\kappa}(\calC), \calD) \rightarrow \Fun(\calC, \calD),$$
where the left hand side denotes the $\infty$-category of all $\kappa$-continuous functors
from $\Ind_{\kappa}(\calC)$ to $\calD$.
\end{proposition}

\begin{proof}
Combine Lemma \ref{waitlong1} with Corollary \ref{leftkanextdef}.
\end{proof}

In other words, if $\calC$ is small and $\calD$ admits $\kappa$-filtered colimits, then any
functor $f: \calC \rightarrow \calD$ determines an essentially unique extension
$F: \Ind_{\kappa}(\calC) \rightarrow \calD$ (such that $f$ is equivalent to $F \circ j$).
We next give a criterion which will allow us to determine when $F$ is an equivalence.

\begin{proposition}\label{uterr}\index{gen}{Ind-object!characterization of}
Let $\calC$ be a small $\infty$-category, $\kappa$ a regular cardinal, and $\calD$ an $\infty$-category which admits $\kappa$-filtered colimits. Let $F: \Ind_{\kappa}(\calC) \rightarrow \calD$
be a $\kappa$-continuous functor, and $f = F \circ j$ its composition with the Yoneda embedding $j: \calC \rightarrow \Ind_{\kappa}(\calC)$. Then:
\begin{itemize}
\item[$(1)$] If $f$ is fully faithful and its essential image consists of $\kappa$-compact objects of $\calD$, then $F$ is fully faithful.
\item[$(2)$] The functor $F$ is an equivalence if and only if the following conditions are satisfied:
\begin{itemize}
\item[$(i)$] The functor $f$ is fully faithful.
\item[$(ii)$] The functor $f$ factors through $\calD^{\kappa}$. 
\item[$(iii)$] The objects $\{ f(C) \}_{C \in \calC}$ generate $\calD$ under $\kappa$-filtered colimits.
\end{itemize}
\end{itemize}
\end{proposition}

\begin{proof}
We first prove $(1)$, using argument of Proposition \ref{trumptow}. Let
$C$ and $D$ be objects of $\Ind_{\kappa}(\calC)$. 
We wish to prove that the map
$$ \eta_{C,D}: \bHom_{\calP(\calC)}(C,D) \rightarrow \bHom_{\calD}(F(C), F(D))$$ is an isomorphism in the homotopy category $\calH$. Suppose first that $C$ belongs to the essential image
of $j$. Let $G: \calP(\calC) \rightarrow \SSet$ be a functor co-represented by $C$, and let
$G': \calD \rightarrow \SSet$ be a functor co-represented by $F(C)$. Then we have a natural transformation of functors $G \rightarrow G' \circ F$. Assumption $(2)$ implies that $G'$ preserves small $\kappa$-filtered colimits, so that $G' \circ F$ preserves small $\kappa$-filtered colimits. 
Proposition \ref{justcut} implies that
$G$ preserves small $\kappa$-filtered colimits. It follows that the collection of objects $D \in \Ind_{\kappa}(\calC)$ such
that $\eta_{C,D}$ is an equivalence is stable under small $\kappa$-filtered colimits colimits. If $D$ belongs to the essential image of $j$, then the assumption that $f$ is fully faithful implies that $\eta_{C,D}$ is a homotopy equivalence. Since the image of the Yoneda embedding generates
$\Ind_{\kappa}(\calC)$ under small $\kappa$-filtered colimits, we conclude that $\eta_{C,D}$ is a homotopy equivalence for every object $D \in \Ind_{\kappa}(\calC)$.

We now drop the assumption that $C$ lies in the essential image of $j$. Fix $D \in \Ind_{\kappa}(\calC)$. Let $H: \Ind_{\kappa}(\calC)^{op} \rightarrow \SSet$
be a functor represented by $D$, and let $H': \calD^{op} \rightarrow \SSet$ be a functor represented by $FD$. Then we have a natural transformation of functors $H \rightarrow H' \circ F^{op}$, which we wish to prove is an equivalence. By assumption, $F^{op}$ preserves small $\kappa$-filtered limits. Proposition \ref{yonedaprop} implies that $H$ and $H'$ preserve small limits. It follows that the collection $P$ of objects $C \in \calP(S)$ such that $\eta_{C,D}$ is an equivalence is stable under small $\kappa$-filtered colimits.
The special case above established that $P$ contains the essential image of the Yoneda embedding. Since $\Ind_{\kappa}(\calC)$ is generated under small $\kappa$-filtered colimits by the image of the Yoneda embedding, we deduce that $\eta_{C,D}$ is an equivalence in general. This completes the proof of $(1)$.

We now prove $(2)$. Suppose first that $F$ is an equivalence. Then
$(i)$ follows from Proposition \ref{fulfaith}, $(ii)$ from Proposition \ref{justcut}, and $(iii)$ from Corollary \ref{indpr}. Conversely, suppose that $(i)$, $(ii)$, and $(iii)$ are satisfied. Using $(1)$, we deduce that $F$ is fully faithful. The essential image of $F$ contains the essential image of $f$ and is stable under small $\kappa$-filtered colimits. Therefore $F$ is essentially surjective, so that $F$ is an equivalence as desired.
\end{proof}

According to Corollary \ref{uterrr}, an $\infty$-category $\calC$ admits small colimits if and only if $\calC$ admits $\kappa$-small colimits and $\kappa$-filtered colimits. Using Proposition \ref{uterr}, we can make a much more precise statement:

\begin{proposition}\label{precst}
Let $\calC$ be a small $\infty$-category and $\kappa$ a regular cardinal. The $\infty$-category
$\calP^{\kappa}(\calC)$ of $\kappa$-compact objects of $\calP(\calC)$ is essentially small: that is, there exists a small $\infty$-category $\calD$ and an equivalence $i: \calD \rightarrow
\calP^{\kappa}(\calC)$. Let $F: \Ind_{\kappa}(\calD) \rightarrow \calP(\calC)$
be a $\kappa$-continuous functor such that the composition of $f$ with the Yoneda embedding
$$ \calD \rightarrow \Ind_{\kappa}(\calD) \rightarrow \calP(\calC)$$
is equivalent to $i$ $($according to Proposition \ref{intprop}, $F$ exists and is unique up
to equivalence$)$. Then $F$ is an equivalence of $\infty$-categories.
\end{proposition}

\begin{proof}
Since $\calP(\calC)$ is locally small, to prove that $\calP^{\kappa}(\calC)$ is small it will suffice to show that the collection of isomorphism classes of objects in the homotopy category
$h \calP^{\kappa}(\calC)$ is small. For this, we invoke Proposition \ref{charsmallpre}:
every $\kappa$-compact object $X$ of $\calP(\calC)$ is a retract of some object $Y$, which is
itself the colimit of some composition
$$ K \stackrel{p}{\rightarrow} \calC \rightarrow \calP(\calC)$$
where $K$ is $\kappa$-small. Since there is a bounded collection of possibilities
for $K$ and $p$ (up to isomorphism in $\sSet$), and a bounded collection of idempotent maps $Y \rightarrow Y$ in $h \calP(\calC)$, there are only a bounded number of possibilities for $X$.

To prove that $F$ is an equivalence, it will suffice to show that $F$ satisfies conditions $(i)$, $(ii)$, and $(iii)$ of Proposition \ref{uterr}. 
Conditions $(i)$ and $(ii)$ are obvious. For $(iii)$, we must prove that every object of
$X \in \calP(\calC)$ can be obtained as a small $\kappa$-filtered colimit of $\kappa$-compact objects of $\calC$. Using Lemma \ref{longwait0}, we can write $X$ as a small colimit
taking values in the essential image of $j: \calC \rightarrow \calP(\calC)$. The proof of Corollary \ref{uterrr} shows that $X$ can be written as a $\kappa$-filtered colimit of a diagram with values
in a full subcategory $\calE \subseteq \calP(\calC)$, where each object of $\calE$ is itself
a $\kappa$-small colimit of some diagram taking values in the essential image of $j$. Using Corollary \ref{tyrmyrr}, we deduce that $\calE \subseteq \calP^{\kappa}(\calC)$, so that $X$ lies in the essential image of $F$ as desired.
\end{proof}

Note that the construction $\calC \mapsto \Ind_{\kappa}(\calC)$ is functorial in $\calC$.
Given a functor $f: \calC \rightarrow \calC'$, Proposition \ref{intprop} implies that the composition of $f$ with the Yoneda embedding $j_{\calC'}: \calC' \rightarrow \Ind_{\kappa} \calC'$ is equivalent to the composition $$ \calC \stackrel{j_{\calC}}{\rightarrow} \Ind_{\kappa} \calC \stackrel{F}{\rightarrow}\Ind_{\kappa} \calC',$$
where $F$ is a $\kappa$-continuous functor. The functor $F$ is well-defined up to equivalence (in fact, up to contractible ambiguity). We will denote $F$ by $\Ind_{\kappa} f$ (though this is perhaps a slight abuse of notation, since $F$ is uniquely determined only up to equivalence).

\begin{proposition}\index{gen}{adjoint functor!between $\Ind$-categories}
Let $f: \calC \rightarrow \calC'$ be a functor between small $\infty$-categories. The following
are equivalent:
\begin{itemize}
\item[$(1)$] The functor $f$ is $\kappa$-right exact.
\item[$(2)$] The map $G: \calP(\calC') \rightarrow \calP(\calC)$ given by composition with $f$ restricts to a functor $g: \Ind_{\kappa}(\calC') \rightarrow \Ind_{\kappa}(\calC)$.
\item[$(3)$] The functor $\Ind_{\kappa} f$ has a right adjoint.
\end{itemize}
Moreover, if these conditions are satisfied, then $g$ is a right adjoint to $\Ind_{\kappa} f$.
\end{proposition}

\begin{proof}
The equivalence $(1) \Leftrightarrow (2)$ is just a reformulation of the definition of $\kappa$-right exactness. Let $\calP(f): \calP(\calC) \rightarrow \calP(\calC')$ be a functor which preserves small colimits such that the diagram of $\infty$-categories
$$ \xymatrix{ \calC \ar[d] \ar[r]^{f} & \calC' \ar[d] \\
\calP(\calC) \ar[r]^{\calP(f)} & \calP(\calC') }$$
is homotopy commutative. Then we may identify $\Ind_{\kappa}(f)$ with the restriction
$\calP(f) | \Ind_{\kappa}(\calC)$. Proposition \ref{adjobs} asserts that $G$ is a right adjoint of
$\calP(f)$. Consequently, if $(2)$ is satisfied, then $g$ is a right adjoint to $\Ind_{\kappa}(f)$. We deduce in particular that $(2) \Rightarrow (3)$. We will complete the proof by showing that $(3)$ implies $(2)$. Suppose that $\Ind_{\kappa}(f)$ admits a right adjoint $g': \Ind_{\kappa}(\calC') \rightarrow \Ind_{\kappa}(\calC)$. Let $X: (\calC')^{op} \rightarrow \SSet$ be an object
of $\Ind_{\kappa}(\calC')$. Then $X^{op}$ is equivalent to the composition
$$ \calC' \stackrel{j}{\rightarrow} \Ind_{\kappa}(\calC') \stackrel{c_X}{\rightarrow} \SSet^{op},$$
where the $c_{X}$ denotes the functor represented by $X$. Since $g'$ is a left adjoint to
$\Ind_{\kappa} f$, the functor $c_{X} \circ \Ind_{\kappa}(f)$ is represented by $g' X$. 
Consequently, we have a homotopy commutative diagram
$$ \xymatrix{ \calC \ar[r]^{j_{\calC}} \ar[d]^{f} & \Ind_{\kappa}(\calC) \ar[r] \ar[d]^{\Ind_{\kappa}(f)} \ar[r]^{c_{g'X}} & \SSet^{op} \ar[d] \\
\calC' \ar[r] & \Ind_{\kappa}(\calC') \ar[r]^{c_{X}} & \SSet^{op} }$$
so that $G(X)^{op} = f \circ X^{op} \simeq c_{g'X} \circ j_{\calC}$, and therefore belongs to
$\Ind_{\kappa}(\calC)$.
\end{proof}

\begin{proposition}\label{turnke}
Let $\calC$ be a small $\infty$-category and $\kappa$ a regular cardinal. 
The Yoneda embedding $j: \calC \rightarrow \Ind_{\kappa}(\calC)$ preserves all
$\kappa$-small colimits which exist in $\calC$.
\end{proposition}

\begin{proof}
Let $K$ be a $\kappa$-small simplicial set, and $\overline{p}: K^{\triangleright} \rightarrow \calC$ a colimit diagram. We wish to show that $j \circ \overline{p}: K^{\triangleright} \rightarrow \Ind_{\kappa}(\calC)$ is also a colimit diagram. Let $C \in \Ind_{\kappa}(\calC)$ be an object, and let
$F: \Ind_{\kappa}(\calC)^{op} \rightarrow \hat{\SSet}$ be the functor represented by $F$. According to Proposition \ref{yonedaprop}, it will suffice to show that 
$F \circ (j \circ \overline{p})^{op}$ is a limit diagram in $\SSet$. We observe that
$F \circ j^{op}$ is equivalent to the object $C \in \Ind_{\kappa}(\calC) \subseteq
\Fun( \calC^{op}, \SSet)$, and therefore $\kappa$-right exact. We now conclude by invoking Proposition \ref{swarmmy}. 
\end{proof}

We conclude this section with a useful result concerning diagrams in $\infty$-categories of $\Ind$-objects:

\begin{proposition}\label{urgh1}
Let $\calC$ be a small $\infty$-category, $\kappa$ a regular cardinal, and 
$j: \calC \rightarrow \Ind_{\kappa}(\calC)$ the Yoneda embedding. Let $A$ be a finite partially ordered set, and let $j': \Fun( \Nerve(A), \calC) \rightarrow \Fun( \Nerve(A), \Ind_{\kappa}(\calC))$ be the induced map. Then $j'$ induces an equivalence
$$ \Ind_{\kappa}( \Fun(\Nerve(A), \calC) ) \rightarrow \Fun( \Nerve(A), \Ind_{\kappa}(\calC)).$$ 
\end{proposition}

In other words, every diagram $\Nerve(A) \rightarrow \Ind_{\kappa}(\calC)$ can be obtained, in
an essentially unique way, as a $\kappa$-filtered colimit of diagrams $\Nerve(A) \rightarrow \calC$.

\begin{warning}
The statement of Proposition \ref{urgh1} fails if we replace $\Nerve(A)$ by an arbitrary finite simplicial set. For example, we may identify the category of abelian groups with the category of $\Ind$-objects of the category of finitely generated abelian groups. If $n > 1$, then the map $q \mapsto \frac{q}{n}$ from the group of rational numbers $\Q$ to itself cannot be obtained as a filtered colimit of endomorphisms finitely generated abelian groups. 
\end{warning}

\begin{proof}[Proof of Proposition \ref{urgh1}]
According to Proposition \ref{uterr}, it will suffice to prove the following:
\begin{itemize}
\item[$(i)$] The functor $j'$ is fully faithful.
\item[$(ii)$] The essential image of $j'$ consists of $\kappa$-compact objects of
$\Fun( \Nerve(A), \Ind_{\kappa}(\calC))$. 
\item[$(iii)$] The essential image of $j'$ generates $\Fun( \Nerve(A), \Ind_{\kappa}(\calC) )$ under small, $\kappa$-filtered colimits. 
\end{itemize}
Since the Yoneda embedding $j: \calC \rightarrow \Ind_{\kappa}(\calC)$ satisfies the analogues of these conditions, $(i)$ is obvious and $(ii)$ follows from 
Proposition \ref{placeabovee}. To prove $(iii)$, we fix an object $F \in \Fun( \Nerve(A),  \Ind_{\kappa}(\calC))$. Let $\calC'$ denote the essential image of $j$, and form a pullback diagram of simplicial sets
$$\xymatrix{ \calD \ar[r] \ar[d] &  \Fun(\Nerve(A), \calC') \ar[d] \\ 
\Fun(\Nerve(A), \Ind_{\kappa}(\calC))_{/F} \ar[r] & \Fun( \Nerve(A), \Ind_{\kappa}(\calC) )}.$$
Since $\calD$ is essentially small, $(iii)$ is a consequence of the following assertions:
\begin{itemize}
\item[$(a)$] The $\infty$-category $\calD$ is $\kappa$-filtered.
\item[$(b)$] The canonical map $\calD^{\triangleright} \rightarrow \Fun( \Nerve(A), \calC)$
is a colimit diagram.
\end{itemize}
To prove $(a)$, it will to show that $\calD$ has the right lifting property with respect
to the inclusion $\Nerve(B) \subseteq \Nerve( B \cup \{\infty\})$, for every $\kappa$-small partially ordered set $B$ (Remark \ref{tweeny}). Regard $B \cup \{ \infty, \infty' \}$ as a partially ordered set with $b < \infty < \infty'$ for each $b \in B$. Unwinding the definitions, we see that
$(a)$ is equivalent to the following assertion:
\begin{itemize}
\item[$(a')$] Let $\overline{F}: \Nerve( A \times ( B \cup \{ \infty' \}) ) \rightarrow \Ind_{\kappa}(\calC)$ 
be such that $\overline{F}| \Nerve(A \times \{ \infty' \}) = F$, and $\overline{F}'| \Nerve(A \times B)$ factors through
$\calC'$. Then there exists a map $\overline{F}': \Nerve(A \times ( B \cup \{ \infty, \infty' \}) ) \rightarrow \Ind_{\kappa}(\calC)$ which extends $\overline{F}$, such that $\overline{F}' | \Nerve(A \times ( B \cup \{\infty\}))$ factors through $\calC'$.
\end{itemize}
To find $\overline{F}'$, we write $A = \{ a_1, \ldots, a_n \}$, where $a_i \leq a_{j}$ implies
$i \leq j$. We will construct a compatible sequence of maps
$$ \overline{F}_{k}: \Nerve( (A \times (B \cup \{ \infty' \})) \cup ( \{ a_1, \ldots, a_k \} \times \{\infty \} )) \rightarrow \calC$$
with $\overline{F}_0 = \overline{F}$ and $\overline{F}_{n} = \overline{F}'$. For each
$a \in A$, we let $ A_{\leq a} = \{ a' \in A: a' \leq a \}$, and we define $A_{< a}$, $A_{\geq a}$, $A_{> a}$ similarly. Supposing that
$\overline{F}_{k-1}$ has been constructed, we observe that constructing $\overline{F}_{k}$
amounts to constructing an object of the $\infty$-category
$$ ( \calC'_{ / F | \Nerve( A_{\geq a_k} )} )_{ \overline{F}_{k-1} | M /},$$
where $M = (A_{\leq a_k} \times B) \cup ( A_{ < a_k } \times \{ \infty \} )$.
The inclusion $\{ a_k \} \subseteq \Nerve( A_{ \geq a_k })$ is left anodyne. It will therefore suffice to construct an object in the equivalent $\infty$-category
$ ( \calC'_{/F(a_k) })_{ \overline{F}_{k-1} |M/ }$. Since $M$ is $\kappa$-small, it suffices to show that
the $\infty$-category $\calC'_{/F(a_k)}$ is $\kappa$-filtered. This is simply a reformulation of the fact that $F(a_k) \in \Ind_{\kappa}(\calC)$. 

We now prove $(b)$. It will suffice to show that for each $a \in A$, the composition
$$ \calD^{\triangleright} \rightarrow \Fun( \Nerve(A), \Ind_{\kappa}(\calC) ) \rightarrow
\Ind_{\kappa}(\calC)$$ is a colimit diagram, where the second map is given by evaluation at $a$. Let $\calD(a) = \calC' \times_{ \Ind_{\kappa}(\calC) } \Ind_{\kappa}(\calC)_{/F(a)}$, so that
$\calD(a)$ is $\kappa$-filtered and the associated map $\calD(a)^{\triangleright} \rightarrow \Ind_\kappa(\calC)$ is a colimit diagram. It will therefore suffice to show that the canonical map
$\calD \rightarrow \calD(a)$ is cofinal. According to Theorem \ref{hollowtt}, it will suffice to show that for each object $D \in \calD(a)$, the fiber product
$\calE = \calD \times_{ \calD(a) } \calD(a)_{D/}$ is weakly contractible. In view of Lemma \ref{stull2}, it will suffice to show that $\calE$ is filtered. This can be established by a minor variation of the argument given above.
\end{proof}

\subsection{Adjoining Colimits to $\infty$-Categories}\label{agileco}

Let $\calC$ be a small $\infty$-category. According to Proposition \ref{intprop}, the 
$\infty$-category $\Ind(\calC)$ enjoys the following properties, which characterize it up to equivalence:
\begin{itemize}
\item[$(1)$] There exists a functor $j: \calC \rightarrow \Ind(\calC)$.
\item[$(2)$] The $\infty$-category $\Ind(\calC)$ admits small filtered colimits.
\item[$(3)$] Let $\calD$ be an $\infty$-category which admits small filtered colimits, and let
$\Fun'( \Ind(\calC), \calD)$ be the full subcategory of $\Fun( \Ind(\calC), \calD)$ spanned by those functors which preserve filtered colimits. Then composition with $j$ induces an equivalence
$\Fun'( \Ind(\calC), \calD) \rightarrow \Fun(\calC, \calD)$.
\end{itemize}
We may informally summarize this characterization as follows: the $\infty$-category $\Ind(\calC)$ is obtained from $\calC$ by freely adjoining the colimits of all small, filtered diagrams. In this section, we will study a generalization of this construction, which allows us to freely adjoin to $\calC$
the colimits of {\em any} collection of diagrams.

\begin{notation}\label{sipser}
Let $\calC$ and $\calD$ be $\infty$-categories, and let $\calR$ be a collection of diagrams
$\{ \overline{p}_{\alpha}: K_{\alpha}^{\triangleright} \rightarrow \calC \}$. We let
$\Fun_{\calR}(\calC, \calD)$ denote the full subcategory of
$\Fun(\calC, \calD)$ spanned by those functors which carry each diagram in $\calR$ to a colimit diagram in $\calD$.\index{not}{FunRCD@$\Fun_{\calR}(\calC, \calD)$}

Let $\calK$ be a collection of simplicial sets. We will say that an $\infty$-category
$\calC$ {\it admits $\calK$-indexed colimits} if it admits $K$-indexed colimits for each $K \in \calK$.
If $f: \calC \rightarrow \calD$ is a functor between $\infty$-categories which admit $\calK$-indexed colimits, then we will say that $f$ {\it preserves $\calK$-indexed colimits} if $f$ preserves $K$-indexed colimits, for each $K \in \calK$. We let $\Fun_{\calK}( \calC, \calD)$ denote the full subcategory of
$\Fun( \calC, \calD)$ spanned by those functors which preserves $\calK$-indexed colimits.\index{not}{FunKCD@$\Fun_{\calK}(\calC,\calD)$}
\end{notation}

\begin{proposition}\label{cupper1}
Let $\calK$ be a collection of simplicial sets, $\calC$ an $\infty$-category, and 
$\calR = \{ \overline{p}_{\alpha}: K^{\triangleright}_{\alpha} \rightarrow \calC \}$ a collection
of diagrams in $\calC$. Assume that each $K_{\alpha}$ belongs to $\calK$. Then there exists
a new $\infty$-category $\calP^{\calK}_{\calR}( \calC)$ and a map $j: \calC \rightarrow \calP^{\calK}_{\calR}(\calC)$ with the following properties:
\begin{itemize}
\item[$(1)$] The $\infty$-category $\calP^{\calK}_{\calR}(\calC)$ admits $\calK$-indexed colimits.
\item[$(2)$] For every $\infty$-category $\calD$ which admits $\calK$-indexed colimits, composition with $j$ induces an equivalence of $\infty$-categories
$$ \Fun_{\calK}( \calP^{\calK}_{\calR}(\calC), \calD) \rightarrow \Fun_{\calR}( \calC, \calD).$$
\end{itemize}
Moreover, if every member of $\calR$ is already a colimit diagram in $\calC$, then we have in addition:
\begin{itemize}
\item[$(3)$] The functor $j$ is fully faithful.
\end{itemize}
\end{proposition}

\begin{remark}
In the situation of Proposition \ref{cupper1}, assertion $(2)$ (applied in the case $\calD = \calP^{\calK}_{\calR}(\calC)$) guarantees that $j$ carries each diagram in $\calR$ to a colimit diagram in $\calP^{\calK}_{\calR}(\calC)$. We can informally summarize conditions $(1)$ and $(2)$ as follows:
the $\infty$-category $\calP^{\calK}_{\calR}(\calC)$ is freely generated by $\calC$ under $\calK$-indexed colimits, subject only to the relations that each diagram in $\calR$ determines a colimit diagram in $\calP^{\calK}_{\calR}(\calC)$. It is clear that this property characterizes $\calP^{\calK}_{\calR}(\calC)$ (and the map $j$) up to equivalence.
\end{remark}

\begin{example}
Suppose that $\calK$ is the collection of all {\em small} simplicial sets, that the $\infty$-category $\calC$ is small, and that the set of diagrams $\calR$ is empty. Then the Yoneda embedding $j: \calC \rightarrow \calP(\calC)$ satisfies the conclusions of Proposition \ref{cupper1}. This is precisely the assertion of Theorem
\ref{charpresheaf} (save for assertion $(3)$, which follows from Proposition \ref{fulfaith}). 
This justifies the notation of Proposition \ref{cupper1}; in the general case we can think of $\calP^{\calK}_{\calR}(\calC)$ as a sort of generalized presheaf category $\calC$, and $j$ as an analogue of the Yoneda embedding.
\end{example}

\begin{proof}[Proof of Proposition \ref{cupper1}:]
We will employ essentially the same argument as in our proof of Proposition \ref{intprop}. 
First, we may enlarge the universe if necessary to reduce to the case where every element of
$\calK$ is a small simplicial set, the $\infty$-category $\calC$ is small, and the collection of diagrams $\calR$ is small. Let $j_0: \calC \rightarrow \calP(\calC)$ denote the Yoneda embedding. For
every diagram $\overline{p}_{\alpha}: K^{\triangleright} \rightarrow \calC$, we let
$p_{\alpha}$ denote the restriction $\overline{p}_{\alpha} | K$, $X_{\alpha} \in \calP(\calC)$ a colimit for
the induced diagram $j \circ p_{\alpha}: K \rightarrow \calP(\calC)$, and
$Y_{\alpha} \in \calC$ the image of the cone point under $\overline{p}_{\alpha}$. The diagram
$j_0 \circ \overline{p}_{\alpha}$ induces a map $s_{\alpha}: X_{\alpha} \rightarrow j_0(Y_{\alpha})$ (well-defined up to homotopy), let $S = \{ s_{\alpha} \}$ be the set of all such morphisms. We let
$S^{-1} \calP(\calC) \subseteq \calP(\calC)$ denote the $\infty$-category of $S$-local objects
of $\calP(\calC)$, and $L: \calP(\calC) \rightarrow S^{-1} \calP(\calC)$ a left adjoint to the inclusion.
We define $\calP^{\calK}_{\calR}(\calC)$ to be the smallest full subcategory of
$S^{-1} \calP(\calC)$ which contains the essential image of the functor $L \circ j_0$ and is closed
under $\calK$-indexed colimits, and let $j = L \circ j_0$ be the induced map. We claim that the
map $j: \calC \rightarrow \calP^{\calK}_{\calR}(\calC)$ has the desired properties.

Assertion $(1)$ is obvious. We now prove $(2)$. Let $\calD$ be an $\infty$-category which admits $\calK$-indexed colimits. In view of Lemma \ref{diverti}, we can assume that there exists a fully faithful inclusion $\calD \subseteq \calD'$, where $\calD'$ admits all small colimits, and $\calD$ is stable under
$\calK$-indexed colimits in $\calD'$. We have a commutative diagram
$$ \xymatrix{ \Fun_{\calK}( \calP^{\calK}_{\calR}(\calC), \calD) \ar[r]^{\phi} \ar[d] & \Fun_{\calR}( \calC, \calD) \ar[d] \\
\Fun_{\calK}( \calP^{\calK}_{\calR}( \calC), \calD') \ar[r]^{\phi'} & \Fun_{\calR}( \calC, \calD'). }$$
We claim that this diagram is a homotopy Cartesian. Unwinding the definitions, this is equivalent to the assertion that a functor $f \in \Fun_{\calK}( \calP^{\calK}_{\calR}(\calC), \calD')$ factors through $\calD$ if and only if $f \circ j: \calC \rightarrow \calD'$ factors through $\calD$. The ``only if'' direction is obvious. Conversely, if $f \circ j$ factors through $\calD$, then $f^{-1} \calD$ is full subcategory of
$\calP^{\calK}_{\calR}(\calC)$ which is stable under $\calK$-indexed limits (since $f$ preserves $\calK$-indexed limits and $\calD$ is stable under $\calK$-indexed limits in $\calD$) and contains the essential image of $j$; by minimality, we conclude that $f^{-1} \calD = \calP^{\calK}_{\calR}(\calC)$. 

Our goal is to prove that the functor $\phi$ is an equivalence of $\infty$-categories. In view of the preceding argument, it will suffice to show that $\phi'$ is an equivalence of $\infty$-categories. In other words, we may replace $\calD$ by $\calD'$ and thereby reduce to the case where
$\calD'$ admit small colimits.

Let $\calE \subseteq \calP(\calC)$ denote the inverse image
$L^{-1} \calP^{\calK}_{\calR}(\calC)$, and let $\overline{S}$ denote the collection of all
morphisms $\alpha$ in $\calE$ such that $L\alpha$ is an equivalence. Composition with
$L$ induces a fully faithful embedding $\Fun( \calP^{\calK}_{\calR}(\calC), \calD) \rightarrow
\Fun( \calE, \calD)$, whose essential image consists of those functors $\calE \rightarrow \calD$
which carry every morphism in $\overline{S}$ to an equivalence in $\calD$. Furthermore, a functor
$f: \calP^{\calK}_{\calR}(\calC) \rightarrow \calD$ preserves $\calK$-indexed colimits if and only if the
composition $f \circ L: \calE \rightarrow \calD$ preserves $\calK$-indexed colimits. The functor $\phi$ factors as a composition
$$ \Fun_{\calK}( \calP^{\calK}_{\calR}(\calC), \calD) \rightarrow \Fun'( \calE, \calD)
\stackrel{\psi}{\rightarrow} \Fun_{\calR}(\calC, \calD),$$
where $\Fun'(\calE, \calD)$ denotes the full subcategory of $\Fun( \calE, \calD)$ spanned by those functors which carry every morphism in $\overline{S}$ to an equivalence and preserve $\calK$-indexed colimits. It will therefore suffice to show that $\psi$ is an equivalence of $\infty$-categories.

In view of Proposition \ref{lklk}, we need only show that if $F: \calE \rightarrow \calD$ is a functor such that $F \circ j_0$ belongs to $\Fun_{\calR}(\calC, \calD)$, then $F$ belongs to
$\Fun'(\calE, \calD)$ if and only if $F$ is a left Kan extension of $F | \calC'$, where
$\calC' \subseteq \calE$ denotes the essential image of the Yoneda embedding $j_0: \calC \rightarrow
\calE$. We first prove the ``if'' direction. Let $F_0 = F | \calC'$. Since $\calD$ admits small colimits, the
functor $F_0$ admits a left Kan extension $\overline{F}: \calP(\calC) \rightarrow \calD$; without loss of generality we may suppose that $F = \overline{F} | \calE$. According to Lemma \ref{longwait1}, 
the functor $\overline{F}$ preserves small colimits. Since $\calE$ is stable under $\calK$-indexed colimits in $\calP(\calC)$, it follows that $\overline{F} | \calE$ preserves $\calK$-indexed colimits.
Furthermore, since $F \circ j_0$ belongs to $\Fun_{\calR}(\calC, \calD)$, the functor
$\overline{F}$ carries each morphism in $S$ to an equivalence in $\calD$. It follows that
$\overline{F}$ factors (up to homotopy) through the localization functor $L$, so that
$\overline{F} | \calE$ carries each morphism in $\overline{S}$ to an equivalence in $\calD$.

For the converse, let us suppose that $F \in \Fun'(\calE, \calD)$; we wish to show that $F$
is a left Kan extension of $F| \calC'$. Let $F'$ denote an arbitrary left Kan extension of
$F | \calC'$, so that the identification $F | \calC' = F' | \calC'$ induces a natural transformation
$\alpha: F' \rightarrow F$. We wish to prove that $\alpha$ is an equivalence. Since $F'$ and
$F$ both carry each morphism in $\overline{S}$ to an equivalence, we may assume without loss of generality that $F = f \circ L$, $F' = f' \circ L$, and $\alpha = \beta \circ L$, where 
$\beta: f' \rightarrow f$ is a natural transformation of functors from $\calP^{\calK}_{\calR}(\calC)$ to
$\calD$. Let $\calX \subseteq \calP^{\calK}_{\calR}(\calC)$ denote the full subcategory spanned by those objects $X$ such that $\beta_{X}: f'(X) \rightarrow f(X)$ is an equivalence. Since both
$f'$ and $f$ preserve $\calK$-indexed colimits, we conclude that $\calX$ is stable under $\calK$-indexed colimits in $\calP^{\calK}_{\calR}(\calC)$. It is clear that $\calX$ contains the essential image
of the functor $j: \calC \rightarrow \calP^{\calK}_{\calR}(\calC)$. It follows by construction that
$\calX = \calP^{\calK}_{\calR}(\calC)$, so that $\beta$ is an equivalence as desired. This completes the proof of $(2)$.

It remains to prove $(3)$. Suppose that every element of $\calR$ is already a colimit diagram in $\calC$. 
We note that the functor $j$ factors as a composition $L \circ j_0$, where the Yoneda embedding $j_0: \calC \rightarrow \calE$ is already known to be fully faithful (Proposition \ref{fulfaith}).
Since the functor $L | S^{-1} \calP(\calC)$ is equivalent to the identity, it will suffice to show that
the essential image of $j_0$ is contained in $S^{-1} \calP(\calC)$. In other words, we must show that
if $s_{\alpha}: X_{\alpha} \rightarrow j_0 Y_{\alpha}$ belongs to $S$, and $C \in \calC$, then the induced map
$$ \bHom_{\calP(\calC)}( j_0 Y_{\alpha}, j_0 C) \rightarrow \bHom_{ \calP(\calC)}( X_{\alpha}, j_0 C)$$
is a homotopy equivalence. Let $\overline{p}: K^{\triangleright}_{\alpha} \rightarrow \calC$
be the corresponding diagram (so that $\overline{p}$ carries the cone point of
$K^{\triangleright}_{\alpha}$ to $Y_{\alpha}$), let $p = \overline{p} | K_{\alpha}$, and let
$\overline{q}: K^{\triangleright}_{\alpha} \rightarrow \calP(\calC)$ be a colimit diagram
extending $q = \overline{q} | K_{\alpha} = j_0 \circ p$. Consider the diagram
$$ \calP(\calC)_{ j_0 Y_{\alpha}/ } \stackrel{g_0}{\leftarrow} \calP(\calC)_{ j_0 \overline{p}/
} \stackrel{g_1}{\rightarrow} \calP(\calC)_{ j_0 p/ } \stackrel{g_2}{\leftarrow} \calP(\calC)_{ \overline{q}/} \stackrel{g_3}{\rightarrow} \calP(\calC)_{X_{\alpha}/ }.$$
The maps $g_0$ and $g_3$ are trivial Kan fibrations (since the inclusion of the cone point into
$K^{\triangleright}_{\alpha}$ is cofinal), and the map $g_2$ is a trivial Kan fibration since
$\overline{q}$ is a colimit diagram. Moreover, for every object $Z \in \calP(\calC)$, the above diagram determines the map $\bHom_{ \calP(\calC)}( j_0 Y_{\alpha}, Z) \rightarrow \bHom_{ \calP(\calC)}( X_{\alpha}, Z)$. Consequently, to prove that this map is an equivalence, it suffices to show that $g_1$
induces a trivial Kan fibration
$$ \calP(\calC)_{ j_0 \overline{p}/ } \times_{ \calP(\calC) } \{ Z \} 
\rightarrow \calP(\calC)_{ j_0 p/} \times_{ \calP(\calC) } \{Z \}.$$ 
Assuming $Z$ belongs to the essential image $\calC'$ of the Yoneda embedding $j_0$, we may
reduce to proving that the induced map
$ \calC'_{ j_0 \overline{p}/ } \rightarrow \calC'_{ j_0 p/}$ is a trivial Kan fibration, which is equivalent to the assertion that $j_0 \circ \overline{p}$ is a colimit diagram in $\calC'$. This is clear, since
$\overline{p}$ is a colimit diagram by assumption and $j_0$ induces an equivalence of
$\infty$-categories from $\calC$ to $\calC'$.
\end{proof}

\begin{definition}\index{not}{PKKC@$\calP^{\calK'}_{\calK}(\calC)$}
Let $\calK \subseteq \calK'$ be collections of simplicial sets, and let $\calC$ be an $\infty$-category which admits $\calK$-indexed limits. We let $\calP^{\calK'}_{\calK}(\calC)$ denote the $\infty$-category $\calP^{\calK'}_{\calR}(\calC)$, where $\calR$ is the set of all colimit diagrams
$\overline{p}: K^{\triangleright} \rightarrow \calC$ such that $K \in \calK$.
\end{definition}

\begin{example}
Let $\calK = \emptyset$, and let $\calK'$ denote the class of {\em all} small simplicial sets.
If $\calC$ is a small $\infty$-category, then we have a canonical equivalence
$\calP^{\calK'}_{\calK}(\calC) \simeq \calP(\calC)$ (Theorem \ref{charpresheaf}).
\end{example}

\begin{example}
Let $\calK = \emptyset$, and let $\calK'$ denote the class of all small $\kappa$-filtered simplicial sets for some regular cardinal $\kappa$. Then for any small $\infty$-category $\calC$, we have a canonical equivalence
$\calP^{\calK'}_{\calK}(\calC) \simeq \Ind_{\kappa}(\calC)$ (Proposition \ref{intprop}).
\end{example}

\begin{example}
Let $\calK$ denote the collection of all $\kappa$-small simplicial sets for some regular cardinal $\kappa$, and let $\calK'$ be the class of all small simplicial sets. Let $\calC$ be a small $\infty$-category which admits $\kappa$-small colimits. Then we have a canonical equivalence
$\calP^{\calK'}_{\calK}(\calC) \simeq \Ind_{\kappa}(\calC)$. This follows from
Theorem \ref{pretop} and Proposition \ref{sumatch}.
\end{example}

\begin{example}
Let $\calK = \emptyset$, and let $\calK' = \{ \Idem \}$, where $\Idem$ is the simplicial set defined
in \S \ref{retrus}. Then, for any $\infty$-category $\calC$, $\calP^{\calK'}_{\calK}(\calC)$ is
an idempotent competion of $\calC$.
\end{example}

\begin{corollary}
Let $\calK \subseteq \calK'$ be classes of simplicial sets. Let $\widehat{\Cat}_{\infty}$ denote the $\infty$-category of $($not necessarily small$)$ $\infty$-categories, let
$\widehat{\Cat}_{\infty}^{\calK}$ denote the subcategory spanned by those $\infty$-categories
which admit $\calK$-indexed colimits and those functors which preserve $\calK$-indexed colimits, and let $\widehat{\Cat}_{\infty}^{\calK'}$ be defined likewise. Then the inclusion
$$ \widehat{ \Cat}_{\infty}^{\calK'} \subseteq \widehat{ \Cat}_{\infty}^{\calK}$$
admits a left adjoint, given by $\calC \mapsto \calP^{\calK'}_{\calK}(\calC)$.
\end{corollary}

\begin{proof}
Combine Proposition \ref{cupper1} with Proposition \ref{sumpytump}.
\end{proof}

We conclude this section by noting the following transitivity property of the construction
$\calC \mapsto \calP^{\calK}_{\calR}(\calC)$:

\begin{proposition}\label{transit1}
Let $\calK \subseteq \calK'$ be collections of simplicial sets and 
$\calC_1, \ldots, \calC_n$ be a sequence of $\infty$-categories. For
$1 \leq i \leq n$, let $\calR_i$ be a collection of diagrams $\{ \overline{p}_{\alpha}: K^{\triangleright}_{\alpha} \rightarrow \calC_i \}$, where each $K_{\alpha}$ belongs to $\calK$, and let
$\calR'_i$ denote the collection of all colimit diagrams $\{ \overline{q}_{\alpha}: K^{\triangleright}_{\alpha} \rightarrow \calP^{\calK}_{\calR_i}( \calC_i)  \}$ such that $K_{\alpha} \in \calK$.
Then the canonical map
$$ \calP^{\calK'}_{ \calR_1 \boxtimes \ldots \boxtimes \calR_n}( \calC_1 \times \ldots \times \calC_n)
\rightarrow \calP^{\calK'}_{ \calR'_1 \boxtimes \ldots \boxtimes \calR'_{n} }( \calP^{\calK}_{\calR_1}( \calC_1) \times \ldots \times \calP^{\calK}_{ \calR_n}(\calC_n) )$$
is an equivalence of $\infty$-categories. Here $\calR_1 \boxtimes \ldots \boxtimes \calR_n$ denotes
the collection of all diagrams of the form
$$ K_{\alpha}^{\triangleright} \stackrel{ \overline{p}_{\alpha}}{\rightarrow}
\calC_i \simeq \{ C_1 \} \times \ldots \times \{C_{i-1} \} \times \calC_i \times \ldots \times \{ C_n \}
\subseteq \calC_1 \times \ldots \times \calC_n$$
where $\overline{p}_{\alpha} \in \calR_i$ and $C_{j}$ is an object of $\calC_j$ for $j \neq i$, and
the collection $\calR'_1 \boxtimes \ldots \boxtimes \calR'_{n}$ is defined likewise.
\end{proposition}

\begin{proof}
Let $\calD$ be an $\infty$-category which admits $\calK'$-indexed colimits. It will suffice to show that
the functor
$$ \Fun_{\calK'} (  \calP^{\calK'}_{ \calR'_1 \boxtimes \ldots \boxtimes \calR'_{n} }( \calP^{\calK}_{\calR_1}( \calC_1) \times \ldots \times \calP^{\calK}_{ \calR_n}(\calC_n) ), \calD) 
\rightarrow \Fun_{\calK'}(  \calP^{\calK'}_{ \calR_1 \boxtimes \ldots \boxtimes \calR_n}( \calC_1 \times \ldots \times \calC_n), \calD)$$ is an equivalence of $\infty$-categories. Unwinding the definitions, we are reduced to proving that the functor
$$ \phi: \Fun_{ \calR'_1 \boxtimes \ldots \boxtimes \calR'_n}(
\calP^{\calK}_{ \calR_1}( \calC_1) \times \ldots \times \calP^{\calK}_{ \calR_n}(\calC_n), \calD)
\rightarrow \Fun_{ \calR \boxtimes \ldots \boxtimes \calR_n}( \calC_1 \times \ldots \times \calC_n, \calD)$$ is an equivalence of $\infty$-categories.
The proof goes by induction
on $n$. If $n = 0$, then both sides are equivalent to $\calD$ and there is nothing to prove.
If $n > 0$, then set $\calD' = \Fun_{\calR_{n}}( \calC_n, \calD)$ and $\calD'' = \Fun_{ \calK}(
\calP^{\calK}_{\calR_1}(\calC_n), \calD)$. Proposition \ref{cupper1} implies that the canonical map
$\calD'' \rightarrow \calD'$ is an equivalence of $\infty$-categories. We can identify $\phi$ with
the functor
$$ 
 \Fun_{\calR'_1 \boxtimes \ldots \boxtimes \calR'_{n-1}} ( \calP^{\calK}_{\calR_1}( \calC)
 \times \ldots \times \calP^{\calK}_{\calR_{n-1}}(\calC), \calD'') \rightarrow
 \Fun_{ \calR_1 \boxtimes \ldots \boxtimes \calR_{n-1}}( \calC_1 \times \ldots \times \calC_{n-1}, \calD').$$
The desired result now follows from the inductive hypothesis.
\end{proof}

\section{Accessible $\infty$-Categories}\label{c5s5}

\setcounter{theorem}{0}

Many of the categories which commonly arise in mathematics can be realized as categories of $\Ind$-objects. For example, the category of
sets is equivalent to $\Ind(\calC)$, where $\calC$ is the category of finite sets; the category
of rings is equivalent to $\Ind(\calC)$, where $\calC$ is the category of finitely presented rings.
The theory of {\it accessible} categories is an axiomatization of this situation. We refer the reader to \cite{adamek} for an exposition of the theory of accessible categories. In this section, we will describe an $\infty$-categorical generalization of the theory of accessible categories.

We will begin in \S \ref{locbrend} by introducing the notion of a {\em locally small} $\infty$-category. A locally small $\infty$-category $\calC$ need not be small, but has small morphism spaces
$\bHom_{\calC}(X,Y)$ for any fixed pair of objects $X,Y \in \calC$. This is analogous to the usual set-theoretic conventions taken in category theory: one allows categories which have a proper class of objects, but requires that morphisms between any pair of objects form a {\em set}. 

In \S \ref{accessible}, we will introduce the definition of an {\em accessible} $\infty$-category.
An $\infty$-category $\calC$ is accessible if it is locally small and has a good supply of filtered colimits and compact objects. Equivalently, $\calC$ is accessible if it is equivalent to
$\Ind_{\kappa}(\calC^{0})$, for some small $\infty$-category $\calC^0$ and some regular cardinal $\kappa$ (Proposition \ref{clear}).

The theory of accessible $\infty$-categories will play an important technical role throughout the remainder of this book. To understand the usefulness of the hypothesis of accessibility, let us consider the following example. Suppose that $\calC$ is an ordinary category, $F: \calC \rightarrow \Set$ is a functor, and we would like to prove that $F$ is representable by an object $C \in \calC$. The functor
$F$ determines a category $\widetilde{\calC} = \{ (C, \eta): C \in \calC, \eta \in F(C) \}$, which
is fibered over $\calC$ in sets. We would like to prove that $\widetilde{\calC}$ is equivalent
to $\calC_{/C}$, for some $C \in \calC$. The object $C$ can then be characterized as the colimit
of the diagram $p: \widetilde{\calC} \rightarrow \calC$. If $\calC$ admits colimits, then we can attempt to construct $C$ by forming the colimit $\varinjlim(p)$. 

We now encounter a set-theoretic difficulty. Suppose that we try to ensure the existence of $\varinjlim(p)$ by assuming that $\calC$ admits {\em all} small colimits. In this case, it is not reasonable to expect $\calC$ itself to be small. The category $\widetilde{\calC}$ is roughly the same size as $\calC$ (or larger), so our assumption will not allow us to construct $\varinjlim(p)$. On
the other hand, if we assume $\calC$ and $\widetilde{\calC}$ are small, then it is not reasonable
to expect $\calC$ to admit colimits of arbitrary small diagrams.

An accessibility hypothesis can be used to circumvent the difficulty described above. 
An accessible category $\calC$ is generally not small, but is ``controlled'' by a small subcategory $\calC^{0} \subseteq \calC$: it therefore enjoys the best features of both the ``small'' and ``large'' worlds. More precisely, the fiber product $\widetilde{\calC} \times_{\calC} \calC^{0}$ is small enough that we might expect the colimit $\varinjlim(p | \widetilde{\calC} \times_{\calC} \calC^{0})$ to exist on general grounds, yet large enough to expect a natural isomorphism
$$ \varinjlim(p) \simeq \varinjlim( p | \widetilde{\calC} \times_{\calC} \calC^{0}).$$
We refer the reader to \S \ref{aftt} for a detailed account of this argument, which we will use to prove an $\infty$-categorical version of the adjoint functor theorem.

The discussion above can be summarized as follows: the theory of accessible $\infty$-categories is a tool which allows us to manipulate large $\infty$-categories as if they were small, without fear of encountering any set-theoretic paradoxes. This theory is quite useful because the condition of accessibility is very robust: the class of accessible $\infty$-categories is stable under most of the basic constructions of higher category theory. To illustrate this, we will prove the following results:

\begin{itemize}
\item[$(1)$] A small $\infty$-category $\calC$ is accessible if and only if $\calC$ is idempotent complete (\S \ref{accessidem}). 
\item[$(2)$] If $\calC$ is an accessible $\infty$-category and $K$ is a small simplicial set, then
$\Fun(K,\calC)$ is accessible (\S \ref{accessfunk}).
\item[$(3)$] If $\calC$ is an accessible $\infty$-category and $p:K \rightarrow \calC$ is a small diagram, then $\calC_{p/}$ and $\calC_{/p}$ are accessible (\S \ref{accessprime} and \S \ref{accessfiber}). 
\item[$(4)$] The collection of accessible $\infty$-categories is stable under homotopy fiber products (\S \ref{accessfiber}).
\end{itemize}

We will apply these facts in \S \ref{accessstable} to deduce a miscellany of further stability results, which will be needed throughout \S \ref{c5s6} and \S \ref{chap6}.

\subsection{Locally Small $\infty$-Categories}\label{locbrend}

In mathematical practice, it is very common to encounter categories $\calC$ for which the collection of all objects is large (too big to be form a set), but the collection of morphisms $\Hom_{\calC}(X,Y)$ is small for every $X,Y \in \calC$. The same situation arises frequently in higher category theory.
However, it is a slightly trickier to describe, because the formalism of $\infty$-categories blurs the distinction between objects and morphisms. Nevertheless, there is an adequate notion of ``local smallness'' in the $\infty$-categorical setting, which we will describe in this section.

Our first step is to give a characterization of the class of essentially small $\infty$-categories. We will need the following lemma.

\begin{lemma}\label{soirt}
Let $\calC$ be a simplicial category, $n$ a positive integer, and
$f_0: \bd \Delta^n \rightarrow \sNerve(\calC)$ a map. Let
$X = f_0 ( \{0 \})$, $Y = f_0 ( \{n\} )$, and $g_0$ denote the induced map
$$\bd (\Delta^1)^{n-1} \rightarrow \bHom_{\calC}(X,Y).$$
Let $f,f': \Delta^n \rightarrow \sNerve(\calC)$ be extensions of $f_0$, and
$g,g': (\Delta^1)^{n-1} \rightarrow \bHom_{\calC}(X,Y)$ the corresponding extensions
of $g_0$. The following conditions are equivalent:
\begin{itemize}
\item[$(1)$] The maps $f$ and $f'$ are homotopic relative to $\bd \Delta^n$.
\item[$(2)$] The maps $g$ and $g'$ are homotopic relative to $\bd (\Delta^1)^{n-1}$.
\end{itemize}
\end{lemma}

\begin{proof}
It is not difficult to show that $(1)$ is equivalent to the assertion that $f$ and $f'$ are left homotopic in the model category $(\sSet)_{\bd \Delta^n /}$ (with the Joyal model structure), and that $(2)$ equivalent to the assertion that $\sCoNerve[f]$ and $\sCoNerve[f']$ are left homotopic in the model category $(\sCat)_{ \sCoNerve[ \bd \Delta^n] /}$. We now invoke the Quillen equivalence of Theorem \ref{biggier} to complete the proof.
\end{proof}

\begin{proposition}\label{grapeape}
Let $\calC$ be an $\infty$-category, and $\kappa$ an uncountable regular cardinal.
The following conditions are equivalent:
\begin{itemize}
\item[$(1)$] The collection of equivalence classes of objects of $\calC$ is $\kappa$-small, and for every morphism $f: C \rightarrow D$ in $\calC$ and every $n \geq 0$, the homotopy set
$\pi_i( \Hom^{\rght}_{\calC}(C,D), f)$ is $\kappa$-small. 

\item[$(2)$] If $\calC' \subseteq \calC$ is a minimal model for $\calC$, then $\calC'$ is $\kappa$-small.

\item[$(3)$] There exists a $\kappa$-small $\infty$-category $\calC'$ and an equivalence
$\calC' \rightarrow \calC$ of $\infty$-categories.

\item[$(4)$] There exists a $\kappa$-small simplicial set $K$ and a categorical equivalences
$K \rightarrow \calC$.

\item[$(5)$] The $\infty$-category $\calC$ is $\kappa$-compact, when regarded as an object of
$\Cat_{\infty}$.

\end{itemize}
\end{proposition}

\begin{proof}
We begin by proving that $(1) \Rightarrow (2)$. Without loss of generality, we may suppose that
$\calC = \sNerve(\calD)$, where $\calD$ is a topological category. Let $\calC' \subseteq \calC$ be a minimal model for $\calC$. We will prove by induction on $n \geq 0$ that the set
$\Hom_{\sSet}(\Delta^n, \calC')$ is $\kappa$-small. If $n=0$, this reduces to the assertion that
$\calC$ has fewer than $\kappa$ equivalence classes of objects. Suppose therefore that $n > 0$. By the inductive hypothesis, the set $\Hom_{\sSet}(\bd \Delta^n, \calC')$ is $\kappa$-small.
Since $\kappa$ is regular, it will suffice to prove that for each map $f_0: \bd \Delta^n \rightarrow \calC'$, the set $S = \{ f \in \Hom_{\sSet}(\Delta^n, \calC'): f|\bd \Delta^n = f_0 \}$ is $\kappa$-small.
Let $C = f_0(\{0\})$, $D = f_0( \{n\})$, and let $g_0: \bd (\Delta^1)^{n-1} \rightarrow \bHom_{\calD}(C,D)$
be the corresponding map. Assumption $(1)$ ensures that there are fewer than $\kappa$
extensions $g: (\Delta^1)^{n-1} \rightarrow \bHom_{\calD}(C,D)$ modulo homotopy relative to $\bd (\Delta^1)^{n-1}$. Invoking Lemma \ref{soirt}, we deduce that there are fewer than $\kappa$ maps
$f: \Delta^n \rightarrow \calC$ modulo homotopy relative to $\bd \Delta^n$. Since $\calC'$ is minimal, no two distinct elements of $S$ are homotopic in $\calC$ relative to $\bd \Delta^n$; therefore $S$ is $\kappa$-small as desired.

It is clear that $(2) \Rightarrow (3) \Rightarrow (4)$. 
We next show that $(4) \Rightarrow (3)$. Let $K \rightarrow \calC$ be a categorical equivalence, where $K$ is $\kappa$-small.
We construct a sequence of inner anodyne inclusions
$$ K = K(0) \subseteq K(1) \subseteq \ldots $$
Supposing that $K(n)$ has been defined, we form a pushout diagram
$$ \xymatrix{ \coprod \Lambda^n_i \ar@{^{(}->}[r] \ar[d] & \coprod \Delta^n \ar[d] \\
K(n) \ar@{^{(}->}[r] & K(n+1) }$$
where the coproduct is taken over all $0 < i < n$ and all maps $\Lambda^n_i \rightarrow K(n)$.
It follows by induction on $n$ that each $K(n)$ is $\kappa$-small. Since $\kappa$ is regular and uncountable, the limit $K(\infty) = \bigcup_{n} K(n)$ is $\kappa$-small.The inclusion
$K \subseteq K(\infty)$ is inner anodyne; therefore the map $K \rightarrow \calC$ factors through an equivalence $K(\infty) \rightarrow \calC$ of $\infty$-categories; thus $(3)$ is satisfied. 

We next show that $(3) \Rightarrow (5)$. Suppose that $(3)$ is satisfied. Without loss of generality, we may replace $\calC$ by $\calC'$ and thereby suppose that $\calC$ is itself $\kappa$-small. Let
$F: \Cat_{\infty} \rightarrow \SSet$ denote the functor co-represented by $\calC$. According
to Lemma \ref{repco}, we may identify $F$ with the simplicial nerve of the functor
$f: \Cat_{\infty}^{\Delta} \rightarrow \Kan$, which carries an $\infty$-category
$\calD$ to the largest Kan complex contained in $\calD^{\calC}$. Let $\calI$ be a $\kappa$-filtered $\infty$-category and $p: \calI \rightarrow \Cat_{\infty}$ a diagram. We wish to prove that
$p$ has a colimit $\overline{p}: \calI^{\triangleright} \rightarrow \Cat_{\infty}$ such that
$F \circ \overline{p}$ is a colimit diagram in $\SSet$. According to Proposition \ref{rot}, we may suppose that $\calI$ is the nerve of a $\kappa$-filtered partially ordered set $A$. Using Proposition \ref{gumby444}, we may further reduce to the case where $p$ is the simplicial nerve of a diagram $P: A \rightarrow \Cat_{\infty}^{\Delta} \subseteq \mSet$ taking values in the {\em ordinary} category of marked simplicial sets. Let $\overline{P}$ be a colimit of $P$. Since the class of weak equivalences in $\mSet$ is stable under filtered colimits, $\overline{P}$ is a homotopy colimit. Theorem \ref{colimcomparee} implies that $\overline{p} = \sNerve(\overline{P})$ is a colimit of $p$.
It therefore suffices to show that $F \circ \overline{p} = \sNerve( f \circ \overline{P})$ is a colimit diagram. Using Theorem \ref{colimcomparee}, it suffices to show that $f \circ \overline{P}$ is a homotopy colimit diagram in $\sSet$. Since the class of weak homotopy equivalences in $\sSet$ is stable under filtered colimits, it will suffice to prove that $f \circ \overline{P}$ is a colimit diagram in the ordinary category $\sSet$. It now suffices to observe that $f$ preserves $\kappa$-filtered colimits, because $\calC$ is $\kappa$-small.

We now complete the proof by showing that $(5) \Rightarrow (1)$. Let $A$ denote the collection
of all $\kappa$-small simplicial subsets $K_{\alpha} \subseteq \calC$, and let $A' \subseteq A$ be the subcollection consisting of indices $\alpha$ such that $K_{\alpha}$ is an $\infty$-category. It is clear that $A$ is a $\kappa$-filtered partially ordered set, and that
$\calC = \bigcup_{ \alpha \in A} K_{\alpha}$. Using the fact that $\kappa > \omega$, it is easy to see that $A'$ is cofinal in $A$, so that $A'$ is also $\kappa$-filtered and $\calC = \bigcup_{ \alpha \in A'} K_{\alpha}$. We may therefore regard $\calC$ as the colimit of a diagram
$P: A' \rightarrow \mSet$ in the ordinary category of fibrant objects of $\mSet$. Since $A'$ is filtered,
we may also regard $\calC$ as a homotopy colimit of $P$. The above argument shows that
$\calC^{\calC} = f \calC$ can be identified with a homotopy colimit of the diagram
$f \circ P: A' \rightarrow \sSet$. In particular, the vertex $\id_{\calC} \in \calC^{\calC}$
must be homotopic to the image of some map $K_{\alpha}^{\calC} \rightarrow \calC^{\calC}$,
for some $\alpha \in A'$. It follows that $\calC$ is a retract of $K_{\alpha}$ in the homotopy category
$\h{\Cat_{\infty}}$. Since $K_{\alpha}$ is $\kappa$-small, we easily deduce that $K_{\alpha}$ satisfies condition $(1)$. Therefore $\calC$, being a retract of $K_{\alpha}$, satisfies condition $(1)$ as well.
\end{proof}

\begin{definition}\index{gen}{essentially $\kappa$-small!$\infty$-category}\index{gen}{$\infty$-category!essentially $\kappa$-small}
An $\infty$-category $\calC$ is {\it essentially $\kappa$-small} if it satisfies the equivalent conditions of Proposition \ref{grapeape}. We will say that $\calC$ is {\it essentially small} if it is essentially $\kappa$-small for some (small) regular cardinal $\kappa$.
\end{definition}

The following criterion for essential smallness is occasionally useful:

\begin{proposition}\label{sumt}
Let $p: \calC \rightarrow \calD$ be a Cartesian fibration of $\infty$-categories and $\kappa$ an uncountable regular cardinal. Suppose that $\calD$ is essentially $\kappa$-small and that, for each
object $D \in \calD$, the fiber $\calC_{D} = \calC \times_{\calD} \{D\}$ is essentially $\kappa$-small. Then $\calC$ is essentially $\kappa$-small.
\end{proposition}

\begin{proof}
We will apply criterion $(1)$ of Proposition \ref{grapeape}. Choose a $\kappa$-small set of representatives $\{ D_{\alpha} \}$ for the equivalence classes of objects of $\calD$. For each $\alpha$, choose a $\kappa$-small set of representatives $\{ C_{\alpha,\beta} \}$ for the equivalence classes of objects of $\calC_{D_{\alpha}}$. The collection of all objects
$C_{\alpha,\beta}$ is $\kappa$-small (since $\kappa$ is regular) and contains representatives for all equivalence classes of objects of $\calC$.

Now suppose that $C$ and $C'$ are objects of $\calC$, having images $D, D' \in \calD$.
Since $\calD$ is essentially $\kappa$-small, the set $\pi_0 \bHom_{\calD}(D,D')$ is $\kappa$-small. Let $f: D \rightarrow D'$ be a morphism, and choose a $p$-Cartesian morphism
$\widetilde{f}: \widetilde{C} \rightarrow D'$ covering $f$. According to Proposition \ref{compspaces}, we have a homotopy fiber sequence
$$ \bHom_{\calC_{D}}( C, \widetilde{C}) \rightarrow \bHom_{\calC}(C,C') \rightarrow
\bHom_{\calD}(D,D')$$
in the homotopy category $\calH$. In particular, we see that $\bHom_{\calC}(C,C')$ contains
fewer than $\kappa$ connected components lying over $f \in \pi_0 \bHom_{\calD}(D,D')$, and therefore fewer than $\kappa$ components in total (since $\kappa$ is regular). Moreover, the long exact sequence of homotopy groups shows that for every $\overline{f}: C \rightarrow C'$ lifting
$f$, the homotopy sets $\pi_{i}( \Hom^r_{\calC}(C,C'), f)$ are $\kappa$-small, as desired. 
\end{proof}

By restricting our attention to {\em Kan complexes}, we obtain an analogue of Proposition \ref{grapeape} for spaces:

\begin{corollary}\label{apegrape}\index{gen}{essentially $\kappa$-small!space}Let $X$ be a Kan complex, and $\kappa$ an {\em uncountable} regular cardinal.
The following conditions are equivalent:
\begin{itemize}

\item[$(1)$] For each vertex $x \in X$ and each $n \geq 0$, the homotopy set
$\pi_n(X,x)$ is $\kappa$-small.

\item[$(2)$] If $X' \subseteq X$ is a minimal model for $X$, then $X'$ is $\kappa$-small.

\item[$(3)$] There exists a $\kappa$-small Kan complex $X'$ and a homotopy equivalence $X' \rightarrow X$.

\item[$(4)$] There exists a $\kappa$-small simplicial set $K$ and a weak homotopy equivalence $K \rightarrow X$.

\item[$(5)$] The $\infty$-category $\calC$ is $\kappa$-compact, when regarded as an object of
$\SSet$.

\item[$(6)$] The Kan complex $X$ is essentially small $($when regarded as an $\infty$-category$)$.

\end{itemize}
\end{corollary}

\begin{proof}
The equivalences $(1) \Leftrightarrow (2) \Leftrightarrow (3) \Leftrightarrow (6)$ follow from Proposition \ref{grapeape}.
The implication $(3) \Rightarrow (4)$ is obvious. We next prove that $(4) \Rightarrow (5)$.
Let $p: K \rightarrow \SSet$ be the constant diagram taking the value $\ast$, let
$\overline{p}: K^{\triangleright} \rightarrow \SSet$ be a colimit of $p$, and let $X' \in \SSet$ be the image under $\overline{p}$ of the cone point of $K^{\triangleright}$. It follows from Proposition \ref{dda} that $\ast$ is a $\kappa$-compact object of $\SSet$. Corollary \ref{tyrmyrr} implies that
$X'$ is a $\kappa$-compact object of $\SSet$. Let $\widetilde{K} \rightarrow K^{\triangleright}$ denote the left fibration associated to $\overline{p}$, and let $X'' \subseteq \widetilde{K}$
denote the fiber lying over the cone point of $K^{\triangleright}$. The inclusion of the cone point in $K^{\triangleright}$ is right anodyne. It follows from Proposition \ref{strokhop} that 
the inclusion $X'' \subseteq \widetilde{K}$ is right anodyne. Since $\overline{p}$ is a colimit diagram, Proposition \ref{charspacecolimit} implies that the inclusion
$K \simeq K \times_{ K^{\triangleright}} \widetilde{K} \subseteq \widetilde{K}$ is a weak homotopy equivalence. We therefore have a chain of weak homotopy equivalences
$$ X \leftarrow K \subseteq \widetilde{K} \leftarrow X'' \leftarrow X',$$
so that $X$ and $X'$ are equivalent objects of $\SSet$. Since $X'$ is $\kappa$-compact, it follows that $X$ is $\kappa$-compact.

To complete the proof, we will show that $(5) \Rightarrow (1)$. 
We employ the argument used in the proof of Proposition \ref{grapeape}. 
Let $F: \SSet \rightarrow \SSet$
be the functor co-represented by $X$. Using Lemma \ref{repco}, we can identify $F$ can be with the simplicial nerve of the functor $f: \Kan \rightarrow \Kan$ given by
$$ Y \mapsto Y^X.$$Let $A$ denote the collection of $\kappa$-small simplicial subsets $X_{\alpha} \subseteq X$ which are Kan complexes. Since $\kappa$ is uncountable,
$A$ is $\kappa$-filtered and $X = \bigcup_{\alpha \in A} K_{\alpha}$. We may regard $X$
as the colimit of a diagram $P: A \rightarrow \sSet$. Since $A$ is filtered, $X$ is also a homotopy colimit of this diagram. Since $F$ preserves $\kappa$-filtered colimits, $f$ preserves $\kappa$-filtered homotopy colimits; therefore $X^X$ is a homotopy colimit of the diagram
$f \circ P$.  In particular, the vertex $\id_{X} \in X^{X}$
must be homotopic to the image of some map $X_{\alpha}^{X} \rightarrow X^{X}$,
for some $\alpha \in A$. It follows that $X$ is a retract of $X_{\alpha}$ in the homotopy category $\calH$. Since $X_{\alpha}$ is $\kappa$-small, we can readily verify that $X_{\alpha}$ satisfies $(1)$. Because $X$ is a retract of $X_{\alpha}$, $X$ satisfies $(1)$ as well.
\end{proof}

\begin{remark}
When $\kappa = \omega$, the situation is quite a bit more complicated. Suppose that
$X$ is a Kan complex representing a compact object of $\SSet$. Then there exists a simplicial set $Y$ with only finitely many nondegenerate simplices, and a map $i: Y \rightarrow X$ which
realizes $X$ as a {\em retract} of $Y$ in the homotopy category $\calH$ of spaces. However, one cannot generally assume that $Y$ is a Kan complex, or that $i$ is a weak homotopy equivalence.
The latter can be achieved if $X$ is connected and simply connected, or more generally if a certain $K$-theoretic invariant of $X$ (the {\em Wall finiteness obstruction} ) vanishes: we refer the reader to \cite{wall} for a discussion.\index{gen}{Wall finiteness obstruction}
\end{remark}

For many applications, it is important to be able to slightly relax the condition that an $\infty$-category be essentiall small.

\begin{proposition}\label{locsm}
Let $\calC$ be an $\infty$-category. The following conditions are equivalent:
\begin{itemize}
\item[$(1)$] For every pair of objects $X,Y \in \calC$, the space
$\bHom_{\calC}(X,Y)$ is essentially small.
\item[$(2)$] For every small collection $S$ of objects of $\calC$, the full subcategory
of $\calC$ spanned by the elements of $S$ is essentially small.
\end{itemize}
\end{proposition}

\begin{proof}
This follows immediately from criterion $(1)$ in Propositions \ref{grapeape} and \ref{apegrape}.
\end{proof}

We will say that an $\infty$-category $\calC$ is {\it locally small} if it satisfies the
equivalent conditions of Proposition \ref{locsm}.\index{gen}{locally small $\infty$-category}\index{gen}{$\infty$-category!locally small}

\begin{example}\label{exlocsm}
Let $\calC$ and $\calD$ be $\infty$-categories. Suppose that $\calC$ is locally small and that $\calD$ is essentially small. Then $\calC^{\calD}$ is essentially small. To prove this, we may assume without loss of generality that $\calC$ and $\calD$ are minimal. Let
$\{ \calC_{\alpha} \}$ denote the collection of all full subcategories of $\calC$, spanned by small collections of objects. Since $\calD$ is small, every finite collection of functors
$\calD \rightarrow \calC$ factors through some small $\calC_{\alpha} \subseteq \calC$. 
It follows that $\Fun(\calD,\calC)$ is the union of it small full subcategories $\Fun(\calD, \calC_{\alpha})$, and is therefore locally small. In particular, for every small $\infty$-category $\calD$, the $\infty$-category $\calP(\calD)$ of presheaves is locally small.
\end{example}

\subsection{Accessibility}\label{accessible}

In this section, we will begin our study of the class of accessible $\infty$-categories.

\begin{definition}\label{kapacc}\index{gen}{accessible!$\infty$-category}\index{gen}{$\kappa$-accessible $\infty$-category}\index{gen}{$\infty$-category!accessible}
Let $\kappa$ be a regular cardinal.
An $\infty$-category $\calC$ is {\it $\kappa$-accessible} if there exists a small
$\infty$-category $\calC^{0}$ and an equivalence $\Ind_{\kappa}(\calC^{0}) \rightarrow \calC$.
We will say that $\calC$ is {\it accessible} if it is $\kappa$-accessible for {\em some} regular cardinal $\kappa$.
\end{definition}

The following result gives a few alternative characterizations of the class of accessible $\infty$-categories.

\begin{proposition}\label{clear}
Let $\calC$ be an $\infty$-category and $\kappa$ a regular cardinal. The
following conditions are equivalent:

\begin{itemize}
\item[$(1)$] The $\infty$-category $\calC$ is $\kappa$-accessible.

\item[$(2)$] The $\infty$-category $\calC$ is locally small, admits $\kappa$-filtered colimits, the
full subcategory $\calC^{\kappa} \subseteq \calC$ of $\kappa$-compact objects is essentially small,  and $\calC^{\kappa}$ generates $\calC$ under small, $\kappa$-filtered colimits.

\item[$(3)$] The $\infty$-category $\calC$ admits small $\kappa$-filtered colimits and contains an essentially small full subcategory $\calC'' \subseteq \calC$ which consists of $\kappa$-compact objects and generates $\calC$ under small $\kappa$-filtered colimits.
\end{itemize}
\end{proposition}

The main obstacle to proving Proposition \ref{clear} is in
verifying that if $\calC_0$ is small, then $\Ind_\kappa(\calC_0)$
has only a bounded number of $\kappa$-compact objects, up to equivalence. It
is tempting to guess that any such object must be equivalent to an
object of $\calC_0$. The following example shows that this is not
necessarily the case.

\begin{example}
Let $R$ be a ring, and let $\calC_0$ denote the (ordinary)
category of finitely generated free $R$-modules. Then $\calC =
\Ind(\calC_0)$ is equivalent to the category of flat $R$-modules
(by Lazard's theorem; see for example the appendix of
\cite{lazard}). The compact objects of $\calC$ are precisely the
finitely generated projective $R$-modules, which need not be free.
\end{example}

Nevertheless, the naive guess is not far off, in virtue of the following result:

\begin{lemma}\label{stylus}\index{gen}{idempotent completion}
Let $\calC$ be a small $\infty$-category, $\kappa$ a regular cardinal, and
$\calC' \subseteq \Ind_{\kappa}(\calC)$ the full subcategory of $\Ind_{\kappa}(\calC)$ spanned by the $\kappa$-compact objects. Then the Yoneda embedding
$j: \calC \rightarrow \calC'$ exhibits $\calC'$ as an idempotent completion of $\calC$. In particular, $\calC'$ is essentially small.
\end{lemma}

\begin{proof}
Corollary \ref{swwe} implies that $\Ind_{\kappa}(\calC)$ is idempotent complete. Since $\calC'$ is stable under retracts in $\Ind_{\kappa}(\calC)$, $\calC'$ is also idempotent complete. Proposition \ref{fulfaith} implies that $j$ is fully faithful. It therefore suffices to prove that every object $C' \in \calC'$ is a retract of $j(C)$, for some $C \in \calC$. 

Let $\calC_{/C'} = \calC \times_{ \Ind_{\kappa}(\calC)} \Ind_{\kappa}(\calC)_{/C'}$. Lemma \ref{longwait0} implies that the diagram
$$ \overline{p}: \calC_{/C'}^{\triangleright} \rightarrow \Ind_{\kappa}(\calC)_{/C'}^{\triangleright} \rightarrow \Ind_{\kappa}(\calC)$$ is a colimit of $p= \overline{p} | \calC_{/C'}$. Let $F: \Ind_{\kappa}(\calC) \rightarrow \SSet$ be the functor co-represented by $C'$; we note that the left fibration associated to $F$ is equivalent to $\Ind_{\kappa}(\calC)_{C'/}$. Since $F$ is $\kappa$-continuous, Proposition \ref{charspacecolimit} implies that the inclusion
$$\calC_{/C'} \times_{ \Ind_{\kappa}(\calC)} \Ind_{\kappa}(\calC)_{C'/}
\subseteq \calC^{\triangleright}_{/C'} \times_{ \Ind_{\kappa}(\calC)} \Ind_{\kappa}(\calC)_{C'/}$$
is a weak homotopy equivalence. The simplicial set on the right has a canonical vertex, corresponding to the identity map $\id_{C'}$. It follows that there exists a vertex on the left hand side belonging to the same path component. Such a vertex classifies a diagram
$$ \xymatrix{ & j(C) \ar[dr] & \\
C' \ar[ur] \ar[rr]^{f} & & C' }$$ where $f$ is homotopic to the identity, which proves that $C'$ is a retract of $j(C)$ in $\Ind_{\kappa}(\calC)$.
\end{proof}

\begin{proof}[Proof of Proposition \ref{clear}]
Suppose that $(1)$ is satisfied. Without loss of generality we may suppose that
$\calC = \Ind_{\kappa} \calC'$, where $\calC'$ is small. Since $\calC$ is a full subcategory of $\calP(\calC')$, it is locally small (see Example \ref{exlocsm}). Proposition \ref{geort}
implies that $\calC$ admits small $\kappa$-filtered colimits. Corollary \ref{indpr} shows that
$\calC$ is generated under $\kappa$-filtered colimits by the essential image of the Yoneda embedding $j: \calC' \rightarrow \calC$, which consists of $\kappa$-compact objects by Proposition \ref{justcut}. Lemma \ref{stylus} implies that full subcategory of $\Ind_{\kappa}(\calC')$ consisting of compact objects is essentially small. We conclude that $(1) \Rightarrow (2)$. 

It is clear that $(2) \Rightarrow (3)$. Suppose that $(3)$ is satisfied. Choose a small $\infty$-category $\calC'$ and an equivalence $i: \calC' \rightarrow \calC''$. Using Proposition \ref{intprop}, we may suppose that $i$ factors as a composition
$$ \calC' \stackrel{j}{\rightarrow} \Ind_{\kappa}(\calC') \stackrel{f}{\rightarrow} \calC$$
where $f$ preserves small $\kappa$-filtered colimits. It follows from Proposition \ref{uterr} that
$f$ is a categorical equivalence. This shows that $(3) \Rightarrow (1)$ and completes the proof.
\end{proof}

\begin{definition}\label{accfun}\index{gen}{accessible!functor}
If $\calC$ is an accessible $\infty$-category, then a functor
$F: \calC \rightarrow \calC'$ is {\it accessible} if it is $\kappa$-continuous for some regular cardinal $\kappa$ (and therefore for all regular cardinals $\tau \geq \kappa$).
\end{definition}

\begin{remark}
Generally we will only speak of the accessibility of a functor $F: \calC \rightarrow \calC'$
in the case where both $\calC$ and $\calC'$ are accessible. However, it is occasionally convenient to use the terminology of Definition \ref{accfun} in the case where $\calC$ is accessible and $\calC'$ is not (or $\calC'$ is not yet known to be accessible).
\end{remark}

\begin{example}\label{spacesareaccessible}
The $\infty$-category $\SSet$ of spaces is accessible. More generally, for any small
$\infty$-category $\calC$, the $\infty$-category $\calP(\calC)$ is accessible: this follows immediately from Proposition \ref{precst}.
\end{example}

If $\calC$ is a $\kappa$-accessible $\infty$-category and $\tau > \kappa$, then $\calC$
is not necessarily $\tau$-accessible. Nevertheless, this is true for many values of $\tau$.

\begin{definition}\label{ineq}\index{not}{kappalltau@$\kappa \ll \tau$}
Let $\kappa$ and $\tau$ be regular cardinals. We write $\tau \ll \kappa$ if the following condition is satisfied: for every $\tau_0 < \tau$ and every $\kappa_0 < \kappa$, we have $\kappa_0^{\tau_0} < \kappa$.
\end{definition}

Note that there exist arbitrarily large regular
cardinals $\kappa'$ with $\kappa' \gg \kappa$: for example, one
may take $\kappa'$ to be the successor of any cardinal having the
form $\tau^{\kappa}$.

\begin{remark}
Every (infinite) regular cardinal $\kappa$ satisfies $\omega \ll \kappa$.
An uncountable regular cardinal $\kappa$ satisfies $\kappa \ll \kappa$ if and only if $\kappa$ is strongly inaccessible.
\end{remark}

\begin{lemma}\label{estimate}
If $\kappa' \gg \kappa$, then any $\kappa'$-filtered partially
ordered set $\calI$ may be written as a union of $\kappa$-filtered
subsets having size $< \kappa'$. Moreover, the family of all such
subsets is $\kappa'$-filtered.
\end{lemma}

\begin{proof}
It will suffice to show that every subset of $S \subseteq \calI$
having cardinality $< \kappa'$ can be included in a larger subset
$S'$, such that $|S'| < \kappa'$, but $S'$ is $\kappa$-filtered.

We define a transfinite sequence of subsets $S_{\alpha} \subseteq
\calI$ by induction. Let $S_{0} = S$, and when $\lambda$ is a
limit ordinal we let $S_{\lambda} = \bigcup_{\alpha < \lambda}
S_{\alpha}$. Finally, we let $S_{\alpha + 1}$ denote a set which
is obtained from $S_{\alpha}$ by adjoining an upper bound for
every subset of $S_{\alpha}$ having size $< \kappa$ (which exists
because $\calI$ is $\kappa'$-filtered). It follows from the
assumption $\kappa' \gg \kappa$ that if $S_{\alpha}$ has size $<
\kappa'$, then so does $S_{\alpha + 1}$. Since $\kappa'$ is
regular, we deduce easily by induction that $|S_{\alpha}| <
\kappa'$ for all $\alpha < \kappa'$. It is easy to check that the
set $S' = S_{\kappa}$ has the desired properties.
\end{proof}

\begin{proposition}\label{enacc}
Let $\calC$ be a $\kappa$-accessible $\infty$-category. Then
$\calC$ is $\kappa'$-accessible for any $\kappa' \gg \kappa$.
\end{proposition}

\begin{proof}
Let $\calC^{\kappa} \subseteq \calC$ denote the full subcategory consisting of $\kappa$-compact objects, and let $\calC' \subseteq \calC$ denote the full subcategory spanned by the colimits of all $\kappa'$-small, $\kappa$-filtered diagrams in $\calC^{\kappa}$. Since $\calC$ is locally small and the collection of all equivalence classes of such diagrams is bounded, we conclude that $\calC'$ is essentially small. Corollary \ref{tyrmyrr} implies that $\calC'$ consists of $\kappa'$-compact objects of $\calC$. According to Proposition \ref{clear}, it will suffice to prove that $\calC'$ generates $\calC$ under small $\kappa'$-filtered colimits. Let $X$ be an object of $\calC$, and let
$p: \calI \rightarrow \calC^{\kappa}$ be a small $\kappa$-filtered diagram with colimit $X$.
Using Proposition \ref{rot}, we may reduce to the case where $\calI$ is the nerve of a $\kappa$-filtered partially ordered set $A$. Lemma \ref{estimate} implies that $A$ can be written as a $\kappa'$-filtered union of $\kappa'$-small, $\kappa$-filtered subsets $\{ A_{\beta} \subseteq A \}_{\beta \in B}$. Using Propositions \ref{extet} and \ref{utl}, we deduce that $X$ can also be obtained as the colimit of a diagram indexed by $\Nerve(B)$, which takes values in $\calC'$.
\end{proof}

\begin{remark}
If $\calC$ is a $\kappa$-accessible $\infty$-category and $\kappa' > \kappa$, then $\calC$ is generally not $\kappa'$-accessible. There are counterexamples even in ordinary category theory: see \cite{adamek}.
\end{remark}

\begin{remark}\label{boundedacc}
Let $\calC$ be an accessible $\infty$-category and $\kappa$ a regular cardinal. Then the full subcategory $\calC^{\kappa} \subseteq \calC$ consisting of $\kappa$-compact objects is essentially small. To prove this, we are free to enlarge $\kappa$ and we may invoke Proposition \ref{enacc} to reduce to the case where $\calC$ is $\kappa$-accessible, in which case the desired result is a consequence of Proposition \ref{clear}.
\end{remark}

\begin{notation}\index{not}{FunAcc@$\aHom(\calC, \calD)$}
If $\calC$ and $\calD$ are accessible $\infty$-categories, we will write
$\aHom(\calC, \calD)$ to denote the full subcategory of $\Fun(\calC, \calD)$
spanned by accessible functors from $\calC$ to $\calD$.
\end{notation}

\begin{remark}
Accessible $\infty$-categories are usually not small. However, they are determined by a ``small'' amount of data: namely, they always have the form $\Ind_{\kappa}(\calC)$ where $\calC$ is a small $\infty$-category. Similarly, an accessible functor $F: \calC \rightarrow \calD$ between accessible categories is determined by a ``small'' amount of data, in the sense that there always exists
a regular cardinal $\kappa$ such that $F$ is $\kappa$-continuous and maps
$\calC^{\kappa}$ into $\calD^{\kappa}$. The restriction $F| \calC^{\kappa}$ then determines
$F$ up to equivalence, by Proposition \ref{intprop}. To prove the existence of $\kappa$, 
we first choose a regular cardinal $\tau$ such that $F$ is $\tau$-continuous. Enlarging $\tau$ if necessary, we may suppose that $\calC$ and $\calD$ are $\tau$-accessible. The collection of equivalence classes of $\tau$-compact objects of $\calC$ is small; consequently, Remark \ref{boundedacc} there exists a (small) regular cardinal $\tau'$ such that
$F$ carries $\calC^{\tau}$ into $\calD^{\tau'}$. We may now choose $\kappa$ to be any regular cardinal such that $\kappa \gg \tau'$. 
\end{remark}

\begin{definition}\index{not}{Acck@$\Acc_{\kappa}$}\index{not}{Acc@$\Acc$}
Let $\kappa$ be a regular cardinal.
We let $\Acc_{\kappa} \subseteq \widehat{\Cat}_{\infty}$ denote the subcategory defined as follows:
\begin{itemize}
\item[$(1)$] The objects of $\Acc_{\kappa}$ are the $\kappa$-accessible $\infty$-categories.
\item[$(2)$] A functor $F: \calC \rightarrow \calD$ between accessible $\infty$-categories
belongs to $\Acc$ if and only if $F$ is $\kappa$-continuous and preserves $\kappa$-compact objects. 
\end{itemize}
Let $\Acc = \bigcup_{\kappa} \Acc_{\kappa}$.
We will refer to $\Acc$ as the {\it $\infty$-category of accessible $\infty$-categories}.
\end{definition}

\begin{proposition}\label{humatch}
Let $\kappa$ be a regular cardinal, and let $\theta: \Acc_{\kappa} \rightarrow \widehat{\Cat}_{\infty}$ be the simplicial nerve of the functor which associates to each $\calC \in \Acc_{\kappa}$ the full subcategory of $\calC$ spanned by the $\kappa$-compact objects. Then:
\begin{itemize}
\item[$(1)$] The functor $\theta$ is fully faithful. 
\item[$(2)$] An $\infty$-category $\calC \in \widehat{\Cat}_{\infty}$ belongs to the essential image of $\theta$ if and only if $\calC$ is essentially small and idempotent complete.
\end{itemize}
\end{proposition}

\begin{proof}
Assertion $(1)$ follows immediately from Proposition \ref{intprop}. If $\calC \in \widehat{\Cat}_{\infty}$ belongs to the essential image of $\theta$, then $\calC$ is
essentially small and idempotent complete (because $\calC$ is stable under retracts in an idempotent complete $\infty$-category). Conversely, suppose that $\calC$ is essentially small and idempotent complete, and choose a minimal model $\calC' \subseteq \calC$. Then $\Ind_{\kappa}(\calC')$ is $\kappa$-accessible. Moreover, the collection of $\kappa$-compact objects of $\Ind_{\kappa}(\calC')$ is an idempotent completion of $\calC'$ (Lemma \ref{stylus}), and therefore equivalent to $\calC$ (since $\calC'$ is already idempotent complete).
\end{proof}

Let $\Cat_{\infty}^{\vee}$ denote the full subcategory of $\Cat_{\infty}$ spanned by the idempotent complete $\infty$-categories.\index{not}{Catinftyvee@$\Cat_{\infty}^{\vee}$}

\begin{proposition}
The inclusion $\Cat_{\infty}^{\vee} \subseteq \Cat_{\infty}$ has a left adjoint.
\end{proposition}

\begin{proof}
Combine Propositions \ref{idmcoo}, \ref{charidemcomp}, and \ref{testreflect}.
\end{proof}

We will refer to a left adjoint to the inclusion $\Cat_{\infty}^{\vee} \subseteq \Cat_{\infty}$ as the
{\it idempotent completion functor}. Proposition \ref{humatch} implies that we have fully faithful embeddings $\Acc_{\kappa} \rightarrow \widehat{\Cat}_{\infty} \hookleftarrow \Cat_{\infty}^{\vee}$ with the same essential image. Consequently, there is a (canonical) equivalence of $\infty$-categories $e: \Cat_{\infty}^{\vee} \simeq \Acc_{\kappa}$, well-defined up to homotopy. We let
$\Ind_{\kappa}: \Cat_{\infty} \rightarrow \Acc_{\kappa}$ denote the composition of $e$ with the idempotent completion functor. In summary:

\begin{proposition}\label{lockap}\index{not}{Indkappa@$\Ind_{\kappa}$}
There is a functor $\Ind_{\kappa}: \Cat_{\infty} \rightarrow \Acc_{\kappa}$ which exhibits
$\Acc_{\kappa}$ as a localization of the $\infty$-category $\Cat_{\infty}$.
\end{proposition}

\begin{remark}
There is a slight danger of confusion with our terminology. The functor $\Ind_{\kappa}: \Cat_{\infty} \rightarrow \Acc_{\kappa}$ is only well-defined up to contractible space of choices, so that
if $\calC$ is an $\infty$-category which admits finite colimits, then the image of $\calC$ under $\Ind_{\kappa}$ is only well-defined up to equivalence. Definition \ref{indsmall} produces a canonical representative for this image.
\end{remark}

\subsection{Accessibility and Idempotent Completeness}\label{accessidem}

Let $\calC$ be an accessible $\infty$-category. Then there exists a regular cardinal
$\kappa$ such that $\calC$ admits $\kappa$-filtered colimits. It follows from Corollary \ref{swwe} that $\calC$ is idempotent complete. Our goal in this section is to prove a converse to this result: if $\calC$ is a small and idempotent complete, then $\calC$ is accessible. 

Let $\calC$ be a small $\infty$-category, and suppose we want to prove that $\calC$ is accessible. 
The main problem is to show that $\calC$ admits $\kappa$-filtered colimits, provided that $\kappa$ is sufficiently large. The idea is that if $\kappa$ is much larger than the size of $\calC$, then any $\kappa$-filtered diagram $\calJ \rightarrow \calC$ is necessarily very ``redundant'' (Proposition \ref{stufenn}). Before we making this precise, we will need a few preliminary results.

\begin{lemma}\label{techycard}
Let $\kappa < \tau$ be uncountable regular cardinals, $A$ a $\tau$-filtered partially ordered set, and $F: A \rightarrow \Kan$ a diagram of Kan complexes indexed by $A$. Suppose that for each $\alpha \in A$, the Kan complex $F(\alpha)$ is essentially $\kappa$-small.
For every $\tau$-small subset $A_0 \subseteq A$, there exists a filtered
$\tau$-small subset $A'_0 \subseteq A$ containing $A_0$, with the property that the map
$$ \colim_{\alpha \in A'_0} F(\alpha) \rightarrow \colim_{\alpha \in A} F(\alpha)$$ is a homotopy equivalence.
\end{lemma}

\begin{proof}
Let $X = \colim_{\alpha \in A} F(\alpha)$. Since $F$ is a filtered diagram, $X$ is also a Kan complex. 
Let $K$ be a simplicial set with only finitely many nondegenerate simplices. Our first claim is that
the set $[K,X]$ of homotopy classes of maps from $K$ into $X$ has cardinality $< \kappa$. For suppose given a collection $\{ g_{\beta}: K \rightarrow X \}$ of pairwise nonhomotopic maps, having cardinality $\kappa$. Since $A$ is $\tau$-filtered, we may suppose that there
is a fixed index $\alpha \in A$ such each $g_{\beta}$ factors as a composition
$$ K \stackrel{g'_{\beta}}{\rightarrow} F(\alpha) \rightarrow X.$$
The maps $g'_{\beta}$ are also pairwise nonhomotopic, which contracts our assumption that
$F(\alpha)$ is weakly homotopy equivalent to a $\kappa$-small simplicial set. 

We now define an increasing sequence $$\alpha_0 \leq \alpha_1 \leq \ldots $$
of elements of $A$. Let $\alpha_0$ be any upper bound for $A_0$. Assuming that $\alpha_i$ has already been selected, choose a representative for every homotopy class of diagrams
$$ \xymatrix{ \bd \Delta^n \ar@{^{(}->}[d] \ar[r] & F(\alpha_i) \ar[d] \\
\Delta^n \ar[r]^{h_{\gamma}} & X. }$$
The argument above proves that we can take the set of all such representatives to be $\kappa$-small, so that there exists $\alpha_{i+1} \geq \alpha_i$ such that each $h_{\gamma}$ factors
as a composition
$$ \Delta^{n} \stackrel{h'_{\gamma}}{\rightarrow} F(\alpha_{i+1}) \rightarrow X$$
and the associated diagram 
$$ \xymatrix{ \bd \Delta^n \ar@{^{(}->}[d] \ar[r] & F(\alpha_i) \ar[d] \\
\Delta^n \ar[r]^-{h'_{\gamma}} & F(\alpha_{i+1}) }$$
is commutative. We now set $A'_0 = A_0 \cup \{\alpha_0, \alpha_1, \ldots \}$; it is easy to check
that this set has the desired properties.
\end{proof}

\begin{lemma}\label{techycardd}
Let $\kappa < \tau$ be uncountable regular cardinals, let $A$ be a $\tau$-filtered partially ordered set, let $\{ F_{\beta} \}_{\beta \in B}$ be a collection of diagrams $A \rightarrow \sSet$ indexed by
a set $B$ of cardinality $< \tau$. Suppose that for each $\alpha \in A$ and each
$\beta \in B$, the Kan complex $F_{\beta}(\alpha)$ is essentially $\kappa$-small.
Then there exists a filtered, $\tau$-small subset
$A' \subseteq A$ such that for each $\beta \in B$, the map
$$ \colim_{A'} F_{\beta}(\alpha) \rightarrow \colim_{A} F_{\beta}(\alpha)$$
is a homotopy equivalence of Kan complexes.
\end{lemma}

\begin{proof}
Without loss of generality, we may suppose that $B = \{ \beta: \beta < \beta_0 \}$ is a set of ordinals.
We will define a sequence of filtered, $\tau$-small subsets $A(n) \subseteq A$ by induction on $n$. For $n = 0$, choose an element $\alpha \in A$ and set $A(0) = \{ \alpha\}$. Suppose next that $A(n)$ has been defined. We define a sequence of enlargements $\{ A(n)_{\beta} \}_{\beta \leq \beta_0}$ by induction on $\beta$. Let $A(n)_{0} = A(n)$, let $A(n)_{\lambda} = \bigcup_{\beta < \lambda} A(n)_{\beta}$ when $\lambda$ is a nonzero limit ordinal, and let
$A(n)_{\beta+1}$ be a $\tau$-small, filtered subset of $A$ such that the map
$$ \colim_{A(n)_{\beta+1} } F_{\beta}(\alpha) \rightarrow \colim_{A} F_{\beta}(\alpha)$$ is a weak homotopy equivalence (such a subset exists in virtue of Lemma 
\ref{techycard}). We now take $A(n+1) = A(n)_{\beta_0}$ and $A' = \bigcup_{n} A(n)$; it is easy to check that $A' \subseteq A$ has the desired properties.
\end{proof}

\begin{lemma}\label{sorens}
Let $\kappa < \tau$ be uncountable regular cardinals. Let $\calC$ be a $\tau$-small $\infty$-category with the property that each of the spaces $\bHom_{\calC}(C,D)$ is essentially $\kappa$-small, and $j: \calC \rightarrow \calP(\calC)$ the Yoneda embedding. Let
$p: \calK \rightarrow \calC$ be a diagram indexed by a $\tau$-filtered $\infty$-category
$\calK$, and $\overline{p}: \calK^{\triangleright} \rightarrow \calP(\calC)$ a colimit of $j \circ p$.
Then there exists a map $i: K \rightarrow \calK$ such that $K$ is $\tau$-small, and the composition
$\overline{p} \circ i^{\triangleright}: K^{\triangleright} \rightarrow \calK^{\triangleright} \rightarrow \calP(\calC)$ is a colimit diagram.
\end{lemma}

\begin{proof}
In view of Proposition \ref{rot}, we may suppose that $\calK$ is the nerve of a $\tau$-filtered partially ordered set $A$. According to Proposition \ref{limiteval}, $\overline{p}$ induces a colimit diagram
$$ \overline{p}_{C}: \calK^{\triangleright} \rightarrow \calP(\calC) \stackrel{e_C}{\rightarrow} \SSet$$
where $e_{C}$ denote the evaluation functor associated to an object $C \in \calC$. We will
identify $\calK^{\triangleright}$ with the nerve of the partially ordered set $A \cup \{\infty\}$. 
Proposition \ref{gumby444} implies that we may replace $\overline{p}_{C}$ with the simplicial nerve
of a functor $F_{C}: A \cup \{ \infty\} \rightarrow \Kan$. Our hypothesis on $\calC$ implies
that $F_{C}|A$ takes values in $\kappa$-small simplicial sets. Applying Theorem \ref{colimcomparee}, we see that the map $\colim_{A} F_{C}(\alpha) \rightarrow F_{C}(\infty)$ is a homotopy equivalence. We now apply Lemma \ref{techycard} to deduce the
existence of a filtered, $\tau$-small subset $A' \subseteq A$ such that each of the maps
$$ \colim_{A'} F_{C}(\alpha) \rightarrow F_{C}(\infty)$$ is a homotopy equivalence.
Let $K = \Nerve(A')$, and let $i: K \rightarrow \calK$ denote the inclusion. Using Theorem \ref{colimcomparee} again, we deduce that the composition
$e_{C} \circ  \overline{p} \circ i^{\triangleright}: K^{\triangleright} \rightarrow \SSet$ is a colimit diagram for each $C \in \calC$. Applying Proposition \ref{limiteval}, we deduce that
$\overline{p} \circ i^{\triangleright}$ is a colimit diagram, as desired.
\end{proof}

\begin{proposition}\label{stufenn}
Let $\kappa < \tau$ be uncountable regular cardinals. Let $\calC$ be an $\infty$-category which is $\tau$-small, such that the morphism spaces $\bHom_{\calC}(C,D)$ are essentially $\kappa$-small.
Let $j: \calC \rightarrow \calP(\calC)$ denote the Yoneda embedding, let
$p: \calK \rightarrow \calC$ be a diagram indexed by a $\tau$-filtered $\infty$-category
$\calK$, and let $X \in \calP(\calC)$ be a colimit of $j \circ p: \calK \rightarrow \calP(\calC)$. 
Then there exists an object $C \in \calC$ such that $X$ is a retract of $j(C)$.
\end{proposition}

\begin{proof}
Let $i: K \rightarrow \calK$ be a map satisfying the conclusions of Lemma \ref{sorens}. Since
$K$ is $\tau$-small and $\calK$ is $\tau$-filtered, there exists an extension
$\overline{i}: K^{\triangleright} \rightarrow \calK$ of $i$. Let $C$ be the image of the cone point
of $K^{\triangleright}$ under $p \circ \overline{i}$, and $\widetilde{C} \in \calC_{p \circ i/}$ the
corresponding lift. Let $\overline{p}: \calK^{\triangleright} \rightarrow \calP(\calC)$
be a colimit of $j \circ p$ carrying the cone point of $\calK^{\triangleright}$ to
$X$. Let $q = j \circ p \circ i: K \rightarrow \calP(\calC)$, $\widetilde{X} \in \calP(\calC)_{q/}$ the corresponding lift of $X$, and $\widetilde{Y} \in \calP(\calC)_{q/}$ a colimit of $q$. 
Since $\widetilde{Y}$ is an initial object of $\calP(\calC)_{q/}$, there is a commutative triangle
$$ \xymatrix{ & j( \widetilde{C}) \ar[dr] & \\
\widetilde{Y} \ar[rr] \ar[ur] & & \widetilde{X} }$$
in the $\infty$-category $\calP(\calC)_{q/}$. Moreover, Lemma \ref{sorens} asserts that
the horizontal map is an equivalence. Thus $\widetilde{X}$ is a retract of
$j(\widetilde{C})$ in the homotopy category of $\calP(\calC)_{q/}$, so that
$X$ is a retract of $j(C)$ in $\calP(\calC)$.
\end{proof}

\begin{corollary}\label{tyrrus}
Let $\kappa < \tau$ be uncountable regular cardinals, and let $\calC$ be a $\tau$-small $\infty$-category whose morphism spaces $\bHom_{\calC}(C,D)$ are essentially $\kappa$-small. 
Then the Yoneda embedding $j: \calC \rightarrow \Ind_{\tau}(\calC)$ exhibits
$\Ind_{\tau}(\calC)$ as an idempotent completion of $\calC$.
\end{corollary}

\begin{proof}
Since $\Ind_{\tau}(\calC)$ admits $\tau$-filtered colimits, it is idempotent complete by Corollary \ref{swwe}. Proposition \ref{stufenn} implies that every object of $\Ind_{\tau}(\calC)$ is a retract of $j(C)$, for some object $C \in \calC$.
\end{proof}

\begin{corollary}\label{sloam}\index{gen}{idempotent complete!and accessibility}
A small $\infty$-category $\calC$ is accessible if and only if it is idempotent complete. Moreover, if these conditions are satisfied and $\calD$ is an any accessible $\infty$-category, then
{\em every} functor $f: \calC \rightarrow \calD$ is accessible.
\end{corollary}

\begin{proof}
The ``only if'' follows from Corollary \ref{swwe}, and the ``if'' direction follows from
Corollary \ref{tyrrus}. Now suppose that $\calC$ is small and accessible, and let
$\calD$ be a $\kappa$-accessible $\infty$-category and $f: \calC \rightarrow \calD$ any functor; we wish to prove that $f$ is accessible. By Proposition \ref{intprop}, we may suppose that
$f = F \circ j$, where $j: \calC \rightarrow \Ind_{\kappa}(\calC)$ is the Yoneda embedding
and $F: \Ind_{\kappa}(\calC) \rightarrow \calD$ is a $\kappa$-continuous functor, and therefore accessible. Enlarging $\kappa$ if necessary, we may suppose that $j$ is an equivalence of $\infty$-categories, so that $f$ is accessible as well.
\end{proof}

\subsection{Accessibility of Functor $\infty$-Categories}\label{accessfunk}

Let $\calC$ be an accessible $\infty$-category, and let $K$ be a small simplicial set. Our goal in this section is to prove that $\Fun(K,\calC)$ is accessible (Proposition \ref{horse1}). In \S \ref{accessstable}, we will prove a much more general stability result of this kind (Corollary \ref{storkus1}), but the proof of that result ultimately rests on the ideas presented here. 

Our proof goes roughly as follows. If $\calC$ is accessible, then $\calC$ has a many $\tau$-compact objects, provided that $\tau$ is sufficiently large. Using Proposition \ref{placeabovee}, we deduce the existence of many $\tau$-compact objects in $\Fun(K,\calC)$. Our main problem is to show that these objects generate $\Fun(K,\calC)$ under $\tau$-filtered colimits. To prove this, we will use a rather technical cofinality result (Lemma \ref{kidav} below). We begin with the following prelimiinary observation:

\begin{lemma}\label{pprekidav}
Let $\tau$ be a regular cardinal, and let $q: Y \rightarrow X$ be a coCartesian fibration with the property that for every vertex $x$ of $X$, the fiber $Y_{x} = Y \times_{X} \{x\}$ is $\tau$-filtered.
Then $q$ has the right lifting property with respect to $K \subseteq K^{\triangleright}$, for every
$\tau$-small simplicial set $K$.
\end{lemma}

\begin{proof}
Using Proposition \ref{princex}, we can reduce to the problem of showing that
$q$ has the right lifting property with respect to the inclusion $K \subseteq K \diamond \Delta^0$. 
In other words, we must show that given any edge $e: C \rightarrow D$ in $X^K$,
where $D$ is a constant map, and any vertex $\widetilde{C}$ of $Y^K$ lifting
$C$, there exists an edge $\widetilde{e}: \widetilde{C} \rightarrow \widetilde{D}$
lifting $\widetilde{e}$, where $\widetilde{D}$ is a constant map from $K$ to $Y$.
We first choose an arbitrary edge $\widetilde{e}': \widetilde{C} \rightarrow \widetilde{D}'$ lifting $e$ (since the map $q^K: Y^K \rightarrow X^K$ is a coCartesian fibration, we can even choose $\widetilde{e}'$ to be $q^K$-coCartesian, though we will not need this). Suppose that
$D$ takes the constant value $x: \Delta^0 \rightarrow X$. Since the fiber $Y_{x}$ is $\tau$-filtered,
there exists an edge $\widetilde{e}'': \widetilde{D}' \rightarrow \widetilde{D}$ in
$Y_{x}^K$, where $\widetilde{D}$ is a constant map from $K$ to $Y_{x}$. We now invoke
the fact that $q^K$ is an inner fibration to supply the dotted arrow in the diagram
$$ \xymatrix{ \Lambda^2_1 \ar[rr]^{ (\widetilde{e}', \bigdot, \widetilde{e}'')} \ar@{^{(}->}[d] & & Y^K \ar[d] \\
\Delta^2 \ar@{-->}[urr]^{\sigma} \ar[rr]^{ s_1 e} & & X^K. }$$
We now define $\widetilde{e} = \sigma | \Delta^{ \{0,2\} }$. 
\end{proof}

\begin{lemma}\label{kidav}
Let $\kappa < \tau$ be regular cardinals. Let $q: Y \rightarrow X$ be a map of simplicial sets with the following properties:
\begin{itemize}
\item[$(i)$] The simplicial set $X$ is $\tau$-small.
\item[$(ii)$] The map $q$ is a coCartesian fibration.
\item[$(iii)$] For every vertex $x \in X$, the fiber $Y_{x} = Y \times_{X} \{x\}$ is $\tau$-filtered
and admits $\tau$-small, $\kappa$-filtered colimits.
\item[$(iv)$] For every edge $e: x \rightarrow y$ in $X$, the associated functor
$Y_{x} \rightarrow Y_{y}$ preserves $\tau$-small, $\kappa$-filtered colimits.
\end{itemize}
Then:
\begin{itemize}
\item[$(1)$] The $\infty$-category $\calC = \bHom_{/X}(X,Y)$ of sections of $q$ is $\tau$-filtered.
\item[$(2)$] For each vertex $x$ of $X$, the evaluation map
$e_{x}: \calC \rightarrow Y_{x}$ is cofinal.
\end{itemize}
\end{lemma}

\begin{proof}
Choose a categorical equivalence $X \rightarrow M$, where $M$ is a minimal $\infty$-category. Since $\tau$ is uncountable, Proposition \ref{grapeape} implies that $M$ is $\tau$-small.
According to Corollary \ref{tttroke}, $Y$ is equivalent to the pullback of a coCartesian fibration $Y' \rightarrow M$. We may therefore replace $X$ by $M$ and thereby reduce to the case where
$X$ is a minimal $\infty$-category. For each ordinal $\alpha$, let $(\alpha)= \{ \beta < \alpha\}$. 

Let $K$ be a $\tau$-small simplicial set equipped with a map $f: K \rightarrow Y$.
We define a new object $K'_{X} \in (\sSet)_{/X}$ as follows. For every finite, nonempty, linearly ordered set $J$, a map $\Delta^J \rightarrow K'_X$ is determined by the following data:

\begin{itemize}
\item A map $\chi: \Delta^{J} \rightarrow X$.

\item A map $\Delta^J \rightarrow \Delta^2$, corresponding to a decomposition
$J = J_0 \coprod J_1 \coprod J_2$.

\item A map $\Delta^{J_0} \rightarrow K$.

\item An order-preserving map $m: J_1 \rightarrow (\kappa)$, having the property that
if $m(i) = m(j)$, then $\chi( \Delta^{ \{i,j\} })$ is a degenerate edge of $X$.
\end{itemize} 

We will prove the existence of a dotted arrow $F'_X$ as indicated in the diagram
$$ \xymatrix{ K \ar[d] \ar[r]^{f} & Y \ar[d]^{q} \\
K'_X \ar@{-->}[ur]^{F'_X} \ar[r] & X. }$$
Let $K'' \subseteq K'_X$ be the simplicial subset corresponding to simplices, as above, where
$J_1 = \emptyset$, and let $F'' = F'_X | K''$. Specializing to the case where $K = Z \times X$, $Z$ a $\tau$-small simplicial set, we will deduce that any diagram $Z \rightarrow \calC$ extends to a map $Z^{\triangleright} \rightarrow \calC$ (given by $F''$), which proves $(1)$. Similarly, by specializing to the case $K = (Z \times X) \coprod_{ Z \times \{x\} } (Z^{\triangleleft} \times \{x\} )$, we will deduce that for every object $y \in Y$ with $q(y) = x$, the $\infty$-category
$\calC \times_{ Y_{x} } (Y_{x})_{y/}$ is $\tau$-filtered, and therefore weakly contractible. 
Applying Theorem \ref{hollowtt}, we deduce $(2)$. 

It remains to construct the map $F'_X$. There
is no harm in enlarging $K$. We may therefore apply the small object argument
to replace $K$ by an $\infty$-category (which we may also suppose is $\tau$-small, since $\tau$ is uncountable). We begin by defining, for each $\alpha \leq \kappa$,
a simplicial subset $K(\alpha) \subseteq K'_X$. The definition is as follows: we will say
that a simplex $\Delta^J \rightarrow K'_X$ factors through $K(\alpha)$ if, in the corresponding decomposition $J = J_0 \coprod J_1 \coprod J_2$, we have $J_2 = \emptyset$, and the map
$J_1 \rightarrow (\kappa)$ factors through $(\alpha)$. Our first task is to construct
$F(\alpha) = F'_X | K(\alpha)$, which we do by induction on $\alpha$.
If $\alpha=0$, $K(\alpha) = K$ and we set $F(0) = f$.
When $\alpha$ is a limit ordinal, we have $K(\alpha) = \bigcup_{ \beta < \alpha } K(\beta)$ and we set $F(\alpha) = \bigcup_{\beta < \alpha} F(\beta)$. 
It therefore suffices to construct $F(\alpha+1)$, assuming that $F(\alpha)$ has already been constructed. For each vertex $x$ of $X$, let $\widetilde{x}=(x, \alpha)$ denote the unique vertex of $K(\alpha+1)$ lying over $x$ which does not belong to $K(\alpha)$. Since $X$ is minimal, Proposition \ref{minstrict} implies that we have a pushout diagram
$$ \xymatrix{ \coprod_{x} K(\alpha)_{/\widetilde{x}} \ar@{^{(}->}[r] \ar[d] & \coprod_{x} (K(\alpha)_{/\widetilde{x}})^{\triangleright} \ar[d]  \\
K(\alpha) \ar@{^{(}->}[r] & K(\alpha+1). }$$
Therefore, to construct $f_{\alpha+1}$, it suffices to prove that $q$ has the right lifting property with respect to each inclusion $K(\alpha)_{/\widetilde{x}} \subseteq (K(\alpha)_{/\widetilde{x}})^{\triangleright}$, which follows
from Lemma \ref{pprekidav}.

We now define, for each simplicial subset $X' \subseteq X$, a corresponding simplicial subset
$K'_{X'} \subseteq K'_{X}$. The definition is as follows: let $\sigma: \Delta^J \rightarrow K'_{X}$ be a simplex corresponding to a decomposition $J = J_0 \coprod J_1 \coprod J_2$. Then $\sigma$
factors through $K'_{X'}$ if and only if the induced map $\Delta^{J_2} \rightarrow X$ factors through $X'$. Our next job is to extend the definition of $F'_{X}$ from $K'_{\emptyset} = K(\kappa)$
to $K'_{X}$, by adjoining simplices to $X$ one at a time.

Let $F'_{\emptyset} = F(\kappa)$, and let $x$ be a vertex of $X$. We begin by defining
a map $F'_{\{x\} } : K'_{\{x\} } \rightarrow Y$ which extends $F'_{\emptyset}$. 
Since $X$ is minimal, there is a pushout diagram
$$ \xymatrix{ K(\kappa)_{/x} \ar@{^{(}->}[r] \ar[d] & K(\kappa)_{/x}^{\triangleright} \ar[d] \\
K_{\emptyset} \ar@{^{(}->}[r] &  K_{\{x\} } } $$
where $K(\kappa)_{/x}$ denotes the fiber product $K(\kappa) \times_{X} X_{/x}$. 
Constructing an extension $F'_{ \{x\} }$ of $F'_{\emptyset}$ is therefore equivalent to providing the dotted arrow indicated in the diagram
$$ \xymatrix{ K(\kappa)_{/x} \ar@{^{(}->}[d] \ar[r]^{p_x} & Y \ar[d] \\
K(\kappa)_{/x}^{\triangleright} \ar[r] \ar@{-->}[ur]^{\overline{p}_x} & X }.$$
We will choose $\overline{p}_x$ to be a relative colimit of $p_x$ over $X$ (see \S \ref{relcol}). To prove that
such a relative colimit exists, we consider the inclusion
$i_{x}: \Nerve (\kappa) \subseteq K(\kappa)_{/x} \times_{X_{/x}} \{ \id_{x} \} \subseteq K(\kappa)_{/x}$. Using Proposition \ref{minstrict}, it is not difficult to see that
$K(\kappa)_{/x}$ is an $\infty$-category. For each 
object $y \in K(\kappa)_{/x}$, the minimality of $X$ implies that $\Nerve (\kappa) \times_{ \calK_{/x} } (\calK_{/x})_{y/}$ is isomorphic to $\Nerve(\{ \alpha: \beta < \alpha < \kappa \})$ for some $\beta < \kappa$, and therefore weakly contractible. Theorem \ref{hollowtt} implies that $i_{x}$ is cofinal. Invoking Proposition \ref{relexist}, it will suffice to prove that 
$p_x \circ i_{x}: \Nerve (\kappa) \rightarrow Y$ admits a relative colimit over $X$. Using conditions
$(ii)$, $(iv)$, and Proposition \ref{relcolfibtest}, we may reduce to producing a colimit
of $p_{x} \circ i_{x}$ in the $\infty$-category $Y_{x}$, which is possible in virtue of assumption $(iii)$.

Applying the above argument separately to each vertex of $X$, we may suppose that
$F'_{X^{(0)}}$ has been constructed, where $X^{(0)}$ denotes the $0$-skeleton of $X$. We now consider the collection of all pairs $( X', F'_{X'})$ where $X'$ is a simplicial subset of $X$ containing all vertices of $X$, and $F'_{X'}: K_{X'} \rightarrow Y$ is a map over $X$ whose restriction to
$K_{X^{0}}$ coincides with $F'_{X^{(0)}}$. This collection is partially ordered, if we write
$(X', F'_{X'}) \leq (X'', F'_{X''})$ to mean that $X'\subseteq X''$ and $F'_{X''} | K_{X'} = F'_{X'}$. 
The hypotheses of Zorn's lemma are satisfied, so that there exists a maximal such pair
$(X', F'_{X'})$. To complete the proof, it suffices to show that $X' = X$. If not, we can choose
$X' \subseteq X'' \subseteq X$, where $X''$ is obtained from $X'$ by adjoining a single nondegenerate simplex $\sigma: \Delta^n \rightarrow X$ whose boundary already belongs to $X'$.
Since $X'$ contains $X^{(0)}$, we may suppose that $n > 0$. Let $K(\kappa)_{/\sigma} = K(\kappa) \times_{X} X_{/\sigma}$, and let
$x = \sigma(0)$. Since $X$ is minimal, we have a pushout diagram
$$ \xymatrix{ K(\kappa)_{/\sigma} \star \bd \Delta^n \ar@{^{(}->}[r] \ar[d] & K(\kappa)_{/\sigma} \star \Delta^n \ar[d] \\
K'_{X'} \ar@{^{(}->}[r] & K'_{X''}. }$$
Let $s: K(\kappa)_{/\sigma} \rightarrow Y$ denote the composition of the projection $K(\kappa)_{/\sigma} \rightarrow K'_{X'}$ with
$F'_{X'}$.  We obtain a commutative diagram
$$ \xymatrix{ \bd \Delta^n \ar[r]^{r} \ar@{^{(}->}[d] & Y_{s/} \ar[d] \\
\Delta^n \ar[r] \ar@{-->}[ur] & X_{q \circ s/}, }$$
and supplying the indicated dotted arrow is tantamount to giving a map $F'_{X''}: K_{X''} \rightarrow Y$ over $X$ which extends $F'_{X'}$. To prove the existence of $F'_{X''}$, it suffices to prove
that the map $\overline{s}: {K'}^{\triangleright} \rightarrow Y$ associated to $r(0)$ is a 
a $q$-colimit diagram. We note that $\overline{s}$ is given as a composition
$$ {K'}^{\triangleright} \rightarrow K_{/x}^{\triangleright} \stackrel{\overline{s}'}{\rightarrow} Y,$$
where $\overline{s}'$ is a $q$-colimit diagram by construction. According to Proposition \ref{relexists}, it will suffice to show that the map $K(\kappa)_{/\sigma} \rightarrow K(\kappa)_{/x}$ is cofinal.
We have a pullback diagram
$$ \xymatrix{ K(\kappa)_{/\sigma} \ar[r] \ar[d] & K(\kappa)_{/x} \ar[d] \\
X_{/\sigma} \ar[r] & X_{/x} }$$ where the lower horizontal map is a trivial fibration of simplicial sets.
It follows that the upper horizontal map is a trivial fibration, and in particular cofinal. Consequently, there exists an extension $F_{X''}$ of $F_{X'}$, which contradicts the maximality of $(X', F_{X'})$ and completes the proof.
\end{proof}

\begin{proposition}\label{horse1}\index{gen}{accessible!functor categories}
Let $\calC$ be an accessible $\infty$-category, and let $K$ be a small simplicial set. Then
$\Fun(K,\calC)$ is accessible.
\end{proposition}

\begin{proof}
Without loss of generality, we may suppose that $K$ is an $\infty$-category.
Choose a regular cardinal $\kappa$ such that $\calC$ admits small $\kappa$-filtered colimits, and choose a second regular cardinal $\tau > \kappa$ such that $\calC$ is also $\tau$-accessible and $K$ is $\tau$-small.
We will prove that $\Fun(K,\calC)$ is $\tau$-accessible. Let $\calC' = \Fun(K,\calC^{\tau}) \subseteq \Fun(K,\calC)$. It is clear
that $\calC'$ is essentially small. Proposition \ref{limiteval} implies that $\Fun(K,\calC)$ admits small $\tau$-filtered colimits, and Proposition \ref{placeabovee} asserts that $\calC'$ consists of $\tau$-compact objects of $\Fun(K,\calC)$. According to Proposition \ref{clear}, it will suffice to prove that $\calC'$ generates $\Fun(K,\calC)$ under small, $\tau$-filtered colimits.

Without loss of generality, we may suppose that $\calC = \Ind_{\tau} \calD'$, where $\calD'$ is a small $\infty$-category. Let $\calD \subseteq \calC$ denote the essential image of the Yoneda embedding. Let $F: K \rightarrow \calC$ be an arbitrary object of $\calC^K$, and let $\Fun(K,\calD)^{/F} = \Fun(K,\calD) \times_{\Fun(K,\calC)} \Fun(K,\calC)^{/F}$. 
Consider the composite diagram
$$ \overline{p}: \Fun(K,\calD)^{/F} \diamond \Delta^0 \rightarrow
\Fun(K,\calC)^{/F} \diamond \Delta^0 \rightarrow \Fun(K,\calC).$$
The $\infty$-category $\Fun(K,\calD)^{/F}$ is equivalent to 
$\Fun(K,\calD') \times_{\Fun(K,\calC)} \Fun(K,\calC)^{/F}$, and therefore essentially small.
To complete the proof, it will suffice to show that $\Fun(K,\calD)^{/F}$ is $\tau$-filtered, and that
$\overline{p}$ is a colimit diagram.

We may identify $F$ with a map $f_K: K \rightarrow \calC \times K$ in $(\sSet)_{/K}$. According to 
Proposition \ref{colimfam}, we obtain a coCartesian fibration $q: (\calC \times K)^{/f_K} \rightarrow K$, 
and the $q$-coCartesian morphisms are precisely those which project to equivalences in $\calC$. Let
$X$ denote the full subcategory of $(\calC \times K)^{/f_K}$ consisting of those objects whose projection to $\calC$ belongs to $\calD$. It follows that $q' = q|X: X \rightarrow K$ is a coCartesian fibration. We may identify the fiber of $q'$ over a vertex
$x \in K$ with $\calD^{/F(x)} = \calD \times_{\calC} \calC^{/F(x)}$. It follows that the fibers
of $q'$ are $\tau$-filtered $\infty$-categories; Lemma \ref{kidav} now guarantees that
$\Fun(K,\calD)^{/F} \simeq \bHom_{/K}(K,X)$ is $\tau$-filtered.

According to Proposition \ref{limiteval}, to prove that $\overline{p}$ is a colimit diagram, it will suffice to prove that for every vertex $x$ of $K$, the composition of $\overline{p}$ with the evaluation
map $e_{x}: \Fun(K,\calC) \rightarrow \calC$ is a colimit diagram. The composition
$e_{x} \circ \overline{p}$ admits a factorization
$$ \Fun(K,\calD)^{/F} \diamond \Delta^0 \rightarrow \calD^{/F(x)} \diamond \Delta^0
\rightarrow \calC$$
where $\calD^{/F(x)} = \calD \times_{\calC} \calC^{/F(x)}$ and the second map is a colimit diagram in $\calC$ by Lemma \ref{longwait0}. It will therefore suffice to prove that the map
$g_x: \Fun(K,\calD)^{/F} \rightarrow \calD^{/F(x)}$ is cofinal, which follows from Lemma \ref{kidav}.
\end{proof}

\subsection{Accessibility of Undercategories}\label{accessprime}

Let $\calC$ be an accessible $\infty$-category, and let $p: K \rightarrow \calC$ be a small diagram. Our goal in this section is to prove that the $\infty$-category $\calC_{p/}$ is accessible (Corollary \ref{horsemn}). 

\begin{remark}
The analogous result for the $\infty$-category $\calC_{/p}$ will be proven in \S \ref{accessfiber}, using Propositions \ref{horse1} and \ref{horse2}. It is possible to use the same argument to give a second proof of Corollary \ref{horsemn}; however, we will {\em need} Corollary \ref{horsemn} in our proof of Proposition \ref{horse2}.
\end{remark}

We begin by studying the behavior of colimits with respect to (homotopy) fiber products of $\infty$-categories.

\begin{lemma}\label{bird1}\index{gen}{initial object!in a homotopy fiber product}
Let $$ \xymatrix{ \calX' \ar[r]^{q'} \ar[d]^{p'} & \calX \ar[d]^{p} \\
\calY' \ar[r]^{q} & \calY }$$
be a diagram of $\infty$-categories which is homotopy Cartesian $($with respect to the Joyal model structure$)$. Suppose that $\calX$ and $\calY$ have initial objects, and that $p$ and $q$ preserve initial objects. An object
$X' \in \calX'$ is initial if and only if $p'(X')$ is an initial object of $\calY'$ and
$q'(X')$ is an initial object of $\calX$. Moreover, there exists an initial object of $\calX'$.
\end{lemma}

\begin{proof}
Without loss of generality, we may suppose that $p$ and $q$ are categorical fibrations, and that
$\calX' = \calX \times_{ \calY} \calY'$. Suppose first that $X'$ is an object of $\calX'$ with the property that $X = q'(X')$ and $Y' = p'(X')$ are initial objects of $\calX$ and $\calY'$. Then
$Y = p(X) = q(Y')$ is an initial object of $\calY$. Let $Z$ be another object of $\calX'$.
We have a pullback diagram of Kan complexes
$$ \xymatrix{ \Hom_{\calX'}^{\rght}(X', Z) \ar[r] \ar[d] & \Hom_{\calX}^{\rght}(X, q'(Z)) \ar[d] \\
\Hom_{\calY'}^{\rght}(Y', p'(Z)) \ar[r] & \Hom_{\calY}^{\rght}(Y, (q \circ p')(Z)). }$$
Since the maps $p$ and $q$ are inner fibrations, Lemma \ref{sharpy} implies that this diagram is homotopy Cartesian (with respect to the usual model structure on $\sSet$). Since
$X$, $Y'$, and $Y$ are initial objects, the Kan complexes $\Hom_{\calX}^{\rght}(X, q'(Z))$,
$\Hom_{\calY'}^{\rght}(Y', p'(Z))$, and $\Hom_{\calY}^{\rght}(Y, (q \circ p')(Z))$ are contractible. It follows that $\Hom_{\calX'}^{\rght}(X', Z)$ is contractible as well, so that $X'$ is an initial object of $\calX'$. 

We now prove that there exists an object $X' \in \calX'$ such that $p'(X')$ and $q'(X')$ are initial.
The above argument shows that $X'$ is an initial obejct of $\calX'$. Since all initial objects of $\calX'$ are equivalent, this will prove that for {\em any} initial object $X'' \in \calX'$, the
objects $p'(X'')$ and $q'(X'')$ are initial.

We begin by selecting arbitrary initial objects $X \in \calX$ and $\overline{Y} \in \calY'$.  Then
$p(X)$ and $q(\overline{Y})$ are both initial objects of $\calY$, so there is an equivalence
$e: p(X) \rightarrow q(\overline{Y})$. Since $q$ is a categorical fibration, there exists an equivalence
$\overline{e}: Y' \rightarrow \overline{Y}$ in $\calY$ such that $q(\overline{e}) = e$. It follows that $Y'$ is an initial object of $\calY'$ with $q(Y') = p(X)$, so that the pair $(X,Y')$ can be identified with an object of $\calX'$ which has the desired properties.
\end{proof}

\begin{lemma}\label{airbirdd}
Let $p: \calX \rightarrow \calY$ be a categorical fibration of $\infty$-categories, and
let $f: K \rightarrow \calX$ be a diagram. Then the induced map
$p': \calX_{f/} \rightarrow \calY_{p  f/}$ is a categorical fibration.
\end{lemma}

\begin{proof}
It suffices to show that $p'$ has the right lifting property with respect to every inclusion
$A \subseteq B$ which is a categorical equivalence. Unwinding the definitions, it suffices
to show that $p$ has the right lifting property with respect to $i: K \star A \subseteq K \star B$.
This is immediate, since $p$ is a categorical fibration and $i$ is a categorical equivalence.
\end{proof}

\begin{lemma}\label{airbird}\index{gen}{undercategory!and homotopy fiber products}
Let $$ \xymatrix{ \calX' \ar[r]^{q'} \ar[d]^{p'} & \calX \ar[d]^{p} \\
\calY' \ar[r]^{q} & \calY }$$
be a diagram of $\infty$-categories which is homotopy Cartesian (with respect to the Joyal model structure), and let $f: K \rightarrow \calX'$ be a diagram in $\calX'$. Then the induced
diagram
$$ \xymatrix{ \calX'_{f/} \ar[r] \ar[d] & \calX_{q'  f/} \ar[d] \\
\calY'_{p'  f/} \ar[r] & \calY_{q  p'  f/}}$$
is also homotopy Cartesian.
\end{lemma}

\begin{proof}
Without loss of generality, we may suppose that $p$ and $q$ are categorical fibrations
and that $\calX' = \calX \times_{\calY} \calY'$. Then $\calX'_{f/} \simeq
\calX_{q'  f/} \times_{ \calY_{q  p'  f/} } \calY'_{p'  f/}$, so the result
follows immediately from Lemma \ref{airbirdd}.
\end{proof}

\begin{lemma}\label{bird3}\index{gen}{colimit!and homotopy fiber products}
Let $$ \xymatrix{ \calX' \ar[r]^{q'} \ar[d]^{p'} & \calX \ar[d]^{p} \\
\calY' \ar[r]^{q} & \calY }$$
be a diagram of $\infty$-categories which is homotopy Cartesian (with respect to the Joyal model structure), and let $K$ be a simplicial set. Suppose that $\calX$ and $\calY'$ admit colimits for all diagrams indexed by $K$, and that $p$ and $q$ preserve colimits of diagrams indexed by $K$. 
Then:
\begin{itemize}
\item[$(1)$] A diagram $\overline{f}: K^{\triangleright} \rightarrow \calX'$ is a colimit of $f=\overline{f}|K$ if and only if
$p' \circ \overline{f}$ and $q' \circ \overline{f}$ are colimit diagrams. In particular, 
$p'$ and $q'$ preserve colimits indexed by $K$.
\item[$(2)$] Every diagram $f: K \rightarrow \calX'$ has a colimit in $\calX'$.
\end{itemize}
\end{lemma}

\begin{proof}
Replacing $\calX'$ by $\calX'_{f/}$,
$\calX$ by $\calX_{q'  f/}$, $\calY'$ by $\calY'_{p'  f/}$, and $\calY$ by
$\calY_{q  p'  f/}$, we may apply Lemma \ref{airbird} to reduce to the case $K = \emptyset$. Now apply Lemma \ref{bird1}.
\end{proof}

\begin{lemma}\label{filterprod}
Let $\calC$ be a small filtered category, and let $\calC^{\triangleright}$ be the category
obtained by adjoining a $($new$)$ final object to $\calC$.
Suppose given a homotopy pullback diagram
$$ \xymatrix{ F' \ar[r] \ar[d] & F \ar[d]^{p} \\
G' \ar[r]^{q} & G }$$
in the diagram category $\Set_{\Delta}^{\calC^{\triangleright}}$ $($which we endow with the {\it projective}
model structure$)$. Suppose further that the diagrams $F, G, G': \calC^{\triangleright} \rightarrow \sSet$ are homotopy colimits. Then $F'$ is also a homotopy colimit diagram.
\end{lemma}

\begin{proof}
Without loss of generality, we may suppose that $G$ is fibrant, $p$ and $q$ are fibrations, and that
$F' = F \times_{G} G'$. Let $\ast$ denote the cone point of $\calC^{\triangleright}$, and let
$F(\infty)$, $G(\infty)$, $F'(\infty)$, and $G'(\infty)$ denote the colimits of the diagrams
$F|\calC$, $G|\calC$, $F'| \calC$, and $G'| \calC$. Since fibrations in $\sSet$ are stable under filtered colimits, the pullback diagram
$$ \xymatrix{ F'(\infty) \ar[r] \ar[d] & F(\infty) \ar[d] \\
G'(\infty) \ar[r] & G(\infty) }$$
exhibits $F'(\infty)$ as a homotopy fiber product of $F(\infty)$ with $G'(\infty)$ over
$G(\infty)$ in $\sSet$. 
Since weak homotopy equivalences are stable under filtered colimits, the natural maps $G(\infty) \rightarrow G(\ast)$, $F'(\infty) \rightarrow F'(\ast)$, and $G'(\infty) \rightarrow G'(\ast)$ are weak homotopy equivalences. Consequently, the diagram
$$ \xymatrix{ F'(\infty) \ar[dr]^{f} & & \\
& F'(\ast) \ar[r] \ar[d] & F(\ast) \ar[d] \\
& G'(\ast) \ar[r] & G(\ast) }$$
exhibits both $F'(\infty)$ and $F'(\ast)$ as homotopy fiber products of $F(\ast)$ with
$G'(\ast)$ over $G(\ast)$. It follows that $f$ is a weak homotopy equivalence, so that $F$ is a homotopy colimit diagram as desired.
\end{proof}

\begin{lemma}\label{yoris}
Let $$ \xymatrix{ \calX' \ar[r]^{q'} \ar[d]^{p'} & \calX \ar[d]^{p} \\
\calY' \ar[r]^{q} & \calY }$$
be a diagram of $\infty$-categories which is homotopy Cartesian $($with respect to the Joyal model structure$)$, and let $\kappa$ be a regular cardinal. Suppose that $\calX$ and $\calY'$ admit small
$\kappa$-filtered colimits, and that $p$ and $q$ preserve small $\kappa$-filtered colimits.
Then:
\begin{itemize}
\item[$(1)$] The $\infty$-category $\calX'$ admits small $\kappa$-filtered colimits.
\item[$(2)$] If $X'$ is an object of $\calX'$ such that $Y' = p'(X')$ and $X = q'(X')$,
and $Y = p(X) = q(Y')$ are $\kappa$-compact, then $X'$ is a $\kappa$-compact object of $\calX'$.
\end{itemize}
\end{lemma}

\begin{proof}
Claim $(1)$ follows immediately from Lemma \ref{bird3}. To prove $(2)$, consider a colimit diagram
$\overline{f}: \calI^{\triangleright} \rightarrow \calX'$. We wish to prove that the composition of $\overline{f}$ with the functor $\calX' \rightarrow \hat{\SSet}$
corepresented by $X'$ is also a colimit diagram. Using Proposition \ref{rot}, we may assume without loss of generality that $\calI$ is the nerve of a $\kappa$-filtered partially ordered set $A$. We may further suppose that $p$ and $q$ are categorical fibrations and that
$\calX' = \calX \times_{\calY} \calY'$. Let $\calI^{\triangleright}_{X'/}$ denote the fiber product
$\calI^{\triangleright} \times_{\calX'} \calX_{X'/}$, and define
$\calI^{\triangleright}_{X/}$, $\calI^{\triangleright}_{Y'/}$, and $\calI^{\triangleright}_{Y/}$ similarly. We have a pullback diagram
$$ \xymatrix{ \calI^{\triangleright}_{X'/} \ar[r] \ar[d] & \calI^{\triangleright}_{X/} \ar[d] \\
\calI^{\triangleright}_{Y'/} \ar[r] & \calI^{\triangleright}_{Y/} }$$
of left fibrations over $\calI^{\triangleright}$. Proposition \ref{sharpen} implies that every arrow in this diagram is a left fibration, so that Corollary \ref{ruy} implies that $\calI^{\triangleright}_{X'/}$
is a homotopy fiber product of $\calI^{\triangleright}_{X/}$ with $\calI^{\triangleright}_{Y'/}$ over
$\calI^{\triangleright}_{Y/}$ in the covariant model category $(\sSet)_{/\calI^{\triangleright}}$. 
Let $G: (\sSet)^{A \cup \{ \infty \} } \rightarrow (\sSet)_{\calI^{\triangleright}}$ denote the 
unstraightening functor of \S \ref{contrasec}. Since $G$ is the right Quillen functor of a Quillen equivalence, the above diagram is weakly equivalent to the image under $G$ of a homotopy pullback diagram
$$ \xymatrix{ F_{X'} \ar[r] \ar[d] & F_{X} \ar[d] \\
F_{Y'} \ar[r] & F_{Y} } $$ 
of (weakly) fibrant objects of $(\sSet)^{A \cup \{ \infty \} }$. Moreover, the simplicial nerve of each $F_{Z}$ can be identified with the composition of $\overline{f}$ with the functor corepresented by $Z$. According to Theorem \ref{colimcomparee}, it will suffice to show that $F_{X'}$ is a homotopy colimit diagram. We now observe that $F_{X}$, $F_{Y'}$, and $F_{Y}$ are homotopy colimit diagrams (since $X$, $Y'$, and $Y$ are assumed to be $\kappa$-compact) and conclude by applying Lemma \ref{filterprod}.
\end{proof}

In some of the arguments below, it will be important to be able to replace colimits of
a diagram $\calJ \rightarrow \calC$ by colimits of some composition
$ \calI \stackrel{f}{\rightarrow} \calJ \rightarrow \calC$.
According to Proposition \ref{gute}, this maneuver is justified provided that $f$ is cofinal. Unfortunately, the class of cofinal morphisms is not sufficiently robust for our purposes. We will therefore introduce a property somewhat stronger than cofinality, which has better stability properties.

\begin{definition}\index{gen}{cofinal!weakly}\index{gen}{weakly cofinal}\index{gen}{$\kappa$-cofinal}
Let $f: \calI \rightarrow \calJ$ be a functor between filtered $\infty$-categories. We will say
that $f$ is {\it weakly cofinal} if, for every object $J \in \calJ$, there exists an object
$I \in \calI$ and a morphism $J \rightarrow f(I)$ in $\calJ$. We will say that $f$ is {\it $\kappa$-cofinal} if, for every diagram $p: K \rightarrow \calI$ where $K$ is $\kappa$-small and weakly contractible, the induced functor $\calI_{p/} \rightarrow \calJ_{f  p/}$ is weakly cofinal.
\end{definition}

\begin{example}\label{easex}
Let $\calI$ be a $\tau$-filtered $\infty$-category, and let $p: K \rightarrow \calI$ be a $\tau$-small diagram. Then the projection $\calI_{p/} \rightarrow \calI$ is $\tau$-cofinal. To prove this, consider
a $\tau$-small diagram $K' \rightarrow \calI_{p/}$ where $K'$ is weakly contractible, corresponding to a map $q: K \star K' \rightarrow \calI$. According to Lemma \ref{chotle2}, the inclusion
$K' \subseteq K \star K'$ is right anodyne, so that the map
$\calI_{q/} \rightarrow \calI_{q|K'/}$ is a trivial fibration (and therefore weakly cofinal).
\end{example}

\begin{lemma}\label{storuse}
Let $A$, $B$, and $C$ be simplicial sets, and suppose that $B$ is weakly contractible. Then the inclusion
$$ (A \star B) \coprod_{B} (B \star C) \subseteq A \star B \star C$$
is a categorical equivalence.
\end{lemma}

\begin{proof}
Let $F(A,B,C) = (A \star B) \coprod_{B} (B \star C)$, and let $G(A,B,C) = A \star B \star C$.
We first observe that both $F$ and $G$ preserve filtered colimits and homotopy pushout squares, separately in each argument. Using standard arguments (see, for example, the proof of Proposition \ref{babyy}), we can reduce to the case where $A$ and $C$ are simplices.

Let us say that a simplicial set $B$ is {\it good} if the inclusion
$F(A,B,C) \subseteq G(A,B,C)$ is a categorical equivalence. We now make the following observations:
\begin{itemize}
\item[$(1)$] Every simplex is good. Unwinding the definitions, this is equivalent to the assertion that
for $0 \leq m \leq n \leq p$, the diagram
$$ \xymatrix{ \Delta^{ \{m, \ldots n \} } \ar@{^{(}->}[r] \ar@{^{(}->}[d] & \Delta^{ \{0, \ldots, n\} } 
\ar@{^{(}->}[d] \\
\Delta^{ \{m, \ldots, p\} } \ar@{^{(}->}[r] & \Delta^{ \{0, \ldots, p\} }}$$ is a homotopy pushout square (with respect to the Joyal model structure). It suffices to check that the equivalent subdiagram
$$ \xymatrix{ \Delta^{ \{m, m+1\}} \coprod_{ \{m+1\} }\ldots \coprod_{ \{n-1\} }
\Delta^{ \{n-1, n\} } \ar@{^{(}->}[r] \ar@{^{(}->}[d] & \Delta^{ \{0,1\} } \coprod_{ \{1\} } \ldots \coprod_{ \{n-1\} }
\Delta^{ \{n-1,,n\} } \ar@{^{(}->}[d] \\
\Delta^{ \{ m, m+1 \}} \coprod_{ \{m+1\} } \ldots \coprod_{ \{n-1\} } \Delta^{ \{n-1, n\} } \ar@{^{(}->}[r] &
\Delta^{ \{0,1\} } \coprod_{ \{1\} } \ldots \coprod_{ \{p-1\} } \Delta^{ \{p-1, p\} } }$$
is a homotopy pushout, which is clear.

\item[$(2)$] Given a pushout diagram of simplicial sets
$$ \xymatrix{ B \ar[r] \ar@{^{(}->}[d] & B' \ar@{^{(}->}[d] \\
B'' \ar[r] & B''' }$$
in which the vertical arrows are cofibrations, if $B$, $B'$, and $B''$ are good, then $B'''$ is good. This follows from the compatibility of the functors $F$ and $G$ with homotopy pushouts in $B$.

\item[$(3)$] Every horn $\Lambda^n_i$ is good. This follows by induction on $n$, using $(1)$ and $(2)$. 

\item[$(4)$] The collection of good simplicial sets is stable under filtered colimits; this follows from the compatibility of $F$ and $G$ with filtered colimits, and the stability of categorical equivalences under filtered colimits.

\item[$(5)$] Every retract of a good simplicial set is good (since the collection of categorical equivalences is stable under the formation of retracts).

\item[$(6)$] If $i: B \rightarrow B'$ is an anodyne map of simplicial sets, and $B$ is good, then $B'$ is good. This follows by combining observations $(1)$ through $(5)$.

\item[$(7)$] If $B$ is weakly contractible, then $B$ is good. To see this, choose a vertex
$b$ of $B$. The simplicial set $\{b\} \simeq \Delta^0$ is good (by $(1)$ ), and the inclusion
$\{b\} \subseteq B$ is anodyne. Now apply $(6)$.
\end{itemize}
\end{proof}

\begin{lemma}\label{wolfpup}
Let $\kappa$ and $\tau$ be regular cardinals, let
$f: \calI \rightarrow \calJ$ be a $\kappa$-cofinal functor between $\tau$-filtered $\infty$-categories, and let $p: K \rightarrow \calJ$ be a $\kappa$-small diagram. Then:
\begin{itemize}
\item[$(1)$] The $\infty$-category $\calI_{p/} = \calI \times_{\calJ} \calJ_{p/}$ is $\tau$-filtered.

\item[$(2)$] The induced functor $\calI_{p/} \rightarrow \calJ_{p/}$ is $\kappa$-cofinal.
\end{itemize}

\end{lemma}

\begin{proof}
We first prove $(1)$. Let $\widetilde{q}: K' \rightarrow \calI_{p/}$ be a $\tau$-small diagram, classifying a compatible pair of maps $q: K' \rightarrow \calI$ and $q': K \star K' \rightarrow \calJ$.
Since $\calI$ is $\tau$-filtered, we can find an extension $\overline{q}: (K')^{\triangleright} \rightarrow \calI$ of $q$. To find a compatible extension of $\widetilde{q}$, it suffices to solve the lifting problem
$$ \xymatrix{ (K \star K') \coprod_{ K'} (K')^{\triangleright} \ar@{^{(}->}[d]^{i} \ar[r] & \calJ \\
(K \star K')^{\triangleright}, \ar@{-->}[ur] & } $$
which is possible since $i$ is a categorical equivalence (Lemma \ref{storuse}) and $\calJ$ is an $\infty$-category.

To prove $(2)$, we consider a map $\widetilde{q}: K' \rightarrow \calI_{p/}$ as above, where
now $K$ is $\kappa$-small and weakly contractible. 
We have a pullback diagram
$$ \xymatrix{ (\calI_{p/})_{\overline{q}/} \ar[r] \ar[d] & \calI_{q/} \ar[d] \\
\calJ_{q'/} \ar[r] & \calJ_{q'|K'/}. }$$
Lemma \ref{chotle2} implies that the inclusion $K' \subseteq K \star K'$ is right anodyne, 
so that the lower horizontal map is a trivial fibration. It follows that the upper horizontal map
is also a trivial fibration. Since $f$ is $\kappa$-cofinal, the right vertical map is weakly cofinal, so that the left vertical map is weakly cofinal as well.
\end{proof}

\begin{lemma}\label{cofinalwolf}
Let $\kappa$ be a regular cardinal, and let $f: \calI \rightarrow \calJ$ be an $\kappa$-cofinal map
of filtered $\infty$-categories. Then $f$ is cofinal.
\end{lemma}

\begin{proof}
According to Theorem \ref{hollowtt}, to prove that $f$ is cofinal it suffices to show that for every
object $J \in \calJ$, the fiber product $\calI_{J/} = \calI \times_{\calJ} \calJ_{J/}$ is weakly contractible. Lemma \ref{wolfpup} asserts that $\calI_{J/}$ is $\kappa$-filtered; now apply Lemma \ref{stull2}.
\end{proof}

\begin{lemma}\label{lemmatp}
Let $\kappa$ be a regular cardinal, let $\calC$ be an $\infty$-category which
admits $\kappa$-filtered colimits, let $\overline{p}: K^{\triangleright} \rightarrow \calC^{\tau}$ be a $\kappa$-small diagram in the $\infty$-category of $\kappa$-compact objects of $\calC$, and let
$p = \overline{p} | K$. Then $\overline{p}$ is a $\kappa$-compact object of $\calC_{p/}$. 
\end{lemma}

\begin{proof}
Let $\overline{p}'$ denote the composition
$$ K \diamond \Delta^0 \rightarrow K^{\triangleright} \stackrel{\overline{p}}{\rightarrow} \calC^{\kappa};$$ it will suffice to prove that $\overline{p}'$ is a $\tau$-compact object
of $\calC^{p/}$. Consider the pullback diagram
$$ \xymatrix{ \calC^{p/} \ar[r] \ar[d] & \Fun(K \times \Delta^1, \calC) \ar[d]^{f} \\
\ast \ar[r]^-{p} & \Fun(K \times \{0\}, \calC). } $$
Corollary \ref{tweezegork} implies that the $f$ is a Cartesian fibration, so we can apply Proposition \ref{basechangefunky} to deduce that the diagram is homotopy Cartesian (with respect to the Joyal model structure). Using Proposition \ref{limiteval}, we deduce that $f$ preserves $\kappa$-filtered colimits, and {\em any} functor $\ast \rightarrow \calD$ preserves filtered colimits (since filtered $\infty$-categories are weakly contractible; see \S \ref{quasilimit7}). Consequently, Lemma \ref{yoris} implies that $\overline{p}'$ is a $\kappa$-compact object of $\calC^{p/}$ provided that its images in $\ast$ and $\Fun(K \times \Delta^1, \calC)$ are $\kappa$-compact. The former condition is obvious, and the latter follows from Proposition \ref{placeabovee}.
\end{proof}

\begin{lemma}\label{sturm}
Let $\calC$ be an $\infty$-category which admits small, $\tau$-filtered colimits, and let
$p: K \rightarrow \calC$ be a small diagram. Then $\calC_{p/}$ admits small, $\tau$-filtered colimits.
\end{lemma}

\begin{proof}
Without loss of generality, we may suppose that $K$ is an $\infty$-category. Let $\calI$ be a $\tau$-filtered $\infty$-category and $q_0: \calI \rightarrow \calC_{p/}$ a diagram, corresponding to
a diagram $q: K \star \calI \rightarrow \calC$. We next observe that $K \star \calI$ is small and $\tau$-filtered, so that $q$ admits a colimit $\overline{q}: (K \star \calI)^{\triangleright}
\rightarrow \calC$. The map $\overline{q}$ can also be identified with a colimit of $q_0$.
\end{proof}

\begin{proposition}\label{accessforwardslice}\index{gen}{undercategory!and compact objects}
Let $\tau \gg \kappa$ be regular cardinals, let $\calC$ be a $\tau$-accessible $\infty$-category, and let $p: K \rightarrow \calC^{\tau}$ be a $\kappa$-small diagram. Then
$\calC_{p/}$ is $\tau$-accessible, and an object of $\calC_{p/}$ is $\tau$-compact if and only if
its image in $\calC$ is $\tau$-compact.
\end{proposition}

\begin{proof}
Let $\calD = \calC_{p/} \times_{\calC} \calC^{\tau}$ be the full subcategory of
$\calC_{p/}$ spanned by those objects whose image in $\calC$ is $\tau$-compact.
Since $\calC_{p/}$ is idempotent complete, and the collection of $\tau$-compact objects
of $\calC$ is stable under the formation of retracts, we conclude that $\calD$ is idempotent complete. We also note that $\calD$ is essentially small; replacing $\calC$ by a minimal model if necessary, we may suppose that $\calD$ is actually small. Proposition \ref{intprop} and
Lemma \ref{sturm} imply that there is an (essentially unique) $\tau$-continuous functor
$F: \Ind_{\tau}(\calD) \rightarrow \calC_{p/}$ such that the composition
$\calD \rightarrow \Ind_{\tau}(\calD) \stackrel{F}{\rightarrow} \calC_{p/}$ is equivalent to the inclusion of $\calD$ in $\calC_{p/}$. To complete the proof, it will suffice to show that
$F$ is an equivalence of $\infty$-categories. According to Proposition \ref{uterr}, it will suffice to show that $\calD$ consists of $\tau$-compact objects of $\calC_{p/}$ and generates
$\calC_{p/}$ under $\tau$-filtered colimits. The first assertion follows from Lemma \ref{lemmatp}. 

To complete the proof, choose an object $\overline{p}: K^{\triangleright} \rightarrow \calC$ of
$\calC_{p/}$, and let $C \in \calC$ denote the image under $\overline{p}$ of the cone point
of $K^{\triangleright}$. Then we may identify $\overline{p}$ with a diagram
$\widetilde{p}: K \rightarrow \calC^{\tau}_{/C}$. Since $\calC$ is $\tau$-accessible, the
$\infty$-category $\calE = \calC^{\tau}_{/C}$ is $\tau$-filtered. It follows that
$\calE_{\widetilde{p}/}$ is $\tau$-filtered and essentially small; to complete the proof, it will suffice to show that the associated map
$$ \calE_{\widetilde{p}/}^{\triangleright} \rightarrow \calC_{p/}$$
is a colimit diagram. Equivalently, we must show that the compositition
$$ K \star \calE_{\widetilde{p}/}^{\triangleright} \stackrel{\theta_0^{\triangleright}}{\rightarrow}
\calE^{\triangleright} \stackrel{\theta_1}{\rightarrow} \calC$$
is a colimit diagram. Since $\theta_1$ is a colimit diagram, it suffices to prove that $\theta_0$ is cofinal. For this, we consider the composition
$$ q: \calE_{\widetilde{p}/} \stackrel{i}{\rightarrow} K \star \calE_{\widetilde{p}/} \stackrel{\theta_0}{\rightarrow} \calE.$$
The $\infty$-category $\calE$ is $\tau$-filtered, so that $\calE_{\widetilde{p}/}$ is also
$\tau$-filtered, and therefore weakly contractible (Lemma \ref{stull2}). It follows that 
$i$ is right anodyne (Lemma \ref{chotle2}), and therefore cofinal. Applying Proposition \ref{cofbasic}, we conclude that $\theta_0$ is cofinal if and only if $q$ is cofinal. We now observe that that $q$ is $\tau$-cofinal (Example \ref{easex}) and therefore cofinal (Lemma \ref{cofinalwolf}).
\end{proof}

\begin{corollary}\label{horsemn}\index{gen}{undercategory!accessibility}\index{gen}{accessible!undercategories}
Let $\calC$ be an accessible $\infty$-category, and let $p:K \rightarrow \calC$ be a diagram indexed by a small simplicial set $K$. Then $\calC_{p/}$ is accessible.
\end{corollary}

\begin{proof}
Choose appropriate cardinals $\tau \gg \kappa$ and apply Proposition \ref{accessforwardslice}.
\end{proof}

\subsection{Accessibility of Fiber Products}\label{accessfiber}

Our goal in this section is to prove that the class of accessible $\infty$-categories is stable under (homotopy) fiber products (Proposition \ref{horse2}). The strategy of proof should now be familiar from \S \ref{accessfunk} and \S \ref{accessprime}. Suppose given a homotopy Cartesian diagram
 $$ \xymatrix{ \calX' \ar[r]^{q'} \ar[d]^{p'} & \calX \ar[d]^{p} \\
\calY' \ar[r]^{q} & \calY }$$
of $\infty$-categories, where $\calX$, $\calY'$, and $\calY$ are accessible $\infty$-categories, and the functors $p$ and $q$ are likewise accessible. If $\kappa$ is a sufficiently large regular cardinal, then we can use Lemma \ref{yoris} to produce a good supply of $\kappa$-compact objects of $\calX'$. Our problem is then to prove that these objects generate $\calX'$ under $\kappa$-filtered colimits. This requires some rather delicate cofinality arguments.

\begin{lemma}\label{supwolf}
Let $\tau \gg \kappa$ be regular cardinals, let $f: \calC \rightarrow \calD$ be a $\tau$-continuous
functor between $\tau$-accessible $\infty$-categories which carries $\tau$-compact objects
of $\calC$ to $\tau$-compact objects of $\calD$. Let $C$ be an object of $\calC$, $\calC^{\tau}_{/C}$ the full subcategory of $\calC_{/C}$ spanned by those objects $C' \rightarrow C$ where
$C'$ is $\tau$-compact, and $\calD^{\tau}_{/f(C)}$ the full subcategory spanned by those
objects $D \rightarrow f(C)$ where $D \in \calD$ is $\tau$-compact. Then $f$ induces
a $\kappa$-cofinal functor $f': \calC^{\tau}_{/C} \rightarrow \calD^{\tau}_{/f(C)}$. 
\end{lemma}

\begin{proof}
Let $\widetilde{p}: K \rightarrow \calC^{\tau}_{/C}$ be a diagram indexed by a $\tau$-small, weakly contractible simplicial set $K$, and let $p: K \rightarrow \calC$ be the underlying map.
We need to show that
the induced functor $(\calC^{\tau}_{/C})_{\widetilde{p}/} \rightarrow (\calD^{\tau}_{/f(C)})_{f'  \widetilde{p}/}$ is weakly cofinal. Using Proposition \ref{accessforwardslice}, we may replace $\calC$ by $\calC_{p/}$ and
$\calD$ by $\calD_{f  p/}$, and thereby reduce to the problem of showing that $f$ is weakly cofinal.
Let $\phi: D \rightarrow f(C)$ be an object of $\calD^{\tau}_{/f(C)}$, and let
$F_{D}: \calD \rightarrow \SSet$ be the functor corepresented by $D$. Since $D$ is
$\tau$-compact, the functor $F_{D}$ is $\tau$-continuous, so that $F_{D} \circ f$
is $\tau$-continuous. Consequently, the space $F_{D}(f(C))$ can be obtained as a colimit
of the $\tau$-filtered diagram
$$p: \calC^{\tau}_{/C} \rightarrow \calD^{\tau}_{/f(C)} \rightarrow \calD \stackrel{F_{D}}{\rightarrow} \SSet.$$
In particular, the path component of $F_{D}(f(C))$ containing $\phi$ lies in the image
of $p(\eta)$, for some $\eta: C' \rightarrow C$ as above. It follows that there exists a commutative diagram
$$ \xymatrix{ D \ar[rr]^{\phi} \ar[dr] & & f(C) \\
& f(C') \ar[ur]^{f(\eta)} & }$$
in $\calD$, which can be identified with a morphism in $\calD^{\tau}_{/f(C)}$ having the desired properties.
\end{proof}

\begin{lemma}\label{pup1}
Let $A = A' \cup \{ \infty\}$ be a linearly ordered set containing a largest element $\infty$, and let
$B \subseteq A'$ be a cofinal subset $($in other words, for every $\alpha \in A'$, there
exists $\beta \in B$ such that $\alpha \leq \beta${}$)$. The inclusion
$$\phi: \Nerve(A') \coprod_{ \Nerve(B) } \Nerve ( B \cup \{ \infty\} ) \subseteq \Nerve(A)$$
is a categorical equivalence.
\end{lemma}

\begin{proof}
For each $\beta \in B$, let $\phi_{\beta}$ denote the inclusion
$$ \Nerve(\{ \alpha \in A': \alpha \leq \beta\}) \coprod_{ \Nerve(\{ \alpha \in B: \alpha \leq \beta \}) }
\Nerve(\{ \alpha \in B: \alpha \leq \beta \} \cup \{ \infty \} ) \subseteq 
\Nerve (\{ \alpha \in A': \alpha \leq \beta \} \cup \{\infty\} ).$$
Since $B$ is cofinal in $A'$, $\phi$ is a filtered colimit of the inclusions 
$\phi_{\beta}$. Replacing $A'$ by $\{ \alpha \in A' : \alpha \leq \beta \}$ and
$B$ by $\{ \alpha \in B: \alpha \leq \beta \}$, we may reduce to the case where
$A'$ has a largest element (which we will continue to denote by $\beta$).

We have a categorical equivalence
$$ \Nerve(B) \coprod_{ \{ \beta \} } \Nerve(\{ \beta, \infty \}) \subseteq \Nerve (B \cup \{\infty\} ).$$
Consequently, to prove that $\phi$ is a categorical equivalence, it will suffice to show that the composition
$$ \Nerve(A') \coprod_{ \{\beta \} } \Nerve(\{ \beta, \infty\}) 
\subseteq \Nerve(A') \coprod_{ \Nerve(B) } \Nerve ( B \cup \{ \infty\} ) \subseteq \Nerve(A)$$
is a categorical equivalence, which is clear. 
\end{proof}

\begin{lemma}\label{remuswolf}
Let $\tau > \kappa$ be regular cardinals, and let
$$ \calX \stackrel{p}{\rightarrow} \calY \stackrel{p'}{\leftarrow} \calX'$$
be functors between $\infty$-categories.
Assume that:
\begin{itemize}
\item[$(1)$] The $\infty$-categories $\calX, \calX'$, and $\calY$ are $\kappa$-filtered, and
admit $\tau$-small, $\kappa$-filtered colimits.
\item[$(2)$] The functors $p$ and $p'$ preserve $\tau$-small, $\kappa$-filtered colimits.
\item[$(3)$] The functors $p$ and $p'$ are $\kappa$-cofinal.
\end{itemize}
Then there exist objects $X \in \calX$, $X' \in \calX'$ such that $p(X)$ and $p'(X')$ are equivalent in $\calY$.
\end{lemma}

\begin{proof}
For every ordinal $\alpha$, we let $[\alpha] = \{ \beta : \beta \leq \alpha \}$ and
$(\alpha) = \{ \beta: \beta < \alpha \}$. 
Let us say that an ordinal $\alpha$ is {\em even} if it is of the form $\lambda + n$, where
$\lambda$ is a limit ordinal and $n$ is an even integer; otherwise we will say that $\alpha$ is {\it odd}. Let $A$ denote the set of all even ordinals smaller than $\kappa$, and $A'$ the set of all odd ordinals smaller than $\kappa$. We regard $A$ and $A'$ as subsets of the linearly ordered set $A \cup A' = (\kappa)$. We will construct a commutative diagram
$$ \xymatrix{ \Nerve(A) \ar[r] \ar[d]^{q} & \Nerve (\kappa) \ar[d]^{Q} & \Nerve(A') \ar[d]^{q'} \ar[l] \\
\calX \ar[r]^{p} & \calY & \calX' \ar[l]^{p'}. }$$
Supposing that this is possible, we choose colimits $X \in \calX$, $X' \in \calX'$, and $Y \in \calY$ for $q$, $q'$, and $Q$, respectively. Since the inclusion $\Nerve(A) \subseteq \Nerve (\kappa)$ is cofinal and $p$-preserves $\kappa$-filtered colimits, we conclude that $p(X)$ and $Y$ are equivalent.
Similarly, $p'(X')$ and $Y$ are equivalent, so that $p(X)$ and $p'(X')$ are equivalent, as desired.

The construction of $q$, $q'$, and $Q$ is given by induction. Let $\alpha < \kappa$, and suppose
that $q| \Nerve(\{ \beta \in A: \beta < \alpha \})$, $q' | \Nerve(\{ \beta \in A': \beta < \alpha \})$
and $Q| \Nerve (\alpha)$ have already been constructed. We will show how to extend the definitions of $q$, $q'$, and $Q$ to include the ordinal $\alpha$. We will suppose that $\alpha$ is even; the case where $\alpha$ is odd is similar (but easier). 

Suppose first that $\alpha$ is a limit ordinal. In this case, define
$q | \Nerve(\{ \beta \in A: \beta \leq \alpha \})$ to be an arbitrary extension of
$q | \Nerve(\{ \beta \in A: \beta < \alpha \})$: such an extension exists in virtue of our assumption that $\calX$ is $\kappa$-filtered. In order to define $Q | \Nerve (\alpha)$
it suffices to verify that $\calY$ has the extension property with respect to the inclusion
$$ \Nerve (\alpha) \coprod_{ \Nerve (\{ \beta \in A: \beta < \alpha \}) }
\Nerve (\{ \beta \in A: \beta \leq \alpha \}) \subseteq \Nerve[\alpha]. $$
Since $\calY$ is an $\infty$-category, this follows immediately from Lemma \ref{pup1}.

We now treat the case where $\alpha = \alpha' + 1$ is a successor ordinal. Let
$q_{< \alpha} = q | \{ \beta \in A: \beta < \alpha\}$, and regard
$Q | \Nerve ( \{ \alpha' \} \cup \{ \beta \in A: \beta < \alpha \} )$ as an object of
$\calY_{f  q_{< \alpha} /}$. We now observe that $\Nerve (\{ \beta \in A: \beta < \alpha\})$ is $\kappa$-small and weakly contractible. Since $p$ is $\kappa$-cofinal, we can construct
$q | \{ \beta \in A: \beta \leq \alpha\}$ extending $q_{< \alpha}$ and a compatible map
$Q | \Nerve ( \{ \alpha' \} \cup \{ \beta \in A: \beta \leq \alpha \} )$. To complete the construction of $Q$, it suffices to show that $\calY$ has the extension property with respect to the inclusion
$$ \Nerve (\alpha)  \coprod_{ \Nerve (\{ \beta \in A: \beta < \alpha \} \cup \{ \alpha' \}) }
\Nerve (\{ \beta \in A: \beta \leq \alpha \} \cup \{\alpha' \}) \subseteq \Nerve [\alpha]. $$ Once again, this follows from Lemma \ref{pup1}.
\end{proof}

\begin{lemma}\label{wolfpup2}
Let $\kappa$ and $\tau$ be regular cardinals, let $f: \calI \rightarrow \calJ$ be a $\kappa$-cofinal functor between $\tau$-filtered $\infty$-categories, and let $p: K \rightarrow \calI$ be a diagram indexed by a $\tau$-small simplicial set $K$. Then the induced functor
$$ \calI_{p/} \rightarrow \calJ_{f  p/}$$ is $\kappa$-cofinal.
\end{lemma}

\begin{proof}
Let $K'$ be a simplicial set which is $\kappa$-small and weakly contractible, and let
$q: K \star K' \rightarrow \calI$ be a diagram. We have a commutative diagram
$$ \xymatrix{ \calI_{q/} \ar[r] \ar[d] & \calJ_{f  q/} \ar[d] \\
\calI_{q|K' / } \ar[r] & \calI_{f  q|K' /}. }$$
Lemma \ref{chotle2} implies that $K' \subseteq K \star K'$ is a right anodyne inclusion, so that the vertical maps are trivial fibrations. Since $f$ is $\kappa$-cofinal, the lower horizontal map is
weakly cofinal; it follows that the upper horizontal map is weakly cofinal as well.
\end{proof}

\begin{lemma}\label{rebuswolf}
Let $\tau > \kappa$ be regular cardinals, and let
$$ \xymatrix{ \calI' \ar[r]^{q'} \ar[d]^{p'} & \calI \ar[d]^{p} \\
\calJ' \ar[r]^{q} & \calJ }$$
a diagram of $\infty$-categories which is homotopy Cartesian $($with respect to the Joyal model structure$)$. Suppose that $\calI$, $\calJ$, and $\calJ'$ are $\tau$-filtered $\infty$-categories which admit $\tau$-small, $\kappa$-filtered colimits. Suppose further that $p$ and $q$ are $\kappa$-cofinal functors which preserve $\tau$-small, $\kappa$-filtered colimits. Then $\calI'$ is $\tau$-filtered, and the functors $p'$ and $q'$ are $\kappa$-cofinal.
\end{lemma}

\begin{proof}
Without loss of generality, we may suppose that $p$ and $q$ are categorical fibrations and that
$\calI' = \calI \times_{\calJ} \calJ'$. To prove that $\calI'$ is $\tau$-filtered, we must show that
$\calI'_{f/}$ is nonempty for every diagram $f: K \rightarrow \calI'$ indexed by a $\tau$-small simplicial set $K$. We have a (homotopy) pullback diagram
$$ \xymatrix{ \calI'_{f/} \ar[r] \ar[d] & \calI_{q'  f/} \ar[d]^{g} \\
\calJ'_{p'  f/} \ar[r]^{h} & \calJ_{p  q'  f/}. }$$
Lemma \ref{forfilt} implies that $\calI_{q'  f/}$, $\calJ'_{p'  f/}$, and
$\calJ_{p  q'  f/}$ are $\tau$-filtered, and Lemma \ref{wolfpup2} implies that
$g$ and $h$ are $\kappa$-cofinal. We may therefore apply Lemma \ref{remuswolf} to deduce
that $\calI'_{f/}$ is nonempty, as desired.

We now prove that $q'$ is $\kappa$-cofinal; the analogous assertion for $p'$ is proven by the same argument. We must show that for every diagram $f: K \rightarrow \calI'$, where $K$ is $\kappa$-small and weakly contractible, the induced map
$\calI'_{f/} \rightarrow \calI_{q'  f/}$ is weakly cofinal. Replacing $\calI'$ by $\calI'_{f/}$ as above, we
may reduce to the problem of showing that $q'$ itself is weakly cofinal. Let $I$ be an object of
$\calI$, let $J = p(I) \in \calJ$, and consider the (homotopy) pullback diagram
$$ \xymatrix{ \calI'_{I/} \ar[r] \ar[d] & \calI_{I/} \ar[d]^{u} \\
\calJ'_{J/} \ar[r]^{v} & \calJ_{J/}. }$$
We wish to show that $\calI'_{I/}$ is nonempty. This follows from Lemma \ref{remuswolf}, since
$u$ and $v$ are $\tau$-cofinal by Lemmas \ref{wolfpup2} and \ref{wolfpup}, respectively.
\end{proof}

\begin{proposition}\label{horse2}\index{gen}{accessible!homotopy fiber products}
Let $$ \xymatrix{ \calX' \ar[r]^{q'} \ar[d]^{p'} & \calX \ar[d]^{p} \\
\calY' \ar[r]^{q} & \calY }$$
be a diagram of $\infty$-categories which is homotopy Cartesian (with respect to the Joyal model structure). Suppose further that $\calX$, $\calY$, and $\calY'$ are accessible, and that
$p$ and $q$ are accessible functors. Then $\calX'$ is accessible. Moreover, for any accessible $\infty$-category $\calC$ and any functor $f: \calC \rightarrow \calX$, $f$ is accessible if and only if the compositions $p' \circ f$ and $q' \circ f$ are accessible. In particular $($taking $f = \id_{\calX}${}$)$, the functors $p'$ and $q'$ are accessible.
\end{proposition}

\begin{proof}
Choose a regular cardinal $\kappa$ such that $\calX$, $\calY'$, and $\calY$ are $\kappa$-accessible. Enlarging $\kappa$ if necessary, we may suppose that $p$ and $q$ are $\kappa$-continuous. It follows from Lemma \ref{bird3} that $\calX'$ admits small $\kappa$-filtered
colimits, and that for any $\kappa' > \kappa$, a functor $f: \calC \rightarrow \calX$ is
$\kappa'$-continuous if and only if $p' \circ f$ and $q' \circ f$ are $\kappa'$-continuous.
This proves the second claim; it now suffices to show that $\calX'$ is accessible. For this, we will use characterization $(3)$ of Proposition \ref{clear}. Without loss of generality, we may suppose
that $p$ and $q$ are categorical fibrations, and that $\calX' = \calX \times_{\calY} \calY'$. It then follows easily that $\calX'$ is locally small. It will therefore suffice to show that there exists
a regular cardinal $\tau$ such that $\calX'$ is generated by a small collection of
$\tau$-compact objects under small, $\tau$-filtered colimits.

Since the $\infty$-categories of $\kappa$-compact objects of $\calX$ and $\calY'$ are essentially small, there exists $\tau > \kappa$ such that $p | \calX^{\kappa} \subseteq
\calY^{\tau}$ and $q| {\calY'}^{\kappa} \subseteq \calY^{\tau}$. Enlarging $\tau$ if necessary, we may suppose that $\tau \gg \kappa$. The proof of Proposition \ref{enacc} shows that
every $\tau$-compact object of $\calX$ can be written as a $\tau$-small, $\kappa$-filtered colimit of objects belonging to $\calX^{\kappa}$. Since $p$ is $\kappa$-continuous, it follows that
$p | \calX^{\tau} \subseteq \calY^{\tau}$ and similarly $q| {\calY'}^{\tau} \subseteq \calY^{\tau}$.
Let $\calX'' = \calX^{\tau} \times_{ \calY^{\tau} } {\calY'}^{\tau}$. Then $\calX''$ is an essentially small, full subcategory of $\calX'$. Lemma \ref{yoris} implies that $\calX''$ consists of $\tau$-compact objects of $\calX'$. To complete the proof, it will suffice to show that $\calX''$ generates $\calX'$ under small $\tau$-filtered colimits.

Let $X' = (X,Y')$ be an object of $\calX'$, and let $Y = pX = qY'$. We have a (homotopy) pullback
diagram
$$ \xymatrix{ \calX''_{/X'} \ar[r]^{f'} \ar[d]^{g'} & \calX^{\tau}_{/X} \ar[d]^{g} \\
{\calY'}^{\tau}_{/Y'} \ar[r]^{f} & \calY^{\tau}_{/Y} }$$
of essentially small $\infty$-categories. Lemma \ref{supwolf} asserts that $f$ and $g$ are $\kappa$-cofinal.
We apply Lemma \ref{rebuswolf} to conclude that $\calX''_{/X'}$ is $\tau$-filtered, and that 
$f'$ and $g'$ are $\kappa$-cofinal. Now consider the diagram
$$ \xymatrix{ (\calX''_{/X'})^{\triangleright} \ar[dd] \ar[dr]^{h} \ar[rr] & & (\calX^{\tau}_{/X})^{\triangleright} \ar[d] \\
& \calX' \ar[r]^{q'} \ar[d]^{p'} & \calX \ar[d]^{p} \\
({\calY'}^{\tau}_{/Y})^{\triangleright} \ar[r]^-{q} & \calY' \ar[r] & \calY. }$$
Lemma \ref{cofinalwolf} allows us to conclude that $f'$ and $g'$ are cofinal, so that
$p' \circ h$ and $q' \circ h$ are colimit diagrams. Lemma \ref{bird3} implies that $h$ is a colimit diagram as well,
so that $X'$ is the colimit of an essentially small, $\tau$-filtered diagram taking values in $\calX''$.
\end{proof}

\begin{corollary}\label{horsemuun}\index{gen}{overcategory!accessible}\index{gen}{accessible!overcategories}
Let $\calC$ be an accessible $\infty$-category, and let $p: K \rightarrow \calC$ be a diagram indexed by a small simplicial set $K$. Then the $\infty$-category $\calC_{/p}$ is accessible. 
\end{corollary}

\begin{proof}
Since the map $\calC_{/p} \rightarrow \calC^{/p}$ is a categorical equivalence, it will suffice to prove that $\calC^{/p}$ is accessible. We have a pullback diagram
$$ \xymatrix{ \calC^{/p} \ar[r] \ar[d] & \Fun(K \times \Delta^1, \calC) \ar[d]^{p} \\
\ast \ar[r]^-{q} & \Fun(K \times \{1\}, \calC) }$$ 
of $\infty$-categories. Since $p$ is a coCartesian fibration, Proposition \ref{basechangefunky} implies that this diagram is homotopy Cartesian. According to Proposition \ref{horse1},
the $\infty$-categories $\calC^{K \times \Delta^1}$ and $\calC^{K \times \{1\} }$ are accessible. Using Proposition \ref{limiteval}, we conclude that for every regular cardinal $\kappa$ such
that $\calC$ admits $\kappa$-filtered colimits, $p$ is $\kappa$-continuous; in particular, $p$ is accessible. Corollary \ref{sloam} implies that $\ast$ is accessible and that $q$ is an accessible functor. Applying Proposition \ref{horse2}, we deduce that $\calC^{/p}$ is accessible.
\end{proof}

\subsection{Applications}\label{accessstable}

In \S \ref{accessfunk} through \S \ref{accessfiber}, we established some of the basic stability properties enjoyed by the class of accessible $\infty$-categories. In this section, we will reap some of the rewards for our hard work.

\begin{lemma}\label{simplehorse}\index{gen}{accessible!coproducts}
Let $\{ \calC_{\alpha} \}_{\alpha \in A}$ be a family of $\infty$-categories indexed by
a small set $A$, and let $\calC = \coprod_{ \alpha \in A} \calC_{\alpha}$ be their coproduct. Then
$\calC$ is an accessible if and only if each $\calC_{\alpha}$ is accessible.
\end{lemma}

\begin{proof}
Immediate from the definitions.
\end{proof}

\begin{lemma}\label{complexhorse}\index{gen}{accessible!products}
Let $\{ \calC_{\alpha} \}_{ \alpha \in A}$ be a family of $\infty$-categories indexed by a small
set $A$, and let $\calC = \prod_{\alpha \in A} \calC_{\alpha}$ be their product. If each
$\calC_{\alpha}$ is accessible, then $\calC$ is accessible. Moreover, if $\calD$ is an accessible $\infty$-category, then a functor $\calD \rightarrow \calC$ is accessible if and only if each of the compositions
$$ \calD \rightarrow \calC \rightarrow \calC_{\alpha}$$
is accessible.
\end{lemma}

\begin{proof}
Let $\calD = \coprod_{\alpha \in A} \calC_{\alpha}$. By Lemma \ref{simplehorse}, $\calD$ is accessible. Let $\Nerve(A)$ denote the constant simplicial set with value $A$.
Proposition \ref{horse1} implies that $\Fun(\Nerve(A), \calD)$ is accessible. We now observe that
$\Fun(\Nerve(A), \calD)$ can be written as a disjoint union of $\calC$ with another $\infty$-category; applying Lemma \ref{simplehorse} again, we deduce that $\calC$ is accessible. The second claim follows immediately from the definitions.
\end{proof}

\begin{proposition}\label{accprop}\index{gen}{limit!of accessible $\infty$-categories}
The $\infty$-category $\Acc$ of accessible $\infty$-categories admits small limits, and the inclusion $i: \Acc \subseteq \widehat{ \Cat}_{\infty}$ preserves small limits.
\end{proposition}

\begin{proof}
By Proposition \ref{alllimits}, it suffices to prove that $\Acc$ admits pullbacks and small products, and that $i$ preserves pullbacks and (small) products. Let $\Acc_{\Delta}$ be the (simplicial) subcategory of $\widehat{\sSet}$ defined as follows:
\begin{itemize}
\item[$(1)$] The objects of $\Acc_{\Delta}$ are the accessible $\infty$-categories.
\item[$(2)$] If $\calC$ and $\calD$ are accessible $\infty$-categories, then
$\bHom_{ \Acc_{\Delta} }( \calC, \calD)$ is the subcategory of $\Fun(\calC,\calD)$ whose
objects are accessible functors, and whose morphisms are {\em equivalences} of functors.
\end{itemize}
The $\infty$-category $\Acc$ is isomorphic to the simplicial nerve
$\sNerve( \Acc_{\Delta} )$. 
In view of Theorem \ref{colimcomparee}, it will suffice to prove that the simplicial category $\Acc_{\Delta}$ admits homotopy fiber products
and (small) homotopy products, and that the inclusion $\Acc_{\Delta} \subseteq
( \widehat{\mSet} )^{\degree}$ preserves
homotopy fiber products and homotopy products. The case of homotopy fiber products follows from Proposition \ref{horse2} and the case of (small) homotopy products follows from Lemma \ref{complexhorse}.
\end{proof}

If $\calC$ is an accessible $\infty$-category, then $\calC$ is the union of full subcategories
$\{ \calC^{\tau} \subseteq \calC \}$, where $\tau$ ranges over all (small) regular cardinals. It seems reasonable to expect that if $\tau$ is sufficiently large, then the properties of $\calC$ are mirrored by properties of $\calC^{\tau}$. The following result provides an illustration of this philosophy:

\begin{proposition}\label{tcoherent}\index{gen}{compact object!limits of}
Let $\calC$ be a $\kappa$-accessible $\infty$-category, and let $\tau \gg \kappa$ be an uncountable regular cardinal such that $\calC^{\kappa}$ is essentially $\tau$-small. Then the full subcategory $\calC^{\tau} \subseteq \calC$ is stable under all $\kappa$-small limits which exist in $\calC$.
\end{proposition}

Before giving the proof, we will need to establish a few lemmas.

\begin{lemma}\label{thinman}
Let $\tau \gg \kappa$ be regular cardinals, and assume $\tau$ is uncountable. 
Let $\calC$ be a $\tau$-small $\infty$-category, and let $D$ be an object of $\Ind_{\kappa}(\calC)$. The following are equivalent:
\begin{itemize}
\item[$(1)$] The object $D$ is $\tau$-compact in $\Ind_{\kappa}(\calC)$. 
\item[$(2)$] For every $C \in \calC$, the space
$\bHom_{\Ind_{\kappa}(\calC)}(j(C),D)$ is essentially $\tau$-small, where
$j: \calC \rightarrow \Ind_{\kappa}(\calC)$ denotes the Yoneda embedding.
\end{itemize}
\end{lemma}

\begin{proof}
Suppose first that $(1)$ is satisfied. Using Lemma \ref{longwait0}, we can write $D$
as the colimit of the $\kappa$-filtered diagram
$$ \calC_{/D} = \calC \times_{ \Ind_{\kappa}(\calC)} \Ind_{\kappa}(\calC)_{/D} \rightarrow \Ind_{\kappa}(\calC).$$
Since $\tau \gg \kappa$, we also write $D$ as a small $\tau$-filtered colimit of objects
$\{ D_{\alpha} \}$, where each $D_{\alpha}$ is the colimit of a $\tau$-small, $\kappa$-filtered
diagram $$ \widetilde{\calC} \rightarrow \calC \rightarrow \Ind_{\kappa}(\calC).$$
Since $D$ is $\tau$-compact, we conclude that $D$ is a retract of $D_{\alpha}$. Let
$F: \Ind_{\kappa}(\calC) \rightarrow \SSet$ denote the functor co-represented by $j(C)$.
According to Proposition \ref{justcut}, $F$ is $\kappa$-continuous. It follows that
$F(D)$ is a retract of $F(D_{\alpha})$, which is itself a $\tau$-small colimit of spaces equivalent
to $\bHom_{\Ind_{\kappa}(\calC)}(j(C), j(C')) \simeq \bHom_{\calC}(C,C')$, which is essentially $\tau$-small by assumption, and therefore a $\tau$-compact object of $\SSet$. It follows that $D$ is also $\tau$-compact object of $\SSet$. 

Now assume $(2)$. Once again, we observe that $D$ can be obtained as the colimit of a diagram
$\calC_{/D} \rightarrow \Ind_{\kappa}(\calC)$. By assumption, $\calC$ is $\tau$-small and the fibers of the right fibration $\calC_{/D} \rightarrow \calC$ are essentially $\tau$-small. Proposition \ref{sumt} implies that $\calC_{/D}$ is essentially $\tau$-small, so that $D$ is a
$\tau$-small colimit of $\kappa$-compact objects of $\Ind_{\kappa}(\calC)$ and therefore $\tau$-compact.
\end{proof}

\begin{lemma}\label{paleman}
Let $\tau \gg \kappa$ be regular cardinals such that $\tau$ is uncountable, and let
$\SSet^{\tau}$ be the full subcategory of $\SSet$ consisting of essentially $\tau$-small spaces.
Then $\SSet^{\tau}$ is stable under $\kappa$-small limits in $\SSet$.
\end{lemma}

\begin{proof}
In view of Proposition \ref{alllimits}, it suffices to prove that $\SSet^{\tau}$ is stable under pullbacks and $\kappa$-small products. Using Theorem \ref{colimcomparee}, it will suffice to show that the full subcategory of $\Kan$ spanned by essentially $\tau$-small spaces is stable under $\kappa$-small products and homotopy fiber products. This follows immediately from characterization $(1)$ given in Proposition \ref{apegrape}.
\end{proof}

\begin{proof}[Proof of Proposition \ref{tcoherent}]
Let $K$ be a $\kappa$-small simplicial set and let $p: K \rightarrow \calC^{\tau}$ be a diagram
which admits a limit $X \in \calC$. We wish to show that $X$ is $\tau$-compact. According to
Lemma \ref{thinman}, it suffices to prove that the space $F(X)$ is essentially $\tau$-small, where
$F: \calC \rightarrow \SSet$ denotes the functor co-represented by a $\kappa$-compact object $C \in \calC$. Since $F$ preserves limits, we note that $F(X)$ is a limit of $F \circ p$. Lemma \ref{thinman} implies that the diagram $F \circ p$ takes values in $\SSet^{\tau} \subseteq \SSet$.
We now conclude by applying Lemma \ref{paleman}.
\end{proof}

We note the following useful criterion for establishing that a functor is accessible:

\begin{proposition}\label{adjoints}\index{gen}{accessible!adjoint functors}
Let $G: \calC \rightarrow \calC'$ be a functor between accessible
$\infty$-categories. If $G$ admits a right or a left adjoint, then $G$ is
accessible.
\end{proposition}

\begin{proof}
If $G$ is a left adjoint, then $G$ commutes with all colimits which exist in $\calC$.
Therefore $G$ is $\kappa$-continuous for any cardinal $\kappa$ having the property that
$\calC$ is $\kappa$-accessible. Let us therefore assume that $G$ is a right
adjoint; choose a left adjoint $F$ for $G$.

Choose a regular cardinal $\kappa$ such that $\calC'$ is $\kappa$-accessible. We may suppose without loss of generality that $\calC' = \Ind_{\kappa} \calD$, where $\calD$ is a small $\infty$-category. Consider the composite functor
$$ \calD \stackrel{j}{\rightarrow} \Ind_{\kappa}(\calD) \stackrel{F}{\rightarrow} \calC.$$
Since $\calD$ is small, there exists a regular cardinal $\tau \gg \kappa$ such that
$\calC$ is $\tau$-accessible and the essential image of $F \circ j$ consists of $\tau$-compact
objects of $\calC$. We will show that $G$ is $\tau$-continuous.

Since $\Ind_{\kappa}(\calD) \subseteq \calP(\calD)$ is stable under small $\tau$-filtered colimits, it will suffice to prove that the composition
$$ G': \calC \stackrel{G}{\rightarrow} \Ind_{\kappa}(\calD) \rightarrow \calP(\calD)$$
is $\tau$-continuous. For each object $D \in \calD$, let $G'_{D}: \calC \rightarrow \hat{\SSet}$
denote the composition of $G'$ with the functor given by evaluation at $D$. According to Proposition \ref{limiteval}, it will suffice to show that each $G'_{D}$ is $\tau$-continuous.
Lemma \ref{repco} implies that $G'_{D}$ is equivalent to the composition of $G$ with
the functor $\calC' \rightarrow \hat{\SSet}$ corepresented by $j(D)$. Since $F$ is left adjoint to $G$, we may identify this with the functor corepresented by $F(j(D))$. Since $F(j(D))$ is $\tau$-compact by construction, this functor is $\tau$-continuous.
\end{proof}

\begin{definition}\label{defaccsub}\index{gen}{accessible!subcategory}
Let $\calC$ be an accessible category. A full subcategory $\calD \subseteq \calC$ is
an {\it accessible subcategory} of $\calC$ if $\calD$ is accessible, and the inclusion
of $\calD$ into $\calC$ is an accessible functor.
\end{definition}

\begin{example}\label{colexam}
Let $\calC$ be an accessible $\infty$-category and $K$ a simplicial set. Suppose that
every diagram $K \rightarrow \calC$ has a limit in $\calC$. Let
$\calD \subseteq \Fun(K^{\triangleleft}, \calC) $ be the full subcategory spanned by the
limit diagrams. Then $\calD$ is equivalent to $\Fun(K,\calC)$, and is therefore accessible (Proposition \ref{horse1}). The inclusion $\calD \subseteq \Fun(K^{\triangleleft}, \calC)$ is a right adjoint, and therefore accessible (Proposition \ref{adjoints}). Thus $\calD$ is an accessible subcategory of
$\Fun(K^{\triangleleft}, \calC)$. 
Similarly, if every diagram $K \rightarrow \calC$
has a colimit, then the full subcategory $\calD' \subseteq 
\Fun(K^{\triangleright}, \calC)$ spanned by the colimit diagrams in an accessible subcategory of $\Fun(K^{\triangleright}, \calC)$. 
\end{example}

\begin{proposition}\label{boundint}
Let $\calC$ be an accessible category, and let $\{ \calD_{\alpha} \subseteq \calC \}_{\alpha \in A}$
be a (small) collection of accessible subcategories of $\calC$. Then
$\bigcap_{\alpha \in A} \calD_{\alpha}$ is an accessible subcategory of $\calC$.
\end{proposition}

\begin{proof}
We have a homotopy Cartesian diagram
$$ \xymatrix{ \bigcap_{\alpha \in A} \calD_{\alpha} \ar[r]^-{i'} \ar[d] & \calC \ar[d]^{f} \\
\prod_{\alpha \in A} \calD_{\alpha} \ar[r]^-{i} & \calC^{A}. }$$
Lemma \ref{complexhorse} implies that $\prod_{\alpha \in A} \calD_{\alpha}$ and
$\calC^{A}$ are accessible, and it is easy to see that $f$ and $i$ are accessible functors. Applying Proposition \ref{horse2}, we conclude that $\bigcap_{\alpha \in A} \calD_{\alpha}$ is accessible, and that $i'$ is an accessible functor, as desired.
\end{proof}

We conclude this chapter by establishing a generalization of Proposition \ref{horse1}.

\begin{proposition}\label{prestorkus}
Let $\calC$ be a subcategory of the $\infty$-category $\widehat{\Cat}_{\infty}$ of $($not necessarily small$)$ $\infty$-categories satisfying the following conditions:
\begin{itemize}
\item[$(a)$] The $\infty$-category $\calC$ admits small limits, and the inclusion
$\calC \subseteq \widehat{\Cat}_{\infty}$ preserves small limits.
\item[$(b)$] If $X$ belongs to $\calC$, then $\Fun(\Delta^1,X)$ belongs to $\calC$.
\item[$(c)$] If $X$ and $Y$ belong to $\calC$, then a functor $X \rightarrow \Fun(\Delta^1,Y)$ is a morphism of $\calC$ if and only if, for every vertex $v$ of $\Delta^1$, the composite functor
$X \rightarrow \Fun(\Delta^1,Y) \rightarrow \Fun( \{v\}, Y) \simeq Y$ is a morphism of $\calC$.
\end{itemize}

Let $p: X \rightarrow S$ be a map of simplicial sets, where $S$ is small.
Assume that:
\begin{itemize}
\item[$(i)$] The map $p$ is a categorical fibration and a locally coCartesian fibration.
\item[$(ii)$] For each vertex $s$ in $S$, the fiber $X_{s}$ belongs to $\calC$.
\item[$(iii)$] For each edge $s \rightarrow s'$ in $S$, the associated functor
$X_{s} \rightarrow X_{s'}$ is a morphism in $\calC$.
\end{itemize}
Let $\calE$ be a set of edges of $S$, and let $Y$ be the full subcategory of
$\bHom_{S}(S, X)$ spanned by those sections $f: S \rightarrow X$ of $p$ which
satisfy the following condition:
\begin{itemize}
\item[$(\ast)$] For every edge $e: \Delta^1 \rightarrow S$ belonging to $\calE$,
$f$ carries $e$ to to a $p_{e}$-coCartesian edge of $\Delta^1 \times_{S} X$, where
$p_{e}: \Delta^1 \times_{S} X \rightarrow \Delta^1$ is the projection.
\end{itemize}
Then $Y$ belongs to to $\calC$. Moreover, if $Z \in \calC$, then a functor
$Z \rightarrow Y$ belongs to $\calC$ if and only if, for every vertex $s$ in $S$, the composite map
$Z \rightarrow Y \rightarrow X_{s}$
belongs to $\calC$.
\end{proposition}

\begin{remark}
Hypotheses $(i)$ through $(iii)$ of Proposition \ref{prestorkus} are satisfied, in particular, if $p: X \rightarrow S$ is a coCartesian fibration classified by a functor $S \rightarrow \calC \subseteq \widehat{\Cat}_{\infty}$.
\end{remark}

\begin{remark}
Hypotheses $(a)$, $(b)$, and $(c)$ of Proposition \ref{prestorkus} are satisfied for the following subcategories $\calC \subseteq \widehat{\Cat}_{\infty}$:
\begin{itemize}
\item Fix a class of simplicial sets $\{ K_{\alpha} \}_{\alpha \in A}$. Then we can take $\calC$ be the subcategory of
$\widehat{\Cat}_{\infty}$ whose objects are $\infty$-categories which admit $K_{\alpha}$-indexed (co)limits, for each $\alpha \in A$, and whose morphisms are functors which preserves $K_{\alpha}$-indexed (co)limits, for each $\alpha \in A$.
\item We can take the objects of $\calC$ to be accessible $\infty$-categories, and the morphisms in $\calC$ to be accessible functors (in view of Propositions \ref{horse1} and \ref{accprop}).
\end{itemize}
We will meet some other examples in \S \ref{c5s6}.
\end{remark}

\begin{remark}
In the situation of Proposition \ref{prestorkus}, we can replace ``coCartesian'' by ``Cartesian'' everywhere to obtain a dual result. This follows by applying Proposition \ref{prestorkus} to the map $X^{op} \rightarrow S^{op}$, after replacing $\calC$ by its preimage under the ``opposition'' involution of $\h{ \widehat{\Cat}_{\infty}}$.
\end{remark}

The proof of Proposition \ref{prestorkus} makes use of the following observation:

\begin{lemma}\label{surgem}
Let $p: \calM \rightarrow \Delta^1$ be a coCartesian fibration, classifying a functor
$F: \calC \rightarrow \calD$, where $\calC = p^{-1} \{0\}$ and $\calD = p^{-1} \{1\}$. 
Let $\calX = \bHom_{\Delta^1}( \Delta^1, \calM)$ be the $\infty$-category of sections of $p$.
Then $\calX$ can be identified with a homotopy limit of the diagram
$$ \calC \stackrel{F}{\rightarrow} \Fun( \{0\}, \calD) \leftarrow \Fun(\Delta^1, \calD).$$
\end{lemma}

\begin{proof}
We first replace the diagram in question by a fibrant one. Let $\calC'$ denote the
$\infty$-category of coCartesian sections of $p$. Then the evaluation map
$e: \calC' \rightarrow \calC$ is a trivial fibration of simplicial sets. Moreover, since $F$ is associated to the correspondence $\calM$, the map $e$ admits a section $s$ such that the composition
$$ \calC \stackrel{s}{\rightarrow} \calC' \rightarrow \calD$$
coincides with $F$. It follows that we have a weak equivalence of diagrams
$$ \xymatrix{ \calC \ar[r]^-{F} \ar[d]^{s} & \Fun( \{0\}, \calD) \ar@{=}[d] & \Fun(\Delta^1, \calD) \ar[l] \ar@{=}[d] \\
\calC' \ar[r]^-{F'} & \Fun( \{0\}, \calD) & \Fun(\Delta^1, \calD) \ar[l] }$$
where $F'$ is given by evaluation at $\{1\}$, and is a categorical fibration. Let
$\calX'$ denote the pullback of the lower diagram, which we can identify with 
the full subcategory of $\bHom_{ \Delta^1}( \Lambda^2_1, \calM )$ spanned
by those functors which carry the first edge of $\Lambda^2_1$ to a coCartesian edge of $\calM$.

Regard $\Delta^2$ as an object of $(\sSet)_{/\Delta^1}$ via the unique retraction
$r: \Delta^2 \rightarrow \Delta^1$ onto the
simplicial subset $\Delta^{ \{0,1\} } \subseteq \Delta^{ \{0,1,2\} }$.
Let $\calX''$ denote the full subcategory of $\bHom_{\Delta^1}(\Delta^2, \calM)$
spanned by those maps $\Delta^2 \rightarrow \calM$ which carry the initial edge of $\Delta^2$ to a $p$-coCartesian edge of $\calM$. 

Let $T$ denote the marked simplicial set whose underlying simplicial set is $\Delta^2$, whose sole nondegenerate marked edge is $\Delta^1 \subseteq \Delta^2$, and let $T' = T \times_{ (\Delta^2)^{\sharp} } ( \Lambda^2_1)^{\sharp}$. Since the opposites of the inclusions
$T' \subseteq T$, $( \Delta^{ \{0,2\} } )^{\flat} \subseteq T$ are marked anodyne, we conclude that the evaluation maps
$$ \calX \leftarrow \calX'' \rightarrow \calX'$$ are trivial fibrations of simplicial sets.
It follows that $\calX$ and $\calX'$ are (canonically) homotopy equivalent, as desired.
\end{proof}

\begin{remark}
In the situation of Lemma \ref{surgem}, the full subcategory of $\calX$ spanned by the
{\em coCartesian} sections of $p$ is equivalent (via evaluation at $\{0\}$) to $\calC$.
\end{remark}

\begin{proof}[Proof of Proposition \ref{prestorkus}]
Let us first suppose that $\calE = \emptyset$.
Let $\sk^n S$ denote the $n$-skeleton of $S$. We observe that
$\bHom_{S}(S,X)$ coincides with the (homotopy) inverse limit
$$ \varprojlim \{ \bHom_{ S }(\sk^n S, X) \}. $$ 
In view of assumption $(a)$, it will suffice to prove the result after replacing $S$ by $\sk^n S$.
In other words, we may reduce to the
case where $S$ is $n$-dimensional. 

We now work by induction on $n$, and observe that
there is a homotopy pushout diagram of simplicial sets
$$ \xymatrix{ S_n \times \bd \Delta^n \ar@{^{(}->}[r] \ar[d] & S_n \times \Delta^n \ar[d] \\
\sk^{n-1} S \ar@{^{(}->}[r] & S. }$$
We therefore obtain a homotopy pullback diagram of $\infty$-categories
$$ \xymatrix{ \bHom_{S}( S, X) \ar[r] \ar[d] & \bHom_{S}( \sk^{n-1} S, X) \ar[d] \\
\bHom_{S}( S_n \times \Delta^n,X) \ar[r] & \bHom_{S}( S_n \times \bd \Delta^n, X). }$$
Invoking assumption $(a)$ again, we are reduced to proving the same result after replacing
$S$ by $\sk^{n-1} S$, $S_n \times \bd \Delta^n$, and $S_{n} \times \Delta^n$. The first two cases follow from the inductive hypothesis; we may therefore assume that $S$ is a disjoint union of copies of $\Delta^n$. Applying $(a)$ once more, we can reduce to the case $S = \Delta^n$.

If $n = 0$, there is nothing to prove. If $n > 1$, then we have a trivial fibration
$$ \bHom_{S}(S, X) \rightarrow \bHom_{S}( \Lambda^n_1, X).$$
Since the horn $\Lambda^n_1$ is of dimension $< n$, we may conclude by applying the inductive hypothesis. We are therefore reduced to the case $S = \Delta^1$.

According to Lemma \ref{surgem}, the $\infty$-category $\bHom_{\Delta^1}(\Delta^1, X)$ can be identified with a homotopy limit of the diagram
$$ X_{ \{ 0\} } \stackrel{F}{\rightarrow} X_{ \{1\} } \leftarrow X_{ \{1\} }^{\Delta^1}.$$
In view of $(a)$, it will suffice to prove that all of the $\infty$-categories and functors in the above diagram belong to $\calC$. This follows immediately from $(b)$ and $(c)$.

We now consider the general case where $\calE$ is not required to be empty. For each
edge $e \in \calE$, let $Y(e)$ denote the full subcategory of $\bHom_{S}(S,X)$ spanned by those sections $f: S \rightarrow X$ which satisfy the condition $(\ast)$ for the edge $e$. We wish to prove:
\begin{itemize}
\item[$(1)$] The intersection $\bigcap_{e \in \calE} Y(e)$ belongs to $\calC$.
\item[$(2)$] If $Z \in \calC$, then a functor $Z \rightarrow \bigcap_{e \in \calE} Y(e)$ is a morphism of $\calC$ if and only if the induced map $Z \rightarrow \bHom_{S}(S,X)$ is a morphism of $\calC$.
\end{itemize}
In view of $(a)$, it will suffice to prove the corresponding results where $\bigcap_{e \in \calE} Y(e)$ is replaced by a single subcategory $Y(e) \subseteq \bHom_{S}(S,X)$. 

Let $e: s \rightarrow s'$ be an edge belonging to $\calE$. Lemma \ref{surgem} implies the existence of a homotopy pullback diagram
We now observe that there is a homotopy pullback diagram 
$$ \xymatrix{ Y(e) \ar[r] \ar[d] & \bHom_{S}(S,X) \ar[d] \\
\Fun'( \Delta^1, X_{s'}) \ar[r] & \Fun( \Delta^1, X_{s'} ), }$$
where $\Fun'( \Delta^1, X_{s'} ) \simeq X_{s'}$ is the full subcategory of $\Fun( \Delta^1, X_{s'} )$ spanned by the equivalences. In view of $(a)$, it suffices to prove the following analogues of $(1)$ and $(2)$:
\begin{itemize}
\item[$(1')$] For each vertex $s' \in S$, the $\infty$-categories
$\Fun'(\Delta^1, X_{s'})$ and $\Fun( \Delta^1, X_{s'})$ belong to $\calC$.
\item[$(2')$] Given an object $Z \in \calC$, a functor 
$Z \rightarrow \Fun'( \Delta^1, X_{s'})$ is a morphism in $\calC$ if and only if the
induced map $Z \rightarrow \Fun( \Delta^1, X_{s'})$ is a morphism of $\calC$.
\end{itemize}
These assertions follow immediately from $(b)$ and $(c)$, respectively.
\end{proof}

\begin{corollary}\label{storkus1}\index{gen}{accessible!$\infty$-category of sections}
Let $p: X \rightarrow S$ be a map of simplicial sets which is a coCartesian fibration $($or a Cartesian fibration$)$. Assume that:

\begin{itemize}
\item[$(1)$] The simplicial set $S$ is small.
\item[$(2)$] For each vertex $s$ of $S$, the $\infty$-category $X_{s} = X \times_{S} \{s\}$ is
accessible.
\item[$(3)$] For each edge $e: s \rightarrow s'$ of $S$, the associated functor
$X_{s} \rightarrow X_{s'}$ $($or $X_{s'} \rightarrow X_{s}${}$)$ is accessible. 
\end{itemize}

Then $\bHom_{S}(S,X)$ is an accessible $\infty$-category. Moreover, if $\calC$ is accessible, then
a functor $$\calC \rightarrow \bHom_{S}(S,X)$$ is accessible if and only if, for every vertex $s$ of $S$, the induced map $\calC \rightarrow X_{s}$ is accessible.
\end{corollary}

\section{Presentable $\infty$-Categories}\label{c5s6}

Our final object of study in this chapter is the theory of {\em presentable} $\infty$-categories.

\begin{definition}\label{presdef}\index{gen}{presentable!$\infty$-category}\index{gen}{$\infty$-category!presentable}
An $\infty$-category $\calC$ is {\it presentable} if $\calC$ is accessible and admits small colimits.
\end{definition}

We will begin in \S \ref{presint} by giving a number of equivalent reformulations of Definition \ref{presdef}. The main result, Theorem \ref{pretop}, is due to Carlos Simpson: an $\infty$-category $\calC$ is presentable if and only if it arises as an (accessible) localization of an $\infty$-category of presheaves. 

Let $\calC$ be an $\infty$-category, and let $F: \calC \rightarrow \SSet^{op}$ be a functor. If $F$ is representable by an object of $\calC$, then $F$ preserves colimits (Proposition \ref{yonedaprop}). In \S \ref{aftt}, we will prove that the converse holds when $\calC$ is presentable. This representability criterion has a number of consequences: it implies that $\calC$ admits (small) limits (Corollary \ref{preslim}), and leads to an $\infty$-categorical analogue of the adjoint functor theorem (Corollary \ref{adjointfunctor}).

In \S \ref{colpres}, we will see that the collection of all presentable $\infty$-categories can be organized into an $\infty$-category $\LPres$. Moreover, we will explain how to compute limits and colimits in $\LPres$. In the course of doing so, we will prove that the class of presentable $\infty$-categories is stable under most of the basic constructions of higher category theory.

In view of Theorem \ref{pretop}, the theory of localizations plays a central role in the study of presentable $\infty$-categories. In \S \ref{invloc}, we will show that the collection of all (accessible) localizations of a presentable $\infty$-category $\calC$ can be parametrized in a very simple way. Moreover, there is a good supply of localizations of $\calC$: given any (small) collection of morphisms $S$ of $\calC$, one can construct a corresponding localization functor
$$\calC \stackrel{L}{\rightarrow} S^{-1} \calC \subseteq \calC,$$
where $S^{-1} \calC$ is a the full subcategory of $\calC$ spanned by the {\it $S$-local} objects.
These ideas are due to Bousfield, who works in the setting of model categories; we will give an exposition here in the language of $\infty$-categories. In \S \ref{factgen2}, we will employ the same techniques to produce examples of factorization systems on the $\infty$-category $\calC$.

Let $\calC$ be an $\infty$-category, and let $C \in \calC$ be an object. We will say that
$C \in \calC$ is {\em discrete} if, for every $D \in \calC$, the nonzero homotopy groups of the mapping space $\bHom_{\calC}(D,C)$ vanish. If we let $\tau_{\leq 0} \calC$ denote the full subcategory of $\calC$ spanned by the discrete objects, then $\tau_{\leq 0} \calC$ is (equivalent to) an ordinary category. If $\calC$ is the $\infty$-category of spaces, then we can identify the discrete objects of $\calC$ with the ordinary category of sets. Moreover, the inclusion $\tau_{\leq 0} \SSet \subseteq \SSet$ has a left adjoint, given by
$$ X \mapsto \pi_0 X.$$
In \S \ref{truncintro}, we will show that the preceding remark generalizes to an arbitrary presentable $\infty$-category $\calC$: the discrete objects of $\calC$ constitute an (accessible) localization of $\calC$. We will also consider a more general condition of $k$-truncatedness (which specializes to the condition of discreteness when $k=0$). The truncation functors which we construct will play an important role throughout \S \ref{chap6}. 

In \S \ref{compactgen}, we will study the theory of {\em compactly generated} $\infty$-categories: $\infty$-categories which are generated (under colimits) by their compact objects. This class of $\infty$-categories includes some of the most important examples, such as $\SSet$ and $\Cat_{\infty}$. 
In fact, the $\infty$-category $\SSet$ satisfies an even stronger condition: it is generated by
compact {\em projective} objects (see Definition \ref{humber}). The presence of enough compact projective objects in an $\infty$-category allows us to construct projective resolutions, which gives rise to the theory of nonabelian homological algebra (or ``homotopical algebra''). We will review the rudiments of this theory in \S \ref{stable11}. Finally, in \S \ref{stable12} we will present the same ideas in a more classical form, following Quillen's manuscript \cite{homotopicalalgebra}. The comparison of these two perspectives is based on a rectification result (Proposition \ref{trent}) which is of some independent interest.

\begin{remark}
We refer the reader to \cite{adamek} for a study of presentability in the setting of ordinary
category theory. Note that \cite{adamek} uses the term {\it locally presentable categories} for what we have chosen to call {\it presentable categories}.
\end{remark}

\setcounter{theorem}{0}
\subsection{Presentability}\label{presint}

Our main goal in this section is to establish the following characterization of presentable $\infty$-categories:

\begin{theorem}[Simpson \cite{simpson}]\label{pretop}\index{gen}{presentable!$\infty$-category}\index{gen}{$\infty$-category!presentable}
Let $\calC$ be an $\infty$-category. The following conditions are
equivalent:

\begin{itemize}
\item[$(1)$] The $\infty$-category $\calC$ is presentable.

\item[$(2)$] The $\infty$-category $\calC$ is accessible, and for every
regular cardinal $\kappa$, the full subcategory $\calC^{\kappa}$ admits $\kappa$-small colimits.

\item[$(3)$] There exists a regular cardinal $\kappa$ such $\calC$ is $\kappa$-accessible
and $\calC^{\kappa}$ admits $\kappa$-small colimits.

\item[$(4)$] There exists a regular cardinal $\kappa$, a small $\infty$-category $\calD$ which
admits $\kappa$-small colimits, and an equivalence $\Ind_{\kappa} \calD \rightarrow \calC$.

\item[$(5)$] There exists a small $\infty$-category $\calD$ such that $\calC$
is an accessible localization of $\calP(\calD)$.

\item[$(6)$] The $\infty$-category $\calC$ is locally small, admits small colimits, and there
exists a regular cardinal $\kappa$ and a (small) set $S$ of $\kappa$-compact objects of $\calC$ such that every object of $\calC$ is a colimit of a small diagram taking values in the full subcategory of $\calC$ spanned by $S$.
\end{itemize}
\end{theorem}

Before giving the proof, we need a few preliminaries remarks. We first observe that condition $(5)$ is potentially ambiguous: it is unclear whether the accessibility hypothesis is on $\calC$ or on the associated localization functor $L: \calP(\calD) \rightarrow \calP(\calD)$. The distinction turns out to be irrelevant, by virtue of the following:

\begin{proposition}\label{accloc}\index{gen}{localization!accessible}\index{gen}{accessible!localization}
Let $\calC$ be an accessible $\infty$-category, and let $L: \calC \rightarrow \calC$ be a functor satisfying the equivalent conditions of Proposition \ref{recloc}. The following conditions are equivalent:
\begin{itemize}
\item[$(1)$] The essential image $L \calC$ of $L$ is accessible.
\item[$(2)$] There exists a localization $f: \calC \rightarrow \calD$, where $\calD$ is accessible,
and an equivalence $L \simeq g \circ f$.
\item[$(3)$] The functor $L$ is accessible $($when regarded as a functor from $\calC$ to itself$)$.
\end{itemize}
\end{proposition}

\begin{proof}
Suppose $(1)$ is satisfied. Then we may take $\calD = L\calC$, $f = L$, and $g$ to be the inclusion $L\calC \subseteq \calC$; this proves $(2)$. If $(2)$ is satisfied, then Proposition \ref{adjoints} shows that $f$ and $g$ are accessible functors, so their composite $g \circ f \simeq L$ is
also accessible; this proves $(3)$. Now suppose that $(3)$ is satisfied. Choose a regular
cardinal $\kappa$ such that $\calC$ is $\kappa$-accessible and $L$ is $\kappa$-continuous.
The full subcategory $\calC^{\kappa}$ consisting of $\kappa$-compact objects of
$\calC$ is essentially small, so there exists a regular cardinal $\tau \gg \kappa$ such that
$LC$ is $\tau$-compact for every $C \in \calC^{\kappa}$. Let $\calC'$ denote the full subcategory of $\calC$ spanned by the colimits of all $\tau$-small, $\kappa$-filtered diagrams in $\calC^{\kappa}$, and let $L\calC'$ denote the essential image of $\calC'$ under $L$. We note that
$L\calC'$ is essentially small. Since $L$ is $\kappa$-continuous, $L \calC$ is stable under
small $\kappa$-filtered colimits in $\calC$. It follows that any $\tau$-compact object
of $\calC$ which belongs to $L \calC$ is also $\tau$-compact when viewed as an object
of $L \calC$, so that $L \calC'$ consists of $\tau$-compact objects of $L \calC$. According to Proposition \ref{clear}, to complete the proof that $L \calC$ is accessible it will suffice to
show that $L \calC'$ generates $L \calC$ under small, $\tau$-filtered colimits.

Let $X$ be an object of $\calC$. Then $X$ can be written as a small $\kappa$-filtered
colimit of objects of $\calC^{\kappa}$. The proof of Proposition \ref{enacc} shows that
we can also write $X$ as the colimit of a small $\tau$-filtered diagram in $\calC'$.
Since $L$ is preserves colimits, it follows that $LX$ can be obtaines as the colimit of a small $\tau$-filtered diagram in $L \calC'$.
\end{proof}

The proof of Theorem \ref{pretop} will require a few easy lemmas:

\begin{lemma}\label{idc}
Let $f: \calC \rightarrow \calD$ be a functor between small $\infty$-category which exhibits
$\calD$ as an idempotent completion of $\calC$, and let $\kappa$ be a regular cardinal. Then
$\Ind_{\kappa}(f): \Ind_{\kappa}(\calC) \rightarrow \Ind_{\kappa}(\calD)$ is an equivalence of $\infty$-categories.
\end{lemma}

\begin{proof}
We first apply Proposition \ref{uterr} to conclude that $\Ind_{\kappa}(f)$ is fully faithful. To prove
that $\Ind_{\kappa}(f)$ is an equivalence, we must show that it generates $\Ind_{\kappa}(\calD)$ under $\kappa$-filtered colimits. Since $\Ind_{\kappa}(\calD)$ is generated under $\kappa$-filtered colimits by the essential image of the Yoneda embedding $j_{\calD}: \calD \rightarrow \Ind_{\kappa}(\calD)$. Let $D$ be an object of $\calD$. Then $D$ is a retract of $f(C)$
for some object $C \in \calC$. Then $j_{\calD}(D)$ is a retract of
$(\Ind_{\kappa}(f) \circ j_{\calC})(C)$. Since $\Ind_{\kappa}(\calC)$ is idempotent complete
(Corollary \ref{swwe}), we conclude that $j_{\calD}(D)$ belongs to the essential image
of $\Ind_{\kappa}(f)$.
\end{proof}

\begin{lemma}\label{easybumb}
Let $F: \calC \rightarrow \calD$ be a functor between $\infty$-categories which admit
small, $\kappa$-filtered colimits, and let $G$ be a right adjoint to $F$. Suppose
that $G$ is $\kappa$-continuous. Then $F$ carries $\kappa$-compact objects of $\calC$ to $\kappa$-compact objects of $\calD$.
\end{lemma}

\begin{proof}
Let $C$ be a $\kappa$-compact object of $\calC$, $e_{C}: \calC \rightarrow \hat{\SSet}$ the functor corepresented by $C$, and $e_{F(C)}: \calD \rightarrow \hat{\SSet}$ the functor corepresented by $F(C)$. Since $G$ is a right adjoint to $F$, we have an equivalence $e_{F(C)} = e_{C} \circ G$.
Since $e_{C}$ and $G$ are both $\kappa$-continuous, $e_{FC}$ is $\kappa$-continuous.
It follows that $F(C)$ is $\kappa$-compact, as desired.
\end{proof}

\begin{proof}[Proof of Theorem \ref{pretop}]
Corollary \ref{tyrmyrr} asserts that the full subcategory $\calC^{\kappa}$ is stable under all
$\kappa$-small colimits which exist in $\calC$. This proves that $(1)$ implies $(2)$.
The implications $(2) \Rightarrow (3) \Rightarrow (4)$ are obvious. We next prove that $(4)$ implies
$(5)$. According to Lemma \ref{idc}, we may suppose without loss of generality that $\calD$ is idempotent complete. Let
$\calP^{\kappa}(\calD)$ denote the full subcategory of $\calP(\calD)$ spanned by the $\kappa$-compact objects, let $\calD'$ be a minimal model for $\calP^{\kappa}(\calD)$, and let $g$ denote the composition
$$ \calD \stackrel{j}{\rightarrow} \calP^{\kappa}(\calD) \rightarrow \calD'$$
where the second map is a homotopy inverse to the inclusion $\calD' \subseteq \calP^{\kappa}(\calD)$. Proposition \ref{fulfaith} implies that $g$ is fully faithful and Proposition \ref{kcolim} implies that $g$ admits a left adjoint $f$. It follows that $F = \Ind_{\kappa}(f)$ and $G = \Ind_{\kappa}(g)$ are adjoint functors, and Proposition \ref{uterr} implies that $G$ is fully faithful. Moreover,
Proposition \ref{precst} implies that $\Ind_{\kappa} \calD'$ is equivalent to $\calP(\calD)$,
so that $\calC$ is equivalent to an accessible localization of $\calP(\calD')$.

We now prove that $(5)$ implies $(6)$. Let $\calD$ be a small $\infty$-category and
$L: \calP(\calD) \rightarrow \calC$ an accessible localization. Remark \ref{localcolim} implies
that $\calC$ admits small colimits and that $\calC$ is generated under colimits by
the essential image of the composition
$$ T: \calD \stackrel{j}{\rightarrow} \calP(\calD) \stackrel{L}{\rightarrow} \calC.$$
To complete the proof of $(6)$, it will suffice to show that there exists a regular cardinal $\kappa$ such that the essential image of $T$ consists of $\kappa$-compact objects.
Let $G$ denote a left adjoint to $L$. By assumption, $G$ is an accessible functor so that there exists a regular cardinal $\kappa$ such that $G$ is $\kappa$-continuous. For each
object $D \in \calD$, the Yoneda image $j(D)$ is a completely compact object of
$\calP(\calD)$, and in particular $\kappa$-compact. Lemma \ref{easybumb} implies that
$T(D)$ is a $\kappa$-compact object of $\calC$. 

We now complete the proof by showing that $(6)
\Rightarrow (1)$.  Assume that there exists a regular cardinal
$\kappa$ and a set $S$ of $\kappa$-compact objects of $\calC$ such
that every object of $\calC$ is a colimit of objects in $S$. Let $\calC' \subseteq \calC$ be the full subcategory of $\calC$ spanned by $S$, and let $\calC'' \subseteq \calC$ be the full subcategory
of $\calC$ spanned by the colimits of all $\kappa$-small diagrams with values in $\calC''$.
Since $\calC'$ is essentially small, there is only a bounded number of such diagrams up to equivalence, so that $\calC''$ is essentially small. Moreover, since every object of $\calC$
is a colimit of a small diagram with values in $\calC'$, the proof of Corollary \ref{uterrr} shows that
every object of $\calC$ can also be obtained as the colimit of a small $\kappa$-filtered diagram with values in $\calC''$. Corollary \ref{tyrmyrr} implies that $\calC''$ consists of $\kappa$-compact objects of $\calC$ (a slightly more refined argument shows that, if $\kappa > \omega$, then $\calC''$ consists of {\em precisely} the $\kappa$-compact objects of $\calC$). We may therefore apply Proposition \ref{clear} to deduce that $\calC$ is accessible.
\end{proof}

\begin{remark}\label{modelcatpresentcat}
The characterization of presentable $\infty$-categories as localizations of presheaf $\infty$-categories was established by Simpson in \cite{simpson} (using a somewhat different language). 
The theory of presentable $\infty$-categories is essentially equivalent to the theory of {\em combinatorial} model categories (see \S \ref{turka} and Proposition \ref{notthereyet}).
Since most of the $\infty$-categories we will meet are presentable, our study
could also be phrased in the language of model categories. However, we will try to avoid this language, since for many purposes the restriction to presentable
$\infty$-categories seems unnatural and is often technically
inconvenient.
\end{remark}

\begin{remark}
Let $\calC$ be a presentable $\infty$-category, and let $\calD$ be an accessible localization of $\calC$. Then $\calD$ is presentable: this follows immediately from characterization $(5)$ of Proposition \ref{pretop}.
\end{remark}

\begin{remark}\label{tensored}
Let $\calC$ be a presentable $\infty$-category. Since $\calC$ admits arbitrary
colimits, it is ``tensored over spaces'', as we explained in 
\S \ref{quasilimit7}. In particular, the homotopy category of $\calC$ is naturally tensored over the homotopy category $\calH$: for each object $C$ of $\calC$ and every simplicial set $S$, there exists an object $C \otimes S$ of $\calC$, well defined up to equivalence, equipped with isomorphisms $$ \bHom_{\calC}(C \otimes S, C') \simeq
\bHom_{\calC}(C,C')^{S}$$
in the homotopy category $\calH$.
\end{remark}

\begin{example}\label{spacesarepresentable}
The $\infty$-category $\SSet$ of spaces is presentable. This follows from characterization
$(1)$ of Theorem \ref{pretop}, since $\SSet$ is accessible (Example \ref{spacesareaccessible}) and admits (small) colimits by Corollary \ref{limitsinmodel}.
\end{example}

According to Theorem \ref{pretop}, if $\calC$ is $\kappa$-accessible, then
$\calC$ admits small colimits if and only if the full subcategory $\calC^{\kappa} \subseteq \calC$
admits $\kappa$-small colimits. Roughly speaking, this is because arbitrary colimits in $\calC$ can be rewritten in terms of $\kappa$-filtered colimits and $\kappa$-small colimits of $\kappa$-compact objects. Our next result is another variation on this idea; it may also be regarded as an 
analogue of Theorem \ref{pretop} (which describes functors, rather than $\infty$-categories):

\begin{proposition}\label{sumatch}
Let $f: \calC \rightarrow \calD$ be a functor between presentable $\infty$-categories. Suppose that
$\calC$ is $\kappa$-accessible. The following conditions are equivalent:
\begin{itemize}
\item[$(1)$] The functor $f$ preserves small colimits.
\item[$(2)$] The functor $f$ is $\kappa$-continuous, and the restriction
$f| \calC^{\kappa}$ preserves $\kappa$-small colimits.
\end{itemize}
\end{proposition}

\begin{proof}
Without loss of generality, we may suppose $\calC = \Ind_{\kappa}( \calC')$, where $\calC'$ is a small, idempotent complete $\infty$-category which admits $\kappa$-small colimits. The proof of Theorem \ref{pretop} shows that the inclusion $\Ind_{\kappa}(\calC') \subseteq \calP(\calC')$ admits a left adjoint $L$.
Let $\alpha: \id_{\calP(\calC')} \rightarrow L$ be a unit for the adjunction, and let
$f': \calC' \rightarrow \calD$ denote the composition of $f$ with the Yoneda embedding
$j: \Ind_{\kappa}(\calC')$. According to Theorem \ref{charpresheaf}, there exists a
colimit-preserving functor $F: \calP(\calC') \rightarrow \calD$ and an equivalence
of $f'$ with $F \circ j$. Proposition \ref{intprop} implies that $f$ and $F| \Ind_{\kappa}(\calC)$
are equivalent; we may therefore assume without loss of generality that $f = F | \Ind_{\kappa}(\calC)$. Let $F' = f \circ L$, so that $\alpha$ induces a natural transformation 
$\beta: F \rightarrow F'$ of functors from $\calP(\calC')$ to $\calD$. We will show that
$\beta$ is an equivalence. Consequently, we deduce that the functor $F'$ is colimit preserving. 
It then follows that $f$ is colimit preserving. To see this, we consider an arbitrary diagram
$p: K \rightarrow \Ind_{\kappa}(\calC')$ and choose a colimit $\overline{p}: K^{\triangleright} \rightarrow \calP(\calC')$. Then $\overline{q} = L \circ \overline{p}$ is a colimit diagram in $\Ind_{\kappa}(\calC')$, and $f \circ \overline{q} = F' \circ \overline{p}$ is a colimit diagram in $\calD$. Since $q = \overline{q} | K$ is equivalent (via $\alpha$) to the original diagram $p$, we conclude that $f$ preserves the colimit of $p$ in $\Ind_{\kappa}(\calC')$, as well.

It remains to prove that $\beta$ is an equivalence of functors. Let $\calE \subseteq \calP(\calC')$ denote the full subcategory spanned by those objects $X \in \calP(\calC')$ for which
$\beta(X): F(X) \rightarrow F'(X)$ is an equivalence in $\calD$. We wish to prove that
$\calE = \calP(\calC')$. Since $F$ and $F'$
are both $\kappa$-continuous functors, $\calE$ is stable under $\kappa$-filtered colimits
in $\calP(\calC')$. It will therefore suffice to prove that $\calE$ contains $\calP^{\kappa}(\calC')$. 

It is clear that $\calE$ contains $\Ind_{\kappa}(\calC')$; in particular, $\calE$ contains
the essential image $\calE'$ of the Yoneda embedding $j: \calC' \rightarrow \calP(\calC')$. 
According to Proposition \ref{charsmallpre}, every object of $\calP^{\kappa}(\calC')$ is a retract of the colimit of a $\kappa$-small diagram $p: K \rightarrow \calE'$. Since $\calC'$ is idempotent complete, we may identify $\calE'$ with the full subcategory of $\Ind_{\kappa}(\calC')$ consisting of $\kappa$-compact objects. In particular, $\calE'$ is stable under $\kappa$-small colimits and retracts in $\Ind_{\kappa}(\calC')$. It follows that $L$ restricts to a functor $L': \calP^{\kappa}(\calC) \rightarrow \calE'$
which preserves $\kappa$-small colimits.

To complete the proof that $\calP^{\kappa}(\calC') \subseteq \calE$, it will suffice to prove that
$F' | \calP^{\kappa}(\calC)$ preserves $\kappa$-small colimits. 
To see this, we write
$F| \calP^{\kappa}(\calC')$ as a composition
$$ \calP^{\kappa}(\calC') \stackrel{L'}{\rightarrow} \calE' \stackrel{F|\calE'}{\rightarrow} \calC,$$
where $L'$ preserves $\kappa$-small colimits (as noted above) and $F|\calE' = f|\calC^{\kappa}$ preserves $\kappa$-small colimits by assumption. 
\end{proof}

\subsection{Representable Functors and the Adjoint Functor Theorem}\label{aftt}

An object $F$ of the $\infty$-category $\calP(\calC)$ of presheaves on $\calC$ is\index{gen}{representable!functor}\index{gen}{functor!representable}
{\it representable} if it lies in the essential image of the Yoneda embedding $j: \calC \rightarrow \calP(\calC)$. If $F: \calC^{op} \rightarrow \SSet$ is representable, then $F$ preserves limits:
this follows from the fact that $F$ is equivalent to the composite map
$$ \calC^{op} \stackrel{j}{\rightarrow} \calP( \calC^{op} ) \rightarrow \SSet $$
where $j$ denotes the Yoneda embedding for $\calC^{op}$ (which is limit-preserving by Proposition \ref{yonedaprop}) and the right map is given by evaluation at $C$ (which is 
limit-preserving by Proposition \ref{limiteval}). If $\calC$ is presentable, then the converse holds.

\begin{lemma}\label{stewgood}
Let $S$ be a small simplicial set, let $f: S \rightarrow \SSet$ be an object
of $\calP(S^{op})$, and let $F: \calP(S^{op}) \rightarrow \hat{\SSet}$ be the functor corepresented by $f$. Then the composition
$$ S \stackrel{j}{\rightarrow} \calP(S^{op}) \stackrel{F}{\rightarrow} \hat{\SSet}$$
is equivalent to $f$.
\end{lemma}

\begin{proof}
According to Corollary \ref{strictify}, we may can choose a (small) fibrant simplicial category $\calC$ and a categorical equivalence $\phi: S \rightarrow \sNerve(\calC^{op})$ such that $f$ is equivalent to the composition of $\psi^{op}$ with the nerve of a simplicial functor $f': \calC \rightarrow \Kan$.
Without loss of generality, we may suppose that $f' \in \Set_{\Delta}^{\calC}$ is projectively cofibrant.
Using Proposition \ref{gumby444}, we have an equivalence of $\infty$-categories
$$ \psi: \sNerve ( \Set_{\Delta}^{\calC} )^{\degree} ) \rightarrow \calP(S).$$
We observe that the composition $F \circ \psi$ can be identified with
the simplicial nerve of the functor $G: (\Set_{\Delta}^{\calC})^{\degree} \rightarrow \Kan$
corepresented by $f'$. The Yoneda embedding factors through $\psi$, via the adjoint of the composition
$$ j': \sCoNerve[S] \rightarrow \calC^{op} \rightarrow (\Set_{\Delta}^{\calC})^{\degree}.$$
It follows that $F \circ j$ can be identified with the adjoint of the composition
$$ \sCoNerve[S] \stackrel{j'}{\rightarrow} ( \Set_{\Delta}^{\calC} )^{\degree} \stackrel{G}{\rightarrow} \Kan.$$
This composition is equal to the functor $f'$, so its simplicial nerve coincides with the original functor $f$.
\end{proof}

\begin{proposition}\label{representable}
Let $\calC$ be a presentable $\infty$-category, and let 
$F: \calC^{op} \rightarrow \SSet$ be a functor. The following are equivalent:
\begin{itemize}
\item[$(1)$] The functor $F$ is representable by an object $C \in \calC$.
\item[$(2)$] The functor $F$ preserves small limits.
\end{itemize}
\end{proposition}

\begin{proof}
The implication $(1) \Rightarrow (2)$ was proven above (for an arbitrary $\infty$-category $\calC$).
For the converse, we first treat the case where $\calC = \calP(\calD)$, for some small $\infty$-category $\calD$. Let $f: \calD^{op} \rightarrow \SSet$ denote the composition
of $F$ with the (opposite) Yoneda embedding $j^{op}: \calD^{op} \rightarrow \calP(\calD)^{op}$, and let $F': \calP(\calD)^{op} \rightarrow \hat{\SSet}$ denote the functor represented by
$f \in \calP(\calD)$. We will prove that $F$ and $F'$ are equivalent. We observe that
$F$ and $F'$ both preserve small limits; consequently, according to Theorem \ref{charpresheaf}, it will suffice to show that the compositions $f = F \circ j^{op}$ and $f' = F' \circ j^{op}$ are equivalent.
This follows immediately from Lemma \ref{stewgood}.

We now consider the case where $\calC$ is an arbitrary presentable $\infty$-category.
According to Theorem \ref{pretop}, we may suppose that $\calC$ is an accessible localization of a presentable $\infty$-category $\calC'$ which has the form $\calP(\calD)$, so that the assertion for $\calC'$ has already been established. Let $L: \calC' \rightarrow \calC$ denote the localization functor. The functor $F \circ L^{op}: (\calC')^{op} \rightarrow \SSet$ preserves small limits, and is therefore representable by an object
$C \in \calC'$. Let $S$ denote the set of all morphisms $\phi$ in $\calC'$ such that
$L(\phi)$ is an equivalence in $\calC$. Without loss of generality, we may identify
$\calC$ with the full subcategory of $\calC'$ consisting of $S$-local objects. By construction,
$C \in \calC'$ is $S$-local and therefore belongs to $\calC$. It follows that
$C$ represents the functor $(F \circ L^{op})| \calC$, which is equivalent to $F$.
\end{proof}

The representability criterion of Proposition \ref{representable}
has many consequences, as now explain.

\begin{lemma}\label{limitscommute}
Let $X$ and $Y$ be simplicial sets, let $\calC$ be an $\infty$-category, and let
$p: X^{\triangleright} \times Y^{\triangleright} \rightarrow \calC$ be a diagram.
Suppose that:
\begin{itemize}
\item[$(1)$] For every vertex $x$ of $X^{\triangleright}$, the associated map
$p_{x}: Y^{\triangleright} \rightarrow \calC$ is a colimit diagram.
\item[$(2)$] For every vertex $y$ of $Y$, the associated map
$p_{y}: X^{\triangleright} \rightarrow \calC$ is a colimit diagram.
\end{itemize}
Let $\infty$ denote the cone point of $Y^{\triangleright}$. Then the restriction
$p_{\infty}: X^{\triangleright} \rightarrow \calC$ is a colimit diagram. 
\end{lemma}

\begin{proof}
Without loss of generality, we can suppose that $X$ and $Y$ are $\infty$-categories.
Since the inclusion $X \times \{\infty\} \subseteq X \times Y^{\triangleright}$ is cofinal, it will suffice to show that the restriction $p| (X \times Y^{\triangleright})^{\triangleright}$ is a colimit diagram.
According to Proposition \ref{stormus} $p|(X \times Y^{\triangleright})$ is a left
Kan extension of $p|(X \times Y)$. By transitivity, it suffices to show that 
$p| (X \times Y)^{\triangleright}$ is a colimit diagram. For this, it will suffice to prove the stronger assertion that $p| (X^{\triangleright} \times Y)^{\triangleright}$ is a left Kan extension of
$p| (X \times Y)$. Since Proposition \ref{stormus} also implies that
$p| (X^{\triangleright} \times Y)$ is a left Kan extension of $p|(X \times Y)$, we may
again apply transitivity and reduce to the problem of showing that
$p| (X^{\triangleright} \times Y)^{\triangleright}$ is a colimit diagram. Let $\infty'$
denote the cone point of $X^{\triangleright}$. Since
the inclusion $\{\infty'\} \times Y \subseteq X^{\triangleright} \times Y$ is cofinal, we
are reduced to proving that $p_{\infty'}: Y^{\triangleright} \rightarrow \calC$ is a colimit diagram, which follows from $(1)$.
\end{proof}

\begin{corollary}\label{preslim}\index{gen}{limit!in a presentable $\infty$-category}
A presentable $\infty$-category $\calC$ admits all (small) limits.
\end{corollary}

\begin{proof}
Let $\hat{\calP}(\calC) = \Fun(\calC^{op}, \hat{ \SSet })$, where $\hat{\SSet}$ denotes the
$\infty$-category of spaces which are not necessarily small, and let $j: \calC \rightarrow \hat{\calP}(\calC)$ be the Yoneda embedding. Since $j$ is fully faithful, it will suffice to show that the essential image of $j$ admits small limits. The $\infty$-category
$\hat{\calP}(\calC)$ admits all small limits (in fact, even limits which are not necessarily small); it therefore suffices to show that the essential image of $j$ is stable under small limits.
This follows immediately from Proposition \ref{representable} and Lemma \ref{limitscommute}.
\end{proof}

\begin{remark}
Let $A$ be a (small) partially ordered set. The $\infty$-category $\Nerve(A)$ is presentable if and only if every subset of $A$ has a least upper bound. Corollary \ref{preslim} can then be regarded as a generalization of the following classical observation: if every subset of $A$ has a least upper bound, then every subset of $A$ has a greatest lower bound (namely, a least upper bound for the collection of all lower bounds). 
\end{remark}

\begin{remark}\label{coten}
Now that we know that every presentable $\infty$-category $\calC$ has arbitrary limits, we can apply
an argument dual to that of Remark \ref{tensored} to show that
$\calC$ is {\it cotensored over $\SSet$}. In other words, for any
$C \in \calC$ and every simplicial set $X$, there exists an object $C^X
\in \calC$ (well defined up to equivalence) together with a collection of
natural isomorphisms
$$\bHom_{\calC}(C',C^X) \simeq
\bHom_{\calC}(C',C)^{X}$$
in the homotopy category $\calH$.
\end{remark}

We can now formulate a ``dual'' version of Proposition \ref{representable}, which requires a slightly stronger hypothesis.

\begin{proposition}\label{representableprime}\index{gen}{corepresentable!functor}\index{gen}{functor!corepresentable}
Let $\calC$ be a presentable $\infty$-category, and let 
$F: \calC \rightarrow \SSet$ be a functor. Then $F$ is {\em corepresentable} by an object
of $\calC$ if and only if $F$ is accessible and preserves small limits.
\end{proposition}

\begin{proof}
The ``only if'' direction is clear, since every object of $\calC$ is $\kappa$-compact for 
$\kappa \gg 0$. We will prove the converse.
Without loss of generality we may suppose that $\calC$ is minimal (this assumption is a technical convenience which will guarantee that various constructions below stay in the world of small $\infty$-categories).
Let $\widetilde{\calC} \rightarrow \calC$ denote the left fibration represented by $F$. 
Choose a regular cardinal $\kappa$ such that $\calC$ is $\kappa$-accessible and $F$ is
$\kappa$-continuous, and let $\widetilde{\calC}^{\kappa}$ denote the fiber product
$\widetilde{\calC} \times_{\calC} \calC^{\kappa}$, where $\calC^{\kappa} \subseteq \calC$ denotes the full subcategory spanned by the $\kappa$-compact objects of $\calC$. The $\infty$-category
$\widetilde{\calC}^{\kappa}$ is small (since $\calC$ is assumed minimal). Corollary \ref{preslim} implies that the diagram $p: \widetilde{\calC}^{\kappa} \rightarrow \calC$ admits a limit
$\overline{p}: (\widetilde{\calC}^{\kappa})^{\triangleleft} \rightarrow \calC$. Since the functor $F$
preserves small limits, Corollary \ref{charspacelimit} implies that there exists a map
$\overline{q}: (\widetilde{\calC}^{\kappa})^{\triangleleft} \rightarrow \widetilde{\calC}$ which extends
the inclusion $q: \widetilde{\calC}^{\kappa} \subseteq \widetilde{\calC}$ and covers $\overline{p}$.
Let $\widetilde{X}_0 \in \widetilde{\calC}$ denote the image of the cone point under $\overline{q}$ and $X_0$ its image in $\calC$. Then $\widetilde{X}_0$ determines a connected component of the space $F(X_0)$. Since $\calC$ is $\kappa$-accessible, we can write
$X_0$ as a $\kappa$-filtered colimit of $\kappa$-compact objects $\{ X_{\alpha} \}$ of $\calC$. Since $F$ is $\kappa$-continuous, there exists a $\kappa$-compact object $X \in \calC$ such that the induced map $F(X) \rightarrow F(X_0)$ has nontrivial image in the connected component classified by $\widetilde{X}_0$. It follows that there exists an object $\widetilde{X} \in \widetilde{\calC}$ lying
over $X_{\alpha}$, and a morphism $f: \widetilde{X} \rightarrow \widetilde{X}_0$ in
$\widetilde{\calC}$. Since $\widetilde{\calC}_{/q} \rightarrow \widetilde{\calC}$ is a right fibration,
we can pull $\overline{q}$ back to obtain a map
$\overline{q}': (\widetilde{\calC}^{\kappa})^{\triangleleft} \rightarrow \widetilde{\calC}$ 
which extends $q$ and carries the cone point to $\widetilde{X}$. It follows that
$\overline{q}'$ factors through $\widetilde{\calC}^{\kappa}$. We have a commutative diagram
$$ \xymatrix{ & \widetilde{\calC} \ar[dr] & \\
\{ \widetilde{X} \} \ar[ur] \ar[rr]^{i} & & \widetilde{\calC}^{\triangleleft} }$$
where $i$ denotes the inclusion of the cone point. The map $i$ is left anodyne, and therefore
a covariant equivalence in $(\sSet)_{/\calC}$. 
It follows that $\widetilde{\calC}^{\kappa}$ is a retract of $\{ X \}$ in
the homotopy category of the covariant model category $(\sSet)_{/\calC^{\kappa}}$. Proposition \ref{othermod} implies that $F | \calC^{\kappa}$ is a retract of the Yoneda image
$j(X)$ in $\calP(\calC^{\kappa})$. Since the $\infty$-category $\calC^{\kappa}$
is idempotent complete and the Yoneda embedding $j: \calC^{\kappa} \rightarrow
\calP(\calC^{\kappa})$ is fully faithful, we deduce that $F| \calC^{\kappa}$ is
equivalent to $j(X')$, where $X' \in \calC^{\kappa}$ is a retract of $X$. Let $F': \calC \rightarrow \SSet$ denote the functor co-represented by $X'$. We note that $F|\calC^{\kappa}$ and
$F'| \calC^{\kappa}$ are equivalent, and that both $F$ and $F'$ are $\kappa$-continuous.
Since $\calC$ is equivalent to $\Ind_{\kappa}(\calC^{\kappa})$, Proposition \ref{intprop} guarantees that $F$ and $F'$ are equivalent, so that $F$ is representable by $X'$.
\end{proof}

\begin{remark}
It is not difficult to adapt our proof of Proposition \ref{representableprime} to obtain an alternative proof of Proposition \ref{representable}.
\end{remark}

From Propositions \ref{representable} and \ref{representableprime} we can deduce a version of the
adjoint functor theorem:

\begin{corollary}[Adjoint Functor Theorem]\label{adjointfunctor}\index{gen}{adjoint functor theorem}
Let $F: \calC \rightarrow \calD$ be functor between presentable $\infty$-categories.
\begin{itemize}
\item[$(1)$] The functor $F$ has a right adjoint if and only if it preserves small colimits.
\item[$(2)$] The functor $F$ has a left adjoint if and only if it is accessible and preserves small limits.
\end{itemize}
\end{corollary}

\begin{proof}
The ``only if'' directions follow from Propositions \ref{adjointcol} and \ref{adjoints}. We now prove the converse direction of $(2)$; the proof of $(1)$ is similar but easier. Suppose that
$F$ is accessible and preserves small limits. Let $F': \calD \rightarrow \SSet$ be a corepresentable functor. Then $F'$ is accessible and preserves small limits, by Proposition \ref{representableprime}. It follows that the composition $F' \circ F: \calC \rightarrow \SSet$ is accessible and preserves small limits. Invoking Proposition \ref{representableprime} again, we deduce that $F' \circ F$ is representable. We now apply Proposition \ref{adjfuncbaby} to deduce that $F$ has a left adjoint.
\end{proof}

\begin{remark}\label{afi}
The proof of $(1)$ in Corollary \ref{adjointfunctor} does not require that $\calD$ is presentable, but only that $\calD$ is (essentially) locally small.
\end{remark}

\subsection{Limits and Colimits of Presentable $\infty$-Categories}\label{colpres}

In this section, we will introduce an $\infty$-category whose objects are presentable $\infty$-categories, and study its properties. In fact, we will introduce two such $\infty$-categories, which are (canonically) anti-equivalent to one another. The basic observation is the following:
given a pair of presentable $\infty$-categories $\calC$ and $\calD$, the proper notion of ``morphism'' between them is a pair of adjoint functors
$$ \Adjoint{F}{\calC}{\calD}{G} $$
Of course, either one of $F$ and $G$ determines the other up to canonical equivalence. We may therefore think of either one as encoding the data of a morphism.

\begin{definition}\index{not}{PresR@$\RPres$}\index{not}{PresL@$\LPres$}
Let $\widehat{\Cat}_{\infty}$ denote the $\infty$-category of (not necessarily small) $\infty$-categories.
We define subcategories $\RPres, \LPres \subseteq \widehat{\Cat}_{\infty}$ as follows:

\begin{itemize}
\item[$(1)$] The objects of both $\RPres$ and $\LPres$ are the presentable $\infty$-categories.
\item[$(2)$] A functor $F: \calC \rightarrow \calD$ between presentable $\infty$-categories is a morphism in $\LPres$ if and only if $F$ preserves small colimits.
\item[$(3)$] A functor $G: \calC \rightarrow \calD$ between presentable $\infty$-categories is a morphism in $\RPres$ if and only if $G$ is accessible and preserves small limits.
\end{itemize}\index{gen}{$\infty$-category!of presentable $\infty$-categories}
\end{definition}

As indicated above, the $\infty$-categories $\RPres$ and $\LPres$ are anti-equivalent to one another. To prove this, it is convenient to introduce the following definition:

\begin{definition}\label{urtus}\index{gen}{presentable!fibration}\index{gen}{fibration!presentable}
A map of simplicial sets $p: X \rightarrow S$ is a {\it presentable fibration}
if it is both a Cartesian fibration and a coCartesian fibration, and each fiber
$X_{s} = X \times_{S} \{s\}$ is a presentable $\infty$-category.
\end{definition}

The following result is simply a reformulation of Corollary \ref{adjointfunctor}:

\begin{proposition}\label{surtog}
\begin{itemize}
\item[$(1)$] Let $p: X \rightarrow S$ be a Cartesian fibration of simplicial sets, classified
by a map $\chi: S^{op} \rightarrow \widehat{\Cat}_{\infty}$. Then $p$ is a presentable fibration if and only if $\chi$ factors through $\RPres \subseteq \widehat{\Cat}_{\infty}$. 

\item[$(2)$] Let $p: X \rightarrow S$ be a coCartesian fibration of simplicial sets, classified by a map $\chi: S \rightarrow \widehat{\Cat}_{\infty}$. Then $p$ is a presentable fibration if and only if
$\chi$ factors through $\LPres \subseteq \widehat{\Cat}_{\infty}$.

\end{itemize}
\end{proposition}

\begin{corollary}\label{warhog}
For every simplicial set $S$, there is a canonical bijection
$$ [ S, \LPres ] \simeq [ S^{op}, \RPres ]$$
where $[S, \calC]$ denotes the collection of equivalence classes of objects of $\Fun(S,\calC)$.
In particular, there is a canonical isomorphism $\LPres \simeq ( \RPres )^{op}$ in the homotopy category of $\infty$-categories.
\end{corollary}

\begin{proof}
According to Proposition \ref{surtog}, both $[S, \LPres]$ and $[S^{op}, \RPres]$ can be identified with the collection of equivalence classes of presentable fibrations $X \rightarrow S$.
\end{proof}

We now commence our study of the $\infty$-category $\LPres$ (or, equivalently, the anti-equivalent $\infty$-category $\RPres$). The next few results express the idea that $\LPres \subseteq \widehat{\Cat}_{\infty}$ is stable under a variety of categorical constructions.

\begin{proposition}\label{complexhorse2}
Let $\{ \calC_{\alpha} \}_{ \alpha \in A}$ be a family of $\infty$-categories indexed by a small
set $A$, and let $\calC = \prod_{\alpha \in A} \calC_{\alpha}$ be their product. If each
$\calC_{\alpha}$ is presentable, then $\calC$ is presentable.
\end{proposition}

\begin{proof}
It follows from Lemma \ref{complexhorse} that $\calC$ is accessible. Let $p: K \rightarrow \calC$ be a diagram indexed by a small simplicial set $K$, corresponding to a family of diagrams
$\{ p_{\alpha}: K \rightarrow \calC_{\alpha} \}_{\alpha \in A}$. Since each $\calC_{\alpha}$ is presentable, each $p_{\alpha}$ has a colimit $\overline{p_{\alpha}}: K^{\triangleright} \rightarrow \calC_{\alpha}$. These colimits determine a map $\overline{p}: K^{\triangleright} \rightarrow \calC$ which is a colimit of $p$.
\end{proof}

\begin{proposition}\label{presexp}\index{gen}{presentable!functor categories}
Let $\calC$ be an presentable $\infty$-category, and let $K$ be a small simplicial set. Then
$\Fun(K,\calC)$ is presentable.
\end{proposition}

\begin{proof}
According to Proposition \ref{horse1}, $\Fun(K,\calC)$ is accessible. It follows from Proposition \ref{limiteval} that if $\calC$ admits small colimits, then $\Fun(K,\calC)$ admits small colimits.
\end{proof}

\begin{remark}
Let $S$ be a (small) simplicial set.
It follows from Example \ref{spacesareaccessible} and Corollary \ref{storum} that $\calP(S)$
is a presentable $\infty$-category. Moreover, Theorem \ref{charpresheaf} has a natural interpretation in the language of presentable $\infty$-categories: informally speaking, it asserts that the construction $$ S \mapsto \calP(S)$$ is left adjoint to the inclusion functor from presentable $\infty$-categories to all $\infty$-categories.
\end{remark}

The following is a variant on Proposition \ref{presexp}:

\begin{proposition}\label{intmap}
Let $\calC$ and $\calD$ be presentable $\infty$-categories. The $\infty$-category
$\LFun(\calC, \calD)$ is presentable.
\end{proposition}

\begin{proof}
Since $\calD$ admits small colimits, the $\infty$-category $\Fun(\calC, \calD)$ admits small colimits (Proposition \ref{limiteval}). Using Lemma \ref{limitscommute}, we conclude that 
$\LFun(\calC, \calD) \subseteq \Fun(\calC, \calD)$ is stable under small colimits. To complete the proof, it will suffice to show that $\LFun(\calC, \calD)$ is accessible. 

Choose a regular cardinal $\kappa$ such that $\calC$ is $\kappa$-accessible, and let
$\calC^{\kappa}$ be the full subcategory of $\calC$ spanned by the $\kappa$-compact objects. 
Propositions \ref{sumatch} and \ref{intprop} imply that the restriction functor $$ \LFun(\calC, \calD) \rightarrow \Fun(\calC^{\kappa}, \calD)$$
is fully faithful, and its essential image is the full subcategory 
$ \calE \subseteq \Fun(\calC^{\kappa}, \calD)$ spanned by those functors which preserve $\kappa$-small colimits. 

Since $\calC^{\kappa}$ is essentially small. the $\infty$-category $\Fun(\calC^{\kappa}, \calD)$
is accessible (Proposition \ref{horse1}). To complete the proof, we will show that $\calE$
is an accessible subcategory of $\Fun(\calC^{\kappa}, \calD)$. For each $\kappa$-small diagram $p: K \rightarrow \calC^{\kappa}$, let $\calE(p)$ denote the full subcategory of
$\Fun(\calC^{\kappa}, \calD)$ which preserve the colimit of $p$. Then $\calE = \bigcap_{p} \calE(p)$, where the intersection is taken over a set of representatives for all equivalence classes of $\kappa$-small diagrams in $\calC^{\kappa}$. According to Proposition \ref{boundint}, it will suffice to show that each $\calE(p)$ is an accessible subcategory of $\Fun(\calC^{\kappa}, \calD)$. 
We now observe that there is a (homotopy) pullback diagram of $\infty$-categories
$$ \xymatrix{ \calE(p) \ar@{^{(}->}[d] \ar[r] & \calE'(p) \ar@{^{(}->}[d] \\
\Fun( \calC^{\kappa}, \calD) \ar[r] & \Fun( K^{\triangleright}, \calD) }$$
where $\calE'$ denotes the full subcategory of $\Fun(K^{\triangleright}, \calD)$ spanned by the colimit diagrams. According to Proposition \ref{horse1}, it will suffice to prove that
$\calE'(p)$ is an accessible subcategory of $\Fun(K^{\triangleright}, \calD)$, which follows from Example \ref{colexam}.
\end{proof}

\begin{remark}
In the situation of Proposition \ref{intmap}, the presentable $\infty$-category
$\LFun(\calC, \calD)$ can be regarded as an {\em internal mapping object} in
$\LPres$. For every presentable $\infty$-category $\calC'$, a colimit-preserving functor
$\calC' \rightarrow \LFun(\calC, \calD)$ can be identified with a bifunctor
$\calC \times \calC' \rightarrow \calD$, which is colimit-preserving separately in each variable. 
There exists a universal recipient for such a bifunctor: a presentable category which we may denote by $\calC \otimes \calC'$. The operation $\otimes$ endows $\LPres$ with the structure of a {\it symmetric monoidal $\infty$-category}. Proposition \ref{intmap} can be interpreted as asserting that this monoidal structure is {\em closed}.
\end{remark}

\begin{proposition}\label{slicstab}\index{gen}{presentable!overcategories}\index{gen}{overcategory!of presentable $\infty$-categories}
Let $\calC$ be an $\infty$-category, and let $p: K \rightarrow \calC$ be a diagram
in $\calC$ indexed by a (small) simplicial set $K$. If $\calC$ is presentable, then the $\infty$-category $\calC_{/p}$ is also presentable.
\end{proposition}

\begin{proof}
According to Corollary \ref{horsemuun}, $\calC_{/p}$ is accessible. The existence of small colimits in $\calC_{/p}$ follows from Proposition \ref{needed17}.
\end{proof}

\begin{proposition}\label{stabslic}\index{gen}{presentable!undercategories}\index{gen}{undercategory!of a presentable $\infty$-categories}
Let $\calC$ be an $\infty$-category, and let $p: K \rightarrow \calC$ be a diagram
in $\calC$ indexed by a small simplicial set $K$. If $\calC$ is presentable, then the $\infty$-category $\calC_{p/}$ is also presentable.
\end{proposition}

\begin{proof}
It follows from Corollary \ref{horsemn} that $\calC_{p/}$ is accessible. It therefore suffices to prove that every diagram $q: K' \rightarrow \calC_{p/}$ has a colimit in $\calC$. We now observe that
$(\calC_{p/})_{q/} \simeq \calC_{q'/}$ where $q': K \star K' \rightarrow \calC$ is the map classified by $q$. Since $\calC$ admits small colimits, $\calC_{q'/}$ has an initial object.
\end{proof}

\begin{proposition}\label{horse22}
Let $$ \xymatrix{ \calX' \ar[r]^{q'} \ar[d]^{p'} & \calX \ar[d]^{p} \\
\calY' \ar[r]^{q} & \calY }$$
be a diagram of $\infty$-categories which is homotopy Cartesian $($with respect to the Joyal model structure$)$. Suppose further that $\calX$, $\calY$, and $\calY'$ are presentable, and that
$p$ and $q$ are presentable functors. Then $\calX'$ is presentable. Moreover, for any presentable $\infty$-category $\calC$ and any functor $f: \calC \rightarrow \calX$, $f$ is presentable if and only if the compositions $p' \circ f$ and $q' \circ f$ are presentable. In particular $($taking $f = \id_{\calX}${}$)$, $p'$ and $q'$ are presentable functors.
\end{proposition}

\begin{proof}
Proposition \ref{horse2} implies that $\calX'$ is accessible. It therefore suffices to prove that
any diagram $f: K \rightarrow \calX'$ indexed by a small simplicial set $K$ has a colimit in $\calX'$. 
Without loss of generality, we may suppose that $p$ and $q$ are categorical fibrations, and that
$\calX' = \calX \times_{ \calY} \calY'$. Let $X$ be an initial object of $\calX_{q' \circ f/}$ and let
$Y'$ be an initial object of $\calY'_{p' f/}$. Since $p$ and $q$ preserve colimits, 
the images $p(X)$ and $q(Y')$ are initial objects in $\calY_{p q' f/}$, and therefore equivalent to one another. Choose an equivalence $\eta: p(X) \rightarrow q(Y')$. Since $q$ is a categorical fibration, $\eta$ lifts to an equivalence $\overline{\eta}: Y \rightarrow Y'$ in $\calY'_{p' f/}$ such that $q(\overline{\eta}) = \eta$. Replacing $Y'$ by $Y$, we may suppose that
$p(X) = q(Y)$ so that the pair $(X,Y)$ may be considered as an object of
$ \calX'_{f/} = \calY'_{ p' f/} \times_{ \calY_{p q f/}} \calX_{q f/}$. 
According to Lemma \ref{bird1}, it is an initial object of $\calX'_{f/}$, so that $f$ has a colimit in $\calX'$. This completes the proof that $\calX'$ is accessible. The second assertion follows immediately from Lemma \ref{bird3}. 
\end{proof}

\begin{proposition}\label{limitpres1}\index{gen}{limit!of presentable $\infty$-categories}
The $\infty$-category $\LPres$ admits all small limits, and the inclusion functor
$\LPres \subseteq \widehat{\Cat}_{\infty}$ preserves all small limits.
\end{proposition}

\begin{proof}
The proof of Proposition \ref{alllimits} shows that it will suffice to consider the case of pullbacks and small products. The desired result now follows by combining Propositions \ref{horse22} and
\ref{complexhorse2}.
\end{proof}

\begin{corollary}\label{kevinet}
Let $p: X \rightarrow S$ be a presentable fibration of simplicial sets, where $S$ is small. Then the $\infty$-category $\calC$ of {\em coCartesian} sections of $p$ is presentable.
\end{corollary}

\begin{proof}
According to Proposition \ref{surtog}, $p$ is classified by a functor
$\chi: S \rightarrow \LPres$. Using Proposition \ref{limitpres1}, we deduce that
the limit of the composite diagram
$$ S \rightarrow \LPres \rightarrow \widehat{\Cat}_{\infty}$$
is presentable. Corollary \ref{blurt} allows us to identify this limit with the $\infty$-category $\calC$.
\end{proof}

Our goal, in the remainder of this section, is to prove the analogue of Proposition \ref{limitpres1} for the $\infty$-category $\RPres$ (which will show that $\LPres$ is equipped with all small {\em colimits} as well as all small limits). The main step is to prove that for every small diagram
$S \rightarrow \RPres$, the limit of the composite functor
$$ S \rightarrow \RPres \rightarrow \widehat{\Cat}_{\infty}$$
is presentable. As in the proof of Corollary \ref{kevinet}, this is equivalent to the assertion
that for any presentable fibration $p: X \rightarrow S$, the $\infty$-category $\calC$ of {\em Cartesian} sections of $p$ is presentable. To prove this, we will show that the $\infty$-category $\bHom_{S}(S,X)$ is presentable, and that $\calC$ is an accessible localization of $\bHom_{S}(S,X)$.

\begin{lemma}\label{doghold}
Let $p: \calM \rightarrow \Delta^1$ be a Cartesian fibration, let $\calC$
denote the $\infty$-category of sections of $p$, and let
$e: X \rightarrow Y$ and $e': X' \rightarrow Y'$ be objects of $\calC$. If
$e'$ is $p$-Cartesian, then the evaluation map
$\bHom_{\calC}( e,e') \rightarrow \bHom_{\calM}(Y,Y')$
is a homotopy equivalence.
\end{lemma}

\begin{proof}
There is a homotopy pullback diagram of simplicial sets, whose image in the homotopy
category $\calH$ is isomorphic to
$$ \xymatrix{ \bHom_{\calC}(e,e') \ar[r] \ar[d] & \bHom_{\calM}(Y,Y') \ar[d] \\
\bHom_{\calM}(X,X') \ar[r] & \bHom_{\calM}(X,Y'). }$$
If $e'$ is $p$-Cartesian, then the lower horizonal map is a homotopy equivalence, so the upper horizonal map is a homotopy equivalence as well.
\end{proof}

\begin{lemma}\label{cathold}
Let $p: \calM \rightarrow \Delta^1$ be a Cartesian fibration. Let
$\calC$ denote the $\infty$-category of sections of $p$, and $\calC' \subseteq \calC$ the full subcategory spanned by Cartesian sections of $p$. Then $\calC'$ is a reflective subcategory of $\calC$.
\end{lemma}

\begin{proof}
Let $e: X \rightarrow Y$ be an arbitrary section of $p$, and choose a Cartesian section
$e': X' \rightarrow Y$ with the same target. Since $e'$ is Cartesian, there exists a diagram
$$ \xymatrix{ X \ar[r] \ar[d] & Y \ar[d]^{\id_{Y}} \\
X' \ar[r] & Y }$$
in $\calM$, which we may regard as a morphism $\phi$ from $e \in \calC$ to $e' \in \calC'$. 
In view of Proposition \ref{testreflect}, it will suffice to show that $\phi$ exhibits
$e'$ as a $\calC'$-localization of $\calC$. In other words, we must show that for {\em any}
Cartesian section $e'': X'' \rightarrow Y''$, composition with $\phi$ induces a
homotopy equivalence
$\bHom_{\calC}( e', e'') \rightarrow \bHom_{\calC}(e,e'').$
This follows immediately from Lemma \ref{doghold}.
\end{proof}

\begin{proposition}\label{seccatdog}
Let $p: X \rightarrow S$ be a presentable fibration, where $S$ is a small simplicial set. 
Then: 
\begin{itemize}
\item[$(1)$] The $\infty$-category $\calC = \bHom_{S}(S,X)$ of sections of $p$ is presentable.
\item[$(2)$] The full subcategory $\calC' \subseteq \calC$ spanned by Cartesian sections of $p$ is an accessible localization of $\calC$.
\end{itemize}
\end{proposition}

\begin{proof}
The accessibility of $\calC$ follows from Corollary \ref{storkus1}. Since $p$ is a Cartesian fibration and the fibers of $p$ admit small colimits, $\calC$ admits small colimits by Proposition \ref{limiteval}. This proves $(1)$.

For each edge $e$ of $S$, let $\calC(e)$ denote the full subcategory of $\calC$ spanned by those maps $S \rightarrow X$ which carry $e$ to a $p$-Cartesian edge of $X$. By definition,
$\calC' = \bigcap \calC(e)$. According to Lemma \ref{stur3}, it will suffice to show that each
$\calC(e)$ is an accessible localization of $\calC$. Consider the map
$$ \theta_e: \calC \rightarrow \bHom_{S}( \Delta^1, X).$$
Proposition \ref{limiteval} implies that $\theta_{e}$ preserves all limits and colimits. 
Moreover, $\calC(e) = \theta_{e}^{-1} \bHom_{S}'(\Delta^1,X)$, where
$\bHom_{S}'(\Delta^1,X)$ denotes the full subcategory of $\bHom_{S}(\Delta^1,X)$ spanned by $p$-Cartesian edges. According to Lemma \ref{stur2}, it will suffice to show that
$\bHom_{S}'(\Delta^1,X) \subseteq \bHom_{S}(\Delta^1,X)$ is an accessible localization
of $\bHom_{S}(\Delta^1,X)$.
In other words, we may suppose $S = \Delta^1$. It then follows that evaluation at 
$\{1\}$ induces a trivial fibration $\calC' \rightarrow X \times_{S} \{1\}$, so that
$\calC'$ is presentable. It therefore suffices to show that $\calC'$ is a reflective subcategory of $\calC$, which follows from Lemma \ref{cathold}.
\end{proof}

\begin{theorem}\label{surbus}\index{gen}{colimit!of presentable $\infty$-categories}
The $\infty$-category $\RPres$ admits small limits, and the inclusion functor
$\RPres \subseteq \widehat{\Cat}_{\infty}$ preserves small limits.
\end{theorem}

\begin{proof}
Let $\chi: S^{op} \rightarrow \RPres$ be a small diagram, and let
$\overline{\chi}: (S^{\triangleright})^{op} \rightarrow \widehat{\Cat}_{\infty}$
be a limit of $\chi$ in $\Cat_{\infty}$. We will show that $\overline{\chi}$ factors
through $\RPres \subseteq \widehat{\Cat}_{\infty}$ and that $\overline{\chi}$ is a limit when regarded as a diagram in $\RPres$.

We first show that $\overline{\chi}$ carries each vertex to a presentable $\infty$-category.
This is clear with the exception of the cone point of $(S^{\triangleright})^{op}$. Let $p: X \rightarrow S$ be a presentable fibration classified by $\chi$. According to Corollary \ref{blurt}, we may identify the image of the cone point under $\overline{\chi}$ with the $\infty$-category $\calC$ of Cartesian sections of $p$. Proposition \ref{seccatdog} implies that this $\infty$-category is presentable.

We next show that $\overline{\chi}$ carries each edge of $(S^{\triangleright})^{op}$ to
an accessible, limit-preserving functor. This is clear for edges which are degenerate or belong
to $S^{op}$. The remaining edges are in bijection with the vertices of $s$, and connect those vertices to the cone point. The corresponding functors can be identified with the composition
$$ \calC \subseteq \bHom_{S}(S,X) \rightarrow X_{s},$$
where the second functor is given by evaluation at $s$. Proposition \ref{seccatdog} implies
that the inclusion $i: \calC \subseteq \bHom_{S}(S,X)$ is accessible and preserves small limits, and
Proposition \ref{limiteval} implies that the evaluation map $\bHom_{S}(S,X) \rightarrow X_{s}$
preserves all limits and colimits. This completes the proof that $\overline{\chi}$ factors through $\RPres$.

We now show that $\overline{\chi}$ is a limit diagram in $\RPres$. Since $\RPres$ is a subcategory of $\widehat{\Cat}_{\infty}$ and $\overline{\chi}$ is already a limit diagram in $\widehat{\Cat}_{\infty}$, it will suffice to verify the following assertion:

\begin{itemize}
\item If $\calD$ is a presentable $\infty$-category, and $F: \calD \rightarrow \calC$ has the property that each of the composite functors
$$ \calD \stackrel{F}{\rightarrow} \calC \stackrel{i}{\subseteq} \bHom_{S}(S,X) \rightarrow X_{s}$$
is accessible and limit-preserving, then $F$ is accessible and limit-preserving.
\end{itemize}

Applying Proposition \ref{seccatdog}, we see that $F$ is accessible and limit preserving if and only if $i \circ F$ is accessible and limit preserving. We now conclude by applying Proposition \ref{limiteval}.
\end{proof}

\subsection{Local Objects}\label{invloc}

According to Theorem \ref{pretop}, every presentable $\infty$-category arises as an (accessible) localization of some presheaf $\infty$-category $\calP(X)$. Consequently, understanding the process of localization is of paramount importance in the study of presentable $\infty$-categories. In this section, we will classify the accessible localizations of an arbitrary presentable $\infty$-category $\calC$. The basic observation is that a localization functor $L: \calC \rightarrow \calC$ is determined, up to equivalence, by the collection $S$ of all morphisms $f$ such that $Lf$ is an equivalence. Moreover, a collection of morphisms $S$ arises from an accessible localization functor if and only if $S$ is {\em strongly saturated} (Definition \ref{saturated2}) and {\em of small generation} (Remark \ref{sat2}). Given any small collection of morphisms $S$ in $\calC$, 
there is a smallest strongly saturated collection containing $S$: this permits us to define a localization
$S^{-1} \calC \subseteq \calC$. The ideas presented in this section go back (at least) to Bousfield; we refer the reader to \cite{bousfield} for a discussion in a more classical setting.

\begin{definition}\index{gen}{$S$-local}\index{gen}{local!object}
Let $\calC$ be an $\infty$-category and $S$ a collection of morphisms of $\calC$. We say that
an object $Z$ of $\calC$ is {\it $S$-local} if, for every morphism $s: X \rightarrow Y$ belonging to $S$, composition with $s$ induces an isomorphism
$\bHom_{\calC}(Y,Z) \rightarrow \bHom_{\calC}(X,Z)$
in the homotopy category $\calH$ of spaces.

A morphism $f: X \rightarrow Y$ of $\calC$ is an {\it $S$-equivalence} if, for every $S$-local
object $Z$, composition with $f$ induces a homotopy equivalence
$\bHom_{\calC}(Y,Z) \rightarrow \bHom_{\calC}(X,Z)$
is an isomorphism in $\calH$.\index{gen}{$S$-equivalence}
\end{definition}

The following result provides a dictionary for relating localization functors to classes of morphisms:

\begin{proposition}\label{localloc}
Let $\calC$ be an $\infty$-category, and let $L: \calC \rightarrow \calC$ be a localization functor.
Let $S$ denote the collection of all morphisms $f$ in $\calC$ such that $Lf$ is an equivalence.
Then:
\begin{itemize}
\item[$(1)$] An object $C$ of $\calC$ is $S$-local if and only if it belongs to $L\calC$.
\item[$(2)$] Every $S$-equivalence in $\calC$ belongs to $S$.
\item[$(3)$] Suppose that $\calC$ is accessible. The following conditions are equivalent:
\begin{itemize}
\item[$(i)$] The $\infty$-category $L \calC$ is accessible.
\item[$(ii)$] The functor $L: \calC \rightarrow \calC$ is accessible.
\item[$(iii)$] There exists a $($small$)$ subset $S_0 \subseteq S$ such that every $S_0$-local object is $S$-local. 

\end{itemize}

\end{itemize}

\end{proposition}

\begin{proof}[Proof of $(1)$ and $(2)$]
By assumption, $L$ is left adjoint to the inclusion $L \calC \subseteq \calC$; let
$\alpha: \id_{\calC} \rightarrow L$ be a unit map for the adjunction. We begin by proving $(1)$.
Suppose that $X \in L \calC$. Let $f: Y \rightarrow Z$ belong to $S$. Then we have a commutative diagram
$$ \xymatrix{ \bHom_{\calC}(LZ,X) \ar[r] \ar[d] & \bHom_{\calC}(LY, X) \ar[d] \\
\bHom_{\calC}(Z,X) \ar[r] & \bHom_{\calC}(Y,X) }$$
in the homotopy category $\calH$, where the vertical maps are given by composition with $\alpha$ and are homotopy equivalences by assumption. Since $Lf$ is an equivalence, the top horizontal map is also a homotopy equivalence. It follows that the bottom horizontal map is a homotopy equivalence, so that $X$ is $S$-local. Conversely, suppose that $X$ is $S$-local. According to Proposition \ref{recloc}, the map $\alpha(X): X \rightarrow LX$ belongs to $S$, so that composition with $\alpha(X)$ induces a homotopy equivalence $\bHom_{\calC}(LX,X) \rightarrow \bHom_{\calC}(X,X)$. In particular, there exists a map $LX \rightarrow X$ whose composition
with $\alpha(X)$ is homotopic to $\id_{X}$. Thus $X$ is a retract of $LX$. Since $\alpha(LX)$
is an equivalence, we conclude that $\alpha(X)$ is an equivalence, so that $X \simeq LX$ and
therefore $X$ belongs to the essential image of $L$, as desired. This proves $(1)$.

Suppose that $f: X \rightarrow Y$ is an $S$-equivalence. 
We have a commutative diagram $$ \xymatrix{ X \ar[r]^{f} \ar[d]^{\alpha(X)} & Y \ar[d]^{\alpha(Y)} \\
LX \ar[r]^{Lf} & LY }$$
where the vertical maps are $S$-equivalences (by Proposition \ref{recloc}), so that $Lf$
is also an $S$-equivalence. Therefore $LX$ and $LY$ corepresent the same functor
on the homotopy category $\h{L \calC}$. Yoneda's lemma implies that $Lf$ is an equivalence, so that $f \in S$. This proves $(2)$.
\end{proof}

The proof of $(3)$ is somewhat more involved, and will require a few preliminaries.

\begin{lemma}\label{neelda}
Let $\tau \gg \kappa$ be regular cardinals, and suppose that $\tau$ is uncountable. Let $A$ be a $\kappa$-filtered partially ordered set, $A' \subseteq A$ a $\tau$-small subset, and
$$ \{ f_{\gamma}: X_{\gamma} \rightarrow Y_{\gamma} \}_{ \gamma \in C}$$ a $\tau$-small collection of natural transformations of diagrams in $\Kan^{A}$. Suppose that for each
$\alpha \in A$, $\gamma \in C$, the Kan complexes $X_{\gamma}(\alpha)$ and
$Y_{\gamma}(\alpha)$ are essentially $\tau$-small.
Suppose further that, for each $\gamma \in C$, the map of Kan complexes $\colim_{A} f_{\gamma}$
is a homotopy equivalence. Then there exists a $\tau$-small, $\kappa$-filtered subset
$A'' \subseteq A$ such that $A' \subseteq A''$, and $\colim_{A''} f_{\gamma}|A''$ 
is a homotopy equivalence for each $\gamma \in C$.
\end{lemma}

\begin{proof}
Replacing each $f_{\gamma}$ by an equivalent transformation if necessary, we may suppose for each $\gamma \in C$, $\alpha \in A$, the map $f_{\gamma}(\alpha)$ is a Kan fibration.

Let $\alpha \in A$, $\gamma \in C$, and let $\sigma(\alpha,\gamma)$ be a diagram
$$ \xymatrix{ \bd \Delta^n \ar[r] \ar@{^{(}->}[d] & X_{\gamma}(\alpha) \ar[d] \\
\Delta^n \ar[r] & Y_{\gamma}(\alpha). }$$
We will say that $\alpha' \geq \alpha$ {\it trivializes} $\sigma(\alpha,\gamma)$ if the lifting problem
depicted in the induced diagram
$$ \xymatrix{ \bd \Delta^n \ar[r] \ar@{^{(}->}[d] & X_{\gamma}(\alpha') \ar[d] \\
\Delta^n \ar[r] \ar@{-->}[ur] & Y_{\gamma'}(\alpha). }$$
admits a solution. Observe that, if $B \subseteq A$ is filtered, then
$\colim_{B} f_{\gamma}|B$ is a Kan fibration, which is trivial if and only if
for every diagram $\sigma(\alpha,\gamma)$ as above, where $\alpha \in B$, there
exists $\alpha' \in B$ such that $\alpha' \geq \alpha$, and $\alpha'$ trivializes
$\sigma(\alpha, \gamma)$. In particular, since $\colim_{A} f_{\gamma}$ is a homotopy equivalence,
every such diagram $\sigma(\alpha,\gamma)$ is trivialized by some $\alpha' \geq \alpha$. 

We now define a sequence of $\tau$-small subsets $A(\lambda) \subseteq A$, indexed by
ordinals $\lambda \leq \kappa$. Let $A(0) = A'$, and let $A(\lambda) = \bigcup_{\lambda' < \lambda} A(\lambda')$ when $\lambda$ is a limit ordinal. Supposing that $\lambda < \kappa$ and that $A(\lambda)$ has been defined, we choose a set of representatives
$\Sigma = \{ \sigma(\alpha, \gamma) \}$ for all homotopy classes of diagrams as above, where
$\alpha \in A(\lambda)$ and $\gamma \in C$. Since the Kan complexes $X_{\gamma}(\alpha)$,
$Y_{\gamma}(\alpha)$ are essentially $\tau$-small, we may choose the set $\Sigma$ to be $\tau$-small. Each $\sigma \in \Sigma$ is trivialized by some $\alpha'_{\sigma} \in A$; let $B = \{ \alpha'_{\sigma} \}_{\sigma \in \Sigma}$. Then $B$ is $\tau$-small. Now choose a $\tau$-small, $\kappa$-filtered subset $A(\lambda+1) \subseteq A$ containing $A(\lambda) \cup B$ (the existence of $A(\lambda+1)$ is guaranteed by Lemma \ref{estimate}).

We now define $A''$ to be $A(\kappa)$; it is easy to see that $A''$ has the desired properties.
\end{proof}

\begin{lemma}\label{neeld}
Let $\tau \gg \kappa$ be regular cardinals, and suppose that $\tau$ is uncountable. Let $A$ be a $\kappa$-filtered partially ordered set, and for every subset $B \subseteq A$, let
$$ \colim_{B}: \Fun(\Nerve(B), \SSet) \rightarrow \SSet$$
denote a left adjoint to the diagonal functor. Let $A' \subseteq A$ be a $\tau$-small subset
and $\{ f_{\gamma}: X_{\gamma} \rightarrow Y_{\gamma} \}_{\gamma \in C}$ a $\tau$-small collection of morphisms in the $\infty$-category $\Fun(\Nerve(A), \SSet^{\tau})$ of
diagrams $\Nerve(A) \rightarrow \SSet^{\tau}$. Suppose that $\colim_A(f_{\gamma})$ is an
equivalence, for each $\gamma \in C$. Then there exists a $\tau$-small, $\kappa$-filtered subset
$A'' \subseteq A$ which contains $A'$, such that each of the morphisms
$\colim_{A''}(f_{\gamma} | \Nerve(A'' ))$ is an equivalence in $\SSet$.
\end{lemma}

\begin{proof}
Using Proposition \ref{gumby444}, we may assume without loss of generality that
each $f_{\gamma}$ is the simplicial nerve of a natural transformation of functors
from $A$ to $\Kan$. According to Theorem \ref{colimcomparee}, we can identify
$\colim_{B}( X_{\gamma} | \Nerve(B))$ and
$\colim_{B}( Y_{\gamma} | \Nerve(B))$ with homotopy colimits in $\Kan$. If $B$ is filtered, then these homotopy colimits reduce to ordinary colimits (since the class of weak homotopy equivalences in $\Kan$ is stable under filtered colimits), and we may apply Lemma \ref{neelda}.
\end{proof}

\begin{proof}[Proof of part $(3)$ of Proposition \ref{localloc}]
If $L \calC$ is accessible, then Proposition \ref{adjoints} implies that both the inclusion $L \calC \rightarrow \calC$ and $L: \calC \rightarrow L \calC$ are accessible functors, so that their composition is accessible. Thus $(i) \Rightarrow (ii)$. Suppose next that $(ii)$ is satisfied.
Let $\alpha: \id_{\calC} \rightarrow L$ denote a unit for the adjunction between $L$ and
the inclusion $L \calC \subseteq \calC$, and let $\kappa$ be a regular cardinal such that
$\calC$ is $\kappa$-accessible and $L$ is $\kappa$-continuous. Without loss of generality,
we may suppose that $\calC$ is minimal, so that $\calC^{\kappa}$ is a small $\infty$-category.
Let $S_0 = \{ \alpha(X): X \in \calC^{\kappa} \}$, and let
$Y \in \calC$ be $S_0$-local. We wish to prove that $Y$ is $S$-local. Let 
$F_Y: \calC \rightarrow \SSet^{op}$ denote the functor represented by $Y$.
Then $\alpha$ induces a natural transformation $F_Y \rightarrow F_Y \circ L$. The functors $F_Y$ and $F_Y \circ L$ are both $\kappa$-continuous, and by assumption $\alpha$ induces
an equivalence of functors $F_Y | \calC^{\kappa} \rightarrow (F_Y \circ L)| \calC^{\kappa}$ when both sides are restricted to $\kappa$-compact objects. Proposition \ref{intprop} now implies that
$F_Y$ and $F_Y \circ L$ are equivalent, so that $Y$ is $S$-local. This proves $(iii)$.

We complete the proof by showing that $(iii)$ implies $(i)$. Let $\kappa$ be a regular cardinal such that $\calC$ is $\kappa$-accessible and $S_0$ is a set of morphisms between $\kappa$-compact objects of $\calC$. We claim that $L \calC$ is stable under $\kappa$-filtered colimits
in $\calC$. To prove this, let $\overline{p}: K^{\triangleright} \rightarrow \calC$ be a colimit diagram, where $K$ is small and $\kappa$-filtered, and $p = \overline{p}|K$ factors through
$L\calC \subseteq \calC$. Let $s: X \rightarrow Y$ be a morphism which belongs to $S_0$, and let
$s': F_X \rightarrow F_Y$ be the corresponding map of co-representable functors
$\calC \rightarrow \hat{\SSet}$. Since $X$ and $Y$ are $\kappa$-compact by assumption, both
$\overline{p}_{X}: F_X \circ \overline{p}$ and $\overline{p}_Y: F_Y \circ \overline{p}$ are colimit diagrams in $\hat{\SSet}$. The map $s'$ induces a transformation
$\overline{p}_{X} \rightarrow \overline{p}_{Y}$, which is an equivalence when restricted to
$K$, and is therefore an equivalence in general. It follows that
$\bHom_{\calC}(Y, \overline{p}(\infty)) \simeq \bHom_{\calC}(X, \overline{p}(\infty))$, where
$\infty$ denotes the cone point of $K^{\triangleright}$. Thus $\overline{p}(\infty)$
is $S_0$-local as desired.

Now choose an uncountable regular cardinal $\tau \gg \kappa$ such that $\calC^{\kappa}$ is essentially $\tau$-small. According to Proposition \ref{clear}, to complete the proof that $\calC$ is accessible it will suffice to show that $L \calC$ is generated by $\tau$-compact objects under $\tau$-filtered colimits. Let $X$ be an object of $L \calC$. Lemma \ref{longwait0} implies
that $X$ can be written as the colimit of a small diagram $p: \calI \rightarrow \calC^{\kappa}$,
where $\calI$ is $\kappa$-filtered. Using Proposition \ref{rot}, we may suppose that
$\calI$ is the nerve of a $\kappa$-filtered partially ordered set $A$. Let $B$ denote the
collection of all $\kappa$-filtered, $\tau$-small subsets $A_{\beta} \subseteq A$ for which
the colimit of $p| \Nerve(A_{\beta})$ is $S_0$-local. Lemma \ref{neeld} asserts that every $\tau$-small subset of $A$ is contained in $A_{\beta}$, for some $\beta \in B$. It follows that
$B$ is $\tau$-filtered, when regarded as partially ordered by inclusion, and that
$A = \bigcup_{\beta \in B} A_{\beta}$. Using Proposition \ref{utl} and Corollary \ref{util}, we can
obtain $X$ as the colimit of a diagram $q: \Nerve(B) \rightarrow \calC$, where each
$q(\beta)$ is a colimit $X_{\beta}$ of $p | \Nerve(A_{\beta})$. The objects $\{ X_{\beta} \}_{\beta \in B}$ are $S_0$-local and $\tau$-compact by construction.
\end{proof}

According to Proposition \ref{localloc}, every localization $L$ of an $\infty$-category $\calC$ is determined by the class $S$ of morphisms $f$ such that $Lf$ is an equivalence. This raises the question: which classes of morphisms $S$ arise in this way? To answer this question, we will begin by isolating some of the most obvious properties enjoyed by $S$.

\begin{definition}\label{saturated2}\index{gen}{strongly saturated}\index{gen}{saturated!strongly}
Let $\calC$ be a $\infty$-category which admits small colimits, and let $S$ be a collection of morphisms of $\calC$. We will say that $S$ is {\it strongly saturated} if it satisfies the following conditions:

\begin{itemize}
\item[$(1)$] Given a pushout diagram 
$$ \xymatrix{ C \ar[r]^{f} \ar[d] & D \ar[d] \\
C' \ar[r]^{f'} \ar[r] & D' }$$
in $\calC$, if $f$ belongs to $S$, then so does $f'$.

\item[$(2)$] The full subcategory of $\Fun(\Delta^1, \calC)$ spanned by $S$ is stable under small  colimits.

\item[$(3)$] Suppose given a $2$-simplex of $\calC$, corresponding to a diagram
$$ \xymatrix{ X \ar[rr]^{f} \ar[dr]^{g} & & Y \ar[dl]_{h} \\
& Z. & }$$
If any two of $f$, $g$, and $h$ belong to $S$, then so does the third.
\end{itemize}
\end{definition}

\begin{remark}\label{simpcons}
Let $\calC$ be an $\infty$-category which admits small colimits, and let
$S$ be a strongly saturated class of morphisms of $\calC$. Let $\emptyset$ be
an initial object of $\calC$. Condition $(2)$ of Definition
\ref{saturated2} implies that $\id_{\emptyset}: \emptyset \rightarrow \emptyset$ belongs to $S$,
since it is an initial object of $\Fun(\Delta^1,\calC)$. Any equivalence in $\calC$ is a pushout of
$\id_{\emptyset}$, so condition $(1)$ implies that $S$ contains all equivalences in $\calC$.
It also follows from condition $(1)$ that if $f: C \rightarrow D$ belongs to $S$ and
$f': C \rightarrow D$ is homotopic to $f$, then $f'$ belongs to $S$ (since $f'$ is a pushout of $f$).
Note also that condition $(2)$ implies that $S$ is stable under retracts, since any retract
of a morphism $f$ can be written as a colimit of copies of $f$ (Proposition \ref{autokan}).
\end{remark}

\begin{remark}\label{sat2}
Let $\calC$ be an $\infty$-category which admits colimits. Given any collection $\{ S_{\alpha} \}_{\alpha \in A}$ of strongly saturated classes of morphisms of $\calC$, the intersection
$S = \bigcap_{\alpha \in A} S_{\alpha}$ is also strongly saturated. It follows that {\em any} collection $S_0$ of morphisms in $\calC$ is contained in a smallest strongly saturated class of morphisms $S$. In this case we will also write $S = \overline{S_0}$; we refer to it as the strongly saturated class of morphisms {\em generated} by $S_0$. We will say that $S$ is {\it of small generation} if $S = \overline{S_0}$, where $S_0 \subseteq S$ is small.\index{gen}{small generation}
\end{remark}

\begin{remark}
Let $\calC$ be an $\infty$-category which admits small colimits. Let $S$ be a strongly saturated
class of morphisms of $\calC$.
If $f: X \rightarrow Y$ lies in $S$ and $K$ is a simplicial set, then
the induced map $X \otimes K \rightarrow Y \otimes K$ (which is well-defined up to equivalence) lies in $S$. This follows from the closure of $S$ under colimits.
We will use this observation in the proof of Proposition \ref{local}, in the case where $K = \bd \Delta^n$ is a (simplicial) sphere.
\end{remark}

\begin{example}
Let $\calC$ be an $\infty$-category which admits small colimits, and let $S$ denote the class of all equivalences in $\calC$. Then $S$ is strongly saturated; it is clearly the smallest strongly saturated class of morphisms of $\calC$.
\end{example}

\begin{remark}\label{prebluse}
Let $F: \calC' \rightarrow \calC$ be a functor between $\infty$-categories. Suppose that $\calC$ and $\calC'$ admit small colimits and that $F$ preserves small colimits. Let $S$ be a strongly saturated class of morphisms in $\calC'$. Then $F^{-1} S$ is a strongly saturated class of morphisms of $\calC$. In particular, if we let $S$ denote the collection of all morphisms $f$ of $\calC'$ such that $F(f)$ is an equivalence, then $S$ is strongly saturated.
\end{remark}

\begin{lemma}\label{bluh}
Let $\calC$ be an $\infty$-category which admits small colimits, let $S_0$ be a class of morphisms
in $\calC$, and let $S$ denote the collection of all $S$-equivalences. Then $S$ is strongly saturated.
\end{lemma}

\begin{proof}
For each object $X \in \calC$, let $F_X: \calC \rightarrow \SSet^{op}$ denote the functor
represented by $X$, and let $S(X)$ denote the collection of all morphisms $f$ such
that $F_X(f)$ is an equivalence. Since $F_X$ preserves small colimits, Remark \ref{prebluse} implies that $S(X)$ is strongly saturated. We now observe that $S$ is the intersection $\bigcap S(X)$, where $X$ ranges over the class of all $S_0$-local objects of $\calC$.
\end{proof}

\begin{lemma}\label{yorkan}
Let $\calC$ be an $\infty$-category which admits small colimits, let $S$ be a strongly saturated collection of morphisms of $\calC$, and let $C \in \calC$ be an object. Let $\calD \subseteq \calC^{C/}$ be the full subcategory of $\calC^{C/}$ spanned by those objects $C \rightarrow C'$ which belong to $S$. Then $\calD$ is stable under small colimits in $\calC^{C/}$.
\end{lemma}

\begin{proof}
The proofs of Corollary \ref{uterrr} and \ref{allfin} show that it will suffice to prove that $\calD$ is stable under filtered colimits, pushouts, and contains the initial objects of $\calC^{C/}$. The last condition is equivalent to the requirement that $S$ contains all equivalences, which follows from Remark \ref{simpcons}. Now suppose that $\overline{p}: K^{\triangleright} \rightarrow \calC^{C/}$
is a colimit of $p = \overline{p} | K$, where $K$ is either filtered or equivalent to $\Lambda^2_0$,
and that $p(K) \subseteq \calD$. We can identify $\overline{p}$ with a map
$P: K^{\triangleright} \times \Delta^1 \rightarrow \calC$ such that
$P | K^{\triangleright} \times \{0\}$ is the constant map taking the value $C \in \calC$.
Since $K$ is weakly contractible, $P | K^{\triangleright} \times \{0\}$ is a colimit diagram in $\calC$.
The map $P | K^{\triangleright} \times \{1\}$ is the image of a colimit diagram under the left
fibration $\calC^{C/} \rightarrow \calC$; since $K$ is weakly contractible, Proposition \ref{goeselse} implies that $P | K^{\triangleright} \times \{1\}$ is a colimit diagram. We now apply
Proposition \ref{limiteval} to deduce that
$P: K^{\triangleright} \rightarrow \calC^{\Delta^1}$ is a colimit diagram. Since $S$ is stable under colimits in $\calC^{\Delta^1}$, we conclude that $P$
carries the cone point of $K^{\triangleright}$ to a morphism belonging to $S$, as we wished to show.
\end{proof}

\begin{lemma}\label{walter}
Let $\calC$ be an $\infty$-category which admits small filtered colimits, let $\kappa$ be an uncountable regular cardinal, let $A$ and $B$ be $\kappa$-filtered partially ordered sets, and
let $p_0: \Nerve(A_0) \rightarrow \calC^{\kappa}$ and $p_1: A_1 \rightarrow \calC^{\kappa}$ be two diagrams which have the same colimit. Let $A''_0 \subseteq A_0$, $A''_1 \subseteq A_1$ be $\kappa$-small subsets. Then
there exist $\kappa$-small, filtered subsets $A'_0 \subseteq A_0$, $A'_1 \subseteq A_1$ such that
$A''_0 \subseteq A'_0$, $A''_1 \subseteq A'_1$, and the diagrams $p_0|\Nerve(A'_0)$, $p_1|\Nerve (A'_1)$ have the same colimit in $\calC$.
\end{lemma}

\begin{proof}
Let $\overline{p}_0$, $\overline{p}_1$ be colimits of $p_0$ and $p_1$, respectively, which carry the cone points to the same object $C \in \calC$. Let
$B = A''_0 \cup A''_1 \cup \Z_{\geq 0} \cup \{ \infty\}$, which we regard as a partially ordered set so that $$\Nerve(B) \simeq (( \Nerve(A''_0) \coprod \Nerve (A''_1)) \star \Nerve(\Z_{\geq 0}))^{\triangleright}.$$
We will construct sequences of elements
$$ \{ a_0 \leq a_2 \leq \ldots \} \subseteq \{ a \in A_0: (\forall a'' \in A''_0) [a'' \leq a] \}$$
$$ \{ a_1 \leq a_3 \leq \ldots \} \subseteq \{ a \in A_1: (\forall a'' \in A''_1) [a'' \leq a] \}$$
and a diagram $\overline{q}: \Nerve(B) \rightarrow \calC$ such that
$$ \overline{q} | (\Nerve(A''_0) \cup \Nerve \{ 0,2,4, \ldots \})^{\triangleright} = \overline{p}_0 | \Nerve ( A''_0 \cup \{ a_0, a_2, \ldots \})^{\triangleright}$$
$$ \overline{q} | (\Nerve(A''_1) \cup \Nerve \{ 1,3,5, \ldots \})^{\triangleright} = \overline{p}_1 | \Nerve ( A''_1 \cup \{ a_1, a_3, \ldots \})^{\triangleright}.$$
Supposing that this has been done, we take $A'_0 = A''_0 \cup \{ a_0, a_2, \ldots \}$,
$A'_1 = A''_1 \cup \{a_1, a_3, \ldots \}$, and observe that the colimits of $p_0 | \Nerve(A'_0)$ 
and $p_1 | \Nerve(A'_1)$ are both equivalent to the colimit of $\overline{q} | \Nerve(\Z_{\geq 0})$.

The construction is by recursion; let us suppose that the sequence 
$a_0, a_1, \ldots, a_{n-1}$ and the map $\overline{q}_{n} = \overline{q}| ((\Nerve(A''_0) \coprod \Nerve(A''_1)) \star (\Nerve\{0, \ldots, n-1\})^{\triangleright}$ have already been constructed (when $n=0$, we observe that $\overline{q}_{0}$ is uniquely determined by $\overline{p}_0$ and $\overline{p}_1$). For simplicity we will treat only the case where $n$ is even; the case where $n$ is odd can be handled by a similar argument.

Let $q_n = \overline{q}_n | ( \Nerve(A''_0) \coprod \Nerve(A''_1)) \star \Nerve \{0, \ldots, n-1\}$ and
$q'_n = \overline{q}_n | \Nerve(A''_0) \star \Nerve \{0, 2, \ldots, n-2\}$.
According to Corollary \ref{jurman}, the left fibrations $\calC_{q_n/} \rightarrow \calC$
and $\calC_{q'_n/} \rightarrow \calC$ are $\kappa$-compact. Let 
$A_0(n) = \{ a \in A_0: (\forall a'' \in A''_0 \cup \{ a_0, \ldots, a_{n-2} \}) [a'' \leq a] \}$, and let
$X = \calC_{q_n/} \times_{\calC} \Nerve(A_0(n))^{\triangleright}$, $X' = \calC_{q'_n/} \times_{\calC} \Nerve(A_0(n))^{\triangleright}$, so that $X$ and $X'$ are left fibrations classified by colimit diagrams $\Nerve(A_0(n))^{\triangleright} \rightarrow \SSet$. Form a pullback diagram
$$ \xymatrix{ Y \ar[r] \ar[d] & X \ar[d] \\
\Nerve(A_0(n))^{\triangleright} \ar[r] & X' }$$
where the left vertical map is a left fibration (by Proposition \ref{sharpen}) and the bottom horizontal map is determined by $\overline{p} |  \Nerve(A''_0 \cup \{0,\ldots, n-2\}) \star \Nerve (A_0(n))^{\triangleright}$. It follows that the diagram is a homotopy pullback, so that 
$Y \rightarrow \Nerve(A_0(n))^{\triangleright}$ is also a left fibration classified by
a colimit diagram $\Nerve(A_0(n))^{\triangleright} \rightarrow \SSet$. The map
$\overline{q}_{n}$ determines a vertex $v$ of $Y$ lying over the cone point of $\Nerve (A_0(n))^{\triangleright}$. According to Proposition \ref{charspacecolimit}, the inclusion $Y \times_{ \Nerve(A_0(n))^{\triangleright} } \Nerve(A_0(n)) \subseteq Y$ is a weak homotopy equivalence of simplicial sets. It follows that there exists an edge $e: v' \rightarrow v$ of $Y$ which joins
$v$ to some vertex $v'$ lying over an element $a \in A_0(n)$. We now define $a_n = a$, and
observe that the edge $e$ corresponds to the desired extension $\overline{q}_{n+1}$ of
$\overline{q}_n$.
\end{proof}

\begin{lemma}\label{perry}
Let $\calC$ be a presentable $\infty$-category, let $S$ be a strongly saturated collection of morphisms in $\calC$, and let $\calD \subseteq \Fun(\Delta^1,\calC)$ be the full subcategory spanned by $S$. The following conditions are equivalent:
\begin{itemize}
\item[$(1)$] The $\infty$-category $\calD$ is accessible.
\item[$(2)$] The $\infty$-category $\calD$ is presentable.
\item[$(3)$] The collection $S$ is of small generation (as a strongly saturated class of morphisms).
\end{itemize}
\end{lemma}

\begin{proof}
We observe that $\calD$ is stable under small colimits in $\Fun(\Delta^1,\calC)$, and therefore admits small colimits; thus $(1) \Rightarrow (2)$. To see that $(2)$ implies $(3)$, we choose a small collection $S_0$ of morphisms in $\calC$ which generates $\calD$ under colimits; it is then obvious that $S_0$ generates $S$ as a strongly saturated class of morphisms. 

Now suppose that $(3)$ is satisfied. 
Choose a small collection of morphisms $\{ f_{\beta} : X_{\beta} \rightarrow Y_{\beta} \}$ which generates $S$ and an uncountable regular cardinal $\kappa$ such that $\calC$ is $\kappa$-accessible and each of the objects $X_{\beta}$, $Y_{\beta}$ is $\kappa$-compact. We will prove that $\calD$ is $\kappa$-accessible.

It is clear that $\calD$ is locally small and admits $\kappa$-filtered colimits. Let $\calD' \subseteq \calD$ be the collection of all morphisms $f: X \rightarrow Y$ such that $f$ belongs to $S$, where both $X$ and $Y$ are $\kappa$-compact. Lemma \ref{hardstuff1} implies that each $f \in \calD'$ is a $\kappa$-compact object of $\Fun(\Delta^1,\calC)$, and in particular a $\kappa$-compact object of $\calD$.
Assume for simplicity that $\calC$ is a minimal $\infty$-category, so that $\calD'$ is small. According to Proposition \ref{intprop}, the inclusion $\calD' \subseteq \calD$ is equivalent to
$j \circ F$, where $j: \calD' \rightarrow \Ind_{\kappa} \calD'$ is the Yoneda embedding and
$F: \Ind_{\kappa} \calD' \rightarrow \calD$ is $\kappa$-continuous. Proposition \ref{uterr} implies that $F$ is fully faithful; let $\calD''$ denote its essential image. To complete the proof, it will suffice to show that $\calD'' = \calD$. Let $S'' \subseteq S$ denote the collection of objects of $\calD''$ (which we may identify with morphisms in $\calC$). By construction, $S''$ contains the collection of
morphisms $\{ f_{\beta} \}$ which generates $S$. Consequently, to prove that $S'' = S$, it will suffice to show that $S''$ is strongly saturated.

It follows from Proposition \ref{sumatch} that $\calD'' \subseteq \Fun(\Delta^1,\calC)$ is stable under small colimits. We next verify that $S''$ is stable under pushouts. Let $K = \Lambda^2_0$, and let
$\overline{p}: K^{\triangleright} \rightarrow \calC$ be a colimit of $p=\overline{p} | K$, 
$$ \xymatrix{ X \ar[r]^{f} \ar[d] & Y \ar[d] \\
X' \ar[r]^{f'} & Y' }$$
such that $f$ belongs to $S''$. The proof of Proposition \ref{horse1} shows that we can write
$p$ as a colimit of a diagram $q: \Nerve(A) \rightarrow \Fun(\Lambda^2_0,\calC^{\kappa})$,
where $A$ is a $\kappa$-filtered partially ordered set. For $\alpha \in A$, we let
$p_{\alpha}$ denote the corresponding diagram, which we may depict as
$$ X'_{\alpha} \leftarrow X_{\alpha} \stackrel{ f_{\alpha} }{\rightarrow} Y_{\alpha}.$$
For each $A' \subseteq A$, we let $p_{A'}$ denote a colimit of $q| \Nerve(A')$, which we will denote by
$$ X'_{A'} \leftarrow X_{A'} \stackrel{ f_{A'} }{\rightarrow} Y_{A'}. $$
Let $B$ denote the collection of $\kappa$-small, filtered subsets $A' \subseteq A$ such that
the $f_{A'}$ belongs to $S''$. Since $f \in S''$, we conclude that $f$ can be 
obtained as the colimit of a $\kappa$-filtered diagram $\Nerve(A') \rightarrow \calD'$,
Applying Lemma \ref{walter}, we deduce that $B$ is $\kappa$-filtered, and that $A = \bigcup_{A' \in B} A'$. Using Proposition \ref{extet} and Corollary \ref{util}, we deduce that $p$ is the colimit of a
diagram $q': \Nerve(B) \rightarrow \Fun(\Lambda^2_0,\calC)$, where $q'(A') = p_{A'}$. Replacing
$A$ by $B$, we may suppose that each $f_{\alpha}$ belongs to $S'$.

Let $\colim: \Fun(\Lambda^2_0,\calC) \rightarrow \Fun(\Delta^1 \times \Delta^1,\calC)$ be a colimit functor
(that is, a left adjoint to the restriction functor). Lemma \ref{limitscommute} implies
that we may identify $\overline{p}$ with a colimit of the diagram $\colim \circ q$. 
Consequently, the morphism $f'$ can be written as a colimit of morphisms $f'_{\alpha}$ which
fit into pushout diagrams
$$ \xymatrix{ X_{\alpha} \ar[r]^{f_{\alpha}} \ar[d] & Y_{\alpha} \ar[d] \\
X'_{\alpha} \ar[r]^{f'_{\alpha}} & Y'_{\alpha}. }$$
Since $f_{\alpha} \in S'' \subseteq S$, we conclude that $f'_{\alpha} \in S$. Since
$X'_{\alpha}$ and $Y'_{\alpha}$ are $\kappa$-compact, we deduce that
$f'_{\alpha} \in S''$. Since $\calD''$ is stable under colimits, we deduce that $f' \in S''$, as desired.

We now complete the proof by showing that $S''$ has the two-out-of-three property, using the same style of argument. 
Let $\sigma: \Delta^2 \rightarrow \calC$ be a simplex corresponding to a diagram
$$ \xymatrix{ X \ar[rr]^{f} \ar[dr]^{g} & & Y \ar[dl]^{h} \\
& Z & }$$
in $\calC$. We will show that if $f,g \in S''$, then $h \in S''$: the argument in the other two cases is the same. The proof of Proposition \ref{horse1} shows that we can write $\sigma$
as the colimit of a diagram $q: \Nerve(A) \rightarrow \Fun(\Delta^2,\calC^{\kappa})$, where
$A$ is a $\kappa$-filtered partially ordered set. For each $\alpha \in A$, we will denote the
corresponding diagram by
$$ \xymatrix{ X_{\alpha} \ar[rr]^{f_{\alpha}} \ar[dr]^{g_{\alpha}} & & Y_{\alpha} \ar[dl]^{h_{\alpha}} \\
& Z_{\alpha}. & }$$
Arguing as above, we may assume (possibly after changing $A$ and $q$) that each
$f_{\alpha}$ belongs to $S''$. Repeating the same argument, we may suppose that
$g_{\alpha}$ belongs to $S''$. Since $S$ has the two-out-of-three property, we conclude that each $h_{\alpha}$ belongs to $S$. Since $X_{\alpha}$ and $Z_{\alpha}$ are $\kappa$-compact, we then have $h_{\alpha} \in S''$. The stability of $\calD''$ under colimits now implies that $h \in S''$, as desired.
\end{proof}

\begin{proposition}\label{local}
Let $\calC$ be a presentable $\infty$-category, and let $S$ be a $($small$)$ collection of morphisms of $\calC$. Let $\overline{S}$ denote the strongly saturated class of morphisms generated by $S$.
Let $\calC' \subseteq \calC$ denote the full subcategory of $\calC$ consisting of $S$-local objects. Then:

\begin{itemize}
\item[$(1)$] For each object $C \in \calC$, there exists a morphism $s: C \rightarrow C'$ such that
$C'$ is $S$-local and $s$ belongs to $\overline{S}$.
\item[$(2)$] The $\infty$-category $\calC'$ is presentable.
\item[$(3)$] The inclusion $\calC' \subseteq \calC$ has a left adjoint $L$.
\item[$(4)$] For every morphism $f$ of $\calC$, the following are equivalent:
\begin{itemize}
\item[$(i)$] The morphism $f$ is an $S$-equivalence.
\item[$(ii)$] The morphism $f$ belongs to $\overline{S}$.
\item[$(iii)$] The induced morphism $Lf$ is an equivalence.
\end{itemize}
\end{itemize}
\end{proposition}

\begin{proof}
Assertion $(1)$ is a consequence of Lemma \ref{superlocal}, which we will prove in \S \ref{factgen2}. The equivalence $(1) \Leftrightarrow (3)$ follows immediately from Proposition \ref{testreflect}. We now prove $(4)$. Lemma \ref{bluh} implies that the collection of $S$-equivalences
is a strongly saturated class of morphisms containing $S$; it therefore contains $\overline{S}$, so that $(ii) \Rightarrow (i)$. Now suppose that $f: X \rightarrow Y$ is such that $Lf$ is an equivalence,
and consider the diagram 
$$ \xymatrix{ X \ar[r] \ar[dr] \ar[d] & Y \ar[d] \\
LX \ar[r] & LY. }$$
Our proof of $(1)$ shows that the vertical morphisms belong to $\overline{S}$, and the
lower horizontal arrow belongs to $\overline{S}$ by Remark \ref{simpcons}. Two applications of the two-out-of-three property now show that $f \in \overline{S}$, so that $(iii) \Rightarrow (ii)$. If $f$ is an $S$-equivalence, then
we may again use the above diagram and the two-out-of-three property to conclude that $Lf$ is an equivalence. It follows that $LX$ and $LY$ co-represent the same functor on the homotopy category $h\calC'$, so that Yoneda's lemma implies that $Lf$ is an equivalence.
Thus $(i) \Rightarrow (iii)$ and the proof of $(4)$ is complete.

It remains to prove $(2)$. Remark \ref{localcolim} implies that $L \calC$ admits small colimits, so
it will suffice to prove that $L \calC$ is accessible. According to Proposition \ref{localloc}, 
this follows from the implication $(iii) \Rightarrow (i)$ of assertion $(4)$.
\end{proof}

Proposition \ref{local} gives a clear picture of the collection of all accessible localizations of a presentable $\infty$-category $\calC$. For any (small) set of morphisms $S$ in $\calC$, the full subcategory $S^{-1} \calC \subseteq \calC$ consisting of $S$-local objects is a localization of $\calC$, and every localization arises in this way. Moreover, the subcategories $S^{-1} \calC$ and $T^{-1} \calC$\index{not}{S-C@$S^{-1} \calC$}
coincide if and only if $S$ and $T$ generate the same strongly saturated class of morphisms.
We will also write $S^{-1} \calC$ for the class of $S$-local objects of $\calC$ in the case where $S$ is {\em not} small; however, this is generally only a well-behaved object in the case where there is a (small) subset $S_0 \subseteq S$ which generates the same strongly saturated class of morphisms.

\begin{proposition}\label{postbluse}
Let $f: \calC \rightarrow \calD$ be a presentable functor between presentable $\infty$-categories, let $S$ be strongly saturated class of morphisms of $\calD$ which is of small generation. Then $f^{-1} S$ is of small generation $($as a strongly saturated class of morphisms of $\calC${}$)$.
\end{proposition}

\begin{proof}
Replacing $\calD$ by $S^{-1} \calD$ if necessary, we may suppose that $S$ consists
of precisely the equivalences in $\calD$. Let $\calE_{\calD} \subseteq \Fun(\Delta^1,\calD)$ denote the full subcategory spanned by those morphisms which are equivalences in $\calD$, and let
$\calE_{\calC} \subseteq \Fun(\Delta^1,\calC)$ denote the full subcategory spanned by
those morphisms which belong to $f^{-1} S$. We have a homotopy Cartesian diagram
of $\infty$-categories
$$ \xymatrix{ \calE_{\calC} \ar[r] \ar[d] & \calE_{\calD} \ar[d] \\
\Fun(\Delta^1,\calC) \ar[r] & \Fun(\Delta^1,\calD). }$$
The $\infty$-category $\calE_{\calD}$ is equivalent to $\calD$, and therefore presentable.
The $\infty$-categories $\Fun(\Delta^1,\calC)$ and $\Fun(\Delta^1,\calD)$ are presentable by 
Proposition \ref{presexp}. It follows from Proposition \ref{horse22} that $\calE_{\calC}$
is presentable. In particular, there is a small collection of objects of $\calE_{\calC}$ which generates $\calE_{\calC}$ under colimits, as desired.
\end{proof}

Let $\calC$ be a presentable $\infty$-category. We will say that a full subcategory
$\calC^0 \subseteq \calC$ is {\it strongly reflective} if it is
the essential image of an accessible localization functor. Equivalently, $\calC^0$ is strongly reflective if it is presentable, stable under equivalence in $\calC$, and the inclusion
$\calC^0 \subseteq \calC$ admits a left adjoint. According to Proposition \ref{local}, $\calC^0$ is strongly reflective if and only if there exists a (small) set $S$ of morphisms of $\calC$ such that $\calC^0$ is the full subcategory of $\calC$ spanned by the $S$-local objects.
For later use, we record a few easy stability properties enjoyed by the collection of strongly reflective subcategories of $\calC$:\index{gen}{strongly reflective}\index{gen}{reflective subcategory!strongly}

\begin{lemma}\label{stur2}
Let $f: \calC \rightarrow \calD$ be a presentable functor between presentable $\infty$-categories, and let $\calC^0 \subseteq \calC$ be a strongly reflective subcategory. Let $f^{\ast}$ denote a right adjoint of $f$, and let
$\calD^{0} \subseteq \calD$ be the full subcategory spanned by those objects
$D \in \calD$ such that $f^{\ast} D \in \calC^0$. Then $\calD^{0}$ is a strongly reflective subcategory of $\calD$.
\end{lemma}

\begin{proof}
Let $S$ be a (small) set of morphisms of $\calC$ such that $\calC^0$ is the full subcategory
of $\calC$ spanned by the $S$-local objects. Then $\calD^0$ is the full subcategory of $\calD$ spanned by the $f(S)$-local objects.
\end{proof}

\begin{lemma}\label{stur3}
Let $\calC$ be a presentable $\infty$-category, and let $\{ \calC_{\alpha} \}_{\alpha \in A}$ be a family of full subcategories of $\calC$ indexed by a $($small$)$ set $A$. Suppose that each
$\calC_{\alpha}$ is strongly reflective. Then
$\bigcap_{\alpha \in A} \calC_{\alpha}$ is strongly reflective.
\end{lemma}

\begin{proof}
For each $\alpha \in A$, choose a (small) set $S(\alpha)$ of morphisms of $\calC$ such that $\calC_{\alpha}$ is the full subcategory of $\calC$ spanned by the $S(\alpha)$-local objects.
Then $\bigcap_{\alpha \in A} \calC_{\alpha}$ is the full subcategory of $\calC$ spanned by the
$\bigcup_{\alpha \in A} S(\alpha)$-local objects.
\end{proof}

\begin{lemma}\label{stur1}
Let $\calC$ be a presentable $\infty$-category and $K$ a small simplicial set.
Let $\calD$ denote the full subcategory of $\Fun(K^{\triangleleft}, \calC)$ spanned by those diagrams $\overline{p}: K^{\triangleleft} \rightarrow \calC$ which are limits of $p = \overline{p}|K$.
Then $\calD$ is a strongly reflective subcategory of $\calC$.
\end{lemma}

\begin{proof}
The restriction functor $\calD \rightarrow \Fun(K,\calC)$ is an equivalence of $\infty$-categories. This proves that $\calD$ is accessible. Let $s: \Fun(K,\calC) \rightarrow \calD$ be a homotopy inverse to the restriction map. Then the composition
$$ \Fun(K^{\triangleright}, \calC) \rightarrow \Fun(K,\calC) \stackrel{s}{\rightarrow} \calD$$
is left adjoint to the inclusion.
\end{proof}

We conclude this section by giving a universal property which characterizes the localization $S^{-1} \calC$.

\begin{proposition}\label{unichar}
Let $\calC$ be a presentable $\infty$-category, and $\calD$ an arbitrary $\infty$-category.
Let $S$ be a $($small$)$ set of morphisms of $\calC$, and $L: \calC \rightarrow S^{-1} \calC \subseteq \calC$ an associated (accessible) localization functor. Composition with $L$ induces
a functor
$$ \eta: \LFun( S^{-1} \calC, \calD) \rightarrow \LFun( \calC, \calD).$$
The functor $\eta$ is fully faithful, and the essential image of $\eta$ consists of those functors
$f: \calC \rightarrow \calD$ such that $f(s)$ is an equivalence in $\calD$, for each $s \in S$.
\end{proposition}

\begin{proof}
Let $\alpha: \id_{\calC} \rightarrow L$ be a unit for the adjunction between $L$ and
the inclusion $S^{-1} \calC \subseteq \calC$.
We first observe that every functor $f_0: S^{-1} \calC \rightarrow \calD$ admits a right
Kan extension $f: \calC \rightarrow \calD$. To prove this, we may first replace $f_0$ by
the equivalent diagram $g_0 = f_0 \circ (L| S^{-1} \calC)$, and define $g= f_0 \circ L$. To prove
that $g$ is a right Kan extension of $g_0$, it suffices to show that for each
object $X \in \calC$, the diagram
$$ \overline{p}: (S^{-1} \calC)_{X/}^{\triangleleft} \rightarrow \calC \stackrel{L}{\rightarrow} 
S^{-1} \calC \stackrel{f_0}{\rightarrow} \calD$$
exhibits $f_0(LX)$ as a limit of $p = \overline{p} | (S^{-1} \calC)_{/X})$. For this, we note
that an $S$-localization map $\alpha(X): X \rightarrow LX$ is an initial object of $(S^{-1} \calC)_{X/}$ (Remark \ref{initrem}), and that $f_0( L \alpha(X))$ is an equivalence by Proposition \ref{recloc}.

Let $\calX$ denote the full subcategory of $\calD^{\calC}$ spanned by those functors
$f: \calC \rightarrow \calD$ which are right Kan extensions of $f| S^{-1} \calC$. According to 
Proposition \ref{lklk}, the restriction map $\calX \rightarrow \Fun(S^{-1} \calC,\calD)$ is a
trivial fibration. Let $\overline{\eta}: \Fun(S^{-1} \calC, \calD) \rightarrow \Fun(\calC, \calD)$ be given
by composition with $L$. The above argument shows that $\overline{\eta}$ factors through $\calX$.
Moreover, the composition of $\overline{\eta}$ with the restriction map is homotopic to the identity on $\Fun(S^{-1} \calC, \calD)$. It follows that $\overline{\eta}$ is an equivalence of $\infty$-categories.

We have a commutative diagram
$$ \xymatrix{ \LFun(S^{-1} \calC, \calD) \ar[r]^{\eta} \ar[d] & \LFun(\calC, \calD) \ar[d] \\
\Fun(S^{-1} \calC,\calD) \ar[r]^{\overline{\eta}} & \Fun(\calC, \calD) }$$
where the vertical maps are inclusions of full subcategories, and the lower horizontal map is fully faithful. It follows that $\eta$ is fully faithful. To complete the proof, we must show that a functor
$f: \calC \rightarrow \calD$ belongs to the essential image of $\eta$ if and only if
$f(s)$ is an equivalence for each $s \in S$. The ``only if'' direction is clear, since
the functor $L$ carries each element of $S$ to an equivalence in $\calC$. Conversely, suppose that
$f$ carries each $s \in S$ to an equivalence. The natural transformation $\alpha$
gives a map of functors $\alpha(f): f \rightarrow f \circ L$; we wish to show that $\alpha(f)$ is an equivalence. Equivalently, we wish to show that for each object $X \in \calC$, $f$ carries
the map $\alpha(X): X \rightarrow LX$ to an equivalence in $\calD$. Let $S'$ denote the
class of all morphisms $\phi$ in $\calC$ such that $f(\phi)$ is an equivalence in $\calD$.
By assumption, $S \subseteq S'$. Lemma \ref{bluh} implies that $S'$ is strongly saturated, so
that Proposition \ref{local} asserts that $\alpha(X) \in S'$, as desired.
\end{proof}

\subsection{Factorization Systems on Presentable $\infty$-Categories}\label{factgen2}

Let $\calC$ be a presentable $\infty$-category. In \S \ref{invloc}, we saw that it is easy to produce localizations of $\calC$: for any small collection of morphisms $S$ in $\calC$, the full subcategory $S^{-1} \calC$ of $S$-local objects of $\calC$ is a presentable localization of $\calC$, which depends only on the strongly saturated class of morphisms $\overline{S}$ generated by $S$. Our goal in this section is to prove a similar result for factorization systems on
$\calC$. The first step is to introduce the analogue of the notion of ``strongly saturated'':

\begin{definition}\index{gen}{saturated}
Let $S$ be a collection of morphisms in a presentable $\infty$-category $\calC$. We will say that
$S$ is {\it saturated} if the following conditions are satisfied:
\begin{itemize}
\item[$(1)$] The collection $S$ is closed under small colimits in $\Fun( \Delta^1, \calC)$.
\item[$(2)$] The collection $S$ contains all equivalences and is stable under composition.
\item[$(3)$] The collection $S$ is closed under the formation of pushouts. That is, given a pushout diagram
$$ \xymatrix{ X \ar[r] \ar[d]^{f} & X' \ar[d]^{f'} \\
Y \ar[r] & Y' }$$
in $\calC$, if $f$ belongs to $S$ then $f'$ also belongs to $S$.
\end{itemize}
\end{definition}

\begin{remark}\index{gen}{of small generation}
Let $\calC$ be a presentable $\infty$-category. Then any intersection of saturated collections of morphisms in $\calC$ is again saturated. It follows that for {\em any} class of morphisms $S$ of $\calC$, there exists a smallest saturated collection of morphisms $\overline{S}$ containing $S$. We will refer to $\overline{S}$ as the {\it saturated collection of morphisms generated by $S$}. We will say that a saturated collection of morphisms $\overline{S}$ is {\it of small generation} if it
is generated by some (small) subset $S \subseteq \overline{S}$.
\end{remark}

\begin{remark}
If $S$ is a saturated collection of morphisms of $\calC$, then $S$ is closed under retracts.
\end{remark}

\begin{remark}
Let $\calC$ be (the nerve of) a presentable category, and let $S$ be a saturated class of morphisms
in $\calC$. Then $S$ is also weakly saturated, in the sense of Definition \ref{saturated}.
\end{remark}

\begin{example}
Let $\calC$ be a presentable $\infty$-category. Then every strongly saturated class of morphisms in $\calC$ is also saturated.
\end{example}

\begin{example}\label{cobblet}
Let $\calC$ be a presentable $\infty$-category, and $S$ any collection of morphisms of $\calC$.
Then $^{\perp}S$ is saturated; this follows immediately from Proposition \ref{swimmm}. 
In particular, if $(S_L, S_R)$ is a factorization system on $\calC$, then $S_L$ is saturated.
\end{example}

The main result of this section is the following converse to Example \ref{cobblet}:

\begin{proposition}\label{nir}
Let $\calC$ be a presentable $\infty$-category, and $S$ a saturated collection of morphisms in $\calC$ which is of small generation. Then $(S, S^{\perp})$ is a factorization system on $\calC$.
\end{proposition}

\begin{corollary}\label{wugg}
Let $\calC$ be a presentable $\infty$-category, let $S$ be a saturated collection of morphisms of $\calC$, and suppose that $S$ is of small generation. Let
$$ \xymatrix{ & Y \ar[dr]^{g} & \\
X \ar[ur]^{f} \ar[rr]^{h} & & Z} $$
be a commutative diagram in $\calC$. If $f$ and $h$ belong to $S$, then $g$ belongs to $S$.
\end{corollary}

\begin{proof}
Combine Propositions \ref{nir}, \ref{swin}, and \ref{swimmm}. 
\end{proof}

In the situation of Proposition \ref{nir}, we will refer to the elements of $S^{\perp}$ as {\it $S$-local morphisms of $\calC$}. Note that an object $X \in \calC$ is $S$-local if and only if a
morphism $X \rightarrow 1_{\calC}$ is $S$-local, where $1_{\calC}$ denotes a final object of $\calC$. 

The proof of Proposition \ref{nir} will be given at the end of this section, after we have established a series of technical lemmas.

\begin{lemma}\label{sweener}
Let $\calC$ be a presentable $\infty$-category, and $S$ a saturated collection of morphisms
of $\calC$. The following conditions are equivalent:
\begin{itemize}
\item[$(1)$] The collection $S$ is of small generation.
\item[$(2)$] The full subcategory $\calD \subseteq \Fun( \Delta^1, \calC)$ spanned by the elements of
$S$ is presentable.
\end{itemize}
\end{lemma}

\begin{proof}
If $\calD$ is presentable, then $\calD$ is generated under small colimits by a small set of objects; these objects clearly generate $S$ as a saturated collection of morphisms. This proves that $(2) \Rightarrow (1)$. To prove the reverse implication, choose a small collection of morphisms $\{ f_{\beta} : X_{\beta} \rightarrow Y_{\beta} \}$ which generates $S$ as a semiaturated class of morphisms and an uncountable regular cardinal $\kappa$ such that $\calC$ is $\kappa$-accessible and each of the objects $X_{\beta}$, $Y_{\beta}$ is $\kappa$-compact. Let $\calD' \subseteq \calD$ be the collection of all morphisms $f: X \rightarrow Y$ such that $f$ belongs to $S$, where both $X$ and $Y$ are $\kappa$-compact. Lemma \ref{hardstuff1} implies that each $f \in \calD'$ is a $\kappa$-compact object of $\Fun(\Delta^1,\calC)$, and in particular a $\kappa$-compact object of $\calD$. 
Assume for simplicity that $\calC$ is a minimal $\infty$-category, so that $\calD'$ is small. According to Proposition \ref{intprop}, the inclusion $\calD' \subseteq \calD$ is equivalent to
$j \circ F$, where $j: \calD' \rightarrow \Ind_{\kappa} \calD'$ is the Yoneda embedding and
$F: \Ind_{\kappa} \calD' \rightarrow \calD$ is $\kappa$-continuous. Proposition \ref{uterr} implies that $F$ is fully faithful; let $\calD''$ denote its essential image. To complete the proof, it will suffice to show that $\calD'' = \calD$. 

Let $S' \subseteq S$ denote the collection of objects of $\calD''$ (which we may identify with morphisms in $\calC$). By construction, $S'$ contains the collection of morphisms $\{ f_{\beta} \}$ which generates $S$. Consequently, to prove that $S' = S$, it will suffice to show that $S'$ is saturated.

It follows from Proposition \ref{sumatch} that $\calD'' \subseteq \Fun(\Delta^1,\calC)$ is stable under small colimits. We next verify that $S'$ is stable under pushouts. Let $K = \Lambda^2_0$, and let
$\overline{p}: K^{\triangleright} \rightarrow \calC$ be a colimit of $p=\overline{p} | K$, 
$$ \xymatrix{ X \ar[r]^{f} \ar[d] & Y \ar[d] \\
X' \ar[r]^{f'} & Y' }$$
such that $f$ belongs to $S'$. Using Proposition \ref{urgh1}, we can write $p$ as the colimit of a diagram $q: \Nerve(A) \rightarrow \Fun(\Lambda^2_0,\calC^{\kappa})$,
where $A$ is a $\kappa$-filtered partially ordered set. For $\alpha \in A$, we let
$p_{\alpha}$ denote the corresponding diagram, which we may depict as
$$ X'_{\alpha} \leftarrow X_{\alpha} \stackrel{ f_{\alpha} }{\rightarrow} Y_{\alpha}.$$
For each $A' \subseteq A$, we let $p_{A'}$ denote a colimit of $q| \Nerve(A')$, which we will denote by
$$ X'_{A'} \leftarrow X_{A'} \stackrel{ f_{A'} }{\rightarrow} Y_{A'}. $$
Let $B$ denote the collection of $\kappa$-small, filtered subsets $A' \subseteq A$ such that
the $f_{A'}$ belongs to $S'$. Since $f \in S'$, we conclude that $f$ can be 
obtained as the colimit of a $\kappa$-filtered diagram in $\calD'$,
Applying Lemma \ref{walter}, we deduce that $B$ is $\kappa$-filtered, and that $A = \bigcup_{A' \in B} A'$. Using Proposition \ref{extet} and Corollary \ref{util}, we deduce that $p$ is the colimit of a diagram $q': \Nerve(B) \rightarrow \Fun(\Lambda^2_0,\calC)$, where $q'(A') = p_{A'}$. Replacing $A$ by $B$, we may suppose that each $f_{\alpha}$ belongs to $S'$.

Let $\colim: \Fun(\Lambda^2_0,\calC) \rightarrow \Fun(\Delta^1 \times \Delta^1,\calC)$ be a colimit functor (that is, a left adjoint to the restriction functor). Lemma \ref{limitscommute} implies
that we may identify $\overline{p}$ with a colimit of the diagram $\colim \circ q$. 
Consequently, the morphism $f'$ can be written as a colimit of morphisms $f'_{\alpha}$ which
fit into pushout diagrams
$$ \xymatrix{ X_{\alpha} \ar[r]^{f_{\alpha}} \ar[d] & Y_{\alpha} \ar[d] \\
X'_{\alpha} \ar[r]^{f'_{\alpha}} & Y'_{\alpha}. }$$
Since $f_{\alpha} \in S' \subseteq S$, we conclude that $f'_{\alpha} \in S$. Since
$X'_{\alpha}$ and $Y'_{\alpha}$ are $\kappa$-compact, we deduce that
$f'_{\alpha} \in S'$. Since $\calD''$ is stable under colimits, we deduce that $f' \in S'$, as desired.

We now complete the proof by showing that $S'$ is stable under composition.
Let $\sigma: \Delta^2 \rightarrow \calC$ be a simplex corresponding to a diagram
$$ \xymatrix{ X \ar[rr]^{f} \ar[dr]^{g} & & Y \ar[dl]^{h} \\
& Z & }$$
in $\calC$. We will show that if $f,g \in S'$, then $h \in S'$. Using Proposition \ref{urgh1}, we can write $\sigma$ as the colimit of a diagram $q: \Nerve(A) \rightarrow \Fun(\Delta^2,\calC^{\kappa})$, where
$A$ is a $\kappa$-filtered partially ordered set. For each $\alpha \in A$, we will denote the
corresponding diagram by
$$ \xymatrix{ X_{\alpha} \ar[rr]^{f_{\alpha}} \ar[dr]^{g_{\alpha}} & & Y_{\alpha} \ar[dl]^{h_{\alpha}} \\
& Z_{\alpha}. & }$$
Arguing as above, we may assume (possibly after changing $A$ and $q$) that each
$f_{\alpha}$ belongs to $S'$. Repeating the same argument, we may suppose that
$g_{\alpha}$ belongs to $S'$. Since $S$ is stable under composition, we conclude that each
$h_{\alpha}$ belongs to $S$. Since each $X_{\alpha}$ and $Z_{\alpha}$ are $\kappa$-compact, we have $h_{\alpha} \in S'$. The stability of $\calD''$ under colimits now implies that $h \in S'$, as desired.
\end{proof}

\begin{lemma}\label{swunl}
Let $\calC$ be a presentable $\infty$-category, and let $S$ be a saturated collection of morphisms in $\calC$. For every object $X \in \calC$, let
$S_X$ denote the collection of all morphisms of $\calC_{/X}$ whose image in $\calC$ belongs to $S$. Then each $S_{X}$ is strongly saturated in $\calC_{/X}$. Moreover, if $S$ is of small generation, then each $S_{X}$ is also of small generation.
\end{lemma}

\begin{proof}
The first assertion follows immediately from the definitions and Proposition \ref{needed17}.
To prove the second, let $\calD$ be the full subcategory of $\Fun( \Delta^1, \calC)$ spanned by the elements of $S$, and $\calD'$ the full subcategory of $\Fun( \Delta^1, \calC_{/X})$ spanned by the elements of $S_X$. We have a (homotopy) pullback diagram of $\infty$-categories
$$ \xymatrix{ \calD' \ar[d] \ar[r] & \Fun( \Delta^1, \calC_{/X} ) \ar[d]^{\psi} \\
\calD \ar[r]^-{\phi} & \Fun(\Delta^1, \calC). }$$
The functors $\phi$ and $\psi$ preserve small colimits, and the $\infty$-categories
$\Fun( \Delta^1, \calC_{/X})$, $\Fun( \Delta^1, \calC)$, and $\calD$ are all presentable
(the last in view of Lemma \ref{sweener}). Using Proposition \ref{horse22}, we deduce that $\calD'$ is presentable, and therefore generated under small colimits by a (small) set of elements of $S_{X}$. This proves that $S_{X}$ is of small generation, as desired.
\end{proof}

\begin{lemma}\label{sweetyork}
Let $\calC$ be a presentable $\infty$-category, $S$ a saturated collection of morphisms of $\calC$, and $X$ an object of $\calC$. Then the full subcategory $\calD \subseteq \calC^{X/}$ spanned by the elements which belong to $S$ is closed under small colimits.
\end{lemma}

\begin{proof}
In view of Corollary \ref{uterrr} and \ref{allfin}, it will suffice to show that
$\calD$ is closed under small filtered colimits, pushouts, and contains the initial objects
of $\calC^{X/}$. The last condition follows from the fact that $S$ contains all equivalences.
Now suppose that $\overline{p}: K^{\triangleright} \rightarrow \calC^{C/}$
is a colimit of $p = \overline{p} | K$, where $K$ is either filtered or equivalent to $\Lambda^2_0$,
and that $p(K) \subseteq \calD$. We can identify $\overline{p}$ with a map
$P: K^{\triangleright} \times \Delta^1 \rightarrow \calC$ such that
$P | K^{\triangleright} \times \{0\}$ is the constant map taking the value $C \in \calC$.
Since $K$ is weakly contractible, $P | K^{\triangleright} \times \{0\}$ is a colimit diagram in $\calC$.
The map $P | K^{\triangleright} \times \{1\}$ is the image of a colimit diagram under the left
fibration $\calC^{C/} \rightarrow \calC$; since $K$ is weakly contractible, Proposition \ref{goeselse} implies that $P | K^{\triangleright} \times \{1\}$ is a colimit diagram. We now apply
Proposition \ref{limiteval} to deduce that
$P: K^{\triangleright} \rightarrow \calC^{\Delta^1}$ is a colimit diagram. Since $S$ is stable under colimits in $\calC^{\Delta^1}$, we conclude that $P$
carries the cone point of $K^{\triangleright}$ to a morphism belonging to $S$, as we wished to show.
\end{proof}

\begin{lemma}\label{sugarplace}
Let $\calC$ be an $\infty$-category, and let $f: C \rightarrow D$, $g: C \rightarrow E$ be morphisms in $\calC$. Then there is a natural identification of $\bHom_{\calC_{C/}}(f,g)$ with
the homotopy fiber of the map
$$ \bHom_{\calC}(D,E) \rightarrow \bHom_{\calC}(C,E)$$
induced by composition with $f$, where the fiber is taken over the point corresponding to $g$.
\end{lemma}

\begin{proof}
We have a commutative diagram of simplicial sets
$$ \xymatrix{ \calC_{f/} \times_{\calC_{C/}} \{g\} \ar[d]^{\phi} \ar[r] & \calC_{f/} \times_{\calC} \{E\} \ar[d]^{\phi'} \ar[r] & \calC_{f/} \ar[d]^{\phi''} \\
\{g\} \ar[r] & \calC_{C/} \times_{\calC} \{E\}  \ar[r] & \calC_{C/} }$$
where both squares are pullbacks. Proposition \ref{sharpen} asserts that $\phi''$ is a left fibration, so that $\phi'$ and $\phi$ are left fibrations as well. 
Since $\calC_{C/} \times_{\calC} \{E\} = \Hom^{\lft}_{\calC}(C,E)$ is a Kan complex, the map
$\phi'$ is actually a Kan fibration (Lemma \ref{toothie2}), so that the square on the left is homotopy pullback, and identifies
$$ \calC_{f/} \times_{ \calC_{C/} } \{g\} \simeq \bHom_{\calC_{C/}}(f,g)$$
with the homotopy fiber of $\phi'$ over $g$; we conclude by observing that $\phi'$ is a model for the map $$ \bHom_{\calC}(D,E) \rightarrow \bHom_{\calC}(C,E)$$ given by composition with $f$. 
\end{proof}

\begin{lemma}\label{sugarplaces}
Let $$ \xymatrix{ X \ar[r]^{f} \ar[d]^{g} & X' \ar[d] \\
Y \ar[r]^{f'} & Y' }$$ be a pushout diagram in an $\infty$-category $\calC$.
Then there exists an isomorphism
$$ \bHom_{\calC_{X/}}(f,g) \simeq \bHom_{\calC_{Y/}}(f', \id_{Y})$$
in the homotopy category $\calH$. 
\end{lemma}

\begin{proof}
According to Corollary \ref{strictify}, we can assume without loss of generality
that $\calC$ is the nerve of a fibrant simplicial category $\calD$, and that the diagram in question is the nerve of a commutative diagram
 $$ \xymatrix{ X \ar[r]^{f} \ar[d]^{g} & X' \ar[d] \\
Y \ar[r]^{f'} & Y' }$$
in $\calD$. Theorem \ref{colimcomparee} implies that this diagram is homotopy coCartesian in $\calD$, so that we have a homotopy pullback diagram
$$ \xymatrix{ \bHom_{\calD}(Y', Y) \ar[r]^{\phi} \ar[d] & \bHom_{\calD}(Y,Y) \ar[d] \\
\bHom_{\calD}(X',Y) \ar[r]^{\phi'} & \bHom_{\calD}(X,Y) }$$
of Kan complexes. Consequently, we obtain an isomorphism in $\calH$ between  the homotopy fiber of $\phi$ over $\id_{Y}$ and the homotopy fiber of $\phi'$ over $g$.
According to Lemma \ref{sugarplace}, these homotopy fibers may be identified with
$\bHom_{\calC_{Y/}}(f', \id_{Y})$ and $ \bHom_{\calC_{X/}}(f,g)$, respectively. 
\end{proof}

\begin{lemma}\label{superlocal}
Let $\calC$ be a presentable $\infty$-category, and let $S$ be a saturated collection of morphisms of $\calC$ which is of small generation. Then, for every object $X \in \calC$, there exists a morphism
$f: X \rightarrow Y$ in $\calC$, such that $f \in S$ and $Y$ is $S$-local.
\end{lemma}

\begin{proof}
Let $\calD$ be the full subcategory of $\Fun(\Delta^1,\calC)$ spanned by the elements of $S$, and form a fiber diagram
$$ \xymatrix{ \calD_{X} \ar[r] \ar[d] & \calD \ar[d] \\
\{X\} \ar[r] & \Fun(\{0\}, \calC). }$$
Since $S$ is stable under pushouts, the right vertical map is a coCartesian fibration, so that the above diagram is homotopy Cartesian by Proposition \ref{basechangefunky}. Lemma 
\ref{sweener} asserts that $\calD$ is accessible, so that $\calD_{X}$ is accessible by Proposition \ref{horse2}. Using Lemma \ref{sweetyork}, we conclude that $\calD_{X}$ is presentable, so that
$\calD_{X}$ has a final object $f: X \rightarrow Y$. To complete the proof, it will suffice to show that $Y$ is $S$-local.

Let $t: A \rightarrow B$ be an arbitrary morphism in $\calC$ which belongs to $S$. We wish to show that composition with $t$ induces a homotopy equivalence $\phi: \bHom_{\calC}(B,Y) \rightarrow \bHom_{\calC}(A,Y)$. Let $g: A \rightarrow Y$ be an arbitrary morphism; using Lemma \ref{sugarplace} we may identify $\bHom_{\calC_{A/}}(t,g)$ with the homotopy fiber of
$\phi$ over the base point $g$ of $\bHom_{\calC}(A,Y)$. We wish to show that this space
is contractible. Form a pushout diagram
$$ \xymatrix{ A \ar[d]^{g} \ar[r]^{t} & B \ar[d] \\
Y \ar[r]^{t'} & Z}$$
in the $\infty$-category $\calC$. 
Lemma \ref{sugarplaces} implies the existence of a homotopy equivalence
$\bHom_{\calC_{A/}}(t,g) \simeq \bHom_{\calC_{Y/}}(t', \id_{Y})$. It will therefore suffice to prove that $\bHom_{\calC_{Y/}}(t', \id_{Y})$ is contractible. Since $t'$ is a pushout of $t$, it belongs to $S$. Let
$\sigma$ be a $2$-simplex of $\calC$ classifying a diagram
$$ \xymatrix{ & Y \ar[dr]^{t'} & \\
X \ar[ur]^{s} \ar[rr]^{s'} & & Z, }$$
so that $s'$ is a composition of the morphisms $s$ and $t'$ in $\calC$, and therefore also belongs to $S$.
Applying Lemma \ref{sugarplace} again, we may identify
$$\bHom_{\calC_{Y/}}(t', \id_{C'}) \simeq \bHom_{\calC_{s/}}(\sigma, s_1(s))$$ with the homotopy fiber of the map
$$ \bHom_{ \calC_{Y/}}(s',s) \rightarrow \bHom_{ \calC_{Y/}}(s,s).$$
given by composition with $\sigma$.
By construction, $\calD_{X}$ is a full subcategory of $\calC^{X/}$ which contains
$s$ and $s'$, and $s$ is a final object of $\calD_{X}$.
In view of the equivalence of $\calC_{X/}$ with $\calC^{X/}$, we conclude that the spaces
$\bHom_{\calC_{X/}}(s',s)$ and $\bHom_{\calC_{X/}}(s,s)$ are contractible, so that
$\phi$ is a homotopy equivalence as desired.
\end{proof}

\begin{proof}[Proof of Proposition \ref{nir}]
Let $h: X \rightarrow Z$ be a morphism in $\calC$; we wish to show that $h$ admits a factorization
$$ \xymatrix{ & Y \ar[dr]^{g} & \\
X \ar[ur]^{f} \ar[rr]^{h} & & Z }$$
where $f \in S$ and $g \in S^{\perp}$. Using Remark \ref{spack}, we deduce that a morphism
$g: Y \rightarrow Z$ belongs to $S^{\perp}$ if and only if it is an $S_{Z}$-local object of
$\calC_{/Z}$, where $S_Z$ is defined as in Lemma \ref{swunl}. The existence of $h$ then follows from
Lemma \ref{superlocal}.
\end{proof}

\subsection{Truncated Objects}\label{truncintro}

Let $X$ be a topological space. The first step in the homotopy-theoretic analysis of the space $X$ is to divide $X$ into path components. The situation can be described as follows: we associate to $X$ a set $\pi_0 X$, which we may view as a discrete topological space. There is a map $f: X \rightarrow \pi_0 X$ which collapses each component of $X$ to a point. If $X$ is a sufficiently nice space (for example, a CW complex), then the path components of $X$ are open, so $f$ is continuous. Moreover, $f$ is universal
among continuous maps from $X$ into a discrete topological space.

The next step in the analysis of $X$ is to consider its fundamental group $\pi_1 X$, which (provided that $X$ is sufficiently nice) may be studied by means of a universal cover $\widetilde{X}$ of $X$.
However, it is important to realize that neither $\pi_1 X$ nor $\widetilde{X}$ is invariantly associated to $X$: both require a choice of base point. The situation can be described more canonically as follows: to $X$ we can associate a {\em fundamental groupoid} $\pi(X)$, and a map $\phi$ from
$X$ to the classifying space $B \pi(X)$. The universal cover $\widetilde{X}$ of $X$ can be identified
(up to homotopy equivalence) with the homotopy fibers of the map $\phi$. The classifying space $B \pi(X)$ can be regarded as a ``quotient'' of $X$, obtained by killing all of the higher homotopy groups of $X$. Like $\pi_0 X$, it can be described by a universal mapping property. 

To continue the analysis, we first recall that a space $Y$ is said to be {\it $k$-truncated} if the homotopy groups of $Y$ vanish in dimensions larger than $k$ (see Definition \ref{trunckan}). 
Every (sufficiently nice) topological space $X$ admits an essentially unique {\em Postnikov tower}
$$ X \rightarrow \ldots \rightarrow \tau_{\leq n} X \rightarrow \ldots \rightarrow \tau_{\leq -1} X$$
where $\tau_{\leq i} X$ is $i$-truncated, and is universal (in a suitable homotopy-theoretic sense) among $i$-truncated spaces which admit a map from $X$. For example, we can take
$\tau_{\leq 0} X = \pi_0 X$, considered as a discrete space, and $\tau_{\leq 1} X = B \pi(X)$.
Moreover, we can recover the space $X$ (up to weak homotopy equivalence) by
taking the homotopy limit of the tower.\index{gen}{Postnikov tower}\index{gen}{tower!Postnikov}

The objective of this section is to construct an analogous theory in the case where $X$ is not a space, but an object of some (abstract) $\infty$-category $\calC$. We begin by observing that the condition that a space $X$ is $k$-truncated can be reformulated in more categorical terms:
a Kan complex $X$ is $k$-truncated if and only if, for every simplicial set $S$, the mapping space
$\bHom_{\sSet}(S,X)$ is $k$-truncated. This motivates the following:

\begin{definition}\label{tooka}
Let $\calC$ be an $\infty$-category and $k \geq -1$ an integer. We will say that an object
$C$ of $\calC$ is {\it $k$-truncated} if, for every object $D \in \calC$, the
space $\bHom_{\calC}(D,C)$ is $k$-truncated. By convention, we will say that $C$ is {\it $(-2)$-truncated} if it is a final object of $\calC$. We will say that an object of $\calC$ is {\it discrete} if it is $0$-truncated. We will generally let $\tau_{\leq k} \calC$ denote the full subcategory of $\calC$ spanned by the $k$-truncated objects.\index{gen}{truncated!object of an $\infty$-category}\index{gen}{$k$-truncated!object of an $\infty$-category}\index{gen}{discrete}
\end{definition}

\begin{notation}
Let $\calC$ be an $\infty$-category. Using Propositions \ref{tokerp} and \ref{huka}, we conclude that the full subcategory $\tau_{\leq 0} \calC$ is equivalent to the nerve of its homotopy category. We will denote this homotopy category by $\Disc(\calC)$, and refer to it as the {\it category of discrete objects of $\calC$}.\index{not}{DiscC@$\Disc(\calC)$}
\end{notation}

\begin{lemma}
Let $C$ be an object of an $\infty$-category $\calC$, and let $k \geq -2$. The following conditions
are equivalent:
\begin{itemize}
\item[$(1)$] The object $C$ is $k$-truncated.
\item[$(2)$] For every $n \geq k+3$ and every diagram
$$ \xymatrix{ \bd \Delta^n \ar@{^{(}->}[d] \ar[r]^{f} & \calC \\
\Delta^n \ar@{-->}[ur], & }$$ for which $f$ carries the final vertex of $\Delta^n$ to
$C$, there exists a dotted arrow rendering the diagram commutative.
\end{itemize}
\end{lemma}

\begin{proof}
Suppose first that $(2)$ is satisfied. Then for every object $D \in \calD$, the Kan complex
$\Hom^{\rght}_{\calC}(D,C)$ and every $n \geq k+3$
has the extension property with respect to $\bd \Delta^{n-1} \subseteq \Delta^{n-1}$, and is therefore $k$-truncated. For the converse, suppose that $(1)$ is satisfied and choose a categorical equivalence $g: \calC \rightarrow \tNerve \calD$, where $\calD$ is a topological category.
According to Proposition \ref{princex}, it will suffice to show that for every $n \geq k+3$ and every diagram $$ \xymatrix{ | \sCoNerve[\bd \Delta^n] | \ar@{^{(}->}[d] \ar[r]^{F} & \calD \\
| \sCoNerve[\Delta^n] | \ar@{-->}[ur], & }$$
having the property that $F$ carries the final object of $| \sCoNerve[\Delta^n] |$ to
$g(C)$, there exists a dotted arrow as indicated, rendering the diagram commutative. Let $D \in \calD$ denote the image of the initial object of $\bd \Delta^n$ under $F$. Then constructing the desired extension is equivalent to extending a map $\bd [0,1]^{n-1} \rightarrow \bHom_{\calD}( D, g(C))$ to a map defined on all of $[0,1]^{n-1}$, which is possibly in virtue of the assumption $(1)$.
\end{proof}

\begin{remark}\label{humpter}
A Kan complex $X$ is $k$-truncated if and only if it is $k$-truncated when regarded as an object in the $\infty$-category $\SSet$ of spaces.
\end{remark}

\begin{proposition}\label{altum}
Let $\calC$ be an $\infty$-category, and $k \geq -2$ an integer. The full subcategory
$\tau_{\leq k} \calC \subseteq \calC$ of $k$-truncated objects is stable under all limits
which exist in $\calC$.\index{not}{taukC@$\tau_{\leq k} \calC$}
\end{proposition}

\begin{proof}
Let $j: \calC \rightarrow \calP(\calC)$ be the Yoneda embedding. By definition,
$\tau_{\leq k} \calC$ is the preimage of $\Fun(\calC^{op}, \tau_{\leq k} \SSet)$ under
$j$. Since $j$ preserves all limits which exist in $\calC$, it will suffice to prove that
$\Fun(\calC^{op},\tau_{\leq k} \SSet) \subseteq \Fun(\calC^{op},\SSet)$ is stable under limits.
Using Proposition \ref{limiteval}, it suffices to prove that the inclusion $i: \tau_{\leq k} \SSet \subseteq \SSet$ is stable under limits. In other words, we must show that $\tau_{\leq k} \SSet$
admits small limits, and that $i$ preserves small limits. According to Propositions \ref{alllimits} and \ref{allimits}, it will suffice to show that $\tau_{\leq k} \SSet \subseteq \SSet$ is stable
under the formation of pullbacks and small products. According to Theorem \ref{colimcomparee}, this is equivalent to the assertion that the full subcategory of $\Kan$ spanned by the $k$-truncated Kan complexes is stable under homotopy products and the formation of homotopy pullback squares. Both assertions can be verified easily by computing homotopy groups.
\end{proof}

\begin{remark}\label{socal}
Let $p: \calC \rightarrow \calD$ be a coCartesian fibration of $\infty$-categories.
Let $C$ and $C'$ be objects of $\calC$, let $f: p(C') \rightarrow p(C)$ be a morphism in $\calC$,
and let $\overline{f}: C' \rightarrow C''$ be a $p$-coCartesian morphism lifting $f$.
According to Proposition \ref{compspaces}, we may identify $\bHom_{\calC_{p(C)}}(C'',C)$ with the homotopy fiber of $\bHom_{\calC}(C',C) \rightarrow \bHom_{\calD}(p(C'), p(C))$ over the base point determined by $f$.
By examining the associated long exact sequences of homotopy groups (as $f$ varies), we conclude that if $C$ is a $k$-truncated object of the fiber $\calC_{p(C)}$ and 
$p(C)$ is a $k$-truncated object of $\calD$, then $C$ is a $k$-truncated object of $\calC$. This can be considered as a generalization of Lemma \ref{sabreto} (which treats the case $k=-2$).
\end{remark}

\begin{remark}\label{trumble}
Let $p: \calM \rightarrow \Delta^1$ be a coCartesian fibration of simplicial sets, which we regard as
a correspondence from the $\infty$-category $\calC = p^{-1} \{0\}$ to $\calD = p^{-1} \{1\}$. 
Suppose that $D$ is a $k$-truncated object of $\calD$. Remark \ref{socal} implies that $D$ is a $k$-truncated object of $\calM$. Let $C,C' \in \calC$, and let $f: C \rightarrow D$ be a $p$-Cartesian morphism of $\calM$. Then composition with $f$ induces a homotopy equivalence
$ \bHom_{\calC}(C',C) \rightarrow \bHom_{\calM}(C',D)$; we conclude that $C$ is a $k$-truncated object of $\calM$.
\end{remark}

\begin{definition}\label{sugarcoat}\index{gen}{$k$-truncated!map between Kan complexes}\index{gen}{truncated!map between Kan complexes}\index{gen}{$k$-truncated!morphism in an $\infty$-category}\index{gen}{truncated!morphism in an $\infty$-category}
We will say that a map $f: X \rightarrow Y$ of Kan complexes is {\it $k$-truncated} if the homotopy
fibers of $f$ (taken over any base point of $Y$) are $k$-truncated. We will say that a morphism
$f: C \rightarrow D$ in an arbitrary $\infty$-category $\calC$ is {\it $k$-truncated}
if composition with $f$ induces a $k$-truncated map
$\bHom_{\calC}(E,C) \rightarrow \bHom_{\calC}(E,D)$
for every object $E \in \calC$.
\end{definition}

\begin{remark}
There is an apparent potential for ambiguity in Definition \ref{sugarcoat} in the case where
$\calC$ is an $\infty$-category whose objects are Kan complexes. However, there is no cause for concern: a map $f: X \rightarrow Y$ of Kan complexes is $k$-truncated if and only if it is $k$-truncated as a morphism in the $\infty$-category $\SSet$.
\end{remark}

\begin{remark}
Let $f: C \rightarrow D$ and $g: E \rightarrow D$ be morphisms in an $\infty$-category $\calC$,
and let $\phi: \bHom_{\calC}(E,C) \rightarrow \bHom_{\calC}(E,D)$ be the map (in the homotopy category $\calH$) given by compostion with $f$. Lemma \ref{sugarplace} implies that the homotopy fiber of $\phi$ over $g$ is homotopy equivalent to 
$\bHom_{\calC_{/D}}(f,g)$. Consequently, we deduce that $f: C \rightarrow D$ is $k$-truncated in the sense of Definition \ref{sugarcoat} if and only if it is $k$-truncated when viewed as an object of the $\infty$-category $\calD_{/D}$.
\end{remark}

\begin{lemma}\label{truncslice}
Let $p: \calC \rightarrow \calD$ be a right fibration of $\infty$-categories, and let
$f: X \rightarrow Y$ be a morphism in $\calC$. Then $f$ is $n$-truncated if and only if
$p(f): p(X) \rightarrow p(Y)$ is $n$-truncated.
\end{lemma}

\begin{proof}
The map $\calC_{/Y} \rightarrow \calD_{/p(Y)}$ is
a trivial fibration, and therefore an equivalence of $\infty$-categories.
\end{proof}

\begin{remark}\label{tunc}
A morphism $f: C \rightarrow D$ in an $\infty$-category $\calC$ is  $k$-truncated if and only if it is $k$-truncated when regarded as an object of the $\infty$-category $\calC^{/D}$ (since the natural
map $\calC_{/D} \rightarrow \calC^{/D}$ is an equivalence of $\infty$-categories). We may identify $\calC^{/D}$ with $p^{-1} \{D \}$, where $p$ denotes the evaluation map
$\calC^{\Delta^1} \rightarrow \calC^{\{1\}}$. Corollary \ref{tweezegork} implies that $p$ is a coCartesian fibration. Consequently, Remark \ref{trumble} translates into the following assertion: if
$$ \xymatrix{ C' \ar[r]^{f'} \ar[d] & D' \ar[d] \\
C \ar[r]^{f} & D }$$
is a pullback diagram in $\calC$, and $f$ is $k$-truncated, then $f'$ is $k$-truncated.
\end{remark}

\begin{example}
A morphism $f: C \rightarrow D$ in an $\infty$-category $\calC$ is $(-2)$-truncated if and only if it is an equivalence.
\end{example}

We will say that a morphism $f: C \rightarrow D$ is a {\it monomorphism} if it is $(-1)$-truncated; this is equivalent to the assertion that the functor $\calC_{/f} \rightarrow \calC_{/D}$ is fully faithful.\index{gen}{monomorphism}

\begin{lemma}\label{trunccomp}
Let $\calC$ be an $\infty$-category and $f: X \rightarrow Y$ a morphism in $\calC$.
Suppose that $Y$ is $n$-truncated. Then $X$ is $n$-truncated if and only if $f$ is $n$-truncated.
\end{lemma}

\begin{proof}
Unwinding the definitions, we reduce immediately to the following statement in classical homotopy theory: given a map $f: X \rightarrow Y$ of Kan complexes, where $Y$ is $n$-truncated, $X$ is $n$-truncated if and only if the homotopy fibers of $f$ are $n$-truncated. This can be established easily, by examining the long exact sequence of homotopy groups associated to $f$.
\end{proof}

The following lemma gives a recursive characterization of the class of $n$-truncated morphisms:

\begin{lemma}\label{trunc}
Let $\calC$ be an $\infty$-category which admits finite limits and let $k
\geq -1$ be an integer. A morphism $f: C \rightarrow C'$ is $k$-truncated if and
only if the diagonal $\delta: C \rightarrow C \times_{C'} C$ (which is well-defined up to homotopy)
is $(k-1)$-truncated.
\end{lemma}

\begin{proof}
For each object $D \in \calC$, let $F_{D}: \calC \rightarrow \SSet$ denote the 
functor co-represented by $D$. Then each $F_{D}$ preserves finite limits, and
a morphism $f$ in $\calC$ is $k$-truncated if and only if each $F_D(f)$ is a $k$-truncated morphism in $\SSet$. We may therefore reduceo to the case where $\calC = \SSet$. Without loss of generality, we may suppose that $f: C \rightarrow C'$ is a Kan fibration. Then Theorem \ref{colimcomparee} allows us to identify the fiber product $C \times_{C'} C$ in $\SSet$ with the same fiber product, formed in the ordinary category $\Kan$. We now reduce to the following assertion in classical homotopy theory (applied to the fibers of $f$): if $X$ is a Kan complex, then $X$ is $k$-truncated if and only if the homotopy fibers of the diagonal map $X \rightarrow X \times X$ are $(k-1)$-truncated. This can be proven readily by examining homotopy groups.
\end{proof}

We immediately deduce the following:

\begin{proposition}\label{eaa}
Let $F: \calC \rightarrow \calC'$ be a left-exact functor between
$\infty$-categories which admit finite limits. Then $F$ carries
$k$-truncated objects into $k$-truncated objects and $k$-truncated
morphisms into $k$-truncated morphisms.
\end{proposition}

\begin{proof}
An object $C$ is $k$-truncated if and only if the morphism $C
\rightarrow 1$ to the final object is $k$-truncated. Since $F$
preserves final object, it suffices to prove the assertion
concerning morphisms. Since $F$ commutes with fiber products,
Lemma \ref{trunc} allows us to use induction on $k$, thereby
reducing to the case where $k=-2$. But the $(-2)$-truncated
morphisms are precisely the equivalences, and these are preserved
by any functor.
\end{proof}

We now specialize to the case of a {\em presentable} $\infty$-category $\calC$. In this setting, we can construct an analogue of the Postnikov tower.

\begin{lemma}\label{slurm}
Let $X$ be a Kan complex, and let $k \geq -2$. The following conditions are equivalent:
\begin{itemize}
\item[$(1)$] The Kan complex $X$ is $k$-truncated.
\item[$(2)$] The diagonal map $\delta: X \rightarrow X^{ \bd \Delta^{k+2} }$ is a homotopy equivalence.
\end{itemize}
\end{lemma}

\begin{proof}
If $k=-2$, then $X^{ \bd \Delta^{k+2} }$ is a point and the assertion is obvious. Assuming $k > -2$,
we can choose a vertex $v$ of $\bd \Delta^{k+2}$, which gives rise to an evaluation map
$e: X^{ \bd \Delta^{k+2} } \rightarrow X$. Since $e \circ \delta = \id_{X}$, $(2)$ is equivalent
to the assertion that $e$ is a homotopy equivalence. We observe that $e$ is a Kan fibration. For each $x$, let $Y_{x} = X^{ \bd \Delta^{k+2} } \times_{X} \{x\}$ denote the fiber of $e$ over the vertex $x$. Then $Y_x$ has a canonical base point, given by the constant map
$\delta(x)$. Moreover, we have a natural isomorphism
$$ \pi_{i}(Y_x, \delta(x) ) \simeq \pi_{i+k+1}(X,x).$$
Condition $(1)$ is equivalent to the assertion that $\pi_{i+k+1}(X,x)$ vanishes for
all $i \geq 0$ and all $x \in X$. In view of the above isomorphism, this is equivalent to the assertion that each $Y_{x}$ is contractible, which is true if and only if the Kan fibration $e$ is trivial.
\end{proof}

\begin{proposition}\label{maketrunc}
Let $\calC$ be a presentable $\infty$-category, $k \geq -2$. Let $\tau_{\leq k} \calC$
denote the full subcategory of $\calC$ spanned by the $k$-truncated objects.
Then the inclusion $\tau_{\leq k} \calC \subseteq \calC$ has an accessible left adjoint, which we will denote by $\tau_{\leq k}$.\index{not}{tauk@$\tau_{\leq k}$}
\end{proposition}

\begin{proof}
Let $f: \bd \Delta^{k+2} \rightarrow \Fun(\calC,\calC)$ denote the constant diagram
taking the value $\id_{\calC}$. Let $\overline{f}: (\bd \Delta^{k+2})^{\triangleright} \rightarrow \Fun(\calC,\calC)$ be a colimit of $f$, and let $F: \calC \rightarrow \calC$ be the image of the cone point under $\overline{f}$. Informally, $F$ is given by the formula
$$ C \mapsto C \otimes S^{k+1},$$
where $S^{k+1}$ denotes the $(k+1)$-sphere and we regard $\calC$ as tensored over spaces (see Remark \ref{tensored}).

Let $\overline{f}': (\bd \Delta^{k+2})^{\triangleright} \rightarrow \calC^{\calC}$ be the constant diagram taking the value $\id_{\calC}$. It follows that there exists an essentially unique map $\overline{f} \rightarrow \overline{f'}$ in $(\calC^{\calC})_{f/}$, which induces a natural transformation of functors $\alpha: F \rightarrow \id_{\calC}$. Let $S = \{ \alpha(C): C \in \calC \}$.
Since $F$ is a colimit of functors which preserve small colimits, $F$ itself preserves small colimits (Lemma \ref{limitscommute}). Applying Proposition \ref{limiteval}, we conclude that
$\alpha: \calC \rightarrow \Fun(\Delta^1,\calC)$ also preserves small colimits. Consequently, there
exists a small subset $S_0 \subseteq S$ which generates $S$ under colimits in $\Fun(\Delta^1,\calC)$.
According to Proposition \ref{local}, the collection of $S$-local objects of $\calC$ is an accessible localization of $\calC$. It therefore suffices to prove that an object $X \in \calC$ is $S$-local if and only if $X$ is $k$-truncated.

According to Proposition \ref{limiteval}, for each $C \in \calC$ we may identify $F(C)$
with the colimit of the constant diagram $\bd \Delta^{k+2} \rightarrow \calC$ taking the value $C$. Corollary \ref{charext} implies that we have a homotopy equivalence
$$ \bHom_{\calC}( F(C), X) \simeq \bHom_{\calC}(C,X)^{ \bd \Delta^{k+2} }.$$
The map $\alpha(C)$ induces a map
$$ \alpha(C)_{X}: \bHom_{\calC}(C,X) \rightarrow \bHom_{\calC}(C,X)^{ \bd \Delta^{k+2} }$$ which can be identified with the inclusion of $\bHom_{\calC}(C,X)$ as the space of constant maps from
$\bd \Delta^{k+2}$ into $\bHom_{\calC}(C,X)$. According to Lemma \ref{slurm}, the map
$\alpha(C)_{X}$ is an equivalence if and only if $\bHom_{\calC}(C,X)$ is $k$-truncated.
Thus $X$ is $k$-truncated if and only if $X$ is $S$-local. 
\end{proof}

\begin{remark}
The notation of Proposition \ref{maketrunc} is self-consistent, in the sense that the existence
of the localization functor $\tau_{\leq k}$ implies that the collection of $k$-truncated objects of $\calC$ may be identified with the essential image of $\tau_{\leq k}$.
\end{remark}

\begin{remark}
If the $\infty$-category $\calC$ is potentially unclear in
context, then we will write $\tau^{\calC}_{\leq k}$ for the truncation
functor in $\calC$. Note also that $\tau^{\calC}_{\leq k}$ is well-defined up to equivalence (in fact, up to a contractible ambiguity).\index{not}{tau^Ck@$\tau^{\calC}_{\leq k}$}
\end{remark}

\begin{remark}
It follows from Proposition \ref{maketrunc} that if $\calC$ is a presentable $\infty$-category, then the full subcategory $\tau_{\leq k} \calC$ of $k$-truncated objects is also presentable. In particular,
the ordinary category $\Disc(\calC)$ of discrete objects of $\calC$ is a presentable
category in the sense of Definition \ref{catpor}.
\end{remark}

Recall that, if $\calC$ and $\calD$ are $\infty$-categories, then $\LFun(\calC, \calD)$ denotes the full subcategory of $\Fun(\calC,\calD)$ spanned by those functors which are left adjoints.
The following result gives a characterization of $\tau_{\leq n} \calC$ by a universal mapping property:

\begin{corollary}\label{truncprop}
Let $\calC$ and $\calD$ be presentable $\infty$-categories. Suppose that $\calD$ is equivalent to an $(n+1)$-category. Then composition with $\tau_{\leq n}$ induces an equivalence
$s: \LFun( \tau_{\leq n} \calC, \calD) \rightarrow \LFun(\calC, \calD).$
\end{corollary}

\begin{proof}
According to Proposition \ref{unichar}, the functor $s$ is fully faithful. A functor $f: \calC \rightarrow \calD$ belongs to the essential image of $s$ if and only if $f$ has a right adjoint $g$ which
factors through $\tau_{\leq n} \calC$. Since $g$ preserves limits, it automatically carries
$\calD = \tau_{\leq n} \calD$ into $\tau_{\leq n} \calC$ (Proposition \ref{eaa}).
\end{proof}

In classical homotopy theory, every space $X$ can be recovered (up to weak homotopy equivalence) as the homotopy inverse limit of its Postnikov tower $\{ \tau_{\leq n} X \}_{n \geq 0}$.
The analogous statement is not true in an arbitrary presentable $\infty$-category, but often holds in specific examples. We now make a general study of this phenomenon.

\begin{definition}\label{postit1}\index{gen}{tower}\index{gen}{pretower}\index{gen}{Postnikov tower}\index{gen}{Postnikov pretower}
Let $\Z_{\geq 0}^{\infty}$ denote the union $\Z_{\geq 0} \cup \{ \infty \}$, regarded as a linearly ordered
set with largest element $\infty$. Let $\calC$ be a presentable $\infty$-category. Recall that
a {\it tower} in $\calC$ is a functor $\Nerve( \Z_{\geq 0}^{\infty})^{op} \rightarrow \calC$, which we view as a diagram
$$ X_{\infty} \rightarrow \ldots \rightarrow X_{2} \rightarrow X_{1} \rightarrow X_{0}.$$
A {\it Postnikov tower} is a tower with the property that for each $n \geq 0$, the map
$X_{\infty} \rightarrow X_{n}$ exhibits $X_{n}$ as an $n$-truncation of $X_{\infty}$. 
We define a {\it pretower} to be a functor from $\Nerve( \Z_{\geq 0})^{op} \rightarrow \calC$.
A {\it Postnikov pretower} is a pretower
$$ \ldots \rightarrow X_{2} \rightarrow X_{1} \rightarrow X_0$$
which exhibits each $X_{n}$ as an $n$-truncation of $X_{n+1}$.
We let $\Post^{+}(\calC)$ denote the full subcategory of $\Fun( \Nerve( \Z_{\geq 0}^{\infty})^{op}, \calC)$\index{not}{Post+C@$\Post^{+}(\calC)$}\index{not}{PostC@$\Post(\calC)$}
spanned by the Postnikov towers, and $\Post(\calC)$ the full subcategory of
$\Fun( \Nerve(\Z_{\geq 0})^{op}, \calC)$ spanned by the Postnikov pretowers.
We have an evident forgetful functor
$\phi: \Post^{+}(\calC) \rightarrow \Post(\calC).$ We will say that
{\it Postnikov towers in $\calC$ are convergent} if $\phi$ is an equivalence of $\infty$-categories.
\end{definition}

\begin{remark}\label{sumptime}
Let $\calC$ be a presentable $\infty$-category, and let
$\calE$ denote the full subcategory of $\calC \times \Nerve( \Z^{\infty}_{\geq 0})^{op}$
spanned by those pairs $(C, n)$ where $C \in \calC$ is $n$-truncated (by convention, we agree that
this condition is always satisfied where $C = \infty$). Then we have a coCartesian fibration
$p: \calE \rightarrow \Nerve( \Z^{\infty}_{\geq 0})^{op}$, which classifies a tower of $\infty$-categories
$$ \calC \rightarrow \ldots \rightarrow \tau_{\leq 2} \calC \stackrel{\tau_{\leq 1}}{\rightarrow} \tau_{\leq 1} \calC
\stackrel{\tau_{\leq 0}}{\rightarrow} \tau_{\leq 0} \calC.$$
We can identify Postnikov towers with coCartesian sections of $p$ (and Postnikov pretowers
with coCartesian sections of the induced fibration
$\calE \times_{ \Nerve( \Z^{\infty}_{\geq 0})^{op}} \Nerve(\Z_{\geq 0})^{op}
\rightarrow \Nerve( \Z_{\geq 0})^{op})$.
According to Proposition \ref{charcatlimit}, Postnikov towers in $\calC$ converge if and only if
the tower above exhibits $\calC$ as the homotopy limit of the sequence of $\infty$-categories
$$ \ldots \rightarrow \tau_{\leq 2} \calC \rightarrow \tau_{\leq 1} \calC \rightarrow
\tau_{\leq 0} \calC.$$
\end{remark}

\begin{remark}\label{swit}
Let $\calC$ be a presentable $\infty$-category, and assume that Postnikov towers in $\calC$
are convergent. Then every Postnikov tower in $\calC$ is a limit diagram. Indeed, given
objects $X, Y \in \calC$, we have natural homotopy equivalences
$$ \bHom_{\calC}( X,Y) \simeq \holim \bHom_{\calC}( \tau_{\leq n} X, \tau_{\leq n} Y)
\simeq \holim \bHom_{\calC}( X, \tau_{\leq n} Y)$$
so that $Y$ is the limit of the pretower $\{ \tau_{\leq n} Y\}$.
\end{remark}

\begin{proposition}\label{slibe1}
Let $\calC$ be a presentable $\infty$-category. Then Postnikov towers in $\calC$ are convergent if and only if, for every tower $X: \Nerve( \Z^{\infty}_{\geq 0})^{op} \rightarrow \calC$, 
the following conditions are equivalent:
\begin{itemize}
\item[$(1)$] The diagram $X$ is a Postnikov tower.
\item[$(2)$] The diagram $X$ is a limit in $\calC$, and the restriction
$X| \Nerve( \Z_{\geq 0})^{op}$ is a Postnikov pretower.
\end{itemize}
\end{proposition}

\begin{proof}
Let $\Post'(\calC)$ be the full subcategory of $\Fun( \Nerve( \Z^{\infty}_{\geq 0})^{op}, \calC)$ spanned by those towers which satisfy condition $(2)$. Using Proposition \ref{lklk}, we deduce that the restriction functor $\Post'(\calC) \rightarrow \Post(\calC)$ is a trivial Kan fibration. If
conditions $(1)$ and $(2)$ are equivalent, then $\Post'(\calC) = \Post^{+}(\calC)$, so that
Postnikov towers in $\calC$ are convergent. Conversely, suppose that
Postnikov towers in $\calC$ are convergent. Using Remark \ref{swit}, we deduce that
$\Post^{+}(\calC) \subseteq \Post'(\calC)$, so we have a commutative diagram
$$ \xymatrix{ \Post^{+}(\calC) \ar[dr] \ar[rr] & & \Post'(\calC) \ar[dl] \\
& \Post(\calC). & }$$
Since both of the vertical arrows are trivial Kan fibrations, we conclude that the inclusion
$\Post^{+}(\calC) \subseteq \Post'(\calC)$ is an equivalence, so that
$\Post^{+}(\calC) = \Post'(\calC)$. This proves that $(1) \Leftrightarrow (2)$.
\end{proof}

\begin{remark}\label{urkan}\index{gen}{tower!highly connected}
Let $\calC$ be a presentable $\infty$-category. We will say that a tower
$X: \Nerve( \Z^{\infty}_{\geq 0})^{op} \rightarrow \calX$ is {\it highly connected} if, for every $n \geq 0$, there exists an integer $k$ such that the induced map $\tau_{\leq n} X(\infty) \rightarrow \tau_{\leq n} X(k')$ is an equivalence for $k' \geq k$. We will say that a pretower 
$Y: \Nerve( \Z_{\geq 0}) \rightarrow
\calX$ is {\it highly connected} if, for every $n \geq 0$, there exists an integer $k$ such that
the map $\tau_{\leq n} Y(k'') \rightarrow \tau_{\leq n} Y(k')$ is an equivalence for $k'' \geq k' \geq k$.\index{gen}{pretower!highly connected}
It is clear that every Postnikov (pre)tower is highly connected. Conversely, if $X$ is a highly connected tower and its underlying pretower is highly connected, then $X$ is a Postnikov tower. Indeed, for each $n \geq 0$ we can choose $k \geq n$ such that the map
$\tau_{\leq n} X(\infty) \rightarrow \tau_{\leq n} X(k)$ is an equivalence. Since $X$ is a Postnikov pretower, this induces an equivalence $\tau_{\leq n} X(\infty) \simeq X(n)$. Consequently, to establish
the implication $(2) \Rightarrow (1)$ in the criterion of Proposition \ref{slibe1}, it suffices to verify the following:
\begin{itemize}
\item[$(\ast)$] Let $X: \Nerve( \Z^{\infty}_{\geq} ) \rightarrow \calC$ be a tower in $\calC$. Assume that $X$ is a limit diagram, and that the underlying pretower is highly connected. Then $X$ is highly connected.
\end{itemize}
In \S \ref{homdim}, we will apply this criterion to prove that Postnikov towers are convergent in a large class of $\infty$-topoi.
\end{remark}

We conclude this section with a useful compatibility property between truncation functors in different $\infty$-categories:

\begin{proposition}\label{compattrunc}
Let $\calC$ and $\calD$ be presentable $\infty$-categories, and let $F: \calC \rightarrow \calD$
be a left-exact presentable functor. Then there is an equivalence of functors
$F \circ \tau_{\leq k}^{\calC} \simeq \tau_{\leq k}^{\calD} \circ F$.
\end{proposition}

\begin{proof}
Since $F$ is left exact, it restricts to a functor from $\tau_{\leq k} \calC$ to $\tau_{\leq k} \calD$ by Proposition \ref{eaa}. We therefore have a diagram
$$ \xymatrix{ \calC \ar[r]^{F} \ar[d]^{\tau_{\leq k}^{\calC} } & \calD \ar[d]^{ \tau_{\leq k}^{\calD}} \\
\tau_{\leq k} \calC \ar[r]^{F} & \tau_{\leq k} \calD }$$
which we wish to prove is commutative up to homotopy. Let $G$ denote a right adjoint to
$F$; then $G$ is left exact and so induces a functor $\tau_{\leq k} \calD
\rightarrow \tau_{\leq k} \calC$. Using Proposition \ref{compadjoint}, we can reduce to proving
that the associated diagram of right adjoints 
$$ \xymatrix{ \calC  & \calD \ar[l]_{G} \\
\tau_{\leq k} \calC \ar[u] & \tau_{\leq k} \calD \ar[u] \ar[l]_{G} }$$
commutes up to homotopy, which is obvious (since the diagram strictly commutes).
\end{proof}

\subsection{Compactly Generated $\infty$-Categories}\label{compactgen}

\begin{definition}\index{gen}{compactly generated}\index{gen}{$\kappa$-compactly generated}\label{compgen}
Let $\kappa$ be a regular cardinal. We will say that an
$\infty$-category $\calC$ is {\it $\kappa$-compactly generated} if it is presentable and $\kappa$-accessible. When $\kappa = \omega$, we will simply say that $\calC$ is {\it compactly generated}.
\end{definition}

The proof of Theorem \ref{pretop} shows that an $\infty$-category $\calC$ is $\kappa$-compactly generated if and only if there exists a small $\infty$-category $\calD$ which admits $\kappa$-small colimits, and an equivalence $\calC \simeq \Ind_{\kappa}(\calD)$. In fact, we can choose $\calD$ to be (a minimal model of)
the $\infty$-category of $\kappa$-compact objects of $\calC$. We would like to assert that this construction establishes an equivalence between two sorts of $\infty$-categories. In order to make this precise, we need to introduce the appropriate notion of functor between $\kappa$-compactly generated $\infty$-categories.

\begin{proposition}\label{comppress}
Let $\kappa$ be a regular cardinal, and let $\Adjoint{F}{\calC}{\calD}{G}$ be a pair of adjoint functors, where $\calC$ and $\calD$ admit small, $\kappa$-filtered colimits.
\begin{itemize}
\item[$(1)$] If $G$ is $\kappa$-continuous, then $F$ carries $\kappa$-compact objects of $\calC$ to $\kappa$-compact objects of $\calD$.
\item[$(2)$] Conversely, if $\calC$ is $\kappa$-accessible and $F$ preserves $\kappa$-compactness, then $G$ is $\kappa$-continuous.
\end{itemize}
\end{proposition}

\begin{proof}
Suppose first that $G$ is $\kappa$-continuous, and let $C \in \calC$ be a $\kappa$-compact object. Let $e: \calC \rightarrow \widehat{\SSet}$ be a functor corepresented by $C$. Then
$e \circ G: \calD \rightarrow \widehat{\SSet}$ is corepresented by $F(C)$. Since $e$ and $G$
are $\kappa$-continuous, so is $e \circ G$; this proves $(1)$.

Conversely, suppose that $F$ preserves $\kappa$-compact objects and that $\calC$ is $\kappa$-accessible. Without loss of generality, we may suppose that there is a small $\infty$-category $\calC'$ such that $\calC = \Ind_{\kappa}(\calC') \subseteq \calP(\calC')$. We wish to prove
that $G$ is $\kappa$-continuous. Since $\Ind_{\kappa}(\calC')$ is stable under $\kappa$-filtered colimits in $\calP(\calC')$, it will suffice to prove that the composite map
$$ \theta: \calD \stackrel{G}{\rightarrow} \calC \subseteq \calP(\calC') $$
is $\kappa$-continuous. In view of Proposition \ref{limiteval}, it will suffice to prove that
for every object $C \in \calC'$, the composition of $\theta$ with evaluation at $C$ is
a $\kappa$-continuous functor. We conclude by observing that this functor is corepresentable
by the image under $F$ of $j(C) \in \calC$ (here $j: \calC' \rightarrow \Ind_{\kappa}(\calC)$ denotes the Yoneda embedding). 
\end{proof}

\begin{corollary}\label{starmin}
Let $\calC$ be a $\kappa$-compactly generated $\infty$-category, and let
$L: \calC \rightarrow \calC$ be a localization functor. The following conditions are equivalent:
\begin{itemize}
\item[$(1)$] The functor $L$ is $\kappa$-continuous.
\item[$(2)$] The full subcategory $L \calC \subseteq \calC$ is stable under $\kappa$-filtered colimits.
\end{itemize}
Suppose that these conditions are satisfied. Then:
\begin{itemize}
\item[$(3)$] The functor $L$ carries $\kappa$-compact objects of $\calC$ to $\kappa$-compact objects
of $L \calC$.

\item[$(4)$] The $\infty$-category $L \calC$ is $\kappa$-compactly generated.

\item[$(5)$] An object $D \in L \calC$ is $\kappa$-compact (in $L \calC$) if and only if
there exists a compact object $C \in \calC$ such that $D$ is a retract of $LC$.
\end{itemize}

\end{corollary}

\begin{proof}
Suppose that $(1)$ is satisfied. Let $p: K \rightarrow L \calC$ be a $\kappa$-filtered diagram. Then the natural transformation $p \rightarrow Lp$ is an equivalence. Using $(1)$, we conclude that
the induced map $\varinjlim(p) \rightarrow L \varinjlim(p)$ is an equivalence, so that
$\varinjlim(p) \in L \calC$. This proves $(2)$.

Conversely, if $(2)$ is satisfied, then the inclusion $L \calC \subseteq \calC$ is $\kappa$-continuous, so that $L: \calC \rightarrow \calC$ is a composition of $\kappa$-continuous functors
$$ \calC \stackrel{L}{\rightarrow} L \calC \rightarrow \calC,$$
which proves $(1)$.

Assume that $(1)$ and $(2)$ are satisfied. Then $L$ is accessible, so that
$L \calC$ is a presentable $\infty$-category. 
Assertion $(3)$ follows from Proposition \ref{comppress}.
Let $D \in L \calC$. Since $\calC$ is $\kappa$-compactly generated, 
$D$ can be written as the colimit of a $\kappa$-filtered diagram $p: K \rightarrow \calC$ taking values in the $\kappa$-compact objects of $\calC$. Then $D \simeq LD$ can be written
as the colimit of $L \circ p$, which takes values $\kappa$-compact objects of $L \calC$. This proves $(4)$. If $D$ is a $\kappa$-compact object of $\calD$, then we deduce that the identity map
$\id_{D}: D \rightarrow D$ factors through $(L \circ p)(k)$ for some vertex $k \in K$, which proves $(5)$.
\end{proof}

\begin{corollary}\label{hunterygreen}
Let $\calC$ be a $\kappa$-compactly generated $\infty$-category, and let
$n \geq -2$. Then:
\begin{itemize}
\item[$(1)$] The full subcategory $\tau_{\leq n} \calC$ is stable under $\kappa$-filtered colimits in $\calC$.
\item[$(2)$] The truncation functor $\tau_{\leq n}: \calC \rightarrow \calC$ is $\kappa$-continuous.
\item[$(3)$] The truncation functor $\tau_{\leq n}$ carries compact objects of $\calC$ to compact objects of $\calC_{\leq n}$. 
\item[$(4)$] The full subcategory $\tau_{\leq n} \calC$ is $\kappa$-compactly generated.
\item[$(5)$] An object $C \in \tau_{\leq n} \calC$ is compact (in $\tau_{ \leq n} \calC$) if and only if there exists a compact object $C' \in \calC$ such that $C$ is a retract of $\tau_{\leq n} C'$.
\end{itemize}
\end{corollary}

\begin{proof}
Corollary \ref{starmin} shows that condition $(1)$ implies $(2)$, $(3)$, $(4)$, and $(5)$.
Consequently, it will suffice to prove that $(1)$ is satisfied. 

Let $C$ be an object of $\calC$. We will show that $C$ is $n$-truncated if and only if the space
$\bHom_{\calC}(D, C)$ is $n$-truncated, for every $\kappa$-compact object $D \in \calC$. 
The ``only if'' direction is obvious. For the converse, let $F_{C}: \calC^{op} \rightarrow \SSet$ be the functor represented by $C$, and let $\calC' \subseteq \calC$ be the full subcategory of $\calC$ spanned by those objects $D$ such that $F_{C}(D)$ is $n$-truncated. Since $F_{C}$ preserves limits, $\calC'$ is stable under colimits in $\calC$. If $\calC'$ contains every $\kappa$-compact object of $\calC$, then $\calC' = \calC$ (since $\calC$ is $\kappa$-compactly generated).

Now suppose that $D$ is a $\kappa$-compact object of $\calC$, let
$G_{D}: \calC \rightarrow \SSet$ be the functor co-represented by $D$, and let
$\calC(D) \subseteq \calC$ be the full subcategory of $\calC$ spanned by those objects
$C$ for which $G_{D}(C)$ is $n$-truncated. Then $\tau_{\leq n} \calC = \bigcap_{D} \calC(D)$. To complete the proof, it will suffice to show that each $\calC(D)$ is stable under $\kappa$-filtered colimits. Since $G_{D}$ is $\kappa$-continuous, it suffices to observe that
$\tau_{\leq n} \SSet$ is stable under $\kappa$-filtered colimits in $\SSet$.
\end{proof}

\begin{definition}\index{not}{PresRkappa@$\RRPres{\kappa}$}
If $\kappa$ is a regular cardinal, we let $\RRPres{\kappa}$ denote the full subcategory
of $\widehat{\Cat}_{\infty}$ whose objects are $\kappa$-compactly generated $\infty$-categories, and whose morphisms are $\kappa$-continuous, limit-preserving functors.
\end{definition}

\begin{proposition}\label{cnote}\index{gen}{limit!of compactly generated $\infty$-categories}
The $\infty$-category $\RRPres{\kappa}$ admits small limits, and the inclusion
$\RRPres{\kappa} \subseteq \widehat{\Cat}_{\infty}$ preserves small limits.
\end{proposition}

\begin{proof}
In view of Theorem \ref{surbus}, the only nontrivial point is to verify that if $p: K \rightarrow \RRPres{\kappa}$ is a diagram of $\kappa$-compactly generated $\infty$-categories $\{ \calC_{\alpha} \}$, then the limit $\calC = \varprojlim(p)$ in $\widehat{\Cat}_{\infty}$ is $\kappa$-compactly generated. In other words, we must show that $\calC$ is generated under colimits by its $\kappa$-compact objects.

For each vertex $\alpha$ of $K$, let 
$$ \Adjoint{ F_{\alpha} }{\calC_{\alpha}}{\calC}{G_{\alpha} }$$ denote the corresponding adjunction. Lemma \ref{steakknife} implies that the identity functor $\id_{\calC}$ can be obtained as the colimit of a diagram $q: K \rightarrow \Fun(\calC, \calC)$, where $q(\alpha) \simeq F_{\alpha} \circ G_{\alpha}$. In particular, $\calC$ is generated (under small colimits) by the essential images of the functors $F_{\alpha}$. Since each $\calC_{\alpha}$ is generated under colimits by $\kappa$-compact objects, and the functors $F_{\alpha}$ preserve colimits and $\kappa$-compact objects (Proposition \ref{comppress}), we conclude that $\calC$ is generated under colimits by its $\kappa$-compact objects, as desired.
\end{proof} 

\begin{notation}\label{funnote}\index{not}{PresLkappa@$\LLPres{\kappa}$}
Let $\kappa$ be a regular cardinal. We let $\LLPres{\kappa}$ denote the full subcategory of
$\widehat{\Cat}_{\infty}$ whose objects are $\kappa$-compactly generated $\infty$-categories, and whose morphisms are functors which preserve small colimits and $\kappa$-compact objects. In view of Proposition \ref{comppress}, the equivalence
$ \LPres \simeq (\RPres)^{op}$ of Corollary \ref{warhog} restricts to an equivalence
$ \LLPres{\kappa} \simeq (\RRPres{\kappa})^{op}$.

Let $\widehat{\Cat}_{\infty}^{\Rex{\kappa}}$ denote the subcategory of $\widehat{\Cat}_{\infty}$ whose objects are $\infty$-categories which admit $\kappa$-small colimits, and whose morphisms are functors which preserve $\kappa$-small colimits, and let $\Cat_{\infty}^{\Rex{\kappa}} = \widehat{\Cat}_{\infty}^{\Rex{\kappa}} \cap \Cat_{\infty}$.\index{not}{catrexhat@$\widehat{\Cat}_{\infty}^{\Rex{\kappa}}$}\index{not}{catrex@$\Cat_{\infty}^{\Rex{\kappa}}$}
\end{notation}

\begin{proposition}\label{suchy}
Let $\kappa$ be a regular cardinal, and let
$$\theta: \LLPres{\kappa} \rightarrow \widehat{\Cat}_{\infty}^{\Rex{\kappa}}$$ be the nerve of the simplicial functor which associates to a $\kappa$-compactly generated $\infty$-category 
$\calC$ the full subcategory $\calC^{\kappa} \subseteq \calC$ spanned by the $\kappa$-compact objects of $\calC$. Then:
\begin{itemize}
\item[$(1)$] The functor $\theta$ is fully faithful. 
\item[$(2)$] The essential image of $\theta$ consists precisely of those objects of $\widehat{\Cat}_{\infty}$ which are essentially small and idempotent complete. 
\end{itemize}
\end{proposition}

\begin{proof}
Combine Propositions \ref{humatch} and \ref{sumatch}.
\end{proof}

\begin{remark}
If $\kappa > \omega$, then Corollary \ref{swwe} shows that the hypothesis of idempotent completeness in $(2)$ is superfluous.
\end{remark}

The proof of Proposition \ref{lockap} yields the following analogue:

\begin{proposition}\label{sumer}
Let $\kappa$ be a regular cardinal. The functor $\Ind_{\kappa}: \Cat_{\infty} \rightarrow \Acc_{\kappa}$ exhibits $\LLPres{\kappa}$ as a localization of $\Cat_{\infty}^{\Rex{\kappa}}$. 
If $\kappa > \omega$, then $\Ind_{\kappa}$ induces an equivalence of $\infty$-categories
$\Cat_{\infty}^{\Rex{\kappa}} \rightarrow \LLPres{\kappa}$.
\end{proposition}

\begin{proof}
The only additional ingredient needed is the following observation: if $\calC$ is an $\infty$-category which admits $\kappa$-small colimits, then the idempotent completion $\calC'$ of $\calC$ also admits $\kappa$-small colimits. To prove this, we observe that $\calC'$ can be identified with the collection of $\kappa$-compact objects of $\Ind_{\kappa}(\calC)$ (Lemma \ref{stylus}). Since $\calC$ admits all small colimits (Theorem \ref{pretop}), we conclude that $\calC'$ admits $\kappa$-small colimits.
\end{proof}

We conclude with a remark about the structure of the $\infty$-category
$\Cat_{\infty}^{\Rex{\kappa}}$.

\begin{proposition}\label{unrose}\index{gen}{filtered colimit!of colimit-preserving functors}
Let $\kappa$ be a regular cardinal. Then the $\infty$-category
$\Cat_{\infty}^{\Rex{\kappa}}$ admits small, $\kappa$-filtered colimits, and the inclusion 
$\Cat_{\infty}^{\Rex{\kappa}} \subseteq \Cat_{\infty}$ preserves small $\kappa$-filtered colimits.
\end{proposition}

\begin{proof}
Let $\calI$ be a small, $\kappa$-filtered $\infty$-category, and let
$p: \calI \rightarrow \Cat_{\infty}^{\Rex{\kappa}}$ be a diagram. Let $\calC$ be a colimit of the induced diagram $\calI \rightarrow \Cat_{\infty}$. 
To complete the proof we must prove the following:
\begin{itemize}
\item[$(i)$] The $\infty$-category $\calC$ admits $\kappa$-small colimits.
\item[$(ii)$] For each $I \in \calI$, the associated functor $p(I) \rightarrow \calC$ preserves $\kappa$-small colimits.
\item[$(iii)$] Let $f: \calC \rightarrow \calD$ be an arbitrary functor. If each of the compositions
$p(I) \rightarrow \calC \rightarrow \calD$ preserves $\kappa$-small colimits, then $f$ preserves $\kappa$-small colimits.
\end{itemize}

Since $\calI$ is $\kappa$-filtered, any $\kappa$-small diagram in $\calC$ factors through one of the maps $p(I) \rightarrow \calC$ ( Proposition \ref{grapeape} ). Thus $(ii) \Rightarrow (i)$ and $(ii) \Rightarrow (iii)$. To prove $(ii)$, we first use Proposition \ref{rot} to reduce to the case where $\calI \simeq \Nerve(A)$, where $A$ is a $\kappa$-filtered partially ordered set. Using Proposition \ref{gumby444}, we can reduce to the case where $p$ is the nerve of a functor from $q: A \rightarrow \sSet$. In view of Theorem \ref{colimcomparee}, we can identify $\calC$ with a homotopy colimit of $q$. Since the collection of categorical equivalences is stable under filtered colimits, we can assume that $\calC$ is actually the filtered colimit of a family of $\infty$-categories $\{ \calC_{\alpha} \}_{\alpha \in A}$. 

Let $K$ be a $\kappa$-small simplicial set, and let $\overline{g}_{\alpha}: K^{\triangleright} \rightarrow \calC_{\alpha}$ be a colimit diagram. We wish to show that the induced map 
$\overline{g}: K^{\triangleright} \rightarrow \calC$ is a colimit diagram. Let
$g = \overline{g} |K$; we need to show that the map $\theta: \calC_{\overline{g}/} \rightarrow \calC_{g/}$ is a trivial Kan fibration. We observe that $\theta$ is a filtered colimit of maps
$\theta_{\beta}: (\calC_{\beta})_{\overline{g}_{\beta}/} \rightarrow (\calC_{\beta})_{g_{\beta}/}$, where $\beta$ ranges over the set $\{ \beta \in A : \beta \geq \alpha \}$. Using the fact that each of the associated maps $\calC_{\alpha} \rightarrow \calC_{\beta}$ preserves $\kappa$-small colimits, we conclude that each $\theta_{\beta}$ is a trivial fibration, so that $\theta$ is a trivial fibration as desired.
\end{proof}

\subsection{Nonabelian Derived Categories}\label{stable11}

According to Corollary \ref{uterrr}, we can analyze arbitrary colimits in an $\infty$-category $\calC$ in terms of finite colimits and filtered colimits. In particular, suppose that $\calC$ admits finite colimits and that we construct new $\infty$-category $\Ind(\calC)$ by formally adjoining filtered colimits to $\calC$. Then $\Ind(\calC)$ admits all small colimits (Theorem \ref{pretop}), and the Yoneda embedding $\calC \rightarrow \Ind(\calC)$ preserves finite colimits (Proposition \ref{turnke}). 
Moreover, we can identify $\Ind(\calC)$ with the $\infty$-category of functors $\calC^{op} \rightarrow \SSet$ which carry finite colimits in $\calC$ to finite limits in $\SSet$. In this section, we will introduce a variation on the same theme. Instead of assuming $\calC$ admits {\em all} finite colimits, we will only assume that $\calC$ admits finite coproducts. We will construct a coproduct-preserving embedding of $\calC$ into a larger $\infty$-category $\calP_{\Sigma}(\calC)$ which admits all small colimits. Moreover, we can characterize $\calP_{\Sigma}(\calC)$ as the $\infty$-category obtained from $\calC$ by formally adjoining colimits of {\em sifted diagrams} (Proposition \ref{surottt}). 

Our first goal in this section is to introduce the notion of a {\em sifted} simplicial set. We begin with a bit of motivation. Let $\calC$ denote the (ordinary) category of groups. Then $\calC$ admits arbitrary colimits. However, colimits of diagrams in $\calC$ can be very difficult to analyze, even if the diagram itself is quite simple. For example, the coproduct of a pair of groups $G$ and $H$ is the {\it amalgamated product} $G \star H$. The group $G \star H$ is typically very complicated, even when $G$ and $H$ are not. For example, the amalgamated product
$\Z / 2 \Z \star \Z / 3 \Z$ is isomorphic to arithmetic group $\PSL_2(\Z)$. 
In general, $G \star H$ is much larger than the coproduct
$G \coprod H$ of the underlying sets of $G$ and $H$. In other words, the forgetful functor $U: \calC \rightarrow \Set$ does not preserve coproducts. However, $U$ does preserve {\em some} colimits: for example, the colimit of a sequence of groups
$$ G_0 \rightarrow G_1 \rightarrow \ldots $$
can be obtained by taking the colimit of the underlying sets, and equipping the result with an appropriate group structure.

The forgetful functor $U$ from groups to sets preserves another important type of colimit: namely, the formation of quotients by equivalence relations. If $G$ is a group, then a subgroup
$R \subseteq G \times G$ is an equivalence relation on $G$ if and only if there exists
a normal subgroup $H \subseteq G$ such that $R = \{ (g,g'): g^{-1} g' \in H \}$. In this case,
the set of $R$-equivalence classes in $G$ is in bijection with the quotient $G/H$, which inherits a group structure from $G$. In other words, the quotient of $G$ by the equivalence relation $R$ can be computed either in the category of groups or the category of sets; the result is the same.

Each of the examples given above admits a generalization: the colimit of a sequence is a special case of a {\it filtered colimit}, and the quotient by an equivalence relation is a special case of a {\it reflexive coequalizer}. The forgetful functor $\calC \rightarrow \Set$ preserves filtered colimits and reflexive coequalizers; moreover, the same is true if we replace the category of groups by any other category of sets with some sort of finitary algebraic structure (for example, abelian groups, or commutative rings). The following definition, which is taken from \cite{homotopyvarieties}, is an attempt to axiomatize the essence of the situation:\index{gen}{coequalizer!reflexive}\index{gen}{reflexive coequalizer}

\begin{definition}[\cite{homotopyvarieties}]\label{siftdef}\index{gen}{sifted}\index{gen}{simplicial set!sifted}
A simplicial set $K$ is {\it sifted} if it satisfies the following conditions:
\begin{itemize}
\item[$(1)$] The simplicial set $K$ is nonempty.
\item[$(2)$] The diagonal map $K \rightarrow K \times K$ is cofinal.
\end{itemize}
\end{definition}

\begin{warning}
In \cite{homotopyvarieties}, Rosicki uses the term {\it homotopy sifted} to describe the analogue of Definition \ref{siftdef} for simplicial categories, and reserves the term {\it sifted} for analogous notion in the setting of ordinary categories. There is some danger of confusion with our terminology:
if $\calC$ is an ordinary category and $\Nerve(\calC)$ is sifted (in the sense of Definition \ref{siftdef}), then $\calC$ is sifted in the sense of \cite{homotopyvarieties}. However, the converse is false in general.
\end{warning}

\begin{example}\label{bin1}
Every filtered $\infty$-category is sifted (Proposition \ref{undertruck}). 
\end{example}

\begin{lemma}\label{bball3}
The simplicial set $\Nerve( \cDelta)^{op}$ is sifted.
\end{lemma}

\begin{proof}
Since $\Nerve(\cDelta)^{op}$ is clearly nonempty, it will suffice to show that the diagonal map
$\Nerve( \cDelta)^{op} \rightarrow \Nerve( \cDelta)^{op} \times \Nerve(\cDelta)^{op}$ is cofinal.
According to Theorem \ref{hollowtt}, this is equivalent to the assertion that for every object
$([m], [n]) \in \cDelta \times \cDelta$, the category
$$ \calC = \cDelta_{/ [m]} \times_{ \cDelta} \cDelta_{/ [n] }$$
has weakly contractible nerve. Let $\calC^{0}$ be the full subcategory of $\calC$
spanned by those objects which correspond to {\em monomorphisms} of partially ordered sets $J \rightarrow [m] \times [n]$. The inclusion of $\calC^{0}$ into $\calC$
has a left adjoint, so the inclusion $\Nerve(\calC^{0}) \subseteq \Nerve(\calC)$ is a weak homotopy equivalence. It will therefore suffice to show that $\Nerve(\calC^{0})$ is weakly contractible. We now observe that $\Nerve(\calC^{0})$ can be identified with the first barycentric subdivision of
$\Delta^m \times \Delta^n$, and is therefore weakly homotopy equivalent to 
$\Delta^m \times \Delta^n$ and so weakly contractible.
\end{proof}

\begin{remark}
The formation of geometric realization of simplicial objects should be thought of as the $\infty$-categorical analogue of the formation of reflexive coequalizers.
\end{remark}

Our next pair of results captures some of the essential features of the theory of sifted simplicial sets:

\begin{proposition}\label{urbil}
Let $K$ be a sifted simplicial set, let
$\calC$, $\calD$, and $\calE$ be $\infty$-categories which admit $K$-indexed colimits, and let
$f: \calC \times \calD \rightarrow \calE$ be a map which preserves $K$-indexed colimits separately in each variable. Then $f$ preserves $K$-indexed colimits. 
\end{proposition}

\begin{proof}
Let $p: K \rightarrow \calC$ and $q: K \rightarrow \calD$ be diagrams indexed by a small simplicial set $K$, and let $\delta: K \rightarrow K \times K$ be the diagonal map. Using the fact that 
$f$ preserves $K$-indexed colimits separately in each variable and Lemma \ref{limitscommute}, we conclude that $\colim( f \circ (p \times q) )$ is a colimit for the diagram $f \circ (p \times q) \circ \delta$. Consequently, $f$ preserves $K$-indexed colimits provided that the diagonal $\delta$ is cofinal. We conclude by invoking the assumption that $K$ is sifted.
\end{proof}

\begin{proposition}\label{siftcont}
Let $K$ be a sifted simplicial set. Then $K$ is weakly contractible.
\end{proposition}

\begin{proof}
Choose a vertex $x$ in $K$. According to Whitehead's theorem, it will suffice to show that for each
$n \geq 0$, the homotopy set $\pi_n( |K|,x)$ consists of a single element. Let $\delta: K \rightarrow K \times K$ be the diagonal map. Since $\delta$ is cofinal, Proposition \ref{cofbasic} implies that
the induced map
$$ \pi_n( |K|, x) \rightarrow \pi_n( |K \times K|, \delta(x) ) \simeq \pi_n(|K|,x) \times \pi_n(|K|,x)$$
is bijective. Since $\pi_n( |K|,x)$ is nonempty, we conclude that it is a singleton.
\end{proof}

We now return to the problem introduced in the beginning of this section.

\begin{definition}\label{vardef}
Let $\calC$ be a small $\infty$-category which admits finite coproducts. We let
$\calP_{\Sigma}(\calC)$ denote the full subcategory of $\calP(\calC)$ spanned by
those functors $\calC^{op} \rightarrow \SSet$ which preserve finite products.
\end{definition}\index{not}{PSigmaC@$\calP_{\Sigma}(\calC)$}

\begin{remark}
The $\infty$-categories of the form $\calP_{\Sigma}(\calC)$ have been studied in
\cite{homotopyvarieties}, where they are called {\it homotopy varieties}. Many of the results proven below can also be found in \cite{homotopyvarieties}.\index{gen}{homotopy varieties}
\end{remark}

\begin{proposition}\label{utut}
Let $\calC$ be a small $\infty$-category which admits finite coproducts. Then:
\begin{itemize}
\item[$(1)$] The $\infty$-category $\calP_{\Sigma}(\calC)$ is an accessible localization
of $\calP(\calC)$. 
\item[$(2)$] The Yoneda embedding $j: \calC \rightarrow \calP(\calC)$ factors
through $\calP_{\Sigma}(\calC)$. Moreover, $j$ carries finite coproducts in $\calC$
to finite coproducts in $\calP_{\Sigma}(\calC)$. 
\item[$(3)$] Let $\calD$ be a presentable $\infty$-category, and let
$$ \Adjoint{ F}{\calP(\calC)}{\calD}{G}$$
be a pair of adjoint functors. Then $G$ factors through $\calP_{\Sigma}(\calC)$ if and only if
$f = F \circ j: \calC \rightarrow \calD$ preserves finite coproducts.
\item[$(4)$] The full subcategory $\calP_{\Sigma}(\calC) \subseteq \calP(\calC)$ is stable
under sifted colimits. 
\item[$(5)$] Let $L: \calP(\calC) \rightarrow \calP_{\Sigma}(\calC)$ be a left adjoint to the inclusion. Then $L$ preserves sifted colimits $($when regarded as a functor from $\calP(\calC)$ to itself$)$.
\item[$(6)$] The $\infty$-category $\calP_{\Sigma}(\calC)$ is compactly generated.
\end{itemize}
\end{proposition}

Before giving the proof, we need a preliminary result concerning the interactions between products sifted colimits.

\begin{lemma}\label{bale2}
Let $K$ be a sifted simplicial set. Let $\calX$ be an $\infty$-category which admits finite products and
$K$-indexed colimits, and suppose that the formation of products in $\calX$ preserves $K$-indexed colimits separately in each variable. Then the colimit functor $\colim: \Fun(K, \calX) \rightarrow \calX$
preserves finite products.
\end{lemma}

\begin{remark}\label{bale3}
The hypotheses of Lemma \ref{bale2} are satisfied when $\calX$ is the $\infty$-category $\SSet$ of spaces: see Lemma \ref{sugartime}. More generally, Lemma \ref{bale2} applies whenever
the $\infty$-category $\calX$ is an $\infty$-topos (see Definition \ref{def1topos}).
\end{remark}

\begin{proof}
Since the simplicial set $K$ is weakly contractible (Proposition \ref{siftcont}), Corollary \ref{charext} implies that the functor $\colim$ preserves final objects. To complete the proof, it will suffice to show that the functor $\colim$ preserves pairwise products. Let $X$ and $Y$ be objects of $\Fun(K, \calX)$. We wish to prove that the canonical map
$$ \colim( X \times Y) \rightarrow \colim(X) \times \colim(Y)$$ is an equivalence.
In other words, we must show that the formation of products commutes with $K$-indexed colimits, which follows immediately by applying Proposition \ref{urbil} to the Cartesian product functor
$\calX \times \calX \rightarrow \calX$.
\end{proof}

\begin{proof}[Proof of Proposition \ref{utut}]
Assertion $(1)$ is an immediate consequence of Lemmas \ref{stur2}, \ref{stur3}, and \ref{stur1}. To prove $(2)$, it will suffice to show that for every representable functor
$e: \calP_{\Sigma}(\calC)^{op} \rightarrow \SSet$, the composition
$$ \calC^{op} \stackrel{j^{op}}{\rightarrow} \calP_{\Sigma}(\calC)^{op} \stackrel{e}{\rightarrow} \SSet$$
preserves finite products (Proposition \ref{yonedaprop}). This is obvious, since the composition can be identified with the object of $\calP_{\Sigma}(\calC) \subseteq
\Fun( \calC^{op}, \SSet)$ representing $e$.

We next prove $(3)$. We note that $f$ preserves finite coproducts if and only if, for
every object $D \in \calD$, the composition
$$ \calC^{op} \stackrel{f^{op}}{\rightarrow} \calD^{op} \stackrel{e_D}{\rightarrow} \SSet$$
preserves finite products, where $e_D$ denotes the functor represented by $D$. 
This composition can be identified with $G(D)$, so that $f$ preserves finite coproducts if and only if
$G$ factors through $\calP_{\Sigma}(\calC)$.

Assertion $(4)$ is an immediate consequence of Lemma \ref{bale2} and Remark \ref{bale3}, and $(5)$ follows formally from $(4)$. To prove $(6)$, we first observe that
$\calP(\calC)$ is compactly generated (Proposition \ref{precst}). Let $\calE
\subseteq \calP(\calC)$ be the full subcategory spanned by the compact objects, and let
$L: \calP(\calC) \rightarrow \calP_{\Sigma}(\calC)$ be a localization functor. Since
$\calE$ generates $\calP(\calC)$ under filtered colimits, $L(\calD)$ generates
$\calP_{\Sigma}(\calC)$ under filtered colimits. Consequently, it will suffice to
show that for each $E \in \calE$, the object $LE \in \calP_{\Sigma}(\calC)$ is compact. Let
$f: \calP_{\Sigma}(\calC) \rightarrow \SSet$ be the functor corepresented by $LE$, and let
$f': \calP(\calC) \rightarrow \SSet$ be the functor corepresented by $E$. Then
the map $E \rightarrow LE$ induces an equivalence
$f \rightarrow f' | \calP_{\Sigma}(\calC)$. Since $f'$ is continuous and
$\calP_{\Sigma}(\calC)$ is stable under filtered colimits in $\calP(\calC)$, we conclude
that $f$ is continuous, so that $LE$ is a compact object of $\calP_{\Sigma}(\calC)$ as desired.
\end{proof}

Our next goal is to prove a converse to part $(4)$ of Proposition \ref{utut}. Namely, we will show that $\calP_{\Sigma}(\calC)$ is generated by the essential image of the Yoneda embedding
under sifted colimits. In fact, we will only need to use special types of sifted colimits: namely, filtered colimits and geometric realizations (Lemma \ref{subato}). The proof is based on the following technical result:

\begin{lemma}\label{presubato}
Let $\calC$ be a small $\infty$-category, and let $X$ be an object of $\calP(\calC)$. Then
there exists a simplicial object $Y_{\bigdot}: \Nerve(\cDelta)^{op} \rightarrow \calP(\calC)$ with the following properties:
\begin{itemize}
\item[$(1)$] The colimit of $Y_{\bigdot}$ is equivalent to $X$.
\item[$(2)$] For each $n \geq 0$, the object $Y_n \in \calP(\calC)$ is equivalent to a small coproduct of objects lying in the image of the Yoneda embedding $j: \calC \rightarrow \calP(\calC)$. 
\end{itemize}
\end{lemma}

We will defer the proof until the end of this section.

\begin{lemma}\label{subato}
Let $\calC$ be a small $\infty$-category which admits finite coproducts, and let
$X \in \calP(\calC)$. The following conditions are equivalent:
\begin{itemize}
\item[$(1)$] The object $X$ belongs to $\calP_{\Sigma}(\calC)$.
\item[$(2)$] There exists a simplicial object $U_{\bigdot}: \Nerve( \cDelta)^{op} \rightarrow \Ind(\calC)$ whose colimit in $\calP(\calC)$ is $X$.
\end{itemize}
\end{lemma}

\begin{proof}
The full subcategory $\calP_{\Sigma}(\calC)$ contains the essential image of the Yoneda embedding and is stable under filtered colimits and geometric realizations (Proposition \ref{utut}); thus $(2) \Rightarrow (1)$. We will prove that $(1) \Rightarrow (2)$. 

We first choose a simplicial object $Y_{\bigdot}$ of $\calP(\calC)$ which satisfies the conclusions of
Lemma \ref{presubato}. Let $L$ be a left adjoint to the inclusion $\calP_{\Sigma}(\calC) \subseteq \calP(\calC)$. Since $X$ is a colimit of $Y_{\bigdot}$, $LX \simeq X$ is a colimit of $LY_{\bigdot}$ (part $(5)$ of Proposition \ref{utut}). It will therefore suffice to prove that each $LY_{n}$ belongs to $\Ind(\calC)$. By hypothesis, each $Y_{n}$ can be written as a small coproduct $\coprod_{\alpha \in A} j(C_{\alpha} )$, where
$j: \calC \rightarrow \calP(\calC)$ denotes the Yoneda embedding. 
Using the results of \S \ref{quasilimit1}, we see that $Y_{n}$ can be obtained also as a filtered colimit of coproducts $\coprod_{ \alpha \in A_0} j(C_{\alpha})$, where $A_0$ ranges over the finite subsets of $A$. Since $L$ preserves filtered colimits (Proposition \ref{utut}), it will suffice to show that each of the objects
$$L(\coprod_{\alpha \in A_0} j(C_{\alpha} ))$$
belongs to $\Ind(\calC)$. We now invoke part $(2)$ of Proposition \ref{utut} to identify this object with $j( \coprod_{\alpha \in A_0} C_{\alpha} )$. 
\end{proof}

\begin{proposition}\label{surottt}
Let $\calC$ be a small $\infty$-category which admits finite coproducts, and let
$\calD$ be an $\infty$-category which admits filtered colimits and geometric realizations.
Let $\Fun_{\Sigma}( \calP_{\Sigma}(\calC), \calD)$ denote the full subcategory spanned by those functors $\calP_{\Sigma}(\calC) \rightarrow \calD$ which preserve filtered colimits and geometric realizations. Then:
\begin{itemize}
\item[$(1)$] Composition with the Yoneda embedding $j: \calC \rightarrow \calP_{\Sigma}(\calC)$ induces an equivalence of categories
$$ \theta: \Fun_{\Sigma}( \calP_{\Sigma}(\calC), \calD) \rightarrow \Fun(\calC, \calD).$$

\item[$(2)$] Any functor $g \in \Fun_{\Sigma}( \calP_{\Sigma}(\calC), \calD)$ preserves sifted colimits.

\item[$(3)$] Assume that $\calD$ admits finite coproducts. A functor $g \in \Fun_{\Sigma}( \calP_{\Sigma}(\calC), \calD)$ preserves small colimits if and only if $g \circ j$ preserves finite coproducts.\index{not}{FunSigma@$\Fun_{\Sigma}( \calC, \calD)$}
\end{itemize}
\end{proposition}

\begin{proof}
Lemma \ref{subato} and Proposition \ref{utut} imply that $\calP_{\Sigma}(\calC)$ is the smallest full subcategory of $\calP(\calC)$ which is closed under filtered colimits, closed under geometric realizations, and contains the essential image of the Yoneda embedding. Consequently, assertion $(1)$ follows from Remark \ref{poweryoga} and Proposition \ref{lklk}. 

We now prove $(2)$. Let $g \in \Fun_{\Sigma}( \calP_{\Sigma}(\calC), \calD)$; we wish to show that $g$ preserves sifted colimits. It will suffice to show that for every representable functor
$e: \calD \rightarrow \SSet^{op}$, the composition $e \circ g$ preserves sifted colimits. In other words, we may replace $\calD$ by $\SSet^{op}$, and thereby reduce to the case where
$\calD$ itself admits sifted colimits. Let $\Fun'_{\Sigma}( \calP_{\Sigma}(\calC), \calD)$ denote the full subcategory of $\Fun_{\Sigma}( \calP_{\Sigma}(\calC), \calD)$ spanned by those functors which preserve sifted colimits. Since $\calP_{\Sigma}(\calC)$ is also the smallest full subcategory of $\calP(\calC)$ which contains the essential image of the Yoneda embedding and is stable under sifted colimits, Remark \ref{poweryoga} implies that $\theta$ induces an equivalence
$$ \Fun'_{\Sigma}( \calP_{\Sigma}(\calC), \calD ) \rightarrow \Fun(\calC, \calD).$$
Combining this observation with $(1)$, we deduce that the inclusion
$\Fun'_{\Sigma}( \calP_{\Sigma}(\calC), \calD) \subseteq \Fun_{\Sigma}( \calP_{\Sigma}(\calC), \calD)$
is an equivalence of $\infty$-categories, and therefore an equality.

The ``only if'' direction of $(3)$ is immediate, since the Yoneda embedding $j: \calC \rightarrow \calP_{\Sigma}(\calC)$ preserves finite coproducts (Proposition \ref{utut}). To prove the converse,
we first apply Lemma \ref{diverti} to reduce to the case where $\calD$ is a full subcategory of an $\infty$-category $\calD'$, with the following properties:
\begin{itemize}
\item[$(i)$] The $\infty$-category $\calD'$ admits small colimits.
\item[$(ii)$] A small diagram $K^{\triangleright} \rightarrow \calD$ is a colimit if and only if the induced diagram $K^{\triangleright} \rightarrow \calD'$ is a colimit.
\end{itemize}
Let $\calC'$ denote the essential image of the Yoneda embedding $j: \calC \rightarrow \calP(\calC)$.
Using Lemma \ref{longwait1}, we conclude that there exists functor $G: \calP(\calC) \rightarrow \calD'$ which is a left Kan extension of $G | \calC' = g | \calC'$, and that $G$ preserves small colimits. 
Let $G_0 = G | \calP_{\Sigma}(\calC)$. Then $G_0$ is a left Kan extension of $g | \calC'$, so there
is a canonical natural transformation $G_0 \rightarrow g$. Let $\calC''$ denote the full subcategory
of $\calP_{\Sigma}(\calC)$ spanned by those objects $C$ for which the map $G_0(C) \rightarrow g(C)$
is an equivalence. Then $\calC''$ contains $\calC'$ and is stable under filtered colimits and geometric realizations, and therefore contains all of $\calP_{\Sigma}(\calC)$. We may therefore replace
$g$ by $G_0$ and thereby assume that $G | \calP_{\Sigma}(\calC) = g$.
Since $G \circ j = g \circ j$ preserves finite coproducts, the right adjoint to $G$ factors through
$\calP_{\Sigma}(\calC)$ (Proposition \ref{utut}), so that $G$ is equivalent to the composition
$$ \calP(\calC) \stackrel{L}{\rightarrow} \calP_{\Sigma}(\calC) \stackrel{G'}{\rightarrow} \calD'$$
for some colimit-preserving functor $G': \calP_{\Sigma}(\calC) \rightarrow \calD'$. Restricting to
the subcategory $\calP_{\Sigma}(\calC) \subseteq \calP(\calC)$, we deduce that
$G'$ is equivalent to $g$, so that $g$ preserves small colimits as desired.
\end{proof}

\begin{remark}\label{spuduse}
Let $\calC$ be a small $\infty$-category which admits finite coproducts. It follows from
Proposition \ref{surottt} that we can identify $\calP_{\Sigma}(\calC)$ with
$\calP_{\calK}^{\calK'}(\calC)$ in each of the following three cases (for
an explanation of this notation, we refer the reader to \S \ref{agileco}):
\begin{itemize}
\item[$(1)$] The collection $\calK$ is empty, and the collection $\calK'$ consists of all
small filtered simplicial sets together with $\Nerve(\cDelta)^{op}$.
\item[$(2)$] The collection $\calK$ is empty, and the collection $\calK'$ consists of all
small sifted simplicial sets.
\item[$(3)$] The collection $\calK$ consists of all finite discrete simplicial sets, and
the collection $\calK'$ consists of all small simplicial sets.
\end{itemize}
\end{remark}

\begin{corollary}\label{swillt}
Let $f: \calC \rightarrow \calD$ be a functor between $\infty$-categories. Assume that $\calC$ admits small colimits. Then $f$ preserves sifted colimits if and only if $f$ preserves filtered colimits and geometric realizations.
\end{corollary}

\begin{proof}
The ``only if'' direction is clear. For the converse, suppose that $f$ preserves filtered colimits and geometric realizations. Let $\calI$ be a small sifted $\infty$-category and $\overline{p}: \calI^{\triangleright} \rightarrow \calC$ a colimit diagram; we wish to prove that $f \circ \overline{p}$ is also a colimit diagram. Let $p = \overline{p} | \calI$. Let $\calJ \subseteq \calP(\calI)$ denote a small full subcategory which contains the essential image of the Yoneda embedding $j: \calI \rightarrow \calP(\calI)$ and is closed under finite coproducts. It follows from Remark \ref{poweryoga} that
the functor $p$ is homotopic to a composition $q \circ j$, where $q: \calJ \rightarrow \calC$ is a functor which preserves finite coproducts. Proposition \ref{surottt} implies that
$q$ is homotopic to a composition
$$ \calJ \stackrel{j'}{\rightarrow} \calP_{\Sigma}(\calJ) \stackrel{q'}{\rightarrow} \calC,$$
where $j'$ denotes the Yoneda embedding and $q'$ preserves small colimits. 
The composition $f \circ q'$ preserves filtered colimits and geometric realization, and therefore
preserves sifted colimits (Proposition \ref{surottt}). 

Let $\overline{p}': \calI^{\triangleright} \rightarrow \calP_{\Sigma}(\calJ)$ be a colimit of the diagram
$j' \circ j$. Since $q'$ preserves colimits, the composition $q' \circ \overline{p}'$ is a colimit of
$q' \circ j' \circ j \simeq p$, and is therefore equivalent to $\overline{p}$. Consequently, it will suffice to show that $f \circ q' \circ \overline{p}'$ is a colimit diagram. Since $\calI$ is sifted, we need only
verify that $f \circ q'$ preserves sifted colimits. By Proposition \ref{surottt}, it will suffice to show that
$f \circ q'$ preserves filtered colimits and geometric realizations. Since $q'$ preserves all colimits, this follows from our assumption that $f$ preserves filtered colimits and geometric realizations.
\end{proof}

In the situation of Proposition \ref{surottt}, every functor $f: \calC \rightarrow \calD$ extends (up to homotopy) to a functor $F: \calP_{\Sigma}(\calC) \rightarrow \calD$, which preserves sifted colimits. We will sometimes refer to $F$ as the {\it left derived functor} of $f$\index{gen}{derived functor}\index{gen}{left derived functor}. In \S \ref{stable12} we will explain the connection of this notion of derived functor with the more classical definition provided by Quillen's theory of homotopical algebra.\index{gen}{derived functor}\index{gen}{functor!derived}

Our next goal is to characterize those $\infty$-categories which have the form $\calP_{\Sigma}(\calC)$.

\begin{definition}\label{humber}
Let $\calC$ be an $\infty$-category which admits geometric realizations of simplicial objects. We will say that an object $P \in \calC$ is {\it projective} if the functor $\calC \rightarrow \SSet$ co-represented by $P$ commutes with geometric realizations.\index{gen}{projective object}\index{gen}{object!projective}
\end{definition}

\begin{remark}\label{untine}
Let $\calC$ be an $\infty$-category which admits geometric realizations of simplicial objects. Then
the collection of projective objects of $\calC$ is stable under all finite coproducts which exist in $\calC$.
This follows immediately from Lemma \ref{bale2} and Remark \ref{bale3}. 
\end{remark}

\begin{remark}\label{comproj}
Let $\calC$ be an $\infty$-category which admits small colimits, and let $X$ be an object of
$\calC$. Then $X$ is compact and projective if and only if $X$ corepresents a functor
$\calC \rightarrow \sSet$ which preserves sifted colimits. The ``only if'' direction is obvious, and
the converse follows from Corollary \ref{swillt}. 
\end{remark}

\begin{example}
Let $\calA$ be an abelian category. Then an object $P \in \calA$ is projective in the sense of classical homological algebra (that is, the functor $\Hom_{\calA}(P, \bigdot)$ is exact) if and only if $P$ corepresents a functor $\calA \rightarrow \Set$ which commutes with geometric realizations of simplicial objects. This is {\em not} equivalent to the condition of Definition \ref{humber}, since the fully faithful embedding $\Set \rightarrow \SSet$ does not preserve geometric realizations. However, it is equivalent to the requirement that $P$ be a projective object (in the sense of Definition \ref{humber}) in the $\infty$-category underlying the homotopy theory of simplicial objects
of $\calA$ (equivalently, the theory of nonpositively graded chain complexes with values in $\calA$;
we will discuss this example in greater detail in \cite{DAG}).
\end{example}

\begin{proposition}\label{smearof}
Let $\calC$ be a small $\infty$-category which admits finite coproducts, 
$\calD$ an $\infty$-category which admits filtered colimits and geometric realizations,
and $F: \calP_{\Sigma}(\calC) \rightarrow \calD$ a left derived functor of
$f = F \circ j: \calC \rightarrow \calD$, where $j: \calC \rightarrow \calP_{\Sigma}(\calC)$ denotes the Yoneda embedding. Consider the following conditions:
\begin{itemize}
\item[$(i)$] The functor $f$ is fully faithful.
\item[$(ii)$] The essential image of $f$ consists of compact projective objects of $\calD$.
\item[$(iii)$] The $\infty$-category $\calD$ is generated by the essential image of $f$ under filtered colimits and geometric realizations.
\end{itemize}
If $(i)$ and $(ii)$ are satisfied, then $F$ is fully faithful. Moreover, $F$ is an equivalence if
and only if $(i)$, $(ii)$, and $(iii)$ are satisfied.
\end{proposition}

\begin{proof}
If $F$ is an equivalence of $\infty$-categories, then $(i)$ follows from Proposition \ref{fulfaith}, and $(iii)$ from Lemma \ref{subato}. To prove $(ii)$, it suffices to show that for
each $C \in \calC$, the functor $e: \calP_{\Sigma}(\calC) \rightarrow \SSet$ corepresented by $C$ preserves filtered colimits and geometric realizations. We can identify $e$ with the composition
$$ \calP_{\Sigma}(\calC) \stackrel{e'}{\subseteq} \calP(\calC) \stackrel{e''}{\rightarrow} \SSet,$$
where $e''$ denotes evaluation at $C$. It now suffices to observe that $e'$ and $e''$
preserve filtered colimits and geometric realizations (Lemma \ref{subato} and Proposition \ref{limiteval}).

For the converse, let us suppose that $(i)$ and $(ii)$ are satisfied. We will show that $F$ is fully faithful. First fix an object $C \in \calC$, and let $\calP'_{\sigma}(\calC)$ be the full subcategory
of $\calP_{\Sigma}(\calC)$ spanned by those objects $M$ for which the map
$$ \bHom_{ \calP_{\Sigma}(\calC)}( j(C), M) \rightarrow \bHom_{\calD}(f(C), F(M) )$$
is an equivalence. Condition $(i)$ implies that $\calP'_{\sigma}(\calC)$ contains the essential image of $j$, and condition $(ii)$ implies that $\calP'_{\sigma}(\calC)$ is stable under filtered colimits and geometric realizations. Lemma \ref{subato} now implies that $\calP'_{\Sigma}(\calC) = \calP_{\Sigma}(\calC)$.

We now define $\calP''_{\Sigma}(\calC)$ to be the full subcategory of $\calP_{\Sigma}(\calC)$ spanned by those objects $M$ such that for all $N \in \calP_{\Sigma}(\calC)$, the map
$$ \bHom_{ \calP_{\Sigma}(\calC)}(M,N) \rightarrow \bHom_{\calD}( F(M), F(N) )$$
is a homotopy equivalence. The above argument shows that $\calP''_{\Sigma}(\calC)$ contains
the essential image of $j$. Since $F$ preserves filtered colimits and geometric realizations,
$\calP''_{\Sigma}(\calC)$ is stable under filtered colimits and geometric realizations. Applying Lemma \ref{subato}, we conclude that $\calP''_{\Sigma}(\calC) = \calP_{\Sigma}(\calC)$; this proves that $F$ is fully faithful.

If $F$ is fully faithful, then the essential image of $F$ contains $f(\calC)$ and is stable under filtered colimits and geometric realizations. If $(iii)$ is satisfied, it follows that $F$ is an equivalence of $\infty$-categories.
\end{proof}

\begin{definition}\label{defpro}
Let $\calC$ be an $\infty$-category which admits small colimits, and let $S$ be a collection of objects of $\calC$. We will say that $S$ is a {\it set of compact projective generators for $\calC$} if the following conditions are satisfied:
\begin{itemize}
\item[$(1)$] Each element of $S$ is a compact projective object of $\calC$.
\item[$(2)$] The full subcategory of $\calC$ spanned by the elements of $S$ is essentially small.
\item[$(3)$] The set $S$ generates $\calC$ under small colimits.
\end{itemize}
We will say that $\calC$ is {\it projectively generated} if there exists a set $S$ of compact projective generators for $\calC$.\index{gen}{projectively generated}\index{gen}{generator!projective} 
\end{definition}

\begin{example}\label{swine}
The $\infty$-category $\SSet$ of spaces is projectively generated. The compact projective objects of $\SSet$ are precisely those spaces which are homotopy equivalent to finite sets (endowed with the discrete topology).
\end{example}

\begin{proposition}\label{protus}
Let $\calC$ be an $\infty$-category which admits small colimits, and let $S$ be a set of compact projective generators for $\calC$. Then:
\begin{itemize}
\item[$(1)$] Let $\calC^{0} \subseteq \calC$ be the full subcategory spanned by finite coproducts of the objects $S$, let $\calD \subseteq \calC^{0}$ be a minimal model for $\calC^{0}$, and let
$F: \calP_{\Sigma}( \calD ) \rightarrow \calC$ be a left derived functor of the inclusion. Then $F$ is an equivalence of $\infty$-categories. In particular, $\calC$ is a compactly generated presentable $\infty$-category.
\item[$(2)$] Let $C \in \calC$ be an object. The following conditions are equivalent:
\begin{itemize}
\item[$(i)$] The object $C$ is compact and projective.
\item[$(ii)$] The functor $e: \calC \rightarrow \widehat{\SSet}$ corepresented by $C$
preserves sifted colimits.
\item[$(iii)$] There exists an object $C' \in \calC^{0}$ such that $C$ is a retract of $C'$.
\end{itemize}
\end{itemize}
\end{proposition}

\begin{proof}
Remark \ref{untine} implies that $\calC^{0}$ consists of compact projective objects of $\calC$.
Assertion $(1)$ now follows immediately from Proposition \ref{smearof}. We now prove $(2)$.
The implications $(iii) \Rightarrow (i)$ and $(ii) \Rightarrow (i)$ are obvious. To complete the proof, we will show that
$(i) \Rightarrow (iii)$. Using $(1)$, we are free to assume $\calC = \calP_{\Sigma}(\calD)$.
Let $C \in \calC$ be a compact projective object.
Using Lemma \ref{subato}, we conclude that there exists a simplicial object
$X_{\bigdot}$ of $\Ind(\calD)$ and an equivalence $C \simeq | X_{\bigdot} |$. Since
$C$ is projective, we deduce $\bHom_{\calC}(C,C)$ is equivalent to the geometric realization of the simplicial space $\bHom_{\calC}(C, X_{\bigdot})$. In particular, $\id_{C} \in \bHom_{\calC}(C,C)$ is homotopic to the image of the some map $f: C \rightarrow X_0$. Using our assumption that
$C$ is compact, we conclude that $f$ factors as a composition
$$ C \stackrel{f_0}{\rightarrow} j(D) \rightarrow X_0,$$
where $j: \calD \rightarrow \Ind(\calD)$ denotes the Yoneda embedding. It follows that
$C$ is a retract of $j(D)$ in $\calC$, as desired.
\end{proof}

\begin{remark}\label{parei}
Let $\calC$ be a small $\infty$-category which admits finite coproducts. Since
the truncation functor $\tau_{\leq n}: \SSet \rightarrow \SSet$ preserves finite products,
it induces a map $\tau: \calP_{\Sigma}(\calC) \rightarrow \calP_{\Sigma}(\calC)$, which is easily seen to be a localization functor. The essential image of $\tau$ consists of those functors
$F \in \calP_{\Sigma}(\calC)$ which take $n$-truncated values. We claim that these
are precisely the $n$-truncated object of $\calP_{\Sigma}(\calC)$. Consequently, we can
identify $\tau$ with the $n$-truncation functor on $\calP_{\Sigma}(\calC)$.

One direction is clear: if $F \in \calP_{\Sigma}(\calC)$ is $n$-truncated, then for each $C \in \calC$ the space $\bHom_{ \calP_{\Sigma}(\calC)}( j(C), F) \simeq F(C)$ must be $n$-truncated. Conversely,
suppose that $F: \calC^{op} \rightarrow \SSet$ takes $n$-truncated values. We wish to prove that
the space $\bHom_{ \calP_{\Sigma}(\calC)}(F', F)$ is $n$-truncated, for each $F' \in \calP_{\Sigma}(\calC)$. The collection of all objects $F'$ which satisfy this condition is stable under small colimits in $\calP_{\Sigma}(\calC)$ and contains the essential image of the Yoneda embedding. It therefore contains the entirety of $\calP_{\Sigma}(\calC)$, as desired.
\end{remark} 

We conclude this section by giving the proof of Lemma \ref{presubato}. Our argument uses some
concepts and results from \S \ref{chap6}, and may be omitted at first reading.

\begin{proof}[Proof of Lemma \ref{presubato}]
For $n \geq 0$, let $\cDelta^{\leq n}$ denote the full subcategory of $\cDelta$ spanned by the objects
$\{ [k] \}_{k \leq n}$. We will construct a compatible sequence of functors
$f_{n}: \Nerve( \cDelta^{\leq n})^{op} \rightarrow \calP(\calC)_{/X}$
with the following properties:
\begin{itemize}
\item[$(A)$] For $n \geq 0$, let $L_n$ denote a colimit of the composite diagram
$$ \Nerve( \cDelta^{\leq n-1})^{op} \times_{ \Nerve(\cDelta)^{op} }
\Nerve( \cDelta_{[n]/})^{op} \rightarrow \Nerve( \cDelta^{\leq n-1})^{op}
\stackrel{ f_{n-1}}{\rightarrow} \calP(\calC)_{/X} \rightarrow \calP(\calC).$$
(the $n$th {\it latching object}). Then there exists an object $Z_{n} \in \calP(\calC)$ which
is a small coproduct of objects in the essential image of the Yoneda embedding
$\calC \rightarrow \calP(\calC)$, and a map $Z_{n} \rightarrow f_{n}([n])$ which,
together with the canonical map $L_n \rightarrow f_{n}([n])$, determines an equivalence
$L_n \coprod Z_{n} \simeq f_{n}([n])$. 

\item[$(B)$] For $n \geq 0$, let $\overline{M}_n$ denote the limit of the diagram
$$ \Nerve( \cDelta^{\leq n-1})^{op} \times_{ \Nerve(\cDelta)^{op} }
\Nerve( \cDelta_{/[n]})^{op} \rightarrow \Nerve( \cDelta^{\leq n-1})^{op}
\stackrel{ f_{n-1}}{\rightarrow} \calP(\calC)_{/X}$$
(the $n$th {\it matching object}), and let $M_n$ denote its image in $\calP(\calC)$.
Then the canonical map $f_{n}([n]) \rightarrow M_n$ is an effective epimorphism
in $\calP(\calC)$ (see \S \ref{surjsurj}).
\end{itemize}
The construction of the functors $f_{n}$ proceeds by induction on $n$,
the case $n < 0$ being trivial. For $n \geq 0$, we invoke Remark \ref{stapler2}:
to extend $f_{n-1}$ to a functor $f_{n}$ satisfying $(A)$ and $(B)$, it suffices to produce
an object $Z_{n}$ and a morphism $\psi: Z_{n} \rightarrow M_n$ in $\calP(\calC)$, such that
the coproduct $L_n \coprod Z_{n} \rightarrow M_n$ is an effective epimorphism.
This is satisfied in particular if $\psi$ itself is an effective epimorphism.

The maps $f_{n}$ together determine a simplicial object $\overline{Y}_{\bigdot}$ of
$\calP(\calC)_{/X}$, which we can identify with a simplicial object $Y_{\bigdot}$
in $\calP(\calC)$ equipped with a map $\theta: \varinjlim Y_{\bigdot} \rightarrow X$.
Assumption $(B)$ guarantees that $\theta$ is a hypercovering of $X$ (see \S \ref{hcovh}),
so that the map $\theta$ is $\infty$-connective (Lemma \ref{fierminus}).
The $\infty$-topos $\calP(\calC)$ has enough points (given by evaluation at objects of $\calC$), and is therefore hypercomplete (Remark \ref{notenough}). It follows that $\theta$ is an equivalence.
We now complete the proof by observing that for $n \geq 0$, we have an equivalence
$Y_{n} \simeq \coprod_{ [n] \rightarrow [k] } Z_{k}$
where the coproduct is taken over all surjective maps of linearly ordered sets $[n] \rightarrow [k]$, so that
$Y_{n}$ is itself a small coproduct of objects lying in the essential image of the Yoneda embedding
$j: \calC \rightarrow \calP(\calC)$.
\end{proof}
  
\subsection{Quillen's Model for $\calP_{\Sigma}(\calC)$}\label{stable12} 
 
Let $\calC$ be a small category which admits finite products. Then $\Nerve(\calC)^{op}$ is an $\infty$-category which admits finite coproducts. In \S \ref{stable11}, we studied the $\infty$-category $\calP_{\Sigma}( \Nerve(\calC)^{op} )$, which we can view as the full subcategory of
$\Fun( \Nerve(\calC), \SSet)$ spanned by those functors which preserve finite products. According to Proposition \ref{gumby444}, $\Fun( \Nerve(\calC), \SSet)$ can be identified with the $\infty$-category underlying the simplicial model category of diagrams $\Set_{\Delta}^{\calC}$ (which we will endow with the {\em projective} model structure described in \S \ref{compp4}). It follows that every functor
$f: \Nerve(\calC) \rightarrow \SSet$ is equivalent to the (simplicial) nerve of a functor $F: \calC \rightarrow \Kan$. Moreover, $f$ belongs to $\calP_{\Sigma}( \Nerve(\calC)^{op})$ if and only if the functor
$F$ is {\em weakly} product preserving, in the sense that for any finite collection of objects $\{ C_i \in \calC \}_{1 \leq i \leq n}$, the natural map $$F( C_1 \times \ldots C_n ) \rightarrow F(C_1) \times \ldots \times F(C_n)$$ is a homotopy equivalence of Kan complexes. Our goal in this section is to prove a refinement of Proposition \ref{gumby444}: if $f$ preserves finite products, then it is possible to arrange that $F$ preserves finite products (up to isomorphism, rather than up to homotopy equivalence). This result is most naturally phrased as an equivalence between model categories (Proposition \ref{trent}), and is due to Bergner (see \cite{bergner3}). We begin by recalling the following result of Quillen (for a proof, we refer the reader to
\cite{homotopicalalgebra}):

\begin{proposition}[Quillen]\label{sutcoat}
Let $\calC$ be a category which admits finite products, and let
$\bfA$ denote the category of functors $F: \calC \rightarrow \sSet$
which preserve finite products. Then $\bfA$ has the structure of a simplicial model category, where:
\begin{itemize}
\item[$(W)$] A natural transformation $\alpha: F \rightarrow F'$ of functors
is a weak equivalence in $\bfA$ if and only if $\alpha(C): F(C) \rightarrow F'(C)$ is a weak
homotopy equivalence of simplicial sets, for every $C \in \calC$. 
\item[$(F)$] A natural transformation $\alpha: F \rightarrow F'$ of functors
is a fibration in $\bfA$ if and only if $\alpha(C): F(C) \rightarrow F'(C)$ is a Kan fibration of simplicial sets, for every $C \in \calC$. 
\end{itemize}
\end{proposition}

Suppose that $\calC$ and $\bfA$ are as in the statement of Proposition \ref{sutcoat}. Then we may regard $\bfA$ as a full subcategory of the category $\Set_{\Delta}^{\calC}$ of {\em all} functors
from $\calC$ to $\sSet$, which we regard as endowed with the projective model structure (so that fibrations and weak equivalences are given pointwise). The inclusion $G:\bfA \subseteq \Set_{\Delta}^{\calC}$ preserves fibrations and trivial fibrations, and therefore determines a Quillen adjunction
$$\Adjoint{F}{ \Set_{\Delta}^{\calC}}{\bfA}{G}.$$
(A more explicit description of the adjoint functor $F$ will be given below.) Our goal in this section is to prove the following result:

\begin{proposition}[Bergner]\label{trent}
Let $\calC$ be a small category which admits finite products, and let
$$ \Adjoint{ F }{ \Set_{\Delta}^{\calC} }{ \bfA }{ G}$$
be as above. Then the right derived functor
$$ RG: \h{\bfA} \rightarrow \h{\Set_{\Delta}^{\calC}}$$
is fully faithful, and an object $f \in \h{\Set_{\Delta}^{\calC}}$ belongs to the essential image of $RG$ if and only if $f$ preserves finite products up to weak homotopy equivalence.
\end{proposition}

\begin{corollary}\label{smokerr}
Let $\calC$ be a small category which admits finite products, and let $\bfA$ be as in Proposition \ref{trent}. Then the natural map $\Nerve( \bfA^{\degree} ) \rightarrow \calP_{\Sigma}( \Nerve(\calC)^{op} )$
is an equivalence of $\infty$-categories.
\end{corollary}

The proof of Proposition \ref{trent} is somewhat technical and will occupy the rest of this section. We begin by introducing some preliminaries.

\begin{notation}\label{bignote}
Let $\calC$ be a small category. We define a pair of categories $\Env(\calC) \subseteq \Env^{+}(\calC)$ as follows:\index{not}{EnvC@$\Env(\calC)$}\index{not}{EnvpC@$\Env^{+}(\calC)$}

\begin{itemize}
\item[$(i)$] An object of $\Env^{+}(\calC)$ is a pair $C = (J, \{ C_j\}_{j \in J} )$, where $J$ is a finite set and each $C_j$ is an object of $\calC$. The object $C$ belongs to $\Env(\calC)$ if and only if $J$ is nonempty. 

\item[$(ii)$] Given objects $C = (J, \{ C_j \}_{j \in J})$ and $C' = (J', \{ C'_{j'} \}_{j' \in J' } )$
of $\Env^{+}(\calC)$, a morphism $C \rightarrow C'$ consists of the following data:
\begin{itemize}
\item[$(a)$] A map $f: J' \rightarrow J$ of finite sets.
\item[$(b)$] For each $j' \in J'$, a morphism $C_{f(j')} \rightarrow C'_{j'}$ in the category $\calC$.
\end{itemize}
Such a morphism belongs to $\Env(\calC)$ if and only if $J$ and $J'$ are nonempty, and $f$ is surjective.
\end{itemize}

There is a fully faithful embedding functor $\theta: \calC \rightarrow \Env(\calC)$, given by
$C \mapsto ( \ast, \{ C \} )$. We can view $\Env^{+}(\calC)$ as the category
obtained from $\calC$ by freely adjoining finite products. In particular, if $\calC$ admits finite products, then $\theta$ admits a (product-preserving) left inverse $\phi^{+}_{\calC}$, given by the formula $( J, \{ C_j \}_{j \in J} ) \mapsto \prod_{j \in J} C_j.$. We let $\phi_{\calC}$ denote the restricton $\phi^{+}_{\calC} | \Env(\calC)$. 

Given a functor $\calF \in \Set_{\Delta}^{\calC}$, we let $E^{+}(\calF) \in
\Set_{\Delta}^{\Env^{+}(\calC)}$ denote the composition
$$ \Env^{+}(\calC) \stackrel{ \Env^{+}(\calF) }{\rightarrow} \Env^{+}(\sSet) \stackrel{ \phi^{+}_{\sSet}}{\rightarrow} \sSet$$
$$ (J, \{ C_j\}_{j \in J} ) \mapsto \prod f(C_j).$$
We let $E(\calF)$ denote the restriction $E^{+}(\calF) | \Env(\calC) \in \Set_{\Delta}^{\Env(\calC)}$.

If the category $\calC$ admits finite products, then we let
$L, L^{+}: \Set_{\Delta}^{\calC} \rightarrow \Set_{\Delta}^{\calC}$ denote the compositions 
$$ \Set_{\Delta}^{\calC} \stackrel{E}{\rightarrow} \Set_{\Delta}^{\Env(\calC)}
\stackrel{ (\phi_{\calC})_{!} }{\rightarrow} \Set_{\Delta}^{\calC}$$
$$ \Set_{\Delta}^{\calC} \stackrel{E^{+}}{\rightarrow} \Set_{\Delta}^{\Env^{+}(\calC)}
\stackrel{ (\phi_{\calC}^{+})_{!} }{\rightarrow} \Set_{\Delta}^{\calC},$$
where $(\phi_{\calC})_{!}$ and $(\phi_{\calC}^{+})_{!}$ indicate left Kan extension functors. 
There is a canonical isomorphism $\theta^{\ast} \circ E \simeq \id$, which induces a natural transformation $\alpha: \id \rightarrow L$. Let $\beta: L \rightarrow L^{+}$ indicate the natural transformation induced by the inclusion $\Env(\calC) \subseteq \Env^{+}(\calC)$.
\end{notation}

\begin{remark}\label{suppper}
Let $\calC$ be a small category. The functor $E^{+}: \Set_{\Delta}^{\calC} \rightarrow
\Set_{\Delta}^{\Env^{+}(\calC)}$ is fully faithful, and has a left adjoint given by $\theta^{\ast}$. 
\end{remark}

We begin by constructing the left adjoint which appears in the statement of Proposition \ref{trent}.

\begin{lemma}\label{toughluff}
Let $\calC$ be a simplicial category which admits finite products, and let $\calF \in \Set_{\Delta}^{\calC}$. Then:
\begin{itemize}
\item[$(1)$] The object $L^{+}(\calF) \in \Set_{\Delta}^{\calC}$ is product-preserving.
\item[$(2)$] If $\calF' \in \Set_{\Delta}^{\calC}$ is product-preserving, then composition with $\beta \circ \alpha$ induces an isomorphism of simplicial sets
$$ \bHom_{ \Set_{\Delta}^{\calC} }(L^{+}(\calF), \calF') \rightarrow \bHom_{\Set_{\Delta}^{\calC}}( \calF, \calF').$$
\end{itemize}
\end{lemma}

\begin{proof}
Suppose given a finite collection of objects $\{ C_1, \ldots, C_n \}$ in $\calC$, and let
$$u: L^{+}(\calF)( C_1 \times \ldots \times C_n) \rightarrow L^{+}(\calF)(C_1) \times \ldots \times L^{+}(\calF)(C_n)$$ be the product of the projection maps. We wish to show that $u$ is an isomorphism of simplicial sets. We will give an explicit construction of an inverse to $u$. For
$C \in \calC$, we let $\Env^{+}(\calC)_{/C}$ denote the fiber product
$\Env^{+}(\calC) \times_{\calD} \calC_{/C}$. For $1 \leq i \leq n$, let $\calG_i$ denote the restriction of $E^{+}(\calF)$ to $\Env^{+}(\calC)_{/C_i}$, and let 
$$ \calG: \prod \Env^{+}(\calD)_{/C_i} \rightarrow \sSet$$
be the product of the functors $\calG_i$. We observe that $L^{+}(\calF)(C_i) \simeq \varinjlim( \calG_i )$, so that the product $\prod L^{+}(\calF)(C_i) \simeq \varinjlim( \calG)$. We now observe that the formation of products in $\calE^{+}(\calC)$ gives an identification of $\calG$ with the composition
$$ \prod \Env^{+}(\calC)_{/C_i} \rightarrow \Env^{+}(\calD)_{/ C_1 \times \ldots \times C_n }
\stackrel{ E^{+}(\calF) }{\rightarrow} \sSet.$$
We thereby obtain a morphism $$v: \varinjlim(\calG) \rightarrow \varinjlim( E^{+}(\calF) | \Env^{+}(\calD)_{/  C_1 \times \ldots \times C_n} \simeq L^{+}(\calF)( C_1 \times \ldots \times C_n).$$
It is not difficult to check that $v$ is an inverse to $u$. 

We observe that $(2)$ is equivalence to the assertion that composition with $\theta^{\ast}$ induces an isomorphism
$$ \bHom_{ \Set_{\Delta}^{\Env^{+}(\calC)} }( E^{+}(\calF), (\phi_{\calC}^{+})^{\ast}(\calF'))
\rightarrow \bHom_{\Set_{\Delta}^{\calC}}( \calF, \calF').$$
Because $\calG$ is product-preserving, there is a natural isomorphism
$( \phi_{\calC}^{+})^{\ast}(\calF') \simeq E^{+}(\calF')$. The desired result now follows from Remark \ref{suppper}.
\end{proof}

It follows that the functor $L^{+}: \Set_{\Delta}^{\calC} \rightarrow \Set_{\Delta}^{\calC}$ factors through $\bfA$, and can be identified with a left adjoint to the inclusion
$\bfA \subseteq \Set_{\Delta}^{\calC}$. In order to prove Proposition \ref{trent}, we need to be able to compute the functor $L^{+}$. We will do this in two steps: first, we show that (under mild hypotheses), the natural transformation $L \rightarrow L^{+}$ is a weak equivalence. Second, we will see that the colimit defining $L$ is actually a homotopy colimit, and therefore has good properties. More precisely, we have the following pair of lemmas, whose proofs will be given at the end of this section.

\begin{lemma}\label{toughstuff}
Let $\calC$ be a small category which admits finite products, and let $\calF \in \Set_{\Delta}^{\calC}$ be a functor which carries the final object of $\calC$ to a contractible Kan complex $K$. Then the canonical map $\beta: L(\calF) \rightarrow L^{+}(\calF)$ is a weak equivalence in $\Set_{\Delta}^{\calC}$.
\end{lemma}

\begin{lemma}\label{trentlem}
Let $\calC$ be a small simplicial category. If $\calF$ is a projectively cofibrant object of $\Set_{\Delta}^{\calC}$, then $E(\calF)$ is a projectively cofibrant object of $\Set_{\Delta}^{\Env(\calC)}$. 
\end{lemma}

We are now almost ready to give the proof of Proposition \ref{trent}. The essential step is contained in the following result:

\begin{lemma}\label{presut}
Let $\calC$ be a simplicial category finite admits finite products, and let
$$ \Adjoint{ F }{ \Set_{\Delta}^{\calC} }{ \bfA }{ G}$$
be as in the statement of Proposition \ref{trent}. Then:
\begin{itemize}
\item[$(1)$] The functors $F$ and $G$ are Quillen adjoints.
\item[$(2)$] If $\calF \in \Set_{\Delta}^{\calC}$ is projectively cofibrant and weakly product preserving, then the unit map $\calF \rightarrow (G \circ F)(\calF)$ is a weak equivalence.
\end{itemize}
\end{lemma}

\begin{proof}
Assertion $(1)$ is obvious, since $G$ preserves fibrations and trivial cofibrations.
It follows that $F$ preserves weak equivalences between projectively cofibrant objects. 
Let $K \in \sSet$ denote the image under $\calF$ of the final object of $\calD$. 
In proving $(2)$, we are free to replace $\calF$ by any weakly equivalent diagram which is also projectively cofibrant. Choosing a fibrant replacement for $\calF$, we may suppose that $K$ is a Kan complex. Since $\calF$ is weakly product preserving, $K$ is contractible.

In view of Lemma \ref{toughluff}, we can identify the composition $G \circ F$ with $L^{+}$ and the unit map with the composition
$$ \calF \stackrel{\alpha}{\rightarrow} L(\calF) \stackrel{\beta}{\rightarrow} L^{+}(\calF).$$
Lemma \ref{toughstuff} implies that $\beta$ is a weak equivalence. Consequently, it will suffice to show that $\alpha$ is a weak equivalence.

We recall the construction of $\alpha$. Let $\theta: \calC \rightarrow \Env(\calC)$ be as in Notation \ref{bignote}, so that there is a canonical isomorphism $\calF \simeq \theta^{\ast} E(\calF)$. This
isomorphism induces a natural transformation $\overline{\alpha}: \theta_{!} \calF \rightarrow E(\calF)$. The functor $\alpha$ is obtained from $\overline{\alpha}$ by applying the functor
$(\phi_{\calC})_{!}$, and identifying $((\phi_{\calC})_{!} \circ \theta_{!} )(\calF)$ with $\calF$. 
We observe that $(\phi_{\calC})_{!}$ preserves weak equivalences between projectively cofibrant objects. Since $\theta_{!}$ preserves projective cofibrations, $\theta_{!} \calF$ is projectively cofibrant. Lemma \ref{trentlem} asserts that $E(\calF)$ is projectively cofibrant. Consequently, it will suffice to prove that $\overline{\alpha}$ is a weak equivalence in $\Set_{\Delta}^{\Env(\calC)}$. Unwinding the definitions, this reduces to the condition that $\calF$ be weakly compatible with (nonempty) products.
\end{proof}

\begin{proof}[Proof of Proposition \ref{trent}]
Lemma \ref{presut} shows that $(F,G)$ is a Quillen adjunction. To complete the proof, we must show:

\begin{itemize}
\item[$(i)$] The counit transformation $LF \circ RG \rightarrow \id$ is an isomorphism of functors
from the homotopy category $\h{\bfA}$ to itself.
\item[$(ii)$] The essential image of $RG: \h{\bfA} \rightarrow \h{ \Set_{\Delta}^{\calC} }$ consists
precisely of those functors which are weakly product-preserving.
\end{itemize}

We observe that $G$ preserves weak equivalences, so we can identify $RG$ with $G$. Since $G$ also detects weak equivalences, $(i)$ will follow if we can show that the induced transformation
$\theta: G \circ LF \circ G \rightarrow G$ is an isomorphism of functors from the homotopy
category $\h{\bfA}$ to itself. This transformation has a right inverse, given by composing the unit transformation $\id \rightarrow G \circ LF$ with $G$. Consequently, $(i)$ follows immediately from Lemma \ref{presut}. 

The image of $G$ consists precisely of the product-preserving diagrams $\calC \rightarrow \sSet$; it follows immediately that every diagram in the essential image of $G$ is weakly product preserving. Lemma \ref{presut} implies the converse: every weakly product preserving functor belongs to the essential image of $G$. This proves $(ii)$.
\end{proof}

\begin{remark}
Proposition \ref{trent} can be generalized to the situation where $\calC$ is a {\em simplicial} category which admits finite products. We leave the necessary modifications to the reader.
\end{remark}

It remains to prove Lemmas \ref{trentlem} and \ref{toughstuff}.

\begin{proof}[Proof of Lemma \ref{trentlem}]
For every object $C \in \calC$ and every simplicial set $K$, we let 
$\calF^{K}_{C} \in \Set_{\Delta}^{\calC}$ denote the functor given by the formula
$\calF^{K}_{C}(D) = \bHom_{\calC}( C, D) \times K$. A cofibration $K \rightarrow K'$ induces
a projective cofibration $\calF^{K}_{C} \rightarrow \calF^{K'}_{C}$. We will refer to a projective cofibration of this form as a {\em generating projective cofibration}.

The small object argument implies that if $\calF \in \Set_{\Delta}^{\calC}$, then there is a transfinite sequence
$$ \calF_0 \subseteq \calF_{1} \subseteq \ldots \subseteq \calF_{\alpha}$$
with the following properties:
\begin{itemize}
\item[$(a)$] The functor $\calF_0: \calD \rightarrow \sSet$ is constant, with value $\emptyset$.
\item[$(b)$] If $\lambda \leq \alpha$ is a limit ordinal, then $\calF_{\lambda} = \bigcup_{ \beta < \lambda} \calF_{\beta}$.
\item[$(c)$] For each $\beta < \alpha$, the inclusion $\calF_{\beta} \subseteq \calF_{\beta+1}$ is a pushout of a generating projective cofibration.
\item[$(d)$] The functor $\calF$ is a retract of $\calF_{\alpha}$.
\end{itemize}

The functor $\calG \mapsto E(\calG)$ preserves initial objects, filtered colimits, and retracts. Consequently, to show that $E(\calF)$ is projectively cofibrant, it will suffice to prove the following 
assertion:
\begin{itemize}
\item[$(\ast)$] Suppose given a cofibration $K \rightarrow K'$ of simplicial sets and a pushout diagram
$$ \xymatrix{ \calF^{K}_{C} \ar[d] \ar[r] & \calG \ar[d] \\
\calF^{K'}_{C} \ar[r] & \calG' }$$
in $\Set_{\Delta}^{\calC}$. If $E(\calG)$ is projectively cofibrant, then $E(\calG')$ is projectively cofibrant.
\end{itemize}

To prove this, we will need to analyze the structure of $E(\calG')$. Given an object
$C' = ( J, \{ C'_{j} \}_{j \in J} )$ of $\Env(\calC)$, we have
$$E(\calG')( C' ) = \prod_{j \in J} ( \calG( C'_{j} ) \coprod_{ K \times \bHom_{\calC}(C,C'_j) }
( K' \times \bHom_{\calC}(C, C'_j) ). $$
Let $\sigma: \Delta^n \rightarrow E(\calG')(C')$ be a simplex, and let $J_{\sigma} \subseteq J$
be the collection of all indices $j$ for which the corresponding simplex 
$\sigma(j): \Delta^n \rightarrow \calG'(C'_j)$ does not factor through $\calG(C'_j)$. 
In this case, we can identify $\sigma(j)$ with an $n$-simplex of $K'$ which does not belong to
$K$. We will say that $\sigma$ is of {\it index $\leq k$} if the set $\{ \sigma(j) : j \in J_{\sigma} \}$ has cardinality $\leq k$. Note that $\sigma$ can be of index smaller than the cardinality of
$J_{\sigma}$, since it is possible for $\sigma(j) = \sigma(j') \in \bHom_{\sSet}(\Delta^n, K')$ even if $j \neq j'$.

Let $E(\calG')^{(k)}(C')$ be the full simplicial subset of $E(\calG')(C')$ spanned by those simplices which are of index $\leq k$. It is easy to see that that $E(\calG')^{(k)}(C')$ depends functorially in $C'$, so we can view $E(\calG')^{(k)}$ as an object of $\Set_{\Delta}^{\Env(\calC)}$. We observe that
$$ E(\calG) \simeq E(\calG')^{(0)} \subseteq E(\calG')^{(1)} \subseteq \ldots $$
and that the union of this sequence is $E(\calG')$. Consequently, it will suffice to prove that
each of the inclusions $E(\calG')^{(k-1)} \subseteq E(\calG')^{(k)}$ is a projective cofibration.

First, we need a bit of notation. Let us say that a simplex of ${K'}^{k}$ is {\it new} if it consists of $k$ distinct simplices of $K'$, none of which belong to $K$. We will say that a simplex of ${K'}^{k}$ is {\it old} if it is not new. The collection of old simplices of ${K'}^{k}$ determine a simplicial subset which we will denote by ${K'}^{(k)}$. We define a functor $\psi: \Env(\calC) \rightarrow \Env(\calC)$ by the formula
$$\psi( J, \{ C'_{j} \}_{j \in J} ) = ( J \cup \{ 1, \ldots, k \}, \{ C'_{j} \}_{j \in J} \cup \{ C \}_{\{1 \ldots k \}} ).$$
Let $\psi^{\ast}: \Set_{\Delta}^{\Env(\calC)} \rightarrow \Set_{\Delta}^{\Env(\calC)}$ be given by composition with $\psi$, and let $\psi_{!}: \Set_{\Delta}^{\Env(\calC)} \rightarrow \Set_{\Delta}^{\Env(\calC)}$ be a left adjoint to $\psi^{\ast}$ (a functor of left Kan extension).
Since $\psi^{\ast}$ preserves projective fibrations and weak equivalences, $\psi_{!}$ preserves projective cofibrations.

Recall that $\Set_{\Delta}^{\Env(\calC)}$ is tensored over the category of simplicial sets:
Given an object $\calM \in \Set_{\Delta}^{\Env(\calC)}$ and a simplicial set $A$, we let $\calM \otimes A \in \Set_{\Delta}^{\Env(\calC)}$ be defined by the formula
$( \calM \otimes A)(D') = \calM(D') \times A$. If $\calM$ is projectively cofibrant, then the operation
$\calM \mapsto \calM \otimes A$ preserves cofibrations in $A$.

There is an obvious map $E(\calG) \otimes {K'}^k \rightarrow \psi^{\ast} E(\calG')^{(k)}$, which restricts to a map $E(\calG) \otimes {K'}^{(k)} \rightarrow \psi^{\ast} E(\calG')^{(k-1)}$. Passing to adjoints, we obtain a commutative diagram
$$ \xymatrix{ \psi_{!}( E(\calG) \otimes {K'}^{(k)} ) \ar[d] \ar[r] & E(\calG')^{(k-1)} \ar[d] \\
\psi_{!}( E(\calG) \otimes {K'}^k ) \ar[r] & E(\calG')^{(k)}. }$$
An easy computation shows that this diagram is coCartesian.
Since $E(\calG)$ is projectively cofibrant, the above remarks imply that the left vertical map is a projective cofibration. It follows that the right vertical map is a projective cofibration as well, which completes the proof.
\end{proof}

The proof of Lemma \ref{toughstuff} is somewhat more difficult, and will require some preliminaries.

\begin{notation}
Let $\Mult: \Set \rightarrow \Set$ be the covariant functor which associates to each set
$S$ the collection of $\Mult(S)$ of nonempty finite subsets of $S$. If $K$ is a simplicial set, we let
$\Multi(K)$ denote the composition of $K$ with $\Mult$, so that an $m$-simplex of $\Multi(K)$ is
a finite nonempty collection of $m$-simplices of $K$.\index{not}{MultiK@$\Multi(K)$}
\end{notation}

\begin{lemma}\label{stuffem}
Let $K$ be a finite simplicial set, let $X \subseteq \Multi(K^{\triangleright}) \times \Delta^n$ be a simplicial subset with the following properties:
\begin{itemize}
\item[$(i)$] The projection $X \rightarrow \Delta^n$ is surjective. 
\item[$(ii)$] If $\tau = (\tau', \tau''): \Delta^m \rightarrow \Multi(K^{\triangleright}) \times \Delta^n$ belongs to $X$, and $\tau' \subseteq \overline{\tau}'$ as subsets of $\Hom_{\sSet}( \Delta^m, K^{\triangleright})$, then $(\overline{\tau}',\tau''): \Delta^m \rightarrow \Multi(K^{\triangleright}) \times \Delta^n$ belongs to $X$.
\end{itemize}
Then $X$ is weakly contractible.
\end{lemma}

\begin{proof}
Let $X' \subseteq X$ be the simplicial subset spanned by those simplices
$\tau = (\tau', \tau''): \Delta^m \rightarrow \Multi(K^{\triangleright}) \times \Delta^n$ which factor through $X$, and for which $\tau' \subseteq \Hom_{\sSet}( \Delta^m, K^{\triangleright})$ includes the constant simplex at the cone point of $K^{\triangleright}$. Our first step is to show that $X'$ is a deformation retract of $X$. More precisely, we will construct a map
$$ h: \Multi(K^{\triangleright}) \times \Delta^n \times \Delta^1 \rightarrow 
\Multi(K^{\triangleright}) \times \Delta^n$$ with the following properties:
\begin{itemize}
\item[$(a)$] The map $h$ carries $X \times \Delta^1$ into $X$ and $X' \times \Delta^1$ into $X'$.
\item[$(b)$] The restriction $h | \Multi(K^{\triangleright}) \times \Delta^n \times \{0\}$ is the identity map.
\item[$(c)$] The restriction $h| X \times \{1\}$ factors through $X'$.
\end{itemize}
The map $h$ will be the product of a map $h': \Multi(K^{\triangleright}) \times \Delta^1 \rightarrow \Multi(K^{\triangleright} )$ with the identity map on $\Delta^n$. To define $h'$, we consider an arbitrary simplex $\tau: \Delta^m \rightarrow \Multi(K^{\triangleright}) \times \Delta^1$, corresponding to a subset $S \subseteq \Hom_{\sSet}( \Delta^m, K^{\triangleright} )$ and a decomposition
$[m] = \{ 0, \ldots, i\} \cup \{i+1, \ldots, m\}$. The subset $h'(\tau) \subseteq \Hom_{\sSet}( \Delta^m, K^{\triangleright})$ is defined as follows: an arbitrary simplex $\sigma: \Delta^m \rightarrow K^{\triangleright}$ belongs to $h'(\tau)$ if there exists $\sigma' \in S$, $i < j \leq n$ such that
$\sigma' | \Delta^{ \{0, \ldots, j-1 \} } = \sigma| \Delta^{ \{0, \ldots, j-1 \} }$, and $\sigma | \Delta^{ \{j, \ldots, m \} }$ is constant at the cone point of $K^{\triangleright}$. It is easy to check that $h'$ has the desired properties.

It remains to prove that $X'$ is weakly contractible. At this point, it is convenient to work in the setting of {\em semisimplicial sets}: that is, we will ignore the degeneracy operations. Let
$X''$ be the semisimplicial subset of $\Multi(K^{\triangleright}) \times \Delta^n$ spanned by those maps $\tau = (\tau', \tau''): \Delta^m \rightarrow \Multi(K^{\triangleright}) \times \Delta^n$ for which
$\tau' = \Hom_{\sSet}( \Delta^m, K^{\triangleright} )$ (we observe that $X''$ is not stable under the degeneracy operators on $\Multi(K^{\triangleright}) \times \Delta^n$). Assumptions $(i)$ and $(ii)$ guarantee that $X'' \subseteq X'$. Moreover, the projection $X \rightarrow \Delta^n$ induces an isomorphism of semisimplicial sets $X'' \rightarrow \Delta^n$. Consequently, it will suffice to prove that $X''$ is a deformation retract of $X'$. 

The proof now proceeds by a variation on our earlier construction. Namely, we will define a map 
of semisimplicial sets
$$ g: \Multi(K^{\triangleright}) \times \Delta^n \times \Delta^1 \rightarrow 
\Multi(K^{\triangleright}) \times \Delta^n$$ with the following properties:
\begin{itemize}
\item[$(a)$] The map $g$ carries $X' \times \Delta^1$ into $X'$ and $X'' \times \Delta^1$ into $X''$.
\item[$(b)$] The restriction $g |X' \times \{1\}$ is the identity map.
\item[$(c)$] The restriction $g| \Multi(K^{\triangleright}) \times \Delta^n \times \{0\}$ factors through $X'$.
\end{itemize}
As before, $g$ is the product of a map $g': \Multi(K^{\triangleright}) \times \Delta^1 \rightarrow \Multi(K^{\triangleright})$ with the identity map on $\Delta^n$. To define $g'$, we consider an arbitrary simplex $\tau: \Delta^m \rightarrow \Multi(K^{\triangleright}) \times \Delta^1$, corresponding to a subset $S \subseteq \Hom_{\sSet}( \Delta^m, K^{\triangleright} )$ and a decomposition
$[m] = \{ 0, \ldots, i\} \cup \{i+1, \ldots, m\}$. We let $g'(\tau) \subseteq \Hom_{\sSet}(\Delta^m, K^{\triangleright} ) = S \cup S'$, where $S'$ is the collection of all simplices $\sigma: \Delta^m \rightarrow K^{\triangleright}$ such that $\sigma | \Delta^{ \{i+1, \ldots, m\} }$ is the constant map at the cone poine of $K^{\triangleright}$. It is readily checked that $g'$ has the desired properties.
\end{proof}

\begin{lemma}\label{toughfluff}
Let $K$ be a contractible Kan complex, and let $X \subseteq \Multi(K) \times \Delta^n$ be a simplicial subset with the following properties:
\begin{itemize}
\item[$(i)$] The projection $X \rightarrow \Delta^n$ is surjective.
\item[$(ii)$] If $\tau = (\tau', \tau''): \Delta^m \rightarrow \Multi(K) \times \Delta^n$ belongs to
$X$, and $\tau' \subseteq \overline{\tau}'$ as subsets of $\Hom_{\sSet}( \Delta^m, K)$, then
$(\overline{\tau}',\tau''): \Delta^m \rightarrow \Multi(K) \times \Delta^n$ belongs to $X$.
\end{itemize}
Then $X$ is weakly contractible.
\end{lemma}

\begin{proof}
It will suffice to show that for every finite simplicial subset $X' \subseteq X$, the inclusion of $X'$ into $X$ is weakly nullhomotopic. Enlarging $X'$ if necessary, we may assume that
$X' = (\Multi(K') \times \Delta^n) \cap X$, where $K'$ is a finite simplicial subset of $K$. By further enlargement, we may suppose that the map $X' \rightarrow \Delta^n$ is surjective. Since $K$ is a contractible Kan complex, the inclusion $K' \subseteq K$ extends to a map $i: {K'}^{\triangleright} \rightarrow K$. Let $\overline{X} \subseteq \Multi( {K'}^{\triangleright} ) \times \Delta^n$ denote
the inverse image of $X$. Then the inclusion $X' \subseteq X$ factors through $\overline{X}$, and Lemma \ref{stuffem} implies that $\overline{X}$ is weakly contractible.
\end{proof}

\begin{proof}[Proof of Lemma \ref{toughstuff}]
Fix an object $C \in \calC$. The simplicial set $L(\calF)(C)$ can be described as follows:
\begin{itemize}
\item[$(\ast)$] For every $n \geq 1$, every map $f: C_1 \times \ldots \times C_n \rightarrow C$ in $\calC$, and every collection of simplices $\{ \sigma_i: \Delta^k \rightarrow \calF(C_i) \}$, there is a simplex $f( \{ \sigma_i \}): \Delta^k \rightarrow L(\calF)(C)$.
\end{itemize}
The simplices $f( \{ \sigma_i \} )$ satisfy relations which are determined by morphisms in the simplicial category $\Env(\calC)$.

To every $k$-simplex $\tau: \Delta^k \rightarrow L(\calF)(D)$, we can associate a nonempty finite subset $S_{\tau} \subseteq \Hom_{\sSet}( \Delta^k, K)$. If $\tau = f( \{ \sigma_i \})$, we assign the set of images of the simplices $\sigma_i$ under the canonical maps $\calF(C_i) \rightarrow \calF(1) = K$. It is easy to see that $S_{\tau}$ is independent of the representation $f( \{ \sigma_i \})$ chosen for $\tau$, and depends functorially on $\tau$. Consequently, we obtain a map of simplicial sets
$L(\calF)(C) \rightarrow \Multi(K)$. Moreover, this map has the following properties:
\begin{itemize}
\item[$(i)$] The product map $\beta': L(\calF)(C) \rightarrow \Multi(K) \times L^{+}(\calF)(C)$ is a monomorphism of simplicial sets. 
\item[$(ii)$] If a $k$-simplex $\tau = (\tau', \tau'') : \Delta^k \rightarrow \Multi(K) \rightarrow L^{+}(\calF)(C)$ belongs to the image of $\beta$, and 
$\tau' \subseteq \overline{\tau}'$ as finite subsets of $\Hom_{\sSet}(\Delta^k, K)$, then
$(\overline{\tau}', \tau''): \Delta^k \rightarrow \Multi(K) \times L^{+}(\calF)(C)$ belongs to the image
of $\beta'$. 
\end{itemize}

We wish to show that $\beta: L(\calF)(C) \rightarrow L^{+}(\calF)(C)$ is a weak homotopy equivalence. It will suffice to show that for every simplex $\Delta^k \rightarrow L^{+}(\calF)(C)$, the fiber product $L(\calF)(C) \times_{ L^{+}(\calF)(C) } \Delta^k$ is weakly contractible. In view of $(i)$, we can identify this fiber product with a simplicial subset $X \subseteq \Delta^k \times \Multi(K)$.
The surjectivity of $\beta$ and condition $(ii)$ imply that $X$ satisfies the hypotheses of Lemma \ref{toughfluff}, so that $X$ is weakly contractible as desired.
\end{proof}

The model category $\bfA$ appearing in Proposition \ref{sutcoat} is very well suited to certain calculations, such as the formation of homotopy colimits of simplicial objects. The following result provides a precise formulation of this idea:

\begin{proposition}\label{eggers}
Let $\calC$ be a category which admits finite products, and $\bfA \subseteq \Set_{\Delta}^{\calC}$,
$\calA \subseteq \Set^{\calC}$ the full subcategories spanned by the product-preserving functors. 
Let $\calF: \cDelta^{op} \rightarrow \bfA$ be a simplicial object of $\bfA$, which we can identify with a {\em bisimplicial} object $F: \cDelta^{op} \times \cDelta^{op} \rightarrow \calA$. Composition with the diagonal $$ \cDelta^{op} \rightarrow \cDelta^{op} \times \cDelta^{op} \stackrel{F}{\rightarrow} \calA$$ gives a simplicial object of $\calA$, which we can identify with an object $| \calF| \in \bfA$. Then the homotopy colimit of $\calF$ is canonically isomorphic to $|\calF|$ in the homotopy category $\h{\bfA}$.
\end{proposition}

The proof requires the following lemma:

\begin{lemma}\label{urplus}
Let $\calC$ be a category which admits finite products, and $\bfA \subseteq \Set_{\Delta}^{\calC}$
be the full subcategory spanned by the product preserving functors. For every object $C \in \calC$, the evaluation map $\bfA \rightarrow \sSet$ preserves homotopy colimits of simplicial objects.
\end{lemma}

\begin{proof}
In view of Corollary \ref{smokerr} and Theorem \ref{colimcomparee}, it will suffice to show that the evaluation functor $\calP_{\Sigma}( \Nerve(\calC)^{op} ) \rightarrow \sSet$ preserves
$\Nerve(\cDelta)^{op}$-indexed colimits. This follows from Proposition \ref{utut}, since
$\Nerve(\cDelta)^{op}$ is sifted (Lemma \ref{bball3}).
\end{proof}

\begin{proof}[Proof of Proposition \ref{eggers}]
Since $\bfA$ is a combinatorial simplicial model category, Corollary \ref{twinner} implies the existence of a canonical map $\gamma: \hocolim \calF \rightarrow | \calF |$ in the homotopy category $\h{\bfA}$; we wish to prove that $\gamma$ is an isomorphism. To prove this, it will suffice to show that
the induced map $\gamma_{C}: ( \hocolim \calF)(C) \rightarrow | \calF|(C)$ is an isomorphism in the homotopy category of simplicial sets, for each object $C \in \calC$. This map fits into a commutative diagram
$$ \xymatrix{ \hocolim( \calF(C) ) \ar[r]^{\gamma'_{C}} \ar[d] & | \calF(C) | \ar[d] \\
\hocolim(\calF)(C) \ar[r] & | \calF|(C).}$$
The left vertical map is an isomorphism in the homotopy category of simplicial sets by Lemma \ref{urplus}, the right vertical map is evidently an isomorphism, and the map $\gamma'_{C}$ is an isomorphism in the homotopy category by Example \ref{swupt}; it follows that $\gamma_{C}$ is also an isomorphism, as desired.
\end{proof}

\chapter{$\infty$-Topoi}\label{chap6}

\setcounter{theorem}{0}
\setcounter{subsection}{0}

In this chapter, we come to the main subject of this book: the theory of $\infty$-topoi.
Roughly speaking, an $\infty$-topos is an $\infty$-category which ``looks like'' the $\infty$-category of spaces, just as an ordinary topos is a category which ``looks like'' the category of sets. As in classical topos theory, there are various ways of making this precise. We will begin in \S \ref{chap6sec1} by reviewing several possible definitions, and proving that they are equivalent to one another.

The main result of \S \ref{chap6sec1} is Theorem \ref{mainchar}, which asserts that an $\infty$-category $\calX$ is an $\infty$-topos if and only if $\calX$ arises as an (accessible) left exact localization of an $\infty$-category of presheaves. In \S \ref{topcomp}, we consider the problem of {\em constructing} left exact localizations. In classical topos theory, there is a bijective correspondence between left exact localizations of a presheaf category $\calP(\calC)$ and Grothendieck topologies on $\calC$. In the $\infty$-category categorical context, one can again use Grothendieck topologies to construct examples of left exact localizations. 
Unfortunately, not every $\infty$-topos arises in this way. Nevertheless, the construction of an $\infty$-category of sheaves $\Shv(\calC)$ from a Grothendieck topology on $\calC$ is an extremely useful construction, which will play an important role throughout \S \ref{chap7}.

In order to understand higher topos theory, we will need to consider $\infty$-topoi not only individually, but in relation to one another. In \S \ref{chap6sec4} we will introduce the notion of a {\em geometric morphism} of $\infty$-topoi. The collection of all $\infty$-topoi and geometric morphisms between them can be organized into an $\infty$-category $\RGeom$. We will study the problem of constructing colimits and (certain) limits in $\RGeom$. In the course of doing so, we will show that the class of $\infty$-topoi is stable under various categorical constructions.

One of our goals in this book is to apply ideas from higher category theory to study more classical mathematical objects, such as topological spaces or ordinary topoi. In order to do so, it is convenient to work in a setting where all of these objects can be considered on the same footing. In \S \ref{chap6sec3}, we will introduce the definition of an {\it $n$-topos} for all $0 \leq n \leq \infty$.
When $n = \infty$, this will reduce to the theory introduced in \S \ref{chap6sec1}. The case $n=1$ will recover classical topos theory, and the case $n=0$ is {\em almost} equivalent to the theory of topological spaces. We will study the theory of $n$-topoi, and introduce constructions which allow us to pass between $n$-topoi and $\infty$-topoi. In particular, we associate an $\infty$-topos $\Shv(X)$ to every topological space $X$, which will be the primary object of study in \S \ref{chap7}.

There are several different ways of thinking about what an $\infty$-topos $\calX$ is.
On the one hand, we can view $\calX$ as a generalized topological space; on the other,
we can think of $\calX$ as an alternative universe in which we can do homotopy theory. 
In \S \ref{chap6sec5}, we will reinforce the second point of view by studying the internal homotopy theory of an $\infty$-topos $\calX$. Just as in classical homotopy theory, one can define homotopy groups, Postnikov towers, Eilenberg MacLane spaces, and so forth. In \S \ref{chap7}, we will bring together these two points of view, by showing that classical geometric properties of a topological space $X$ are reflected in internal homotopy of the $\infty$-topos $\Shv(X)$ of sheaves on $X$.

There are several papers on higher topos theory in the literature. The
papers \cite{street} and \cite{ditopoi} both discuss a notion of
$2$-topos (the second from an elementary point of view). However,
the basic model for these $2$-topoi is the $2$-category of
(small) categories, rather than the $2$-category of (small)
groupoids. Jardine (\cite{jardine}) has exhibited a model
structure on the category of simplicial presheaves on a Grothendieck site, and
the $\infty$-category associated to this model category is an
$\infty$-topos in our sense. This construction is generalized
from ordinary categories with a Grothendieck topology to
simplicial categories with a Grothendieck topology in \cite{toen} (and again produces $\infty$-topoi).
However, not every $\infty$-topos arises in this way: one can
construct only $\infty$-topoi which are {\it hypercomplete} (called {\it $t$-complete} in
\cite{toen}); we will summarize the situation in Section
\ref{hyperstacks}. Our notion of an $\infty$-topos is {\em essentially} equivalent to the notion of
a {\it Segal topos} introduced in \cite{toen}, and to Charles Rezk's notion of a {\it model topos}. 
We note also that the paper \cite{toen} has considerable overlap with the ideas discussed here.

\section{$\infty$-Topoi: Definitions and Characterizations}\label{chap6sec1}
 
\setcounter{theorem}{0}

Before we study the $\infty$-categorical version of topos theory, it seems appropriate to briefly review the classical theory. Recall that a {\it topos} is a category $\calC$ which behaves like the category of sets, or (more generally) the category of sheaves of sets on a topological space. There are several (equivalent) ways of making this idea precise. The following result is proved (in a slightly different form) in \cite{SGA}:

\begin{proposition}\label{toposdefined}\index{gen}{topos}
Let $\calC$ be a category. The following conditions are equivalent:
\begin{itemize}
\item[$(A)$] The category $\calC$ is $($equivalent to$)$ the category of
sheaves of sets on some Grothendieck site. \item[$(B)$] The category $\calC$ is
$($equivalent to$)$ a left-exact localization of the category of presheaves of sets
on some small category $\calC_0$. \item[$(C)$]\index{gen}{Giraud's axioms!for ordinary topoi} Giraud's axioms are satisfied: \begin{itemize}
\item[$(i)$] The category $\calC$ is presentable $($that is, $\calC$ has small colimits and a set of small generators$)$.
\item[$(ii)$] Colimits in $\calC$ are universal. 
\item[$(iii)$] Coproducts in $\calC$ are disjoint. 
\item[$(iv)$] Equivalence relations in $\calC$ are effective.
\end{itemize}
\end{itemize}
\end{proposition}

\begin{definition}\label{def1topos}
A category $\calC$ is called a {\it topos} if it satisfies the equivalent conditions of
Proposition \ref{toposdefined}.
\end{definition}

\begin{remark}
A reader who is unfamiliar with some of the terminology used in the statement of Proposition \ref{toposdefined} should not worry: we will review the meaning of each condition in \S \ref{axgir} as we search for $\infty$-categorical generalizations of axioms $(i)$ through $(iv)$.
\end{remark}

Our goal in this section is to introduce the $\infty$-categorical analogue of Definition \ref{def1topos}.
Proposition \ref{toposdefined} suggests several possible approaches. We begin with the simplest of these:

\begin{definition}\label{itoposdef}\index{gen}{$\infty$-topos}
Let $\calX$ be an $\infty$-category. We will say that $\calX$ is an {\it $\infty$-topos} if there
exists a small $\infty$-category $\calC$ and an accessible left exact localization functor
$\calP(\calC) \rightarrow \calX$. 
\end{definition}

\begin{remark}
Definition \ref{itoposdef} involves an accessibility condition which was not
mentioned in Proposition \ref{toposdefined}. This is because every left exact localization of a
category of {\em set}-valued presheaves is automatically accessible (see Proposition \ref{alltoploc}). We do not know if the corresponding result holds for $\SSet$-valued presheaves. However,
it is true under a suitable hypercompleteness assumption:
see \cite{toenvezz}.
\end{remark}

Adopting Definition \ref{itoposdef} amounts to selecting an {\em extrinsic} approach to higher topos theory: the class of $\infty$-topoi is defined to be the smallest collection of $\infty$-categories which contains $\SSet$ and is stable under certain constructions (left exact localizations and the formation of functor categories). The main objective of this section is to give several reformulations of Definition \ref{def1topos} which have a more intrinsic flavor. Our results may be summarized in the following statement (all our our terminology will be explained later in this section):

\begin{theorem}\label{mainchar}\index{gen}{Giraud's theorem!for $\infty$-topoi}
Let $\calX$ be an $\infty$-category. The following conditions are equivalent:

\begin{itemize}
\item[$(1)$] The $\infty$-category $\calX$ is an $\infty$-topos.

\item[$(2)$] The $\infty$-category $\calX$ is presentable, and for every small simplicial
set $K$ and every natural transformation
$\overline{\alpha}: \overline{p} \rightarrow \overline{q}$ of diagrams
$\overline{p}, \overline{q}: K^{\triangleright} \rightarrow \calX$, the following condition is satisfied:
\begin{itemize} 
\item
If $\overline{q}$ is a colimit diagram and $\alpha = \overline{\alpha} | K$ is a Cartesian transformation, then $\overline{p}$ is a colimit diagram if and only if $\overline{\alpha}$ is a Cartesian transformation.
\end{itemize}

\item[$(3)$] The $\infty$-category $\calX$ satisfies the following $\infty$-categorical analogues of Giraud's axioms:
\begin{itemize}
\item[$(i)$] The $\infty$-category $\calX$ is presentable.
\item[$(ii)$] Colimits in $\calX$ are universal.
\item[$(iii)$] Coproducts in $\calX$ are disjoint.
\item[$(iv)$] Every groupoid object of $\calX$ is effective.
\end{itemize}

\end{itemize}
\end{theorem}

We will review the meanings of conditions $(i)$ through $(iv)$ in \S \ref{axgir} and \S \ref{gengroup}.  In \S \ref{magnet} we will give several equivalent formulations of $(2)$, and prove the implications
$(1) \Rightarrow (2) \Rightarrow (3)$. The implication $(3) \Rightarrow (1)$ is the most difficult; we will give the proof in \S \ref{proofgiraud} after establishing a crucial technical lemma in \S \ref{freegroup}. Finally, in \S \ref{rezk2} we will establish yet another characterization of $\infty$-topoi, based on the theory of classifying objects.

\begin{remark}
The characterization of the class of $\infty$-topoi given by part $(2)$ of Theorem \ref{mainchar} is due to Rezk, as are many of the ideas presented in \S \ref{magnet}.
\end{remark}

The equivalence $(1) \Leftrightarrow (3)$ of Theorem \ref{mainchar} can be viewed as an $\infty$-categorical analogue of the equivalence $(B) \Leftrightarrow (C)$ in Proposition \ref{toposdefined}. 
It is natural to ask if there is also some equivalent of the characterization $(A)$. To put the question another way: given a small $\infty$-category $\calC$, does there exist some natural description of the class of all left-exact localizations of $\calC$? Experience with classical topos theory suggests that we might try to characterize such localizations in terms of {\em Grothendieck topologies} on $\calC$. We will introduce a theory of Grothendieck topologies on $\infty$-categories in \S \ref{cough}, and show that every Grothendieck topology on $\calC$ determines a left-exact localization of $\calP(\calC)$. However, it turns out that not every $\infty$-topos arises via this construction. This raises a natural question: is it possible to give an explicit description of {\em all} left-exact localizations of $\calP(\calC)$, perhaps in terms of some more refined theory of Grothendieck topologies? We will give a partial answer to this question in \S \ref{chap6sec5}.

\subsection{Giraud's Axioms in the $\infty$-Categorical Setting}\label{axgir}\index{gen}{Giraud's axioms!for $\infty$-topoi}

Our goal in this section is to formulate higher-categorical analogues of the conditions $(i)$ through $(iv)$ which appear in Proposition \ref{toposdefined}. 
We consider each axiom in turn. In each case, our objective is to find an analogous axiom which makes sense in the setting of $\infty$-categories, and is satisfied by the $\infty$-category $\SSet$
of spaces. 

\begin{itemize}
\item[$(i)$] The category $\calC$ is presentable.
\end{itemize}

The generalization to the case where $\calC$ is a $\infty$-category is obvious: we should merely require $\calC$ to be a presentable $\infty$-category in the sense of Definition \ref{presdef}.
According to Example \ref{spacesarepresentable}, this condition is satisfied when $\calC$ is the $\infty$-category of spaces.

\begin{itemize}\index{gen}{colimit!universal}
\item[$(ii)$] Colimits in $\calC$ are universal.
\end{itemize}

Let us first recall the meaning of this condition in classical category theory. If the axiom $(i)$ is satisfied, then $\calC$ is presentable and therefore admits all (small) limits and colimits. In particular, every diagram
$$ X \rightarrow S \stackrel{f}{\leftarrow} T$$
has a limit $X_{T} = X \times_{S} T$.
This construction determines
a functor
$$ f^{\ast}: \calC_{/S} \rightarrow \calC_{/T}$$
$$X \mapsto X_{T},$$
which is a right adjoint to the functor given by composition with $f$.

We say that colimits in $\calC$ are {\it universal} if the functor $f^{\ast}$ is colimit-preserving, for
every map $f: T \rightarrow S$ in $\calC$. (In other words, colimits are universal in $\calC$ if any colimit in $\calC$ {\em remains} a colimit in $\calC$ after pulling back along a morphism $T \rightarrow S$.)

Let us now attempt to make this notion precise in the setting of an arbitrary $\infty$-category
$\calC$. Let $\calO_{\calC} = \Fun(\Delta^1, \calC)$\index{not}{OC@$\calO_{\calC}$}, and let $p: \calO_{\calC} \rightarrow \calC$ be given by evaluation at
$\{1\} \subseteq \Delta^1$. Corollary \ref{tweezegork} implies that $p$ is a coCartesian fibration.

\begin{lemma}\label{charpull}
Let $\calX$ be an $\infty$-category and let $p: \calO_{\calX} \rightarrow \calX$ be defined as above. Let
$F$ be a morphism in $\calO_{\calX}$, corresponding to a diagram $\sigma: \Delta^1 \times \Delta^1 \simeq (\Lambda^2_2)^{\triangleleft} \rightarrow \calX$, which we will denote by
$$ \xymatrix{ X' \ar[r]^{f'} \ar[d] & Y' \ar[d]^{g} \\
X \ar[r]^{f} & Y}$$
Then $F$ is $p$-Cartesian if and only if the above diagram is a pullback in $\calX$.
In particular, $p$ is a Cartesian fibration if and only if the $\infty$-category $\calX$ admits pullbacks.
\end{lemma}

\begin{proof}
For every simplicial set $K$, let $K^{+}$ denote the full simplicial subset
of $(K \star \{x\} \star \{y\}) \times \Delta^1$ spanned by all of the vertices except
$(x,0)$, and define a simplicial set $\calC$ by setting 
$$\Fun(K, \calC) = \{ m: K^{+} \rightarrow \calX : m| (\{x\} \star \{y\}) \times \{1\} = f,
m| \{y\} \times \Delta^1 = g \}.$$
We observe that we have a commutative diagram
$$ \xymatrix{ \calC \ar[r] \ar[d] & (\calO_{\calX})_{/g} \ar[d] \\
\calX_{/f} \ar[r] & \calX_{/Y'} }$$
which induces a map $q: \calC \rightarrow (\calO_{\calX})_{/g} \times_{ \calX_{/Y'} } \calX_{/f}$.
We first claim that $q$ is a trivial fibration. Unwinding the definitions, we observe that the right lifting property of $q$ with respect to an inclusion $\bd \Delta^n \subseteq \Delta^n$ follows from the extension property of $\calX$ with respect to $\Lambda^{n+2}_{n+1}$, which follows in turn from our assumption that $\calX$ is an $\infty$-category.

The inclusion $K^{+} \subseteq K \times \Delta^1$ induces a projection $q': ( \calO_{\calX})_{/F} \rightarrow \calC$ which fits into a pullback diagram
$$ \xymatrix{ ( \calO_{\calX})_{/F} \ar[r] \ar[d] & \calC \ar[d]^{g} \\
\calX_{/ \sigma } \ar[r]^{q''} & \calX_{/ \sigma | \Lambda^2_2}. }$$
It follows that $q'$ is a right fibration, and that $q'$ is trivial if $\sigma$ is a pullback
diagram. Conversely, we observe that $( \Lambda^2_2)^{\triangleleft}$ is a retract of
$(\Delta^{0})^{+}$, so that the map $g$ is surjective on vertices. Consequently, if
$q'$ is a trivial fibration, then the fibers of $q''$ are contractible, so that $q''$ is a trivial fibration
(Lemma \ref{toothie}) and $\sigma$ is a pullback diagram.

By definition, $F$ is $p$-Cartesian if and only if the composition
$$q \circ q': (\calO_{\calX})_{/F} \rightarrow (\calO_{\calX})_{/g} \times_{ \calX_{/Y'} } \calX_{/f}$$
is a trivial fibration. Since $q$ is a trivial fibration and $q'$ is a right fibration, this is also equivalent to the assertion that $q'$ is a trivial fibration (Lemma \ref{toothie}).
\end{proof}

Now suppose that $\calX$ is an $\infty$-category which admits pullbacks, so that the projection
$p: \calO_{\calX} \rightarrow \calX$ is both a Cartesian fibration and a coCartesian fibration.
Let $f: S \rightarrow T$ be a morphism in $\calX$. Taking the pullback of $p$ along
the corresponding map $\Delta^1 \rightarrow \calX$, we obtain a correspondence from
$p^{-1}(S) = \calX^{/S}$ to $p^{-1}(T) = \calX^{/T}$, associated to a pair of adjoint functors
$$ f_{!}: \calX^{/S} \rightarrow \calX^{/T}$$
$$ f^{\ast}: \calX^{/T} \rightarrow \calX^{/S}.$$
The functors $f_{!}$ and $f^{\ast}$ are well-defined up to homotopy (in fact, up to a contractible space of choices. We may think of $f_{!}$ as the functor given by composition with $f$, and $f^{\ast}$ as the functor given by pullback along $f$ (in view of Lemma \ref{charpull}).\index{gen}{pullback functor} 

We can now formulate the $\infty$-categorical analogue of $(ii)$:

\begin{definition}\label{colu}\index{gen}{colimit!universal}
Let $\calC$ be a presentable $\infty$-category. We will say that {\it colimits in $\calC$
are universal} if, for any morphism $f: T \rightarrow S$ in $\calC$, the associated pullback functor
$$ f^{\ast}: \calC^{/S} \rightarrow \calC^{/T}$$
preserves (small) colimits.
\end{definition}

Assume that $\calC$ is a presentable $\infty$-category, and let $f: T \rightarrow S$ be a morphism in $\calC$. By the adjoint functor theorem, $f^{\ast}: \calC^{/S} \rightarrow \calC^{/T}$ preserves all colimits if and only if it has a right adjoint $f_{\ast}$. Since the existence of adjoint functors can be tested inside the enriched homotopy category, this gives a convenient criterion which allows us to test whether or not colimits in $\calC$ are universal.

\begin{remark}
Let $\calX$ be an $\infty$-category. The assumption that colimits in $\calX$ are universal can be viewed as a kind of distributive law. We have the following table of vague analogies:
$$
\begin{array}{cc}
\text{ Higher Category Theory } & \text{ Algebra } \\
\hline
\\
\text{ $\infty$-Category } & \text{ Set } \\ \\
\text{ Presentable $\infty$-Category } & \text{ Abelian Group } \\ \\
\text{ Colimits } & \text{ Sums } \\ \\
\text{ Limits } & \text{ Products } \\ \\
\colim( X_{\alpha}) \times_{S} T \simeq \colim( X_{\alpha} \times_{S} T) & (x+y)z=xz+yz \\ \\
\text{ $\infty$-Topos } & \text{ Commutative Ring } \\ \\
\end{array}$$
\end{remark}

Definition \ref{colu} has a reformulation in the language of classifying functors (\S \ref{universalfib}):

\begin{proposition}\label{gentur}
Let $\calX$ be an $\infty$-category which admits finite limits. The following conditions are equivalent:
\begin{itemize}
\item[$(1)$] The $\infty$-category $\calX$ is presentable, and colimits in $\calX$ are universal.
\item[$(2)$] The Cartesian fibration $p: \calO_{\calX} \rightarrow \calX$ is classified by a functor
$\calX^{op} \rightarrow \LPres$.
\end{itemize}
\end{proposition}

\begin{proof}
We can restate condition $(2)$ as follows: each fiber $\calX^{/U}$ of $p$ is presentable, and
each of the pullback functors $f^{\ast}: \calX^{/V} \rightarrow \calX^{/U}$ preserves small colimits.
It is clear that $(1) \Rightarrow (2)$, and that $(2)$ implies that colimits in $\calX$ are universal.
Since $\calX$ admits finite limits, it has a final object $1$; condition $(2)$ implies that
$\calX \simeq \calX^{/1}$ is presentable, which proves $(1)$.
\end{proof}

\begin{itemize}\index{gen}{coproduct!disjoint}
\item[$(iii)$] Coproducts in $\calC$ are disjoint.
\end{itemize}

If $\calC$ is an $\infty$-category which admits finite coproducts, then we will say that
{\it coproducts in $\calC$ are disjoint} if every coCartesian diagram
$$ 
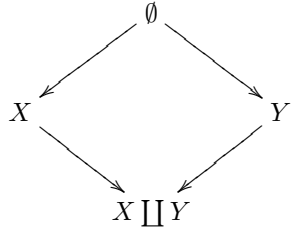
\xymatrix{ & \emptyset \ar[dr] \ar[dl] & \\ X \ar[dr] & & Y \ar[dl] \\ & X \amalg Y & }$$
is also Cartesian, provided that $\emptyset$ is an initial object of $\calC$.
More informally, to say that coproducts are disjoint is to say that the intersection of $X$ and $Y$ inside the union $X \amalg Y$ is empty. 

We now come to the most subtle and interesting of Giraud's axioms:

\begin{itemize}\index{gen}{effective equivalence relation}
\item[$(iv)$] Every equivalence relation in $\calC$ is effective.
\end{itemize}
 
Recall that if $X$ is an object in an (ordinary) category $\calC$,
then an {\it equivalence relation} $R$ on $X$ is an object of
$\calC$ equipped with a map $p: R \rightarrow X \times X$ such
that for any $S$, the induced map $$\Hom_{\calC}(S,R) \rightarrow
\Hom_{\calC}(S,X) \times \Hom_{\calC}(S,X)$$ exhibits
$\Hom_{\calC}(S,R)$ as an equivalence relation on
$\Hom_{\calC}(S,X)$.

If $\calC$ admits finite limits, then it is easy to construct
equivalence relations in $\calC$: given any map $X \rightarrow Y$
in $\calC$, the induced map $X \times_Y X \rightarrow X \times X$
is an equivalence relation on $X$. If the category $\calC$ admits
finite colimits, then one can attempt to invert this process:
given an equivalence relation $R$ on $X$, one can form the
coequalizer of the two projections $R \rightarrow X$ to obtain an
object which we will denote by $X/R$. In the category of sets, one
can recover $R$ as the fiber product $X \times_{X/R} X$. In
general, this need not occur: one always has $R \subseteq X
\times_{X/R} X$, but the inclusion may be strict (as subobjects of
$X \times X$). If equality holds, then $R$ is said to be an {\it
effective equivalence relation}, and the map $X \rightarrow X/R$
is said to be an {\it effective epimorphism}.\index{gen}{effective epimorphism}

\begin{remark}
Recall that a map $f: X \rightarrow Y$ in a category $\calC$ is said to be
a {\it categorical epimorphism} if the natural map $\Hom_{\calC}(Y,Z) \rightarrow \Hom_{\calC}(X,Z)$ is {\em injective} for every object $Z \in \calC$, so that we may identify
$\Hom_{\calC}(Y,Z)$ with a subset of $\Hom_{\calC}(X,Z)$. 
To say that $f$ is an {\em effective} epimorphism is to say that we can characterize this subset: it is the collection of all maps $g: X \rightarrow Z$ such that the diagram
$$ 
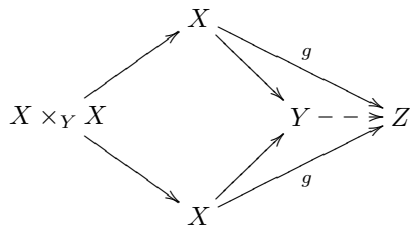
\xymatrix{ & X \ar[dr] \ar[drr]^{g} & & \\
X \times_{Y} X \ar[ur] \ar[dr] & & Y \ar@{-->}[r] & Z \\
& X \ar[ur] \ar[urr]_{g} & & }$$
commutes (which is obviously a necessary condition for the indicated dotted arrow to exist).
\end{remark}

Using the terminology introduced above, we can neatly summarize some of the fundamental properties of the category of sets:

\begin{fact}\label{factoid}
In the category of sets, every equivalence relation is
effective and the effective epimorphisms are precisely the
surjective maps.
\end{fact}

The first assertion of Fact \ref{factoid} remains valid in any
topos, and according to the axiomatic point of view it is one of
the defining features of a topos.

If $\calC$ is a category with finite limits and colimits in which
all equivalence relations are effective, then we obtain a
one-to-one correspondence between equivalence relations on an
object $X$ and {\em quotients} of $X$ (that is, isomorphism classes
of effective epimorphisms $X \rightarrow Y$). This correspondence
is extremely useful because it allows us to make elementary
descent arguments: one can deduce statements about quotients of
$X$ from statements about $X$ and about equivalence relations on
$X$ (which live over $X$). We would like to formulate an $\infty$-categorical analogue of this condition which will allow us to make similar arguments.

In the $\infty$-category $\SSet$ of spaces, the situation is more complicated.
The correct notion of surjection of spaces $X \rightarrow Y$ is a map which
induces a surjection on path components $\pi_0 X \rightarrow \pi_0 Y$.
However, in this case, the (homotopy) fiber product $R= X \times_Y X$ does
not give an equivalence relation on $X$, because the map $R
\rightarrow X \times X$ is not necessarily injective in any
reasonable sense. However, it does retain some of the pleasant
features of an equivalence relation: instead of transitivity, we
have a {\it coherently associative} composition law $R \times_X R
\rightarrow R$ (this is perhaps most familiar in the situation where $X$ is a point: in this case, $R$ can be identified with the based loop space of $Y$, which is endowed with a multiplication given by concatenation of loops). In \S \ref{gengroup} we will make this idea precise, and define {\it groupoid objects} and {\it effective groupoid objects} in an arbitrary $\infty$-category. Granting these notions
for the moment, we have a natural candidate for the $\infty$-categorical generalization
of condition $(iv)$:

\begin{itemize}\index{gen}{groupoid object}\index{gen}{groupoid object!effective}
\item[$(iv)'$] Every groupoid object of $\calC$ is effective.
\end{itemize}

\subsection{Groupoid Objects}\label{gengroup}

Let $\calC$ be a category which admits finite limits. A {\it groupoid object}\index{gen}{groupoid object!of a category} of $\calC$ is a functor $F$ from $\calC$ to the category $\Cat$ of (small) groupoids, which has the following properties:

\begin{itemize}
\item[$(1)$] There exists an object $X_0 \in \calC$ and a (functorial) identification of
$\Hom_{\calC}(C,X_0)$ with the set of objects in the groupoid $F(C)$, for each $C \in \calC$.
\item[$(2)$] There exists an object $X_1 \in \calC$ and a (functorial) identification of
$\Hom_{\calC}(C,X_1)$ with the set of morphisms in groupoid $F(C)$, for each $C \in \calC$.
\end{itemize}

\begin{example}
Let $\calC$ be the category $\Set$ of sets. Then a groupoid object of $\calC$ is simply a (small) groupoid.
\end{example}

Giving a groupoid object of a category $\calC$ is equivalent to giving a pair of objects $X_0 \in \calC$ (the ``object classifier'') and $X_1 \in \calC$ (the ``morphism classifier''), together with a collection of maps which relate $X_0$ to $X_1$ and satisfy appropriate identities, which imitate the usual axiomatics of category theory. These identities can be very efficiently encoded using the formalism of simplicial objects. For every $n \geq 0$, let
$[n]$ denote the category associated to the linearly ordered set $\{0, \ldots, n \}$, and consider
the functor $F_n: \calC \rightarrow \Set$ defined so that
$$F_n(C) = \Hom_{\Cat}( [n], F(C) ).$$
By assumption, $F_0$ and $F_1$ are representable by objects $X_0, X_1 \in \calC$.
Since $\calC$ is stable under finite limits, it follows that 
$$ F_n = F_1 \times_{F_0} \ldots \times_{F_0} F_1$$
is representable by an object $X_n = X_1 \times_{X_0} \ldots \times_{X_0} X_1$.
The objects $X_{n}$ can be assembled into a simplicial object $X_{\bigdot}$ of
$\calC$. We can think of this construction as a generalization of the process which associates
to every groupoid $\calD$ its nerve $\Nerve(\calD)$ (a simplicial set). Moreover, as in the classical case, the association $F \mapsto X_{\bigdot}$ is fully faihtful. In other words, we can identify groupoid objects of $\calC$ with the corresponding simplicial objects. Of course, not every simplicial object $X_{\bigdot}$ of $\calC$ arises via this construction. This is true if and only if certain additional conditions are met: for instance, the diagram
$$ \xymatrix{ X_2 \ar[r]^{d_0} \ar[d]^{d_2} & X_1 \ar[d]^{d_1} \\
X_1 \ar[r]^{d_0}  & X_0 }$$
must be Cartesian. 

The purpose of this section is to generalize the notion of a groupoid object to the setting where $\calC$ is an $\infty$-category. We begin by introducing the class of {\it simplicial objects} of $\calC$; we then define groupoid objects to be simplicial objects which satisfy additional conditions.

\begin{definition}\label{siminf}\index{gen}{simplicial object!of an $\infty$-category}
Let $\cDelta_{+}$\index{not}{Delta+@$\cDelta_{+}$} denote the category of finite (possibly empty) linearly ordered sets.
A {\it simplicial object} of an $\infty$-category $\calC$ is a map of $\infty$-categories
$$U_{\bigdot}: \Nerve(\cDelta)^{op} \rightarrow \calC.$$
An {\it augmented simplicial object}\index{gen}{simplicial object!augmented} of $\calC$ is a map
$$U_{\bigdot}^{+}: \Nerve(\cDelta_{+})^{op} \rightarrow \calC.$$

We let $\calC_{\Delta}$\index{not}{Delta@$\calC_{\Delta}$} denote the $\infty$-category $\Fun( \Nerve(\cDelta)^{op}, \calC)$; we will refer
to $\calC_{\Delta}$ as the {\it $\infty$-category of simplicial objects of $\calC$}. Similarly, we 
will refer to $\Fun(\Nerve(\cDelta_{+})^{op}, \calC)$ 
as the {\it $\infty$-category of augmented simplicial objects of $\calC$} and we will denote it by $\calC_{\Delta_{+}}$\index{not}{DeltaPlus@$\calC_{\Delta_{+}}$}.

If $U_{\bigdot}$ is an (augmented) simplicial object of $\calC$ and $n \geq 0$ ($n \geq -1$), we will write $U_{n}$ for the object $U([n]) \in \calC$. 
\end{definition}

\begin{remark}
In the case where $\calC$ is the nerve of an ordinary category $\calD$, Definition \ref{siminf} recovers the usual notion of a simplicial object of $\calD$. More precisely, the $\infty$-category $\calC_{\Delta}$ of simplicial objects of $\calC$ is naturally isomorphic to the nerve of the category of simplicial objects of $\calD$.
\end{remark}

\begin{lemma}\label{silling}
Let $f: X \rightarrow Y$ be a map of simplicial sets. Suppose that:

\begin{itemize}
\item[$(1)$] The map $f$ induces a bijection $X_0 \rightarrow Y_0$ on vertex sets.
\item[$(2)$] The simplicial set $Y$ is a Kan complex.
\item[$(3)$] The map $f$ has the right lifting property with respect to every horn inclusion
$\Lambda^n_i \subseteq \Delta^n$, for $n \geq 2$.
\item[$(4)$] The map $f$ is a weak homotopy equivalence.
\end{itemize}

Then $f$ is a trivial Kan fibration.
\end{lemma}

\begin{proof}
In view of condition $(4)$, it suffices to prove that $f$ is a Kan fibration. In other words, we must show that $p$ has the right lifting property with respect to every horn inclusion $\Lambda^n_i \subseteq \Delta^n$. If $n > 1$, this follows from $(3)$. We may therefore reduce to the case where $n=1$; by symmetry, we may suppose that $i=0$.

Let $e: y \rightarrow y'$ be an edge of $Y$. Condition $(1)$ implies that there is a (unique) pair of vertices $x,x' \in X_0$ with $y=f(x)$, $y' = f(x')$.
Since $f$ is a homotopy equivalence, there is a path $p$ from $x$ to $x'$ in the topological space $|X|$, such that the induced path $|f| \circ p$ in $|Y|$ is homotopic to $e$ via a homotopy which keeps the endpoints fixed. By cellular approximation we may suppose that this path is contained in the $1$-skeleton of $|X|$. 
Consequently, there is a positive integer $k$, a sequence
of vertices $\{ z_0, \ldots, z_k \}$ with $z_0 = x$, $z_k = x'$ such
each adjacent pair $(z_i, z_{i+1})$ is joined by an edge $p_{i}$ (running in either direction), such that $p$ is homotopic (relative to its boundary) to the path obtained by concatenating the edges $p_i$. Using conditions $(2)$ and $(3)$, we note that $X$ has the extension property with respect to the inclusion $\Lambda^n_i \subseteq \Delta^n$ for each $n \geq 2$.
It follows that we may assume that
$p_i$ runs from $z_i$ to $z_{i+1}$: if it runs in the opposite direction, then we can extend the map $$(p_i, s_0 z_i, \bigdot): \Lambda^2_2 \rightarrow X$$ to a $2$-simplex
$\sigma: \Delta^2 \rightarrow X$, and then replace $p_i$ by $d_2 \sigma$.

Without loss of generality, we may suppose that $k > 0$ is chosen as small as possible.
We claim that $k = 1$. Otherwise, choose an extension $\tau: \Delta^2 \rightarrow X$ of the map
$$ (p_1, \bigdot, p_0 ): \Lambda^2_1 \rightarrow X.$$
We can then replace the initial segment
$$ z_0 \stackrel{p_0}{\rightarrow} z_1 \stackrel{p_1}{\rightarrow} z_2$$
by the edge $d_1(\tau): z_0 \rightarrow z_2$ and obtain a shorter path from $x$ to $x'$, contradicting the minimality of $k$. 

The edges $e$ and $f(p_0)$ are homotopic in $Y$ relative to their endpoints. Using $(3)$, we see that $p_0$ is homotopic (relative to its endpoints) to an edge $\overline{e}$ which satisfies $f(\overline{e}) = e$. This completes the proof that $f$ is a Kan fibration.
\end{proof}

\begin{notation}
Let $K$ be a simplicial set. We let $\cDelta_{/K}$ denote the {\it category of simplices of $K$}\index{gen}{category!of simplices} defined in \S \ref{quasilimit1}. The objects of $\cDelta_{/K}$\index{not}{DeltaK@$\cDelta_{/K}$} are pairs $(J, \eta)$ where $J$ is an object of $\cDelta$ and $\eta \in \Hom_{\sSet}( \Delta^J, K)$. A morphism from $(J, \eta)$ to $(J', \eta')$ is a commutative diagram
$$ \xymatrix{ \Delta^{J} \ar[rr] \ar[dr] & & \Delta^{J'} \ar[dl] \\
& K. & }$$
Equivalently, we can describe $\cDelta_{/K}$ as the fiber product 
$\cDelta \times_{ \sSet } (\sSet)_{/K}$.

If $\calC$ is an $\infty$-category, $U: \Nerve(\cDelta)^{op} \rightarrow \calC$ is a simplicial object of $\calC$, and $K$ is a simplicial set, then we let $U[K]$ denote the composite map
$$ \Nerve(\cDelta_{/K})^{op} \rightarrow \Nerve(\cDelta)^{op} \rightarrow \calC.$$
\end{notation}

\begin{proposition}\label{grpobjdef}\index{gen}{groupoid object!of an $\infty$-category}
Let $\calC$ be an $\infty$-category and $U: \Nerve(\cDelta)^{op} \rightarrow \calC$ a simplicial
object of $\calC$. The following conditions are equivalent:
\begin{itemize}

\item[$(1)$] For every weak homotopy equivalence $f: K \rightarrow K'$ of simplicial
sets which induces a bijection $K_0 \rightarrow K'_0$ on vertices, the induced map
$\calC_{/U[K']} \rightarrow \calC_{/U[K]}$ is a categorical equivalence.

\item[$(2)$] For every cofibration $f: K \rightarrow K'$ of simplicial sets which is a weak homotopy equivalence and bijective on vertices, the induced map
$\calC_{/U[K']} \rightarrow \calC_{/U[K]}$ is a categorical equivalence.

\item[$(2')$] For every cofibration $f: K \rightarrow K'$ of simplicial sets which is a weak homotopy equivalence and bijective on vertices, the induced map
$\calC_{/U[K']} \rightarrow \calC_{/U[K]}$ is a trivial fibration.

\item[$(3)$] For every $n \geq 2$ and every $0 \leq i \leq n$, the induced map
$\calC_{/U[\Delta^n]} \rightarrow \calC_{/U[\Lambda^n_i]}$ is a categorical equivalence.

\item[$(3')$] For every $n \geq 2$ and every $0 \leq i \leq n$, the induced map
$\calC_{/U[\Delta^n]} \rightarrow \calC_{/U[\Lambda^n_i]}$ is a trivial fibration.

\item[$(4)$] For every $n \geq 0$ and every partition $[n] = S \cup S'$ such that
$S \cap S'$ consists of a single element $s$, the induced map
$\calC_{/U[\Delta^n]} \rightarrow \calC_{ /U[K]}$
is a categorical equivalence, where $K = \Delta^{S} \amalg_{ \{s\} } \Delta^{ S'} \subseteq \Delta^n$.

\item[$(4')$] For every $n \geq 0$ and every partition $[n] = S \cup S'$ such that
$S \cap S'$ consists of a single element $s$, the induced map
$\calC_{/U[\Delta^n]} \rightarrow \calC_{ /U[K]}$
is a trivial fibration, where $K = \Delta^{S} \amalg_{ \{s\} } \Delta^{ S'} \subseteq \Delta^n$.

\item[$(4'')$] For every $n \geq 0$ and every partition $[n] = S \cup S'$ such that
$S \cap S'$ consists of a single element $s$, the diagram
$$ \xymatrix{ U([n]) \ar[r] \ar[d] & U(S) \ar[d] \\
U(S') \ar[r] & U( \{s\}) }$$
is a pullback square in the $\infty$-category $\calC$.
\end{itemize}
\end{proposition}

\begin{proof}
The dual of Proposition \ref{sharpen} implies that any monomorphism $K \rightarrow K'$ of simplicial sets induces a right fibration $\calC_{/U[K]} \rightarrow \calC_{/U[K]}$. By Corollary \ref{heath}, a right fibration is a trivial fibration if and only if it is a categorical equivalence.
This proves that $(2) \Leftrightarrow (2')$, $(3) \Leftrightarrow (3')$, and $(4) \Leftrightarrow (4')$. The implications $(1) \Rightarrow (2) \Rightarrow (3)$ are obvious. 

We now prove that $(3)$ implies $(1)$. Let $A$ denote the class of all morphisms $f: K' \rightarrow K$ which induce a categorical equivalence $\calC_{/U[K]} \rightarrow \calC_{/U[K']}$. 
Let $A'$ denote the class of all {\em cofibrations} which have the same property; equivalently, $A'$ is the class of all cofibrations which induce a trivial fibration $\calC_{/U[K]} \rightarrow \calC_{/U[K']}$. From this characterization it is easy to see that $A'$ is weakly saturated. Let $A''$ be the weakly saturated class of morphisms generated by the inclusions $\Lambda^n_i \subseteq \Delta^n$ for $n > 1$.
If we assume $(3)$, then we have the inclusions $A'' \subseteq A' \subseteq A$.

Let $f: K \rightarrow K'$ be an arbitrary morphism of simplicial sets. By Proposition \ref{quillobj}, we can choose a map $h': K' \rightarrow M'$ which belongs to $A''$, where $M'$ has the extension property with respect to $\Lambda^n_i \subseteq \Delta^n$ for $n > 1$ and is therefore a Kan complex. Applying Proposition \ref{quillobj} again, we can construct a commutative diagram
$$ \xymatrix{ K \ar[d]^{f} \ar[r]^{h} & M \ar[d]^{g} \\
K' \ar[r]^{h'} & M' }$$
where the horizontal maps belong to $A''$ and $g$ has the right lifting property with
respect to every morphism in $A''$. If $f$ is a weak homotopy equivalence which is bijective on vertices, then $g$ has the same properties, so that $g$ is a trivial fibration by Lemma \ref{silling}. It follows that $g$ has the right lifting property with respect to the cofibration $g \circ h: K \rightarrow M'$, so that $g \circ h$ is a retract of $h$ and therefore belongs to $A''$. Since $g \circ h = h' \circ f$ and $h'$ belong to $A'' \subseteq A$, it follows
that $f$ belongs to $A$.

It is clear that $(1) \Rightarrow (4)$. We next prove that $(4') \Rightarrow (3)$. We must show
that if $n > 1$, then every inclusion $\Lambda^n_i \subseteq \Delta^n$ belongs to the class $A$ defined above. The proof is by induction on $n$. Replacing $i$ by $n-i$ if necessary, we may suppose that $i < n$. If $(n,i) \neq (2,0)$, we consider the composition
$$ \Delta^{n-1} \amalg_{ \{n-1\} } \Delta^{ \{n-1, n\} } \stackrel{f}{\hookrightarrow} \Lambda^n_i \stackrel{f'}{\hookrightarrow} \Delta^n.$$
Here $f$ belongs to $A'$ by the inductive hypothesis and $f' \circ f$ belongs to $A'$ by virtue of the assumption $(4')$; therefore $f'$ also belongs to $A$. If $n=2$ and $i=0$, then we observe
that the inclusion $\Lambda^2_1 \subseteq \Delta^2$ is of the form
$\Delta^S \amalg_{ \{s\} } \Delta^{S'} \subseteq \Delta^2$, where $S = \{0,1\}$ and
$S' = \{0,2\}$. 

To complete the proof, we show that $(4)$ is equivalent to $(4'')$.  
Fix $n \geq 0$, let $S \cup S' = [n]$ be such that $S \cap S' = \{s\}$, and let
$K = \Delta^{S} \amalg_{ \{s\} } \Delta^{S'} \subseteq \Delta^n$.
Let
$\calI'$ denote the full subcategory of $\cDelta_{/\Delta^n}$ spanned by the objects
$[n]$, $S$, $S'$, and $\{s\}$. Let $\calI \subseteq \calI'$ be the full subcategory obtained by omitting the object $[n]$. Let $p'$ denote the composition
$$ \Nerve(\calI')^{op} \rightarrow \Nerve(\cDelta)^{op} \stackrel{U}{\rightarrow} \calC$$ and let
$p = p' | \Nerve(\calI)^{op}$. 
Consider the diagram
$$ \xymatrix{ \calC_{/U[\Delta^n]} \ar[r] \ar[d]^{u} & \calC_{/U[K]} \ar[d]^{v} \\
\calC_{/p'} \ar[r] & \calC_{/p}. }$$ 
Condition $(4)$ asserts that the upper horizontal map is a categorical equivalence, and condition $(4'')$ asserts that the lower horizontal map is a categorical euivalence. To prove that
$(4) \Leftrightarrow (4'')$, it suffices to show that the vertical maps $u$ and $v$ are categorical equivalences.

We have a commutative diagram
$$ \xymatrix{ \calC_{/U[\Delta^n]} \ar[rr]^{u} \ar[dr] & & \calC_{/p'} \ar[dl] \\
& \calC_{/U[\Delta^n]}.& }$$
Since $\Delta^n$ is a final object of both $\cDelta_{/\Delta^n}$ and $\calI'$, the unlabelled maps are trivial fibrations. It follows that $u$ is a categorical equivalence.

To prove that $v$ is a categorical equivalence, it suffices to show that the inclusion
$g: \calI \subseteq \cDelta_{/K}$ induces a right anodyne map
$$\Nerve(g): \Nerve(\calI) \subseteq \Nerve(\cDelta_{/K})$$ of simplicial sets. We observe that the functor $g$ has a left adjoint $f$, which associates to each simplex $\sigma: \Delta^m \rightarrow
K$ the smallest simplex in $\calI$ which contains the image of $\sigma$. 
The map $\Nerve(g)$ is a section of $\Nerve(f)$, and there is a (fiberwise) simplicial homotopy from
$\id_{ \Nerve(\cDelta_{/K}) }$ to $\Nerve(g) \circ \Nerve(f)$. We now invoke Proposition \ref{trull11} to deduce that $\Nerve(g)$ is right anodyne, as desired.
\end{proof}

\begin{definition}\index{gen}{groupoid object!of an $\infty$-category}
Let $\calC$ be an $\infty$-category. We will often denote simplicial objects of $\calC$ by
$U_{\bigdot}$, and write $U_{n}$ for $U_{\bigdot}([n]) \in \calC$. We will say that a simplicial object $U_{\bigdot} \in \calC_{\Delta}$ is a {\em groupoid object} of $\calC$ if it satisfies the equivalent conditions of Proposition \ref{grpobjdef}. We will let $\Grp(\calC)$\index{not}{Grp@$\Grp(\calC)$} denote the full subcategory of $\calC_{\Delta}$ spanned by the groupoid objects of $\calC$.
\end{definition}

\begin{remark}
It follows from the proof of Proposition \ref{grpobjdef} that to verify that a simplicial object
$X_{\bigdot} \in \calC_{\Delta}$ is a groupoid object, we need only verify condition $(4'')$
in a small class of specific examples, but we will not need this observation.
\end{remark}

\begin{proposition}\label{refle}
Let $\calC$ be a presentable $\infty$-category. The full subcategory $\Grp(\calC) \subseteq \calC_{\Delta}$ is strongly reflective.
\end{proposition}

\begin{proof}
Let $n \geq 0$ and $[n] = S \cup S'$ be as in the statement of $(4'')$ of Proposition \ref{grpobjdef}. 
Let $\calD(S,S') \subseteq \calC_{\Delta}$
be the full subcategory consisting of those simplicial objects $U \in \calC_{\Delta}$
for which the associated diagram
$$ \xymatrix{ U([n]) \ar[r] \ar[d] & U(S) \ar[d] \\
U(S') \ar[r] & U( \{s\}) }$$
is Cartesian. Lemmas \ref{stur1} and \ref{stur2} imply that $\calD(S,S')$ is a strongly reflective subcategory of $\calC_{\Delta}$. 
Let $\calD$ denote the intersection of all these subcategories, taken over all $n \geq 0$ and all such decompositions $[n] = S \cup S'$. Lemma \ref{stur3} implies that
$\calD \subseteq \calC_{\Delta}$ is strongly reflective, and Proposition \ref{grpobjdef}
implies that $\calD = \Grp(\calC)$.
\end{proof}

Our next step is to exhibit a large class of examples of groupoid objects. We first sketch the idea.
Suppose that $\calC$ is an $\infty$-category which admits finite limits, and let $u: U \rightarrow X$ be a morphism in $\calC$. Using this data, we can construct a simplicial object $U_{\bigdot}$ of $\calC$, where $U_{n}$ is given by the $(n+1)$-fold fiber power of $U$ over $X$. In order to describe this construction more precisely, we need to introduce a bit of notation.

\begin{notation}
Let $\cDelta^{\leq n}_{+}$\index{not}{DeltaLeq0Plus@$\cDelta^{\leq n}_{+}$} denote the full subcategory of $\cDelta_{+}$ spanned by the objects $\{ [k] \}_{ -1 \leq k \leq n}$.
\end{notation}

\begin{proposition}\label{strump}
Let $\calC$ be an $\infty$-category, and let $U_{\bigdot}^{+}: \Nerve(\cDelta_{+})^{op} \rightarrow \calC$ be an augmented simplicial object of $\calC$. The following conditions are equivalent:
\begin{itemize}
\item[$(1)$] The augmented simplicial object $U_{\bigdot}^{+}$ is a right Kan extension
of $U_{\bigdot}^{+}| \Nerve(\cDelta^{\leq 0}_{+})^{op}$.
\item[$(2)$] The underlying simplicial object $U_{\bigdot}$ is a groupoid object of $\calC$, and the diagram $U_{\bigdot}^{+} | \Nerve( \cDelta^{\leq 1}_{+})^{op}$ is a pullback square
$$ \xymatrix{ U_{1} \ar[r] \ar[d] & U_0 \ar[d] \\
U_{0} \ar[r] & U_{-1} }$$ 
in the $\infty$-category $\calC$.
\end{itemize}
\end{proposition}

\begin{proof}
Suppose first that $(1)$ is satisfied. It follows immediately from the definition of right Kan extensions that the diagram $$ \xymatrix{ U_{1} \ar[r] \ar[d] & U_0 \ar[d] \\
U_{0} \ar[r] & U_{-1} }$$ 
is a pullback. To prove that $U_{\bigdot}$ is a groupoid, we show that $U_{\bigdot}$ satisfies criterion $(4'')$ of Proposition \ref{grpobjdef}. Let $S$ and $S'$ be sets with union $[n]$ and intersection $S \cap S' = \{s\}$. Let $\calI$ be the nerve of the category
$(\cDelta_{+})_{/\Delta^n}$. For each subset $J \subseteq [n]$, let 
$\calI(J)$ denote the full subcategory of $\calI$ spanned by the initial object
together with the inclusions $\{j\} \rightarrow \Delta^n$, $j \in S$. By assumption,
$U_{\bigdot}^{+}$ exhibits $U_{\bigdot}(S)$ as a limit of
$U_{\bigdot}^{+} | \Nerve(\calI(S))$, $U_{\bigdot}(S')$ as a limit of
$U_{\bigdot}^{+} | \Nerve(\calI(S'))$, $U_{\bigdot}([n])$ as a limit of
$U_{\bigdot}^{+} | \Nerve(\calI( [n] ))$, and $U_{\bigdot}( \{s\} )$ as a limit of
$U_{\bigdot}^{+} | \Nerve(\calI(\{s\}))$. It follows from Corollary \ref{util}
that the diagram
$$ \xymatrix{ U_{\bigdot}( [n] ) \ar[r] \ar[d] & U_{\bigdot}(S) \ar[d] \\
U_{\bigdot}(S') \ar[r] & U_{\bigdot}( \{s\} ) }$$
is a pullback.

We now prove that $(2)$ implies $(1)$. Using the above notation, we must show that for each
$n \geq -1$, $U_{\bigdot}^{+}$ exhibits $U_{\bigdot}^{+}([n])$ as a limit of
$U_{\bigdot}^{+} | \calI( [n] )$. For $n \leq 0$, this is obvious; for $n=1$ it is equivalent to the assumption that
$$ \xymatrix{ U_{1} \ar[r] \ar[d] & U_0 \ar[d] \\
U_{0} \ar[r] & U_{-1} }$$ is a pullback diagram. We prove the general case by induction on $n$.
Using the inductive hypothesis, we conclude that $U_{\bigdot}(\Delta^{S})$ is a limit of
$U_{\bigdot}^{+} | \calI(S)$ for all {\em proper} subsets $S \subset [n]$. Choose
a decomposition $\{0, \ldots, n\} = S \cup S'$, where $S \cap S' = \{s\}$. According to
Proposition \ref{train}, the desired result is equivalent to the assertion that
$$ \xymatrix{ U_{\bigdot}([n]) \ar[r] \ar[d] & U_{\bigdot}(S) \ar[d] \\
U_{\bigdot}(S') \ar[r] & U_{\bigdot}( \{s\} ) }$$
is a pullback diagram, which follows from our assumption that $U_{\bigdot}$ is a groupoid object of $\calC$.
\end{proof}

We will say that augmented simplicial object $U_{\bigdot}^{+}$ in an $\infty$-category $\calC$ is a {\it \Cech nerve} if it satisfies the equivalent conditions of Proposition \ref{strump}. In this case, $U_{\bigdot}^{+}$ is determined up to equivalence by the map $u: U_0 \rightarrow U_{-1}$; we will also say that $U_{\bigdot}^{+}$ is the {\it \Cech nerve of $u$}\index{gen}{Cech nerve@\Cech nerve}. 

\begin{notation}
Let $U_{\bigdot}$ be a simplicial object in an $\infty$-category $\calC$. We may regard
$U_{\bigdot}$ as a diagram in $\calC$ indexed
by $\Nerve(\cDelta)^{op}$. We let $|U_{\bigdot}|: \Nerve(\cDelta_{+})^{op} \rightarrow \calC$
denote a colimit for $U_{\bigdot}$ (if such a colimit exists). We will refer to any such colimit as a {\it geometric realization}\index{gen}{geometric realization} of $U_{\bigdot}$. 
\end{notation}

\begin{remark}
Note that we are regarding $|U_{\bigdot}|$ as a colimit diagram in $\calC$, not as an object of $\calC$. We also note that our notation is somewhat abusive, since $|U_{\bigdot}|$ is not uniquely determined by $U_{\bigdot}$. However, if a colimit of $U_{\bigdot}$ exists, then it is determined up to contractible ambiguity.
\end{remark}

\begin{definition}
Let $U_{\bigdot}$ be a simplicial object of an $\infty$-category $\calC$. We will say that
$U_{\bigdot}$ is an {\it effective groupoid}\index{gen}{groupoid object!effective} if can be extended to a colimit diagram $U^{+}_{\bigdot}: \Nerve( \cDelta_{+})^{op} \rightarrow \calC$ such that
$U^{+}_{\bigdot}$ is a \Cech nerve.
\end{definition}

\begin{remark}
It follows immediately from characterization $(3)$ of Proposition \ref{strump} that any effective groupoid $U_{\bigdot}$ is a groupoid.
\end{remark}

We can now state the $\infty$-categorical counterpart of Fact \ref{factoid}: every groupoid object in $\SSet$ is effective. This statement is somewhat less trivial than its classical analogue. For example, a groupoid object $U_{\bigdot}$ in $\SSet$ with $U_0 = \ast$ can be thought of
as a space $U_{1}$ equipped with a coherently associative multiplication operation. If $U_{\bigdot}$ is effective, then there exists a fiber diagram
$$ \xymatrix{ U_1 \ar[r] \ar[d] & \ast \ar[d] \\
\ast \ar[r] & U_{-1} }$$
so that $U_{1}$ is homotopy equivalent to a loop space.
This is an classical result (see, for example, \cite{stasheff}). We will give a somewhat indirect proof in the next section.

\subsection{$\infty$-Topoi and Descent}\label{magnet}

In this section, we will describe an elegant characterization of the notion of an $\infty$-topos, based on the theory of {\em descent}. We begin by explaining the idea in informal terms.
Let $\calX$ be an $\infty$-category. To each object $U$ of $\calX$ we can associate
the overcategory $\calX^{/U}$. If $\calX$ admits finite limits, then this construction gives a contravariant functor from $\calX$ to the $\infty$-category
$\widehat{ \Cat}_{\infty}$ of (not necessarily small) $\infty$-categories. If $\calX$ is an $\infty$-topos, then this functor carries colimits in $\calX$ to limits of $\infty$-categories. In other words, if an object $X \in \calX$ is obtained as the colimit of some diagram $\{ X_{\alpha} \}$ in $\calX$, then giving a morphism $Y \rightarrow X$ is equivalent to a suitably compatible diagram of morphisms
$\{ Y_{\alpha} \rightarrow X_{\alpha} \}$. Moreover, we will eventually show that this property {\em characterizes} the class of $\infty$-topoi. The ideas presented in this section are due to Charles Rezk.

\begin{definition}\index{gen}{Cartesian transformation}
Let $\calX$ be an $\infty$-category, $K$ a simplicial set, and $p,q: K \rightarrow \calX$
two diagrams. We will say that a natural transformation $\alpha: p \rightarrow q$ is
{\it Cartesian} if, for each edge $\phi: x \rightarrow y$ in $K$, the associated diagram
$$ \xymatrix{ p(x) \ar[r]^{p(\phi)} \ar[d]^{\alpha(x)} & p(y) \ar[d]^{\alpha(y)} \\
q(x) \ar[r]^{q(\phi)} & q(y) }$$
is a pullback in $\calX$. 
\end{definition}

\begin{lemma}\label{ib0}
Let $\calX$ be an $\infty$-category, and 
let $\alpha: p \rightarrow q$ be a natural transformation of diagrams
$p,q: K \diamond \Delta^0 \rightarrow \calX$. Suppose that, for every vertex
$x$ of $K$, the associated transformation
$$ p | \{x\} \diamond \Delta^0 \rightarrow q| \{x\} \diamond \Delta^0 $$
is Cartesian. Then $\alpha$ is Cartesian.
\end{lemma}

\begin{proof}
Let $z$ be the ``cone point'' of $K \diamond \Delta^0$. We note that to every
edge $e: x \rightarrow y$ in $K \diamond \Delta^0$ we can associate a diagram
$$ \xymatrix{ x \ar[d]^{e} \ar[r] \ar[dr]^{g} & z \ar[d]^{\id_{z}} \\
y \ar[r] & z. }$$
The transformation $\alpha$ restricts to a Cartesian transformation on the horizontal edges and the right vertical edge, either by assumption or because they are degenerate. Applying Lemma \ref{transplantt}, we deduce first that $\alpha(g)$ is a Cartesian transformation, then that
$\alpha(e)$ is a Cartesian transformation.
\end{proof}

The condition that an $\infty$-category has universal colimits can be formulated in the language of Cartesian transformations:

\begin{lemma}\label{ib1}
Let $\calX$ be a presentable $\infty$-category. The following conditions are equivalent:
\begin{itemize} 
\item[$(1)$] Colimits in $\calX$ are universal.

\item[$(2)$] Let $p, q: (K^{\triangleright} \diamond \Delta^0)
\rightarrow \calX$ be diagrams which carry $\Delta^0$ to vertices $X, Y \in \calX$, and let
$\alpha: p \rightarrow q$ be a Cartesian transformation. 
If the map $q': K^{\triangleright} \rightarrow \calX^{/Y}$ associated to
$q$ is a colimit diagram, then the map $p': K^{\triangleright} \rightarrow \calX^{/X}$ associated to $p$ is a colimit diagram.

\item[$(3)$] Let $p, q: K \star \Delta^1
\rightarrow \calX$ be diagrams which carry $\{1\}$ to vertices $X, Y \in \calX$, and let
$\alpha: p \rightarrow q$ be a Cartesian transformation. 
If the map $K^{\triangleright} \rightarrow \calX_{/Y}$ associated to
$q$ is a colimit diagram, then the map $K^{\triangleright} \rightarrow \calX_{/X}$ associated to $p$ is a colimit diagram.

\item[$(4)$] Let $p, q: K \star \Delta^1
\rightarrow \calX$ be diagrams which carry $\{1\}$ to vertices $X, Y \in \calX$, and let
$\alpha: p \rightarrow q$ be a Cartesian transformation. 
If $q | K \star \{0\}$ is a colimit diagram, then $p| K \star \{0\}$ is a colimit diagram.

\item[$(5)$] Let $\alpha: p \rightarrow q$ be a Cartesian transformation of diagrams
$K^{\triangleright} \rightarrow \calX$. If $q$ is a colimit diagram, then $p$ is a colimit diagram.
\end{itemize}
\end{lemma}

\begin{proof}
Assume that $(1)$ is satisfied; we will prove $(2)$. The transformation $\alpha$ induces a map $f: X \rightarrow Y$.  Consider the map
$$ \phi: \Fun(K^{\triangleright}, \calO_{\calX}) \rightarrow \Fun(K^{\triangleright}, \calX )$$
given by evaluation at the final vertex of $\Delta^1$.
Let $\delta(f)$ denote the image of $f$ under the diagonal map $\delta: \calX \rightarrow
\Fun(K^{\triangleright}, \calX)$. Then we may identify $\alpha$ with an edge $e$ of
$\Fun(K^{\triangleright}, \calO_{\calX})$ which covers $\delta(f)$. Since $\alpha$ is Cartesian, we can apply Lemma \ref{charpull} and Proposition \ref{doog} to deduce that $e$ is $\phi$-Cartesian.
The composition $f^{\ast} \circ q'$ is the origin of a $\phi$-Cartesian edge $e': f^{\ast} \circ q' \rightarrow q'$ of $\Fun(K^{\triangleright}, \calO_{\calX})$ covering $\delta(f)$, so we conclude that
$f^{\ast} \circ q'$ and $p'$ are equivalent in 
$\Fun(K^{\triangleright}, \calX^{/X})$. Since $q'$ is a colimit diagram and $f^{\ast}$ preserves colimits, $f^{\ast} \circ q'$ is a colimit diagram. It follows that $p'$ is a colimit diagram, as desired.

We now prove that $(2) \Rightarrow (1)$. Let $f: X \rightarrow Y$ be a morphism in $\calX$, and let
$q': K^{\triangleright} \rightarrow \calX^{/Y}$ be a colimit diagram. Choose a $\phi$-Cartesian
edge $e': f^{\ast} \circ q' \rightarrow q'$ as above, corresponding to a natural transformation
$\alpha: p \rightarrow q$ of diagrams $p,q: (K^{\triangleright} \diamond \Delta^0) \rightarrow \calX$.
Since $e$ is $\phi$-Cartesian, we may invoke Proposition \ref{doog} and
Lemma \ref{charpull} to deduce that $\alpha$ restricts to a Cartesian transformation
$p|(\{x\} \diamond \Delta^0) \rightarrow q|( \{x\} \diamond \Delta^0 )$ for every vertex $x$ of
$K^{\triangleright}$. It follows from Lemma \ref{ib0} that $\alpha$ is Cartesian. Invoking $(2)$, we conclude that $f^{\ast} \circ q'$ is a colimit diagram, as desired. 

The equivalence $(2) \Leftrightarrow (3)$ follows from Proposition \ref{rub3}, and the equivalence $(3) \Leftrightarrow (4)$ follows from Proposition \ref{needed17}. The implication $(5) \Rightarrow (4)$ is obvious. The converse implication $(4) \Rightarrow (5)$ follows from the observation that $K \star \{0\}$ is a retract of $K \star \Delta^1$.
\end{proof}

\begin{notation}\label{ugaboo}\index{gen}{stable under pullbacks}
Let $\calX$ be an $\infty$-category which admits pullbacks, and let $S$ be a class of morphisms in $\calX$. We will say that $S$ is {\it stable under pullback} if for any pullback diagram
$$ \xymatrix{ X' \ar[r] \ar[d]^{f'} & X \ar[d]^{f} \\
Y' \ar[r] & Y }$$
in $\calX$ such that $f$ belongs to $S$, $f'$ also belongs to $S$. We let
$\calO_{\calX}^{S}$ denote the full subcategory of $\calO_{\calX}$ spanned by $S$, and
$\calO_{\calX}^{(S)}$ the subcategory of $\calO_{\calX}$ whose objects are 
are elements of $S$, and whose morphisms are pullback diagrams as above.
We observe that evaluation at $\{1\} \subseteq \Delta^1$ induces a Cartesian fibration $\calO_{\calX}^{S} \rightarrow \calX$, which restricts to a right fibration
$\calO_{\calX}^{(S)} \rightarrow \calX$ (Corollary \ref{relativeKan}).\index{not}{OXS@$\calO_{\calX}^{S}$}\index{not}{OX(S)@$\calO_{\calX}^{(S)}$} 
\end{notation}

\begin{lemma}\label{ib2}
Let $\calX$ be a presentable $\infty$-category, and suppose that colimits in $\calX$ are universal. 
Let $S$ be a class of morphisms of $\calX$ which is stable under pullback, $K$ a small simplicial
set, and $\overline{q}: K^{\triangleright} \rightarrow \calX$ a colimit diagram.
The following conditions are equivalent:
\begin{itemize}
\item[$(1)$] The composition $f \circ \overline{q}: K^{\triangleright} \rightarrow \hat{\Cat}_{\infty}^{op}$ is a colimit diagram, where $f: \calX \rightarrow \hat{\Cat}_{\infty}^{op}$
classifies the Cartesian fibration $\calO_{\calX}^{S} \rightarrow \calX$. 

\item[$(2)$] The composition $f' \circ \overline{q}: K^{\triangleright} \rightarrow \hat{\SSet}^{op}$ is a colimit diagram, where $f: \calX \rightarrow \hat{\SSet}^{op}$
classifies the right fibration $\calO_{\calX}^{(S)} \rightarrow \calX$. 

\item[$(3)$] For every natural transformation $\overline{\alpha}: \overline{p} \rightarrow \overline{q}$ of colimit diagrams $K^{\triangleright} \rightarrow \calX$, if $\alpha = \overline{\alpha} | K$ is a Cartesian transformation and $\alpha(x) \in S$ for each vertex $x \in K$, then $\overline{\alpha}$ is a Cartesian transformation and $\overline{\alpha}(\infty) \in S$, where $\infty$ denotes the cone point of $K^{\triangleright}$.
\end{itemize}
\end{lemma}

\begin{proof}
Let $\overline{\calC} = \Fun(K^{\triangleright}, \calX)^{/\overline{q}}$ and
$\calC = \Fun(K,\calX)^{/q}$. Let $\overline{\calC}^{0}$ denote the full subcategory of
$\overline{\calC}$ spanned by {\em Cartesian} natural tranformations $\overline{\alpha}: \overline{p} \rightarrow \overline{q}$ with the property that
$\overline{\alpha}(x)$ belongs to $S$ for each vertex $x \in K^{\triangleright}$, and let $\calC^0$ be defined similarly. Finally, let
$\overline{\calC}^{1}$ denote the full subcategory of $\overline{\calC}$ spanned by those natural transformations $\overline{\alpha}: \overline{p} \rightarrow \overline{q}$ such that $\overline{p}$ is a colimit diagram, $\alpha = \overline{\alpha} | K$ is a Cartesian transformation, and
$\overline{\alpha}(x)$ belongs to $S$ for each vertex $x \in K$. Lemma \ref{ib1} implies that $\overline{\calC}^{0} \subseteq \overline{\calC}^{1}$. 

Let $\calD$ denote the full subcategory of $\Fun(K^{\triangleright}, \calC)$ spanned by the colimit diagrams. Proposition \ref{lklk} asserts that the restriction map $\calD \rightarrow \Fun(K,\calC)$ is a trivial fibration. It follows that the associated map $\calD^{/\overline{q}} \rightarrow \Fun(K,\calC)^{/q}$ is also a trivial fibration, and therefore restricts to a trivial fibration
$\overline{\calC}^{1} \rightarrow \calC^{0}$.

According to Proposition \ref{charcatlimit}, condition $(1)$ is equivalent to the assertion that
the projection $\overline{\calC}^{0} \rightarrow \calC^0$ is an equivalence of $\infty$-categories.
In view of the above argument, this is equivalent to the assertion that the fully faithful inclusion
$\overline{\calC}^{0} \subseteq \overline{\calC}^{1}$ is essentially surjective. Since
$\overline{\calC}^{0}$ is clearly stable under equivalence in $\overline{\calC}$, $(1)$
holds if and only if $\overline{\calC}^0 = \overline{\calC}^1$, which is manifestly equivalent to $(3)$.
The proof that $(2) \Leftrightarrow (3)$ is similar, using Proposition \ref{charspacelimit} in place
of Proposition \ref{charcatlimit} and the maximal Kan complexes contained in
$\overline{\calC}^{0}$, $\overline{\calC}^{1}$, and $\calC^0$.
\end{proof}

\begin{lemma}\label{sumoto}
Let $\calX$ be a presentable category in which colimits are universal.
Let $f: X \rightarrow \emptyset$ be a morphism in $\calX$, where $\emptyset$ is an initial
object of $\calX$. Then $X$ is also initial.
\end{lemma}

\begin{proof}
Observe that $\id_{\emptyset}$ is both an initial object of $\calX^{/ \emptyset}$ (Proposition \ref{needed17}) and a final object of $\calX^{/\emptyset}$. Let $f^{\ast}: \calX^{/\emptyset} \rightarrow \calX^{/X}$ be a pullback functor. Then $f^{\ast}$ preserves limits (since it is a right adjoint) and colimits (since colimits in $\calX$ are universal). Therefore $f^{\ast} \id_{\emptyset}$
is both initial and final in $\calX^{/X}$. It follows that $\id_{X}: X \rightarrow X$, being another final object of $\calX^{/X}$, is also initial. Applying Proposition \ref{needed17}, we deduce that
$X$ is an initial object of $\calX$, as desired.
\end{proof}

\begin{lemma}\label{ib3}
Let $\calX$ be a presentable $\infty$-category in which colimits are universal, and let
$S$ be a class of morphisms in $\calX$ which is stable under pullback. The following conditions are equivalent:
\begin{itemize}
\item[$(1)$] The Cartesian fibration $\calO^{S}_{\calX} \rightarrow \calX$ is classified by a 
colimit-preserving functor $\calX \rightarrow \hat{\Cat}^{op}_{\infty}$. 
\item[$(2)$] The right fibration $\calO^{(S)}_{\calX} \rightarrow \calX$ is classified by a colimit-preserving functor $\calX \rightarrow \hat{\SSet}^{op}$.  
\item[$(3)$] The class $S$ is stable under $($arbitrary$)$ coproducts, and 
for every pushout diagram
$$ \xymatrix{ f \ar[r]^{\alpha} \ar[d]^{\beta} & g \ar[d]^{\beta'} \\
f' \ar[r]^{\alpha'} & g' }$$
in $\calO_{\calX}$, if $\alpha$ and $\beta$ are Cartesian transformations and
$f,f',g \in S$, then $\alpha'$ and $\beta'$ are also Cartesian transformations and $g' \in S$.
\end{itemize}
\end{lemma}

\begin{proof}
The equivalence of $(1)$ and $(2)$ follows easily from Lemma \ref{ib2}.
Let $s: \calX \rightarrow \hat{\Cat}_{\infty}^{op}$ be a functor which classifies $\calO_{\calX}$.
Then $(1)$ is equivalent to the assertion that $s$ preserves small colimits. Supposing that $(1)$ is satisfied, we deduce $(3)$ by applying Lemma \ref{ib2} in the special cases of sums and coproducts. For the converse, let us suppose that $(3)$ is satisfied. Let $\emptyset$ denote an initial object of $\calX$. Since colimits in $\calX$ are universal, Lemma \ref{sumoto} implies that $\calX^{/\emptyset}$ is equivalent to final $\infty$-category $\Delta^0$.
Since the morphism $\id_{\emptyset}$ belongs to $S$ (since $S$ is stable under {\em empty} coproducts), we conclude that $s(\emptyset)$ is a final $\infty$-category, so that
$s$ preserves initial objects. It follows from Corollary \ref{allfinn} that $s$ preserves finite coproducts. According to Proposition \ref{alllimits}, it will suffice to prove that $s$ preserves arbitrary coproducts. To handle the case of infinite coproducts we apply Lemma \ref{ib2} again: we must show that if
$\{ f_{\alpha} \}_{\alpha \in A}$ is a collection of elements of $S$ having a coproduct
$f = \coprod_{\alpha \in A} f_{\alpha}$, then $f \in S$ and each of the maps $f_{\alpha} \rightarrow f$ is a Cartesian transformation. The first condition is true by assumption; for the second we 
let $f'$ be a coproduct of the family $\{ f_{\beta} \}_{\beta \in A, \beta \neq \alpha}$, so that $f \simeq f' \amalg f_{\alpha}$ and $f' \in S$. Applying Lemma \ref{ib2} (and the fact that $s$ preserves finite coproducts), we deduce that $f_{\alpha} \rightarrow f$ is a Cartesian transformation as desired.
\end{proof}

\begin{definition}\label{localitie}\index{gen}{local!class of morphisms}
Let $\calX$ be a presentable $\infty$-category in which colimits are universal, and let $S$ be a class of morphisms in $\calX$. We will say that $S$ is {\it local} if it is stable under pullbacks
and satisfies the equivalent conditions of Lemma \ref{ib3}.
\end{definition}

\begin{theorem}\label{charleschar}
Let $\calX$ be a presentable $\infty$-category. The following conditions are equivalent:
\begin{itemize}
\item[$(1)$] Colimits in $\calX$ are universal, and for every 
pushout diagram
$$ \xymatrix{ f \ar[r]^{\alpha} \ar[d]^{\beta} & g \ar[d]^{\beta'} \\
f' \ar[r]^{\alpha'} & g' }$$
in $\calO_{\calX}$, if $\alpha$ and $\beta$ are Cartesian transformations, then
$\alpha'$ and $\beta'$ are also Cartesian transformations.

\item[$(2)$] Colimits in $\calX$ are universal, and the class of {\em all} morphisms
in $\calX$ is local.

\item[$(3)$] The Cartesian fibration $\calO_{\calX} \rightarrow \calX$ is classified by a limit-preserving functor $\calX^{op} \rightarrow \LPres$.

\item[$(4)$] Let $K$ be a small simplicial set and 
$\overline{\alpha}: \overline{p} \rightarrow \overline{q}$ a natural transformation of diagrams
$\overline{p}, \overline{q}: K^{\triangleright} \rightarrow \calX$. Suppose that
$\overline{q}$ is a colimit diagram, and that $\alpha = \overline{\alpha} | K$ is a Cartesian transformation. Then $\overline{p}$ is a colimit diagram if and only if $\overline{\alpha}$ is 
a Cartesian transformation.
\end{itemize}
\end{theorem}

\begin{proof}
The equivalences $(1) \Leftrightarrow (2) \Leftrightarrow (3)$ follow from Lemma \ref{ib3} and Proposition \ref{gentur}. The equivalence $(3) \Leftrightarrow (4)$ follows from Lemmas \ref{ib1} and \ref{ib2}. 
\end{proof}

We now have most of the tools required to establish the implication $(1) \Rightarrow (2)$ of Theorem \ref{mainchar}. In view of Theorem \ref{charleschar}, it will suffice to prove the following:

\begin{proposition}\label{lemonade}
Let $\calX$ be an $\infty$-topos. Then:
\begin{itemize}
\item[$(1)$] Colimits in $\calX$ are universal.
\item[$(2)$] For every pushout diagram $$ \xymatrix{ f \ar[r]^{\alpha} \ar[d]^{\beta} & g \ar[d]^{\beta'} \\
f' \ar[r]^{\alpha'} & g' }$$
in $\calO_{\calX}$, if $\alpha$ and $\beta$ are Cartesian transformations, then
$\alpha'$ and $\beta'$ are also Cartesian transformations.
\end{itemize}
\end{proposition}

\begin{remark}
Once we have established Theorem \ref{mainchar} in its entirety, it will follow from Theorem \ref{charleschar} that the converse of Proposition \ref{lemonade} is also valid: a presentable $\infty$-category $\calX$ is an $\infty$-topos if and only if it satisfies conditions $(1)$ and $(2)$ as above. Condition $(1)$ is equivalent to the requirement that
for every morphism $f: X \rightarrow Y$ in $\calX$, the pullback functor $f^{\ast}: \calX_{/Y} \rightarrow \calX_{/X}$ has a right adjoint (in the case where $Y$ is a final object of $\calX$, this simply amounts to the requirement that every object $Z \in \calX$ admits an exponential
$Z^X$; in other words, the requirement that $\calX$ be {\it Cartesian closed}), and condition $(2)$ involves only finite diagrams in the $\infty$-category $\calX$. One could conceivably obtain a theory of {\em elementary $\infty$-topoi} by dropping the requirement that $\calX$ be presentable (or replacing it by weaker conditions which are also finite in nature). We will not pursue this idea further.\index{gen}{$\infty$-topos!elementary}
\end{remark}

Before giving the proof of Proposition \ref{lemonade}, we need to establish a few easy lemmas.

\begin{lemma}\label{stugart}
Let 
$$ \xymatrix{ \phi \ar[r]^{p} \ar[d]^{q} & \psi \ar[d]^{q'} \\
\phi' \ar[r]^{p'} & \psi' }$$
be a coCartesian square in the category of arrows of $\sSet$. Suppose that
$p$ and $q$ are homotopy Cartesian and that $q$ is a cofibration. Then:
\begin{itemize}
\item[$(1)$] The maps $p'$ and $q'$ are homotopy Cartesian.
\item[$(2)$] Given any map of arrows $r: \psi' \rightarrow \theta$ such that
$r \circ p'$ and $r \circ q'$ are homotopy Cartesian, the map $r$ is itself homotopy Cartesian.
\end{itemize}
\end{lemma}

\begin{proof}
Let $r: \psi' \rightarrow \theta$ be as in $(2)$. We must show that $r$ is homotopy Cartesian if and only if $r \circ p'$ and $r \circ q'$ are homotopy Cartesian (taking $r = \id_{\psi'}$, we will deduce $(1)$). Without loss of generality, we may replace $\phi$, $\psi$, $\phi'$ and $\theta$ with minimal Kan fibrations. We now observe that $r$, $r \circ p'$, and $r \circ q'$ are homotopy Cartesian if and only if they are Cartesian; the desired result now follows immediately.
\end{proof}

\begin{lemma}\label{tulmand}
Let $\bfA$ be a simplicial model category containing an object $Z$ which is both fibrant and cofibrant, and let $\bfA_{/Z}$ be endowed with the induced model structure. Then
the natural map $\theta: \sNerve ( \bfA_{/Z}^{\degree} ) \rightarrow \sNerve(\bfA^{\degree})_{/Z}$
is an equivalence of $\infty$-categories.
\end{lemma}

\begin{proof}
Let $\phi: Z' \rightarrow Z$ be an object of $\sNerve(\bfA^{\degree})_{/Z}$. Then we can choose a factorization
$$ Z' \stackrel{i}{\rightarrow} Z'' \stackrel{\psi}{\rightarrow} Z$$
where $i$ is a trivial cofibration and $\psi$ is a fibration, corresponding to a fibrant-cofibrant object
of $\bfA_{/Z}$. The above diagram classifies an equivalence between $\phi$ and $\psi$
in $\sNerve(\bfA^{\degree})_{/Z}$, so that $\theta$ is essentially surjective.

Recall that for any simplicial category $\calC$ containing a pair of objects $X$ and $Y$, there is a natural isomorphism of simplicial sets
$$ \Hom^{\rght}_{\sNerve(\calC)}( X, Y ) \simeq \Sing_{Q^{\bigdot}}( \bHom_{\calC}(X,Y) ),$$ 
where $Q^{\bigdot}$ is the cosimplicial object of $\sSet$ introduced in \S \ref{twistt}. The same calculation shows that if $\phi: X \rightarrow Z$, $\psi: Y \rightarrow Z$ are two morphisms
in $\calC$, then
$$ \Hom^{\rght}_{ \sNerve(\calC)_{/Z} }(\phi, \psi) \simeq \Sing_{Q^{\bigdot}}(P),$$
where $P$ denotes the path space 
$$ \bHom_{\calC}(X,Y) \times_{ \bHom_{\calC}(X,Z)^{ \{0\} } }
\bHom_{\calC}(X,Z)^{\Delta^1} \times_{ \bHom_{\calC}(X,Z)^{ \{1\} }} \{ \phi \}.$$
If $\calC$ is fibrant, then we may identify $M$ with the homotopy fiber of the map
$$ f: \bHom_{\calC}(X,Y) \stackrel{\psi}{\rightarrow} \bHom_{\calC}(X,Z)$$
over the vertex $\phi$. Consequently, we may identify the natural map
$$ \Hom^{\rght}_{ \sNerve( \calC_{/Z} )}(\phi, \psi) \rightarrow
\Hom^{\rght}_{ \sNerve(\calC)_{/Z} }(\phi, \psi)$$
with $\Sing_{Q^{\bigdot}}( \theta)$, where $\theta$ denotes the inclusion of
the fiber of $f$ into the homotopy fiber of $f$. Consequently, to show that
$\Sing_{Q^{\bigdot}}(\theta)$ is a homotopy equivalence, it suffices to prove that
$f$ is a Kan fibration. In the special case where $\calC = \bfA^{\degree}$ and
$\psi$ is a fibration, this follows from the definition of a simplicial model category.
\end{proof}

\begin{lemma}\label{sugartime}
Let $\SSet$ denote the $\infty$-category of spaces. Then:
\begin{itemize}
\item[$(1)$] Colimits in $\SSet$ are universal.
\item[$(2)$] For every pushout diagram $$ \xymatrix{ f \ar[r]^{\alpha} \ar[d]^{\beta} & g \ar[d]^{\beta'} \\
f' \ar[r]^{\alpha'} & g' }$$
in $\calO_{\SSet}$, if $\alpha$ and $\beta$ are Cartesian transformations, then
$\alpha'$ and $\beta'$ are also Cartesian transformations.
\end{itemize}
\end{lemma}

\begin{proof}
We first prove $(1)$.
Let $f: X \rightarrow Y$ be a morphism in $\SSet$. Without loss of generality, we may
suppose that $f$ is a Kan fibration. We wish to show that the projection
$$ \SSet_{/f} \rightarrow \SSet_{/Y}$$ has a right adjoint which preserves colimits.
We have a commutative diagram of $\infty$-categories
$$ \xymatrix{ \sNerve( (\sSet)^{\degree}_{/X} ) \ar[r]^{F} \ar[d]^{\phi} & \Nerve((\sSet)^{\degree}_{/Y}) \ar[d]^{\psi} \\
\SSet_{/f} \ar[r] \ar[d]^{\phi'} & \SSet_{/Y} \\
\SSet_{/X} & } $$
Lemma \ref{tulmand} asserts that $\psi$ and $\phi' \circ \phi$ are categorical equivalences, and $\phi'$ is a trivial fibration. It follows that $\phi$ is also a categorical equivalence. Consequently,
it will suffice to show that the functor $F$ has a right adjoint $G$ which preserves colimits.

We observe that $F$ is obtained by restricting the simplicial nerve of the functor $f_{!}: (\sSet)_{/X} \rightarrow (\sSet)_{/Y}$, given by composition with $f$. The functor $f_{!}$ is a left Quillen functor: it has a right adjoint $f^{\ast}$, given by the formula $f^{\ast}(Y') = Y' \times_{Y} X.$
According to Proposition \ref{quiladj}, $F$ admits a right adjoint $G$, which is given by
the restricting the simplicial nerve of the functor $f^{\ast}$. To prove that $G$ preserves colimits, it will suffice to show that $G$ itself admits a right adjoint. Using Proposition \ref{quiladj} again, we are reduced to proving that $f^{\ast}$ is a {\em left} Quillen functor. We observe that
$f^{\ast}$ admits a right adjoint $f_{\ast}$, given by the formula
$f_{\ast}(X') = \bHom_{Y}(X,X')$. It is clear that $f^{\ast}$ preserves cofibrations; it also preserves weak equivalences, since $f$ is a fibration and $\sSet$ is a {\em right proper} model category (with its usual model structure). 

To prove $(2)$, we first apply Proposition \ref{gumby444} to reduce to the case where
the pushout diagram in question arises from a strictly commutative square
$$ \xymatrix{ f \ar[r]^{\alpha} \ar[d]^{\beta} & g \ar[d]^{\beta'} \\
f' \ar[r]^{\alpha'} & g' }$$
of morphisms in the category $\Kan$. We now complete the proof by applying 
Lemma \ref{stugart} and Theorem \ref{colimcomparee}. 
\end{proof}

\begin{lemma}\label{tryme}
Let $\calX$ be a presentable $\infty$-category, and let $L: \calX \rightarrow \calY$ be an accessible left exact localization. If colimits in $\calX$ are universal, then colimits in $\calY$ are universal.
\end{lemma}

\begin{proof}
We will use characterization $(5)$ of Lemma \ref{ib0}. Let $G$ be a right adjoint to $L$, and let
$\alpha: \overline{p} \rightarrow \overline{q}$ be a Cartesian transformation of diagrams $K^{\triangleright} \rightarrow \calY$. Suppose that $\overline{q}$ is a colimit of $q = \overline{q}|K$.
Choose a colimit $\overline{q}'$ of $G \circ q$, so that there exists a morphism
$\overline{q}' \rightarrow G \circ \overline{q}$ in $\calX_{G \circ q/}$ which determines a natural transformation $\beta: \overline{q}' \rightarrow \overline{q}$ in $\Fun(K^{\triangleright}, \calX)$. 
Form a pullback diagram
$$ \xymatrix{ \overline{p}' \ar[r]^{\alpha'} \ar[d] & \overline{q}' \ar[d]^{\beta} \\
G \circ \overline{p} \ar[r]^{G \circ \alpha} & G \circ \overline{q}. }$$
in $\calX^{K^{\triangleright}}$. Since $G$ is left exact, $G \circ \alpha$ is a Cartesian transformation. It follows that $\alpha'$, being a pullback of $G \circ \alpha$, is also a Cartesian transformation. Since colimits in $\calX$ are universal, we conclude that $\overline{p}'$ is a colimit diagram.
Since $L$ is left exact, we obtain a pullback diagram
$$ \xymatrix{ L \circ \overline{p}' \ar[r] \ar[d] & L \circ \overline{q}' \ar[d]^{L \circ \beta} \\
L \circ G \circ \overline{p} \ar[r] & L \circ G \circ \overline{q}. }$$
Since $L$ preserves colimits, $L \circ \overline{q}'$ and $L \circ \overline{p'}$ are colimit diagrams. 
The diagram $L \circ G \circ \overline{q}$ is equivalent to $q$, and therefore also a colimit diagram. We deduce that $L \circ \beta$ is an equivalence. Since the diagram is a pullback, the left vertical arrow is an equivalence as well, so that $L \circ G \circ \overline{p}$ is a colimit diagram. We finally conclude that $\overline{p}$ is a colimit diagram, as desired.
\end{proof}

We are now ready to give the proof of Proposition \ref{lemonade}.

\begin{proof}[Proof of Proposition \ref{lemonade}]
Let us say that a presentable $\infty$-category $\calX$ is {\em good} if it satisfies conditions $(1)$ and $(2)$.
Lemma \ref{sugartime} asserts that $\SSet$ is good. Using Proposition \ref{limiteval}, it is easy to see that if $\calX$ is good then so is $\Fun(K,\calX)$, for every small simplicial set $K$. It follows that
every $\infty$-category $\calP(\calC)$ of presheaves is good. To complete the proof, it will suffice to show that if $\calX$ is good and $L: \calX \rightarrow \calY$ is an accessible left exact localization functor, then $\calY$ is good. Lemma \ref{tryme} shows that colimits in $\calY$ are universal. 
Consider a diagram $\sigma: \Lambda^2_0 \rightarrow \calO_{\calY}$, denoted by 
$$ g \stackrel{\alpha}{\leftarrow} f \stackrel{\beta}{\rightarrow} h$$
where $\alpha$ and $\beta$ are Cartesian transformations. We wish to show that if
$\overline{\sigma}$ is a colimit of $\sigma$ in $\calO_{\calY}$, then $\overline{\sigma}$ carries
each edge to a Cartesian transformation. Without loss of generality, we may suppose that
$\sigma = L \circ \sigma'$ for some $\sigma': \Lambda^2_0 \rightarrow \calO_{\calX}$ which
is equivalent to $G \circ \sigma$. Since $G$ is left exact, $G(\alpha)$ and $G(\beta)$ are
Cartesian transformations. Because $\calX$
satisfies $(2)$, there exists a colimit $\overline{\sigma}'$ of $\sigma'$ which carries each edge to a Cartesian transformation. Then $L \circ \overline{\sigma}'$ is a colimit of $\sigma$. Since $L$ is left exact, $L \circ \overline{\sigma}'$ carries each edge to a Cartesian transformation in 
$\calO_{\calY}$. 
\end{proof}

Our final objective in this section is to prove the implication $(2) \Rightarrow (3)$ of Theorem \ref{mainchar} (Proposition \ref{lemonade2} below). 

\begin{lemma}\label{aclock}
Let $\calX$ be an $\infty$-category and $U^{+}_{\bigdot}: \Nerve(\cDelta_{+})^{op} \rightarrow \calX$
an augmented simplicial object of $\calX$. Let $\cDelta_{\infty}$ denote the category
whose objects are finite, linearly ordered sets $J$, where $\Hom_{\cDelta_{\infty}}(J,J')$ is the collection of all order-preserving maps $J \cup \{\infty\} \rightarrow J' \cup \{\infty\}$ which carry $\infty$ to $\infty$ $($here $\infty$ is regarded as a maximal element of $J \cup \{ \infty \}$ and
$J' \cup \{\infty\}${}$)$. Suppose that $U^{+}_{\bigdot}$ extends to a functor
$F: \Nerve(\cDelta_{\infty})^{op} \rightarrow \calX$. Then $U^{+}_{\bigdot}$ is a colimit diagram in $\calX$.
\end{lemma}

\begin{proof}
Let $\overline{\calC}$ denote the category whose objects are triples $(J, J_{+})$, where $J$ is a finite, linearly ordered set, and $J_{+}$ is an upward-closed subset of $J$. We define
$\Hom_{\overline{\calC}}( (J, J_{+}), (J', J'_{+})$ to be the set of all order-preserving maps
from $J$ into $J'$ that carry $J_{+}$ into $J'_{+}$. Observe that we have a functor
$\overline{\calC} \rightarrow \cDelta_{\infty}$, given by
$$ (J, J_{+}) \mapsto J - J_{+}. $$
Let $F'$ denote the composite functor
$$ \Nerve(\overline{\calC})^{op} \rightarrow \Nerve(\cDelta_{\infty})^{op} \rightarrow \calX.$$

Let $\calC$ be the full subcategory of $\overline{\calC}$ spanned by those pairs $(J, J_+)$ where
$J \neq \emptyset$. Let $\overline{\calC}^{0}$ denote the full subcategory spanned by those pairs $(J, J_{+})$ where $J_{+} = \emptyset$, and let $\calC^{0} = \overline{\calC}^0 \cap \calC$.
We observe that $\overline{\calC}^{0}$ can be identified with $\cDelta_{+}$ and
that $\calC^{0}$ can be identified with $\cDelta$, in such a way that $U^{+}_{\bigdot}$ is identified with $F' | \Nerve (\overline{\calC}^{0})^{op}$.

Our first claim is that the inclusion $\Nerve(\calC^{0})^{op} \subseteq \Nerve(\calC)^{op}$ is cofinal. According to Theorem \ref{hollowtt}, it will suffice
to show that for every object $X = (J, J_{+})$ of $\calC$, the category
$\calC^{0}_{/X}$ has a contractible nerve. This is clear, since the relevant category has a final object: namely, the map $(J, \emptyset) \rightarrow (J, J_{+})$. As a consequence, we conclude that 
$U_{\bigdot}^{+}$ is a colimit diagram if and only if $F'$ is a colimit diagram.

We now define $\calC^{1}$ to be the full subcategory of $\calC$ spanned by those pairs
$(J,J_{+})$ such that $J_{+}$ is nonempty. We claim that $F' | \Nerve(\calC)^{op}$ is a left Kan extension of $F' | \Nerve (\calC^{1})^{op}$. To prove this, we must show that for every
$(J, \emptyset) \in \calC^{0}$, the induced map
$$ (\Nerve (\calC^{1}_{(J, \emptyset)/} )^{op} )^{\triangleright} \rightarrow \Nerve(\calC)^{op}
\rightarrow \calX$$
is a colimit diagram. Let $\calD$ denote the full subcategory of
$\calC^{1}_{(J, \emptyset)/}$ spanned by those morphisms $(J, \emptyset) \rightarrow
(J' , J'_{+} )$ which induce isomorphisms $J \simeq J' - J'_{+}$. We claim that the inclusion
$\Nerve(\calD)^{op} \subseteq \Nerve ( \calC^{1}_{(J, \emptyset)/} )^{op}$ is cofinal.
To prove this, we once again invoke Theorem \ref{hollowtt}, to reduce to the following assertion:
for every morphism $\phi: (J, \emptyset) \rightarrow (J'', J''_{+})$, if $J''_{+} \neq \emptyset$, then the category $\calD_{/\phi}$ of all factorizations
$$ (J, \emptyset) \rightarrow (J', J'_{+}) \rightarrow (J'', J''_{+}) $$ 
such that $J'_{+} \neq \emptyset$ and $J \simeq J' - J'_{+}$, has weakly contractible nerve.
This is clear, since $\calD_{/\phi}$ has a final object $(J \amalg J'''_{+}, J'''_{+})$, where
$J'''_{+} = \{ j \in J''_{+}: (\forall i \in J) [ j \geq \phi(i) ] \}$. Consequently, it will suffice to prove that the induced functor
$$ \Nerve(\calD^{op})^{\triangleright} \rightarrow \calX$$
is a colimit diagram. This diagram can be identified with the constant diagram
$$ \Nerve(\cDelta_{+})^{op} \rightarrow \calX$$ taking the value $U_{\bigdot}(J)$, and
is a colimit diagram because the category $\cDelta$ has weakly contractible nerve (Corollary \ref{silt}).

We now apply Lemma \ref{kan0}, which asserts that $F'$ is a colimit diagram if and only if
$F' | (\Nerve (\calC^{1})^{op})^{\triangleright}$ is a colimit diagram.
Let $\calC^{2} \subseteq \calC^{1}$ be the full subcategory spanned by those objects
$(J, J_{+})$ such that $J = J_{+}$. We claim that the inclusion
$\Nerve(\calC^{2})^{op} \subseteq \Nerve(\calC^{1})^{op}$ is cofinal. According to Theorem
\ref{hollowtt}, it will suffice to show that, for every object $(J, J_+) \in \calC^{1}$, the
category $\calC^{2}_{/(J, J_{+})}$ has weakly contractible nerve. This is clear,
since the map $(J_{+}, J_{+}) \rightarrow (J, J_{+})$ is a final object of 
the category $\calC^{(2)}_{/(J,J_{+})}$. Consequently, to prove that
$F' | (\Nerve (\calC^{1})^{op})^{\triangleright}$ is a colimit diagram, it will suffice to prove that
$F' | ( \Nerve (\calC^{2})^{op})^{\triangleright}$ is a colimit diagram. But this diagram can be identified with the constant map $\Nerve(\cDelta_{+})^{op} \rightarrow \calX$ taking the value
$U_{\bigdot}( \Delta^{-1})$, which is a colimit diagram because the simplicial set
$\Nerve(\cDelta)^{op}$ is weakly contractible (Corollary \ref{silt}). 
\end{proof}

\begin{lemma}\label{bclock}
Let $\calX$ be an $\infty$-category, and let $U_{\bigdot}: \Nerve(\cDelta)^{op} \rightarrow \calX$ be a simplicial object of $\calX$. Let $U'_{\bigdot}$ be the augmented simplicial object
given by composing $U_{\bigdot}$ with the functor
$$ \cDelta_{+} \rightarrow \cDelta$$
$$ J \rightarrow J \amalg \{ \infty \}.$$
Then:
\begin{itemize}
\item[$(1)$] The augmented simplicial object $U'_{\bigdot}$ is a colimit diagram.
\item[$(2)$] If $U_{\bigdot}$ is a groupoid object of $\calX$, then the evident natural transformation of simplicial objects $\alpha: U'_{\bigdot} | \Nerve(\cDelta)^{op} \rightarrow U_{\bigdot}$ is Cartesian.
\end{itemize}
\end{lemma}

\begin{proof}
Assertion $(1)$ follows immediately from Lemma \ref{aclock}. To prove $(2)$, let us consider
the collection $S$ of all morphisms $f: J \rightarrow J'$ in $\cDelta$ such that
$\alpha(f)$ is a pullback square
$$ \xymatrix{ U'_{\bigdot}(J') \ar[r] \ar[d] & U_{\bigdot}(J') \ar[d] \\
U'_{\bigdot}(J) \ar[r] & U_{\bigdot}(J) }$$
in $\calX$. We wish to prove that every morphism of $\cDelta$ belongs to $S$. Using Lemma \ref{transplantt}, we deduce that if $f' \in S$, then $f \in S \Leftrightarrow (f \circ f' \in S)$. Consequently, it will suffice to prove that every inclusion $\{j\} \subseteq J$ belongs to $S$.
Unwinding the definition, this amounts the the requirement that the diagram
$$ \xymatrix{ U_{\bigdot}( J \cup \{\infty\} ) \ar[d]  \ar[r] & U_{\bigdot}( J) \ar[d] \\
U_{\bigdot} ( \{j, \infty\} ) \ar[r] & U_{\bigdot}( \{j\} )}$$
is Cartesian, which follows immediately from Criterion $(4'')$ of Proposition \ref{grpobjdef}.
\end{proof}

\begin{remark}
Assertion $(2)$ of Lemma \ref{bclock} has a converse: if $\alpha$ is a Cartesian transformation, then $U_{\bigdot}$ is a groupoid object of $\calX$. This can be deduced easily by examining the proof of Proposition \ref{grpobjdef}, but we will not have need of it.
\end{remark}

\begin{proposition}\label{lemonade2}
Let $\calX$ be an $\infty$-category satisfying the equivalent conditions of Theorem \ref{charleschar}. Then $\calX$ satisfies the $\infty$-categorical Giraud axioms:
\begin{itemize}
\item[$(i)$] The $\infty$-category $\calX$ is presentable.
\item[$(ii)$] Colimits in $\calX$ are universal.
\item[$(iii)$] Coproducts in $\calX$ are disjoint.
\item[$(iv)$] Every groupoid object of $\calX$ is effective.
\end{itemize}
\end{proposition}

\begin{proof}
Axioms $(i)$ and $(ii)$ are obvious. To prove $(iii)$, let us consider an arbitrary pair of objects
$X, Y \in \calX$, and let $\emptyset$ denote an initial object of $\calX$. Let $f: \emptyset \rightarrow X$ be a morphism (unique up to homotopy, since $\emptyset$ is initial). We observe that
$\id_{\emptyset}$ is an initial object of $\calO_{\calX}$. Form a pushout diagram
$$ \xymatrix{ \id_{\emptyset} \ar[r]^{\alpha} \ar[d]^{\beta} & \id_{Y} \ar[d]^{\beta'}  \\
f \ar[r]^{\alpha'} & g }$$
in $\calO_{\calX}$. It is clear that $\alpha$ is a Cartesian transformation, and Lemma \ref{sumoto} implies that $\beta$ is Cartesian as well. Invoking condition $(2)$ of Theorem \ref{charleschar}, we deduce that $\alpha'$ is
a Cartesian transformation. But $\alpha'$ can be identified with a pushout diagram
$$ \xymatrix{ \emptyset \ar[r] \ar[d] & Y \ar[d] \\
X \ar[r] & X \amalg Y. }$$

It remains to prove that every groupoid object in $\calX$ is effective. Let $U_{\bigdot}$
be a groupoid object of $\calX$, and let $\overline{U}_{\bigdot}  : \Nerve(\cDelta_{+})^{op} \rightarrow \calX$ be a colimit of $U_{\bigdot}$. Let $U'_{\bigdot}: \Nerve(\cDelta_{+})^{op} \rightarrow \calX$
be the result of composing $\overline{U}_{\bigdot}$ with the ``shift'' functor
$$ \cDelta_{+} \rightarrow \cDelta_{+}$$
$$ J \mapsto J \amalg \{ \infty \}.$$
(In other words, $U'_{\bigdot}$ is the shifted simplicial object given by
$U'_{n} = U_{n+1}$.)
Lemma \ref{bclock} implies that $U'_{\bigdot}$ is a colimit diagram in $\calX$.
We have a transformation $\overline{\alpha}: U'_{\bigdot} \rightarrow \overline{U}_{\bigdot}$.
Since $U_{\bigdot}$ is a groupoid, $\alpha = \overline{\alpha} | \Nerve(\cDelta)^{op}$
is a Cartesian transformation (Lemma \ref{bclock} again). Applying $(4)$, we deduce
that $\overline{\alpha}$ is a Cartesian transformation. In particular, we conclude that
$$ \xymatrix{ U'_0 \ar[r] \ar[d] & U'_{-1} \ar[d] \\
\overline{U}_0 \ar[r] & \overline{U}_{-1} }$$
is a pullback diagram in $\calX$. But this diagram can be identified with
$$ \xymatrix{ U_1 \ar[r] \ar[d] & U_0 \ar[d] \\
U_0 \ar[r] & \overline{U}_{-1}, }$$
so that $U_{\bigdot}$ is effective by Proposition \ref{strump}.
\end{proof}

\begin{corollary}\label{alleff}
Every groupoid object of $\SSet$ is effective.
\end{corollary}

\subsection{Free Groupoids}\label{freegroup}

Let $\calX$ be an $\infty$-category which satisfies the $\infty$-categorical Giraud axioms
$(i)-(iv)$ of Theorem \ref{mainchar}. We wish to prove that $\calX$ is an $\infty$-topos. It is clear that any proof will need to make use of the full strength of axioms $(i)$ through $(iv)$; in particular, we will need to apply $(iv)$ to a class of groupoid objects of $\calX$ which are not obviously effective.
The purpose of this section is to describe a construction which will yields nontrivial examples of groupoid objects, and to deduce a consequence (Proposition \ref{storytell})
which we will use in the proof of Theorem \ref{mainchar}.

\begin{definition}\label{swarpy}\index{gen}{left exact!at an object $Z$}
Let $f: \calX \rightarrow \calY$ be a functor between $\infty$-categories which admit finite limits. Let
$Z$ be an object of $\calX$. We will say that $f$ is {\it left exact at $Z$} if, for every pullback square
$$ \xymatrix{ W \ar[r] \ar[d] & Y \ar[d] \\
X \ar[r] & Z }$$
in $\calX$, the induced square
$$ \xymatrix{ f(W) \ar[r] \ar[d] & f(Y) \ar[d] \\
f(X) \ar[r] & f(Z) }$$
is a pullback in $\calY$.
\end{definition}

We can now state the main result of this section:

\begin{proposition}\label{storytell}
Let $\calX$ and $\calY$ be presentable $\infty$-categories, and let $f: \calX \rightarrow \calY$ be a functor which preserves small colimits. Suppose that every groupoid object in either $\calX$ or $\calY$ is effective.
Let
$$ \xymatrix{ U_{1} \ar@<.5ex>[r] \ar@<-.5ex>[r] & U_0 \ar[r]^{s} & U_{-1}. } $$ be a coequalizer diagram in $\calX$, and let
$$ \xymatrix{ X \ar[r] \ar[d] & U_0 \ar[d]^{s} \\
U_0 \ar[r]^{s} & U_{-1} }$$ be a pullback diagram in $\calX$. Suppose that $f$ is left exact
at $U_0$. Then the associated diagram 
$$ \xymatrix{ f(X) \ar[r] \ar[d] & f(U_0) \ar[d]^{s} \\
f(U_0) \ar[r]^{s} & f(U_{-1}) }$$
is a pullback square in $\calY$.
\end{proposition}

Before giving the proof, we must establish some preliminary results.

\begin{lemma}\label{drumb}
Let $\calX$ and $\calY$ be $\infty$-categories which admit finite limits, let
$f: \calX \rightarrow \calY$ be a functor, and let $U_{\bigdot}$ be a groupoid object
of $\calX$. Suppose that $f$ is left exact at $U_0$. Then $f \circ U_{\bigdot}$ is a groupoid
object of $\calY$. 
\end{lemma}

\begin{proof}
This follows immediately from characterization $(4'')$ given in Proposition \ref{grpobjdef}.
\end{proof}

Let $\calX$ be a presentable $\infty$-category. We define a {\it simplicial resolution}\index{gen}{resolution!simplicial} in $\calX$ to be an augmented simplicial object $U_{\bigdot}^{+} : \Nerve(\cDelta_{+})^{op} \rightarrow \calX$
which is a colimit of the underlying simplicial object $U_{\bigdot} = U_{\bigdot}^{+} | \Nerve(\cDelta)^{op}$. We let $\Res(\calX)$\index{not}{ResX@$\Res(\calX)$} denote the full subcategory of $\calX_{\Delta_{+}}$ spanned by the simplicial resolutions. Note that since every simplicial object of $\calX$ has a colimit, 
the restriction functor $\Res(\calX) \rightarrow \calX_{\Delta}$ is a trivial fibration, and therefore an equivalence of $\infty$-categories. We will say that a simplicial resolution $U_{\bigdot}^{+}$
is a {\it groupoid resolution}\index{gen}{resolution!groupoid} if the underlying simplicial object
$U_{\bigdot}$ is a groupoid object of $\calX$.

We will say that a map $f: U_{\bigdot}^{+} \rightarrow V_{\bigdot}^{+}$ of simplicial resolutions {\it exhibits $V_{\bigdot}^{+}$ as the groupoid resolution generated by $U_{\bigdot}^{+}$} if $V_{\bigdot}^{+}$ is a groupoid resolution and the induced map $$ \bHom_{\Res(\calX)}( V_{\bigdot}^{+}, W_{\bigdot}^{+} )
\rightarrow \bHom_{\Res(\calX)}( U_{\bigdot}^{+}, W_{\bigdot}^{+} )$$
is a homotopy equivalence for every groupoid resolution $W_{\bigdot}^{+} \in \Res(\calX)$.

\begin{remark}\label{swather}
Let $\calX$ be a presentable $\infty$-category. Then for every simplicial resolution
$U_{\bigdot}^{+}$ in $\calX$, there is a map $f: U_{\bigdot}^{+} \rightarrow V_{\bigdot}^{+}$ which exhibits $V_{\bigdot}^{+}$ as the groupoid resolution generated by $U_{\bigdot}^{+}$. In view of
the equivalence $\Res(\calX) \rightarrow \calX_{\Delta}$, this is equivalent to the assertion that
$\Grp(\calX)$ is a localization of $\calX_{\Delta}$. This follows from Proposition \ref{grpobjdef} together with Lemmas \ref{stur3} and \ref{stur1}.
\end{remark}

\begin{lemma}\label{step1}
Let $\calX$ be a presentable $\infty$-category and let $f: U_{\bigdot}^{+} \rightarrow V_{\bigdot}^{+}$ be a map of simplicial resolutions which exhibits $V_{\bigdot}^{+}$ as the groupoid resolution generated by $U_{\bigdot}^{+}$. Let $W_{\bigdot}^{+}$ be an augmented simplicial object of $\calX$ such that the underlying simplicial object $W_{\bigdot} \in \calX_{\Delta}$
is a groupoid. Composition with $f$ induces a homotopy equivalence
$$ \bHom_{\calX_{\Delta_{+}}}( V_{\bigdot}^{+}, W_{\bigdot}^{+} )
\rightarrow \bHom_{\calX_{\Delta_{+}}}( U_{\bigdot}^{+}, W_{\bigdot}^{+} ).$$
\end{lemma}

\begin{proof}
Let $|W_{\bigdot}|$ be a colimit of $W_{\bigdot}$. Then we have a commutative diagram
$$ \xymatrix{ \bHom_{\calX_{\Delta_{+}}}( V_{\bigdot}^{+}, |W_{\bigdot}| )
\ar[r] \ar[d] &  \bHom_{\calX_{\Delta_{+}}}( U_{\bigdot}^{+}, |W_{\bigdot}| ) \ar[d] \\
\bHom_{\calX_{\Delta_{+}}}( V_{\bigdot}^{+}, W_{\bigdot}^{+} )
\ar[r] & \bHom_{\calX_{\Delta_{+}}}( U_{\bigdot}^{+}, W_{\bigdot}^{+} )} $$
where the vertical maps are homotopy equivalences (since $U_{\bigdot}^{+}$ and $V_{\bigdot}^{+}$ are resolutions) and the upper horizontal map is a homotopy equivalence (since $|W_{\bigdot}|$ is a groupoid resolution).
\end{proof}

\begin{lemma}\label{step9}
Let $\calX$ be a presentable $\infty$-category. Suppose that $f: U_{\bigdot}^{+} \rightarrow V_{\bigdot}^{+}$ be a map in $\Res(\calX)$ which exhibits $V_{\bigdot}^{+}$ as the groupoid resolution generated by $U_{\bigdot}^{+}$. Then $f$ induces equivalences $U_{-1} \rightarrow V_{-1}$ and
$U_{0} \rightarrow V_{0}$.
\end{lemma}

\begin{proof}
Let $\cDelta^{ \leq 0}_{+}$ be the full subcategory of $\cDelta_{+}$ spanned by the objects $\Delta^{-1}$ and $\Delta^0$. Let $j: \cDelta^{\leq 0}_{+} \rightarrow \cDelta_{+}$ denote the inclusion functor, let
$j^{\ast}: \calX_{\Delta_{+}} \rightarrow \calO_{\calX}$ be the associated restriction functor. We wish to show that $j^{\ast}(f)$ is an equivalence. Equivalently, we show that for every $W \in \calO_{\calX}$, 
composition with $j^{\ast}(f)$ induces a homotopy equivalence
$$ \bHom_{\calO_{\calX}}( j^{\ast} V_{\bigdot}^{+}, W )
\rightarrow \bHom_{\calO_{\calX}}( j^{\ast} U_{\bigdot}^{+}, W ).$$

Let $j_{\ast}$ be a right adjoint to $j^{\ast}$ (a right Kan extension functor). It will suffice to prove that composition with $f$ induces a homotopy equivalence
$$ \bHom_{\calX_{\Delta_{+}}}( V_{\bigdot}^{+}, j_{\ast} W )
\rightarrow \bHom_{\calX_{\Delta_{+}}}( U_{\bigdot}^{+}, j_{\ast} W ).$$
The augmented simplicial object $j_{\ast} W$ is a \Cech nerve, so that the underlying simplicial object of $j_{\ast} W$ is a groupoid by Proposition \ref{strump}. We now conclude by applying Lemma \ref{step1}.
\end{proof}

Let $\calI$ denote the subcategory of $\cDelta_{+}$ spanned by the objects
$\emptyset$, $[0]$, and $[1]$, where the morphisms are given by {\em injective} maps of linearly ordered sets.
This category may be depicted as follows: 
$$ \xymatrix{ \emptyset \ar[r] & [0] \ar@<.5ex>[r] \ar@<-.5ex>[r] & [1]} $$
We let $\calI_0$ denote the full subcategory of $\calI$ spanned by the objects $[0]$ and $[1]$. We will say that a diagram $\Nerve(\calI)^{op} \rightarrow \calX$ is a {\it coequalizer diagram} if it is a colimit of its restriction to $\Nerve(\calI_0)^{op} \rightarrow \calX$.

Let $i$ denote the inclusion $\calI \subseteq \cDelta_{+}$, and let $i^{\ast}$ denote the
restriction functor $\calX_{\Delta_{+}} \rightarrow \Fun( \Nerve(\calI)^{op}, \calX)$. If $\calX$ is a presentable $\infty$-category, then $i^{\ast}$ has a left adjoint $i_{!}$ (a left Kan extension).

\begin{lemma}\label{step7}
Let $\calX$ be a presentable $\infty$-category. The left Kan extension
$i_{!}: \Fun( \Nerve(\calI)^{op}, \calX) \rightarrow \calX_{\Delta_{+}}$ carries
coequalizer diagrams to simplicial resolutions.
\end{lemma}

\begin{proof}
We have a commutative diagram of inclusions of subcategories
$$ \xymatrix{ \calI_0 \ar[r]^{j'} \ar[d]^{i'} & \calI \ar[d]^{i} \\
\cDelta \ar[r]^{j} & \cDelta_{+} }$$
which gives rise to a homotopy commutative diagram of $\infty$-categories
$$ \xymatrix{ \Fun( \Nerve(\calI_0)^{op}, \calX) \ar[r]^{j'_{!}} \ar[d]^{i'_{!}} & 
\Fun( \Nerve(\calI)^{op}, \calX) \ar[d]^{i_{!}} \\
\calX_{\Delta} \ar[r]^{j_{!}} & \calX_{\Delta_{+}} }$$
in which the morphisms are given by left Kan extensions. An object $U \in \Fun( \Nerve(\calI)^{op} \calX)$ is a coequalizer diagram if and only if it lies in the essential image of $j'_{!}$. In
this case, $i_{!} U$ lies in the essential image of $i_{!} \circ j'_{!} \simeq j_{!} \circ i'_{!}$, 
which is contained in the essential image of $j_{!}$: namely, the resolutions. 
\end{proof}

\begin{lemma}\label{step8}
Let $\calX$ be a presentable $\infty$-category and suppose given a diagram
$U: \Nerve(\calI)^{op} \rightarrow \calC$, which we may depict as
$$ \xymatrix{ U_{1} \ar@<.5ex>[r] \ar@<-.5ex>[r] & U_0 \ar[r] & U_{-1}. } $$
Let $V_{\bigdot} = i_{!} U \in \calX_{\Delta_{+}}$ be a left Kan extension of $U$ along
$i: \Nerve(\calJ)^{op} \rightarrow \cDelta_{+}^{op}$. 
Then the augmentation maps $V_0 \rightarrow V_{-1}$ and $U_0 \rightarrow U_{-1}$
are equivalent in the $\infty$-category $\calO_{\calX}$.
\end{lemma}

\begin{proof}
This follows from Proposition \ref{timeless}, since $\Hom_{\calI}(\Delta^i, \bigdot)
\simeq \Hom_{\cDelta_{+}}( \Delta^i, \bigdot )$ for $i \leq 0$.
\end{proof}

\begin{proof}[Proof of Proposition \ref{storytell}]
Let $U: \Nerve(\calI)^{op} \rightarrow \calC$ be a coequalizer diagram in $\calX$, which we denote by
$$ \xymatrix{ U_{1} \ar@<.5ex>[r] \ar@<-.5ex>[r] & U_0 \ar[r]^{s} & U_{-1}, } $$
and form a pullback square
$$ \xymatrix{ X \ar[r] \ar[d] & U_0 \ar[d]^{s} \\
U_0 \ar[r]^{s} & U_{-1}. }$$ 

Let $V_{\bigdot} = i_{!} U \in \calX_{\Delta_{+}}$ be a left Kan extension of $U$. According to Lemma \ref{step7}, $V_{\bigdot}$ is a simplicial resolution. We may therefore choose a map
$V_{\bigdot} \rightarrow W_{\bigdot}$ which exhibits $W_{\bigdot}$ as the groupoid resolution generated by $V_{\bigdot}$ (Remark \ref{swather}). Since every groupoid object in $\calX$ is effective, $W_{\bigdot}$ is a \Cech nerve. It follows from the characterization given in Proposition \ref{strump} that there is a pullback diagram
$$ \xymatrix{ W_{1} \ar[r] \ar[d] & W_{0} \ar[d] \\
W_0 \ar[r] & W_{-1} }$$
in $\calX$. Using Lemma \ref{step8} and Lemma \ref{step9}, we see that this diagram is
equivalent to the pullback diagram $$ \xymatrix{ X \ar[r] \ar[d] & U_0 \ar[d]^{s} \\
U_0 \ar[r]^{s} & U_{-1}. }$$ It therefore suffices to prove that the induced diagram
$$ \xymatrix{ f(W_{1}) \ar[r] \ar[d] & f(W_{0}) \ar[d] \\
f(W_0) \ar[r] & f(W_{-1}) }$$
is a pullback. We make a slightly stronger claim: the augmented simplicial object
$f \circ W_{\bigdot}$ is a \Cech nerve.
Since every groupoid object in $\calY$ is effective, it will suffice to prove that $f \circ W_{\bigdot}$ is a groupoid resolution. Since $f$ preserves colimits, it is clear that $f \circ W_{\bigdot}$ is a simplicial resolution. It follows from Lemma \ref{drumb} that the underlying simplicial object of $f \circ W_{\bigdot}$ is a groupoid.
\end{proof}

\subsection{Giraud's Theorem for $\infty$-Topoi}\label{proofgiraud}

In this section, we will complete the proof of Theorem \ref{mainchar} by showing that every $\infty$-category $\calX$ which satisfies the $\infty$-categorical Giraud axioms $(i)$ through $(iv)$ arises as a left exact localization of an $\infty$-category of presheaves. Our strategy is simple: we choose a small category $\calC$ equipped with a functor $f: \calC \rightarrow \calX$. According to Theorem \ref{charpresheaf}, we obtain a colimit-preserving functor $F: \calP(\calC) \rightarrow \calX$ which extends $f$, up to homotopy. We will apply Proposition \ref{storytell} to show that, under suitable hypotheses, $F$ is a left exact localization functor (Proposition \ref{natash}). 

\begin{lemma}\label{sumdescription}
Let $\calX$ be a presentable $\infty$-category in which colimits are universal and coproducts are disjoint.

Let $\{ \phi_i: Z_i \rightarrow Z\}_{i \in I}$ be a family of morphisms in $\calX$ which exhibit $Z$ as a coproduct of the family of objects $\{Z_i\}_{i \in I}$. Let
$$ \xymatrix{ W \ar[r]^{\alpha} \ar[d] & Z_i \ar[d]^{\phi_i} \\
Z_j \ar[r]^{\phi_j} & Z}$$ be a square diagram in $\calX$. Then:
\begin{itemize}
\item[$(1)$] If $i \neq j$, then the diagram is a pullback square if and only if $W$ is an initial object of $\calX$.
\item[$(2)$] If $i = j$, then the diagram is a pullback square if and only if $\alpha$
is an equivalence.
\end{itemize}
\end{lemma}

\begin{proof}
Let $Z_{i}^{\vee}$ be a coproduct for the objects $\{ Z_{k} \}_{k \in I, k \neq i}$, and
let $\psi: Z_{i}^{\vee} \rightarrow Z$ be a morphism such that each of the compositions
$$ Z_{k} \rightarrow Z_{i}^{\vee} \stackrel{\psi}{\rightarrow} Z$$
is equivalent to $Z$. Then there is a pushout square
$$ \xymatrix{ \emptyset \ar[r]^{\beta} \ar[d] & Z_i \ar[d]^{\phi_i} \\
Z_i^{\vee} \ar[r]^{\psi} & Z }$$
where $\emptyset$ denotes an initial object of $\calX$.
Since coproducts in $\calX$ are disjoint, this pushout square is also a pullback.

Let $\phi_{i}^{\ast}: \calX^{/Z} \rightarrow \calX^{/Z_i}$ denote a pullback functor. The above
argument shows that $\phi_{i}^{\ast}(\psi)$ is an initial object of $\calX^{/Z_i}$.
If $j \neq i$, then there is a map of arrows $\phi_j \rightarrow \psi$ in $\calX^{/Z}$, and therefore a map $\phi_i^{\ast}(\phi_j) \rightarrow \phi_i^{\ast}(\psi)$ in $\calX^{/Z_i}$. Consequently,
if $\alpha \simeq \phi_i^{\ast}(\phi_j)$, then $W$ admits a map to an initial object of $\calX$, and is therefore itself initial by Lemma \ref{sumoto}. This proves the ``only if'' direction of $(1)$. The converse follows from the uniqueness of initial objects.

Now suppose that $i = j$. We observe that $\id_{/Z}$ is a coproduct of $\phi_i$ and $\psi$ in the
$\infty$-category $\calX^{/Z}$. Since $\phi_i^{\ast}$ preserves coproducts, we deduce
that $\id_{Z_i}$ is a coproduct of $\phi^{\ast}(\phi_j): X \rightarrow Z_i$ and $\beta: \emptyset \rightarrow Z_i$ in $\calX^{/Z_i}$. Since $\beta$ is an initial object of $\calX^{/Z_i}$, we see that
$\phi^{\ast}(\phi_j)$ is an equivalence. The natural map $\gamma: \alpha \rightarrow \phi_i^{\ast}(\phi_i)$
corresponds to a commutative diagram
$$ \xymatrix{ W \ar[r]^{\alpha} \ar[d]^{\gamma_0} & Z_i \ar[d]^{\id_{Z_i} } \ar[d] \\
X \ar[r]^{ \phi_i^{\ast}(\phi_i)} &  Z_i }$$
in the $\infty$-category $\calX$. Consequently, $\alpha$ is an equivalence if and only if
$\gamma_0$ is an equivalence, if and only if $\gamma$ is an equivalence in
$\calX^{/Z_i}$. This proves $(2)$.
\end{proof}

\begin{proposition}\label{natash}
Let $\calC$ be a small $\infty$-category which admits finite limits, and let
$\calX$ be an $\infty$-category which satisfies the $\infty$-categorical Giraud axioms
$(i)-(iv)$ of Theorem \ref{mainchar}. Let $F: \calP(\calC) \rightarrow \calX$ be
a colimit-preserving functor. Suppose that the composition
$F \circ j: \calC \rightarrow \calX$ is left exact, where $j: \calC \rightarrow \calP(\calC)$ denotes the Yoneda embedding. Then $F$ is left exact.
\end{proposition}

\begin{proof}
According to Corollary \ref{allfinn}, to prove that $F$ is left exact, it will suffice to prove that $F$ preserves pullbacks and final objects. Since all final objects are equivalent, to prove that $F$ preserves final objects it suffices to exhibit a single final object $Z$ of $\calP(\calC)$ such that $FY \in \calX$ is final. Let $z$ be a final object of $\calC$ (which exists in virtue of our assumption that $\calC$ admits finite limits). Then $Z = j(z)$ is a final object of $\calP(\calC)$, since $j$ preserves limits by Proposition \ref{yonedaprop}. Consequently $F(Z) = f(z)$ is final, since $f$ is left-exact.

Let $\alpha: Y \rightarrow Z$ be a morphism in $\calP(\calC)$. We will say that $\alpha$ is {\it good} if for every pullback square
$$ \xymatrix{ W \ar[r] \ar[d] & Y \ar[d]^{\alpha} \\
X \ar[r] & Z}$$ 
in $\calP(\calC)$, the induced square
$$ \xymatrix{ F(W) \ar[r] \ar[d] & F(Y) \ar[d]^{F(\alpha)} \\
F(X) \ar[r]^{\beta} & F(Z)}$$ 
is a pullback in $\calX$. Note that Lemma \ref{transplantt} implies that the class of good morphisms in $\calP(\calC)$ is stable under composition.

We rephrase this condition that a morphism $\alpha$ be good in terms of the pullback functors
$\alpha^{\ast}: \calP(\calC)^{/Z} \rightarrow \calP(\calC)^{/Y}$, $F(\alpha)^{\ast}: \calX^{/F(Z)} \rightarrow \calP(\calC)^{/F(Y)}$. Application of the functor $F$ gives a map
$$ t: F \circ \alpha^{\ast} \rightarrow F(\alpha)^{\ast} \circ F$$ in the $\infty$-category
of functors from $\calP(\calC)^{/Z}$ to $\calX^{/F(Z)}$, and $\alpha$ is good if and only if $t$ is an equivalence. Note that $t$ is a natural transformation of colimit-preserving functors. Since the image of the Yoneda embedding $j: \calC \rightarrow \calP(\calC)$ generates $\calP(\calC)$ under colimits, it will suffice to prove that $t$ is an equivalence when evaluated on objects
of the form $\beta: j(x) \rightarrow Z$, where $x$ is an object of $\calC$.

Let us say that an object $Z \in \calP(\calC)$ is {\it good} if every morphism $\alpha: Y \rightarrow Z$ is good. In other words, an object $Z \in \calP(\calC)$ is good if $F$
is left exact at $Z$ in the sense of Definition \ref{swarpy}.
By repeating the above argument, we deduce that $Z$ is good if and only if every morphism of the form $\alpha: j(y) \rightarrow Z$ is good for $y \in \calC$. 

We next claim that for every object $z \in \calC$, the Yoneda image $j(z) \in \calP(\calC)$ is good.
In other words, we must show that for every pullback square
$$ \xymatrix{ W \ar[r] \ar[d] & j(y) \ar[d]^{\alpha} \\
j(x) \ar[r]^{\beta} & j(z) }$$
in $\calP(\calC)$, the induced square
$$ \xymatrix{ F(W) \ar[r] \ar[d] & f(x) \ar[d] \\
f(y) \ar[r] & f(z)} $$
is a pullback in $\calX$. Since the Yoneda embedding is fully faithful, we may suppose
that $\alpha$ and $\beta$ are the Yoneda images of morphisms $x \rightarrow z$,
$y \rightarrow z$. Since $j$ preserves limits, we may reduce to the case where the first
diagram is the Yoneda image of a pullback diagram in $\calC$. The desired result then follows from the assumption that $f$ is left exact.

To complete the proof that $F$ is left exact, it will suffice to prove that every object
of $\calP(\calC)$ is good. Because the Yoneda embedding $j: \calC \rightarrow \calP(\calC)$
generates $\calP(\calC)$ under colimits, it will suffice to prove that the collection of good objects of $\calP(\calC)$ is stable under colimits. According to Proposition \ref{appendicites}, it will suffice to prove that the collection of good objects of $\calP(\calC)$ is stable under coequalizers and small coproducts. 

We first consider the case of coproducts. Let $\{ Z_{i} \}_{i \in I}$ be a family of good objects
of $\calP(\calC)$ indexed by a (small) set $I$, and $\{ \phi_i: Z_i \rightarrow Z\}_{i \in I}$ be a family
of morphisms which exhibit $Z$ as a coproduct of the family $\{ Z_i \}_{i \in I}$.
Suppose given a pullback diagram
$$ \xymatrix{ W \ar[r] \ar[d] & j(y) \ar[d]^{\alpha} \\
j(x) \ar[r]^{\beta} & Z }$$
in $\calP(\calC)$. According to Proposition \ref{limiteval}, evaluation at the object
$y$ induces a colimit-preserving functor $\calP(\calC) \rightarrow \SSet$. Consequently,
we have a homotopy equivalence
$$ \bHom_{\calP(\calC)}( j(y), Z ) \simeq 
\coprod_{i \in I} \bHom_{\calP(\calC)}( j(y), Z_i) $$
in the homotopy category $\calH$. Therefore we may assume that
$\alpha$ factors as a composition
$$ j(y) \stackrel{\alpha'}{\rightarrow} Z_i \stackrel{\phi_i}{\rightarrow} Z$$
for some $i \in I$. By assumption, the morphism $\alpha'$ is good; it therefore suffices
to prove that $\phi_i$ is good. By a similar argument, we can replace
$\beta$ by a map $\phi_j: Z_j \rightarrow Z$, for some $j \in I$. We are now required to show
that if $$ \xymatrix{ W' \ar[r] \ar[d] & Z_i \ar[d]^{\phi_i} \\
Z_j \ar[r]^{\phi_j} & Z }$$ is a pullback diagram in $\calP(\calC)$, then
$$ \xymatrix{ F(W') \ar[r] \ar[d] & F(Z_i) \ar[d]^{\phi_i} \\
Z_j \ar[r]^{\phi_j} & Z }$$
is a pullback diagram in $\calX$. Since $F$ preserves initial objects, this follows immediately
from Lemma \ref{sumdescription}.

We now complete the proof by showing that the collection of good objects of $\calP(\calC)$ is stable under the formation of coequalizers. Let
$$ \xymatrix{ Z_{1} \ar@<.5ex>[r] \ar@<-.5ex>[r] & Z_0 \ar[r]^{s} & Z_{-1}. } $$ 
be a coequalizer diagram in $\calP(\calC)$, and suppose that $Z_0$ and $Z_1$ are good. We must show that any pullback diagram
$$ \xymatrix{ W \ar[r] \ar[d] & j(y) \ar[d]^{\alpha} \\
j(x) \ar[r]^{\beta} & Z_{-1} }$$
remains a pullback diagram after applying the functor $F$. The functor
$$ \calP(\calC) \rightarrow \Nerve(\Set)$$
$$ T \mapsto \Hom_{\h{\calP(\calC)}}( j(x), T)$$
can be written as a composition
$$ \calP(\calC) \rightarrow \SSet \stackrel{\pi_0}{\rightarrow} \Nerve(\Set)$$
where the first functor is given by evaluation at $x$. Both of these functors commute
with colimits. Consequently, we have a coequalizer diagram
$$ \xymatrix{ \Hom_{\h{\calP(\calC)}}(j(x), Z_{1}) \ar@<.5ex>[r] \ar@<-.5ex>[r] & 
\Hom_{\h{\calP(\calC)}}(j(x),Z_0) \ar[r] & \Hom_{\h{\calP(\calC)}}(j(x),Z_{-1}). } $$ 
in the category of sets. In particular, the map $\beta$ factors as a composition
$$ j(x) \stackrel{\beta'}{\rightarrow} Z_0 \stackrel{s}{\rightarrow} Z_{-1}. $$
Since we have already assumed that $\beta'$ is good, we can replace
$\beta$ by the map $s: Z_0 \rightarrow Z_{-1}$ in the above diagram. By a similar argument, we can replace $\alpha: Y \rightarrow Z_{-1}$ by the map $s: Z_0 \rightarrow Z_{-1}$. We now obtain the desired result by applying Proposition \ref{storytell}.
\end{proof}

We are now ready to complete the proof of Theorem \ref{mainchar}:

\begin{proposition}\label{precisechar}
Let $\calX$ be an $\infty$-category. Suppose that $\calX$ satisfies the $\infty$-categorical
Giraud axioms:
\begin{itemize}
\item[$(i)$] The $\infty$-category $\calX$ is presentable.
\item[$(ii)$] Colimits in $\calX$ are universal.
\item[$(iii)$] Coproducts in $\calX$ are disjoint.
\item[$(iv)$] Every groupoid object of $\calX$ is effective.
\end{itemize}
Then there exists a small $\infty$-category $\calC$ which admits finite limits, and an accessible left exact localization functor $\calP(\calC) \rightarrow \calX$. In particular, $\calX$ is an $\infty$-topos.
\end{proposition}

\begin{proof}
Let $\calX$ be an $\infty$-topos. According to Proposition \ref{tcoherent}, there exists a regular cardinal $\tau$ such that $\calX$ is $\tau$-accessible, and the full subcategory $\calX^{\tau}$ spanned by the $\tau$-compact objects of $\calX$ is stable under finite limits.
Let $\calC$ be a minimal model for $\calX^{\tau}$, so that there is an equivalence $\Ind_{\tau}(\calC) \rightarrow \calX$. The proof of Theorem \ref{pretop} shows that 
the inclusion $\Ind_{\tau}(\calC) \subseteq \calP(\calC)$ has a left adjoint $L$. 
The composition of $L$ with the Yoneda embedding $\calC \rightarrow \calP(\calC)$ can be identified with the Yoneda embedding $\calC \rightarrow \Ind_{\tau}(\calC)$, therefore preserves all limits which exist in $\calC$ (Proposition \ref{yonedaprop}). Applying Proposition \ref{natash}, we deduce that $L$ is left exact, so that $\Ind_{\tau}(\calC)$ is a left exact localization (automatically accessible) of $\calP(\calC)$. Since $\calX$ is equivalent to $\Ind_{\tau}(\calC)$, we conclude that
$\calX$ is also an accessible left exact localization of $\calP(\calC)$.
\end{proof}

\subsection{$\infty$-Topoi and Classifying Objects}\label{rezk2}

Let $\calX$ be an ordinary category, and let $X$ be an object of $\calX$. Let $\Sub(X)$\index{not}{SubX@$\Sub(X)$}
denote the partially ordered collection of {\it subobjects} of $X$: an object of $\Sub(X)$ is an equivalence class of monomorphisms $Y \rightarrow X$\index{gen}{subobject}. If $\calC$ is accessible, then $\Sub(X)$ is actually a set. If $\calX$ admits finite limits, then $\Sub(X)$ is contravariantly functorial in $X$:
given a subobject $Y \rightarrow X$ and any map $X' \rightarrow X$, the fiber product
$Y' = X' \times_{X} Y$ is a subobject of $X'$. A {\it subobject classifier}\index{gen}{classifying map!for subobjects} is an object $\Omega$ of $\calX$ which {\em represents} the functor $\Sub$. In other words, $\Omega$ has a universal subobject 
$\Omega_0 \subseteq \Omega$ such that every monomorphism $Y \rightarrow X$ fits into a {\em unique} Cartesian diagram
$$ \xymatrix{ Y \ar[r] \ar@{^{(}->}[d] & \Omega_0 \ar@{^{(}->}[d] \\
X \ar[r] & \Omega.}$$
(In this case, $\Omega_0$ is automatically a final object of $\calC$.)

Every topos has a subobject classifier. In fact, in the theory of {\it elementary topoi}, the existence of a subobject classifier is taken as one of the axioms. Thus, the existence of a subobject classifier is one of the defining characteristics of a topos. We would like to discuss the appropriate $\infty$-categorical generalization of the theory of subobject classifiers. The ideas presented here are due to Charles Rezk.

\begin{definition}\label{ugatoo}
Let $\calX$ be an $\infty$-category which admits pullbacks, and $S$ a collection of morphisms
of $\calX$ which is stable under pullback. We will say that a morphism $f: X \rightarrow Y$ {\it classifies $S$} if it is a final object of $\calO_{\calX}^{(S)}$ (see Notation \ref{ugaboo}). 
In this situation, we will also say that the object $Y \in \calX$ {\it classifies $S$}\index{gen}{classifying map!for a collection of morphisms}. 
A {\it subobject classifier} for $\calX$ is an object which classifies the collection of all monomorphisms in $\calX$.
\end{definition}

\begin{example}
The $\infty$-category $\SSet$ of spaces has a subobject classifier: namely, the discrete space $\{ 0, 1\}$ with two elements.
\end{example}

The following result provides a necessary and sufficient condition for the existence of a classifying object for $S$: 

\begin{proposition}\label{classexist}
Let $\calX$ be a presentable $\infty$-category in which colimits are universal, and let $S$
be a class of morphisms in $\calX$ which is stable under pullbacks. There exists a classifying object for $S$ if and only if the following conditions are satisfied:
\begin{itemize}
\item[$(1)$] The class $S$ is local $($Definition \ref{localitie}$)$.
\item[$(2)$] For every object $X \in \calX$, the full subcategory of
$\calX_{/X}$ spanned by the elements of $S$ is essentially small.
\end{itemize}
\end{proposition}

\begin{proof}
Let $s: \calX^{op} \rightarrow \hat{\SSet}$ be a functor which classifies the right fibration
$\calO_{\calX}^{(S)} \rightarrow \calX$. Then $S$ has a classifying object if and only if
$s$ is a representable functor. According to the representability criterion of Proposition \ref{representable}, this is equivalent to the assertion that $s$ preserves small limits, and the essential image of $s$ consists of essentially small spaces. According to Lemma \ref{ib3}, $s$ preserves small limits if and only if $(1)$ is satisfied. It now suffices to observe that for each
$X \in \calX$, the space $s(X)$ is essentially small if and only if the full subcategory of
$\calX_{/X}$ spanned by $S$ is essentially small. 
\end{proof}

Using Proposition \ref{classexist}, one can show that every $\infty$-topos has a subobject classifier. However, in the $\infty$-categorical context, the emphasis on {\em subobjects} misses the point. To see why, let us return to considering an ordinary category $\calX$ with a subobject classifier $\Omega$. By definition, for every object $X \in \calX$, we may identify maps $X \rightarrow \Omega$ with subobjects of $X$: that is, isomorphism classes of maps $Y \rightarrow X$ which happen to be monomorphisms. Even better would be an {\it object classifier}\index{gen}{classifying map!for objects}: that is, an object $\widetilde{\Omega}$ such that $\Hom_{\calX}(X, \widetilde{\Omega})$ could be identified with {\em arbitrary} maps $Y \rightarrow X$. But this is an unreasonable demand: if $Y \rightarrow X$ is not an monomorphism, then there may be automorphisms of $Y$ as an object of $\calX_{/X}$. It would be unnatural to ignore these automorphisms. However, it is also not possible to take them into account, since $\Hom_{\calX}(X, \widetilde{\Omega})$ must be a set rather than a groupoid.

If we allow $\calX$ to be an $\infty$-category, this objection loses its force. Informally speaking, we can consider the functor which associates to each $X \in \calX$ the maximal $\infty$-groupoid contained in $\calX_{/X}$ (this is contravariantly functorial in $X$, provided that $\calX$ has finite limits). We might hope that this functor is representable by some $\Omega_{\infty} \in \calX$, which we would then call an {\it object classifier}.

Unfortunately, a new problem arises: it is generally unreasonable to ask for the collection of {\em all} morphisms in $\calX$ to be classified by an object of $\calX$, since this would require each slice $\calX_{/X}$ to be essentially small (Proposition \ref{classexist}). This is essentially a technical difficulty, which we will circumvent by introducing a cardinality bound.

\begin{definition}\index{gen}{relatively $\kappa$-compact morphism}
Let $\calX$ be a presentable $\infty$-category. We will say that a morphism
$f: X \rightarrow Y$ is {\it relatively $\kappa$-compact} if, for every pullback diagram
$$ \xymatrix{ X' \ar[r] \ar[d]^{f'} & X \ar[d]^{f} \\
Y' \ar[r] & Y }$$
such that $Y'$ is $\kappa$-compact, $X'$ is also $\kappa$-compact.
\end{definition}

\begin{lemma}\label{sumarus}
Let $\calX$ be a presentable $\infty$-category, $\kappa$ a regular cardinal, $\calJ$ a $\kappa$-filtered $\infty$-category, and $\overline{p}: \calJ^{\triangleright} \rightarrow \calX$ a colimit diagram. Let $f: X \rightarrow Y$ be a morphism in $\calX$, where $Y$ is the image
under $\overline{p}$ of the cone point of $\calJ^{\triangleright}$. For each $\alpha$ in $\calJ$, let $Y_{\alpha} = \overline{p}(\alpha)$ and form a pullback diagram
$$ \xymatrix{ X_{\alpha} \ar[r] \ar[d]^{f_{\alpha}} & X \ar[d]^{f} \\
Y_{\alpha} \ar[r]^{g_{\alpha}} & Y. }$$
Suppose that each $f_{\alpha}$ is relatively $\kappa$-compact. Then $f$ is relatively $\kappa$-compact.
\end{lemma}

\begin{proof}
Let $Z$ be a $\kappa$-compact object of $\calX$, and $g: Z \rightarrow Y$ a morphism.
Since $Z$ is $\kappa$-compact and $\calJ$ is $\kappa$-filtered, there exists a $2$-simplex of $\calX$, corresponding to a diagram
$$ \xymatrix{ & Y_{\alpha} \ar[dr]^{g_{\alpha}} & \\
Z \ar[ur] \ar[rr]^{g} & & Y. }$$
Form a Cartesian rectangle $\Delta^2 \times \Delta^1 \rightarrow \calX$, which we will depict
as
$$ \xymatrix{ Z' \ar[r] \ar[d]^{f'} & X_{\alpha} \ar[r] \ar[d]^{f_{\alpha}} & X \ar[d]^{f} \\
Z \ar[r] & Y_{\alpha} \ar[r] & Y. }$$
Since $f'$ is a pullback of $f_{\alpha}$, we conclude that $Z'$ is $\kappa$-compact. 
Lemma \ref{transplantt} implies that $f'$ is also a pullback of $f$ along $g$, so that $f$ is relatively $\kappa$-compact as desired.
\end{proof}

\begin{lemma}\label{sumaris}
Let $\calX$ be a presentable $\infty$-category in which colimits are universal. Let $\tau > \kappa$ be regular cardinals such that
$\calX$ is $\kappa$-accessible and the full subcategory $\calX^{\tau}$ consisting of $\tau$-compact objects of $\calX$ is stable under pullbacks in $\calX$. Let $\alpha: \sigma \rightarrow \sigma'$ be a Cartesian transformation between pushout squares $\sigma, \sigma': \Delta^1 \times \Delta^1 \rightarrow \calX$, which we may view as a pushout square
$$ \xymatrix{ f \ar[r]^{\alpha} \ar[d]^{\beta} & g \ar[d]^{\beta'} \\
f' \ar[r]^{\alpha'} & g' }$$
in $\Fun(\Delta^1,\calX)$. Suppose that $f$, $g$, and $f'$ are relatively $\tau$-compact. Then $g'$ is relatively $\tau$-compact. 
\end{lemma}

\begin{proof}
Let $\calC$ denote the full subcategory of $\Fun(\Delta^1 \times \Delta^1, \calX)$ spanned by the pushout squares, and let $\calC^{\tau} = \calC \cap \Fun( \Delta^1 \times \Delta^1, \calX^{\tau})$.
Since the class of $\tau$-compact objects of $\calX$ is stable under pushouts (Corollary \ref{tyrmyrr}), we have a commutative diagram
$$ \xymatrix{ \calC^{\tau} \ar[r] \ar[d] & \Fun( \Lambda^2_0, \calX^{\tau})  \ar[d] \\
\calC \ar[r] & \Fun(\Lambda^2_0, \calX) }$$
where the horizontal arrows are trivial fibrations (Proposition \ref{lklk}). The proof of Proposition
\ref{horse1} shows that every object of $\Fun(\Lambda^2_0, \calX)$ can be written as the colimit of a $\tau$-filtered diagram in $\Fun(\Lambda^2_0, \calX^{\tau})$. It follows that $\sigma' \in \calC$
can be obtained as the colimit of a $\tau$-filtered diagram in $\calC^{\tau}$. Since
colimits in $\calX$ are universal, we conclude that the natural transformation $\alpha$
can be obtained as a $\tau$-filtered colimit of natural transformations
$ \alpha_{i}: \sigma_{i} \rightarrow \sigma'_{i}$ in $\calC^{\tau}$. Lemma \ref{limitscommute} implies that the inclusion $\calC \subseteq \Fun(\Delta^1 \times \Delta^1, \calX)$ is colimit-preserving. Consequently, we deduce that $g'$ can be written as a $\tau$-filtered colimit of morphisms
$\{ g'_{i} \}$ determined by restricting $\{ \alpha_i \}$. According to Lemma \ref{sumaris}, it will suffice to prove that each morphism $g'_{i}$ is relatively $\tau$-compact. In other words, we may replace $\sigma'$ by $\sigma'_i$ and thereby reduce to the case where $\sigma'$ belongs to $\calC^{\tau}$. Since $f,g,$ and $f'$ are relatively $\tau$-compact, we conclude that
$\sigma | \Lambda^2_0$ takes values in $\calX^{\tau}$. Since $\sigma$ is a pushout diagram, Corollary \ref{tyrmyrr} implies that $\sigma$ takes values in $\calX^{\tau}$. Now we observe that $g'$ is a morphism between $\tau$-compact objects of $\calX$, and therefore automatically relatively $\tau$-compact in virtue of our assumption that $\calX^{\tau}$ is stable under pullbacks in $\calX$.
\end{proof}

\begin{proposition}\label{cardyp}  
Let $\calX$ be a presentable $\infty$-category in which colimits are universal, and let
$S$ be a local class of morphisms in $\calX$. For each regular cardinal $\kappa$, 
let $S_{\kappa}$ denote the collection of all morphisms $f$ which belong to $S$ and are relatively $\kappa$-compact. If $\kappa$ is sufficiently large, then $S_{\kappa}$ has a classifying object.
\end{proposition}

\begin{proof}
Choose $\kappa'$ such that $\calX$ is $\kappa'$-accessible.
The restriction functor $r: 
\Fun( (\Lambda^2_2)^{\triangleleft}, \calX) \rightarrow \Fun( \Lambda^2_2, \calX)$
is accessible: in fact, it preserves all colimits (Proposition \ref{limiteval}). Let $g$ be a right adjoint to $r$ (a limit functor); Proposition \ref{adjoints} implies that $g$ is also accessible. Choose a regular cardinal $\kappa'' > \kappa'$ such that $g$ is $\kappa''$-continuous, and choose $\kappa \geq \kappa''$ such that $g$ carries $\kappa''$-compact objects of 
$\Fun(\Lambda^2_2, \calX)$ into $\Fun( (\Lambda^2_2)^{\triangleleft}, \calX^{\kappa})$. It follows that the class of $\kappa$-compact objects of $\calX$ is stable under pullbacks. We will show that $S_{\kappa}$ has a classifying object.

We will verify the hypotheses of Proposition \ref{classexist}. First, we must show that $S_{\kappa}$ is local. For this, we will verify condition $(3)$ of Lemma \ref{ib3}. We begin by showing that
$S_{\kappa}$ is stable under small coproducts. Let $\{ f_{\alpha}: X_{\alpha} \rightarrow Y_{\alpha} \}_{\alpha \in A}$ be a small collection of morphisms belonging to $S_{\kappa}$, and let $f: X \rightarrow Y$
be a coproduct $\coprod_{\alpha \in A} f_{\alpha}$ in $\Fun(\Delta^1,\calX)$. We wish to show that $f \in S_{\kappa}$. Since $S$ is local, we conclude that $f \in S$ (using Lemma \ref{ib3}). It therefore suffices to show that $f$ is relatively $\kappa$-compact. Suppose given a $\kappa$-compact
object $Z \in \calX$ and a morphism $g: Z \rightarrow Y$. Using Proposition \ref{extet} and Corollary \ref{util}, we conclude that $Y$ can be obtained as a $\kappa$-filtered colimit of objects
$Y_{A_0} = \coprod_{ \alpha \in A_0} Y_{\alpha}$, where $A_0$ ranges over the $\kappa$-small subsets of $A$. Since $Z$ is $\kappa$-compact, we conclude that there exists a factorization
$$ Z \stackrel{g'}{\rightarrow} Y_{A_0} \stackrel{g''}{\rightarrow} Y$$
of $g$. Form a Cartesian rectangle $\Delta^2 \times \Delta^1 \rightarrow \calX$,
$$ \xymatrix{ Z' \ar[r] \ar[d] & X_{A_0} \ar[r] \ar[d] & X \ar[d] \\
Z \ar[r] & Y_{A_0} \ar[r] & Y. }$$
Since $S$ is local, we can identify $X_{A_0}$ with the coproduct $\coprod_{\alpha \in A_0} X_{\alpha}$. Since colimits are universal, we conclude that $Z'$ is a coproduct of
objects $Z'_{\alpha} = X_{\alpha} \times_{Y_{\alpha} } Z$, where $\alpha$ ranges over
$A_0$. Since each $f_{\alpha}$ is relatively $\kappa$-compact, we conclude that each $Z'_{\alpha}$ is $\kappa$-compact. Thus $Z'$, as a $\kappa$-small colimit of $\kappa$-compact objects, is also $\kappa$-compact (Corollary \ref{tyrmyrr}). 

We must now show that for every pushout diagram
$$ \xymatrix{ f \ar[r]^{\alpha} \ar[d]^{\beta} & g \ar[d]^{\beta'} \\
f' \ar[r]^{\alpha'} & g' }$$
in $\calO_{\calX}$, if $\alpha$ and $\beta$ are Cartesian transformations and
$f,f',g \in S_{\kappa}$, then $\alpha'$ and $\beta'$ are also Cartesian transformations and $g' \in S_{\kappa}$. The first assertion follows immediately from Lemma \ref{ib3} (since $S$ is local), and we deduce also that $g' \in S$. It therefore suffices to show that $g$ is relatively $\kappa$-compact, which follows from Lemma \ref{sumaris}.

It remains to show that, for each $X \in \calX$, the full subcategory of $\calX_{/X}$ spanned by
the elements of $S$ is essentially small. Equivalently, we must show that the right fibration
$p: \calO_{\calX}^{(S)} \rightarrow \calX$ has essentially small fibers. Let
$F: \calX^{op} \rightarrow \hat{\SSet}$ classify $p$. Since $S$ is local, $F$ preserves limits.
The full subcategory of $\hat{\SSet}$ spanned by the essentially small Kan complexes is stable under small limits, and $\calX$ is generated by $\calX^{\kappa}$ under small ($\kappa$-filtered) colimits. Consequently, it will suffice to show that $F(X)$ is essentially small, when $X$ is $\kappa$-compact. In other words, we must show that there are only a bounded number equivalence classes of morphisms $f: Y \rightarrow X$ such that $f \in S_{\kappa}$. We now observe that if
$f \in S_{\kappa}$, then $f$ is relatively $\kappa$-compact, so that $Y$ also belongs to $\calX^{\kappa}$. We now conclude by observing that the $\infty$-category $\calX^{\kappa}$ is essentially small.
\end{proof}

We now give a characterization of $\infty$-topoi based on the existence of object classifiers.

\begin{theorem}[Rezk]\label{colimsurt}\index{gen}{classifying map!for relatively $\kappa$-compact morphisms}
Let $\calX$ be a presentable $\infty$-category. Then
$\calX$ is an $\infty$-topos if and only if the following conditions are satisfied:
\begin{itemize}
\item[$(1)$] Colimits in $\calX$ are universal.
\item[$(2)$] For all sufficiently large regular cardinals $\kappa$, there exists
a classifying object for the class of all relatively $\kappa$-compact morphisms
in $\calX$.
\end{itemize}
\end{theorem}

\begin{proof}
Assume that colimits in $\calX$ are universal. According to Theorems \ref{mainchar} and \ref{charleschar}, $\calX$ is an $\infty$-topos if and only if the class $S$ consisting of all morphisms of $\calX$ is local. This clearly implies $(2)$, in view of Proposition \ref{cardyp}. Conversely,
suppose that $(2)$ is satisfied, and let $S_{\kappa}$ be defined as in the statement of Proposition \ref{cardyp}. Proposition \ref{classexist} ensures that $S_{\kappa}$ is local for all sufficiently large regular cardinals $\kappa$. We note that $S = \bigcup S_{\kappa}$. It follows from Criterion $(3)$ of Lemma \ref{ib3} that $S$ is also local, so that $\calX$ is an $\infty$-topos.
\end{proof}
 
\section{Constructions of $\infty$-Topoi}\label{topcomp}

\setcounter{theorem}{0}

According to Definition \ref{itoposdef}, an $\infty$-category $\calX$ is an $\infty$-topos if and only if $\calX$ arises as an (accessible) left exact localization of a presheaf $\infty$-category $\calP(\calC)$. To complete the analogy with classical topos theory, we would like to have some concrete description of the collection of left exact localizations of $\calP(\calC)$. In \S \ref{leloc}, we will study left exact localization functors in general, and single out a special class which we call {\it topological} localizations. In \S \ref{cough}, we will study topological localizations of $\calP(\calC)$, and show that they are in bijection with {\it Grothendieck topologies} on the $\infty$-category $\calC$, in exact analogy with classical topos theory. In particular, given a Grothendieck topology on $\calC$, one can define an $\infty$-topos $\Shv(\calC) \subseteq \calP(\calC)$ of {\it sheaves on $\calC$}. In \S \ref{surjsurj}, we will characterize $\Shv(\calC)$ by a universal mapping property. 
Unfortunately, not every $\infty$-topos $\calX$ can be obtained as topological localization of an $\infty$-category of presheaves. Nevertheless, in \S \ref{cantopp} we will construct $\infty$-categories of sheaves which closely approximate $\calX$, using the formalism of {\it canonical topologies}. These ideas will be applied in \S \ref{chap6sec3}, to obtain a classification theorem for $n$-topoi.

\subsection{Left Exact Localizations}\label{leloc}

Let $\calX$ be an $\infty$-category. Up to equivalence, a localization $L: \calX \rightarrow \calY$
is determined by the collection $S$ of all morphisms $f: X \rightarrow Y$ in $\calX$ such that
$Lf$ is an equivalence in $\calY$ (Proposition \ref{localloc}). Our first result provides a useful criterion for testing the left-exactness of $L$.

\begin{proposition}\label{charleftloc}\index{gen}{localization!left exact}
Let $L: \calX \rightarrow \calY$ be a localization of $\infty$-categories. 
Suppose that $\calX$ admits finite limits. The following conditions are equivalent:
\begin{itemize}
\item[$(1)$] The functor $L$ is left exact.
\item[$(2)$] For every pullback diagram
$$ \xymatrix{ X' \ar[r] \ar[d]^{f'} & X \ar[d]^{f} \\
Y' \ar[r] & Y }$$ in $\calX$ such that $Lf$ is an equivalence in $\calY$,
$Lf'$ is also an equivalence in $\calY$.
\end{itemize}
\end{proposition}

\begin{proof}
It is clear that $(1)$ implies $(2)$. Suppose that $(2)$ is satisfied. We wish to show that $L$ is
left exact. Let $S$ be the collection of morphisms $f$ in $\calX$ such that $Lf$ is an equivalence.
Without loss of generality, we may identify $\calY$ with the full subcategory of $\calX$ spanned by the $S$-local objects. Since the final object $1 \in
\calX$ is obviously $S$-local, we have $L1 \simeq 1$. Thus it will
suffice to show that $L$ commutes with pullbacks. We observe that given
any diagram $X \rightarrow Y \leftarrow Z$, the pullback $LX  \times_{LY} LZ$ is a limit
of $S$-local objects of $\calX$, and therefore $S$-local. To complete the proof, it will suffice
to show that the natural map $f: X \times_Y Z \rightarrow LX \times_{LY} LZ$ belongs to $S$.
We can write $f$ as a composition
of maps $$X \times_Y Z \rightarrow X \times_{LY} Z \rightarrow LX
\times_{LY} Z \rightarrow LX \times_{LY} LZ.$$ The last two maps
are obtained from $X \rightarrow LX$ and $Z \rightarrow LZ$ by
base change. Assumption $(2)$ implies that they belong to $S$. Thus, it will
suffice to show that $f': X \times_Y Z \rightarrow X \times_{LY}
Z$ belongs to $S$. This map is a pullback of the diagonal $f'': Y \rightarrow Y \times_{LY} Y$, so it will suffice to prove that $f'' \in S$.
Projection to the first factor gives a left homotopy inverse $g: Y \times_{LY} Y
\rightarrow Y$ of $f''$, so it suffices to prove that $g \in S$. 
But $g$ is a base change of the morphism $Y \rightarrow LY$.
\end{proof}

\begin{proposition}\label{swimmer}
Let $\calX$ be a presentable $\infty$-category in which colimits are universal. Let $S$ be a class of morphisms in $\calX$, and let $\overline{S}$ be the strongly saturated class of morphisms generated by $S$. Suppose
that $S$ has the following property: for every pullback diagram
$$ \xymatrix{ X' \ar[r] \ar[d]^{f'} & X \ar[d]^{f} \\
Y' \ar[r] & Y }$$ in $\calX$, if $f \in S$, then $f' \in \overline{S}$. Then $\overline{S}$ is stable under pullbacks.
\end{proposition}

\begin{proof}
Let $S'$ be the set of all morphisms $f$ in $\calX$ with the property that for any pullback diagram
$$ \xymatrix{ X' \ar[r] \ar[d]^{f'} & X \ar[d]^{f} \\
Y' \ar[r] & Y, }$$
the morphism $f'$ belongs to $\overline{S}$. By assumption, $S \subseteq S'$. Using the fact
that colimits are universal, we deduce that $S'$ is strongly saturated. Consequently
$\overline{S} \subseteq S'$, as desired.
\end{proof}

\begin{corollary}\label{sweetums}
Let $\calX$ be a presentable $\infty$-category in which colimits are universal, let $S$ be a $($small$)$ set of morphisms in $\calX$, and let $\overline{S}$ denote the smallest strongly saturated class of morphisms which contains $S$ and is stable under pullbacks. Then $\overline{S}$ is generated $($as a strongly saturated class of morphisms$)$ by a $($small$)$ set.
\end{corollary}

\begin{proof}
Choose a (small) set $U$ of objects of $\calX$ which generates $\calX$ under colimits. Enlarging $U$ if necessary, we may suppose that $U$ contains the codomain of every morphism belonging to $S$.
Let $S'$ be the set of all morphisms $f'$ which fit into a pullback diagram
$$ \xymatrix{ X' \ar[r] \ar[d]^{f'} & X \ar[d]^{f} \\
Y' \ar[r] & Y }$$
where $f \in S$ and $Y' \in U$, and let $\overline{S}'$ denote the strongly saturated class of morphisms generated by $S'$. To complete the proof it will suffice to show that $\overline{S}' = \overline{S}$.
The inclusions $S \subseteq S' \subseteq \overline{S}' \subseteq \overline{S}$ are obvious.
To show that $\overline{S} \subseteq \overline{S}'$, it will suffice to show that $\overline{S}'$ is stable under pullbacks. In view of Proposition \ref{swimmer}, it will suffice to show that for every pullback diagram
$$ \xymatrix{ X'' \ar[r] \ar[d]^{f''} & X' \ar[d]^{f'} \\
Y'' \ar[r] & Y' }$$
such that $f' \in S'$, the morphism $f''$ belongs to $\overline{S}'$. Using our assumption that colimits in $\calX$ are universal and that $U$ generates $\calX$ under colimits, we can reduce to the case where
$Y'' \in U$. In this case, $f'' \in S'$ by construction.
\end{proof}

Recall that a morphism $f: Y \rightarrow Z$ in an $\infty$-category $\calX$ is a {\it monomorphism}\index{gen}{monomorphism}
if it is a $(-1)$-truncated object of the $\infty$-category $\calX_{/Z}$. Equivalently, $f$ is a monomorphism if for every object $X \in \calX$, the induced map
$$ \bHom_{\calX}(X,Y) \rightarrow \bHom_{\calX}(X,Z)$$
exhibits $\bHom_{\calX}(X,Y) \in \calH$ as a summand of $\bHom_{\calX}(X,Z)$ in the homotopy category $\calH$. If we fix $Z \in \calX$, then the collection of equivalence classes of monomorphisms $Y \rightarrow Z$ are partially ordered under inclusion. We will denote this partially ordered collection by $\Sub(Z)$.

\begin{proposition}\label{subobjset}\index{not}{SubX@$\Sub(X)$}
Let $\calX$ be a presentable $\infty$-category, and let $X$ be an object of $\calX$. Then
$\Sub(X)$ is a $($small$)$ partially ordered set.
\end{proposition}

\begin{proof}
By definition, the partially ordered set $\Sub(X)$ is characterized by the existence of
an equivalence 
$$ \tau_{ \leq -1} \calX_{/X} \rightarrow \Nerve( \Sub(X)) .$$ 
Propositions \ref{slicstab} and \ref{maketrunc} imply that $\Nerve(\Sub(X))$ is presentable.
Consequently, there exists a small subset $S \subseteq \Sub(X)$ which generates $\Nerve (\Sub(X))$ under colimits. It follows that every element of $\Sub(X)$ can be written as the supremum of a subset of $S$, so that $\Sub(X)$ is also small.
\end{proof}

\begin{definition}\label{deftoploc}\index{gen}{localization!topological}\index{gen}{topological!localization}\index{gen}{topological!class of morphisms}
Let $\calX$ be a presentable $\infty$-category, and let $\overline{S}$ be a strongly saturated class of morphisms of $\calX$. We will say that $\overline{S}$ is {\it topological} if the following conditions are satisfied:
\begin{itemize}
\item[$(1)$] There exists $S \subseteq \overline{S}$ consisting of {\em monomorphisms} such that
$S$ generates $\overline{S}$ as a strongly saturated class of morphisms.
\item[$(2)$] Given a pullback diagram
$$ \xymatrix{ X' \ar[r] \ar[d]^{f'} & X \ar[d]^{f} \\
Y' \ar[r] & Y }$$
in $\calX$ such that $f$ belongs to $\overline{S}$, the morphism $f'$ also belongs to $\overline{S}$.
\end{itemize}

We will say that a localization $L: \calX \rightarrow \calY$ is {\it topological} if the collection $\overline{S}$ of all morphisms $f: X \rightarrow Y$ in $\calX$ such that $Lf$ is an equivalence
is topological.
\end{definition}

\begin{proposition}\label{toplocsmall}
Let $\calX$ be a presentable $\infty$-category in which colimits are universal, and let
$\overline{S}$ be a strongly saturated class of morphisms of $\calX$ which is topological.
Then there exists a $($small$)$ subset $S_0 \subseteq \overline{S}$ which consists of monomorphisms and generates $\overline{S}$ as a strongly saturated class of morphisms.
\end{proposition}

\begin{proof}
For every object $U \in \calX$, let $\Sub'(U) \subseteq \Sub(U)$ denote the collection of equivalence classes of monomorphisms $U' \rightarrow U$ which belong to $\overline{S}$.
Choose a small collection of objects $\{ U_{\alpha} \}_{\alpha \in A}$ which generates $\calX$ under colimits. For each $\alpha \in A$ and each element
$\widetilde{\alpha} \in \Sub'(U_{\alpha})$, choose a representative monomorphism $f_{\widetilde{\alpha}}: V_{\widetilde{\alpha}} \rightarrow U_{\alpha}$ which belongs to $\overline{S}$. Let 
$$S_0 = \{ f_{\widetilde{\alpha}} | \alpha \in A, \widetilde{\alpha} \in \Sub'(U_{\alpha}) \}.$$
It follows from Proposition \ref{subobjset} that $S_0$ is a (small) set. Let $\overline{S}_0$ denote the strongly saturated class of morphisms generated by $S_0$. We will show that $\overline{S}_0 = \overline{S}$.

Let $\calX^{0}$ be the full subcategory of $\calX$ spanned by objects $U$ with the following property: if $f: V \rightarrow U$ is a monomorphism and $f \in \overline{S}$, then $f \in \overline{S_0}$. By construction, for each $\alpha \in A$, $U_{\alpha} \in \calX^{0}$. 
Since colimits in $\calX$ are universal, it is easy to see that $\calX^{0}$ is stable under colimits in $\calX$. It follows that $\calX^{0} = \calX$, so that every monomorphism which belongs to $\overline{S}$ also belongs to $\overline{S_0}$. Since $\overline{S}$ is generated by monomorphisms, we conclude that $\overline{S} = \overline{S_0}$, as desired.
\end{proof}

\begin{corollary}\label{topaccess}
Let $\calX$ be a presentable $\infty$-category. Every topological localization $L: \calX \rightarrow \calY$ is accessible and left exact. 
\end{corollary}

\subsection{Grothendieck Topologies and Sheaves in Higher Category Theory}\label{cough}

Every ordinary topos is equivalent to the category of sheaves on some Grothendieck site. This can be deduced from the following pair of statements:
\begin{itemize}
\item[$(i)$] Every topos is equivalent to a left exact localization of the some presheaf category
$\Set^{\calC^{op}}$.
\item[$(ii)$] There is a bijective correspondence between left exact localizations of
$\Set^{\calC^{op}}$ and Grothendieck topologies on $\calC$.
\end{itemize}
In \S \ref{chap6sec1}, we proved the $\infty$-categorical analogue of assertion $(i)$. Unfortunately, $(ii)$ is not quite true in the $\infty$-categorical setting. In this section, we will establish a slightly weaker statement: for every $\infty$-category $\calC$, there is a bijective correspondence between Grothendieck topologies on $\calC$ and {\em topological} localizations of $\calP(\calC)$ (Proposition \ref{toprole}).
Our first step is to introduce the $\infty$-categorical analogue of a Grothendieck site. The following definition is taken from \cite{toen}:

\begin{definition}\label{grtop}\index{gen}{sieve}\index{gen}{sieve!covering}
Let $\calC$ be a $\infty$-category. A {\it sieve} on $\calC$ is a full subcategory of $\calC^{(0)} \subseteq \calC$ having the property that if $f: C \rightarrow D$ is a morphism in $\calC$, and $D$ belongs to $\calC^{(0)}$, then $C$ also belongs to $\calC^{(0)}$.\index{not}{calC0@$\calC^{(0)}$}

Observe that if $f: \calC \rightarrow \calD$ is a functor between $\infty$-categories and
$\calD^{(0)} \subseteq \calD$ is a sieve on $\calD$, then $f^{-1} \calD^{(0)} = \calD^{(0)} \times_{\calD} \calC$ is a sieve on $\calC$.
Moreover, if $f$ is an equivalence, then $f^{-1}$ induces a bijection between sieves on
$\calD$ and sieves on $\calC$.

If $C \in \calC$ is an object, then a {\it sieve on $C$} is a sieve on the $\infty$-category $\calC_{/C}$. Given a morphism $f: D \rightarrow C$ and a sieve $\calC^{(0)}_{/C}$ on $C$, we let $f^{\ast} \calC^{(0)}_{/C}$ denote the unique sieve on $D$ such that $f^{\ast} \calC^{(0)}_{/C}
\subseteq \calC_{/D}$ and
$\calC^{(0)}_{/C}$ determine same sieve on $\calC_{/f}$.

A {\it Grothendieck topology}\index{gen}{Grothendieck topology} on an $\infty$-category $\calC$ consists of a specification, for each
object $C$ of $\calC$, of a collection of sieves on $C$, which we will refer to as {\it covering sieves}. The collections of covering sieves are required to possess the following properties:

\begin{itemize}
\item[$(1)$] If $C$ is an object of $\calC$, then the sieve $\calC_{/C} \subseteq \calC_{/C}$ on $C$ is a covering sieve.

\item[$(2)$] If $f: C \rightarrow D$ is a morphism in $\calC$, and $\calC^{(0)}_{/D}$ is a covering sieve on $D$, then $f^{\ast} \calC^{(0)}_{/C}$ is a covering sieve on $C$.

\item[$(3)$] Let $C$ be an object of $\calC$, $\calC^{(0)}_{/C}$ a covering sieve on $C$, and $\calC^{(1)}_{/C}$ an arbitrary sieve on $C$. Suppose that, for each $f: D \rightarrow C$ belonging to the sieve $\calC^{(0)}_{/C}$, the pullback
$f^{\ast} \calC^{(1)}_{/C}$ is a covering sieve on $D$. Then $\calC^{(1)}_{/C}$ is a covering sieve on $C$.
\end{itemize}
\end{definition}

\begin{example}\label{trivvtop}
Any $\infty$-category $\calC$ may be equipped with the {\it trivial topology}, in which a sieve
$\calC^{(0)}_{/C}$ on an object $C$ of $\calC$ is covering if and only if $\calC^{(0)}_{/C} = \calC_{/C}$.
\end{example}

\begin{remark}
In the case where $\calC$ is (the nerve of) an ordinary category, the definition given above reduces to the usual notion of a Grothendieck topology on $\calC$. Even in the general case, a Grothendieck
topology on $\calC$ is just a Grothendieck topology on the homotopy category $h \calC$.
This is not completely obvious, since for an object $C$ of $\calC$, the functor
$$ \eta: \h{(\calC_{/C})} \rightarrow (\h{\calC})_{/C}$$ is usually not an equivalence of categories.
A morphism from on the left hand side corresponds to a commutative triangle
$$ \xymatrix{ D \ar[rr] \ar[dr] & & D' \ar[dl] \\
& C & }$$ given by a {\it specified} $2$-simplex $\sigma: \Delta^2 \rightarrow \calC$ (taken modulo homotopy), while on the right hand side one requires only that the above diagram commutes up to homotopy: this amounts to requiring the existence of $\sigma$, but $\sigma$ itself is not taken as part of the data.

Although $\eta$ need not be an equivalence of categories, $\eta^{\ast}$ does induce
a bijection from the set of sieves on $(\h{\calC})_{/C}$ to the set of sieves on $\h{(\calC_{/C})}$: for this, it suffices to observe that $\eta$ induces surjective maps
$$ \Hom_{\h(\calC_{/C})} (D,D') \rightarrow \Hom_{ (\h{\calC})_{/C} }(D,D')$$
on morphism sets, which is obvious from the description given above.
\end{remark}

The main objective of this section is to prove that for any (small) $\infty$-category $\calC$, there is a bijective correspondence between Grothendieck topologies on $\calC$ and (equivalence classes of) topological localizations of $\calP(\calC)$. We begin by establishing a correspondence between sieves on $\calC$ and $(-1)$-truncated objects of $\calP(\calC)$. 
For each object $U \in \calP(\calC)$, let $\calC^{(0)}(U) \subseteq \calC$ be the full subcategory spanned by those objects $C \in \calC$ such that $U(C) \neq \emptyset$. It is easy to see that
$\calC^{(0)}(U)$ is a sieve on $\calC$. Conversely, given a sieve $\calC^{(0)} \subseteq \calC$, there is a unique map $\calC \rightarrow \Delta^1$ such that $\calC^{(0)}$ is the preimage of
$\{0\}$. This construction determines a bijection between sieves on $\calC$ and functors
$f: \calC \rightarrow \Delta^1$, and we may identify $\Delta^1$ with the full subcategory of
$\SSet^{op}$ spanned by the objects $\emptyset, \Delta^0 \in \Kan$. Since every $(-1)$-truncated Kan complex is equivalent to either $\emptyset$ or $\Delta^0$, we conclude:

\begin{lemma}\label{siffer}
For every small $\infty$-category $\calC$, the construction
$U \mapsto \calC^{(0)}(U)$ determines a bijection between the set of equivalence classes of
$(-1)$-truncated objects of $\calP(\calC)$ and the set of all sieves on $\calC$.
\end{lemma}

We now introduce a relative version of the above construction. Let $\calC$ be a small $\infty$-category as above, and let $j: \calC \rightarrow \calP(\calC)$ be
the Yoneda embedding. Let $C \in \calC$ be an object, and let $i: U \rightarrow j(C)$ be a monomorphism in $\calP(\calC)$. Let $\calC_{/C}(U)$ denote the full subcategory of $\calC$
spanned by those objects $f: D \rightarrow C$ of $\calC_{/C}$ such that there exists a commutative triangle
$$ \xymatrix{ j(D) \ar[rr]^{j(f)} \ar[dr] & & j(C) \\
& U \ar[ur]^{i} & }$$
It is easy to see that $\calC_{/C}(U)$ is a sieve on $C$. Moreover, it clear that if $i: U \rightarrow j(C)$ and $i': U' \rightarrow j(C)$ are equivalent subobjects of $j(C)$, then $\calC_{/C}(U) = \calC_{/C}(U')$. 

\begin{proposition}\label{surry}
Let $\calC$ be a small $\infty$-category containing an object $C$, and let $j: \calC \rightarrow \calP(\calC)$ be the Yoneda embedding. The construction described above yields a bijection
$$ ( i: U \rightarrow j(C) ) \mapsto \calC_{/C}(U)$$
from $\Sub(j(C))$ to the set of all sieves on $C$.
\end{proposition}

\begin{proof}
Use Corollary \ref{swapKK} to reduce to Lemma \ref{siffer}.
\end{proof}

\begin{definition}\label{defsheaff}\index{gen}{sheaf}\index{not}{ShvC@$\Shv(\calC)$}
Let $\calC$ be a $($small$)$ $\infty$-category equipped with a Grothendieck topology.
Let $S$ be the collection of all monomorphisms $U \rightarrow j(C)$ which correspond
to covering sieves $\calC^{(0)}_{/C} \subseteq \calC_{/C}$. An object
$\calF \in \calP(\calC)$ is a {\it sheaf} if it is $S$-local. We let $\Shv(\calC)$ denote the full subcategory of $\calP(\calC)$ spanned by $S$-local objects. 
\end{definition}

\begin{lemma}\label{stokeworth}
Let $\calC$ be a $($small$)$ $\infty$-category equipped with a Grothendieck topology. Then
$\Shv(\calC)$ is a topological localization of $\calP(\calC)$. In particular,
$\Shv(\calC)$ is an $\infty$-topos. 
\end{lemma}

\begin{proof}
By definition, $\Shv(\calC) = S^{-1} \calP(\calC)$, where $S$ is the class of all monomorphisms 
$i: U \rightarrow j(C)$ which correspond to covering sieves on $C \in \calC$. Let $\overline{S}$ be the strongly saturated class of morphisms generated by $S$; we wish to show that $\overline{S}$ is stable under pullback.

Let $S'$ denote the collection of all morphisms $f: X \rightarrow Y$ such that for any
pullback diagram $\sigma: \Delta^1 \times \Delta^1 \rightarrow \calP(\calC)$ depicted as follows:
$$ \xymatrix{ X' \ar[d]^{f'} \ar[r] & X \ar[d]^{f} \\
Y' \ar[r]^{g} & Y, }$$
the morphism $f'$ belongs to $\overline{S}$. Since colimits in $\calP(\calC)$ are universal, it is easy to prove that $S'$ is strongly saturated. We wish to prove that $\overline{S} \subseteq S'$. Since $\overline{S}$ is the smallest saturated class containing $S$, it will suffice to prove that $S \subseteq S'$. We may therefore suppose that $Y = j(C)$ in the diagram above, and that $f: X \rightarrow j(C)$ is the monomorphism corresponding to a covering sieve
$\calC^{(0)}_{/C}$ on $C$. 

Since $\calP(\calC)_{/j(C)} \simeq \calP(\calC_{/C})$ is generated under colimits by the Yoneda embedding, there exists a diagram $p: K \rightarrow \calC_{/C}$ such that the composite
map $j \circ p: K \rightarrow \calP(\calC)_{/j(C)}$ has $g: Y' \rightarrow j(C)$ as a colimit.
Because colimits in $\calP(\calC)$ are universal, we can extend $j \circ p$ to a 
diagram $P: K \rightarrow (\calP(\calC)^{\Delta^1})_{/f}$ which carries each vertex 
$k \in K$ to a pullback diagram,
$$ \xymatrix{ X_k \ar[d]^{f_k} \ar[r] & X \ar[d] \\
j(D_k) \ar[r]^{j(g_k)} & j(C) }$$
such that $\sigma$ is a colimit of $P$. Each $f_k$ is a monomorphism associated to the
covering sieve $g_k^{\ast} \calC^{(0)}_{/C}$, and therefore belongs to $S \subseteq \overline{S}$. It follows that $f'$ is a colimit in $\calP(\calC)^{\Delta^1}$ of morphisms belonging to $\overline{S}$, and therefore itself belongs to $\overline{S}$.
\end{proof}

The next lemma ensures us that we can recover a Grothendieck topology on $\calC$ from
its $\infty$-category of sheaves $\Shv(\calC) \subseteq \calP(\calC)$.

\begin{lemma}\label{recloose}
Let $\calC$ be a $($small $)$ $\infty$-category equipped with a Grothendieck topology, and
let $L: \calP(\calC) \rightarrow \Shv(\calC)$ denote a left adjoint to the inclusion. Let
$j: \calC \rightarrow \calP(\calC)$ denote the Yoneda embedding, and let
$i: U \rightarrow j(C)$ be a monomorphism corresponding to a sieve $\calC^{(0)}_{/C}$
on $C$. Then $Li$ is an equivalence if and only if
$\calC^{(0)}_{/C}$ is a covering sieve.
\end{lemma}

\begin{proof}
It is clear that if $\calC^{(0)}_{/C}$ is a covering sieve, then $Li$ is an equivalence. Conversely, suppose that $Li$ is an equivalence. Then $\tau_{\leq 0}(L i)$ is an equivalence. In view of
Proposition \ref{compattrunc}, we can identify $\tau_{\leq 0}(Li)$ with $L (\tau_{\leq 0} i)$.
The morphism $\tau_{\leq 0} i$ can be identified with a monomorphism 
$\eta: \calF \subseteq \Hom_{\h{\calC}}( \bigdot, C)$ in the ordinary category of presheaves of sets on $\h{\calC}$, where $$\calF(D) = \{ f \in \Hom_{\h{\calC}}(D,C): f \in \calC^{(0)}_{/C} \}.$$ 
If $\eta$ becomes an equivalence after sheafification, then the identity map
$\id_{C}: C \rightarrow C$ belongs to $\calF(C)$ locally; in other words, 
there exists a collection of morphisms $\{ f_{\alpha}: C_{\alpha} \rightarrow C \}$ which
generated a covering sieve on $C$ such that each $f_{\alpha}$ belongs to
$\calF(C_{\alpha})$, and therefore to $\calC^{(0)}_{/C}$. It follows that
$\calC^{(0)}_{/C}$ contains a covering sieve on $C$ and is therefore itself covering.
\end{proof}

We may summarize the results of this section as follows:

\begin{proposition}\label{toprole}
Let $\calC$ be a small $\infty$-category. There is a bijective correspondence
between Grothendieck topologies on $\calC$ and $($equivalence classes of$)$ topological localizations of $\calP(\calC)$.
\end{proposition}

\begin{proof}
According to Lemma \ref{stokeworth}, every Grothendieck topology on $\calC$ determines
a topological localization $\Shv(\calC) \subseteq \calP(\calC)$. Lemma \ref{recloose} shows that two Grothendieck topologies which determine the same $\infty$-categories of sheaves must coincide. To complete the proof, it will suffice to show that {\em every} topological localization of
$\calP(\calC)$ arises in this way. Let $\overline{S}$ be a strongly saturated collection of morphisms in $\calP(\calC)$, and suppose that $\overline{S}$ is topological. Let $S \subseteq \overline{S}$
be the collection of all monomorphisms $U \rightarrow j(C)$ which belong to $\overline{S}$, where
$j: \calC \rightarrow \calP(\calC)$ denotes the Yoneda embedding. Since the objects $\{ j(C) \}_{C \in \calC}$ generate $\calP(\calC)$ under colimits, and colimits in $\calP(\calC)$ are universal, 
we conclude that every monomorphism in $\overline{S}$ is a colimit of elements of $S$. 
Since $\overline{S}$ is generated by monomorphisms, we conclude that $\overline{S}$ is generated by $S$. 

Let us say that a sieve $\calC^{(0)}_{/C} \subseteq \calC_{/C}$ on an object $C \in \calC$
is {\it covering} if the corresponding monomorphism $U \rightarrow j(C)$ belongs to $S$.
We will show that the collection of covering sieves determines a Grothendieck topology on $\calC$. Granting this, we observe that $\overline{S}^{-1} \calP(\calC)$ is the $\infty$-category $\Shv(\calC) \subseteq \calP(\calC)$ of sheaves with respect to this Grothendieck topology, which will complete the proof.

We now verify the axioms $(1)$ through $(3)$ of Definition \ref{grtop}:

\begin{itemize}
\item[$(1)$] Every sieve of the form $\calC_{/C} \subseteq \calC_{/C}$ is covering, since every identity map $\id_{j(C)}: j(C) \rightarrow j(C)$ belongs to $S$.
\item[$(2)$] Let $f: C \rightarrow D$ be a morphism in $\calC$, and let
$\calC^{(0)}_{/D} \subseteq \calC_{/D}$ be a covering sieve, corresponding to a monomorphism
$i: U \rightarrow j(D)$ which belongs to $S$. Then $f^{\ast} \calC^{(0)}_{/C} \subseteq \calC_{/C}$
corresponds to a monomorphism $u: U' \rightarrow j(C)$ which is a pullback of $i$ along
$j(f)$, and therefore belongs to $S$ (since $\overline{S}$ is stable under pullbacks).
\item[$(3)$] Let $C$ be an object of $\calC$, $\calC^{(0)}_{/C}$ a covering sieve on $C$ corresponding to a monomorphism $i: U \rightarrow j(C)$ which belongs to $S$, 
and $\calC^{(1)}_{/C}$ an arbitrary sieve on $C$ corresponding to a monomorphism
$v: U' \rightarrow j(C)$. Suppose that, for each $f: D \rightarrow C$ belonging to the sieve $\calC^{(0)}_{/C}$, the pullback $f^{\ast} \calC^{(1)}_{/C}$ is a covering sieve on $D$.
Since $j': \calC_{/C} \rightarrow \calP(\calC)_{/j(C)}$ is a fully faithful embedding which generates $\calP(\calC)_{/j(C)}$ under colimits (see the proof of Corollary \ref{swapKK}), we conclude
there is a diagram $K \rightarrow \calC_{/C}$ such that $j' \circ K$ has colimit $i'$.
Since colimits in $\calP(\calC)$ are universal, we conclude that the map
$v': U \times_{j(C)} U' \rightarrow U$ is a colimit of morphisms of the form
$j(D) \times_{j(C)} U' \rightarrow j(D)$, which belong to $\overline{S}$ by assumption.
Since $\overline{S}$ is stable under colimits, we conclude that $i''$ belongs to $\overline{S}$.
We now have a pullback diagram
$$ \xymatrix{ U \times_{j(C)} U' \ar[r]^{v'} \ar[d]^{u'} & U \ar[d]^{u} \\
U' \ar[r]^{v} & j(C). }$$
By assumption, $u \in S$. 
Thus $v \circ u' \sim u \circ v' \in \overline{S}$. Since $u'$ is a pullback of
$u$, we conclude that $u' \in \overline{S}$, so that $v \in \overline{S}$.
This implies that $\calC^{(1)}_{/C} \subseteq \calC_{/C}$ is a covering sieve, as we wished to prove.
\end{itemize}
\end{proof}

For later use, we record the following characterization of initial objects in $\infty$-categories of sheaves:

\begin{proposition}\label{suture}\index{gen}{initial object!in an $\infty$-category of sheaves}
Let $\calC$ be a small $\infty$-category equipped with a Grothendieck topology, and let
$\calC' \subseteq \calC$ denote the full subcategory spanned by those objects
$C \in \calC$ such that $\emptyset \subseteq \calC_{/C}$ is a covering sieve on $C$.
An object $\calF \in \Shv(\calC)$ is initial if and only if it satisfies the following conditions:
\begin{itemize}
\item[$(1)$] If $C \in \calC'$, then $\calF(C)$ is contractible.
\item[$(2)$] If $C \notin \calC'$, then $\calF(C)$ is empty.
\end{itemize}
\end{proposition}

\begin{proof}
Let $L: \calP(\calC) \rightarrow \Shv(\calC)$ be a left adjoint to the inclusion, and
let $\emptyset$ be an initial object of $\calP(\calC)$. Then $L \emptyset$ is an initial object of $\Shv(\calC)$. Since $L$ is left exact, it preserves $(-1)$-truncated objects, as does the inclusion
$\Shv(\calC) \subseteq \calP(\calC)$. Thus $L \emptyset$ is $(-1)$-truncated, and corresponds to some sieve $\calC^{(0)} \subseteq \calC$ (Lemma \ref{siffer}). As we saw in the proof of Lemma \ref{recloose}, a sieve $\calC^{(0)}$ classifies an object of
$\Shv(\calC)$ is and only if $\calC^{(0)}$ is saturated in the following sense:
if $C \in \calC$ and the induced sieve $\calC^{(0)} \times_{\calC} \calC_{/C}$ is covering,
then $C \in \calC^{(0)}$. An initial object
of $\Shv(\calC)$ is an initial object of $\tau_{\leq -1} \Shv(\calC)$, and must therefore
correspond to the {\em smallest} saturated sieve on $\calC$. An easy argument shows that this sieve is $\calC'$, and that $\calF \in \calP(\calC)$ is a $(-1)$-truncated object classified by
$\calC'$ if and only if conditions $(1)$ and $(2)$ are satisfied.
\end{proof}

\subsection{Effective Epimorphisms}\label{surjsurj}

In classical topos theory, the assumption that every equivalence relation is effective
leads to a bijective correspondence between equivalence relations on an object $X$
and {\em effective epimorphisms} $X \rightarrow Y$. The purpose of this section is to
generalize the notion of an effective epimorphism to the $\infty$-categorical setting.

Our primary interest is studying the class of effective epimorphisms in an $\infty$-topos $\calX$.
However, we will later need to employ the same ideas when $\calX$ is an $n$-topos, for $n < \infty$. It is therefore convenient to work in a slightly more general setting.

\begin{definition}\index{gen}{semitopos}
An $\infty$-category $\calX$ is a {\it semitopos} if it satisfies the following conditions:
\begin{itemize}
\item[$(1)$] The $\infty$-category $\calX$ is presentable.
\item[$(2)$] Colimits in $\calX$ are universal.
\item[$(3)$] For every morphism $f: U \rightarrow X$, the underlying groupoid
of the \Cech nerve $\mCech(f)$ is effective (see \S \ref{gengroup}). 
\end{itemize}
\end{definition}

\begin{remark}
Every $\infty$-topos is a semitopos; this follows immediately from Theorem \ref{mainchar}.
\end{remark}

\begin{remark}
If $\calX$ is a semitopos, then so is $\calX_{/X}$ for every object $X \in \calX$.
\end{remark}

\begin{proposition}\label{slurpme}
Let $\calX$ be a semitopos. Let $p: U \rightarrow X$ be a morphism
in $\calX$, let $U_{\bigdot}$ be the underlying simplicial object of the \Cech nerve $\mCech(p)$, let
$V \in \calX$ be a colimit of $U_{\bigdot}$. The induced diagram
$$ \xymatrix{ U \ar[rr] \ar[dr]^{p} & & V \ar[dl]^{p'} \\
& X & }$$
identifies $p'$ with a $(-1)$-truncation of $p$ in $\calX_{/X}$.
\end{proposition}

\begin{proof}
We first show that $V$ is $(-1)$-truncated. It suffices to show
that the diagonal map $V \rightarrow V \times_{X} V$ is an
equivalence. We may identify $V$ with $V \times_{V} V$.
Since colimits in $\calX$ are universal, it will suffice to prove
that for each $m,n \geq 0$, the natural map
$$ p_{n.m}: U_m \times_{V} U_n \rightarrow U_{m} \times_{X} U_{n}$$
is an equivalence. We next observe that each $p_{n,m}$ is a pullback of
$$ p_{0,0}: U \times_{V} U \rightarrow U \times_{X} U. $$
Because $U_{\bigdot}$ is an effective groupoid, both sides may be identified with $U_{1}$.

To complete the proof, it suffices to show that the natural map $\bHom_{\calX_{/X}}(p',q)
\rightarrow \bHom_{\calX_{/X}}(p,q)$ is an equivalence whenever $q: E \rightarrow X$
is a monomorphism. Note that both sides are either empty or contractible. We must show that
if $\bHom_{\calX_{/X}}(p,q)$ is nonempty, then so is $\bHom_{\calX_{/X}}(p',q)$. 
We observe that the map $\calX_{/q} \rightarrow \calX_{/X}$ is fully faithful, 
and that its essential image is a sieve on $\calX_{/X}$. If that sieve contains $p$, then
it contains the entire groupoid $U_{\bigdot}$ (viewed as a groupoid in $\calX_{/X}$).
We conclude that there exists a groupoid object $W_{\bigdot}: \Nerve(\cDelta)^{op} \rightarrow \calX_{/q}$ lifting $U_{\bigdot}$. Let $\widetilde{V} \in \calX_{/q}$ be a colimit of
$V_{\bigdot}$. According to Proposition \ref{needed17}, the image of $\widetilde{V}$ in $\calX_{/X}$ can be identified with the map $p': V \rightarrow X$. The existence of $\widetilde{V}$ 
proves that $\bHom_{\calX{/X}}(p',q)$ is nonempty, as desired.
\end{proof}

\begin{corollary}\label{subobj}
Let $\calX$ be a semitopos, and let $f: U \rightarrow X$
be a morphism in $\calX$. The following conditions are equivalent:
\begin{itemize}
\item[$(1)$] If we regard $f$ as an object of the $\infty$-category $\calX_{/X}$, then
$\tau_{\leq -1}(f)$ is a final object of $\calX_{/X}$.
\item[$(2)$] The \Cech nerve $\mCech(f)$ is a simplicial resolution of $X$.
\end{itemize}
\end{corollary}

We will say that a morphism $f: U \rightarrow X$ in an semitopos $\calX$ is
an {\it effective epimorphism} if it satisfies the equivalent conditions of Corollary \ref{subobj}.\index{gen}{effective epimorphism} There is a one-to-one correspondence between effective epimorphisms and effective groupoids. More precisely, let $\Res_{\Eff}(\calX)$ denote the full subcategory of the $\infty$-category
$\calX_{\Delta_{+}}$ spanned by those augmented simplicial objects $U_{\bigdot}$ which
are both \Cech nerves and simplicial resolutions. The restriction functors
$$ \xymatrix{ & \calX_{\Delta_{+}} \ar[dr] \ar[dl] & \\
\calX_{\Delta} & & \Fun(\Delta^1, \calX) }$$
induce equivalences of $\infty$-categories from $\Res_{\Eff}(\calX)$ to the full subcategory
of $\calX_{\Delta}$ spanned by the effective groupoids, and from $\Res_{\Eff}(\calX)$ to the full subcategory of $\Fun(\Delta^1, \calX)$ spanned by the effective epimorphisms.

\begin{remark}\label{geoeff}
Let $f_{\ast}: \calX \rightarrow \calY$ be a geometric morphism of $\infty$-topoi, and let
$u: U \rightarrow Y$ be an effective epimorphism in $\calY$. Then $f^{\ast}(u)$ is an effective epimorphism in $\calX$. To see this, choose a \Cech nerve $U_{\bigdot}$ of $u$. Since
$u$ is an effective epimorphism, $U_{\bigdot}$ is a colimit diagram. The left exactness of $f^{\ast}$ implies that $f^{\ast} \circ U_{\bigdot}$ is a \Cech nerve of $f^{\ast}(u)$. Since $f^{\ast}$ is a left adjoint, we conclude that $f^{\ast} \circ U_{\bigdot}$ is a colimit diagram so that $f^{\ast}(u)$ is an effective epimorphism.
\end{remark}

The following result summarizes a few basic properties of effective epimorphisms:

\begin{proposition}\label{sinn}
Let $\calX$ be a semitopos.

\begin{itemize}
\item[$(1)$] Any equivalence $X \rightarrow Y$ in $\calX$ is an effective epimorphism.

\item[$(2)$] If $f,g: X \rightarrow Y$ are homotopic morphisms in $\calX$, then
$f$ is an effective epimorphism if and only if $g$ is an effective epimorphism.

\item[$(3)$] If $F: \calY \rightarrow \calX$ is a left exact presentable functor between semitopoi, and $f: U \rightarrow X$ is an effective epimorphism in $\calX$, then $F(f)$ is an effective epimorphism
in $\calY$.
\end{itemize}
\end{proposition}

\begin{proof}
Assertions $(1)$ and $(2)$ are obvious. To prove $(3)$, we observe that $f$ is an effective epimorphism if and only if it can be extended to an augmented simplicial object $U_{\bigdot}$ which is both a simplicial resolution and a \Cech nerve. Since $F$ is left exact, it preserves the property of being a \Cech nerve; since $F$ preserves colimits, it preserves the property of being a simplicial resolution.
\end{proof}

\begin{remark}
Let $\calX$ be a semitopos, and let $f: X \rightarrow T$ be an effective epimorphism in $\calX$
Applying part $(3)$ of Proposition \ref{sinn} to the geometric morphism $f: \calX_{/S} \rightarrow \calX_{/T}$ induced by a morphism $S \rightarrow T$ in $\calX$, we deduce that any base change
$X \times_{T} S \rightarrow X$ of $f$ is also an effective epimorphism.
\end{remark}

In order to verify other basic properties of the class of effective epimorphisms, such as stability under composition, we will need to reformulate the property of surjectivity in terms of subobjects.
Let $\calX$ be presentable $\infty$-category. For each $X \in \calX$, the $\infty$-category
$\tau_{\leq -1} \calX_{/X}$ of subobjects of $X$ is equivalent to the nerve of a partially
ordered set which we will denote by $\Sub(X)$; we may identify $\Sub(X)$ with the 
set of equivalence classes of monomorphisms $U \rightarrow X$.
A morphism $f: X \rightarrow Y$ in $\calX$ induces a left exact pullback functor $\calX_{/X} \rightarrow \calX_{/Y}$. This functor preserves $(-1)$-truncated objects by Proposition \ref{eaa}, and therefore induces a map $f^{\ast}: \Sub(Y) \rightarrow \Sub(X)$ of partially ordered sets. 

\begin{remark}\label{summep}
Let $\calX$ be a presentable $\infty$-category in which colimits are universal. Then any monomorphism $u: U \rightarrow \coprod X_{\alpha}$ can be obtained as a coproduct
of maps $u_{\alpha}: U_{\alpha} \rightarrow X_{\alpha}$, where each $u_{\alpha}$ is a pullback of $u$ and therefore also a monomorphism. It follows that the natural map
$$ \theta: \Sub( \coprod X_{\alpha} ) \rightarrow \prod \Sub(X_{\alpha})$$
is a monomorphism of partially ordered sets. If coproducts in $\calX$ are disjoint, then
$\theta$ is bijective.
\end{remark}

\begin{proposition}\label{charsurj}
Let $\calX$ be a semitopos.
A morphism $f: U \rightarrow X$ in $\calX$ is an effective epimorphism if and
only if $f^{\ast}: \Sub(X) \rightarrow \Sub(U)$ is injective.
\end{proposition}

\begin{proof}
Suppose first that $f^{\ast}$ is injective. Let $U_{\bigdot}$
be the underlying groupoid of a \Cech nerve of $f$, let
$V$ be a colimit of $U_{\bigdot}$, let $u: V \rightarrow X$
be the corresponding monomorphism, and let $[V]$ denote the corresponding element
of $\Sub(X)$. Since $f$ factors through $u$, we conclude that $f^{\ast}[V] = f^{\ast} [X] = [U] \in \Sub(U)$. Invoking the injectivity of $f^{\ast}$, we conclude that $[V] = [X]$ so that $u$ is an equivalence.

For the converse, let us suppose that $f$ is an effective epimorphism.
Let $[V]$ and $[V']$ be elements of $\Sub(X)$, represented by monomorphisms
$u: V \rightarrow X$ and $u': V' \rightarrow X$, and suppose that
$f^{\ast} [V] = f^{\ast} [V']$. We wish to prove that $[V] = [V']$. 
Since $f^{\ast}$ is a left exact functor, we have
$f^{\ast}( [V] \cap [V'] ) = f^{\ast} [ V \times_{X} V' ]$. It will suffice
to prove that $[V'] = [V \times_{X} V']$; the same argument will then establish
that $[V] = [V \times_{X} V']$ and the proof will be complete. In other words, we may assume without loss of generality that $[V] \subseteq [V']$ so that there is a commutative diagram
$$ \xymatrix{ V_0 \ar[dr]^{u} \ar[rr]^{g} & & V' \ar[dl]^{u'} \\
& X. & }$$
We wish to show that $g$ is an equivalence. The map $g$ induces
a natural transformation of augmented simplicial objects
$$ \alpha_{\bigdot}: u^{\ast} \circ \mCech(f) \rightarrow {u'}^{\ast} \circ \mCech(f).$$
We observe that $g$ can be identified with $\alpha_{-1}$. Since
$f$ is an effective epimorphism, $\mCech(f)$ is a colimit diagram.
Since colimits in $\calX$ are universal, we conclude that
$\alpha_{-1}$ is a colimit of $\alpha | \Nerve(\cDelta)^{op}$. Consequently, to
prove that $\alpha_{-1}$ is an equivalence, it will suffice to prove that $\alpha_n$
is an equivalence for $n \geq 0$. Since each $\alpha_n$ is a pullback of $\alpha_0$, it will
suffice to prove that $\alpha_0$ is an equivalence. But this is simply a reformulation
of the condition that $f^{\ast} [V] = f^{\ast} [V']$.
\end{proof}

From this we immediately deduce some corollaries.

\begin{corollary}\label{sumepi}
Let $\calX$ be a semitopos, and let $\{ f_{\alpha}: X_{\alpha} \rightarrow Y_{\alpha} \}$ be a $($small$)$ collection of effective epimorphisms in $\calX$. Then the induced map
$$ f: \coprod_{\alpha} X_{\alpha} \rightarrow \coprod_{\alpha} Y_{\alpha} $$
is an effective epimorphism.
\end{corollary}

\begin{proof}
Combine Proposition \ref{charsurj} with Remark \ref{summep}.
\end{proof}

\begin{corollary}\label{composite}
Let $\calX$ be a semitopos containing a diagram
$$ \xymatrix{ & Y \ar[dr]^{g} & \\
X \ar[ur]^{f} \ar[rr]^{h} & & Z. }$$
\begin{itemize}
\item[$(1)$] If $f$ and $g$ are effective epimorphisms, then so is $h$.
\item[$(2)$] If $h$ is an effective epimorphism, then so is $g$.
\end{itemize}
\end{corollary}

\begin{proof}
This follows immediately from Proposition \ref{charsurj}, and the observation
that we have an equality $f^{\ast} \circ g^{\ast} = h^{\ast}$ of functions
$\Sub(Z) \rightarrow \Sub(X)$.
\end{proof}

The theory of effective epimorphism is a mechanism for
proving theorems by descent.

\begin{lemma}\label{epie}
Let $\calX$ be a semitopos, let $\overline{p}: K^{\triangleright} \rightarrow \calX$ be
a colimit diagram, let $\infty$ denote the cone point of $K^{\triangleright}$. Then
the associated map
$$ \coprod_{x \in K_0} p(x) \rightarrow p(\infty)$$
$($which is well-defined up to homotopy$)$ is an effective epimorphism.
\end{lemma}

\begin{proof}
For each vertex $x$ of $K^{\triangleright}$, let $Z_x =p(x)$, and if $x$ belongs to $K$ we will 
denote the corresponding map $Z_{x} \rightarrow Z_{\infty}$ by $f_{x}$.
Let $E'' \subseteq E' \in \Sub(Z_{\infty})$ be such that
$f_{x}^{\ast}E'' = f_{x}^{\ast} E'$ for each vertex $x$ of $K$; we wish to show that $E'' = E'$. 
We can represent $E''$ and $E'$ by a $2$-simplex $\sigma_{\infty}: \Delta^2 \rightarrow \calX$, which we depict as
$$ \xymatrix{ & Z'_{\infty} \ar[dr] & \\ 
Z''_{\infty} \ar[ur] \ar[rr] & & Z_{\infty}. }$$
Lift the above diagram to a $2$-simplex $\sigma: \Delta^2 \rightarrow 
\Fun( K^{\triangleright}, \calX)$
$$ \xymatrix{ & p' \ar[dr]^{g''} & \\
p'' \ar[ur]^{g'} \ar[rr]^{g} & & p }$$
where $g$, $g'$, and $g''$ are Cartesian transformations. Our assumption guarantees that
the restriction of $g'$ induces an equivalence $p'' | K \rightarrow p' | K$. Since colimits in
$\calX$ are universal, $g'$ is itself an equivalence, so that $E'' = E'$ as desired.
\end{proof}

\begin{proposition}\label{torque}
Let $\calX$ be an $\infty$-topos, and let $S$ be a collection of morphisms of $\calX$ which is stable under pullbacks and coproducts. The following conditions are equivalent:
\begin{itemize}
\item[$(1)$] The class $S$ is local $($Definition \ref{localitie}$)$. 

\item[$(2)$] Given a pullback diagram
$$ \xymatrix{ X' \ar[r] \ar[d]^{f'} & X \ar[d]^{f} \\
Y' \ar[r]^{g} & Y }$$
where $g$ is an effective epimorphism and $f' \in S$, we have $f \in S$.
\end{itemize}
\end{proposition}

\begin{proof}
We first show that $(1) \Rightarrow (2)$. 
Let $Y_{\bigdot}: \Nerve(\cDelta_{+})^{op} \rightarrow \calX$ be a \Cech nerve of
the map $g$, and choose a Cartesian transformation $f_{\bigdot}: X_{\bigdot} \rightarrow Y_{\bigdot}$ of augmented simplicial objects which extends $f$. Then we can identify
$f'$ with $f_{0}: X_0 \rightarrow Y_0$. Each $f_n$ is a pullback of $f_0$, and therefore belongs to $S$. Applying Lemma \ref{ib2}, we deduce that $f$ belongs to $S$ as well.

Conversely, suppose that $(2)$ is satisfied. We will show that $S$ satisfies criterion $(3)$ of Lemma \ref{ib3}. Let $$ \xymatrix{ u \ar[r]^{\alpha} \ar[d]^{\beta} & v \ar[d]^{\beta'} \\
u' \ar[r]^{\alpha'} & v' }$$
be a pushout diagram in $\calO_{\calX}$, where $\alpha$ and $\beta$ are Cartesian
and $u,v,u' \in S$. Since $\calX$ is an $\infty$-topos, we conclude that $\alpha'$ and $\beta'$ are also Cartesian. To complete the proof, it will suffice to show that $v' \in S$. For this, we observe that
there is a pullback diagram
$$ \xymatrix{ X \amalg X' \ar[d]^{v \amalg u'} \ar[r] & X'' \ar[d]^{v'} \\
Y \amalg Y' \ar[r]^{g} & Y'' }$$
where $g$ is an effective epimorphism (Lemma \ref{epie}) and apply hypothesis $(2)$.
\end{proof}

\begin{proposition}\label{hintdescent1}
Let $\calX$ be a semitopos, and suppose given a pullback square
$$ \xymatrix{ X' \ar[r]^{g'} \ar[d]^{f'} & X \ar[d]^{f} \\
S' \ar[r]^{g} & S }$$
in $\calX$. If $f$ is an effective epimorphism, then so is $f'$. 
The converse holds if $g$ is an effective epimorphism.
\end{proposition}

\begin{proof}
Let $g^{\ast}: \calX^{/S} \rightarrow \calX^{/S'}$ be a pullback functor. Without loss of 
generality we may suppose that $f' = g^{\ast} f$. 
Let $U_{\bigdot}: \Nerve(\cDelta_{+})^{op} \rightarrow \calX$ be a \Cech nerve of $f$.
Since $g^{\ast}$ is left exact (being a right adjoint), we conclude that
$g^{\ast} \circ U_{\bigdot}$ is a \Cech nerve of $f'$. If $f$ is an effective epimorphism, then
$U_{\bigdot}$ is a colimit diagram. Because colimits in $\calX$ are universal, $g^{\ast} \circ U_{\bigdot}$ is also a colimit diagram, so that $f'$ is an effective epimorphism.

Conversely, suppose that $f'$ and $g$ are effective epimorphisms. Corollary \ref{composite} implies that $g \circ f'$ is an effective epimorphism. The commutativity of the diagram
implies that $f \circ g'$ is an effective epimorphism, so that $f$ is an effective epimorphism (Corollary \ref{composite} again).
\end{proof}

\begin{lemma}\label{hint0}
Let $\calX$ be a semitopos, and suppose given a pullback square
$$ \xymatrix{ X' \ar[r]^{g'} \ar[d]^{f'} & X \ar[d]^{f} \\
S' \ar[r]^{g} & S }$$
in $\calX$. Suppose that $f'$ is an equivalence that $g$ is an effective epimorphism.
Then $f$ is an equivalence.
\end{lemma}

\begin{proof}
Let $U_{\bigdot}$ be a \Cech nerve of $g'$, and let $V_{\bigdot}$ be a \Cech nerve of $g$.
The above diagram induces a transformation $\alpha_{\bigdot}: U_{\bigdot} \rightarrow V_{\bigdot}$. The map $\alpha_0$ can be identified with $f'$, and is therefore an equivalence.
For $n \geq 0$, $\alpha_n: U_n \rightarrow V_n$ is a pullback of $\alpha_0$, and therefore also an equivalence. Since $g$ is an effective epimorphism, $V_{\bigdot}$ is a colimit diagram.
Applying Proposition \ref{hintdescent1}, we conclude that $g'$ is also an effective epimorphism so that $U_{\bigdot}$ is a colimit diagram. It follows that $f = \alpha_{-1}$ is a colimit of
equivalences, and is therefore an equivalence.
\end{proof}

\begin{proposition}\label{hintdescent0}
Let $\calX$ be a semitopos, and suppose given a pullback square
$$ \xymatrix{ X' \ar[r] \ar[d]^{f'} & X \ar[d]^{f} \\
S' \ar[r]^{g} & S }$$
in $\calX$. If $f$ is $n$-truncated, then so is $f'$. The converse holds if $g$ is an effective epimorphism.
\end{proposition}

\begin{proof}
Let $g^{\ast}: \calX^{/S} \rightarrow \calX^{/S'}$ be a pullback functor. The first part of
$(1)$ asserts that $g^{\ast}$ carries $n$-truncated objects to $n$-truncated objects.
This follows immediately from Proposition \ref{eaa}, since $g^{\ast}$ is a right adjoint and therefore left exact. We will prove the converse in a slightly stronger form: if $i: U \rightarrow V$
is a morphism in $\calX^{/S}$ such that then $g^{\ast}(i)$ is an $n$-truncated morphism in
$\calX^{/S'}$, then $i$ is $n$-truncated. The proof is by induction on $n$. If $n \geq -1$, we
can use Lemma \ref{trunc} to reduce to the problem of showing that the diagonal map
$\delta: U \rightarrow U \times_{V} U$ is $(n-1)$-truncated. Since $g^{\ast}$ is left exact, 
we can identify $g^{\ast}(\delta)$ with the diagonal map $g^{\ast} U \rightarrow g^{\ast} U \times_{ g^{\ast} V} g^{\ast} U$, which is $(n-1)$-truncated according to Lemma \ref{trunc}; the desired result then follows from the inductive hypothesis. In the case $n=-2$, we have a pullback
diagram
$$ \xymatrix{ g^{\ast} U \ar[r] \ar[d]^{g^{\ast} i} & U \ar[d]^{i} \\
g^{\ast} V \ar[r]^{g'} \ar[d] & V \ar[d] \\
S' \ar[r]^{g} & S. }$$
Proposition \ref{hintdescent1} implies that $g'$ is an effective epimorphism, and
$g^{\ast} i$ is an equivalence, so that $i$ is also an equivalence by Lemma \ref{hint0}.
\end{proof}

Let $\calC$ be a small $\infty$-category equipped with a Grothendieck topology. 
Our final goal in this section is to use the language of effective epimorphisms to characterize
the $\infty$-topos $\Shv(\calC)$ by a universal property.

\begin{lemma}\label{stubba}
Let $\calC$ be a $($small$)$ $\infty$-category containing an object $C$, let
$\{ f_{\alpha}: C_{\alpha} \rightarrow C \}_{\alpha \in A}$ be a collection of morphisms indexed by a set $A$ and let $\calC^{(0)}_{/C} \subseteq \calC_{/C}$ be the sieve on $C$ that they generate.
Let $j: \calC \rightarrow \calP(\calC)$ denote the Yoneda embedding and $i: U \rightarrow j(C)$
a monomorphism corresponding to the sieve $\calC^{(0)}_{/C}$. Then $i$ can be identified with a $(-1)$-truncation of the induced map $\coprod_{\alpha \in A} j(C_{\alpha}) \rightarrow j(C)$
in the $\infty$-topos $\calP(\calC)_{/C}$.
\end{lemma}

\begin{proof}
Using Proposition \ref{surry}, we can identify the equivalence classes of $(-1)$-truncated object $U \in \calP(\calC)_{/j(C)}$ with sieves $\calC^{(0)}_{/C} \subseteq \calC_{/C}$. It is not difficult to see that $j( f_{\alpha})$ factors through $U$ if and only if $f_{\alpha} \in \calC^{(0)}_{/C}$. 
Consequently, the $(-1)$-truncation of $\coprod_{\alpha \in A} j(C_{\alpha}) \rightarrow j(C)$
is associated to the {\em smallest} sieve on $\calC$ which contains each $f_{\alpha}$.
\end{proof}

\begin{lemma}\label{pregrute}
Let $\calX$ be an $\infty$-topos, $\calC$ a small $\infty$-category equipped with a Grothendieck topology, and $f_{\ast}: \calX \rightarrow \calP(\calC)$ a functor with a left exact
left adjoint $f^{\ast}: \calP(\calC) \rightarrow \calX$.

The following conditions are equivalent:
\begin{itemize}
\item[$(1)$] The functor $f_{\ast}$ factors through $\Shv(\calC) \subseteq \calP(\calC)$.

\item[$(2)$] For every collection of morphisms $\{ v_{\alpha}: C_{\alpha} \rightarrow C\}$ which generate
a covering sieve in $\calC$, the induced map
$$ \coprod f^{\ast}( j(C_{\alpha}) ) \rightarrow f^{\ast} ( j(C) )$$
is an effective epimorphism in $\calX$, where $j: \calC \rightarrow \calP(\calC)$ denotes the Yoneda embedding.
\end{itemize}
\end{lemma}

\begin{proof}
Suppose first that $(1)$ is satisfied, and let $\{ v_{\alpha}: C_{\alpha} \rightarrow C\}$ be a collection of morphisms as in the statement of $(2)$. Let $L: \calP(\calC) \rightarrow \Shv(\calC)$ be a left adjoint to the inclusion. Then we have an equivalence of functors $f^{\ast} \simeq (f^{\ast}|\Shv(\calC)) \circ L$. Applying Remark \ref{geoeff}, we are reduced to showing that if
$$ u: \coprod j(C_{\alpha}) \rightarrow j(C)$$
is the natural map, then $Lu$ is an effective epimorphism in $\calP(\calC)$. 
Factor $u$ as a composition
$$ \coprod j(C_{\alpha}) \stackrel{u'}{\rightarrow} U \stackrel{u''}{\rightarrow} j(C)$$
where $u'$ is an effective epimorphism and $u''$ is a monomorphism. We wish to show that
$Lu''$ is an equivalence. Lemma \ref{stubba} allows us to identify $u''$ with the monomorphism
associated to the sieve $\calC^{(0)}_{/C}$ on $C$ generated by the maps $v_{\alpha}$. By assumption, this is a covering sieve, so that $L u''$ is an equivalence in $\Shv(\calC)$ by construction.

Conversely, suppose that $(2)$ is satisfied. Let $C \in \calC$ and let $\calC^{(0)}_{/C} \subseteq \calC_{/C}$ be a covering sieve on $C$ associated to a monomorphism
$u'': U \rightarrow j(C)$. We wish to show that $f^{\ast} u''$ is an equivalence.
According to Lemma \ref{stubba}, we have a factorization
$$ \coprod_{\alpha} j(C_{\alpha}) \stackrel{u'}{\rightarrow} U \stackrel{u''}{\rightarrow} j(C), $$
where the maps $v_{\alpha}: C_{\alpha} \rightarrow C$ are chosen to generate the sieve
$\calC^{(0)}_{/C}$, and $u'$ is an effective epimorphism. Let $u$ be a composition of $u'$ and $u''$. Then $f^{\ast} u'$ is an effective epimorphism (Remark \ref{geoeff}), and $f^{\ast} u$ is an effective epimorphism by assumption $(2)$. Corollary \ref{composite} now shows that
$f^{\ast} u''$ is an effective epimorphism. Since $f^{\ast} u''$ is also a monomorphism, we conclude that $f^{\ast} u''$ is an equivalence as desired.
\end{proof}

\begin{proposition}\label{igrute}
Let $\calX$ be an $\infty$-topos, and let $\calC$ be a small $\infty$-category equipped with a Grothendieck topology. Let
$L: \calP(\calC) \rightarrow \Shv(\calC)$ denote a left adjoint to the inclusion, and
$j: \calC \rightarrow \calP(\calC)$ the Yoneda embedding. Let
$\Fun^{\ast}( \Shv(\calC), \calX )$ denote the $\infty$-category of left exact, colimit-preserving functors from $\Shv(\calC)$ to $\calX$ (Definition \ref{defhomst}). The composition
$$ J: \Fun^{\ast}(\Shv(\calC), \calX) \stackrel{L}{\rightarrow} \Fun^{\ast}( \calP(\calC), \calX)
\stackrel{j}{\rightarrow} \Fun( \calC, \calX)$$
is fully faithful. Suppose furthermore that $\calC$ admits finite limits. Then
a functor $f: \calC \rightarrow \calX$ belongs to the essential image of $J$ if and only if
the following conditions are satisfied:
\begin{itemize}
\item[$(1)$] The functor $f$ is left exact.
\item[$(2)$] For every collection of morphisms $\{ C_{\alpha} \rightarrow C\}_{\alpha \in A}$ which generate a covering sieve on $C$, the associated morphism
$$ \coprod_{\alpha \in A} f(C_{\alpha}) \rightarrow f(C)$$ is an effective epimorphism in $\calX$.
\end{itemize}
\end{proposition}

\begin{proof}
If the topology on $\calC$ is trivial, then Theorem \ref{charpresheaf} implies that
$J$ is fully faithful, and the description of the essential image of $J$ follows from Proposition \ref{natash}. In the general case, Proposition \ref{unichar} implies that
composition with $L$ induces a fully faithful embedding
$$ J': \Fun^{\ast}( \Shv(\calC) ,\calX) \rightarrow \Fun^{\ast}( \calP(\calC), \calX),$$
so that $J$ is a composition of $J'$ with a fully faithful functor
$$ J'': \Fun^{\ast}(\calP(\calC),\calX) \rightarrow \Fun(\calC, \calX).$$ 
Suppose that $\calC$ admits finite limits and that $f$ satisfies $(1)$, so that
$f$ is equivalent to $J''(u^{\ast})$ for some left exact, colimit preserving
$u^{\ast}: \calP(\calC) \rightarrow \calX$. The functor $u^{\ast}$ is unique up to equivalence,
and Lemma \ref{pregrute} ensures that $u^{\ast}$ belongs to the essential image of $J'$ if
and only if condition $(2)$ is satisfied.
\end{proof}

\begin{remark}
It is possible to formulate a generalization of Proposition \ref{igrute} which describes the essential image of $J$ even when $\calC$ does not admit finite limits. The present version will be sufficient for the applications in this book. 
\end{remark}

\subsection{Canonical Topologies}\label{cantopp}

Let $\calX$ be an $\infty$-topos. Suppose that we wish to identify $\calX$ with an $\infty$-category of sheaves. The first step is to choose a pair of adjoint functors
$$ \Adjoint{F}{\calP(\calC)}{\calX}{G}$$
where $F$ is left exact. According to Theorem \ref{charpresheaf}, $F$ is determined up
to equivalence by the composition
$$ f: \calC \stackrel{j}{\rightarrow} \calP(\calC) \stackrel{F}{\rightarrow} \calX.$$
We might then try to choose a topology on $\calC$ such that $G$ factors as a composition
$$ \calX \stackrel{G'}{\rightarrow} \Shv(\calC) \subseteq \calP(\calC).$$
Though it is not always possible to guarantee that $G'$ is an equivalence, we will show that for an appropriately chosen topology (Definition \ref{defcantop}), the $\infty$-topos $\Shv(\calC)$ is a close approximation to $\calX$ (Proposition \ref{preciselate}). 

\begin{definition}\label{defcantop}\index{gen}{canonical covering}
Let $\calX$ be a semitopos, $\calC$ a small $\infty$-category which admits finite limits, and
$f: \calC \rightarrow \calX$ a left exact functor. We will say that a sieve $\calC^{(0)}_{/C} \subseteq \calC_{/C}$ on an object $C \in \calC$ is a {\it canonical covering relative to $f$} if there exists
a collection of morphisms $\{ u_{\alpha}: C_{\alpha} \rightarrow C \}$ belonging to $\calC^{(0)}_{/C}$ such that the induced map $ \coprod f(C_{\alpha}) \rightarrow f(C)$ is an effective epimorphism in $\calX$.
\end{definition}

Our first goal is to verify that the canonical topology is actually a Grothendieck topology on $\calC$.

\begin{proposition}\label{cantop}\index{gen}{canonical topology}
Let $f: \calC \rightarrow \calX$ be as in Definition \ref{defcantop}. 
The collection of canonical coverings relative to $f$ determine a Grothendieck topology on $\calC$.
\end{proposition}

\begin{proof}
Since any identity map $\id_{f(C)}: f(C) \rightarrow f(C)$ is an effective epimorphism, it is clear that
the sieve $\calC_{/C}$ is a canonical covering of $C$ for every $C \in \calC$. Suppose
that $\calC_{/C}^{(0)} \subseteq \calC_{/C}$ is a canonical covering of $C$, and that
$g: D \rightarrow C$ is a morphism in $\calC$. We wish to prove that
the induced sieve $g^{\ast} \calC_{/C}^{(0)}$ is a canonical covering. Choose a 
collection of objects $u_{\alpha}: C_{\alpha} \rightarrow C$ of $\calC_{/C}^{(0)}$ 
such that the induced map $\coprod_{\alpha} f(C_{\alpha}) \rightarrow f(C)$ is an
effective epimorphism, and form pullback diagrams
$$ \xymatrix{ D_{\alpha} \ar[r]^{v_{\alpha}} \ar[d] & D \ar[d]^{g} \\
C_{\alpha} \ar[r]^{u_{\alpha}} & C }$$
in $\calC$. Using the fact that $f$ is left exact and that colimits in $\calX$ are universal, we conclude that the diagram
$$ \xymatrix{ \coprod f(D_{\alpha}) \ar[r] \ar[d] & f(D) \ar[d] \\
\coprod f(C_{\alpha}) \ar[r] & f(C) }$$
is a pullback, so that the upper  horizontal map is an effective epimorphism by Proposition \ref{hintdescent1}. Since each $v_{\alpha}$ belongs to $g^{\ast} \calC_{/C}^{(0)}$, it follows that
$g^{\ast} \calC_{/C}^{(0)}$ is a canonical covering.

Now suppose that $\calC_{/C}^{(0)}$ and $\calC_{/C}^{(1)}$ are sieves on $C \in \calC$, where
$\calC_{/C}^{(0)}$ is a canonical covering, and for each $g: D \rightarrow C$ in $\calC_{/C}^{(0)}$, the covering $g^{\ast} \calC_{/C}^{(1)}$ is a canonical covering of $D$. Choose a collection
of morphisms $g_{\alpha}: D_{\alpha} \rightarrow C$ belonging to $\calC_{/C}^{(0)}$ with the property that $\coprod f(D_{\alpha}) \rightarrow f(C)$ is an effective epimorphism. For each
$D_{\alpha}$, choose a collection of morphisms $h_{\alpha,\beta}: E_{\alpha,\beta} \rightarrow D_{\alpha}$ belonging to $g_{\alpha}^{\ast} \calC_{/C}^{(1)}$ such that the map
$\coprod_{\beta} f(E_{\alpha,\beta}) \rightarrow f(D_{\alpha})$ is an effective epimorphism.
Using Corollary \ref{sumepi}, we conclude that the map
$$ \coprod_{\alpha,\beta} f(E_{\alpha,\beta}) \rightarrow \coprod_{\alpha} f(D_{\alpha})$$ is an effective epimorphism. Since effective epimorphisms are stable under composition (Corollary \ref{composite}), we have an effective epimorphism $\coprod_{\alpha,\beta} f(E_{\alpha,\beta}) \rightarrow f(C)$, induced by the collection of compositions $g_{\alpha} \circ h_{\alpha,\beta}: E_{\alpha,\beta} \rightarrow C$. Each of these compositions belong to $\calC_{/C}^{(1)}$, so
that $\calC_{/C}^{(1)}$ is a canonical covering of $C$.
\end{proof}

For later use, we record a few features of the canonical topology:

\begin{lemma}\label{caninit}
Let $f: \calC \rightarrow \calX$ be as in Definition \ref{defcantop}, and regard
$\calC$ as endowed with the canonical topology relative to $f$. Let $j: \calC \rightarrow \calP(\calC)$ denote the Yoneda embedding and let $L: \calP(\calC) \rightarrow \Shv(\calC)$
be a left adjoint to the inclusion. Suppose that $C \in \calC$ is such that $f(C)$ is an initial object of $\calX$. Then $L j(C)$ is an initial object of $\Shv(\calC)$.
\end{lemma}

\begin{proof}
If $f(C)$ is an initial object of $\calX$, then the empty sieve $\emptyset \subseteq \calC_{/C}$
is a covering sieve with respect to the canonical topology. By construction, the associated monomorphism $\emptyset \rightarrow j(C)$ becomes an equivalence after applying $L$, so that $L j(C)$ is initial in $\Shv(\calC)$.
\end{proof}

\begin{lemma}\label{canonicalcoproducts}
Let $f: \calC \rightarrow \calX$ be as in Definition \ref{defcantop}. Suppose that $f$ is fully faithful, coproducts in $\calX$ are disjoint, and let $\{ u_{\alpha}: C_{\alpha} \rightarrow C \}$ be a small collection of morphisms in $\calC$ such that the morphisms $f(u_{\alpha})$ exhibit
$f(C)$ as a coproduct of the family $\{ f(C_{\alpha}) \}$. Let $\calF: \calC^{op} \rightarrow \SSet$
be a sheaf on $\calC$ $($with respect to the canonical topology induced by $f${}$)$. Then the morphisms
$\{ \calF(u_{\alpha}) \}$ exhibit $\calF(C)$ as a product of $\{ \calF(C_{\alpha}) \}$ in $\SSet$.
\end{lemma}

\begin{proof}
We wish to show that the natural map $\calF(C) \rightarrow \prod \calF(C_{\alpha} )$ is an isomorphism in the homotopy category $\calH$. We may identify the left hand side with
$\bHom_{\calP(\calC)}( j(C), \calF)$, and the right hand side with
$\bHom_{\calP(\calC)}( \coprod j(C_{\alpha}), \calF)$. Consequently, it will suffice to show that the natural map
$$ v: \coprod j(C_{\alpha}) \rightarrow j(C)$$
becomes an equivalence after applying the localization functor $L: \calP(\calC) \rightarrow \Shv(\calC)$. Choose a factorization of $v$ as a composite
$$ \coprod j(C_{\alpha}) \stackrel{v'}{\rightarrow} U \stackrel{v''}{\rightarrow} j(C)$$
where $v'$ is an effective epimorphism, and $v''$ is a monomorphism. We observe that 
$v''$ is the monomorphism associated to the
sieve $\calC^{(0)}_{/C} \rightarrow \calC$ generated by the morphisms $u_{\alpha}$.
This is clearly a covering sieve with respect to the canonical topology, so that $Lv''$
is an equivalence in $\Shv(\calC)$. It follows that $Lv$ is equivalent to $Lv'$, and is therefore an effective epimorphism (Remark \ref{geoeff}). 
Form a pullback diagram
$$ \xymatrix{ V \ar[r]^{\overline{v}} \ar[d] & \coprod j(C_{\beta}) \ar[d]^{v} \\
\coprod j(C_{\alpha}) \ar[r]^{v} & j(C) }$$
We wish to prove that $Lv$ is an equivalence. According to Lemma \ref{hint0}, it will suffice
to show that $L \overline{v}$ is an equivalence. Since colimits in $\calP(\calC)$ are universal, we may identify $\overline{v}$ with a coproduct of morphisms 
$$\overline{v}_{\beta}: V_{\beta} \rightarrow j(C_{\beta}),$$ 
where $V_{\beta}$ can be written as a coproduct $\coprod_{\alpha} j( C_{\alpha} \times_{C} C_{\beta})$. Using Lemma \ref{sumdescription}, we can identify the summand
$ j( C_{\beta} \times_{C} C_{\beta})$ of $V_{\beta}$ with $j(C_{\beta})$, and the restriction
of $\overline{v}_{\beta}$ to this summand is an equivalence. To complete the proof, it will suffice to show that for every other summand $D_{\alpha, \beta} = j( C_{\alpha} \times_{C} C_{\beta} )$,
the localization $LD$ is an initial object of $\Shv(\calC)$. To prove this, we observe Lemma \ref{sumdescription} implies that $f( C_{\alpha} \times_{C} C_{\beta} )$ is an initial object
of $\calX$, and apply Lemma \ref{caninit}.
\end{proof}

\begin{lemma}\label{charepii}
Let $\calC$ be a small $\infty$-category equipped with a Grothendieck topology, and let
$u: \calF' \rightarrow \calF$ be a morphism in $\Shv(\calC)$. Suppose that, for each
$C \in \calC$ and each $\eta \in \pi_0 \calF(C)$, there exists a collection of morphisms
$\{ C_{\alpha} \rightarrow C \}$ which generates a covering sieve on $C$ and
a collection of $\eta_{\alpha} \in \pi_0 \calF'(C_{\alpha})$ such that $\eta$ and
$\eta_{\alpha}$ have the same image in $\pi_0 \calF(C_{\alpha})$. Then $u$ is an effective epimorphism.
\end{lemma}

\begin{proof}
Replacing $\calF$ by its image in $\calF'$ if necessary, we may suppose that
$u$ is a monomorphism. Let $L: \calP(\calC) \rightarrow \Shv(\calC)$ be a left adjoint to the inclusion, and let $\calD$ be the full subcategory of $\calP(\calC)$ spanned by those objects
$\calG$ such that, for every pullback diagram
$$ \xymatrix{ \calG' \ar[r]^{u'} \ar[d] & \calG \ar[d] \\
\calF' \ar[r]^{u} & \calF }$$
in $\calP(\calC)$, $Lu'$ is an equivalence in $\Shv(\calC)$. To prove that $u$ is an equivalence, it will suffice to show that the equivalent morphism $Lu$ is an equivalence. For this, it will suffice to prove that $\calF \in \calD$. We will in fact prove that $\calD = \calP(\calC)$. We first observe
that, since colimits in $\calP(\calC)$ are universal and $L$ commutes with colimits, 
$\calD$ is stable under colimits in $\calP(\calC)$. Since $\calP(\calC)$ is generated under colimits by the image of the Yoneda embedding, it will suffice to prove that $j(C) \in \calD$, for each $C \in \calC$. Choose a map $j(C) \rightarrow \calF$, classified up to homotopy by $\eta \in \pi_0 \calF(C)$, and form a pullback diagram 
$$ \xymatrix{ U \ar[r]^{u'} \ar[d] & j(C) \ar[d] \\
\calF' \ar[r]^{u} & \calF }$$
as above. Then $u'$ is a monomorphism; according to Proposition \ref{surry} it is classified
by a sieve $\calC^{(0)}_{/C}$ on $\calC$. Our hypothesis guarantees that $\calC^{(0)}_{/C}$ contains a collection of morphisms $\{ C_{\alpha} \rightarrow C \}$ which generate a covering sieve, so that $\calC^{(0)}_{/C}$ is itself covering. It follows immediately from the construction of $\Shv(\calC)$ that $Lu'$ is an equivalence.
\end{proof}

We close with the following result, which implies that any $\infty$-topos is closely approximated by an $\infty$-category of sheaves.

\begin{proposition}\label{preciselate}
Let $\calX$ be a semitopos, $\calC$ a small $\infty$-category which admits finite limits, and
$$ \Adjoint{F}{\calP(\calC)}{\calX}{G}$$ a pair of adjoint functors. Suppose that the composition
$$ f: \calC \stackrel{j}{\rightarrow} \calP(\calC) \stackrel{F}{\rightarrow} \calX$$
is left exact, and regard $\calC$ as endowed with the canonical topology relative to $f$.
Then:

\begin{itemize}
\item[$(1)$] The functor $G$ factors through $\Shv(\calC)$.

\item[$(2)$] Suppose that $f$ is fully faithful and generates $\calX$ under colimits. Then $G$ carries effective epimorphisms in $\calX$ to effective epimorphisms in $\Shv(\calC)$. 

\end{itemize}
\end{proposition}

\begin{proof}
In view of the definition of the canonical topology, $(1)$ is equivalent to the following assertion: given a collection of morphisms $\{ u_{\alpha}: C_{\alpha} \rightarrow C \}$ in
$\calC$ such that the induced map $u: \coprod_{\alpha} C_{\alpha} \rightarrow C$ is an effective epimorphism in $\calX$, if $i: U \rightarrow j(C)$ is the monomorphism in $\calP(\calC)$
corresponding to the sieve $\calC^{(0)}_{/C} \subseteq \calC_{/C}$ generated by
the collection $\{ u_{\alpha} \}$, then $F(i)$ is an equivalence in $\calX$.
Let $u': \coprod_{\alpha} j(C_{\alpha}) \rightarrow j(C)$ be the coproduct of the family
$\{ j( u_{\alpha} ) \}$, and let $V_{\bigdot}: \cDelta_{+}^{op} \rightarrow \calP(\calC)$ be a \Cech nerve of $u'$. Then $i$ can be identified with the induced map from the colimit of 
$V_{\bigdot} | \Nerve(\cDelta)^{op}$ to $V_{-1}$. Since $F$ preserves colimits, to show that
$F(i)$ is an equivalence, it will suffice to show that $F \circ V_{\bigdot}$ is a colimit diagram.
Since $u$ is an effective epimorphism, it suffices to observe that $F \circ V_{\bigdot}$ is equivalent to the \Cech nerve of $u$.

We now prove $(2)$. Suppose that $u: Y \rightarrow Z$ is an effective epimorphism in $\calX$.
We wish to prove that $Gu$ is an effective epimorphism in $\Shv(\calC)$. We will show that the criterion of Lemma \ref{charepii} is satisfied. Choose an object $C \in \calC$ and a point
$\eta \in \pi_0 \bHom_{\calP(\calC)}( j(C), GZ) \simeq \pi_0 \bHom_{\calX}(f(C), Z)$. Form a pullback diagram
$$ \xymatrix{ Y' \ar[r]^{u'} \ar[d]^{s} & f(C) \ar[d] \\
Y \ar[r]^{u} & Z }$$
so that $u'$ is an effective epimorphism. Since $f(\calC)$ generates $\calX$ under colimits, there
exists an effective epimorphism $u'': \coprod_{\alpha} f(C_{\alpha}) \rightarrow Y$. The composition
$u' \circ u''$ is an effective epimorphism, and corresponds to a family of maps
$w_{\alpha}: f(C_{\alpha}) \rightarrow f(C)$ in $\calX$. Since $f$ is fully faithful, we may suppose
that each $w_{\alpha} = fv_{\alpha}$ for some map $v_{\alpha}: C_{\alpha} \rightarrow C$ in $\calC$. It follows that the collection of maps $\{ v_{\alpha} \}$ generate a covering sieve on $C$ with respect to the canonical topology. Moreover, each of the compositions
$$f(C_{\alpha}) \rightarrow \coprod_{\alpha} f(C_{\alpha}) \rightarrow Y $$
gives rise to a point $\eta_{\alpha} \in \pi_0 \bHom_{\calX}( f(C_{\alpha}), Y) \simeq
\pi_0 \bHom_{\calP(\calC)}( j( C_{\alpha}), G(Y))$ with the desired properties.
\end{proof}
\section{The $\infty$-Category of $\infty$-Topoi}\label{chap6sec4}
\setcounter{theorem}{0}

In this section, we will show that the collection of all $\infty$-topoi can be organized into an $\infty$-category $\RGeom$. The objects of $\RGeom$ are $\infty$-topoi, and the morphisms are called {\it geometric morphisms}; we will give a definition in \S \ref{gemor1}. In \S \ref{colimtop}, we will show that $\RGeom$ admits (small) colimits. In \S \ref{inftyfiltlim}, we will show that $\RGeom$ admits 
(small) {\em filtered} limits; we will treat the case of general limits in \S \ref{genlim}.

Let $\calX$ be an $\infty$-topos containing an object $U$. In \S \ref{gemor2}, we will show that the 
$\infty$-category $\calX_{/U}$ is an $\infty$-topos. Moreover, this $\infty$-topos is equipped with
a canonical geometric morphism $\calX_{/U} \rightarrow \calX$. Geometric morphisms which arise via this construction are said to be {\em \'{e}tale}. In \S \ref{structuretheor}, we will define a more general notion of {\em algebraic} morphism between $\infty$-topoi. We will also prove a structure theorem which implies that, under certain hypotheses, every $\infty$-topos $\calX$ admits an algebraic morphism to an $\infty$-category of sheaves on a $2$-category.

\subsection{Geometric Morphisms}\label{gemor1}

In classical topos theory, the correct notion of {\em morphism} between two topoi is an adjunction
$$ \Adjoint{f^{\ast}}{\calX}{\calY}{f_{\ast}}$$
where the functor $f^{\ast}$ is left exact. We will introduce the same ideas in the $\infty$-categorical setting.

\begin{definition}\label{geomorph}\index{gen}{geometric morphism}
Let $\calX$ and $\calY$ be $\infty$-topoi.
A {\it geometric morphism} from $\calX$ to $\calY$ is a functor
$f_{\ast}: \calX \rightarrow \calY$ which admits a left exact left adjoint (which we will typically denote by $f^{\ast}$).
\end{definition}

\begin{remark}\index{gen}{pullback functor}\index{gen}{pushforward functor}
Let $f_{\ast}: \calX \rightarrow \calY$ be a geometric morphism from an $\infty$-topos $\calX$ to another $\infty$-topos $\calY$, so that $f_{\ast}$ admits a left adjoint $f^{\ast}$. 
Either functor $f_{\ast}$ and $f^{\ast}$ determines the other up to equivalence (in fact, up to contractible ambiguity). We will often abuse terminology by referring to $f^{\ast}$ as a geometric morphism from $\calX$ to $\calY$. We will always indicate in our notation whether the left or right adjoint is being considered: a superscripted asterisk indicates a left adjoint (pullback functor), and a subscripted asterisk indicates a right adjoint (pushforward functor).
\end{remark}

\begin{remark}
Any equivalence of $\infty$-topoi is a geometric morphism. If $f_{\ast}, g_{\ast}: \calX \rightarrow \calY$ are homotopic, then $f_{\ast}$ is a geometric morphism if and only if $g_{\ast}$ is a geometric morphism (because we can identify left adjoints of $f_{\ast}$ with left adjoints of $g_{\ast}$). 
\end{remark}

\begin{remark}
Let $f_{\ast}: \calX \rightarrow \calY$ and $g_{\ast}: \calY \rightarrow \calZ$ be geometric morphisms. Then $f_{\ast}$ and $g_{\ast}$ admit left exact left adjoints, which we will denote by
$f^{\ast}$ and $g^{\ast}$, respectively. The composite functor $f^{\ast} \circ g^{\ast}$ is left
exact, and is a left adjoint to $g_{\ast} \circ f_{\ast}$ by Proposition \ref{compadjoint}. We conclude that $g_{\ast} \circ f_{\ast}$ is a geometric morphism, so the class of geometric morphisms is stable under composition.
\end{remark}

\begin{definition}\index{not}{LTop@$\LGeom$}\index{not}{RTop@$\RGeom$}
Let $\widehat{\Cat}_{\infty}$ denote the $\infty$-category of (not necessarily small) $\infty$-categories. We define subcategories $\LGeom, \RGeom \subseteq \widehat{\Cat}_{\infty}$ as follows:
\begin{itemize}
\item[$(1)$] The objects of $\LGeom$ and $\RGeom$ are the $\infty$-topoi.
\item[$(2)$] A functor $f^{\ast}: \calX \rightarrow \calY$ between $\infty$-topoi belongs to
$\LGeom$ if and only if $f^{\ast}$ preserves small colimits and finite limits.
\item[$(3)$] A functor $f_{\ast}: \calX \rightarrow \calY$ between $\infty$-topoi belongs to
$\RGeom$ if and only if $f_{\ast}$ has a left adjoint which is left exact.
\end{itemize}
\end{definition}

The $\infty$-categories $\LGeom$ and $\RGeom$ are canonically anti-equivalent. To prove this, we will use the argument of Corollary \ref{warhog}. First, we need a definition.

\begin{definition}\index{gen}{fibration!topos}\label{skuzz}
A map $p: X \rightarrow S$ of simplicial sets is a {\it topos fibration} if the following conditions are satisfied:
\begin{itemize}
\item[$(1)$] The map $p$ is both a Cartesian fibration and a coCartesian fibration.
\item[$(2)$] For every vertex $s$ of $S$, the corresponding fiber $X_{s} = X \times_{S} \{s\}$ is an $\infty$-topos.
\item[$(3)$] For every edge $e: s \rightarrow s'$ in $S$, the associated functor
$X_{s} \rightarrow X_{s'}$ is left exact.
\end{itemize}
\end{definition}

The following analogue of Proposition \ref{surtog} follows immediately from the definitions:

\begin{proposition}\label{surtog2}
\begin{itemize}
\item[$(1)$] Let $p: X \rightarrow S$ be a Cartesian fibration of simplicial sets, classified
by a map $\chi: S^{op} \rightarrow \widehat{\Cat}_{\infty}$. Then $p$ is a topos fibration if and only if $\chi$ factors through $\RGeom \subseteq \widehat{\Cat}_{\infty}$. 

\item[$(2)$] Let $p: X \rightarrow S$ be a coCartesian fibration of simplicial sets, classified by a map $\chi: S \rightarrow \widehat{\Cat}_{\infty}$. Then $p$ is a topos fibration if and only if
$\chi$ factors through $\LGeom \subseteq \widehat{\Cat}_{\infty}$.
\end{itemize}
\end{proposition}

\begin{corollary}\label{suytoy}
For every simplicial set $S$, there is a canonical bijection
$$ [ S, \LGeom ] \simeq [ S^{op}, \RGeom ]$$
where $[K, \calC]$ denotes the collection of equivalence classes of objects of $\Fun(K,\calC)$.
In particular, $\LGeom$ and $\RGeom^{op}$ are canonically isomorphic
in the homotopy category of $\infty$-categories.
\end{corollary}

\begin{proof}
According to Proposition \ref{surtog2}, both $[S, \LGeom]$ and $[S^{op}, \RGeom]$ can be identified with the collection of equivalence classes of topos fibrations $X \rightarrow S$.
\end{proof}

The following proposition is a simple reformulation of some of the results of \S \ref{truncintro}.

\begin{proposition}
Let $f_{\ast}: \calX \rightarrow \calY$ be a geometric morphism between $\infty$-topoi, having a left adjoint $f^{\ast}$. Then $f^{\ast}$ and
$f_{\ast}$ carry $k$-truncated objects to $k$-truncated objects and $k$-truncated morphisms to $k$-truncated morphisms, for any integer $k \geq -2$. Moreover, there is a $($canonical$)$ equivalence
of functors $f^{\ast} \tau^{\calY}_{\leq k} \simeq \tau^{\calX}_{\leq k} f^{\ast}$.
\end{proposition}

\begin{proof}
The first assertion follows immediately from Lemma \ref{trunc}, since $f_{\ast}$ and $f^{\ast}$ are both left-exact. The second follows from Proposition \ref{compattrunc}.
\end{proof}

\begin{definition}\index{not}{FunLast@$\Fun_{\ast}(\calX,\calY)$}\index{not}{FunUast@$\Fun^{\ast}(\calX,\calY)$}\label{defhomst}
Let $\calX$ and $\calY$ be $\infty$-topoi. We let
$\Fun_{\ast}(\calX, \calY)$ denote the full subcategory of $\Fun(\calX,\calY)$ spanned by
geometric morphisms $f_{\ast}: \calX \rightarrow \calY$, and
$\Fun^{\ast}(\calY, \calX)$ the full subcategory of $\Fun(\calY,\calX)$ spanned by
their left adjoints.
\end{definition}

\begin{remark}
According to Proposition \ref{switcheroo}, the $\infty$-categories
$\Fun_{\ast}(\calX, \calY)$ and $\Fun^{\ast}(\calY, \calX)$ are canonically anti-equivalent to one another.
\end{remark}

\begin{warning}\label{tooobig}
If $\calX$ and $\calY$ are $\infty$-topoi, then the $\infty$-category $\Fun_{\ast}(\calX, \calY)$ of
geometric morphisms from $\calX$ to $\calY$ is {\em not} necessarily small, or even equivalent to a small $\infty$-category. This phenomenon is familiar in classical topos theory. For example,
there is a classifying topos $\calA$ for abelian groups, having the property that for {\em any} topos
$\calX$, the category $\calC$ of geometric morphisms $\calX \rightarrow \calA$ is equivalent to the category of abelian group objects of $\calX$. This category is almost never small (for example, when $\calX$ is the topos of sets, $\calC$ is equivalent to the category of abelian groups).
\end{warning}

In spite of Warning \ref{tooobig}, the $\infty$-category of geometric morphisms between two $\infty$-topoi can be reasonably controlled:

\begin{proposition}\label{nottoobig}
Let $\calX$ and $\calY$ be $\infty$-topoi. Then the
$\infty$-category $\Fun^{\ast}(\calY, \calX)$ of geometric morphisms from
$\calX$ to $\calY$ is accessible.
\end{proposition}

\begin{proof}
For each regular cardinal $\kappa$, let $\calY^{\kappa}$ denote the full subcategory of
$\calY$ spanned by $\kappa$-compact objects. Choose a regular cardinal $\kappa$ such that $\calY$ is $\kappa$-accessible and $\calY^{\kappa}$ is stable under finite limits in $\calY$. We may therefore identify $\calY$ with $\Ind^{\kappa}(\calC)$, where
$\calC$ is a minimal model for $\calY^{\kappa}$. According to Proposition \ref{intprop},
composition with the Yoneda embedding $j: \calC \rightarrow \calY$ induces an equivalence
from the $\infty$-category of $\kappa$-continuous functors $\Fun_{\kappa}(\calY,\calX)$
to the $\infty$-category $\Fun(\calC,\calX)$. We now make the following observations:

\begin{itemize}
\item[$(1)$] A functor $F: \calY \rightarrow \calX$ preserves all small colimits if and only if
$F \circ j: \calC \rightarrow \calX$ preserves $\kappa$-small colimits (Proposition \ref{sumatch}).

\item[$(2)$] A colimit-preserving functor $F: \calY \rightarrow \calX$ is left exact if and only if the composition $F \circ j: \calC \rightarrow \calX$ is left exact (Proposition \ref{natash}). 

\end{itemize}

Invoking Proposition \ref{switcheroo}, we deduce that the $\infty$-category $\Fun^{\ast}(\calY, \calX)$ is equivalent to the full subcategory $\calM \subseteq \calX^{\calC}$ consisting of functors which preserve $\kappa$-small colimits and finite limits. Proposition \ref{horse1} implies that $\Fun(\calC,\calX)$ is accessible. For every $\kappa$-small (finite) diagram $p: K \rightarrow \calC$, the full subcategory of $\Fun(\calC,\calX)$ which
preserve colimits (limits) of $p$ is an accessible subcategory of $\Fun(\calC,\calX)$
(Example \ref{colexam}). Up to isomorphism, there are only a bounded number of $\kappa$-small (finite) diagrams in $\calC$. Consequently, $\calM$ is an intersection of a bounded number of accessible subcategories of $\Fun(\calC,\calX)$, and therefore accessible by (Proposition \ref{boundint}).
\end{proof}

\subsection{Colimits of $\infty$-Topoi}\label{colimtop}

Our goal in this section is to construct colimits in the $\infty$-category $\RGeom$ of $\infty$-topoi. 
According to Corollary \ref{suytoy}, it will suffice construct {\em limits} in the $\infty$-category $\LGeom$. 

\begin{proposition}\label{quathorse}\index{gen}{coproduct!of $\infty$-topoi}
Let $\{ \calX_{\alpha} \}_{\alpha \in A}$ be a collection of $\infty$-topoi, parametrized by a $($small$)$
set $A$. Then the product $\calX = \prod_{\alpha \in A} \calX_{\alpha}$ is an $\infty$-topos.
Moreover, each projection $\pi^{\ast}_{\alpha}: \calX \rightarrow \calX_{\alpha}$ is left exact and colimit preserving. The corresponding geometric morphisms exhibit $\calX$ as a product of the family $\{ \calX_{\alpha} \}_{\alpha \in A}$ in the $\infty$-category $\LGeom$.
\end{proposition}

\begin{proof}
Proposition \ref{complexhorse2} implies that $\calX$ is presentable.
It is clear that a diagram $\overline{p}: K^{\triangleright} \rightarrow \calX$ is a colimit if and only if each composition
$\pi^{\ast}_{\alpha} \circ \overline{p}: K^{\triangleright} \rightarrow \calX_{\alpha}$
is a colimit diagram in $\calX_{\alpha}$. 
Similarly, a diagram $\overline{q}: K^{\triangleleft} \rightarrow \calX$ is a limit if and only if each composition
$\pi^{\ast}_{\alpha} \circ \overline{q}: K^{\triangleleft} \rightarrow \calX_{\alpha}$
is a limit diagram in $\calX_{\alpha}$. Using criterion $(2)$ of Theorem \ref{mainchar}, we deduce that $\calX$ is an $\infty$-topos, and that each $\pi^{\ast}_{\alpha}$ preserves {\em all} limits and colimits that exist in $\calX$. Choose a right adjoint
$\pi_{\ast}^{\alpha}: \calX_{\alpha} \rightarrow \calX$ to each $\pi^{\ast}_{\alpha}$. 

According to Theorem \ref{colimcomparee}, the $\infty$-category $\calX$ is a product
of the family $\{ \calX_{\alpha} \}_{\alpha \in A}$ in the $\infty$-category $\widehat{\Cat}_{\infty}$. Since $\LGeom$ is a subcategory of $\widehat{\Cat}_{\infty}$, it will suffice to prove the following assertion:

\begin{itemize}
\item For every $\infty$-topos $\calY$ and every functor $f^{\ast}: \calY \rightarrow \calX$
such that each of the composite functors $\calY \rightarrow \calX_{\alpha}$ is left exact and colimit preserving, $f^{\ast}$ is itself left exact and colimit preserving.
\end{itemize}

This follows immediately from the fact that limits and colimits are computed pointwise.
\end{proof}

\begin{proposition}\label{horse32}\index{gen}{pushout!of $\infty$-topoi}
Let $$ \xymatrix{ \calX' \ar[r]^{{q'}^{\ast}} \ar[d]^{{p'}^{\ast}} & \calX \ar[d]^{p^{\ast}} \\
\calY' \ar[r]^{q^{\ast}} & \calY }$$
be a diagram of $\infty$-categories which is homotopy Cartesian (with respect to the Joyal model structure). Suppose further that $\calX$, $\calY$, and $\calY'$ are $\infty$-topoi, and that
$p^{\ast}$ and $q^{\ast}$ are left exact and colimit preserving. Then $\calX'$ is an $\infty$-topos. Moreover, for any $\infty$-topos $\calZ$ and any functor $f^{\ast}: \calZ \rightarrow \calX$, $f^{\ast}$ is left exact and colimit preserving if and only if the compositions ${p'}^{\ast} \circ f^{\ast}$ and ${q'}^{\ast} \circ f^{\ast}$ are left exact and colimit preserving. In particular $($taking $f^{\ast} = \id_{\calX}${}$)$, the functors
${p'}^{\ast}$ and ${q'}^{\ast}$ are left exact and colimit preserving.
\end{proposition}

\begin{proof}
The second claim follows immediately from Lemma \ref{bird3} and the dual result concerning limits.
To prove the first, we observe that $\calX'$ is presentable by Proposition \ref{horse22}. To show that
$\calX$ is an $\infty$-topos, it will suffice to show that it satisfies criterion $(2)$ of Theorem \ref{mainchar}. This follows immediately from Lemma \ref{bird3}, given that
$\calX$ and $\calY'$ satisfy criterion $(2)$ of Theorem \ref{mainchar}.
\end{proof}

\begin{proposition}\label{colimtopoi}\index{gen}{colimit!of $\infty$-topoi}
The $\infty$-category $\LGeom$ admits small limits, and the inclusion functor
$\LGeom \subseteq \widehat{\Cat}_{\infty}$ preserves small limits.
\end{proposition}

\begin{proof}
According to Proposition \ref{alllimits}, it suffices to prove this result for pullbacks and small products. In the case of products, we apply Proposition \ref{quathorse}. For pullbacks, we
use Proposition \ref{horse32} and Theorem \ref{colimcomparee}.
\end{proof}

\subsection{Filtered Limits of $\infty$-Topoi}\label{inftyfiltlim}

We now consider the problem of computing limits in the $\infty$-category $\RGeom$ of $\infty$-topoi. This is quite a bit more difficult than the analogous problem for colimits, because the inclusion functor
$i: \RGeom \subseteq \widehat{\Cat}_{\infty}$ does not commute with limits in general.
However, the inclusion $i$ does commute with {\em filtered} limits:

\begin{theorem}\label{sutcar}\index{gen}{filtered limit!of $\infty$-topoi}
The $\infty$-category $\RGeom$ admits small filtered limits $($that is, limits indexed by diagrams
$\calC^{op} \rightarrow \RGeom$ where $\calC$ is a small, filtered $\infty$-category$)$. Moreover, the inclusion $\RGeom \subseteq \widehat{\Cat}_{\infty}$ preserves small, filtered limits.
\end{theorem}

The remainder of this section is devoted to the proof of Theorem \ref{sutcar}. Our basic strategy is to mimic the proof of Theorem \ref{surbus}. Our first step is to show that the limit (in $\widehat{\Cat}_{\infty}$) of a filtered diagram of $\infty$-topoi is itself an $\infty$-topos. This is equivalent to a more concrete assertion: if $p: X \rightarrow S$ is a topos fibration, and
$S^{op}$ is a small, filtered $\infty$-category, then the $\infty$-category $\calC$ of Cartesian sections of $p$ is an $\infty$-topos. We saw in Proposition \ref{seccatdog} that
$\calC$ is an accessible localization of the $\infty$-category $\bHom_{S}(S,X)$ spanned by
{\em all} sections of $p$. Our first step will be to show that $\bHom_{S}(S,X)$ is an $\infty$-topos. For this, the hypothesis that $S^{op}$ is filtered is irrelevant.

\begin{lemma}
Let $p: X \rightarrow S$ be a topos fibration, where $S$ is a small simplicial set. The
$\infty$-category $\bHom_{S}(S,X)$ of sections of $p$ is an $\infty$-topos.
\end{lemma}

\begin{proof}
This is a special case of Proposition \ref{prestorkus}.
\end{proof}

\begin{proposition}\label{steak1}
Let $A$ be a $($small$)$ filtered partially ordered set, and let $p: X \rightarrow \Nerve(A)$
be a topos fibration. Let $\calC = \bHom_{ \Nerve(A) }( \Nerve(A), X)$ be the $\infty$-category of sections of $p$, and let $\calC' \subseteq \calC$ be the full subcategory of $\calC$
spanned by the Cartesian sections of $p$. Then $\calC'$ is a topological localization of $\calC$.
\end{proposition}

\begin{proof}
Let us say that a subset $A' \subseteq A$ is {\it dense} if there exists $\alpha \in A$ such that
$$ \{ \beta \in A: \beta \geq \alpha \} \subseteq A'.$$
For each morphism $f$ in $\calC$, let $A(f) \subseteq A$ be the collection of all
$\alpha \in A$ such that the image of $f$ in $X_{\alpha}$ is an equivalence. Let
$S$ be the collection of all monomorphisms $f$ in $\calC$ such that $A(f)$ is dense.
It is clear that $S$ is stable under pullbacks, so that $S^{-1} \calC$ is a topological localization of $\calC$. To complete the proof, it will suffice to show that $\calC' = S^{-1} \calC$.

We first claim that each object of $\calC'$ is $S$-local. Let $f: C \rightarrow C'$
belong to $S$, and let $D \in \calC'$. Choose $\alpha_0$ such that $A(f)$ contains
$A' = \{ \beta \in A: \beta \geq \alpha_0 \}$, and let $R^{\ast}$ denote a right adjoint to
the restriction functor 
$$ R: \bHom_{ \Nerve(A)}( \Nerve(A), X) \rightarrow \bHom_{ \Nerve(A)}( \Nerve(A'), X).$$
According to Proposition \ref{leftkanadj}, the essential image of $R^{\ast}$ consists of those functors
$E: \Nerve(A) \rightarrow X$ which are $p$-right Kan extensions of $E | \Nerve(A')$. We claim that $D$ satisfies this condition. In other words, we claim that
for each $\alpha \in A$, the map
$$ \overline{q}: \Nerve(A'')^{\triangleleft}
\rightarrow \Nerve(A) \stackrel{D}{\rightarrow} X$$
is a $p$-limit, where $A'' = \{ \beta \in A: \beta \geq \alpha, \beta \geq \alpha_0 \}$. Since $\overline{q}$ carries each edge of
$\Nerve(A'')^{\triangleleft}$
to a $p$-Cartesian edge of $X$, it suffices to verify that the simplicial set
$\Nerve(A'')$ is weakly contractible (Proposition \ref{timal}). This follows immediately from the observation that $A''$ is a filtered partially ordered set.

We may therefore suppose that $D = R^{\ast} \overline{D}$, 
where $\overline{D} = D | \Nerve(A')$ is a Cartesian section of the induced
map $p': X \times_{ \Nerve(A)} \Nerve(A') \rightarrow \Nerve(A')$. 
We wish to prove that composition with $f$ induces a homotopy equivalence
$$ \bHom_{\calC}( C', R^{\ast} \overline{D}) \rightarrow \bHom_{\calC}(C, R^{\ast} \overline{D}).$$
This follows immediately from the fact that $R$ and $R^{\ast}$ are adjoint, since $R(f)$ is an equivalence.

We now show that every $S$-local object of $\calC$ belongs to $\calC'$. Let
$C \in \calC$ be a section of $p$ which is $S$-local. Choose $\alpha \leq \beta$ in $A$, and let 
$$ \Adjoint{F}{X_{\alpha}}{X_{\beta}}{G}$$
denote the (adjoint) functors associated to the (co)Cartesian fibration $p: X \rightarrow \Nerve(A)$.
The section $C$ gives rise to a pair of objects 
$C_{\alpha} \in X_{\alpha}$, $C_{\beta} \in X_{\beta}$, and a morphism 
$\phi: C_{\alpha} \rightarrow C_{\beta}$ in the $\infty$-category $X$. The map $\phi$ induces a morphism $u: C_{\alpha} \rightarrow G C_{\beta}$ in $X_{\alpha}$, which is well-defined up to equivalence. We wish to show that $\phi$ is $p$-Cartesian, which is equivalent to the assertion that $u$ is an equivalence in $X_{\alpha}$. Equivalently, we wish to show that for each object
$P \in X_{\alpha}$, composition with $u$ induces a homotopy equivalence
$$ \bHom_{X_{\alpha}}( P, C_{\alpha}) \rightarrow \bHom_{X_{\alpha}}(P, G C_{\beta}).$$

We may identify $P$ with a diagram
$$ \xymatrix{ \{ \alpha \} \ar[r]^{P} \ar@{^{(}->}[d] & X \ar[d] \\
\Nerve(A) \ar@{=}[r] \ar@{-->}[ur]^{D} & \Nerve(A). }$$
Using Corollary \ref{kanexistleft}, choose an extension $D$ as indicated in the diagram above, so that $D$ is a left Kan extension of $D | \{ \alpha \}$ over $\Nerve(A)$.
Similarly, we have a diagram
$$ \xymatrix{ \{ \beta \} \ar[r]^{F(P)} \ar@{^{(}->}[d] & X \ar[d] \\
\Nerve(A) \ar@{=}[r] \ar@{-->}[ur]^{D'} & \Nerve(A). }$$
and we can choose $D'$ to be a $p$-left Kan extension of
$D' | \{ \beta \}$.

Proposition \ref{leftkanadj} implies that for every object $E \in \calC$, the restriction maps
$$ \bHom_{\calC}(D,E) \rightarrow \bHom_{X_{\alpha}}( P, E(\alpha) )$$
$$ \bHom_{\calC}(D',E) \rightarrow \bHom_{X_{\beta}}( F(P), E(\beta) )$$
are equivalences. In particular, the equivalence between $D(\beta)$ and $F(P)$ induces
a morphism $\theta: D' \rightarrow D$. 

We have a commutative diagram in the homotopy category $\calH$:
$$ \xymatrix{ \bHom_{\calC}(D,C) \ar[rr]^{ \circ \theta} \ar[d] & & \bHom_{\calC}(D',C) \ar[d] \\
\bHom_{X_{\alpha}}(P, C_{\alpha}) \ar[r]^{\circ u} & \bHom_{X_{\alpha}}(P, G(C_{\beta})) &
\bHom_{X_{\beta}}( F(P), C_{\beta} ). \ar[l] }$$
The vertical maps are homotopy equivalences, and the 
the horizontal map on the lower right is a homotopy equivalence because $F$ and $G$ are adjoint.
To complete the proof, it will suffice to show that the upper horizontal map is an equivalence.
Since $C$ is $S$-local, it will suffice to show that $\theta \in S$.

Let $\beta \leq \beta'$, and consider the diagram
$$ \xymatrix{ & D'(\beta) \ar[r]^{w'} \ar[d]^{\theta(\beta)} & D'(\beta') \ar[d]^{\theta(\beta')} \\
D(\alpha) \ar[r]^{v} & D(\beta) \ar[r]^{w} & D(\beta'). }$$
in the $\infty$-category $X$. Since $D'$ is a $p$-left Kan extension of $D' | \{ \beta \}$, we conclude that $w'$ is $p$-coCartesian. Similarly, since $D$ is a $p$-left Kan
extension of $D | \{ \alpha \}$, we conclude that $v$ and $w \circ v$ are $p$-coCartesian.
It follows that $w$ is $p$-coCartesian as well (Proposition \ref{protohermes}). Since $\theta(\beta)$ is an equivalence by construction, we conclude that $\theta(\beta')$ is an equivalence. Thus
$A(\theta) \subseteq A$ is dense.

It remains only to show that $\theta$ is a monomorphism. For this, it suffices to show that
$\theta(\gamma)$ is an monomorphism in $X_{\gamma}$ for each $\gamma \in A$.
If $\gamma \geq \beta$, this follows from the above argument. Suppose $\gamma \ngeq \beta$.
Since $D'$ is a $p$-left Kan extension of $D' | \{ \beta \}$ over $\Nerve(A)$, we conclude
that $D'(\gamma)$ is a $p$-colimit of the empty diagram, and therefore an initial object
of $X_{\gamma}$. It follows that any map $D'(\gamma) \rightarrow D(\gamma)$ is a monomorphism.
\end{proof}

\begin{proposition}\label{steak2}
Let $A$ be a $($small$)$ filtered partially ordered set, let $p: X \rightarrow \Nerve(A)^{op}$, and let
$\calY \subseteq \bHom_{ \Nerve(A)}(\Nerve(A), X)$ be the full subcategory spanned by the Cartesian sections of $p$. For each $\alpha \in A$, the evaluation map
$\pi_{\ast}: \calY \rightarrow X_{\alpha}$ is a geometric morphism of $\infty$-topoi.
\end{proposition}

\begin{proof}
Let $A' = \{ \beta \in A: \alpha \leq \beta \}$. Using Theorem \ref{hollowtt}, we conclude that
the inclusion $\Nerve(A') \subseteq \Nerve(A)$ is cofinal. Corollary \ref{blurt} implies that
the restriction map
$$ \bHom_{ \Nerve(A)}( \Nerve(A), X) \rightarrow \bHom_{\Nerve(A)}(\Nerve(A'), X)$$
induces an equivalence on the full subcategories spanned by Cartesian sections. Consequently, we are free to replace $A$ by $A'$ and thereby assume that $\alpha$ is a least element of $A$.

The functor $\pi_{\ast}$ factors as a composition
$$ \calY \stackrel{\phi_{\ast}}{\rightarrow} \bHom_{ \Nerve(A)}( \Nerve(A), X) \stackrel{\psi_{\ast}}{\rightarrow}
X_{\alpha}$$
where $\phi_{\ast}$ denotes the inclusion functor and $\psi_{\ast}$ the evaluation functor.
Proposition \ref{steak1} implies that $\phi_{\ast}$ is a geometric morphism; it therefore suffices to show that $\psi_{\ast}$ is a geometric morphism as well. 

Let $\psi^{\ast}$ be a left adjoint to $\psi_{\ast}$ (the existence of $\psi^{\ast}$ follows from Proposition \ref{leftkanadj}, as indicated below). We wish to show that $\psi^{\ast}$ is left exact. According to Proposition \ref{limiteval}, it will suffice to show that the composition
$$ \theta: X_{\alpha} \stackrel{\psi^{\ast}}{\rightarrow} \bHom_{ \Nerve(A)}( \Nerve(A), X) 
\stackrel{e_{\beta}}{\rightarrow} X_{\beta},$$
is left exact, where $e_{\beta}$ denotes the functor given by evaluation at $\beta$.

Let $f: \Delta^1 \rightarrow \Nerve(A)$ be the edge joining $\alpha$ to $\beta$, let
$\calC$ be the $\infty$-category of coCartesian sections of $p$, and let
$\calC'$ be the $\infty$-category of coCartesian sections of the induced map
$p': X \times_{ \Nerve(A)} \Delta^1 \rightarrow \Delta^1$. We observe that
$\calC$ consists precisely of those sections $s: \Nerve(A) \rightarrow X$ of $p$ which
are $p$-left Kan extensions of $s | \{ \alpha \}$. Applying Proposition \ref{lklk},
we conclude that the evaluation map $e_{\alpha}: \calC \rightarrow X_{\alpha}$ is a trivial fibration, and that (by Proposition \ref{leftkanadj}) we may identify $\psi^{\ast}$ with the composition
$$ X_{\alpha} \stackrel{q}{\rightarrow} \calC \subseteq \bHom_{\Nerve(A)}(\Nerve(A),X),$$
where $q$ is a section of $e_{\alpha}|\calC$. Let $q': X_{\alpha} \rightarrow \calC'$
be the composition of $q$ with the restriction map $\calC \rightarrow \calC'$. Then
$\theta$ can be identified with the composition
$$ X_{\alpha} \stackrel{q'}{\rightarrow} \calC' \stackrel{e_{\beta}}{\rightarrow} X_{\beta},$$
which is the functor $X_{\alpha} \rightarrow X_{\beta}$ associated to $f: \alpha \rightarrow \beta$ by the coCartesian fibration $p$. Since $p$ is a topos fibration, $\theta$ is left exact as desired.
\end{proof}

Let $G$ be a profinite group and let $X$ be a set with a continuous action of $G$.
Then we can recover $X$ as the direct limit of the fixed point sets $X^{U}$, where $U$ ranges over the collection of open subgroups of $G$. Our next result is an $\infty$-categorical analogue of this observation.

\begin{lemma}\label{suruti}
Let $p: X \rightarrow S^{\triangleright}$ be a Cartesian fibration of simplicial sets, which
is classified by a colimit diagram $S^{\triangleright} \rightarrow \Cat_{\infty}^{op}$, and let
$\overline{s}: S^{\triangleright} \rightarrow X$ be a Cartesian section of $p$. Then $\overline{s}$ is a $p$-colimit diagram.
\end{lemma}

\begin{proof}
In virtue of Corollary \ref{tttroke}, we may suppose that $S$ is an $\infty$-category. Unwinding the definitions, we must show that the map $X_{\overline{s}/} \rightarrow X_{s/}$ induces an equivalence of $\infty$-categories when restricted to the inverse image of the cone point of
$S^{\triangleright}$. Fix an object $x \in X$ lying over the cone point of $S^{\triangleright}$.
Let $\overline{f}: S^{\triangleright} \rightarrow X$ be the constant map with value $x$, and let
$f = \overline{f} | S$. To complete the proof, it will suffice to show that the restriction map
$$ \theta: \bHom_{ \Fun(S^{\triangleright}, X)}( \overline{s}, \overline{f} )
\rightarrow \bHom_{ \Fun(S,X) }( s,f)$$ is a homotopy equivalence. To prove this, choose a
$p$-Cartesian transformation $\overline{\alpha}: \overline{g} \rightarrow \overline{f}$, where
$\overline{g}: S^{\triangleright} \rightarrow X$ is a section of $p$ (automatically Cartesian).
Let $g = \overline{g} | S$ and let $\alpha: g \rightarrow f$ be the associated transformation. Let
$\overline{\calC}$ be the full subcategory of $\bHom_{S^{\triangleright}}(S^{\triangleright},X)$
spanned by the Cartesian sections of $p$, and let $\calC \subseteq \bHom_{S^{\triangleright}}(S,X)$ be defined similarly.
We have a commutative diagram in the homotopy category $\calH$
$$ \xymatrix{ \bHom_{\overline{\calC} }( \overline{s}, \overline{g} )
\ar[r]^{\theta'} \ar[d]^{\overline{\alpha}} & \bHom_{\calC}( s, g) \ar[d]^{\alpha} \\
\bHom_{ \Fun(S^{\triangleright}, X)}( \overline{s}, \overline{f} ) \ar[r]^{\theta} &
\bHom_{ \Fun(S, X)}(s,f). }$$
Proposition \ref{compspaces} implies that the vertical maps are homotopy equivalences, and
Proposition \ref{charcatlimit} implies that $\theta'$ is a homotopy equivalence (since the restriction map $\overline{\calC} \rightarrow \calC$ is an equivalence of $\infty$-categories). It follows that
$\theta$ is a homotopy equivalence as well.
\end{proof}

\begin{lemma}\label{steakknife}
Let $p: X \rightarrow S$ be a presentable fibration, let $\calC$ 
be the full subcategory of $\bHom_{S}(S,X)$ spanned by the Cartesian sections of $p$. For
each vertex $s$ of $S$, let $\psi(s)_{\ast}: \calC \rightarrow X_{s}$ be the functor given by evaluation at $s$, and let $\psi(s)^{\ast}$ be a left adjoint to $\psi(s)_{\ast}$. There
exists a diagram $\theta: S \rightarrow \Fun(\calC, \calC)$ with the following properties:
\begin{itemize}
\item[$(1)$] For each vertex $s$ of $S$, $\theta(s)$ is equivalent to the composition
$\psi(s)^{\ast} \circ \psi(s)_{\ast}$.
\item[$(2)$] The identity functor $\id_{\calC}$ is a colimit of $\theta$ in the $\infty$-category of functors $\Fun(\calC, \calC)$. 
\end{itemize}
\end{lemma}

\begin{proof}
Without loss of generality, we may suppose that $p$ extends to a presentable fibration
$\overline{p}: \overline{X} \rightarrow S^{\triangleright}$, which is classified by colimit
diagram $S^{\triangleright} \rightarrow \LPres$ (and therefore by a colimit diagram
$S^{\triangleright} \rightarrow \Cat^{op}_{\infty}$, in virtue of Theorem \ref{surbus}). Let
$\overline{\calC}$ be the $\infty$-category of Cartesian sections of $\overline{p}$, so that
we have trivial fibrations
$$ \calC \leftarrow \overline{\calC} \rightarrow X_{\infty},$$
where $X_{\infty} = \overline{X} \times_{S^{\triangleright}} \{\infty\}$, and $\infty$
denotes the cone point of $S^{\triangleright}$. For each vertex $s$ of
$S^{\triangleright}$, we let $\overline{\psi}(s)_{\ast}: \overline{\calC} \rightarrow \overline{X}_s$
be the functor given by evaluation at $s$, and $\overline{\psi}(s)^{\ast}$ a left adjoint to
$\overline{\psi}(s)_{\ast}$. To complete the proof, it will suffice to construct a map 
$\theta': S \rightarrow \Fun( \overline{\calC}, X_{\infty})$ with the following properties:
\begin{itemize}
\item[$(1')$] For each vertex $s$ of $S$, $\theta'(s)$ is equivalent to the composition
$\overline{\psi}(\infty)_{\ast} \circ \overline{\psi}(s)^{\ast} \circ \overline{\psi}(s)_{\ast}$. 
\item[$(2')$] The functor $\overline{\psi}(\infty)_{\ast}$ is a colimit of $\theta'$.
\end{itemize}

Let $e: \overline{\calC} \times S^{\triangleright} \rightarrow \overline{X}$ be the evaluation map. 
Choose an $\overline{p}$-coCartesian natural transformation $e \rightarrow e'$, where
$e'$ is a map from $\overline{\calC} \times S^{\triangleright}$ to $X_{\infty}$. Lemma \ref{suruti} implies that for each object $X \in \overline{\calC}$, the restriction 
$e| \{X\} \times S^{\triangleright}$ is an $\overline{p}$-colimit diagram in $\overline{X}$. 
Applying Proposition \ref{chocolatelast}, we deduce that 
$e' | \{X\} \times S^{\triangleright}$ is a colimit diagram in $X_{\infty}$. According to Proposition
\ref{limiteval}, $e'$ determines a colimit diagram $S^{\triangleright} \rightarrow \Fun(\overline{\calC}, X_{\infty})$. Let $\theta'$ be the restriction of this diagram to $S$. Then the colimit of $\theta'$
can be identified with $e' | \overline{\calC} \times \{ \infty \}$, which is equivalent
to $e | \overline{\calC} \times \{\infty \} = \overline{\psi}(\infty)_{\ast}$. This proves $(2')$. 
To verify $(1')$, we observe that $e' | \overline{\calC} \times \{s\}$ can be identified with
the composition of $\overline{\psi}(s)_{\ast} = e | \overline{\calC} \times \{s\}$ with the
functor $X_{s} \rightarrow X_{\infty}$ associated to the coCartesian fibration
$\overline{p}$, which can in turn be identified with $\overline{\psi}(\infty)_{\ast} \circ \overline{\psi}(s)^{\ast}$ (both are left adjoints to the pullback functor
$X_{\infty} \rightarrow X_{s}$ associated to $\overline{p}$). 
\end{proof}

\begin{proposition}\label{steak3}
Let $A$ be a $($small$)$ filtered partially ordered set, let $p: X \rightarrow \Nerve(A)$, and let
$\calY \subseteq \bHom_{ \Nerve(A)}(\Nerve(A), X)$ be the full subcategory spanned by the Cartesian sections of $p$. Let $\calZ$ be an $\infty$-topos, and $\pi_{\ast}: \calZ \rightarrow \calY$ an arbitrary functor. Suppose that, for each $\alpha \in A$, the composition
$$ \calZ \stackrel{\pi_{\ast}}{\rightarrow} \calY \rightarrow X_{\alpha}$$
is a geometric morphism of $\infty$-topoi. Then $\pi_{\ast}$ is a geometric morphism of $\infty$-topoi.
\end{proposition}

\begin{proof}
Let $\pi^{\ast}$ denote a left adjoint to $\pi_{\ast}$. Since $\pi^{\ast}$ commutes with colimits,
Lemma \ref{steakknife} implies that $\pi^{\ast}$ can be written as the colimit of a diagram
$q: \Nerve(A) \rightarrow \calZ^{\calY}$ having the property that for each
$\alpha \in A$, $q(\alpha)$ is equivalent to $\pi^{\ast} \psi(\alpha)^{\ast} \psi(\alpha)_{\ast}$,
where $\psi(\alpha)_{\ast}$ denotes the evaluation functor at $\alpha$ and $\psi(\alpha)^{\ast}$ its left adjoint. Each composition $\pi^{\ast} \psi(\alpha)^{\ast}$ is left adjoint to the geometric morphism
$\psi(\alpha)_{\ast} \pi_{\ast}$, and therefore left exact. It follows that $q(\alpha)$ is left exact.
Since filtered colimits in $\calZ$ are left exact (Example \ref{tucka}), we conclude that
the functor $\pi^{\ast}$ is left exact, as desired.
\end{proof}

\begin{proof}[Proof of Theorem \ref{sutcar}]
Let $\calC$ be a small, filtered $\infty$-category, and let $q: \calC^{op} \rightarrow \RGeom$
be an arbitrary diagram. Choose a limit $\overline{q}: (\calC^{\triangleright})^{op} \rightarrow \widehat{\Cat}_{\infty}$ of $q$ in the $\infty$-category $\widehat{\Cat}_{\infty}$. We must show that
$\overline{q}$ factors through $\RGeom$, and is a limit diagram in $\RGeom$.

Using Proposition \ref{rot}, we may assume without loss of generality that $\calC$ is the nerve of a filtered partially ordered set $A$. Let $p: X \rightarrow \Nerve(A)^{op}$ be the topos fibration classified by $q$ (Proposition \ref{surtog2}). Then the image of the cone point of $( \calC^{\triangleright})^{op}$ under $\overline{q}$ is equivalent to the $\infty$-category $\calX$ of Cartesian sections of $p$ (Corollary \ref{blurt}). It follows from Proposition \ref{steak1} that
$\calX$ is an $\infty$-topos. Moreover, Proposition \ref{steak2} ensures that for each $\alpha \in A$, the evaluation map $\calX \rightarrow X_{\alpha}$ is a geometric morphism. This proves that 
$\overline{q}$ factors through $\RGeom$. To complete the proof, we must show that $\overline{q}$ is a limit diagram in $\RGeom$. Since $\RGeom$ is a subcategory of $\widehat{\Cat}_{\infty}$, and $\overline{q}$ is a limit diagram in $\widehat{\Cat}_{\infty}$, this reduces immediately to the statement of Proposition \ref{steak3}.
\end{proof}

\subsection{General Limits of $\infty$-Topoi}\label{genlim}

Our goal in this section is to construct general limits in the $\infty$-category $\RGeom$. Our strategy is necessarily rather different from that of \S \ref{inftyfiltlim}, because the inclusion
$i: \RGeom \rightarrow \widehat{\Cat}$ does not preserve limits in general.
In fact, $i$ does not even preserve the final object:

\begin{proposition}\label{spacefinall}\index{gen}{final object!of the $\infty$-category of $\infty$-topoi}
Let $\calX$ be an $\infty$-topos. Then $\Fun^{\ast}(\SSet, \calX)$ is a contractible Kan complex.
In particular, $\SSet$ is a final object in the $\infty$-category $\RGeom$ of $\infty$-topoi.
\end{proposition}

\begin{proof}
We observe that $\SSet \simeq \Shv(\Delta^0)$ where the $\infty$-category $\Delta^0$ is endowed with the ``discrete'' topology (so that the empty sieve does not constitute a cover
of the unique object). According to Proposition \ref{igrute}, the $\infty$-category
$\Fun^{\ast}(\SSet, \calX)$ is equivalent to the full subcategory of $\calX \simeq \Fun(\Delta^0, \calX)$ spanned
by those objects $X \in \calX$ which correspond to left exact functors $\Delta^0 \rightarrow \calX$.
It is clear that these are precisely the final objects of $\calX$, which form a contractible Kan complex
(Proposition \ref{initunique}). 
\end{proof}

To construct limits in general, we first develop some tools for describing $\infty$-topoi via
``generators and relations''. This will allow us to reduce the construction of limits in $\RGeom$ to the problem of constructing colimits in $\Cat_{\infty}$.

\begin{lemma}\label{sleepyswine}
Let $\calC$ be a small $\infty$-category, $\kappa$ a regular cardinal, and suppose given
a $($small$)$ collection of $\kappa$-small diagrams $\{ \overline{f}_{\alpha}: K_{\alpha}^{\triangleright} \rightarrow \calC \}_{\alpha \in A}$. Then there exists a functor $F: \calC \rightarrow \calD$ with the following properties:
\begin{itemize}
\item[$(1)$] The $\infty$-category $\calD$ is small and admits $\kappa$-small colimits.
\item[$(2)$] For each $\alpha \in A$, the induced map $F \circ \overline{f}_{\alpha}: K_{\alpha}^{\triangleright} \rightarrow \calC$ is a colimit diagram.
\item[$(3)$] Let $\calE$ be an arbitrary $\infty$-category which admits $\kappa$-small colimits. Let
$\Fun'(\calD, \calE)$ denote the full subcategory of $\Fun(\calD, \calE)$ spanned by those functors which preserve $\kappa$-small colimits. Then composition with $F$ induces a fully faithful embedding
$$ \theta: \Fun'(\calD, \calE) \rightarrow \Fun( \calC, \calE).$$
The essential image of $\theta$ consists those functors $F': \calC \rightarrow \calE$ such that
each $F' \circ \overline{f}_{\alpha}$ is a colimit diagram in $\calE$.
\end{itemize}
\end{lemma}

\begin{proof}
Let $j: \calC \rightarrow \calP(\calC)$ denote the Yoneda embedding. For each $\alpha \in A$, let
$f_{\alpha} = \overline{f}_{\alpha} | K_{\alpha}$, and let $C_{\alpha} \in \calC$ denote the image of the cone point under $\overline{f}_{\alpha}$. Let $D_{\alpha} \in \calP(\calC)$ denote a colimit of the induced diagram $j \circ f_{\alpha}$, so that $j \circ \overline{f}_{\alpha}$ induces a map
$s_{\alpha} = D_{\alpha} \rightarrow j( C_{\alpha} )$. Let $S = \{ s_{\alpha} \}_{\alpha \in A}$, let $\calX$ denote the localization $S^{-1} \calP(\calC)$, and let $L: \calP(\calC) \rightarrow \calX$ denote a left adjoint to the inclusion. Let $\calD'$ be the smallest full subcategory of $\calX$ that contains the essential image of the functor $L \circ j$ and is stable under $\kappa$-small colimits, let
$\calD$ be a minimal model for $\calD'$, and let $F: \calC \rightarrow \calD$ denote the composition of
$L \circ j$ with a retraction of $\calD'$ onto $\calD$. It follows immediately from the construction that $\calD$ satisfies conditions $(1)$ and $(2)$.

We observe that for each $\alpha \in A$, the domain and codomain of
$s_{\alpha}$ are both $\kappa$-compact objects of $\calP(\calC)$. It follows that
$\calX$ is stable under $\kappa$-filtered colimits in $\calP(\calC)$. Corollary \ref{starmin} implies that $L$ carries $\kappa$-compact objects of $\calP(\calC)$ to $\kappa$-compact objects of $\calX$. 
Since the collection of $\kappa$-compact objects of $\calX$ is stable under $\kappa$-small colimits, we conclude that $\calD'$ consists of $\kappa$-compact objects of $\calX$. Invoking Proposition \ref{uterr}, we deduce that the inclusion $\calD \subseteq \calX$ determines an equivalence
$\Ind_{\kappa}(\calD) \simeq \calX$. 

We now prove $(3)$. We observe that there exists a fully faithful embedding $i: \calE \rightarrow \calE'$ which preserves $\kappa$-small colimits, where $\calE'$ admits arbitrary small colimits (for example, we can take $\calE' = \Fun( \calE, \widehat{\SSet})^{op}$ and $i$ to be the Yoneda embedding).
Replacing $\calE$ by $\calE'$ if necessary, we may assume that $\calE$ itself admits arbitrary small colimits. We have a homotopy commutative diagram
$$ \xymatrix{ \Fun^{L}( \calX, \calE) \ar[r]^{\theta'} \ar[d] & \Fun^{L}( \calP(\calC), \calE) \ar[d] \\
\Fun'(\calD, \calE) \ar[r]^{\theta} & \Fun( \calC, \calE), }$$
where $\Fun^{L}( \calY, \calE)$ denotes the full subcategory of $\Fun(\calY, \calE)$ spanned by those functors which preserve small colimits. Propositions \ref{intprop} and \ref{sumatch} imply that the left vertical arrow is an equivalence, while Theorem \ref{charpresheaf} implies that the right vertical arrow is an equivalence. It will therefore suffice to show that $\theta'$ is fully faithful, and that the essential image of $\theta'$ consists of those colimit-preserving functors $F'$ from $\calP(\calC)$ to $\calE$ such that $F' \circ j \circ \overline{f}_{\alpha}$ is a colimit diagram, for each $\alpha \in A$. This follows immediately from Proposition \ref{unichar}.
\end{proof}

\begin{definition}\index{not}{Catinftylex@$\Cat_{\infty}^{\lex}$}
Let $\Cat_{\infty}^{\lex}$ denote the subcategory of $\Cat_{\infty}$ defined as follows:
\begin{itemize}
\item[$(1)$] A small $\infty$-category $\calC$ belongs to $\Cat_{\infty}^{\lex}$ if and only if
$\calC$ admits finite limits.
\item[$(2)$] Let $f: \calC \rightarrow \calD$ be a functor between small $\infty$-categories which admit finite limits. Then $f$ is a morphism in $\Cat_{\infty}^{\lex}$ if and only if $f$ preserves finite limits.
\end{itemize}
\end{definition}

\begin{lemma}\label{pugswell}
The $\infty$-category $\Cat_{\infty}^{\lex}$ admits small colimits. 
\end{lemma}

\begin{proof}
Let $p: \calJ \rightarrow \Cat_{\infty}^{\lex}$ be a small diagram, which carries each
vertex $j \in \calJ$ to an $\infty$-category $\calC_{j}$.
Let $\calC$ be a colimit of the diagram $p$ in $\Cat_{\infty}$, and for each
$j \in J$, let $\phi_{j}: \calC_{j} \rightarrow \calC$ be the associated functor.
Consider the collection of all isomorphism classes of diagrams $\{ f: K^{\triangleleft} \rightarrow \calC \}$, where $K$ is a finite simplicial set and the map $f$ admits a factorization
$$ K^{\triangleleft} \stackrel{f_0}{\rightarrow} \calC_{j} \stackrel{\phi_{j}}{\rightarrow} \calC,$$
where $f_0$ is a limit diagram in $\calC_{j}$. Invoking the dual of Lemma \ref{sleepyswine}, we deduce the existence of a functor $F: \calC \rightarrow \calD$ with the following properties:
\begin{itemize}
\item[$(1)$] The $\infty$-category $\calD$ is small and admits finite limits.
\item[$(2)$] Each of the compositions $F \circ \phi_{j}$ is left exact.
\item[$(3)$] For every $\infty$-category $\calE$ which admits finite limits, composition with
$F$ induces an equivalence from the full subcategory of $\Fun(\calD, \calE)$ spanned by the left exact functors to the full subcategory of $\Fun( \calC, \calE)$ spanned by those functors $F': \calC \rightarrow \calE$ such that each $F' \circ \phi_{j}$ is left exact.
\end{itemize}
It follows that that $\calD$ can be identified with a colimit of the diagram $p$ in
the $\infty$-category $\Cat_{\infty}^{\lex}$.
\end{proof}

\begin{lemma}\label{squarepeg}
Let $\calC$ be a small $\infty$-category which admits finite limits, and let
$f_{\ast}: \calX \rightarrow \calP(\calC)$ be a geometric morphism of $\infty$-topoi.
Then there exists a small $\infty$-category $\calD$ which admits finite limits, a left
exact functor $f'': \calC \rightarrow \calD$ such that $f_{\ast}$ is equivalent to the composition,
$$ \calX \stackrel{f'_{\ast}}{\rightarrow} \calP( \calD) \stackrel{ f''_{\ast}}{\rightarrow} \calP(\calC)$$
where $f'_{\ast}$ is a fully faithful geometric morphism and $f''_{\ast}$ is given by composition with $f''$.
\end{lemma}

\begin{proof}
Without loss of generality we may assume that $\calX$ is minimal. Let $f^{\ast}$ be a left adjoint to $f_{\ast}$. Choose a regular cardinal $\kappa$ large enough that the composition
$$ \calC \stackrel{j_{\calC}}{\rightarrow} \calP(\calC) \stackrel{f^{\ast}}{\rightarrow} \calX$$
carries each object $C \in \calC$ to a $\kappa$-compact object of $\calX$. Enlarging
$\kappa$ if necessary, we may assume that $\calX$ is $\kappa$-accessible and that the collection
of $\kappa$-compact objects is stable under finite limits (Proposition \ref{tcoherent}). Let $\calD$ be the collection of $\kappa$-compact objects of $\calX$. The proof of Proposition \ref{precisechar} shows that the inclusion $\calD \subseteq \calX$ can be extended to a left exact localization functor ${f'}^{\ast}: \calP(\calD) \rightarrow \calX$. 

Using Theorem \ref{charpresheaf}, we conclude that the composition
$j_{\calD} \circ f^{\ast} \circ j_{\calC}: \calC \rightarrow \calP(\calD)$ can be extended
to a colimit-preserving functor ${f''}^{\ast}: \calP(\calC) \rightarrow \calP(\calD)$, and that
${f'}^{\ast} \circ {f''}^{\ast}$ is homotopic to $f^{\ast}$. Proposition \ref{natash} implies that ${f''}^{\ast}$ is left exact. It follows that ${f'}^{\ast}$ and ${f''}^{\ast}$ admit right adjoints $f'_{\ast}$ and
$f''_{\ast}$ with the desired properties.
\end{proof}

\begin{proposition}\label{swunder}
The $\infty$-category $\RGeom$ of $\infty$-topoi admits pullbacks.
\end{proposition}

\begin{proof}
Suppose first that we are given a pullback square
$$ \xymatrix{ \calW \ar[r]^{f'_{\ast}} \ar[d]^{g'_{\ast}} & \calX \ar[d]^{g_{\ast}} \\
\calY \ar[r]^{f_{\ast}} & \calZ }$$
in the $\infty$-category of $\RGeom$. We make the following observations:

\begin{itemize}
\item[$(a)$] Suppose that $\calZ$ is a left exact localization of another $\infty$-topos $\calZ'$. Then the induced diagram
$$ \xymatrix{ \calW \ar[r] \ar[d] & \calX \ar[d] \\
\calY \ar[r] & \calZ }$$
is also a pullback square.

\item[$(b)$] Let $S^{-1} \calX$ and $T^{-1} \calY$ be left exact localizations of $\calX$ and $\calY$, respectively. Let $U$ be the smallest strongly saturated collection of morphisms in
$\calW$ which contains ${f'}^{\ast} S$ and ${g'}^{\ast} T$, and is closed under pullbacks. Using Corollary \ref{sweetums}, we deduce that $U$ is generated by a (small) set of morphisms in $\calW$. It follows that the diagram
$$ \xymatrix{ U^{-1} \calW \ar[r] \ar[d] & S^{-1} \calX \ar[d] \\
T^{-1} \calY \ar[r] & \calZ }$$
is again a pullback in $\RGeom$.
\end{itemize}

Now suppose given an arbitrary diagram
$$ \calX \stackrel{g_{\ast} }{\rightarrow} \calZ \stackrel{f_{\ast}}{\leftarrow} \calY$$
in $\RGeom$. We wish to prove that there exists a fiber product $\calX \times_{ \calZ} \calY$
in $\RGeom$. The proof of Proposition \ref{precisechar} implies that there exists
a small $\infty$-category $\calC$ which admits finite limits, such that
$\calZ$ is a left exact localization of $\calP(\calC)$. Using $(a)$, we can reduce to the case where $\calZ = \calP(\calC)$. Using $(b)$ and Lemma \ref{squarepeg}, we can reduce to the case where
$\calX = \calP(\calD)$ for some small $\infty$-category $\calD$ which admits finite limits, and $g_{\ast}$ is induced by composition with a left exact functor $g: \calC \rightarrow \calD$. Similarly, we can assume that $f_{\ast}$ is determined by a left exact functor $f: \calC \rightarrow \calD'$.
Using Lemma \ref{pugswell}, we can form a pushout diagram
$$ \xymatrix{ \calE & \calD \ar[l] \\
\calD' \ar[u] & \calC \ar[u]^{g} \ar[l]^{f} }$$
in the $\infty$-category $\Cat_{\infty}^{\lex}$. Using Proposition \ref{natash} and Theorem \ref{charpresheaf}, it is not difficult to see that the induced diagram
$$ \xymatrix{ \calP(\calE) \ar[r] \ar[d] & \calP(\calD) \ar[d]^{g_{\ast}} \\
\calP( \calD') \ar[r]^{f_{\ast}} & \calP( \calC) }$$
is the desired pullback square in $\RGeom$.
\end{proof}

\begin{corollary}\label{geolit}
The $\infty$-category $\RGeom$ admits small limits.
\end{corollary}

\begin{proof}
Using Corollaries \ref{uterrr} and \ref{allfin}, it suffices to show that $\RGeom$ admits filtered limits, a final object, and pullbacks. The existence of filtered limits follows from Theorem \ref{sutcar}, the existence of a final object follows from Proposition \ref{spacefinall}, and the existence of pullbacks follows from Proposition \ref{swunder}.
\end{proof}

\begin{remark}
Our construction of fiber products in $\RGeom$ is somewhat inexplicit. We will later give a more concrete construction in the case of ordinary products; see \S \ref{products}.
\end{remark}

We conclude this section by proving a companion result to Corollary \ref{geolit}. First, a few general remarks. The $\infty$-category $\RGeom$ is most naturally viewed as an {\it $\infty$-bicategory}, since we can consider also noninvertible natural transformations between geometric morphisms. Correspondingly, we can consider a more general theory of {\it $\infty$-bicategorical} limits in $\RGeom$. While we do not want to give any precise definitions, we would like to point out Corollary \ref{geolit} can be generalized to show that $\RGeom$ admits all (small) $\infty$-bicategorical limits. In more concrete terms, this just means that $\RGeom$ is {\it cotensored} over $\Cat_{\infty}$ in the following sense:

\begin{proposition}\label{cotens}
Let $\calX$ be an $\infty$-topos, and let $\calD$ be a small $\infty$-category. Then there exists an $\infty$-topos $\Mor(\calC, \calX)$ and a functor $e: \calC \rightarrow \Fun_{\ast}( \Mor(\calC, \calX), \calX)$
with the following universal property:
\begin{itemize}
\item[$(\ast)$] For every $\infty$-topos $\calY$, composition with $e$ induces an equivalence of $\infty$-categories
$$ \Fun_{\ast}( \calY, \Mor(\calC, \calX) ) \rightarrow \Fun( \calC, \Fun_{\ast}( \calY, \calX) ).$$
\end{itemize}
\end{proposition}

\begin{proof}
We first treat the case where $\calX = \calP(\calD)$, where $\calD$ is a small $\infty$-category which admits finite limits. Using Lemma \ref{sleepyswine}, we conclude that there exists a functor
$e_0: \calC^{op} \times \calD \rightarrow \calD'$ with the following properties:
\begin{itemize}
\item[$(1)$] The $\infty$-category $\calD'$ is small and admits finite limits.
\item[$(2)$] For each object $C \in \calC$, the induced functor
$$ \calD \simeq \{ C \} \times \calD \subseteq \calC^{op} \times \calD \stackrel{e_0}{\rightarrow} \calD'$$
is left exact.
\item[$(3)$] Let $\calE$ be an arbitrary $\infty$-category which admits finite limits. Then
composition with $e_0$ induces an equivalence from full subcategory of $\Fun( \calD', \calE)$ spanned by the left exact functors to the full subcategory of $\Fun( \calC^{op} \times \calD, \calE)$ spanned by those functors which restrict to left-exact functors $\{ C \} \times \calD \rightarrow \calE$, for each $C \in \calC$.
\end{itemize}
In this case, we can define $\Mor(\calC, \calX)$ to be $\calP( \calD')$, and $e: \calC \rightarrow
\Fun_{\ast}( \calP( \calD'), \calP(\calD) )$ to be given by composition with $e_0$; the universal property
$(\ast)$ follows immediately from Theorem \ref{charpresheaf} and Proposition \ref{natash}.

In the general case, we invoke Proposition \ref{precisechar} to reduce to the case where
$\calX = S^{-1} \calX'$ is an accessible left-exact localizaton of an $\infty$-topos
$\calX'$, where $\calX' \simeq \calP(\calD)$ is as above so that we can
construct an $\infty$-topos $\Mor( \calC, \calX')$ and a map 
$e': \calC \rightarrow \Fun_{\ast}( \Mor(\calC, \calX'), \calX')$ satisfying the condition $(\ast)$.
For each $C \in \calC$, let $e'(C)_{\ast}$ denote the corresponding geometric morphism
from $\Mor(\calC, \calX')$ to $\calX'$, let $e'(C)^{\ast}$ denote a left adjoint to
$e'(C)_{\ast}$, and let $S(C) = e'(C)^{\ast} S$ be the image of $S$ in the collection of morphisms
of $\Mor(\calC, \calX')$. Since each $e'(C)^{\ast}$ is a colimit-preserving functor, each of the sets 
$S(C)$ is generated under colimits by a small collection of morphisms in $\Mor(\calC, \calX')$. Let $T$ be the smallest collection of morphisms in
$\Mor(\calC, \calX')$ which is strongly saturated, stable under pullbacks, and contains each of the sets $S_{C}$. Using Corollary \ref{sweetums}, we conclude that $T$ is generated (as a strongly saturated class of morphisms) by a small collection of morphisms in $\Mor(\calC, \calX')$. It follows that
$\Mor(\calC, \calX) = T^{-1} \Mor( \calC, \calX')$ is an $\infty$-topos. By construction,
the map $e'$ restricts to give a functor
$e: \calC \rightarrow \Fun_{\ast}( \Mor( \calC, \calX), \calX)$. Unwinding the definitions, we see that
$e$ has the desired properties.
\end{proof}

\begin{remark}\label{sablewise}
Let $\calX$ be an $\infty$-topos, and let $\RGeom^{\Delta}$ denote the simplicial
subcategory of $\widehat{\Cat}_{\infty}^{\Delta}$ corresponding to the subcategory
$\RGeom \subseteq \widehat{\Cat}_{\infty}$, so that
$\RGeom \simeq \Nerve( \RGeom^{\Delta})$. The construction 
$\calY \mapsto \Fun_{\ast}( \calX, \calY)$ determines a simplicial functor from
$\RGeom^{\Delta}$ to $\widehat{\Cat}_{\infty}^{\Delta}$, which in turn induces a functor
$$\theta_{\calX}: \RGeom \rightarrow \widehat{\Cat}_{\infty}.$$
We claim that $\theta_{\calX}$ preserves small limits (this translates
into the condition that limits in $\RGeom$ really give {\em $\infty$-bicategorical limits} in the
$\infty$-bicategory of $\infty$-topoi).

To prove this, fix an arbitrary $\infty$-category $\calC$, and let $e_{\calC}: \widehat{\Cat}_{\infty} \rightarrow \widehat{\SSet}$ be the functor corepresented by $\calC$. It will suffice to show
that $e_{\calC} \circ \theta_{\calX}$ preserves small limits. The collection of all $\infty$-categories $\calC$ which satisfy this condition is stable under all colimits, so we may assume without loss of generality that $\calC$ is small. It now suffices to observe that $e_{\calC} \circ \theta_{\calX}$ is
equivalent to the functor corepresented by the $\infty$-topos $\Fun(\calC, \calX)$.
\end{remark}

\subsection{\'{E}tale Morphisms of $\infty$-Topoi}\label{gemor2}

Let $f: X \rightarrow Y$ be a continuous map of topological spaces.
We say that $f$ is {\it \'{e}tale} (or a {\it local homeomorphism}) if, for
every point $x \in X$, there exist open sets $U \subseteq X$
containing $x$ and $V \subseteq Y$ containing $f(x)$ such that
$f$ induces a homeomorphism $U \rightarrow V$.  
Let $\calF$ denote the sheaf of sections of $f$: that is, $\calF$ is a sheaf of
sets on $Y$ such that for every open set $V \subseteq Y$,
$\calF(V)$ is the set of all continuous maps $s: V \rightarrow X$ such that
$f \circ s = \id: V \rightarrow Y$. The construction
$(f: X \rightarrow Y) \mapsto \calF$ determines an equivalence of
categories, from the category of topological spaces which
are \'{e}tale over $Y$ to the category of sheaves (of sets) on $Y$.
In particular, we can recover the topological space $X$ (up to homeomorphism) from the sheaf of sets $\calF$ on $Y$. For example, we
can reconstruct the category $\Shv_{\Set}(X)$ of sheaves on $X$
as the overcategory $\Shv_{\Set}(Y)_{/ \calF}$.

Our goal in this section is to develop an analogous theory of \'{e}tale morphisms in the setting of $\infty$-topoi. Suppose given a geometric morphism
$f_{\ast}: \calX \rightarrow \calY$. Under what circumstances should we say that
$f_{\ast}$ is \'{e}tale? By analogy with the case of topological spaces, we should expect that an \'{e}tale morphism determines a ``sheaf'' on $\calY$: that is, an object $U$ of the $\infty$-category $\calY$. Moreover, we should then be able to recover the $\infty$-category $\calX$ as an overcategory $\calY_{/U}$. 
The following result guarantees that this expectation is somewhat reasonable:

\begin{proposition}\label{generalslice}\index{gen}{overcategory!of an $\infty$-topos}
Let $\calX$ be an $\infty$-topos, and let $U$ be an object of $\calX$.
\begin{itemize}
\item[$(1)$] The $\infty$-category $\calX_{/U}$ is an $\infty$-topos.
\item[$(2)$] The projection $\pi_{!}: \calX_{/U} \rightarrow \calX$ has a right adjoint $\pi^{\ast}$ which
commutes with colimits. Consequently, $\pi^{\ast}$ itself has a right adjoint $\pi_{\ast}: \calX_{/U} \rightarrow \calX$, which is a geometric morphism of $\infty$-topoi.
\end{itemize}
\end{proposition}

\begin{proof}
The existence of a right adjoint $\pi^{\ast}$ to the projection $\pi_{!}: \calX_{/U} \rightarrow \calX$ follows from the assumption that $\calX$ admits finite limits. Moreover, the assertion that $\pi^{\ast}$ preserves colimits is a special case of the assumption that colimits in $\calX$ are universal. This proves $(2)$.

To prove $(1)$, we will show that $\calX_{/U}$ satisfies criterion $(2)$ of Theorem \ref{mainchar}. 
We first observe that $\calX_{/U}$ is presentable (Proposition \ref{slicstab}).
Let $K$ be a small simplicial set, and let $\overline{\alpha}: \overline{p} \rightarrow \overline{q}$
be a natural transformation of diagrams $\overline{p}, \overline{q}: K^{\triangleright}
\rightarrow \calX_{/U}$. Suppose that $\overline{q}$ is a colimit diagram, and that
$\alpha = \overline{\alpha} | K$ is a Cartesian transformation. The projection
$\pi_{!}$ preserves all colimits (since it is a left adjoint), so that $\pi_{!} \circ \overline{q}$ is a colimit
diagram in $\calX$. Since $\pi_{!}$ preserves pullback squares (Proposition \ref{goeselse}), $\pi_{!} \circ \alpha$
is a Cartesian transformation. By assumption, $\calX$ is an $\infty$-topos, so that
Theorem \ref{mainchar} implies that $\pi_{!} \circ \overline{p}$ is a colimit diagram
if and only if $\pi_{!} \circ \overline{\alpha}$ is a Cartesian transformation. Using
Propositions \ref{goeselse} and \ref{needed17}, we conclude
that $\overline{p}$ is a colimit diagram if and only if $\overline{\alpha}$ is a Cartesian transformation, as desired.
\end{proof}

A geometric morphism $f_{\ast}: \calX \rightarrow \calY$ of $\infty$-topoi is said to be {\it \'{e}tale}\index{gen}{\'{e}tale morphism} if it arises via the construction of Proposition \ref{generalslice}; that is, if $f$ admits a factorization
$$ \calX \stackrel{f'_{\ast}}{\rightarrow} \calY_{/U} \stackrel{f''_{\ast}}{\rightarrow} \calY$$ where
$U$ is an object of $\calY$, $f'_{\ast}$ is a categorical equivalence, and $f''_{\ast}$ is a right adjoint to the pullback functor ${f''}^{\ast}: \calY \rightarrow \calY_{/U}$. We note that in this case, $f^{\ast}$ has a {\em left adjoint} $f_{!} = f''_{!} \circ f'_{\ast}$. Consequently, $f^{\ast}$ preserves {\em all} limits, not just finite limits. 

\begin{remark}\label{mark}
Given an \'{e}tale geometric morphism $f: \calX_{/U} \rightarrow \calX$ of $\infty$-topoi, the
description of the pushforward functor $f_{\ast}$ is slightly more complicated than that of $f_{!}$ (which is merely the forgetful functor) or $f^{\ast}$ (which is given by taking products with $U$).
Given an object $p: X \rightarrow U$ of $\calX_{/U}$, the pushforward $f_{\ast} X$ is an object of
$\calX$ which represents the functor ``sections of $p$''.
\end{remark}

The collection of \'{e}tale geometric morphisms contains all equivalences and is stable under composition. Consequently, we can consider the subcategory $\RGeom_{\mathet} \subseteq \RGeom$\index{not}{RGeomet@$\RGeom_{\mathet}$} containing all objects of $\RGeom$, whose morphisms are precisely the \'{e}tale geometric morphisms. Our goal in this section is to study the $\infty$-category $\RGeom_{\mathet}$. Our main results are the following:

\begin{itemize}
\item[$(a)$] If $\calX$ is an $\infty$-topos containing an object $U$, then the associated
\'{e}tale geometric morphism $\pi_{\ast}: \calX_{/U} \rightarrow \calX$ can be described by a
universal property. Namely, $\calX_{/U}$ is universal among $\infty$-topoi $\calY$
with a geometric morphism $\phi_{\ast}: \calY \rightarrow \calX$ such that
$\phi^{\ast} U$ admits a global section (Proposition \ref{goodking}). 

\item[$(b)$] There is a simple criterion for testing whether a geometric morphism
$\pi_{\ast}: \calX \rightarrow \calY$ is \'{e}tale. Namely, $\pi_{\ast}$ is \'{e}tale if and only if
the functor $\pi^{\ast}$ admits a left adjoint $\pi_{!}$, the functor $\pi_{!}$ is conservative, and an appropriate push-pull formula holds in the the $\infty$-category $\calY$ (Proposition \ref{pushpo}).

\item[$(c)$] Given a pair of topological spaces $X_0$ and $X_1$ and a homeomorphism
$\phi: U_0 \simeq U_1$ between open subsets $U_0 \subseteq X_0$ and $U_1 \subseteq X_1$,
we can ``glue'' $X_0$ to $X_1$ along $\phi$ to obtain a new topological space
$X = X_0 \coprod_{U_0} X_1$. Moreover, the topological space $X$ contains
open subsets homeomorphic to $X_0$ and $X_1$. In the setting of $\infty$-topoi, it is possible to make much more general ``gluing'' constructions of the same type. We can formulate this idea more precisely as follows: given any diagram $\{ \calX_{\alpha} \}$ in the $\infty$-category $\RGeom_{\mathet}$ having a colimit $\calX$ in $\RGeom$, each of the associated geometric morphisms $\calX_{\alpha} \rightarrow \calX$ is \'{e}tale (Theorem \ref{prescan}). Using this fact, we will show that the $\infty$-category $\RGeom_{\mathet}$ admits small colimits.
\end{itemize}

\begin{remark}
We will say that a geometric morphism of $\infty$-topoi $f^{\ast}: \calY \rightarrow \calX$ is
\'{e}tale if and only if its right adjoint $f_{\ast}: \calX \rightarrow \calY$ is \'{e}tale. We let
$\LGeom_{\mathet}$ denote the subcategory of $\LGeom$ spanned by the \'{e}tale geometric morphisms, so that there is a canonical equivalence $\RGeom_{\mathet} \simeq \LGeom_{\mathet}^{op}$.\index{not}{Ltopet@$\LGeom_{\mathet}$}
\end{remark}

Our first step is to obtain a more precise formulation of the universal property described in $(a)$:

\begin{definition}
Let $f^{\ast}: \calX \rightarrow \calY$ be a geometric morphism of $\infty$-topoi. Let
$U$ be an object of $\calX$ and $\alpha: 1_{\calY} \rightarrow f^{\ast} U$ a morphism in
$\calY$, where $1_{\calY}$ denotes a final object of $\calY$. We will say that
$\alpha$ {\it exhibits $\calY$ as a classifying $\infty$-topos for sections of $U$}
if, for every $\infty$-topos $\calZ$, the diagram
$$ \xymatrix{ \Fun^{\ast}( \calY, \calZ) \ar[r]^{\circ f^{\ast}} \ar[d]^{\phi} & \Fun^{\ast}(\calX, \calZ) \ar[d]^{\phi_0} \\
\calZ_{\ast} \ar[r] & \calZ }$$
is a homotopy pullback square of $\infty$-categories. Here $\calZ_{\ast}$ denotes the
$\infty$-category of pointed objects of $\calZ$ (that is, the full subcategory of
$\Fun(\Delta^1, \calZ)$ spanned by morphisms $f: Z \rightarrow Z'$ where $Z$ is a final object
of $\calZ$), and the morphisms $\phi$ and $\phi_0$ are given by evaluation on $\alpha$ and $U$, respectively.
\end{definition}

Let $\calX$ be an $\infty$-topos containing an object $U$. It follows immediately from the definition that
a classifying $\infty$-topos for sections of $U$ is uniquely determined up to equivalence, provided that it exists. For the existence, we have the following result:

\begin{proposition}\label{goodking}
Let $\calX$ be an $\infty$-topos containing an object $U$, let
$\pi_{!}: \calX_{/U} \rightarrow \calX$ be the projection map, and let
$\pi^{\ast}: \calX \rightarrow \calX_{/U}$ be a right adjoint to $\pi_{!}$. 
Let $1_{U}$ denote the identity map from $U$ to itself, regarded as a
$($final$)$ object of $\calX_{/U}$, and let $\alpha: 1_{U} \rightarrow \pi^{\ast} U$ be
adjoint to the identity map $\pi_{!} 1_{U} \simeq U$. Then $\alpha$ exhibits
$\calX_{/U}$ as a classifying $\infty$-topos for sections of $U$.
\end{proposition}

Before giving the proof of Proposition \ref{goodking}, we summarize some of the pleasant consequences.

\begin{corollary}\label{goodelk}
Let $\calX$ be an $\infty$-topos containing an object $U$, and let
$\pi^{\ast}: \calX \rightarrow \calX_{/U}$ be the corresponding \'{e}tale geometric morphism.
For every $\infty$-topos $\calZ$, composition with $\pi^{\ast}$ induces a left fibration
$$ \Fun^{\ast}( \calX_{/U}, \calZ ) \rightarrow \Fun^{\ast}( \calX, \calZ).$$
Moreover, the fiber over a geometric morphism $\phi^{\ast}: \calX \rightarrow \calZ$
is homotopy equivalent to the mapping space $\bHom_{\calZ}( 1_{\calZ}, \phi^{\ast} U )$. 
\end{corollary}

\begin{remark}\label{goodilk}
Corollary \ref{goodelk} implies in particular the existence of homotopy fiber sequences
$$ \bHom_{\calZ}( 1_{\calZ}, \phi^{\ast} U ) \rightarrow \bHom_{ \LGeom }( \calX_{/U}, \calZ)
\rightarrow \bHom_{\LGeom}(\calX, \calZ)$$
(where the fiber is taken over a geometric morphism $\phi^{\ast} \in \bHom_{\LGeom}(\calX, \calZ)$).

Suppose that $\calZ = \calX_{/V}$, and that $\phi^{\ast}$ is a right adjoint to the projection
$\calX_{/V} \rightarrow \calZ$. We then deduce the existence of a canonical homotopy equivalence
$$ \bHom_{\calX}(V,U) \simeq \bHom_{\calZ}( 1_{\calZ}, \phi^{\ast} U ) \simeq
\bHom_{ \LGeom_{\calX/}}( \calX_{/U}, \calX_{/V} ).$$
\end{remark}

\begin{remark}\label{pusha}
It follows from Remark \ref{goodilk} that if $f^{\ast}: \calX \rightarrow \calY$ is a geometric morphism
of $\infty$-topoi and $U \in \calX$ is an object, then the induced diagram
$$ \xymatrix{ \calX \ar[r] \ar[d] & \calY \ar[d] \\
\calX_{/U} \ar[r] & \calY_{/ f^{\ast} U } }$$
is a pushout square in $\LGeom$.
\end{remark}

\begin{corollary}\label{toadscan}
Suppose given a commutative diagram
$$ \xymatrix{ & \calY \ar[dr]^{g_{\ast}} & \\
\calX \ar[ur]^{f_{\ast}} \ar[rr]^{h_{\ast}} & & \calZ }$$
in $\LGeom^{op}$, where $g_{\ast}$ is \'{e}tale. Then $f_{\ast}$ is \'{e}tale if and only if $h_{\ast}$ is \'{e}tale.
\end{corollary}

\begin{proof}
The ``only if'' direction is obvious. To prove the converse, let us suppose
that $g_{\ast}$ and $h_{\ast}$ are both \'{e}tale, so that we have
equivalences $\calX \simeq \calZ_{/U}$ and $\calU \simeq \calZ_{/V}$ for some pair
of objects $U, V \in \calZ$. 
Using Remark \ref{goodilk}, we deduce that the morphism $f_{\ast}$ is determined by
a map $U \rightarrow V$ in $\calZ$, which we can identify with an object $\overline{V} \in \calY$
such that $\calX \simeq \calY_{/ \overline{V} }$. 
\end{proof}

\begin{remark}\label{postit}
Let $\calX$ be an $\infty$-topos. The projection map
$$p: \Fun( \Delta^1, \calX) \rightarrow \Fun( \{1\}, \calX) \simeq \calX$$
is a Cartesian fibration. Moreover, for every morphism $\alpha: U \rightarrow V$ in
$\calX$, the associated functor $\alpha^{\ast}: \calX^{/V} \rightarrow \calX^{/U}$ is an \'{e}tale 
geometric morphism of $\infty$-topoi, so that $p$ is classified by a functor
$\chi_0: \calX^{op} \rightarrow \LGeom_{\mathet}$. The functor $\chi_0$ carries the final object
of $\calX$ to an $\infty$-topos equivalent to $\calX$, and therefore factors as a composition
$$ \calX^{op} \stackrel{\chi}{\rightarrow} (\LGeom_{\mathet})_{\calX/} \rightarrow \LGeom_{\mathet}.$$
The argument of Remark \ref{goodilk} shows that $\chi$ is fully faithful, and
it follows immediately from the definitions that $\chi$ is essentially surjective.
Corollary \ref{toadscan} allows us to identify $(\LGeom_{\mathet})_{\calX/}$ with
the full subcategory of $\LGeom_{\calX/}$ spanned by the \'{e}tale geometric morphisms
$f^{\ast}: \calX \rightarrow \calY$. Consequently, we can regard $\chi$ as a fully faithful embedding
of $\calX$ into the $\infty$-category $(\LGeom^{op})_{/\calX}$ of $\infty$-topoi over $\calX$, whose essential image consists of those $\infty$-topoi which are \'{e}tale over $\calX$.
\end{remark}



\begin{proof}[Proof of Proposition \ref{goodking}]
Let $p: \calM \rightarrow \Delta^1$ be a correspondence from $\calX_{/U} \simeq \calM \times_{ \Delta^1} \{0\}$ to $\calX \simeq \calM \times_{ \Delta^1} \{1\}$ associated to the adjoint functors
$$ \Adjoint{ \pi_{!} }{ \calX_{/U} }{\calX}{\pi^{\ast}.}$$
Let $\alpha_0$ be a morphism from $1_{U} \in \calX_{/U}$ to $1_{\calX} \in \calX$
in $\calM$ (so that $\alpha_0$ is determined uniquely up to homotopy). We observe
that there is a retraction $r: \calM \rightarrow \calX_{/U}$ which restricts to $\pi^{\ast}$
on $\calX \subseteq \calM$, and we can identify $\alpha$ with $r( \alpha_0 )$. 

Let $\calZ$ be an arbitrary $\infty$-topos. Let $\calC$ be the full subcategory of
$\Fun( \calM, \calZ)$ spanned by those functors $F: \calM \rightarrow \calZ$ with the following properties:
\begin{itemize}
\item[$(a)$] The restriction $F | \calX_{/U}: \calX_{/U} \rightarrow \calZ$ preserves small colimits and finite limits.
\item[$(b)$] The functor $F$ is a left Kan extension of $F| \calX_{/U}$. In other words, 
$F$ carries $p$-Cartesian morphisms in $\calM$ to equivalences in $\calZ$.
\end{itemize}
Proposition \ref{lklk} implies that the restriction map
$\calC \rightarrow \Fun^{\ast}( \calX_{/U}, \calZ)$ is a trivial Kan fibration. Moreover,
this trivial Kan fibration has a section given by composition with $r$. It will therefore suffice
to show that the diagram
$$ \xymatrix{ \calC \ar[r] \ar[d] & \Fun^{\ast}( \calX, \calZ) \ar[d] \\
\calZ_{\ast} \ar[r] & \calZ }$$
is a homotopy pullback square. In other words, we wish to show that restriction along
$\alpha_0$ and the inclusion $\calX \subseteq \calM$ induce a categorical equivalence
$\calC \rightarrow \calZ_{\ast} \times_{ \calZ} \Fun^{\ast}( \calX, \calZ)$.

We define simplicial subsets $\calM'' \subseteq \calM' \subseteq \calM$ as follows:
\begin{itemize}
\item[$(i)$] Let $\calM'' \simeq \calX \coprod_{ \{1\} } \Delta^1$ be the union of
$\calX$ with the $1$-simplex of $\calM$ corresponding to the morphism $\alpha_0$.
\item[$(ii)$] Let $\calM'$ be the full subcategory of $\calM$ spanned by $\calX$
together with the object $1_U$.
\end{itemize}

We can identify $\calZ_{\ast} \times_{ \calZ} \Fun^{\ast}(\calX, \calZ)$ with the full subcategory
$\calC'' \subseteq \Fun( \calM'', \calZ)$ spanned by those functors $F$ satisfying the following conditions:
\begin{itemize}
\item[$(a')$] The restriction $F| \calX$ preserves small colimits and finite limits.
\item[$(b')$] The object $F( 1_{U} )$ is final in $\calZ$.  
\end{itemize}
Let $\calC'$ be the full subcategory of $\Fun( \calM', \calZ)$ spanned by those functors which
satisfy $(a')$ and $(b')$. To complete the proof, it will suffice to show that the restriction maps
$$ \calC \stackrel{u}{\rightarrow} \calC' \stackrel{v}{\rightarrow} \calC''$$
are trivial Kan fibrations.

We first show that $u$ is a trivial Kan fibration. In view of Proposition \ref{lklk}, it will suffice to prove the following:
\begin{itemize}
\item[$(\ast)$] A functor $F: \calM \rightarrow \calZ$ satisfies $(a)$ and $(b)$ if and only if
it satisfies $(a')$ and $(b')$, and $F$ is a right Kan extension of $F | \calM'$.
\end{itemize}
To prove the ``only if'' direction, let us suppose that $F$ satisfies $(a)$ and $(b)$. Without loss
of generality, we may suppose $F = F_0 \circ r$, where $F_0 = F | \calX_{/U}$. Then 
$F | \calX = F_0 \circ \pi^{\ast}$. Since $F_0$ and $\pi^{\ast}$ both 
preserve small colimits and finite limits, we deduce $(a')$. Condition $(b')$ is an immediate
consequence of $(a)$. We must show that $F$ is a right Kan extension of $F | \calX$.
Unwinding the definitions (and applying Corollary \ref{hollowtt}), we are reduced to showing that for every object $\overline{V} \in \calX_{/U}$ corresponding to a morphism $V \rightarrow U$ in $\calX$, the diagram
$$ \xymatrix{ F(\overline{V}) \ar[r] \ar[d] & F( V) \ar[d] \\
F(1_U) \ar[r] & F(U) }$$
is a pullback square in $\calZ$. Since $F = F_0 \circ r$ and $F_0$ preserves finite limits,
it suffices to show that the square
$$ \xymatrix{ \overline{V} \ar[r] \ar[d] & \pi^{\ast} V \ar[d] \\
1_{U} \ar[r] & \pi^{\ast} U }$$
is a pullback square in $\calX_{/U}$. In view of Proposition \ref{needed17}, it suffices
to observe that the square
$$ \xymatrix{ V \ar[r] \ar[d] & V \times U \ar[d] \\
U \ar[r] & U \times U }$$
is a pullback in $\calX$.

Now let us suppose that $F$ is a right Kan extension of $F_1 = F | \calM'$, and that
$F_1$ satisfies conditions $(a')$ and $(b')$. We first claim that $F$ satisfies $(b)$.
In other words, we claim that for every object $V \in \calX$, the canonical map
$F( \pi^{\ast} V) \rightarrow F(V)$ is an equivalence. Consider the diagram
$$ \xymatrix{ F( \pi^{\ast} V) \ar[r] \ar[d] & F(V \times U) \ar[d] \ar[r] & F(V) \ar[d] \\
F(1_U) \ar[r] & F(U) \ar[r] & F(1_{\calX}). }$$
Since $F$ is a right Kan extension of $F_1$, the left square is a pullback.
Since $F_1$ satisfies $(a)$, the right square is a pullback. Therefore the outer square is a pullback. Condition $(b')$ implies that the lower horizontal composition is an equivalence, so the upper horizontal composition is an equivalence as well.

We now prove that $F$ satisfies $(a)$. Condition $(b')$ implies that the functor
$F_0 = F | \calX_{/U}$ preserves final objects. It will therefore suffice to show that
$F_0$ preserves small colimits and pullback squares.
Since $F$ is a right Kan extension of $F_1$, the
functor $F_0$ can be described by the formula
$$ V \mapsto F( \pi_{!} V) \times_{ F(U) } F(1_U).$$
It therefore suffices to show that the functors $\pi_{!}$, $F|\calX$, and
$\bigdot \times_{ F(U)} F(1_U)$ preserve small colimits and pullback squares.
For $\pi_{!}$, this follows from Propositions \ref{needed17} and \ref{goeselse}.
For $F|\calX$, we invoke assumption $(a')$. For the functor $\bigdot \times_{ F(U)} F(1_{U})$, we
invoke our assumption that $\calZ$ is an $\infty$-topos (so that colimits in $\calZ$ are universal).
This completes the verification that $u$ is a trivial Kan fibration.

To complete the proof, we must show that the functor $v$ is a trivial Kan fibration.
We note that $v$ fits into a pullback diagram
$$ \xymatrix{ \calC' \ar[r]^{v} \ar[d] & \calC'' \ar[d] \\
\Fun( \calM', \calZ) \ar[r]^{v'} & \Fun( \calM'', \calZ ). }$$
It will therefore suffice to show that $v'$ is a trivial Kan fibration. Since
$\calZ$ is an $\infty$-category, we need only show that the inclusion
$\calM'' \subseteq \calM'$ is a categorical equivalence of simplicial sets. This is a special case of
Proposition \ref{simplexplay}.
\end{proof}


We next establish the recognition principle promised in $(b)$:

\begin{proposition}\label{pushpo}
Let $f^{\ast}: \calX \rightarrow \calY$ be a geometric morphism of $\infty$-topoi. Then
$f^{\ast}$ is \'{e}tale if and only if the following conditions are satisfied:
\begin{itemize}
\item[$(1)$] The functor $f^{\ast}$ admits a left adjoint $f_{!}$ $($in view of Corollary \ref{adjointfunctor}, this is equivalent to the assumption that $f^{\ast}$ preserves small limits$)$.
\item[$(2)$] The functor $f_{!}$ is conservative. That is, if $\alpha$ is a morphism in $\calY$ such
that $f_{!} \alpha$ is an equivalence in $\calX$, then $\alpha$ is an equivalence in $\calY$.
\item[$(3)$] For every morphism $X \rightarrow Y$ in $\calX$, every object $Z \in \calY$, and
every morphism $f_{!} Z \rightarrow Y$, the induced diagram
$$ \xymatrix{ f_{!} ( f^{\ast} X \times_{ f^{\ast} Y } Z) \ar[r] \ar[d] & f_{!} Z \ar[d] \\
X \ar[r] & Y }$$
is a pullback square in $\calX$.
\end{itemize}
\end{proposition}

\begin{remark}
Condition $(3)$ of Proposition \ref{pushpo} can be regarded as a push-pull formula: it provides a canonical equivalence
$$ f_{!}( f^{\ast} X \times_{ f^{\ast} Y} Z) \simeq X \times_{Y} f_{!} Z.$$
In particular, when $Y$ is final in $\calX$, we have an equivalence
$f_{!}( f^{\ast} X \times Z) \simeq X \times f_{!} Z$: in other words, the functor
$f_{!}$ is ``linear'' with respect to the action of $\calX$ on $\calY$.
\end{remark}

\begin{proof}[Proof of Proposition \ref{pushpo}]
Suppose first that $f^{\ast}$ is an \'{e}tale geometric morphism. Without loss of generality, we may suppose that $\calY = \calX_{/U}$, and that $f^{\ast}$ is right adjoint to the forgetful functor
$f_{!}: \calX_{/U} \rightarrow \calX$. Assertions $(1)$ and $(2)$ are obvious, and assertion
$(3)$ follows from the observation that, for every diagram
$$ X \rightarrow Y \leftarrow Z \rightarrow U,$$
the induced map $(X \times U) \times_{ Y \times U } Z \rightarrow X \times_{Y} Z$
is an equivalence in $\calX$.

For the converse, let us suppose that $(1)$, $(2)$ and $(3)$ are satisfied. We wish to show that
$f^{\ast}$ is \'{e}tale. Let $U = f_{!} 1_{\calY}$. Let 
$F$ denote the composition $\calY \simeq \calY_{/1_{\calY}} \stackrel{f_{!}}{\rightarrow} \calX_{/U}$. 
To complete the proof, it will suffice to show that $F$ is an equivalence of $\infty$-categories. Proposition \ref{curpse} implies that $F$ admits a right adjoint $G$, given by the formula
$$(X \rightarrow U) \mapsto f^{\ast} X \times_{ f^{\ast} U} 1_{\calY}.$$
Assumption $(3)$ guarantees that the counit map $v: FG \rightarrow \id_{\calX_{/U}}$ is
an equivalence. To complete the proof, it suffices to show that for each $Y \in \calY$, the
unit map $u_Y: Y \rightarrow GFY$ is an equivalence. The map $Fu_Y: FY \rightarrow FGFY$
has a left homotopy inverse (given by $v_{FY}$) which is an equivalence, so that $F u_{Y}$ is
an equivalence. It follows that $f_{!} u_Y$ is an equivalence, so that $u_Y$ is an equivalence
by virtue of assumption $(2)$. Thus $G$ is a homotopy inverse to $F$, so that $F$ is an equivalence of $\infty$-categories as desired.
\end{proof}

Our final goal in this section is to prove the following result:

\begin{theorem}\label{prescan}
The $\infty$-category $\RGeom_{\mathet}$ admits small colimits, and the inclusion
$\RGeom_{\mathet} \subseteq \RGeom$ preserves small colimits.
\end{theorem}

The proof of Theorem \ref{prescan} is rather technical and will occupy our attention for the remainder of this section. However, the analogous result
is elementary if we work with $\infty$-topoi which are assume to be \'{e}tale over a fixed
base $\calX$. In this case, Theorem \ref{prescan} can be reduced (with the aid of Remark \ref{postit}) to the following assertion:

\begin{proposition}\label{toadsteal}
Let $\calX$ be an $\infty$-topos, and let $\chi: \calX \rightarrow \LGeom^{op}_{/ \calX}$ be
the functor of Remark \ref{postit}. Then $\chi$ preserves small colimits.
\end{proposition}

\begin{proof}
Combine Propositions \ref{needed17}, \ref{colimtopoi}, and Theorem \ref{charleschar}.
\end{proof}

\begin{proof}[Proof of Theorem \ref{prescan}]
As a first step, we establish the following:
\begin{itemize}
\item[$(\ast)$] Suppose given a small diagram $p: K \rightarrow \LGeom^{op}$ 
and a geometric morphism of $\infty$-topoi $\phi_{\ast}: \colim(p) \rightarrow \calY$. Suppose
further that for each vertex $v$ in $K$, the induced geometric morphism
$\phi(v)_{\ast}: p(v) \rightarrow \calY$ is \'{e}tale. Then $\phi_{\ast}$ is \'{e}tale.
\end{itemize}
To prove $(\ast)$, we note that $\phi_{\ast}$ determines a functor
$\overline{p}: K \rightarrow \LGeom^{op}_{/\calY}$ lifting $p$. Since
each $\phi(v)_{\ast}$ is \'{e}tale, Remark \ref{postit} implies that 
$\overline{p}$ factors as a composition
$$ K \stackrel{q}{\rightarrow} \calY \stackrel{\chi}{\rightarrow} \LGeom^{op}_{/\calY}.$$
Let $U \in \calY$ be a colimit of the diagram $q$. Then Corollary \ref{toadscan} implies that
$\colim(p) \simeq \calY_{/U}$, so that $\phi_{\ast}$ is \'{e}tale as desired.

We now return to the proof of Theorem \ref{prescan}. Using Proposition \ref{appendixdiagram} and its proof, we can reduce the proof to the following special cases:
\begin{itemize}
\item[$(a)$] The $\infty$-category $\LGeom^{op}_{\mathet}$ admits small coproducts, and
the inclusion $\LGeom^{op}_{\mathet} \subseteq \LGeom^{op}$ preserves small coproducts.
\item[$(b)$] The $\infty$-category $\LGeom^{op}_{\mathet}$ admits coequalizers, and the inclusion
$\LGeom^{op}_{\mathet} \subseteq \LGeom^{op}$ preserves coequalizer diagrams.
\end{itemize}

We first prove $(a)$. In view of $(\ast)$, it will suffice to prove the following:
\begin{itemize}
\item[$(a')$] Let $\{ \calX_{\alpha} \}$ be a small collection of $\infty$-topoi, and let
$\calX$ be their coproduct in $\LGeom^{op}$ (so that we have an equivalence of $\infty$-categories $\calX \simeq \prod_{\alpha} \calX_{\alpha}$). Then each of the associated geometric morphisms
$\phi_{\alpha}^{\ast}: \calX \rightarrow \calX_{\alpha}$ is \'{e}tale.
\end{itemize}
To prove $(a')$, we may assume without loss of generality that $\calX = \prod_{\alpha} \calX_{\alpha}$
and that $\phi_{\alpha}^{\ast}$ is given by projection onto the corresponding factor. The desired
result then follows from the criterion of Proposition \ref{pushpo} (more concretely: let
let $U \in \calX$ be an object whose image in $\calX_{\alpha}$ is
a final object $U_{\alpha} \in \calX_{\alpha}$, and whose image in $\calX_{\beta}$ is an initial object 
$U_{\beta} \in \calX_{\beta}$ for $\beta \neq \alpha$. Then $\calX_{/U} \simeq \prod_{\beta} (\calX_{\beta})_{/U_{\beta}} \simeq \calX_{\alpha}$.) 

To prove $(b)$, we can again invoke $(\ast)$ to reduce to the following assertion:
\begin{itemize}
\item[$(b')$] Suppose given a diagram
$$\xymatrix{ \calY \ar@<.4ex>[r] \ar@<-.4ex>[r] & \calX_0}.$$
in $\LGeom_{\mathet}^{op}$, having colimit $\calX$ in $\LGeom^{op}$. Then
the induced geometric morphism $\phi_{\ast}: \calX_0 \rightarrow \calX$ is \'{e}tale.
\end{itemize}

To prove $(b')$, we identify the diagram in question with a functor
$p: \calI \rightarrow \LGeom^{op}$, and $\calI$ with the subcategory of
$\Nerve(\cDelta)^{op}$ spanned by the objects $\{ [0], [1] \}$ and injective maps of linearly ordered sets.
Let $\calX_{\bigdot}$ be the simplicial object of $\LGeom^{op}$ given by left Kan extension along
the inclusion $\calI \subseteq \Nerve(\cDelta)^{op}$, so that each $\calX_{n}$ is equivalent to a coproduct (in $\LGeom^{op}$) of $\calX_0$ with $n$ copies of $\calY$. Using $(\ast)$ and Corollary \ref{toadscan}, we deduce that $\calX_{\bigdot}$ is a simplicial object in $\LGeom^{op}_{\mathet}$.
Consequently, assertion $(b')$ is an immediate consequence of Lemma \ref{kan0} and
the following:
\begin{itemize}
\item[$(b'')$] Let $\calX_{\bigdot}$ be a simplicial object of $\LGeom^{op}_{\mathet}$, and let
$\calX$ be its geometric realization in $\LGeom^{op}$. Then the induced geometric morphism
$\phi_{\ast}: \calX_0 \rightarrow \calX$ is \'{e}tale.
\end{itemize}

The proof of $(b'')$ is based on the following Lemma, whose proof we defer until the end of this section:

\begin{lemma}\label{santan}
Suppose given a simplicial object $\calX_{\bigdot}$ in $\LGeom^{op}_{\mathet}$. Then
there exists a morphism of simplicial objects $\calX_{\bigdot} \rightarrow \calX'_{\bigdot}$
of $\LGeom^{op}_{\mathet}$ with the following properties:
\begin{itemize}
\item[$(1)$] The induced map $\calX_0 \rightarrow \calX'_0$ is an equivalence of $\infty$-topoi.
\item[$(2)$] The simplicial object $\calX'_{\bigdot}$ is a groupoid object in $\LGeom^{op}$.
\item[$(3)$] The induced map of geometric realizations $($in $\LGeom^{op}${}$)$ is an equivalence
$| \calX_{\bigdot} | \rightarrow | \calX'_{\bigdot} |$.
\end{itemize}
\end{lemma}

Using Lemma \ref{santan}, we can reduce the proof of $(b'')$ to the special case where
$\calX_{\bigdot}$ is a groupoid object of $\LGeom^{op}$. The diagram
$$ \Nerve( \cDelta )^{op} \stackrel{\calX_{\bigdot}}{\rightarrow} \LGeom^{op}_{\mathet} 
\subseteq \widehat{\Cat}_{\infty}^{op}$$
is classified by a Cartesian fibration $q: \calZ \rightarrow \Nerve( \cDelta)^{op}$.
Here we can identify $\calX_{n}$ with the fiber $\calZ_{[n]} = \calZ \times_{ \Nerve(\cDelta)^{op} } \{ [n] \}$, and every map of linearly ordered sets $\alpha: [m] \rightarrow [n]$ induces a geometric morphism $\alpha^{\ast}: \calZ_{[m]} \rightarrow \calZ_{[n]}$. Since the geometric morphism $\alpha^{\ast}$ is \'{e}tale, it admits a left adjoint $\alpha_{!}$, so that $q$ is also a coCartesian fibration (Corollary \ref{grutt1}).

It follows from Propositions \ref{colimtopoi} and \ref{charcatlimit} that we can identify
$\calX$ with the full subcategory of $\Fun_{ \Nerve(\cDelta)^{op} } ( \Nerve(\cDelta)^{op}, \calZ)$ spanned by the Cartesian sections of $q$; under this identification, the pullback functor
$\phi^{\ast}$ corresponds to the functor $\calX \rightarrow \calZ_{[0]} \simeq \calX_0$ given by evaluation at $[0]$. 

Let $1_{\calX}$ denote a final object of $\calX$, which we regard as a section of $q$. 
Let $T: \Nerve(\cDelta)^{op} \rightarrow \Nerve(\cDelta)^{op}$ denote the shift functor
$[n] \mapsto [n] \star [0]$, and let $\beta_0: T \rightarrow \id_{ \Nerve(\cDelta)^{op} }$ denote the evident natural transformation. Let $\beta: (1_{\calX} \circ T) \rightarrow U_{\bigdot}$ be a
natural transformation in $\Fun( \Nerve(\cDelta)^{op}, \calZ)$ lifting $\beta$ which is
$q$-coCartesian. Since $\calX_{\bigdot}$ is a groupoid object of $\LGeom^{op}_{\mathet}$, we deduce that $U_{\bigdot}$ is a Cartesian section of $q$, which we can identify with an object of
$\calX$.

Let $S = \Nerve(\cDelta)^{op} \times \Delta^1$, so that $\beta_0$ defines a map
$S \rightarrow \Nerve(\cDelta)^{op}$. Let $\calZ' = \calZ \times_{ \Nerve(\cDelta)^{op} } S$
and let $\beta_{S} = \beta$, regarded as a section of of the projection $q': \calZ' \rightarrow S$.
Let $\calZ'' = {\calZ'}^{/\beta_S}$ (see \S \ref{consweet} for an explanation of this notation). 
Let $q'': \calZ'' \rightarrow S$. The fibers of $q''$ can be described as follows:
\begin{itemize}
\item The fiber of $q''$ over $([n], 0)$ can be identified with
$\calZ_{[n+1]}^{/ 1_{\calZ_{[n+1]}}} \simeq \calZ_{[n+1]}$.
\item The fiber of $q''$ over $([n], 1)$ can be identified with
$\calZ_{[n]}^{/U_{n}} \simeq \calZ_{[n+1]}$. 
\end{itemize}
Proposition \ref{colimfam} implies that the projection $q'': \calZ'' \rightarrow S$ is a coCartesian fibration, classified by a map $\chi: S \rightarrow \widehat{\Cat}_{\infty}$. The above description
shows that $\chi$ can be regarded as an equivalence from $\chi^0 = \chi | \Nerve(\cDelta)^{op} \times \{0\}$ to $\chi^1 = \chi | \Nerve(\cDelta)^{op} \times \{1\}$ in the $\infty$-category of simplicial objects
of $\widehat{\Cat}_{\infty}$. Moreover, the functor $\chi^0$ classifies the pullback of
the coCartesian fibration $q$ by the translation map $T: \Nerve(\cDelta)^{op} \rightarrow
\Nerve(\cDelta^{op})$, so that $\chi^0$ and $\chi^1$ factor through $\LGeom^{op}_{\mathet}$. 
Lemma \ref{bclock} implies that the colimit of $\chi^0$ (hence also of $\chi^1$) in
$\LGeom^{op}$ is canonically equivalent to $\calX_0$. On the other hand, Propositions \ref{charcatlimit} and \ref{colimtopoi} allow us to identify $\colim( \chi^1)$ with the $\infty$-category
of Cartesian sections of the projection $\calZ'' \times_{ S} ( \Nerve(\cDelta)^{op} \times \{1\} )
\rightarrow \Nerve(\cDelta)^{op}$, which is isomorphic to $\calX^{/U_{\bigdot}}$ as a simplicial set.
We now complete the proof by observing that the resulting identification
$\calX_0 \simeq \calX^{/U_{\bigdot}}$ is compatible with the projection $\phi_{\ast}: \calX_0 \rightarrow \calX$.
\end{proof}

The remainder of this section is devoted to the proof of Lemma \ref{santan}. We first need to introduce a bit of notation. We begin with a few remarks about the behavior of $\infty$-topoi under change of universe.

\begin{notation}
Let $\calX$ be an $\infty$-topos and $\calC$ an arbitrary $\infty$-category. We let
$\Shv_{\calC}(\calX)$ denote the full subcategory of $\Fun( \calX^{op}, \calC)$ spanned by those functors which preserve small limits.\index{not}{ShvCX@$\Shv_{\calC}(\calX)$} We will refer to
$\Shv_{\calC}(\calX)$ as the {\it $\infty$-category of $\calC$-valued sheaves on $\calX$}.
\end{notation}

\begin{remark}\label{quest}
Let $\calX$ be an $\infty$-topos. Proposition \ref{representable} implies that the Yoneda
embedding $\calX \rightarrow \Shv_{\SSet}(\calX)$ is an equivalence; in other words, we can identify $\calX$ with the $\infty$-category of sheaves of (small) spaces on itself. Let $\widehat{\SSet}$ denote the $\infty$-category of spaces which belong to some larger universe $\calU$. We claim the following:
\begin{itemize}
\item[$(a)$] The $\infty$-category $\Shv_{ \widehat{\SSet} }( \calX)$ can be regarded
as an $\infty$-topos in $\calU$.
\item[$(b)$] The inclusion $\Shv_{\SSet}(\calX) \subseteq \Shv_{\widehat{\SSet}}(\calX)$ preserves small colimits.
\end{itemize}
To prove $(a)$, let us suppose that $\calX = S^{-1} \calP(\calC)$, where
$\calC$ is a small $\infty$-category and $S$ is a strongly saturated class of morphisms
in $\calP(\calC)$, which is stable under pullbacks and of small generation.
Theorem \ref{charpresheaf} and Proposition \ref{unichar} allow us to identify
$\Shv_{\widehat{\SSet}}(\calX)$ with $S^{-1} \widehat{\calP}(\calC)$, where
$\widehat{\calP}(\calC)$ denotes the presheaf $\infty$-category
$\Fun( \calC^{op}, \widehat{\SSet})$. Let $\widehat{S}$ denote the strongly saturated class of morphisms of $\widehat{\calP}(\calC)$ generated by $S$. Then $\widehat{S}$ is of small generation (and therefore of $\calU$-small generation); to complete the proof of $(a)$ it will suffice to show that $\widehat{S}$ is stable under pullbacks. 

Let $\widehat{\calP}(\calC)^{0}$ denote the full subcategory of $\widehat{\calP}(\calC)$ spanned by those objects $X$ with the following property: 
\begin{itemize}
\item[$(\ast)$] Let
$$ \xymatrix{ Y \ar[d]^{f} \ar[r] & Y' \ar[d]^{f'} \\
X \ar[r] & X' }$$
be a pullback diagram in $\widehat{\calP}(\calC)$. If $f' \in S$, then $f \in \widehat{S'}$.
\end{itemize}
Since colimits in $\widehat{\calP}(\calC)$ are universal, the subcategory
$\widehat{\calP}(\calC)^{0}$ is stable under $\calU$-small colimits in
$\widehat{\calP}(\calC)$. Moreover, since $S$ is stable under pullbacks
in $\calP(\calC)$ (and since the inclusion $\calP(\calC) \subseteq \widehat{\calP}(\calC)$ is fully faithful), the $\infty$-category $\widehat{\calP}(\calC)^{0}$ contains $\calP(\calC)$.
Since $\widehat{\calP}(\calC)$ is generated (under $\calU$-small colimits) by the essential image of the Yoneda embedding $\calC \rightarrow \calP(\calC)$, we conclude that
$\widehat{\calP}(\calC)^{0} = \widehat{\calP}(\calC)$. 

We now let $S'$ denote the collection of all morphisms in $\widehat{\calP}(\calC)$ such that, for every pullback diagram
$$ \xymatrix{ Y \ar[d]^{f} \ar[r] & Y' \ar[d]^{f'} \\
X \ar[r] & X' }$$
in $\calP(\calC)$, if $f' \in S'$ then $f \in \widehat{S}$. The above argument shows that
$S \subseteq S'$. Since $S'$ is strongly saturated, we conclude that $\widehat{S} \subseteq S'$, so that $\widehat{S}$ is stable under pullbacks as desired. This completes the proof of $(a)$.

To prove $(b)$, it will suffice to show that the composite map
$$ \calP(\calC) \rightarrow S^{-1} \calP(\calC) \rightarrow \widehat{S}^{-1} \widehat{\calP}(\calC)$$
preserves small colimits. We can rewrite this as the composition of a pair of functors
$$ \calP(\calC) \stackrel{i}{\rightarrow} \widehat{\calP}(\calC)
\stackrel{L}{\rightarrow} \widehat{S}^{-1} \widehat{\calP}(\calC).$$
The functor $L$ is left adjoint to the inclusion of
$\widehat{S}^{-1} \widehat{\calP}(\calC)$ into $\widehat{\calP}(\calC)$, and therefore
preserves all $\calU$-small colimits. It therefore suffices to show that the inclusion
$i: \Fun( \calC^{op}, \SSet) \rightarrow \Fun( \calC^{op}, \widehat{\SSet})$ preserves small colimits.
In view of Proposition \ref{limiteval}, it will suffice to prove the inclusion
$i_0: \SSet \rightarrow \widehat{\SSet}$ preserves small colimits. We note that
$i_0$ is an equivalence from $\SSet$ to the full subcategory
$\widehat{\SSet}^0 \subseteq \widehat{\SSet}$ spanned by the essentially small spaces.
It now suffices to observe that the collection of essentially small spaces is stable under small colimits
(this follows from Corollaries \ref{apegrape} and \ref{tyrmyrr}).
\end{remark}

\begin{remark}\label{postquest}
Let $\calU$ be a universe as in Example \ref{quest}, let $f^{\ast}: \calX \rightarrow \calY$ be a geometric morphism of $\infty$-topoi, and let
$\widehat{f}_{\ast}: \Shv_{ \widehat{\SSet}}(\calX) \rightarrow \Shv_{ \widehat{\SSet}}(\calY)$ be
given by composition with $f^{\ast}$. Then $\widehat{f}_{\ast}$ can be identified with a geometric morphism in the universe $\calU$. To prove this, let $\kappa$ denote the regular cardinal
in the universe $\calU$ such that small sets (in our original universe) can be identified with
$\kappa$-small sets in $\calU$. It follows from Corollary \ref{indpr} that we can identify
$\Shv_{\widehat{\SSet}}(\calX)$ and $\Shv_{\widehat{\SSet}}(\calX)$ with
$\widehat{\Ind}_{\kappa}(\calX)$ and $\widehat{\Ind}_{\kappa}(\calY)$, respectively.
Proposition \ref{adjobs} implies that $\widehat{f}_{\ast}$ admits a left adjoint
$\widehat{f}^{\ast}$ which fits into a commutative diagram
$$ \xymatrix{ \calX \ar[r]^{f^{\ast}} \ar[r] \ar[d] & \calY \ar[d] \\
\Shv_{\widehat{\SSet}}(\calX) \ar[r]^{ \widehat{f}^{\ast} } & \Shv_{ \widehat{\SSet}}(\calY). }$$
To complete the proof, it will suffice to show that $\widehat{f}^{\ast}$ is left exact.
Since $f^{\ast}$ preserves final objects, the functor $\widehat{f}^{\ast}$ preserves final objects as well.
It therefore suffices to show that $\widehat{f}^{\ast}$ preserves pullback diagrams. Using Proposition \ref{urgh1} and Example \ref{tucka}, we conclude that every pullback diagram
in $\Shv_{ \widehat{\SSet}}(\calX)$ can be obtained as a $\calU$-small, $\kappa$-filtered colimit of pullback diagrams in $\calX$. The desired result now follows from the assumption that $f^{\ast}$ is left exact, and the observation that the class of pullback diagrams in $\Shv_{ \widehat{\SSet}}(\calY)$ is stable under $\calU$-small filtered colimits (Example \ref{tucka}).
\end{remark}

For the remainder of this section, we fix a larger universe $\calU$. Let
$\widehat{\SSet}$ denote the $\infty$-category of $\calU$-small spaces.

\begin{notation}
Let $F: \LGeom \rightarrow \widehat{\SSet}$ be a functor. For every $\infty$-topos
$\calX$, we let $F_{\calX}: \calX^{op} \rightarrow \widehat{\SSet}$ denote the composition
$$ \calX^{op} \simeq \LGeom_{\mathet}^{\calX/} \rightarrow \LGeom \stackrel{F}{\rightarrow} \widehat{\SSet}.$$
We will say that $F_{\calX}$ is a {\it sheaf} if, for every $\infty$-topos $\calX$, the functor
$F_{\calX}$ preserves small limits. We let $\widehat{\Shv}( \LGeom^{op} )$ denote the full subcategory of $\Fun( \LGeom, \widehat{\SSet} )$ spanned by the sheaves.
\end{notation}

\begin{example}\label{stup}
Let $\calX$ be an $\infty$-topos, and let $e_{\calX}: \LGeom \rightarrow \widehat{\SSet}$ be the functor
represented by $\calX$. Proposition \ref{toadsteal} implies that $e_{\calX}$ belongs to
$\widehat{\Shv}( \LGeom^{op} )$. We will say that a sheaf $F \in \widehat{\Shv}( \LGeom^{op} )$ is
{\it representable} if $F \simeq e_{\calX}$ for some $\infty$-topos $\calX$.
\end{example}

\begin{lemma}\label{kumba}
The $\infty$-category $\widehat{\Shv}( \LGeom^{op} )$ is an $\infty$-topos in the universe $\calU$.
Moreover, for every $\infty$-topos $\calX$, the restriction functor
$F \mapsto F_{\calX}$ determines a functor $\widehat{\Shv}(\LGeom^{op}) \rightarrow
\Shv_{ \widehat{\SSet} }(\calX)$ which preserves $\calU$-small colimits and finite limits.
\end{lemma}

\begin{proof}
Let $\Fun^{\mathet}( \Delta^1, \LGeom)$ denote the full subcategory $\Fun( \Delta^1, \LGeom)$ spanned by the \'{e}tale morphisms, and let
$e: \Fun^{\mathet}( \Delta^1, \LGeom) \rightarrow \LGeom$ be given by evaluation at the vertex
$\{0\} \in \Delta^1$. Since the collection of \'{e}tale morphisms in $\LGeom$ is stable under pushouts (Remark \ref{pusha}), the map $e$ is a coCartesian fibration.

We define a simplicial set $\calK$ equipped with a projection
$p: \calK \rightarrow \LGeom$ so that the following universal property is satisfied: for every
simplicial set $K$, we have a natural bijection
$$ \Hom_{\Ind(\calG^{op})}(K, \calK) = \Hom_{\sSet}(
K \times_{ \LGeom } \Fun^{\mathet}(\Delta^1, \LGeom), \widehat{\SSet} ).$$
Then $\calK$ is an $\infty$-category, whose objects can be identified with pairs
$(\calX, F_{\calX})$, where $\calX$ is an $\infty$-topos
and $F_{\calX}: \LGeom^{\calX/}_{\mathet} \rightarrow \widehat{\SSet}$ is a functor.
It follows from Corollary \ref{skinnysalad} that the projection $p$ is a Cartesian fibration, and that
a morphism $(\calX,F_{\calX}) \rightarrow (\calY, F_{\calY})$ is $p$-Cartesian if and only if, for every object $U \in \calX$, the canonical map
$F_{\calX}( \calX_{/U}) \rightarrow F_{\calY}( \calY_{/f^{\ast} U})$ is a homotopy equivalence, where $f^{\ast}$ denotes the underlying geometric morphism from $\calX$ to $\calY$.

Let $\calK_0$ denote the full subcategory of $\calK$ spanned by pairs $(\calX, F_{\calX})$ where
the functor $F_{\calX}$ preserves small limits. It follows from the above that the Cartesian
fibration $p$ restricts to a Cartesian fibration $p_0: \calK_0 \rightarrow \LGeom$ (with the same class of Cartesian morphisms). The fiber of $\calK_0$ over an object $\calX \in \LGeom$ can be identified with
$\Shv_{\widehat{\SSet}}( \calX)$, which is an $\infty$-topos in the universe $\calU$ (Remark \ref{quest}). Moreover, to every geometric morphism $f^{\ast}: \calX \rightarrow \calY$ in
$\LGeom$, the Cartesian fibration $p_0$ associates the pushforward functor $\widehat{f}_{\ast}: \Shv_{ \widehat{\SSet} }(\calY) \rightarrow \Shv_{ \widehat{\SSet} }(\calX)$ given by composition with $f^{\ast}$. It follows from Remark \ref{postquest} that $\widehat{f}_{\ast}$ admits a left adjoint
$\widehat{f}^{\ast}$, and that $\widehat{f}^{\ast}$ is left exact. We may summarize the situation
by saying that $p_0$ is a topos fibration (see Definition \ref{skuzz}); in particular, $p_0$ is a coCartesian fibration.

Let $\calY = \Fun_{ \LGeom}( \LGeom, \calK_0)$ denote the $\infty$-category of sections of $p_0$. Unwinding the definition, we can identify $\calY$ with the $\infty$-category
$\Fun( \Fun^{\mathet}( \Delta^1, \LGeom),\widehat{\SSet})$.
Let $\LGeom'$ denote the essential image of the (fully faithful) diagonal embedding
$\LGeom \rightarrow \Fun( \Delta^1, \LGeom)$. Consider the following conditions on a section $s: \LGeom \rightarrow \calK_0$ of $p_0$:
\begin{itemize}
\item[$(a)$] The functor $s$ carries \'{e}tale morphisms in $\LGeom$ to $p$-coCartesian morphisms in $\calX$.
\item[$(b)$] Let $S: \Fun^{\mathet}( \Delta^1, \LGeom) \rightarrow \widehat{\SSet}$ be the functor
corresponding to $s$. Then, for every commutative diagram
$$ \xymatrix{ & \calY \ar[dr] & \\
\calX \ar[ur] \ar[rr] & & \calZ }$$
of \'{e}tale morphisms in $\LGeom$, the induced map
$S(\calX \rightarrow \calZ) \rightarrow S(\calY \rightarrow \calZ)$ is an equivalence in $\widehat{\SSet}$.
\item[$(c)$] For every \'{e}tale morphism $f^{\ast}: \calX \rightarrow \calY$ in $\LGeom$, the canonical
map $S( \id_{\calX} ) \rightarrow S(f^{\ast})$ is an equivalence in $\widehat{\SSet}$.
\item[$(d)$] The functor $S$ is a left Kan extension of $S | \LGeom'$.
\end{itemize}
Unwinding the definitions, we see that $(a) \Leftrightarrow (b) \Rightarrow (c) \Leftrightarrow (d)$.
Moreover, the implication $(c) \Rightarrow (b)$ follows by a two-out-of-three argument.
Let $\calY'$ denote the full subcategory of $\calY$ spanned by those sections which
satisfy the equivalent conditions $(a)$ through $(d)$; it follows from Proposition \ref{lklk}
that composition with the diagonal embedding $\LGeom \rightarrow \Fun^{\mathet}( \Delta^1, \LGeom)$ induces an equivalence $\theta: \calY' \rightarrow \widehat{\Shv}(\LGeom^{op})$. The desired result now follows from Proposition \ref{prestorkus}.
\end{proof}

To prove Lemma \ref{santan}, we need a criterion which will allow us to detect \'{e}tale geometric morphisms of $\infty$-topoi in terms of the functors that they represent. To formulate this criterion, we introduce a bit of temporary terminology.

\begin{definition}\label{squee}
Let $\alpha: F \rightarrow G$ be a morphism in $\widehat{\Shv}(\LGeom^{op})$. We will say that
$\alpha$ is {\it universal} if, for every geometric morphism of $\infty$-topoi
$f^{\ast}: \calX \rightarrow \calY$, the induced diagram
$$ \xymatrix{ \widehat{f}^{\ast} F_{\calX} \ar[r] \ar[d] & \widehat{f}^{\ast} G_{\calX} \ar[d] \\
F_{\calY} \ar[r] & G_{\calY} }$$
is a pullback square in $\Shv_{ \widehat{\SSet}}( \calY)$. (Here $\widehat{f}^{\ast}$ denotes
the geometric morphism described in Remark \ref{postquest}).
\end{definition}

\begin{remark}
The collection of universal morphisms in $\widehat{\Shv}(\LGeom^{op})$ is stable under pullbacks
and composition, and contains every equivalence in $\widehat{\Shv}(\LGeom^{op})$.
\end{remark}

\begin{remark}\label{sapper}
Let $p: K \rightarrow \widehat{\Shv}(\LGeom^{op})$ be a small diagram having a colimit $F$.
Assume that for every edge $v \rightarrow v'$ of $K$, the induced map $p(v) \rightarrow p(v')$
is universal. Then each of the induced maps $p(v) \rightarrow F$ is universal. This follows immediately from Theorem \ref{charleschar}.
\end{remark}

\begin{lemma}\label{sabus}
Let $\alpha: F \rightarrow G$ be a morphism in $\widehat{\Shv}(\LGeom^{op})$, and assume
that $G$ is representable by an $\infty$-topos $\calX$. Then $F$ is representable by an
$\infty$-topos \'{e}tale over $\calX$ if the following conditions are satisfied:
\begin{itemize}
\item[$(1)$] The morphism $\alpha$ is universal $($in the sense of Definition \ref{squee}$)$.
\item[$(2)$] For every $\infty$-topos $\calY$, the homotopy fibers of the induced map
$F(\calY) \rightarrow G(\calY)$ are essentially small.
\end{itemize}
\end{lemma}

\begin{remark}
In fact, the converse to Lemma \ref{sabus} is true as well, but we will not need this fact.
\end{remark}

\begin{proof}
Choose a point $\eta \in G(\calX)$ which induces an equivalence $e_{\calX} \rightarrow G$.
The map $\eta$ induces a global section $\overline{\eta}$ of $G_{\calX}$ in the $\infty$-topos
$\Shv_{ \widehat{\SSet}}(\calX)$. Let $F_0$ denote the fiber product $F_{\calX} \times_{ G_{\calX} } 1_{ \Shv_{ \widehat{\SSet}}(\calX)}.$
Assumption $(2)$ implies that the functor $F_0: \calX^{op} \rightarrow \widehat{\SSet}$
takes values which are essentially small. It follows from Proposition \ref{representable} that
$F_0$ is representable by an object $U \in \calX$. In particular, we have a tautological point
$\eta' \in F_0(U)$, which determines a commutative diagram
$$ \xymatrix{ e_{ \calX_{/U}} \ar[r] \ar[d] & F \ar[d]^{\alpha} \\
e_{\calX} \ar[r] & G }$$
in $\widehat{\Shv}(\LGeom^{op})$. To complete the proof, it will suffice to show that the upper horizontal map is an equivalence. 

Fix an $\infty$-topos $\calY$, and let $R: \widehat{\Shv}( \LGeom^{op} )
\rightarrow \Shv_{ \widehat{\SSet}}(\calY)$ be the restriction map; we will show that the induced map
$R( e_{\calX_{/U}}) \rightarrow R(F)$ is an equivalence in $\Shv_{ \widehat{\SSet}}(\calY)$. It
will suffice to show that for every $V \in \calY$, the map
$R(e_{\calX_{/U}}(V) \rightarrow R(F)(V)$ induces a homotopy equivalence after passing
to the homotopy fibers over any point $\eta' \in R(G)(V)$. Replacing $\calY$ by
$\calY_{/V}$, we may assume that $\eta'$ is induced by a geometric morphism
$f^{\ast}: \calX \rightarrow \calY$, which determines a map
$1 \rightarrow R(G)$, where $1$ denotes the final object of $\Shv_{\widehat{\SSet}}(\calY)$.
Let $F' = R( e_{\calX_{/U}}) \times_{ R(G)} 1$ and $F'' = R(F) \times_{ R(G) } 1$; to complete the proof it will suffice to show that the induced map $F' \rightarrow F''$ is an equivalence.

We have a commutative diagram
$$ \xymatrix{ F'' \ar[r] \ar[d] & \widehat{f}^{\ast} F_{\calX} \ar[r] \ar[d] & R(F) \ar[d] \\
1 \ar[r]^{ \widehat{f}^{\ast} \overline{\eta} }& \widehat{f}^{\ast} G_{\calX} \ar[r] & R(G) }$$
in the $\infty$-category $\Shv_{ \widehat{\SSet}}(\calY)$. Here the right square is a pullback
since $\alpha$ is universal, and the outer square is a pullback by construction.
It follows that the left square is also a pullback, so that 
$$F'' \simeq \widehat{f}^{\ast} (1_{\calX} \times_{ G_{\calX} } F_{\calX})
\simeq \widehat{f}^{\ast} F_0.$$
We note that $\widehat{f}^{\ast} F_0$ can be identified with the functor represented by
the object $f^{\ast} U \in \calY$, which (by virtue of Remark \ref{goodilk}) is equivalent
to $F'$ as desired. 
\end{proof}

\begin{lemma}\label{guesstimate}
Let $\kappa$ be an uncountable regular cardinal, and let $X_{\bigdot}$ be a simplicial
object of $\SSet$ with the following properties:
\begin{itemize}
\item[$(a)$] For each $n \geq 0$, the connected components of $X_{\bigdot}$ are essentially $\kappa$-small.
\item[$(b)$] For every morphism $[m] \rightarrow [n]$ in $\cDelta$, the induced map
$X_{n} \rightarrow X_{m}$ has essentially $\kappa$-small homotopy fibers.
\end{itemize}
Let $X$ be the geometric realization of $X_{\bigdot}$. Then the induced map
$X_0 \rightarrow X$ has essentially $\kappa$-small homotopy fibers.
\end{lemma}

\begin{proof}
Replacing $X$ by one of its connected components $X'$ (and each $X_{n}$ by the
inverse image $X' \times_{X} X_n$), we may suppose that $X$ is connected.

Let $R \subseteq \pi_0 X_0 \times \pi_0 X_0$ denote the image of $\pi_0 X_1$, and let
$\sim$ denote the equivalence relation on $\pi_0 X_0$ generated by $R$. It follows from
assumption $(b)$ that for every $\kappa$-small subset $A \subseteq \pi_0 X_0$, the
intersections $R \cap (A \times \pi_0 X_0)$ and $R \cap (\pi_0 X_0 \times A)$ are
again $\kappa$-small. Since $\kappa$ is uncountable, it follows that the every $\sim$-equivalence class is $\kappa$-small. Since $(\pi_0 X_0)/ \sim$ is isomorphic to $\pi_0 X \simeq \ast$, 
we conclude that $\pi_0 X$ is itself $\kappa$-small. Combining this with $(a)$, we conclude that $X_0$ is essentially $\kappa$-small. Invoking $(b)$, we deduce that
each $X_{n}$ is essentially $\kappa$-small, so that $X$ is essentially $\kappa$-small. The desired conclusion now follows from the long exact sequences associated to the fibration sequences
$X_0 \times_{X} \{\ast\} \rightarrow X_0 \rightarrow X.$
\end{proof}

\begin{lemma}\label{carbcount}
Let $\calX$ be an $\infty$-topos. Then:
\begin{itemize}
\item[$(1)$] The inclusion $\Shv_{ \widehat{\SSet}}(\calX) \subseteq \Fun( \calX^{op}, \widehat{\SSet})$ admits a left exact left adjoint $L$.
\item[$(2)$] Let $F \in \Fun( \calX^{op}, \widehat{\SSet})$ be a functor such that each of the spaces $F(X)$ is essentially small. Then each of the spaces $LF(X)$ is essentially small. 
\item[$(3)$] Let $\alpha: F \rightarrow G$ be a morphism in $\Fun( \calX^{op}, \widehat{\SSet})$ such that, for each $X \in \calX$, the homotopy fibers of the induced map $F(X) \rightarrow G(X)$ are essentially small. Then for each $X \in \calX$, the homotopy fibers of the map $LF(X) \rightarrow LG(X)$ are also essentially small.
\end{itemize}
\end{lemma}

\begin{proof}
The existence of the left adjoint $L$ follows from Lemma \ref{stur1}. Since
$\Shv_{ \widehat{\SSet}}(\calX)$ contains the essential image of the Yoneda embedding
$j: \calX \rightarrow \Fun( \calX^{op}, \widehat{\SSet})$, we can identify $L \circ j$ with $j$.
Since $j$ is left exact, Proposition \ref{natash} implies that $L$ is also left exact. This proves $(1)$.

We now prove $(2)$. Choose a (small) regular cardinal $\kappa$ such that $\calX$ is $\kappa$-accessible, and let $\calX^{\kappa}$ denote the full subcategory of $\calX$ spanned by the $\kappa$-compact objects. Let $T$ denote the composition
$$ \Fun( \calX^{op}, \widehat{\SSet}) \stackrel{T'}{\rightarrow}
\Fun( (\calX^{\kappa})^{op}, \widehat{\SSet}) \stackrel{T''}{\rightarrow} \Fun( \calX^{op}, \widehat{\SSet}),$$ where $T'$ is the restriction functor and $T''$ is given by the right Kan extension. We have an evident natural transformation $\id \rightarrow T$, which exhibits $T$ as a localization functor
on $\Fun( \calX^{op}, \widehat{\SSet})$. Proposition \ref{sumoto} implies that
every $\widehat{\SSet}$-valued sheaf on $\calX$ is $T$-local. It follows that the canonical map
$L \rightarrow LT$ is an equivalence of functors. In particular, to prove that $LF(X)$ is locally small, 
we may assume without loss of generality that $F$ is $T$-local.

Let $\Fun'( \calX^{op}, \widehat{\SSet})$ denote the full subcategory of $\Fun(\calX^{op}, \widehat{\SSet})$ spanned by the $T$-local functors (in other words, those functors which are right Kan extensions of their restriction to $(\calX^{\kappa})^{op}$; by Proposition \ref{lklk} this
$\infty$-category is equivalent to $\Fun( (\calX^{\kappa})^{op}, \widehat{\SSet} )$). We can identify
$\Fun'( \calX^{op}, \widehat{\SSet} )$ with the $\infty$-category of $\widehat{\SSet}$-valued sheaves
$\Shv_{ \widehat{\SSet}}( \calP( \calX^{\kappa}) )$ on the $\infty$-topos $\calP( \calX^{\kappa})$. Let $F'$
be the image of $F$ under this identification; we observe that the functor
$F': \calP( \calX^{\kappa})^{op} \rightarrow \widehat{\SSet}$ takes essentially small values.
In Remark \ref{quest}, we saw that this $\infty$-category contains $\Shv_{\widehat{\SSet}}( \calX)$ as a left exact localization, and that the localization functor $L': \Shv_{\widehat{\SSet}}(\calP(\calX^{\kappa})) \rightarrow \Shv_{\widehat{\SSet}}( \calX)$. Since $F'$ belongs to the essential image of
the inclusion $\Shv_{ \SSet}( \calP(\calX^{\kappa})) \subseteq \Shv_{ \widehat{\SSet}}( \calP( \calX^{\kappa}))$, the argument given there proves that $L'F'$ belongs to the essential image of
the inclusion $\Shv_{\SSet}(\calX) \subseteq \Shv_{ \widehat{\SSet}}(\calX)$, so that
$LF(X)$ is essentially small as desired.

To prove $(3)$, let us fix a point $\eta \in LG(X)$. We wish to prove the following stronger version of
$(3)$:
\begin{itemize}
\item[$(3')$] For every map $U \rightarrow X$ in $\calX$, the homotopy fiber of the induced map
$LF \rightarrow LG$ is essentially small (here the homotopy fiber is taken over the point determined by $\eta$).
\end{itemize}

Let $\calX^{0}_{/X}$ denote the full subcategory of $\calX_{/X}$ spanned by those morphisms
$U \rightarrow X$ for which condition $(2')$ is satisfied. Since $LF$ and $LG$ belong to
$\Shv_{\widehat{\SSet}}(\calX)$ (and since the collection of essentially small spaces is stable under small limits), we conclude that $\calX^{0}_{/X}$ is stable under small colimits in $\calX_{/X}$. 

Let $\calX^{1}_{/X}$ be the largest sieve contained in $\calX^{0}_{/X}$ (in other words, a morphism $U \rightarrow X$ belongs to $\calX^{1}_{/X}$ if and only if, for every morphism $V \rightarrow U$ in $\calX$, the composite map $V \rightarrow X$ beongs to $\calX^{0}_{/X}$). Since colimits in $\calX$ are universal, we conclude that $\calX^{1}_{/X}$ is stable under small colimits in $\calX_{/X}$. It
follows that $\calX^{1}_{/X} \simeq \calX_{/X_0}$ for some monomorphism
$i: X_0 \rightarrow X$ in $\calX$. We wish to show that $i$ is an equivalence.

Since $L$ is left exact, we have $L( G \times_{LG} j(X) ) \simeq L j(X) \simeq j(X)$.
In particular, the map $G \times_{LG} j(X) \rightarrow j(X_0)$ cannot factor through $j(X_0)$ 
unless $i$ is an equivalence. It will therefore suffice to show that $G \times_{ j(X) } j(X_0) \simeq G$.
In other words, it will suffice to show that if $U \in \calX_{/X}$ and $\eta' \in G(U)$ is a point such that
the images of $\eta$ and $\eta'$ lie in the same connected component of $LG(U)$, then
$U \in \calX^{1}_{/X}$. Since the existence of $\eta'$ is stable under the process of
replacing $U$ by some further refinement $V \rightarrow U$, it will suffice to show that
$U \in \calX^{0}_{/X}$. Replacing $X$ by $U$, we obtain the following reformulation of $(3')$:

\begin{itemize}
\item[$(3'')$] Let $\eta' \in G(X)$. Then the homotopy fiber $Z$ of the induced map
$LF(X) \rightarrow LG(X)$ (over the point determined by $\eta$) is essentially small.
\end{itemize}

Since $L$ is left exact, we can identify $Z$ with $LF_0(X)$, where
$F_0 = F \times_{G} j(X)$. Since the homotopy fibers of the maps $F(Y) \rightarrow G(Y)$ are
essentially small, we may assume without loss of generality that $F_0 \in \Fun( \calX^{op}, \SSet)$. 
Invoking $(2)$, we deduce that the values of $LF_0$ are essentially small as desired.
\end{proof}

\begin{proof}[Proof of Lemma \ref{santan}]
Let $\calX_{\bigdot}$ be a simplicial object of $\LGeom^{op}_{\mathet}$, and let
$F_{\bigdot}$ be its image under the Yoneda embedding 
$j: \LGeom^{op} \rightarrow \widehat{\Shv}( \LGeom^{op} )$. Let $F$ be a geometric realization
of $|F_{\bigdot}|$. We will prove the following:
\begin{itemize}
\item[$(\ast)$] The map $\beta: F_0 \rightarrow F$ satisfies conditions $(1)$ and $(2)$ of Lemma \ref{sabus}.
\end{itemize}
Assuming $(\ast)$ for the moment, we will complete the proof of Lemma \ref{santan}. Let
$F'_{\bigdot}$ be a \Cech nerve of the induced map
$F_0 \rightarrow F$ (so that $F'_{n} \simeq F_0 \times_{F} F_0 \times \ldots \times_{F} F_0$; in particular $F'_0 \simeq F_0$). We first claim that each $F'_{n}$ is representable by an $\infty$-topos $\calX'_{n}$, and that each inclusion $[0] \hookrightarrow [n]$ induces an \'{e}tale map of $\infty$-topos $\calX'_{n} \rightarrow \calX'_0 \simeq \calX_0$. Since $F'_{\bigdot}$ is a groupoid object of $\widehat{\Shv}( \LGeom^{op} )$ (and $F'_0 \simeq F_0$ is representable by the $\infty$-topos $\calX_0$), it will suffice to prove this result when $n=1$. Consider the pullback diagram
$$ \xymatrix{ F'_1 \ar[r]^{\beta'} \ar[d] & F'_0 \ar[d] \\
F_0 \ar[r]^{\beta} & F. }$$
It follows from condition $(\ast)$ that $\beta'$ satisfies conditions $(1)$ and $(2)$ of Lemma \ref{sabus},
so that $F'_1$ is representable by an $\infty$-topos \'{e}tale over $\calX_0$ as desired.

Since the Yoneda embedding $j$ is fully faithful, we may assume without loss of generality
that $F'_{\bigdot}$ is the image under $j$ of a groupoid object $\calX'_{\bigdot}$ of
$\LGeom^{op}$. Using Corollary \ref{toadscan}, we deduce that $\calX'_{\bigdot}$ 
defines a simplicial object of the subcategory $\LGeom^{op}_{\mathet}$. The evident
natural transformation $F_{\bigdot} \rightarrow F'_{\bigdot}$ induces a map of
simplicial objects $\alpha: \calX_{\bigdot} \rightarrow \calX'_{\bigdot}$; we claim that $\alpha$ has the desired properties. The only nontrivial point is to verify that the induced map of geometric realizations
$| \calX_{\bigdot} | \rightarrow | \calX'_{\bigdot} |$ is an equivalence of $\infty$-topoi. For this, it suffices to show that for every $\infty$-topos $\calY$, the upper horizontal map in the diagram
$$\xymatrix{ \bHom_{ \LGeom^{op} }( | \calX'_{\bigdot} |, \calY) \ar[r] \ar[d] & \bHom_{ \LGeom^{op} }( | \calX_{\bigdot } |, \calY) \ar[d] \\
\varprojlim \bHom_{ \LGeom^{op} }( \calX'_{n}, \calY) \ar[r] & \varprojlim \bHom_{ \LGeom^{op} }( \calX_{n}, \calY ) }$$
is a homotopy equivalence. Since the vertical maps are homotopy equivalences, it suffices to show that the lower horizontal map is a homotopy equivalence. Since $j$ is fully faithful, it suffices to show that the lower horizontal map in the analogous diagram
$$\xymatrix{ \bHom_{ \widehat{\Shv}(\LGeom^{op}) }( | F'_{\bigdot} |, e_{\calY}) \ar[r] \ar[d] & \bHom_{ \widehat{\Shv}(\LGeom^{op})}( | F_{\bigdot } |, e_{\calY}) \ar[d] \\
\varprojlim \bHom_{ \widehat{\Shv}(\LGeom^{op})}( F'_{n}, e_{\calY}) \ar[r] & \varprojlim \bHom_{ \widehat{\Shv}(\LGeom^{op}) }( F_{n}, e_{\calY} ) }$$
is a homotopy equivalence. Again, the vertical maps are homotopy equivalences, so we are reduced to showing that the upper horizontal map is a homotopy equivalence. This follows from the fact that we have an equivalence $|F_{\bigdot}| \simeq F \simeq |F'_{\bigdot}|$ in $\widehat{\Shv}( \LGeom^{op} )$, since groupoid objects in $\widehat{\Shv}(\LGeom^{op})$ are effective (Lemma \ref{kumba}).

It remains to prove $(\ast)$. Remark \ref{sapper} implies that $\beta: F_0 \rightarrow F$ is universal.
To complete the proof, we must show that for every $\infty$-topos $\calY$, the homotopy fibers
of the induced map $F_0( \calY) \rightarrow F(\calY)$ are essentially small. Let
$R: \widehat{\Shv}( \LGeom^{op}) \rightarrow \Shv_{ \widehat{\SSet}}( \calY)$ denote the restriction map,
let $G_{\bigdot}$ denote the image of $F_{\bigdot}$ under $R$, and let $G = | G_{\bigdot} | \simeq R(F)$. Let $G'$ denote the geometric realization of $G_{\bigdot}$ in the larger $\infty$-category
$\Fun( \calY^{op}, \widehat{\SSet})$, and let $L: \Fun( \calY^{op}, \widehat{\SSet} ) \rightarrow \Shv_{ \widehat{\SSet} }(\calY)$ be a left adjoint to the inclusion (see Lemma \ref{carbcount}). Then we can identify the map $G_0 \rightarrow G$ with the image under $L$ of the map $u: G_0 \rightarrow G'$. In view of Lemma \ref{carbcount}, it will suffice to show that for every object $U \in \calY$, the induced map $G_0(U) \rightarrow G'(U)$ has essentially small homotopy fibers.

For each $n \geq 0$ and each object $U \in \calY$, we can identify $G_{n}(U)$ with the
maximal Kan complex contained in $\Fun^{\ast}( \calX_{n}, \calY_{/U})$. Since the $\infty$-category $\Fun^{\ast}( \calX_{n}, \calY_{/U})$ is locally small (Proposition \ref{nottoobig}), we conclude that each connected component of $G_{n}(U)$ is essentially small.
Moreover, for every morphism $[m] \rightarrow [n]$ in $\cDelta$, the induced map
$\beta: G_{n}(U) \rightarrow G_{m}(U)$ is induced by composition with an \'{e}tale geometric morphism
$g^{\ast}: \calX_{m} \rightarrow \calX_{n}$, so that the homotopy fibers of $\beta$ are essentially
small by Remark \ref{goodilk}. The desired result now follows from Lemma \ref{guesstimate}.
\end{proof}

\subsection{Structure Theory for $\infty$-Topoi}\label{structuretheor}

In this section we will analyze the following question: given a geometric morphism
$f: \calX \rightarrow \calY$ of $\infty$-topoi, when is $f$ an equivalence? Clearly, this is true if and only if the pullback functor $f^{\ast}$ is both fully faithful and essentially surjective. It is useful to isolate and study these conditions individually.

\begin{definition}\label{nutro}\index{gen}{image of a geometric morphism}
Let $f: \calX \rightarrow \calY$ be a geometric morphism of $\infty$-topoi. The
{\it image} of $f$ is defined to be the smallest full subcategory of $\calX$ which contains
$f^{\ast} \calY$ and is stable under small colimits and finite limits. 
We will say that $f$ is {\it algebraic} if the image of $f$ coincides with $\calX$.\index{gen}{algebraic morphism}
\end{definition}

Our first goal is to prove that the image of a geometric morphism is itself an $\infty$-topos.

\begin{proposition}\label{proet}
Let $f: \calX \rightarrow \calZ$ be a geometric morphism of $\infty$-topoi, and let $\calY$
be the image of $f$. Then $\calY$ is an $\infty$-topos. Moreover, the inclusion
$\calY \subseteq \calX$ is left exact and colimit-preserving, so we have obtain a factorization of $f$ as a composition of geometric morphisms
$$ \calX \stackrel{g}{\rightarrow} \calY \stackrel{h}{\rightarrow} \calZ$$
where $h$ is algebraic and $g^{\ast}$ is fully faithful. 
\end{proposition}

\begin{proof}
We will show that $\calY$ satisfies the $\infty$-categorical versions of Giraud's axioms (see
Theorem \ref{mainchar}). Axioms $(ii)$, $(iii)$, and $(iv)$ are concerned with the interaction between colimits and finite limits. Since $\calX$ satisfies these axioms, and $\calY \subseteq \calX$ is stable under the relevant constructions, $\calY$ automatically satisfies these axioms as well. The only nontrivial point is to verify $(i)$, which asserts that $\calY$ is presentable. 

Choose a small collection of objects $\{ Z_{\alpha} \}$ which generate $\calZ$ under colimits.
Now choose an uncountable regular cardinal $\tau$ with the following properties:

\begin{itemize}
\item[$(1)$] Each $f^{\ast}(Z_{\alpha})$ is a $\tau$-compact object of $\calX$.
\item[$(2)$] The final object $1_{\calX}$ is $\tau$-compact.
\item[$(3)$] The limits functor $\Fun(\Lambda^2_2, \calX) \rightarrow \calX$
(a right adjoint to the diagonal functor) is $\tau$-continuous and preserves $\tau$-compact objects.
\end{itemize}

Let $\calY'$ be the collection of all objects of $\calY$ which are $\tau$-compact when considered as objects of $\calX$. Clearly, each object of $\calY'$ is also $\tau$-compact when regarded as an object of $\calY$. Moreover, because $\calX$ is accessible, $\calY'$ is essentially small. It will therefore suffice to prove that $\calY'$ generates $\calY$ under colimits.

Choose a minimal model $\calY'_0$ for $\calY$. Since $\calX$ is accessible, the full subcategory $\calX^{\kappa}$ spanned by the $\kappa$-compact objects is essentially small, so that $\calY'_0$ is small. According to Proposition \ref{intprop}, there exists a $\tau$-continuous functor
$F: \Ind_{\tau}( \calY'_{0} ) \rightarrow \calX$ whose composition with the Yoneda embedding
is equivalent to the inclusion $\calY'_0 \subseteq \calX$. Since $\calY'_0$ admits $\tau$-small colimits, $\Ind_{\tau}( \calY'_{0} )$ is presentable.
Proposition \ref{uterr} implies that $F$ is fully faithful; let $\calY''$ be its essential image. To complete the proof, it will suffice to show that $\calY'' = \calY$. 

Since $\calY$ is stable under colimits in $\calX$, we have $\calY'' \subseteq \calY$. 
According to Proposition \ref{sumatch}, $F$ preserves small colimits, so that $\calY''$ is stable under small colimits in $\calX$. By construction, $\calY''$ contains each $f^{\ast}(Z_{\alpha})$. Since $f^{\ast}$ preserves colimits, we conclude that $\calY''$ contains $f^{\ast} \calZ$. By definition $\calY$ is the smallest full subcategory of $\calX$ which contains $f^{\ast} \calZ$ and is stable under small colimits and finite limits. It remains only to show that $\calY''$ is stable under finite limits. Assumption $(2)$ guarantees that $\calY''$ contains the final object of $\calX$, so we
need only show that $\calY''$ is stable under pullbacks. Consider a diagram
$p: \Lambda^2_2 \rightarrow \calY''$. The proof of Proposition \ref{horse1} (applied
with $K = \Lambda^2_2$ and $\kappa = \omega$) shows that $p$ can be written
as a $\tau$-filtered colimit of diagrams $p_{\alpha}: \Lambda^2_2 \rightarrow \calY''$. 
Since filtered colimits in $\calX$ are left exact (Example \ref{tucka}), we conclude that the limit of
$p$ can be obtained as a $\tau$-filtered colimit of limits of the diagrams $p_{\beta}$. In view of assumption $(3)$, each of these limits lies in $\calY'$, so that the limit of $p$ lies in $\calY''$ as desired.
\end{proof}

\begin{remark}
The factorization of Proposition \ref{proet} is unique up to (canonical) equivalence.
\end{remark}

The terminology of Definition \ref{nutro} is partially justified by the following observations: 

\begin{proposition}\label{charproet}
\begin{itemize}
\item[$(1)$] Every \'{e}tale geometric morphism between $\infty$-topoi is algebraic.
\item[$(2)$] The collection of algebraic geometric morphisms of $\infty$-topoi is stable
under filtered limits $($in $\RGeom${}$)$.
\end{itemize}
\end{proposition}

\begin{proof}
We first prove $(1)$. Let $\calX$ be an $\infty$-topos, let $U$ be an object of $\calX$,
let $\pi_{!}: \calX^{/U} \rightarrow \calX$ be the projection functor, and let $\pi^{\ast}$ be a 
left adjoint to $\pi_{!}$. Let $f: X \rightarrow U$ be an object of $\calX^{/U}$, and let
$F: f \rightarrow \id_{U}$ be a morphism in $\calX^{/U}$ (uniquely determined up to equivalence; for example, we can take $F$ to be the composition of $f$ with a retraction $\Delta^1 \times \Delta^1 \rightarrow \Delta^1$). Let $g: F \rightarrow \pi^{\ast} \pi_{!} F$ be the unit map for the adjunction
between $\pi^{\ast}$ and $\pi_{!}$. We claim that $g$ is a pullback square in $\calX_{/U}$. According to Proposition \ref{goeselse}, it will suffice to verify that the image of $g$ under $\pi_{!}$ is a pullback square in $\calX$. But this square can be identified with
$$ \xymatrix{ X \ar[r] \ar[d] & X \times U \ar[d] \\
U \ar[r]^-{\delta} & U \times U, }$$
which is easily shown to be Cartesian. It follows that, in $\calX_{/U}$, $f$ can be obtained
as a fiber product of the final object with objects that lie in the essential image of $\pi^{\ast}$.
It follows that $\pi^{\ast} \calX$ generates $\calX_{/U}$ under finite limits, so that
$\pi$ is algebraic.

To prove $(2)$, we consider a geometric morphism $f: \calX \rightarrow \calY$ which
is a filtered limit of algebraic geometric morphisms $\{ f_{\alpha}: \calX_{\alpha} \rightarrow \calY_{\alpha} \}$ in the $\infty$-category $\Fun(\Delta^1, \RGeom)$. Let $\calX' \subseteq \calX$
be a full subcategory which is stable under finite limits, small colimits, and contains
$f^{\ast} \calY$. We wish to prove that $\calX' = \calX$. For each $\alpha$, we have
a diagram of $\infty$-topoi
$$ \xymatrix{ \calX \ar[r]^{f} \ar[d]^{\psi(\alpha)} & \calY \ar[d] \\
\calX_{\alpha} \ar[r]^{f_{\alpha}} & \calY_{\alpha}. }$$
Let $\calX'_{\alpha}$ be the preimage of $\calX'$ under $\psi(\alpha)^{\ast}$. Then
$\calX'_{\alpha} \subseteq \calX_{\alpha}$ is stable under finite limits, small colimits, and
contains the essential image of $f_{\alpha}^{\ast}$. Since $f_{\alpha}$ is algebraic, we conclude
that $\calX'_{\alpha} = \calX_{\alpha}$. In other words, $\calX'$ contains the essential image
of each $\psi(\alpha)^{\ast}$. Lemma \ref{steakknife} implies that every object of $\calX$ can be realized as a filtered colimit of objects, each of which belongs to the essential image of
$f_{\alpha}^{\ast}$ for $\alpha$ appropriately chosen. Since $\calX'$ is stable under small colimits, we conclude that $\calX' = \calX$. It follows that $f$ is algebraic, as desired.
\end{proof}

\begin{remark}\label{weakcon}
It is possible to formulate a converse to Proposition \ref{charproet}. Namely, one can characterize the class of algebraic morphisms as the smallest class of geometric morphisms which contains all \'{e}tale morphisms and is stable under certain kinds of filtered limits. However, it is necessarily to allow limits which are parametrized not just by filtered $\infty$-categories, but filtered {\em stacks} over $\infty$-topoi. The precise statement requires ideas which lie outside the scope of this book.
\end{remark}

Having achieved a rudimentary understanding of the class of algebraic geometric morphisms, we now turn our attention to the opposite extreme: namely, geometric morphisms $f: \calX \rightarrow \calY$ where $f^{\ast}$ is fully faithful. 

\begin{proposition}\label{unterware}
Let $f: \calX \rightarrow \calY$ be a geometric morphism of $\infty$-topoi. Suppose that
$f^{\ast}$ is fully faithful and essentially surjective on $1$-truncated objects. Then $f^{\ast}$ is essentially surjective on $n$-truncated objects for all $n$.
\end{proposition}

The proof uses ideas which will be introduced in \S \ref{homotopysheaves} and \S \ref{chmdim}.

\begin{proof}
Without loss of generality, we may identify $\calY$ with the essential image of $f^{\ast}$.
We use induction on $n$. The result is obvious for $n = 1$. Assume that $n > 1$, and let $X$ be an $n$-truncated object of $\calX$. By the inductive hypothesis, $U= \tau_{\leq n-1} X$ belongs to $\calY$.
Replacing $\calX$ and $\calY$ by $\calX_{/ U}$ and $\calY_{ / U}$, we may suppose that $X$ is $n$-connective.

We observe that $\pi_{n} X$ is an abelian group object of the ordinary topos
$\Disc(\calX_{/X})$. Since $X$ is $2$-connective, Proposition \ref{nicelemma}
implies that the pullback functor $\Disc(\calX) \rightarrow \Disc( \calX_{/X})$ is an equivalence of categories. We may therefore identify $\pi_n X$ with an abelian group object
$A \in \Disc(\calX)$. Since $A$ is discrete, it belongs to $\calY$. It follows that
the Eilenberg-MacLane object $K(A,n+1)$ belongs to $\calY$. Since $X$ is an $n$-gerb banded by $A$, Theorem \ref{starthm} implies the existence of a pullback diagram
$$ \xymatrix{ X \ar[r] \ar[d] & 1_{\calX} \ar[d] \\
1_{\calX} \ar[r] & K(A,n+1). }$$
Since $\calY$ is stable under pullbacks in $\calX$, we conclude that $X \in \calY$ as desired.
\end{proof}

\begin{corollary}\label{unwhere}
Let $f: \calX \rightarrow \calY$ be a geometric morphism of $\infty$-topoi. Suppose that
$f^{\ast}$ is fully faithful and essentially surjective on $1$-truncated objects, and that
$\calX$ is $n$-localic $($see \S \ref{nlocalic}$)$. Then $f$ is an equivalence of $\infty$-topoi.
\end{corollary}

\begin{remark}
In the situation of Corollary \ref{unwhere}, one can eliminate the hypothesis that $\calX$ is $n$-localic in the presence of suitable finite-dimensionality assumptions on $\calX$ and $\calY$; see \S \ref{homdim}.
\end{remark}

\begin{remark}
Let $\calX$ be an $n$-localic $\infty$-topos, and let $\calY$ be the $2$-localic $\infty$-topos
associated to the $2$-topos $\tau_{\leq 1} \calX$, so that we have a geometric morphism
$f: \calX \rightarrow \calY$. It follows from Corollary \ref{unwhere} that $f$ is algebraic.
Roughly speaking, this tells us that there is only a very superficial interaction between the theory of $k$-categories and ``topology'', for $k > 2$. On the other hand, this statement fails dramatically if
$k=1$: the relationship between an ordinary topos and its underlying locale is typically very complicated, and not algebraic in any reasonable sense. It is natural to ask what happens when $k=2$. In other words, does Proposition \ref{unterware} remain valid if $f^{\ast}$ is only assumed to be essentially surjective on discrete objects? An affirmative answer would indicate that our theory of $\infty$-topoi is a relatively modest extension of classical topos theory. A counterexample could be equally interesting, if it were to illustrate a nontrivial interaction between higher category theory and geometry.
\end{remark}
\section{$n$-Topoi}\label{chap6sec3}

\setcounter{theorem}{0}

Roughly speaking, an ordinary topos is a category which resembles the category of sheaves of {\em sets} on a topological space $X$. In \S \ref{chap6sec1}, we introduced the definition of an $\infty$-topos. In the same rough terms, we can think of an $\infty$-topos as an $\infty$-category which resembles the $\infty$-category of sheaves of $\infty$-groupoids on a topological space $X$. Phrased in this way, it is natural to guess that these two notions have a common generalization. In \S \ref{c631}, we will introduce the notion of an $n$-topos, for every $0 \leq n \leq \infty$. The idea is that an $n$-topos should be an $n$-category which resembles the $n$-category of sheaves of $(n-1)$-groupoids on a topological space $\calX$. Of course, there are many approaches to making this idea precise. Our main result, Theorem \ref{nchar}, asserts that several candidate definitions are equivalent to one another. The proof of Theorem \ref{nchar} will occupy our attention for most of this section. In \S \ref{provengiraud}, we study an axiomatization of the class of $n$-topoi in the spirit of Giraud's theorem, and in \S \ref{provengiraudeasy} we will give a characterization of $n$-topoi based on their descent properties. The case of $n=0$ is somewhat exceptional, and merits special treatement. In \S \ref{0topoi} we will show that a $0$-topos is essentially the same thing as a {\em locale} (a mild generalization of the notion of a topological space). 

Our main motivation for introducing the definition of an $n$-topos is that it allows us to study $\infty$-topoi and topological spaces (or, more generally, $0$-topoi) in the same setting. In \S \ref{nlocalic}, we will introduce constructions which allow us to pass back and forth between $m$-topoi and $n$-topoi, for any $0 \leq m \leq n \leq \infty$. We introduce an $\infty$-category
$\Geo_{n}^{\GeoR}$ of $n$-topoi, for each $n \leq 0$, and show that each
$\Geo_{n}^{\GeoR}$ can be regarded as a {\em localization} of the $\infty$-category
$\RGeom$. In other words, the study of $n$-topoi for $n < \infty$ can be regarded as a special case of the theory of $\infty$-topoi.

\subsection{Characterizations of $n$-Topoi}\label{c631}

In this section, we will introduce the definition of $n$-topos for $0 \leq n < \infty$. In view of Theorem \ref{mainchar}, there are several reasonable approaches to the subject. We will begin with extrinsic approach.

\begin{definition}\label{ntopdefff}\index{gen}{$n$-topos}
Let $0 \leq n < \infty$. An $\infty$-category $\calX$ is an {\it $n$-topos} if there
exists a small $\infty$-category $\calC$ and an (accessible) left exact localization
$$ L: \calP_{\leq n-1}(\calC) \rightarrow \calX,$$
where $\calP_{\leq n-1}(\calC)$ denotes the full subcategory of $\calP(\calC)$ spanned by
the $(n-1)$-truncated objects.\index{not}{Pcalleq@$\calP_{\leq n}(\calC)$}
\end{definition}

\begin{remark}
The accessibility condition on the localization functor $L: \calP_{\leq n-1}(\calC) \rightarrow \calX$ of Definition \ref{ntopdefff} is superfluous: we will show that such left exact localization of $\calP_{\leq n-1}(\calC)$ is automatically accessible (combine Proposition \ref{alltoploc} with Corollary \ref{topaccess}).
\end{remark}

\begin{remark}
An $\infty$-category $\calX$ is a $1$-topos if and only if it is equivalent to the nerve of an ordinary (Grothendieck) topos; this follows immediately from characterization $(B)$ of Proposition \ref{toposdefined}.
\end{remark}

\begin{remark}
Definition \ref{ntopdefff} makes sense also in the case $n=-1$, but is not very interesting. Up to equivalence, there is precisely one $(-1)$-topos: the final $\infty$-category $\ast$.
\end{remark}

Our main goal is to prove the following result:

\begin{theorem}\label{nchar}\index{gen}{Giraud's theorem!for $n$-topoi}
Let $\calX$ be a presentable $\infty$-category and let $0 \leq n < \infty$. The following conditions are equivalent:
\begin{itemize}
\item[$(1)$] There exists a small $n$-category $\calC$ which admits finite limits, a Grothendieck topology on $\calC$, and an equivalence of $\calX$ with the full subcategory of $\Shv_{\leq n-1}(\calC) \subseteq \Shv(\calC)$ consisting of $(n-1)$-truncated objects of $\Shv(\calC)$.

\item[$(2)$] There exists an $\infty$-topos $\calY$ and an equivalence
$\calX \rightarrow \tau_{\leq n-1} \calY$.

\item[$(3)$] The $\infty$-category $\calX$ is an $n$-topos.

\item[$(4)$] Colimits in $\calX$ are universal, $\calX$ is equivalent to an $n$-category, and the class of $(n-2)$-truncated morphisms
in $\calX$ is local $($ see \S \ref{magnet} $)$.

\item[$(5)$] Colimits in $\calX$ are universal, $\calX$ is equivalent to an $n$-category, and for all sufficiently large regular cardinals $\kappa$, there exists an object of $\calX$ which classifies
$(n-2)$-truncated, relatively $\kappa$-compact morphisms in $\calX$.

\item[$(6)$] The $\infty$-category $\calX$ satisfies the following $n$-categorical versions
of Giraud's axioms:

\begin{itemize}\index{gen}{Giraud's axioms!for $n$-topoi}
\item[$(i)$] The $\infty$-category $\calX$ is equivalent to a presentable $n$-category.
\item[$(ii)$] Colimits in $\calX$ are universal.
\item[$(iii)$] If $n > 0$, then coproducts in $\calX$ are disjoint.
\item[$(iv)$] Every $n$-efficient $($see \S \ref{provengiraud}$)$ groupoid object of $\calX$ is effective.
\end{itemize}
\end{itemize}
\end{theorem}

\begin{proof}
The case $n=0$ will be analyzed very explicitly in \S \ref{0topoi}; let us therefore restrict our attention to the case $n > 0$.
The implication $(1) \Rightarrow (2)$ is obvious (take $\calY = \Shv(\calC)$). Suppose
that $(2)$ is satisfied. Without loss of generality, we may suppose that $\calY$ is an (accessible) left exact localization of $\calP(\calC)$ for some small $\infty$-category $\calC$. Then
$\calX$ is a left-exact localization of $\calP_{\leq n-1}(\calC)$, which proves $(3)$.

We next prove the converse $(3) \Rightarrow (2)$. We first observe that
$\calP_{\leq n-1}(\calC) = \Fun( \calC^{op}, \tau_{\leq n-1} \SSet)$. Let $h_n \calC$ be the underlying $n$-category of $\calC$, as in Proposition \ref{undern}. Since $\tau_{\leq n-1} \SSet$
is equivalent to an $n$-category, we conclude that composition with the projection
$\calC \rightarrow \hn{n}{\calC}$ induces an equivalence $\calP_{\leq n-1}(\hn{n}{\calC}) \rightarrow \calP_{\leq n-1}(\calC)$. Consequently, we may assume without loss of generality (replacing $\calC$ by $\hn{n}{\calC}$ if necessary) that there is an accessible left exact localization $L: \calP_{\leq n-1}(\calC) \rightarrow \calX$, where $\calC$ is an $n$-category. Let $S$ be the collection of all morphisms $u$ in $\calP_{\leq n-1}(\calC)$ such that $Lu$ is an equivalence, so that $S$ is of small generation. Let
$\overline{S}$ be the strongly saturated class of morphisms in $\calP(\calC)$ generated by $S$.
We observe that $\tau_{\leq n-1}^{-1}(S)$ is a strongly saturated class of morphisms containing
$S$, so that $\overline{S} \subseteq \tau_{\leq n-1}^{-1}(S)$. It follows that $S^{-1} \calP_{\leq n-1}(\calC)$
is contained in $\calY = \overline{S}^{-1} \calP(\calC)$, and may therefore be identified with the collection of $(n-1)$-truncated objects of $\calY$. To complete the proof, it will suffice to show that $\calY$ is an $\infty$-topos. For this, it will suffice to show that $\overline{S}$ is stable under pullbacks. Let $T$ be the collection of all morphisms $f: X \rightarrow Y$ in $\calP(\calC)$ such that
for every pullback diagram
$$ \xymatrix{ X' \ar[r] \ar[d]^{f'} & X \ar[d]^{f} \\
Y' \ar[r] & Y }$$
the morphism $f'$ belongs to $\overline{S}$. It is easy to see that $T$ is strongly saturated; we wish to show that $T \subseteq \overline{S}$. It will therefore suffice to prove that $S \subseteq T$. Let us therefore fix $f:X \rightarrow Y$ belonging to $S$, and let $\calD$ be the full subcategory of $\calP(\calC)$ spanned by those objects $Y'$ such that for {\em any} pullback diagram
$$ \xymatrix{ X' \ar[r] \ar[d]^{f'} & X \ar[d]^{f} \\
Y' \ar[r] & Y, }$$
$f'$ belongs to $\overline{S}$. Since colimits in $\calP(\calC)$ are universal and $\overline{S}$ is stable under colimits, we conclude that $\calD$ is stable under colimits in $\calP(\calC)$. 
Since $\calP(\calC)$ is generated under colimits by the essential image of the Yoneda embedding
$j: \calC \rightarrow \calP(\calC)$, it will suffice to show that $j(C) \in \calD$ for each $C \in \calC$.
We now observe that $\calP_{\leq n-1}(\calC) \subseteq \calD$ (since $S$ is stable under pullbacks in $\calP_{\leq n-1}(\calC)$), and that $j(C) \in \calP_{\leq n-1}(\calC)$ in virtue of our assumption that $\calC$ is an $n$-category.

The implication $(2) \Rightarrow (4)$ will be established in \S \ref{provengiraudeasy} (Propositions \ref{tigre} and \ref{tigress}).
The proof of Theorem \ref{colimsurt} adapts without change to show that $(4) \Leftrightarrow (5)$.
The implication $(4) \Rightarrow (6)$ will be proven in \S \ref{provengiraudeasy} (Proposition \ref{ncharles}). Finally, the ``difficult'' implication $(6) \Rightarrow (1)$ will be proven in \S \ref{provengiraud} (Proposition \ref{diamondstep}), using an inductive argument quite similar to the proof of Giraud's original result.
\end{proof}

\begin{remark}
Theorem \ref{nchar} is slightly stronger than its $\infty$-categorical analogue, Theorem \ref{mainchar}: it asserts that every $n$-topos arises as an $n$-category of sheaves on some
$n$-category $\calC$ equipped with a Grothendieck topology.
\end{remark}

\begin{remark}
Let $\calX$ be a presentable $\infty$-category in which colimits are universal. Then
there exists a regular cardinal $\kappa$ such that every monomorphism is relatively $\kappa$-compact. In this case, characterization $(5)$ of Theorem \ref{nchar} recovers a classical
description of ordinary topos theory: a category $\calX$ is a topos if and only if it is presentable, colimits in $\calX$ are universal, and $\calX$ has a subobject classifier.
\end{remark}

\subsection{$0$-Topoi and Locales}\label{0topoi}

Our goal in this section is to prove Theorem \ref{nchar} in the special case $n=0$. A byproduct of our proof is a classification result (Corollary \ref{charlocale}), which identifies the theory of $0$-topoi with the classical theory of {\em locales} (Definition \ref{deflocale}).

We begin by observing that when $n=0$, a morphism in an $\infty$-category $\calX$ is
$(n-2)$-truncated if and only if it is an equivalence. Consequently, any final object of $\calX$ is an $(n-2)$-truncated morphism classifier, and the class of $(n-2)$-truncated morphisms is automatically local (in the sense of Definition \ref{localitie}). 
Moreover, if $\calX$ is a $0$-category then every groupoid object in $\calX$ is equivalent to a constant groupoid, and therefore automatically effective. Consequently, characterizations $(4)$ through $(6)$ in Theorem \ref{nchar} all reduce to the same condition on $\calX$, and we may restate the desired result as follows:

\begin{theorem}\label{sumatc}\index{gen}{Giraud's theorem!for $0$-topoi}
Let $\calX$ be a presentable $0$-category. The following conditions are equivalent:
\begin{itemize}
\item[$(1)$] There exists a small $0$-category $\calC$ which admits finite limits, 
a Grothendieck topology on $\calC$, and an equivalence $\calX \rightarrow \Shv_{\leq -1}(\calC)$.
\item[$(2)$] There exists an $\infty$-topos $\calY$ and an equivalence
$\calX \rightarrow \tau_{\leq -1} \calY$.
\item[$(3)$] The $\infty$-category $\calX$ is a $0$-topos.
\item[$(4)$] Colimits in $\calX$ are universal.
\end{itemize}
\end{theorem}

Before giving a proof of Theorem \ref{sumatc}, it is convenient to reformulate
condition $(4)$. Recall that any $0$-category $\calX$ is equivalent to
$\Nerve(\calU)$, where $\calU$ is a partially ordered set which is well-defined up to canonical isomorphism (see Example \ref{0catdef}). The presentability of $\calX$ is equivalent to the assertion that $\calU$ is a {\em complete lattice}: that is, every subset of $\calU$ has a least upper bound in $\calU$ (this condition formally implies the existence of greater lower bounds, as well).\index{gen}{complete lattice}

\begin{remark}
If $n=0$, then every presentable $n$-category is essentially small. This is typically not true for $n > 0$.
\end{remark}

We note that the condition that colimits in $\calX$ be universal can also be formulated in terms of the partially ordered set $\calU$: it is equivalent to the assertion that meets in $\calU$ commute with infinite joins in the following sense:

\begin{definition}\label{deflocale}
Let $\calU$ be a partially ordered set. We will say that $\calU$ is a {\em locale} if the following conditions are satisfied:\index{gen}{locale}
\begin{itemize}
\item[$(1)$] Every subset $\{ U_{\alpha} \}$ of elements of $\calU$ has a least upper bound
$\bigcup_{\alpha} U_{\alpha}$ in $\calU$.

\item[$(2)$] The formation of least upper bounds commutes with meets, in the sense that
$$ \bigcup (U_{\alpha} \cap V) = ( \bigcup U_{\alpha} ) \cap V.$$ (Here $(U \cap V)$ denotes the greatest lower bound of $U$ and $V$, which exists in virtue of assumption $(1)$.)
\end{itemize}
\end{definition}

\begin{example}
For every topological space $X$, the collection $\calU(X)$ of open subsets of $X$ forms a locale. 
Conversely, if $\calU$ is a locale, then there is a natural topology on the collection of prime filters of $\calU$ which allows us to extract from $\calU$ a topological space. These two constructions are adjoint to one another, and in good cases they are actually inverse equivalences. More precisely, the adjunction gives rise to an equivalence between the category of {\em spatial} locales and the category of {\em sober} topological spaces\label{sober}. In general, a locale can be regarded as a sort of generalized topological space, in which one may speak of open sets but one does not generally have a sufficient supply of points. We refer the reader to \cite{johnstone} for details.
\end{example}

We can summarize the above discussion as follows:

\begin{proposition}\label{colim0univ}
Let $\calX$ be a presentable $0$-category. Then colimits in $\calX$ are universal if and only if $\calX$ is equivalent to $\Nerve(\calU)$, where $\calU$ is a locale.
\end{proposition}

We are now ready to give the proof of Theorem \ref{sumatc}.

\begin{proof}
The implications $(1) \Rightarrow (2) \Rightarrow (3)$ are easy. Suppose that $(3)$ is satisfied, so that $\calX$ is a left-exact localization of $\calP_{\leq -1}(\calC)$, for some small $\infty$-category $\calC$. Up to equivalence, there are precisely two $(-1)$-truncated spaces: $\emptyset$ and $\ast$. Consequently, $\tau_{\leq -1} \SSet$ is equivalent to the 
two-object $\infty$-category $\Delta^1$. It follows that $\calP_{\leq -1}(\calC)$ is equivalent
to $\Fun( \calC^{op},\Delta^1)$.

Let $\widetilde{X}$ denote the collection of sieves on $\calC$, ordered by inclusion. Then, identifying a
a functor $f: \calC \rightarrow \Delta^1$
with the sieve $f^{-1} \{0\} \subseteq \calC$, we deduce that
$\Fun(\calC, \Delta^1)$ is isomorphic to the nerve $\Nerve(\widetilde{X})$.

Without loss of generality, we may identify $\calX$ with the essential image of a localization
functor $L: \Nerve(\widetilde{X}) \rightarrow \Nerve(\widetilde{X})$. The map $L$ may be identified with a map
of partially ordered sets from $\widetilde{X}$ to itself. Unwinding the definitions, we find that the condition that $L$ be a left exact localization is equivalent to the following three properties:

\begin{itemize}
\item[$(A)$] The map $L: \widetilde{X} \rightarrow \widetilde{X}$ is idempotent.
\item[$(B)$] For each $U \in \widetilde{X}$, $U \subseteq L(U)$.
\item[$(C)$] The map $L: \widetilde{X} \rightarrow \widetilde{X}$ preserves
finite intersections (since $\calX$ is a {\em left exact} localization of $\Nerve(\widetilde{X})$.)
\end{itemize}

Let $\calU = \{ U \in \widetilde{X}: LU = U \}$.
Then it is easy to see that $\calX$ is equivalent to the nerve $\Nerve(\calU)$, and that the partially ordered set $X$ satisfies conditions $(1)$ and $(2)$ of Definition \ref{deflocale}. Therefore $\calU$ is a locale, so that colimits in $\Nerve( \calU)$ are universal by Proposition \ref{colim0univ}. This proves that $(3) \Rightarrow (4)$. 

Now suppose that $(4)$ is satisfied. Using Proposition \ref{colim0univ}, we may suppose without loss of generality that $\calX = \Nerve(\calU)$, where $\calU$ is a locale. We observe that
$\calX$ is itself small. Let us say that
a sieve $\{ U_{\alpha} \rightarrow U \}$ on an object $U \in \calX$ is {\it covering} if 
$$ U = \bigcup_{\alpha} U_{\alpha} $$ in $\calU$. Using the assumption that $\calU$ is a locale, it is easy to see that the collection of covering sieves determines a Grothendieck topology on
$\calX$. The $\infty$-category $\calP_{\leq -1}(\calX)$ can be identified with the nerve of
the partially ordered set of all downward-closed subsets $\calU_0 \subseteq \calU$. Moreover, an object of $\calP_{\leq -1}(\calX)$ belongs to $\Shv_{\leq -1}(\calX)$ if and only if the corresponding subset $\calU_0 \subseteq \calU$ is stable under joins. Every such subset $\calU_0 \subseteq \calU$ has a largest element $U \in \calU$, and we then have an identification
$\calU_0 = \{ V \in \calU: V \leq U\}$. It follows that $\Shv_{\leq -1}(\calX)$ is equivalent
to the nerve of the partially ordered set $\calU$, which is $\calX$. This proves $(1)$, and concludes the argument.
\end{proof}

We may summarize the results of this section as follows:

\begin{corollary}\label{charlocale}
An $\infty$-category $\calX$ is a $0$-topos if and only if it is equivalent to
$\Nerve(\calU)$, where $\calU$ is a locale.
\end{corollary}

\begin{remark}
Coproducts in a $0$-topos are typically {\em not} disjoint.
\end{remark}

In classical topos theory, there are functorial constructions for passing back and forth between topoi and locales. Given a locale $\calU$ (such as the locale $\calU(X)$ of open subsets of a topological space $X$), one 
may define a topos $\calX$ of {\it sheaves $($ of sets $)$ on $\calU$}. The original locale $\calU$ may then be recovered as the partially ordered set of subobjects of the final object of $\calX$. In fact, for {\em any} topos $\calX$, the partially ordered set $\calU$ of subobjects of the final object forms a locale. In general, $\calX$ cannot be recovered as the category of sheaves on $\calU$; this is true if and only if
$\calX$ is a {\it localic} topos\index{gen}{localic topos}\index{gen}{topos!localic}: that is, if and only if $\calX$ is generated under colimits
by the collection of subobjects of the final object ${1}_{\calX}$. In \S \ref{chap6sec4} we will discuss a generalization of this picture, which will allow us to pass between $m$-topoi and $n$-topoi for any $m \leq n$.

\subsection{Giraud's Axioms for $n$-Topoi}\label{provengiraud}

In \S \ref{axgir}, we sketched an axiomatic approach to the theory of $\infty$-topoi. The axioms
we introduced were closely parallel to Giraud's axioms for ordinary topoi, with one important difference. If $\calX$ is an $\infty$-topos, then {\em every} groupoid object of $\calX$ is effective. If $\calX$ is an ordinary topos, then a groupoid $U_{\bigdot}$ is effective only if the diagram
$$\xymatrix{ U_{1} \ar@<.5ex>[r] \ar@<-.5ex>[r] & U_0 }$$
exhibits $U_{1}$ as an equivalence relation on $U_0$. Our first goal in this section is to formulate an analogue of this condition, which will lead us to an axiomatic description of $n$-topoi for all
$0 \leq n \leq \infty$.

\begin{definition}\index{gen}{groupoid object!$n$-efficient}
Let $\calX$ be an $\infty$-category and $U_{\bigdot}$ a groupoid object of $\calX$. We will say that
$U_{\bigdot}$ is {\it $n$-efficient} if the natural map $$ U_{1} \rightarrow U_0 \times U_0 $$
(which is well-defined up to equivalence) is $(n-2)$-truncated.
\end{definition}

\begin{remark}
By convention, we regard every groupoid object as $\infty$-efficient.
\end{remark}

\begin{example}
If $\calC$ is (the nerve of) an ordinary category, then giving a $1$-efficient groupoid object $U_{\bigdot}$ of $\calC$ is equivalent to giving an object $U_0$ of $\calC$ and an equivalence relation $U_1$ on $U_0$.
\end{example}

\begin{proposition}\label{storage}
An $\infty$-category $\calX$ is equivalent to an $n$-category if and only if every
effective groupoid in $\calX$ is $n$-efficient.
\end{proposition}

\begin{proof}
Suppose first that $\calC$ is equivalent to an $n$-category. Let $U_{\bigdot}$ be an effective groupoid in $\calX$. Then $U_{\bigdot}$ has a colimit $U_{-1}$. The existence of a pullback diagram
$$ \xymatrix{ U_1 \ar[r] \ar[d] & U_0 \ar[d] \\
U_0 \ar[r] & U_{-1} }$$
implies that the map $f': U_{1} \rightarrow U_0 \times U_0$ is a pullback of the diagonal map
$f: U_{-1} \rightarrow U_{-1} \times U_{-1}$. We wish to show that $f'$ is $(n-2)$-truncated.
By Lemma \ref{tunc} it suffices to show that $f$ is $(n-2)$-truncated. By Lemma \ref{trunc}, this is equivalent to the assertion that $U_{-1}$ is $(n-1)$-truncated. Since $\calC$ is equivalent to an $n$-category, every object of $\calC$ is $(n-1)$-truncated.

Now suppose that every effective groupoid in $\calX$ is $n$-efficient. Let $U \in \calX$ be
an object; we wish to show that $U$ is $(n-1)$-truncated.
The constant simplicial object $U_{\bigdot}$ taking the value $U$ is an effective groupoid, and therefore $n$-efficient. It follows that the diagonal map $U \rightarrow U \times U$ is $(n-2)$-truncated. Lemma \ref{trunc} implies that $U$ is $(n-1)$-truncated as desired.
\end{proof}

We are now almost ready to supply the ``hard'' step in the proof of Theorem \ref{nchar} (namely, the implication $(6) \Rightarrow (1)$). We first need a slightly technical lemma, whose proof requires routine cardinality estimates.

\begin{lemma}\label{nocake}
Let $\calX$ be a presentable $\infty$-category in which colimits are universal. There exists a regular cardinal $\tau$ such that $\calX$ is $\tau$-accessible, and the full subcategory of $\calX^{\tau} \subseteq \calX$ spanned by the $\tau$-compact objects is stable under the formation of subobjects and finite limits.
\end{lemma}

\begin{proof}
Choose a regular cardinal $\kappa$ such that $\calX$ is $\kappa$-accessible. We observe that, up to equivalence, there are a bounded number of $\kappa$-compact objects of $\calX$, and therefore a bounded number of {\em subobjects} of $\kappa$-compact objects of $\calX$. Now choose an uncountable regular cardinal $\tau \gg \kappa$ such that:
\begin{itemize}
\item[$(1)$] The $\infty$-category $\calX^{\kappa}$ is essentially $\tau$-small.
\item[$(2)$] For each $X \in \calX^{\kappa}$ and each monomorphism $i: U \rightarrow X$, $U$ is $\tau$-compact.
\end{itemize}
It is clear that $\calX$ is $\tau$-accessible, and $\calX^{\tau}$ is stable under finite limits (in fact, $\kappa$-small limits) by Proposition \ref{tcoherent}. To complete the proof, we must show that $\calX^{\tau}$ is stable under the formation of subobjects. Let $i: U \rightarrow X$ be a monomorphism, where $X$ is $\tau$-compact. Since $\calX$ is $\kappa$-accessible, we can write
$X$ as the colimit of $\kappa$-filtered diagram $p: \calJ \rightarrow \calX^{\kappa}$. Since $X$ is $\tau$-compact, it is a retract of the colimit $X'$ of some $\tau$-small subdiagram $p| \calJ'$.
Since $\tau$ is uncountable, we can use Proposition \ref{autokan} to write
$X$ as the colimit of a $\tau$-small diagram $\Idem \rightarrow \calX$, which carries
the unique object of $\Idem$ to $X'$. Since colimits in $\calX$ are universal, it follows that
$U$ can be written as a $\tau$-small colimit of a diagram $\Idem \rightarrow \calX$ which takes the value $U' = U \times_{X} X'$. It will therefore suffice to prove that $U'$ is $\tau$-compact.
Invoking the universality of colimits once more, we observe that $U'$ is a $\tau$-small colimit of objects of the form $U'' = U' \times_{X'} p(J)$, where $J$ is an object of $\calJ'$. We now observe
that $U''$ is a subobject of $p(J) \in \calX^{\kappa}$, and is therefore $\tau$-compact by assumption $(2)$. It follows that $U'$, being a $\tau$-small colimit of $\tau$-compact objects of $\calX$, is also $\tau$-compact.
\end{proof}

\begin{proposition}\label{diamondstep}
Let $0 < n < \infty$, and let $\calX$ be an $\infty$-category satisfying the following conditions:
\begin{itemize}
\item[$(i)$] The $\infty$-category $\calX$ is presentable.
\item[$(ii)$] Colimits in $\calX$ are universal.
\item[$(iii)$] Coproducts in $\calX$ are disjoint.
\item[$(iv)$] The effective groupoid objects of $\calX$ are precisely the $n$-efficient groupoids.
\end{itemize}
Then there exists a small $n$-category $\calC$ which admits finite limits, a Grothendieck topology on $\calC$, and an equivalence $\calX \rightarrow \Shv_{\leq n-1}(\calC)$.
\end{proposition}

\begin{proof}
Without loss of generality we may suppose that $\calX$ is minimal. Choose a regular cardinal $\kappa$ such that $\calX$ is $\kappa$-accessible, and the full subcategory $\calC \subseteq \calX$ spanned by the $\kappa$-compact objects of $\calX$ is stable under the formation of subobjects and finite limits (Lemma \ref{nocake}).
We endow $\calC$ with the canonical topology induced by the inclusion
$\calC \subseteq \calX$. According to Theorem \ref{charpresheaf}, there is an (essentially unique) colimit preserving functor $F: \calP(\calC) \rightarrow \calX$ such that $F \circ j$ is equivalent
to the inclusion $\calC \subseteq \calX$, where $j: \calC \rightarrow \calP(\calC)$ denotes the Yoneda embedding. The proof of Theorem \ref{pretop} shows that $F$ has a fully faithful right adjoint $G: \calX \rightarrow \calP(\calC)$. We will complete the proof by showing that the essential image of $G$ is precisely $\Shv_{\leq n-1}(\calC)$. 

Since $\calX$ is equivalent to an $n$-category (Proposition \ref{storage}) and $G$ is left exact, we conclude that $G$ factors through $\calP_{\leq n-1}(\calC)$. It follows from Proposition \ref{preciselate} that $G$ factors through $\Shv_{\leq n-1}(\calC)$. Let $\calX' \subseteq \Shv_{\leq n-1}(\calC)$ denote the essential image of $G$. To complete the proof, it will suffice to show that $\calX' = \Shv_{\leq n-1}(\calC)$. Let $\emptyset$ be an initial object of $\calX$. The space $\bHom_{\calX}(X, \emptyset)$ is contractible if $X$ is an initial object of $\calX$, and empty otherwise (Lemma \ref{sumoto}). It follows from Proposition \ref{suture} that $G(\emptyset)$ is an initial object of $\Shv_{\leq n-1}(\calC)$. 

We next claim that $\calX'$ is stable under small coproducts in $\Shv_{\leq n-1}(\calC)$. It will suffice to show that the map $G$ preserves coproducts. Let $\{ U_{\alpha} \}$ be a small collection of objects of $\calX$ and $U$ their coproduct in $\calX$. According to Lemma \ref{sumdescription}, we have a pullback diagram
$$ \xymatrix{ V_{\alpha,\beta} \ar[r]^{\phi} \ar[d]^{\phi'} & U_{\alpha} \ar[d] \\
U_{\beta} \ar[r] & U, }$$
where $V_{\alpha,\beta}$ is an initial object of $\calX$ if $\alpha \neq \beta$, while
$\phi$ and $\phi'$ are equivalences if $\alpha = \beta$.
The functor $G$ preserves all limits, so that the diagram
$$ \xymatrix{ G(V_{\alpha,\beta}) \ar[r] \ar[d] & G(U_{\alpha}) \ar[d] \\
G( U_{\beta} ) \ar[r] & G(U) }$$
is a pullback in $\Shv_{\leq n-1}(\calC)$. Let $U'$ denote a coproduct
of the objects $i(U_{\alpha})$ in $\Shv_{\leq n-1}(\calC)$, and let $g: U' \rightarrow G(U)$ be the induced map. Since colimits in $\calX$ are universal, we obtain a natural identification of
$U' \times_{G(U)} U'$ with the coproduct
$$ \coprod_{\alpha, \beta} ( G(U_{\alpha}) \times_{ G(U) } G(U_{\beta}) )
\simeq \coprod_{\alpha} U_{\alpha} \simeq U',$$
where the second equivalence follows from our observation that $G$ preserves initial objects.
Applying Lemma \ref{trunc}, we deduce that $g$ is a monomorphism.

To prove that $g$ is an equivalence, it will suffice to show that the map 
$$\pi_0 U'(C) \rightarrow \pi_0 G(U)(C) = \pi_0 \bHom_{\calX}(C,U)$$ is surjective for every object $C \in \calC$.
Since colimits in $\calX$ are universal,
every map $h: C \rightarrow U$ can be written as a coproduct of maps $h_{\alpha}: C_{\alpha} \rightarrow U_{\alpha}$. Each $C_{\alpha}$ is a subobject of $C$ (Lemma \ref{corsumoto}) and therefore belongs to $\calC$. Let $h'_{\alpha} \in \pi_0 U'(C_{\alpha})$ denote the homotopy class
of the composition $G(C_{\alpha}) \stackrel{h_{\alpha}}{\rightarrow} G(U_{\alpha}) \rightarrow U'$.
Since the topology on $\calC$ is canonical, Lemma \ref{canonicalcoproducts} implies that
$\pi_0 U'(C) \simeq \prod_{\alpha} \pi_0 U'(C_{\alpha})$ contains an element $h'$ which
restricts to each $h'_{\alpha}$. It is now clear that $h$ is the image of $h'$ under
the map $\pi_0 U'(C) \rightarrow \pi_0 \bHom_{\calX}(C,U)$. 

We will prove the following result by induction on $k$: if there exists a $k$-truncated morphism
$f: X \rightarrow Y$, where $Y \in \calX'$ and $X \in \Shv_{\leq n-1}(\calC)$, then $X \in \calX'$.
Taking $k = n-1$ and $Y$ to be a final object of $\Shv_{\leq n-1}(\calC)$ (which belongs to
$\calX'$ because $\calC$ contains a final object), we conclude that every object of
$\Shv_{\leq n-1}(\calC)$ belongs to $\calX'$, which completes the proof.

If $k = -2$, then $f$ is an equivalence so that $X \in \calX'$ as desired.
Assume now that $k \geq -1$. Since $\calX'$ contains the essential image of the Yoneda embedding and is stable under coproducts, there exists an effective epimorphism $p: U \rightarrow X$ in $\Shv_{\leq n-1}(\calC)$, where $U \in \calX'$. Let $\overline{U}_{\bigdot}$ be a \Cech nerve of $p$ in 
$\Shv_{\leq n-1}(\calC)$, and $U_{\bigdot}$ the associated groupoid object. We claim that
$U_{\bigdot}$ is a groupoid object of $\calX'$. Since $\calX'$ is stable under limits in
$\Shv_{\leq n-1}(\calC)$, it suffices to prove that $U_0 = U$ and $U_1 = U \times_{X} U$ belong to $\calX'$. We now observe that there exists a pullback diagram
$$ \xymatrix{ U \times_{X} U \ar[r]^-{\delta'} \ar[d] & U \times_{Y} U \ar[d] \\
X \ar[r]^-{\delta} & X \times_{Y} X. }$$
Since $f$ is $k$-truncated, $\delta$ is $(k-1)$-truncated (Lemma \ref{trunc}), so 
that $\delta'$ is $(k-1)$-truncated. Since $U \times_{Y} U$ belongs to $\calX'$ (because
$\calX'$ is stable under limits), our inductive hypothesis allows us to conclude that $U \times_{X} U
\in \calX'$, as desired.  

We observe that $U_{\bigdot}$ is an $n$-efficient groupoid object of $\calX'$. Invoking assumption $(iv)$, we conclude that $U_{\bigdot}$ is effective in $\calX'$. Let $X' \in \calX'$ be a colimit of $U_{\bigdot}$ in $\calX'$, so that we have a morphism $u: X \rightarrow X'$ in $\Shv_{\leq n-1}(\calC)_{U_{\bigdot}/}$. To complete the proof that $X \in \calX'$, it will suffice to show that $u$ is an equivalence. Since $u$ induces an equivalence
$$ U \times_{X} U \rightarrow U \times_{X'} U,$$
it is a monomorphism (Lemma \ref{trunc}). It will therefore suffice to show that $u$ is an effective epimorphism. We have a commutative diagram
$$ \xymatrix{ U \ar[dr]^{p} \ar[rr]^{p'} & & X' \\
& X \ar[ur]^{u} & }$$
where $p$ is an effective epimorphism; it therefore suffices to show that $p'$ is an effective epimorphism which follows immediately from Proposition \ref{preciselate}.
\end{proof}

\begin{remark}
Proposition \ref{diamondstep} is valid also for $n=0$, but is almost vacuous: coproducts
in a $0$-topos $\calX$ are never disjoint unless $\calX$ is trivial (equivalent to the final $\infty$-category $\ast$).
\end{remark}

\begin{remark}
In a certain respect, the theory of $\infty$-topoi is {\em simpler than} the theory of ordinary topoi: in an $\infty$-topos, {\em every} groupoid object is effective; it is not necessary to impose any additional conditions like $n$-efficiency. The absense of this condition gives the theory
of $\infty$-topoi a slightly different flavor than ordinary topos theory. In an $\infty$-topos, we are free
to form quotients of objects not only by equivalence relations, but by arbitrary groupoid actions.
In geometry, this extra flexibility allows the construction of useful objects such as orbifolds and algebraic stacks, which are useful in a variety of mathematical situations.

One can imagine weakening the gluing conditions even further, and
considering axioms having the form ``every category object is
effective''. This seems like a natural approach to a theory
of topos-like $(\infty,\infty)$-categories. However, we will not pursue the matter any further here.
\end{remark}

It follows from Proposition \ref{diamondstep} (and arguments to be given in \S \ref{provengiraudeasy}) that every left-exact localization of a presheaf $n$-category
$\calP_{\leq n-1}(\calC)$ can also be obtained as an $n$-category of sheaves. According to the next two results, this is no accident: every left exact localization of $\calP_{\leq n-1}(\calC)$ is topological, and the topological localizations of $\calP_{\leq n-1}(\calC)$ correspond precisely to the Grothendieck topologies on $\calC$ (provided that $\calC$ is an $n$-category).

\begin{proposition}\label{alltoploc}
Let $\calX$ be a presentable $n$-category, $0 \leq n < \infty$, and suppose that colimits in $\calX$ are universal. Let $L: \calX \rightarrow \calY$ be a left exact localization. Then $L$ is a topological localization.
\end{proposition}

\begin{proof}
Let $S$ denote the collection of all monomorphisms $f: U \rightarrow V$ in $\calX$ such that $Lf$ is an equivalence. Since $L$ is left exact, it is clear that $S$ is stable under pullback. Let $\overline{S}$ be the strongly saturated class of morphisms generated by $S$. Proposition \ref{swimmer} implies that $\overline{S}$ is stable under pullback, and therefore topological.
Proposition \ref{toplocsmall} implies that $\overline{S}$ is generated by a (small) set of morphisms. 
Let $\calX' \subseteq \calX$ denote the full subcategory spanned by $\overline{S}$-local objects. According to Proposition \ref{local}, $\calX'$ is an accessible localization of $\calX$; let $L'$ denote the associated localization functor. Since $Lf$ is an equivalence for each $f \in \overline{S}$, the localization $L$ is equivalent to the composition
$$ \calX \stackrel{L'}{\rightarrow} \calX' \stackrel{L|\calX'}{\rightarrow} \calY.$$
We may therefore replace $\calX$ by $\calX'$ and thereby reduce to the case where
$S$ consists precisely of the equivalences in $\calX$; we wish to prove that $L$ is an equivalence.

We now prove the following claim: if $f: X \rightarrow Y$ is a $k$-truncated morphism in $\calC$
such that $Lf$ is an equivalence, then $f$ is an equivalence. The proof goes by induction on 
$k$. If $k = -1$, then $f$ is a monomorphism, and so belongs to $S$, and is therefore an equivalence. Suppose that $k \geq 0$. Let
$\delta: X \rightarrow X \times_{Y} X$ be the diagonal map (which is well-defined up to equivalence). According to Lemma \ref{trunc}, $\delta$ is $(k-1)$-truncated. Since
$L$ is left exact, $L(\delta)$ can be identified with a diagonal map $LX \rightarrow LX \times_{LY} LX$, which is therefore an equivalence. The inductive hypothesis implies that $\delta$ is an equivalence. Applying Lemma \ref{trunc} again, we deduce that $f$ is a monomorphism, so that $f \in S$ and is therefore an equivalence as noted above. 

Since $\calX$ is an $n$-category, every morphism in $\calX$ is $(n-1)$-truncated. We conclude that for {\em every} morphism $f $ in $\calX$, $f$ is an equivalence if and only if $Lf$ is an equivalence. Since $L$ is a localization functor, it must be an equivalence.
\end{proof}

\subsection{$n$-Topoi and Descent}\label{provengiraudeasy}

Let $\calX$ be an $\infty$-category which admits finite limits, and let
$\calO_{\calX}$ denote the functor $\infty$-category $\Fun(\Delta^1, \calX)$ equipped
with the Cartesian fibration $e: \calO_{\calX} \rightarrow \calX$ (given by evaluation
at $\{1\} \subseteq \Delta^1$), as in \S \ref{axgir}. Let $F: \calX^{op} \rightarrow \widehat{\Cat}_{\infty}$ be a functor which classifies $e$; informally, $F$ associates to each object $U \in \calX$ the $\infty$-category $\calX_{/U}$. According to Theorem \ref{charleschar}, 
$\calX$ is an $\infty$-topos if and only if $F$ the functor
$F$ preserves limits, and factors through $\LPres \subseteq \widehat{\Cat}_{\infty}$.
The assumption that $F$ preserves limits can be viewed as a descent condition: it asserts that if $X \rightarrow U$ is a morphism of $\calX$, and $U$ is decomposed into ``pieces'' $U_{\alpha}$, then $X$ can be canonically reconstructed from the ``pieces'' $X \times_{U} U_{\alpha}$.
The goal of this section is to obtain a similar characterization of the class of $n$-topoi, for $0 \leq n < \infty$. 

We begin by considering the case where $\calX$ is the (nerve of) the category of sets. In this case, we can think of $F$ as a contravariant functor from sets to categories,
which carries a set $U$ to the category $\Set_{/U}$. This functor does not preserve pullbacks: given a pushout square $$ \xymatrix{ & X \ar[dl] \ar[dr] &  \\
Y \ar[dr] & & Z \ar[dl] \\
& Y \amalg_{X} Z & }$$
in the category $\Set$, there is an associated functor
$$ \theta: \Set_{/ Y \amalg_{X} Z} \rightarrow \Set_{/Y} \times_{ \Set_{/X} } \Set_{/Z}$$
(here the right hand side indicates a {\em homotopy} fiber product of categories). The functor $\theta$ is generally not an equivalence of categories: for example, $\theta$ fails to be an 
equivalence if $Y = Z = \ast$, provided that $X$ has cardinality at least $2$. However, $\theta$ is always fully faithful. Moreover, we have the following result:

\begin{fact}\label{stiro}
The functor $\theta$ induces an isomorphism of partially ordered sets
$$ \Sub(Y \amalg_{X} Z) \rightarrow \Sub(Y) \times_{ \Sub(X) } \Sub(Z)$$
where $\Sub(M)$ denotes the partially ordered set of subsets of $M$.
\end{fact}

In this section, we will show that an appropriate generalization of Fact \ref{stiro} can be used to characterize the class of $n$-topoi, for all $0 \leq n \leq \infty$. First, we need to introduce some terminology.

\begin{notation}\label{sumahum}
Let $\calX$ be an $\infty$-category which admits pullbacks, and let $0 \leq n \leq \infty$. We let
$\calO_{\calX}^{n}$ denote the full subcategory of $\calO_{\calX}$ spanned by morphisms
$f: U \rightarrow X$ which are $(n-2)$-truncated, and $\calO_{\calX}^{(n)} \subseteq
\calO^{n}_{\calX}$ the subcategory whose objects are $(n-2)$-truncated morphisms
in $\calX$, and whose morphisms are Cartesian transformations (see Notation
\ref{ugaboo}).\index{not}{OcalXn@$\calO_{\calX}^{n}$}\index{not}{OcalX(n)@$\calO_{\calX}^{(n)}$}
\end{notation}

\begin{example}\label{trivexa}
Let $\calX$ be an $\infty$-category which admits pullbacks. Then $\calO^{0}_{\calX}$
is the full subcategory of $\calO_{\calX}$ spanned by the final objects in each fiber of the morphism $p: \calO_{\calX} \rightarrow \calX$. Since $p$ is a coCartesian fibration (Corollary \ref{tweezegork}), Proposition \ref{topaz} asserts that the restriction $p | \calO_{\calX}^{0}$ is a trivial fibration of simplicial sets.
\end{example}

\begin{lemma}\label{nstab}
Let $\calX$ be a presentable $\infty$-category in which colimits are universal and coproducts are disjoint, and let $n \geq -2$. Then the class of $n$-truncated morphisms in $\calX$ is stable under small coproducts.
\end{lemma}

\begin{proof}
The proof is by induction on $n$, where the case $n=-2$ is obvious. Suppose that the
$\{ f_{\alpha}: X_{\alpha} \rightarrow Y_{\alpha} \}$ is a family of $n$-truncated morphisms
in $\calX$ having coproduct $f: X \rightarrow Y$. Since colimits in $\calX$ are universal, we conclude that $X \times_{Y} X$ can be written as a coproduct 
$$\coprod_{ \alpha, \beta} (X_{\alpha} \times_{Y} X_{\beta})
\simeq \coprod_{\alpha,\beta} (X_{\alpha} \times_{Y_{\alpha} } (Y_{\alpha} \times_{Y} Y_{\beta} )
\times_{Y_{\beta}} X_{\beta}). $$
Applying Lemma \ref{sumdescription}, we can rewrite this coproduct as
$$ \coprod_{\alpha} (X_{\alpha} \times_{ Y_{\alpha} } X_{\alpha}). $$
Consequently, the diagonal map $\delta: X \rightarrow X \times_{Y} X$ is a coproduct
of diagonal maps $\{ \delta_{\alpha}: X_{\alpha} \rightarrow X_{\alpha} \times_{Y_{\alpha} } X_{\alpha} \}$. Applying Lemma \ref{trunc}, we deduce that each $\delta_{\alpha}$ is $(n-1)$-truncated, so that $\delta$ is $(n-1)$-truncated by the inductive hypothesis. We now apply Lemma \ref{trunc} again to deduce that $f$ is $n$-truncated, as desired.
\end{proof}

Combining Lemmas \ref{ib1}, \ref{ib2}, \ref{ib3}, and \ref{nstab}, we deduce the following analogue of Theorem \ref{charleschar}.

\begin{theorem}\label{countercon}
Let $\calX$ be a presentable $\infty$-category in which colimits are universal and coproducts are disjoint. The following conditions are equivalent:
\begin{itemize}
\item[$(1)$] For every 
pushout diagram
$$ \xymatrix{ f \ar[r]^{\alpha} \ar[d]^{\beta} & g \ar[d]^{\beta'} \\
f' \ar[r]^{\alpha'} & g' }$$
in $\calO^{n}_{\calX}$, if $\alpha$ and $\beta$ are Cartesian transformations, then
$\alpha'$ and $\beta'$ are also Cartesian transformations.

\item[$(2)$] The class of $(n-2)$-truncated morphisms
in $\calX$ is local.

\item[$(3)$] The Cartesian fibration $\calO^{n}_{\calX} \rightarrow \calX$ is classified
by a limit-preserving functor $\calX^{op} \rightarrow \hat{\Cat}_{\infty}$

\item[$(4)$] The right fibration $\calO^{(n)}_{\calX} \rightarrow \calX$ is classified by a limit-preserving functor $\calX^{op} \rightarrow \hat{\SSet}$.

\item[$(5)$] Let $K$ be a small simplicial set and 
$\overline{\alpha}: \overline{p} \rightarrow \overline{q}$ a natural transformation of colimit diagrams
$\overline{p}, \overline{q}: K^{\triangleright} \rightarrow \calX$. Suppose that
$\alpha = \overline{\alpha} | K$ is a Cartesian transformation, and that $\alpha(x)$ is $(n-2)$-truncated for every vertex $x \in K$. Then $\overline{\alpha}$ is a Cartesian transformation, and
$\overline{\alpha}(\infty)$ is $(n-2)$-truncated, where $\infty$ denotes the cone point of
$K^{\triangleright}$.
\end{itemize}
\end{theorem}

Our next goal is to establish the implication $(2) \Rightarrow (4)$ of Theorem \ref{nchar}.
We will deduce this from the equivalence $(2) \Leftrightarrow (3)$ (which we have already established) together with Propositions \ref{tigre} and \ref{tigress} below.

\begin{proposition}\label{tigre}
Let $\calX$ be an $n$-topos, $0 \leq n \leq \infty$. Then colimits in $\calX$ are universal.
\end{proposition}

\begin{proof}
Using Lemma \ref{tryme}, we may reduce to the case $\calX = \calP_{\leq n-1}(\calC)$ for some small $\infty$-category $\calC$. Using Proposition \ref{limiteval}, we may further reduce to the case where $\calX = \tau_{\leq n-1} \SSet$. 

Let $f: X \rightarrow Y$ be a map of $(n-1)$-truncated spaces, and let $f^{\ast}: \SSet^{/Y} \rightarrow \SSet^{/X}$ be a pullback functor.
Since $\calX$ is stable under limits in $\SSet$, $f^{\ast}$ restricts to give a functor $\calX^{/Y} \rightarrow \calX^{/X}$; we wish to prove that this restricted functor commutes with colimits. 
We observe
that $\calX^{/X}$ and $\calX^{/Y}$ can be identified with the full subcategories of
$\SSet^{/X}$ and $\SSet^{/Y}$ spanned by the $(n-1)$-truncated objects, by Lemma \ref{trunccomp}. Let $\tau_{X}: \SSet^{/X} \rightarrow \calX^{/X}$ and $\tau_{Y}: \SSet^{/Y} \rightarrow \calX^{/Y}$ denote left adjoints to the inclusions. The functor $f^{\ast}$ preserves all colimits (Lemma \ref{sugartime}) and all limits (since
$f^{\ast}$ has a left adjoint). Consequently, Proposition \ref{compattrunc} implies that
$\tau_{X} \circ f^{\ast} \simeq f^{\ast} \circ \tau_{Y}$.

Let $p: K^{\triangleright} \rightarrow \calX^{/Y}$ be a colimit diagram. We wish to show that
$f^{\ast} \circ p$ is a colimit diagram. According to Remark \ref{localcolim}, we may assume that $p = \tau_{Y} \circ p'$, for some colimit diagram $p': K^{\triangleright} \rightarrow \SSet^{/Y}$. Since colimits in $\SSet$ are universal (Lemma \ref{sugartime}), the composition
$f^{\ast} \circ p': K^{\triangleright} \rightarrow \SSet^{/X}$ is a colimit diagram. Since $\tau_{X}$ preserves colimits, we conclude that $\tau_{X} \circ f^{\ast} \circ p': K^{\triangleright} \rightarrow \calX^{/X}$ is a colimit diagram, so that $f^{\ast} \circ \tau_{Y} \circ p' = f^{\ast} \circ p$ is also a colimit diagram, as desired.
\end{proof}

\begin{proposition}\label{tigress}
Let $\calY$ be an $\infty$-topos and let $\calX = \tau_{\leq n} \calY$, $0 \leq n \leq \infty$. Then the class of $(n-2)$-truncated morphisms in $\calX$ is local.
\end{proposition}

\begin{proof}
Combining Propositions \ref{hintdescent0}, \ref{torque}, and Lemma \ref{nstab}, we conclude
that the class of $(n-2)$-truncated morphisms in $\calY$ is local. Consequently, the Cartesian fibration
$\calO_{\calY}^{n} \rightarrow \calY$ is classified by a colimit preserving functor
$ F: \calY \rightarrow \hat{\Cat}_{\infty}^{op}$. It follows that $\calO_{\calX}^{(n)} \rightarrow \calX$
is classified by $F | \calX$. To prove that $F | \calX$ is colimit-preserving, it will suffice to show that 
$F$ is equivalent to $F \circ \tau_{\leq n}$; in other words, that $F$ carries
each $n$-truncation $Y \rightarrow \tau_{\leq n} Y$ to an equivalence in $\hat{\Cat}_{\infty}^{op}$.
Replacing $\calY$ by $\calY_{/ \tau_{\leq n} Y}$, we reduce to Lemma \ref{nicelemma}.
\end{proof}

We conclude this section by proving the following generalization of Proposition \ref{lemonade2}, which also establishes the implication $(4) \Rightarrow (6)$ of Theorem \ref{nchar}. We will assume $n > 0$; the case $n = 0$ was analyzed in \S \ref{0topoi}.

\begin{lemma}\label{corsumoto}
Let $\calX$ be a presentable $\infty$-category in which colimits are universal, and let
$f: \emptyset \rightarrow X$ be a morphism in $\calX$, where $\emptyset$ is an initial object of $\calX$. Then $f$ is a monomorphism.
\end{lemma}

\begin{proof}
Let $Y$ be an arbitrary object of $\calX$, we wish to show that composition with $f$
induces a $(-1)$-truncated map
$$ \bHom_{\calX}( Y, \emptyset) \rightarrow \bHom_{\calX}(Y, X).$$
If $Y$ is an initial object of $\calX$, then both sides are contractible; otherwise
the left side is empty (Lemma \ref{sumoto}).
\end{proof}

\begin{proposition}\label{ncharles}
Let $1 \leq n \leq \infty$, and let $\calX$ be a presentable $n$-category. Suppose that
colimits in $\calX$ are universal, and that the class of $(n-2)$-truncated morphisms in $\calX$ is local. Then $\calX$ satisfies the $n$-categorical Giraud axioms:
\begin{itemize}
\item[$(i)$] The $\infty$-category $\calX$ is equivalent to a presentable $n$-category.
\item[$(ii)$] Colimits in $\calX$ are universal.
\item[$(iii)$] Coproducts in $\calX$ are disjoint.
\item[$(iv)$] Every $n$-efficient groupoid object of $\calX$ is effective.
\end{itemize}
\end{proposition}

\begin{proof}
Axioms $(i)$ and $(ii)$ hold by assumption.
To show that coproducts in $\calX$ are disjoint, let us consider an arbitrary pair of objects
$X, Y \in \calX$, and let $\emptyset$ denote an initial object of $\calX$. Let $f: \emptyset \rightarrow X$ be a morphism (unique up to homotopy, since $\emptyset$ is initial). Since
colimits in $\calX$ are universal, $f$ is a monomorphism (Lemma \ref{corsumoto}) and therefore belongs to
$\calO_{\calX}^{n}$, since $n \geq 1$. We observe that
$\id_{\emptyset}$ is an initial object of $\calO_{\calX}$, so we can form a pushout diagram
$$ \xymatrix{ \id_{\emptyset} \ar[r]^{\alpha} \ar[d]^{\beta} & \id_{Y} \ar[d]^{\beta'}  \\
f \ar[r]^{\alpha'} & g }$$
in $\calO^{n}_{\calX}$. It is clear that $\alpha$ is a Cartesian transformation, and Lemma \ref{sumoto} implies that $\beta$ is Cartesian as well. Invoking Theorem \ref{countercon}, we deduce that $\alpha'$ is
a Cartesian transformation. But $\alpha'$ can be identified with a pushout diagram
$$ \xymatrix{ \emptyset \ar[r] \ar[d] & Y \ar[d] \\
X \ar[r] & X \amalg Y. }$$
This proves $(iii)$. 

Now suppose that $U_{\bigdot}$ is an $n$-efficient groupoid object of $\calX$; we wish to prove that $U_{\bigdot}$ is effective. Let $\overline{U}_{\bigdot}  : \Nerve(\cDelta_{+})^{op} \rightarrow \calX$ be a colimit of $U_{\bigdot}$. Let $U'_{\bigdot}: \Nerve(\cDelta_{+})^{op} \rightarrow \calX$
be the result of composing $\overline{U}_{\bigdot}$ with the shift functor
$$ \cDelta_{+} \rightarrow \cDelta_{+}$$
$$ J \mapsto J \amalg \{ \infty \}.$$
(In other words, $U'_{\bigdot}$ is the shifted simplicial object given by
$U'_{n} = U_{n+1}$.)
Lemma \ref{bclock} asserts that $U'_{\bigdot}$ is a colimit diagram in $\calX$.
We have a transformation $\overline{\alpha}: U'_{\bigdot} \rightarrow \overline{U}_{\bigdot}$.
Let $\overline{V}_{\bigdot}$ denote the constant augmented simplicial object of $\calX$ taking the value $U_0$, so that we have a natural transformation $\overline{\beta}: U'_{\bigdot} \rightarrow \overline{V}_{\bigdot}$. Let $\overline{W}_{\bigdot}$ denote a product of $\overline{U}_{\bigdot}$ and $\overline{V}_{\bigdot}$ in the $\infty$-category $\calX_{\Delta_{+}}$ of augmented simplicial objects, and let $\overline{\gamma}: U'_{\bigdot} \rightarrow \overline{W}_{\bigdot}$ be
the induced map. We observe that for each $n \geq 0$, the map 
$\overline{\gamma}(\Delta^n): U_{n+1} \rightarrow \overline{W}_n$ is a pullback of
$U_1 \rightarrow U_0 \times U_0$, and therefore $(n-2)$-truncated (since $U_{\bigdot}$ is
assumed to be $n$-efficient). Since $U_{\bigdot}$ is a groupoid, we conclude that
$\gamma = \overline{\gamma} | \Nerve(\cDelta)^{op}$ is a Cartesian transformation. Invoking Theorem \ref{countercon}, we deduce that $\overline{\gamma}$ is also a Cartesian transformation, so that the diagram
$$ \xymatrix{ U_1 \ar[r] \ar[d] &  U_0 \ar[d] \\
\overline{W}_0 \ar[r] & U_0 \times \overline{W}_{-1} }$$
is Cartesian. Combining this with the Cartesian diagram
$$ \xymatrix{ \overline{W}_0 \ar[r] \ar[d] & \overline{W}_{-1} \ar[d] \\
U_0 \ar[r] & \overline{U}_{-1} }$$
we deduce that $U$ is effective, as desired.
\end{proof}

\subsection{Localic $\infty$-Topoi}\label{nlocalic}

The standard example of an ordinary topos is the category $\Shv(X; \Set)$ of sheaves (of sets) on a topological space $X$. Of course, not every topos is of this form: the category $\Shv(X; \Set)$ 
is generated under colimits by subobjects of its final object (which can be identified with open subsets of $X$). A topos $\calX$ with this property is said to be {\it localic}, and is determined up to equivalence by the locale $\Sub(1_{\calX})$, which we may view as a $0$-topos. The objective of this section is to obtain an $\infty$-categorical analogue of this picture, which will allow us to relate the theory of $n$-topoi to that of $m$-topoi, for all $0 \leq m \leq n \leq \infty$.

\begin{definition}\label{geomorphn}\index{gen}{geometric morphism!of $n$-topoi}
Let $\calX$ and $\calY$ be $n$-topoi, $0 \leq n \leq \infty$
A {\it geometric morphism} from $\calX$ to $\calY$ is a functor $f_{\ast}: \calX \rightarrow \calY$ which admits a left exact left adjoint (which we will typically denote by $f^{\ast}$). 

We let $\Fun_{\ast}(\calX, \calY)$ denote the full subcategory
of the $\infty$-category $\Fun(\calX,\calY)$ spanned by the geometric morphisms, and
$\Geo^{\GeoR}_{n}$ denote the subcategory of $\widehat{\Cat}_{\infty}$ whose objects are $n$-topoi and whose morphisms are geometric morphisms.\index{not}{FunLast@$\Fun_{\ast}(\calX, \calY)$}\index{not}{TopRn@$\Geo^{\GeoR}_{n}$}
\end{definition}

\begin{remark}
In the case where $n=1$, the $\infty$-category of geometric morphisms $\Fun_{\ast}(\calX, \calY)$ between two $1$-topoi is equivalent to (the nerve of) the category of geometric morphisms
between the ordinary topoi $\h{\calX}$ and $\h{\calY}$.
\end{remark}

\begin{remark}
In the case where $n=0$, the $\infty$-category of geometric morphisms $\Fun_{\ast}(\calX, \calY)$
between two $0$-topoi is equivalent to the nerve of the partially ordered set of homomorphisms
from the underlying locale of $\calY$ to the underlying locale of $\calX$. (A {\it homomorphism} between locales is a map of partially ordered sets which preserve finite meets and arbitrary joins.) In the case where $\calX$ and $\calY$ are associated to (sober) topological spaces $X$ and $Y$, this is simply the set of continuous maps from $X$ to $Y$, partially ordered by specialization.
\end{remark}

If $m \leq n$, then the $\infty$-categories 
$\Geo^{\GeoR}_{m}$ and $\Geo^{\GeoR}_{n}$ are related by the following observation:

\begin{proposition}\label{goto1topos}
Let $\calX$ be an $n$-topos, and let $0 \leq m \leq n$. Then the full subcategory
$\tau_{\leq m-1} \calX$ spanned by the $(m-1)$-truncated objects is an $m$-topos.
\end{proposition}

\begin{proof}
If $m = n = \infty$, the result is obvious. Otherwise, it follows immediately from
$(2)$ of Theorem \ref{nchar}.
\end{proof}

\begin{lemma}\label{eoi1}
Let $\calC$ be a small $n$-category which admits finite limits, and let $\calY$ be an $\infty$-topos.
Then the restriction map
$$ \Fun_{\ast}(\calY, \calP(\calC) ) \rightarrow \Fun_{\ast}( \tau_{\leq n-1} \calY, \calP_{\leq n-1}(\calC) )$$
is an equivalence of $\infty$-categories.
\end{lemma}

\begin{proof}
Let $\calM \subseteq \Fun(\calP(\calC), \calY)$ and
$\calM' \subseteq \Fun(\calP_{\leq n-1}(\calC) , \tau_{\leq n-1} \calY)$
denote the full subcategories spanned by left exact, colimit preserving functors. 
In view of Proposition \ref{switcheroo}, it will suffice to prove that the restriction map
$\theta: \calM \rightarrow \calM'$ is an equivalence of $\infty$-categories.

Let $\calM''$ denote full subcategory of $\Fun(\calP(\calC), \tau_{\leq n-1} \calY)$ spanned by
colimit preserving functors whose restriction to $\calP_{\leq n-1}(\calC)$ is left exact. 
Corollary \ref{truncprop} implies that the restriction map $\theta': \calM'' \rightarrow \calM'$
is an equivalence of $\infty$-categories.

Let $j: \calC \rightarrow \calP_{\leq n-1}(\calC) \subseteq \calP(\calC)$ denote the Yoneda embedding. Composition with $j$ yields a commutative diagram
$$ \xymatrix{ \calM \ar[r]^{\theta} \ar[d]^{\psi} & \calM' \ar[d]^{\psi'} \\
\Fun(\calC,\tau_{\leq n-1} \calY) \ar@{=}[r] & \Fun(\calC,\tau_{\leq n-1} \calY). }$$
Theorem \ref{charpresheaf} implies that $\psi$ and $\psi' \circ \theta'$ are fully faithful. Since $\theta'$ is an equivalence of $\infty$-categories, we deduce that $\psi'$ is fully faithful.
Thus $\theta$ is fully faithful; to complete the proof, we must show that $\psi$ and $\psi'$ have the same essential image. Suppose that $f: \calC \rightarrow \tau_{\leq n-1} \calY$ belongs to the essential image of $\psi'$. Without loss of generality, we may suppose that $f$ is a composition
$$ \calC \stackrel{j}{\rightarrow} \calP_{\leq n-1}(\calC) \stackrel{g^{\ast}}{\rightarrow} \tau_{\leq n-1} \calY. $$
As a composition of left exact functors, $f$ is left exact. We may now invoke Proposition \ref{natash} to deduce that $f$ belongs to the essential image of $\psi$.
\end{proof}

\begin{lemma}\label{eoi2}
Let $\calC$ be a small $n$-category which admits finite limits and is equipped with a Grothendieck topology, and let $\calY$ be an $\infty$-topos.
Then the restriction map
$$ \theta: \Fun_{\ast}(\calY, \Shv(\calC) ) \rightarrow \Fun_{\ast}( \tau_{\leq n-1} \calY, \Shv_{\leq n-1}(\calC) )$$
is an equivalence of $\infty$-categories.
\end{lemma}

\begin{proof}
We have a commutative diagram
$$ \xymatrix{ \Fun_{\ast}( \calY, \Shv(\calC) ) \ar[r]^-{\theta} \ar[d] & \Fun_{\ast}( \tau_{\leq n-1} \calY, \Shv_{\leq n-1}(\calC) ) \ar[d] \\
\Fun_{\ast}(\calY, \calP(\calC) ) \ar[r]^-{\theta'} & \Fun_{\ast}( \tau_{\leq n-1} \calY, \calP_{\leq n-1}(\calC))}$$
where the vertical arrows are inclusions of full subcategories, and $\theta'$ is an equivalence of $\infty$-categories (Lemma \ref{eoi1}). To complete the proof, it will suffice to show that if
$f_{\ast}: \calY \rightarrow \calP(\calC)$ is a geometric morphism such that $f_{\ast} | \tau_{\leq n-1} \calY$ factors through $\Shv_{\leq n-1}(\calC)$, then $f_{\ast}$ factors through $\Shv(\calC)$. 

Let $f^{\ast}$ be a left adjoint to $f_{\ast}$, and let $\overline{S}$ denote the collection of all
morphisms in $\calP(\calC)$ which localize to equivalences in $\Shv(\calC)$. We must show that $f^{\ast}\overline{S} $ consists of equivalences in $\calY$. Let $S \subseteq \overline{S}$ be the collection of monomorphisms which belong to $\overline{S}$. Since $\Shv(\calC)$ is a topological localization of $\calP(\calC)$, it will suffice to show that $f^{\ast} S$ consists of equivalences in $\calY$. Let $g: X \rightarrow Y$ belong to $S$. Since $\calP(\calC)$ is generated under colimits by the essential image of the Yoneda embedding, we can write $Y$ as a colimit of a diagram
$K \rightarrow \calP_{\leq n-1}(\calC)$. Since colimits in $\calP(\calC)$ are universal, we obtain a corresponding expression of $g$ as a colimit of morphisms $\{ g_{\alpha}: X_{\alpha} \rightarrow Y_{\alpha} \}$ which are pullbacks of $g$, where $Y_{\alpha} \in \calP_{\leq n-1}(\calC)$. In this case, $g_{\alpha}$ is again a monomorphism, so that $X_{\alpha}$ is also $(n-1)$-truncated.
Since $f^{\ast}$ commutes with colimits, it will suffice to show that each $f^{\ast}(g_{\alpha})$ is
an equivalence. But this follows immediately from our assumption that $f_{\ast} | \tau_{\leq n-1} \calY$ factors through $\Shv_{\leq n-1}(\calY)$.
\end{proof}

\begin{proposition}\label{swmad}
Let $0 \leq m \leq n \leq \infty$, and let $\calY$ be an $m$-topos. There exists an $n$-topos
$\calX$ and a categorical equivalence 
$f_{\ast}: \tau_{\leq m-1} \calX \rightarrow \calY$ with the following universal property: for any $n$-topos $\calZ$, composition with $f_{\ast}$ induces
an equivalence of $\infty$-categories
$$ \theta: \Fun_{\ast}( \calZ, \calX ) \rightarrow \Fun_{\ast}( \tau_{\leq m-1} \calZ, \calY).$$
\end{proposition}

\begin{proof}
If $m = \infty$, then also $n= \infty$ and we may take $\calX = \calY$. Otherwise, we may apply
Theorem \ref{nchar} to reduce to the case where $\calY = \Shv_{\leq m-1}(\calC)$, where $\calC$ is a small $m$-category which admits finite limits and is equipped with a Grothendieck topology.
In this case, we
let $\calX = \Shv_{\leq n-1}(\calC)$ and define $f_{\ast}$ to be the identity. Let $\calZ$ be an arbitrary $n$-topos. According to Theorem \ref{nchar}, we may assume without loss of generality that $\calZ = \tau_{\leq n-1} \calZ'$, where $\calZ'$ is an $\infty$-topos. We have a commutative diagram
$$ \xymatrix{ & \Fun_{\ast}( \calZ, \calX) \ar[dr]^{\theta} & \\
\Fun_{\ast}(\calZ', \Shv(\calC)) \ar[ur]^{\theta'} \ar[rr]^{\theta''} & & \Fun_{\ast}(\tau_{\leq m-1} \calZ, \calY). }$$
Lemma \ref{eoi2} implies that $\theta'$ and $\theta''$ are equivalences of $\infty$-categories,
so that $\theta$ is also an equivalence of $\infty$-categories.
\end{proof}

\begin{definition}\index{gen}{localic!$n$-topos}\label{stuffera}
Let $0 \leq m \leq n \leq \infty$, and let $\calX$ be an $n$-topos. We will say that
$\calX$ is {\it $m$-localic} if, for any $n$-topos $\calY$, the natural map
$$ \Fun_{\ast}( \calY, \calX) \rightarrow \Fun_{\ast}( \tau_{\leq m-1} \calY, \tau_{\leq m-1} \calX )$$
is an equivalence of $\infty$-categories.
\end{definition}

According to Proposition \ref{swmad}, every
$m$-topos $\calX$ is equivalent to the subcategory of $(m-1)$-truncated objects in an $m$-localic $n$-topos $\calX'$, and $\calX'$ is determined up to equivalence. More precisely, the truncation functor
$$ \Geo^{\GeoR}_{n}  \stackrel{\tau_{\leq m-1}}{\rightarrow} \Geo^{\GeoR}_{m}$$
induces an equivalence $ \calC \rightarrow \Geo^{\GeoR}_{m}$, where
$\calC \subseteq \Geo^{\GeoR}_{n}$ is the full subcategory spanned by the $m$-localic $n$-topoi.
In other words, we may view the $\infty$-category of
$m$-topoi as a {\it localization} of the $\infty$-category of $n$-topoi. In particular, the theory of
$m$-topoi for $m < \infty$ can be regarded as a special case of the theory of $\infty$-topoi. For this reason, we will focus our attention on the case $n = \infty$ for most of the remainder of this book.

\begin{proposition}\label{useiron}
Let $\calX$ be an $n$-localic $\infty$-topos. Then any topological localization of $\calX$ is also $n$-localic.
\end{proposition}

\begin{proof}
The proof of Proposition \ref{swmad} shows that $\calX$ is $n$-localic if and only if there
exists a small $n$-category $\calC$ which admits finite limits, a Grothendieck topology on $\calC$, and an equivalence $\calX \rightarrow \Shv(\calC)$. In other words, $\calX$ is $n$-localic if and only if it is equivalent to a topological localization of $\calP(\calC)$, where $\calC$ is a small $n$-category which admits finite limits. It is clear that any topological localization of $\calX$ has the same property.
\end{proof}


Let $\calX$ be an $\infty$-topos. One should think of the $\infty$-categories $\tau_{\leq n-1} \calX$ as  ``Postnikov sections'' of $\calX$. The classical $1$-truncation $\tau_{\leq 1} X$ of a
homotopy type $X$ remembers only the fundamental groupoid of $X$.
It therefore knows all about local systems of sets on $X$, but
nothing about fibrations over $X$ with non-discrete fibers. The
relationship between $\calX$ and $\tau_{\leq 0} \calX$ is analogous:
$\tau_{\leq 0} \calX$ knows about the sheaves of sets on $\calX$, but
has forgotten about sheaves with nondiscrete stalks.

\begin{remark}
In view of the above discussion, the notation $\tau_{\leq 0} \calX$ is unfortunate
because the analogous notation for the $1$-truncation of a
homotopy type $X$ is $\tau_{\leq 1} X$. We caution the reader not to
regard $\tau_{\leq 0} \calX$ as the result of applying an operation
$\tau_{\leq 0}$ to $\calX$; it instead denotes the essential image of the truncation
functor  $\tau_{\leq 0}: \calX \rightarrow \calX$.
\end{remark}

\section{Homotopy Theory in an $\infty$-Topos}\label{chap6sec5}

\setcounter{theorem}{0}

In classical homotopy theory, the most important invariants of a (pointed) space
$X$ are its homotopy groups $\pi_{i}(X,x)$. Our first objective in this section is to define analogous invariants in the case where $X$ is an object of an arbitrary $\infty$-topos $\calX$. In this setting, the homotopy groups are not ordinary groups but are instead {\em sheaves} of groups on the underlying topos $\Disc(\calX)$. In \S \ref{homotopysheaves}, we will study these homotopy groups
and the closely related theory of {\it $n$-connectivity}. The main theme is that the internal homotopy theory of a general $\infty$-topos $\calX$ behaves much like the classical case $\calX = \SSet$.

One important classical fact which does {\em not} hold in general for an $\infty$-topos is Whitehead's theorem. If $f: X \rightarrow Y$ is a map of CW-complexes, then $f$ is a homotopy equivalence if and only if $f$ induces bijective maps $\pi_{i}(X,x) \rightarrow \pi_i(Y, f(x))$ for
any $i \geq 0$ and any base point $x \in X$. If $f: X \rightarrow Y$ is a map in an arbitrary $\infty$-topos $\calX$ satisfying an analogous condition on (sheaves of) homotopy groups, then we say that
$f$ is {\it $\infty$-connective}. We will say that an $\infty$-topos $\calX$ is {\it hypercomplete}
if every $\infty$-connective morphism in $\calX$ is an equivalence. Whitehead's theorem may be interpreted as saying that the $\infty$-topos $\SSet$ is hypercomplete. An arbitrary $\infty$-topos $\calX$ need not be hypercomplete. We will survey the situation in \S \ref{hyperstacks}, where we also give some reformulations of the notion of hypercompleteness and show that every topos $\calX$ has a {\it hypercompletion} $\calX^{\hyp}$. In \S \ref{hcovh}, we will show that an $\infty$-topos $\calX$ is hypercomplete if and only if $\calX$ satisfies a descent condition with respect to hypercoverings (other versions of this result can be found in \cite{hollander} and \cite{toen}).

\begin{remark}
The Brown-Joyal-Jardine theory of simplicial (pre)sheaves on a topological space $X$ is a model for the hypercomplete $\infty$-topos $\Shv(X)^{\hyp}$. In many respects, the $\infty$-topos $\Shv(X)$ of sheaves of spaces on $X$ is better behaved {\em before} hypercompletion. We will outline some of the advantages of $\Shv(X)$ in \S \ref{versus} and in \S \ref{chap7}.
\end{remark}

\subsection{Homotopy Groups}\label{homotopysheaves}

Let $\calX$ be an $\infty$-topos, and let $X$ be an object of
$\calX$. We will refer to a discrete object of $\calX_{/X}$ as a {\it sheaf of sets on $X$}.
Since $\calX$ is presentable, it is automatically {\em cotensored} over
spaces, as explained in Remark \ref{coten}. Consequently, for any object $X$ of
$\calX$ and any simplicial set $K$, there exists an object $X^K$ of $\calX$ equipped with natural isomorphisms
$$ \bHom_{\calX}( Y, X^K ) \rightarrow \bHom_{\calH}( K, \bHom_{\calX}(Y,X) )$$
in the homotopy category $\calH$ of spaces.

\begin{definition}\index{not}{pinX@$\pi_n(X)$}\index{gen}{homotopy groups!in an $\infty$-topos}
Let $S^n  = \bd \Delta^{n+1} \in \calH$ denote the (simplicial) $n$-sphere, and fix a base point $\ast \in S^n$. Then evaluation at $\ast$ induces a morphism $s: X^{S^n} \rightarrow X$ in $\calX$.
We may regard $s$ as an object of $\calX_{/X}$, and we define $\pi_n(X) = \tau_{\leq 0} s \in \calX_{/X}$ to be the associated discrete object of $\calX_{/X}$.
\end{definition}

 We will generally identify $\pi_n(X)$ with its image in the underlying topos $\Disc(\calX_{/X})$ (where it is well-defined up to canonical isomorphism). The constant map $S^n \rightarrow \ast$ induces
a map $X \rightarrow X^{S^n}$, which determines a base point of $\pi_n(X)$. 

Suppose that $K$ and $K'$ are pointed simplicial sets, and let $K \vee K'$ denote the coproduct $K \amalg_{\ast} K'$. There is a pullback diagram
$$ \xymatrix{ & X^{K \vee K'} \ar[dr] \ar[dl] & \\
X^K \ar[dr] & & X^{K'} \ar[dl] \\
& X & }$$
in $\calX$, so that $X^{K \vee K'}$ may be identified with a product
of $X^K$ and $X^{K'}$ in the $\infty$-topos $\calX_{/X}$. We now make the following
general observation:

\begin{lemma}\label{slurpy}
Let $\calX$ be an $\infty$-topos. The truncation functor $\tau_{\leq n}: \calX \rightarrow \calX$
preserves finite products.
\end{lemma}

\begin{proof}
We must show that for any finite collection of objects $\{ X_{\alpha} \}_{\alpha \in A}$
having product $X$, the induced map
$$ \tau_{\leq n} X \rightarrow \prod_{\alpha \in A} \tau_{\leq n} X_{\alpha} $$
is an equivalence. If $\calX$ is the $\infty$-category of spaces, then this follows from Whitehead's theorem: simply compute homotopy groups (and sets) on both sides. If $\calX = \calP(\calC)$, then to prove that a map in $\calX$ is an equivalence, it suffices to show that it remains an equivalence after evaluation at any object $C \in \calC$; thus we may reduce to the case where $\calX = \SSet$
considered above. In the general case, $\calX$ is equivalent to the essential image of a left-exact localization functor $L: \calP(\calC) \rightarrow \calP(\calC)$ for some small $\infty$-category $\calC$. Without loss of generality, we may identify $\calX$ with a full subcategory of $\calP(\calC)$. Then
$\calX \subseteq \calX' = \calP(\calC)$ is stable under limits, so that $X$ may be identified with
a product of the family $\{ X_{\alpha} \}_{\alpha \in A}$ in $\calX'$. It follows from the case treated above that the natural map
$$ \tau_{\leq n}^{\calX'} X \rightarrow \prod_{\alpha \in A} \tau_{\leq n}^{\calX'} X_{\alpha}$$
is an equivalence. But Proposition \ref{compattrunc} implies that $L \circ \tau_{\leq n}^{\calX'} | \calX$ is an $n$-truncation functor for $\calX$. The desired result now follows by applying the functor $L$ to both sides of the above equivalence, and invoking the assumption that $L$ is left-exact (here we must require the finiteness of $A$).
\end{proof}

It follows from Lemma \ref{slurpy} that there is a canonical isomorphism
$$\tau^{\calX_{/X}}_{\leq 0} (X^{K \vee K'}) \simeq \tau^{\calX_{/X}}_{\leq 0}(X^K) \times \tau^{\calX_{/X}}_{\leq 0}(X^{K'})$$ in the topos $\Disc(\calX_{/X})$. 
In particular, for $n > 0$, the usual comultiplication $S^n \rightarrow S^n \vee S^n$
(a well-defined map in the homotopy category $\calH$)
induces a multiplication map $\pi_n(X) \times \pi_n(X) \rightarrow \pi_n(X)$. As in ordinary homotopy theory, we conclude that $\pi_n(X)$ is a group object of $\Disc(\calX_{/X})$ for $n > 0$, which is commutative for $n > 1$.

In order to work effectively with homotopy sets, it is convenient
to define the homotopy sets $\pi_n(f)$ of a morphism $f: X
\rightarrow Y$ to be the homotopy sets of $f$ considered as an object
of the $\infty$-topos $\calX_{/Y}$. In view of the equivalences
$\calX_{/f} \rightarrow \calX_{/X}$, we may identify $\pi_n(f)$ with an object
of $\Disc(\calX_{/X})$, which is again a sheaf of groups if
$n \geq 1$, and abelian groups if $n\geq 2$.
The intuition is that the stalk of these sheaves at a point $p$ of $X$ is the $n$th homotopy group of the homotopy fiber of $f$, taken with respect to the base point $p$.

\begin{remark}\label{recgroup}
It is useful to have the following recursive definition for homotopy groups.
Let $f: X \rightarrow Y$ be a morphism in an $\infty$-topos $\calX$. 
Regarding $f$ as an object of the topos
$\calX_{/Y}$, we may take its $0$th truncation
$\tau_{\leq 0}^{\calX_{/Y}} f$. This is a discrete object of $\calX_{/Y}$, and by definition
we have $\pi_0(f) \simeq f^{\ast} \tau_{\leq 0}^{\calX_{/Y}}(X) \simeq
X \times_{Y} \tau_{\leq 0}^{\calX_{/Y}}(f)$. The natural map $X
\rightarrow \tau_{\leq 0}^{\calX_{/Y}}(f)$ gives a global section of
$\pi_0(f)$. Note that in this case, $\pi_0(f)$ is the pullback of
a discrete object of $\calX_{/Y}$: this is because the definition of $\pi_0$ does not
require a base point. 

If $n >0$, then we have a natural isomorphism $\pi_{n}(f) \simeq
\pi_{n-1}(\delta)$ in $\Disc(\calX_{/X})$, where $\delta: X \rightarrow X \times_{Y} X$ is the associated diagonal map. \end{remark}

\begin{remark}\label{emmy}
Let $f: \calX \rightarrow \calY$ be a geometric morphism of
$\infty$-topoi, and let $g: Y \rightarrow Y'$ be a morphism in
$\calY$. Then there is a canonical isomorphism $f^{\ast}( \pi_n(g) ) \simeq \pi_n(f^{\ast}(g))$
in $\Disc(\calX_{/f^{\ast} Y})$. This follows immediately from Proposition \ref{compattrunc}.
\end{remark}

\begin{remark}\label{sequence}\index{gen}{long exact sequence of homotopy groups}
Given a pair of composable morphisms $X \stackrel{f}{\rightarrow}
Y \stackrel{g}{\rightarrow} Z$, there is an associated sequence of pointed objects
$$\ldots \rightarrow f^{\ast} \pi_{n+1}(g) \stackrel{\delta_{n+1}}{\rightarrow} \pi_n(f) \rightarrow \pi_n(g \circ f) \rightarrow f^{\ast}
\pi_n(g) \stackrel{\delta_{n}}{\rightarrow} \pi_{n-1}(f) \rightarrow \ldots$$
in the ordinary topos $\Disc(\calX_{/X})$, with the
usual exactness properties. To construct the boundary map $\delta_n$, we observe
that the $n$-sphere $S^n$ can be written as a (homotopy) pushout
$D^{-} \amalg_{ S^{n-1} } D^{+}$ of two hemispheres along the equator. By construction,
$f^{\ast} \pi_n(g)$ can be identified with the $0$-truncation of
$$X \times_{Y} Y^{S^n} \times_{ Z^{S^n} } Z \simeq
X^{D^{-}} \times_{ Y^{D^{-}} } Y^{S^n} \times_{ Z^{S^{n}} } Z,$$
which maps by restriction to
$$ X^{S^{n-1}} \times_{ Y^{S^{n-1}} } Y^{D^{+}} \simeq X^{ S^{n-1} } \times_{Y^{S^{n-1}} } Y.$$
We now observe that the $0$-truncation of the latter object is naturally isomorphic to
$\pi_{n-1}(f) \in \Disc(\calX_{/X})$.

To prove the exactness of the above sequence in an $\infty$-topos $\calX$, we
first choose an accessible left exact localization $L: \calP(\calC) \rightarrow \calX$.
Without loss of generality, we may suppose that the diagram
$X \stackrel{f}{\rightarrow} Y \stackrel{g}{\rightarrow} Z$ is the image under $L$ of
a diagram in $\calP(\calC)$. Using Remark \ref{emmy}, we conclude that the sequence constructed above is equivalent to the image under $L$ of an analogous sequence in the $\infty$-topos
$\calP(\calC)$. Since $L$ is left exact, it will suffice to prove that this second sequence is exact; in other words, we may reduce to the case $\calX = \calP(\calC)$. Working componentwise, we
can reduce further to the case where $\calX = \SSet$. The desired result now follows from classical homotopy theory. (Special care should be taken regarding the exactness of the above sequence
at $\pi_0(f)$: really this should be interpreted in terms of an action of the group $f^{\ast} \pi_1(g)$ on
$\pi_0(f)$. We leave the details of the construction of this action to the reader. )
\end{remark}

\begin{remark}
If $\calX = \SSet$, and $\eta: \ast \rightarrow X$ is a pointed
space, then $\eta^{\ast} \pi_{n}(X)$ can be identified with the $n$th homotopy group
of $X$ with base point $\eta$.
\end{remark}

We now study the implications of the vanishing of homotopy groups.

\begin{proposition}\label{ditz}\index{gen}{truncated!and homotopy groups}
Let $f: X \rightarrow Y$ be an $n$-truncated morphism in an $\infty$-topos $\calX$. Then
$\pi_{k}(f) \simeq \ast$ for all $k > n$. If $n \geq 0$ and $\pi_{n}(f)
\simeq \ast$, then $f$ is $(n-1)$-truncated.
\end{proposition}

\begin{proof}
The proof goes by induction on $n$. If $n = -2$, then $f$ is an
equivalence and there is nothing to prove. Otherwise, the diagonal map
$\delta: X \rightarrow X \times_Y X$ is $(n-1)$-truncated (Lemma \ref{trunc}). The
inductive hypothesis and Remark \ref{recgroup}
allow us to deduce that $\pi_{k}(f) \simeq \pi_{k-1}(\delta) \simeq \ast$ whenever $k > n$ and $k > 0$. Similarly, if $n \geq 1$ and $\pi_{n}(f) \simeq \pi_{n-1}(\delta) \simeq \ast$, then $\delta$ is
$(n-2)$-truncated by the inductive hypothesis, so that $f$ is
$(n-1)$-truncated (Lemma \ref{trunc}). 

The case of small $k$ and $n$ requires special attention: we must
show that if $f$ is $0$-truncated, then $f$ is $(-1)$-truncated if
and only if $\pi_0(f) \simeq \ast$. Because $f$ is $0$-truncated, we have
an equivalence $\tau_{\leq 0}^{\calX_{/Y}}(f) \simeq f$, so that $\pi_0(f) \simeq X
\times_Y X$. To say $\pi_0(f) \simeq \ast$ is to assert that the diagonal map
$\delta: X \rightarrow X \times_{Y} X$ is an equivalence, which is equivalent
to the assertion that $f$ is $(-1)$-truncated (Lemma \ref{trunc}). 
\end{proof}

\begin{remark}
The Proposition \ref{ditz} implies that if $f$ is $n$-truncated
for {\em some} $n \gg 0$, then we can test whether or not $f$ is
$m$-truncated for any particular value of $m$ by computing the
homotopy groups of $f$. In contrast to the classical situation, it
is not possible to drop the assumption that $f$ is
$n$-truncated for $n \gg 0$. 
\end{remark}

\begin{lemma}\label{truncatepin}
Let $X$ be an object in an $\infty$-topos $\calX$, and let
$p: X \rightarrow Y$ be an $n$-truncation of $X$. Then $p$ induces isomorphisms
$\pi_{k}(X) \simeq p^{\ast} \pi_{k}(Y)$ for all $k
\leq n$.
\end{lemma}

\begin{proof}
Let $\phi: \calX \rightarrow \calY$ be a geometric morphism such
that $\phi_{\ast}$ is fully faithful. By Proposition
\ref{compattrunc} and Remark \ref{emmy}, it will suffice to prove
the lemma in the case where $\calX = \calY$. We may therefore assume that $\calY$ is an $\infty$-category of presheaves. In this case, homotopy groups and truncations are
computed pointwise. Thus we may reduce to the case $\calX =
\SSet$, where the conclusion follows from classical homotopy theory. 
\end{proof}

\begin{definition}\label{stooog}\index{gen}{connective!$n$-connective object}\index{gen}{connective!$n$-connective morphism}
Let $f: X \rightarrow Y$ be a morphism in an $\infty$-topos
$\calX$, and let $0 \leq n \leq \infty$. We will say that $f$ is
{\it $n$-connective} if it is an effective epimorphism and $\pi_k(f) = \ast$ for
$0 \leq k < n$. We shall say that the object $X$ is {\it
$n$-connective} if $f: X \rightarrow 1_{\calX}$ is $n$-connective, where $1_{\calX}$ denotes the final object of $\calX$. By convention, we will say that every morphism $f$ in $\calX$ is
$(-1)$-connective.
\end{definition}

\begin{definition}\index{gen}{connected!object of an $\infty$-topos}
Let $X$ be an object of an $\infty$-topos $\calX$. We will say that $X$ is {\it connected} if
it is $1$-connective: that is, if the truncation $\tau_{\leq 0} X$ is a final object in $\calX$.
\end{definition}


\begin{proposition}\label{goober}
Let $X$ be an object in an $\infty$-topos $\calX$ and let $n \geq -1$. Then $X$ is
$n$-connective if and only if $\tau_{\leq n-1} X$ is a final object of $\calX$.
\end{proposition}

\begin{proof}
The case $n=-1$ is trivial. The proof in general goes by induction on $n \geq 0$. If $n=0$, then the
conclusion follows from Proposition \ref{slurpme}. 
Suppose $n > 0$. Let $p: X \rightarrow \tau_{n-1} X$ be an $(n-1)$-truncation of
$X$. If $\tau_{\leq n-1} X$ is a final object of $\calX$, then 
$$\pi_{k} X \simeq p^{\ast} \pi_k(\tau_{n-1} X) \simeq \ast$$ 
for $k < n$ by Lemma \ref{truncatepin}. Since 
the map $p: X \rightarrow \tau_{\leq n-1} X \simeq 1_{\calX}$ is an effective epimorphism
(Proposition \ref{pi00detects}), it follows that $X$ is $n$-connective.

Conversely, suppose that $X$ is $n$-connective. Then $p^{\ast}
\pi_{n-1}(\tau_{\leq n-1} X) \simeq \ast$. Since $p$ is an effective epimorphism, Lemma \ref{hint0} implies that
$\pi_{n-1}(\tau_{\leq n-1} X) = \ast$. Using Proposition \ref{ditz}, we conclude that $\tau_{\leq n-1} X$ is
$(n-2)$-truncated, so that $\tau_{\leq n-1} X \simeq \tau_{\leq n-2} X$.
Repeating this argument, we reduce to the case where $n=0$ which was handled above.
\end{proof}

\begin{corollary}\label{togoto}
The class of $n$-connective objects of an $\infty$-topos $\calX$ is stable under finite products.
\end{corollary}

\begin{proof}
Combine Proposition \ref{goober} with Lemma \ref{slurpy}. 
\end{proof}

Let $\calX$ be an $\infty$-topos and $X$ an object of $\calX$. 
Since $\bHom_{\calX}(X,Y) \simeq \bHom_{\calX}(\tau_{\leq n} X, Y)$ whenever
$Y$ is $n$-truncated, we deduce that $X$ is $(n+1)$-connective if and
only if the natural map $\bHom_{\calX}(1_{\calX},Y) \rightarrow \bHom_{\calX}(X,Y)$ is an
equivalence for all $n$-truncated $Y$. From this, we can
immediately deduce the following relative version of Proposition
\ref{goober}:

\begin{corollary}\label{goober2}
Let $f: X \rightarrow X'$ be a morphism in an $\infty$-topos
$\calX$. Then $f$ is $(n+1)$-connective if and only if composition with $f$ induces a homotopy equivalence $$\bHom_{\calX_{/X'}}(\id_{X'},Y)
\rightarrow \bHom_{\calX_{/X'}}(f,Y)$$ 
for every $n$-truncated object $Y \in \calX_{/X'}$.
\end{corollary}

\begin{remark}\label{nconn}
Let $L: \calX \rightarrow \calY$ be a left exact localization of $\infty$-topoi, and let
$f: Y \rightarrow Y'$ be an $n$-connective morphism in $\calY$. Then $f$ is equivalent (in $\Fun(\Delta^1,\calY)$) to $Lf_0$, where $f_0$ is an $n$-connective morphism in $\calX$.
To see this, we choose a (fully faithful) right adjoint $G$ to $L$, and a factorization
$$ \xymatrix{ & X \ar[dr]^{f''} & \\
G(Y) \ar[ur]^{f'} \ar[rr]^{ G(f_0)} & & G(Y') }$$
where $f'$ is $n$-connective and $f''$ is $(n-1)$-truncated. Then
$Lf'' \circ Lf'$ is equivalent to $f$, and is therefore $n$-connective. It follows that $Lf''$ is
an equivalence, so that $Lf'$ is equivalent to $f$.
\end{remark}

We conclude by noting the following stability properties of the
class of $n$-connective morphisms:

\begin{proposition}\label{inftychange}
Let $\calX$ be an $\infty$-topos.
\begin{itemize}
\item[$(1)$] Let $f: X \rightarrow Y$ be a morphism in $\calX$. If $f$
is $n$-connective, then it is $m$-connective for all $m \leq n$. Conversely, if
$f$ is $n$-connective for all $n < \infty$, then $f$ is $\infty$-connective.

\item[$(2)$] Any equivalence in $\calX$ is $\infty$-connective.

\item[$(3)$] Let $f,g: X \rightarrow Y$ be a homotopic morphisms in $\calX$. Then
$f$ is $n$-connective if and only if $g$ is $n$-connective.

\item[$(4)$] Let $\pi^{\ast}: \calX \rightarrow \calY$ be left adjoint to a geometric morphism
from $\pi_{\ast}: \calY \rightarrow \calX$, and let $f: X \rightarrow X'$ be an $n$-connective
morphism in $\calX$. Then $\pi^{\ast} f$ is an $n$-connective morphism in $\calY$.

\item[$(5)$] Suppose given a diagram
$$ \xymatrix{ & Y \ar[dr]^{g} & \\
X \ar[ur]^{f} \ar[rr]^{h} & & Z }$$
in $\calX$, where $f$ is $n$-connective. Then $g$ is $n$-connective if and only if
$h$ is $n$-connective. 

\item[$(6)$] Suppose given a pullback diagram
$$ \xymatrix{ X' \ar[d]^{f'} \ar[r]^{q'} & X \ar[d]^{f} \\
Y' \ar[r]^{q} & Y }$$
in $\calX$. If $f$ is $n$-connective, then so is $f'$. The converse holds if $q$ is an effective epimorphism.
\end{itemize}
\end{proposition}

\begin{proof}
The first three assertions are obvious. Claim $(4)$ follows from Propositions \ref{goober} and \ref{compattrunc}. To prove $(5)$, we first observe that Corollary \ref{composite} implies that
$g$ is an effective epimorphism if and only if $h$ is an effective epimorphism. According to
Remark \ref{sequence}, we have a long exact sequence
$$\ldots \rightarrow f^{\ast} \pi_{i+1}(g) {\rightarrow} \pi_i(f) \rightarrow \pi_i(h) \rightarrow f^{\ast}
\pi_i(g) \rightarrow \pi_{i-1}(f) \rightarrow \ldots$$ 
of pointed objects in the topos $\Disc(\calX_{/X})$. It is then clear that if $f$ and $g$
are $n$-connective, then so is $h$. Conversely, if $f$ and $h$ are $n$-connective, then
$f^{\ast} \pi_i(g) \simeq \ast$ for $i \leq n$. Since $f$ is an effective epimorphism, Lemma \ref{hint0} implies that $\pi_i(g) \simeq \ast$ for $i \leq n$, so that $g$ is also $n$-connective.

The first assertion of $(6)$ follows from $(4)$, since a pullback functor
$q^{\ast}: \calX_{/Y} \rightarrow \calX_{/Y'}$ is left adjoint to a geometric morphism.
For the converse, let us suppose that $q$ is an effective epimorphism and that $f'$ is $n$-connective. According to Lemma \ref{hintdescent1}, the maps $f$ and $q'$ are effective epimorphisms. Applying Remark \ref{emmy}, we conclude that there are canonical isomorphisms
${q'}^{\ast} \pi_k(f) \simeq \pi_k(f')$ in the topos $\Disc( \calX_{/X'})$, so that
${q'}^{\ast} \pi_k(f) \simeq \ast$ for $k < n$. Applying Lemma \ref{hint0}, we conclude
that $\pi_k(f) \simeq \ast$ for  $k < n$, so that $f$ is $n$-connective as desired.
\end{proof}

\begin{corollary}\label{pusherr}
Let 
$$\xymatrix{ X' \ar[r]^{g} \ar[d]^{f'} & X \ar[d]^{f} \\
Y' \ar[r] & Y }$$
be a pushout diagram in an $\infty$-topos $\calX$. Suppose that $f'$ is
$n$-connective. Then $f$ is $n$-connective.
\end{corollary}

\begin{proof}
Choose an accessible, left exact localization functor $L: \calP(\calC) \rightarrow \calX$.
Using Remark \ref{nconn}, we can assume without loss of generality that
$f' = Lf'_0$, where $f'_0: X'_0 \rightarrow Y'_0$ is a morphism in $\calP(\calC)$.
Similarly, we may assume $g = Lg_0$, for some morphism $g_0: X'_0 \rightarrow X_0$.
Form a pushout diagram
$$\xymatrix{ X'_0 \ar[r]^{g_0} \ar[d]^{f'_0} & X_0 \ar[d]^{f_0} \\
Y'_0 \ar[r] & Y_0 }$$
in $\calP(\calC)$. Then the original diagram is equivalent to the image (under $L$) of the diagram above. In view of Proposition \ref{inftychange}, it will suffice to show that $f_0$ is $n$-connective.
Using Propositions \ref{goober} and \ref{compattrunc}, we see that $f_0$ is $n$-connective if and only if its image under the evaluation map $\calP(\calC) \rightarrow \SSet$ associated to any object $C \in \calC$ is $n$-connective. In other words, we can reduce to the case where $\calX = \SSet$, and the result now follows from classical homotopy theory.
\end{proof}

We conclude by establishing a few results which will be needed in \S \ref{dimension}:

\begin{proposition}\label{trowler}
Let $f: X \rightarrow Y$ be a morphism in an $\infty$-topos $\calX$, 
$\delta: X \rightarrow X \times_{Y} X$ the associated diagonal morphism, and $n \geq 0$.
The following conditions are equivalent:
\begin{itemize}
\item[$(1)$] The morphism $f$ is $n$-connective.

\item[$(2)$] The diagonal map $\delta: X \rightarrow X \times_{Y} X$ is $(n-1)$-connective, and
$f$ is an effective epimorphism.
\end{itemize}
\end{proposition}

\begin{proof}
Immediate from Definition \ref{stooog} and Remark \ref{recgroup}.
\end{proof}

\begin{proposition}\label{conslice}
Let $\calX$ be an $\infty$-topos containing an object $X$, and let $\sigma: \Delta^2 \rightarrow \calX$ be a $2$-simplex corresponding to a diagram
$$ \xymatrix{ Y \ar[dr] \ar[rr]^{f} & & Z \ar[dl]^{g} \\
& X. & }$$
Then $f$ is an $n$-connective morphism in $\calX$ if and only if $\sigma$ is an
$n$-connective morphism in $\calX_{/X}$.
\end{proposition}

\begin{proof}
We observe that $\calX_{/g} \rightarrow \calX_{/Z}$ is a
trivial fibration, so that an object of $\calX_{/g}$ is $n$-connective if and only if its
image in $\calX_{/Z}$ is $n$-connective. 
\end{proof}

\begin{proposition}\label{sectcon}
Let $f: X \rightarrow Y$ be a morphism in an $\infty$-topos $\calX$, let
$s:Y \rightarrow X$ be a section of $f$ (so that $f \circ s$ is homotopic to $\id_{Y}$), and
let $n \geq 0$. Then $f$ is $n$-connective if and only if $s$ is $(n-1)$-connective.
\end{proposition}

\begin{proof}
We have a $2$-simplex $\sigma: \Delta^2 \rightarrow \calX$ which we may depict as follows:
$$ \xymatrix{ & X \ar[dr]^{f} & \\
Y \ar[ur]^{s} \ar[rr]^{\id_Y} & & Y. }$$
Corollary \ref{composite} implies that $f$ is an effective epimorphism; this completes the proof in the case $n = 0$. Suppose that $n > 0$, and that $s$ is $(n-1)$-connective.
In particular, $s$ is an effective epimorphism. The long exact sequence of Remark \ref{sequence} gives an isomorphism $\pi_{i}(s) \simeq s^{\ast} \pi_{i+1}(f)$, so that $s^{\ast} \pi_k(f)$ vanishes for
$1 \leq k < n$. Applying Lemma \ref{hint0}, we conclude that $\pi_k(f) \simeq \ast$ for $1 \leq k < n$. Moreover, since $s$ is an effective epimorphism it induces an effective epimorphism
$\pi_0( \id_{Y} ) \rightarrow \pi_0(f)$ in the ordinary topos $\Disc( \calX_{/Y} )$, so that
$\pi_0(f) \simeq \ast$ as well. This proves that $f$ is $n$-connective.

Conversely, if $f$ is $n$-connective, then $\pi_i(s) \simeq \ast$ for
$i < n-1$; the only nontrivial point is to verify that $s$ is an effective epimorphism. According to Proposition \ref{conslice}, it will suffice to prove that $\sigma$ is an effective epimorphism when viewed as a morphism in $\calX_{/Y}$. Using Proposition \ref{pi00detects}, we may reduce to proving that $\sigma' = \tau_{\leq 0}^{\calX_{/Y}}(\sigma)$ is an equivalence in $\calX_{/Y}$.
This is clear, since the source and target of $\sigma'$ are both final objects of $\calX_{/Y}$ (in virtue of our assumption that $f$ is $1$-connective).
\end{proof}

\subsection{$\infty$-Connectedness}\label{hyperstacks}

Let $\calC$ be an ordinary category equipped with a Grothendieck topology, and let
$\bfA = \Set_{\Delta}^{\calC^{op}}$ be the category of simplicial presheaves on $\calC$.

\begin{proposition}[Jardine \cite{jardine}]\label{jardinesardine}\index{gen}{model category!local}
There exists a left proper, combinatorial, simplicial model structure on the category $\bfA$, which admits the following description:

\begin{itemize}
\item[$(C)$] A map $f: F_{\bigdot} \rightarrow G_{\bigdot}$ of simplicial presheaves on $\calC$
is a {\it local cofibration} if it is a injective cofibration: that is, if and only if the induced map
$F_{\bigdot}(C) \rightarrow G_{\bigdot}(C)$ is a cofibration of simplicial sets for each object $C \in \calC$.

\item[$(W)$] A map $f: F_{\bigdot} \rightarrow G_{\bigdot}$ of simplicial presheaves on $\calC$
is a {\it local equivalence} if and only if, for any object $C \in \calC$ and any commutative diagram of topological spaces
$$ \xymatrix{ S^{n-1} \ar[r] \ar@{^{(}->}[d] & |F_{\bigdot}(C)| \ar[d] \\
D^n \ar[r] & |G_{\bigdot}(C)|, } $$
there exists a collection of morphisms $\{ C_{\alpha} \rightarrow C\}$ which generate a covering
sieve on $C$, such that in each of the induced diagrams
$$ \xymatrix{ S^{n-1} \ar[r] \ar@{^{(}->}[d] & |F_{\bigdot}(C_{\alpha})| \ar[d] \\
D^n \ar[r] \ar@{-->}[ur] & |G_{\bigdot}(C_{\alpha})|, } $$
one can produce a dotted arrow so that the upper triangle commutes and the lower triangle
commutes up to a homotopy which is fixed on $S^{n-1}$.
\end{itemize}
\end{proposition}

We refer the reader to \cite{jardine} for a proof (one can also deduce Proposition \ref{jardinesardine} from Proposition \ref{goot}). We will refer to the model structure of Proposition \ref{jardinesardine} as the {\it local} model structure on $\bfA$.

\begin{remark}\label{pointeddesc}
In the case where the topos $\toposX$ of sheaves of sets on $\calC$ has enough points, there is a simpler description of the class $(W)$ of local equivalences: a map $F \rightarrow G$ of simplicial presheaves is a local equivalence if and only if it induces weak homotopy equivalences 
$F_{x} \rightarrow G_{x}$ of simplicial sets, after passing to the stalk at any point $x$ of $\toposX$.
We refer the reader to \cite{jardine} for details.
\end{remark}

Let $\bfA^{\degree}$ denote the full subcategory of $\bfA$ consisting of fibrant-cofibrant objects (with respect to the local model structure), and let $\calX = \sNerve(\bfA^{\degree})$ be the associated $\infty$-category. We observe that the local model structure on $\bfA$ is a localization of the injective model structure on $\bfA$. Consequently, the $\infty$-category $\calX$ is a localization of the $\infty$-category associated to the injective model structure on $\bfA$, which (in view of 
Proposition \ref{othermod}) is equivalent to $\calP( \Nerve(\calC))$. It is tempting to guess
that $\calX$ is equivalent to the left exact localization $\Shv( \Nerve(\calC))$ constructed in \S \ref{cough}. This is not true in general; however, as we will explain below, we can always recover $\calX$ as an accessible left-exact localization of $\Shv( \Nerve(\calC))$. In particular, $\calX$
is itself an $\infty$-topos.

In general, the difference between $\calX$ and $\Shv( \Nerve(\calC))$ is measured by the failure of Whitehead's theorem. Essentially by construction, the equivalences in $\bfA$ are those maps which induce isomorphisms on homotopy sheaves. In general, this assumption is not strong enough to guarantee that a morphism in $\Shv( \Nerve(\calC) )$ is an equivalence. However, this is the only difference: the $\infty$-category $\calX$ can be obtained from
$\Shv(\Nerve(\calC))$ by inverting the class of $\infty$-connective morphisms (Proposition \ref{suga}). Before proving this, we study the class of $\infty$-connective morphisms in an arbitrary
$\infty$-topos.

\begin{lemma}\label{swarp}
Let $p: \calC \rightarrow \calD$ be a Cartesian fibration of $\infty$-categories, let
$\calC'$ be a full subcategory of $\calC$, and suppose that for every $p$-Cartesian morphism
$f: C \rightarrow C'$ in $\calC$, if $C' \in \calC'$, then $C \in \calC'$. Let
$D$ be an object of $\calD$, and let $f: C \rightarrow C'$ be a morphism in the fiber
$\calC_{D} = \calC \times_{\calD} \{D\}$ which exhibits $C'$ as a $\calC^{0}_{D}$-localization
of $C$ (see Definition \ref{locaobj}). Then $f$ exhibits $C'$ as a $\calC$-localization of $\calC$.
\end{lemma}

\begin{proof}
According to Proposition \ref{verylonger}, $p$ induces a Cartesian fibration
$\calC_{C/} \rightarrow \calD_{D/}$, which restricts to give a Cartesian fibration
$p': \calC'_{C/} \rightarrow \calD_{D/}$. We observe that $f$ is an object of
$\calC'_{C/}$ which is an initial object of $(p')^{-1} \{ \id_{D} \}$ (Remark \ref{initrem}), 
and that $\id_{D}$ is an initial object of $\calD_{D/}$. Lemma \ref{sabreto} implies that $f$ is an initial object of $\calC'_{C/}$, so that $f$ exhibits $C'$ as $\calC'$-localization of $C$ (Remark \ref{initrem}) as desired.
\end{proof}

\begin{lemma}\label{swarp2}
Let $p: \calC \rightarrow \calD$ be a Cartesian fibration of $\infty$-categories, let
$\calC'$ be a full subcategory of $\calC$, and suppose that for every $p$-Cartesian morphism
$f: C \rightarrow C'$ in $\calC$, if $C' \in \calC'$, then $C \in \calC'$. Suppose that for each
object $D \in \calD$, the fiber $\calC'_{D} = \calC' \times_{\calD} \{D\}$ is a reflective subcategory
of $\calC_{D} = \calC \times_{\calD} \{D\}$ (see Remark \ref{reflective}). Then $\calC'$ is a reflective subcategory of $\calC$.
\end{lemma}

\begin{proof}
Combine Lemma \ref{swarp} with Proposition \ref{testreflect}.
\end{proof}

\begin{lemma}\label{swarp3}
Let $\calX$ be a presentable $\infty$-category, let $\calC$ be an accessible $\infty$-category,
and let $\alpha: F \rightarrow G$ be a natural transformation between accessible functors
$F,G: \calC \rightarrow \calX$. Let $\calC(n)$ be the full subcategory of $\calC$ spanned by those objects $C$ such that $\alpha(C): F(C) \rightarrow G(C)$ is $n$-truncated. Then $\calC(n)$ is an accessible subcategory of $\calC$ (see Definition \ref{defaccsub}).
\end{lemma}

\begin{proof}
We work by induction on $n$. If $n=-2$, then we have a (homotopy) pullback diagram
$$ \xymatrix{ \calC(n) \ar[r] \ar[d] & \calC \ar[d]^{\alpha} \\
\calE \ar[r] & \Fun(\Delta^1,\calX) }$$
where $\calE$ is the full subcategory of $\Fun(\Delta^1,\calX)$ spanned by
equivalences. The inclusion of $\calE$ into $\Fun(\Delta^1,\calX)$ is equivalent
to the diagonal map $\calX \rightarrow \Fun(\Delta^1, \calX)$, and therefore accessible.
Proposition \ref{horse2} implies that $\calC(n)$ is an accessible subcategory of $\calC$, as desired.

If $n \geq -1$, we apply the the inductive hypothesis to the diagonal functor
$\delta: F \rightarrow F \times_{G} F$, using Lemma \ref{trunc}.
\end{proof}

\begin{lemma}\label{tur}
Let $\calX$ be a presentable $\infty$-category, and let $-2 \leq n < \infty$. Let
$\calC$ be the full subcategory of $\Fun(\Delta^1,\calX)$ spanned by the $n$-truncated
morphisms. Then $\calC$ is a strongly reflective subcategory of $\Fun(\Delta^1,\calX)$.
\end{lemma}

\begin{proof}
Applying Lemma \ref{swarp2} to the projection $\Fun(\Delta^1,\calX) \rightarrow \Fun( \{1\}, \calX)$, we conclude that $\calC$ is a reflective subcategory of $\Fun(\Delta^1,\calX)$. The accessibility
of $\calC$ follows from Lemma \ref{swarp3}.
\end{proof}

\begin{lemma}\label{swarp4}
Let $\calX$ be an $\infty$-topos, let $0 \leq n \leq \infty$, and let $\calD(n)$ be the full subcategory of $\Fun(\Delta^1,\calX)$ spanned by the $n$-connective morphisms of $\calX$. Then
$\calD(n)$ is an accessible subcategory of $\calX$ and is stable under colimits in 
$\calX$.
\end{lemma}

\begin{proof}
Suppose first that $n < \infty$. Let $\calC(n) \subseteq \Fun(\Delta^1,\calX)$ be the full subcategory spanned by the $n$-truncated morphisms in $\calX$. According to Lemma \ref{tur}, the
inclusion $\calC(n) \subseteq \Fun(\Delta^1,\calX)$ has a left adjoint $L: \Fun(\Delta^1,\calX) \rightarrow \calC(n)$. Moreover, the proof of Lemma \ref{swarp} shows that $f$ is $n$-connective if and only if $Lf$ is an equivalence. It is easy to see that the full subcategory $\calE \subseteq \calC(n)$ spanned by the equivalences is stable under colimits in $\calC(n)$, so that $\calD(n)$ is stable under colimits in $\Fun(\Delta^1,\calX)$. The accessibility of $\calD(n)$ follows from the existence of the (homotopy) pullback diagram
$$ \xymatrix{ \calD(n) \ar[r] \ar[d] & \Fun(\Delta^1,\calX) \ar[d]^{L} \\
\calE \ar[r] & \calC(n) }$$
and Proposition \ref{horse2}.

If $n=\infty$, we observe that $\calD(n) = \cup_{m < \infty} \calD(m)$, which is manifestly stable under colimits, and is an accessible subcategory of $\calX^{\Delta^1}$ by Proposition \ref{boundint}.
\end{proof}

\begin{proposition}\label{goober3}
Let $\calX$ be an $\infty$-topos, and let $S$ denote the
collection of $\infty$-connective morphisms of $\calX$. Then $S$ is strongly 
saturated and of small generation $($see Definition \ref{saturated2}$)$. 
\end{proposition}

\begin{proof}
Lemma \ref{swarp4} implies that $S$ is stable under colimits in $\Fun(\Delta^1,\calX)$, and
Corollary \ref{pusherr} shows that $S$ is stable under pushouts. To prove that $S$ has the two-out-of-three property, we consider a diagram $\sigma: \Delta^2 \rightarrow \calX$, which we depict as
$$ \xymatrix{ & Y \ar[dr]^{g} & \\
X \ar[ur]^{f} \ar[rr]^{h} & & Z. }$$
If $f$ is $\infty$-connective, then Proposition \ref{inftychange} implies that $g$
is $\infty$-connective if and only if $h$ is $\infty$-connective. Suppose that $g$ and $h$
are $\infty$-connective. The long exact sequence
$$\ldots \rightarrow f^{\ast} \pi_{n+1}(g) \rightarrow \pi_n(f) \rightarrow \pi_n(h) \rightarrow f^{\ast}
\pi_n(g) \rightarrow \pi_{n-1}(f) \rightarrow \ldots$$ 
of Remark \ref{sequence} shows that $\pi_n(f) \simeq \ast$ for all $n \geq 0$. It
will therefore suffice to prove that $f$ is an effective epimorphism. According to Proposition \ref{conslice}, it will suffice to show that $\sigma$ is an effective epimorphism in 
$\calX_{/Z}$. According to Proposition \ref{pi00detects}, it suffices to show that
$\tau_{\leq 0}^{\calX_{/Z}}(h)$ and $\tau_{\leq 0}^{\calX_{/Z}}(g)$ are both final objects of
$\calX_{/Z}$, which follows from the $0$-connectivity of $g$ and $h$ (Proposition \ref{goober}). 

To show that $S$ is of small generation, it suffices (in view of Lemma \ref{perry}) to show that
the full subcategory of $\Fun(\Delta^1,\calX)$ spanned by $S$ is accessible. This follows from Lemma \ref{swarp4}.
\end{proof}

Let $\calX$ be an $\infty$-topos. We will say that an object $X$ of $\calX$ is
{\it hypercomplete}\index{gen}{hypercomplete!object} if it is local with respect to the class of $\infty$-connective morphisms. 
Let $\calX^{\hyp}$\index{gen}{hypercomplete!$\infty$-topos}\index{not}{Xhyp@$\calX^{\hyp}$} denote the full subcategory of $\calX$ spanned by the hypercomplete objects of $\calX$. Combining Propositions \ref{goober3} and \ref{local}, we deduce that
$\calX^{\hyp}$ is an accessible localization of $\calX$. Moreover, since Proposition \ref{inftychange} implies that the class of $\infty$-connective morphisms is stable under pullback, we deduce from Proposition \ref{charleftloc} that $\calX^{\hyp}$ is a {\em left exact} localization
of $\calX$. It follows that $\calX^{\hyp}$ is itself an $\infty$-topos. We will show in a moment that
$\calX^{\hyp}$ can be described by a universal property.

\begin{lemma}\label{sshock1}
Let $\calX$ be an $\infty$-topos, and let $n < \infty$. Then
$\tau_{\leq n} \calX \subseteq \calX^{\hyp}$.
\end{lemma}

\begin{proof}
Corollary \ref{goober2} implies that an $n$-truncated object of $\calX$ is local with respect to every $n$-connective morphism of $\calX$, and therefore with respect to every $\infty$-connective morphism of $\calX$. 
\end{proof}

\begin{lemma}\label{sshock2}
Let $\calX$ be an $\infty$-topos, let $L: \calX \rightarrow \calX^{\hyp}$ be a left adjoint to the inclusion, and let $X \in \calX$ be such that $LX$ is an $\infty$-connective object of
$\calX^{\hyp}$. Then $LX$ is a final object of $\calX^{\hyp}$.
\end{lemma}

\begin{proof}
For each $n  < \infty$, we have equivalences
$$1_{\calX} \simeq \tau_{\leq n}^{\calX^{\hyp}} LX \simeq L \tau_{\leq n}^{\calX} X\simeq \tau_{\leq n}^{\calX} X$$
where the first is because of our hypothesis that $LX$ is $\infty$-connective, the second
is given by Proposition \ref{compattrunc}, and the third by Lemma \ref{sshock1}. It follows
that $X$ is an $\infty$-connective object of $\calX$, so that $LX$ is a final object of
$\calX^{\hyp}$ by construction.
\end{proof}

We will say that an $\infty$-topos $\calX$ is {\it hypercomplete} if $\calX^{\hyp} = \calX$; in
other words, $\calX$ is hypercomplete if every $\infty$-connective morphism of $\calX$ is an equivalence, so that Whitehead's theorem holds in $\calX$.

\begin{remark}
In \cite{toen}, the authors use the term {\it $t$-completeness} to refer to the property that we have called hypercompleteness.
\end{remark}

\begin{lemma}
Let $\calX$ be an $\infty$-topos. Then the $\infty$-topos $\calX^{\hyp}$ is hypercomplete.
\end{lemma}

\begin{proof}
Let $f: X \rightarrow Y$ be an $\infty$-connective morphism in $\calX^{\hyp}$. Applying
Lemma \ref{sshock2} to the $\infty$-topos $(\calX^{\hyp})_{/Y} \simeq (\calX_{/Y})^{\hyp}$, we deduce that $f$ is an equivalence.
\end{proof}

We are now prepared to characterize $\calX^{\hyp}$ by a universal property:

\begin{proposition}
Let $\calX$ and $\calY$ be $\infty$-topoi. Suppose that $\calY$ is hypercomplete.
Then composition with the inclusion $\calX^{\hyp} \subseteq \calX$ induces an isomorphism
$$ \Fun_{\ast}(\calY, \calX^{\hyp} ) \rightarrow \Fun_{\ast}(\calY, \calX).$$
\end{proposition}

\begin{proof}
Let $f_{\ast}: \calY \rightarrow \calX$ be a geometric morphism; we wish to prove that
$f_{\ast}$ factors through $\calX^{\hyp}$. Let $f^{\ast}$ denote a left adjoint to $f_{\ast}$; it will suffice to show that $f^{\ast}$ carries each $\infty$-connective morphism $u$ of
$\calX$ to an equivalence in $\calY$. Proposition \ref{inftychange} implies that
$f^{\ast}(u)$ is $\infty$-connective, and the hypothesis that $\calY$ is hypercomplete guarantees that $u$ is an equivalence.
\end{proof}

The following result establishes the relationship between our theory of hypercompleteness
and the Brown-Joyal-Jardine theory of simplicial presheaves.

\begin{proposition}\label{suga}
Let $\calC$ be a small category equipped with a Grothendieck topology, and let $\bfA$
denote the category of simplicial presheaves on $\calC$, endowed with the local model structure $($ see Proposition \ref{jardinesardine} $)$.
Let $\bfA^{\degree}$ denote the full subcategory consisting of fibrant-cofibrant objects, and let
$\calA = \sNerve( \bfA^{\degree} )$ be the corresponding $\infty$-category. Then $\calA$
is equivalent to $\Shv(\calC)^{\hyp}$; in particular, it is a hypercomplete $\infty$-topos.
\end{proposition}

\begin{proof}
Let $\calP(\calC)$ denote the $\infty$-category $\calP(\Nerve(\calC))$ of presheaves on
$\Nerve(\calC)$, and let
$\bfA'$ denote the model category of simplicial presheaves on $\calC$, endowed with the {\em injective} model structure of \S \ref{quasilimit3}. According to Proposition \ref{gumby444}, 
the simplicial nerve functor induces an equivalence
$$ \theta: \sNerve ({\bfA'}^{\degree}) \rightarrow \calP(\calC).$$
We may identify $\sNerve(\bfA^{\degree})$ with the full subcategory of 
$\sNerve({\bfA'}^{\degree})$ spanned by the $S$-local objects, where $S$ is the class
of local equivalences (Proposition \ref{suritu}). 

We first claim that $\theta | \sNerve(\bfA^{\degree})$ factors
through $\Shv(\calC)$. Consider an object $C \in \calC$ and a sieve
$\calC^{(0)}_{/C} \subseteq \calC_{/C}$. Let $\chi_{C}: \calC \rightarrow \Set$ be the functor
$D \mapsto \Hom_{\calC}(D,C)$ represented by $\calC$, let $\chi_{C}^{(0)}$ be the subfunctor of $\chi_{C}$ determined by the sieve $\calC^{(0)}_{/C}$, and let $i: \chi_{C}^{(0)} \rightarrow \chi_{C}$ be the inclusion. We regard $\chi_{C}$ and $\chi_{C}^{(0)}$ as objects of
simplicial presheaves on $\calC$, which take values in the full subcategory 
of $\sSet$ spanned by the {\em constant} simplicial sets. We observe that
every simplicial presheaf on $\calC$ which is valued in constant simplicial sets
is automatically fibrant, and every object of $\bfA'$ is cofibrant. Consequently, we
may regard $i$ as a morphism in the $\infty$-category $\sNerve (\bfA')^{\degree}$.
It is easy to see that $\theta(i)$ represents the monomorphism
$U \rightarrow j(C)$ classified by the sieve $\calC^{(0)}_{/C}$. If $\calC^{(0)}_{/C}$ is a covering sieve on $C$, then $i$ is a local equivalence. Consequently, every object $X \in \sNerve (\bfA^{\degree})$ is $i$-local, so that $\theta(X)$ is $\theta(i)$-local. By construction, $\Shv(\calC)$ 
is the full subcategory of $\calP(\calC)$ spanned by those objects which are $\theta(i)$-local
for every covering sieve $\calC^{(0)}_{/C}$ on every object $C \in \calC$. We conclude
that $\theta | \sNerve(\bfA^{\degree})$ factors through $\Shv(\calC)$. 

Let $\calX = \theta^{-1} \Shv(\calC)$, so that $\sNerve(\bfA^{\degree})$ can be identified
with the collection of $S'$-local objects of $\calX$, where $S'$ is the collection of all morphisms
in $\calX$ which belong to $S$. Then $\theta$ induces an equivalence
$\sNerve (\bfA^{\degree}) \rightarrow \theta(S')^{-1} \Shv(\calC)$. We now observe
that a morphism $f$ in $\calX$ belongs to $S'$ if and only if $\theta(f)$ is an $\infty$-connective morphism in $\Shv(\calC)$ (since the condition of being a local equivalence can be tested on homotopy sheaves). It follows that $\theta(S')^{-1} \Shv(\calC) = \Shv(\calC)^{\hyp}$, as desired. 
\end{proof}

\begin{remark}
In \cite{toen}, the authors discuss a generalization of Jardine's construction, in which
the category $\calC$ is replaced by a simplicial category. Proposition \ref{suga} holds in this more general situation as well.
\end{remark}

We conclude this section with a few remarks about localizations of an $\infty$-topos
$\calX$. In \S \ref{leloc} we introduced the class of topological localizations of $\calX$, consisting of those left exact localizations which can be obtained by inverting monomorphisms in $\calX$.
The hypercompletion $\calX^{\hyp}$ is, in some sense, at the other extreme: it is obtained by
inverting the $\infty$-connective morphisms in $\calX$, which are never monomorphisms unless they are already equivalences. In fact, $\calX^{\hyp}$ is the {\em maximal} (left exact) localization of $\calX$ which can be obtained without inverting monomorphisms:

\begin{proposition}\label{antitopchar}
Let $\calX$ and $\calY$ be $\infty$-topoi, and let $f^{\ast}: \calX \rightarrow \calY$ be a left exact, colimit preserving functor.
The following conditions are equivalent:
\begin{itemize}
\item[$(1)$] For every monomorphism $u$ in $\calX$, if $f^{\ast} u$ is an equivalence in $\calY$, then
$u$ is an equivalence in $\calX$.
\item[$(2)$] For every morphism $u \in \calX$, if $f^{\ast} u$ is an equivalence in $\calY$, then $f$ is $\infty$-connective.
\end{itemize}
\end{proposition}

\begin{proof}
Suppose first that $(2)$ is satisfied. If $u$ is a monomorphism and $f^{\ast} u$ is an equivalence in $\calY$, then $u$ is $\infty$-connective. In particular, $u$ is both a monomorphism and an effective epimorphism, and therefore an equivalence in $\calX$. This proves $(1)$. Conversely, suppose that $(1)$ is satisfied, and let $u: X \rightarrow Z$ be an arbitrary morphism in $\calX$ such
that $f^{\ast}(u)$ is an equivalence. We will prove by induction on $n$ that $u$ is $n$-connective.

We first consider the case $n=0$. Choose a factorization
$$ \xymatrix{ & Y \ar[dr]^{u''} & \\
X \ar[ur]^{u'} \ar[rr]^{u} & & Z }$$
where $u'$ is an effective epimorphism, and $u''$ is a monomorphism. Since
$f^{\ast} u$ is an equivalence, Corollary \ref{composite} implies that $f^{\ast} u''$
is an effective epimorphism. Since $f^{\ast} u''$ is also a monomorphism (in virtue of our
assumption that $f$ is left exact), we conclude that $f^{\ast} u''$ is an equivalence.
Applying $(1)$, we deduce that $u''$ is an equivalence, so that $u$ is an effective epimorphism as desired.

Now suppose $n > 0$. According to Proposition \ref{trowler}, it will suffice to show that the
diagonal map $\delta: X \rightarrow X \times_{Z} X$ is $(n-1)$-connective. By the inductive hypothesis, it will suffice to prove that $f^{\ast}(\delta)$ is an equivalence in $\calY$. We conclude by observing that 
$f^{\ast}$ is left exact, so we can identify $\delta$ with the diagonal map associated to
the equivalence $f^{\ast}(u): f^{\ast} X \rightarrow f^{\ast} Z$.
\end{proof}

\begin{definition}\index{gen}{localization!cotopological}
Let $\calX$ be an $\infty$-topos, and let $\calY \subseteq \calX$ be an accessible
left exact localization of $\calX$. We will say that $\calY$ is an {\em cotopological} localization
of $\calX$ if the left adjoint $L: \calX \rightarrow \calY$ to the inclusion of $\calY$ in $\calX$
satisfies the equivalent conditions of Proposition \ref{antitopchar}.
\end{definition}

\begin{remark}\label{sorkum}
Let $f^{\ast}: \calX \rightarrow \calY$ be the left adjoint of a geometric morphism between $\infty$-topoi, and suppose that the equivalent conditions of Proposition \ref{antitopchar} are satisfied.
Let $u: X \rightarrow Z$ be a morphism in $\calX$, and choose a factorization
$$ \xymatrix{ & Y \ar[dr]^{u''} & \\
X \ar[ur]^{u'} \ar[rr]^{u} & & Z }$$
where $u'$ is an effective epimorphism and $u''$ is a monomorphism. Then
$u''$ is an equivalence if and only if $f^{\ast}(u'')$ is an equivalence. Applying
Corollary \ref{composite}, we conclude that $u$ is an effective epimorphism if and only if
$f^{\ast}(u)$ is an effective epimorphism.
\end{remark}

The hypercompletion $\calX^{\hyp}$ of an $\infty$-topos $\calX$ can be characterized as the {\em maximal} cotopological localization of $\calX$ (that is, the cotopological localization which
is obtained by inverting as many morphisms as possible). According to our next result, every localization can be obtained by combining topological and cotopological localizations:

\begin{proposition}\label{factanti}
Let $\calX$ be an $\infty$-topos, and let $\calX'' \subseteq \calX$ be an accessible, left
exact localization of $\calX$. Then there exists a topological localization
$\calX' \subseteq \calX$ such that $\calX'' \subseteq \calX'$ is a cotopological localization
of $\calX'$.
\end{proposition}

\begin{proof}
Let $L: \calX \rightarrow \calX''$ be a left adjoint to the inclusion, let
$S$ be the collection of all monomorphisms $u$ in $\calX$ such that $Lu$ is an equivalence,
and let $\calX' = S^{-1} \calX$ be the collection of $S$-local objects of $\calX$. Since $L$ is left exact, $S$ is stable under pullbacks and therefore determines a topological localization of $\calX$. By construction, we have $\calX'' \subseteq \calX'$. The restriction $L| \calX'$ exhibits
$\calX''$ as an accessible left exact localization of $\calX'$. Let $u$ be a monomorphism
in $\calX'$ such that $Lu$ is an equivalence. Then $u$ is a monomorphism in $\calX$, so
that $u \in S$. Since $\calX'$ consists of $S$-local objects, we conclude that $u$ is an equivalence. It follows that $\calX''$ is a cotopological localization of $\calX'$, as desired.
\end{proof}

\begin{remark}
It is easy to see that the factorization of Proposition \ref{factanti} is essentially uniquely determined: more precisely, $\calX'$ is unique provided that we assume that it is stable under equivalences in $\calX$.
\end{remark}

Combining Proposition \ref{factanti} with Remark \ref{charnice}, we see that every $\infty$-topos $\calX$ can be obtained in following way:
\begin{itemize}
\item[$(1)$] Begin with the $\infty$-category $\calP(\calC)$ of presheaves on some small $\infty$-category $\calC$.

\item[$(2)$] Choose a Grothendieck topology on $\calC$: this is equivalent to choosing a 
left exact localization of the underlying topos $\Disc(\calP(\calC)) = \Set^{ \h{\calC^{op}}}$.

\item[$(3)$] Form the associated topological localization $\Shv(\calC) \subseteq \calP(\calC)$, which can be described as the pullback 
$$ \calP(\calC) \times_{ \calP( \Nerve ( \h{\calC}) )} \Shv( \Nerve ( \h{\calC} ) )$$
in $\RGeom$.

\item[$(4)$] Form a cotopological localization of $\Shv(\calC)$ by inverting some class
of $\infty$-connective morphisms of $\Shv(\calC)$.
\end{itemize}

\begin{remark}\label{comk1}
Let $\calX$ be an $\infty$-topos. The collection of all $\infty$-connective morphisms
in $\calX$ is saturated. It follows from Proposition \ref{nir} that there exists a factorization system
$(S_L, S_R)$ on $\calX$, where $S_L$ is the collection of all $\infty$-connective morphisms in
$\calX$. We will say that a morphism in $\calX$ is {\it hypercomplete} if it belongs to $S_R$.
Unwinding the definitions (and using the fact that a morphism in $\calX_{/Y}$ is $\infty$-connective if and only if its image in $\calX$ is $\infty$-connective), we conclude that a morphism
$f: X \rightarrow Y$ is hypercomplete if and only if it is hypercomplete when viewed as an object of the $\infty$-topos $\calX_{/Y}$ (see \S \ref{hyperstacks}). 

Using Proposition \ref{swimmm}, we deduce that the collection of hypercomplete morphisms in $\calX$ is stable under limits and the formation of pullback squares.
\end{remark}

\begin{remark}\label{suchlike}
Let $\calX$ be an $\infty$-topos. The condition that a morphism $f: X \rightarrow Y$ be hypercomplete is {\em local}: that is, if $\{ Y_{\alpha} \rightarrow Y \}$ is a collection of morphisms which
determine an effective epimorphism $\coprod Y_{\alpha} \rightarrow Y$, and each of the induced maps
$f_{\alpha}: X \times_{Y} Y_{\alpha} \rightarrow Y_{\alpha}$ is hypercomplete, then $f$ is hypercomplete. To prove this, we set $Y_0 = \coprod_{\alpha} Y_{\alpha}$; then
$\calX_{/Y_0} \simeq \prod_{ \alpha} \calX_{/Y_{\alpha}}$ (since coproducts in $\calX$ are disjoint),
so it is easy to see that the induced map $f': X \times_{Y} Y_0 \rightarrow Y_0$ is hypercomplete.
Let $Y_{\bigdot}$ be the simplicial object of $\calX$ given by the \Cech nerve of the effective epimorphism $Y_0 \rightarrow Y$. For every map $Z \rightarrow Y$, let $Z_{\bigdot}$ be
the simplicial object described by the formula
$Z_{n} = Y_{n} \times_{Y} Z$
(equivalently, $Z_{\bigdot}$ is the \Cech nerve of the effective epimorphism $Z \times_{Y} Y_0 \rightarrow Z$). Using Remark \ref{comk1}, we conclude that each of the maps $X_{n} \rightarrow Y_{n}$ is hypercomplete.

For every map $A \rightarrow Y$, the mapping space
$\bHom_{ \calX_{/Y}}( A, X)$ can be obtained as the totalization of a cosimplicial space
$$ n \mapsto \bHom_{ \calX_{/Y_n} }( A_{n}, X_n).$$
If $g: A \rightarrow B$ is an $\infty$-connective morphism in $\calX_{/Y}$, then
each of the induced maps $A_{n} \rightarrow B_{n}$ is $\infty$-connective, so the induced map
$$ \bHom_{ \calX_{/Y_n}}( B_{n}, X_{n} ) \rightarrow \bHom_{ \calX_{/Y_{n}}}( A_{n}, X_n)$$
is a homotopy equivalence. Passing to the totalization, we obtain a homotopy equivalence
$\bHom_{\calX_{/Y}}(B,X) \rightarrow \bHom_{\calX_{/Y}}( A, X)$. Thus $f$ is hypercomplete, as desired.
\end{remark}

\subsection{Hypercoverings}\label{hcovh}

Let $\calX$ be an $\infty$-topos. In \S \ref{hyperstacks}, we defined the {\em hypercompletion}
$\calX^{\hyp} \subseteq \calX$ to be the left exact localization of $\calX$ obtained by inverting the $\infty$-connective morphisms. In this section, we will give an alternative description of the hypercomplete objects $X \in \calX^{\hyp}$: they are precisely those objects of $\calX$ which satisfy a descent condition with respect to hypercoverings (Theorem \ref{surp}). We begin by reviewing the definition of a hypercovering.

Let $X$ be a topological space, and let $\calF$ be a presheaf of sets on $X$. To construct the sheaf associated to $\calF$, it is natural to consider the presheaf $\calF^{+}$, defined by
$$ \calF^{+} = \varinjlim_{ \calU } \varprojlim_{ V \in \calU } \calF(V).$$ 
Here the direct limit is taken over all sieves $\calU$ which cover $U$. There is an obvious map
$ \calF \rightarrow \calF^{+}$, which is an isomorphism whenever $\calF$ is a sheaf. Moreover, $\calF^{+}$ is ``closer'' to being a sheaf than $\calF$ is. More precisely, $\calF^{+}$ is always a separated presheaf: two sections of $\calF^{+}$ which agree locally automatically coincide.
If $\calF$ is itself a separated presheaf, then $\calF^{+}$ is a sheaf.\index{not}{Fcal+@$\calF^{+}$}

For a general presheaf $\calF$, we need to apply the above construction twice to construct the associated sheaf $(\calF^{+})^{+}$. To understand the problem, let us try to prove that $\calF^{+}$ is a sheaf (to see where the argument breaks down). Suppose given an open covering
$X = \bigcup U_{\alpha}$, and a collection of sections $s_{\alpha} \in \calF^{+}(U_{\alpha})$ such that
$$ s_{\alpha} | U_{\alpha} \cap U_{\beta} = s_{\beta} | U_{\alpha} \cap U_{\beta}. $$
Refining the covering $U_{\alpha}$ if necessary, we may assume that each $s_{\alpha}$
is the image of some section $t_{\alpha} \in \calF(U_{\alpha})$. However, the equation
$$ t_{\alpha} | U_{\alpha} \cap U_{\beta} = t_{\beta} | U_{\alpha} \cap U_{\beta} $$
only holds locally on $U_{\alpha} \cap U_{\beta}$, so the sections $t_{\alpha}$ do not necessarily determine a global section of $\calF^{+}$. To summarize: the freedom to consider
arbitrarily fine open covers $\calU = \{ U_{\alpha} \}$ is not enough; we also need to be able to refine the intersections $U_{\alpha} \cap U_{\beta}$. This leads very naturally to the notion of a {\it hypercovering}. Roughly speaking, a hypercovering of $X$ consists of an open covering
$\{ U_{\alpha} \}$ of $X$, an open covering of $\{ V_{\alpha\beta\gamma} \}$ of each intersection $U_{\alpha} \cap U_{\beta}$, and analogous data associated to more complicated intersections (see Definition \ref{worum} for a more precise formulation).\index{gen}{hypercovering}

In classical sheaf theory, there are two ways to construct the sheaf associated to a presheaf $\calF$:
\begin{itemize}
\item[$(1)$] One can apply the construction $\calF \mapsto \calF^{+}$ twice.
\item[$(2)$] Using the theory of hypercoverings, one can proceed directly by defining
$$\calF^{\dagger}(U) = \varinjlim_{\calU} \varprojlim \calF(V)$$
where the direct limit is now taken over arbitrary {\em hypercoverings} $\calU$.
\end{itemize}

In higher category theory, the difference between these two approaches becomes more prominent. For example, suppose that $\calF$ is not a presheaf of sets, but a presheaf of {\em groupoids} on
$X$. In this case, one can construct the associated sheaf of groupoids using either approach. However, in the case of approach $(1)$, it is necessary to apply the construction
$\calF \mapsto \calF^{+}$ {\em three} times: the first application guarantees that the automorphism groups of sections of $\calF$ are separated presheaves, the second guarantees that they are sheaves, and the third guarantees that $\calF$ itself satisfies descent. More generally, if $\calF$ is a sheaf of $n$-truncated spaces, then the sheafification of $\calF$ via approach $(1)$ takes place in $(n+2)$-stages. 

When we pass to the case $n= \infty$, the situation becomes more complicated. If $\calF$ is a presheaf of spaces on $X$, then it is not reasonable to expect to obtain a sheaf by applying the construction $\calF \mapsto \calF^{+}$ any finite number of times. In fact, it is not obvious that $\calF^{+}$ is any closer than $\calF$ to being a sheaf. Nevertheless, this is true: we can construct the sheafification of $\calF$ via a {\em transfinite iteration} of the construction $\calF \mapsto \calF^{+}$. More precisely, we define a transfinite sequence of presheaves
$$ \calF(0) \rightarrow \calF(1) \rightarrow \ldots $$
as follows:
\begin{itemize}
\item[$(i)$] Let $\calF(0) = \calF$.
\item[$(ii)$] For every ordinary $\alpha$, let $\calF(\alpha+1) = \calF(\alpha)^{+}$.
\item[$(iii)$] For every limit ordinal $\lambda$, let $\calF(\lambda) = \colim_{\alpha} \calF(\alpha)$, where $\alpha$ ranges over ordinals less than $\lambda$.
\end{itemize}

One can show that the above construction {\em converges}, in the sense that $\calF(\alpha)$
is a sheaf for $\alpha \gg 0$ (and therefore $\calF(\alpha) \simeq \calF(\beta)$ for $\beta \geq \alpha)$. Moreover, $\calF(\alpha)$ is universal among sheaves of spaces which admits a map from $\calF$.

Alternatively, one use the construction $\calF \mapsto \calF^{\dagger}$ to construct a sheaf of spaces from $\calF$ in a single step. The universal property asserted above guarantees the existence of a morphism of sheaves $\theta: \calF(\alpha) \rightarrow \calF^{\dagger}$. However, the morphism $\theta$ is generally {\em not} an equivalence. Instead, $\theta$ realizes
$\calF^{\dagger}$ as the {\em hypercompletion} of $\calF(\alpha)$ in the $\infty$-topos $\Shv(X)$. 
We will not prove this statement directly, but will instead establish a reformulation (Corollary \ref{charhyp}) which does not make reference to the sheafification constructions outlined above.

Before we can introduce the definition of a hypercovering, we need to review some simplicial terminology.

\begin{notation}\index{not}{Deltaleqn@$\cDelta^{\leq n}$}
For each $n \geq 0$, let $\cDelta^{\leq n}$ denote the full
subcategory of $\cDelta$ spanned by the set of objects
$\{ [0], \ldots, [n] \}$. 
If $\calX$ is a presentable $\infty$-category, the 
restriction functor
$$ \sk_n: \calX_{\Delta} \rightarrow \Fun(\Nerve (\cDelta^{\leq n})^{op}, \calX)$$
has a right adjoint, given by right Kan extension along the inclusion $\Nerve (\cDelta^{\leq n})^{op} \subseteq \Nerve(\cDelta)^{op}$. Let $\cosk_n: \calX_{\Delta} \rightarrow \calX_{\Delta}$ be the composition of $\sk_n$ with its right adjoint. We will refer to $\cosk_{n}$ as the {\it $n$-coskeleton functor}.\index{gen}{skeleton}\index{gen}{coskeleton}\index{not}{skn@$\sk_n$}\index{not}{coskn@$\cosk_n$}
\end{notation}

\begin{definition}\label{worum}\index{gen}{hypercovering}\index{gen}{hypercovering!effective}
Let $\calX$ be an $\infty$-topos. A simplicial object $U_{\bigdot} \in \calX_{\Delta}$ is a {\it hypercovering of $\calX$} if, for each $n \geq 0$, the unit map
$$ U_{n} \rightarrow ( \cosk_{n-1} U_{\bigdot} )_n $$
is an effective epimorphism. We will say that $U_{\bigdot}$ is an {\em effective hypercovering of $\calX$}
if the colimit of $U_{\bigdot}$ is a final object of $\calX$.
\end{definition}

\begin{remark}
More informally, a simplicial object $U_{\bigdot} \in \calX_{\Delta}$ is a hypercovering of $\calX$ if each of the associated maps
$$ U_0 \rightarrow 1_{\calX}$$
$$ U_1 \rightarrow U_0 \times U_0 $$
$$ U_2 \rightarrow \ldots $$
is an effective epimorphism.
\end{remark}

\begin{lemma}\label{fier1}
Let $\calX$ be an $\infty$-topos, and let $U_{\bigdot}$ be a simplicial object
in $\calX$.  Let $L: \calX \rightarrow \calX^{\hyp}$ be a left adjoint to the inclusion. The following conditions are equivalent:
\begin{itemize}
\item[$(1)$] The simplicial object $U_{\bigdot}$ is a hypercovering of $\calX$.
\item[$(2)$] The simplicial object $L \circ U_{\bigdot}$ is a hypercovering
of $\calX^{\hyp}$.
\end{itemize}
\end{lemma}

\begin{proof}
Since $L$ is left exact, we can identify $L \circ \cosk_{n} U_{\bigdot}$ with
$\cosk_{n} (L \circ U_{\bigdot})$. The desired result now follows from Remark \ref{sorkum}.
\end{proof}

\begin{lemma}\label{fier0}
Let $\calX$ be an $\infty$-topos, and let $U$ be an $\infty$-connective object of $\calX$.
Let $U_{\bigdot}$ be the constant simplicial object with value $U$.
Then $U_{\bigdot}$ is a hypercovering of $\calX$.
\end{lemma}

\begin{proof}
Using Lemma \ref{fier1}, we can reduce to the case where $\calX$ is hypercomplete.
Then $U \simeq 1_{\calX}$, so that $U_{\bigdot}$ is equivalent to the constant functor with value $1_{\calX}$, and is therefore a final object of $\calX_{\Delta}$. For each $n \geq 0$, the coskeleton functor $\cosk_{n-1}$ preserves small limits, so $\cosk_{n-1} U_{\bigdot}$ is also a final object of $U_{\bigdot}$. It follows that the unit map $U_{\bigdot} \rightarrow \cosk_{n-1} U_{\bigdot}$ is an equivalence.
\end{proof}

\begin{notation}
Let $\cDelta_{s}$ be the subcategory of $\cDelta$ with the same objects, but where the morphisms are given by {\em injective} order preserving maps between nonempty linearly ordered sets.
If $\calX$ is an $\infty$-category, we will refer to a diagram $\Nerve(\cDelta_{s})^{op} \rightarrow \calX$ as a {\it semisimplicial object of $\calX$}.
\end{notation}

\begin{lemma}\label{bball4}
The inclusion $\Nerve(\cDelta^{op}_{s}) \subseteq \Nerve(\cDelta^{op})$ is cofinal. 
\end{lemma}

\begin{proof}
According to Theorem \ref{hollowtt}, it will suffice to prove that for every $n \geq 0$, the 
category $\calC = \cDelta_{s} \times_{\cDelta} \cDelta_{/ [n]}$ has weakly contractible nerve.
To prove this, we let $F: \calC \rightarrow \calC$ be the constant functor taking value given by the inclusion $[0] \subseteq [n]$, and $G: \calC \rightarrow \calC$ the functor which carries
an arbitrary map $[m] \rightarrow [n]$ to the induced map
$[0] \amalg [m] \rightarrow [n]$. We have natural transformations of functors
$$ F \rightarrow G \leftarrow \id_{\calC}. $$
Let $X$ be the topological space $| \Nerve(\calC) |$. The natural transformations above show that the identity map $\id_{X}$ is homotopic to a constant, so that $X$ is contractible as desired.
\end{proof}

Consequently, if $U_{\bigdot}$ is a simplicial object in an $\infty$-category $\calX$, and $U_{\bigdot}^{s} = U_{\bigdot} | \Nerve(\cDelta^{op}_{s})$ is the associated semisimplicial object, then we can identify colimits of $U_{\bigdot}$ with colimits of $U_{\bigdot}^{s}$.

We will say that a simplicial object $U_{\bigdot}$ in an $\infty$-category
$\calX$ is {\it $n$-coskeletal}\index{gen}{coskeletal} if it is a right Kan extension of its restriction to
$\Nerve(\cDelta^{op}_{\leq n})$. Similarly, we will say that a semisimplicial
object of $U_{\bigdot}$ of $\calX$ is {\it $n$-coskeletal} if it is a right Kan extension
of its restriction to $\Nerve(\cDelta^{op}_{s,\leq n})$, where
$\cDelta_{s,\leq n} = \cDelta_{s} \times_{\cDelta} \cDelta_{\leq n}$.

\begin{lemma}\label{bball5}
Let $\calX$ be an $\infty$-category, let $U_{\bigdot}$ be a simplicial object of $\calX$, and let $U^{s}_{\bigdot} = U_{\bigdot} | \Nerve(\cDelta_{s}^{op})$ the associated semisimplicial object. Then $U_{\bigdot}$ is $n$-coskeletal if and only if $U^{s}_{\bigdot}$ is $n$-coskeletal.
\end{lemma}

\begin{proof}
It will suffice to show that, for each $\Delta^m \in \cDelta$, the nerve of the inclusion
$$ (\cDelta_{s})_{/[m]} \times_{ \cDelta_{s} } \cDelta_{s, \leq n} \subseteq 
\cDelta_{ / [m] } \times_{\cDelta} \cDelta_{\leq n}$$
is cofinal. Let $\theta: [m'] \rightarrow [m]$ be an object of 
$\cDelta_{ / [m] } \times_{\cDelta} \cDelta_{\leq n}$. We let
$\calC$ denote the category of all factorizations 
$$ [m'] \stackrel{\theta'}{\rightarrow} [m''] \stackrel{\theta''}{\rightarrow} [m]$$
for $\theta$ such that $\theta''$ is a monomorphism and $m'' \leq n$. According to Theorem \ref{hollowtt}, it will suffice to prove that $\Nerve(\calC)$ is weakly contractible (for every choice of $\theta$). We now
simply observe that $\calC$ has an initial object (given by the unique a factorization where
$\theta'$ is an epimorphism). 
\end{proof}

\begin{lemma}[\cite{hollander}]\label{bball1}
Let $\calX$ be an $\infty$-topos, and let $U_{\bigdot}$ be an $n$-coskeletal hypercovering of $\calX$. Then $U_{\bigdot}$ is effective.
\end{lemma}

\begin{proof}
We will prove this result by induction on $n$. If $n=0$, then $U_{\bigdot}$ can be identified with the underlying groupoid of the \Cech nerve of the map $\theta: U_0 \rightarrow 1_{\calX}$, where $1_{\calX}$ is a final object of $\calX$. Since $U_{\bigdot}$ is a hypercovering, $\theta$ is an effective epimorphism, so the \Cech nerve of $\theta$ is a colimit diagram and the desired result follows. Let us therefore assume that $n > 0$. Let $V_{\bigdot} = \cosk_{n-1} U_{\bigdot}$, and let
$f_{\bigdot}: U_{\bigdot} \rightarrow V_{\bigdot}$ be the adjunction map. 
For each $m \geq 0$, the map $f_{m}: U_{m} \rightarrow V_{m}$ is a composition of finitely many pullbacks of $f_{n}$. Since $U_{\bigdot}$ is a hypercovering, $f_n$ is an effective epimorphism, so each $f_m$ is also an effective epimorphism. We also observe that $f_m$ is an equivalence for $m < n$.

Let $W_{+}: \Nerve (\cDelta_{+} \times \cDelta)^{op} \rightarrow \calX$
be a \Cech nerve of $f_{\bigdot}$ (formed in the $\infty$-category $\calX_{\Delta}$ of simplicial objects of $\calX$). We observe that
$W_{+} | \Nerve ( \{ \emptyset \} \times \cDelta)^{op}$ can be identified with
$V_{\bigdot}$. Since $V_{\bigdot}$
is an $(n-1)$-coskeletal hypercovering of $\calX$, the inductive hypothesis implies that any colimit $|V_{\bigdot}|$ is a final object of $\calX$. The inclusion $\Nerve ( \{ \emptyset \} \times \cDelta)^{op}
\subseteq \Nerve( \cDelta_{+} \times \cDelta)^{op}$ is cofinal (being a product
of $\Nerve(\cDelta)^{op}$ with the inclusion of a final object into $\Nerve(\cDelta_{+})^{op}$ ), so
we may identify colimits of $W_{+}$ with colimits of $V_{\bigdot}$. It follows that
any colimit of $W_{+}$ is a final object of $\calX$.
We next observe that each of the augmented simplicial objects
$W_{+}| \Nerve( \cDelta_{+} \times \{ [m] \})^{op}$ is a \Cech nerve of $f_{m}$, and therefore a colimit diagram (since $f_m$ is an effective epimorphism). Applying Lemma \ref{longerwait}, we conclude that $W_{+}$ is a left Kan extension of the bisimplicial object
$W = W_{+}| \Nerve( \cDelta \times \cDelta)^{op}$. According to Lemma \ref{kan0}, we can
identify colimits of $W_{+}$ with colimits of $W$, so any colimit of $W$ is a final object of $\calX$.

Let $D_{\bigdot}: \Nerve(\cDelta^{op}) \rightarrow \calX$ be the simplicial object of $\calX$
obtained by composing $W$ with the diagonal map $\delta: \Nerve(\cDelta^{op}) \rightarrow \Nerve (\cDelta \times \cDelta)^{op}$. According to Lemma \ref{bball3}, $\delta$ is cofinal.
We may therefore identify colimits of $W$ with colimits of $D_{\bigdot}$, so that
any colimit $|D_{\bigdot}|$ of $D_{\bigdot}$ is a final object of $\calX$.

Let $U_{\bigdot}^{s} = U_{\bigdot} | \Nerve(\cDelta_{s}^{op})$, and let
$D_{\bigdot}^{s} = D_{\bigdot} | \Nerve(\cDelta_{s}^{op})$. We will prove that
$U_{\bigdot}^{s}$ is a retract of $D_{\bigdot}^{s}$ in the $\infty$-category of semisimplicial
objects of $\calX$. According to Lemma \ref{bball4}, we can identify colimits of
$D_{\bigdot}^{s}$ with colimits of $D_{\bigdot}$. It will follow that any colimit of
$U_{\bigdot}^{s}$ is a retract of a final object of $\calX$, and therefore itself final.
Applying Lemma \ref{bball4} again, we will conclude that any colimit of $U_{\bigdot}$ is
a final object of $\calX$, and the proof will be complete.

We observe that $D^{s}_{\bigdot}$ is the result of composing $W$ with the (opposite of the nerve of the) diagonal functor
$$\delta^{s}: \cDelta_{s} \rightarrow \cDelta \times \cDelta.$$
Similarly, the semisimplicial object $U^{s}_{\bigdot}$ is obtained from $W$ via the composition
$$ \epsilon: \cDelta_{s} \subseteq \cDelta \simeq \{ [0] \} \times \cDelta \subseteq
\cDelta \times \cDelta.$$
There is a obvious natural transformation of functors $\delta^{s} \rightarrow \epsilon$, which
yields a map of semisimplicial objects $\theta: U^{s}_{\bigdot} \rightarrow D^{s}_{\bigdot}$.
To complete the proof, it will suffice to show that there exists a map
$$ \theta': D^{s}_{\bigdot} \rightarrow U^{s}_{\bigdot}$$
such that $\theta' \circ \theta$ is homotopic to the identity on $U^{s}_{\bigdot}$.

According to Lemma \ref{bball5}, $U^{s}_{\bigdot}$ is $n$-coskeletal as a {\em semisimplicial} object of $\calX$. Let $D^{s}_{\leq n}$ and $U^{s}_{\leq n}$ denote
restrictions of $D^{s}_{\bigdot}$ and $U^{s}_{\bigdot}$ to $\Nerve(\cDelta_{s, \leq n}^{op})$, and $\theta_{\leq n}:
U^{s}_{\leq n} \rightarrow D^{s}_{\leq n}$ the morphism induced by $\theta$.
We have canonical homotopy equivalences
$$\bHom_{ \Fun(\Nerve(\cDelta^{op}_{s}), \calX)} ( D^{s}_{\bigdot}, U^{s}_{\bigdot})
\simeq \bHom_{ \Fun(\Nerve(\cDelta^{op}_{s, \leq n}),\calX)}( D^{s}_{\leq n}, U^{s}_{\leq n})$$
$$\bHom_{ \Fun(\Nerve(\cDelta^{op}_{s}),\calX)} ( U^{s}_{\bigdot}, U^{s}_{\bigdot})
\simeq \bHom_{ \Fun(\Nerve(\cDelta^{op}_{s, \leq n}),\calX)}( U^{s}_{\leq n}, U^{s}_{\leq n}).$$
It will therefore suffice to prove that there exists a map
$$ \theta'_{\leq n}: D^{s}_{\leq n} \rightarrow U^{s}_{\leq n}$$
such that $\theta'_{\leq n} \circ \theta_{\leq n}$ is homotopic to the identity on $U^{s}_{\leq n}$.

Consider the functors 
$$\overline{\delta}^{s}: \cDelta_{s, \leq n} \rightarrow \cDelta_{+} \times \cDelta$$
$$\overline{\epsilon}: \cDelta_{s, \leq n} \rightarrow \cDelta_{+} \times \cDelta$$
defined as follows:
$$ \overline{\delta}^{s}( [m]) = 
\begin{cases} ( \emptyset, [m] ) & \text{if } m < n \\
([n], [n] ) & \text{if } m = n \end{cases}$$
$$ \overline{\epsilon}( [m] ) = 
\begin{cases} ( \emptyset, [m] ) & \text{if } m < n \\
( [0] , [n] ) & \text{if } m = n. \end{cases}$$
We have a commutative diagram of natural transformations
$$ \xymatrix{ \overline{\delta}^{s} \ar[r] \ar[d] & \delta^{s} \ar[d] \\
\overline{\epsilon} \ar[r] & \epsilon }$$
which gives rise to a diagram
$$ \xymatrix{ \overline{D}^{s}_{\leq n} & D^{s}_{\leq n} \ar[l] \\
\overline{U}^{s}_{\leq n} \ar[u]^{\psi_{\leq n}} & U^{s}_{\leq n} \ar[u]^{\theta_{\leq n}} \ar[l] }$$
in the $\infty$-category $\Fun(\Nerve(\cDelta^{op}_{s, \leq n}),\calX)$. The
vertical arrows are equivalences. Consequently, it will suffice to produce
a (homotopy) left inverse to $\psi_{\leq n}$. 

For $m \geq 0$, let $V^{s}_{\leq m} = V_{\bigdot} | \cDelta_{s, \leq m}$. We can identify
$\overline{D}^{s}_{\leq n}$ and $\overline{U}^{s}_{\leq n}$ with objects
$X,Y \in \calX_{/V^{s}_{\leq n-1}}$, and $\psi_{\leq n}$ with a morphism
$f: X \rightarrow Y$. To complete the proof, it will suffice to produce a left
inverse to $f$ in the $\infty$-category $\calX_{/V^{s}_{\leq n-1}}$. 
We observe that, since $V_{\bigdot}$ is $(n-1)$-coskeletal, we have a diagram of trivial fibrations
$$ \calX_{/V_{n}} \leftarrow \calX_{/ V^{s}_{\leq n} } \rightarrow \calX_{ / V^{s}_{\leq n-1}}.$$
Using this diagram (and the construction of $W$), we conclude that $Y$ can be identified
with a product of $(n+1)$ copies of $X$ in $\calX_{/V^{s}_{\leq n-1}}$, and that
$f$ can be identified with the identity map. The existence of a left homotopy inverse to $f$ is now obvious (choose any of the $(n+1)$-projections from $Y$ onto $X$).
\end{proof}

\begin{lemma}\label{bball2}
Let $\calX$ be an $\infty$-topos, and let $f_{\bigdot}: U_{\bigdot} \rightarrow V_{\bigdot}$ be a natural transformation between simplicial objects of $\calX$. Suppose that, for each
$k \leq n$, the map $f_{k}: U_k \rightarrow V_{k}$ is an equivalence. Then the induced map
$|f_{\bigdot}| : | U_{\bigdot} | \rightarrow | V_{\bigdot} |$ of colimits is $n$-connective.
\end{lemma}

\begin{proof}
Choose a left exact localization functor $L: \calP(\calC) \rightarrow \calX$. Without loss of generality, we may suppose that $f_{\bigdot} = L \circ \overline{f}_{\bigdot}$, where
$\overline{f}_{\bigdot}: \overline{U}_{\bigdot} \rightarrow \overline{V}_{\bigdot}$ is
a transformation between simplicial objects of $\calP(\calC)$, where $\overline{f}_{k}$ is an equivalence for $k \leq n$. Since $L$ preserves colimits and $n$-connectivity (Proposition \ref{inftychange}), it will suffice to prove that $|f_{\bigdot}|$ is $n$-connective. Using Propositions \ref{goober} and \ref{compattrunc}, we see that $|f_{\bigdot}|$ is $n$-connective if and only if, for each object $C \in \calC$, the induced morphism in $\SSet$ is $n$-connective. In other words, we may assume without loss of generality that $\calX = \SSet$.

According to Proposition \ref{gumby444}, we may assume that
$f_{\bigdot}$ is obtained by taking the simplicial nerve of a map $f'_{\bigdot}: U'_{\bigdot} \rightarrow V'_{\bigdot}$ between simplicial objects in the ordinary category $\Kan$.
Without loss of generality, we may suppose that
 $U'_{\bigdot}$ and $V'_{\bigdot}$ are projectively cofibrant (as diagrams in the model category $\sSet$). According to Theorem \ref{colimcomparee}, it will suffice to prove that the induced map from the (homotopy) colimit of $U'_{\bigdot}$ to the (homotopy) colimit of $V'_{\bigdot}$ has $n$-connective homotopy fibers, which follows from classical homotopy theory.
\end{proof}

\begin{lemma}\label{fierminus}
Let $\calX$ be an $\infty$-topos, let $U_{\bigdot}$ be a hypercovering of $\calX$.
Then the colimit $| U_{\bigdot} |$ is $\infty$-connective.
\end{lemma}

\begin{proof}
We will prove that $\theta$ is $n$-connective for every $n \geq 0$. Let
$V_{\bigdot} = \cosk_{n+1} U_{\bigdot}$, and let
$u: U_{\bigdot} \rightarrow V_{\bigdot}$ be the adjunction map.
Lemma \ref{bball2} asserts that the induced map $|U_{\bigdot}| \rightarrow |V_{\bigdot}|$
is $n$-connective, and Lemma \ref{bball1} asserts that $| V_{\bigdot} |$ is a final object of $\calX$. It follows that $| U_{\bigdot} | \in \calX$ is $n$-connective, as desired.
\end{proof}

The preceding results lead to an easy characterization of the class of hypercomplete $\infty$-topoi:

\begin{theorem}\label{surp}
Let $\calX$ be an $\infty$-topos. The following conditions are equivalent:
\begin{itemize}
\item[$(1)$] For every $X \in \calX$, every hypercovering $U_{\bigdot}$ of
$\calX_{/X}$ is effective.
\item[$(2)$] The $\infty$-topos $\calX$ is hypercomplete.
\end{itemize}
\end{theorem}

\begin{proof}
Suppose that $(1)$ is satisfied. Let $f: U \rightarrow X$ be an $\infty$-connective morphism in $\calX$, and let $f_{\bigdot}$ be the constant simplicial object of $\calX_{/X}$ with value $f$.
According to Lemma \ref{fier0}, $f$ is a hypercovering of $\calX_{/X}$. Invoking $(1)$, we conclude that $f \simeq | f_{\bigdot} |$ is a final object of $\calX_{/X}$; in other words, $f$ is an equivalence. This proves that $(1) \Rightarrow (2)$.

Conversely, suppose that $\calX$ is hypercomplete. Let $X \in \calX$ be an object and
$U_{\bigdot}$ a hypercovering of $\calX_{/X}$. Then Lemma \ref{fierminus} implies that
$| U_{\bigdot} |$ is an $\infty$-connective object of $\calX_{/X}$. Since $\calX$ is hypercomplete, we conclude that $|U_{\bigdot}|$ is a final object of $\calX_{/X}$, so that $U_{\bigdot}$ is effective.
\end{proof}

\begin{corollary}[Dugger-Hollander-Isaksen \cite{hollander}, To\"{e}n-Vezzosi \cite{toen}]\label{charhyp}
Let $\calX$ be an $\infty$-topos. For each $X \in \calX$ and each hypercovering
$U_{\bigdot}$ of $\calX_{/X}$, let $|U_{\bigdot}|$ be the associated morphism of $\calX$
$(${}which has target $X${}$)$. 
Let $S$ denote the collection of all such morphisms $|U_{\bigdot}|$. Then
$\calX^{\hyp} = S^{-1} \calX$. In other words, an object of
$\calX$ is hypercomplete if and only if it is $S$-local.
\end{corollary}

\begin{remark}
One can generalize Corollary \ref{charhyp} as follows: let $L: \calX \rightarrow \calY$ be an arbitrary left exact localization of $\infty$-topoi, and let $S$ be the collection of all morphisms of the form
$|U_{\bigdot}|$, where $U_{\bigdot}$ is a simplicial object of $\calX_{/X}$ such that
$L \circ U_{\bigdot}$ is an effective hypercovering of $\calY_{/LX}$. 
Then $L$ induces an equivalence $S^{-1} \calX \rightarrow \calY$.

It follows that every $\infty$-topos can be obtained by starting with an $\infty$-category of presheaves $\calP(\calC)$, selecting a collection of augmented simplicial objects
$U^{+}_{\bigdot}$, and inverting the corresponding
maps $|U_{\bigdot}| \rightarrow U_{-1}$. The specification of the desired class of augmented simplicial objects can be viewed as a kind of ``generalized topology'' on $\calC$, in which one specifies not only the covering sieves but also the collection of hypercoverings which are to become effective after localization. It seems plausible that this notion of topology can be described more directly in terms of the $\infty$-category $\calC$, but we will not pursue the matter further.
\end{remark}

\subsection{Descent versus Hyperdescent}\label{versus}

Let $X$ be a topological space, and let $\calU(X)$ denote the category of open subsets of $X$.\index{not}{Ucal(X)@$\calU(X)$} The category $\calU(X)$ is equipped with a Grothendieck topology in which the covering sieves
on $U$ are those sieves $\{ U_{\alpha} \subseteq U \}$ such that $U = \bigcup_{\alpha} U_{\alpha}$. We may therefore consider the $\infty$-topos $\Shv(\Nerve(\calU(X)))$, which we will 
call the $\infty$-topos of {\it sheaves on $X$}\index{gen}{$\infty$-topos!of sheaves on a topological space} and denote by $\Shv(X)$\index{not}{ShvX@$\Shv(X)$}
In \S \ref{hyperstacks} we discussed an alternative theory of sheaves on $X$, which can be obtained either through Jardine's local model structure on the category of simplicial presheaves or by passing to the hypercompletion $\Shv(X)^{\hyp}$ of $\Shv(X)$. According to Theorem \ref{surp}, $\Shv(X)^{\hyp}$ is distinguished from $\Shv(X)$ in that objects of $\Shv(X)^{\hyp}$ are required to satisfy a descent condition for arbitrary hypercoverings of $X$, while objects of $\Shv(X)$ are required to satisfy a descent condition only for ordinary coverings.

\begin{warning}
The notation $\Shv(X)$ will always represent the $\infty$-category of $\SSet$-valued
sheaves on $X$, rather than the ordinary category of set-valued sheaves. If we need to
indicate the latter, we will denote it by $\Shv_{\Set}(X)$.
\end{warning}

The $\infty$-topos $\Shv(X)^{\hyp}$ seems to have received more attention than $\Shv(X)$ in the literature (though there is some discussion of $\Shv(X)$ in \cite{hollander} and \cite{toen}).
We would like to make the case that for most purposes, $\Shv(X)$ has better properties.
A large part of \S \ref{chap7} will be devoted to justifying some of the claims made below.

\begin{itemize}

\item[$(1)$] In \S \ref{nlocalic}, we saw that the construction
$$ X \mapsto \Shv(X)$$
could be interpreted as a {\em right} adjoint to the functor which associates to every
$\infty$-topos $\calY$ the underlying locale of subobjects of the final object of $\calY$.
In other words, $\Shv(X)$ occupies a universal position among $\infty$-topoi which are related to the original space $X$.

\item[$(2)$] Suppose given a Cartesian square
$$ \xymatrix{ X' \ar[r]^{\psi'} \ar[d]^{\pi'} & X \ar[d]^{\pi} \\
S' \ar[r]^{\psi} & S }$$
in the category of locally compact topological spaces.
In classical sheaf theory, there is a {\it base change} transformation
$$\psi^{\ast} \pi_{\ast} \rightarrow \pi'_{\ast} \psi'^{\ast}$$
of functors between the derived categories of (left-bounded) complexes of (abelian) sheaves
on $X$ and on $S'$. The proper base change theorem
asserts that this transformation is an equivalence whenever the
map $\pi$ is proper.\index{gen}{proper!base change theorem}

The functors $\psi^{\ast}$, ${\psi'}^{\ast}$, $\pi_{\ast}$, and $\pi'_{\ast}$
can be defined on the $\infty$-topoi $\Shv(X), \Shv(X'), \Shv(S),$ and $\Shv(S')$, and on their hypercompletions. Moreover, one has a base change map
$$\psi^{\ast} \pi_{\ast} \rightarrow \pi'^{\ast} {\psi'}^{\ast}$$
in this nonabelian situation as well. 

It is natural to ask if the
base change transformation is an equivalence when $\pi$ is proper.
It turns out that this is {\em false} if we work with hypercomplete $\infty$-topoi.
Let us sketch a counterexample:

\begin{counterexample}\label{thrust}\index{gen}{Hilbert cube}
Let $Q$ denote the Hilbert cube $[0,1] \times [0,1] \times
\ldots$. For each $i$, we let $Q_i \simeq Q$ denote ``all but the
first $i$'' factors of $Q$, so that $Q = [0,1]^i \times Q_i$.

We construct a sheaf of spaces $\calF$ on $X = Q \times [0,1]$ as follows.
Begin with the empty stack. Adjoin to it two sections, defined
over the open sets $[0,1) \times Q_1 \times [0,1)$ and $(0,1]
\times Q_1 \times [0,1)$. These sections both restrict to give
sections of $\calF$ over the open set $(0,1) \times Q_1 \times
[0,1)$. We next adjoin paths between these sections, defined over
the smaller open sets $(0,1) \times [0,1) \times Q_2 \times
[0,\frac{1}{2})$ and $(0,1) \times (0,1] \times Q_2 \times [0,
\frac{1}{2})$. These paths are both defined on the smaller open
set $(0,1) \times (0,1) \times Q_2 \times [0, \frac{1}{2})$, so we
next adjoin two homotopies between these paths over the open sets
$(0,1) \times (0,1) \times [0,1) \times Q_3 \times [0,
\frac{1}{3})$ and $(0,1) \times (0,1) \times (0,1] \times Q_3
\times [0, \frac{1}{3})$. Continuing in this way, we obtain a
sheaf $\calF$. On the closed subset $Q \times \{0\} \subset X$,
the sheaf $\calF$ is $\infty$-connective by construction, and
therefore its hypercompletion admits a global section.
However, the hypercompletion of $\calF$ does not admit a
global section in any neighborhood of $Q \times \{0\}$, since such
a neighborhood must contain $Q \times [ 0, \frac{1}{n})$ for $n
\gg 0$ and the higher homotopies required for the construction of
a section are eventually not globally defined.
\end{counterexample}

However, in the case where $\pi$ is a proper map, the base-change map
$$\psi^{\ast} \pi_{\ast} \rightarrow \pi'_{\ast} \psi'^{\ast}$$
{\em is} an equivalence of functors from $\Shv(X)$ to $\Shv(S')$. One may regard this
fact as a nonabelian generalization of the classical proper-base change theorem.
We refer the reader to \S \ref{chap7sec3} for a precise statement and proof.

\begin{remark}
A similar issue arises in classical sheaf theory if one chooses to
work with unbounded complexes. In \cite{spaltenstein},
Spaltenstein defines a derived category of unbounded complexes of
sheaves on $X$, where $X$ is a topological space. His definition
forces all quasi-isomorphisms to become invertible, which is
analogous to procedure of obtaining $\calX^{\hyp}$ from $\calX$ by
inverting the $\infty$-connective morphisms. Spaltenstein's work
shows that one can extend the {\it definitions} of all of the
basic objects and functors. However, it turns out that the {\it
theorems} do not all extend: in particular, one does not have the
proper base change theorem in Spaltenstein's setting
(Counterexample \ref{thrust} can be adapted to the setting
of complexes of abelian sheaves). The problem may be rectified by
imposing weaker descent conditions, which do not invert all
quasi-isomorphisms; we will give a more detailed discussion in \cite{DAG}.
\end{remark}

\item[$(3)$] The $\infty$-topos $\Shv(X)$ often has better finiteness properties
than $\Shv(X)^{\hyp}$. Recall that a topological space $X$ is {\it coherent}
if the collection of compact open subsets of $X$ is stable under finite intersections, and forms a basis for the topology of $X$.\index{gen}{coherent topological space}

\begin{proposition}\label{cohcomp}
Let $X$ be a coherent topological space. Then the $\infty$-category
$\Shv(X)$ is compactly generated: that is, $\Shv(X)$ is generated under filtered colimits
by its compact objects.\index{gen}{compactly generated}
\end{proposition}

\begin{proof}
Let $\calU_c(X)$ be the partially ordered set of {\em compact} open subsets of $X$, let
$\calP_c(X) = \calP( \Nerve(\calU_c(X)) )$, and let $\Shv_c(X)$ be the full subcategory of 
$\calP_c(X)$ spanned by those presheaves $\calF$ with the following properties:
\begin{itemize}
\item[$(1)$] The object $\calF(\emptyset) \in \calC$ is final.
\item[$(2)$] For every pair of compact open sets $U, V \subseteq X$, the associated diagram
$$ \xymatrix{ \calF( U \cap V) \ar[r] \ar[d] & \calF(U) \ar[d] \\
\calF(V) \ar[r] & \calF(U \cup V) }$$
is a pullback.
\end{itemize}

In \S \ref{cohthm}, we will prove that the restriction functor 
$\Shv(X) \rightarrow \Shv_c(X)$ is an equivalence of $\infty$-categories (Theorem \ref{surm}). 
It will therefore suffice to prove that $\Shv_c(X)$ is compactly generated. 

Using Lemmas \ref{stur1}, \ref{stur2}, and \ref{stur3}, we conclude that $\Shv_c(X)$ is an accessible localization of $\calP_c(X)$. 
Let $X$ be a compact object of $\calP_{c}(X)$. We observe that $X$ and $LX$ co-represent the same functor on $\Shv_{c}(X)$. Proposition \ref{frent} implies that the subcategory $\Shv_c(X) \subseteq \calP_c(X)$ is stable under filtered colimits in $\calP_c(X)$. It follows that
$LX$ is a compact object of $\Shv_0(X)$. Since $\calP_{c}(X)$ is generated under filtered colimits by its compact objects (Proposition \ref{precst}), we conclude that $\Shv_{c}(X)$ has the same property.
\end{proof}

\begin{remark}
In the situation of Proposition \ref{cohcomp}, we can give an explicit description of the class of compact objects of $\Shv(X)$. Namely, they are precisely those sheaves $\calF$ whose stalks are compact objects of $\SSet$, and which are {\em locally constant} along a suitable stratification of $X$. In other words, we may interpret Proposition \ref{cohcomp} as asserting that there is a good theory of {\em constructible} sheaves on $X$.
\end{remark}

It is not possible to replace to $\Shv(X)$ by $\Shv(X)^{\hyp}$ in the statement of Proposition \ref{cohcomp}.

\begin{counterexample}\label{trust}
Let $S = \{ x,y,z \}$ be a topological space consisting of three points, with topology generated
by the open subsets $S^{+} = \{x,y\} \subset S$ and $S^{-} = \{x,z\} \subset S$.
Let $X = S \times S \times \ldots$ be a product of infinitely many copies of $S$.
Then $X$ is a coherent topological space. We will show that the global sections functor
$\Gamma: \Shv(X)^{\hyp} \rightarrow \SSet$ does not commute with filtered colimits, so that
the final object of $\Shv(X)^{\hyp}$ is not compact. A more elaborate version of the same argument shows that $\Shv(X)^{\hyp}$ contains no compact objects other than its initial object.

To show that $\Gamma$ does not commute with filtered colimits, we use a variant on the
construction of Counterexample \ref{thrust}. We define a sequence of
sheaves $$\calF_0 \rightarrow \calF_1 \rightarrow \ldots $$
as follows. Let $\calF_0$ be generated by sections
$$ \eta^0_{+} \in \calF(S^{+} \times S \times \ldots) $$
$$ \eta^0_{-} \in \calF(S^{-} \times S \times \ldots). $$
Let $\calF_1$ be the sheaf obtained from $\calF_0$ by adjoining
paths 
$$ \eta^{1}_{+}: \Delta^1 \rightarrow \calF( \{x\} \times S^{+} \times S \times \ldots )$$
$$ \eta^{1}_{-}: \Delta^1 \rightarrow \calF( \{x\} \times S^{-} \times S \times \ldots )$$
from $\eta^0_{+}$ to $\eta^{0}_{-}$.
Similarly, let $\calF_2$ be obtained from $\calF_1$ by adjoining homotopies
$$ \eta^2_{+}: (\Delta^1)^2 \rightarrow \calF( \{x\} \times \{x\} \times S^{+} \times S \times \ldots )$$
$$ \eta^2_{-}: (\Delta^1)^2 \rightarrow \calF( \{x \} \times \{x\} \times S^{-} \times S \times \ldots ),$$
from $\eta^1_{+}$ to $\eta^{1}_{-}$. Continuing this procedure, we obtain a
sequence of sheaves $$\calF_0 \rightarrow \calF_1 \rightarrow \calF_2 \rightarrow \ldots $$
whose colimit $\calF_{\infty} \in \Shv(X)^{\hyp}$ admits a section (since we allow descent
with respect to hypercoverings). However, none of the individual sheaves $\calF_n$ admits a global section.
\end{counterexample}

\begin{remark}
The analogue of Proposition \ref{cohcomp} fails, in general, if we replace the coherent topological space $X$ by a coherent topos. For example, we cannot take $X$ to be the topos of \'{e}tale sheaves on an algebraic variety. However, it turns out that the analogue Proposition \ref{cohcomp} {\em is} true for the topos of {\em Nisnevich} sheaves on an algebraic variety; we refer the reader to \cite{DAG} for details.
\end{remark}

\begin{remark}\label{notenough}
A {\it point}\index{gen}{point!of an $\infty$-topos} of an $\infty$-topos $\calX$ is a geometric morphism $p_{\ast}: \SSet \rightarrow \calX$, where $\SSet$ denotes the $\infty$-category of spaces (which is a final object of $\RGeom$, in virtue of Proposition \ref{spacefinall}). We say that $\calX$ has {\it enough points}\index{gen}{enough points} if, for every morphism $f: X \rightarrow Y$ in $\calX$ having the property that $p^{\ast}(f)$ is an equivalence for {\em every} point $p$ of $\calX$, $f$ is itself an equivalence in $\calX$. If $f$ is $\infty$-connective, then every stalk $p^{\ast}(f)$ is $\infty$-connective, hence an equivalence by Whitehead's theorem. Consequently, if $\calX$ has enough points, then it is hypercomplete.

In classical topos theory, Deligne's version of the G\"{o}del completeness theorem (see \cite{where}) asserts that every coherent topos has enough points.  Counterexample \ref{trust} shows that there exist coherent topological spaces with
$\Shv(X)^{\hyp} \neq \Shv(X)$, so that $\Shv(X)$ does not necessarily have enough points. 
Consequently, Deligne's theorem does not hold in the $\infty$-categorical context.
\end{remark}

\item[$(4)$] Let $k$ be a field, and let $\calC$ denote the category of chain complexes of $k$-vector spaces. Via the Dold-Kan correspondence we may regard $\calC$ as a simplicial category. We let $\Mod(k) = \sNerve(\calC)$ denote the simplicial nerve. We will refer to $\Mod(k)$ as the {\it $\infty$-category of $k$-modules}; it is a presentable $\infty$-category which we will discuss at greater length in \cite{DAG}.\index{not}{Modk@$\Mod(k)$}

Let $X$ be a compact topological space, and choose a functorial injective resolution
$$ \calF \rightarrow I^0( \calF) \rightarrow I^1( \calF ) \rightarrow \ldots $$
on the category of sheaves $\calF$ of $k$-vector spaces on $X$. For every open subset
$U$ on $X$, we let $k_U$ denote the constant sheaf on $U$ with value $k$, extended by zero to $X$. Let $\HH^{BM}(U) = \Gamma(X, I^{\bigdot}( k_U) )^{\vee}$, the {\it dual} of the complex
of global sections of the injective resolution $I^{\bigdot}(k_U)$. Then $\HH^{BM}(U)$ is a 
complex of $k$-vector spaces, whose homologies are precisely the Borel-Moore homology of $U$ with coefficients in $k$ (in other words, they are the dual spaces of the compactly supported cohomology groups of $U$). The assignment
$$ U \mapsto \HH^{BM}(U)$$
determines a presheaf on $X$ with values in the $\infty$-category $\Mod(k)$.

In view of the existence of excision exact sequences for Borel-Moore homology, it is natural to suppose that $\HH^{BM}(U)$ is actually a {\it sheaf} on $X$ with values in $\Mod(k)$. This is true provided that the notion of ``sheaf'' is suitably interpreted: namely, 
$\HH^{BM}$ extends (in an essentially unique fashion) to a colimit-preserving functor
$$ \phi: \Shv(X) \rightarrow \Mod(k)^{op}.$$
(In other words, the functor $U \mapsto \HH^{BM}(U)$ determines a $\Mod(k)$-valued sheaf
on $X$ in the sense of Definition \ref{valsheaf}.)
However, the sheaf $\HH^{BM}$ is not necessarily hypercomplete, in the sense that $\phi$ does not necessarily factor through $\Shv(X)^{\hyp}$.

\begin{counterexample}\label{spacerk}
There exists a compact Hausdorff space $X$ and a hypercovering $U_{\bigdot}$ of $X$ such that the natural map $\HH^{BM}(X) \rightarrow \varprojlim \HH^{BM}(U_{\bigdot})$ is not an equivalence. Let
$X$ be the Hilbert cube $Q = [0,1] \times [0,1] \times \ldots$ (more generally, we could take $X$ to be any nonempty Hilbert cube manifold). It is proven in \cite{chapman} that every point of $X$ has arbitrarily small neighborhoods which are homeomorphic to $Q \times [0,1)$. Consequently, 
there exists a hypercovering $U_{\bigdot}$ of $X$, where each $U_n$ is a disjoint union
of open subsets of $X$ homeomorphic to $Q \times [0,1)$. The Borel-Moore homology of
every $U_n$ vanishes; consequently, $\varprojlim \HH^{BM}(U_{\bigdot})$ is zero. However,
the (degree zero) Borel-Moore homology of $X$ itself does not vanish, since $X$ is nonempty and compact.
\end{counterexample}

Borel-Moore homology is a very useful tool in the study of a locally compact space $X$, and its descent properties (in other words, the existence of various Mayer-Vietoris sequences) is very naturally encoded in the statement that $\HH^{BM}$ is a {\it $k$-module in the $\infty$-topos $\Shv(X)$} (in other words, a sheaf on $X$ with values in $\Mod(k)$); however, this $k$-module generally does not lie in $\Shv(X)^{\hyp}$. We see from this example that non-hypercomplete sheaves (with values in $\Mod(k)$, in this case) on $X$ often arise naturally in the study of infinite-dimensional spaces.

\item[$(5)$] Let $X$ be a topological space, and $f: \Shv(X) \rightarrow \Shv(\ast) \simeq \SSet$ the geometric morphism induced by the projection $X \rightarrow \ast$. Let $K$ be a Kan complex, regarded as an object of $\SSet$. Then $\pi_0 f_{\ast} f^{\ast} K$
is a natural definition of the sheaf cohomology of
$X$ with coefficients in $K$. If $X$ is paracompact, then the cohomology set defined above is naturally isomorphic to the set $[X,|K|]$ of homotopy classes of maps from $X$ into the geometric realization $|K|$; we will give a proof of this statement in \S \ref{paracompactness}.
The analogous statement fails if we replace $\Shv(X)$ by $\Shv(X)^{\hyp}$.

\item[$(6)$] Let $X$ be a topological space. Combining Remark \ref{pointeddesc} with Proposition \ref{suga}, we deduce that $\Shv(X)^{\hyp}$ has enough points, and that $\Shv(X)^{\hyp} = \Shv(X)$ if and only if $\Shv(X)$ has enough points. The possible failure of Whitehead's theorem in $\Shv(X)$ may be viewed either as a bug or a feature. The existence of enough points for $\Shv(X)$ is extremely convenient; it allows us to reduce many statements about the $\infty$-topos $\Shv(X)$
to statements about the $\infty$-topos $\SSet$ of spaces, where we can apply classical homotopy theory. On the other hand, if $\Shv(X)$ does {\em not} have enough points, then there is the possibility that it detects certain global phenomena 
which cannot be properly understood by restricting to points. Let us consider an example from geometric topology. A map $f: X \rightarrow Y$ of compact metric spaces is called {\it cell-like} if each fiber $X_{y} = X \times_{Y} \{y\}$ has trivial shape (see \cite{cellmap})\index{gen}{cell-like}. This notion has good formal properties provided that we restrict our attention to metric spaces which are {\em absolute neighborhood retracts}. In the general case, the theory of cell-like maps can be badly behaved: for example, a composition of cell-like maps need not be cell-like. 

The language of $\infty$-topoi provides a convenient formalism for discussing the problem. 
In \S \ref{celluj}, we will introduce the notion of a {\em cell-like} morphism $p_{\ast}: \calX \rightarrow \calY$ between $\infty$-topoi. By definition, $p_{\ast}$ is cell-like if it is proper and if the unit map $u: \calF \rightarrow p_{\ast} p^{\ast} \calF$ is an equivalence for each $\calF \in \calY$. A cell-like map $p: X \rightarrow Y$ of compact metric spaces {\em need not} give rise to a cell-like morphism $p_{\ast}: \Shv(X) \rightarrow \Shv(Y)$. The hypothesis
that each fiber $X_{y}$ has trivial shape ensures that the unit $u: \calF \rightarrow p_{\ast} p^{\ast} \calF$ is an equivalence after passing to stalks at each point $y \in Y$. This implies only that
$u$ is $\infty$-connective, and in general $u$ need not be an equivalence.

\begin{remark}
It is tempting to try to evade the problem described above by working instead with the hypercomplete $\infty$-topoi $\Shv(X)^{\hyp}$ and $\Shv(Y)^{\hyp}$. In this case, we {\em can} test whether or not $u: \calF \rightarrow p_{\ast} p^{\ast} \calF$ is an equivalence by passing to stalks. However, since the proper base change theorem does not hold in the hypercomplete context, the stalk $(p_{\ast} p^{\ast} \calF)_{y}$ is not generally equivalent to the global sections of
$p^{\ast} \calF | X_{y}$. Thus, we still encounter difficulties if we want to deduce global consequences from information about the individual fibers $X_{y}$. 
\end{remark}

\item[$(7)$] The counterexamples described in this section have one feature in common: the underlying space $X$ is infinite-dimensional. In fact, this is necessary: if the space $X$ is finite-dimensional (in a suitable sense), then the $\infty$-topos $\Shv(X)$ is hypercomplete (Corollary \ref{fdfd}). This finite-dimensionality condition on $X$ is satisfied in many of the situations to which the theory of simplicial presheaves is commonly applied, such as the Nisnevich
topology on a scheme of finite Krull dimension.

\end{itemize}

\chapter{Higher Topos Theory in Topology}\label{chap7}

\setcounter{theorem}{0}
\setcounter{subsection}{0}

In this chapter, we will sketch three applications of the theory of $\infty$-topoi to the study classical topology. We begin in \S \ref{paracompactness}, by showing that if $X$ is a paracompact topological space, then the $\infty$-topos $\Shv(X)$ of sheaves on $X$ can be interpreted as a homotopy theory of topological spaces $Y$ equipped with a map to $X$. We will deduce, as an application, that if $p_{\ast}: \Shv(X) \rightarrow \Shv(\ast)$ is the geometric morphism induced by the projection $X \rightarrow \ast$, then the composition $p_{\ast} p^{\ast}$ is equivalent to the functor $$K \mapsto K^X$$ from (compactly generated) topological spaces to itself.

Our  second application is to the dimension theory of topological spaces. There are many different notions of {\it dimension} for a topological space $X$, including the notion of covering dimension (when $X$ is paracompact), Krull dimension (when $X$ is Noetherian), and cohomological dimension. We will define the {\it homotopy dimension} of an $\infty$-topos $\calX$, which specializes to the covering dimension when $\calX = \Shv(X)$ for a paracompact space $X$, and is closely related to both cohomological dimension and Krull dimension. We will show that any $\infty$-topos which is (locally) finite-dimensional is hypercomplete, thereby justifying assertion $(7)$ of \S \ref{versus}. We will conclude by proving a bound on the homotopy dimension of $\Shv(X)$ where $X$ is a {\it Heyting space} (see \S \ref{heyt} for a definition); this may be regarded as a generalization of Grothendieck's vanishing theorem, which applies to non-abelian cohomology and to (certain) non-Noetherian spaces $X$.

Our third application is a generalization of the {\em proper base change theorem.}
Suppose given a Cartesian diagram
$$ \xymatrix{ X' \ar[r]^{p'} \ar[d]^{q'} & X \ar[d]^{q} \\
Y' \ar[r]^{p} & Y}$$
of locally compact topological spaces.
There is a natural transformation
$$ \eta: p^{\ast} q_{\ast} \rightarrow q'_{\ast} {p'}^{\ast}$$
of functors from the derived category of abelian sheaves on $X$ to the derived category of abelian sheaves on $Y'$. The proper base change theorem asserts that $\eta$ is an isomorphism whenever $q$ is a proper map. In \S \ref{chap7sec3}, we will generalize this statement to allow nonabelian coefficient systems. To give the proof, we will develop a theory of {\it proper morphisms} between $\infty$-topoi, which is of some interest in itself.

\section{Paracompact Spaces}\label{paracompactness}

\setcounter{theorem}{0}

Let $X$ be a topological space and $G$ an abelian group. There are
many different definitions for the cohomology group $\HH^n(X;G)$;
we will single out three of them for discussion here. First of
all, we have the singular cohomology groups $\HH^n_{\text{sing}}(X;G)$,
which are defined to be cohomology of a chain complex of $G$-valued
singular cochains on $X$. An alternative is to regard $\HH^n( \bigdot,
G)$ as a representable functor on the homotopy category of
topological spaces, so that $\HH^n_{\text{rep}}(X;G)$ can be identified with the
set of homotopy classes of maps from $X$ into an Eilenberg-MacLane
space $K(G,n)$. A third possibility is to use the sheaf cohomology
$\HH^n_{\text{sheaf}}(X; \underline{G})$ of $X$ with coefficients
in the constant sheaf $\underline{G}$ on $X$.

If $X$ is a sufficiently nice space (for example, a CW complex),
then these three definitions give the same result. In general, however,
all three give different answers. The singular cohomology of $X$ is defined using continuous maps from $\Delta^k$ into $X$, and is useful only when there is a good supply of such maps.
Similarly, the cohomology group $\HH^{n}_{\text{rep}}(X;G)$ is defined using continuous maps
from $X$ to a simplicial complex, and is useful only when there is a good supply of
real-valued functions on $X$. However, the sheaf cohomology of $X$ seems to be
a good invariant for arbitrary spaces: it has excellent formal
properties and gives sensible answers in situations where the other definitions break down (such as the \'{e}tale topology of algebraic varieties).

We will take the position that the sheaf cohomology of a
space $X$ is the correct answer in all cases. It is then natural to ask
for conditions under which the other definitions of cohomology
give the same answer. We should expect this to be true for
singular cohomology when there are many continuous functions {\em
into $X$}, and for Eilenberg-MacLane cohomology when there are
many continuous functions {\em out of} $X$. It seems that the
latter class of spaces is much larger than the former: it
includes, for example, all paracompact spaces, and consequently
for paracompact spaces one can show that the sheaf cohomology
$\HH^n_{\text{sheaf}}(X;G)$ coincides with the Eilenberg-MacLane
cohomology $\HH^n_{\text{rep}}(X;G)$. Our goal in this section is to prove a generalization of the preceding statement to the setting of nonabelian cohomology (Theorem \ref{nice} below; see also Theorem \ref{main} for the case where the coefficient system $G$ it not assumed to be constant).

As we saw in \S \ref{versus}, we can associate to every topological space $X$ an $\infty$-topos $\Shv(X)$ of sheaves (of spaces) on $X$. Moreover, given a continuous map $p: X \rightarrow Y$ of topological spaces, $p^{-1}$ induces a map from the category of open subsets of $Y$ to the category of open subsets of $X$. Composition with $p^{-1}$ induces a geometric morphism $p_{\ast}: \Shv(X) \rightarrow \Shv(Y)$. 

Fix now a topological space $X$ and let $p: X \rightarrow \ast$ denote the projection from $X$ to a point. Let $K$ be a Kan complex, which we may identify with an object of $\SSet \simeq \Shv(\ast)$. 
Then $p^{\ast} K \in \Shv(X)$ may be regarded as the constant sheaf on $X$ having
value $K$, and $p_{\ast} p^{\ast} K \in \SSet$ as the space of global sections of $p^{\ast} K$. 
Let $|K|$ denote the geometric realization of $K$ (a topological space), and let
$[ X, |K| ]$ denote the {\em set} of homotopy classes of maps from $X$ into $|K|$. The main goal of this section is to prove the following:

\begin{theorem}\label{nice}
If $X$ is paracompact, then there is a canonical bijection
$$ \phi: [X, |K| ] \rightarrow \pi_0(p_{\ast} p^{\ast} K).$$
\end{theorem}

\begin{remark}
In fact, the map $\phi$ exists without the assumption that $X$ is paracompact: the construction in general can be formally reduced to the paracompact case, since the universal example $X = |K|$ is paracompact. However, in the case where $X$ is not paracompact, the map $\phi$ is not necessarily bijective.
\end{remark}

Our first step in proving Theorem \ref{nice} is to realize the space of maps
from $X$ into $|K|$ as a mapping space in an appropriate simplicial category
of {\it spaces over $X$}. In \S \ref{para1}, we define this category
and endow it with a (simplicial) model structure. We may therefore extract an underlying $\infty$-category $\Nerve(\Top^{\degree}_{/X})$. 

Our next goal is to construct an equivalence between $\Nerve(\Top_{/X}^{\degree})$ and the $\infty$-topos $\Shv(X)$ of sheaves of spaces on $X$ (a very similar comparison result has been obtained by 
To\"{e}n; see \cite{toen2}). To prove this, we will attempt to realize $\Nerve(\Top_{/X}^{\degree})$ as a localization of a certain $\infty$-category of presheaves. We will give an explicit description of the relevant localization in \S \ref{para2}, and show that it is equivalent to $\Nerve(\Top_{/X}^{\degree})$ in \S \ref{para3}. 
In \S \ref{dooky}, we will deduce Theorem \ref{nice} as a corollary of this more general comparison result. We conclude with \S \ref{shapesec}, in which we apply our results to obtain a reformulation of classical shape theory in the language of $\infty$-topoi.

\subsection{Some Point-Set Topology}

Let $X$ be a paracompact topological space. In order to prove
Theorem \ref{nice}, we will need to understand the homotopy
theory of presheaves on $X$. We then encounter the following
technical obstacle: an open subset of a paracompact space need not
be paracompact. Because we wish to deal only with paracompact
spaces, it will be convenient to restrict our attention to
presheaves which are defined only with respect to a particular
basis $\calB$ for $X$ consisting of paracompact open sets. The existence of a well-behaved basis is guaranteed by the following result:

\begin{proposition}\label{gooffy}
Let $X$ be a paracompact topological space and $U$ an open subset of $X$.
The following conditions are equivalent:
\begin{itemize}
\item[$(i)$] There exists a continuous function $f: X \rightarrow [0,1]$ such that
$U = \{ x \in X: f(x) > 0 \}$.
\item[$(ii)$] There exists a sequence of closed subsets 
$\{ K_{n} \subseteq X \}_{n \geq 0}$ such that each $K_{n+1}$ contains an open neighborhood of
$K_{n}$, and $U = \bigcup_{n \geq 0} K_n$.
\item[$(iii)$] There exists a sequence of closed subsets $\{ K_{n} \subseteq X \}_{n \geq 0}$ such that
$U = \bigcup_{n \geq 0} K_n$.
\end{itemize}
Let $\calB$ denote the collection of all open subsets of $X$ which satisfy these conditions. Then:
\begin{itemize}
\item[$(1)$] The elements of $\calB$ form a basis for the topology of
$X$. 

\item[$(2)$] Each element of $\calB$ is paracompact. 

\item[$(3)$] The
collection $\calB$ is stable under finite intersections (in particular, $X \in \calB$).

\item[$(4)$] The empty set $\emptyset$ belongs to $\calB$.
\end{itemize}
\end{proposition}

\begin{remark}
A subset of $X$ which can be written as a countable union of closed subsets of $X$ is called an\index{not}{Fsigma@$F_{\sigma}$}
{\it $F_{\sigma}$}-subset of $X$. Consequently, the basis $\calB$ for the topology of $X$ appearing in Proposition \ref{gooffy} can be characterized as the collection of open $F_{\sigma}$-subsets of $X$.
\end{remark}

\begin{remark}
If the topological space $X$ admits a metric $d$, then {\em every} open subset $U \subseteq X$
belongs to the basis $\calB$ of Proposition \ref{gooffy}. Indeed, we may assume without loss of generality that the diameter of $X$ is at most $1$ (adjusting the metric if necessary), in which case the function
$$ f(x) = d(x, X-U) = \inf_{y \notin U} d(x,y)$$
satisfies condition $(i)$.
\end{remark}

\begin{proof}
We first show that $(i)$ and $(ii)$ are equivalent. If $(i)$ is satisfied, then the closed subsets
$K_n = \{ x \in X: f(x) \geq \frac{1}{n} \}$ satisfy the demands of $(ii)$. Suppose next that $(ii)$ is satisfied.
For each $n \geq 0$, let $G_{n}$ denote the closure of $X - K_{n+1}$, so that $G_n \cap K_{n} = \emptyset$. It follows that there exists a continuous function $f_{n}: X \rightarrow [0, 1]$
such that that $f_n$ vanishes on $G_n$ and the restriction of $f$ to $K_{n}$ is the constant function taking the value $1$. Then the function $f = \Sum_{n > 0} \frac{f_n}{2^n}$ has the property required by $(i)$.

We now prove that $(ii) \Leftrightarrow (iii)$. The implication $(ii) \Rightarrow (iii)$ is obvious.
For the converse, suppose that $U = \bigcup_{n} K_n$, where the $K_n$ are closed subsets
of $X$. We define a new sequence of closed subsets $\{ K'_{n} \}_{n \geq 0}$ by induction as follows.
Let $K'_0 = K_0$. Assuming that $K'_n$ has already been defined, let $V$ and $W$ be
disjoint open neighborhoods of the closed sets $K'_{n} \cup K_{n+1}$ and $X-U$, respectively (the existence of such neighborhoods follows from the assumption that $X$ is paracompact; in fact, it would suffice to assume that $X$ is normal), and define $K'_{n+1}$ to be the closure of $V$. It is then
easy to see that the sequence of closed sets $\{ K'_{n} \}_{n \geq 0}$ satisfies the requirements of $(ii)$.

We now verify properties $(1)$ through $(4)$ of the collection of open sets $\calB$. Assertions $(3)$ and $(4)$ are obvious. To prove $(1)$, consider an arbitrary point $x \in X$ and an open set $U$ containing $x$. Then the closed sets $\{x\}$ and $X - U$ are disjoint, so there exists a continuous function
$f: X \rightarrow [0,1]$ supported on $U$ such that $f(x) = 1$. Then $U' = \{ y \in X: f(y) > 0$ is
an open neighborhood of $x$ contained in $U$, and $U' \in \calB$.

It remains to prove $(2)$. Let $U \in \calB$; we wish to prove that $U$ is paracompact. Write
$U = \bigcup_{n \geq 0} K_n$, where each $K_{n}$ is a closed subset of $X$ containing a neighborhood of $K_{n-1}$ (by convention, we set $K_{n} = \emptyset$ for $n < 0$). 
Let $\{ U_{\alpha} \}$ be an open covering of $X$. Since
each $K_n$ is paracompact, we can choose a locally finite covering
$\{ V_{\alpha, n} \}$ of $K_n$ which refines $\{ U_{\alpha} \cap K_n \}$. Let
$V^{0}_{\alpha,n}$ denote the intersection of $V_{\alpha,n}$ with 
the interior of $K_n$, and let $W_{\alpha,n} = V^{0}_{\alpha,n} \cap (X - K_{n-2})$. 
Then $\{ W_{\alpha, n} \}$ is a locally finite open covering of $X$ which refines
$\{ U_{\alpha} \}$.
\end{proof}

Let $X$ be a paracompact topological space, and let $\calB$ be the basis constructed in
Proposition \ref{gooffy}. Then $\calB$ can be viewed as a category with finite limits, and
is equipped with a natural Grothendieck topology. To simplify the notation, we will let
$\Shv(\calB)$ denote the $\infty$-topos $\Shv( \Nerve(\calB))$. Note that because
$\Nerve(\calB)$ is the nerve of a partially ordered set, the $\infty$-topos $\Shv(\calB)$ is
$0$-localic. Moreover, the corresponding locale $\Sub(\bf{1})$ of subobjects of
the final object ${\bf 1} \in \Shv(\calB)$ is isomorphic to the lattice of open subsets of $X$. It follows that the restriction map $\Shv(X) \rightarrow \Shv(\calB)$ is an equivalence of $\infty$-topoi.

\begin{warning}
Let $X$ be a topological space and $\calB$ a basis of $X$, regarded as a partially ordered set
with respect to inclusions. Then $\calB$ inherits a Grothendieck topology, and we can define
$\Shv(\calB)$ as above. However, the induced map $\Shv(X) \rightarrow \Shv(\calB)$ is generally not an equivalence of $\infty$-categories: this requires the assumption that $\calB$ is stable under finite intersections. In other words, a sheaf (of spaces) on $X$ generally {\em cannot} be recovered by knowing its sections on a basis for the topology of $X$; see Counterexample \ref{spacerk}.
\end{warning}

\subsection{Spaces over $X$}\label{para1}

Let $X$ be a topological space with a specified basis $\calB$, fixed throughout this section.
We wish to study the homotopy theory of {\it spaces over $X$}; that is, spaces $Y$ equipped with a map $p: Y \rightarrow X$.
We should emphasize that we do not wish to assume that the map $p$ is a fibration, or that $p$ is equivalent to a fibration in any reasonable sense: we are imagining that $p$ encodes a {\it sheaf of spaces} on $X$, and we do not wish to impose any condition of local triviality on this sheaf.

Let $\Top$\index{not}{Top@$\Top$} denote the category of topological spaces, and $\Top_{/X}$ the category of topological spaces mapping to $X$. For each $p: Y \rightarrow X$ and every open subset $U \subseteq X$,
we define a simplicial set $\Sing_{X}(Y,U)$ by the formula
$$ \Sing_{X}(Y,U)_{n} = \Hom_{X} ( U \times |\Delta^n|, Y )\index{not}{Sing_X@$\Sing_X$}.$$
Face and degeneracy maps are defined in the obvious way. We note the simplicial set $\Sing_X(Y,U)$ is {\em always} a Kan complex. We will simply write $\Sing_{X}(Y)$ to denote the simplicial presheaf on $X$ given by
$$ U \mapsto \Sing_{X}(Y,U).$$

\begin{proposition}\label{qur}\index{gen}{model category!of spaces over $X$}
There exists a model structure on the category $\Top_{/X}$, uniquely determined by the following properties:
\begin{itemize}
\item[$(W)$] A morphism $$ \xymatrix{ Y \ar[rr] \ar[dr]^{p} & & Z \ar[dl]^{q} \\
& X & }$$ is a {\it weak equivalence} if and only if, for every $U \subseteq X$ belonging to $\calB$, the induced map $\Sing_X(Z,U)_{\bigdot} \rightarrow \Sing_X(Y,U)_{\bigdot}$ is a homotopy equivalence of Kan complexes.
\item[$(F)$] A morphism $$ \xymatrix{ Y \ar[rr] \ar[dr]^{p} & & Z \ar[dl]^{q} \\
& X & }$$ is a {\it fibration} if and only if, for every $U \subseteq X$ belonging to $\calB$, the induced map $\Sing_X(Z,U)_{\bigdot} \rightarrow \Sing_X(Y,U)_{\bigdot}$ is a Kan fibration.
\end{itemize}
\end{proposition}

\begin{remark}
The model structure on $\Top_{/X}$ described in Proposition \ref{qur} depends on the chosen basis $\calB$ for $X$, and not only on the topological space $X$ itself.
\end{remark}

\begin{proof}[Proof of Proposition \ref{qur}]
The proof uses the theory of {\em cofibrantly generated} model categories; we give a sketch and refer the reader to \cite{hirschhorn} for more details. We will say that a morphism $Y \rightarrow Z$ in $\Top_{/X}$ is a cofibration if it has the left lifting property with respect to every trivial fibration in $\Top_{/X}$. 

We begin by observing that a map $Y \rightarrow Z$ in $\Top_{/X}$ is a fibration if and only if it has the right lifting property with respect to every inclusion $U \times \Lambda^n_i \subseteq U \times \Delta^n$, where $0 \leq i \leq n$ and $U$ is in $\calB$. Let $\calI$ denote the weakly saturated class of morphisms in $\Top_{/X}$ generated by these inclusions. Using the small object argument, one can show that every morphism $Y \rightarrow Z$ in $\Top_{/X}$ admits a factorization
$$ Y \stackrel{f}{\rightarrow} Y' \stackrel{g}{\rightarrow} Z$$
where $f$ belongs to $\calI$ and $g$ is a fibration. (Although the objects in $\Top_{/X}$ are not generally small, one can still apply the small object argument since they are small {\it relative} to the class $\calI$ of morphisms: see \cite{hirschhorn}).

Similarly, a map $Y \rightarrow Z$ is a trivial fibration if and only if it has the right lifting property with respect to every inclusion $U \times |\bd \Delta^n| \subseteq U \times |\Delta^n|$, where $U \in \calB$. 
Let $\calJ$ denote the weakly saturated class of morphisms generated by these inclusions: then every
morphism $Y \rightarrow Z$ admits a factorization 
$$ Y \stackrel{f}{\rightarrow} Y' \stackrel{g}{\rightarrow} Z$$
where $f$ belongs to $\calJ$ and $g$ is a trivial fibration.

The only nontrivial point to verify is that every morphism which belongs to $\calI$ is a trivial cofibration; once this is established, the axioms for a model category follow formally. Since it is clear that $\calI$ is contained in $\calJ$, and that $\calJ$ consists of cofibrations, it suffices to show that every morphism in $\calI$ is a weak equivalence. To prove this, let us consider the class $\calK$ of all closed immersions $k: Y \rightarrow Z$ in $\Top_{/X}$ such that there exist functions $\lambda: Z \rightarrow [0, \infty)$ and $h: Z \times [0, \infty) \rightarrow Z$ such that $k(Y) = \lambda^{-1} \{0\}$, $h(z,0)=z$, and $h(z, \lambda(z)) \in k(Y)$. Now we make the following observations:
\begin{itemize}
\item[$(1)$] Every inclusion $U \times |\Lambda^n_i| \subseteq U \times |\Delta^n|$ belongs to $\calK$.
\item[$(2)$] The class $\calK$ is weakly saturated; consequently, $\calI \subseteq \calK$.
\item[$(3)$] Every morphism $k: Y \rightarrow Z$ which belongs to $\calK$ is a homotopy equivalence in $\Top_{/X}$, and is therefore a weak equivalence.
\end{itemize}
\end{proof}

The category $\Top_{/X}$ is naturally tensored over simplicial sets, if we define $Y \otimes \Delta^n = Y \times |\Delta^n|$ for $Y \in \Top_{/X}$. This induces a simplicial structure on $\Top_{/X}$, which is obviously compatible with the model structure of Proposition \ref{qur}.

We note that $\Sing_X$ is a (simplicial) functor from $\Top_{/X}$ to the category of simplicial presheaves on $\calB$ (here we regard $\calB$ as a category whose morphisms are given by inclusions of open subsets of $X$). We regard $\Set_{\Delta}^{\calB^{op}}$ as a simplicial model category, via the {\it projective} model structure described in \S \ref{quasilimit3}.
By construction, $\Sing_X$ preserves fibrations and trivial fibrations. Moreover, the functor $\Sing_X$ has a left adjoint $$ F \mapsto |F|_{X};$$
we will refer to this left adjoint as {\it geometric realization} (in the case where $X$ is a point, it coincides with the usual geometric realization functor from $\sSet$ to the category of topological spaces).
The functor $|F|_{X}$ is determined by the property that $|F_U|_{X} \simeq U$ if $F_U$ denotes the presheaf (of sets) represented by $U$, and the requirement that geometric realization commutes with colimits and with tensor products by simplicial sets.

We may summarize the situation as follows:

\begin{proposition}\label{exquill}
The adjoint functors $(||_{X}, \Sing_X)$ determine a $($simplicial$)$ Quillen adjunction between
$\Top_{/X}$ $($with the model structure of Proposition \ref{qur} $)$ and $\Set_{\Delta}^{\calB^{op}}$ $($with the projective model structure$)$.
\end{proposition}

\subsection{The Sheaf Condition}\label{para2}

Let $X$ be a topological space and $\calB$ a basis for the topology
of $X$ which is stable under finite intersections. Let $\bfA$ denote the category
of $\Set_{\Delta}^{\calB^{op}}$ of simplicial presheaves on $\calB$; we regard
$\bfA$ as a model category with respect to the {\em projective} model structure defined in \S \ref{quasilimit3}. According to Proposition \ref{othermod}, the $\infty$-category $\sNerve (\bfA^{\degree})$ associated to $\bfA$ is equivalent to the $\infty$-category 
$\calP(\calB) = \calP( \Nerve(\calB))$\index{not}{PcalcalB@$\calP(\calB)$} of presheaves on $\calB$. In particular, the homotopy category
$\h{\calP( \Nerve(\calB))}$ is equivalent to the homotopy category $\h{\bfA}$ ( the category obtained from $\bfA$ by formally inverting all weak equivalences of simplicial presheaves). The $\infty$-category $\Shv(\calB)$ is a reflective subcategory of $\calP(\Nerve(\calB))$. Consequently, we may identify the homotopy category $\h{\Shv(\calB)}$ with a reflective subcategory of $\h{\bfA}$. We will say that a simplicial presheaf $F: \calB^{op} \rightarrow \sSet$ is a {\it sheaf} if it belongs to this reflective subcategory. The purpose of this section is to obtain an explicit criterion which will allow us to test whether or not a given simplicial presheaf $F: \calB^{op} \rightarrow \sSet$ is a sheaf.

\begin{warning}
The condition that a simplicial presheaf $F: \calB^{op} \rightarrow \sSet$ be a sheaf, in the sense defined above, is generally unrelated to the condition that $F$ be a simplicial object in the category of sheaves of sets on $X$ (though these two notions do agree in the special case where the simplicial presheaf $F$ takes values in {\em constant} simplicial sets).
\end{warning}

Let $j: \Nerve(\calB) \rightarrow \calP(\calB)$ be the Yoneda embedding. By definition, an object
$F \in \calP(\calB)$ belongs to $\Shv(\calB)$ if and only if, for every $U \in \calB$ and every
monomorphism $i: U^{0} \rightarrow j(U)$ which corresponds to a {\em covering} sieve $\calU$ on $U$, the induced map
$$ \bHom_{ \calP(\calB)}( j(U), F) \rightarrow \bHom_{ \calP(\calB)}( U^0, F)$$
is an isomorphism in the homotopy category $\calH$. In order to make this condition explicit
in terms of simplicial presheaves, we note that $i: U^{0} \rightarrow j(U)$ can be identified with the inclusion $\chi_{\calU} \subseteq \chi_{U}$ of simplicial presheaves, where
$$\chi_{U}(V) = \begin{cases} \ast & \text{if } V \subseteq U \\
\emptyset & \text{otherwise.} \end{cases}$$
$$\chi_{\calU}(V) = \begin{cases} \ast & \text{if } V \in \calU \\
\emptyset & \text{otherwise.} \end{cases}$$
However, we encounter a technical issue: in order to extract the correct space of maps
$\bHom_{\calP(\calB)}( U^0, F)$, we need to select a {\em projectively cofibrant} model for $U^0$ in $\bfA$. In general, the simplicial presheaf $\chi_{\calU}$ defined above is not projectively cofibrant.
To address this problem, we will construct a new simplicial presheaf, equivalent to $\chi_{\calU}$, which has better mapping properties.

\begin{definition}
Let $\calU$ be a linearly ordered set equipped with a map $s: \calU \rightarrow \calB$.\index{not}{NCalU@$N_{\calU}$}
We define a simplicial presheaf $N_{\calU}: \calB^{op} \rightarrow \sSet$ as follows: for each $V \in \calB$, let $N_{\calU}(V)$ be the nerve of the linearly ordered set $\{ U \in \calU: V \subseteq s(U) \}$. 
$N_{\calU}$ may be viewed as a subobject of the constant presheaf $\underline{ \Delta^{\calU} }$ taking the value $\Nerve(\calU) = \Delta^{\calU}$. 
\end{definition}

\begin{remark}
The above notation is slightly abusive, in that $N_{\calU}$ depends not only on $\calU$, but on the map $s$ and on the linear ordering of $\calU$.  If the map $s$ is injective (as it will be in most applications), we will frequently simply identify $\calU$ with its image in $\calB$. In practice, $\calU$ will usually be a covering sieve on some object $U \in \calB$.
\end{remark}

\begin{remark}
The linear ordering of $\calU$ is {\em unrelated} to the partial ordering of $\calB$ by inclusion.
We will write the former as $\leq$ and the latter as $\subseteq$.
\end{remark}

\begin{example}
Let $\calU = \emptyset$. Then $N_{\calU} = \emptyset$.
\end{example}

\begin{example}
Let $\calU = \{ U\}$ for some $U \in \calB$, and let $s: \calU \rightarrow \calB$ be the inclusion. Then $N_{\calU} \simeq \chi_{U}$.
\end{example}

\begin{proposition}\label{murkminus}
Let $$ \xymatrix{ \calU \ar[dd]^{p} \ar[dr]^{s} & \\
& \calB \\
\calU' \ar[ur]^{s'} & }$$
be a commutative diagram, where $p$ is an order-preserving injection between linearly ordered sets. Then the induced map $N_{\calU} \rightarrow N_{\calU'}$ is a projective cofibration of simplicial presheaves.
\end{proposition}

\begin{proof}
Without loss of generality, we may identify $\calU$ with a linearly ordered subset of $\calU'$ via $p$. Choose a transfinite sequence of simplicial subsets of $N \calU'$ 
$$ K_0 \subseteq K_1 \subseteq \ldots$$
where $K_0 = N \calU$, $K_{\lambda} = \bigcup_{\alpha < \lambda} K_{\alpha}$ if $\lambda$ is a nonzero limit ordinal, and $K_{\alpha+1}$ is obtained from $K_{\alpha}$ by adjoining a single nondegenerate simplex (if such a simplex exists). For each ordinal $\alpha$, let
$F_{\alpha} \subseteq N_{\calU'}$ be defined by
$$ F_{\alpha}(V) = N_{\calU'}(V) \cap K_{\alpha} \subseteq \Nerve(\calU'). $$
Then $F_0 = N_{\calU}$, $F_{\lambda} = \colim_{\alpha < \lambda} F_{\alpha}$ when $\alpha$ is a nonzero limit ordinal, and $F_{\alpha} \simeq N_{\calU'}$ for $\alpha \gg 0$. It therefore suffices to show that each map $F_{\alpha} \rightarrow F_{\alpha+1}$ is a projective cofibration. If
$K_{\alpha} = K_{\alpha+1}$, this is clear; otherwise, we may suppose that $K_{\alpha+1}$ is obtained from $K_{\alpha}$ by adjoining a single nondegenerate simplex $\{ U_0 < U_1 < \ldots U_{n} \}$ of $\Nerve(\calU')$. Let $U = s'(U_0) \cap \ldots \cap s'(U_n) \in \calB$. 

Then there is a coCartesian square
$$ \xymatrix{ \chi_{U} \otimes \bd \Delta^n \ar[r] \ar[d] & \chi_U \otimes \Delta^n \ar[d] \\
F_{\alpha} \ar[r] & F_{\alpha+1}, }$$
The desired result now follows, since the upper horizontal arrow is clearly a projective cofibration.
\end{proof}

\begin{corollary}
Let $\calU$ be a linearly ordered set and $s: \calU \rightarrow \calB$ a map. Then the simplicial presheaf $N_{\calU} \in \Set_{\Delta}^{\calB^{op}}$ is projectively cofibrant.
\end{corollary}

Note that $N_{\calU}(V)$ is contractible if $V \subseteq s(U)$ for some $U \in \calU$, and empty otherwise. Consequently, we deduce:

\begin{corollary}
Let $\calU \subseteq \calB$ be a sieve, equipped with a linear ordering.
The unique map $N_{\calU} \rightarrow \chi_{\calU}$ is a weak equivalence of simplicial presheaves.
\end{corollary}

\begin{notation}
Let $\calU$ be a linearly ordered set equipped with a map $s: \calU \rightarrow \calB$.
For any simplicial presheaf $F: \calB^{op} \rightarrow \sSet$, we let
$F(\calU)$ denote the simplicial set $\bHom_{\bfA}( N_{\calU}, F)$.
\end{notation}

\begin{remark}
Let $U \in \calB$, $\calU = \{ U\}$ and $s: \calU \rightarrow \calB$ is the inclusion. Then $F( \calU ) = F(U)$. In general, we can think of $F( \calU)$ as a {\it homotopy limit} of $F(V)$ taken over
$V$ in the sieve generated by $s: \calU \rightarrow \calB$. To give a vertex of $F(\calU)$, we must give for each $U \in \calU$ a point of $F( sU)$; for every pair of objects $U, V \in \calU$ a path between the corresponding points in $F( sU \cap sV)$, and so forth.
\end{remark}

\begin{corollary}\label{critsheaf}
Let $F: \calB^{op} \rightarrow \Kan$ be a $($projectively fibrant$)$ simplicial presheaf on $\calB$. Then $F$ is a sheaf if and only if, for every $U \in \calB$ and every sieve $\calU$ that covers $U$, there
exists a linearly ordered set $\calU_0$ equipped with a map $\calU_0 \rightarrow \calU$, which generates $\calU$ as a sieve, such that the induced map
$F(U) \rightarrow F(\calU_0)$ is a weak homotopy equivalence of simplicial sets.
\end{corollary}

\begin{lemma}\label{partit}
Suppose that $U \subseteq X$ is paracompact, and let $\calU \subseteq \calB$ be a covering
of $U$. Choose a linear ordering of $\calU$. Then the natural map 
$\pi: |N_{\calU}|_{X} \rightarrow U$ is a homotopy equivalence in $\Top_{/X}$.
$($In other words, there exists a section $s: U \rightarrow N_{\calU}$ of $\pi$, such that
$s \circ \pi$ is fiberwise homotopic to the identity.$)$
\end{lemma}

\begin{proof}
Any partition of unity subordinate to the open cover $\calU$ gives rise to a section of $\pi$. 
To check that $s \circ \pi$ is fiberwise homotopic to the identity, use a ``straight line'' homotopy.
\end{proof}

\begin{proposition}\label{aese}
Let $X$ be a topological space, $\calB$ a basis for the topology of $X$. Assume that $\calB$ is stable under finite intersections and that each element of $\calB$ is paracompact. For
every continuous map of topological spaces $p: Y \rightarrow X$, the simplicial presheaf $\Sing_{X}(Y)$ of sections of $p$ is sheaf.
\end{proposition}

\begin{proof}
Let $F = \Sing_{X}(Y)$. We note that $F$ is a projectively fibrant simplicial presheaf on $\calB$.
By Corollary \ref{critsheaf}, it suffices to show that for every $U \in \calB$,
every covering $\calU$ of $U$, and every linear ordering on $\calU$, 
the natural map $F(U) \rightarrow F(\calU)$ is a homotopy equivalence of simplicial sets. In other words, it suffices to show that composition with the projection $\pi: N_{\calU} \rightarrow U$ induces a homotopy equivalence
$$ \bHom_{/X}( U, Y ) \rightarrow \bHom_{/X}( |N_{\calU}|_{X}, Y)$$
of simplicial sets. This follows immediately from Lemma \ref{partit}.
\end{proof}

\begin{remark}
Under the hypotheses of Proposition \ref{aese}, the object of $\Shv(X)$ corresponding
to the simplicial presheaf $\Sing_{X}(Y)$ is {\em not necessarily hypercomplete}. 
\end{remark}

\subsection{The Main Result}\label{para3}

Suppose that $X$ is a paracompact topological space and $\calB$ is the basis for the topology of $X$ described in Proposition \ref{gooffy}. Our main goal is to show that the composition of the adjoint functors
$$ F \mapsto \Sing_{X} |F|_{X}$$ may be identified with a ``sheafification'' of $F$, at least in the case where $F$ is a projectively cofibrant simplicial presheaf on $\calB$. 

In proving this, we have some flexibility regarding the choice of $F$: it will suffice to treat the question after replacing $F$ by a weakly equivalent simplicial presheaf $F'$, provided that $F'$ is also projectively cofibrant. Our first step is to make a particularly convenient choice for $F'$.

\begin{lemma}\label{grumppp}
Let $\calB$ be a partially ordered set $($via $\subseteq${}$)$ with a least element $\emptyset$, and $F: \calB^{op} \rightarrow \sSet$ be an arbitrary simplicial presheaf such that $F(\emptyset)$ is weakly contractible.

There
exists a $($linearly ordered$)$ set $V$ and a simplicial presheaf $F': \calB^{op} \rightarrow \sSet$ with the following properties:
\begin{itemize}
\item[$(1)$] There exists a monomorphism $F' \rightarrow \underline{\Delta^{V}}$ from $F'$ to the (constant)
simplicial presheaf $\underline{\Delta^{V}}$ on $\calB$ taking the value $\Delta^V$. (Recall that
$\Delta^V$ denotes the nerve of the linearly ordered set $V$.)

\item[$(2)$] For every finite subset $V_0 \subseteq V$, there exists $U \in \calB$ such that 
$U' \subseteq U$ if and only if $\Delta^{V_0} \subseteq F'(U') \subseteq \Delta^{V}$.

\item[$(3)$] As a simplicial presheaf on $\calB$, $F'$ is projectively cofibrant.

\item[$(4)$] In the homotopy category of $\Set_{\Delta}^{\calB^{op}}$, $F'$ and $F$ are equivalent to one another.

\end{itemize}

\end{lemma}

\begin{proof}
Without loss of generality, we may suppose that $F$ is (weakly) fibrant. We now build a ``cellular model'' of $F$. More precisely, we construct the following data:

\begin{itemize}
\item[$(A)$] A transfinite sequence of simplicial sets
$$ Y_0 \rightarrow Y_1 \rightarrow \ldots, $$
where $Y_{\alpha}$ is defined for all ordinals $< \alpha_0$.

\item[$(B)$] For each $\alpha < \alpha_0$, a subsheaf $F_{\alpha}$ of the constant presheaf
on $\calB$ taking the value $Y_{\alpha}$.

\item[$(C)$] A compatible family of maps $F_{\alpha} \rightarrow F$, so that we may regard $\{F_{\alpha} \}$ as a functor from the linearly ordered set $\{ \alpha: \alpha < \alpha_0 \}$ to
$(\sSet)^{\calB^{op}}_{/F}$. 

\item[$(D)$] For each $\alpha < \alpha_0$, there exists $U \in \calB$, $n \geq 0$, and compatible pushout diagrams
$$ \xymatrix{ \bd \Delta^n \ar[d] \ar@{^{(}->}[r] & \Delta^n \ar[d] \\
\colim_{ \beta < \alpha} Y_{\beta} \ar@{^{(}->}[r] & Y_{\alpha}, }$$
$$ \xymatrix{ \chi_{U} \times \bd \Delta^n \ar[d] \ar@{^{(}->}[r] & \chi_{U} \times \Delta^n \ar[d] \\
\colim_{ \beta < \alpha} F_{\beta} \ar@{^{(}->}[r] & F_{\alpha}, }$$
where $$ \chi_U(W) = \begin{cases} \ast & \text{if } W \subseteq U \\
\emptyset & \text{otherwise.} \end{cases} $$

\item[$(E)$] The canonical map $\colim_{\beta < \alpha_0} F_{\beta} \rightarrow F$
is a weak equivalence in $\Set_{\Delta}^{\calB^{op}}$.

\end{itemize}

The construction of this data is reasonably standard, and left to the reader.
Let $Y = \colim_{ \beta < \alpha_0} Y_{\beta}$. Let $Y''$ be the second barycentric subdivision of $Y$, so that $Y''$ may be identified with a simplicial complex (that is, a simplicial subset of $\Delta^{V}$ for some linearly ordered set $V$. For each $\alpha$, let $F''_{\alpha}$ denote
the result of applying the second barycentric subdivision functor to $F_{\alpha}$ termwise.
Let $F'' = \colim_{ \beta < \alpha_0 } F''_{\beta}$. Finally, we define $F'$ by the coCartesian square
$$ \xymatrix{ F''(\emptyset) \times \chi_{\emptyset} \ar[r] \ar[d] & \Delta^V \times \chi_{\emptyset} \ar[d] \\ F'' \ar[r] & F'. }$$

The simplicial presheaf $F'$ satisfies $(1)$ by construction. Properties $(2)$ and $(3)$ are reasonably clear (in fact, $(3)$ is a formal consequence of $(2)$). Condition $(4)$ holds for the simplicial presheaf $F''$ as a consequence of $(E)$, and the fact that there is a canonical
weak homotopy equivalence $K'' \rightarrow K$, for any simplicial set $K$ (see Variant \ref{baryvar}).
Moreover, the assumption that $F(\emptyset)$ is weakly contractible ensures that $(4)$ remains valid for the pushout $F'$.
\end{proof}

Before we can state the next lemma, let us introduce a bit of notation. Let $F: \calB^{op} \rightarrow \sSet$ be a simplicial presheaf. Then we let $|F|$ denote the presheaf of topological spaces on $\calB$ obtained by composing $F$ with the geometric realization functor; similarly, if 
$G$ is a presheaf of topological spaces on $\calB$, then we let $\Sing G$ denote the presheaf of simplicial sets obtained by composing $G$ with the functor $\Sing$. We note that there is a natural transformation
$$ \Sing |F| \rightarrow \Sing_{X} |F|_X.$$

\begin{lemma}\label{hardd}
Let $X$ be a topological space, $\calB$ the collection of open $F_{\sigma}$ subsets of 
$X$ (see Proposition \ref{gooffy}). Let $F: \calB^{op} \rightarrow \sSet$ be a 
projectively cofibrant simplicial presheaf, which is a sheaf $($that is, a fibrant model for $F$ satisfies
the criterion of Corollary \ref{critsheaf}$)$. 
Then the unit map $F \rightarrow \Sing_{X} |F|_{X}$ is an equivalence.
\end{lemma}

\begin{proof}
We note that the functor $F \mapsto \Sing |F|$ preserves weak equivalences in $F$, and the functor $F \mapsto \Sing_{X} |F|_{X}$ preserves weak equivalences between {\em projectively cofibrant}
presheaves $F$. Consequently, by Lemma \ref{grumppp}, we may suppose without loss of generality that there is a linearly ordered set $V$ and that $F$ is a subsheaf of the constant simplicial presheaf taking the value $\Delta^{V}$, such that $F(\emptyset) = \Delta^V$.

It will be sufficient to prove that
$$\Sing |F| \rightarrow \Sing_{X} |F|_X$$ is an equivalence: in other words, we wish to show that
$(\Sing |F|)(U) \rightarrow (\Sing_{X} |F|_X)(U)$ is a homotopy equivalence of Kan complexes, for every $U \in \calB$. Replacing $X$ by $U$, we can reduce to the problem of showing that
$$ p: (\Sing |F|)(X) \rightarrow (\Sing_X |F|_X)(X)$$
is a homotopy equivalence. 
It now suffices to show that
for every inclusion $K' \subseteq K$ of {\em finite} simplicial sets (that is, simplicial sets with only finitely many nondegenerate simplices), a commutative diagram
$$ \xymatrix{ K' \times \{0\}  \ar@{^{(}->}[d] \ar[r] & (\Sing |F|)(X) \ar[d] \\
K \times \{0\} \ar[r]^-{g} & ( \Sing_{X} |F|_{X} )(X) }$$
can be expanded to a commutative diagram
$$ \xymatrix{ (K' \times \Delta^1) \coprod_{ K' \times \{1\} } (K \times \{1\})  \ar@{^{(}->}[d] \ar[rr] & & (\Sing |F|)(X) \ar[d] \\
K \times \Delta^1 \ar[rr] & & ( \Sing_{X} |F|_{X} )(X). }$$
(In fact, it suffices to treat the case where $K' \subseteq K$ is the inclusion $\bd \Delta^n \subseteq \Delta^n$; however, this will result in no simplification in the following arguments.) 

Now let $\calB = \{ U_{\alpha} \}_{\alpha \in A}$, where $A$ is a linearly ordered set.
Since $F$ is assumed to be a sheaf on $\calB$, the equivalent presheaf $\Sing |F|$ is also a sheaf. 
Consequently, for any covering $\calU \subseteq \calB$ (and any linear ordering of $\calU$), the natural map $(\Sing |F|)(X) \rightarrow (\Sing |F|)(\calU)$ is an equivalence.
Likewise, by Proposition \ref{aese}, the map
$(\Sing_{X} |F|_X)(X) \rightarrow ( \Sing_{X} |F|_X )(\calU)$ is an equivalence. Consequently, it suffices to find a covering $\calU \subseteq \calB$ of $X$ and a diagram
$$ \xymatrix{ (K' \times \Delta^1) \coprod_{ K' \times \{1\} } (K \times \{1\})  \ar[d] \ar[rr] & & (\Sing |F|)(\calU) \ar[d] \\
K \times \Delta^1 \ar[rr]^{G} & & ( \Sing_{X} |F|_{X} )(\calU) }$$
which extends $g$.

Since $K$ is finite, the map $g: K \rightarrow (\Sing_X |F|_X)(X)$ may be identified with a continuous, fiber-preserving map $X \times |K| \rightarrow |F|_X$, which we will also denote by $g$.
By assumption, $F$ is a subsheaf of the constant presheaf taking the value $\Delta^{V}$; constantly, we may identify $|F|_X$ with a subspace of $\Delta^V_X = X \otimes \Delta^V$.
(We may identify $\Delta^V_{X}$ with the product
$X \times | \Delta^V |$ {\em as a set}, though it generally has a finer topology.) We may represent a point of $\Delta^V_{X}$ by an ordered pair $(x, q)$, where $x \in X$ and $q: V \rightarrow [0,1]$
has the property that $\{ v \in V: q(v) \neq 0\}$ is finite and 
$\Sum_{v \in V} q(v) = 1$. For each $v \in V$, we let $\Delta^V_{X,v}$ denote the {\em open} subset of $\Delta^V_{X}$ consisting of all pairs $(x,g)$ such that $q(v) > 0$; note that the sets
$\{ \Delta^V_{X,v}\}_{ v \in V}$ form an open cover of cover $\Delta^V_{X}$.
Consequently, the open sets $\{ g^{-1} \Delta^V_{X,v} \}_{v \in V}$ form an open cover of $X \times |K|$. Let $x$ be a point of $X$. The compactness of $|K|$ implies that there is a finite subset
$V_0 \subseteq V$, an open neighborhood $U_x$ of $X$ containing $x$, and an open covering
$\{ W_{x,v}: v \in V_0 \}$ of $|K|$, such that $g( U_x \times W_{x,v} ) \subseteq \Delta^V_{X,v}$.
Choose a partition of unity subordinate to the covering $\{ W_{x,v} \}$, thereby determining a map $f_{x}: |K| \rightarrow |\Delta^{V_0}|$. The open sets $\{ U_x \}$ cover $X$; since $X$ is paracompact, this covering has a locally finite refinement. Shrinking the $U_x$ if necessary, we may suppose that this refinement is given by $\{ U_{x} \}_{x \in X_0}$, and that each $U_{x}$ belongs to $\calB$. Let $\calU = \{ U_{x} \}_{x \in X_0}$, and choose a linear ordering of $\calU$.

We now define a new map
$g': K \rightarrow (\Sing |F|)(\calU)$. To do so, we must give, for every finite
$\calU_0 = \{ U_{x_0} < \ldots < U_{x_n} \} \subseteq \calU$, a map
$$ g'_{A_0}: |\Delta^{ \{x_0, \ldots, x_n\} }| \times |K| \rightarrow |F(U_{x_0} \cap \ldots \cap U_{x_n})| \subseteq |\Delta^V|,$$
which are required to satisfy some obvious compatibilities. 
Define $g'_{\calU_0}$ by the formula
$$ g'_{\calU_0}( \Sum \lambda_i x_i, z) = \Sum \lambda_i f_{x_i}(z).$$
It is clear that $g'_{\calU_0}$ is well-defined as a map from $| \Delta^{ \{x_0, \ldots, x_n \} } | \times |K|$ to
$| \Delta^{V} |$. We claim that, in fact, this map factors through $| F( U_{x_0} \cap \ldots \cal U_{x_n} )|$. Let $z \in |K|$, and consider the set $V' = \{ v \in V: (\exists 0 \leq i \leq n) [f_{x_i}(z)(v) \neq 0] \} \subseteq V$. Condition $(2)$ of Lemma \ref{grumppp} ensures that there exists $U \in \calB$ such that $\Delta^{V'} \subseteq F(U')$ if and only if
$U' \subseteq U$. We note that, for each $y \in U_{x_0} \cap \ldots \cap U_{x_n}$, we have $g(y,z)(v) \neq 0$ for $v \in V'$; it follows that $y \in U$. Consequently, we deduce that
$U_{\alpha_0} \cap \ldots \cap U_{\alpha_n} \subseteq U$, so that
$\Delta^{V'} \subseteq F( U_{x_0} \cap \ldots \cap U_{x_n} ).$ It follows that
$g'_{\calU_0}| \{z\} \times |\Delta^{\calU_0}|$ factors through $| F(U_{x_0} \cap \ldots \cap U_{x_n} )|$. Since this holds for every $z \in |K|$, it follows that $g'_{A_0}$ is well-defined; evidently these maps are compatible with one another and give the desired map $g': K \rightarrow (\Sing |F|)(\calU)$.

We now observe that the composite maps
$$ K \stackrel{g'}{\rightarrow } (\Sing |F|)(\calU) \rightarrow (\Sing_{X} |F|_X)(\calU)$$
$$ K \stackrel{g}{\rightarrow } (\Sing_X |F|_X)(X) \rightarrow (\Sing_X |F|_X)(\calU)$$
are homotopic via a ``straight-line'' homotopy 
$G: K \times \Delta^1 \rightarrow (\Sing_X |F|_X)(\calU)$, which has the desired properties.
\end{proof}

Now, the hard work is done and we are ready to enjoy the fruits of our labors.

\begin{theorem}\label{main}
Let $X$ be a paracompact topological space and $\calB$ the collection of open
$F_{\sigma}$ subsets of $X$ (see Proposition \ref{gooffy}). Then, for any
projectively cofibrant $F: \calB^{op} \rightarrow \sSet$, the natural map
$$ F \rightarrow \Sing_{X} |F|_X$$ exhibits $\Sing_{X} |F|_X$ as a sheafification of $F$.
\end{theorem}

\begin{proof}
Let $\h{\Top_{/X}}$ be the homotopy category of the model category $\Top_{/X}$ (the category obtained by inverting all of the weak equivalences defined in Proposition \ref{qur}) and
$\h{\Set_{\Delta}^{\calB^{op}}}$ the homotopy category of the category of simplicial presheaves on $\calB$. It follows from Proposition \ref{exquill} that 
the adjoint functors $\Sing_X$ and $||_X$ induce adjoint functors
$$ \Adjoint{||_X^L}{\h{\Set_{\Delta}^{\calB^{op}}}}{\h{\Top_{/X}}}{\Sing_{X}}$$
Here $||^{L}_{X}$ denotes the left-derived functor of the geometric realization (since
every object of $\Top_{/X}$ is fibrant, $\Sing_X$ may be identified with its right-derived functor).

We first claim that for any $Y \in \Top_{/X}$, the counit map
$ | \Sing_X Y |_{X}^{L} \rightarrow Y$ is a weak equivalence. To see this, choose a projectively cofibrant model $F \rightarrow \Sing_X Y$ for $\Sing_X Y$; we wish to show that
the induced map $|F|_{X} \rightarrow Y$ is a weak equivalence. By definition, this is
equivalent to the assertion that $\Sing_X |F|_{X} \rightarrow \Sing_X Y$ is a weak equivalence. But we have a commutative triangle
$$ \xymatrix{ F \ar[rr] \ar[dr] & & \Sing_X Y \ar[dl] \\
& \Sing_X |F|_X & }$$
where left-diagonal map is a weak equivalence by Lemma \ref{hardd} and Proposition \ref{aese}, and the top horizontal map is a weak equivalence by construction; the desired result now follows from the two-out-of-three property.

It follows that we may identify $\h{\Top_{/X}}$ with a full subcategory $\calC \subseteq \h{\Set_{\Delta}^{\calB^{op}}}$. By Proposition \ref{aese}, the objects of this subcategory are sheaves on $\calB$; by Lemma \ref{hardd}, every sheaf on $\calB$ is equivalent to $\Sing_X Y$ for an appropriately chosen $Y$; thus $\calC$ consists of precisely the sheaves on $\calB$.

The composite functor $F \mapsto \Sing_{X} |F|^{L}_{X}$ may be identified with a localization
functor from $\h{\Set_{\Delta}^{\calB^{op}}}$ to the subcategory $\calC$. In particular, when $F$ is projectively cofibrant, the unit of the adjunction $F \rightarrow \Sing_{X} |F|_{X}$ is a localization of $F$.
\end{proof}

\begin{corollary}\label{wamain}
Under the hypotheses of Theorem \ref{main}, the functor $\Sing_{X}$ induces an equivalence of $\infty$-categories
$$ \sNerve( \Top_{/X}^{\degree} ) \rightarrow \Shv(X).$$
In particular, the $\infty$-category $\sNerve ( \Top_{/X}^{\degree} )$ is an $\infty$-topos.
\end{corollary}

\begin{remark}\label{qurk}
In the language of model categories, we may interpret Corollary \ref{wamain} as asserting
that $\Sing_X$, $||_X$ furnish a Quillen equivalence between $\Top_{/X}$ (with the model structure of Proposition \ref{qur}) and $\Set_{\Delta}^{\calB^{op}}$ where the latter is equipped with the following {\em localization} of the projective model structure:
\begin{itemize}
\item[$(1)$] A map $F \rightarrow F'$ in $\Set_{\Delta}^{\calB^{op}}$ is a {\it cofibration} if it is a projective cofibration (in the sense of Definition \ref{projinj}).
\item[$(2)$] A map $F \rightarrow F'$ in $\Set_{\Delta}^{\calB^{op}}$ is a {\it weak equivalence} if it induces an
equivalence in the $\infty$-category $\Shv(X)$.
\end{itemize}
\end{remark}

\subsection{Base Change}\label{dooky}

With Corollary \ref{wamain} in hand, we are {\em almost} ready to deduce Theorem \ref{nice}.
Suppose given a paracompact space $X$, and let $\calB$ denote the collection of all
open $F_{\sigma}$ subsets of $X$. Let $p: \Shv(X) \rightarrow \Shv(\ast) \simeq \SSet$ be the geometric morphism induced by the projection $X \rightarrow \ast$.

For any simplicial set $K$, let $F_{K}$ denote the constant simplicial presheaf
on $\calB$ taking the value $K$. Then, if we endow $\Set_{\Delta}^{\calB^{op}}$ with the localized model structure of Remark \ref{qurk}, then $F_K$ is a model for the sheaf $p^{\ast} K$. Consequently, the space $p_{\ast} p^{\ast} K$ may be identified up to homotopy with the mapping space
$$ \bHom_{\Shv(X)}( F_{\ast}, F_{K} )$$
which, in virtue of Corollary \ref{wamain}, is equivalent to
$$ \bHom_{\Top_{/X}}( X, X \otimes K ) = (\Sing_X (X \otimes K))(X)$$
However, at this point, a technical wrinkle appears: $X \otimes K$ agrees with $X \times |K|$ as a set, but it is equipped with a finer topology (given by the direct limit of the product topologies $X \times |K_0|$, where $K_0 \subseteq K$ is a finite simplicial subset). In general, we have only
an {\em inclusion} of simplicial presheaves
$$ \eta: \Sing_X (X \otimes K) \subseteq \Sing_{X} (X \times |K|),$$
which need not be an isomorphism. However, we will complete the proof of Theorem \ref{nice}
by showing that $\eta$ is an equivalence of simplicial presheaves.

We consider a slightly more general situation. Let $p: X \rightarrow Y$ be a continuous map between paracompact spaces, and let $\calB_{X}$ and $\calB_Y$ the collections of open $F_{\sigma}$ subsets in $X$ and $Y$, respectively. Note that the inverse image along $p$ determines a map
$q: \calB_Y \rightarrow \calB_X$. Composition with $q$ induces a
``pushforward'' functor $q_{\ast}: \Set_{\Delta}^{\calB_X^{op}} \rightarrow \Set_{\Delta}^{\calB_Y^{op}}$, which
has a left adjoint which we will denote by $q^{\ast}$. Similarly, there is a ``pullback'' functor
$p^{\ast}: \Top_{/Y} \rightarrow \Top_{/X}$; however, $p^{\ast}$ generally does not possess a right adjoint. Consider the square
$$ \xymatrix{ \Set_{\Delta}^{\calB_Y^{op}} \ar[r]^{||_Y} \ar[d]^{q^{\ast}} & \Top_{/Y} \ar[d]^{p^{\ast}} \\
\Set_{\Delta}^{\calB_X^{op}} \ar[r]^{||_X} & \Top_{/X}. }$$
This square is ``lax commutative'', in the sense that there exists a natural transformation of functors
$$\eta_{F}: | q^{\ast} F |_{X} \rightarrow p^{\ast} |F|_Y = |F|_Y \times_Y X.$$
The map $\eta_{F}$ is always a bijection of topological spaces, but is generally not a homeomorphism. Nevertheless, we have the following:

\begin{proposition}\label{basechang}
Under the hypotheses above, if $F: \calB_Y^{op} \rightarrow \sSet$ is a projectively cofibrant simplicial presheaf on $Y$, then the map $\eta_{F}: | q^{\ast} F|_{X} \rightarrow |F|_Y \times_Y X$ is a weak equivalence in $\Top_{/X}$.
\end{proposition}

The proof is based on the following lemma:

\begin{lemma}\label{prechange}
Let $Y$ be a paracompact topological space and let $\calB$ be the collection of open $F_{\sigma}$ subsets of $Y$ (see Proposition \ref{gooffy}). Let $V$ be a linearly ordered set. Suppose that for
every nonempty finite subset $V_0 \subseteq V$, we are given an basic open set
$U(V_0) \in \calB$ satisfying the following conditions:
\begin{itemize}
\item[$(a)$] If $V_0 \subseteq V_1$, then $U(V_1) \subseteq U(V_0)$.
\item[$(b)$] The open set $U( \emptyset)$ concides with $X$.
\end{itemize}
Let $F: \calB^{op} \rightarrow \sSet$ be the simplicial presheaf which assigns to
each $U \in \calB$ the simplicial subset $F(U) \subseteq \Delta^{V}$ spanned by those nondegenerate simplices $\sigma$ corresponding to finite subsets $V_0 \subseteq V$ such that $U \subseteq U(V_0)$
(see Lemma \ref{grumppp}). 

For every object $X \in \Top_{/Y}$, an $n$-simplex $\tau$ of 
$\bHom_{ \Top_{/Y}}(Y, |F|_Y)$ determines a map of topological spaces from
$X \times | \Delta^n |$ to $| \Delta^{V} |$, which in turn determines a collection of maps
$\phi_{v}: X \times | \Delta^n | \rightarrow [0,1]$ such that for every $x \in X \times | \Delta^n |$, the 
sum $\Sigma_{v \in V} \phi_v(x)$ is equal to $1$. Let $\bHom^{0}_{\Top_{/Y}}(X, |F|_Y)$ denote the simplicial subset of $\bHom_{ \Top_{/Y}}(X, |F|_Y)$ spanned by those simplices $\tau$ which satisfy the following condition, where $K = \Delta^n$: 
\begin{itemize}
\item[$(\ast)$] There exists a locally finite collection of open sets $\{ U_{v} \subseteq X \times |K| \}_{v \in V}$
such that each $U_{v}$ contains the closure of the support of the function $\phi_{v}$, and
$\bigcap_{v \in V_0} U_{v}$ is contained in the inverse image of $U(V_0)$ for every finite subset $V_0 \subseteq V$.
\end{itemize}
If the topological space $X$ is paracompact, then the inclusion
$i: \bHom^{0}_{ \Top_{/Y}}( X, |F|_Y) \subseteq \bHom_{ \Top_{/Y}}( X, |F|_Y)$ is a homotopy equivalence of Kan complexes.
\end{lemma}

\begin{proof}
Note that for any finite simplicial set $K$, we can identify
$\Hom_{\sSet}(K, \bHom^{0}_{ \Top_{/Y}}( X, |F|_Y)$ with the set of all collections
of continuous maps $\{ \phi_{v}: X \times |K| \rightarrow [0,1] \}$ satisfying the condition $(\ast)$. Composition with a retraction of $| \Delta^n |$ onto a horn $| \Lambda^n_i |$ determines a section of the restriction map 
$$\Hom_{\sSet}( \Delta^n, \bHom^{0}_{ \Top_{/Y}}( X, |F|_Y))
\rightarrow \Hom_{\sSet}( \Lambda^n_i, \bHom^{0}_{ \Top_{/Y}}( X, |F|_Y)),$$
from which it follows that $\bHom^{0}_{ \Top_{/Y}}( X, |F|_Y)$ is a Kan complex. 

To prove that $i$ is a homotopy equivalence, we argue as in the proof of Lemma \ref{hardd}: it will suffice to show that for every  inclusion $K' \subseteq K$ of finite simplicial sets, every commutative diagram $$ \xymatrix{ K' \times \{0\}  \ar@{^{(}->}[d] \ar[r] & \bHom^{0}_{\Top_{/Y}}(X, |F|_Y) \ar[d] \\
K \times \{0\} \ar[r]^-{g} & \bHom_{ \Top_{/Y}}(X, |F|_Y) }$$
can be expanded to a commutative diagram
$$ \xymatrix{ (K' \times \Delta^1) \coprod_{ K' \times \{1\} } (K \times \{1\})  \ar@{^{(}->}[d] \ar[rr] & &  \bHom^{0}_{\Top_{/Y}}(X, |F|_Y) \ar[d] \\
K \times \Delta^1 \ar[rr]^{G} & & \bHom_{\Top_{/Y}}(X, |F|_Y). }$$
The map $g$ is classified by a collection of continuous maps
$\{ g_v: X \times |K| \rightarrow [0,1] \}_{v \in V}$ such that
$\Sigma_{v \in V} g_v(x) = 1$. Let $\{ U_{v} \}_{v \in V}$ be a collection of open subsets of
$X \times |K'|$ satisfying condition $(\ast)$ for the functions $\{ g_{v} | X \times |K'| \}$.
For each $v \in V$, let $W_v = \{ x \in X \times |K|: g_v(x) \neq 0 \}$. Choose a locally finite open
covering $\{ U'_{v} \}_{v \in V}$ of $X \times |K|$ which refines $\{ W_{v} \}$. Let
$\{ g'_{v} \}_{v \in V}$ be a partition of unity such that the closure of the support of each
$g'_{v}$ is contained in $U'_{v}$. We define maps
$\{ G_{v}: X \times |K| \times [0,1] \rightarrow [0,1] \}_{v \in V}$ by the formula
$$ G_{v}(x, t) = \begin{cases} (2t) g'_{v}(x) + (1-2t) g_{v}(x) & \text{ if } t \leq \frac{1}{2} \\
g'_{v}(x) & \text{ if  } t \geq \frac{1}{2}. \end{cases}$$
Then the maps $\{ G_v \}$ determine a continuous map
$G: X \times |K| \times [0,1] \rightarrow |F|_Y$, which we can identify with a map of simplicial
sets $K \times \Delta^1 \rightarrow \bHom_{ \Top_{/Y}}(X, |F|_Y)$. The restriction of
this map to $(K' \times \Delta^1) \coprod_{ K' \times \{1\} } (K \times \{1\})$ factors
through $\bHom^{0}_{ \Top_{/Y}}(X, |F|_Y)$, since the open subsets
$\{ (U_{v} \times [0,1]) \cup (U'_{v} \times \{1\} \}$ satisfy condition $(\ast)$.
\end{proof}

\begin{remark}\label{postchan}
In the situation of Lemma \ref{prechange}, suppose we are given a Kan complex
$$\bHom^{1}_{ \Top_{/Y}}( X, |F|_Y) \subseteq \bHom_{ \Top_{/Y}}( X, |F|_Y)$$ which contains $\bHom^{0}_{ \Top_{/Y}}(X, |F|_Y)$ and is closed
under the formation of ``straight line'' homotopies. More precisely, suppose that any map
$G: \Delta^n \times \Delta^1 \rightarrow \bHom_{ \Top_{/Y}}( X, |F|_Y)$ factors through
$\bHom_{ \Top_{/Y}}^{1}(X, |F|_Y)$, provided that it satisfies the properties listed below:
\begin{itemize} 

\item[$(i)$] The map $G$ is classified by a collection of continuous functions $\{ G_{v}: X \times | \Delta^n| \times [0,1] \rightarrow [0,1] \}_{v \in V}$.

\item[$(ii)$] Each $G_{v}$ can be described by the formula
$$ G_{v}(x, t) = \begin{cases} (2t) g'_{v}(x) + (1-2t) g_{v}(x) & \text{ if } t \leq \frac{1}{2} \\
g'_{v}(x) & \text{ if  } t \geq \frac{1}{2}. \end{cases}$$

\item[$(iii)$] The closure of the support of each $g'_v$ is contained in the open set
$\{ x \in X \times | \Delta^n |: g_v(x) \neq 0 \}$.

\item[$(iv)$] The restriction of $G$ to $\Delta^n \times \{0\}$ belongs to
$\bHom_{ \Top_{/Y}}^{1}(X, |F|_Y)$ (by virtue of $(iii)$, this implies that
$G| \Delta^n \times \{1\}$ factors through
$ \bHom_{ \Top_{/Y}}^{0}(X, |F|_Y) \subseteq \bHom_{ \Top_{/Y}}^{1}(X, |F|_Y)$). 
\end{itemize}
Then the proof of Lemma \ref{prechange} shows that the inclusions
$$ \bHom^{0}_{ \Top_{/Y}}( X, |F|_Y)
\subseteq \bHom^{1}_{ \Top_{/Y}}(X, |F|_Y) \subseteq \bHom_{ \Top_{/Y}}( X, |F|_Y)$$
are homotopy equivalences.
\end{remark}

\begin{proof}[Proof of Proposition \ref{basechang}]
Suppose given a weak equivalence $F \rightarrow F'$ between projectively cofibrant simplicial presheaves $F,F': \calB_Y^{op} \rightarrow \sSet$. Both $q^{\ast}$ and $||_X$ are left Qullen functors, and therefore preserve weak equivalences between cofibrant objects; it follows that $| q^{\ast} F|_{X} \rightarrow | q^{\ast} F'|_{X}$ is a weak equivalence. Similarly, 
$|F|_{Y} \rightarrow |F'|_{Y}$ is a weak equivalence between cofibrant objects of $\Top_{/Y}$.
Since {\em every} object of $\Top_{/Y}$ is fibrant, we conclude that $|F|_{Y} \rightarrow |F'|_{Y}$
is a homotopy equivalence in $\Top_{/Y}$; thus $|F|_{Y} \times_{Y} X \rightarrow |F'|_{Y} \times_{Y} X$ is a homotopy equivalence in $\Top_{/X}$. Consequently, we deduce that $\eta_{F}$ is a weak equivalence if and only if $\eta_{F'}$ is a weak equivalence.

Let $F$ be an arbitrary projectively cofibrant simplicial presheaf; we wish to show that $\eta_{F}$
is a weak equivalence. There exists a trivial projective cofibration $F \rightarrow F'$, where $F'$ is projectively fibrant. It now suffices to show that $\eta_{F'}$ is a weak equivalence. Replacing $F$ by $F'$, we reduce to the case where $F$ is projectively fibrant.

Let $F'$ be a simplicial presheaf on $\calB_Y$ satisfying the conditions of Lemma \ref{grumppp}.
Then $F'$ and $F$ are equivalent in the homotopy category of simplicial presheaves on $\calB_Y$.
Since $F'$ is projectively cofibrant and $F$ is projectively fibrant, there exists a weak equivalence
$F' \rightarrow F$. We may therefore once again reduce to proving that $\eta_{F'}$ is a weak equivalence. Replacing $F$ by $F'$, we may suppose that $F$ satisfies the conditions of
Lemma \ref{grumppp}, for some linearly ordered set $V$.

For each $U \in \calB_Y$, we have a commutative diagram of Kan complexes
$$ \xymatrix{  \bHom^{0}_{ \Top_{/X} }( U, | q^{\ast} F|_{X} ) \ar[r]^{\phi_0} \ar[d] & 
\bHom^{0}_{ \Top_{/Y} }( U, |F|_Y) \ar[d] \\
\bHom_{ \Top_{/X} }( U, | q^{\ast} F|_{X} ) \ar[r]^{\phi} & 
\bHom_{ \Top_{/Y} }( U, |F|_Y) }$$ 
where the vertical maps are defined as in Lemma \ref{prechange}. We wish to show that $\phi$ is a homotopy equivalence. This follows from the observation that $\phi_0$ is an isomorphism, and the vertical arrows are homotopy equivalences by Lemma \ref{prechange}.
\end{proof}

Theorem \ref{nice} now follows immediately from Proposition \ref{basechang}, applied in the case where $Y = \ast$ and $F$ is the constant simplicial presheaf $\calB_{X} \rightarrow \sSet$ taking the value $K$.

\begin{remark}
There is another solution to the technical difficulty presented by the fact that the bijection 
$X \otimes K \rightarrow X \times |K|$ is not necessarily a homeomorphism: one can work in a suitable category of compactly generated topological spaces, where the base change functor
$Z \mapsto X \times_Y Z$ has a right adjoint and therefore automatically commutes with all colimits.
This is perhaps a more conceptually satisfying approach; however, it leads to a proof of
Theorem \ref{nice} only in the special case where the space $X$ is itself compactly generated.
\end{remark}

We close this section by describing a few applications of Proposition \ref{basechang} and its proof to the theory of sheaves (of spaces) on a paracompact topological space $X$.

\begin{corollary}\label{hidebound}
Let $X$ be a paracompact topological space, $Y$ a closed subset of $X$, and $i: Y \rightarrow X$
the inclusion map. Let $\calF$ be an object of $\Shv(X)$, and let $\eta_0$ be a global section of
$i^{\ast} \calF$. Then there exists an open subset $U$ of $X$ which contains $Y$ and a section
$\eta \in \calF(U)$ whose image under the restriction map $\calF( U ) \rightarrow (i^{\ast} \calF)(U \cap Y) = (i^{\ast} \calF)(Y)$ lies in the path component of $\eta_0$.
\end{corollary}

\begin{proof}
Let $\calB$ denote the collection of open $F_{\sigma}$ subsets of $X$. Without loss of generality, we may assume that $\calF$ is represented by a projectively cofibrant simplicial presheaf
$F \subseteq \underline{ \Delta^{V} }$ satisfying the conditions of Lemma \ref{grumppp},
where $V$ is a linearly ordered set. Using Proposition \ref{basechang}, Corollary \ref{wamain}, and 
Lemma \ref{prechange}, it will suffice to prove the following assertion:
\begin{itemize}
\item[$(a)$] For every vertex $\eta_0$ of $\bHom^{0}_{\Top_{/X}}( Y, |F|_X)$ can be
lifted to a vertex of $\bHom^{0}_{ \Top_{/X}}(U, |F|_X)$, for some sufficiently small
paracompact open neighborhood $U$ of $Y$.
\end{itemize}
To prove $(a)$, suppose we are given a vertex of
$\eta_0 \in \bHom^{0}_{\Top_{/X}}( Y, |F|_X)$, corresponding to a collection of functions
$\{ \phi_v: Y \rightarrow [0,1] \}$. Since $\eta_0$ belongs to
$\bHom^{0}_{\Top_{/X}}(Y, |F|_X)$, there exist open sets
$U_v \subseteq Y$ satisfying condition $(\ast)$ of Lemma \ref{prechange}.

For each $y \in Y$, there exists an open neighborhood $W_y$ of $y$ in $X$ for which the
set $V(y) = \{ v \in V: W_y \cap U_v \}$ is finite. Let $W = \bigcup_{y \in Y} W_y$.
Shrinking $W$ if necessary, we may suppose that $W$ is a paracompact open neighborhood of $Y$
in $X$. Since $W$ is paracompact, there exists a locally finite open covering
$\{ W'_{\alpha} \}_{\alpha}$ of $W$, so that for each index $\alpha$ there exists a point
$y_{\alpha} \in Y$ such that $W'_{\alpha} \subseteq W_{y}$. For $v \in V$, let
$U'_{v} = \bigcup_{v \in V(y_{\alpha})} W'_{\alpha}$. The open sets
$U'_{v}$ form a locally finite open covering of $W$, and each intersection
$U'_{v} \cap Y$ is an open subset of $U_{v}$ which contains the closure of the support
of $\phi_{v}$.

For each $v \in V$, choose a continuous function $\phi'_{v}: X \rightarrow [0,1]$ such that
$\phi'_{v} | Y = \phi_{v}$, and the closure of the support of $\phi'_{v}$ is contained in
$U'_{v}$. There exists another open set $U''_{v}$ whose closure is contained in $U'_{v}$, which again contains the closure of the support of $\phi'_{v}$. For every finite subset $V_0 \subseteq V$, let
$K_{V_0}$ denote the intersection $\bigcap_{v \in V_0} \overline{ U}''_{v}$, and let
$K^{0}_{V_0}$ denote the open subset of $K_{V_0}$ given by the inverse image of the
open set $U(V_0) \subseteq X$ (the largest open subset for which
$\Delta^{V_0}$ belongs to $F( U(V_0) ) \subseteq \Delta^V$). Then
$\{ K_{V_0} - K^{0}_{V_0} \}$ is a locally finite collection of closed subsets of $X$, none of which intersects $Y$. Let $K = \bigcup_{V_0}(K_{V_0} - K^{0}_{V_0})$; then $K$ is a closed subset
of $X$. Let $W'$ be an open $F_{\sigma}$-subset of $W$ which contains $Y$ and does not intersect $K$. Replacing $X$ by $W'$, we may assume that
$W = X$ and that $K = \emptyset$.

Since the collection of functions $\{ \phi'_{v} \}_{v \in V}$ have locally finite support,
the function $\phi' = \Sigma_{v \in V} \phi'_{v}$ is well-defined, and takes the value
$1$ on $Y$. The open set $\{ x \in X: \phi'(x) > 0 \}$ is a paracompact open subset of $X$ (Proposition \ref{gooffy}). Shrinking $X$ further, we may suppose
that $\phi'$ is everywhere nonzero on $X$. Set $\phi''_{v} = \frac{ \phi_{v} }{\phi}$ for each
$v \in V$. Then the functions $\phi''_{v}$ determine a vertex
$\eta \in \bHom_{ \Top_{/X}}(X, |F|_X)$. Moreover, the open sets
$\{ U''_{v} \}$ satisfy condition $(\ast)$ appearing in the statement of Lemma \ref{prechange}, so that
$\eta$ belongs to $\bHom^{0}_{\Top_{/X}}( X, |F|_X)$ as desired.
\end{proof}

Corollary \ref{hidebound} admits the following refinement:

\begin{corollary}\label{snottle}
Let $X$ be a paracompact topological space, $Y$ a closed subset of $X$, and $i: Y \rightarrow X$
the inclusion map. Let $\calF$ be an object of $\Shv(X)$. Then the canonical map
$$ \alpha_{\calF}: \varinjlim_{Y \subseteq U} \calF(U) \rightarrow
\varinjlim_{Y \subseteq U} (i^{\ast} \calF)(U \cap Y) \simeq (i^{\ast} \calF)(Y)$$
is a homotopy equivalence. Here the colimit is taken over the filtered partially ordered set of
all open subsets of $X$ which contain $Y$.
\end{corollary}

\begin{proof}
We will prove by induction on $n \geq 0$ that the map $\alpha_{\calF}$ is $n$-connective. The case
$n=0$ follows from Corollary \ref{hidebound}. Suppose that $n > 0$. We must show that,
for every pair of points $\eta, \eta' \in \varinjlim_{Y \subseteq U} \calF(U)$, the induced map
of fiber products
$$ \alpha'_{\calF}: \ast \times_{ \varinjlim_{ Y \subseteq U} \calF(U) } \ast
\rightarrow \ast \times_{ (i^{\ast} \calF)(Y)} \ast$$
is $(n-1)$-connective. Without loss of generality, we may assume that $\eta$ and $\eta'$
arise from sections of $\calF$ over some $U \subseteq X$ containing $Y$. Shrinking $U$ if necessary, we may assume that $U$ is paracompact. Replacing $X$ by $U$, we may assume that
$\eta$ and $\eta'$ arise from global sections $f,f': 1 \rightarrow \calF$, where $1$ denotes the final object of $\Shv(X)$. Let $\calG = 1 \times_{\calF} 1 \in \Shv(X)$. Using the left exactness of
$i^{\ast}$ and Proposition \ref{frent}, we can identify $\alpha'_{\calF}$ with $\alpha_{\calG}$.
We now invoke the inductive hypothesis to deduce that $\alpha_{\calG}$ is $(n-1)$-connective, as desired.
\end{proof}

\begin{lemma}\label{agint}
Let $Y$ be a paracompact topological and $\calB$ the collection of open $F_{\sigma}$ subsets of
$Y$ (see Proposition \ref{gooffy}). Let $V$ be a linearly ordered set, and let $F: \calB^{op} \rightarrow \sSet$ be as in the statement of Lemma \ref{prechange}. Suppose given a paracompact space
$X \in \Top_{/Y}$ and a closed subspace $X' \subseteq X$. Then the map
$$ \bHom^{0}_{ \Top_{/Y}}( X, |F|_Y) \rightarrow \bHom^{0}_{ \Top_{/Y}}( X', |F|_{Y} )$$
is a Kan fibration.
\end{lemma}

\begin{proof}
We must show that every lifting problem of the form
$$ \xymatrix{ \Lambda^{n}_i \ar[r] \ar@{^{(}->}[d] &  \bHom^{0}_{ \Top_{/Y}}( X, |F|_Y ) \ar[d] \\
\Delta^n \ar[r] & \bHom^{0}_{ \Top_{/Y}}( X', |F|_Y ) }$$
admits a solution. Since the pair $( |\Delta^m|, |\Lambda^m_i|)$ is homeomorphic to
$( |\Delta^{m-1}| \times [0,1], | \Delta^{m-1} | \times \{0\})$, we can replace $X$ by
$X \times | \Delta^{m-1} |$ and thereby reduce to the case $m=1$.

Let $Z = ( X \times \{0\} ) \coprod_{ X' \times \{0\} } (X' \times [0,1] )$, and let
$\eta_0 \in \bHom^{0}_{ \Top_{/Y}}( Z, |F|_Y)$; we wish to show that 
$\eta_0$ can be lifted to a point in $\bHom^{0}_{ \Top_{/Y}}( X \times [0,1], |F|_Y)$.
The proof of Corollary \ref{hidebound} shows that we can lift
$\eta_0$ to a point $\eta_1 \in \bHom^{0}_{ \Top_{/Y}}( U, |F|_Y)$, for some open set
$U \subseteq X \times [0,1]$ containing $Z$. 

For each $x \in X$, there exists a real number $\epsilon_x > 0$ and an open neighborhood
$V_x \subseteq X$ such that $V_x \times [0, \epsilon_x) \subseteq U$. 
Since $X$ is paracompact, the open covering $\{ V_x \}_{x \in X}$ admits a locally
finite refinement $\{ W_{\alpha} \}$, so that for each index $\alpha$ there
exists a point $x(\alpha) \in X$ such that $W_{\alpha} \subseteq V_{x(\alpha)}$.
Let $\{ \phi_{\alpha} \}$ be a partition of unity subordinate the covering $W_{\alpha}$, and let
$$ \psi = \Sigma_{ \alpha} \epsilon_{x(\alpha)} \phi_{\alpha}.$$
Since the interval $[0,1]$ is compact, there exists an open neighborhood
$V \subseteq X$ containing $X'$ such that $V \times [0,1] \subseteq U$.
Choose a function $\psi': X \rightarrow [0,1]$ such that
$\psi' | (X-V) = \psi|(X-V)$ and $\psi' | X'$ is equal to $1$. Set
$K = \{ (x,t) \in X \times [0,1]: t \leq \phi(x) \}$, so that
$Z \subseteq K \subseteq U$, and let $\eta_2 \in \bHom^{0}_{ \Top_{/Y}}(K, |F|_X)$ be the 
restriction of $\eta_1$. Since $K$ is a retract of $X \times [0,1]$ in the category
$\Top_{/Y}$, we can lift $\eta_2$ to a point $\eta in \bHom^{0}_{ \Top_{/Y}}( X \times [0,1], |F|_X)$ as desired.
\end{proof}

\begin{proposition}\label{siegland}
Let $X$ be a paracompact topological space. Suppose given a sequence of closed subspaces
$$ X_0 \subseteq X_1 \subseteq X_2 \subseteq \ldots \subseteq X$$
with the following properties:
\begin{itemize}
\item[$(1)$] The union $\bigcup X_{i}$ coincides with $X$.
\item[$(2)$] A subset $U \subseteq X$ is open if and only if each of the intersections
$U \cap X_i$ is an open subset of $X_i$.
\end{itemize}
Then the induced diagram
$$ \Shv(X_0) \rightarrow
\Shv(X_1) \rightarrow \ldots
\rightarrow \Shv(X)$$
exhibits $\Shv(X)$ as the colimit of the sequence $\{ \Shv(X_i) \}_{i \geq 0}$ in the 
$\infty$-category $\RGeom$ of $\infty$-topoi.
\end{proposition}

\begin{remark}\label{stasis}
Hypotheses $(1)$ and $(2)$ of Proposition \ref{siegland} can be summarized by saying that
$X$ is the direct limit of the sequence $\{ X_i \}$ in the category of topological spaces.
It follows from this condition that for any locally compact space $Y$, the product
$X \times Y$ is also the direct limit of the sequence $\{ X_i \times Y \}$.
To prove this, we observe that for any topological space $Z$, we have bijections
$$ \Hom(X \times Y, Z)
\simeq \Hom(X, Z^Y) \simeq \varprojlim \Hom(X_i, Z^Y) \simeq \varprojlim \Hom(X_i \times Y, Z),$$
where $Z^Y$ is endowed with the compact-open topology. In particular, we deduce that
$X \times \Delta^n$ is the direct limit of the topological spaces $X_i \times \Delta^n$, for each $n \geq 0$.
\end{remark}

\begin{proof}
For each nonnegative integer $n$, let $i(n)$ denote the inclusion from
$X_n$ to $X_{n+1}$, and $j(n)$ the inclusion of $X_{n}$ into $X$. 
These functors induce geometric morphisms
$$ \Adjoint{ i(n)^{\ast} }{ \Shv(X_{n+1}) }{ \Shv(X_n)}{i(n)_{\ast}} $$
$$ \Adjoint{ j(n)^{\ast} }{ \Shv(X)}{ \Shv(X_n)}{ j(n)_{\ast}}.$$
Let $\calC$ denote a homotopy inverse limit of the tower of $\infty$-categories
$$ \ldots \rightarrow \Shv( X_2) \stackrel{ i(1)^{\ast}}{\rightarrow} \Shv(X_1)
\stackrel{ i(0)^{\ast}}{\rightarrow} \Shv(X_0).$$
In view of Proposition \ref{colimtopoi}, we can also identify
$\calC$ with the direct limit of the sequence $\{ \Shv(X_i) \}_{i \geq 0}$ in $\RGeom$.
The maps $j(n)$ determine a geometric morphism
$$ \Adjoint{ j^{\ast} }{ \Shv(X)}{\calC}{j_{\ast}}.$$
To complete the proof, it will suffice to show that the functor $j^{\ast}$ is an equivalence of
$\infty$-categories.

We first show that the unit map $u: \id_{\Shv(X)} \rightarrow j_{\ast} j^{\ast}$ is an equivalence of functors.
Let $\calF \in \Shv(X)$; we wish to show that the map
$$u_{\calF}: \calF \rightarrow j_{\ast} j^{\ast} \calF
\simeq \varprojlim j(n)_{\ast} j(n)^{\ast} \calF$$
is an equivalence in $\Shv(X)$. It will suffice to prove the analogous assertion after evaluating
both sides on every open $F_{\sigma}$ subset $U \subseteq X$. Replacing $X$ by $U$, we are reduced to proving that the induced map
$$ \alpha_{\calF}: \calF(X) \rightarrow (j_{\ast} j^{\ast} \calF)(X) \simeq \varprojlim (j(n)^{\ast} \calF)(X_n)$$
is a homotopy equivalence. Let $\calB$ be the collection of all open $F_{\sigma}$ subsets of $X$. Without loss of generality, we may assume that $\calF$ is represented by a projectively cofibrant simplicial presheaf $F: \calB^{op} \rightarrow \SSet$ satisfying the conditions
of Lemma \ref{grumppp}. Using Theorem \ref{main}, Corollary \ref{wamain}, and Proposition \ref{basechang}, we can identify
$\calF(X)$ with the Kan complex of sections $K = \bHom^{0}_{ \Top_{/X}}( X, \widetilde{X})$ and
each $( j(n)^{\ast} \calF)(X_n)$ with the Kan complex of sections $K(n) = \bHom_{ \Top_{/X}}(X_n, \widetilde{X})$. It will therefore suffice to show that the canonical map
$K \rightarrow \varprojlim K(n)$ exhibits $K$ as a homotopy inverse limit of the tower
$\{ K(n) \}$. 

It follows from Remark \ref{stasis} that the map $K \rightarrow \varprojlim K(n)$ is
an isomorphism of simplicial sets. For each $n \geq 0$, let 
$K(n)^{0} = \bHom^{0}_{ \Top_{/X}}( X_n, \widetilde{X} ) \subseteq K(n)$
(with notation as in Lemma \ref{prechange}), and let $K^{0} = \varprojlim K(n)^{0} \subseteq K$.
Lemma \ref{prechange} implies that each inclusion $K(n)^{0} \subseteq K(n)$ is a homotopy equivalence. Lemma \ref{agint} implies that the restriction maps
$K(n+1)^{0} \rightarrow K(n)^{0}$ are Kan fibrations. It follows that the inverse limit $K^0$ of the tower
$\{ K(n)^{0} \}$ is a Kan complex, and that the map
$K^0 \simeq \varprojlim \{ K(n)^{0} \}$ exhibits $K^0$ as the homotopy inverse limit of
$\{ K(n)^{0} \}$. Invoking Remark \ref{postchan}, we deduce that the inclusion $K^0 \subseteq K$
is a homotopy equivalence, so that the equivalent diagram $K \simeq \varprojlim \{ K(n) \}$ exhibits
$K$ as a homotopy inverse limit of $\{ K(n) \}$ as desired.

We now argue that the counit map $v: j^{\ast} j_{\ast} \rightarrow \id$ is an equivalence of functors.
Unwinding the definitions, we must prove the following: given a collection of sheaves
$\calF_{n} \in \Shv(X_n)$ and equivalences $\calF_{n} \simeq i(n)^{\ast} \calF_{n+1}$, 
the canonical map
$$ j(n)^{\ast}( \varprojlim j(n+k)_{\ast} \calF_{n+k}) \rightarrow \calF_n$$
is an equivalence of sheaves on $X_{n}$, for each $n \geq 0$. It will suffice to show that this map induces a homotopy equivalence after passing to the global sections over every open
$F_{\sigma}$ subset $U \subseteq X_n$. There exists a function $\phi_0: X_n \rightarrow [0,1]$ such that
$U = \{ x \in X_n: \phi_0(x) > 0 \}$. Choose a map $\phi: X \rightarrow [0,1]$ such that
$\phi_0 = \phi| X_n$. Replacing $X$ by the paracompact open subset
$\{ x \in X: \phi(x) > 0 \}$, we can reduce to the case where $U = X_n$. 

We will prove by induction on $k$ that, for any compatible collection of sheaves
$\{ \calF_{n} \in \Shv(X_n), \calF_n \simeq i(n)^{\ast} \calF_{n+1} \}$, the map 
$$ \psi: (j(n)^{\ast} \calF)(X_n) \rightarrow \calF_n(X_n)$$
is $k$-connective, where $\calF = \varprojlim j(m)_{\ast} \calF_m$.
If $k > 0$, then it will suffice to show that for any pair of points
$\eta, \eta' \in (j(n)^{\ast} \calF)(X_n)$, the induced map
$$ \psi': \ast \times_{ (j(n)^{\ast} \calF)(X_n) } \ast \rightarrow
\ast \times_{ \calF_n(X_n)} \ast$$
is $(k-1)$-connective. Using Corollary \ref{hidebound}, we may assume that
$\eta$ and $\eta'$ arise from sections $\overline{\eta}, \overline{\eta}' \in
\calF(U)$, for some open neighborhood $U$ of $X_n$. Shrinking $U$ if necessary, we may assume that $U$ is paracompact. Replacing $X$ by $U$, we may assume that $\overline{\eta}$ and
$\overline{\eta}'$ are global sections of $\calF$. Since $j(n)^{\ast}$ is left exact, we can identify
$\psi'$ with the map
$$ j(n)^{\ast} ( \ast \times_{\calF} \ast)(X_n) \rightarrow ( \ast \times_{\calF_n} \ast)(X_n).$$
The $(k-1)$-connectivity of this map now follows from the inductive hypothesis.

It remains to treat the case $k=0$. Fix an element $\eta_n \in \calF_n(X_n)$; we wish to show that
$\eta_n$ lies in the image of $\pi_0 \psi$. For every open set $U \subseteq X$,
the map
$$ \pi_0 \calF(U) = \pi_0 \varprojlim (j(m)_{\ast} \calF_m)(U)
\rightarrow \varprojlim (\pi_0 j(m)_{\ast} \calF_m)(U)
\simeq \varprojlim \pi_0 \calF_m( U \cap X_m)$$
is surjective. Consequently, to prove that $\eta_n$ lies in the image of
$\pi_0 \psi$, it will suffice to show that there exists an open set $U$ containing $X_n$
such that $\eta_n$ can be lifted to $\varprojlim \pi_0 \calF_m(U \cap X_m)$. 
By virtue of assumption $(2)$, it will suffice to construct a sequence of
open $F_{\sigma}$ subsets $\{ U_{m} \subseteq X_m \}_{m \geq n}$ and a sequence of compatible
sections $\gamma_{m} \in \pi_0 \calF_m(U_m)$, such that $U_n = X_n$ and
$\gamma_m = \eta_m$. The construction goes by induction on $m$. Assuming that
$(U_m, \eta_m)$ has already been constructed, we invoke the assumption that
$U_m$ is an $F_{\sigma}$ to choose a continuous function $f: X_m \rightarrow [0,1]$ such that
$U_m = \{ x \in X_m: f(x) > 0 \}$. Let $f': X_{m+1} \rightarrow [0,1]$ be a continuous extension of
$f$, and let $V = \{ x \in X_{m+1}: f'(x) > 0 \}$. Then $V$ is a paracompact open subset of $X_{m+1}$, and $U_{m}$ can be identified with a closed subset of $V$. Applying Corollary \ref{hidebound} to the restriction $\calF_{n+1} | V$, we deduce the existence of an open set $U_{m+1} \subseteq V$
such that $\eta_{k}$ can be extended to a section $\eta_{k+1} \in \pi_0 \calF_{m+1}(U_{m+1})$. 
Shrinking $U_{m+1}$ if necessary, we may assume that $U_{m+1}$ is itself an $F_{\sigma}$, which completes the induction.
\end{proof}

\begin{remark}
Suppose given a given a sequence of closed embeddings of topological spaces
$$ X_0 \subseteq X_1 \subseteq X_2 \subseteq \ldots,$$
and let $X$ be the direct limit of the sequence. Suppose further that:
\begin{itemize}
\item[$(a)$] For each $n \geq 0$, the space $X_n$ is paracompact.
\item[$(b)$] For each $n \geq 0$, there exists an open neighborhood
$Y_n$ of $X_n$ in $X_{n+1}$ and a retraction $r_n$ of $Y_{n}$ onto $X_n$.
\end{itemize}
Then $X$ is itself paracompact, so that the hypotheses of Proposition \ref{siegland} are satisfied and $\Shv(X)$ is the direct limit of the sequence of $\infty$-topoi $\{ \Shv(X_n) \}_{n \geq 0}$. To prove this, it will suffice every open covering $\{ U_{\alpha} \}_{\alpha \in A}$ of $X$ admits a refinement
$\{ V_{\beta} \}_{\beta \in B}$ which is {\em countably locally finite}: that is, there exists a decomposition
$B = \bigcup_{n \geq 0} B_{n}$ such that each of the collections
$\{ V_{\beta} \}_{\beta \in B_n }$ is a locally finite collection of open sets, each of which is contained in some $U_{\alpha}$ (see \cite{munkres}). To construct this locally finite open covering, we choose for
each $n \geq 0$ a locally finite open covering $\{ W_{\beta} \}_{\beta \in B_n}$ of $X_n$ which refines the covering $\{ U_{\alpha} \cap X_n \}_{\alpha \in A}$. For each $\beta \in B_n$, we have
$W_{\beta} \subseteq U_{\alpha}$ for some $\alpha \in A$. We now define
$V_{\beta}$ to be the union of a collection of open subsets 
$\{ V_{\beta}(m) \subseteq X_m \}_{m \geq n}$, which are constructed as follows:
\begin{itemize}
\item If $m = n$, we set $V_{\beta}(m) = W_{\beta}$.
\item Let $m > n$, and let $Z_{m-1}$ be an open neighborhood of $X_{m-1}$ in
$X_m$ whose closure is contained in $Y_{m-1}$. We then set
$V_{\beta}(m)= \{ z \in Z_{m-1}: r_{m-1}(z) \in V_{\beta}(m-1) \} \cap U_{\alpha}.$ 
\end{itemize}
It is clear that each $V_{\beta}(m)$ is an open subset of $X_m$ contained in $U_{\alpha}$, and that
$V_{\beta}(m+1) \cap X_m = V_{\beta}(m)$. Since $X$ is equipped with the direct limit topology,
the union $V_{\beta} = \bigcup_{m} V_{\beta}(m)$ is open in $X$. The only nontrivial point is to
verify that the collection $\{ V_{\beta} \}_{\beta \in B_n}$ is locally finite. 

Pick a point $x \in X$; we wish to prove the existence of a neighborhood $S_x$ of $x$ such that
$\{ \beta \in B_n: S \cap V_{\beta} \neq \emptyset \}$ is finite. Then there exists some
$m \geq n$ such that $x \in X_{m}$; we will construct $S_{x}$ using induction on $m$. 
If $m > n$ and $x \in \overline{Z}_{m-1}$, then let $x' = r_{m-1}(x)$, and set
$S_{x} = S_{x'}$. If $m > n$ and $x \notin \overline{Z}_{m-1}$, or if
$m=n$, then we define $S_{x} = \bigcup_{k \geq m} S_{x}(k)$, where
$S_x(k)$ is an open subset of $X_{k}$ containing $x$, defined as follows.
If $m > n$, let $S_x(m) = X_{m} - \overline{Z}_{m-1}$, and if
$m=n$ let $S_x(m)$ be an open subset of $X_n$ which intersects only finitely many of the sets
$\{ W_{\beta} \}_{ \beta \in B_n }$. If $k > m$, we let 
$S_{x}(k) =  \{ z \in Y_{k-1}: r_{k-1}(z) \in S_x(k-1) \}.$
It is not difficult to verify that the open set $S_x$ has the desired properties.
\end{remark}

\subsection{Higher Topoi and Shape Theory}\label{shapesec}

If $X$ is a sufficiently nice topological space (for example, an absolute neighborhood retract), then there exists a homotopy equivalence $Y \rightarrow X$, where $Y$ is a CW complex.
If $X$ is merely assumed to be paracompact, then it is generally not possible to approximate $X$ well by means of a CW-complex $Y$ equipped with a map to $X$. However, in view of Theorem \ref{main}, one can still extract a substantial amount of information by considering maps
from $X$ to CW complexes. {\it Shape theory} is an attempt to summarize all of this information in a single invariant, called the {\it shape} of $X$. In this section, we will sketch a generalization of shape theory to the setting of $\infty$-topoi.

\begin{definition}\label{defshape}\index{gen}{pro-space}\index{gen}{shape}\index{not}{ProSSet@$\Pro(\SSet)$}
We let $\Pro(\SSet)$ denote the full subcategory of $\Fun(\SSet, \SSet)^{op}$ spanned by
left exact, accessible functors $f: \SSet \rightarrow \SSet$. We will refer to $\Pro(\SSet)$
as the $\infty$-category of {\em pro-spaces}, or as the $\infty$-category of {\it shapes}.
\end{definition}

\begin{remark}
If $\calC$ is a small $\infty$-category which admits finite limits, then any functor
$f: \calC \rightarrow \SSet$ is accessible and may be viewed as an object of $\calP(\calC^{op})$. 
The left exactness of $f$ is then equivalent to the condition that $f$ belongs to
$\Ind(\calC^{op}) = \Pro(\calC)^{op}$. Definition \ref{defshape} constitutes a natural extension
of this terminology to a case where $\calC$ is not necessarily small; here it is convenient to add a hypothesis of accessibility for technical reasons (which will not play any role in the discussion below).
\end{remark}

\begin{definition}\label{deshape}
Let $\calX$ be an $\infty$-topos. According to Proposition \ref{spacefinall}, there exists a
geometric morphism $q_{\ast}: \calX \rightarrow \SSet$, which is unique up to homotopy.
Let $q^{\ast}$ be a left adjoint to $q_{\ast}$ (also unique up to homotopy).
The composition $q_{\ast} q^{\ast}: \SSet \rightarrow \SSet$ is an accessible left-exact functor, which we will refer to as the {\it shape of $\calX$} and denote by $\Sh(\calX) \in \Pro(\SSet)$.
\index{gen}{shape!of an $\infty$-topos}\index{not}{ShcalX@$\Sh(\calX)$}
\end{definition}

\begin{remark}
This definition of the shape of an $\infty$-topos appears also in \cite{toen}.
\end{remark}

\begin{remark}
Let $p_{\ast}: \calY \rightarrow \calX$ be a geometric morphism of $\infty$-topoi and $p^{\ast}$ a left adjoint to $p_{\ast}$. Let $q_{\ast}: \calX \rightarrow \SSet$ and $q^{\ast}$ be as in Definition \ref{deshape}. The unit map $\id_{\calX} \rightarrow p_{\ast} p^{\ast}$ induces a transformation 
$$q_{\ast} q^{\ast} \rightarrow q_{\ast} p_{\ast} p^{\ast} q^{\ast} \simeq (q \circ p)_{\ast} (q \circ p)^{\ast},$$
which we may view as a map $\Sh(\calX) \rightarrow \Sh(\calY)$ in $\Pro(\SSet)$. Via this construction, we may view $\Sh$ as a functor from the homotopy category $\h{\RGeom}$ of $\infty$-topoi to the homotopy category $\h{\Pro(\SSet)}$. We will say that a geometric morphism $p_{\ast}: \calY \rightarrow \calX$ is a {\it shape equivalence} if it induces an equivalence $\Sh(\calY) \rightarrow \Sh(\calX)$ of pro-spaces.\index{gen}{shape!equivalence}
\end{remark}

\begin{remark}
By construction, the shape of an $\infty$-topos $\calX$ is well-defined up to equivalence
in $\Pro(\SSet)$. By refining the above construction, it is possible construct a shape functor from $\RGeom$ to the $\infty$-category $\Pro(\SSet)$, rather than on the level of homotopy.
\end{remark}

\begin{remark}\label{struke}
Our terminology does not quite conform to the usage in classical topology. Recall that if $X$ is a compact metric space, the {\it shape} of $X$ is defined as a pro-object in the {\em homotopy} category of spaces. There is a refinement of shape, known as {\it strong shape}, which takes values in the homotopy category of pro-spaces. Definition \ref{deshape} is a generalization of strong shape, rather than shape. We refer the reader to \cite{shapetheory} for a discussion of classical shape theory.
\end{remark}

\begin{proposition}\label{parashape}
Let $p: X \rightarrow Y$ be a continuous map of paracompact topological spaces.
Then $p_{\ast}: \Shv(X) \rightarrow \Shv(Y)$ is a shape equivalence if and only if,
for every Kan complex $K$, the induced map of Kan complexes 
$\bHom_{\Top}(Y, |K|)  \rightarrow \bHom_{\Top}(X, |K|)$ is a homotopy equivalence.
$($Here $\bHom_{\Top}(Y, |K| )$ denotes the simplicial set whose $n$-simplices are
given by continuous maps $Y \times | \Delta^n | \rightarrow |K|$, and
$\bHom_{\Top}(X, |K| )$ is defined likewise.$)$
\end{proposition}

\begin{proof}
Corollary \ref{wamain} and Proposition \ref{basechang} imply that for any paracompact topological space $Z$ and any Kan complex $K$, there is a natural isomorphism
$$ \Sh( \Shv(Z) ) ( K ) \simeq \bHom_{\Top}(Z, |K|)$$ in the homotopy category $\calH$.
\end{proof}

\begin{example}
Let $X$ be a scheme, let $\toposX$ be the topos of \'{e}tale sheaves on $X$, and let
$\calX$ be the associated $1$-localic $\infty$-topos (see \S \ref{nlocalic}). The shape $\Sh(\calX)$ defined above is closely related to the \'{e}tale homotopy type introduced by Artin and Mazur (see \cite{artinmazur}). There are three important differences:
\begin{itemize}
\item[$(1)$] Artin and Mazur work with pro-objects in the homotopy category $\calH$, rather than with actual pro-objects of $\SSet$. Our definition is closer in spirit to that of Friedlander, who works instead in the homotopy category of pro-objects in $\sSet$ (see \cite{fried}).
\item[$(2)$] The \'{e}tale homotopy type of \cite{artinmazur} is constructed by considering \'{e}tale hypercoverings of $X$; it is therefore more closely related to the shape of the hypercompletion
$\calX^{\hyp}$.
\item[$(3)$] Artin and Mazur generally study a certain completion of $\Sh(\calX^{\hyp})$ with respect to the class of truncated spaces, which has the effect of erasing the distinction between $\calX$ and $\calX^{\hyp}$ and discarding a bit of (generally irrelevant) information.
\end{itemize}
\end{example}

\begin{remark}\label{surety}
Let $\ast$ denote a topological space consisting of a single point. By definition,
$\Shv(\ast)$ is the full subcategory of $\Fun(\Delta^1,\SSet)$ spanned by those morphisms
$f: X \rightarrow Y$ where $Y$ is a final object of $\SSet$. We observe that $\Shv(\ast)$ is equivalent to the full subcategory spanned by those morphisms $f$ as above where
$Y = \Delta^0 \in \SSet$, and that this full subcategory is {\em isomorphic} to $\SSet$.
\end{remark}
\begin{definition}
We will say that an $\infty$-topos $\calX$ has {\it trivial shape} if $\Sh(\calX)$ is equivalent
to the identity functor $\SSet \rightarrow \SSet$.\index{gen}{shape!trivial}
\end{definition}

\begin{remark}
Let $q_{\ast}: \calX \rightarrow \SSet$ be a geometric morphism. Then the unit map
$u: \id_{\SSet} \rightarrow q_{\ast} q^{\ast}$ induces a map of pro-spaces
$\Sh(X) \rightarrow \id_{\SSet}$. Since $\id_{\SSet}$ is a final object in $\Pro(\SSet)$, 
we observe that $\calX$ has trivial shape if and only if $u$ is an equivalence; in other words, if and only if the pullback functor $q^{\ast}$ is fully faithful.
\end{remark}

We now sketch another interpretation of shape theory, based on the $\infty$-topoi associated
to {\em pro-spaces}. Let $X = \SSet$, let $\pi: \SSet \times \SSet \rightarrow \SSet$ be the projection onto the first factor, let $\delta: \SSet \rightarrow \SSet \times \SSet$ denote the diagonal map, and let $\phi: (\SSet \times \SSet)^{/\delta_{\SSet}} \rightarrow \SSet$ be defined as in \S \ref{consweet}. Proposition \ref{colimfam} implies that $\phi$ is a coCartesian fibration. We may identify
the fiber of $\phi$ over an object $X \in \SSet$ with the $\infty$-category $\SSet^{/X}$. For each morphism $f: X \rightarrow Y$ in $\SSet$, $\phi$ associates a functor
$f_{!}: \SSet^{/X} \rightarrow \SSet^{/Y}$, given by composition with $f$.
Since $\SSet$ admits pullbacks, each $f_{!}$ admits a right adjoint $f^{\ast}$, so that
$\phi$ is also a Cartesian fibration, associated to some functor
$\psi: \SSet^{\op} \rightarrow \LGeom$.

Let $\hat{X}: \SSet \rightarrow \SSet$ be a pro-space. Then $\hat{X}$ classifies a left fibration 
$M^{op} \rightarrow \SSet$, where $M$ is a filtered $\infty$-category. 
Let $\theta$ denote the composition
$$ M^{op} \rightarrow \SSet \stackrel{\psi^{op}}{\rightarrow} (\LGeom)^{op}.$$
Although $M$ is generally not small, the accessibility condition on $F$ guarantees the existence of a cofinal map $M' \rightarrow M$, where $M'$ is a small, filtered $\infty$-category. Theorem \ref{sutcar} implies that the diagram $\theta$ has a limit, which we will denote by
$$ \SSet_{/ \hat{X} }$$ and refer to as the {\it $\infty$-topos of local systems on $\hat{X}$}.\index{gen}{$\infty$-topos!of local systems}

\begin{remark}
If $\hat{X}$ is a pro-space, then Proposition \ref{charproet} implies that
the associated geometric morphism $\SSet_{/ \hat{X} } \rightarrow \SSet$ is pro-\'{e}tale. However, the converse is false in general.
\end{remark}

\begin{remark}
Let $G$ be a profinite group, which we may identify with a $\Pro$-object in the category of finite groups. We let $BG$ denote the corresponding $\Pro$-object of $\SSet$, obtained by applying the classifying space functor objectwise. Then $\SSet_{/ BG}$ can be identified with
the $1$-localic $\infty$-topos associated to the ordinary topos of sets with a continuous $G$-action.
It follows from the construction of filtered limits in $\RGeom$ (see \S \ref{inftyfiltlim}) that
we can describe objects $Y \in \SSet_{/BG}$ informally as follows: $Y$ associates to
each open subgroup $U \subseteq G$ a space $Y^{U}$ of {\it $U$-fixed points}, which depends functorially on the finite $G$-space $G/U$. Moreover, if $U$ is a normal subgroup of $V$, then
the natural map from $Y^{V}$ to the (homotopy) fixed point space $(Y^{U})^{V/U}$ should be a homotopy equivalence.
\end{remark}

\begin{remark}
By refining the construction above, it is possible to construct a functor
$$ \Pro(\SSet) \rightarrow \RGeom$$
$$ \hat{X} \mapsto \SSet_{/ \hat{X} }.$$
This functor has a left adjoint, given by
$$ \calX \mapsto \Sh( \calX ).$$
\end{remark}

\begin{warning}
If $\hat{X}$ is a pro-space, then the shape of $\SSet_{/\hat{X}}$ is not necessarily equivalent to $\hat{X}$. In general we have only a counit morphism
$$ \Sh ( \SSet_{/ \hat{X} } ) \rightarrow \hat{X}.$$
\end{warning}

\section{Dimension Theory}\label{dimension}

\setcounter{theorem}{0}

In this section, we will discuss the dimension theory of
topological spaces from the point of view of higher topos theory. We
begin in \S \ref{homdim} by introducing the {\it homotopy dimension} of an $\infty$-topos. We will show that the finiteness
of the homotopy dimension of an $\infty$-topos $\calX$ has pleasant
consequences: it implies that every object is the inverse limit of
its Postnikov tower, and in particular that $\calX$ is hypercomplete.

In \S \ref{chmdim}, we define the {\it cohomology groups} of an $\infty$-topos $\calX$. These cohomology groups have a natural interpretation in terms of the classification of higher gerbes on $\calX$. Using this interpretation, we will show that the cohomology dimension of an $\infty$-topos $\calX$ {\em almost} coincides with its homotopy dimension.

In \S \ref{covdim}, we review the classical theory of {\it covering dimension} for paracompact topological spaces. Using the results of \S \ref{paracompactness}, we will show that the covering dimension of a paracompact space $X$ coincides with the homotopy dimension of the $\infty$-topos $\Shv(X)$.

We conclude in \S \ref{heyt} by introducing a dimension theory for {\em Heyting spaces}, which generalizes the classical theory of Krull dimension for Noetherian topological spaces. Using this theory, we will prove an upper bound for the homotopy dimension of $\Shv(X)$, for suitable Heyting spaces $X$. This result can be regarded as a generalization of Grothendieck's vanishing theorem for the cohomology of Noetherian topological spaces.

\subsection{Homotopy Dimension}\label{homdim}

Throughout this section, we will use the symbol $1_{\calX}$ to denote the final object of an $\infty$-topos
$\calX$.

\begin{definition}\label{sugarpie}\index{gen}{dimension!homotopy}\index{gen}{homotopy dimension}
\index{gen}{homotopy dimension!finite}
Let $\calX$ be an $\infty$-topos. We shall say that $\calX$ has
{\it homotopy dimension $\leq n$} if every $n$-connective
object $U \in \calX$ admits a global section $1_{\calX} \rightarrow U$. We
say that $\calX$ has {\it finite homotopy dimension} if there
exists $n \geq 0$ such that $\calX$ has homotopy dimension $\leq
n$.
\end{definition}

\begin{example}
An $\infty$-topos $\calX$ is of homotopy
dimension $\leq -1$ if and only if $\calX$ is equivalent to
the trivial $\infty$-category $\ast$ (the $\infty$-topos of sheaves on the empty space
$\emptyset$). The ``if'' direction is obvious. Conversely, if $\calX$ has homotopy dimension
$\leq -1$, then the initial object $\emptyset$ of $\calX$ admits a global section
$1_{\calX} \rightarrow \emptyset$. For every object $X \in \calX$, we have a map
$X \rightarrow 1_{\calX} \rightarrow \emptyset$, so that $X$ is also initial (Lemma \ref{sumoto}).
Since the collection of initial objects of $\calX$ span a contractible Kan complex (Proposition \ref{initunique}), we deduce that $\calX$ is itself a contractible Kan complex.
\end{example}

\begin{example}\label{honeypie}
The $\infty$-topos $\SSet$ has homotopy dimension $0$. 
More generally, if $\calC$ is an $\infty$-category with a final object $1_{\calC}$, then
$\calP(\calC)$ has homotopy dimension $\leq 0$. To see this, we first observe that the Yoneda embedding $j: \calC \rightarrow \calP(\calC)$ preserves limits, so that $j(1_{\calC})$ is a final object of $\calP(\calC)$. To prove that $\calP(\calC)$ has homotopy dimension $\leq 0$, we need to show that the functor $\calP(\calC) \rightarrow \SSet$ corepresented by $j(1_{\calC})$ preserves
effective epimorphisms. This functor can be identified with evaluation at $1_{\calC}$. It therefore preserves all limits and colimits, and so carries effective epimorphisms to effective epimorphisms by Proposition \ref{sinn}.
\end{example}

\begin{example}\label{honepie}
Let $X$ be a Kan complex, and let $n \geq -1$. The following conditions are equivalent:
\begin{itemize}
\item[$(1)$] The $\infty$-topos $\SSet_{/X}$ has homotopy dimension $\leq n$.
\item[$(2)$] The geometric realization $|X|$ is a retract (in the homotopy category $\calH$) of a CW complex $K$ of dimension $\leq n$.
\end{itemize}
To prove that $(2) \Rightarrow (1)$, let us choose an $n$-connective
object of $\calX_{/X}$ corresponding to a Kan fibration $p: Y \rightarrow X$, whose homotopy fibers are $n$-connective. Choose a map $K \rightarrow |X|$ which admits a right homotopy inverse. 
To prove that $p$ admits a section up to homotopy, tt will suffice to show that there exists a dotted arrow $$ \xymatrix{ & |Y| \ar[d]^{p} \\
K \ar@{-->}[ur]^{f} \ar[r] & |X| }$$
in the category of topological spaces, rendering the diagram commutative. The construction
of $f$ proceeds cell-by-cell on $K$, using the $n$-connectivity of $p$ to solve lifting problems of the form
$$ \xymatrix{ S^{k-1} \ar@{^{(}->}[d] \ar[r] & |Y| \ar[d]^{p} \\
D^{k} \ar[r] \ar@{-->}[ur] & |X| }$$
for $k \leq n$.

To prove that $(1) \Rightarrow (2)$, we choose any $n$-connective map $q: K \rightarrow |X|$, where
$K$ is an $n$-dimensional CW complex. Condition $(1)$ guarantees that $q$ admits a right homotopy inverse, so that $|X|$ is a retract of $K$ in the homotopy category $\calH$.
\end{example}

\begin{remark}\label{inter}
If $\calX$ is a coproduct (in the $\infty$-category $\RGeom$) of $\infty$-topoi $\calX_{\alpha}$, then
$\calX$ is of homotopy dimension $\leq n$ if and
only if each $\calX_{\alpha}$ is of homotopy dimension $\leq n$.
\end{remark}

It is convenient to introduce a relative version of Definition \ref{sugarpie}.

\begin{definition}\index{gen}{homotopy dimension!of a geometric morphism}
Let $f: \calX \rightarrow \calY$ be a geometric morphism of $\infty$-topoi.
We will say that $f$ is of {\it homotopy dimension $\leq n$} if, for every $k \geq n$ and every
$k$-connective morphism $X \rightarrow X'$ in $\calX$, the induced map
$f_{\ast} X \rightarrow f_{\ast} X'$ is a $(k-n)$-connective morphism in $\calY$ (since $f_{\ast}$ is well-defined up to equivalence, this condition is independent of the choice of $f_{\ast}$).
\end{definition}

\begin{lemma}\label{pie}
Let $\calX$ be an $\infty$-topos, and let $F_{\ast}: \calX \rightarrow \SSet$ be a geometric morphism (which is unique up to equivalence). The following are equivalent:
\begin{itemize}
\item[$(1)$] The $\infty$-topos $\calX$ is of homotopy dimension $\leq n$.
\item[$(2)$] The geometric morphism $F_{\ast}$ is of homotopy dimension $\leq n$.
\end{itemize}
\end{lemma}

\begin{proof}
Suppose first that $(2)$ is satisfied, and let $X$ be an $n$-connective object
of $\calX$. Then $F_{\ast} X$ is a $0$-connective object of $\SSet$: that is, it is a nonempty Kan complex. It therefore has a point $1_{\SSet} \rightarrow F_{\ast} X$. By adjointness, we see that there exists a map $1_{\calX} \rightarrow X$ in $\calX$, where
$1_{\calX} = F^{\ast} 1_{\SSet}$ is a final object of $\calX$ because $F^{\ast}$ is left exact. This proves $(1)$.

Now assume $(1)$, and let
$s: X \rightarrow Y$ be an $k$-connective morphism in $\calX$; we wish to show that
$F_{\ast}s$ is $(k-n)$-connective. The proof goes by induction on $k \geq n$. If $k = n$, then 
are reduced to proving the surjectivity of the horizontal maps in the diagram
$$ \xymatrix{ \pi_0 \bHom_{\calX}( 1_{\calX} , X) \ar[r] \ar@{=}[d] & \pi_0 \bHom_{\calX}( 1_{\calX}, Y) \ar@{=}[d] \\
\pi_0 \bHom_{\SSet}( 1_{\SSet}, F_{\ast} X) \ar[r] & \pi_0 \bHom_{\SSet}(1_{\SSet}, F_{\ast} Y) }$$
of sets. Let $p: 1_{\calX} \rightarrow Y$ be any morphism in $\calX$, and form a pullback diagram
$$ \xymatrix{ Z \ar[d]^-{s'} \ar[r] & X \ar[d]^-{s} \\
1_{\calX} \ar[r]^-{p} & Y. }$$
The map $s'$ is a pullback of $s$, and therefore $n$-connective by Proposition \ref{inftychange}. Using $(1)$, we deduce the existence of a map $1_{\calX} \rightarrow Z$, and a composite map
$$1_{\calX} \rightarrow Z \rightarrow X$$ is a lifting of $p$ up to homotopy.

We now treat the case where $k > n$. Form a diagram
$$ \xymatrix{ X \ar[drr]^-{s'} & & & \\
& & X \times_{Y} X \ar[r] \ar[d] & X \ar[d]^-{s} \\
& & X \ar[r]^-{s} & Y}$$ 
where the square on the bottom-right is a pullback in $\calX$. 
According to Proposition \ref{trowler}, $s'$ is $(k-1)$-connective. Using the inductive hypothesis, we deduce that $F_{\ast}(s')$ is $(k-n-1)$-connective. We now invoke Proposition \ref{trowler} in the $\infty$-topos $\SSet$ deduce that $F_{\ast}(s)$ is $(k-n)$-connective, as desired.
\end{proof}

\begin{definition}\index{gen}{homotopy dimension!locally finite}
We will say that an $\infty$-topos $\calX$ is {\it locally of homotopy dimension $\leq n$} if there exists a collection $\{U_{\alpha} \}$ of objects of $\calX$ which generate $\calX$ under colimits, such that each $\calX_{/U_{\alpha} }$ is of homotopy dimension $\leq n$.
\end{definition}

\begin{example}
Let $\calC$ be a small $\infty$-category. Then $\calP(\calC)$ is locally of homotopy dimension $\leq 0$. To prove this, we first observe that $\calP(\calC)$ is generated under colimits by the Yoneda embedding $j: \calC \rightarrow \calP(\calC)$. It therefore suffices to prove that
each of the $\infty$-topoi $\calP(\calC)_{/j(C)}$ has finite homotopy dimension. According to Corollary \ref{swapKK}, the $\infty$-topos $\calP(\calC)_{/j(C)}$ is equivalent to $\calP( \calC_{/C})$, which is of homotopy dimension $0$ (see Example \ref{honeypie}). 
\end{example}

Our next goal is to prove the following result:

\begin{proposition}\label{cumba}
Let $\calX$ be an $\infty$-topos which is locally of homotopy dimension $\leq n$ for some integer $n$. Then Postnikov towers in $\calX$ are convergent.
\end{proposition}

\begin{proof}
We will show that $\calX$ satisfies the criterion of Remark \ref{urkan}. Let $X: \Nerve( \Z^{\infty}_{\geq 0} )^{op} \rightarrow \calX$ be a limit tower, and assume that the underlying pretower is highly connected. We wish to show that $X$ is highly connected. Choose $m \geq -1$; we wish to show that the map $X(\infty) \rightarrow X(k)$ is $m$-connective for $k \gg 0$.
Reindexing the tower if necessary, we may suppose that for every $p \geq q$, the map
$X(p) \rightarrow X(q)$ is $(m+q)$-connective. We claim that, in this case, we can take $k = 0$. 
The proof goes by induction on $m$. If $m > 0$, we can deduce the desired result by applying the inductive hypothesis to the tower
$$ X(\infty) \rightarrow \ldots X(\infty) \times_{ X(2) } X(\infty) \rightarrow X(\infty) \times_{ X(1) } X(\infty)
\rightarrow X(\infty) \times_{ X(0) } X(\infty). $$
Let us therefore assume that $m = 0$; we wish to show that the map $X(\infty) \rightarrow X(0)$ is an effective epimorphism. Since the objects $\{ U_{\alpha} \}$ generate $\calX$ under colimits, there
is an effective epimorphism $\phi: U \rightarrow X(0)$, where $U$ is a coproduct of objects of the form
$\{ U_{\alpha} \}$. Using Remark \ref{inter}, we deduce that $\calX_{/U}$ has homotopy dimension $\leq n$. Let $F: \calX \rightarrow \SSet$ denote the functor corepresented by $U$. Then
$F$ factors as a composition
$$ \calX \stackrel{f^{\ast}}{\rightarrow} \calX_{/U} \stackrel{\Gamma}{\rightarrow} \SSet,$$
where $f^{\ast}$ is the left adjoint to the geometric morphism $\calX_{/U} \rightarrow \calX$
and $\Gamma$ is the global sections functor. It follows that $F$ carries $n$-connective morphisms to effective epimorphisms (Lemma \ref{pie}). The map $\phi$ determines a point of
$F( X(0) )$. Since each of the maps $F( X(k+1)  ) \rightarrow F( X(k) )$ induces a surjection on
connected components, we can lift this point successively to each $F( X(k) )$ and thereby obtain
a point in $F( X(\infty) ) \simeq \holim \{ F( X(n) ) \}$. This point determines a diagram
$$ \xymatrix{ & X(\infty) \ar[dr]^{\psi} & \\
U \ar[ur] \ar[rr]^{\phi} & & X(0) }$$
which commutes up to homotopy. Since $\phi$ is an effective epimorphism, we deduce that
the map $\psi$ is an effective epimorphism, as desired.
\end{proof}

\begin{lemma}\label{tow}
Let $\calX$ be a presentable $\infty$-category, let $\Fun(\Nerve (\Z_{\geq 0}^{\infty})^{op}, \calX)$ be the $\infty$-category of towers in $\calX$, and let $\calX_{\tau} \subseteq 
\Fun( \Nerve (\Z_{\geq 0}^{\infty})^{op}, \calX)$ denote the full subcategory spanned by the Postnikov towers. Evaluation at $\infty$ induces a trivial fibration of simplicial sets
$ \calX_{\tau} \rightarrow \calX$. In particular, every object $X(\infty) \in \calX$ can be extended to a Postnikov tower 
$$ X(\infty) \rightarrow \ldots \rightarrow X(1) \rightarrow X(0).$$
\end{lemma}

\begin{proof}
Let $\calC$ be the full subcategory of $\calX \times \Nerve(\Z_{\geq 0}^{\infty})^{op}$
spanned by the pairs $(X, n)$ where $X$ is an object of $\calX$, $n \in \Z_{\geq 0}^{\infty}$,
and $X$ is $n$-truncated, and let $p: \calC \rightarrow \calX$ denote the natural projection.
Since every $m$-truncated object of $\calX$ is also $n$-truncated for $m \geq n$, it is easy to see that $p$ is a Cartesian fibration. Proposition \ref{maketrunc} implies that each of the inclusion
functors $\tau_{\leq m} \calX \subseteq \tau_{\leq n} \calX$ has a left adjoint, so that $p$ is also a coCartesian fibration (Corollary \ref{getcocart}). By definition, $\calX_{\tau}$ can be identified
with the simplicial set
$$ \bHom^{\flat}_{\Nerve(\Z_{\geq 0}^{\infty})}( \Nerve(\Z_{\geq 0}^{\infty})^{\sharp}, 
(\calC^{op})^{\natural} )^{op}$$
and $\calX$ itself can be identified with
$$ \bHom^{\flat}_{\Nerve(\Z_{\geq 0}^{\infty})}( \{\infty\}^{\sharp}, 
(\calC^{op})^{\natural} )^{op}.$$
It now suffices to observe that the inclusion
$\{ \infty\}^{\sharp} \subseteq \Nerve(\Z_{\geq 0}^{\infty})^{\sharp}$ is marked anodyne.
\end{proof}

\begin{corollary}[Jardine]\label{fdfd}
Let $\calX$ be an $\infty$-topos which is locally of homotopy dimension $\leq n$ for some integer $n$. Then $\calX$ is hypercomplete.
\end{corollary}

\begin{proof}
Let $X(\infty)$ be an arbitrary object of $\calX$. By Lemma \ref{tow} we can find a Postnikov tower
$$ X(\infty) \rightarrow \ldots \rightarrow X(1) \rightarrow X(0).$$
Since $X(n)$ is $n$-truncated, it belongs to $\calX^{\hyp}$ by Corollary \ref{goober2}.
By Proposition \ref{cumba}, the tower exhibits $X(\infty)$ as a limit of objects
of $\calX^{\hyp}$, so that $X(\infty)$ belongs to $\calX^{\hyp}$ as well since
the full subcategory $\calX^{\hyp} \subseteq \calX$ is stable under limits.
\end{proof}

\begin{lemma}\label{nicelemma}
Let $\calX$ be an $\infty$-topos, $n \geq 0$, $X$ an
$(n+1)$-connective object of $\calX$, and $f^{\ast}: \calX \rightarrow \calX_{/X}$ a right adjoint to the projection $\calX_{/X} \rightarrow \calX$.
Then $f^{\ast}$ induces a fully faithful functor
$\tau_{\leq n} \calX \rightarrow \tau_{\leq n} \calX_{/X}$ which restricts to an equivalence from $\tau_{\leq n-1} \calX$ to $\tau_{\leq n-1} \calX_{/X}$
\end{lemma}

\begin{proof}
We first prove that $f^{\ast}$ is fully faithful when restricted to the
$\infty$-category of $n$-truncated objects of $\calX$. Let $Y,
Z \in \calX$ be objects, where $Y$ is $n$-truncated. We have a commutative
diagram 
$$ \xymatrix{ \bHom_{\calX_{/X}}(f^{\ast} Y, f^{\ast} Z) \ar@{=}[r] & 
 \bHom_{\calX}( X \times Y, Z) &  \bHom_{\calX}( \tau_{\leq n}(X \times Y), Z) \ar[l] \\
 \bHom_{\calX}(Y,Z) \ar[u] & & \bHom_{\calX}( \tau_{\leq n} Y, Z) \ar[u] \ar[ll] }$$
in the homotopy category $\calH$, where the horizontal arrows are homotopy equivalences.
Consequently, to prove that the left vertical map is a homotopy equivalence, it suffices to show that the projection $\tau_{\leq n} (X \times Y) \rightarrow \tau_{\leq n} Y$ is an equivalence.
This follows immediately from Lemma \ref{slurpy} and our assumption that
$X$ is $(n+1)$-connective.

Now suppose that $\overline{Y}$ is an $(n-1)$-truncated object of $\calX_{/X}$. We wish to show that $\overline{Y}$ lies in the essential image of $f^{\ast} | \tau_{\leq n-1} \calX$. Let $Y$
denote the image of $\overline{Y}$ in $\calX$, and let $Y \rightarrow Z$ exhibit $Z$ as
an $(n-1)$-truncation of $Y$ in $\calX$. To complete the proof, it will suffice to show that the composition
$$ u: \overline{Y} \stackrel{u'}{\rightarrow} f^{\ast} Y \stackrel{u''}{\rightarrow} f^{\ast} Z$$
is an equivalence in $\calX_{/X}$. Since both $\overline{Y}$ and $f^{\ast} Z$
are $(n-1)$-truncated, it suffices to prove that $u$ is $n$-connective. According to Proposition \ref{inftychange},
it suffices to prove that $u'$ and $u''$ are $n$-connective. Proposition \ref{compattrunc}
implies that $u''$ exhibits $f^{\ast} Z$ as an $(n-1)$-truncation of $f^{\ast} Y$, and is therefore
$n$-connective. 

We now complete the proof by showing that $u'$ is $n$-connective.
Let $v'$ denote the image of image of $u'$ in the $\infty$-topos $\calX$. According to
Proposition \ref{conslice}, it will suffice to show that $v'$ is $n$-connective.
We observe that $v'$ is a section of the projection $q: Y \times X \rightarrow Y$.
$q: Y \times X \rightarrow Y$. According to Proposition \ref{sectcon}, it will suffice to prove that
$q$ is $(n+1)$-connective. Since $q$ is a pullback of the projection $X \rightarrow 1_{\calX}$, Proposition \ref{inftychange} allows us to conclude the proof (since $X$ is $(n+1)$-connective, by assumption).
\end{proof}

Lemma \ref{nicelemma} has some pleasant consequences.

\begin{proposition}\label{pi00detects}
Let $\calX$ be an $\infty$-topos and let $\tau_{\leq 0} : \calX \rightarrow \tau_{\leq 0} \calX$
denote a left adjoint to the inclusion. A morphism $\phi: U \rightarrow X$ in $\calX$ is an effective epimorphism if and only if $\tau_{\leq 0}(\phi)$ is an effective epimorphism in the ordinary topos
$\h{(\tau_{\leq 0} \calX)}$.
\end{proposition}

\begin{proof}
Suppose first that $\phi$ is an effective epimorphism. Let $U_{\bigdot}: \Nerve \cDelta_{+}^{op} \rightarrow \calX$ be a \Cech nerve of $\phi$, so that $U_{\bigdot}$ is a colimit diagram.
Since $\tau_{\leq 0}$ is a left adjoint,
$\tau_{\leq 0} U_{\bigdot}$ is a colimit diagram in $\tau_{\leq 0} \calX$. 
Using Proposition \ref{charsurj}, we deduce easily that $\tau_{\leq 0} \phi$ is an effective epimorphism.

For the converse, choose a factorization of $\phi$ as a composition
$$ U \stackrel{\phi'}{\rightarrow} V \stackrel{\phi''}{\rightarrow} X$$
where $\phi'$ is an effective epimorphism and $\phi''$ is a monomorphism. Applying Lemma \ref{nicelemma} to the $\infty$-topos $\calX_{/ \tau_{\leq 0} X}$, we conclude that
$\phi''$ is the pullback of a monomorphism $i: \overline{V} \rightarrow \tau_{\leq 0} X$.
Since the effective epimorphism $\tau_{\leq 0}(\phi)$ factors through $i$, we conclude
that $i$ is an equivalence, so that $\phi''$ is likewise an equivalence. It follows that
$\phi$ is an effective epimorphism as desired.
\end{proof}

Proposition \ref{pi00detects} can be regarded as a generalization of the following well-known property of the $\infty$-category of spaces, which can itself be regarded as the $\infty$-categorical analogue of the second part of Fact \ref{factoid}:

\begin{corollary}
Let $f: X \rightarrow Y$ be a map of Kan complexes. Then $f$ is an effective epimorphism
in the $\infty$-category $\SSet$ if and only if the induced map $\pi_0 X \rightarrow \pi_0 Y$ is surjective.
\end{corollary}

\begin{remark}\label{charnice}
It follows from Proposition \ref{pi00detects} that the class of $\infty$-topoi having the form
$\Shv(\calC)$, where $\calC$ is a small $\infty$-category, is not substantially larger than the class of ordinary topoi. More precisely, every topological localization of $\calP(\calC)$ can be obtained by inverting
morphisms between {\em discrete} objects of $\calP(\calC)$. It follows that there exists a pullback diagram of $\infty$-topoi
$$ \xymatrix{ \Shv(\calC) \ar[r] \ar[d] & \calP(\calC) \ar[d] \\
\Shv(\Nerve( \h{\calC})) \ar[r] & \calP(\Nerve(\h{\calC})) }$$
where the $\infty$-topoi on the bottom line are $1$-localic, and therefore determined by the
ordinary topoi of presheaves of sets on the homotopy category $h\calC$ and sheaves of sets on $h\calC$, respectively.
\end{remark}

\begin{corollary}\label{enuff}
Let $X$ be a topological space. Suppose that $\Shv(X)$ is locally of homotopy dimension $\leq n$ for some integer $n$. Then $\Shv(X)$ has enough points.
\end{corollary}

\begin{proof}
Note that every point $x \in X$ gives rise to a point $x_{\ast}: \Shv(\ast) \rightarrow \Shv(X)$
of the $\infty$-topos $\Shv(X)$. Let $f: \calF \rightarrow \calF'$ be a morphism in
$\Shv(X)$ such that $x^{\ast}(f)$ is an equivalence in $\SSet$ for each $x \in X$. We wish to prove that $f$ is an equivalence. According to Corollary \ref{fdfd}, it will suffice to prove that $f$ is $\infty$-connective. We will prove by induction on $n$ that $f$ is $n$-connective. If
$n > 0$, we simply apply the inductive hypothesis to the diagonal morphism
$\delta: \calF \rightarrow \calF \times_{\calF'} \calF$. We may therefore reduce to the case
$n=0$; we wish to show that $f$ is an effective epimorphism. Since $\Shv(X)$ is generated under colimits by the sheaves $\chi_U$ associated to open subsets $U \subseteq X$,  we may assume without loss of generality that $\calF' = \chi_{U}$. We may now invoke Proposition \ref{pi00detects} to reduce to the case where $\calF$ is an object of $\tau_{\leq 0} \Shv(X)_{/\chi_U}$. This $\infty$-category is equivalent to the nerve of the category of sheaves of {\em sets} on $U$. We are therefore reduced to proving that if $\calF$ is a sheaf of sets on $U$ whose stalk
$\calF_{x}$ is a singleton at each point $x \in U$, then $\calF$ has a global section, which is clear.
\end{proof}

\subsection{Cohomological Dimension}\label{chmdim}

In classical homotopy theory, one can analyze a space $X$ by means of its Postnikov tower
$$ \ldots \tau_{\leq n} X \stackrel{\phi_n}{\rightarrow} \tau_{\leq n-1} X \rightarrow \ldots.$$
In this diagram, the homotopy fiber $F$ of $\phi_n$ $(n \geq 1)$ is a space which has only a single nonvanishing homotopy group, in dimension $n$. The space $F$ is determined up to homotopy equivalence by $\pi_n F$: in fact $F$ is homotopy equivalent to an Eilenberg-MacLane space
$K( \pi_n F, n)$ which can be functorially constructed from the group $\pi_n F$. The study of these Eilenberg-MacLane spaces is of central interest, because (according to the above analysis) they constitute basic building blocks out of which any arbitrary space can be constructed.
Our goal in this section is to generalize the theory of Eilenberg-MacLane spaces to the setting of 
an arbitrary $\infty$-topos $\calX$.

\begin{definition}\label{gropab}\index{gen}{object!pointed}\index{gen}{pointed object}\index{not}{Xast@$\calX_{\ast}$}
Let $\calX$ be an $\infty$-category. A {\it pointed object} is a morphism
$X_{\ast}: 1 \rightarrow X$ in $\calX$, where $1$ is a final object of $\calX$.
We let $\calX_{\ast}$ denote the full subcategory of $\Fun(\Delta^1,\calX)$ spanned by the pointed objects of
$\calX$.

A {\it group object} of $\calX$ is a groupoid object $U_{\bigdot}: \Nerve \cDelta^{op} \rightarrow \calX$ for which $U_0$ is a final object of $\calX$. Let $\Group(\calX)$ denote the full subcategory
of $\calX_{\Delta}$ spanned by the group objects of $\calX$.\index{gen}{object!group}\index{gen}{group object}\index{not}{groupcalX@$\Group(\calX)$}

We will say that a pointed object $1 \rightarrow X$ of an $\infty$-topos $\calX$ is an {\it Eilenberg-MacLane object of degree $n$} if $X$ is $n$-truncated and $n$-connective. We let $\EM_n(\calX)$ denote the full subcategory of $\calX_{\ast}$ spanned by the Eilenberg-MacLane objects of degree $n$.\index{gen}{object!Eilenberg-MacLane}\index{gen}{Eilenberg-MacLane object}\index{not}{EMncalX@$\EM_n(\calX)$}
\end{definition}

\begin{example}
Let $\calC$ be an ordinary category which admits finite limits. A {\it group object}
of $\calC$ is an object $X \in \calC$ which is equipped with an identity section
$1_{\calC} \rightarrow X$, an inversion map $X \rightarrow X$, and
a multiplication $m: X \times X \rightarrow X$, which satisfy the usual group axioms. Equivalently, a group object of $\calC$ is an object $X$ together with a group structure on each morphism
space $\Hom_{\calC}(Y,X)$, which depends functorially on $Y$. We will denote the category
of group objects of $\calC$ by $\Group(\calC)$. The $\infty$-category
$\Nerve(\Group(\calC))$ is equivalent to the $\infty$-category of group objects of
$\Nerve(\calC)$, in the sense of Definition \ref{gropab}. Thus, the notion of a group object
of an $\infty$-category can be regarded as a generalization of the notion of a group object of an ordinary category.
\end{example}

\begin{remark}\label{prodem}
Let $\calX$ be an $\infty$-topos and $n \geq 0$ an integer. Then the full subcategory of
$\Fun(\Delta^1,\calX)$ consisting of Eilenberg-MacLane objects $p: 1 \rightarrow X$
is stable under finite products. 
This is clear, since:
\begin{itemize}
\item[$(1)$] A finite product of $n$-connective objects of $\calX$ is $n$-connective (Corollary \ref{togoto}).
\item[$(2)$] {\em Any} limit of $n$-truncated objects of $\calX$ is $n$-truncated (since 
$\tau_{\leq n} \calX$ is a localization of $\calX$).
\end{itemize}
\end{remark}

\begin{proposition}\label{tinner}
Let $\calX$ be an $\infty$-category, and let $U_{\bigdot}$ be a simplicial object of $\calX$. 
Then $U_{\bigdot}$ is a group object of $\calX$ if and only if the following conditions are satisfied:
\begin{itemize}
\item[$(1)$] The object $U_0$ is final in $\calX$.
\item[$(2)$] For every decomposition $[n] = S \cup S'$, where
$S \cap S' = \{s\}$, the maps
$$ U(S) \leftarrow U_n \rightarrow U(S')$$
exhibit $U_n$ as a product of $U(S)$ and $U(S)$ in $\calX$.
\end{itemize}
\end{proposition}

\begin{proof}
This follows immediately from characterization $(4'')$ of Proposition \ref{grpobjdef}.
\end{proof}

\begin{corollary}\label{grpstable}
Let $\calX$ and $\calY$ be $\infty$-categories which admit finite products, and let
$f: \calX \rightarrow \calY$ be a functor which preserves finite products. Then
the induced functor $\calX_{\Delta} \rightarrow \calY_{\Delta}$ carries group objects
of $\calX$ to group objects of $\calY$.
\end{corollary}

\begin{corollary}
Let $\calX$ be an $\infty$-category which admits finite products, and let $\calY \subseteq \calX$ be a full subcategory which is stable under finite products. Let $Y_{\bigdot}$ be a simplicial
object of $\calY$. Then $Y_{\bigdot}$ is a group object of $\calY$ if and only if it is a group object of $\calX$. 
\end{corollary}

\begin{definition}
Let $\calX$ be an $\infty$-category. A {\it zero object} of $\calX$ is an object which is both initial and final.
\end{definition}

\begin{lemma}\label{pointer}
Let $\calX$ be an $\infty$-category with a final object $1_{\calX}$. Then the inclusion
$i: \calX^{1_{\calX}/} \subseteq \calX_{\ast}$ is an equivalence of $\infty$-categories.
\end{lemma}

\begin{proof}
Let $K$ be the full subcategory of $\calX$ spanned by the final objects, and let
$1_{\calX}$ be an object of $K$. Proposition \ref{initunique} implies that $K$ is a contractible Kan complex, so that the inclusion $\{1_{\calX} \} \subseteq K$ is an equivalence of $\infty$-categories.
Corollary \ref{tweezegork} implies that the projection 
$\calX_{\ast} \rightarrow K$ is a coCartesian fibration. We now apply Proposition \ref{basechangefunky} to deduce the desired result.
\end{proof}

\begin{lemma}\label{pointerprime}
Let $\calX$ be an $\infty$-category with a final object. Then the $\infty$-category
$\calX_{\ast}$ has a zero object. If $\calX$ already has a zero object, then the forgetful functor
$\calX_{\ast} \rightarrow \calX$ is an equivalence of $\infty$-categories.
\end{lemma}

\begin{proof}
Let $1_{\calX}$ be a final object of $\calX$, and let $U = \id_{1_{\calX}} \in \calX_{\ast}$. 
We wish to show that $U$ is a zero object of $\calX_{\ast}$. According to Lemma \ref{pointer}, 
it will suffice to show that $U$ is a zero object
of $\calX_{1_{\calX}/}$. It is clear that $U$ is initial, and the finality of $U$ follows from
Proposition \ref{needed17}. 

For the second assertion, let us suppose that $1_{\calX}$ is also an initial object of $\calX$.
We wish to show that the forgetful functor $\calX_{\ast} \rightarrow \calX$ is an equivalence of $\infty$-categories. Applying Lemma \ref{pointer}, it will suffice to show that the projection
$f: \calX^{1_{\calX}/} \rightarrow \calX$ is an equivalence of $\infty$-categories. 
But $f$ is a trivial fibration of simplicial sets.
\end{proof}

\begin{lemma}\label{postEM}
Let $\calX$ be an $\infty$-category, and let $f: \calX_{\ast} \rightarrow \calX$
be the forgetful functor $($ which carries a pointed object $1 \rightarrow X$ to $X$ $)$.
Then $f$ induces an equivalence of $\infty$-categories
$$ \Group( \calX_{\ast}) \rightarrow \Group( \calX).$$
\end{lemma}

\begin{proof}
The functor $f$ factors as a composition
$$ \calX_{\ast} \subseteq \Fun(\Delta^1,\calX) \rightarrow \calX$$
where the first map is the inclusion of a full subcategory which is stable under limits,
and the second map preserves all limits (Proposition \ref{limiteval}). It follows that $f$ preserves limits, and therefore composition with $f$ induces a functor $F: \Group(\calX_{\ast}) \rightarrow \Group(\calX)$ by Corollary \ref{grpstable}. 

Observe that the $0$-simplex $\Delta^0$ is an initial object of $\cDelta^{op}$. Consequently, there exists a functor $T: \Delta^1 \times \Nerve(\cDelta)^{op} \rightarrow \Nerve(\cDelta)^{op}$, which is a natural transformation from the constant functor taking the value $\Delta^0$ to the identity functor. Composition with $T$ induces a functor
$$ \calX_{\Delta} \rightarrow \Fun(\Delta^1,\calX)_{\Delta}.$$
Restricting to group objects, we get a functor $s: \Group(\calX) \rightarrow \Group(\calX_{\ast})$.
It is clear that $F \circ s$ is the identity. 

We observe that if $\calX$ has a zero object, then
$f$ is an equivalence of $\infty$-categories (Lemma \ref{pointerprime}). It follows immediately that $F$ is an equivalence of $\infty$-categories. Since $s$ is a right inverse to $F$, we conclude that $s$ is an equivalence of $\infty$-categories as well.

To complete the proof in the general case, it will suffice to show that the composition $s \circ F$ is an equivalence of $\infty$-categories. To prove this, we set $\calY = \calX_{\ast}$, and let
$F': \Group(\calY_{\ast}) \rightarrow \Group(\calX_{\ast})$ and $s': \Group(\calY) \rightarrow \Group(\calY_{\ast})$ be defined as above. We then have a commutative diagram
$$ \xymatrix{ \Group(\calY) \ar[r]^-{F} \ar[d]^-{s'} & \Group(\calX) \ar[d]^-{s} \\
\Group(\calY_{\ast}) \ar[r]^-{F'} & \Group(\calX_{\ast}) }$$ 
so that $s \circ F = F' \circ s'$. Lemma \ref{pointerprime} implies that $\calY$ has a zero object, so that $F'$ and $s'$ are equivalences of $\infty$-categories. Therefore $F' \circ s' = s \circ F$ is an equivalence of $\infty$-categories, and the proof is complete.
\end{proof}

The following Proposition guarantees a good supply of Eilenberg-MacLane objects in an $\infty$-topos $\calX$.

\begin{lemma}\label{preEM}
Let $\calX$ be an $\infty$-topos containing a final object $1_{\calX}$ and let $n \geq 1$.
Let $p$ denote the composition
$$ \Fun(\Delta^1,\calX) \stackrel{\mCech}{\rightarrow} \calX_{\cDelta_{+}} \rightarrow \calX_{\cDelta}$$
which associates to each morphism $U \rightarrow X$ the underlying groupoid of its \Cech nerve.
Then:
\begin{itemize}
\item[$(1)$] Let $\calX'$ denote the full subcategory of $\Fun(\Delta^1,\calX)$ consisting
of {\em connected} pointed objects of $\calX$. Then the restriction of $p$ induces an equivalence of $\infty$-categories from $\calX'$ to the $\infty$-category $\Group(\calX)$.

\item[$(2)$] The essential image of $p| \EM_{n}(\calX)$ coincides with
the essential image of the composition
$$ \Group( \EM_{n-1}(\calX)) \subseteq \Group( \calX_{\ast}) \rightarrow \Group(\calX).$$
\end{itemize}
\end{lemma}

\begin{proof}
Let $\calX''$ be the full subcategory of $\Fun(\Delta^1,\calX)$ spanned by the effective epimorphisms
$u: U \rightarrow X$. Since $\calX$ is an $\infty$-topos, $p$ induces an equivalence
from $\calX''$ to the $\infty$-category of groupoid objects of $\calX$. Consequently, to prove
$(1)$, it will suffice to show that if $u: 1_{\calX} \rightarrow X$ is a morphism in $\calX$ and
$1_{\calX}$ is a final object, then $u$ is an effective epimorphism if and only if $X$ is connected.
We note that $X$ is connected if and only if the map $\tau_{\leq 0}(u):
\tau_{\leq 0} 1_{\calX} \rightarrow \tau_{\leq 0} X$ is an isomorphism in the ordinary
topos $\Disc(\calX)$. According to Proposition \ref{pi00detects}, $u$
is an effective epimorphism if and only if $\tau_{\leq 0}(u)$ is an effective epimorphism.
We now observe that in any ordinary category $\calC$, an effective epimorphism
$u': 1_{\calC} \rightarrow X'$ whose source is a final object of $\calC$ is automatically an isomorphism, since the equivalence relation $1_{\calC} \times_{X'} 1_{\calC} \subseteq 1_{\calC} \times 1_{\calC}$ automatically consists of the whole of $1_{\calC} \times 1_{\calC} \simeq 1_{\calC}$.

To prove $(2)$, we consider an augmented simplicial object $X_{\bigdot}$ of $\calX$
which is a \Cech nerve, having the property that $X_0$ is a final object of $\calX$.
We wish to show that the pointed object $X_0 \rightarrow X_{-1}$ belongs to
$\EM_{n}(\calX)$ if and only if each $X_{k}$ is $(n-1)$-truncated and $(n-1)$-connective,
for $k \geq 0$. We conclude by making the following observations:
\begin{itemize}
\item[$(a)$] Since $X_{k}$ is equivalent to a $k$-fold product of copies of $X_{1}$, the
objects $X_{k}$ are $(n-1)$-truncated ($(n-1)$-connective) for all $k \geq 0$ if and only if
$X_{1}$ is $(n-1)$-truncated ($(n-1)$-connective).

\item[$(b)$] We have a pullback diagram
$$ \xymatrix{ X_{1} \ar[r]^-{f} \ar[d] & X_0 \ar[d]^-{g} \\
X_0 \ar[r]^-{g} & X_{-1}. }$$
The object $X_{1}$ is $(n-1)$-truncated if and only if $f$ is $(n-1)$-truncated.
Since $g$ is an effective epimorphism, $f$ is $(n-1)$-truncated if and only if $g$ is $(n-1)$-truncated (Proposition \ref{hintdescent0}). Using the long exact sequence of Remark \ref{sequence}, we conclude that this is equivalent to the vanishing of $g^{\ast} \pi_k X_{-1}$ for $k > n$.
Since $g$ is an effective epimorphism, this is equivalent to the vanishing of
$\pi_k X_{-1}$ for $k > n$; in other words, to the requirement that $X_{-1}$ is $n$-truncated.

\item[$(c)$] The object $X_{1}$ is $(n-1)$-connective if and only if
$f$ is $(n-1)$-connective. Arguing as above, we conclude that $f$ is $(n-1)$-connective if and only if $g$ is $(n-1)$-connective (Proposition \ref{inftychange}). Using the long exact sequence of Remark \ref{sequence}, this is equivalent to the vanishing of the homotopy sheaf $g^{\ast} \pi_k X_{-1}$ for $k < n$. Since $g$ is an effective epimorphism, this is equivalent to the vanishing of $\pi_k X_{-1}$ for $k < n$; in other words, to the condition that $X_{-1}$ is $(n-1)$-truncated.
\end{itemize}

\end{proof}

\begin{proposition}\label{EM}
Let $\calX$ be an $\infty$-topos and $n \geq 0$ a nonnegative integer, and let 
$\pi_n: \calX_{\ast} \rightarrow \Nerve(\Disc(\calX))$ denote the associated homotopy group functor.

Then:
\begin{itemize}
\item[$(1)$] If $n = 0$, then $\pi_n$ determines an equivalence from the $\infty$-category $\EM_{0}(\calX)$ to the $($nerve of the$)$ category of pointed objects of $\Disc(\calX)$.
\item[$(2)$] If $n = 1$, then $\pi_n$ determines an equivalence from the $\infty$-category $\EM_{1}(\calX)$ to the $($nerve of the$)$ category of group objects of $\Disc(\calX)$.
\item[$(3)$] If $n \geq 2$, then $\pi_n$ determines an equivalence from the $\infty$-category $\EM_{n}(\calX)$ to the $($nerve of the$)$ category of commutative group objects of $\Disc(\calX)$.
\end{itemize}
\end{proposition}

\begin{proof}
We use induction on $n$. The case $n=0$ follows immediately from the definitions. The
case $n=1$ follow from the case $n=0$, by combining Lemmas \ref{preEM} and \ref{postEM}.
If $n=2$, we apply the inductive hypothesis, together with Lemma \ref{preEM} and the observation
that if $\calC$ is an ordinary category which admits finite products, then $\Group( \Group(\calC))$ is equivalent to category $\Ab(\calC)$ of {\em commutative} group objects of $\calC$. The argument
in the case $n > 2$ makes use of the inductive hypothesis, Lemma \ref{preEM}, and the observation that $\Group( \Ab(\calC) )$ is equivalent to $\Ab(\calC)$ for any ordinary category $\calC$ which admits finite products.
\end{proof}

Fix an $\infty$-topos $\calX$, a final object $1_{\calX} \in \calX$, and an integer $n \geq 0$. According to Proposition \ref{EM}, there exists a homotopy inverse to the functor $\pi$. We will denote this functor by
$$ A \mapsto (p: 1_{\calX} \rightarrow K(A,n) )$$
where $A$ is a pointed object of the topos $\Disc(\calX)$ if $n = 0$, a group object
if $n =1$, and an abelian group object if $n \geq 2$.\index{not}{K(A,n)@$K(A,n)$}

\begin{remark}\label{produm}
The functor $A \mapsto K(A,n)$ preserves finite products. This is clear, since the class of Eilenberg-MacLane objects is stable under finite products (Remark \ref{prodem}) and the homotopy inverse functor $\pi$ commutes with finite products (since homotopy groups are constructed using pullback and truncation functors, each of which commutes with finite products).
\end{remark}

\begin{definition}
Let $\calX$ be an $\infty$-topos, $n \geq 0$ an integer, and $A$ an abelian group object of the topos $\Disc(\calX)$.
We define 
$$ \HH^{n}(\calX; A) = \pi_0 \bHom_{\calX}( 1_{\calX}, K(A,n) ); $$
we refer to $\HH^n(\calX; A)$ as the {\it $n$th cohomology group of $\calX$ with coefficients in $A$}.\index{not}{HncalXA@$\HH^{n}(\calX;A)$}\index{gen}{cohomology group!of an $\infty$-topos}
\end{definition}

\begin{remark}
It is clear that we can also make sense of $\HH^{1}(\calX; G)$ when $G$ is a sheaf of nonabelian groups, or $\HH^0(\calX; E)$ when $E$ is only a sheaf of (pointed) sets.
\end{remark}

\begin{remark}
It is clear from the definition that $\HH^{n}(\calX; A)$ is functorial in $A$. Moreover, this functor
commutes with finite products by Remark \ref{produm} (and the fact that products in $\calX$
are products in the homotopy category $h \calX$). If $A$ is an abelian group, then the multiplication
map $A \times A \rightarrow A$ induces a (commutative) group structure on $\HH^{n}(\calX;A)$. This justifies our terminology in referring to $\HH^{n}(\calX;A)$ as a cohomology {\em group}.
\end{remark}

\begin{remark}\label{compfood}
Let $\calC$ be a small category equipped with a Grothendieck topology, and let
$\calX$ be the $\infty$-topos $\Shv( \Nerve \calC)$ of sheaves of {\em spaces} on
$\calC$, so that the underlying topos $\Disc(\calX)$ is equivalent to the category of sheaves of {\em sets} on $\calC$. Let $A$ be a sheaf of abelian groups on $\calC$.
Then $\HH^{n}(\calX;A)$ may be identified with the $n$th cohomology group of $\Disc(\calX)$ with coefficients in $A$, in the sense of ordinary sheaf theory. To see this, choose a resolution
$$ A \rightarrow I^0 \rightarrow I^1 \rightarrow \ldots \rightarrow I^{n-1} \rightarrow J$$
of $A$ by abelian group objects of $\Disc(\calX)$, where each $I^{k}$ is injective. The complex
$$ I^0 \rightarrow \ldots \rightarrow J$$
may be identified, via the Dold-Kan correspondence, with a simplicial abelian group object
$C_{\bigdot}$ of $\Disc(\calX)$. Regard $C_{\bigdot}$ as a presheaf on $\calC$ with values in $\sSet$. Then:
\begin{itemize}
\item[$(1)$] The induced presheaf $F: \Nerve(\calC)^{op} \rightarrow \SSet$ belongs to
$\calX = \Shv(\Nerve(\calC)) \subseteq \calP( \Nerve(\calC))$ (this uses the injectivity
of the objects $I^k$) and is equipped with a canonical basepoint $p: 1_{\calX} \rightarrow F$.
\item[$(2)$] The pointed object $p: 1_{\calX} \rightarrow F$ is an Eilenberg-MacLane object of $\calX$, and there is a canonical identification $A \simeq p^{\ast}(\pi_n F)$. We may therefore identify
$F$ with $K(A,n)$.
\item[$(3)$] The set of homotopy classes of maps from $1_{\calX}$ to $F$ in $\calX$ may be identified with
the cokernel of the map $\Gamma( \Disc(\calX); I^{n-1} ) \rightarrow \Gamma( \Disc(\calX); J)$, which is
also the $n$th cohomology group of $\Disc(\calX)$ with coefficients in $A$ in the sense of classical sheaf theory.
\end{itemize} 
For further discussion of this point, we refer the reader to \cite{jardine}.
\end{remark}

We are ready to define the cohomological dimension of an $\infty$-topos.

\begin{definition}\label{codimmm}\index{gen}{dimension!cohomological}\index{gen}{cohomoogical dimension}
Let $\calX$ be an $\infty$-topos. We will say that $\calX$ has
{\it cohomological dimension $\leq n$} if, for any sheaf of
abelian groups $A$ on $\calX$, the cohomology group $\HH^k(\calX,A)$ vanishes for
$k > n$.
\end{definition}

\begin{remark}
For small values of $n$, some authors prefer to require a stronger
vanishing condition which applies also when $A$ is a non-abelian
coefficient system. The appropriate definition requires the
vanishing of cohomology for coefficient systems which are defined
only up to inner automorphisms, as in \cite{giraud}. With the
appropriate modifications, Theorem \ref{cohdim} below remains
valid for $n < 2$.
\end{remark}

The cohomological dimension of an $\infty$-topos $\calX$ is closely related to the homotopy dimension of $\calX$. If $\calX$ has homotopy dimension $\leq n$, then
$$ \HH^m (\calX; A)  = \pi_0 \bHom_{\calX}(1_{\calX}, K(A,m)) = \ast $$
for $m > n$ by Lemma \ref{pie}, so that $\calX$ is also of cohomological dimension
$\leq n$. We will establish a partial converse to this result. 

\begin{definition}\index{gen}{gerbe}\index{gen}{$n$-gerbe}
Let $\calX$ be an $\infty$-topos. An {\it $n$-gerbe} on $\calX$ is an object
$X \in \calX$ which is $n$-connective and $n$-truncated.
\end{definition}

Let $\calX$ be an $\infty$-topos containing an $n$-gerbe $X$, and let $f: \calX_{/X} \rightarrow \calX$ denote the associated geometric morphism. If $X$ is equipped with a base point $p: 1_{\calX} \rightarrow X$, then $X$
is canonically determined (as a pointed object) by $p^{\ast} \pi_n X$, by Proposition \ref{EM}.
We now wish to consider the case in which $X$ is {\em not} pointed. If $n \geq 2$, 
then $\pi_n X$ can be regarded as an abelian group object in the topos
$\Disc(\calX_{/X})$. 
Proposition \ref{nicelemma} implies that $\pi_n X \simeq f^{\ast} A$, where $A$ is a sheaf
of abelian groups on $\calX$, which is determined up to canonical isomorphism.
(In concrete terms, this boils down the observation that the $1$-connectivity of $X$ allows us
to extract higher homotopy groups without specifying a basepoint on $X$. )
In this situation, we will say that $X$ is {\it banded by $A$}.\index{gen}{gerbe!banded}

\begin{remark}
For $n < 2$, the situation is more complicated. We refer the reader to \cite{giraud} for a discussion.
\end{remark}

Our next goal is to show that the cohomology groups of an $\infty$-topos $\calX$ can be interpreted as classifying equivalence classes of $n$-gerbes over $\calX$. Before we can prove this, we need to establish some terminology.

\begin{notation}\index{not}{BandcalX@$\Band(\calX)$}
Let $\calX$ be an $\infty$-topos. We define a category $\Band(\calX)$ as follows:
\begin{itemize}
\item[$(1)$] The objects of $\Band(\calX)$ are pairs $(U,A)$, where $U$ is an object of $\calX$
and $A$ is an abelian group object of the homotopy category $\Disc( \calX_{/U} )$. 
\item[$(2)$] Morphisms from $(U,A)$ to $(U',A')$ are given by pairs
$(\eta, f)$, where $\eta \in \pi_0 \bHom_{\calX}(U,U')$ and $f: A \rightarrow A'$ is a map
which induces an isomorphism $A \simeq \eta^{\ast} A'$ of abelian group objects. Composition of morphisms is defined in the obvious way.
\end{itemize}

For $n \geq 2$, let $\Gerb_n(\calX)$ denote the subcategory of $\Fun(\Delta^1,\calX)$ spanned by those objects $f: X \rightarrow S$ which are $n$-gerbes in $\calX_{/S}$ and those morphisms which correspond to pullback diagrams
$$ \xymatrix{ X' \ar[r] \ar[d]^-{f} & X \ar[d]^-{f} \\
S' \ar[r] & S. }$$\index{not}{GerbncalX@$\Gerb_n(\calX)$}

\begin{remark}\label{sumh}
Since the class of morphisms $f: X \rightarrow S$ which belong to $\calX^{\Delta^1}$ is stable under pullback, we can apply Corollary \ref{tweezegork} (which asserts that
$p: \Fun(\Delta^1,\calX) \rightarrow \Fun( \{1\}, \calX)$ is a Cartesian fibration), Lemma \ref{charpull} (which characterizes the $p$-Cartesian morphisms of $\Fun(\Delta^1,\calX)$), and Corollary \ref{relativeKan} to deduce that the projection $\Gerb_{n}(\calX) \rightarrow \calX$ is a right fibration.
\end{remark}

If $f: X \rightarrow U$ belongs to $\Gerb_n(\calX)$, then 
there exists an abelian group object $A$ of $\Disc( \calX_{/U} )$ such that $X$ is banded by $A$. The construction
$$ (f: X \rightarrow U) \mapsto (U,A)$$
determines a functor
$$ \chi: \Gerb_n(\calX) \rightarrow \Nerve(\Band(\calX)).$$

Let $A$ be an abelian group object of $\Disc( \calX)$. We let $\Band^{A}(\calX)$\index{not}{BandAcalX@$\Band^{A}(\calX)$}
be the category whose objects are triples $(X, A_{X}, \phi)$, where $X \in \h{\calX}$,
$A_{X}$ is an abelian group object of $\Disc( \calX_{/X})$, and $\phi$ is a map
$A_{X} \rightarrow A$ which induces an isomorphism $A_{X} \simeq A \times X$ of abelian
group objects of $\Disc( \calX_{/X})$. We have forgetful functors
$$ \Band^{A}(\calX) \stackrel{\phi}{\rightarrow} \Band(\calX) \rightarrow \h{\calX},$$
both of which are Grothendieck fibrations and whose composition is an equivalence of categories.
We define $\Gerb^{A}_n(\calX)$ by the following pullback diagram:
$$ \xymatrix{ \Gerb_n^{A}(\calX) \ar[r] \ar[d] & \Gerb_n(\calX) \ar[d]^-{\chi} \\
\Nerve(\Band^{A}(\calX)) \ar[r] & \Nerve(\Band(\calX)). }$$
Note that since $\phi$ is a Grothendieck fibration, $\Nerve \phi$ is 
a Cartesian fibration (Remark \ref{gcart}), so that the diagram above is homotopy Cartesian
( Proposition \ref{basechangefunky} ). We will refer to $\Gerb_n^{A}(\calX)$ as the {\it sheaf of gerbes over $\calX$ banded by $A$}.\index{not}{GerbnAcalX@$\Gerb_n^{A}(\calX)$}
\end{notation}

More informally: an object of $\Gerb_n^{A}(\calX)$ is an $n$-gerbe $f: X \rightarrow U$
in $\calX_{/U}$ {\em together with} an isomorphism $\phi_{X}: \pi_{n} X \simeq X \times A$
of abelian group objects of $\Disc( \calX_{/X})$. Morphisms in $\Gerb_n^{A}$
are given by pullback squares
$$ \xymatrix{ X' \ar[d] \ar[r]^-{f} & X \ar[d] \\
U' \ar[r] & U }$$
such that the associated diagram of abelian group objects of $\Disc( \calX_{/X'})$
$$ \xymatrix{ & f^{\ast}(\pi_n X) \ar[dr]^-{ f^{\ast} \phi_{X} } & \\
\pi_n X' \ar[ur]^-{\pi_n f} \ar[rr]^-{\phi_{X'}} & & A \times X' }$$
is commutative.

\begin{lemma}\label{stareye}
Let $\calX$ be an $\infty$-topos, $n \geq 1$, and $A$ an abelian group object in the topos
$\Disc( \calX)$. Let $X$ be an $n$-gerbe in $\calX$ equipped with a fixed isomorphism
$\phi: \pi_n X \simeq X \times A$ of abelian group objects of $\Disc( \calX_{/X})$,
and let $u: 1_{\calX} \rightarrow K(A,n)$ be an Eilenberg-MacLane object of $\calX$ classified by $A$. Let $\bHom_{\calX}^{\phi}( K(A,n), X)$ be the summand of $\bHom_{\calX}( K(A,n), X)$
corresponding to those maps $f: K(A,n) \rightarrow X$ for which the composition
$$ A \times K(A,n) \simeq \pi_n K(A,n) \rightarrow f^{\ast}(\pi_n X) \stackrel{f^{\ast} \phi}{\rightarrow}
A \times K(A,n)$$
is the identity $($ in the category of abelian group objects of $\calX_{/K(A,n)}$ $)$. Then
composition with $u$ induces a homotopy equivalence
$$ \theta^{\phi}: \bHom_{\calX}^{\phi}( K(A,n), X) \rightarrow \bHom_{\calX}( 1_{\calX}, X).$$
\end{lemma}

\begin{proof}
Let $\theta: \bHom_{\calX}( K(A,n), X) \rightarrow \bHom_{\calX}(1_{\calX}, X)$, and let
$f: 1_{\calX} \rightarrow X$ be any map (which we may identify with an Eilenberg-MacLane object of $\calX$. The homotopy fiber of $\theta$ over the point represented by $f$ can be identified
with $\bHom_{\calX_{1_{\calX}/}}( u, f)$. In view of the equivalence between 
$\calX_{1_{\calX}/}$ and $\calX_{\ast}$, we can identify this mapping space with
$\bHom_{\calX_{\ast}}( u, f)$. Applying Proposition \ref{EM}, we deduce that
the homotopy fiber of $\theta$ is equivalent to the (discrete) set of all endomorphisms $v: A \rightarrow A$ (in the category of group objects of $\Disc( \calX)$). We now observe that
the homotopy fiber of $\theta^{\phi}$ over $f$ is a summand of the homotopy fiber of $\theta$ over $f$, corresponding to those components for which $v = \id_{A}$. It follows that the homotopy fibers
of $\theta^{\phi}$ are contractible, so that $\theta^{\phi}$ is a homotopy equivalence as desired.
\end{proof}

\begin{lemma}\label{starsky}
Let $\calX$ be an $\infty$-topos, $n \geq 1$, and $A$ an abelian group object of 
$\Disc( \calX)$. Let $f: K(A,n) \times X \rightarrow X$ be a trivial $n$-gerbe over
$X$ banded by $A$, and $g: \widetilde{Y} \rightarrow Y$ any $n$-gerbe over $Y$ banded by $A$.
Then there is a canonical homotopy equivalence
$$ \bHom_{ \Gerb^A_{n} }(f,g) \simeq \bHom_{ \calX } (X, \widetilde{Y} ).$$
\end{lemma}

\begin{proof}
Choose a morphism $\alpha: \id_{X} \rightarrow f$ in $\calX_{/X}$ as depicted below:
$$ \xymatrix{ X \ar[d] \ar[r]^-{s} & X \times K(A,n) \ar[d]^-{f} \\
X \ar[r]^-{\id_{X} } & X }$$
which exhibits $f$ as an Eilenberg-MacLane object of $\calX_{/X}$. We observe that
evaluation at $\{0\} \subseteq \Delta^1$ induces a trivial fibration
$$ \Hom^{\lft}_{\calX^{\Delta^1}}( \id_{X}, g) \rightarrow \Hom^{\lft}_{\calX}(X, \widetilde{Y}).$$
Consequently, we may identify $\bHom_{\calX}(X, \widetilde{Y})$ with the Kan complex
$$ Z = \Fun(\Delta^1,\calX)_{\id_{X}/} \times_{ \Fun(\Delta^1,\calX) } \{ g \}. $$
Similarly, the trivial fibration $\Fun(\Delta^1,\calX)_{\alpha/} \rightarrow \Fun(\Delta^1,\calX)_{f/}$
allows us to identify $\bHom_{\Gerb_n}(f,g)$ with the Kan complex
$$Z' = \Fun(\Delta^1,\calX)_{\alpha/} \times_{\Fun(\Delta^1,\calX) } \{ g\},$$
and $\bHom_{\Gerb_n}(f,g)$ with the summand $Z''$ of $Z'$ corresponding to those maps
which induce the identity isomorphism of $A \times (K(A,n) \times X)$ (in the category of group
objects of $\Disc( \calX_{/ K(A,n) \times X})$). We now observe that evaluation at
$\{1\} \subseteq \Delta^1$ gives a commutative diagram
$$ \xymatrix{ Z'' \ar[r] \ar[dr]^-{\psi''} & Z' \ar[d]^-{\psi'} \ar[r] & Z \ar[d]^-{\psi} \\
& \calX_{\id_{X}/} \times_{\calX} \{Y \} \ar[r] & \calX_{X/} \times_{\calX} \{Y\} }.$$
where the vertical maps are Kan fibrations. If we fix a pullback square
$$ \xymatrix{ \widetilde{X} \ar[r] \ar[d]^-{g'} & \widetilde{Y} \ar[d] \\
X \ar[r]^-{h} & Y,} $$
then we can identify $\psi^{-1} \{h\}$ with $\bHom_{\calX^{/X}}( \id_{X}, g')$,
${\psi'}^{-1} \{ s^0 h \}$ with $\bHom_{\calX^{/X}}( X \times K(A,n), g')$, 
${\psi'})^{-1} \{ s^0 h\}$ with the summand of $\bHom_{\calX^{/X}}( X \times K(A,n), g')$
corresponding to those maps which induce the identity on $A \times (K(A,n) \times X)$ (in the category of group objects of $\Disc( \calX_{/ K(A,n) \times X})$), and $\theta$ with
the map given by composition with $s$. Invoking Lemma \ref{stareye} in the $\infty$-topos
$\calX^{/X}$, we deduce that the map $\theta$ in the diagram
$$ \xymatrix{
Z'' \ar[r]^-{\theta} \ar[d]^-{\psi''} & Z \ar[d]^-{\psi} \\
 \calX_{\id_{X}/} \times_{\calX} \{Y \} \ar[r] & \calX_{X/} \times_{\calX} \{Y\} }$$
induces homotopy equivalences from the fibers of $\psi''$ to the fibers of $\psi$. Since the 
lower horizontal map is a trivial fibration of simplicial sets, we conclude that $\theta$ is itself a homotopy equivalence, as desired.
\end{proof}

\begin{theorem}\label{starthm}
Let $\calX$ be an $\infty$-topos, $n \geq 1$, and $A$ an abelian group object of
$\Disc( \calX)$. Then:

\begin{itemize}
\item[$(1)$] The composite map
$$ \theta: \Gerb_n^{A}(\calX) \rightarrow \Gerb_n(\calX) \subseteq \Fun(\Delta^1,\calX) \rightarrow
\Fun(\{1\}, \calX) \simeq \calX $$ is a right fibration.

\item[$(2)$] The right fibration $\theta$ is representable by an Eilenberg-MacLane object
$K(A,n+1)$. 
\end{itemize}
\end{theorem}

\begin{proof}
For each object $X \in \calX$, we let $A_{X}$ denote the projection $A \times X \rightarrow X$, viewed as an abelian group object of $\Disc( \calX_{/X})$.  
The functor $\phi: \Band^{A}(\calX) \rightarrow \Band(\calX)$
is a fibration in groupoids, so that $\Nerve \phi$ is a right fibration (Proposition \ref{stinkyer}). 
The functor $\theta$ admits a factorization
$$ \Gerb_n^{A}(\calX) \stackrel{\theta'}{\rightarrow} \Gerb_n(\calX) \stackrel{\theta''}{\rightarrow} \calX$$
where $\theta''$ is a right fibration (Remark \ref{sumh}) and $\theta'$ is a pullback of
$\Nerve \phi$, and therefore also a right fibration. It follows that $\theta$, being a composition of right fibrations, is a right fibration; this proves $(1)$.

To prove $(2)$, we consider an Eilenberg-MacLane object $u: 1_{\calX} \rightarrow K(A,n+1)$.
Since $K(A,n+1)$ is $(n+1)$-truncated and $1_{\calX}$ is $n$-truncated (in fact, $(-2)$-truncated), 
Lemma \ref{trunccomp} implies that $u$ is $n$-truncated. 
The long exact sequence
$$\ldots \rightarrow u^{\ast} \pi_{i+1} K(A,n+1) \rightarrow \pi_i u \rightarrow \pi_i(1_{\calX}) \rightarrow i^{\ast}
\pi_i( K(A,n+1) ) \rightarrow \pi_{i-1}(u) \rightarrow \ldots$$
of Remark \ref{sequence} shows that $u$ is $n$-connective, and provides an
isomorphism $\phi: A \simeq \pi_n(u)$ in the category of group objects of $\Disc( \calX)$,
so that we may view the pair $(u,\phi)$ as an object of $\Gerb_n^{A}(\calX)$. Since
$1_{\calX}$ is a final object of $\calX$, Lemma \ref{starsky} implies that $(u,\phi)$ 
is a final object of $\Gerb_n^{A}(\calX)$, so that the right fibration $\theta$ is representable
by $\theta(u,\phi) = K(A,n+1)$.
\end{proof}

\begin{corollary}\label{coclassify}
Let $\calX$ be an $\infty$-topos, $n \geq 2$, and $A$ an abelian group object
of $\Disc( \calX)$. There is a canonical bijection of
$\HH^{n+1}(\calX; A)$ with the set of equivalence classes of $n$-gerbes on $\calX$ banded by $A$.
\end{corollary}

\begin{remark}
Under the correspondence of Proposition \ref{coclassify}, an $n$-gerbe $X$ on $\calX$
admits a global section $1_{\calX} \rightarrow X$ if and only if the associated cohomology
class in $\HH^{n+1}(\calX;A)$ vanishes.
\end{remark}

\begin{theorem}\label{cohdim}
Let $\calX$ be an $\infty$-topos and $n \geq 2$. Then $\calX$ has
cohomological dimension $\leq n$ if and only if it satisfies the
following condition: any $n$-connective, truncated object of
$\calX$ admits a global section.
\end{theorem}

\begin{proof}
Suppose that $\calX$ has the property that every
$n$-connective, truncated object $X \in \calX$ admits a global
section. As in the proof of Lemma \ref{pie}, we deduce that for any
truncated, $(n+1)$-connective object $X \in \calX$, the space of global sections 
$\bHom_{\calX}(1,X)$ is connected.
Let $k > n$, and let $G$ be a sheaf of abelian groups on $\calX$. Then $K(G,k)$ is
$(n+1)$-connective, so that $\HH^k(\calX,G) = \ast$. Thus $\calX$ has
cohomological dimension $\leq n$.

For the converse, let us assume that $\calX$ has cohomological
dimension $\leq n$ and let $X$ denote an $n$-connective,
$k$-truncated object of $\calX$. We will show that $X$ admits a
global section by descending induction on $k$. If $k \leq n-1$, then $X$ is a final
object of $\calX$, so there is nothing to prove. In the general case, choose a truncation
$X \rightarrow \tau_{\leq k-1} X$; we may assume
by the inductive hypothesis that $\tau_{\leq k-1} X$ has a global section $s: 1
\rightarrow \tau_{\leq k-1} X$.
Form a pullback square
$$ \xymatrix{ X' \ar[r] \ar[d] & X \ar[d] \\
1 \ar[r]^-{s} & \tau_{\leq k-1} X. }$$ 
It now suffices to prove that $X'$ has a global section. We note that $X'$ is $k$-connective, where $k \geq n \geq 2$. It follows that $X'$ is a $k$-gerbe on $\calX$; suppose it is banded by
an abelian group object $A \in \Disc(\calX)$. 
According to Corollary \ref{coclassify}, $X'$ is classified up to equivalence by an element in $\HH^{k+1}(\calX, A)$, which vanishes in virtue of the fact that
$k+1 > n$ and the cohomological dimension of $\calX$ is $\leq n$. Consequently, $X'$ is equivalent to $K(A,k)$ and therefore admits a global section.
\end{proof}

\begin{corollary}\label{confusion}
Let $\calX$ be an $\infty$-topos. If $\calX$ has homotopy
dimension $\leq n$, then $\calX$ has cohomological dimension $\leq
n$. The converse holds provided that $\calX$ has finite homotopy
dimension and $n \geq 2$.
\end{corollary}

\begin{proof}
Only the last claim requires proof. Suppose that $\calX$ has
cohomological dimension $\leq n$ and homotopy dimension $\leq k$.
We must show that every $n$-connective object $X$ of $\calX$
has a global section. Choose a truncation $X \rightarrow \tau_{\leq k-1} X$.
Then $\tau_{\leq k-1} X$ is truncated and
$n$-connective, so it admits a global section by Theorem
\ref{cohdim}. Form a pullback square
$$ \xymatrix{ X' \ar[r] \ar[d] & X \ar[d] \\
1 \ar[r] & \tau_{\leq k-1} X.}$$
It now suffices to prove that $X'$ has a global section. But $X'$ is
$k$-connective, and therefore has a global section in virtue of the assumption that
$\calX$ has homotopy dimension $\leq k$.
\end{proof}

\begin{warning}\label{injur}[Weiland]
The converse to Corollary \ref{confusion} is false if we do not assume that
$\calX$ has finite homotopy dimension. To see this, we discuss the following example, which
we learned from Ben Wieland. Let $G$ denote the group $\Z_{p}$ of $p$-adic integers
(viewed as a profinite group). Let $\calC$ denote the category whose objects are
the finite quotients $\{ \Z_{p} / p^{n} \Z_{p} \}_{ n \geq 0}$, and whose morphisms
are given by $G$-equivariant maps. We regard $\calC$ as endowed with a Grothendieck topology in which every nonempty sieve is a covering. The $\infty$-topos $\Shv(\Nerve \calC)$
is $1$-localic, and the underlying ordinary topos $\h{ \tau_{\leq 0} \Shv(\Nerve \calC) }$ can be identified
with the category $BG$ of continuous $G$-sets (that is, sets $C$ equipped with an action of $G$
such that the stabilizer of each element $x \in C$ is an open subgroup of $G$). Since
the profinite group $G$ has cohomology dimension $2$ (see \cite{serre}), we deduce that $\calX$ is
of cohomological dimension $2$. However, we will show that $\calX$ is not hypercomplete, and therefore cannot be of finite homotopy dimension.

Let $K$ be a finite CW complex whose homotopy groups consist entirely of $p$-torsion
(for example, we could take $K$ to be a Moore space $M( \Z/p\Z)$), and let
$X = \Sing K \in \SSet$. Let $F: \Nerve(\calC)^{op} \rightarrow \SSet$ denote the constant functor taking the value $X$. We claim that $F$ belongs to
$\Shv(\calC)$. Unwinding the definitions, we must show that for each $m \leq n$, 
the diagram $F$ exhibits $F( \Z_{p} / p^{m} \Z_p)$ as equivalent to the homotopy
invariants for the trivial action of $p^{m} \Z_{p}/ p^{n} \Z_p$ on $F( \Z_{p} / p^{n} \Z_p)$. In
other words, we must show that the diagonal embedding
$$ \alpha: X \rightarrow \Fun( BH, X)$$
is a homotopy equivalence, where $H$ denotes the quotient group
$p^{m} \Z_{p}/ p^{n} \Z_p$. Since both sides are $p$-adically complete, it will suffice
to show that $\alpha$ is a $p$-adic homotopy equivalence, which follows from a suitable version of the Sullivan conjecture (see, for example, \cite{schwartz}). 

We define another functor $F': \Nerve(\calC)^{op} \rightarrow \SSet$, which is obtained as the simplicial nerve of the functor described by the formula
$$ \Z_{p}/ p^{n} \Z_{p} \mapsto \Sing (K^{\R / p^{n} \Z}).$$
For $m \leq n$, the loop space $K^{ \R/ p^{m} \Z}$ can be identified with the homotopy
fixed points of the (nontrivial) action of $H = p^{m} \Z_{p} / p^{n} \Z_{p} \simeq p^{m} \Z / p^{n} \Z$
on the loop space $K^{ \R/ p^{n} \Z}$: this follows from the observation that $H$ acts freely on
$\R / p^{n} \Z$, with quotient $\R / p^{m} \Z$. Consequently, $F'$ belongs to $\Shv( \Nerve(\calC))$.

The inclusion of $K$ into each loop space $K^{ \R / p^{n} \Z}$ induces a morphism
$\alpha: F \rightarrow F'$ in the $\infty$-topos $\Shv( \Nerve(\calC) )$. Using the fact
that the homotopy groups of $K$ are $p$-torsion, we deduce that the morphism
$\alpha$ is $\infty$-connective (this follows from the observation that the map
$$X \simeq \varinjlim F( \Z_p / p^{n} \Z_p) \rightarrow
\varinjlim F'( \Z / p^{n} \Z_P)$$ is a homotopy equivalence). However,
the morphism $\alpha$ is not an equivalence in $\Shv( \Nerve(\calC) )$ unless $K$ is essentially discrete. Consequently, $\Shv( \Nerve(\calC) )$ is not hypercomplete, and therefore cannot be of finite homotopy dimension.
\end{warning}

In spite of Warning \ref{injur}, many situations which guarantee that a topological space (or topos)
$X$ is of bounded cohomological dimension also guarantee that the associated $\infty$-topos is
of bounded homotopy dimension. We will see some examples in the next two sections.

\subsection{Covering Dimension}\label{covdim}

In this section, we will review the classical theory of covering
dimension for paracompact spaces, and then show that the covering
dimension of a paracompact space $X$ coincides with its homotopy
dimension.

\begin{definition}\label{paradim}\index{gen}{dimension!covering}\index{gen}{covering dimension}
A paracompact topological space $X$ has {\it covering
dimension $\leq n$} if the following condition is satisfied: for
any open covering $\{ U_{\alpha} \}$ of $X$, there exists an open
refinement $\{ V_{\alpha} \}$ of $X$ such that each intersection
$V_{\alpha_0} \cap \ldots \cap V_{\alpha_{n+1}} = \emptyset$
provided the $\alpha_i$ are pairwise distinct.
\end{definition}

\begin{remark}
When $X$ is paracompact, the condition of Definition \ref{paradim}
is equivalent to the (a priori weaker) requirement that such a refinement exist whenever
$\{U_{\alpha} \}$ is a finite covering of $X$. This weaker
condition gives a good notion whenever $X$ is a normal topological
space. Moreover, if $X$ is normal, then the covering dimension of
$X$ (by this second definition) coincides with the covering dimension of
the Stone-\Cech compactification of $X$. Thus, the dimension
theory of normal spaces is controlled by the dimension theory of
compact Hausdorff spaces.
\end{remark}

\begin{remark}
Suppose that $X$ is a compact Hausdorff space, which is written as
a filtered inverse limit of compact Hausdorff spaces $\{
X_{\alpha} \}$, each of which has dimension $\leq n$. Then $X$ has
dimension $\leq n$. Conversely, any compact Hausdorff space of
dimension $\leq n$ can be written as a filtered inverse limit of
finite simplicial complexes having dimension $\leq n$. Thus, the
dimension theory of compact Hausdorff spaces is controlled by the
(completely straightforward) dimension theory of finite simplicial
complexes.
\end{remark}

\begin{remark}
There are other approaches to classical dimension theory. For
example, a topological space $X$ is said to have {\it small
$($ large $)$ inductive dimension $\leq n$} if every point of $X$ (every
closed subset of $X$) has arbitrarily small open neighborhoods $U$
such that $\bd U$ has small inductive dimension $\leq n-1$. These
notions are well-behaved for separable metric spaces, where they
coincide with the covering dimension (and with each other). In
general, the covering dimension has better formal properties.
\end{remark}

Our goal in this section is to prove that the covering dimension of a paracompact topological space $X$ coincides with the homotopy dimension of $\Shv(X)$.
First, we need a technical
lemma.

\begin{lemma}\label{core}
Let $X$ be a paracompact space, $k \geq 0$, $\{U_{\alpha}
\}_{\alpha \in A}$ be a covering of $X$. Suppose that for every
$A_0 \subseteq A$ of size $k+1$, we are given a covering
$\{V_{\beta} \}_{\beta \in B(A_0)}$ of the intersection $U_{A_0} =
\bigcap_{\alpha \in A_0} U_{\alpha}$. Then there exists a covering
$\{ W_{\alpha} \}_{\alpha \in \widetilde{A}}$ of $X$ and a map
$\pi: \widetilde{A} \rightarrow A$ with the following properties:
\begin{itemize}
\item[$(1)$] For $\widetilde{\alpha} \in \widetilde{A}$ with $\pi(\widetilde{\alpha}) = \alpha$, we have $W_{\widetilde{\alpha}} \subseteq U_{\alpha}$. 

\item[$(2)$] Suppose that $\widetilde{\alpha}_0, \ldots, \widetilde{\alpha}_k$ is a collection
of elements of $\widetilde{A}$, with $\pi( \widetilde{\alpha}_i ) = \alpha_i$. Suppose
further that $A_0 = \{ \alpha_0, \ldots, \alpha_k \}$ has cardinality $(k+1)$ (in other words, the $\alpha_i$ are all disjoint from one another). 
Then there exists
$\beta \in B(A_0)$ such that $W_{\widetilde{\alpha}_0} \cap \ldots \cap
W_{\widetilde{\alpha}_k} \subseteq V_{\beta}$.
\end{itemize}
\end{lemma}

\begin{proof}
Since $X$ is paracompact, we may find a locally finite covering
$\{ U'_{\alpha} \}_{\alpha \in A}$ of $X$, such that the
each closure $\overline{ U'_{\alpha} }$ is contained in
$U_{\alpha}$. Let $S$ denote the set of all subsets $A_0 \subseteq
A$ having size $k+1$. For $A_0 \in S$, let $K(A_0) = \bigcap_{\alpha
\in A_0} \overline{U_{\alpha}}$. Now set 
$$\widetilde{A} = \{ (\alpha, A_0, \beta): \alpha \in A_0 \in S, \beta \in
B(A_0) \} \cup A.$$ 
For $\widetilde{\alpha} = (\alpha, A_0, \beta) \in \widetilde{A}$, we
set $\pi(\widetilde{\alpha}) = \alpha$ and 
$$W_{\widetilde{\alpha}} =
(U'_{\alpha} - \bigcup_{\alpha \in A'_0 \in S} K(A'_0) ) \cup (V_{\beta}
\cap U'_{\alpha}).$$
If $\alpha \in A \subseteq \widetilde{A}$, we let $\pi(\alpha) =
\alpha$ and $W_{\alpha} = U'_{\alpha} - \bigcup_{\alpha \in A_0 \in S}
K(A_0)$. The local finiteness of the cover $\{ U'_{\alpha} \}$ ensures that each
$W_{\widetilde{\alpha}}$ is an open set. It is now easy to check that the covering $\{ W_{ \widetilde{\alpha} } \}_{\widetilde{\alpha} \in \widetilde{A}}$ has the desired properties.
\end{proof}

\begin{theorem}\label{paradimension}
Let $X$ be a paracompact topological space of covering dimension
$\leq n$. Then the $\infty$-topos $\Shv(X)$ of sheaves on $X$ has
homotopy dimension $\leq n$.
\end{theorem}

\begin{proof}
We make use of the results and notations of \S \ref{paracompactness}.
Let $\calB$ denote the collection of all open $F_{\sigma}$ subsets of $X$, and fix a linear ordering on $\calB$. We may identify $\Shv(X)$ with the simplicial nerve of the category
of all functors $F: \calB^{op} \rightarrow \Kan$ which have the property that for any
$\calU \subseteq \calB$ with $U = \bigcup_{V \in \calU} V$, the natural map
$F(U) \rightarrow F(\calU)$ is a homotopy equivalence.

Suppose that $F: \calB^{op} \rightarrow \sSet$ represents an $n$-connective
sheaf; we wish to show that the simplicial set $F(X)$ is nonempty. It suffices to prove that
$F(\calU)$ is nonempty, for some covering $\calU$ of $X$; in other words, it suffices to produce a map $N_{\calU} \rightarrow F$. The idea is that since $X$ has finite covering dimension, we can choose arbitrarily fine covers $\calU$ such that $N_{\calU}$ is {\it $n$-dimensional}; that is, equal to its $n$-skeleton. 

For every simplicial set $K$, let $K^{(i)}$ denote the {\it $i$-skeleton} of $K$ (the union of all
nondegenerate simplices of $K$ of dimension $\leq i$). If $G: \calB^{op} \rightarrow \sSet$ is a simplicial presheaf, we let $G^{(i)}$ denote the simplicial presheaf given by the formula
$$ G^{(i)}(U) = (G(U))^{(i)}.$$

We will prove the following statement by induction on $i$, $-1
\leq i \leq n$:

\begin{itemize}
\item There exists an open cover $\calU_{i} \subseteq \calB$ of $X$
and a map $\eta_i: N_{\calU_i}^{(i)} \rightarrow
F$.
\end{itemize}

Assume that this statement holds for $i = n$. Passing to a
refinement, we may assume that the cover $\calU_{n}$ has the
property that no more than $n+1$ of its members intersect (this is
the step where we shall use the assumption on the covering
dimension of $X$). It follows that $N_{\calU_{n}}^{(n)} = N_{\calU_{n}}$, and the proof
will be complete.

To begin the induction in the case $i = -1$, we let $\calU_{-1} = \{ X \}$; the $(-1)$-skeleton of 
$N_{\calU_{-1}}$ is empty, so that $\eta_{-1}$ exists (and is unique).

Now suppose that $\calU_{i} = \{ U_{\alpha} \}_{ \alpha \in A} $ and $\eta_i$ have been constructed, $i < n$. Let $A_0 \subseteq A$ have cardinality $(i+2)$, and set $U(A_0) = \bigcap_{\alpha \in A_0} U_{\alpha}$; then $A_0$ determines an $n$-simplex of $N_{\calU_{i}}(U(A_0))$, so that
$\eta_i$ restricts to give a map
$$\eta_{i, A_0}: \bd \Delta^{i+1} \rightarrow F(U(A_0)).$$
By assumption, $F$ is $n$-connective; it follows that there is an open covering
$$ \{ V_{\beta} \}_{ \beta \in B(A_0) }$$
of $U(A_0)$, such for each $V_{\beta}$ there is a commutative diagram
$$ \xymatrix{ \bd \Delta^{i+1} \ar@{^{(}->}[d] \ar[r] & F(U(A_0)) \ar[d] \\
\Delta^{i+1} \ar[r] & F(V_{\beta} ). } $$

We apply Lemma \ref{core} to this data, to
obtain an new open cover $\calU_{i+1} = \{ W_{\widetilde{\alpha}} \}_{\widetilde{\alpha} \in
\widetilde{A} }$ which refines $\{ U_{\alpha} \}_{\alpha \in A}$.
Refining the cover further if necessary, we may assume that each
of its members belongs to $\calB$. By functoriality, we obtain a
map
$$ N^{(i)}_{\calU_{i+1} } \rightarrow F.$$
To complete the proof,
it will suffice to extend $f$ to the $(i+1)$-skeleton of the nerve
of $\{ W_{\alpha} \}_{\alpha \in \widetilde{A}}$. Let $\widetilde{A_0} \subseteq \widetilde{A}$
have cardinality $i+2$, and let $W( \widetilde{A_0}) = \bigcap_{\widetilde{\alpha} \in \widetilde{A_0}} W_{\widetilde{\alpha}}$; then we must solve a lifting problem
$$ \xymatrix{ \bd \Delta^{i+1} \ar@{^{(}->}[d] \ar[r] & F(W) \\
\Delta^{i+1}. \ar@{-->}[ur] & } $$
Let $\pi: \widetilde{A} \rightarrow A$ denote the map of Lemma \ref{core}.
If $A_0 = \pi( \widetilde{A}_0 )$ has cardinality smaller than $i+2$, then there is a canonical extension, given by applying $\pi$ and using $\eta_{i}$. Otherwise, Lemma \ref{core} guarantees that $W( \widetilde{A_0} ) \subseteq V_{\beta}$ for some $\beta \in B(A_0)$, so that the
desired extension exists by construction.
\end{proof}

\begin{corollary}\label{corub}
Let $X$ be a paracompact topological space. The following conditions are equivalent:
\begin{itemize}
\item[$(1)$] The covering dimension of $X$ is $\leq n$.
\item[$(2)$] The homotopy dimension of $\Shv(X)$ is $\leq n$.
\item[$(3)$] For every closed subset $A \subseteq X$, every $m \geq n$, and every
continuous map $f_0: A \rightarrow S^m$, there exists $f: X \rightarrow S^m$ extending
$f_0$.
\end{itemize}
\end{corollary}

\begin{proof}
The implication $(1) \Rightarrow (2)$ is Theorem \ref{paradimension}.
The equivalence $(1) \Leftrightarrow (3)$ follows from classical dimension theory (see, for example,  \cite{dimtheory}).
We will complete the proof by showing that $(2) \Rightarrow (3)$.
Let $A$ be a closed subset of $X$, $m \geq n$, and $f_0: A \rightarrow S^m$ a continuous map.
Let $\calB$ be the collection of all open $F_{\sigma}$ subsets of $X$.
We define a simplicial presheaf $F: \calB \rightarrow \Kan$, so that an $n$-simplex of $F(U)$
is a map $f$ rendering the diagram
$$ \xymatrix{ (U \cap A) \times | \Delta^n | \ar[r] \ar[d] & A \ar[d]^-{f_0} \\
U \times | \Delta^n | \ar[r]^-{f} \ar[r] & S^m}$$ commutative. 
To prove $(3)$, it will suffice to show that $F(X)$ is nonempty. In virtue of the assumption
that $\Shv(X)$ has homotopy dimension $\leq n$, it will suffice to show that $\calF$ is an $n$-connective sheaf on $X$.

We first show that $F$ is a sheaf. Choose a linear ordering on $\calB$. We must show that
for every open covering $\calU$ of $U \in \calB$, the natural map
$\calF(U) \rightarrow \calF(\calU)$ is a homotopy equivalence. The proof is similar to that of Proposition \ref{aese}. Let $\pi: |N_{\calU}|_{X} \rightarrow U$ be the projection; then
we may identify $F(\calU)$ with the simplicial set parametrizing
continuous maps $|N_{\calU}|_{X} \rightarrow S^m$, whose restriction to $\pi^{-1}(A)$ is
given by $f_0$. The desired equivalence now follows from the fact that $|N_{\calU}|_{X}$ is fiberwise homotopy equivalent to $U$ (Lemma \ref{partit}).

Now we claim that $\calF$ is $n$-connective as an object of $\Shv(X)$. In other words, we must show that for any $U \in \calB$, any $k \leq n$, and any map $g: \bd \Delta^k \rightarrow F(U)$, there is an open covering $\{ U_{\alpha} \}$ of $U$ and a family of commutative diagrams
$$ \xymatrix{ \bd \Delta^k \ar@{^{(}->}[d] \ar[r]^-{g} & F(U) \ar[d] \\
\Delta^k \ar[r]^-{g_{\alpha}} & F(U_{\alpha} ).} $$

We may identify $g$ with a continuous map
$$g: S^{k-1} \times U \rightarrow S^{m}$$
such that $g( z, a) = f_0(a)$ for $a \in A$. Choose a point $x \in U$.
Consider the map $g| S^{k-1} \times \{x\}$. Since $k-1 < n \leq m$, this map is nullhomotopic;
therefore it admits an extension $g'_{x}: D^k \times \{x\} \rightarrow S^m$. Moreover, if
$x \in A$, then we may choose $g'_{x}$ to be the constant map with value $f_0(x)$.
Amalgamating $g$, $g'_{x}$, and $f_0$, we obtain a continuous map
$$ g'_{0}: (S^{k-1} \times U) \cup (D^k \times ( A \cup \{x\} ) ) \rightarrow S^m.$$
Since $(S^{k-1} \times U) \cup (D^k \times (A \cup \{x\}))$ is a closed subset of
the paracompact space $U \times D^k$, and the sphere $S^m$ is an absolute
neighborhood retract, the map $g'_0$ extends continuously to a map $g'': W \rightarrow S^m$,
where $W$ is an open neighborhood of $(S^{k-1} \times U) \cup (D^k \times (A \cup \{x\}))$
in $U \times D^k$. The compactness of $D^k$ implies that $W$ contains
$D^k \times U_x$, where $U_x \subseteq U$ is an open neighborhood of $x$.
Shrinking $U_x$ if necessary, we may suppose that $U_x$ belongs to $\calB$;
these open sets $U_x$ form an open cover of $U$, with the required
extension $\Delta^k \rightarrow F(U_x)$ supplied by the map $g'' | D^k \times U_x$.
\end{proof}

\subsection{Heyting Dimension}\label{heyt}

For the purposes of studying paracompact topological spaces, Definition \ref{paradim}
gives a perfectly adequate theory of dimension. However, there are other situations in which
Definition \ref{paradim} is not really appropriate. For example, in algebraic geometry one often considers the Zariski topology on an algebraic variety $X$. This topology is generally not Hausdorff, and is typically of infinite covering dimension. In this setting, there is a better dimension theory: the theory of Krull dimension. In this section, we will introduce a mild generalization of the theory of Krull dimension, which we will call the {\it Heyting dimension} of a topological space $X$. We will then study the relationship between the Heyting dimension of $X$ and the homotopy dimension of the associated $\infty$-topos $\Shv(X)$.

Recall that a topological space $X$ is said to be {\it Noetherian}\index{gen}{topological space!Noetherian} if the collection of closed subsets of $X$ satisfies the descending chain condition. A closed subset $K \subseteq X$ is said to be
{\it irreducible}\index{gen}{irreducible!closed set} if it cannot be written as a finite union of proper closed subsets of $K$ (in particular, the empty set is {\em not} irreducible, since it can be written as an empty union).
The collection of irreducible closed subsets of $X$ forms a well-founded partially ordered set, therefore it has a unique ordinal rank function $\rk$, which may be characterized as follows:

\begin{itemize}
\item If $K$ is an irreducible closed subset of $X$, then $\rk(K)$ is the smallest ordinal which is larger than $\rk(K')$, for all proper irreducible closed subsets $K' \subset K$.
\end{itemize}

We call $\rk(K)$ the {\it Krull dimension} of $K$; the {\it Krull dimension} of $X$ is the supremum of $\rk(K)$, as $K$ ranges over all irreducible closed subsets of $X$.\index{gen}{dimension!Krull}\index{gen}{Krull dimension}

We next introduce a generalization of the Krull dimension to a suitable class of non-Noetherian spaces. We shall say that a topological space $X$ is a {\it
Heyting space} if satisfies the following conditions:\index{gen}{topological space!Heyting}\index{gen}{Heyting!space}

\begin{itemize}
\item[$(1)$] The compact open subsets of $X$ form a basis for the topology of $X$.

\item[$(2)$] A finite intersection of compact open subsets of $X$ is compact (in particular, $X$ is compact).

\item[$(3)$] If $U$ and $V$ are compact open subsets of
$X$, then the interior of $U \cup (X-V)$ is compact. 

\end{itemize}

\begin{remark}\index{gen}{Heyting!algebra}
Recall that a {\it Heyting algebra} is a distributive lattice $L$
with the property that for any $x,y \in L$, there exists a maximal
element $z$ with the property that $x \wedge z \subseteq y$. It
follows immediately from our definition that the lattice of
compact open subsets of a Heyting space forms a Heyting algebra.
Conversely, given any Heyting algebra one may form its spectrum,
which is a Heyting space. This sets up a duality between the
category of {\em sober} Heyting spaces (Heyting spaces in which every irreducible
closed subset has a unique generic point) and the category of Heyting algebras. This
duality is a special case of a more general duality between coherent
topological spaces and distributive lattices. We refer the reader
to \cite{johnstone} for further details.
\end{remark}

\begin{remark}
Suppose that $X$ is a Noetherian topological space. Then $X$ is
a Heyting space, since every open subset of $X$ is compact.
\end{remark}

\begin{remark}
If $X$ is a Heyting space and $U \subseteq X$ is a compact open
subset, then $X$ and $X-U$ are also Heyting spaces. In this case,
we say that $X-U$ is a {\it cocompact} closed subset of $X$.
\end{remark}

We next define the dimension of a Heyting space. The definition is
recursive. Let $\alpha$ be an ordinal. A Heyting space $X$ has
{\it Heyting dimension $\leq \alpha$} if and only if, for any
compact open subset $U \subseteq X$, the boundary of $U$ has
Heyting dimension $< \alpha$ (we note that the boundary of $U$ is
also a Heyting space); a Heyting space has {\it Heyting dimension} $< 0$ if and
only if it is empty.\index{gen}{Heyting!dimension}\index{gen}{dimension!Heyting}

\begin{remark}
A Heyting space has dimension $\leq 0$ if and only if it is
Hausdorff. The Heyting spaces of dimension $\leq 0$ are precisely
the compact, totally disconnected Hausdorff spaces. In particular,
they are also paracompact spaces and their Heyting dimension
coincides with their covering dimension.
\end{remark}

\begin{proposition}\label{closs}

\begin{itemize}
\item[$(1)$] Let $X$ be a Heyting space of dimension $\leq \alpha$. Then
for any compact open subset $U \subseteq X$, both $U$ and $X-U$
have Heyting dimension $\leq \alpha$.

\item[$(2)$] Let $X$ be a Heyting space which is a union of finitely many
compact open subsets $U_{\alpha}$ of dimension $\leq \alpha$. Then
$X$ has dimension $\leq \alpha$.

\item[$(3)$] Let $X$ be a Heyting space which is a union of finitely many
cocompact closed subsets $K_{\alpha}$ of Heyting dimension $\leq
\alpha$. Then $X$ has Heyting dimension $\leq \alpha$.
\end{itemize}
\end{proposition}

\begin{proof}
All three assertions are proven by induction on $\alpha$. The
first two are easy, so we restrict our attention to $(3)$. Let $U$
be a compact open subset of $X$, having boundary $B$. Then $U \cap
K_{\alpha}$ is a compact open subset of $K_{\alpha}$, so that the
boundary $B_{\alpha}$ of $U \cap K_{\alpha}$ in $K_{\alpha}$ has
dimension $\leq \alpha$. We see immediately that $B_{\alpha}
\subseteq B \cap K_{\alpha}$, so that $\bigcup B_{\alpha}
\subseteq B$. Conversely, if $b \notin \bigcup B_{\alpha}$ then,
for every $\beta$ such that $b \in K_{\beta}$, there exists a
neighborhood $V_{\beta}$ containing $b$ such that $V_{\beta} \cap
K_{\beta} \cap U = \emptyset$. Let $V$ be the intersection of the
$V_{\beta}$, and let $W = V - \bigcup_{b \notin K_{\gamma}}
K_{\gamma}$. Then by construction, $b \in W$ and $W \cap U =
\emptyset$, so that $b \in B$. Consequently, $B = \bigcup
B_{\alpha}$. Each $B_{\alpha}$ is closed in $K_{\alpha}$, thus in
$X$ and also in $B$. The hypothesis implies that $B_{\alpha}$ has
dimension $< \alpha$. Thus the inductive hypothesis guarantees
that $B$ has dimension $< \alpha$, as desired.
\end{proof}

\begin{remark}\label{DVR}
It is not necessarily true that a Heyting space which is a union
of finitely many {\em locally closed} subsets of dimension $\leq
\alpha$ is also of dimension $\leq \alpha$. For example, a
topological space with $2$ points and a nondiscrete, nontrivial
topology has Heyting dimension $1$, but is a union of two locally
closed subsets of Heyting dimension $0$.
\end{remark}

\begin{proposition}\label{krullheyt}
If $X$ is a Noetherian topological space, then the Krull
dimension of $X$ coincides with the Heyting dimension of $X$.
\end{proposition}

\begin{proof}
We first prove, by induction on $\alpha$, that if the Krull
dimension of a Noetherian space $X$ is $\leq \alpha$, then the
Heyting dimension of $X$ is $\leq \alpha$. Since $X$ is Noetherian,
it is a union of finitely many closed irreducible subspaces, each of which automatically has Krull
dimension $\leq \alpha$. Using Proposition \ref{closs}, we may
reduce to the case where $X$ is irreducible. Consider any open subset $U \subseteq X$, and let $Y$ be its boundary. We must show that $Y$ has Heyting dimension $\leq
\alpha$. Using Proposition \ref{closs} again, it suffices to prove
this for each irreducible component of $Y$. Now we simply apply
the inductive hypothesis and the definition of the Krull
dimension.

For the reverse inequality, we again use induction on $\alpha$.
Assume that $X$ has Heyting dimension $\leq \alpha$. To show that
$X$ has Krull dimension $\leq \alpha$, we must show that every
irreducible closed subset of $X$ has Krull dimension $\leq
\alpha$. Without loss of generality we may assume that $X$ is
irreducible. Now, to show that $X$ has Krull dimension $\leq
\alpha$, it will suffice to show that any {\em proper} closed
subset $K \subseteq X$ has Krull dimension $< \alpha$. By the
inductive hypothesis, it will suffice to show that $K$ has Heyting
dimension $< \alpha$. By the definition of the Heyting dimension,
it will suffice to show that $K$ is the boundary of $X - K$. In
other words, we must show that $X - K$ is dense in $X$. This
follows immediately from the irreducibility of $X$.
\end{proof}

We now prepare the way for our vanishing theorem. First, we
introduce a modified notion of connectivity:

\begin{definition}\label{strongcon}\index{gen}{connective!strongly}\index{gen}{strongly $k$-connective}
Let $X$ be a Heyting space and $k$ any integer. Let $\calF \in \Shv(V)$ 
be a sheaf of spaces
on a compact open subset $V \subseteq X$. We will say that $\calF$ is {\it strongly $k$-connective}
if the following condition is satisfied: for every compact open subset $U \subseteq V$
and every map $\zeta: \bd \Delta^m \rightarrow \calF(U)$, there exists a cocompact
closed subset $K \subseteq U$ such that $\overline{K} \subseteq X$ has Heyting dimension
$< m-k$, an open cover $\{ U_{\alpha} \}$ of $U-K$, and a collection of commutative diagrams
$$ \xymatrix{ \bd \Delta^m \ar[r]^-{\zeta} \ar@{^{(}->}[d] & \calF(U) \ar[d] \\
\Delta^m \ar[r]^-{\eta_{\alpha}} \ar[r] & \calF( U_{\alpha} ). } $$
\end{definition}

\begin{remark}
There is a slight risk of confusion with the terminology of Definition \ref{strongcon}.
The condition that a sheaf $\calF$ on $V \subseteq X$ be strongly $k$-connective
depends not only on $V$ and $\calF$, but also on $X$: this is because the Heyting dimension of 
a cocompact closed subset $K \subseteq U$ can increase when we take its closure $\overline{K}$ in $X$.
\end{remark}

\begin{remark}
Strong $k$-connectivity is an unstable analogue of the
connectivity conditions on complexes of sheaves, associated to the dual of the standard perversity. For a discussion of perverse sheaves in the
abelian context we refer the reader to \cite{deligne}.
\end{remark}

\begin{remark}
It is clear from the definition that a strongly $k$-connective sheaf $\calF$ on $V \subseteq X$
is $k$-connective. Conversely, suppose that $X$ has Heyting dimension $\leq n$ and that
$\calF$ is $k$-connective, then $\calF$ is strongly $(k-n)$-connected (if $\bd \Delta^{m} \rightarrow \calF(U)$ is any map, then we may take $K=U$ for $m > n$ and $K = \emptyset$ for $m \leq n$).
\end{remark}

The strong $k$-connectivity of a sheaf $\calF$ is, by definition, a
local property. The key to our vanishing result is that this is
equivalent to an apparently stronger {\it global} property.

\begin{lemma}\label{precheesit}
Let $X$ be a Heyting space, $V$ a compact open subset of $X$, and
$\calF: \calU(V)^{op} \rightarrow \Kan$ a strongly $k$-connective sheaf on $V$.
Let $A \subseteq B$ be an inclusion of finite simplicial sets of dimension $\leq m$, 
let $U \subseteq V$, and let $\zeta: A \rightarrow \calF(U)$ be a map of simplicial sets.

There exists a cocompact closed subset $K \subseteq U$ whose closure
$\overline{K} \subseteq X$ has Heyting dimension $< m-1-k$, an open covering
$\{ U_{\alpha} \} $ of $U - K$, and a collection of commutative diagrams
$$ \xymatrix{ A \ar@{^{(}->}[d] \ar[r]^-{\zeta} & \calF(U) \ar[d] \\
B \ar[r]^-{\eta_{\alpha}} & \calF( U_{\alpha} ). }$$
\end{lemma}

\begin{proof}
Induct on the number of simplices of $B$ which do not belong to $A$, and invoke Definition \ref{strongcon}.
\end{proof}

\begin{lemma}\label{cheezit}
Let $X$ be a Heyting space, $V$ a compact open subset of $X$, let 
$\calF: \calU(V)^{op} \rightarrow \Kan$ be a sheaf on $X$, let
$\eta: \bd \Delta^m \rightarrow \calF(V)$ be a map, and form a pullback square
$$ \xymatrix{ \calF' \ar[r] \ar[d] & \calF^{\Delta^m} \ar[d] \\
\ast \ar[r]^-{\eta} & \calF^{\bd \Delta^m}. }$$
Suppose that $\calF$ is strongly $k$-connective. Then $\calF'$ is strongly $(k-m)$-connective.
\end{lemma}

\begin{proof}
Unwinding the definitions, we must show that for every compact $U \subset V$ and every
map 
$$\zeta: ( \bd \Delta^m \times \Delta^n ) \coprod_{ \bd \Delta^m \times \bd \Delta^n }
( \Delta^m \times \bd \Delta^n ) \rightarrow \calF(U)$$
whose restriction $\zeta | \bd \Delta^m \times \Delta^n$ is given by $\eta$, there
exists a cocompact closed subset $K \subseteq U$ such that $\overline{K} \subseteq X$ has Heyting dimension $< n+m -k$, an open covering $\{ U_{\alpha} \}$ of $U - K$, and a collection of maps
$$ \zeta_{\alpha}: \Delta^m \times \Delta^n \rightarrow \calF( U_{\alpha})$$ which extend $\zeta$. 
This follows immediately from Lemma \ref{precheesit}.
\end{proof}

\begin{theorem}\label{vanishing}
Let $X$ be a Heyting space of dimension $\leq n$, let $W \subseteq
X$ be a compact open set, and let $\calF \in \Shv(W)$. The
following conditions are equivalent:

\begin{itemize}

\item[$(1)$] For any compact open sets $U \subseteq V \subseteq W$
and any commutative diagram
$$ \xymatrix{ \bd \Delta^{m} \ar[r]^-{\zeta} \ar@{^{(}->}[d] & \calF(V) \ar[d] \\
\Delta^m \ar[r]^-{\eta} & \calF(U), }$$
there exists a cocompact closed subset $K \subseteq V-U$
such that $\overline{K} \subseteq X$ has dimension $< m-k$
and a commutative diagram
$$ \xymatrix{ \bd \Delta^{m} \ar[r]^-{\zeta} \ar@{^{(}->}[d] & \calF(V) \ar[d] \\
\Delta^m \ar[r]^-{\eta'} & \calF(V-K), }$$
such that the composition
$ \Delta^m \stackrel{\eta'}{\rightarrow} \calF(V-K) \rightarrow \calF(U)$
is homotopic to $\eta$ relative to $\bd \Delta^m$.

\item[$(2)$] For any compact open sets $V \subseteq W$ and
any map $\zeta: \bd \Delta^m \rightarrow \calF(V)$, there
exists a commutative diagram
$$ \xymatrix{ \bd \Delta^{m} \ar[r]^-{\zeta} \ar@{^{(}->}[d] & \calF(V) \ar[d] \\
\Delta^m \ar[r]^-{\eta'} & \calF(V-K), }$$
where $K \subseteq V$ is a cocompact closed subset
and $\overline{K} \subseteq X$ has dimension $< m-k$.

\item[$(3)$] The sheaf $\calF$ is strongly $k$-connective.
\end{itemize}
\end{theorem}

\begin{proof}
It is clear that $(1)$ implies $(2)$ (take $U$ to be empty) and
that $(2)$ implies $(3)$ (by definition). We must show that $(3)$
implies $(1)$. So let $\calF$ be a strongly $k$-connective sheaf on $W$
and $$ \xymatrix{ \bd \Delta^{m} \ar[r]^-{\zeta} \ar@{^{(}->}[d]  & \calF(V) \ar[d] \\
\Delta^m \ar[r]^-{\eta} & \calF(U) }$$
a commutative diagram as above. Without loss of generality, we may replace $W$ by $V$ and $\calF$ by $\calF|V$.

We may identify $\calF$ with a functor from $\calU(V)^{op}$ into the category $\Kan$ of Kan complexes. Form a pullback square 
$$ \xymatrix{ \calF' \ar[r] \ar[d] & \calF^{\Delta^m} \ar[d] \\
\ast \ar[r]^-{\zeta} & \calF^{\bd \Delta^m} }$$
in $\Set_{\Delta}^{ \calU(V)^{op} }$. The right vertical map is a projective fibration, so that
the diagram is homotopy Cartesian (with respect to the projective model structure).
It follows that $\calF'$ is also a sheaf on $V$, which is strongly $(k-m)$-connective by Lemma \ref{cheezit}. Replacing $\calF$ by $\calF'$, we may reduce to the case $m=0$.

The proof now proceeds by induction on $k$. For our base case, we take
$k=-n-1$, so that there is no connectivity assumption on the stack
$\calF$. We are then free to choose $K = V-U$ (it is clear that
$\overline{K}$ has dimension $\leq n$).

Now suppose that the theorem is known for strongly
$(k-1)$-connective stacks on any compact open subset of $X$; we
must show that for any strongly $k$-connective $\calF$ on $V$ and
any $\eta \in \calF(U)$, there exists a closed subset
$K \subseteq V-U$ such that $\overline{K} \subseteq X$ has Heyting dimension
$ < -k$, and a point $\eta' \in \calF(V-K)$ whose restriction to $U$ lies in
the same component of $\calF(U)$ as $\eta$.

Since $\calF$ is strongly $k$-connective, we deduce that there
exists an open cover $\{V_{\alpha} \}$ of some open subset
$V-K_0$, where $K_0$ has dimension $< -k$ in $X$, together with
points $\psi_{\alpha} \in \calF( V_{\alpha})$. Adjoining
the open set $U$ and the point $\eta$ if necessary, we may suppose
that $K_0 \cap U = \emptyset$. Replacing $V$ by
$V- K_0$ we may reduce to the case $K_0 = \emptyset$.

Since $V$ is compact, we may assume that there exist only finitely
many indices $\alpha$. Proceeding by induction on the number of
indices, we may reduce to the case where $V = U \cup V_{\alpha}$ for some $\alpha$. 
Let $\eta'$ and $\psi'$ denote the images of $\eta$ and $\psi$ in $U \cap V_{\alpha}$,
and form a pullback
diagram $$ \xymatrix{ \calF' \ar[rr] \ar[d] & & (\calF|(U \cap V_{\alpha}))^{\Delta^1} \ar[d] \\
\ast \ar[rr]^-{(\eta', \psi')} & & (\calF|(U \cap V_{\alpha}))^{\bd \Delta^1}.}$$
Again, this diagram is a homotopy pullback, so that $\calF'$ is a sheaf on
$U \cap V_{\alpha}$ which is strongly $(k-1)$-connective by Lemma \ref{cheezit}.
According to the inductive hypothesis, there exists a closed subset 
$K \subset U \cap V_{\alpha}$ such that $\overline{K} \subseteq X$ has dimension
$< -k+1$, such that the images of $\psi_{\alpha}$ and $\eta$ belong to the same
component of $\calF( (U \cap V_{\alpha}) - K )$. Replacing $V_{\alpha}$ by
Since $\overline{K}$ has dimension $< -k+1$ in $X$, the
boundary $\bd K$ of $K$ has codimension $< -k$ in $X$. Let $V' =
V_{\alpha} - (V_{\alpha} \cap \overline{K})$. Since $\calF$ is a sheaf, we have a 
homotopy pullback diagram
$$ \xymatrix{ \calF( V' \cup U) \ar[r] \ar[d] & \calF(U) \ar[d] \\
\calF(V') \ar[r] & \calF( V' \cap U ). }$$
We observe that there is a path joining the images of $\eta$ and
$\psi_{\alpha}$ in $\calF(V' \cap U) = \calF( (U \cap V_{\alpha}) - K)$, so that there
is a vertex $\widetilde{\eta} \in \calF(V' \cup U)$ whose image in $\calF(U)$ lies
in the same component as $\eta$. We now observe that $V' \cup U = V - (V \cap \bd K)$,
and that $\overline{ V \cap \bd K }$ is contained in $\bd K$ and therefore has Heyting dimension
$\leq -k$.
\end{proof}

\begin{corollary}\label{hphp}
Let $\pi: X \rightarrow Y$ be a continuous map between Heyting
spaces of finite dimension. Suppose that $\pi$ has the property
that for any cocompact closed subset $K \subseteq X$ of dimension
$\leq n$, $\pi(K)$ is contained in a cocompact closed subset of
dimension $\leq n$. Then the functor $\pi_{\ast}: \Shv(X) \rightarrow \Shv(Y)$ 
carries strongly $k$-connective sheaves on $X$ to strongly $k$-connective sheaves on $Y$.
\end{corollary}

\begin{proof}
This is clear from the characterization $(2)$ of Theorem
\ref{vanishing}.
\end{proof}

\begin{corollary}\label{gra}
Let $X$ be a Heyting space of finite Heyting dimension, and let $\calF$ be
a strongly $k$-connective sheaf on $X$. Then $\calF(X)$ is
$k$-connective.
\end{corollary}

\begin{proof}
Apply Corollary \ref{hphp} in the case where $Y$ is a point.
\end{proof}

\begin{corollary}\label{turkant}\index{gen}{Grothendieck's vanishing theorem!nonabelian version}
Let $X$ be a Heyting space of Heyting dimension $\leq n$, and let $\calF$ be an
$n$-connective sheaf on $X$. Then for any compact open $U \subseteq X$, the
map $\pi_0 \calF(X) \rightarrow \pi_0 \calF(U)$ is surjective. In particular, $\Shv(X)$ has homotopy dimension $\leq n$.
\end{corollary}

\begin{proof}
Suppose first that $(1)$ is satisfied. Let $\calF$ be an $n$-connective sheaf on $X$.
Then $\calF$ is strongly $0$-connective; by characterization $(2)$ of Theorem \ref{vanishing}, we deduce that $\calF(X) \rightarrow \calF(U)$ is surjective. The last claim follows by taking
$U = \emptyset$.
\end{proof}

\begin{remark}
Let $X$ be a Heyting space of Heyting dimension $\leq n$. Then any compact open subset
of $X$ also has Heyting dimension $\leq n$. It follows that $\Shv(X)$ is locally of homotopy dimension $\leq n$, and therefore hypercomplete by Corollary \ref{fdfd}.
\end{remark}

\begin{remark}
It is not necessarily true that a Heyting space $X$ such that $\Shv(X)$ has homotopy dimension $\leq n$ is itself of Heyting dimension $\leq n$. For example, if $X$ is the Zariski spectrum of a discrete valuation ring (that is, a two point space with a nontrivial topology), then $X$ has homotopy dimension zero (see Example \ref{honeypie}).
\end{remark}

In particular, we obtain Grothendieck's vanishing theorem (see \cite{tohoku} for the original, quite different proof):

\begin{corollary}\index{gen}{Grothendieck's vanishing theorem}
Let $X$ be a Noetherian topological space of Krull dimension $\leq
n$. Then $X$ has cohomological dimension $\leq n$.
\end{corollary}

\begin{proof}
Combine Proposition \ref{krullheyt}, Corollary \ref{turkant}, and Corollary \ref{confusion}.
\end{proof}

\begin{example}
Let $V$ be a real algebraic variety (defined over the real
numbers, say). Then the lattice of open subsets of $V$ that can be
defined by polynomial equations and inequalities is a Heyting
algebra, and the spectrum of this Heyting algebra is a Heyting
space $X$ having dimension at most equal to the dimension of $V$.
The results of this section therefore apply to $X$.

More generally, let $T$ be an o-minimal theory (see for example
\cite{lou}), and let $S_n$ denote the set of complete $n$-types of
$T$. We endow $S_n$ with the topology generated by the
sets $U_{\phi} = \{p: \phi \in p\}$, where $\phi$ ranges over
formula with $n$ free variables such that the openness of the set of points satisfying
$\phi$ is provable in $T$. Then $S_n$ is a Heyting
space of Heyting dimension $\leq n$.
\end{example}

\begin{remark}
The methods of this section can be adapted to slightly more
general situations, such as the Nisnevich topology on a Noetherian
scheme of finite Krull dimension. It follows that the
$\infty$-topoi associated to such sites have (locally) finite homotopy
dimension and are therefore hypercomplete.
We will discuss this matter in more detail in \cite{DAG}.
\end{remark}
\section{The Proper Base Change Theorem}\label{chap7sec3}

\setcounter{theorem}{0}

Let 
$$ \xymatrix{ X' \ar[r]^{q'} \ar[d]^{p'} & X \ar[d]^{p} \\
Y' \ar[r]^{q} & Y }$$
be a pullback diagram in the category of locally compact Hausdorff spaces. One has a natural isomorphism of pushforward functors $$  q_{\ast} p'_{\ast} \simeq p_{\ast} q'_{\ast}$$
from the category of sheaves of sets on $Y$ to the category of sheaves of sets on $X'$. This isomorphism induces a natural transformation
$$ \eta: q^{\ast} p_{\ast} \rightarrow p'_{\ast} {q'}^{\ast}.$$
If $p$ (and therefore also $p'$) is a proper map, then $\eta$ is an isomorphism: this is a simple version of the classical {\it proper base change theorem}.\index{gen}{proper base change theorem!for sheaves of sets}

The purpose of this section is to generalize the above result, allowing sheaves which take values in the $\infty$-category $\SSet$ of spaces rather than in the ordinary category of sets. Our generalization can be viewed as a proper base change theorem for nonabelian cohomology. 

We will begin in \S \ref{propertopoi} by defining the notion of a {\em proper morphism} of $\infty$-topoi. Roughly speaking, a geometric morphism $\pi_{\ast}: \calX \rightarrow \calY$ of $\infty$-topoi is proper if and only if it satisfies the conclusion of the proper base change theorem. Using this language, our job is to prove that a proper map of topological spaces $p: X \rightarrow Y$ induces a proper morphism $p_{\ast}: \Shv(X) \rightarrow \Shv(Y)$ of $\infty$-topoi. We will outline the proof of this result in \S \ref{propertopoi} by reducing to two special cases: the case where $p$ is a closed embedding, and the case where $Y$ is a point. We will treat the first case in \S \ref{closedsub},
after introducing a general theory of {\em closed immersions} of $\infty$-topoi. This allows us to reduce to the case where $Y$ is a point and $X$ a compact Hausdorff space. Our approach is now in two parts:

\begin{itemize}
\item[$(1)$] In \S \ref{products}, we will show that we can identify the $\infty$-category 
$\Shv(X') = \Shv(X \times Y')$ with an $\infty$-category of sheaves on $X$, taking values in 
$\Shv(Y')$.
\item[$(2)$] In \S \ref{properproper}, we give an analysis of the category of sheaves
on a compact Hausdorff space $X$, taking values in a general $\infty$-category
$\calC$. Combining this analysis with $(1)$, we will deduce the desired base change theorem.
\end{itemize}

The techniques used in \S \ref{properproper} to analyze $\Shv(X)$ can be applied also in the (easier) setting of coherent topological spaces, as we explain in \S \ref{cohthm}. Finally, we conclude in \S \ref{celluj} by reformulating the classical theory of {\em cell-like} maps in the language of $\infty$-topoi.

\subsection{Proper Maps of $\infty$-Topoi}\label{propertopoi}

In this section, we introduce the notion of a {\em proper} geometric morphism between $\infty$-topoi. Here we follow the ideas of \cite{moerdijk}, and turn the conclusion of the proper base change theorem into a definition. First, we require a bit of terminology.

Suppose given a diagram of categories and functors
$$ \xymatrix{ \calC' \ar[r]^{q'_{\ast}} \ar[d]^{p'_{\ast}} & \calD' \ar[d]^{p_{\ast}} \\
\calC \ar[r]^{q_{\ast}} & \calD }$$
which commutes up to a specified isomorphism 
$\eta: p_{\ast} q'_{\ast} \rightarrow
q_{\ast} p'_{\ast}$. Suppose furthermore that the functors $q_{\ast}$ and $q'_{\ast}$ admit
left adjoints, which we will denote by $q^{\ast}$ and ${q'}^{\ast}$. Consider the composition
$$ \gamma: q^{\ast} p_{\ast} \stackrel{u}{\rightarrow} q^{\ast} p_{\ast} q'_{\ast} {q'}^{\ast}
\stackrel{\eta}{\rightarrow} q^{\ast} q_{\ast} p'_{\ast} {q'}^{\ast} \stackrel{v}{\rightarrow}
p'_{\ast} {q'}^{\ast},$$
where $u$ denotes a unit for the adjunction $( {q'}^{\ast}, q'_{\ast})$ and $v$ a counit for the adjunction $( q^{\ast}, q_{\ast})$. We will refer to $\gamma$ as the {\it push-pull} transformation
associated to the above diagram.\index{gen}{push-pull transformation}

\begin{definition}
A diagram of categories $$ \xymatrix{ \calC' \ar[r]^{q'_{\ast}} \ar[d]^{p'_{\ast}} & \calD' \ar[d]^{p_{\ast}} \\ \calC \ar[r]^{q_{\ast}} & \calD }$$ which commutes up to specified isomorphism is
{\it left adjointable} if the functors $q_{\ast}$ and $q'_{\ast}$ admit left adjoints
$q^{\ast}$ and ${q'}^{\ast}$, and the associated push-pull transformation
$$ \gamma: q^{\ast} p_{\ast} \rightarrow p'_{\ast} {q'}^{\ast}$$
is an isomorphism of functors.\index{gen}{left adjointable}
\end{definition}

\begin{definition}
A diagram of $\infty$-categories
$$ \xymatrix{ \calC' \ar[r]^{q'_{\ast}} \ar[d]^{p'_{\ast}} & \calD' \ar[d]^{p_{\ast}} \\ \calC \ar[r]^{q_{\ast}} & \calD }$$ which commutes up to (specified) homotopy is {\it left adjointable} if the associated
diagram of homotopy categories is left adjointable.
\end{definition}

\begin{remark}\label{toadcatcher}
Suppose given a diagram of simplicial sets
$$ \calM' \stackrel{P}{\rightarrow} \calM \stackrel{f}{\rightarrow} \Delta^1,$$
where both $f$ and $f \circ P$ are Cartesian fibrations. Then we may view
$\calM$ as a correspondence from $\calD = f^{-1} \{0\}$ to $\calC = f^{-1} \{1\}$, associated
to some functor $q_{\ast}: \calC \rightarrow \calD$. Similarly, we may view
$\calM'$ as a correspondence from $\calD' = (f \circ P)^{-1} \{0\}$ to
$\calC' = (f \circ P)^{-1} \{1\}$, associated to some functor $q'_{\ast}: \calC' \rightarrow \calD'$.
The map $P$ determines functors $p'_{\ast}: \calC' \rightarrow \calC$, $q'_{\ast}: \calD' \rightarrow \calD$, and (up to homotopy) a natural transformation $\alpha: p_{\ast} q'_{\ast} \rightarrow
q_{\ast} p'_{\ast}$, which is an equivalence if and only if the map $P$ carries
$(f \circ P)$-Cartesian edges of $\calM'$ to $f$-Cartesian edges of $\calM$. In this case,
we obtain a diagram of homotopy categories
$$ \xymatrix{ \h{\calC'} \ar[r]^{q'_{\ast}} \ar[d]^{p'_{\ast}} & \h{\calD'} \ar[d]^{p_{\ast}} \\ \h{\calC} \ar[r]^{q_{\ast}} & \h{\calD} }$$
which commutes up to canonical isomorphism. 

Now suppose that the functors $q_{\ast}$ and $q'_{\ast}$ admit left adjoints, which we will denote by $q^{\ast}$ and ${q'}^{\ast}$, respectively. Then the maps $f$ and $f \circ P$ are coCartesian fibrations. Moreover, the associated push-pull transformation can be described as follows. 
Choose an object $D' \in \calD'$, and a $(f \circ P)$-coCartesian morphism
$\phi: D' \rightarrow C'$, where $C' \in \calC$. Let $D = P(D')$, and choose an $f$-coCartesian morphism $\psi: D \rightarrow C$ in $\calM$, where $C \in \calC$. Using the fact that
$\psi$ is $f$-coCartesian, we can choose a $2$-simplex in $\calM$ depicted as follows:
$$ \xymatrix{ & C \ar[dr]^{\theta} & \\
D \ar[ur]^{\psi} \ar[rr]^{P(\phi)} & & P(C'). }$$ 
We may then identify $C$ with $q^{\ast} p_{\ast} D'$, $P(C')$ with $p'_{\ast} {q'}^{\ast} D'$, and
$\theta$ with the value of the push-pull transformation $q^{\ast} p_{\ast} \rightarrow p'_{\ast} {q'}^{\ast} D'$ on the object $D' \in \calD'$. The morphism $\theta$ is an equivalence if and only if $P(\phi)$ is $f$-coCartesian. Consequently, we deduce that the original diagram
$$ \xymatrix{ \h{\calC'} \ar[r]^{q'_{\ast}} \ar[d]^{p'_{\ast}} & \h{\calD'} \ar[d]^{p_{\ast}} \\ \h{\calC} \ar[r]^{q_{\ast}} & \h{\calD} }$$
is left adjointable if and only if $P$ carries $(f \circ P)$-coCartesian edges to $f$-coCartesian edges. We will make use of this criterion in \S \ref{properproper}.
\end{remark}

\begin{definition}\label{smurk}\index{gen}{proper!morphism of $\infty$-topoi}
Let $p_{\ast}: \calX \rightarrow \calY$ be a geometric morphism of $\infty$-topoi. We will say that
$p_{\ast}$ is {\em proper} if the following condition is satisfied:
\begin{itemize}
\item[$(\ast)$] For every Cartesian rectangle 
$$ \xymatrix{ \calX'' \ar[d] \ar[r] & \calX' \ar[d] \ar[r] & \calX \ar[d]^{p_{\ast}} \\
\calY'' \ar[r] & \calY' \ar[r] & \calY }$$
of $\infty$-topoi, the left square is left adjointable.
\end{itemize}
\end{definition}


\begin{remark}\label{swurk}
Let $\calX$ be an $\infty$-topos, and let $\calJ$ be a small $\infty$-category.
The diagonal functor $\delta: \calX \rightarrow \Fun(\calJ, \calX)$ preserves all (small) limits and colimits, by Proposition \ref{limiteval}, and therefore admits both a left adjoint $\delta_{!}$ and a right adjoint $\delta_{\ast}$. If $\calJ$ is filtered, then $\delta_{!}$ is left exact (Proposition \ref{frent}). Consequently, we have a diagram of geometric morphisms
$$ \calX \stackrel{ \delta }{\rightarrow} \Fun(\calJ,\calX) \stackrel{ \delta_{\ast} }{\rightarrow} \calX. $$

Now suppose that $p_{\ast}: \calX \rightarrow \calY$ is a proper geometric morphism of $\infty$-topoi. We obtain a rectangle
$$ \xymatrix{ \calX \ar[r] \ar[d]^{p_{\ast}} & \Fun(\calJ, \calX) \ar[d]^{p_{\ast}^{\calJ}} \ar[r] & \calX \ar[d] \\
\calY \ar[r] & \Fun(\calJ, \calY) \ar[r] & \calY }$$ 
which commutes up to (specified) homotopy. One can show that this is a Cartesian
rectangle in $\RGeom$, so that the square on the left is left adjointable.
Unwinding the definitions, we conclude that $p_{\ast}$ commutes with filtered colimits. 
Conversely, if
$p_{\ast}: \calX \rightarrow \calY$ is an arbitrary geometric morphism of $\infty$-topoi which
commutes with colimits indexed by {\em filtered $\calY$-stacks} (over each object of $\calY$), then $p_{\ast}$ is proper. To give a proof (or even a precise formulation) of this statement would require ideas from relative category theory which we will not develop in this book. We refer the reader to \cite{moerdijk}, where the analogous result is established for proper maps between ordinary topoi.
\end{remark}

The following properties of the class of proper morphisms follow immediately from Definition \ref{smurk}:

\begin{proposition}\label{properties}
\begin{itemize}
\item[$(1)$] Every equivalence of $\infty$-topoi is proper.
\item[$(2)$] If $p_{\ast}$ and $p'_{\ast}$ are equivalent geometric morphisms from
an $\infty$-topos $\calX$ to another $\infty$-topos $\calY$, then $p_{\ast}$ is proper if and only if $p'_{\ast}$ is proper.
\item[$(3)$] Let 
$$ \xymatrix{ \calX' \ar[d]^{p'_{\ast}} \ar[r] & \calX \ar[d]^{p_{\ast}} \\
\calY' \ar[r] & \calY }$$ be a pullback diagram of $\infty$-topoi. If
$p_{\ast}$ is proper, then so is $p'_{\ast}$.
\item[$(4)$] Let 
$$\calX \stackrel{p_{\ast}}{\rightarrow} \calY \stackrel{q_{\ast}}{\rightarrow} \calZ$$
be proper geometric morphisms between $\infty$-topoi. Then $q_{\ast} \circ p_{\ast}$ is a 
proper geometric morphism.
\end{itemize}
\end{proposition}

In order to relate Definition \ref{smurk} to the classical statement of the proper base change theorem, we need to understand the relationship between products in the category of topological spaces and products in the $\infty$-category of $\infty$-topoi. A basic result asserts that these are compatible, provided that a certain local compactness condition is met.

\begin{definition}\index{gen}{locally compact}\index{gen}{topological space!locally compact}
Let $X$ be a topological space which is not assumed to be Hausdorff. We say that $X$ is {\it locally compact} if, for every open set $U \subseteq X$ and every point $x \in U$, there
exists a (not necessarily closed) compact set $K \subseteq U$, where $K$ contains an open neighborhood of $x$.
\end{definition}

\begin{example}
If $X$ is Hausdorff space, then $X$ is locally compact in the sense defined above if and only if
$X$ is locally compact in the usual sense.
\end{example}

\begin{example}
Let $X$ be a topological space for which the compact open subsets of $X$ form a basis for the topology of $X$. Then $X$ is locally compact.
\end{example}

\begin{remark}
Local compactness of $X$ is precisely the condition which is needed for function
spaces $Y^X$, endowed with the compact-open topology, to represent the functor
$Z \mapsto \Hom( Z \times X, Y)$.
\end{remark}

\begin{proposition}\label{cartmun}
Let $X$ and $Y$ be topological spaces, and assume that $X$ is locally compact.
The diagram
$$ \xymatrix{ \Shv(X \times Y) \ar[r] \ar[d] & \Shv(X) \ar[d] \\
\Shv(Y) \ar[r] & \Shv(\ast) }$$
is a pullback square in the $\infty$-category $\RGeom$ of $\infty$-topoi.
\end{proposition}

\begin{proof}
Let $\calC \subseteq \RGeom$ be the full subcategory spanned by the $0$-localic $\infty$-topoi. Since $\calC$ is a localization of $\RGeom$, the inclusion $\calC \subseteq \RGeom$ preserves limits. It therefore suffices to prove that
$$ \xymatrix{ \Shv(X \times Y) \ar[r] \ar[d] & \Shv(X) \ar[d] \\
\Shv(Y) \ar[r] & \Shv(\ast) }$$
gives a pullback diagram in $\calC$. Note that $\calC^{op}$ is equivalent to the (nerve of the) ordinary category of locales. For each topological space $M$, let $\calU(M)$ denote the locale of open subsets of $M$. Let
$$ \calU(X) \stackrel{\psi_X}{\rightarrow} \calP \stackrel{\psi_{Y}}{\leftarrow} \calU(Y)$$
be a diagram which exhibits $\calP$ as a coproduct of $\calU(X)$ and $\calU(Y)$ in the category of locales, and let $\phi: \calP \rightarrow \calU(X \times Y)$ be the induced map. We wish to prove that $\phi$ is an isomorphism. This is a standard result in the theory of locales; we will include a proof for completeness.

Given open subsets $U \subseteq X$ and $V \subseteq Y$, let
$U \otimes V = ( \psi_{X} U) \cap (\psi_Y V) \in \calP$, so that
$\phi(U \otimes V) = U \times V \in \calU(X \times Y)$. We define a map
$\theta: \calU(X \times Y) \rightarrow \calP$ by the formula
$$ \theta(W) = \bigcup_{ U \times V \subseteq W } U \otimes V.$$
Since every open subset of $X \times Y$ can be written as a union of products $U \times V$, where $U$ is an open subset of $X$ and $V$ an open subset of $Y$, it is clear that
$\phi \circ \theta: \calU( X \times Y) \rightarrow \calU(X \times Y)$ is the identity.
To complete the proof, it will suffice to show that $\theta \circ \phi: \calP \rightarrow \calP$ is the identity. Every element of $\calP$ can be written as $\bigcup_{\alpha} U_{\alpha} \otimes V_{\alpha}$ for $U_{\alpha} \subseteq X$ and $V_{\alpha} \subseteq Y$ appropriately chosen.
We therefore wish to show that
$$ \bigcup_{ U \times V \subseteq \bigcup_{\alpha} U_{\alpha} \otimes V_{\alpha}} U \times V = \bigcup_{\alpha} U_{\alpha} \otimes V_{\alpha}.$$
It is clear that the right hand side is contained in the left hand side. The reverse containment is equivalent to the assertion that if $U \times V \subseteq \bigcup_{\alpha} U_{\alpha} \times V_{\alpha}$, then $U \otimes V \subseteq \bigcup_{\alpha} U_{\alpha} \otimes V_{\alpha}$.

We now invoke the local compactness of $X$.
Write $U = \bigcup K_{\beta}$, where each $K_{\beta}$ is a compact subset of $U$
and the interiors $\{ K^{\degree}_{\beta} \}$ cover $U$.
Then $U \otimes V = \bigcup_{\beta} K^{\degree}_{\beta} \otimes V$; it therefore suffices to prove
that $K^{\degree}_{\beta} \otimes V \subseteq \bigcup_{\alpha} U_{\alpha} \otimes V_{\alpha}$.
Let $v$ be a point of $V$. Then $K_{\beta} \times \{v\}$ is a compact subset of
$\bigcup_{ \alpha} U_{\alpha} \times V_{\alpha}$. Consequently, there exists a finite
set of indices $\{ \alpha_1, \ldots, \alpha_n \}$ such that $v \in V_{v,\beta} = V_{\alpha_1} \cap \ldots \cap V_{\alpha_n}$ and $K_{\beta} \subseteq U_{\alpha_1} \cup \ldots \cup U_{\alpha_n}$. 
It follows that $K^{\degree}_{\beta} \otimes V_{v,\beta} \subseteq \bigcup_{\alpha} U_{\alpha} \otimes V_{\alpha}$. Taking a union over all $v \in V$, we deduce the desired result.
\end{proof}

Let us now return to the subject of the proper base change theorem. We have essentially defined a proper morphism of $\infty$-topoi to be one for which the proper base change theorem holds. The challenge, then, is to produce examples of proper geometric morphisms.
The following results will be proven in \S \ref{closedsub} and \S \ref{properproper}, respectively:

\begin{itemize}
\item[$(1)$] If $p: X \rightarrow Y$ is a closed embedding of topological spaces, then
$p_{\ast}: \Shv(X) \rightarrow \Shv(Y)$ is proper.
\item[$(2)$] If $X$ is a compact Hausdorff space, then the global sections functor
$\Gamma: \Shv(X) \rightarrow \Shv(\ast)$ is proper.
\end{itemize}

Granting these statements for the moment, we can deduce the main result of this section.
First, we must recall a bit of point-set topology:

\begin{definition}\index{gen}{completely regular}\index{gen}{topological space!completely regular}
A topological space $X$ is said to be {\it completely regular} if every point of $X$ is closed in $X$, and for every closed subset $Y \subseteq X$ and every point $x \in X-Y$ there is a continuous function $f: X \rightarrow [0,1]$ such that $f(x) = 0$ and $f|Y$ takes the constant value $1$. 
\end{definition}

\begin{remark}
A topological space $X$ is completely regular if and only if it is homeomorphic to a subspace of a compact Hausdorff space $\overline{X}$ (see \cite{munkres}).
\end{remark}

\begin{definition}\index{gen}{proper!map of topological spaces}
A map $p: X \rightarrow Y$ of (arbitrary) topological spaces is said to be {\it proper} if it is universally closed. In other words, $p$ is proper if and only if for every pullback diagram of topological spaces
$$ \xymatrix{ X' \ar[d]^{p'} \ar[r] & X \ar[d]^{p} \\
Y' \ar[r] & Y }$$
the map $p'$ is closed.
\end{definition}

\begin{remark}
A map $p: X \rightarrow Y$ of topological spaces is proper 
if and only if it is closed and each of the fibers of $p$ is compact (though not necessarily Hausdorff).
\end{remark}

\begin{theorem}\label{pbct}
Let $p: X \rightarrow Y$ be a proper map of topological spaces, where $X$ is completely regular. Then $p_{\ast}: \Shv(X) \rightarrow \Shv(Y)$ is proper.
\end{theorem}

\begin{proof}
Let $q: X \rightarrow \overline{X}$ be an identification of $X$ with a subspace of a compact Hausdorff space $\overline{X}$. The map $p$ admits a factorization
$$ X \stackrel{q \times p}{\rightarrow}\overline{X} \times Y \stackrel{\pi_{Y}}{\rightarrow} Y.$$
Using Proposition \ref{properties}, we can reduce to proving that $(q \times p)_{\ast}$ and $(\pi_{Y})_{\ast}$ are proper.

Because $q$ identifies $X$ with a subspace of $\overline{X}$, $q \times p$ identifies
$X$ with a subspace over $\overline{X} \times Y$. Moreover, $q \times p$ factors
as a composition 
$$ X \rightarrow \overline{X} \times X \rightarrow \overline{X} \times Y$$
where the first map is a closed immersion (since $\overline{X}$ is Hausdorff) and
the second map is closed (since $p$ is proper). It follows that $q \times p$ is a closed immersion,
so that $(q \times p)_{\ast}$ is a proper geometric morphism by Proposition \ref{closeduse2}.

Proposition \ref{cartmun} implies that the geometric morphism
$(\pi_{Y})_{\ast}$ is a pullback of the global sections functor $\Gamma: \Shv( \overline{X} ) \rightarrow \Shv(\ast)$ in the $\infty$-category $\RGeom$. Using Proposition \ref{properties}, we may reduce to proving that $\Gamma$ is proper, which follows from Corollary \ref{compactprop}.
\end{proof}

\begin{remark}\label{pbct2}
The converse to Theorem \ref{pbct} holds as well (and does not require the assumption that $X$ is completely regular): if $p_{\ast}: \Shv(X) \rightarrow \Shv(Y)$ is a proper geometric morphism, then
$p$ is a proper map of topological spaces. This can be proven easily, using the characterization of properness described in Remark \ref{swurk}.
\end{remark}

\begin{corollary}[Nonabelian Proper Base Change Theorem]\index{gen}{proper base change theorem!nonabelian}\label{suman}
Let
$$ \xymatrix{ X' \ar[r]^{q'} \ar[d]^{p'} & X \ar[d]^{p} \\
Y' \ar[r]^{q} & Y }$$
be a pullback diagram of locally compact Hausdorff spaces, and suppose that
$p$ is proper. Then the associated diagram
$$ \xymatrix{ \Shv(X') \ar[r]^{q'_{\ast}} \ar[d]^{p'_{\ast}} & \Shv(X) \ar[d]^{p_{\ast}} \\
\Shv(Y') \ar[r]^{q_{\ast}} & \Shv(Y) }$$ is left adjointable.
\end{corollary}

\begin{proof}
In view of Theorem \ref{pbct}, it suffices to show that 
$$ \xymatrix{ \Shv(X') \ar[r]^{q'_{\ast}} \ar[d]^{p'_{\ast}} & \Shv(X) \ar[d]^{p_{\ast}} \\
\Shv(Y') \ar[r]^{q_{\ast}} & \Shv(Y) }$$
is a pullback diagram of $\infty$-topoi. Let $\overline{X}$ denote a compactification of $X$ (for example, the one-point compactification) and consider the larger diagram of $\infty$-topoi
$$ \xymatrix{ \Shv(X') \ar[r] \ar[d] & \Shv(X) \ar[d] & \\
\Shv(\overline{X} \times Y' ) \ar[r] \ar[d] & \Shv( \overline{X} \times Y) \ar[r] \ar[d] & \Shv( \overline{X} ) \ar[d] \\
\Shv(Y') \ar[r] & \Shv(Y) \ar[r] & \Shv(\ast). }$$
The upper square is a (homotopy) pullback by Proposition \ref{closeduse2} and Corollary \ref{closeduse1}. The lower right square and the lower rectangle are (homotopy) Cartesian
by Proposition \ref{cartmun}, so that the lower left square is (homotopy) Cartesian as well.
It follows that the vertical rectangle is also (homotopy) Cartesian, as desired.
\end{proof}

\begin{remark}\label{classicpbct}
The classical proper base change theorem, for sheaves of abelian groups on locally compact
topological spaces, is a formal consequence of Corollary \ref{suman}. We give a brief sketch.
The usual formulation of the proper base change theorem (see, for example, \cite{kashiwara}) is equivalent to the statement that
if
$$ \xymatrix{ X' \ar[r]^{q'} \ar[d]^{p'} & X \ar[d]^{p} \\
Y' \ar[r]^{q} & Y }$$
is a pullback diagram of locally compact topological spaces, and $p$ is proper, then
the associated diagram
$$ \xymatrix{ D^{-}(X') \ar[r]^{q'_{\ast}} \ar[d]^{p'_{\ast}} & D^{-}(X) \ar[d]^{p_{\ast}} \\
D^{-}(Y') \ar[r]^{q_{\ast}} & D^{-}(Y) }$$
is left adjointable. Here $D^{-}(Z)$ denotes the (bounded below) derived category
of abelian sheaves on a topological space $Z$. 

Let $\bfA$ denote the category whose objects are chain complexes
$$ \ldots \rightarrow A^{-1} \rightarrow A^0 \rightarrow A^1 \rightarrow \ldots $$
of abelian groups. Then $\bfA$ admits the structure of a combinatorial model category, in which the weak equivalences are given by quasi-isomorphisms. Let $\calC = \sNerve( \bfA^{\degree})$ be the underlying $\infty$-category.
For any topological space $Z$, one can define an $\infty$-category
$\Shv(Z;\calC)$ of sheaves on $Z$ with values in $\calC$; see \S \ref{products}. The homotopy
category $\h{\Shv(Z;\calC)}$ is an unbounded version of the derived category
$D^{-}(Z)$; in particular, it contains $D^{-}(Z)$ as a full subcategory. Consequently, we obtain a natural generalization of the proper base change theorem where boundedness hypotheses have been removed, which asserts that the diagram
$$ \xymatrix{ \Shv(X'; \calC) \ar[r]^{q'_{\ast}} \ar[d]^{p'_{\ast}} & \Shv(X;\calC) \ar[d]^{p_{\ast}} \\
\Shv(Y'; \calC) \ar[r]^{q_{\ast}} & \Shv(Y;\calC) }$$
is left adjointable. Using the fact that $\calC$ has enough compact objects, one can deduce this statement formally from Corollary \ref{suman}.
\end{remark}

\subsection{Closed Subtopoi}\label{closedsub}

If $X$ is a topological space and $U \subseteq X$ is an open subset, then we may view
the closed complement $X-U \subseteq X$ as a topological space in its own right. Moreover, the inclusion $(X-U) \hookrightarrow X$ is a proper map of topological spaces (that is, a closed map whose fibers are compact). The purpose of this section is to present an analogous construction in the case where $\calX$ is an $\infty$-topos.

\begin{lemma}\label{cst}
Let $\calX$ be an $\infty$-topos and $\emptyset$ an initial object of $\calX$. Then $\emptyset$ is $(-1)$-truncated.
\end{lemma}

\begin{proof}
Let $X$ be an object of $\calX$. The space $\bHom_{\calX}(X, \emptyset)$ is contractible if $X$ is an initial object of $\calX$, and empty otherwise (by Lemma \ref{sumoto}). In either case,
$\bHom_{\calX}(X, \emptyset)$ is $(-1)$-truncated.
\end{proof}

\begin{lemma}\label{drupe}
Let $\calX$ be an $\infty$-topos and let $f: \emptyset \rightarrow X$ be a morphism in $\calX$, where $\emptyset$ is an initial object. Then $f$ is a monomorphism.
\end{lemma}

\begin{proof}
Apply Lemma \ref{cst} to the $\infty$-topos $\calX_{/X}$.
\end{proof}

\begin{proposition}\label{turb}
Let $\calX$ be an $\infty$-topos and let $U$ be an object of $\calX$. Let
$S_{U}$ be the smallest strongly saturated class of morphisms of $\calX$ which is stable under pullbacks and contains a morphism $f: \emptyset \rightarrow U$, where $\emptyset$ is an initial object of $\calX$. Then $S_{U}$ is topological (in the sense of Definition \ref{deftoploc}).
\end{proposition}

\begin{proof}
For each morphism $g: X \rightarrow U$ in $\calC$, form a pullback square
$$ \xymatrix{ \emptyset' \ar[r]^{f_{Y}} \ar[d] & Y \ar[d]^{g} \\
\emptyset \ar[r]^{f} & U. }$$
Let $S = \{ f_{X} \}_{ g: X \rightarrow U}$ and let $\overline{S}$ be the strongly saturated class of morphisms generated by $S$. We note that each $f_{X}$ is a pullback of $f$, and therefore a monomorphism (by Lemma \ref{drupe}). Let $S'$ be the collection of all morphisms $h: V \rightarrow W$ with the property that for every pullback diagram
$$ \xymatrix{ V' \ar[r] \ar[d]^{h'} & V \ar[d]^{h} \\
W' \ar[r] & W }$$
in $\calX$, the morphism $h$ belongs to $\overline{S}$. Since colimits in $\calX$ are universal, we deduce that $S'$ is strongly saturated, and $S \subseteq S' \subseteq \overline{S}$ by construction. Therefore $S' = \overline{S}$, so that $\overline{S}$ is stable under pullbacks. Since
$f \in \overline{S}$, we deduce that $S_{U} \subseteq \overline{S}$. On the other hand,
$S \subseteq S_{U}$ and $S_{U}$ is strongly saturated, so $\overline{S} \subseteq S_U$.
Therefore $S_{U} = \overline{S}$. Since $S$ consists of monomorphisms, we conclude that $S_{U}$ is topological.
\end{proof}

In the situation of Proposition \ref{turb}, we will say that a morphism of $\calX$ is an {\it equivalence away from $U$} if it belongs to $S_{U}$.\index{gen}{equivalence!away from $U$}

\begin{lemma}\label{charclosedtub}
Let $\calX$ be an $\infty$-topos containing a pair of objects $U,X \in \calX$, and let $S_{U}$ denote the class of morphism in $\calX$ which are equivalences away from $U$. The following are equivalent:
\begin{itemize}
\item[$(1)$] The object $X$ is $S_{U}$-local.
\item[$(2)$] For every map $\widetilde{U} \rightarrow U$ in $\calX$, the space
$\bHom_{\calX}(\widetilde{U}, X)$ is contractible.
\item[$(3)$] There exists a morphism $g: U \rightarrow X$ such that the diagram
$$ \xymatrix{ & U \ar[dl]^{\id_{U}} \ar[dr]^{g} & \\
U & & X }$$
exhibits $U$ as a product of $U$ with $X$ in $\calX$.
\end{itemize}
\end{lemma}

\begin{proof}
Let $S$ be the collection of all morphisms $f_{\widetilde{U}}$ which come from pullback diagrams
$$ \xymatrix{ \emptyset' \ar[r]^{f_{\widetilde{U}}} \ar[d] & \widetilde{U} \ar[d] \\
\emptyset \ar[r] & U }$$
where $\emptyset$ and therefore also $\emptyset'$ are initial objects of $\calX$. We saw in the proof of Proposition \ref{turb} that $S$ generates $S_{U}$ as a strongly saturated class of morphisms. Therefore, $X$ is $S_{U}$-local if and only if each $f_{\widetilde{U}}$ induces
an isomorphism
$$ \bHom_{\calX}( \widetilde{U} , X) \rightarrow \bHom_{\calX}( \emptyset', X) \simeq \ast$$
in the homotopy category $\calH$. This proves that $(1) \Leftrightarrow (2)$.

Now suppose that $(2)$ is satisfied. Taking $\widetilde{U}  = U$, we deduce that there exists a morphism $g: U \rightarrow X$. We will prove that $g$ and $\id_{X}$ exhibit $U$ as a product of $U$ with $X$. As explained in \S \ref{quasilimit5}, this is equivalent to the assertion that for every $Z \in \calX$, the map
$$ \bHom_{\calX}(Z,U) \rightarrow \bHom_{\calX}(Z,U) \times \bHom_{\calX}(Z,X)$$ is an isomorphism in $\calH$. If there are no morphisms from $Z$ to $U$ in $\calX$, then both sides are empty and the result is obvious. Otherwise, we may invoke $(2)$ to deduce that
$\bHom_{\calX}(Z,X)$ is contractible, and the desired result follows. This completes the proof that $(2) \Rightarrow (3)$.

Suppose now that $(3)$ is satisfied for some morphism $g: U \rightarrow X$. For any object
$Z \in \calX$, we have a homotopy equivalence
$$ \bHom_{\calX}(Z,U) \rightarrow \bHom_{\calX}(Z,U) \times \bHom_{\calX}(Z,X).$$
If $\bHom_{\calX}(Z,U)$ is nonempty, then we may pass to the fiber over a point
of $\bHom_{\calX}(Z,U)$ to obtain a homotopy equivalence $\ast \rightarrow \bHom_{\calX}(Z,X)$, so that $\bHom_{\calX}(Z,X)$ is contractible. This proves $(2)$.
\end{proof}

If $\calX$ is an $\infty$-topos and $U \in \calX$, then we will say that an object
$X \in \calX$ is {\it trivial on $U$}\index{gen}{trivial on $U$}\index{not}{X/U@$\calX/U$}
if it satisfies the equivalent conditions of Lemma \ref{charclosedtub}. We let $\calX/U$ denote the full subcategory of $\calX$ spanned by the objects $X$ which are trivial on $U$. It follows from Proposition \ref{turb} that $\calX/U$ is a topological localization of $\calX$, and in particular $\calX/U$ is an $\infty$-topos. We next show that
$\calX/U$ depends only on the support of $U$.

\begin{lemma}\label{yurba}
Let $\calX$ be an $\infty$-topos, and let $g: U \rightarrow V$ be a morphism in
$\calX$. Then $\calX/V \subseteq \calX/U$. Moreover, if $g$ is an effective epimorphism, then
$\calX/U = \calX/V$. 
\end{lemma}

\begin{proof}
The first assertion follows immediately from Lemma \ref{charclosedtub}. To prove the
second, it will suffice to prove that if $g$ is strongly saturated then $S_{V} \subseteq S_{U}$. Since $S_{U}$ is strongly saturated and stable under pullbacks, it will suffice to prove that $S_{U}$ contains a morphism
$f: \emptyset \rightarrow V$, where $\emptyset$ is an initial object of $\calX$.

Form a pullback diagram $\sigma: \Delta^1 \times \Delta^1 \rightarrow \calX$:
$$ \xymatrix{ \emptyset' \ar[r]^{f'} \ar[d] & U \ar[d]^{g} \\
\emptyset \ar[r]^{f} & V.}$$ We may view $\sigma$ as an effective epimorphism
from $f'$ to $f$ in the $\infty$-topos $\calX^{\Delta^1}$. Let $f_{\bigdot} = \mCech(\sigma): \cDelta_{+} \rightarrow \calX^{\Delta^1}$ be a \Cech nerve of $\sigma: f' \rightarrow f$. We note that
for $n \geq 0$, the map $f_n$ is a pullback of $f'$, and therefore belongs to $S_{U}$. Since
$f_{\bigdot}$ is a colimit diagram, we deduce that $f$ belongs to $S_{U}$ as desired.
\end{proof}

If $\calX$ is an $\infty$-topos, we let $\Sub(1_{\calX})$ denote the partially ordered set of equivalence classes of $(-1)$-truncated objects of $\calX$. We note that this set is independent of the choice of a final object $1_{\calX} \in \calX$, up to canonical isomorphism.
Any $U \in \Sub(1_{\calX})$ can be represented by a $(-1)$-truncated object $\widetilde{U} \in \calX$. We define $\calX/U = \calX/ \widetilde{U} \subseteq \calX$. It follows from Lemma \ref{yurba} that $\calX/U$ is independent of the choice of $\widetilde{U}$ representing $U$, and that for any object $X \in \calX$, we
have $\calX/X = \calX/U$ where $U \in \Sub(1_{\calX})$ is the ``support'' of $X$ (namely, the equivalence class of the truncation $\tau_{-1} X$).\index{gen}{support}\index{not}{X/U@$\calX/U$}

\begin{definition}\index{gen}{closed!subtopos}\index{gen}{closed!immersion}
If $\calX$ is an $\infty$-topos and $U \in \Sub(1_{\calX})$, then we will refer to 
 $\calX/U$ as the {\it closed subtopos of $\calX$ complementary to $U$}. More generally, we will say that a geometric morphism $\pi: \calY \rightarrow \calX$ is a {\it closed immersion} if
there exists $U \in \Sub(1_{\calX})$ such that $\pi_{\ast}$ induces an equivalence of $\infty$-categories from $\calY$ to $\calX/U$.
\end{definition}

\begin{proposition}\label{slurppp}
Let $\calX$ be an $\infty$-topos, and let $U \in \Sub(1_{\calX})$. Then the closed immersion
$$ \pi: \calX/U \rightarrow \calX$$ induces an isomorphism of partially ordered sets from
$\Sub( 1_{\calX/U})$ to $\{ V \in \Sub(1_{\calX}): U \subseteq V\} )$.
\end{proposition}

\begin{proof}
Choose a $(-1)$-truncated object $\widetilde{U} \in \calX$ representing $U$.
Since $\pi^{\ast}$ is left exact, an object $X$ of $\calX/U$ is $(-1)$-truncated as an object of $\calX/U$ if and only if it is $(-1)$-truncated as an object of $\calX$. It therefore suffices to prove that
if $\widetilde{V}$ is a $(-1)$-truncated object of $\calX$ representing an element $V \in \Sub(1_{\calX})$, then $\widetilde{V}$ is $S_{U}$-local if and only if $U \subseteq V$. One direction is clear: if $\widetilde{V}$ is $S_U$-local, then we have an isomorphism
$$ \bHom_{\calX}( \widetilde{U}, \widetilde{V} ) \rightarrow \bHom_{\calX}( \emptyset, \widetilde{V}) = \ast$$
in the homotopy category $\calH$, so that $U \subseteq V$. The converse follows from characterization $(3)$ given in Lemma \ref{charclosedtub}.
\end{proof}

\begin{corollary}
Let $\calX$ be an $\infty$-topos, and let $U,V \in \Sub(1_{\calX})$. Then $S_{U} \subseteq S_{V}$ if and only if $U \subseteq V$.
\end{corollary}

\begin{proof}
The ``if'' direction follows from Lemma \ref{yurba} and the converse from Proposition \ref{slurppp}.
\end{proof}

\begin{corollary}\label{glad1}
Let $\calX$ be a $0$-localic $\infty$-topos, associated to the locale $\calU$, and let
$U \in \calU$. Then $\calX/U$ is a $0$-localic $\infty$-topos
associated to the locale $\{ V \in \calU: U \subseteq V \}$.
\end{corollary}

\begin{proof}
The $\infty$-topos $\calX/U$ is a topological localization of a $0$-localic $\infty$-topos, and therefore also $0$-localic (Proposition \ref{useiron}). The identification of the underlying locale follows from Proposition \ref{slurppp}.
\end{proof}

\begin{corollary}\label{closeduse1}
Let $X$ be a topological space, $U \subseteq X$ an open subset and $Y = X - U$. The inclusion
of $Y$ in $X$ induces a closed immersion of $\infty$-topoi $\Shv(Y) \rightarrow \Shv(X)$ and
an equivalence $\Shv(Y) \rightarrow \Shv(X)/U$.
\end{corollary}

\begin{lemma}\label{unipropclose}
Let $\calX$ and $\calY$ be $\infty$-topoi, and let $U \in \calY$ be an object. The map
$$ \Fun_{\ast}(\calX, \calY/U) \rightarrow \Fun_{\ast}(\calX, \calY)$$ identifies
$\Fun_{\ast}(\calX, \calY/U)$ with the full subcategory of $\Fun_{\ast}(\calX,\calY)$ spanned by those geometric morphisms $\pi_{\ast}: \calX \rightarrow \calY$ such that $\pi^{\ast} U$ is an initial object of $\calX$ $($here $\pi^{\ast}$ denotes a left adjoint to $\pi_{\ast}${}$)$.
\end{lemma}

\begin{proof}
Let $\pi_{\ast}: \calX \rightarrow \calY$ be a geometric morphism. Using the adjointness of
$\pi_{\ast}$ and $\pi^{\ast}$, it is easy to see that $\pi_{\ast} X$ is $S_U$-local if and onyl if
$X$ is $\pi^{\ast}(S_U)$-local. In particular, $\pi_{\ast}$ factors through $\calY/U$ if and only if
$\pi^{\ast}(S_U)$ consists of equivalences in $\calX$. 
Choosing $f \in S_{U}$ of the form $f: \emptyset \rightarrow U$, where
$\emptyset$ is an initial object of $\calX$, we deduce that $\pi^{\ast} f$ is an equivalence
so that $\pi^{\ast} U \simeq \pi^{\ast} \emptyset$ is an initial object of $\calX$. 
Conversely, suppose that $\pi^{\ast} U$ is an initial object of $\calX$. Then
$\pi^{\ast} f$ is a morphism between two initial objects of $\calX$, and therefore an equivalence.
Since $\pi^{\ast}$ is left exact and colimit-preserving, the collection of all morphisms
$g$ such that $\pi^{\ast} g$ is an equivalence is strongly saturated, stable under pullbacks, and contains $f$; it therefore contains $S_{U}$, so that $\pi_{\ast}$ factors through $\calY/U$ as desired.
\end{proof}

\begin{proposition}\label{closeduse2}
Let $\pi_{\ast}: \calX \rightarrow \calY$ be a geometric morphism of $\infty$-topoi, and 
let $\pi^{\ast}: \Sub(1_{\calX}) \rightarrow \Sub(1_{\calY})$ denote the induced map of
partially ordered sets. Let $U \in \Sub(1_{\calX})$. There is a commutative diagram
$$ \xymatrix{ \calX/\pi^{\ast}U \ar[rrr]^{\pi_{\ast} | (\calX/\pi^{\ast}U)} \ar[d] & & & \calY/U \ar[d] \\
\calX \ar[rrr]^{\pi_{\ast}} & & &  \calY }$$
of $\infty$-topoi and geometric morphisms, where the vertical maps are given by the natural inclusions. This diagram is left adjointable, and exhibits $\calX/ (\pi^{\ast}U)$ as a fiber product of
$\calX$ and $\calY/U$ over $\calY$ in the $\infty$-category $\RGeom$.
\end{proposition}

\begin{proof}
Let $\pi^{\ast}$ denote a left adjoint to $\pi_{\ast}$.
Our first step is to show that the upper horizontal map $\pi_{\ast} |(\calX/\pi^{\ast}U)$ is well-defined. In other words, we must show that if $X \in \calX$ is trivial on $\pi^{\ast} U$, then
$\pi_{\ast} X \in \calY$ is trivial on $U$. Suppose that $Y \in \calY$ has support contained in $U$; we must show that $\bHom_{\calY}(Y, \pi_{\ast} X)$ is contractible. But this space is homotopy equivalent to $\bHom_{\calX}( \pi^{\ast} Y, X) \simeq \ast$, since $\pi^{\ast} Y$ has support contained in $\pi^{\ast} U$ and $X$ is trivial on $\pi^{\ast} U$.

We note also that $\pi^{\ast}$ carries $\calY/U$ into $\calX/\pi^{\ast}U$. This follows immediately from characterization $(3)$ of Lemma \ref{charclosedtub}, since $\pi^{\ast}$ is left exact. Therefore $\pi^{\ast}| \calY/U$ is a left adjoint of $\pi_{\ast} | \calX/\pi^{\ast}U$. From the fact that $\pi^{\ast}$ is left-exact we easily deduce that $\pi^{\ast}| \calY/U$ is left exact. It follows that $\pi_{\ast}|\calX/\pi^{\ast}U$ has a left-exact left adjoint, and is therefore a geometric morphism of $\infty$-topoi.
Moreover, the diagram
$$ \xymatrix{ \calX/\pi^{\ast}U \ar[d] & & & \calY/U \ar[d] \ar[lll]_{\pi^{\ast} | \calY/Y} \\
\calX & & & \calY \ar[lll]_{\pi^{\ast}} }$$
is (strictly) commutative, which proves that the diagram of pushforward functors is left adjointable.

We now claim that the diagram 
$$ \xymatrix{ \calX/\pi^{\ast}U \ar[rrr]^{\pi_{\ast} | \calX/\pi^{\ast}U} \ar[d] & &  & \calY/U \ar[d] \\
\calX \ar[rrr] & &  & \calY }$$
is a pullback diagram of $\infty$-topoi. For every pair of $\infty$-topoi $\calA$ and $\calB$,
let $[\calA, \calB]$ denote the largest Kan complex contained in $\Fun_{\ast}(\calA, \calB)$. According to Theorem \ref{colimcomparee}, it will suffice to show that for any $\infty$-topos $\calZ$, the associated diagram of Kan complexes
$$ \xymatrix{ [\calZ, \calX/\pi^{\ast}U] \ar[r] \ar[d] & [\calZ, \calY/U] \ar[d] \\
[\calZ, \calX] \ar[r] & [\calZ, \calY] }$$
is homotopy Cartesian. Lemma \ref{unipropclose} implies that the vertical maps are inclusions of full simplicial subsets. 
It therefore suffices to show that if $\phi_{\ast}: \calZ \rightarrow \calY$ is a geometric morphism
such that $\pi_{\ast} \circ \phi_{\ast}$ factors through $\calY/U$, then
$\phi_{\ast}$ factors through $\calX/\pi^{\ast}U$. This follows immediately from the characterization given in Lemma \ref{unipropclose}.
\end{proof}

\begin{corollary}
Let $$ \xymatrix{ \calX' \ar[d]^{p'_{\ast}} \ar[r] & \calX \ar[d]^{p_{\ast}} \\
\calY' \ar[r] & \calY }$$
be a pullback diagram in the $\infty$-category $\RGeom$ of $\infty$-topoi. If
$p_{\ast}$ is a closed immersion, then $p'_{\ast}$ is a closed immersion.
\end{corollary}

\subsection{Products of $\infty$-Topoi}\label{products}\index{gen}{product!of $\infty$-topoi}

In \S \ref{genlim}, we showed that the $\infty$-category $\RGeom$ of $\infty$-topoi admits all
(small) limits. Unfortunately, the construction of general limits was rather inexplicit. Our goal in this section is to give a very concrete description of the product of two $\infty$-topoi, at least in a special case.

\begin{definition}\label{valsheaf}\index{gen}{presheaf!with values in $\calC$}
Let $X$ be a topological space, and $\calC$ an $\infty$-category. We let
$\calU(X)$ denote the collection of all open subsets of $X$, partially ordered by inclusion.
A {\it presheaf on $X$ with values in $\calC$} is a functor
$\calU(X)^{op} \rightarrow \calC$. 

Let $\calF: \calU(X)^{op} \rightarrow \calC$ be a presheaf on $X$ with values in $\calC$. We will say that $\calF$ is a {\it sheaf} with values in $\calC$ if, for every $U \subseteq X$ and every
covering sieve $\calU(X)_{/U}^{(0)} \subseteq \calU(X)_{/U}$, the composition
$$ \Nerve(\calU(X)_{/U}^{(0)})^{\triangleright}
\subseteq \Nerve(\calU(X)_{/U})^{\triangleright} \rightarrow
\Nerve(\calU(X)) \stackrel{\calF}{\rightarrow} \calC^{op}$$
is a colimit diagram.

We let $\calP(X; \calC)$ denote the $\infty$-category
$\Fun( \calU(X)^{op}, \calC)$ consisting of all presheaves on $X$ with values in $\calC$, and
$\Shv(X;\calC)$ the full subcategory of $\calP(X;\calC)$ spanned by the sheaves on
$X$ with values in $\calC$.
\end{definition}

\begin{remark}
We can phrase the sheaf condition informally as follows: a $\calC$-valued presheaf $\calF$ on
a topological space $X$ is a sheaf if, for every open subset $U \subseteq X$ and every
covering sieve $\{ U_{\alpha} \subseteq U \}$, the natural map
$\calF(U) \rightarrow \projlim_{\alpha} \calF(U_{\alpha})$
is an equivalence in $\calC$.
\end{remark}

\begin{remark}
If $X$ is a topological space, then $\Shv(X) = \Shv(X, \SSet)$, where $\SSet$ denotes the $\infty$-category of spaces.
\end{remark}

\begin{lemma}\label{strippus}
Let $\calC$, $\calD$, and $\calE$ be $\infty$-categories which admit finite limits, 
let $\calC^{0} \subseteq \calC$ and $\calD^{0} \subseteq \calD$ be the full subcategories of $\calC$ and $\calD$ consisting of final objects. Let $F: \calC \times \calD \rightarrow \calE$ be a functor.
The following conditions are equivalent:

\begin{itemize}
\item[$(1)$] The functor $F$ preserves finite limits.

\item[$(2)$] 
The functors $F | \calC^{0} \times \calD$ and
$F| \calC \times \calD^{0}$ preserve finite limits, and for every pair of morphisms
$C \rightarrow 1_{\calC}$, $D \rightarrow 1_{\calD}$ where $1_{\calC} \in \calC$
and $1_{\calD} \in \calD$ are final objects, the associated diagram
$$ F(1_{\calC},D) \leftarrow F(C,D) \rightarrow F(C, 1_{\calD})$$ exhibits
$F(C,D)$ as a product of $F(1_{\calC},D)$ and $F(C,1_{\calD})$ in $\calE$.

\item[$(3)$] The functors $F | \calC^{0} \times \calD$ and
$F| \calC \times \calD^{0}$ preserve finite limits, and $F$ is a right
Kan extension of the restriction
$$F^0=  F | (\calC \times \calD^{0}) \coprod_{ \calC^{0} \times \calD^{0} } (\calC^{0} \times \calD).$$
\end{itemize}
\end{lemma}

\begin{proof}
The implication $(1) \Rightarrow (2)$ is obvious. To see that $(2) \Rightarrow (1)$, we
choose final objects $1_{\calC} \in \calC$, $1_{\calD} \in \calD$, and natural transformations
$\alpha: \id_{\calC} \rightarrow \underline{1}_{\calC}$, $\beta: \id_{\calD} \rightarrow \underline{1}_{\calD}$ (where $\underline{X}$ denotes the constant functor with value $X$).
Let $F_{\calC}: \calC \rightarrow \calE$ denote the composition
$$ \calC \simeq \calC \times \{1_{\calD} \} \subseteq \calC \times \calD \stackrel{F}{\rightarrow} \calE,$$ and define $F_{\calD}$ similarly. Then $\alpha$ and $\beta$ induce natural transformations
$$ F_{\calC} \circ \pi_{\calC} \leftarrow F \rightarrow F_{\calD} \circ \pi_{\calD}.$$
Assumption $(2)$ implies that the functors $F_{\calC}$, $F_{\calD}$ preserve finite limits, and that the above diagram exhibits $F$ as a product of $F_{\calC} \circ \pi_{\calC}$ with
$F_{\calD} \circ \pi_{\calD}$ in the $\infty$-category $\calE^{\calC \times \calD}$. We now apply Lemma \ref{limitscommute} to deduce that $F$ preserves finite limits as well.

We now show that $(2) \Leftrightarrow (3)$. Assume that $F| \calC^{0} \times \calD$
and $F| \calC \times \calD^{0}$ preserve finite limits, so that in particular 
$F | \calC^{0} \times \calD^{0}$ takes values in the full subcategory $\calE^0 \subseteq \calE$ spanned by the final objects. Fix morphisms $u: C \rightarrow 1_{\calC}$, $v: D \rightarrow 1_{\calD}$, where $1_{\calC} \in \calC$ and $1_{\calD} \in \calD$ are final obejcts. We will show that
the diagram $$ F(1_{\calC},D) \leftarrow F(C,D) \rightarrow F(C, 1_{\calD})$$ exhibits
$F(C,D)$ as a product of $F(1_{\calC},D)$ and $F(C,1_{\calD})$ if and only if
$F$ is a right Kan extension of $F^0$ at $(C,D)$.

The morphisms $u$ and $v$ determine a map $u \times v: \Delta^1 \times \Delta^1 \rightarrow \calC \times \calD$, which we may identify with a map
$$ w: \Lambda^2_2 \rightarrow ((\calC^{0} \times \calD) \coprod_{ \calC^0 \times \calD^0} (\calC \times \calD^{0}))_{(C,D)/}.$$
Using Theorem \ref{hollowtt}, it is easy to see that $w^{op}$ is cofinal. Consequently,
$F$ is a right Kan extension of $F^{0}$ at $(C,D)$ if and only if the diagram
$$ \xymatrix{ F(C,D) \ar[r] \ar[d] & F(C, 1_{\calD}) \ar[d] \\
F(1_{\calC}, D) \ar[r] & F(1_{\calC}, 1_{\calD}) }$$
is a pullback square. Since $F(1_{\calC}, 1_{\calD})$ is a final object of $\calE$, this is
equivalent to assertion $(2)$.
\end{proof}

\begin{lemma}\label{sablesilk}
Let $\calC$ and $\calD$ be small $\infty$-categories which admit finite limits, and let
$1_{\calC} \in \calC$, $1_{\calD} \in \calD$ be final objects, and let $\calX$ be
an $\infty$-topos. The projections
$$ \calP(\calC \times \{1_{\calD} \}) \stackrel{p_{\ast}}{\leftarrow} \calP(\calC \times \calD) \stackrel{q_{\ast}}{\rightarrow} \calP(
\{1_{\calC} \} \times \calD)$$
induce a categorical equivalence
$$ \Fun_{\ast}(\calX, \calP(\calC \times \calD)) \rightarrow \Fun_{\ast}(\calX, \calP(\calC)) \times
\Fun_{\ast}(\calX, \calP(\calD)).$$
In particular, $\calP(\calC \times \calD)$ is a product of $\calP(\calC)$ with
$\calP(\calD)$ in the $\infty$-category $\RGeom$ of $\infty$-topoi.
\end{lemma}

\begin{proof}
For every $\infty$-category $\calY$ which admits finite limits, let
$[\calY,\calX]$ denote the full subcategory of $\Fun(\calY,\calX)$ spanned by the left
exact functors $\calY \rightarrow \calX$. If $\calY$ is an $\infty$-topos, we let
$[\calY, \calX]_0$ denote the full subcategory of $[\calY, \calX]$ spanned by
the {\em colimit-preserving} left exact functors $\calY \rightarrow \calX$. In view of
Proposition \ref{switcheroo} and Remark \ref{switcheroo2}, it will suffice to prove that composition with the left adjoints to
$p_{\ast}$ and $q_{\ast}$ induces an equivalence of $\infty$-categories
$$ [\calP(\calC \times \calD),\calX]_0 \rightarrow [\calP(C),\calX]_0 \times [\calP(\calD),\calX]_0.$$
Applying Proposition \ref{igrute}, we may reduce to the problem of showing that the map
$$ [ \calC \times \calD, \calX] \rightarrow [\calC, \calX] \times [\calD, \calX]$$
is an equivalence of $\infty$-categories. 

Let $\calC^{0} \subseteq \calC$ and $\calD^{0} \subseteq \calD$ denote the full subcategories
consisting of final objects of $\calC$ and $\calD$, respectively. Proposition \ref{initunique}
implies that $\calC^{0}$ and $\calD^{0}$ are contractible. It will therefore suffice to prove
that the restriction map
$$ \phi: [ \calC \times \calD, \calX] \rightarrow [\calC \times \calD^{0}, \calX]
\times_{[ \calC^0 \times \calD^{0}, \calX] } [\calC^{0} \times \calD, \calX]$$
is a trivial fibration of simplicial sets. This follows immediately from
Lemma \ref{strippus} and Proposition \ref{lklk}.
\end{proof}

\begin{notation}\index{not}{timesO@$\otimes$}\index{not}{timesOC@$\otimes^{\calC}$}
Let $\calX$ be an $\infty$-topos, and $p^{\ast}: \SSet \rightarrow \calX$ a geometric
morphism (essentially unique in view of Proposition \ref{spacefinall}). Let $\pi_{\calX}: \calX \times \SSet \rightarrow \calX$ and $\pi_{\SSet}: \calX \times \SSet \rightarrow \SSet$ denote the projection functors. Let $\otimes$ be a product of $\pi_{\calX}$ with $p^{\ast} \circ \pi_{\SSet}$ in the
$\infty$-category of functors from $\calX \times \SSet$ to $\calX$. Then $\otimes$
is uniquely defined up to equivalence, and we have natural transformations
$$ X \leftarrow X \otimes S \rightarrow p^{\ast} S$$
which exhibit $X \otimes S$ as product of $X$ with $p^{\ast} S$ for all $X \in \calX$, $S \in \SSet$.
We observe that $\otimes$ preserves colimits separately in each variable.

If $\calC$ is a small $\infty$-category, we let 
$\otimes^{\calC}$ denote the composition
$$ \calP(\calC; \calX) \times \calP(\calC) \simeq \calP(\calC; \calX \times \SSet)
\stackrel{\circ \otimes}{\rightarrow} \calP(\calC, \calX).$$
We observe that if $F \in \calP(\calC; \calX)$ and $G \in \calP(\calC)$, then
$F \otimes^{\calC} G$ can be identified with a product of $F$ with $p^{\ast} \circ G$ in
$\calP(\calC; \calX)$.
\end{notation}

\begin{lemma}\label{goldenboy}
Let $\calC$ be a small $\infty$-category, $\calX$ an $\infty$-topos. Let
$g: \calX \rightarrow \SSet$ a functor corepresented by an object $X \in \calX$,
and $G: \calP(\calC; \calX) \rightarrow \calP(\calC)$ the induced functor.
Let $\underline{X} \in \calP(\calC; \calX)$ denote the constant functor with the value $X$. Then
the functor
$$ F = \underline{X} \otimes^{\calC} \id_{ \calP(\calC)}.$$
is a left adjoint to $G$.
\end{lemma}

\begin{proof}
Since adjoints and $\otimes^{\calC}$ can both be computed pointwise on $\calC$, it suffices
to treat the case where $\calC = \Delta^0$. In this case, we deduce the existence of a left adjoint
$F'$ to $G$ using Corollary \ref{adjointfunctor} (the accessibility of $G$ follows from the fact
that $X$ is $\kappa$-compact for sufficiently large $\kappa$, since $\calX$ is accessible). Now
$F$ and $F'$ are both colimit-preserving functors $\SSet \rightarrow \calX$. In virtue of Theorem \ref{charpresheaf}, to prove that $F$ and $F'$ are equivalent, it will suffice to show that the objects
$F(\ast), F'(\ast) \in \calX$ are equivalent. In other words, we must prove that $F'(\ast) \simeq X$. By adjointness, we have natural isomorphisms
$$ \bHom_{\calX}( F'(\ast), Y) \simeq \bHom_{\calH}( \ast, G(Y)) \simeq \bHom_{\calX}(X,Y)$$
in $\calH$ for each $Y \in \calX$, so that $F'(\ast)$ and $X$ corepresent the same functor on
the homotopy category $\h{\calX}$, and are therefore equivalent by Yoneda's lemma.
\end{proof}

\begin{lemma}\label{goldenrod}
Let $\calC$ be a small $\infty$-category which admits finite limits and contains
a final object $1_{\calC}$, let $\calX$ and $\calY$ be $\infty$-topoi, and let 
$p_{\ast}: \calX \rightarrow \SSet$ be a geometric morphism $($essentially unique, in virtue of Proposition \ref{spacefinall}{}$)$. Then the maps 
$$ \calP(\calC) \stackrel{P_{\ast}}{\leftarrow} \calP(\calC; \calX) \stackrel{e_{1_{\calC}}}{\rightarrow} \calX$$
induce equivalences of $\infty$-categories
$$ \Fun_{\ast}( \calY, \calP(\calC; \calX)) \rightarrow \Fun_{\ast}(\calY, \calX) \times \Fun_{\ast}(\calY, \calP(\calC)).$$
In particular, $\calP(\calC; \calX)$ is a product of $\calP(\calC)$ and $\calX$ in the 
$\infty$-category $\RGeom$ of $\infty$-topoi.
Here $e_{1_{\calC}}$ denote the evaluation map at the object $1_{\calC} \in \calC$, 
and $P_{\ast}: \calP(\calC; \calX) \rightarrow \calP(\calC)$ is given by composition with
$p_{\ast}$. 
\end{lemma}

\begin{proof}
According to Proposition \ref{precisechar}, we may assume without loss of generality
that there exists a small $\infty$-category $\calD$ such that $\calX$ is the essential
image of an accessible left exact localization functor $L: \calP(\calD) \rightarrow \calP(\calD)$, 
and that $p_{\ast}$ is given by evaluation at a final object $1_{\calD} \in \calD$.
We have a commutative diagram
$$ \xymatrix{ \Fun_{\ast}(\calY, \calP(\calC; \calX)) \ar[r] \ar[d] & \Fun_{\ast}(\calY, \calP(\calC))
\times \Fun_{\ast}(\calY, \calX) \ar[d] \\
\Fun_{\ast}(\calY, \calP(\calC \times \calD)) \ar[r] & \Fun_{\ast}(\calY, \calP(\calC))
\times \Fun_{\ast}(\calY, \calP(\calD))}$$
where the vertical arrows are inclusions of full subcategories, and the bottom arrow is an equivalence of $\infty$-categories by Lemma \ref{sablesilk}. Consequently, it will suffice
to show that if $q_{\ast}: \calY \rightarrow \calP(\calC \times \calD)$ is a geometric morphism with the property that the composition
$$ r_{\ast}: \calY \rightarrow \calP(\calC \times \calD) \rightarrow \calP(\calD)$$
factors through $\calX$, then $q_{\ast}$ factors through $\calP(\calC; \calX)$. 

Let $Y \in \calY$ and $C \in \calC$; we wish to show that
$q_{\ast}(Y)(C) \in \calX$. It will suffice to show that if $s: D \rightarrow D'$ is
a morphism in $\calP(\calD)$ such that $L(s)$ is an equivalence in $\calX$, then
$q_{\ast}(Y)(C)$ is $s$-local. Let $F: \calP(\calD) \rightarrow \calP(\calC \times \calD)$
be a left adjoint to the functor given by evaluation at $C$. 
We have a commutative diagram
$$ \xymatrix{  \bHom_{\calP(\calD)}( D', q_{\ast}(Y)(C)) \ar[r] \ar[d] & \bHom_{\calY}(q^{\ast} F(D'), Y) \ar[d] \\
\bHom_{\calP(\calD)}( D, q_{\ast}(Y)(C) ) \ar[r] & \bHom_{\calY}(
q^{\ast} F(D) , Y) }$$
where the horizontal arrows are homotopy equivalences. Consequently, to prove that the left vertical map is an equivalence, it will suffice to prove that $q^{\ast} F(s)$ is an equivalence in $\calY$. According to Lemma \ref{goldenboy}, the functor
$F$ can be identified with a product of a left adjoint $r^{\ast}$ to the projection
$r_{\ast}: \calP(\calC \times \calD) \rightarrow \calP(\calD)$ with a constant functor. Since $q^{\ast}$ preserves finite products,
it will suffice to show that $(q^{\ast} \circ r^{\ast})(s)$ is an equivalence in $\calY$. This follows immediately
from our assumption that $r_{\ast} \circ q_{\ast}: \calY \rightarrow \calP(\calD)$ factors through
$\calX$.
\end{proof}

The main result of this section is the following:

\begin{theorem}\label{conprod}
Let $X$ be a topological space, $\calX$ an $\infty$-topos, and
$\pi_{\ast}: \calX \rightarrow \SSet$ a geometric morphism $($which is essentially unique,
by virtue of Proposition \ref{spacefinall}{}$)$. Then $\Shv(X; \calX)$ is an $\infty$-topos, and the
diagram
$$ \calX \stackrel{\Gamma}{\leftarrow} \Shv(X; \calX) \stackrel{\pi_{\ast}}{\rightarrow} \Shv(X)$$
exhibits $\Shv(X;\calX)$ as a product of $\Shv(X)$ and $\calX$ in the $\infty$-category
$\RGeom$ of $\infty$-topoi. Here $\Gamma$ denotes the global sections functor, given by evaluation at $X \in \calU(X)$.
\end{theorem}

\begin{proof}
We first show that $\Shv(X; \calX)$ is an $\infty$-topos. Let $\calP(X; \calX)$ be the $\infty$-category
$\Fun(\Nerve(\calU(X))^{op},\calX)$ of $\calX$-valued presheaves on $\calX$. For each object
$Y \in \calX$, choose a morphism $e_Y: \emptyset_{\calX} \rightarrow Y$ in $\calX$, whose source is
an initial object of $\calX$.
For each sieve $\calV$ on $\calX$, let $\chi^Y_{\calV}: \calU(X)^{op} \rightarrow \calX$
be the composition
$$ \calU(X)^{op} \rightarrow \Delta^1 \stackrel{e_Y}{\rightarrow} \calX,$$
so that $$\chi^Y_{\calV}(U) \begin{cases} Y & \text{if } U \in \calV \\
\emptyset_{\calX} & \text{if } U \notin \calV.
\end{cases}$$
so that we have a natural map $\chi^Y_{\calV} \rightarrow \chi^Y_{\calV'}$ if
$\calV \subseteq \calV'$. For each open subset $U \subseteq X$, let
$\chi^{Y}_{U}= \chi^{Y}_{\calV}$, where $\calV = \{ V \subseteq U \}$. 
Let $S$ be the set of all morphisms
$f^{Y}_{\calV}: \chi^{Y}_{\calV} \rightarrow \chi^{Y}_{U}$, where $\calV$ is a sieve covering $U$,
and let $\overline{S}$ be the strongly saturated class of morphisms generated by $\calX$. We first claim that $\overline{S}$ is setwise generated. To see this, we observe that the passage from
$Y$ to $f^{Y}_{\calV}$ is a colimit-preserving functor of $Y$, so it suffices to consider a set of objects $Y \in \calX$ which generates $\calX$ under colimits.

We next claim that $\overline{S}$ is topological, in the sense of Definition \ref{deftoploc}.
By a standard argument, it will suffice to show that there is a class of objects
$F_{\alpha} \in \calP(X; \calX)$ which generates $\calP(X; \calX)$ under colimits, such that
for every pullback diagram
$$ \xymatrix{ F'_{\alpha} \ar[d]^{f'} \ar[r] & \chi^{Y}_{\calV} \ar[d]^{f^{Y}_{\calV}} \\
F_{\alpha} \ar[r] & \chi^{Y}_{U}, }$$
the morphism $f'$ belongs to $\overline{S}$. We observe that if
$\calX$ is a left exact localization of $\calP(\calD)$, then $\calP(X; \calX)$ is a left
exact localization of $\calP( \calU(X) \times \calD)$ and is therefore generated under colimits
by the Yoneda image of $\calU(X) \times \calD$. In other words, it will suffice to consider 
$F_{\alpha}$ of the form $\chi^{Y'}_{U'}$, where $Y' \in \calX$ and $U' \subseteq X$.
If $Y'$ is an initial object of $\calX$, then $g$ is an equivalence and there is nothing to prove. Otherwise, the existence of the lower horizontal map implies that $U' \subseteq U$. 
Let $\calV' = \{ V \in \calV: V \subseteq U' \}$; then it is easy to see that $f'$ is
equivalent to $\chi^{Y'}_{\calV'}$, and therefore belongs to $\overline{S}$.

We next claim that $\Shv(X; \calX)$ consists precisely of the $S$-local objects of
$\calP(X; \calX)$. To see this, let $Y \in \calX$ be an arbitrary object, and consider
the functor $G_Y: \calX \rightarrow \SSet$ corepresented by $Y$. It follows from Proposition \ref{yonedaprop} that an arbitrary $F \in \calP(X; \calX)$ is a $\calX$-valued sheaf on $X$ if and only if, for each $Y \in \calX$, the composition $G_Y \circ F \in \calP(X)$ belongs to $\Shv(X)$.
This is equivalent to the assertion that, for every sieve $\calV$ which covers $U \subseteq X$,
the presheaf $G_Y \circ F$ is $s_{\calV}$-local, where $s_{\calV}: \chi_{\calV} \rightarrow \chi_{U}$ is the associated monomorphism of presheaves. Let $G^{\ast}_Y$ denote a left adjoint to $G_{Y}$; then $G_Y \circ F$ is $s_{\calV}$-local if and only if $F$ is $G^{\ast}_{Y}(s_{\calV})$-local. We now apply Lemma \ref{goldenboy} to identify $G^{\ast}_{Y}(s_{\calV})$ with $f_{\calV}^{Y}$.

We now have an identification $\Shv(X; \calX) \simeq \overline{S}^{-1} \calP(X; \calX)$, so that
$\Shv(X; \calX)$ is a topological localization of $\calP(X; \calX)$ and in particular an $\infty$-topos. 
We now consider an arbitrary $\infty$-topos $\calY$. We have a commutative diagram
$$ \xymatrix{ \Fun_{\ast}(\calY, \Shv(X; \calX)) \ar[r] \ar[d] & \Fun_{\ast}(\calY, \Shv(X)) \times \Fun_{\ast}(\calY; \calX) \ar[d] \\
\Fun_{\ast}(\calY, \calP(X; \calX)) \ar[r] & \Fun_{\ast}(\calY,\calP(X)) \times \Fun_{\ast}(\calY, \calX), } $$
where the vertical arrows are inclusions of full subcategories and the lower horizontal arrow is an equivalence by Lemma \ref{goldenrod}. To complete the proof, it will suffice to show that the upper horizontal arrow is also an equivalence. In other words, we must show that if
$g_{\ast}: \calY \rightarrow \calP(X; \calX)$ is a geometric morphism with the property that
the composition 
$$ \calY \stackrel{g_{\ast}}{\rightarrow} \calP(X; \calX) \stackrel{h_{\ast}}{\rightarrow} \calP(X) $$
factors through $\Shv(X)$, then $g_{\ast}$ factors through $\Shv(X; \calX)$.
Let $g^{\ast}$ and $h^{\ast}$ denote left adjoints to $g_{\ast}$ and $h_{\ast}$, respectively. It will suffice to show that for every
morphism $f_{\calV}^{Y} \in S$, the pullback $g^{\ast} f_{\calV}^{Y}$ is an equivalence in $\calY$.
We now observe that $f_{\calV}^{Y}$ is a pullback of $f_{\calV}^{1_{\calX}}$; since
$g^{\ast}$ is left exact, it will suffice to show that $g^{\ast} f_{\calV}^{1_{\calX}}$
is an equivalence in $\calY$. We have an equivalence $f_{\calV}^{1_{\calX}} \simeq h^{\ast} s_{\calV}$, where $s_{\calV}$ is the monomorphism in $\calP(X)$ associated to the sieve $\calV$.
The composition $(g^{\ast} \circ h^{\ast})( s_{\calV} )$ is an equivalence because 
$h_{\ast} \circ g_{\ast}$ factors through $\Shv(X)$, which consists of $s_{\calV}$-local objects
of $\calP(X)$.
\end{proof}

\begin{remark}\index{gen}{fiber product!of $\infty$-topoi}
It is not difficult to extend the proof of Theorem \ref{conprod} to the case where
$\Shv(X)$ is replaced by an arbitrary $\infty$-topos $\calY$. In this case, one must replace
$\Shv(X; \calX)$ by the $\infty$-category of all limit-preserving functors $\calY^{op} \rightarrow \calX$. Using these ideas, one can construct the fiber product
$$ \calX \times_{ \calZ} \calY$$ in $\RGeom$ where
$\calZ = \SSet$ is the final object in $\RGeom$. To give a construction which works in general,  one needs to repeat all of the above arguments in a relative setting over the $\infty$-topos $\calZ$. We will not pursue the subject any further in this book.
\end{remark}

\subsection{Sheaves on Locally Compact Spaces}\label{properproper}

By definition, a sheaf of sets $\calF$ on a topological space $X$ is determined by the sets
$\calF(U)$ as $U$ ranges over the open subsets of $X$. If $X$ is a locally compact Hausdorff space, then there is an alternative collection of data which determines $X$: the values
$\calF(K)$, where $K$ ranges over the compact subsets of $X$. Here $\calF(K)$ denotes the direct limit $\colim_{K \subseteq U} \calF(U)$ taken over all open neighborhoods of $K$ (or, equivalently, the collection of global sections of the restriction $\calF|K$). The goal of this section is to prove a generalization of this result, where the sheaf $\calF$ is allowed to take values in a more general $\infty$-category $\calC$.

\begin{definition}\index{not}{Kcal(X)@$\calK(X)$}
Let $X$ be a locally compact Hausdorff space. We let $\calK(X)$ denote the collection of all compact subsets of $X$.
If $K, K' \subseteq X$, we write $K \Subset K'$ if there exists an open subset $U \subseteq X$ such that $K \subseteq U \subseteq K'$. If $K \in \calK(X)$, we let $\calK_{K \Subset}(X) = \{ K' \in \calK(X): K \Subset K' \}$.\index{not}{KSubsetK'@$K \Subset K'$}\index{not}{KcalKSubsetX@$\calK_{K \Subset}(X)$}

Let $\calF: \Nerve(\calK(X))^{op} \rightarrow \calC$ be a presheaf on $\Nerve(\calK(X))$ (here $\calK(X)$ is viewed as a partially ordered set with respect to inclusion) with values in $\calC$.
We will say that
$\calF$ is a {\it $\calK$-sheaf} if the following conditions are satisfied:\index{gen}{$\calK$-sheaf}
\begin{itemize}
\item[$(1)$] The object $\calF(\emptyset) \in \calC$ is final.
\item[$(2)$] For every pair $K, K' \in \calK(X)$, the associated diagram
$$ \xymatrix{ \calF( K \cup K' ) \ar[r] \ar[d] & \calF(K) \ar[d] \\
\calF(K') \ar[r] & \calF(K \cap K') }$$
is a pullback square in $\calC$.
\item[$(3)$] For each $K \in \calK(X)$, the restriction of $\calF$ exhibits
$\calF(K)$ as a colimit of $\calF | \Nerve(\calK_{K \Subset}(X))^{op}$.
\end{itemize}

We let $\Shv_{\calK}(X;\calC)$ denote the full subcategory of $\Fun(\Nerve(\calK(X))^{op},\calC)$ spanned by the $\calK$-sheaves. In the case where $\calC = \SSet$, we will write
$\Shv_{\calK}(X)$ instead of $\Shv_{\calK}(X;\calC)$.\index{not}{ShvcalKXC@$\Shv_{\calK}(X;\calC)$}\index{not}{ShvcalKX@$\Shv_{\calK}(X)$}
\end{definition}

\begin{definition}\label{leftexactcolim}\index{gen}{filtered colimit!left exactness of}
Let $\calC$ be a presentable $\infty$-category. We will say that {\it filtered colimits
in $\calC$ are left exact} if the following condition is satisfied: for every small filtered $\infty$-category $\calI$, the colimit functor $\Fun(\calI,\calC) \rightarrow \calC$
is left exact.
\end{definition}

\begin{example}\index{gen}{Grothendieck abelian category}\index{gen}{abelian category!Grothendieck}
A {\it Grothendieck abelian category} is an abelian category $\calA$ whose nerve 
$\Nerve(\calA)$ is a presentable $\infty$-category with left exact filtered colimits, in the sense of Definition \ref{leftexactcolim}. We refer the reader to \cite{tohoku} for further discussion.
\end{example}

\begin{example}\label{sumta1}
Filtered colimits are left exact in the $\infty$-category $\SSet$ of spaces; this follows immediately from Proposition \ref{frent}. It follows that filtered colimits in $\tau_{\leq n} \SSet$ are left exact for each $n \geq -2$, since the full subcategory $\tau_{ \leq n} \SSet \subseteq \SSet$ is stable under filtered colimits and finite limits (in fact, under all limits).
\end{example}

\begin{example}\label{sumta2}
Let $\calC$ be a presentable $\infty$-category in which filtered colimits are left exact, and let $X$ be an arbitrary simplicial set. Then filtered colimits are left exact in $\Fun(X,\calC)$. This follows
immediately from Proposition \ref{limiteval}, which asserts that the relevant limits and colimits can be computed pointwise.
\end{example}

\begin{example}\label{sumta3}
Let $\calC$ be a presentable $\infty$-category in which filtered colimits are left exact, and let
$\calD \subseteq \calC$ be the essential image of an (accessible) left exact localization functor $L$. Then filtered colimits in $\calD$ are left exact. To prove this, we consider an arbitrary filtered $\infty$-category $\calI$, and observe that the colimit functor $\varinjlim: \Fun(\calI, \calD) \rightarrow \calD$ is equivalent to the composition
$$ \Fun(\calI, \calD) \subseteq \Fun(\calI,\calC) \rightarrow \calC \stackrel{L}{\rightarrow} \calD,$$
where the second arrow is given by the colimit functor $\varinjlim \Fun(\calI,\calC) \rightarrow \calC$.
\end{example}

\begin{example}\label{tucka}
Let $\calX$ be an $n$-topos, $0 \leq n \leq \infty$. Then filtered colimits in $\calX$ are left exact. This follows immediately from Examples \ref{sumta1}, \ref{sumta2}, and \ref{sumta3}.
\end{example}

Our goal is to prove that if $X$ is a locally compact Hausdorff space and $\calC$ is a presentable $\infty$-category, then the $\infty$-categories $\Shv(X)$ and $\Shv_{\calK}(X)$ are equivalent. As a first step, we prove that a $\calK$-sheaf on $X$ is determined ``locally''.

\begin{lemma}\label{noodlesoup}
Let $X$ be a locally compact Hausdorff space and $\calC$ a presentable $\infty$-category in which filtered colimits are left exact.
Let $\calW$ be a collection of open
subsets of $X$ which covers $X$, and let $\calK_{\calW}(X) = \{ K \in \calK(X): (\exists W \in \calW) [ K \subseteq W] \}$. Suppose that $\calF \in \Shv_{\calK}(X;\calC)$. Then
$\calF$ is a right Kan extension of $\calF | \Nerve(\calK_{\calW}(X))^{op}$. 
\end{lemma}

\begin{proof}
Let us say that an open covering $\calW$ of a locally compact Hausdorff space $X$ is {\it good} if it satisfies the conclusion of the Lemma. 
Note that $\calW$ is a good covering of $X$ if and only if, for every compact subset $K \subseteq X$, the open sets $\{ K \cap W: W \in \calW \}$ 
form a good covering of $K$. We wish to prove that {\em every} covering $\calW$ of a locally compact topological space $X$ is good. In virtue of the preceding remarks, we can reduce to the case where $X$ is compact, and thereby assume that $\calW$ has a finite subcover.

We will prove, by induction on $n \geq 0$, that if $\calW$ is collection of open subsets of a locally compact Hausdorff space $X$ such that there exist $W_1, \ldots, W_n \in \calW$ with
$W_1 \cup \ldots \cup W_n = X$, then $\calW$ is a good covering of $X$. If $n = 0$, then
$X = \emptyset$. In this case, we must prove that $\calF(\emptyset)$ is final, which is part of the definition of $\calK$-sheaf.

Suppose that $\calW \subseteq \calW'$ are coverings of $X$, and that for every $W' \in \calW'$ the induced covering $\{ W \cap W': W \in \calW \}$ is a good covering of $W'$. It then follows from Proposition \ref{acekan} that $\calW'$ is a good covering of $X$ if and only if $\calW$ is a good covering of $X$.

Now suppose $n > 0$. Let $V = W_2 \cup \ldots \cup W_n$, and let $\calW' = \calW \cup \{V\}$. Using the above remark and the inductive hypothesis, it will suffice to show that $\calW'$ is a good covering of $X$. Now $\calW'$ contains a pair of open sets $W_1$ and $V$ which cover $X$. We thereby reduce to the case $n=2$; using the above remark we can furthermore suppose that
$\calW = \{ W_1, W_2 \}$. 

We now wish to show that for every compact $K \subseteq X$, $\calF$ exhibits
$\calF(K)$ as the limit of $\calF | \Nerve(\calK_{\calW}(X))^{op}$. Let $P$ be the collection
of all pairs $K_1, K_2 \in \calK(X)$ such that $K_1 \subseteq W_1$, $K_2 \subseteq W_2$, and $K_1 \cup K_2 = K$. We observe that $P$ is a filtered when ordered by inclusion.
For $\alpha = (K_1, K_2) \in P$, let $\calK_{\alpha} = \{ K' \in \calK(X): (K' \subseteq K_1) \vee (K' \subseteq K_2) \}$. We note that $\calK_{\calW}(X) = \bigcup_{\alpha \in P} \calK_{\alpha}$.
Moreover, Theorem \ref{hollowtt} implies that for $\alpha = (K_1, K_2) \in P$, the inclusion
$ \Nerve \{ K_1, K_2, K_1 \cap K_2 \} \subseteq \Nerve(\calK_{\alpha})$ is cofinal.
Since $\calF$ is a $\calK$-sheaf, we deduce that $\calF$ exhibits
$\calF(K)$ as a limit of the diagram $\calF | \Nerve(\calK_{\alpha})^{op}$ for each
$\alpha \in P$. Using Proposition \ref{extet}, we deduce that $\calF(K)$ is a limit of
$\calF | \Nerve(\calK_{\calW}(X))^{op}$ if and only if $\calF(K)$ is a limit of the constant
diagram $\Nerve(P)^{op} \rightarrow \SSet$ taking the value $\calF(K)$. This is clear, since
$P$ is filtered so that the map $\Nerve(P) \rightarrow \Delta^0$ is cofinal by Theorem \ref{hollowtt}.
\end{proof}

\begin{theorem}\label{kuku}
Let $X$ be a locally compact Hausdorff space and $\calC$ a presentable $\infty$-category in which filtered colimits are left exact.
Let $\calF: \Nerve (\calK(X) \cup \calU(X))^{op} \rightarrow \calC$ be a presheaf on the partially ordered set $\calK(X) \cup \calU(X)$.  
The following conditions are equivalent:
\begin{itemize}
\item[$(1)$] The presheaf $\calF_{\calK} = \calF| \Nerve(\calK(X))^{op}$ is a $\calK$-sheaf and $\calF$ is a right Kan extension of $\calF_{\calK}$.
\item[$(2)$] The presheaf $\calF_{\calU} = \calF | \Nerve(\calU(X))^{op}$ is a sheaf and $\calF$ is a left Kan extension of $\calF_{\calU}$.
\end{itemize}
\end{theorem}

\begin{proof}
Suppose first that $(1)$ is satisfied. We first prove that $\calF$ is a left Kan extension of $\calF_{\calU}$. Let $K$ be a compact subset of $X$, and let $\calU_{K \subseteq}(X) = \{ U \in \calU(X): K \subseteq U \}$. Consider the diagram
$$ \xymatrix{ \Nerve(\calU_{K \subseteq}(X))^{op} \ar[r]^-{p} \ar[d] & \Nerve (\calU_{K \subseteq}(X) \cup \calK_{K \Subset}(X))^{op} \ar[d] & \Nerve(\calK_{K \Subset}(X))^{op} \ar[d] \ar[l]_-{p'} \\
\Nerve(\calU_{K \subseteq}(X)^{op})^{\triangleright} \ar[r] \ar[ddr]^{\psi} & 
\Nerve (\calU_{K \subseteq}(X)) \cup \calK_{K \Subset}(X))^{op})^{\triangleright} \ar[d] & \Nerve (\calK_{K \Subset}^{op})^{\triangleright} \ar[ddl]_{\psi'} \ar[l] \\
& \Nerve (\calU(X) \cup \calK(X))^{op} \ar[d]^{\calF} & \\
& \calC. & }$$
We wish to prove that $\psi$ is a colimit diagram. Since $\calF_{\calK}$ is a $\calK$-sheaf, we deduce that $\psi'$ is a colimit diagram. It therefore suffices to check that $p$ and $p'$ are cofinal. 
According to Theorem \ref{hollowtt}, it suffices to show that for every $Y \in \calU_{K \subseteq}(X) \cup \calK_{K \Subset}(X)$, the partially ordered sets $\{ K' \in \calK(X): K \Subset K' \subseteq Y \}$ and $\{ U \in \calU(X): K \subseteq U \subseteq Y \}$ have contractible nerves. We now observe that both of these partially ordered sets is filtered, since they are nonempty and stable under finite unions.

We now show that $\calF_{\calU}$ is a sheaf. Let $U$ be an open subset of $X$ and let $\calW$ be a sieve which covers $U$. Let $\calK_{\subseteq U}(X) = \{ K \in \calK(X): K \subseteq U \}$ and let
$\calK_{\calW}(X) = \{ K \in \calK(X): (\exists W \in \calW) [ K \subseteq W ] \}$. 
We wish to prove that the diagram 
$$ \Nerve(\calW^{op})^{\triangleleft} \rightarrow \Nerve(\calU(X))^{op} \stackrel{\calF_{\calU}} \rightarrow \SSet$$ is a limit. Using Theorem \ref{hollowtt}, we deduce that the inclusion
$\Nerve(\calW) \subseteq \Nerve (\calW \cup \calK_{\calW}(X) )$ is cofinal. It therefore suffices to prove that $ \calF | (\calW \cup \calK_{\calW}(X) \cup \{U \})^{op}$ is a right Kan extension of
$\calF | (\calW \cup \calK_{\calW}(X))^{op}$. Since $\calF|(\calW \cup \calK_{\calW}(X))^{op}$ is a right Kan extension of $\calF|\calK_{\calW}(X)^{op}$ by assumption, it suffices to prove that $\calF| (\calW \cup \calK_{\calW}(X) \cup \{U\})^{op}$ is a right Kan extension of $\calF| \calK_{\calW}(X)^{op}$. This is clear at every object distinct from $U$; it will therefore suffice
to prove that $\calF| (\calK_{\calW}(X) \cup \{U\})^{op}$ is a right Kan extension of
$\calF| \calK_{\calW}(X)^{op}$.

By assumption, $\calF | \Nerve (\calK_{\subseteq U}(X) \cup \{U\})^{op}$ is a right Kan extension of
$\calF| \Nerve(\calK_{\subseteq U}(X))^{op}$ and Lemma \ref{noodlesoup} implies that
$\calF| \Nerve(\calK_{\subseteq U}(X))^{op}$ is a right Kan extension of $\calF| \Nerve (\calK_{\calW}(X))^{op}$. Invoking Proposition \ref{acekan}, we deduce that
$\calF| \Nerve (\calK_{\calW}(X) \cup \{U\})^{op}$ is a right Kan extension of 
$\calF| \Nerve (\calK_{\calW}(X))^{op}$. This shows that $\calF_{\calU}$ is a sheaf, and completes the proof that $(1) \Rightarrow (2)$.

Now suppose that $\calF$ satisfies $(2)$. We first verify that $\calF_{\calK}$ is a $\calK$-sheaf. 
The space $\calF_{\calK}(\emptyset) = \calF_{\calU}(\emptyset)$ is contractible because $\calF_{\calU}$ is a sheaf (and because the empty sieve is a covering sieve on $\emptyset \subseteq X$). Suppose next that $K$ and $K'$ are compact subsets of $X$. We wish to prove that the diagram 
$$ \xymatrix{ \calF( K \cup K' ) \ar[r] \ar[d] & \calF(K) \ar[d] \\
\calF(K') \ar[r] & \calF(K \cap K') }$$
is a pullback in $\SSet$. Let us denote this diagram by $\sigma: \Delta^1 \times \Delta^1 \rightarrow \SSet$. Let $P$ be the set of all pairs $U,U' \in \calU(X)$ such that $K \subseteq U$ and $K' \subseteq U'$. The functor $\calF$ induces a map $\sigma_{P}: \Nerve(P^{op})^{\triangleright} \rightarrow \SSet^{\Delta^1 \times \Delta^1}$, which carries each pair $(U,U')$ to the diagram
$$ \xymatrix{ \calF( U \cup U' ) \ar[r] \ar[d] & \calF(U) \ar[d] \\
\calF(U') \ar[r] & \calF(U \cap U') }$$
and carries the cone point to $\sigma$. Since $\calF_{U}$ is a sheaf, each $\sigma_P(U,U')$ is
a pullback diagram in $\calC$. Since filtered colimits in $\calC$ are left exact, it will suffice to show that $\sigma_P$ is a colimit diagram. By Proposition \ref{limiteval}, it suffices to show that each of the four maps
$$ \Nerve(P^{op})^{\triangleright} \rightarrow \SSet$$, given by evaluating $\sigma_P$ at the four vertices of $\Delta^1 \times \Delta^1$, is a colimit diagram. We will treat the case of the final vertex; the other cases are handled in the same way. Let $Q = \{ U \in \calU(X): K \cap K' \subseteq U\}$. T
We are given a map $g: \Nerve(P^{op})^{\triangleright} \rightarrow \SSet$ which admits a factorization
$$ \Nerve(P^{op})^{\triangleright} \stackrel{g''}{\rightarrow} \Nerve(Q^{op})^{\triangleright} \stackrel{g'}{\rightarrow}
\Nerve (\calU(X) \cup \calK(X))^{op} \stackrel{\calF}{\rightarrow} \calC.$$
Since $\calF$ is a left Kan extension of $\calF_{\calU}$, the diagram $\calF \circ g''$ is a colimit.
It therefore suffices to show that $g''$ induces a cofinal map $\Nerve(P)^{op} \rightarrow \Nerve (Q)^{op}$. Using Theorem \ref{hollowtt}, it suffices to prove that for every $U'' \in Q$, the
partially ordered set $P_{U''} = \{ (U,U') \in P: U \cap U' \subseteq U'' \}$ has contractible nerve. It now suffices to observe that $P^{op}_{U''}$ is filtered (since $P_{U''}$ is nonempty and stable under intersections). 

We next show that for any compact subset $K \subseteq X$, the map
$$ \Nerve(\calK_{K \Subset}(X)^{op})^{\triangleright} \rightarrow \Nerve (\calK(X) \cup \calU(X))^{op} \stackrel{\calF}{\rightarrow} \calC$$
is a colimit diagram. Let $\calV = \calU(X) \cup \calK_{K \Subset}(X)$, and let $\calV' = \calV \cup \{K\}$. It follows from Proposition \ref{acekan} that $\calF | \Nerve(\calV)^{op}$ and
$\calF | \Nerve(\calV')^{op}$ are left Kan extensions of $\calF | \Nerve (\calU(X))^{op}$, so
that $\calF | \Nerve (\calV')^{op}$ is a left Kan extension of $\calF| \Nerve(\calV)^{op}$. Therefore the diagram
$$ ( \Nerve ( \calK_{K \Subset}(X) \cup \{ U \in \calU(X): K \subseteq U \} )^{op})^{\triangleright}
\rightarrow \Nerve( \calK(X) \cup \calU(X))^{op} \stackrel{\calF}{\rightarrow} \calC $$
is a colimit. It therefore suffices to show that the inclusion
$$ \Nerve (\calK_{K \Subset}(X))^{op} \subseteq \Nerve (\calK_{K \Subset}(X) \cup \{ U \in \calU(X): K\subseteq U\} )^{op}$$
is cofinal. Using Theorem \ref{hollowtt}, we are reduced to showing that if 
$Y \in \calK_{K \Subset}(X) \cup \{ U \in \calU(X): K \subseteq U\}$, then the nerve of the partially ordered set $R = \{ K' \in \calK(X): K \Subset K' \subset Y \}$ is weakly contractible. It now suffices to observe that $R^{op}$ is filtered, since $R$ is nonempty and stable under intersections. This completes the proof that $\calF_{\calK}$ is a $\calK$-sheaf.

We now show that $\calF$ is a right Kan extension of $\calF_{\calK}$. Let $U$ be an open subset of $X$, and for $V \in \calU(X)$ write $V \Subset U$ if the closure $\overline{V}$ is compact and contained in $U$. Let $\calU_{\Subset U}(X) = \{ V \in \calU(X): V \Subset U\}$,
and consider the diagram
$$ \xymatrix{ \Nerve (\calU_{\Subset U}(X))^{op} \ar[r]^-{f} \ar[d] & \Nerve ( \calU_{\Subset U}(X) \cup \calK_{\subseteq U}(X))^{op} \ar[d] & \Nerve (\calK_{\subseteq U}(X))^{op} \ar[d] \ar[l]_-{f'} \\
\Nerve (\calU_{\Subset U}(X)^{op})^{\triangleleft} \ar[ddr]^{\phi} & 
\Nerve ( \calU_{\Subset U}(X) \cup \calK_{\subseteq U}(X))^{op})^{\triangleleft} \ar[d] & 
\Nerve (\calK_{\subseteq U}(X)^{op})^{\triangleleft} \ar[ddl]_{\phi'} \\
& \Nerve (\calK(X) \cup \calU(X))^{op} \ar[d]^{\calF} & \\
& \calC. & }$$
We wish to prove that $\phi'$ is a limit diagram. Since the sieve $\calU_{\Subset U}(X)$ covers $U$ and $\calF_{\calU}$ is a sheaf, we conclude that $\phi$ is a limit diagram. It therefore suffices to prove that $f^{op}$ and $(f')^{op}$ are cofinal maps of simplicial sets. According to Theorem \ref{hollowtt}, it suffices to prove that if $Y \in \calK_{\subseteq U}(X) \cup \calU_{\Subset U}(X)$, then the partially ordered sets $\{ V \in \calU(X): Y \subseteq V \Subset U \}$ and
$\{ K \in \calK(X): Y \subseteq K \subseteq U \}$ have weakly contractible nerves. We now observe that both of these partially ordered sets are filtered (since they are nonempty and stable under unions). This completes the proof that $\calF$ is a right Kan extension of $\calF_{\calK}$.
\end{proof}

\begin{corollary}\label{streeem}
Let $X$ be a locally compact topological space and $\calC$ a presentable $\infty$-category in which filtered colimits are left exact. Let
$$\Shv_{\calKU}(X;\calC) \subseteq \Fun(\Nerve (\calK(X) \cup \calU(X))^{op}, \calC)$$
be the full subcategory spanned by those presheaves which satisfy the equivalent conditions of Theorem \ref{kuku}. Then the restriction functors $$ \Shv(X;\calC) \leftarrow \Shv_{\calKU}(X;\calC) \rightarrow \Shv_{\calK}(X;\calC)$$\index{not}{ShvcalKUX@$\Shv_{\calKU}(X;\calC)$}
are equivalences of $\infty$-categories.
\end{corollary}

\begin{corollary}\label{compactprop}
Let $X$ be a compact Hausdorff space. Then the global sections functor
$\Gamma: \Shv(X) \rightarrow \SSet$ is a proper morphism of $\infty$-topoi.
\end{corollary}

\begin{proof}
The existence of fiber products $\Shv(X) \times_{\SSet} \calY$ in
$\RGeom$ follows from Theorem \ref{conprod}. It will therefore suffice to prove that for any (homotopy) Cartesian rectangle
$$ \xymatrix{ \calX'' \ar[r] \ar[d] & \calX' \ar[r] \ar[d] & \Shv(X) \ar[d] \\
\calY'' \ar[r]^{f_{\ast}} & \calY' \ar[r] & \SSet,}$$
the square on the left is left adjointable. Using Theorem \ref{conprod}, we can
identify the square on the left with
$$ \xymatrix{ \Shv(X; \calY'') \ar[r] \ar[d] & \Shv(X; \calY') \ar[d] \\
\calY'' \ar[r]^{f_{\ast}} & \calY', }$$ where the vertical morphisms are given by taking global sections.

Choose a correspondence $\calM$ from $\calY'$ to $\calY''$ which is associated to the functor $f_{\ast}$. Since $f_{\ast}$ admits a left adjoint $f^{\ast}$, the projection $\calM \rightarrow \Delta^1$ is both a Cartesian fibration and a coCartesian fibration. For every simplicial
set $K$, let $\calM_{K} = \Fun(K,\calM) \times_{ \Fun(K,\Delta^1) } \Delta^1$. Then $\calM_{K}$ determines a correspondence from $\Fun(K,\calY')$ to $\Fun(K,\calY'')$. Using Proposition \ref{doog}, we
conclude that $\calM_{K} \rightarrow \Delta^1$ is both a Cartesian and a coCartesian fibration,
and that it is associated to the functors given by composition with $f_{\ast}$ and $f^{\ast}$.

Before proceeding further, let us adopt the following convention for the remainder of the proof:
given a simplicial set $Z$ with a map $q: Z \rightarrow \Delta^1$, we will say that an edge
of $Z$ is {\it Cartesian} or {\it coCartesian} if it is $q$-Cartesian or $q$-coCartesian, respectively.
The map $q$ to which we are referring should be clear from context.

Let $\calM_{\calU}$ denote the full subcategory of $\calM_{ \Nerve(\calU(X))^{op}}$ whose objects
correspond to {\em sheaves} on $X$ (with values in either $\calY'$ or $\calY''$). Since
$f_{\ast}$ preserves limits, composition with $f_{\ast}$ carries $\Shv(X; \calY'')$ into
$\Shv(X; \calY')$. We conclude that the projection $\calM_{\calU} \rightarrow \Delta^1$
is a Cartesian fibration, and that the inclusion $\calM_{\calU} \subseteq \calM_{ \Nerve (\calU(X))^{op}}$
preserves Cartesian edges.

Similarly, we define $\calM_{\calK}$ to be the full subcategory of $\calM_{ \Nerve(\calK(X))^{op} }$
whose objects correspond to $\calK$-sheaves on $X$ (with values in either $\calY'$ or
$\calY''$). Since $f^{\ast}$ preserves finite limits and filtered colimits, composition with $f^{\ast}$ carries $\Shv_{\calK}(X; \calY')$ into $\Shv_{\calK}(X; \calY'')$. It follows that the projection
$\calM_{\calK} \rightarrow \Delta^1$ is a coCartesian fibration, and that the inclusion
$\calM_{\calK} \subseteq \calM_{ \Nerve(\calU(X))^{op}}$ preserves coCartesian edges.

Now let $\calM'_{\calKU} = \calM_{ \Nerve (\calK(X) \cup \calU(X))^{op}}$ and let
$\calM_{\calKU}$ be the full subcategory of $\calM'_{\calKU}$
spanned by the objects of $\Shv_{\calKU}(X; \calY')$ and $\Shv_{\calKU}(X; \calY'')$. 
We have a commutative diagram
$$ \xymatrix{ & \calM_{\calKU} \ar[dr]^{\phi_{\calU}} \ar[dl]^{\phi_{\calK}} & \\
\calM_{\calU} \ar[dr]^{\Gamma_{\calU}} & & \calM_{\calK} \ar[dl]^{\Gamma_{\calK}} \\
& \calM & }$$
where $\Gamma_{\calU}$ and $\Gamma_{\calK}$ denote the global sections functors (given by evaluation at $X \in \calU(X) \cap \calK(X)$). According to Remark \ref{toadcatcher}, to complete the proof it will suffice to show that $\calM_{\calU} \rightarrow \Delta^1$ is a coCartesian fibration, and that $\Gamma_{\calU}$ preserves both Cartesian and coCartesian edges. It is clear that
$\Gamma_{\calU}$ preserves Cartesian edges, since it is a composition of maps
$$ \calM_{\calU} \subseteq \calM_{ \Nerve( \calU(X))^{op} } \rightarrow \calM$$
which preserve Cartesian edges. Similarly, we already know that $\calM_{\calK} \rightarrow \Delta^1$ is a coCartesian fibration, and that $\Gamma_{\calK}$ preserves coCartesian edges.
To complete the proof, it will therefore suffice to show that $\phi_{\calU}$ and $\phi_{\calK}$ are
equivalences of $\infty$-categories. We will give the argument for $\phi_{\calU}$; the proof in the case of $\phi_{\calK}$ is identical and left to the reader.

According to Corollary \ref{streeem}, the map $\phi_{\calU}$ induces equivalences
$$ \Shv_{\calKU}(X; \calY') \rightarrow \Shv(X; \calY')$$
$$ \Shv_{\calKU}(X; \calY'') \rightarrow \Shv(X; \calY'')$$ 
after passing to the fibers over either vertex of $\Delta^1$.
We will complete the proof by applying Corollary \ref{usefir}. In order to do so, we must verify that
$p: \calM_{\calKU} \rightarrow \Delta^1$ is a Cartesian fibration, and that $\phi_{\calU}$ preserves Cartesian edges.

To show that $p$ is a Cartesian fibration, we begin with an arbitrary
$\calF \in \Shv_{\calKU}(X; \calY'')$. Using Proposition \ref{doog}, we conclude the existence of a $p'$-Cartesian morphism $\alpha: \calF' \rightarrow \calF$, where $p'$ denotes the projection
$\calM'_{\calKU}$ and $\calF' = \calF \circ p_{\ast} \in
\Fun( \Nerve(\calK(X) \cup \calU(X))^{op}, \calY')$. Since $p_{\ast}$ preserves limits, we conclude that
$\calF' | \Nerve(\calU(X))^{op}$ is a sheaf on $\calX$ with values in $\calY'$; however,
$\calF'$ is not necessarily a left Kan extension of $\calF' | \Nerve(\calU(X))^{op}$. 
Let $\calC$ denote the full subcategory of $\Fun(\Nerve (\calK(X) \cup \calU(X))^{op}, \calY')$
spanned by those functors $\calG: \Nerve (\calK(X) \cup \calU(X))^{op}$ which are left Kan extensions of $\calG | \Nerve(\calU(X))^{op}$, and $s$ a section of the trivial fibration
$\calC \rightarrow (\calY')^{ \Nerve(\calU(X))^{op}}$, so that $s$ is a left adjoint to the
restriction map $r: \calM'_{\calKU} \rightarrow (\calY')^{ \Nerve(\calU(X))^{op} }$. 
Let $\calF'' = (s \circ r) \calF'$ be a left Kan extension of $\calF' | \Nerve(\calU(X))^{op}$. Then
$\calF''$ is an initial object of the fiber $\calM'_{\calKU} \times_{ \Fun(\Nerve(\calU(X))^{op}, \calY')}
\{ \calF' | \Nerve(\calU(X))^{op} \}$, so that there exists a map $\beta: \calF'' \rightarrow \calF'$
which induces the identity on $\calF'' | \Nerve(\calU(X))^{op} = \calF' | \Nerve(\calU(X))^{op}$. 

Let
$\sigma: \Delta^2 \rightarrow \calM'_{\calKU}$ classify
a diagram
$$ \xymatrix{ & \calF' \ar[dr]^{\alpha} & \\
\calF'' \ar[ur]^{\beta} \ar[rr]^{\gamma} & & \calF, }$$
so that $\gamma$ is a composition of $\alpha$ and $\beta$. It is easy
to see that $\phi_{\calU}(\gamma)$ is a Cartesian edge of $\calM_{\calU}$ (since it is a composition of a Cartesian edge with an equivalence in $\Shv(X; \calY')$). We claim that $\gamma$ is $p$-Cartesian. To prove this, consider the diagram
$$ \xymatrix{ \Shv_{\calKU}(X; \calY') \times_{ \calM'_{\calKU}} (\calM'_{\calKU})_{/\sigma} \ar[d]^{\theta_0} \ar[r]^-{\eta'} & (\calM_{\calKU})_{\gamma} \ar[dd]^{\eta} \\
(\Shv_{\calKU}(X; \calY'))_{/\beta}
\times_{ \Shv_{\calKU}(X; \calY')_{/ \calF'}} (\calM'_{\calKU})_{/\alpha}
\ar[d]^{\theta_1} & \\ 
\Shv_{\calKU}(X; \calY') \times_{ \calM'_{\calKU} } (\calM'_{\calKU})_{/\alpha} \ar[r]^-{\theta_2} &
\Shv_{\calKU}(X; \calY') \times_{ \calM_{\calKU} }
(\calM_{\calKU})_{/\calF}. }$$
We wish to show that $\eta$ is a trivial fibration. Since $\eta$ is a right fibration, it suffices
to show that the fibers of $\eta$ are contractible. The map $\eta'$ is a trivial fibration
(since the inclusion $\Delta^{ \{0,2\} } \subseteq \Delta^2$ is right anodyne), so it will suffice to prove that $\eta \circ \eta'$ is a trivial fibration. In view of the commutativity of the diagram, it will suffice to show that $\theta_0$, $\theta_1$, and $\theta_2$ are trivial fibrations. The triviality of $\theta_0$ follows from the fact that the horn inclusion $\Lambda^2_1 \subseteq \Delta^2$ is right anodyne.
The triviality of $\theta_2$ follows from the fact that $\alpha$ is $p'$-Cartesian. Finally, we observe that $\theta_1$ is a pullback of the map $\theta'_1: \Shv_{\calKU}(X; \calY')_{/\beta} \rightarrow
\Shv_{\calKU}(X; \calY')_{/\calF'}$. Let $\calC = (\calY')^{ \Nerve (\calK(X) \cup \calU(X))^{op} }$.
To prove that $\theta'_1$ is a trivial fibration, we must show that
for every $\calG \in \Shv_{\calKU}$, composition with $\beta$ induces a homotopy equivalence
$$ \bHom_{ \calC}( \calG, \calF'' ) \rightarrow \bHom_{\calC}(\calG, \calF').$$
Without loss of generality, we may suppose that $\calG = s( \calG')$, where
$\calG' \in \Shv(X; \calY')$; now we simply invoke the adjointness of $s$ with the
restriction functor $r$ and the observation that $r(\beta)$ is an equivalence.
\end{proof}

\begin{corollary}\label{compactcase}
Let $X$ be a compact Hausdorff space.
The global sections functor $\Gamma: \Shv(X) \rightarrow \SSet$ preserves filtered colimits.
\end{corollary}

\begin{proof}
Applying Theorem \ref{kuku}, we can replace $\Shv(X)$ by $\Shv_{\calK}(X)$. Now observe that the full subcategory $\Shv_{\calK}(X) \subseteq \calP( \Nerve(\calK(X))^{op})$ is stable under filtered colimits. We thereby reduce to proving that the evaluation functor $\calP( \Nerve(\calK(X))^{op} ) \rightarrow \SSet$ commutes with filtered colimits, which follows from Proposition \ref{limiteval}.
Alternatively, one can apply Corollary \ref{streeem} and Remark \ref{swurk}.
\end{proof}

\begin{remark}
One can also deduce Corollary \ref{compactcase} using the geometric model for $\Shv(X)$ introduced in \S \ref{paracompactness}. Using the characterization of properness in terms of filtered colimits described in Remark \ref{swurk}, one can formally deduce Corollary \ref{compactprop} from
Corollary \ref{compactcase}. This leads to another proof of the proper base change theorem, which does not make use of Theorem \ref{kuku} or the other ideas of this section. However, this alternative proof is considerably more difficult than the one described here, since it requires a rigorous justification of Remark \ref{swurk}. We also note that Theorem \ref{kuku} and Corollary \ref{streeem} are interesting in their own right, and could conceivably be applied in other contexts.
\end{remark}

\subsection{Sheaves on Coherent Spaces}\label{cohthm}\index{gen}{topological space!coherent}\index{gen}{coherent topological space}

Theorem \ref{kuku} has an analogue in the setting of coherent topological spaces which is somewhat easier to prove. First, we need the analogue of Lemma \ref{noodlesoup}:

\begin{lemma}\label{oldtime}\index{not}{Ucal0X@$\calU_0(X)$}
Let $X$ be a coherent topological space, let $\calU_0(X)$ denote the collection
of compact open subsets of $X$, and let $\calF: \Nerve(\calU_0(X))^{op} \rightarrow \calC$
be a presheaf taking values in an $\infty$-category $\calC$, having the following properties:
\begin{itemize}
\item[$(1)$] The object $\calF(\emptyset) \in \calC$ is final.
\item[$(2)$] For every pair of compact open sets $U, V \subseteq X$, the diagram
$$ \xymatrix{ \calF( U \cap V) \ar[r] \ar[d] & \calF(U) \ar[d] \\
\calF(V) \ar[r] & \calF(U \cup V) }$$
is a pullback.
\end{itemize}
Let $\calW$ be a covering of $X$ by compact open subsets, and let $\calU_{1}(X) \subseteq \calU_{0}(X)$ be collection of all compact open subsets of $X$ which are contained in some element of $\calW$. Then $\calF$ is a right Kan extension of $\calF | \Nerve(\calU_1(X))^{op}$. 
\end{lemma}

\begin{proof}
The proof is similar to that of Lemma \ref{noodlesoup}, but slightly easier.
Let us say that a covering $\calW$ of a coherent topological space $X$ by compact open subsets
is {\it good} if it satisfies the conclusions of the Lemma. We observe that $\calW$ automatically has a finite subcover. We will prove, by induction on $n \geq 0$, that if $\calW$ is collection of open subsets of a locally coherent topological space $X$ such that there exist $W_1, \ldots, W_n \in \calW$ with
$W_1 \cup \ldots \cup W_n = X$, then $\calW$ is a good covering of $X$. If $n = 0$, then
$X = \emptyset$. In this case, we must prove that 
$\calF(\emptyset)$ is final, which is one of our assumptions.

Suppose that $\calW \subseteq \calW'$ are coverings of $X$ by compact open sets, and that for every $W' \in \calW'$ the induced covering $\{ W \cap W': W \in \calW \}$ is a good covering of $W'$. It then follows from Proposition \ref{acekan} that $\calW'$ is a good covering of $X$ if and only if $\calW$ is a good covering of $X$.

Now suppose $n > 0$. Let $V = W_2 \cup \ldots \cup W_n$, and let $\calW' = \calW \cup \{V\}$. Using the above remark and the inductive hypothesis, it will suffice to show that $\calW'$ is a good covering of $X$. Now $\calW'$ contains a pair of open sets $W_1$ and $V$ which cover $X$. 
We thereby reduce to the case $n=2$; using the above remark we can furthermore suppose that
$\calW = \{ W_1, W_2 \}$. 

We now wish to show that for every compact $U \subseteq X$, $\calF$ exhibits
$\calF(U)$ as the limit of $\calF | \Nerve(\calU_1(X)_{/U})^{op}$. Without loss of generality,
we may replace $X$ by $U$ and thereby reduce to the case $U=X$. Let
$\calU_2(X) = \{ W_1, W_2, W_1 \cap W_2 \} \subseteq \calU_1(X)$. Using Theorem
\ref{hollowtt}, we deduce that the inclusion $\Nerve(\calU_2(X)) \subseteq \Nerve(\calU_1(X))$ is cofinal. Consequently, it suffices to prove that $\calF(X)$ is the limit of the diagram
$\calF | \Nerve(\calU_2(X) )^{op}$. In other words, we must show that the diagram
$$ \xymatrix{ \calF(X) \ar[r] \ar[d] & \calF(W_1) \ar[d] \\
\calF(W_2) \ar[r] & \calF(W_1 \cap W_2) }$$
is a pullback in $\calC$, which is true by assumption.
\end{proof}

\begin{theorem}\label{surm}
Let $X$ be a coherent topological space, and let $\calU_0(X) \subseteq \calU(X)$ denote the collection of compact open subsets of $X$. Let $\calC$ be an $\infty$-category which admits small limits.
The restriction map 
$$ \Shv(\calX; \calC) \rightarrow \Fun( \Nerve (\calU_0(X))^{op}, \calC)$$
is fully faithful, and its essential image consists of precisely those functors
$\calF_0: \Nerve( \calU_0(X))^{op} \rightarrow \calC$ satisfying the following conditions:
\begin{itemize}
\item[$(1)$] The object $\calF_0(\emptyset) \in \calC$ is final.
\item[$(2)$] For every pair of compact open sets $U, V \subseteq X$, the diagram
$$ \xymatrix{ \calF_0( U \cap V) \ar[r] \ar[d] & \calF_0(U) \ar[d] \\
\calF_0(V) \ar[r] & \calF_0(U \cup V) }$$
is a pullback.
\end{itemize}
\end{theorem}

\begin{proof}
Let $\calD \subseteq \calC^{ \Nerve(\calU(X))^{op} }$ be the full subcategory spanned by those presheaves $\calF: \Nerve( \calU(X))^{op} \rightarrow \calC$ which are right Kan extensions of
$\calF_0 = \calF | \Nerve(\calU_0(X))^{op}$, and such that $\calF_0$ satisfies conditions $(1)$ and $(2)$. According to Proposition \ref{lklk}, it will suffice to show that $\calD$ coincides with
$\Shv(X; \calC)$. 

Suppose that $\calF: \Nerve(\calU(X))^{op} \rightarrow \calC$ is a sheaf. We first show that $\calF$ is a right Kan extension of $\calF_0 = \calF| \Nerve(\calU_0(X))^{op}$. Let $U$ be an open subset of $X$, let $\calU(X)_{/U}^{(0)}$ denote the collection of compact open subsets of $U$, and let
$\calU(X)_{/U}^{(1)}$ denote the sieve generated by $\calU(X)_{/U}^{(0)}$. Consider the diagram

$$\xymatrix{ \Nerve(\calU(X)_{/U}^{(0)})^{\triangleright} \ar[drrr]^{f} \ar[r]^{i} &
\Nerve(\calU(X)_{/U}^{(1)})^{\triangleright} \ar[drr]^{f'} \ar[r] & 
\Nerve(\calU(X)_{/U})^{\triangleright} \ar[dr] \ar[r] & 
\Nerve(\calU(X)) \ar[d]^{\calF} \\
& & & \calC^{op}. }$$
We wish to prove that $f$ is a colimit diagram. Using Theorem \ref{hollowtt}, we deduce that the inclusion $\Nerve(\calU(X))^{(0)}_{/U} \subseteq \Nerve(\calU(X))^{(1)}_{/U}$ is cofinal.
It therefore suffices to prove that $f'$ is a colimit diagram. Since $\calF$ is a sheaf, it suffices to prove that $\calU(X)_{/U}^{(1)}$ is a covering sieve. In other words, we need to prove that $U$ is
a union of compact open subsets of $X$, which follows immediately from our assumption that
$X$ is coherent.

We next prove that $\calF_0$ satisfies $(1)$ and $(2)$. To prove $(1)$, we simply observe that
the empty sieve is a cover of $\emptyset$ and apply the sheaf condition. To prove $(2)$, we may assume without loss of generality that neither $U$ nor $V$ is contained in the other (otherwise the result is obvious). 
Let $\calU(X)^{(0)}_{/U \cup V}$ be the full subcategory spanned by $U$, $V$, and $U \cap V$, and let $\calU(X)^{(1)}_{/U \cup V}$ be the sieve on $U \cup V$ generated by $\calU(X)^{(0)}_{/U \cup V}$. As above, we have a diagram 
$$\xymatrix{ (\Nerve(\calU(X))_{/U \cup V}^{(0)})^{\triangleright} \ar[drrr]^{f} \ar[r]^{i} &
\Nerve(\calU(X)_{/U \cup V}^{(1)})^{\triangleright} \ar[drr]^{f'} \ar[r] & 
\Nerve(\calU(X)_{/U \cup V})^{\triangleright} \ar[dr] \ar[r] & 
\Nerve(\calU(X)) \ar[d]^{\calF} \\
& & & \calC^{op}, }$$
and we wish to show that $f$ is a colimit diagram. Theorem \ref{hollowtt} implies that
the inclusion $\Nerve(\calU(X))_{/U \cup V}^{(0)} \subseteq \Nerve(\calU(X))_{/U \cup V}^{(1)}$ is cofinal. It therefore suffices to prove that $f'$ is a colimit diagram, which follows from the sheaf condition since $\calU(X)^{(1)}_{/U \cup V}$ is a covering sieve. This completes the proof 
that $\Shv(X;\calC) \subseteq \calD$.

It remains to prove that $\calD \subseteq \Shv(X;\calC)$. In other words, we must show that
if $\calF$ is a right Kan extension of $\calF_0 = \calF | \Nerve(\calU_0(X))^{op}$, and
$\calF_0$ satisfies conditions $(1)$ and $(2)$, then $\calF$ is a sheaf. Let $U$ be an open subset of $X$, $\calU(X)_{/U}^{(0)}$ a sieve which covers $U$. Let $\calU_0(X)_{/U}$ denote the category of compact open subsets of $U$ and $\calU_0(X)^{(0)}_{/U}$ the category of compact
open subsets of $U$ which belong to the sieve $\calU(X)_{/U}^{(0)}$. We wish to prove
that $\calF(U)$ is a limit of $\calF | \Nerve(\calU(X)_{/U}^{(0)})^{op}$. We will in fact prove the slightly stronger assertion that $\calF | \Nerve (\calU(X)_{/U})^{op}$ is a right Kan extension of 
$\calF | \Nerve( \calU(X)_{/U}^{(0)})^{op}$.

We have a commutative diagram
$$ \xymatrix{ \calU_0(X)^{(0)}_{/U} \ar[r] \ar[d] & \calU_0(X)_{/U} \ar[d] \\
\calU(X)_{/U}^{(0)} \ar[r] & \calU(X)_{/U}. }$$
By assumption, $\calF$ is a right Kan extension of $\calF_0$. It follows that
$\calF | \Nerve (\calU(X)_{/U}^{(0)})^{op}$ is a right Kan extension of
$\calF | \Nerve (\calU_0(X)_{/U}^{(0)})^{op}$ and that
$\calF | \Nerve (\calU(X)_{/U})^{op}$ is a right Kan extension of
$\calF | \Nerve (\calU_0(X)_{/U})^{op}$. By the transitivity of Kan extensions
(Proposition \ref{acekan}), it will suffice to prove that $\calF | \Nerve(\calU_0(X)_{/U})^{op}$ is a right
Kan extension of $\calF | \Nerve (\calU_0(X)^{(0)}_{/U})^{op}$. This follows immediately
from Lemma \ref{oldtime}.
\end{proof}

\begin{corollary}\label{applyZR}
Let $X$ be a coherent topological space. Then the global sections functor
$\Gamma: \Shv(X) \rightarrow \SSet$ is a proper map of $\infty$-topoi.
\end{corollary}

\begin{proof}
Identical to the proof of Corollary \ref{compactprop}, using
Theorem \ref{surm} in place of Corollary \ref{streeem}.
\end{proof} 

\begin{corollary}
Let $X$ be a coherent topological space. Then the global sections functor
$$ \Gamma: \Shv(X) \rightarrow \SSet$$ commutes with filtered colimits.
\end{corollary}

\subsection{Cell-Like Maps}\label{celluj}

Recall that a topological space $X$ is an {\it absolute neighborhood retract} if $X$ is\index{gen}{topological space!absolute neighborhood retract}\index{gen}{absolute neighborhood retract}
metrizable and if for any closed immersion $X \hookrightarrow Y$ of $X$ in a metric space $Y$, there exists an open set $U \subseteq Y$ containing the image of $X$, such that the inclusion
$X \hookrightarrow U$ has a left inverse (in other words, $X$ is a {\em retract} of $U$).

Let $p: X \rightarrow Y$ be a continuous map between locally compact absolute neighborhood retracts. The map $p$ is said to be {\it cell-like} if $p$ is proper and each fiber $X_{y} = X \times_{Y} \{y\}$ has trivial shape (in the sense of Borsuk; see \cite{shapetheory} and \S \ref{shapesec}). The theory of cell-like maps plays an important role in geometric topology: we refer the reader to \cite{cellmap} for a discussion (and for several equivalent formulations of the condition that a map be cell-like).\index{gen}{cell-like!map of topological spaces}

The purpose of this section is to describe a class of geometric morphisms between $\infty$-topoi, which we will call {\it cell-like} morphisms. We will then compare our theory of cell-like morphisms with the classical theory of cell-like maps. We will also give a ``nonclassical'' example which arises in the theory of rigid analytic geometry.

\begin{definition}\label{celldef}\index{gen}{cell-like!map of $\infty$-topoi}
Let $p_{\ast}: \calX \rightarrow \calY$ be a geometric morphism of $\infty$-topoi. We will say that
$p_{\ast}$ is {\it cell-like} if it is proper and if the right adjoint $p^{\ast}$ (which is well-defined up to equivalence) is fully faithful.
\end{definition}

\begin{warning}
Many authors refer to a map $p: X \rightarrow Y$ of {\em arbitrary} compact metric spaces
as {\it cell-like} if each fiber $X_{y} = X \times_{Y} \{y\}$ has trivial shape. This condition
is generally {\em weaker} than the condition that $p_{\ast}: \Shv(X) \rightarrow \Shv(Y)$ be cell-like in the sense of Definition \ref{celldef}. However, the two definitions are equivalent provided that $X$ and $Y$ are sufficiently nice (for example, if they are locally compact absolute neighborhood retracts). Our departure from the classical terminology is perhaps justified by the 
fact that the class of morphisms introduced in Definition \ref{celldef} has good formal properties: for example, stability under composition.
\end{warning}

\begin{remark}\label{urjk}
Let $p_{\ast}: \calX \rightarrow \calY$ be a cell-like geometric morphism between $\infty$-topoi.
Then the unit map $\id_{\calY} \rightarrow p_{\ast} p^{\ast}$ is an equivalence of functors. It follows immediately that $p_{\ast}$ induces an equivalence of shapes $\Sh(\calX) \rightarrow \Sh(\calY)$ (see \S \ref{shapesec}.
\end{remark}

\begin{proposition}\label{fibersilly}
Let $p_{\ast}: \calX \rightarrow \calY$ be a proper morphism of $\infty$-topoi. Suppose that
$\calY$ has enough points. Then $p_{\ast}$ is cell-like if and only if, for every pullback diagram
$$ \xymatrix{ \calX' \ar[r] \ar[d] & \calX \ar[d]^{p_{\ast}} \\
\SSet \ar[r] & \calY }$$
in $\RGeom$, the $\infty$-topos $\calX'$ has trivial shape.
\end{proposition}

\begin{proof}
Suppose first that each fiber $\calX'$ has trivial shape. Let $\calF \in \calY$. We wish to show that the unit map $u: \calF \rightarrow p_{\ast} p^{\ast} \calF$ is an equivalence. Since $\calY$ has enough points, it suffices to show that for each point $q_{\ast}: \SSet \rightarrow \calY$, the map $q^{\ast} u$ is an equivalence in $\SSet$, where $q^{\ast}$ denotes a left adjoint to $q_{\ast}$.
Form a pullback diagram of $\infty$-topoi
$$ \xymatrix{ \calX' \ar[r] \ar[d]^{s_{\ast}} & \calX \ar[d]^{p_{\ast}} \\
\SSet \ar[r]^{q_{\ast}} & \calY. }$$
Since $p_{\ast}$ is proper, this diagram is left-adjointable. Consequently,
$q^{\ast} u$ can be identified with the unit map
$$ K \rightarrow s_{\ast} s^{\ast} K,$$
where $K = q^{\ast} \calF \in \SSet$. If $\calX'$ has trivial shape, then this map is an equivalence.

Conversely, if $p_{\ast}$ is cell-like, then the above argument shows that for every diagram
$$ \xymatrix{ \calX' \ar[r] \ar[d]^{s_{\ast}} & \calX \ar[d]^{p_{\ast}} \\
\SSet \ar[r]^{q_{\ast}} & \calY }$$
as above and every $\calF \in \calY$, the adjunction map
$$ K \rightarrow s_{\ast} s^{\ast} K$$ is an equivalence, where
$K = q^{\ast} \calF$. To prove that $\calX'$ has trivial shape, it will suffice to show that $q^{\ast}$ is essentially surjective. For this, we observe that since $\SSet$ is a final object in the $\infty$-category of $\infty$-topoi, there exists a geometric morphism
$r_{\ast}: \calY \rightarrow \SSet$ such that $r_{\ast} \circ q_{\ast}$ is homotopic to $\id_{\SSet}$.
It follows that $q^{\ast} \circ r^{\ast} \simeq \id_{\SSet}$. Since $\id_{\SSet}$ is essentially surjective, we conclude that $q^{\ast}$ is essentially surjective.
\end{proof}

\begin{corollary}\label{comb0}
Let $p: X \rightarrow Y$ be a map of paracompact topological spaces. Assume that $p_{\ast}$ is proper, and that $Y$ has finite covering dimension. Then $p_{\ast}: \Shv(X) \rightarrow \Shv(Y)$ is cell-like if and only if each fiber $X_y = X \times_{Y} \{y\}$ has trivial shape.
\end{corollary}

\begin{proof}
Combine Proposition \ref{fibersilly} with Corollary \ref{enuff}.
\end{proof}

\begin{proposition}\label{cell}
Let $p: X \rightarrow Y$ be a proper map of locally compact ANRs. The following conditions are equivalent:

\begin{itemize}
\item[$(1)$] The geometric morphism $p_{\ast}: \Shv(X) \rightarrow \Shv(Y)$ is cell-like.
\item[$(2)$] For every open subset $U \subseteq Y$, the restriction map
$X \times_{Y} U \rightarrow U$ is a homotopy equivalence.
\item[$(3)$] Each fiber $X_{y} = X \times_{Y} \{y\}$ has trivial shape.
\end{itemize}
\end{proposition}

\begin{proof}
It is easy to see that if $p_{\ast}$ is cell-like, then each of the restrictions $p': X \times_{Y} U \rightarrow U$ induces a cell-like geometric morphism. According to Remark \ref{urjk}, 
$p'_{\ast}$ is a shape equivalence, and therefore a homotopy equivalence by Proposition \ref{parashape}. Thus $(1) \Rightarrow (2)$.

We next prove that $(2) \Rightarrow (1)$. 
Let $\calF \in \Shv(Y)$, and let $u: \calF \rightarrow p_{\ast} p^{\ast} \calF$ be a unit map; we wish to show that $u$ is an equivalence. It will suffice to show that the induced map
$\calF(U) \rightarrow (p_{\ast} p^{\ast} \calF)(U)$ is an equivalence in $\SSet$ for each paracompact open subset $U \subseteq Y$. Replacing $Y$ by $u$, we may reduce to the problem of showing that the map $\calF(Y) \rightarrow (p^{\ast} \calF)(X)$ is a homotopy equivalence. According
to Corollary \ref{wamain}, we may assume that $\calF$ is the simplicial nerve of $\Sing_{Y} \widetilde{Y}$, where $\widetilde{Y}$ is a fibrant-cofibrant object of $\Top_{/Y}$. According to
Proposition \ref{basechang}, we may identify $p^{\ast} \calF$ with $\Sing_{X} \widetilde{X}$, where
$\widetilde{X} = X \times_{Y} \widetilde{Y}$. It therefore suffices to prove that the induced map of simplicial function spaces
$$ \bHom_{Y}( Y, \widetilde{Y}) \rightarrow \bHom_{X}(X, \widetilde{X}) \simeq \bHom_{Y}(X, \widetilde{Y})$$
is a homotopy equivalence, which follows immediately from $(2)$. 

The implication $(1) \Rightarrow (3)$ follows from the proof of Proposition \ref{cell}, and the implication $(3) \Rightarrow (2)$ is classical (see \cite{haver}).
\end{proof}

\begin{remark}
It is possible to prove the following generalization of Proposition \ref{cell}: a proper geometric morphism $p_{\ast}: \calX \rightarrow \calY$ is cell-like if and only if, for each
object $U \in \calY$, the associated geometric morphism
$\calX_{/p^{\ast} U} \rightarrow \calY_{/U}$ is a shape equivalence (and, in fact, it is only necessary to check this on a collection of objects $U \in \calY$ which generates $\calY$ under colimits). 
\end{remark}

\begin{remark}
Another useful property of the class of cell-like morphisms, which we will not prove here, is stability under base change: given a pullback diagram
$$ \xymatrix{ \calX' \ar[d]^{p'_{\ast}} \ar[r] & \calX \ar[d]^{p_{\ast}} \\
\calY' \ar[r] & \calY }$$
where $p_{\ast}$ is cell-like, $p'_{\ast}$ is also cell-like.
\end{remark}

If $p_{\ast}: \calX \rightarrow \calY$ is a cell-like morphism of $\infty$-topoi, then many properties of $\calY$ are controlled by the analogous properties of $\calX$. For example:

\begin{proposition}
Let $p_{\ast}: \calX \rightarrow \calY$ be a cell-like morphism of $\infty$-topoi. If
$\calX$ has homotopy dimension $\leq n$, then $\calY$ also has homotopy dimension $\leq n$.
\end{proposition}

\begin{proof}
Let $1_{\calY}$ be a final object of $\calY$, $U$ an $n$-connective object of $\calY$, and
$p^{\ast}$ a left adjoint to $p_{\ast}$. We wish to prove that $\Hom_{h \calY}(1_{\calY},U)$ is nonempty. Since $p^{\ast}$ is fully faithful, it will suffice to prove that
$\Hom_{\h{\calX}}( p^{\ast} 1_{\calY}, p^{\ast} U)$. We now observe that $p^{\ast} 1_{\calY}$ is a final object of $\calX$ (since $p$ is left exact), $p^{\ast}U$ is $n$-connective (Proposition \ref{inftychange}), and $\calX$ has homotopy dimension $\leq n$, so that
$\Hom_{\h{\calX}}( p^{\ast} 1_{\calY}, p^{\ast} U)$ is nonempty as desired.
\end{proof}

We conclude with a different example of a class of cell-like maps. We will assume in the following discussion that the reader is familiar with the basic ideas of rigid analytic geometry; for an account of this theory we refer the reader to \cite{rigidgeom}. Let
$K$ be field which is complete with respect to a non-Archimedean absolute value
$||_K: K \rightarrow \R$. Let $A$ be an affinoid algebra over $K$: that is, a quotient of an algebra of convergent power series (in several variables) with values in $K$. Let $X$ be
the rigid space associated to $A$. One can associate to $X$ two different ``underlying'' topological spaces:

\begin{itemize}
\item[$(ZR1)$] The category $\calC$ of rational open subsets of $X$ has a Grothendieck topology, given by admissible affine covers. The topos of sheaves of sets on $\calC$ is localic, and the underlying locale has enough points: it is therefore isomorphic to the locale of open subsets of a (canonically determined) topological space $X_{ZR}$, the {\it Zariski-Riemann} space of $X$.
 
\item[$(ZR2)$] In the case where $K$ is a {\em discretely} valued field with ring of integers $R$, one may define $X_{ZR}$ to be the inverse limit of the underlying spaces of all formal schemes
$\hat{X} \rightarrow \Spf R$ which have generic fiber $X$.
\item[$(ZR3)$] Concretely, $X_{ZR}$ can be identified with the set of all isomorphism classes of continuous multiplicative seminorms $||_A: A \rightarrow M \cup \{\infty\}$, where $M$ is an ordered abelian group containing the value group $|K^{\ast}|_{K} \subseteq \R^{\ast}$, and the restriction
of $||_A$ to $K$ is $||_K$.

\item[$(B1)$] The category of sheaves of sets on $\calC$ contains a full subcategory, consisting of {\em overconvergent} sheaves. This category is also a localic topos, and the underyling locale is isomorphic to the lattice of open subsets of a (canonically determined) topological space $X_{B}$, the {\it Berkovich space} of $X$. The category of overconvergent sheaves is a localization of the category of all sheaves on $\calC$, and there is an associated map of topological spaces
$p: X_{ZR} \rightarrow X_B$.\index{gen}{Berkovich space}

\item[$(B2)$] Concretely, $X_{B}$ can be identified with the set of all continuous multiplicative seminorms $||_A: A \rightarrow \R \cup \{\infty\}$ which extend $||_K$. It is equipped with the topology of pointwise convergence, and is a compact Hausdorff space. 
\end{itemize}

The relationship between the Zariski-Riemann space $X_{ZR}$ and the Berkovich space $X_{B}$ (or, more conceptually, the relationship between the category of {\em all} sheaves on $X$ and the category of {\em overconvergent} sheaves on $X$) is neatly summarized by the following result.

\begin{proposition}\label{rigidex}
Let $K$ be a field which is complete with respect to a non-Archimedean absolute value $||_K$, let $A$ be an affinoid algebra over $K$, let $X$ be the associated rigid space, and $p: X_{ZR} \rightarrow X_{B}$ the natural map. Then $p$ induces a
cell-like morphism of $\infty$-topoi $p_{\ast}: \Shv(X_{ZR}) \rightarrow \Shv(X_{B})$. 
\end{proposition}

Before giving the proof, we need an easy lemma.
Recall that a topological space $X$ is {\it irreducible} if every finite collection of nonempty open subsets of $X$ has nonempty intersections.\index{gen}{irreducible!topological space}\index{gen}{topological space!irreducible}

\begin{lemma}\label{silshape}
Let $X$ be an irreducible topological space. Then $\Shv(X)$ has trivial shape.
\end{lemma}

\begin{proof}
Let $\pi: X \rightarrow \ast$ be the projection from $X$ to a point, $\pi_{\ast}: \Shv(X) \rightarrow \Shv(\ast)$ the induced geometric morphism. We will construct a left adjoint $\pi^{\ast}$ to
$\pi_{\ast}$ such that the unit map $\id \rightarrow \pi_{\ast} \pi^{\ast}$ is an equivalence. 

We begin by defining $G: \calP(X) \rightarrow \calP(\ast)$ to be the functor
given by composition with $\pi^{-1}$, so that $G | \Shv(X) = \pi_{\ast}$. Let
$$i: \Nerve (\calU(X))^{op} \rightarrow \Nerve(\calU(\ast))^{op}$$
be defined so that
$$ i(U) = \begin{cases} \emptyset & \text{if } U = \emptyset \\
\{ \ast \} & \text{if } U \neq \emptyset, \end{cases}$$
and let $F: \calP(\ast) \rightarrow \calP(U)$ be given by composition with $i$.
We observe that $F$ is a left Kan extension functor, so that the identity map
$$ \id_{ \calP(\ast)} \rightarrow G \circ F$$
exhibits $F$ as a left adjoint to $G$. We will show that $F( \Shv(\ast) ) \subseteq \Shv(X)$.
Setting $\pi^{\ast} = F | \Shv(\ast)$, we conclude that the identity map
$$ \id_{ \Shv(\ast)} \rightarrow \pi_{\ast} \pi^{\ast} $$
is the unit of an adjunction between $\pi_{\ast}$ and $\pi^{\ast}$, which will complete the proof.

Let $\calU \subseteq \calU(X)$ be a sieve which covers the open set $U \subseteq X$.
We wish to prove that the diagram
$$p: \Nerve(\calU^{op})^{\triangleleft} \rightarrow \Nerve(\calU(X))^{op}
\stackrel{i}{\rightarrow} \Nerve(\calU(\ast))^{op} \stackrel{\calF}{\rightarrow} \SSet$$
is a limit. Let $\calU_0 = \{ V \in \calU: V \neq \emptyset \}$. Since
$\calF(\emptyset)$ is a final object of $\SSet$, $p$ is a limit if and only if
$p | \Nerve(\calU_0^{op})^{\triangleleft}$ is a limit diagram. If $U = \emptyset$,
then this follows from the fact that $\calF(\emptyset)$ is final in $\SSet$.
If $U \neq \emptyset$, then $p | \Nerve(\calU_0^{op})^{\triangleleft}$ is
a constant diagram, so it will suffice to prove that the simplicial set
$\Nerve(\calU_0)^{op}$ is weakly contractible. This follows from the observation
that $\calU_0^{op}$ is a filtered partially ordered set, since $\calU_0$ is nonempty and stable under finite intersections (because $X$ is irreducible). 
\end{proof}

\begin{proof}[Proof of Proposition \ref{rigidex}]
We first show that $p_{\ast}$ is a proper map of $\infty$-topoi. We note that $p$ factors as a composition
$$ X_{ZR} \stackrel{p'}{\rightarrow} X_{ZR} \times X_{B} \stackrel{p''}{\rightarrow} X_{B}.$$
The map $p'$ is a pullback of the diagonal map $X_{B} \rightarrow X_{B} \times X_{B}$. Since
$X_{B}$ is Hausdorff, $p'$ is a closed immersion. It follows $p'_{\ast}$ is a closed immersion of $\infty$-topoi (Corollary \ref{glad1}) and therefore a proper morphism (Proposition \ref{closeduse2}). It therefore suffices to prove that $p''$ is a proper map of $\infty$-topoi. We note the existence of a commutative diagram
$$ \xymatrix{ \Shv(X_{ZR} \times X_{B}) \ar[d]^{p''_{\ast}} \ar[r] & \Shv(X_{ZR}) \ar[d]^{g_{\ast}} \\
\Shv( X_{B}) \ar[r] & \Shv(\ast). }$$
Using Proposition \ref{cartmun}, we deduce that this is a homotopy Cartesian diagram of $\infty$-topoi. It therefore suffices to show that the global sections functor $g_{\ast}: \Shv(X_{ZR}) \rightarrow \Shv(\ast)$ is proper, which follows from Corollary \ref{applyZR}.

We now observe that the topological space $X_{B}$ is paracompact and has finite covering dimension (\cite{berkovich}, Corollary $3.2.8$), so that $\Shv(X_{B})$ has enough points (Corollary \ref{enuff}). According to Proposition \ref{fibersilly}, it suffices to show that for every fiber diagram
$$ \xymatrix{ \calX' \ar[r] \ar[d] & \Shv(X_{ZR}) \ar[d] \\
\Shv(\ast) \ar[r]^-{q_\ast} & \Shv(X_{B}), }$$
the $\infty$-topos $\calX'$ has trivial shape. Using Lemma \ref{eoi2}, we conclude that $q_{\ast}$ is necessarily induced by a homomorphism of locales $\calU(X_{B}) \rightarrow \calU(\ast)$, which corresponds to an irreducible closed subset of $X_{B}$. Since $X_{B}$ is Hausdorff, this subset consists of a single (closed) point $x$. Using Proposition \ref{closeduse2} and Corollary \ref{glad1}, we can identify $\calX'$ with the $\infty$-topos $\Shv(Y)$, where $Y = X_{ZR} \times_{X_B} \{x\}$.
We now observe that the topological space $Y$ is coherent and irreducible (it contains a unique ``generic'' point), so that $\Shv(Y)$ has trivial shape by Lemma \ref{silshape}.
\end{proof}

\begin{remark}
Let $p_{\ast}: \Shv(X_{ZR}) \rightarrow \Shv(X_{B})$ be as in Proposition \ref{rigidex}.
Then $p_{\ast}$ has a fully faithful left adjoint $p^{\ast}$. We might say that an object of
$\Shv(X_{ZR})$ is {\it overconvergent} if it belongs to the essential image of $p^{\ast}$; for sheaves of sets, this agrees with the classical terminology.\index{gen}{overconvergent sheaf}
\end{remark}

\begin{remark}
One can generalize Proposition \ref{rigidex} to rigid spaces which are not affinoid; we leave the details to the reader.
\end{remark}

\appendix

\chapter{Appendix}

\setcounter{theorem}{0}
\setcounter{subsection}{0}

This appendix is comprised of three parts. In \S \ref{catreview}, we will review some ideas from classical category theory, such as monoidal structures, enriched categories, and Quillen's small object argument. We give a brief overview of the theory of model categories in \S \ref{appmodelcat}.
The main result here is Proposition \ref{goot}, which will allow us to establish the existence of model category structures in a variety of situations with a minimal amount of effort. In \S \ref{techapp} we will use this result to make detailed study of the theory of simplicial categories.
Our exposition is rather dense; for a more leisurely account of the theory of model categories we refer the reader to one of the standard texts, such as \cite{hovey}.

$\{W(t)\}\equiv \{W(t),t\geq 0\}\glossary{Wt@$W(t)$ & State or window
size of the source at time $t$}$

\section{Category Theory}\label{catreview}

\setcounter{theorem}{0}

Familiarity with classical category theory is the main prerequisite for reading this book.
In this section, we will fix some of the notation that we will use when discussing categories, and summarize (generally without proofs) some of the concepts which we will use in the body of the text.\index{gen}{category}

If $\calC$ is a category, we let $\Ob(\calC)$ denote the set of objects of $\calC$.\index{not}{ObC@$\Ob(\calC)$} We will write $X \in \calC$ to
mean that $X$ is an object of $\calC$. For $X,Y \in \calC$, we write $\Hom_{\calC}(X,Y)$ for the set of morphisms from $X$ to $Y$ in $\calC$. We also write
$\id_X$ for the identity automorphism of $X \in \calC$ (regarded as an element of $\Hom_{\calC}(X,X)$).\index{not}{HomC@$\Hom_{\calC}(X,Y)$}

If $Z$ is an object in a category $\calC$, then the {\it overcategory}\index{gen}{overcategory}
$\calC_{/Z}$ of {\it objects over $Z$} is defined as follows: the
objects of $\calC_{/Z}$ are diagrams $X \rightarrow Z$ in $\calC$. A morphism
from $f: X \rightarrow Z$ to $g: Y \rightarrow Z$ is a commutative triangle
$$ \xymatrix{ X \ar[rr] \ar[dr]_{f} & & Y \ar[dl]^{g} \\
& Z.}$$
Dually, we have an {\it undercategory} $\calC_{Z/} = ((\calC^{op})_{Z/})^{op}$ of {\em objects under $Z$}\index{gen}{undercategory}.\index{not}{calC/X@$\calC_{/X}$}\index{not}{calCX/@$\calX_{X/}$}

If $f: X \rightarrow Z$ and $g: Y \rightarrow Z$ are objects in
$\calC_{/Z}$, then we will often write $\Hom_{Z}(X,Y)$ rather than
$\Hom_{\calC_{/Z}}(f,g)$.\index{not}{HomZXY@$\Hom_{Z}(X,Y)$}

We let $\Set$ denote the category of sets, and $\Cat$ the category of (small) categories (where the morphisms are given by functors).\index{not}{Cat@$\Cat$}\index{not}{Set@$\Set$}

If $\kappa$ is a regular cardinal, we will say that a set $S$ is {\it $\kappa$-small}\index{gen}{$\kappa$-small} if it has cardinality less than $\kappa$. We will also use this terminology when discussing mathematical objects other than sets, which are built out of sets. For example, we will say that a category $\calC$ is {\it $\kappa$-small} if the set of all objects of $\calC$ is $\kappa$-small, and the set of all morphisms in $\calC$ is likewise $\kappa$-small.

We will need to discuss categories which are not small. In order to minimize the effort spent dealing with set-theoretic complications, we will adopt the usual device of ``Grothendieck universes''. We fix a strongly inaccessible cardinal $\kappa$, and refer to a mathematical object (such as a set or category) as {\it small} if it is $\kappa$-small, and {\it large} otherwise.
It should be emphasized that this is primarily a linguistic device, and that none of our results depend in an essential way on the existence of a strongly inaccessible cardinal $\kappa$.\index{gen}{small}

Throughout this book, the word ``topos'' will always mean {\em Grothendieck topos}. Strictly speaking, a knowledge of classical topos theory is not required to read this paper: all of the relevant classical concepts will be introduced (though sometimes in a hurried fashion) in the course of our search for suitable $\infty$-categorical analogues. 

\subsection{Compactness and Presentability}

Let $\kappa$ be a regular cardinal.

\begin{definition}\index{gen}{filtered!partially ordered set}
A partially ordered set $\calI$ is {\it $\kappa$-filtered} if, for any subset $\calI_0 \subseteq \calI$ having cardinality $< \kappa$, there exists an upper bound for $\calI_0$ in $\calI$.
\end{definition}

Let $\calC$ be a category which admits (small) colimits, and 
let $X$ be an object of $\calC$.
Suppose given a $\kappa$-filtered partially ordered set $\calI$ and a diagram $\{Y_{\alpha} \}_{\alpha \in \calI}$ in $\calC$, indexed by $\calI$. Let $Y$ denote a colimit of this diagram.
There there is an associated map of sets
$$ \psi: \colim \Hom_{\calC}(X, Y_{\alpha}) \rightarrow \Hom_{\calC}(X,Y).$$
We say that $X$ is {\it $\kappa$-compact} if $\psi$ is bijective, for {\em every} $\kappa$-filtered
partially ordered set $\calI$ and {\em every} diagram $\{ Y_{\alpha} \}$ indexed by $\calI$.
We say that $X$ is {\it small} if it is $\kappa$-compact for some (small) regular cardinal $\kappa$.
In this case, $X$ is $\kappa$-compact for all sufficiently large regular cardinals $\kappa$.\index{gen}{compact!object of a category}\index{gen}{small!object of a category}

\begin{definition}\label{catpor}\index{gen}{presentable!category}
A category $\calC$ is {\it presentable} if it satisfies the following conditions:
\begin{itemize}
\item[$(1)$] The category $\calC$ admits all (small) colimits.

\item[$(2)$] There exists a (small) set $S$ of objects of $\calC$ which generates
$\calC$ under colimits; in other words, every object of $\calC$ may be obtained as the colimit of a (small) diagram taking values in $S$.

\item[$(3)$] Every object in $\calC$ is small. (Assuming $(2)$, this is equivalent to the assertion
that every object which belongs to $S$ is small.)

\item[$(4)$] For any pair of objects $X,Y \in \calC$, the set $\Hom_{\calC}(X,Y)$ is small.
\end{itemize}
\end{definition}

\begin{remark}
In \S \ref{c5s6}, we describe an $\infty$-categorical generalization of Definition \ref{catpor}.
\end{remark}

\begin{remark}
For more details of the theory of presentable categories, we refer the reader to \cite{adamek}. Note that our terminology differs slightly from that of \cite{adamek}, in which our presentable categories are called {\it locally presentable} categories.
\end{remark}

\subsection{Lifting Problems and the Small Object Argument}\label{liftingprobs}

Let $\calC$ be a category, and let $p: A \rightarrow B$ and $q: X
\rightarrow Y$ be morphisms in $\calC$. Recall that $p$ is said to
have the {\it left lifting property} with respect to $q$, and $q$
the {\it right lifting property} with respect to $p$, if given any diagram
$$ \xymatrix{ A \ar[d]^{p} \ar[r] & X \ar[d]^q \\
B \ar@{-->}[ur] \ar[r] & Y \\} $$
there exists a dotted arrow as indicated, rendering the diagram commutative.\index{gen}{left lifting property}\index{gen}{right lifting property}

\begin{remark}
In the case where $Y$ is a final object of $\calC$, we will instead say that
$X$ has the {\it extension property} with respect to $p: A \rightarrow B$.\index{gen}{extension property}
\end{remark}

Let $S$ be any collection of morphisms in $\calC$. We define
$_{\perp} S$ to be the class of all morphisms which have the right
lifting property with respect to all morphisms in $S$, and $S_{\perp}$
to be the class of all morphisms which have the left lifting
property with respect to all morphisms in $S$. We observe that
$$S \subseteq (_{\perp}S)_{\perp}.$$\index{not}{Sperp@$_{\perp}S$}\index{not}{Stop@$S_{\perp}$}

The class of morphisms $(_{\perp} S)_{\perp}$ enjoys several
stability properties which we axiomatize in the following definition.

\begin{definition}\label{saturated}\index{gen}{weakly saturated}\index{gen}{saturated!weakly}
Let $\calC$ be a category with all (small) colimits, and let $S$ be a class of morphisms of $\calC$. We will say that $S$ is {\it weakly saturated} if it has the following properties:

\begin{itemize}
\item[$(1)$] (Closure under the formation of pushouts) Given a pushout diagram
$$ \xymatrix{ C \ar[r]^{f} \ar[d] & D \ar[d] \\
C' \ar[r]^{f'} & D' }$$
such that $f$ belongs to $S$, the morphism $f'$ also belongs to $S$. 

\item[$(2)$] (Closure under transfinite composition) Let $C \in \calC$ be an object, $\alpha$
an ordinal, and let $\{ D_{\beta} \}_{\beta < \alpha} $ be a system of objects of $\calC_{C/}$
indexed by $\alpha$: in other words, for each $\beta < \alpha$, we are supplied with a morphism $C \rightarrow D_{\beta}$, and for each $\gamma \leq \beta < \alpha$ a commutative diagram
$$ \xymatrix{ & D_{\gamma} \ar[dd]^{\phi_{\gamma,\beta}} \\
C \ar[ur] \ar[dr] \\
& D_{\beta} }$$
satisfying $\phi_{\beta,\gamma} \circ \phi_{\gamma, \delta} = \phi_{\beta, \delta}.$
For $\beta \leq \alpha$, we let $D_{<\beta}$ be a colimit of the system
$\{ D_{\gamma} \}_{\gamma < \beta}$, taken in the category $\calC_{C/}$. 

Suppose that, for each $\beta < \alpha$, the natural map $D_{< \beta} \rightarrow D_{\beta}$
belongs to $S$. Then the induced map $C \rightarrow D_{<\alpha}$ belongs to $S$.

\item[$(3)$] (Closure under the formation of retracts) Given a commutative diagram
$$ \xymatrix{ C \ar[d]^{f} \ar[r] & C' \ar[d]^{g} \ar[r] & C \ar[d]^{f} \\
D \ar[r] & D' \ar[r] & D }$$
in which both horizontal compositions are the identity, if $g$ belongs to $S$, then so does $f$.
\end{itemize}
\end{definition}

It is worth noting that saturation has the following consequences:

\begin{proposition}
Let $\calC$ be a category which admits all $($small$)$ colimits, and let $S$ be a weakly saturated class of morphism in $\calC$. Then:
\begin{itemize}
\item[$(1)$] Every isomorphism belongs to $S$.
\item[$(2)$] The class $S$ is stable under composition: if $f: X \rightarrow Y$ and $g: Y \rightarrow Z$ belong to $S$, then so does $g \circ f$.
\end{itemize}
\end{proposition}

\begin{proof}
Assertion $(1)$ is equivalent to the closure of $S$ under transfinite composition, in the special case where $\alpha=0$; $(2)$ is equivalent to the special case where $\alpha=2$.
\end{proof}

\begin{remark}
A reader who is ill-at-ease with the style of the preceding argument should feel free to take
the asserted properties as part of the definition of a weakly saturated class of morphisms.
\end{remark}

The intersection of any collection of weakly saturated classes of morphisms is itself weakly saturated. Consequently, for any category $\calC$ which admits small colimits, and any collection
$A$ of morphisms of $\calC$, there exists a {\em smallest} weakly saturated class of morphisms containing $A$: we will call this the weakly saturated class of morphisms {\it generated} by $A$.
We note that $(_{\perp} A)_{\perp}$ is weakly saturated. Under appropriate set-theoretic assumptions, Quillen's ``small object'' argument can be used to establish that $(_{\perp} A)_{\perp}$ is the weakly saturated class generated by $A$:

\begin{proposition}[Small Object Argument]\label{quillobj}\index{gen}{small object argument}
Let $\calC$ be a presentable category and $A_0 = \{ \phi_i: C_i \rightarrow D_i \}_{i \in I}$ a collection
of morphisms in $\calC$ indexed by a $($small$)$ set $I$. For each $n \geq 0$, let
$\calC^{[n]}$ denote the category of functors from the linearly ordered set $[n] = \{0, \ldots, n\}$ into $\calC$. There exists a functor
$T: \calC^{[1]} \rightarrow \calC^{[2]}$ with the following properties:
\begin{itemize}
\item[$(1)$] The functor $T$ carries a morphism $f: X \rightarrow Z$ to a diagram
$$ \xymatrix{ & Y \ar[dr]^{f''} & \\
X \ar[ur]^{f'} \ar[rr]^{f} & & Z }$$
where $f'$ belongs to the weakly saturated class of morphisms generated by $A_0$ and $f''$
has the right lifting property with respect to each morphism in $A_0$.
\item[$(2)$] If $\kappa$ is a regular cardinal such that each of the objects $C_i$, $D_i$ is $\kappa$-compact, then $T$ commutes with $\kappa$-filtered colimits.
\end{itemize}
\end{proposition}

\begin{proof}
Fix a regular cardinal $\kappa$ as in $(2)$, and fix a morphism $f: X \rightarrow Z$ in $\calC$.
We will give a functorial construction of the desired diagram
$$ \xymatrix{ & Y \ar[dr]^{f''} & \\
X \ar[ur]^{f'} \ar[rr]^{f} & & Z }$$
We define a transfinite sequence of objects
$$ Y_0 \rightarrow Y_1 \rightarrow \ldots $$
in $\calC_{/Z}$, indexed by ordinals smaller than $\kappa$. Let $Y_0 = X$, and let $Y_{\lambda} = \colim_{ \alpha < \lambda} Y_{\alpha}$ when $\lambda$ is a nonzero limit ordinal. For $i \in I$, let
$F_i: \calC_{/Z} \rightarrow \Set$ be the functor
$$(T \rightarrow Z ) \mapsto \Hom_{\calC}( D_i, Z) \times_{ \Hom_{\calC}(C_i, Z)}
\Hom_{\calC}(C_i, T).$$
Supposing that $Y_{\alpha}$ has been defined, we define $Y_{\alpha+1}$ by the following pushout diagram
$$ \xymatrix{ \coprod_{i \in I, \eta \in F_i(Y_{\alpha})} C_i \ar[r] \ar[d] & Y_{\alpha} \ar[d] \\
\coprod_{i \in I, \eta \in F_i(Y_{\alpha})} D_i \ar[r] & Y_{\alpha+1}. }$$
We conclude by defining $Y$ to be $\colim_{\alpha < \kappa} Y_{\alpha}$. It is easy to check that this construction has the desired properties.
\end{proof}

\begin{remark}
If $\calC$ is enriched, tensored and cotensored over another presentable monoidal category
$\bfS$ (see \S \ref{enrichcat}), then a similar construction shows that we can choose
$T$ to be a $\bfS$-enriched functor.
\end{remark}

\begin{corollary}\label{tilobj}
Let $\calC$ be a presentable category, and let $A$ be a {\em set} of morphisms of $\calC$.
Then $(_{\perp} A)_{\perp}$ is the smallest weakly saturated class of morphisms containing $A$.
\end{corollary}

\begin{proof}
Let $\overline{A}$ be the smallest weakly saturated class of morphisms containing $A$, so that
$\overline{A} \subseteq (_{\perp} A)_{\perp}$. To establish the reverse inclusion, 
For the reverse inclusion, let us suppose that $f: X \rightarrow Z$ belongs to
$(_{\perp} A)_{\perp}$. Proposition \ref{quillobj} implies the existence of a factorization
$$ X \stackrel{f'}{\rightarrow} Y \stackrel{f''}{\rightarrow} Z$$
where $f' \in \overline{A}$ and $f''$ belongs to $_{\perp} A$. It follows that $f$ has the left lifting property with respect to $f''$, so that $f$ is a retract of $f'$ and therefore belongs to $\overline{A}$.
\end{proof}

\begin{remark}\label{easyprest}
Let $\calC$ be a presentable category, let $S$ be a (small) set of morphisms in $\calC$, and suppose that $f: X \rightarrow Y$ belongs to the weakly saturated class of morphisms generated by $S$.
The proofs of Proposition \ref{quillobj} and Corollary \ref{tilobj} show that there exists a transfinite sequence 
$$ Y_0 \rightarrow Y_1 \rightarrow \ldots $$
of objects of $\calC_{X/}$, indexed by a set of ordinals $\{ \beta | \beta < \alpha \}$, with the following properties:
\begin{itemize}
\item[$(i)$] For each $\beta < \alpha$, there is a pushout diagram 
$$ \xymatrix{ C \ar[r]^{g} \ar[d] & D \ar[d] \\
\colim_{ \gamma < \beta} Y_{\gamma} \ar[r] & Y_{\beta},}$$
where the colimit is formed in $\calC_{X/}$ and $g \in S$.
\item[$(ii)$] The object $Y$ is a retract of $\colim_{\gamma < \alpha} Y_{\gamma}$
in the category $\calC_{X/}$.
\end{itemize}
\end{remark}

\subsection{Monoidal Categories}\label{monoidaldef}

A {\it monoidal category}\index{gen}{monoidal category}\index{gen}{category!monoidal} is a category $\calC$ equipped with a (coherently) associative ``product''
functor $\otimes: \calC \times \calC \rightarrow \calC$ and a unit object ${\bf 1}$. 
The associativity is expressed by demanding isomorphisms
$$ \eta_{A,B,C}: (A \otimes B) \otimes C \rightarrow A \otimes (B \otimes C),$$
and the requirement that ${\bf 1}$ be unital is expressed by demanding isomorphisms
$$ \alpha_{A}: A \otimes {\bf 1} \rightarrow A$$
$$ \beta_{A}: {\bf 1} \otimes A \rightarrow A.$$
We do not merely require the existence of these isomorphisms: they are part of the structure of a monoidal category. Moreover, these isomorphisms are required to satisfy the following conditions:

\begin{itemize}
\item The isomorphism $\eta_{A,B,C}$ depends {\em functorially} on the triple $(A,B,C)$; in other words, $\eta$ may be regarded as a natural isomorphism between the functors
$$ \calC \times \calC \times \calC \rightarrow \calC.$$
$$ (A,B,C) \mapsto (A \otimes B) \otimes C$$
$$ (A,B,C) \mapsto A \otimes (B \otimes C).$$
Similarly $\alpha_A$ and $\beta_A$ depend functorially on $A$. 

\item Given any quadruple $(A,B,C,D)$ of objects of $\calC$, the {\it MacLane pentagon}\index{gen}{MacLane pentagon}\index{gen}{pentagon axiom}
$$ \xymatrix{ & ((A \otimes B) \otimes C) \otimes D \ar[dl]^{\eta_{A,B,C} \otimes \id_D}
\ar[dr]^{\eta_{ A \otimes B,C,D}} \\
(A \otimes (B \otimes C)) \otimes D \ar[d]^{ \eta_{A, B \otimes C, D}} & & (A \otimes B) \otimes (C \otimes D) \ar[d]^{\eta_{A,B, C \otimes D}} \\
A \otimes (( B \otimes C) \otimes D) \ar[rr]^{\id_A \otimes \eta_{B,C,D}} & & A \otimes (B \otimes (C \otimes D))}$$
is commutative. 

\item For any pair $(A,B)$ of objects of $\calC$, the triangle
$$ \xymatrix{ (A \otimes {\bf 1}) \otimes B \ar[rr]^{ \eta_{A, {\bf 1}, B}} \ar[dr]^{ \alpha_A \otimes \id_B} & & A \otimes ( {\bf 1} \otimes B) \ar[dl]^{\id_A \otimes \beta_{B}} \\
& A \otimes B}$$ is commutative.
\end{itemize}

MacLane's coherence theorem asserts that the commutativity of this pair of diagrams implies the commutativity of {\em all} diagrams that can be written using only the isomorphisms $\eta_{A,B,C}$, $\alpha_A$, and $\beta_A$. More precisely, any monoidal category is equivalent (as a monoidal category) to a {\em strict} monoidal category: that is, a monoidal category in which $\otimes$ is literally associative, ${\bf 1}$ is literally a unit with respect to $\otimes$, and the isomorphisms $\eta_{A,B,C}$, $\alpha_A$, $\beta_A$ are the identity maps.\index{gen}{monoidal category!strict}
\index{gen}{MacLane's coherence theorem}

\begin{example}
Let $\calC$ be a category which admits finite products. Then $\calC$ admits the structure of a monoidal category, where the operation $\otimes$ is given by Cartesian product
$$ A \otimes B \simeq A \times B$$
and the isomorphisms $\eta_{A,B,C}$ are induced from the evident associativity of the Cartesian product. The identity ${\bf 1}$ is defined to be the final object of $\calC$, and the isomorphisms
$\alpha_A$ and $\beta_A$ are determined in the obvious way. We refer to this monoidal structure on $\calC$ as the {\em Cartesian monoidal structure}.

We remark that the Cartesian product $A \times B$ is only well-defined up to (unique) isomorphism (as is the final object ${\bf 1}$), so that strictly speaking the Cartesian monoidal structure on $\calC$ depends on various choices; however, all such choices lead to (canonically) equivalent monoidal categories.\index{gen}{monoidal category!Cartesian}
\end{example}

\begin{remark}
Let $(\calC, \otimes, {\bf 1}, \eta, \alpha, \beta)$ be a monoidal category. We will generally abuse notation by simply saying that $\calC$ is a monoidal category, or that $(\calC, \otimes)$ is a monoidal category, or that $\otimes$ is a {\it monoidal structure} on $\calC$; the other structure is implicitly understood to be present as well.
\end{remark}

\begin{remark}
Let $\calC$ be a category equipped with a monoidal structure $\otimes$. Then we may define a new monoidal structure on $\calC$, by setting $A \otimes^{op} B = B \otimes A$. We refer to this monoidal structure $\otimes^{op}$ as the {\it opposite} of the monoidal structure $\otimes$.
\end{remark}

\begin{definition}\label{tukerdef}
A monoidal category $(\calC,\otimes)$ is said to be {\it left closed} if, for each $A \in \calC$, the functor $$ N \mapsto A \otimes N$$
admits a right adjoint
$$ Y \mapsto {}^A\!Y.$$
We say that $(\calC, \otimes)$ is {\it right-closed} if the opposite monoidal structure $(\calC, \otimes^{op})$ is left-closed; in other words, if every functor
$$ N \mapsto N \otimes A$$ has a right adjoint
$$ Y \mapsto Y^A.$$
Finally, we say that $(\calC, \otimes)$ is {\it closed} if it is both right-closed and left-closed.\index{gen}{monoidal category!left closed}\index{gen}{monoidal category!right closed}\index{gen}{monoidal category!closed}
\end{definition}

In the setting of monoidal categories, it is appropriate to consider only those functors which
are compatible with the monoidal structures in the following sense:

\begin{definition}\index{gen}{functor!lax monoidal}
Let $(\calC, \otimes)$ and $( \calD, \otimes)$ be monoidal categories. A {\it right-lax monoidal functor} from $\calC$ to $\calD$ consists of the following data:
\begin{itemize}
\item A functor $G: \calC \rightarrow \calD$.
\item A natural transformation $\gamma_{A,B}: G(A) \otimes G(B) \rightarrow G(A \otimes B)$ rendering commutative the diagram
$$ \xymatrix{ (G(A) \otimes G(B)) \otimes G(C) \ar[r]^{\gamma_{A,B}}
\ar[d]^{ \eta_{G(A),G(B),G(C)} } 
& G(A \otimes B) \otimes G(C) \ar[r]^{ \gamma_{A \otimes B, C} } & G((A \otimes B) \otimes C)
\ar[d]^{ G(\eta_{A,B,C})} \\
G(A) \otimes( G(B) \otimes G(C) ) \ar[r]^{\gamma_{B,C}}
& G(A) \otimes G(B \otimes C) \ar[r]^{ \gamma_{A, B \otimes C} } & G(A \otimes (B \otimes C)). }$$
\item A map $e: {\bf 1_{\calD}} \rightarrow G( {\bf 1_{\calC}} )$ rendering commutative the diagrams
$$ \xymatrix{ G(A) \otimes {\bf 1_{\calD} } \ar[r]^{ \id \otimes e} \ar[dr]^{\alpha_{G(A)}} & G(A) \otimes G( \bf{1_{\calC}}) \ar[r]^{\gamma_{A, {\bf 1_{\calC}}}} & G(A \otimes {\bf 1_{\calC}}) \ar[dl]^{G(\alpha_A)} \\ 
& G(A) }$$
$$ \xymatrix{  {\bf 1}_{\calD} \otimes G(B) \ar[r]^{ e \otimes \id} 
\ar[dr]^{\beta_{G(B)}} & G( {\bf 1_{\calC}}) \otimes G(B) \ar[r]^{\gamma_{{\bf 1_{\calC}},A}} & G({\bf 1_{\calC}} \otimes B) \ar[dl]^{G(\alpha_B)} \\ 
& G(B) & }.$$
\end{itemize}

A natural transformation between right-lax monoidal functors is {\it monoidal} if it commutes with
the maps $\gamma_{A,B}$, $e$.\index{gen}{functor!monoidal}
\end{definition}

Dually, a {\it left-lax monoidal functor} from $\calC$ to $\calD$ consists of a right-lax monoidal functor from $\calC^{op}$ to $\calD^{op}$; it is determined by giving a functor
$F: \calC \rightarrow \calD$ together with a map $e': F( {\bf 1_{\calC}}) \rightarrow {\bf 1}_{\calD}$
and a natural transformation
$$ \gamma'_{A,B}: F(A \otimes B) \rightarrow F(A) \otimes F(B)$$
satisfying the appropriate analogues of the conditions listed above.

If $F$ is a right-lax monoidal functor via {\em isomorphisms}
$$ e: {\bf 1}_{\calD} \rightarrow F( { \bf 1_{\calC} })$$
$$ \gamma_{A,B}: F(A) \otimes F(B) \rightarrow F(A \otimes B),$$
then $F$ may be regarded as a left-lax monoidal functor by setting $e' = e^{-1}$,
$\gamma'_{A,B} = \gamma_{A,B}^{-1}$. In this case, we simply say that $F$ is a {\em monoidal} functor.

\begin{remark}
Let 
$$ \Adjoint{F}{\calC}{\calD}{G}$$ be an adjunction between categories $\calC$ and $\calD$. Suppose that $\calC$ and $\calD$ are equipped with monoidal structures. Then endowing $G$ with the structure of a right-lax monoidal functor is equivalent to endowing $F$ with the structure of a left-lax monoidal functor.
\end{remark}

\begin{example}
Let $\calC$ and $\calD$ be categories which admit finite products, and let
$F: \calC \rightarrow \calD$ be a functor between them. Then, if we regard $\calC$ and
$\calD$ as endowed with the Cartesian monoidal structure, then $F$ acquires the structure
of a left lax-monoidal functor in a canonical way, via the maps
$F(A \times B) \rightarrow F(A) \times F(B)$ induced from the functoriality of $F$. In this case,
$F$ is a monoidal functor if and only if it commutes with finite products.
\end{example}

\subsection{Enriched Category Theory}\label{enrichcat}

One frequently encounters categories $\calD$ in which the collections of morphisms
$\Hom_{\calD}(X,Y)$ between two objects $X,Y \in \calD$ has additional structure: for example, a topology, or a group structure, or the structure of a vector space. These situations may all be efficiently described using the language of {\it enriched category theory}, which we now introduce.

Let $(\calC, \otimes)$ be a monoidal category. A {\it $\calC$-enriched category} $\calD$ consists of the following data:\index{gen}{enriched category}\index{gen}{category!enriched}

\begin{itemize}
\item[$(1)$] A collection of objects.

\item[$(2)$] For every pair of objects $X,Y \in \calD$, a mapping object
$\bHom_{\calD}(X,Y)$ of $\calC$.

\item[$(3)$] For every triple of objects $X,Y,Z \in \calD$, a composition map
$$ \bHom_{\calD}(Y,Z) \otimes \bHom_{\calD}(X,Y) \rightarrow \bHom_{\calD}(X,Z).$$
Composition is required to be associative in the sense that for any $W,X,Y,Z \in \calC$, the diagram
$$ \xymatrix{ \bHom_{\calD}(Z,Y) \otimes \bHom_{\calD}(Y,X) \otimes \bHom_{\calD}(X,W)
\ar[r] \ar[d] & 
\bHom_{\calD}(Z,X) \otimes \bHom_{\calD}(X,W) \ar[d] \\ 
\bHom_{\calD}(Z,Y) \otimes \bHom_{\calD}(Y,W) \ar[r] & \bHom_{\calD}(Z,W)}$$
is commutative.

\item[$(4)$] For every object $X \in \calD$, a unit map ${\bf 1} \rightarrow \bHom_{\calD}(X,X)$
rendering commutative the diagrams
$$ \xymatrix{ { \bf 1} \otimes \bHom_{\calD}(Y,X) \ar[rr] \ar[dr] & &  \bHom_{\calD}(X,X) \otimes \bHom_{\calD}(Y,X) \ar[dl] \\
& \bHom_{\calD}(Y,X) }$$

$$ \xymatrix{ \bHom_{\calD}(X,Y) \otimes {\bf 1} \ar[rr] \ar[dr] & &  \bHom_{\calD}(X,Y) \otimes \bHom_{\calD}(X,X) \ar[dl] \\
& \bHom_{\calD}(X,Y).}$$
\end{itemize}

\begin{example}
Suppose that $(\calC, \otimes)$ is a {\em right-closed} monoidal category. Then $\calC$ is enriched over itself in a natural way, if one defines $\bHom_{\calC}(X,Y) = Y^{X}$.
\end{example}

\begin{example}
Let $\calC$ be the category of sets, with the Cartesian monoidal structure. Then a $\calC$-enriched category is simply a category in the usual sense. 
\end{example}

\begin{remark}\label{laxcon}
Let $G: \calC \rightarrow \calC'$ be a right-lax monoidal functor between monoidal categories.
Suppose that $\calD$ is a category enriched over $\calC$. We may define a category
$G(\calD)$, enriched over $\calC'$, as follows:

\begin{itemize}
\item[$(1)$] The objects of $G(\calD)$ are the objects of $\calD$.
\item[$(2)$] Given objects $X,Y \in \calD$, we set $$\bHom_{G(\calD)}(X,Y) = G( \bHom_{\calD}(X,Y) ).$$
\item[$(3)$] The composition in $G(\calD)$ is given by the map
$$ G( \bHom_{\calD}(Y,Z) ) \otimes G( \bHom_{\calD}(X,Y) ) \rightarrow
G( \bHom_{\calD}(Y,Z) \otimes \bHom_{\calD}(X,Y) ) \rightarrow G( \bHom_{\calD}(X,Z)).$$
Here the first map is determined by the right-weakly monoidal structure on the functor $G$, and the second is obtained by applying $G$ to the composition law in the category $\calD$.
\item[$(4)$] For every object $X \in \calD$, the associated unit $G(\calD)$ is given by the composition
$$ {\bf 1_{\calC'} } \rightarrow G( { \bf 1_{\calC} } ) \rightarrow G( \bHom_{\calD}(X,X)).$$
\end{itemize}
\end{remark}

\begin{remark}\index{gen}{functor!enriched}
If $\calD$ and $\calD'$ are categories enriched over the same monoidal category $\calC$, then one can define a category of {\em $\calC$-enriched} functors from $\calD$ to $\calD'$ in the evident way. Namely, an enriched functor $F: \calD \rightarrow \calD'$ consists of a map from the objects of $\calD$ to the objects of $\calD'$ and a collection of morphisms
$$ \eta_{X,Y}: \bHom_{\calD}(X,Y) \rightarrow \bHom_{\calD'}(FX, FY)$$ with the following properties:
\begin{itemize}
\item[$(i)$] For each object $X \in \calD$, the composition
$$ {\bf 1_{\calC}} \rightarrow \bHom_{\calD}(X,X) \stackrel{\eta_{X.X}}{\rightarrow}
\bHom_{\calD'}(FX,FX) $$
coincides with the unit map for $FX \in \calD'$.
\item[$(ii)$] For every triple of objects $X, Y, Z \in \calD$, the diagram
$$ \xymatrix{ 
\bHom_{\calD}(X,Y) \otimes \bHom_{\calD}(Y,Z) \ar[r] \ar[d] & \bHom_{\calD}(X,Z) \ar[d] \\
\bHom_{\calD'}(FX,FY) \otimes \bHom_{\calD}(FY,FZ) \ar[r] & \bHom_{\calD}(FX,FZ) }$$
is commutative.
\end{itemize}
If $F$ and $F'$ are enriched functors, an {\em enriched natural transformation $\alpha$} from
$F$ to $F'$ consists of specifying, for each object $X \in \calD$, a morphism
$\alpha_{X} \in \Hom_{\calD'}( FX, F'X)$ which renders commutative the diagram
$$ \xymatrix{ \bHom_{\calD}(X,Y) \ar[r] \ar[d] & \bHom_{\calD'}(FX,FY) \ar[d]^{\alpha_Y} \\
\bHom_{\calD'}(F'X, F'Y) \ar[r]^{\alpha_X} & \bHom_{\calD'}(FX, F'Y). }$$\index{gen}{natural transformation!enriched}
\end{remark}

Suppose that $\calC$ is any monoidal category. Consider the functor $\calC \rightarrow \Set$ given by
$$ X \mapsto \Hom_{\calC}( { \bf 1}, X).$$
This is a right-lax monoidal functor from $(\calC, \otimes)$ to $\Set$, where the latter is equipped with the Cartesian monoidal structure. By the above remarks, we see that we may equip any $\calC$-enriched category $\calD$ with the structure of an ordinary category by setting
$$ \Hom_{\calD}(X,Y) = \Hom_{\calC}( {\bf 1}, \bHom_{\calD}(X,Y) ).$$ 
We will generally not distinguish notationally between $\calD$ as a $\calC$-enriched category
and this (underlying) category having the same objects. However, to avoid confusion, we use different notations for the morphisms: $\bHom_{\calD}(X,Y)$ is an object of $\calC$, while $\Hom_{\calD}(X,Y)$ is a set.

Let $\calC$ be a right-closed monoidal category, and $\calD$ a category enriched over $\calC$.
Fix objects $C \in \calC$, $X \in \calD$, and consider the functor
$$ \calD \rightarrow \calC$$
$$ Y \mapsto \bHom_{\calD}(X,Y)^{C}.$$
This functor may or may not be {\em corepresentable}\index{gen}{corepresentable!functor}, in the sense that there exists an object
$Z \in \calD$ and an isomorphism of functors
$$ \eta: \bHom_{\calD}(X, \bigdot)^{C} \simeq \bHom_{\calD}(Z, \bigdot).$$
If such an object $Z$ exists, we will denote it by $X \otimes C$. The natural isomorphism $\eta$ is determined by specifying a single map $\eta(X): C \rightarrow \bHom_{\calD}(X, X \otimes C)$. By general nonsense, the map $\eta(X)$ determines $X \otimes C$ up to (unique) isomorphism, provided that $X \otimes C$ exists. If the object $X \otimes C$ exists for every $C \in \calC$, $X \in \calD$, then we say that $\calD$ is {\it tensored over $\calC$}.\index{gen}{tensored} In this case, we may regard $$(X,C) \mapsto X \otimes C$$
as determining a functor $\calD \otimes \calC \rightarrow \calD$. Moreover, one has canonical
isomorphisms $$X \otimes (C \otimes D) \simeq (X \otimes C) \otimes D$$
which express the idea that $\calD$ may be regarded as equipped with an ``action'' of $\calC$. Here we imagine $\calC$ as a kind of generalized monoid (via its monoidal structure).

Dually, if $\calC$ is right-closed, then an object of $\calD$ which represents the functor
$$ Y \mapsto ^{C}\!\bHom_{\calD}(Y,X)$$
will be denoted by $^{C}\!X$; the object $^{C}\!X$ (if it exists) is determined up to (unique) isomorphism by a map $C \rightarrow \bHom_{\calD}(^{C}\!X,X)$. 
If this object exists for all $C \in \calC$, $X \in \calD$, then we say that
$\calD$ is {\it cotensored over $\calC$}.\index{gen}{tensored}\index{gen}{cotensored}

\begin{example}
Let $\calC$ be a right-closed monoidal category. Then $\calC$ may be regarded as enriched over itself in a natural way. It is automatically tensored over itself; it is cotensored over itself if and only if it is left-closed.
\end{example}

\subsection{Trees}

Let $\calC$ be a presentable category and $S$ a small collection of morphisms in $\calC$.
According to Remark \ref{easyprest}, the smallest weakly saturated class of morphisms $\overline{S}$
containing $S$ can be obtained from $S$ using pushouts, retracts, and transfinite composition.
It is natural to ask if the formation of retracts is necessary: that is, does the weakly saturated class of morphisms generated by $S$ coincide with the class of morphisms which are given by transfinite compositions of pushouts of morphisms of $S$? Our goal for the remainder of this section is to give an affirmative answer, at least after $S$ has been suitably enlarged (Proposition \ref{easycrust}). This result is of a somewhat technical nature, and will be needed only during our discussion of combinatorial model categories in \S \ref{combimod}.

We begin by introducing a generalization of the notion of a transfinite chain of morphisms.

\begin{definition}\index{gen}{tree}\index{gen}{$S$-tree}\index{gen}{root of an $S$-tree}
Let $\calC$ be a presentable category, and let $S$ be a collection of morphisms in $\calC$. 
An {\it $S$-tree} in $\calC$ consists of the following data:
\begin{itemize}
\item[$(1)$] An object $X \in \calC$, called the {\it root} of the $S$-tree.
\item[$(2)$] A partially ordered set $A$ which is {\it well-founded} (that is, every nonempty subset of $P$ has a minimal element).
\item[$(3)$] A diagram $A \rightarrow \calC_{X/}$, which we will denote by $\alpha \mapsto Y_{\alpha}$.
\item[$(4)$] For each $\alpha \in A$, a pushout diagram
$$ \xymatrix{ C \ar[r]^{f} \ar[d] & D \ar[d] \\
\colim_{\beta < \alpha} X_{\beta} \ar[r] & X_{\alpha}, }$$
where $f \in S$.
\end{itemize}
Let $\kappa$ be a regular cardinal. We will say that an $S$-tree in $\calC$ is {\it $\kappa$-good} if each of the objects $C$ and $D$ appearing above is $\kappa$-compact, and for each
$\alpha \in A$, the set $\{ \beta \in A: \beta < \alpha \}$ is $\kappa$-small.\index{gen}{$S$-tree!$\kappa$-good}\index{gen}{tree!$\kappa$-good}
\end{definition}

\begin{notation}
Let $\calC$ be a presentable category and $S$ a collection of morphisms in $\calC$. We will
indicate an $S$-tree by writing $\{ Y_{\alpha} \}_{\alpha \in A}$. Here the root $X \in \calC$ and the relevant pushout diagrams are understood implicitly to be part of the data.

Suppose given an $S$-tree $\{ Y_{\alpha} \}_{\alpha \in A}$, and a subset $B \subseteq A$ which is closed downwards in the following sense: if $\alpha \in B$ and $\beta \leq \alpha$, then $\beta \in B$. Then $\{ Y_{\alpha} \}_{\alpha \in B}$ is an $S$-tree. We let $Y_{B}$ denote the colimit
$\colim_{\alpha \in B} Y_{\alpha}$, formed in the category $\calC_{X/}$. In particular we have a canonical isomorhism
$Y_{\emptyset} \simeq X$. If $B = \{ \alpha \in A| \alpha \leq \beta \}$, then $Y_{B} \simeq Y_{\alpha}$. 
\end{notation}

\begin{remark}\index{gen}{$S$-tree!associated}\label{asstree}
Let $\calC$ be a presentable category, $S$ a collection of morphisms in $\calC$, and
$\{ Y_{\alpha} \}_{ \alpha in A}$ an $S$-tree in $\calC$ with root $X$. Given a map
$f: X \rightarrow X'$, we can form an {\it associated $S$-tree}
$ \{ Y_{\alpha} \coprod_{X} X' \}_{\alpha \in A}$, having root $X'$.
\end{remark}

\begin{example}
Let $\calC$ be a presentable category, $S$ a collection of morphisms in $\calC$, and
$\{ Y_{\alpha} \}_{\alpha \in A}$ an $S$-tree in $\calC$ with root $X$. If $A$ is linearly ordered, then we may identify $\{ Y_{\alpha} \}_{\alpha \in A}$ with a (possibly transfinite) sequence of morphisms belonging to $S$,
$$ X \rightarrow Y_0 \rightarrow Y_1 \rightarrow \ldots, $$
as in the statement of $(2)$ in Definition \ref{saturated}.
\end{example}

\begin{remark}\label{relci}
Let $\calC$ be a presentable category, $S$ a collection morphisms in $\calC$, and
$\{Y_{\alpha} \}_{\alpha \in A}$ an $S$-tree in $\calC$. Let $B \subseteq A$ be closed downward. 
For $\alpha \in A - B$, let $B_{\alpha} = B \cup \{ \beta \in A: \beta \leq \alpha \}$, and let
$Z_{\alpha} = Y_{B_{\alpha}}$. Then $ \{ Z_{\alpha} \}_{ \alpha \in A-B}$ is an $S$-tree in
$\calC$ with root $Y_{B}$.
\end{remark}

\begin{lemma}\label{uper}
Let $\calC$ be a presentable category and let $S$ be a collection of morphisms
in $\calC$. Let $\{ Y_{\alpha} \}_{\alpha \in A}$ be an $S$-tree in $\calC$, and let
$A'' \subseteq A' \subseteq A$ be subsets which are closed downward in $A$. Then
the induced map
$ Y_{A''} \rightarrow Y_{A'}$ belongs to the weakly saturated class of morphisms generated by $S$.
In particular, the canonical map
$Y_{\emptyset} \rightarrow Y_{A}$ belongs to the weakly saturated class of morphisms generated by $S$.
\end{lemma}

\begin{proof} Using Remarks \ref{relci} and \ref{asstree}, we can assume without loss of enerality that $A'' = \emptyset$ and $A' = A$. Using the assumption that $A$ is well-founded, we can write $A$ as the union of a transfinite sequence (downward closed) subsets $\{ B( \gamma ) \subseteq A \}_{\gamma < \beta }$ with the following property:
\begin{itemize}
\item[$(\ast)$] For each $\gamma < \beta$, the set $B(\gamma)$ is obtained from
$B'(\gamma) = \bigcup_{\gamma' < \gamma} B(\gamma')$ by adjoining a minimal element of
$\alpha_{\gamma}$ of $A - B'(\gamma)$. 
\end{itemize}
For $\gamma < \beta$, let $Z_{\gamma} = Y_{B(\gamma)}$. We now observe that
$Y_{A} \simeq \colim_{\gamma < \beta} Z_{\gamma}$, and that for each
$\gamma < \beta$ there is a pushout diagram
$$ \xymatrix{ \colim_{\alpha < \alpha_{\gamma}} Y_{\alpha} \ar[r] \ar[d] & Y_{\alpha} \ar[d] \\
\colim_{\gamma' < \gamma} Z_{\gamma'} \ar[r]^{f} & Z_{\gamma},  }$$
so that $f$ is the pushout of a morphism belonging to $S$.
\end{proof}

\begin{lemma}\label{humber2}
Let $\calC$ be a presentable category, $\kappa$ a regular cardinal, and let $S = \{ f_{s}: C_{s} \rightarrow D_{s} \}$ be a collection of morphisms in $\calC$, where each of the objects $C_{s}$ and $D_{s}$ is $\kappa$-compact. Suppose that $\{ Y_{\alpha} \}_{\alpha \in A}$ is an $S$-tree in $\calC$, indexed by a partially ordered set $(A, \leq)$. Then there exists the following:
\begin{itemize}
\item[$(1)$] A new ordering $\preceq$ on $A$, which refines $\leq$ in the following sense:
if $\alpha \preceq \beta$, then $\alpha \leq \beta$. Let $A'$ denote the partially ordered set $A$, with this new partial ordering.
\item[$(2)$] A $\kappa$-good $S$-tree $\{ Y'_{\alpha} \}_{\alpha \in A'}$, having the same root
$X$ as $\{ Y_{\alpha} \}_{\alpha \in A}$.
\item[$(3)$] A collection of maps $f_{\alpha}: Y'_{\alpha} \rightarrow Y_{\alpha}$, which form a commutative diagram
$$ \xymatrix{ Y'_{\alpha'} \ar[r] \ar[d]^{f_{\alpha'}} & Y'_{\alpha} \ar[d]^{f_{\alpha}} \\
Y_{\alpha'} \ar[r] & Y_{\alpha} }$$
when $\alpha' \preceq \alpha$.
\item[$(4)$] For every subset $B \subseteq A$ which is closed downwards with respect to $\preceq$, the induced map $f_{B}: Y'_{B} \rightarrow Y_{B}$ is an isomorphism.
\end{itemize}
\end{lemma}

\begin{proof}
Choose a transfinite sequence of downward-closed subsets $\{ A(\gamma) \subseteq A \}_{\gamma \leq \beta}$ so that the following conditions are satisfied:
\begin{itemize}
\item[$(i)$] If $\gamma' \leq \gamma \leq \beta$, then $A(\gamma') \subseteq A(\gamma)$.
\item[$(ii)$] If $\lambda \leq \beta$ is a limit ordinal (possibly zero), then
$A(\lambda) = \bigcup_{ \gamma < \lambda} A(\gamma)$.
\item[$(iii)$] If $\gamma + 1 \leq \beta$, then $A(\gamma+1) = A(\gamma) \cup \{ \alpha_{\gamma} \}$, where $\alpha_{\gamma}$ is a minimal element of $A - A(\gamma)$.
\item[$(iv)$] The subset $A(\beta)$ coincides with $A$.
\end{itemize}

We will construct a compatible family of orderings $A'(\gamma) = (A(\gamma), \preceq)$, $S$-trees $\{ Y'_{\alpha} \}_{ \alpha \in A'(\gamma)} \}$, and collections of morphisms 
$\{ Y'_{\alpha} \rightarrow Y_{\alpha} \}_{\alpha \in A(\gamma)}$ by induction on $\gamma$, so that the analogues of conditions $(1)$ through $(4)$ are satisfied. If $\gamma$ is a limit ordinal, there is nothing to do; let us assume therefore that $\gamma < \beta$ and that the data
$( A'(\gamma), \{ Y'_{\alpha} \}_{\alpha \in A'(\gamma)}, \{ f_{\alpha} \}_{\alpha \in A(\gamma)} )$ has already been constructed. Let $B = \{ \alpha \in A: \alpha < \alpha_{\gamma} \}$, so that we have a pushout diagram
$$ \xymatrix{ C \ar[r]^{f} \ar[d]^{i} & D \ar[d] \\
Y_B \ar[r] & Y_{\alpha} }$$
where $f \in S$. By the inductive hypothesis, we may identify $Y_{B}$ with $Y'_{B}$. 
Since $C$ is $\kappa$-compact, the map $i$ admits a factorization
$$C \stackrel{i'}{\rightarrow} Y'_{B'} \stackrel{i''}{\rightarrow} Y'_{B}$$
where $B'$ is $\kappa$-small. Enlarging $B'$ if necessary, we may suppose that
$B'$ is closed downwards under $\preceq$. We now extend the partial ordering
$\preceq$ to $A'(\gamma+1) = A'(\gamma) \cup \{ \alpha_{\gamma} \}$ by declaring that
$\alpha \leq \alpha_{\gamma}$ if and only if $\alpha \in B'$. We define
$Y'_{\alpha_{\gamma}}$ by forming a pushout diagram
$$ \xymatrix{ C \ar[r]^{f} \ar[d]^{i'} & D \ar[d] \\
Y'_{B'} \ar[r] & Y'_{\alpha_{\gamma}}, }$$
and we define $f_{\alpha_{\gamma}}: Y'_{\alpha_{\gamma}} \rightarrow Y_{\alpha_{\gamma}}$
to be the map induced by $i''$. It is readily verified that these data satisfy the desired conditions.
\end{proof}

\begin{lemma}\label{turkteck}
Let $\calC$ be a presentable category, $\kappa$ an uncountable regular cardinal, and 
$S$ a collection of morphisms in $\calC$. Let $\{ Y_{\alpha} \}_{\alpha \in A}$ be a $\kappa$-good $S$-tree with root $X$, and $T_A: Y_{A} \rightarrow Y_{A}$ an idempotent endomorphism of
$Y_{A}$ in the category $\calC_{X/}$. Let $B_0$ be an arbitrary $\kappa$-small subset of
$A$. Then there exists a $\kappa$-small subset $B \subseteq A$ which is downward closed and contains $B_0$ and an idempotent endomorphism $T_{B}: Y_{B} \rightarrow Y_{B}$
such that the following diagram commutes:
$$ \xymatrix{ X \ar[r] \ar[d]^{=} & Y_B \ar[d]^{T_B} \ar[r] & Y_{A} \ar[d]^{T_A} \ar[d] \\
X \ar[r] & Y_B \ar[r] & Y_A. }$$
\end{lemma}

\begin{proof}
Enlarging $B_0$ if necessary, we may assume that $B_0$ is closed downwards. 
For every pair of downward closed subsets $A'' \subseteq A' \subseteq A$, let
$i_{A'',A'}$ denote the canonical map from $Y_{A''}$ to $Y_{A'}$.
Note that because $\{ Y_{\alpha} \}_{\alpha \in A}$ is a $\kappa$-good $S$-tree, if
$A' \subseteq A$ is closed downward and $\kappa$-small, $Y_{A'}$ is $\kappa$-compact when 
viewed as an object of $\calC_{X/}$. In particular, $Y_{B_0}$ is a $\kappa$-compact object of $\calC_{X/}$. It follows that the composition
$$ Y_{B_0} \stackrel{i_{B_0,A}}{\rightarrow} Y_{A} \stackrel{T_A}{\rightarrow} Y_{A}$$
can also be factored as a composition
$$ Y_{B_0} \stackrel{T_0}{\rightarrow} Y_{B_1} \stackrel{ i_{B_1,A}}{\rightarrow} Y_A,$$
where $B_1 \subseteq A$ is closed downwards and $\kappa$-small. Enlarging $B_1$ if necessary, we may suppose that $B_1$ contains $B_0$.

We now proceed to define a sequence of $\kappa$-small, downward closed subsets
$$ B_0 \subseteq B_1 \subseteq B_2 \subseteq \ldots $$
of $A$, and maps $T_i: Y_{B_i} \rightarrow Y_{B_{i+1}}$. Suppose that $i > 0$, and that
$B_{i}$ and $T_{i-1}$ have already been constructed. By compactness again, we conclude that the composite map
$$ Y_{B_i} \stackrel{ i_{B_i, A}}{\rightarrow} Y_{A} \stackrel{T_A}{\rightarrow} Y_A $$
can be factored as
$$ Y_{B_i} \stackrel{ T_{i} }{\rightarrow} Y_{B_{i+1} } \stackrel{ i_{B_{i+1}, A}}{\rightarrow} Y_A,$$
where $B_{i+1}$ is $\kappa$-small. Enlarging $B_{i+1}$ if necessary, we may assume that
$B_{i+1}$ contains $B_{i}$ and that the following diagrams commute:
$$ \xymatrix{ Y_{B_{i-1}} \ar[r]^{T_{i-1}} \ar[d]^{i_{B_{i-1}, B_{i}}} & Y_{B_{i}} \ar[d]^{i_{B_{i}, B_{i+1}}} \\
Y_{B_i} \ar[r]^{T_i} & Y_{B_{i+1}} }$$
$$ \xymatrix{ Y_{B_{i-1}} \ar[r]^{T_{i-1}} \ar[d]^{T_{i-1}} & Y_{B_i} \ar[d]^{T_{i}} \\
Y_{B_i} \ar[r]^{i_{B_i, B_{i+1}}} & Y_{B_{i+1}}. }$$
Let $B = \bigcup B_i$; then $B$ is $\kappa$-small in virtue of our assumption that $\kappa$ is uncountable. The collection of maps $\{ T_i \}$ assemble to a map $T_{B}: Y_{B} \rightarrow Y_{B}$ with the desired properties.
\end{proof}

\begin{lemma}\label{superturk}
Let $\calC$ be a presentable category, $\kappa$ an uncountable regular cardinal, and 
$S$ a collection of morphisms in $\calC$. Let $\{ Y_{\alpha} \}_{\alpha \in A}$ be a $\kappa$-good $S$-tree with root $X$, let $B \subseteq A$ be downward closed, and suppose given a commutative diagram
$$ \xymatrix{ Y_{B} \ar[r] \ar[d]^{T_B} \ar[r] & Y_A \ar[d]^{T_A} \\
Y_{B} \ar[r] & Y_A }$$
in $\calC_{X/}$, where $T_A$ and $T_B$ are idempotent. Let $C_0 \subseteq A$ be a $\kappa$-small subset. Then there exists a downward closed $\kappa$-small subset $C \subseteq A$ containing $C_0$ and a pair idempotent maps $$T_{C}: Y_{C} \rightarrow Y_{C}$$
$$T_{B \cap C}: Y_{B \cap C} \rightarrow Y_{B \cap C}$$ such that the following diagram
commutes $($in $\calC_{X/}${}$)$:
$$ \xymatrix{ Y_{B} \ar[d]^{T_B} & Y_{B \cap C} \ar[l] \ar[r] \ar[d]^{T_{B \cap C}} & Y_{C} \ar[d]^{T_{C}} \ar[r] & Y_{A} \ar[d]^{T_A} \\
Y_{B} & Y_{B \cap C} \ar[l] \ar[r] & Y_{C} \ar[r] & Y_A.}$$
\end{lemma}

\begin{proof}
Enlarging $C_0$ if necessary, we may suppose that $C_0$ is downward closed.
We will define sequences of $\kappa$-small, downward closed subsets
$$C_0 \subseteq C_1 \subseteq \ldots \subseteq A$$
$$D_1 \subseteq D_2 \subseteq \ldots \subseteq B$$
and idempotent maps $\{ T_{C_i}: Y_{C_i} \rightarrow Y_{C_{i}} \}_{i \geq 1}$,
$\{ T_{D_i}: Y_{D_i} \rightarrow Y_{D_{i}} \}_{i \geq 1}$. Moreover, we will guarantee that the following conditions are satisfied:
\begin{itemize}
\item[$(i)$] For each $i > 0$, the set $D_i$ contains the intersection $B \cap C_{i-1}$.
\item[$(ii)$] For each $i > 0$, the set $C_{i}$ contains $D_i$.
\item[$(iii)$] For each $i > 0$, the diagrams 
$$ \xymatrix{ Y_{D_{i}} \ar[r] \ar[d]^{T_{D_i}} & Y_{B} \ar[d] & Y_{C_{i}} \ar[d]^{T_{C_i}} \ar[r] & Y_{A} \ar[d]^{T_A} \\
Y_{D_i} \ar[r] & Y_{B} & Y_{C_i} \ar[r] & Y_A}$$
are commutative.
\item[$(iv)$] For each $i > 2$, the diagrams
$$ \xymatrix{ Y_{D_{i-2}} \ar[d]^{T_{D_{i-2}}} \ar[r] & Y_{D_{i-1}} \ar[d]^{T_{D_{i-2}}} & Y_{C_{i-2}} \ar[d]^{T_{C_{i-2}}} \ar[r] & Y_{C_{i-1}} \ar[d]^{T_{C_{i-1}}} \\
Y_{D_{i-2}} \ar[d] & Y_{D_{i-1}} \ar[d] & Y_{C_{i-2}} \ar[d] & Y_{C_{i-1}} \ar[d] \\
Y_{D_{i}} \ar[r]^{=} & Y_{D_{i}} & Y_{C_{i}} \ar[r]^{=} & Y_{C_{i}} }$$
commute.
\item[$(v)$] For each $i > 1$, the diagram
$$ \xymatrix{ Y_{D_{i-1}} \ar[r] \ar[d]^{T_{D_{i-1}}} & Y_{C_{i-1}} \ar[d]^{T_{C_{i-1}}} \\
Y_{D_{i-1}} \ar[d] & Y_{C_{i-1}} \ar[d] \\
Y_{D_{i}} \ar[r] & Y_{C_{i}}}$$ is commutative.
\end{itemize}
The construction goes by induction on $i$. Using a compactness argument, we see that conditions $(iv)$ and $(v)$ are satisfied provided that we choose $C_i$ and $D_i$ to be sufficiently large. The existence of the desired idempotent maps satisfying $(iii)$ then follows from Lemma \ref{turkteck}, applied to the roots $\{ Y_{\alpha} \}_{\alpha \in A}$ and $\{ Y_{\alpha} \}_{\alpha \in B}$.
We now take $C = \bigcup C_{i}$. Conditions $(i)$ and $(ii)$ guarantee that
$B \cap C = \bigcup D_{i}$. Using $(iv)$, it follows that the maps $\{ T_{C_{i}} \}$ and
$\{ T_{D_{i}} \}$ glue to give idempotent endomorphisms $T_{C}: Y_{C} \rightarrow Y_{C}$,
$T_{B \cap C}: Y_{B \cap C} \rightarrow Y_{B \cap C}$. Using $(iii)$ and $(v)$, we deduce that all of the desired diagrams are commutative.
\end{proof}

\begin{lemma}\label{tirun}
Let $\calC$ be a presentable category, $\kappa$ a regular cardinal, and suppose that $\calC$ is $\kappa$-accessible: that is, $\calC$ is generated under $\kappa$-filtered colimits by $\kappa$-compact objects $($Definition \ref{kapacc}$)$. Let $f: C \rightarrow D$ be a morphism between $\kappa$-compact object of $\calC$, let
$g: X \rightarrow Y$ be a pushout of $f$ $($so that $Y \simeq X \coprod_{C} D${}$)$, and let
$g': X' \rightarrow Y'$ be a retract of $g$ in the category of morphisms of $\calC$.
Then there exists a morphism $f': C' \rightarrow D'$ with the following properties:
\begin{itemize}
\item[$(1)$] The objects $C', D' \in \calC$ are $\kappa$-compact.
\item[$(2)$] The morphism $g'$ is a pushout of $f'$.
\item[$(3)$] The morphism $f'$ belongs to the weakly saturated class of morphisms generated by $f$.
\end{itemize}
\end{lemma}

\begin{proof}
Since $g'$ is a retract of $g$, there exists a commutative diagram
$$ \xymatrix{ X' \ar[r] \ar[d]^{g'} & X \ar[d]^{g} \ar[r] & X' \ar[d]^{g'} \\
Y' \ar[r] & Y \ar[r] & Y'. }$$
Replacing $g$ by the induced map $X' \rightarrow X' \coprod_{X} Y$, we can reduce to the case
where $X = X'$, and $Y'$ is a retract of $Y$ in $\calC_{X/}$. Then $Y'$ can be identified with
the image of some idempotent $i: Y \rightarrow Y$. 

Since $\calC$ is $\kappa$-accessible, we can write $X$ as the colimit of a $\kappa$-filtered
diagram $\{ X_{\lambda} \}$. The object $C$ is $\kappa$-compact by assumption. Refining our diagram if necessary, we may assume that it takes values in $\calC_{C/}$, and that $Y$ is
given as the colimit of the $\kappa$-filtered diagram $\{ X_{\lambda} \coprod_{C} D \}$. 

Because $D$ is $\kappa$-compact, the composition
$D \rightarrow Y \stackrel{i}{\rightarrow} Y$ admits a factorization
$$ D \stackrel{j}{\rightarrow} X_{\lambda} \coprod_{C} D \rightarrow Y.$$
The $\kappa$-compactness of $C$ implies that, after enlarging $\lambda$ if necessary, we may suppose that the composition $j \circ f$ coincides with the canonical map from $C$
to $X_{\lambda} \coprod_{C} D$. Consequently, $j$ and the $\id_{X_{\lambda}}$ determine
a map $i'$ from $Y_{\lambda} = X_{\lambda} \coprod_{C} D$ to itself. Enlarging $\lambda$ once more, we may suppose that $i'$ is idempotent, and that the diagram
$$ \xymatrix{ Y_{\lambda} \ar[d]^{i'} \ar[r] & Y \ar[d]^{i} \\
Y_{\lambda} \ar[r] & Y }$$
is commutative. Let $Y'_{\lambda}$ be the image of the idempotent $i'$, and let
$f': X_{\lambda} \rightarrow Y'_{\lambda}$ be the canonical map. Then $f'$ is a retract of the map
$X_{\lambda} \rightarrow Y_{\lambda}$, which is a pushout of $f$. This proves $(3)$. The objects
$X_{\lambda}$ and $Y'_{\lambda}$ are $\kappa$-compact by construction, so that $(1)$ is satisfied. We now observe that the diagram
$$ \xymatrix{ X_{\lambda} \ar[r] \ar[d] & Y'_{\lambda} \ar[d] \\
X \ar[r] & Y' }$$
is a retract of the pushout diagram
$$ \xymatrix{ X_{\lambda} \ar[r] \ar[d] & Y_{\lambda} \ar[d] \\
X \ar[r] & Y, }$$
and therefore itself a pushout diagram. This proves $(2)$ and completes the proof. 
\end{proof}

\begin{lemma}\label{tiura}
Let $\calC$ be a presentable category, $\kappa$ a regular cardinal such that $\calC$ is $\kappa$-accessible, and  $S = \{ f_s: C_s \rightarrow D_s \}$ a collection of morphisms $\calC$ such that each $C_s$ is $\kappa$-compact. Let $\{ Y_{\alpha} \}_{ \alpha \in A}$ be an $S$-tree in $\calC$,
with root $X$, and suppose that $A$ is $\kappa$-small. Then there exists a map
$X' \rightarrow X$, where $X$ is $\kappa$-compact, an $S$-tree
$\{ Y'_{\alpha} \}_{\alpha \in A}$ with root $X'$, and an isomorphism of $S$-trees
$$ \{ Y'_{\alpha} \coprod_{X'} X \}_{\alpha \in A} \simeq \{Y_{\alpha} \}_{\alpha \in A}$$
$($see Remark \ref{asstree}$)$. 
\end{lemma}

\begin{proof}
Since $\calC$ is $\kappa$-accessible, we can write $X$ as the colimit of 
diagram $\{ X_{i} \}_{i \in I}$ indexed by a $\kappa$-filtered partially ordered set $I$, where each $X_{i}$ is $\kappa$-compact. 
Choose a transfinite sequence of downward-closed subsets $\{ A(\gamma) \subseteq A \}_{\gamma \leq \beta}$ so that the following conditions are satisfied:
\begin{itemize}
\item[$(i)$] If $\gamma' \leq \gamma \leq \beta$, then $A(\gamma') \subseteq A(\gamma)$.
\item[$(ii)$] If $\lambda \leq \beta$ is a limit ordinal (possibly zero), then
$A(\lambda) = \bigcup_{ \gamma < \lambda} A(\gamma)$.
\item[$(iii)$] If $\gamma + 1 \leq \beta$, then $A(\gamma+1) = A(\gamma) \cup \{ \alpha_{\gamma} \}$, where $\alpha_{\gamma}$ is a minimal element of $A - A(\gamma)$.
\item[$(iv)$] The subset $A(\beta)$ coincides with $A$.
\end{itemize}
Note that, since $A$ is $\kappa$-small, we have $\beta < \kappa$.

We will construct:
\begin{itemize}
\item[$(a)$] A transfinite sequence of elements $\{ i_{\gamma} \in I \}_{\gamma \leq \beta}$, such that $i_{\gamma} \leq i_{\gamma'}$ for $\gamma \leq \gamma'$. 
\item[$(b)$] A sequence of $S$-trees $\{ Y^{\gamma}_{\alpha} \}_{\alpha \in A(\gamma)} \}$, having roots $X_{i_{\gamma}}$. 
\item[$(c)$] A collection of isomorphisms of $S$-trees
$$ \{ Y^{\gamma}_{\alpha} \coprod_{X_{i_{\gamma}}} X_{i_{\gamma'}} \}_{\alpha \in A(\gamma)}
\simeq \{ Y^{\gamma'}_{\alpha} \}_{\alpha \in A(\gamma)}$$
$$ \{ Y^{\gamma}_{\alpha} \coprod_{ X_{i_{\gamma}}} X \}_{\alpha \in A(\gamma)}
\simeq \{ Y_{\alpha} \}_{ \alpha \in A(\gamma) }$$
which are compatible with one another in the obvious sense.
\end{itemize}
If $\gamma$ is a limit ordinal (or zero), we simply choose $i_{\gamma}$ to be
any upper bound for $\{ i_{\gamma'} \}_{\gamma' < \gamma}$ in $I$. The rest of the data is uniquely determined. The existence of such an upper bound is guaranteed by our assumption that $I$ is $\kappa$-filtered, since $\gamma \leq \beta < \kappa$. Let us therefore suppose that the above data has been constructed for all ordinals $\leq \gamma$, and proceed to define $i_{\gamma+1}$.
Let $i = i_{\gamma}$, $\alpha = \alpha_{\gamma}$, and let $B = \{ \beta \in A: \beta < \alpha \}$. Then
we have canonical isomorphisms
$$Y_{B} \simeq Y^{\gamma}_{B} \coprod_{X_i} X
\simeq \colim \{ Y^{\gamma}_{B} \coprod_{X_i} X_j \}_{j \geq i},$$
and a pushout diagram
$$ \xymatrix{ C_{s} \ar[r]^{f_{s}} \ar[d]^{g} & D_{s} \ar[d] \\
Y_B \ar[r] & Y_{\alpha}. }$$
The $\kappa$-compactness of $C_{s}$ implies that $g$ factors as a composition
$$ C_{s} \stackrel{g'}{\rightarrow} Y^{\gamma}_{B} \coprod_{ X_i} X_j $$
for some $j \geq i$. We now define $i_{\gamma+1} = j$, and
$Y^{\gamma+1}_{\alpha}$ by forming a pushout diagram
$$ \xymatrix{ C_{s} \ar[d]^{g_{s}} \ar[r] & D_{s} \ar[d] \\
Y^{\gamma}_{B} \coprod_{X_i} X_j \ar[r] & Y^{\gamma+1}_{\alpha}. }$$
\end{proof}

\begin{proposition}\label{easycrust}
Let $\calC$ be presentable $\infty$-category, $\kappa$ a regular cardinal, $\overline{S}$ a weakly saturated class of morphisms in $\calC$. Let $S \subseteq \overline{S}$ be the subset consisting of those morphisms $f: X \rightarrow Y$ in $\overline{S}$ such that $X$ and $Y$ are $\kappa$-compact. Assume that:
\begin{itemize}
\item[$(i)$] The regular cardinal $\kappa$ is uncountable, and $\calC$ is $\kappa$-accessible.
\item[$(ii)$] The set $S$ generates $\overline{S}$ as a weakly saturated class of morphisms.
\end{itemize}
Then, for every morphism $f: X \rightarrow Y$ belonging to $\overline{S}$, there
exists an transfinite sequence of objects $\{ Y_{\gamma} \}_{\gamma < \beta}$ of
$\calC_{X/}$ with the following properties:
\begin{itemize}
\item[$(1)$] For every ordinal $\gamma < \beta$, the natural map
$\colim_{\gamma' < \gamma} Z_{\gamma'} \rightarrow Z_{\gamma}$ is the pushout of a morphism in  $S$.
\item[$(2)$] The colimit $\colim_{\gamma < \beta} Z_{\gamma}$ is isomorphic to $Y$
(as objects of $\calC_{X/}$). 
\end{itemize}
\end{proposition}

\begin{proof}
Remark \ref{easyprest} implies the existence of a transfinite sequence of objects
$$ Y_0 \rightarrow Y_1 \rightarrow \ldots $$
in $\calC_{X/}$ indexed by a set of ordinals $A = \{ \alpha | \alpha < \lambda \}$, 
satisfying condition $(1)$, such that
$Y$ is a {\em retract} of $\colim_{\alpha < \lambda} Y_{\alpha}$ in $\calC_{X/}$. We may view
the sequence $\{ Y_{\alpha} \}_{ \alpha \in A}$ as an $S$-tree in $\calC$, having root $X$. According to Lemma \ref{humber2}, we can choose a new $S$-tree $\{ Y'_{\alpha} \}_{ \alpha \in A'}$
which is $\kappa$-good, where $Y'_{A'} \simeq Y_A$, so that $Y$ is a retract of
$Y'_{A'}$. Choose an idempotent map $T_{A'}: Y'_{A'} \rightarrow Y'_{A'}$ in $\calC_{X/}$, whose
image is isomorphic to $Y$.

We now define a transfinite sequence
$$ B(0) \subseteq B(1) \subseteq B(2) \subseteq \ldots, $$
indexed by ordinals $\gamma < \beta$, and a compatible system of idempotent maps $T_{B(\gamma)}: Y'_{B_\gamma} \rightarrow Y'_{B_{\gamma}}.$
Fix an ordinal $\gamma$, and suppose that $B(\gamma')$ and $T_{B(\gamma')}$ have been defined for $\gamma' < \gamma$. Let $B'(\gamma) = \bigcup_{ \gamma' < \gamma} B(\gamma')$, and let $T_{B'(\gamma)}$ be the result of amalgamating the maps $\{ T_{B(\gamma')} \}_{\gamma' < \gamma}$. If $B'(\gamma) = A'$, we set $\beta = \gamma$ and conclude the construction; otherwise, choose a minimal element $a \in A' - B'(\gamma)$. Applying Lemma \ref{superturk}, we deduce the existence of a downward closed subset $C(\gamma) \subseteq A'$, and a compatible collection of idempotent maps 
$$ T_{C(\gamma)} : Y'_{C(\gamma)} \rightarrow Y'_{C(\gamma)} $$
$$ T_{C(\gamma) \cap B'(\gamma)}: Y'_{C(\gamma) \cap B'(\gamma)} \rightarrow Y'_{C(\gamma) \cap B'(\gamma)}.$$
We then define $B(\gamma) = B'(\gamma) \cup C(\gamma)$, and
$T_{B(\gamma)}$ to be the result of amalgamating $T_{B'(\gamma)}$ and $T_{C(\gamma)}$. 

We observe that, for every ordinal $\gamma$, there is a $\kappa$-good $S$-tree 
$\{ Y''_{\alpha} \}_{ \alpha \in B(\gamma) - B'(\gamma)}$
with root $Y'_{B(\gamma)}$, such that $Y''_{B(\gamma)-B'(\gamma)} \simeq Y'_{B(\gamma)}$ (Remark \ref{relci}). Combining Lemma \ref{tiura} with the observation that 
$B(\gamma)-B'(\gamma)$ is $\kappa$-small, we deduce that the map
$$ Y'_{B'(\gamma)} \rightarrow Y'_{B(\gamma)}$$ is the pushout of a morphism in $S$.

For each ordinal $\gamma < \beta$, let $Z_{\gamma}$ denote the image of the idempotent map
$T_{B(\gamma)}$. Then $\colim_{\gamma < \beta} Z_{\gamma} \simeq Y$, so that $(2)$ is satisfied. Condition $(1)$ follows from Lemma \ref{tirun}.
\end{proof}

\begin{corollary}\label{unitape}
Under the hypotheses of Proposition \ref{easycrust}, there exists a $\kappa$-good
$S$-tree $\{ Y_{\alpha} \}_{\alpha \in A}$ such that $Y_{A} \simeq Y$ in $\calC_{X/}$.
\end{corollary}

\begin{proof}
Combine Proposition \ref{easycrust} with Lemma \ref{humber2}.
\end{proof}

\section{Model Categories}\label{appmodelcat}

One of the oldest and most successful approaches to the study of
$\infty$-categorical phenomena is Quillen's theory of model
categories. In this book, Quillen's theory will play two (related) roles:

\begin{itemize}
\item[$(1)$] The structures that we use to describe higher categories are naturally organized into model categories. For example, $\infty$-categories are precisely those simplicial sets which are fibrant with respect to the Joyal model structure (Theorem \ref{joyalcharacterization}).
The theory of model categories provides a convenient framework for phrasing certain results and for comparing different models of higher category theory (see, for example, \S \ref{compp3}).

\item[$(2)$] The theory of model categories can itself be regarded as an approach to higher
category theory. If $\bfA$ is a simplicial model category, then the subcategory $\bfA^{\degree} \subseteq \bfA$ of fibrant-cofibrant objects forms a fibrant simplicial category. Proposition \ref{toothy} implies that the simplicial nerve $\sNerve(\bfA^{\degree})$ is an $\infty$-category. We will refer to $\sNerve(\bfA^{\degree})$ as the {\it underlying $\infty$-category} of\index{not}{Adegree@$\bfA^{\degree}$}
$\bfA$. Of course, not every $\infty$-category arises in this way, even up to equivalence: for example, the existence of homotopy limits and homotopy colimits in $\bfA$ implies the existence of various limits and colimits in $\sNerve(\bfA^{\degree})$ (Corollary \ref{limitsinmodel}). Nevertheless, we can often use the theory of model categories to prove theorems about general $\infty$-categories, by reducing to the situation of $\infty$-categories which arise via the above construction (every
$\infty$-category $\calC$ admits a fully faithful embedding into $\sNerve( \bfA^{\degree})$, for
an appropriately chosen simplicial model category $\bfA$). For example, our proof of the $\infty$-categorical Yoneda lemma (Proposition \ref{fulfaith}) uses this strategy.
\end{itemize}


The purpose of this section is to review the theory of model categories, with an eye towards the sort of applications described above. Our exposition is somewhat terse and we will omit many proofs.
For a more detailed account, we refer the reader to \cite{hovey} (or any other text on the theory of model categories).

\subsection{The Model Category Axioms}

\begin{definition}\label{modelcatdef}\index{gen}{model category}\index{gen}{category!model}
A {\it model category} is a category $\calC$ which is equipped with three distinguished classes of morphisms in $\calC$, called {\it cofibrations, fibrations,} and {\it weak equivalences}, in which the following axioms are satisfied:
\begin{itemize}
\item[$(1)$] The category $\calC$ admits (small) limits and colimits.
\item[$(2)$] Given a composable pair of maps $X \stackrel{f}{\rightarrow} Y \stackrel{g}{\rightarrow} Z$, if any two of $g \circ f$, $f$, and $g$ are weak equivalences, then so is the third.
\item[$(3)$] Suppose $f: X \rightarrow Y$ is a retract of $g: X' \rightarrow Y'$: that is, suppose
there exists a commutative diagram
$$ \xymatrix{ X \ar[r]^{i} \ar[d]^{f} & X' \ar[d]^{g} \ar[r]^{r} & X \ar[d]^{f} \\
Y \ar[r]^{i'} & Y' \ar[r]^{r'} & Y } $$
where $r \circ i = \id_X$ and $r' \circ i' = \id_Y$. Then
\begin{itemize}
\item[$(i)$] If $g$ is a fibration, so is $f$.\index{gen}{fibration}
\item[$(ii)$] If $g$ is a cofibration, then so is $f$.\index{gen}{cofibration}
\item[$(iii)$] If $g$ is a weak equivalence, then so is $f$.\index{gen}{weak equivalence}
\end{itemize}

\item[$(4)$] Given a diagram
$$ \xymatrix{ A \ar[d]^{i} \ar[r] & X \ar[d]^{p} \\
B \ar[r] \ar@{-->}[ur] & Y,}$$
a dotted arrow can be found rendering the diagram commutative if either
\begin{itemize}
\item[$(i)$] The map $i$ is a cofibration, and the map $p$ is both a fibration and a weak equivalence.
\item[$(ii)$] The map $i$ is both a cofibration and a weak equivalence, and the map $p$ is a fibration.
\end{itemize}
\item[$(5)$] Any map $X \rightarrow Z$ in $\calC$ admits factorizations
$$ X \stackrel{f}{\rightarrow} Y \stackrel{g}{\rightarrow} Z$$
$$ X \stackrel{f'}{\rightarrow} Y' \stackrel{g'}{\rightarrow} Z$$
where $f$ is a cofibration, $g$ is a fibration and a weak equivalence, $f'$ is a cofibration and a weak equivalence, and $g'$ is a fibration.
\end{itemize}
\end{definition}

A map $f$ in a model category $\calC$ is called a {\it trivial cofibration} if it is both a cofibration and a weak equivalence; similarly $f$ is called a {\it trivial fibration} if it is both a fibration and a weak equivalence. By axiom $(1)$, any model category $\calC$ has an initial object $\emptyset$
and a final object $\ast$. An object $X \in \calC$ is said to be {\it fibrant} if the unique map
$X \rightarrow \ast$ is a fibration, and {\it cofibrant} if the unique map $\emptyset \rightarrow X$ is a cofibration.\index{gen}{fibrant}\index{gen}{cofibrant}\index{gen}{trivial!fibration}\index{gen}{trivial!cofibration}

\begin{example}\label{trivmodel}
Let $\calC$ be any category which admits small limits and colimits. Then
$\calC$ can be endowed with the {\em trivial} model structure:
\begin{itemize}
\item[$(W)$] The weak equivalences in $\calC$ are the isomorphisms.
\item[$(C)$] Every morphism in $\calC$ is a cofibration.
\item[$(F)$] Every morphism in $\calC$ is a fibration. 
\end{itemize}
\end{example}

\subsection{The Homotopy Category of a Model Category}

Let $\calC$ be a model category containing an object $X$. A {\it cylinder object} for $X$\index{gen}{cylinder object}\index{gen}{object!cylinder}
is an object $C$ together with a diagram
$ X \coprod X \stackrel{i}{\rightarrow} C \stackrel{j}{\rightarrow} X$
where $i$ is a cofibration and $j$ is a weak equivalence, and the composition $j \circ i$
is the ``fold map'' $X \coprod X \rightarrow X$. 
Dually, a {\it path object} for $Y \in \calC$ is an object $P$ together with a diagram $$Y \stackrel{q}{\rightarrow} P \stackrel{p}{\rightarrow} Y \times Y$$ such that $q$ is a weak equivalence, $p$ is a fibration, and $p \circ q$ is the diagonal map $Y \rightarrow Y \times Y$.
The existence of cylinder and path objects follows from the factorization axiom $(5)$ of Definition \ref{modelcatdef} (factor the ``fold map'' $X \coprod X \rightarrow X$ as a cofibration followed by a trivial fibration and the diagonal map $Y \rightarrow Y \times Y$ as a trivial cofibration followed by a fibration).\index{gen}{path object}\index{gen}{object!path}

\begin{proposition}\label{homotopy}
Let $\calC$ be a model category. Let $X$ be a cofibrant object of $\calC$, $Y$ a fibrant object of $\calC$, and $f,g: X \rightarrow Y$ two maps. The following conditions are equivalent:
\begin{itemize}
\item[$(1)$] For every cylinder object $X \coprod X \stackrel{j}{\rightarrow} C$, there exists a commutative diagram
$$ \xymatrix{ X \coprod X \ar[rr]^{j} \ar[dr]^{(f,g)} & & C \ar[dl] \\
 & Y}$$

\item[$(2)$] There exists a cylinder object $X \coprod X \stackrel{j}{\rightarrow} C$ and a commutative diagram
$$ \xymatrix{ X \coprod X \ar[rr]^{j} \ar[dr]^{(f,g)} & & C \ar[dl] \\
 & Y}$$

\item[$(3)$] For every path object $P \stackrel{p}{\rightarrow} Y \times Y$, there exists a commutative diagram $$ \xymatrix{ X  \ar[rr] \ar[dr]^{(f,g)} & & P \ar[dl]^{p} \\
 & Y \times Y}$$

\item[$(4)$] There exists a path object $P \stackrel{p}{\rightarrow} Y \times Y$ and a commutative diagram $$ \xymatrix{ X  \ar[rr] \ar[dr]^{(f,g)} & & P \ar[dl]^{p} \\
 & Y \times Y}$$
\end{itemize}
\end{proposition}

If $\calC$ is a model category containing a cofibrant object $X$ and a fibrant object $Y$, we say two maps $f,g: X \rightarrow Y$ are {\it homotopic} if the hypotheses of Proposition \ref{homotopy} are satisfied, and write $f \simeq g$. The relation $\simeq$ is an equivalence relation on\index{gen}{homotopy!between morphisms in a model category} $\Hom_{\calC}(X,Y)$. The {\it homotopy category} $\h{\calC}$ may be defined as follows:\index{gen}{homotopy category!of a model category}\index{gen}{category!homotopy}\index{not}{hcalC@$\h{\calC}$}

\begin{itemize}
\item The objects of $\h{\calC}$ are the fibrant-cofibrant objects of $\calC$.
\item For $X,Y \in \h{\calC}$, the set $\Hom_{\h{\calC}}(X,Y)$ is the set of $\simeq$-equivalence classes
of $\Hom_{\calC}(X,Y)$.
\end{itemize}

Composition is well-defined in $\h{\calC}$, in virtue of the fact that if $f \simeq g$, then
$f \circ h \simeq g \circ h$ (this is clear from characterization $(2)$ of Proposition \ref{homotopy}) and $h' \circ f \simeq h' \circ g$ (this is clear from characterization $(4)$ of Proposition \ref{homotopy}), for any maps $h,h'$ such that the compositions are defined in $\calC$.

There is another way of defining $\h{\calC}$ (or at least, a category equivalent to $\h{\calC}$): one begins with all of $\calC$ and formally adjoins inverses to all weak equivalences. Let $H(\calC)$ denote the category so-obtained. If $X \in \calC$ is cofibrant and $Y \in \calC$ is fibrant, then homotopic maps $f,g: X \rightarrow Y$ have the same image in $H(\calC)$; consequently we obtain a functor $\h{\calC} \rightarrow H(\calC)$ which can be shown to be an equivalence. We will generally ignore the distinction between these two categories, employing whichever description is more useful for the problem at hand.

\begin{remark}
Since $\calC$ is (generally) not a small category, it is not immediately clear that $H(\calC)$ has small morphism sets; however, this follows from the equivalence between $H(\calC)$ and $\h{\calC}$.
\end{remark}

\subsection{A Lifting Criterion}

The following basic principle will be used many times throughout this book:

\begin{proposition}\label{princex}
Let $\calC$ be model category containing cofibrant objects $A$ and $B$, and a fibrant object $X$.
Suppose given a cofibration $i: A \rightarrow B$ and any map $f: A \rightarrow X$. Suppose
moreover that there exists a commutative diagram
$$ \xymatrix{ A \ar[dd]^{[i]} \ar[dr]^{[f]} & \\
& X \\
B \ar[ur]^{\overline{g}} } $$
in the homotopy category $h \calC$. Then there exists a commutative diagram
$$ \xymatrix{ A \ar[dd]^{i} \ar[dr]^{f} & \\
& X \\
B \ar[ur]^{g} } $$ in $\calC$, with $[g] = \overline{g}$. $($Here we let $[p]$ denote the homotopy class in $\h{ \calC}$ of a morphism $p$ in $\calC$.$)$ 
\end{proposition}

\begin{proof}
Choose a map $g': B \rightarrow X$ representing the homotopy class $\overline{g}$.
Choose a cylinder object
$$A \coprod A \rightarrow C(A) \rightarrow A,$$
and a factorization
$$ C(A) \coprod_{A \coprod A} (B \coprod B) \rightarrow C(B) \rightarrow B$$
where the first map is a cofibration and the second a trivial fibration. We observe that
$C(B)$ is a cylinder object for $B$.

Since $g' \circ i$ is homotopic to $f$, there exists a map $h_0: C(A) \coprod_{A} B \rightarrow X$
with $h|B = g'$ and $h|A = f$. The inclusion $C(A) \coprod_{A} B \rightarrow C(B)$ is a trivial cofibration, so $h_0$ extends to a map $h: C(B) \rightarrow X$. We may regard $h$ as a homotopy from $g'$ to $g$, where $g \circ i = f$.
\end{proof}

Proposition \ref{princex} will often be applied in the following way. Suppose given a diagram
$$ \xymatrix{ A' \ar[r] \ar[dd] & A \ar[dd]^{i} \ar[dr]^{f} & \\
& & X \\
B' \ar[r] & B \ar@{-->}[ur] } $$
which we would like to extend as indicated by the dotted arrow. If $X$ is fibrant, $i$ is a cofibration between cofibrant objects, and the horizontal arrows are weak equivalences, then it suffices
to solve the (frequently easier) problem of constructing the dotted arrow in the diagram
$$ \xymatrix{ A' \ar[dd] \ar[drr] & \\
& & X \\
B' \ar@{-->}[urr] }.$$

\subsection{Left Properness and Homotopy Pushout Squares}\label{hopush}

\begin{definition}
A model category $\calC$ is {\it left proper} if, for any pushout square\index{gen}{left proper}\index{gen}{model category!left proper}
$$ \xymatrix{ A \ar[r]^{i} \ar[d]^{j} & B \ar[d]^{j'} \\
A' \ar[r]^{i'} & B'}$$
in which $i$ is a cofibration and $j$ is a weak equivalence, the map $j'$ is also a weak equivalence. Dually, $\calC$ is {\it right proper} if, for any pullback square\index{gen}{model category!right proper}\index{gen}{right proper}
$$ \xymatrix{ X' \ar[r]^{p'} \ar[d]^{q'} & Y' \ar[d]^{q} \\
X \ar[r]^{q} & Y}$$
in which $p$ is a fibration and $q$ is a weak equivalence, the map $q'$ is also a weak equivalence.
\end{definition}

In this book, we will deal almost exclusively with left proper model categories. The following provides a useful criterion for establishing left-properness.

\begin{proposition}\label{propob}
Let $\calC$ be a model category in which every object is cofibrant. Then $\calC$ is left proper.
\end{proposition}

Proposition \ref{propob} is an immediate consequence of the following basic lemma:

\begin{lemma}
Let $\calC$ be a model category containing a pushout diagram
$$ \xymatrix{ A \ar[r]^{i} \ar[d]^{j} & B \ar[d]^{j'} \\
A' \ar[r]^{i'} & B'.}$$
Suppose that $A$ and $A'$ are cofibrant, $i$ is a cofibration, and $j$ is a weak equivalence.
Then $j'$ is a weak equivalence.
\end{lemma}

\begin{proof}
We wish to show that $j'$ is an isomorphism in the homotopy category $h\calC$. In other words, we need to show that for every fibrant object
$Z$ of $\calC$, composition with $j'$ induces a bijection $\Hom_{\h{\calC}}(B',Z) \rightarrow \Hom_{\h{ \calC}}(B,Z)$.

We first show that composition with $j'$ is surjective on homotopy classes. Suppose given
a map $f: B \rightarrow Z$. Since $j$ is a weak equivalence, the composition $f \circ i$ is homotopic to $g \circ j$, for some $g: A' \rightarrow B$. According to Proposition \ref{princex}, there
is a map $f': B \rightarrow Z$ such that $f' \circ i = g \circ j$, and such that $f'$ is homotopic to $f$. The amalgamation of $f'$ and $g$ determines a map $B' \rightarrow Z$ which lifts $f'$.

We now show that $j'$ is injective on homotopy classes. Suppose given a pair of maps
$s,s': B' \rightarrow Z$. Let $P$ be a path object for $Z$. If $s \circ j'$ and $s' \circ j'$ are homotopic, then there exists a commutative diagram
$$ \xymatrix{ B \ar[r]^{h} \ar[d]^{j'} & P \ar[d] \\
B' \ar[r]^{s \times s'} & Z \times Z.}$$
We now replace $\calC$ by $\calC_{/Z \times Z}$ and apply the surjectivity statement above
to deduce that there is a map $h': B' \rightarrow P$ such that $h$ is homotopic to $h' \circ j'$. The existence of $h'$ shows that $s$ and $s'$ are homotopic, as desired.
\end{proof}

Suppose given a diagram
$$ A_0 \leftarrow A \rightarrow A_1$$ in a model category $\calC$. In general, the pushout
$ A_0 \coprod_A A_1$ is poorly behaved, in the sense that a map of diagrams
$$ \xymatrix{ A_0 \ar[d] & A \ar[l] \ar[r] \ar[d] & A_1 \ar[d]\\
B_0 & B \ar[l] \ar[r] & B_1 }$$
need not induce a weak equivalence $A_0 \coprod_A A_1 \rightarrow B_0 \coprod_B B_1$, even if
each of the vertical arrows in the diagram is individually a weak equivalence. To correct this difficulty, it is convenient to introduce the left-derived functor of ``pushout''. The {\it homotopy pushout}\index{gen}{homotopy pushout}\index{gen}{pushout!homotopy} of the diagram
$$ \xymatrix{ A_0 & A \ar[l] \ar[r] & A_1 } $$
is defined to be the pushout $A'_0 \coprod_{ A' } A'_1$, where we have chosen a commutative diagram
$$ \xymatrix{ A'_0 \ar[d] & A' \ar[d] \ar[r]^{i} \ar[l]_{j} & A'_1 \ar[d]\\
A_0 & A \ar[l] \ar[r] & A_1 } $$
in which the top row is a {\em cofibrant} diagram, in the sense that $A'$ is cofibrant and the maps
$i$ and $j$ are both cofibrations. One can show that such a diagram exists, and that the pushout $A_0' \coprod_{A'} A_1'$ depends on the choice of diagram only up to weak equivalence. (For a more systematic approach which includes a definition of ``cofibrant'' for more complicated diagrams, we refer the reader to \S \ref{quasilimit3}.) 

More generally, we will say that a diagram
$$ \xymatrix{ & A \ar[dr] \ar[dl] & \\
A_0 \ar[dr] & & A_1 \ar[dl] \\
& M }$$
is a {\it homotopy pushout square} if the composite map
$$ A'_0 \coprod_{A'} A'_1 \rightarrow A_0 \coprod_{A} A_1 \rightarrow M$$
is a weak equivalence. In this case we will also say that $M$ is a {\it homotopy pushout} of
$A_0$ and $A_1$ over $A$. One can show that this condition is independent of the choice of
``cofibrant resolution'' $$ \xymatrix{ A'_0 & A' \ar[l] \ar[r] & A'_1}$$ of the original diagram.
In particular, we note that if the diagram
$$ \xymatrix{ A_0 & A \ar[r] \ar[l] & A_1 }$$
is {\em already} cofibrant, then the ordinary pushout $A_0 \coprod_A A_1$ is a homotopy pushout. However, the condition that the diagram be cofibrant is quite strong; in good situations we can get away with quite a bit less:

\begin{proposition}\label{leftpropsquare}
Let $\calC$ be a model category, and let
$$ \xymatrix{ & A \ar[dl]^{i} \ar[dr]^{j} & \\
A_0 \ar[dr] & & A_1 \ar[dl] \\
& A_0 \coprod_A A_1 } $$
be a pushout square in $\calC$. This diagram is also a homotopy pushout square if
either of the following conditions is satisfied:
\begin{itemize}
\item[$(i)$] The objects $A$ and $A_0$ are cofibrant, and $j$ is a cofibration.
\item[$(ii)$] The map $j$ is a cofibration, and $\calC$ is left proper.
\end{itemize}
\end{proposition}

\begin{remark}
The above discussion of homotopy pushouts can be dualized; one obtains the notion of {\it homotopy pullbacks}, and the analogue of Proposition \ref{leftpropsquare} requires either that $\calC$ be a {\em right} proper model category or that the objects in the diagram be fibrant.\index{gen}{homotopy pullback}\index{gen}{pullback!homotopy}
\end{remark}

\subsection{Quillen Adjunctions and Quillen Equivalences}\label{quilladj}

Let $\calC$ and $\calD$ be model categories, and suppose given a pair of adjoint functors
$$ \Adjoint{F}{\calC}{\calD}{G} $$
(here $F$ is the left adjoint and $G$ is the right adjoint). The following conditions are equivalent:

\begin{itemize}
\item[$(1)$] The functor $F$ preserves cofibrations and trivial cofibrations.

\item[$(2)$] The functor $G$ preserves fibrations and trivial fibrations.

\item[$(3)$] The functor $F$ preserves cofibrations and the functor $G$ preserves fibrations.

\item[$(4)$] The functor $F$ preserves trivial cofibrations and the functor $G$ preserves trivial fibrations.
\end{itemize}

If any of these equivalent conditions is satisfied, then we say that the pair $(F,G)$ is a {\it Quillen adjunction} between $\calC$ and $\calD$. We also say that $F$ is a {\it left Quillen functor} and that $G$ is a {\it right Quillen functor}. In this case, one can show that $F$ preserves weak equivalences between cofibrant objects, and $G$ preserves weak equivalences between fibrant objects.\index{gen}{adjunction!Quillen}\index{gen}{Quillen adjunction}

Suppose that $\Adjoint{F}{\calC}{\calD}{G} $ is a Quillen adjunction.
We may view the homotopy category $\h{\calC}$ as obtained from $\calC$ by first passing to the full subcategory consisting of cofibrant objects, and then inverting all weak equivalences. Applying a similar procedure with $\calD$, we see that because $F$ preserves weak equivalence between cofibrant objects, it induces a functor $\h{\calC} \rightarrow \h{\calD}$; this functor is called the {\it left derived functor of $F$} and denoted $LF$.\index{gen}{left derived functor}\index{gen}{derived functor!left} Similarly, one may define the {\it right derived functor} $RG$ of $G$\index{gen}{right derived functor}\index{gen}{derived functor!right}. One can show that $LF$ and $RG$ determine an adjunction between the homotopy categories $\h{\calC}$ and $\h{\calD}$. 

\begin{proposition}\label{quilleq}
Let $\calC$ and $\calD$ be model categories, and let
$$ \Adjoint{F}{\calC}{\calD}{G} $$
be a Quillen adjunction. The following
are equivalent:
\begin{itemize}
\item[$(1)$] The left derived functor $LF: \h{\calC} \rightarrow \h{\calD}$ is an equivalence of categories.
\item[$(2)$] The right derived functor $RG: \h{\calD} \rightarrow \h{\calC}$ is an equivalence of categories.
\item[$(3)$] For every cofibrant object $C \in \calC$ and every fibrant object $D \in \calD$, a map
$C \rightarrow G(D)$ is a weak equivalence in $\calC$ if and only if the adjoint map $F(C) \rightarrow D$ is a weak equivalence in $\calD$.
\end{itemize}
\end{proposition}

\begin{proof}
Since the derived functors $LF$ and $RG$ are adjoint to one another, it is clear that $(1)$ is equivalent to $(2)$. Moreover, $(1)$ and $(2)$ are equivalent to the assertion that the unit and counit of the adjunction
$$ u: \id_{\calC} \rightarrow RG \circ LF$$
$$ v: LF \circ RG \rightarrow \id_{\calD}$$
are weak equivalences. Let us consider the unit $u$. Choose a fibrant object $C$ of $\calC$.
The composite functor $(RG \circ LF)(C)$ is defined to be $G(D)$, where $F(C) \rightarrow D$ is a weak equivalence in $\calD$, and $D$ is a fibrant object of $\calD$. Thus, $u$ is a weak equivalence when evaluated on $C$ if and only if for any weak equivalence
$F(C) \rightarrow D$, the adjoint map $C \rightarrow G(D)$ is a weak equivalence. Similarly,
the counit $v$ is a weak equivalence if and only if the converse holds. Thus $(1)$ and $(2)$ are equivalent to $(3)$.
\end{proof}

If the equivalent conditions of Proposition \ref{quilleq} are satisfied, then we say that the adjunction $(F,G)$ gives a {\it Quillen equivalence} between the model categories $\calC$ and $\calD$.\index{gen}{Quillen equivalence}

\subsection{Combinatorial Model Categories}\label{combimod}

In this section, we give an overview of Jeff Smith's theory of {\it combinatorial model categories}. Our main goal is to prove Proposition \ref{goot}, which allows us to construct model structures on a  
category $\calC$ by specifying the class of weak equivalences, together with a small amount of additional data.

\begin{definition}[Smith]\index{gen}{combinatorial model category}\index{gen}{model category!combinatorial}\label{sittu}
Let $\bfA$ be model category. We say that $\bfA$ is {\it combinatorial} if the following conditions are satisfied:
\begin{itemize}
\item[$(1)$] The category $\bfA$ is presentable.
\item[$(2)$] There exists a set $I$ of {\it generating cofibrations}, such that the collection of all cofibrations in $\bfA$ is the smallest weakly saturated class of morphisms containing $I$ (see Definition \ref{saturated}). 
\item[$(3)$] There exists a set $J$ of {\it generating trivial cofibrations}, such that the collection of all trivial cofibrations in $\bfA$ is the smallest weakly saturated class of morphisms containing $J$.
\end{itemize}
\end{definition}

If $\calC$ is a combinatorial model category, then the model structure on $\calC$ is uniquely determined by the generating cofibrations and generating trivial cofibrations. However, in practice these generators might be difficult to find. Our goal in this section is to reformulate Definition \ref{sittu} in a manner which puts more emphasis on the category of weak equivalences in $\bfA$. 

In practice, it is often easier to describe the class of {\em all} weak equivalences than it is to describe a class of generating trivial cofibrations. 

\begin{definition}\index{gen}{$\kappa$-accessible subcategory}\index{gen}{accessible!subcategory}\label{kappar}
Let $\calC$ be a presentable category and $\kappa$ a regular cardinal. We will say that a full subcategory $\calC_0 \subseteq \calC$ is an {\it $\kappa$-accessible subcategory} of
$\calC$ if the following conditions are satisfied:
\begin{itemize}
\item[$(1)$] The full subcategory $\calC_0 \subseteq \calC$ is stable under $\kappa$-filtered colimits.
\item[$(2)$] There exists a (small) set of objects of $\calC_0$ which generates $\calC_0$ under $\kappa$-filtered colimits. 
\end{itemize}
We will say that $\calC_0 \subseteq \calC$ is an {\it accessible subcategory} if $\calC_0$ is a $\kappa$-accessible subcategory of $\calC$, for some regular cardinal $\kappa$.
\end{definition}

Condition $(2)$ of Definition \ref{kappar} admits the following reformulation:

\begin{proposition}\label{reefa} Let $\kappa$ be a regular cardinal, let $\calC$ be a presentable category, and let $\calC_0 \subseteq \calC$ be a full subcategory which is stable under $\kappa$-filtered colimits. Then $\calC_0$ satisfies condition $(2)$ of Definition \ref{kappar} if and only if the following condition is satisfied, for all sufficiently large regular cardinals $\tau \gg \kappa$
\begin{itemize}
\item[$(2'_{\tau})$] Let $A$ be a $\tau$-filtered partially ordered set and
$\{ X_{\alpha} \}_{\alpha \in A}$ a diagram of $\tau$-compact objects of $\calC$ indexed by $A$.
For every $\kappa$-filtered subset $B \subseteq A$, we let
$X_B$ denote $(${}$\kappa$-filtered$)$ colimit of the diagram $\{ X_{\alpha} \}_{\alpha \in B}$.
Suppose that $X_{A}$ belongs to $\calC_0$. Then for every $\tau$-small subset $C \subseteq A$,
there exists a $\tau$-small, $\kappa$-filtered subset $B \subseteq A$ which contains $C$, such that
$X_B$ belongs to $\calC_0$. 
\end{itemize}
\end{proposition}

First, we need the following preliminary result:

\begin{lemma}\label{constunt}
Let $\tau \gg \kappa$ be regular cardinals such that $\tau > \kappa$, let $\calD$ be presentable $\infty$-category, let $\{ C_{a} \}_{a \in A}$ and $\{ D_{b} \}_{b \in B}$ be families of $\tau$-compact objects in $\calD$ indexed by $\tau$-filtered partially ordered sets $A$ and $B$, such that
$$ \colim_{a \in A} C_{a} \simeq \colim_{b \in B} D_{b}.$$
Then, for every pair of $\tau$-small subsets $A_0 \subseteq A$, $B_0 \subseteq B$, there
exist $\tau$-small, $\kappa$-filtered subsets $A' \subseteq A$, $B' \subseteq B$ such that
$A_0 \subseteq A'$, $B_0 \subseteq B'$, and $\colim_{a \in A'} C_a \simeq \colim_{b \in B'} D_b$.
\end{lemma}

\begin{proof}
Let $\calA$ be the partially ordered set of all $\tau$-small, $\kappa$-filtered subsets of
$A$ which contain $A_0$, let $\calB$ be the partially ordered set of all $\tau$-small, $\kappa$-filtered subsets of $B$ which contain $B_0$, let $X \in \calD$ be the common colimit
$ \colim_{a \in A} C_{a} \simeq \colim_{b \in B} D_{b}$, and let $\calC$ be the full subcategory
of $\calD_{/X}$ spanned by those morphisms $Y \rightarrow X$ where $Y$ is a $\tau$-compact object of $\calD$. Let $f: \calA \rightarrow \calC$ and $g: \calB \rightarrow \calC$ be the functors described by the formulas
$$f( A' ) = ( \colim_{a \in A'} C_{a} \rightarrow \colim_{a \in A} C_a )$$
$$g( B' ) = ( \colim_{b \in B'} D_{b} \rightarrow \colim_{b \in B} D_b ).$$
The desired result now follows by applying Lemma \ref{remuswolf} to the associated diagram
$$ \Nerve(\calA) \rightarrow \Nerve(\calC) \leftarrow \Nerve(\calB).$$
\end{proof}

\begin{proof}[Proof of Proposition \ref{reefa}]
First suppose that $(2'_{\tau})$ is satisfied for all sufficiently large $\tau \gg \kappa$. Choose
$\tau \gg \kappa$ large enough that $\calC$ is generated under colimits by its full subcategory $\calC^{\tau}$ of $\tau$-compact objects, and such that $(2'_{\tau})$ is satisfied. Let 
$\calD = \calC^{\tau} \cap \calC_0$, so that $\calD$ is essentially small. We will show that
$\calD$ generates $\calC_0$ under $\tau$-filtered colimits.
By assumption, every object $X \in \calC$ can be obtained as a $\tau$-filtered colimit of $\tau$-compact objects $\{ X_{\alpha} \}_{\alpha \in A}$.
Let $A'$ denote the collection of all
$\tau$-small, $\kappa$-filtered subsets $B \subseteq A$ such that $X_{B} \in \calC_0$. 
We regard $A'$ as partially-ordered via inclusions. Invoking condition $(2'_{\tau})$, we deduce that $X_{A}$ is the colimit of the $\tau$-filtered collection of objects $\{ X_{A'} \}_{A' \in B}$. We now observe that each $X_{A'}$ belongs to $\calD$.

Now suppose that condition $(2)$ is satisfied, so that $\calC_0$ is generated under $\kappa$-filtered colimits by a small subcategory $\calD \subseteq \calC_0$. Choose
$\tau \gg \kappa$ large enough that every object of $\calD$ is $\tau$-compact. Enlarging $\tau$ if necessary, we may suppose that $\tau > \kappa$. We claim that $(2'_{\tau})$ is satisfied.
To prove this, we consider any system of morphisms $\{ X \}_{\alpha \in A}$ satisfying the hypotheses of $(5'_{\tau})$. In particular, $X_A$ belongs to
$\calC_0$, so that $X_A$ may be obtained in some {\em other} way as a $\kappa$-filtered colimit of a system $\{ Y_{\beta}  \}_{\beta \in B}$, where each of the objects $Y_{\beta}$ belongs to $\calD$ and is therefore $\tau$-compact. Let $C'$ denote the family of all $\tau$-small, $\kappa$-filtered subsets $B_0 \subseteq B$. Replacing $B$ by $B'$ and the family
$\{ Y_{\beta} \}_{\beta \in B}$ by
$\{ Y_{B_0} \}_{ B_0 \in B'}$, we may assume that $B$ is $\tau$-filtered.

Let $A_0 \subseteq A$ be a $\tau$-small subset. Applying Lemma \ref{constunt} to the diagram category $\calC$, we deduce that $A_0 \subseteq A'$, where $A'$ is a $\tau$-small, $\kappa$-filtered subset of $A$, and there is an isomorphism
$X_{A'} \simeq Y_{B'}$; here $B'$ is a $\kappa$-filtered subset of $B$, so that $Y_{B'} \in \calC_0$ in virtue of our assumption that $\calC_0$ is stable under $\kappa$-filtered colimits.
\end{proof}

\begin{corollary}\label{sundert}
Let $f: \calC \rightarrow \calD$ be a functor between presentable categories which preserves $\kappa$-filtered colimits, and let $\calD_0 \subseteq \calD$ be a $\kappa$-accessible subcategory. Then $f^{-1} \calD_0 \subseteq \calC$ is a $\kappa$-accessible subcategory.
\end{corollary}

\begin{corollary}[Smith]\label{smitty}
Let $\bfA$ be a combinatorial model category, let $\bfA^{[1]}$ be the category of morphisms in $\bfA$, let $W \subseteq \bfA^{[1]}$ be the full subcategory spanned by the weak equivalences, and let
$F \subseteq \bfA^{[1]}$ be the full subcategory spanned by the fibrations. Then
$F$, $W$, and $F \cap W$ are accessible subcategories of $\bfA^{[1]}$. 
\end{corollary}

\begin{proof}
For every morphism $i: A \rightarrow B$, let $F_i: \bfA^{[1]} \rightarrow \Set^{[1]}$ be the functor which carries a morphism $f: X \rightarrow Y$ to the induced map of sets
$$ \Hom_{\bfA}(B, X) \rightarrow \Hom_{\bfA}(B,Y) \times_{ \Hom_{\bfA}(A,Y) } \Hom_{\bfA}(A,X).$$
We observe that if $A$ and $B$ are $\kappa$-compact objects of $\bfA$, then $F_i$ preserves $\kappa$-filtered colimits.

Let $\calC_0$ be the full subcategory of $\Set^{[1]}$ spanned by the collection of {\em surjective} maps between sets. It is easy to see that $\calC_0$ is an accessible category of $\Set^{[1]}$. It follows that the full subcategories $R(i) = F_i^{-1} \calC_0 \subseteq \bfA^{[1]}$ are accessible subcategories of $\bfA^{[1]}$ (Corollary \ref{sundert}). 

Let $I$ be a set of generating cofibrations for $\bfA$, and $J$ a set of generating trivial cofibrations. Then Proposition \ref{boundint} implies that the subcategories
$$ F = \bigcap_{j \in J} R(j)$$
$$ W \cap F = \bigcap_{i \in I} R(i)$$
are accessible subcategories of $\bfA^{[1]}$. 

Applying Proposition \ref{quillobj}, we deduce that there exists a pair of functors
$T', T'': \bfA^{[1]} \rightarrow \bfA^{[1]}$, which carry an arbitrary morphism $f: X \rightarrow Z$ to a factorization
$$ X \stackrel{ T'(f) }{\rightarrow} Y \stackrel{T''(f)}{\rightarrow} Z$$
where $T'(f)$ is a trivial cofibration, and $T''(f)$ is a fibration. Moreover, the functor $T''$ can be chosen to commute with $\kappa$-filtered colimits, for a sufficiently large regular cardinal $\kappa$. We now observe that $W$ is the inverse image of $F \cap W$ under the functor
$T'': \bfA^{[1]} \rightarrow \bfA^{[1]}$, and is therefore an accessible subcategory of
$\bfA^{[1]}$ by Corollary \ref{sundert}.
\end{proof}

Our next goal is to prove a converse to Corollary \ref{smitty}, which will allow us to construct examples of combinatorial model categories. First, we need the following preliminary result.

\begin{lemma}\label{seeva}
Let $\bfA$ be a presentable category. Suppose $W$ and $C$ are collections of morphisms
of $\bfA$ with the following properties:
\begin{itemize}
\item[$(1)$] The collection $C$ is a weakly saturated class of morphisms of $\bfA$, and there exists a $($small$)$ subset $C_0 \subseteq C$ which generates $C$ as a weakly saturated class of morphisms.
\item[$(2)$] The intersection $C \cap W$ is a weakly saturated class of morphisms of $\bfA$.
\item[$(3)$] The full subcategory $W \subseteq \bfA^{[1]}$ is an accessible subcategory
of $\bfA^{[1]}$.
\item[$(4)$] The class $W$ has the two-out-of-three property.
\end{itemize}
Then $C \cap W$ is generated, as a weakly saturated class of morphisms, by a $($small$)$ subset
$S \subseteq C \cap W$.
\end{lemma}

\begin{proof}
Let $\kappa$ be a regular cardinal such that $W$ is $\kappa$-accessible. Choose a regular cardinal $\tau \gg \kappa$ such that $W$ satisfies condition $(2'_{\tau})$ of Proposition \ref{reefa}. Enlarging $\tau$ if necessary, we may assume that $\tau > \kappa$ (so that $\tau$ is uncountable), that $\calC$ is $\tau$-accessible, and that the source and target of every morphism in $C_0$ is 
$\tau$-compact. Enlarging $C_0$ if necessary, we may suppose that $C_0$ consists of {\em all} morphisms $f: X \rightarrow Y$ in $C$ such that $X$ and $Y$ are $\tau$-compact. Let $S = C_0 \cap W$. We will show that $S$ generates $C \cap W$ as a weakly saturated class of morphisms.

Let $\overline{S}$ be the weakly saturated class of morphisms generated by $S$, and let
$f: X \rightarrow Y$ be a morphism which belongs to $C \cap W$. We wish to show that $f \in \overline{S}$. Corollary \ref{unitape} implies that there exists a $\tau$-good $C_0$-tree $\{ Y_{\alpha} \}_{\alpha \in A}$ with root $X$, such that
$Y$ is isomorphic to $Y_{A}$ as objects of $\calC_{X/}$. Let us say that a subset $B \subseteq A$ is {\it good} if it is closed downwards and the canonical map $i: X \rightarrow Y_{B}$ belongs to $W$ (we note that $i$ automatically belongs to $C$, in virtue of Lemma \ref{uper}). 

We now make the following observations:
\begin{itemize}
\item[$(i)$] Given an increasing transfinite sequence of good subsets $\{ A_{\gamma} \}_{\gamma < \beta}$,
the union $\bigcup A_{\gamma}$ is good. This follows from the assumption that $C \cap W$
is weakly saturated.
\item[$(ii)$] Let $B \subseteq A$ be good, and let $B_0 \subseteq B$ be $\tau$-small. Then there exists a $\tau$-small subset $B' \subseteq B$ containing $B_0$. This follows from our assumption
that $W$ satisfies condition $(2'_{\tau})$ of Proposition \ref{reefa}. 
\item[$(iii)$] Suppose that $B,B' \subseteq A$ are such that $B$, $B'$, and $B \cap B'$ are good.
Then $B \cup B'$ is good. To prove this, we consider the pushout diagram
$$ \xymatrix{ Y_{B \cap B'} \ar[r] \ar[d] & Y_{B} \ar[d] \\
Y_{B'} \ar[r] & Y_{B \cup B'}. }$$
Every morphism in this diagram belongs to $C$ (Lemma \ref{uper}), and the upper horizontal map belongs to $W$ in virtue of assumption $(4)$. Since $C \cap W$ is stable under pushouts, we conclude that the lower vertical map belongs to $W$. Assumption $(4)$ now implies that the composite map $X \rightarrow Y_{B'} \rightarrow Y_{B \cup B'}$ belongs to $W$, as desired.
\end{itemize}

The next step is to prove the following claim:
\begin{itemize}
\item[$(\ast)$] Let $A'$ be a good subset of $A$, and let $B_0 \subseteq A$ be $\tau$-small. Then
there exists a $\tau$-small subset $B \subseteq A$ such that $B_0 \subseteq B$, $B$ is good, and $B \cap A'$ is good.
\end{itemize}

To prove $(\ast)$, we begin by setting $B'_0 = A' \cap B_0$. We now define sequences of $\tau$-small subsets
$$ B_0 \subseteq B_1 \subseteq B_2 \subseteq \ldots $$
$$ B'_0 \subseteq B'_1 \subseteq B'_2 \subseteq \ldots$$
as follows. Suppose that $B_i$ and $B'_{i}$ have been defined. Applying $(ii)$, we choose
$B_{i+1}$ to be any $\tau$-small good subset of $A$ which contains $B_i \cup B'_{i}$.
Applying $(ii)$ again, we select $B'_{i+1}$ to be any $\tau$-small good subset of $A'$ which contains $A' \cap B_{i+1}$. Let $B = \bigcup B_{i}$. It follows from $(i)$ that $B$ and $A' \cap B = \bigcup_{i} B'_{i}$ are both good.

We now choose a transfinite sequence of good subsets $\{ A(\gamma) \subseteq A \}_{\gamma < \beta}$. Suppose that $A(\gamma')$ has been defined for $\gamma' < \gamma$, and let
$A'(\gamma) = \bigcup_{\gamma' < \gamma} A(\gamma')$. It follows from $(i)$ that
$A'(\gamma)$ is good. If $A'(\gamma) = A$, we set $\beta = \gamma$ and conclude the construction. Otherwise, choose a minimal element 
$a \in A - A'(\gamma)$. Applying $(\ast)$, we deduce that there exists a $\tau$-small good subset
$B(\gamma) \subseteq A$ containing $a$, such that $A'(\gamma) \cap B(\gamma)$ is good. Let $A(\gamma) = A'(\gamma) \cup B(\gamma)$. It follows from $(iii)$ that $A(\gamma)$ is good. 

We observe that $\{ Y_{A(\gamma)} \}_{\gamma < \beta}$ is a transfinite sequence of objects of
$\calC_{X/}$ having colimit $Y$. To prove that $f: X \rightarrow Y$ belongs to $\overline{S}$, it will suffice to show that for each $\gamma < \beta$, the map $g: Y_{A'(\gamma)} \rightarrow Y_{A(\gamma)}$ belongs to $\overline{S}$. Remark \ref{relci} implies the existence of a $C_0$-tree
$\{ Z_{\alpha} \}_{\alpha \in A(\gamma) - A'(\gamma)}$ with root $Y_{A'(\gamma)}$ and colimit $Y_{A(\gamma)}$. Since $A(\gamma) - A'(\gamma)$ is $\tau$-small, Lemma \ref{tiura} implies the existence of a pushout diagram 
$$ \xymatrix{ M \ar[r] \ar[d] & N \ar[d] \\
Y_{A'(\gamma)} \ar[r] & Y_{A(\gamma)}}$$
where $g \in C_0$. 

Since $\calC$ is $\tau$-accessible, we can write $Y_{A'(\gamma)}$ as the colimit of a family of $\tau$-compact objects $\{ Z_{\lambda} \}_{\lambda \in P}$, indexed by a $\tau$-filtered partially ordered set $P$. Since $M$ is $\tau$-compact, we can assume (reindexing the colimit if necessary) that we have a compatible family of maps $\{ M \rightarrow Z_{\lambda} \}$. For each $\lambda$, let $g_{\lambda}: Z_{\lambda} \rightarrow Z_{\lambda} \coprod_{M} N$ be the induced map. Then $g$ is the filtered colimit of the family $\{ g_{\lambda} \}_{\lambda \in P}$. Since $W$ satisfies condition
$(2'_{\tau})$ of Proposition \ref{reefa}, we conclude that there exists a $\tau$-small, $\kappa$-filtered subset $P_0 \subseteq P$, such that $g' = \colim_{\lambda \in P_0} g_{\lambda}$ belongs to $W$. We now observe that $g' \in S$, and that $g$ is a pushout of $g'$, so that $g \in \overline{S}$ as desired.
\end{proof}

\begin{proposition}\label{bigmaker}
Let $\bfA$ be a presentable category, and let $W$ and $C$ be classes of morphisms
in $\bfA$ with the following properties:
\begin{itemize}
\item[$(1)$] The collection $C$ is a weakly saturated class of morphisms of $\bfA$, and there exists a $($small$)$ subset $C_0 \subseteq C$ which generates $C$ as a weakly saturated class of morphisms.
\item[$(2)$] The intersection $C \cap W$ is a weakly saturated class of morphisms of $\bfA$.
\item[$(3)$] The full subcategory $W \subseteq \bfA^{[1]}$ is an accessible subcategory
of $\bfA^{[1]}$.
\item[$(4)$] The class $W$ has the two-out-of-three property.
\item[$(5)$] If $f$ is a morphism in $\bfA$ which has the right lifting property with respect to each element of $C$, then $f \in W$.
\end{itemize}
Then $\bfA$ admits a combinatorial model structure, which may be described as follows:
\begin{itemize}
\item[$(C)$] The cofibrations in $\bfA$ are the elements of $C$.
\item[$(W)$] The weak equivalences in $\bfA$ are the elements of $W$.
\item[$(F)$] A morphism in $\bfA$ is a fibration if it has the right lifting property with respect to every morphism in $C \cap W$.
\end{itemize}
\end{proposition}

\begin{proof}
The category $\bfA$ has all (small) limits and colimits, since it is presentable. The two-out-three property for $W$ is among our assumptions, and the stability of $W$ under retracts follows from the accessibility of $W \subseteq \bfA^{[1]}$ (Corollary \ref{swwe}). The class of cofibrations is stable under retracts by $(1)$, and the class of fibrations is stable under retracts by definition.
The classes of fibrations and cofibrations are stable under retracts by definition.

We next establish the factorization axioms. By the small object argument, any morphism
$X \rightarrow Z$ admits a factorization
$$X \stackrel{f}{\rightarrow} Y \stackrel{g}{\rightarrow} Z$$
where $f \in C$ and $g$ has the right lifting property with respect to every morphism in $C$.
In particular, $g$ has the right lifting property with respect to every morphism in $C \cap W$, so that $g$ is a fibration; assumption $(5)$ then implies that $g$ is a trivial fibration.
Similarly, using Lemma \ref{seeva} we may choose a factorization as above where $f \in C \cap W$ and $g$ has the right lifting property with respect to $C \cap W$; $g$ is then a fibration by definition.

To complete the proof, it suffices to show that cofibrations have the left lifting property with respect to trivial fibrations, and trivial cofibrations have the left lifting property with respect to fibrations. The second of these statements is clear (it is the definition of a fibration). For the first statement, let us consider an arbitrary trivial fibration $p: X \rightarrow Z$. By the small object argument,
there exists a factorization of $p$
$$ X \stackrel{q}{\rightarrow} Y \stackrel{r}{\rightarrow} Z$$
where $q$ is a cofibration, and $r$ has the right lifting property with respect to all cofibrations.
Then $r$ is a weak equivalence by $(3)$, so that $q$ is a weak equivalence by the two-out-of-three property. Considering the diagram
$$ \xymatrix{ X \ar@{=}[r] \ar[d]^{q} & X \ar[d]^{p} \\
Y \ar[r]^{r} \ar@{-->}[ur] & Z,}$$
we deduce the existence of the dotted arrow from the fact that $p$ is a fibration and $q$ is a trivial
cofibration. It follows that $p$ is a retract of $r$, and therefore $p$ also has the right lifting property with respect to all cofibrations. This completes the proof that $\bfA$ is a model category. The assertion that $\bfA$ is combinatorial follows immediately from $(1)$ and from Lemma \ref{seeva}.
\end{proof}

\begin{corollary}\label{uryt}
Let $\bfA$ be a presentable category equipped with a model structure. Suppose that there exists
a $($small$)$ set which generates the collection of cofibrations in $\bfA$ $($as a weakly saturated class of morphisms$)$. Then the following are equivalent:
\begin{itemize}
\item[$(1)$] The model category $\bfA$ is combinatorial; in other words, there exists a $($small$)$ set which generates the collection of trivial cofibrations in $\bfA$ $($as a weakly saturated class of morphisms$)$.
\item[$(2)$] The collection of weak equivalences in $\bfA$ determines an accessible subcategory
of $\bfA^{[1]}$. 
\end{itemize}
\end{corollary}

\begin{proof}
The implication $(1) \Rightarrow (2)$ follows from Corollary \ref{smitty}, and the reverse implication follows from Proposition \ref{bigmaker}.
\end{proof}

Our next goal is to prove a weaker version of Proposition \ref{bigmaker} which is somewhat easier to apply in practice.

\begin{definition}\label{perfequiv}\index{gen}{perfect!class of morphisms}
Let $\bfA$ be a presentable category. A class $W$ of morphisms in $\calC$ is {\it perfect}
if it satisfies the following conditions:

\begin{itemize}
\item[$(1)$] Every isomorphism belongs to $W$.
\item[$(2)$] Given a pair of composable morphisms $X \stackrel{f}{\rightarrow} Y \stackrel{g}{\rightarrow} Z$, if any two of the morphisms $f$, $g$, and $g \circ f$ belong to $W$, then so does the third.
\item[$(3)$] The class $W$ is stable under filtered colimits. More precisely, suppose given a family of morphisms $\{ f_{\alpha}: X_{\alpha} \rightarrow Y_{\alpha} \}$ which is indexed by a filtered partially ordered set. Let $X$ denote a colimit of $\{ X_{\alpha} \}$ and $Y$ a colimit of
$\{ Y_{\alpha} \}$, and $f: X \rightarrow Y$ the induced map. If each $f_{\alpha}$ belongs to $W$, then so does $f$. 
\item[$(4)$] There exists a (small) subset $W_0 \subseteq W$ such that every morphism belonging to $W$ can be obtained as a filtered colimit of morphisms belonging to $W_0$.
\end{itemize}
\end{definition}

\begin{example}
If $\calC$ is a presentable category, then the class $W$ consisting of all isomorphisms in $\calC$
is perfect.
\end{example}

The following is an immediate consequence of Corollary \ref{sundert}:

\begin{corollary}\label{perfpull}
Let $F: \calC \rightarrow \calC'$ be a functor between presentable categories which preserves filtered colimits, and let $W_{\calC'}$ be a perfect class of morphisms in $\calC'$. Then
$W_{\calC} = F^{-1} W_{\calC'}$ is a perfect class of morphisms in $\calC$.
\end{corollary}

\begin{proposition}\label{goot}
Let $\bfA$ be a presentable category. Suppose given a class $W$ of morphisms of $\calC$, which we will call {\it weak equivalences}, and a $($small$)$ {\em set} $C_0$ of morphisms of $\calC$, which we will call {\it generating cofibrations}. Suppose furthermore that the following assumptions are satisfied:

\begin{itemize}
\item[$(1)$] The class $W$ of weak equivalences is perfect $($Definition \ref{perfequiv}$)$.
\item[$(2)$] For any diagram
$$ \xymatrix{ X \ar[r]^{f} \ar[d] & Y \ar[d] \\
	X' \ar[r] \ar[d]^{g} & Y' \ar[d]^{g'} \\
	X'' \ar[r] & Y'' } $$
in which both squares are coCartesian, $f$ belongs to $C_0$, and $g$ belongs to $W$, 
the map $g'$ also belongs to $W$.
\item[$(3)$] If $g: X \rightarrow Y$ is a morphism in $\bfA$ which has the right lifting property with respect to every morphism in $C_0$, then $g$ belongs to $W$. 
\end{itemize}

Then there exists a left proper, combinatorial model structure on $\calC$ which may be described as follows:

\begin{itemize}
\item[$(C)$] A morphism $f: X \rightarrow Y$ in $\bfA$ is a {\it cofibration} if it belongs to the weakly saturated class of morphisms generated by $C_0$.
\item[$(W)$] A morphism $f: X \rightarrow Y$ in $\calC$ is a {\it weak equivalence} if it belongs to $W$.
\item[$(F)$] A morphism $f: X \rightarrow Y$ in $\calC$ is a {\it fibration} if it has the right lifting property with respect to every map which is both a cofibration and a weak equivalence.
\end{itemize}\index{gen}{perfect!model category}\index{gen}{model category!perfect}
\end{proposition}

\begin{proof}
We first show that the class of weak equivalences is stable under pushouts by cofibrations.
Let $P$ denote the collection of all morphisms $f$ in $\bfA$ with the following property:
for coCartesian diagram
$$ \xymatrix{ X \ar[r]^{f} \ar[d] & Y \ar[d] \\
	X' \ar[r] \ar[d]^{g} & Y' \ar[d]^{g'} \\
	X'' \ar[r] & Y'' } $$
where $g$ belongs to $W$, the map $g'$ also belongs to $W$. By assumption, 
$C_0 \subseteq P$. It is easy to see that $P$ is weakly saturated (using the stability of $W$ under filtered colimits), so that every cofibration belongs to $P$. 

It remains only to show that $\bfA$ is a model category. In view of Proposition \ref{bigmaker}, it will suffice to show that $C \cap W$ is a weakly saturated class of morphisms. It is clear that $C \cap W$ is stable under retracts. It will therefore suffice to verify the stability of $C \cap W$ under pushouts and transfinite composition. The case of transfinite composition is easy: $C$ is stable under transfinite composition because $C$ is weakly saturated, and $W$ is stable under transfinite composition because
it is stable under finite composition and filtered colimits.

It remains to show that $C \cap W$ is stable under pushouts. Suppose given a coCartesian diagram
$$ \xymatrix{ X \ar[d]^{f} \ar[r] & X'' \ar[d]^{f''} \\
Y \ar[r] & Y'' }$$
in which $f$ belongs to $C \cap W$; we wish to show that $f''$ also belongs to $C \cap W$. Since
$C$ is weakly saturated, it will suffice to show that $f''$ belongs to $W$. Using the small object argument, we can factor the top horizontal map to produce a coCartesian rectangle
$$ \xymatrix{ X \ar[d]^{f} \ar[r]^{g} & X' \ar[d]^{f'} \ar[r]^{h} & X'' \ar[d]^{f''} \\
Y \ar[r] & Y' \ar[r]^{h'} & Y'' }$$
in which $g$ is a cofibration and $h$ has the right lifting property with respect to all the morphisms in $C_0$. Since $W$ is stable under the formation of pushouts by cofibrations, we deduce that $f'$
belongs to $W$. Moreover, by assumption $(3)$, $h$ belongs to $W$. Since $h'$ is a pushout
of $h$ by the cofibration $f'$, we deduce that $h'$ belongs to $W$ as well. Applying the two-out-of-three property (twice), we deduce that $f''$ belongs to $W$. 
\end{proof}

\begin{remark} 
Let $\bfA$ be a model category. Then $\bfA$ arises via the construction of Proposition \ref{goot} if and only if it is combinatorial, left proper, and the collection of weak equivalences in $\bfA$ is stable under filtered colimits. 
\end{remark}

\subsection{Simplicial Sets}\label{simpset}

The formalism of simplicial sets plays a prominent role throughout this book. In this section, we will review the definition of a simplicial set, and establish some notation.

For each $n \geq 0$, we let
$[n]$ denote the linearly ordered set $\{ 0, \ldots, n \}$.\index{not}{[n]@$[n]$}
We let $\cDelta$ denote\index{not}{Delta@$\cDelta$}
the category of {\it combinatorial simplices}: the objects of
$\cDelta$ are the linearly ordered sets $[n]$, and morphisms
in $\cDelta$ are given by (nonstrictly) order-preserving maps.

If $\calC$ is any category, a {\it simplicial object} of $\calC$\index{gen}{simplicial object!of a category}
is a functor $\cDelta^{op} \rightarrow \calC$. Dually, a {\it
cosimplicial object} of $\calC$ is a functor $\cDelta \rightarrow
\calC$. A {\it simplicial set} is a simplicial object in the
category of sets. More explicitly, a simplicial set $S$ is
determined by the following data:\index{gen}{simplicial set}

\begin{itemize}
\item A set $S_{n}$ for each $n \geq 0$ (the value of $S$
on the object $[n] \in \cDelta$).

\item A map $p^{\ast}: S_n \rightarrow S_m$ for each
order-preserving map $[m] \rightarrow [n]$, the formation of which is compatible with composition
(including empty composition, so that
$( \id_{[n]} )^{\ast} = \id_{S_n}$).
\end{itemize}

Let recall a bit of standard notation for working with a
simplicial set $S$. For each $0 \leq j \leq n$, the {\it face map}
$d_j: S_n \rightarrow S_{n-1}$ is defined to
be the pullback $p^{\ast}$, where $p: [n-1] \rightarrow [n]$ is given by
$$p(i) =
\begin{cases} i & \text{if } i<j \\
i+1 & \text{if } i \geq j. \end{cases} $$ 
Similarly, the {\it degeneracy map} $s_j: S_{n} \rightarrow S_{n+1}$
is defined to be the pullback $q^{\ast}$, where $q: [n+1] \rightarrow [n]$ is defined by the formula
$$q(i) = \begin{cases} i & \text{if } i \leq j \\
i-1 & \text{if } i > j. \end{cases}$$\index{not}{d_i@$d_i$}\index{not}{s_i@$s_i$}\index{gen}{face map}\index{gen}{degeneracy map}
Because every order-preserving map from 
$[n]$ to $[m]$ can be factored as a composition of face and degeneracy maps, 
the structure of a simplicial set $S$ is completely determined by the sets
$S_n$ for $n \geq 0$, together with the face and degeneracy operations
defined above. These operations are required to satisfy certain identities, which we will not make explicit here.

\begin{remark}
The category $\cDelta$ is equivalent to the (larger) category of all finite, nonempty linearly ordered sets. We will sometimes abuse notation by identifying $\cDelta$ with this larger subcategory, and regarding simplicial sets (or more general simplicial objects) as functors which are defined on all nonempty linearly ordered sets.
\end{remark}

\begin{notation}
The category of simplicial sets will be denoted by $\sSet$.\index{not}{sSet@$\sSet$}
If $J$ is a linearly ordered set, we let $\Delta^J \in \sSet$ denote the
representable functor 
$[n] \mapsto \Hom( [n], J)$, where the morphisms are taken in the category of linearly ordered sets. For each $n \geq 0$, we will write
$\Delta^{n}$ in place of $\Delta^{[n]}$. We observe that, for any simplicial set $S$, there
is a natural identification of sets $S_{n} \simeq \Hom_{\sSet}(\Delta^n, S)$.\index{not}{DeltaJ@$\Delta^J$}\index{not}{Deltan@$\Delta^n$}
\end{notation}

\begin{example}
For $0 \leq j \leq n$, we let $\Lambda^n_j \subset \Delta^n$\index{not}{Lambdanj@$\Lambda^n_{k}$}
denote the ``$j$-th horn''. It is determined by the following
property: an element of $(\Lambda^n_j)_m$ is given by an
order-preserving map $p: [m] \rightarrow [n]$ satisfying the condition that 
$\{j \} \cup p( [m] ) \neq [n]$. Geometrically,
$\Lambda^n_j$ corresponds to the subset of an $n$-simplex
$\Delta^n$ in which the $j$th face and the interior have been
removed.\index{gen}{horn}

More generally, if $J$ is any linearly ordered set containing an element $j$, we let
$\Lambda^J_j$ denote the simplicial subset of $\Delta^J$ obtained by removing the interior
and the ``opposite face'' to the vertex $j$.\index{not}{LambdaJj@$\Lambda^J_{j}$}
\end{example}

The category $\sSet$ of simplicial sets has a (combinatorial, left and right proper) model structure, which we will refer to as the {\it Kan model structure}. It may be described as follows:
\index{gen}{simplicial set!Kan model structure}\index{gen}{Kan!model structure}

\begin{itemize}
\item A map of simplicial sets $f: X \rightarrow Y$ is a {\it cofibration} if it is a monomorphism; that is, if the induced map $X_n \rightarrow Y_n$ is injective for all $n \geq 0$.\index{gen}{cofibration!of simplicial sets}
\item A map of simplicial sets $f: X \rightarrow Y$ is a {\it fibration} if it is a Kan fibration; that is, if for any diagram
$$ \xymatrix{ \Lambda^n_i \ar@{^{(}->}[d] \ar[r] & X \ar[d]^{f} \\
\Delta^n \ar[r] \ar@{-->}[ur] & Y}$$
it is possible to supply the dotted arrow rendering the diagram commutative.\index{gen}{fibration!Kan}\index{gen}{Kan fibration}
\item A map of simplicial sets $f: X \rightarrow Y$ is a {\it weak equivalence} if the induced map of geometric realizations $|X| \rightarrow |Y|$ is a homotopy equivalence of topological spaces.\index{gen}{weak homotopy equivalence!of simplicial sets}
\end{itemize}

To prove this, we observe that the class of all cofibrations is generated by the collection of all inclusions $\bd \Delta^n \subseteq \Delta^n$; it is then easy to see that the conditions of
Proposition \ref{goot} are satisfied. The nontrivial point
is to verify that the fibrations for the resulting model structure are precisely the Kan fibrations, and
that $\sSet$ is right proper; these facts ultimately rely on a delicate analysis due to Quillen (see \cite{goerssjardine}).

\begin{remark}
In \S \ref{compp3}, we introduce another model structure on $\sSet$, the {\it Joyal model structure}. This model structure has the same class of cofibrations, but the fibrations and the weak equivalences differ from those defined in this section. To avoid confusion, we will refer to the fibrations and weak equivalences for the usual model structure on simplicial sets as {\it Kan fibrations} and {\it weak homotopy equivalences}, respectively.
\end{remark}


\subsection{Diagram Categories and Homotopy (Co)limits}\label{qlim7}

Let $\bfA$ be a combinatorial model category and $\calC$ a small category.
We let $\Fun(\calC, \bfA)$ denote the category of all functors from $\calC$ to $\bfA$.
In this section, we will see that $\Fun(\calC, \bfA)$ again admits the structure of a combinatorial model category: in fact, it admits two such structures. Moreover,
by considering the functoriality of this construction in the category $\calC$, we will obtain
the theory of {\it homotopy limits} and {\it homotopy colimits}.

\begin{definition}\label{injproj}\index{gen}{cofibration!projective}\index{gen}{cofibration!injective}\index{gen}{fibration!injective}\index{gen}{fibration!projective}\index{gen}{injective!cofibration}\index{gen}{injective!fibration}\index{gen}{projective!cofibration}\index{gen}{projective!fibration}\label{cooper}
Let $\calC$ be a small category, and $\bfA$ a model category.
A natural transformation $\alpha: F \rightarrow G$ in $\Fun(\calC, \bfA)$ is a:

\begin{itemize}
\item {\it injective cofibration} if the induced map $F(C) \rightarrow
G(C)$ is a cofibration in $\bfA$, for each $C \in \calC$.

\item {\it projective fibration} if the induced map $F(C) \rightarrow
G(C)$ is a fibration in $\bfA$, for each $C \in \calC$.

\item {\it weak equivalence} if the induced map $F(C) \rightarrow
G(C)$ is a weak equivalence in $\bfA$, for each $C \in \calC$.

\item {\it injective fibration} if it has the right lifting property
with respect to every morphism $\beta$ in $\Fun(\calC, \bfA)$ which is
simultaneously a weak equivalence and a injective cofibration.

\item {\it projective cofibration} if it has the left lifting property
with respect to every morphism $\beta$ in $\Fun(\calC, \bfA)$ which is
simultaneously a weak equivalence and a projective fibration.
\end{itemize}
\end{definition}

\begin{proposition}\label{smurff}\index{gen}{model category!projective}\index{gen}{model category!injective}
Let $\bfA$ be a combinatorial model category and $\calC$
be a small category. Then there exist two combinatorial model structures on $\Fun(\calC, \bfA)$:

\begin{itemize}
\item The {\it projective model structure}, determined by the strong
cofibrations, weak equivalences, and projective fibrations.

\item The {\it injective model structure}, determined by the weak
cofibrations, weak equivalences, and injective fibrations.
\end{itemize}
\end{proposition}

The following is the key step in the proof of Proposition \ref{smurf}:

\begin{lemma}\label{mainerthyme}
Let $\bfA$ be a presentable category, and $\calC$ a small category. Let $S_0$ be a (small) set of morphisms of $\bfA$, and let $\overline{S}_0$ be the weakly saturated class of morphisms generated by $S_0$. Let $\widetilde{S}$ be the collection of all morphisms $F \rightarrow G$ in $\Fun(\calC, \bfA)$ with the following property: for every $C \in \calC$, the map $F(C) \rightarrow G(C)$ belongs to $\overline{S}_0$. Then there exists a $($small$)$ set of morphisms $S$ of $\Fun(\calC, \bfA)$ which generates $\widetilde{S}$ as a weakly saturated class of morphisms.
\end{lemma}

We will prove a generalization of Lemma \ref{mainerthyme} in \S \ref{quasilimit3}
(Lemma \ref{mainertime}).

\begin{proof}[Proof of Proposition \ref{smurff}]
We first treat the case of the projective model structure. For each
For each object $C \in \calC$ and each $A \in \bfA$, we define 
$$\calF^C_A: \calC \rightarrow \bfA$$
by the formula\index{not}{FcalCA@$\calF^{C}_{A}$}
$$\calF^C_A(C') = \coprod_{ \alpha \in \bHom_{\calC}(C,C') } A.$$
We note that if $i: A \rightarrow A'$ is a (trivial) cofibration in $\bfA$, then the induced map
$\calF^C_{A} \rightarrow \calF^C_{A'}$ is a strong (trivial) cofibration in $\Fun(\calC, \bfA)$. 

Let $I_0$ be a set of generating cofibrations $i: A \rightarrow B$ for $\bfA$, and let $I$ be the set of all induced maps $\calF^{C}_{A} \rightarrow \calF^{C}_{B}$ (where $C$ ranges over $\calC$. Let $J_0$ be a set of generating trivial cofibrations for $\bfA$, and define 
$J$ likewise. It follows immediately from the definitions that a morphism in $\Fun(\calC, \bfA)$ is a projective fibration if and only if it has the right lifting property with respect to every morphism in $J$, and a weak trivial fibration if and only if it has the right lifting property with respect to every morphism in $I$. Let $\overline{I}$ and $\overline{J}$ be the weakly saturated classes of morphisms of $\Fun(\calC, \bfA)$ generated by $I$ and $J$, respectively. Using the small object argument, we deduce:
\begin{itemize}
\item[$(i)$] Every morphism $f: X \rightarrow Z$ in $\Fun(\calC, \bfA)$ admits a factorization
$$ X \stackrel{f'}{\rightarrow} Y \stackrel{f''}{\rightarrow} Z$$
where $f' \in \overline{I}$ and $f''$ is a weak trivial fibration.
\item[$(ii)$] Every morphism $f: X \rightarrow Z$ in $\Fun(\calC, \bfA)$ admits a factorization
$$ X \stackrel{f'}{\rightarrow} Y \stackrel{f''}{\rightarrow} Z$$
where $f' \in \overline{J}$ and $f''$ is a projective fibration. 
\item[$(iii)$] The class $\overline{I}$ coincides with the class of projective cofibrations in $\bfA$.
\end{itemize}
Furthermore, since the class of trivial projective cofibrations in $\Fun(\calC, \bfA)$ is weakly saturated and contains $J$, it contains $\overline{J}$. This proves that $\Fun(\calC, \bfA)$ satisfies the factorization axioms. The only other nontrivial point to check is that $\Fun(\calC, \bfA)$ satisfies the lifting axioms. Consider a diagram
$$ \xymatrix{ A \ar[d]^{i} \ar[r] & X \ar[d]^{p} \\
C \ar[r] \ar@{-->}[ur] & Y }$$
in $\Fun(\calC, \bfA)$, where $i$ is a projective cofibration and $p$ is a projective fibration. We wish to show that there exists a dotted arrow as indicated, provided that either $i$ or $p$ is a weak equivalence.
If $p$ is a weak equivalence then this follows immediately from the definition of a injective fibration.
Suppose instead that $i$ is a trivial projective cofibration. We wish to show that $i$ has the left lifting property with respect to every projective fibration. It will suffice to show every trivial injective fibration belongs to $\overline{J}$ (this will also show that $J$ is a set of generating trivial cofibrations for $\Fun(\calC, \bfA)$, which will show that the projective model structure on $\Fun(\calC, \bfA)$ is combinatorial). Suppose then that $i$ is a trivial weak coibration, and choose a factorization
$$ A \stackrel{i'}{\rightarrow} B \stackrel{i''}{\rightarrow} C$$
where $i' \in \overline{J}$ and $i''$ is a projective fibration. Then $i'$ is a weak equivalence, so
that $i''$ is a weak equivalence by the two-out-of-three property. Consider the diagram
$$ \xymatrix{ A \ar[d]^{i} \ar[r]^{i'} & B \ar[d]^{i''} \\
C \ar[r]^{=} \ar@{-->}[ur] & C. }$$
Since $i$ is a cofibration, there exists a dotted arrow as indicated. This proves that $i$ is a retract of $i'$, and therefore belongs to $\overline{J}$ as desired.

We now prove the existence of the injective model structure on $\Fun(\calC, \bfA)$. Here
it is difficult to proceed directly, so we will instead apply Proposition \ref{bigmaker}. It will suffice to check each of the hypotheses in turn:
\begin{itemize}
\item[$(1)$] The collection of injective cofibrations in $\Fun(\calC, \bfA)$ is generated (as a weakly saturated class) by some small set of morphisms. This follows from Lemma \ref{mainertime}.
\item[$(2)$] The collection of trivial injective cofibrations in $\Fun(\calC, \bfA)$ is weakly saturated: this follows immediately from the fact that the class of injective cofibrations in $\bfA$ is weakly saturated.
\item[$(3)$] The collection of weak equivalences in $\Fun(\calC, \bfA)$ is an accessible subcategory
of $\Fun(\calC, \bfA)^{[1]}$: this follows from the proof of Proposition \ref{horse1}, since the collection of weak equivalences in $\bfA$ form an accessible subcategory of $\bfA^{[1]}$.  

\item[$(4)$] The collection of weak equivalences in $\Fun(\calC, \bfA)$ satisfy the two-out-of-three property: this follows immediately from the fact that the weak equivalences in $\bfA$ satisfy the two-out-of-three property.

\item[$(5)$] Let $f: X \rightarrow Y$ be a morphism in $\bfA$ which has the right lifting property with respect to every injective cofibration. In particular, $f$ has the right lifting property with respect to each of the morphisms in the class $I$ defined above, so that $f$ is a trivial projective fibration, and in particular a weak equivalence.
\end{itemize}
\end{proof}

\begin{remark}\label{postsm}
In the situation of Proposition \ref{smurff}, if $\bfA$ is assumed to be right or left proper, then $\Fun(\calC, \bfA)$ is likewise right or left proper (with respect to either the projective or the injective model structures). 
\end{remark}

\begin{remark}\label{postsmurff}
It follows from the proof of Proposition \ref{smurf} that the class of projective cofibrations
is generated (as a weakly saturated class of morphisms) by the maps $j: \calF^{C}_{A} \rightarrow \calF^{C}_{A'}$, where $C \in \calC$ and $A \rightarrow A'$ is a cofibration in $\bfA$. 
We observe that $j$ is a injective cofibration. It follows that every projective cofibration is a injective cofibration; dually, every injective fibration is a projective fibration.
\end{remark}

\begin{remark}\label{twofus}
The construction of Proposition \ref{smurff} is functorial in the following sense:
given a Quillen adjunction of combinatorial model categories
$\Adjoint{F}{\bfA}{\bfB}{G}$ and
a small category $\calC$, composition with $F$ and $G$ determines a Quillen adjunction
$$ \Adjoint{ F^{\calC}}{\Fun(\calC, \bfA)}{\Fun(\calC, \bfB)}{G^{\calC}}$$
(with respect to either the injective or the projective model structures).
Moreover, if $(F,G)$ is a Quillen equivalence, then so is $( F^{\calC}, G^{\calC} )$.
\end{remark}

Because the projective and injective model structures on
$\Fun(\calC, \bfA)$ have the same weak equivalences, the identity
functor $\id_{\Fun(\calC, \bfA)}$ is a Quillen equivalence between them. However, it is
important to keep distinguish these two model structures, because
they have different variance properties as we now explain.

Let $f: \calC \rightarrow \calC'$ be a functor between small categories. Then
composition with $f$ yields a pullback functor $f^{\ast}:
\Fun(\calC', \bfA) \rightarrow \Fun(\calC, \bfA)$. Since $\bfA$ admits small limits and colimits, $f^{\ast}$ has a right adjoint which we shall denote by $f_{\ast}$ and a left
adjoint which we shall denote by $f_{!}$.

\begin{proposition}\label{colbinn}
Let $\bfA$ be a combinatorial model category,
and let $f: \calC \rightarrow \calC'$ be a functor between small categories.
Then:

\begin{itemize}
\item[$(1)$] The pair $( f_{!}, f^{\ast} )$ determines a Quillen
adjunction between the {\em projective} model structures on
$\Fun(\calC, \bfA)$ and $\Fun(\calC', \bfA)$.

\item[$(2)$] The pair $( f^{\ast}, f_{\ast} )$ determines a
Quillen adjunction between the {\em injective} model structures on
$\Fun(\calC, \bfA)$ and $\Fun(\calC', \bfA)$.
\end{itemize}
\end{proposition}

\begin{proof}
This follows immediately from the simple observation that
$f^{\ast}$ preserves weak equivalences, projective fibrations, and weak
cofibrations.
\end{proof}

We now review the theory of homotopy limits and colimits in a combinatorial model category $\bfA$. For simplicity, we will discuss homotopy limits and leave the analogous theory of homotopy colimits
to the reader. Let $\bfA$ be a combinatorial model category, and let 
$f: \calC \rightarrow \calC'$ be functor betweeen (small) categories. 
We wish to consider the right-derived functor $Rf_{\ast}$ of the right Kan extension $f_{\ast}: \Fun(\calC, \bfA) \rightarrow
\Fun(\calC', \bfA)$. This derived functor is called the {\it homotopy right Kan extension} functor.
The usual way of defining it involves choosing a ``fibrant replacement functor'' $Q: \Fun(\calC, \bfA) \rightarrow \Fun(\calC, \bfA)$, and setting $Rf_{\ast} = f_{\ast} \circ Q$. The assumption that $\bfA$ is combinatorial guarantees that such a fibrant replacement functor exists. However, for our purposes it is more convenient to address the indeterminacy in the definition of $Rf_{\ast}$ in another way.

Let $F \in \Fun(\calC, \bfA)$, $G \in \Fun(\calC', \bfA)$, and let $\eta: G \rightarrow f_{\ast} F$ be a map in $\Fun(\calC', \bfA)$. We will say that $\eta$ {\it exhibits $G$ as the homotopy right Kan extension of $F$}
if, for some weak equivalence $F \rightarrow F'$ where $F'$ is injectively fibrant in $\Fun(\calC, \bfA)$, the composite map $G \rightarrow f_{\ast} F \rightarrow f_{\ast} F'$ is a weak equivalence in
$\Fun(\calC', \bfA)$. Since $f_{\ast}$ preserves weak equivalences between injectively fibrant objects, this condition is independent of the choice of $F'$.\index{gen}{Kan extension!homotopy}

\begin{remark}
Given an object $F \in \Fun(\calC, \bfA)$, it is not necessarily the case that there exists a map
$\eta: G \rightarrow f_{\ast} F$ which exhibits $G$ as a homotopy right Kan extension of $F$. 
However, such a map can always be found after replacing $F$ by a weakly equivalent object; for example, if $F$ is injectively fibrant, we may take $G = f_{\ast} F$ and $\eta$ to be the identity map.
\end{remark}

Let $[0]$ denote the final object of $\Cat$: that is, the category with one object and only the identity morphism. For {\em any} category $\calC$, there is a unique
functor $f: \calC \rightarrow [0]$. If $\bfA$ is a combinatorial
model category, $F: \calC \rightarrow \bfA$ a functor, and $A \in \bfA \simeq \Fun( [0], \bfA)$
is an object, then we will say that a natural transformation
$\alpha: f^{\ast} A \rightarrow F$ {\it exhibits $A$ as a homotopy limit of $F$} if it exhibits $A$
as a homotopy right Kan extension of $F$. Note that we can identify
$\alpha$ with a map $A \rightarrow \lim_{C \in \calC} F(C)$ in the model category $\bfA$.

The theory of homotopy right Kan extensions in general can be reduced to the theory
of homotopy limits, in view of the following result:

\begin{proposition}\label{sabke}
Let $\bfA$ be a combinatorial model category, let $f: \calC \rightarrow \calD$ be a functor
between small categories, and let $F: \calC \rightarrow \bfA$ and $G: \calD \rightarrow \bfA$
be diagrams. A natural transformation $\alpha: f^{\ast} G \rightarrow F$ exhibits
$G$ as a homotopy right Kan extension of $F$ if and only if, for each object
$D \in \calD$, $\alpha$ exhibits $G(D)$ as a homotopy limit of the composite diagram
$$ F_{D/}: \calC \times_{ \calD } \calD_{D/} \rightarrow \calC \stackrel{F}{\rightarrow} \bfA. $$
\end{proposition}

To prove Proposition \ref{sabke}, we can immediately reduce to the case where
$F$ is a injectively fibrant diagram. In this case, $\alpha$ exhibits $G$ as a homotopy
right Kan extension of $F$ if and only if it induces a weak homotopy equivalence
$G(D) \rightarrow \lim F_{D/}$, for each $D \in \calD$. It will therefore suffice to prove the following result
(in the case $\calC' = \calC \times_{\calD} \calD_{D/}$):

\begin{lemma}\label{sumtuous}
Let $\bfA$ be a combinatorial model category and $g: \calC' \rightarrow \calC$
a functor which exhibits $\calC'$ as cofibered in sets over $\calC$.
Then the pullback functor $g^{\ast}: \Fun(\calC, \bfA) \rightarrow \Fun( \calC', \bfA)$
preserves injective fibrations.
\end{lemma}

\begin{proof}
It will suffice to show that the left adjoint $g_{!}$ preserves weak trivial cofibrations. 
Let $\alpha: F \rightarrow F'$ be a map in $\Fun( \calC', \bfA)$.
We observe that for each object $C \in \calC$, the map 
$(q_{!} \alpha)(C): (q_{!} F)(C) \rightarrow (q_{!} F')(C)$ can be identified
with the coproduct of the maps $\{ \alpha(C'): F(C') \rightarrow F'(C') \}_{ C' \in g^{-1} \{C\} }$.
If $\alpha$ is a weak trivial cofibration, then each of these maps is a trivial cofibration in $\bfA$, so that
$q_{!} \alpha$ is again a weak trivial cofibration as desired.
\end{proof}

\begin{remark}
In the preceding discussion, we considered injective model structures,
$Rf_{\ast}$, and homotopy limits. An entirely dual discussion may
be carried out with projective model structures and $Lf_{!}$; one obtains a
notion of {\it homotopy colimit} which is the dual of the notion
of homotopy limit.
\end{remark}

\begin{example}
Let $\bfA$ be a combinatorial model category, and consider a diagram
$$ X' \stackrel{f}{\leftarrow} X \stackrel{g}{\rightarrow} X''.$$
This diagram is projectively cofibrant if and only if the object $X$ is cofibrant, and the maps
$f$ and $g$ are both cofibrations. Consequently, the definition of homotopy colimits given above recovers, as a special case, the theory of homotopy pushouts presented in \S \ref{hopush}.
\end{example}

\subsection{Reedy Model Structures}\label{coreed}

Let $\bfA$ be a combinatorial model category and $\calJ$ a small category. In
\S \ref{qlim7}, we saw that the diagram category $\Fun( \calJ, \bfA)$ can again be regarded as a combinatorial model category, via either the projective or injective model structure of
Proposition \ref{smurff}. In the special case where $\calJ$ is a {\em Reedy category}
(see Definition \ref{quellreed}), it is often useful to consider still another model structure on $\Fun(\calJ, \bfA)$: the {\it Reedy model structure}. We will sketch the definition and some of the basic properties of Reedy model categories below; we refer the reader to \cite{hirschhorn} for a more detailed treatment.

\begin{definition}\label{quellreed}\index{gen}{Reedy!category}\index{gen}{category!Reedy}
A {\it Reedy category} is a small category $\calJ$ equipped with a factorization
system $\calJ^{L}, \calJ^{R} \subseteq \calJ$ satisfying the following conditions:
\begin{itemize}
\item[$(1)$] Every isomorphism in $\calJ$ is an identity map.
\item[$(2)$] Given a pair of object $X, Y \in \calJ$, let us write
$X \preceq_0 Y$ if either there exists a morphism $f: X \rightarrow Y$ belonging
to $\calJ^{R}$ or there exists a morphism $g: Y \rightarrow X$ belonging to $\calJ^{L}$.
We will write $X \prec_0 Y$ if $X \preceq_0 Y$ and $X \neq Y$.
Then there are no infinite descending chains $$ \ldots \prec_0 X_2 \prec_0 X_1 \prec_0 X_0.$$
\end{itemize}
\end{definition}

\begin{remark}
Let $\calJ$ be a category equipped with a factorization system $(\calJ^{L}, \calJ^{R})$, and let
$\preceq_0$ be the relation described in Definition \ref{quellreed}. This relation is generally not transitive. We will denote its transitive closure by $\preceq$. Then condition
$(2)$ of Definition \ref{quellreed} guarantees that $\preceq$ is a well-founded partial ordering on the set of objects of $\calJ$. In other words, every nonempty set $S$ of objects
of $\calJ$ contains a $\preceq$-minimal element.
\end{remark}

\begin{remark}
In the situation of Definition \ref{quellreed}, we will often abuse terminology and simply refer to $\calJ$ as a Reedy category, implicitly assuming that a factorization system on $\calJ$ has been specified as well.
\end{remark}

\begin{warning}
Condition $(1)$ of Definition \ref{quellreed} is not stable under equivalence of categories.
Suppose that $\calJ$ is equivalent to a Reedy category. Then $\calJ$ can itself be regarded as a Reedy category if and only if every isomorphism class of objects in $\calJ$ contains
a unique representative. (Definition \ref{quellreed} can easily be modified so as to be invariant under equivalence, but it is slightly more convenient not to do so.)
\end{warning}

\begin{example}\label{onep}
The category $\cDelta$ of combinatorial simplices is a Reedy category with respect to
the factorization system $( \cDelta^{L}, \cDelta^{R} )$; here a morphism
$f: [m] \rightarrow [n]$ belongs to $\cDelta^{L}$ if and only if $f$ is surjective, and
to $\cDelta^{R}$ if and only if $f$ is injective.
\end{example}

\begin{example}\label{twop}
Let $\calJ$ be a Reedy category with respect to the factorization system
$( \calJ^{L}, \calJ^{R})$. Then $\calJ^{op}$ is a Reedy category with respect to the
factorization system $( (\calJ^R)^{op}, (\calJ^{L})^{op} )$. 
\end{example}

\begin{notation}\index{gen}{latching object}\index{gen}{object!latching}\index{not}{LsubJ@$L_{J}(X)$}\label{bugga}
Let $\calJ$ be a Reedy category, $\calC$ a category which admits small limits and colimits, and $X: \calJ \rightarrow \calC$ a functor. For every object $J \in \calJ$, we define the
{\it latching object} $L_{J}(X)$ to be the colimit
$$ \varinjlim_{J' \in \calJ^{R}_{/J}, J' \neq J} X(J').$$
Similarly, we define the {\it matching object}\index{gen}{matching object}\index{gen}{object!matching}
\index{not}{MsubJ@$M_{J}(X)$} to be the limit
$$ \varprojlim_{J' \in \calJ^{L}_{J/}, J' \neq J} X(J').$$
We then have canonical maps
$L_{J}(X) \rightarrow X(J) \rightarrow M_{J}(X).$
\end{notation}

\begin{example}\label{croupus}
Let $X: \cDelta^{op} \rightarrow \Set$ be a simplicial set, and regard
$\cDelta^{op}$ as a Reedy category using Examples \ref{onep} and \ref{twop}.
For every nonnegative integer $n$, the latching object $L_{[n]} X$ can be identified
with the collection of all degenerate simplices of $X$. In particular, the map
$L_{[n]}(X) \rightarrow X([n])$ is always a monomorphism.

More generally, we observe that a map of simplicial sets $f: X \rightarrow Y$ is a monomorphism if and only if, for every $n \geq 0$, the map
$$ L_{[n]}(Y) \coprod_{ L_{[n]}(X) } X([n]) \rightarrow Y([n])$$
is a monomorphism of sets. The ``if'' direction is obvious. For the converse,
let us suppose that $f$ is a monomorphism; we must show that if
$\sigma$ is an $n$-simplex of $X$ such that $f(\sigma)$ is degenerate, then
$\sigma$ is already degenerate. If $f(\sigma)$ is degenerate, then
$f(\sigma) = \alpha^{\ast} f(\sigma) = f( \alpha^{\ast} \sigma)$ where $\alpha: [n] \rightarrow [n]$ is a map of linearly ordered sets other than the identity. Since $f$ is a monomorphism, we deduce that $\sigma = \alpha^{\ast} \sigma$, so that $\sigma$ is degenerate as desired.
\end{example}

\begin{remark}
Let $X: \calJ \rightarrow \calC$ be as in Notation \ref{bugga}. Then the
$J$th matching object $M_{J}(X)$ can be identified with the $J$th latching object
of the induced functor $X^{op}: \calJ^{op} \rightarrow \calC^{op}$. 
\end{remark}

\begin{remark}\label{surpose}
Let $X: \calJ \rightarrow \calC$ be as in Notation \ref{bugga}. Then the
$J$th matching object $M_{J}(X)$ can also be identified with the colimit
$$ \varinjlim_{(f: J' \rightarrow J) \in S} X(J')$$
where $S$ is any full subcategory of $\calJ_{/J}$ with the following properties:
\begin{itemize}
\item[$(1)$] Every morphism $f: J' \rightarrow J$ which belongs to
$\calJ^{R}$ and is not an isomorphism also belongs to $S$.
\item[$(2)$] If $f: J' \rightarrow J$ belongs to $S$, then
$J \npreceq J'$.
\end{itemize}
This follows from a cofinality argument, since every morphism
$f: J' \rightarrow J$ in $S$ admits a canonical factorization
$$ J' \stackrel{f'}{\rightarrow} J'' \stackrel{f''}{\rightarrow} J,$$
where $f'$ belongs to $\calJ^{L}$ and $f''$ belongs to $\calJ^{R}$. Assumption
$(2)$ guarantees that the map $f''$ is not an isomorphism.

Similarly, we can replace the limit $\varprojlim_{ f: J \rightarrow J'} X(J')$
defining the matching object $M_{J}(X)$ by a limit over a slightly larger category, when convenient.
\end{remark}

\begin{notation}
Let $\calJ$ be a Reedy category. A {\it good filtration} of $\calJ$ is a transfinite sequence
$$ \{ \calJ_{\beta} \}_{\beta < \alpha}$$ of full subcategories of $\calJ$ with the following properties:
\begin{itemize}
\item[$(a)$] The filtration is exhaustive in the following sense: every object of
$\calJ$ belongs to $\calJ_{\beta}$ for sufficiently large $\beta < \alpha$.
\item[$(b)$] For each ordinal $\beta < \alpha$, the category $\calJ_{\beta}$ is obtained
from the subcategory $\calJ_{< \beta} = \bigcup_{ \gamma < \beta} \calJ_{\gamma}$
by adjoining a single new object $J_{\beta}$ satisfying the following condition:
if $J \in \calJ$ satisfies $J \prec J_{\beta}$, then $J \in \calJ_{< \beta}$.
\end{itemize}
\end{notation}

\begin{remark}
Let $\calJ$ be a Reedy category. Then there exists a good filtration of $\calJ$. In fact, the
existence of a good filtration is {\em equivalent} to the second assumption of Definition \ref{quellreed}.  
\end{remark}

\begin{remark}
Let $\calJ$ be a Reedy category with respect to the factorization system
$(\calJ^{L}, \calJ^{R})$, and let $\{ \calJ_{\beta} \}_{\beta < \alpha}$ be a good filtration of $\calJ$. Then each $\calJ_{\beta}$ admits a factorization system
$( \calJ^{L} \cap \calJ_{\beta}, \calJ^{R} \cap \calJ_{\beta})$. In other words, if
$f: I \rightarrow K$ is a morphism in $\calJ_{\beta}$ which admits a factorization
$$ I \stackrel{f'}{\rightarrow} J \stackrel{f''}{\rightarrow} K$$
where $f'$ belongs to $\calJ^{L}$ and $f''$ belongs to $\calJ^{R}$, then
the object $J$ also belongs to $\calJ_{\beta}$. This is clear: either
$f''$ is an isomorphism, in which case $J = K \in \calJ_{\beta}$, or
$f''$ is not an isomorphism so that $J \prec K$ implies that $J \in \calJ_{< \beta}$.
\end{remark}

The following result summarizes the essential features of a good filtration:

\begin{proposition}\label{twingood}
Let $\calJ$ be a Reedy category with a good filtration $\{ \calJ_{\beta} \}_{ \beta < \alpha}$,
let $\beta < \alpha$ be an ordinal, so that $\calJ_{\beta}$ is obtained from $\calJ_{< \beta}$ by 
adjoining a single new object $J$. Then we have a homotopy pushout square
(with respect to the Joyal model structure)
$$ \xymatrix{ \Nerve (\calJ_{< \beta})_{/J} \star \Nerve (\calJ_{< \beta})_{J/}
\ar@{^{(}->}[d] \ar[r] & \Nerve(\calJ_{< \beta}) \ar@{^{(}->}[d] \\
\Nerve( \calJ_{< \beta})_{/J} \star \{J\} \star \Nerve( \calJ_{< \beta})_{J/} \ar[r]
& \Nerve( \calJ_{\beta}). }$$
\end{proposition}

\begin{corollary}\label{stapler}
Let $\calJ$ be a Reedy category with a good filtration $\{ \calJ_{\beta} \}_{\beta < \alpha}$,
let $\beta < \alpha$ be an ordinal, so that $\calJ_{\beta}$ is obtained from
$\calJ_{< \beta}$ by adjoining a single new object $J$. Let
$\calC$ be a category which admits small limits and colimits, and let
$X: \calJ_{< \beta} \rightarrow \calC$ be a functor, and let the latching and matching
objects $L_J(X)$ and $M_{J}(X)$ be defined as in Notation \ref{bugga}
(note that this does not require that the functor $X$ be defined on the object $J$), so that we have a canonical map $\alpha: L_{J}(X) \rightarrow M_{J}(X)$. The following data are equivalent:
\begin{itemize}
\item[$(1)$] A functor $\overline{X}: \calJ_{\beta} \rightarrow \calC$ extending $X$.
\item[$(2)$] A commutative diagram
$$ \xymatrix{ & C \ar[dr] & \\
L_{J}(X) \ar[ur] \ar[rr]^{\alpha} & & M_{J}(X) }$$
in the category $\calC$.
\end{itemize}
The equivalence carries a functor $\overline{X}$ to the evident diagram with
$C = \overline{X}(J)$.
\end{corollary}

\begin{proof}
Using Proposition \ref{twingood}, we see that giving an extension $\overline{X}: \calJ_{\beta} \rightarrow \calC$ of $X$ is equivalent to giving an extension $\overline{Y}: ( \calJ_{< \beta})_{/J} \star \{J\} \star ( \calJ_{< \beta})_{J/} \rightarrow \calC$ of the composite functor
$$Y: (\calJ_{< \beta})_{/J} \star (\calJ_{< \beta})_{J/} \rightarrow \calJ_{< \beta} \stackrel{X}{\rightarrow} \calC.$$
This, in turn, is equivalent to giving a commutative diagram
$$ \xymatrix{ & C \ar[dr] & \\
\varinjlim Y| (\calJ_{< \beta})_{/J} \ar[ur] \ar[rr]^{\alpha'} & & \varprojlim Y|(\calJ_{< \beta})_{J/}) }$$
where $\alpha'$ is the map induced by the diagram $Y$. The equivalence of this with the
data $(2)$ follows immediately from Remark \ref{surpose}.
\end{proof}

\begin{remark}\label{stapler2}
The proof of Corollary \ref{stapler} carries over without essential change to the case where
$\calC$ is an $\infty$-category which admits small limits and colimits. In this case, to
extend a functor $X: \Nerve( \calJ_{< \beta}) \rightarrow \calC$ to a functor $\overline{X}$ defined on the whole of $\calJ_{\beta}$, it will suffice to specify the object
$$ \overline{X}(J) \in \calC_{ X | (\calJ_{< \beta})_{/J} / \, X | (\calJ_{< \beta})_{J/}}
\simeq \calC_{ L_{J}(X) / \, / M_J(X) },$$
where the latching and matching objects $L_J(X), M_J(X) \in \calC$ are defined in the obvious way.
\end{remark}

The proof of Proposition \ref{twingood} will require a few preliminaries.

\begin{lemma}\label{gump1}
Let $\calJ$ be a Reedy category equipped with a good filtration
$\{ \calJ_{\beta} \}_{\beta < \alpha}$. Fix $\beta < \alpha$, and let
$\calJ_{\beta}$ be obtained from $\calJ_{< \beta}$ by adjoining the object $J$.
Let $f: J \rightarrow J$ be a map which is not the identity, let
$\calI$ denote the category $(\calJ_{J/})_{/f} \simeq (\calJ_{/J})_{f/}$ of
factorizations of the morphism $f$, and let $\calI_0 = \calI \times_{\calJ} \calJ_{< \beta}$. Then
the nerve $\Nerve \calI_0$ is weakly contractible.
\end{lemma}

\begin{proof}
Let $\calI_1$ denote the full subcategory of $\calI_0$ spanned by those diagrams
$$ \xymatrix{ & I \ar[dr]^{f''} & \\
J \ar[rr]^{f} \ar[ur]^{f'} & & J }$$
where $I \in \calJ_{< \beta}$ and $f''$ is a morphism in $\calJ^{R}$. The inclusion
$\calI_1 \subseteq \calI_0$ admits a left adjoint, so that $\Nerve \calI_1$ is a deformation retract
of $\Nerve \calI_0$. It will therefore suffice to show that $\Nerve \calI_1$ is weakly contractible.
Let $\calI_2$ denote the full subcategory of $\calI_1$ spanned by those diagrams as above
where, in addition, the morphism $f'$ belongs to $\calJ^{L}$. Then the inclusion
$\calI_2 \subseteq \calI_1$ admits a right adjoint, so that $\Nerve \calI_2$ is a deformation retract of $\Nerve \calI_1$. It will therefore suffice to show that $\Nerve \calI_2$ is weakly contractible. This is clear, since the category $\calI_2$ consists of a single object (with no nontrivial endomorphisms).
\end{proof}

\begin{lemma}\label{gump2}
Let $n \geq 1$, and suppose given a sequence of weakly contractible simplicial sets
$\{ A_i \}_{1 \leq i \leq n}$. Let $L$ denote the iterated join
$$ \{ J_0 \} \star A_1 \star \{J_1 \} \star A_2 \star \ldots \star A_n \star \{J_n \},$$
and let $K$ denote the simplicial subset of $L$ spanned by those simplices which do not contain
all of the vertices $\{ J_i \}_{0 \leq i \leq n}$. Then the inclusion $K \subseteq L$ is a categorical equivalence of simplicial sets.
\end{lemma}

\begin{proof}
If $n=1$, then this follows immediately from Lemma \ref{storuse}.
Suppose that $n > 1$. Let $X$ denote the iterated join
$A_1 \star \{J_1 \} \star A_2 \star \ldots \star \{J_{n-1} \} \star A_n$.
For every subset $S \subseteq \{ 1, \ldots, n-1\}$, let
$X(S)$ denote the simplicial subset of $X$ spanned by those simplices
which do not contain any vertex $J_i$ for $i \in S$.
Let $X' = \bigcup_{ S \neq \emptyset} X(S) \subseteq X(\emptyset) = X$.
Then $X'$ is the homotopy colimit of the diagram of simplicial sets
$\{ X(S) \}_{S \neq \emptyset}$. Each $X(S)$ is a join of weakly contractible simplicial sets,
and is therefore weakly contractible. Since $n > 1$, the partially ordered set
$\{ S \subseteq \{1, \ldots, n-1\}: S \neq \emptyset \}$ has a largest element, and is therefore
weakly contractible. It follows that the simplicial set $X'$ is weakly contractible.

The assertion that the inclusion $K \subseteq L$ is a categorical equivalence is
equivalent to the assertion that the diagram
$$ \xymatrix{ (\{J_0\} \star X') \coprod_{X'} ( X' \star \{J_n\}) \ar@{^{(}->}[d] \ar@{^{(}->}[r] & 
( \{J_0\} \star X) \coprod_{ X} (X \star \{J_0\}) \ar@{^{(}->}[d] \\
\{J_0\} \star X' \star \{J_n\} \ar@{^{(}->}[r] & \{J_0\} \star X \star \{J_n \} }$$
is a homotopy pushout square (with respect to the Joyal model structure).
To prove this, it suffices to observe that the vertical maps are both categorical equivalences (Lemma \ref{storuse}).
\end{proof}

\begin{proof}[Proof of Proposition \ref{twingood}]
Let $S$ denote the collection of all composable chains of morphisms
$$\overline{f}: J \stackrel{ f_1}{\rightarrow} J \stackrel{ f_2} \ldots \stackrel{f_n}{\rightarrow} J.$$
where $n \geq 1$ and each $f_{i} \neq \id_{J}$. For every subset $S' \subseteq S$, let
$X(S')$ denote the simplicial subset of $\Nerve( \calJ_{\beta})$ spanned by those simplices
$\sigma$ satisfying the following condition:
\begin{itemize}
\item[$(\ast)$] For every nondegenerate face $\tau$ of $\sigma$ of positive dimension, if every
vertex of $\tau$ coincides with $J$, then $\tau$ belongs to $S'$.
\end{itemize}
Note that $X(S)$ coincides with $\Nerve( \calJ_{\beta})$, while $X( \emptyset)$ coincides with
the pushout 
$$(\Nerve (\calJ_{< \beta})_{/J} \star \{ J \} \star \Nerve (\calJ_{< \beta})_{J/})
\coprod_{ \Nerve( \calJ_{< \beta})_{/J} \star \Nerve( \calJ_{< \beta})_{J/}}
\Nerve( \calJ_{< \beta} ). $$
It will therefore suffice to show that the inclusion $X(\emptyset) \subseteq X(S)$ is a
categorical equivalence of simplicial sets.

Choose a well-ordering
$$ S = \{ \overline{f}_0 < \overline{f}_1 < \overline{f}_2 < \ldots \}$$
with the following property: if $\overline{f}$ has length shorter than $\overline{g}$
(when regarded as a chain of morphisms), then $\overline{f} < \overline{g}$.
For every ordinal $\alpha$, let $S_{\alpha} = \{ \overline{f}_{\beta} \}_{\beta < \alpha}$.
We will prove that for every ordinal $\alpha$, the inclusion
$X( \emptyset) \subseteq X( S_{\alpha})$ is a categorical equivalence.
The proof proceeds by induction on $\alpha$. If $\alpha = 0$ there is nothing to prove, and
if $\alpha$ is a limit ordinal then the desired result follows from the inductive hypothesis and the fact that the class of categorical equivalences is stable under filtered colimits. We may therefore assume that
$\alpha = \beta + 1$ is a successor ordinal. The inductive hypothesis guarantees that
$X(\emptyset) \subseteq X( S_{\beta})$ is a categorical equivalence. It will therefore suffice to show that
the inclusion $j: X( S_{\beta} ) \subseteq X( S_{\alpha})$ is a categorical equivalence.
We may also suppose that $\beta$ is smaller than the order type of $S$, so that
$\overline{f}_{\beta}$ is well-defined (otherwise, the inclusion $j$ is an isomorphism and the result is obvious). 

Let $\overline{f} = \overline{f}_{\beta}$ be the composable chain of morphisms
$$\overline{f}: J \stackrel{ f_1}{\rightarrow} J \stackrel{ f_2}{\rightarrow} \ldots \stackrel{f_n}{\rightarrow} J.$$
For $1 \leq i \leq n$, let $A_i$ denote the nerve of the category
$$\calJ_{< \beta} \times_{\calJ} ( \calJ_{J/})_{/ f_i} \simeq \calJ_{< \beta} \times_{\calJ}
( \calJ_{/J})_{f_i/ }.$$
Let $K$ denote the simplicial subset of
$$ \{ J_0 \} \star A_1 \star \{J_1 \} \star A_2 \star \ldots \star A_n \star \{J_n \}$$
spanned by those simplices which do not contain every vertex $J_{n}$.
We then have a homotopy pushout diagram
$$ \xymatrix{ \Nerve(\calJ_{< \beta})_{/J} \star K \star \Nerve(\calJ_{< \beta})_{J/} \ar@{^{(}->}[d] \ar[r] & 
X( S_{\beta}) \ar[d] \\
\Nerve(\calJ_{< \beta})_{/J} \star \{J_0\} \star A_1 \star \ldots \star A_n \star \{J_n\} 
\star \Nerve(\calJ_{< \beta})_{J/} \ar[r] & X( S_{\alpha}).}$$
It will therefore suffice to prove that the left vertical map is a categorical equivalence.
In view of Corollary \ref{gyyyt}, it will suffice to show that the inclusion
$$ K \subseteq \{ J_0 \} \star A_1 \star \{J_1 \} \star A_2 \star \ldots \star A_n \star \{J_n \}$$
is a categorical equivalence. Since each $A_i$ is weakly contractible
(Lemma \ref{gump1}), this follows immediately from Lemma \ref{gump2}.
\end{proof}

\begin{proposition}\index{gen}{Reedy!model structure}\index{gen}{model category!Reedy}\label{reedmod}
Let $\calJ$ be a Reedy category, and let $\bfA$ be a model category.
Then there exists a model structure on the category of functors $\Fun(\calJ, \bfA)$
with the following properties:
\begin{itemize}\index{gen}{cofibration!Reedy}\index{gen}{Reedy!cofibration}
\item[$(C)$] A morphism $X \rightarrow Y$ in $\Fun(\calJ, \bfA)$ is a 
{\it Reedy cofibration} if and only if, for every object $J \in \calJ$, the induced map
$X(J) \coprod_{ L_{J}(X) } L_{J}(Y) \rightarrow Y(J)$ is a cofibration in $\bfA$.
\item[$(F)$] A morphism $X \rightarrow Y$ in $\Fun(\calJ, \bfA)$ is a\index{gen}{fibration!Reedy}
{\it Reedy fibration} if and only if, for every object $J \in \calJ$, the induced map
$X(J) \rightarrow Y(J) \times_{ M_{J}(Y) } M_{J}(X)$ is a fibration in $\bfA$.\index{gen}{Reedy!fibration}
\item[$(W)$] A morphism $X \rightarrow Y$ in $\Fun(\calJ, \bfA)$ is a weak equivalence if and only if, for every $J \in \calJ$, the map $X(J) \rightarrow Y(J)$ is a weak equivalence.
\end{itemize} 
Moreover, a morphism $f: X \rightarrow Y$ in $\Fun( \calJ, \bfA)$ is a trivial cofibration if
and only if the following condition is satisfied:
\begin{itemize}
\item[$(WC)$] For every object $J \in \calJ$, the map $X(J) \coprod_{ L_{J}(X) } L_{J}(Y) \rightarrow Y(J)$ is a trivial cofibration in $\bfA$.
\end{itemize}
Similarly, $f$ is a fibration if and only if it satisfies the dual condition:
\begin{itemize}
\item[$(WF)$] For every object $J \in \calJ$, the map $X(J) \rightarrow Y(J) \times_{ M_J(Y)} M_{J}(X)$ is a trivial fibration in $\bfA$.
\end{itemize}
\end{proposition}

The model structure of Proposition \ref{reedmod} is called the {\it Reedy model structure} on 
$\Fun(\calJ, \bfA)$. Note that Proposition \ref{reedmod} does not require that the model category $\bfA$ is combinatorial.

\begin{lemma}\label{jackal}
Let $\calJ$ be a Reedy category containing an object $J$, let $\bfA$ be a model category, and let $f: F \rightarrow G$ be a natural transformation in $\Fun(\calJ, \bfA)$.
Let $\calI \subseteq \calJ^{R}_{/J}$ be a sieve: that is, $\calI$ is a full subcategory
of $\calJ^{R}_{/J}$ with the property that if $I \rightarrow I'$ is a morphism in
$\calJ^{R}_{/J}$ such that $I' \in \calI$, then $I \in \calI$. Let $\calI' \subseteq \calI$ be another sieve. Then:
\begin{itemize}
\item[$(a)$] If the map $f$ satisfies condition $(C)$ of Proposition \ref{reedmod} for every object $I \in \calI$, then
the induced map
$$ \chi_{\calI', \calI}: \varinjlim( F | \calI) \coprod_{ \varinjlim(F | \calI') } \varinjlim( G | \calI')
\rightarrow \varinjlim( G | \calI )$$
is a cofibration in $\bfA$.
\item[$(b)$] If the map $f$ satisfies condition $(WC)$ of Proposition \ref{reedmod} for every object $I \in \calI$, then the map $\chi_{\calI', \calI}$ is a trivial cofibration in $\bfA$.
\end{itemize}
\end{lemma}

\begin{proof}
We will prove $(a)$; the proof of $(b)$ is identical. Choose a transfinite sequence of sieves
$\{ \calI_{\beta} \subseteq \calI \}_{ \beta < \alpha }$ with the following properties:
\begin{itemize}
\item[$(i)$] The union $\bigcup_{\beta < \alpha} \calI_{\beta}$ coincides with $\calI$.
\item[$(ii)$] For each $\beta < \alpha$, the sieve $\calI_{\beta}$ is obtained
from $\calI_{< \beta} = \calI' \cup (\bigcup_{ \gamma < \beta} \calI_{\gamma})$
by adjoining a single new object $J_{\beta} \in \calJ^{R}_{/J}$. 
\end{itemize}
For every triple $\delta \leq \gamma \leq \beta \leq \alpha$, let
$\chi_{ \delta, \gamma, \beta}$ denote the induced map
$$ \varinjlim( F| \calI_{< \beta}) \coprod_{ \varinjlim(F | \calI_{< \delta})}
\varinjlim( G| \calI_{< \delta}) \rightarrow
\varinjlim( F | \calI_{< \beta} ) \coprod_{ \varinjlim( F| \calI_{< \gamma}) }
\varinjlim( G| \calI_{< \gamma} ).$$
We wish to prove that $\chi_{0, \alpha, \alpha}$ is a cofibration.
We will prove more generally that $\chi_{\delta, \gamma, \beta}$ is an
equivalence for every $\delta \leq \gamma \leq \beta \leq \alpha$.
The proof uses induction on $\gamma$. If $\gamma$ is a limit ordinal, 
then we can write $\chi_{\delta, \gamma, \beta}$ as the transfinite composition
of the maps $\{ \chi_{ \epsilon, \epsilon + 1, \beta} \}_{\delta \leq \epsilon < \gamma}$, which are cofibrations by the inductive hypothesis. We may therefore assume that
$\gamma = \gamma_0 + 1$ is a successor ordinal. If $\delta = \gamma$, then
$\chi_{\delta, \gamma, \beta}$ is an isomorphism; otherwise, we have
$\delta \leq \gamma_0$. In this case, we have
$$ \chi_{ \delta, \gamma, \beta} = \chi_{ \gamma_0, \gamma, \beta} \circ \chi_{ \delta, \gamma_0, \beta}.$$
Using the inductive hypothesis, we can reduce to the case $\delta = \gamma_0$. The map
$\chi_{\gamma_0, \gamma, \beta}$ is a pushout of the map $\chi_{\gamma_0, \gamma, \gamma}$. We are therefore reduced to proving that $\chi_{\gamma_0, \gamma, \gamma}$
is a cofibration. But $\chi_{\gamma_0, \gamma, \gamma}$ is a pushout of the map
$L_{I}(G) \coprod_{ L_{I}(F)} F(I) \rightarrow G(I)$ for $I = J_{\gamma_0}$.
This map is a cofibration by virtue of our assumption that $f$ satisfies $(C)$.
\end{proof}

\begin{proof}[Proof of Proposition \ref{reedmod}]
Let $f: X \rightarrow Z$ be a morphism in $\Fun(\calJ, \bfA)$. We will prove that
$f$ admits a factorization
$$ X \stackrel{f'}{\rightarrow} Y \stackrel{f''}{\rightarrow} Z$$
where:
\begin{itemize}
\item[$(i)$] The map $f''$ is a fibration, and $f'$ satisfies $(WC)$.
\item[$(ii)$] The map $f'$ is a cofibration, and $f''$ satisfies $(WF)$.
\end{itemize}
By symmetry, it will suffice to consider case $(i)$. Choose a good filtration
$\{ \calJ_{\beta} \}_{ \beta < \alpha}$ of $\calJ$. For $\beta < \alpha$, let
$X_{\beta} = X | \calJ_{\beta}$, $Z_{\beta} = Z | \calJ_{\beta}$, and let
$f_{\beta}: X_{\beta} \rightarrow Z_{\beta}$ be the restriction of $f$. We will construct
a compatible family of factorizations of $f_{\beta}$ as a composition
$$ X_{\beta} \stackrel{f'_{\beta}}{\rightarrow} Y_{\beta} \stackrel{f''_{\beta}}{\rightarrow} Z_{\beta}.$$
Suppose that $\calJ_{\beta}$ is obtained from $\calJ_{< \beta}$ by adjoining a single new object $J$. Assuming that $(f'_{\gamma}, f''_{\gamma})$ has been constructed for all
$\gamma < \beta$, we note that constructing $(f'_{\beta}, f''_{\beta})$ is equivalent
(by virtue of Corollary \ref{stapler}) to giving a commutative diagram
$$ \xymatrix{ L_{J}(X) \ar[d] \ar[r] & L_{J}(Y_{< \beta}) \ar[d] & \\
X(J) \ar[r] & Y_{\beta}(J) \ar[r] \ar[d] & Z(J) \ar[d] \\
& M_J(Y_{< \beta}) \ar[r] & M_{J}(Z). }$$
In other words, we must factor a certain map
$$g: L_{J}( Y_{< \beta}) \coprod_{ L_{J}(X) } X(J) \rightarrow
M_{J}(Y_{< \beta} ) \times_{ M_{J}(Z)} Z(J)$$
as a composition
$$ L_{J}( Y_{< \beta}) \coprod_{ L_{J}(X) } X(J) \stackrel{g'}{\rightarrow}
Y_{\beta}(J) \stackrel{g''}{\rightarrow}
M_{J}(Y_{< \beta} ) \times_{ M_{J}(Z)} Z(J).$$
Using the fact that $\bfA$ is a model category, we can choose a factorization
where $g'$ is a trivial cofibration and $g''$ a fibration. It is readily verified that this construction has the desired properties.

We now prove the following:
\begin{itemize}
\item[$(i')$] A morphism $f: X \rightarrow Y$ in $\Fun(\calJ, \bfA)$ satisfies $(WC)$ if and only if
$f$ is both a fibration and a weak equivalence.
\item[$(ii')$] A morphism $f: X \rightarrow Y$ in $\Fun(\calJ, \bfA)$ satisfies $(WF)$ if and only if $f$ is both a cofibration and a weak equivalence.
\end{itemize}
By symmetry, it will suffice to prove $(i')$. The ``only if'' direction follows from
Lemma \ref{jackal}. For the ``if'' direction, it will suffice to show that for
each $\beta < \alpha$, the induced transformation
$f_{\beta}: X_{\beta} \rightarrow Y_{\beta}$ satisfies $(WC)$ when regarded as a morphism
of $\Fun( \calJ_{\beta}, \bfA)$. Suppose that $\calJ_{\beta}$ is obtained from $\calJ_{< \beta}$ by adjoining a single new element $J$. We have a commutative diagram
$$ \xymatrix{ & L_{J}(Y) \coprod_{ L_{J}(X) } X(J) \ar[dr]^{q} & \\
X(J) \ar[ur]^{p} \ar[rr]^{r} & & Y(J). }$$
We wish to prove that $q$ is a trivial cofibration in $\bfA$. Since $f$ is a cofibration in $\Fun(\calJ, \bfA)$, the map $q$ is a cofibration in $\bfA$. It will therefore suffice to show that
$q$ is a weak equivalence. By the two-out-of-three property, it will suffice to show that
$p$ and $r$ are weak equivalences. For $r$, this follows from our assumption that
$f$ is a weak equivalence in $\Fun(\calJ, \bfA)$. The map $p$ is a pushout of the map of latching objects $L_{J}(X) \rightarrow L_{J}(Y)$, which is a cofibration in $\bfA$ by virtue of the inductive hypothesis and Lemma \ref{jackal}.

Combining $(i)$ and $(i')$ (and the analogous assertions $(ii)$ and $(ii')$), we deduce that
$\Fun(\calJ, \bfA)$ satisfies the factorization axioms for a model category. To complete the proof, it will suffice to verify the lifting axioms:
\begin{itemize}
\item[$(i'')$] Every fibration in $\Fun(\calJ, \bfA)$ has the right lifting property with respect to
morphisms in $\Fun(\calJ, \bfA)$ which satisfy $(WC)$.
\item[$(ii'')$] Every cofibration in $\Fun(\calJ, \bfA)$ has the left lifting property with respect to morphisms in $\Fun(\calJ, \bfA)$ which satisfy $(WF)$.
\end{itemize}
Again, by symmetry it will suffice to prove $(i'')$. Consider a diagram
$$ \xymatrix{ A \ar[r] \ar[d]^{f} & X \ar[d]^{g} \\
B \ar[r] \ar@{-->}[ur]^{h} & Y, }$$
where $f$ satisfies $(WC)$ and $g$ satisfies $(F)$; we wish to prove that there
exists a dotted arrow $h$ as indicated, rendering the diagram commutative.
To prove this, we will construct a compatible family of natural transformations
$\{ h_{\beta}: B| \calJ_{\beta} \rightarrow X| \calJ_{\beta} \}_{ \beta < \alpha}$
which render the diagrams
$$ \xymatrix{ A | \calJ_{\beta} \ar[r] \ar[d] & X| \calJ_{\beta} \ar[d]^{g} \\
B| \calJ_{\beta} \ar[r] \ar[ur]^{h_{\beta}} & Y|\calJ_{\beta} }$$
commutative. Suppose that $\calJ_{\beta}$ is obtained from
$\calJ_{< \beta}$ by adjoining a single new object $J$. Assume that the maps
$\{ h_{\gamma} \}_{\gamma < \beta}$ have already been constructed, and can be amalgamated to a single natural transformation $h_{< \beta}: B | \calJ_{< \beta}
\rightarrow X| \calJ_{< \beta}$. Using Corollary \ref{stapler}, we see that extending
$h_{< \beta}$ to a map $h_{\beta}$ with the desired properties is equivalent to solving a lifting problem of the kind depicted in the following diagram:
$$ \xymatrix{ L_{J}(B) \coprod_{ L_{J}(A)} A(J) \ar[d]^{f'} \ar[r] & X(J) \ar[d]^{g'} \\
B(J) \ar@{-->}[ur] \ar[r] & Y(J) \times_{ M_J(Y)} M_J(X). }$$
Since our assumptions guarantee that $f'$ is a trivial cofibration and that $g'$ is a fibration, 
this lifting problem has a solution as desired.
\end{proof}

\begin{example}\label{tetsu}
Let $\bfA$ be the category of {\em bisimplicial} sets, which we will identify
with $\Fun( \cDelta^{op}, \sSet)$ and endow with the Reedy model structure.
It follows from Example \ref{croupus} that a morphism $f: X \rightarrow Y$ of bisimplicial sets is a Reedy cofibration if and only if it is a monomorphism. Consequently, the
Reedy model structure on $\bfA$ coincides with the injective model structure on $\bfA$.
\end{example}

\begin{example}\label{sued}
Let $\calJ$ be a Reedy category with $\calJ^{L} = \calJ$, and let
$\bfA$ be a model category. Then the weak equivalences and cofibrations of the
Reedy model structure (Proposition \ref{reedmod}) are the injective cofibrations and the
weak equivalences appearing in Definition \ref{cooper}. It follows that the Reedy model structure on $\Fun(\calJ, \bfA)$ coincides with the injective model structure of Proposition \ref{smurff} (in particular, the injective model structure is well defined in this case even without the assumption that $\bfA$ is combinatorial). Similarly, if $\calJ^{R} = \calJ$, then we can identify the Reedy model structure on $\Fun(\calJ, \bfA)$ with the projective model structure of Proposition \ref{smurff}.

In the general case, we can regard the Reedy model structure on $\Fun(\calJ, \bfA)$ as
a mixture of the projective and injective model structures. More precisely, we have the following:
\begin{itemize}
\item[$(i)$] A natural transformation $F \rightarrow G$ in $\Fun(\calJ, \bfA)$ satisfies
condition $(C)$ of Proposition \ref{reedmod} if and only if the induced transformation
$F|\calJ^{R} \rightarrow G| \calJ^{R}$ is a projective cofibration in $\Fun(\calJ^{R}, \bfA)$.
\item[$(ii)$] A natural transformation $F \rightarrow G$ in $\Fun(\calJ, \bfA)$ satisfies
condition $(F)$ of Proposition \ref{reedmod} if and only if the induced transformation
$F|\calJ^{L} \rightarrow G| \calJ^{L}$ is a injective fibration in $\Fun(\calJ^{L}, \bfA)$.
\item[$(iii)$] A natural transformation $F \rightarrow G$ in $\Fun(\calJ, \bfA)$ satisfies
condition $(WC)$ of Proposition \ref{reedmod} if and only if the induced transformation
$F|\calJ^{R} \rightarrow G| \calJ^{R}$ is a trivial projective cofibration in $\Fun(\calJ^{R}, \bfA)$.
\item[$(iv)$] A natural transformation $F \rightarrow G$ in $\Fun(\calJ, \bfA)$ satisfies
condition $(WF)$ of Proposition \ref{reedmod} if and only if the induced transformation
$F|\calJ^{L} \rightarrow G| \calJ^{L}$ is a trivial injective fibration in $\Fun(\calJ^{L}, \bfA)$.
\end{itemize}
\end{example}

\begin{remark}
Let $\calJ$ be a Reedy category and $\bfA$ a combinatorial model category, so that
the injective and projective model structures on $\Fun( \calJ, \bfA)$ are well-defined.
The identity functor from $\Fun(\calJ, \bfA)$ to itself can be regarded as a left Quillen equivalence from the projective model structure to the Reedy model structure, and from the
Reedy model structure to the injective model structure.
\end{remark}

\begin{corollary}\label{dimbu}
Let $\calC$ be a small category. Suppose that there exists a well-ordering
$\leq$ on the collection of objects of $\calC$ satisfying the following condition:
for every pair of objects $X,Y \in \calC$, we have
$$ \Hom_{\calC}(X,Y) = \begin{cases} \emptyset & \text{ if } $X > Y$ \\
\{ \id_X \} & \text{ if } X = Y. \end{cases}$$
Let $\bfA$ be a model category. Then:
\begin{itemize}
\item[$(i)$] A natural transformation $F \rightarrow G$ in $\Fun(\calC, \bfA)$ is a (trivial) projective cofibration if and only if, for every object $C \in \calC$, the induced map
$$ F(C) \coprod_{ \varinjlim_{ D \rightarrow C, D \neq C} F(D) }
\varinjlim_{ D \rightarrow C, D \neq C} G(D) \rightarrow G(C)$$
is a (trivial) cofibration in $\bfA$.
\item[$(ii)$] A natural transformation $F \rightarrow G$ in $\Fun(\calC^{op}, \bfA)$ is a
(trivial) injective fibration if and only if, for every object $C \in \calC$, the induced map
$$ F(C) \rightarrow G(C) \times_{ \varprojlim_{D \rightarrow C, D \neq C} G(D) }
\varprojlim_{D \rightarrow C, D \neq C} F(D) $$
is a (trivial) fibration in $\bfA$.
\end{itemize}
\end{corollary}

\begin{proof}
Combine Example \ref{sued} with Proposition \ref{reedmod}.
\end{proof}

\begin{corollary}\label{jonnyt}
Let $\bfA$ be a model category, let $\alpha$ be an ordinal, and let
$(\alpha)$ denote the linearly ordered set $\{ \beta < \alpha \}$, regarded as a category. Then:
\begin{itemize}
\item[$(1)$] Let $F \rightarrow F'$ be a natural transformation of diagrams
$(\alpha) \rightarrow \bfA$. Suppose that, for each $\beta < \alpha$, the maps
$$ \colim_{\gamma < \beta} F(\gamma) \rightarrow F(\beta)$$
$$ \colim_{\gamma < \beta} F'(\gamma) \rightarrow F'(\beta)$$
are cofibrations, while the map $F(\beta) \rightarrow F'(\beta)$ is a weak equivalence.
Then the induced map $$\colim_{\gamma < \alpha} F(\gamma) \rightarrow \colim_{\gamma < \alpha} F'(\gamma)$$ is a weak equivalence.
\item[$(2)$] Let $G \rightarrow G'$ be a natural transformation of diagrams
$(\alpha)^{op} \rightarrow \bfA$. Suppose that, for each $\beta < \alpha$, the maps
$$ G(\beta) \rightarrow \lim_{\gamma < \beta} G(\gamma)$$
$$ G'(\beta) \rightarrow \lim_{\gamma < \beta} G'(\gamma)$$
are fibrations, while the map $G(\beta) \rightarrow G'(\beta)$ is a weak equivalence.
Then the induced map $$\lim_{\gamma < \alpha} G(\gamma) \rightarrow \lim_{\gamma < \alpha} G'(\gamma)$$ is a weak equivalence.
\end{itemize}
\end{corollary}

\begin{proof}
We will prove $(1)$; $(2)$ follows by the same argument. Let $p: (\alpha) \rightarrow \ast$ be
the unique map, let $p^{\ast}: \bfA \rightarrow \bfA^{(\alpha)}$ be the diagonal map, and let
$p_{!}: \bfA^{(\alpha)} \rightarrow \bfA$ be a left adjoint to $p^{\ast}$. Then $p_{!}$ can be identified with the functor $F \mapsto \colim_{\gamma < \alpha} F(\gamma)$. We observe that $(p_!, p^{\ast})$ is a Quillen adjunction (where
$\bfA^{(\alpha)}$ is endowed with the projective model structure) so that
$p_{!}$ preserves weak equivalence between projectively cofibrant objects. The desired result now follows from Corollary \ref{dimbu}.
\end{proof}

Suppose that we are given a bifunctor
$$ \otimes: \bfA \times \bfB \rightarrow \bfC, $$
where $\bfC$ is a category which admits small limits. For any smallcategory $\calJ$, we define the
{\em coend} functor $\int_{\calJ}: \Fun(\calJ, \bfA) \times \Fun( \calJ^{op}, \bfB) \rightarrow \bfC$
so that the integral $\int_{\calJ}(F,G)$ is defined to be the coequalizer of the diagram
$$\xymatrix{ \coprod_{ J \rightarrow J'} F(J) \otimes G(J') 
\ar@<.4ex>[r] \ar@<-.4ex>[r] & \coprod_{J} F(J) \otimes G(J)}.$$
We then have the following result:

\begin{proposition}\label{intreed}
Let $\otimes: \bfA \times \bfB \rightarrow \bfC$ be a left Quillen bifunctor
(see Proposition \ref{biquill}), and let $\calJ$ be a Reedy category. Then the
coend functor
$$ \int_{\calJ}: \Fun( \calJ, \bfA) \times \Fun( \calJ^{op}, \bfB) \rightarrow \bfC$$
is also a left Quillen bifunctor, where we regard $\Fun(\calJ, \bfA)$ and
$\Fun(\calJ^{op}, \bfB)$ as endowed with the Reedy model structure.
\end{proposition}

\begin{proof}
Let $f: F \rightarrow F'$ be a Reedy cofibration in
$\Fun(\calJ, \bfA)$ and $g: G \rightarrow G'$ a Reedy cofibration in
$\Fun(\calJ^{op}, \bfB)$. Set $C = \int_{\calJ}(F,G') \coprod_{ \int_{\calJ}(F,G) } \int_{\calJ}(F',G) \in \bfC$, and $C' = \int_{\calJ}(F',G')$. We wish to show that the induced map
$C \rightarrow \int_{\calJ}(F',G')$ is a cofibration, which is trivial if either
$f$ or $g$ is trivial. 

Choose a good filtration $\{ \calJ_{\beta} \}_{\beta < \alpha}$ of $\calJ$.
For $\beta \leq \alpha$, we define
$$ C_{\beta} = \int_{ \calJ_{\beta}}( F|\calJ_{<\beta}, G'|\calJ_{<\beta})
\coprod_{ \int{\calJ_{\beta}}(F| \calJ_{<\beta}, G|\calJ_{<beta})} \int_{\calJ_{\beta}}( F'|\calJ_{<\beta}, G| \calJ_{<\beta})$$
$$ C'_{\beta} = \int_{\calJ_{\beta}}( F'| \calJ_{<\beta}, G'| \calJ_{<\beta}).$$
We wish to show that the map 
$$C_{\alpha} \simeq C_{\alpha} \coprod_{ C_0} C'_0 \rightarrow C_{\alpha}
\coprod_{ C_{\alpha} } C'_{\alpha}$$
is a cofibration (which is trivial if either $f$ or $g$ is trivial). We will prove more generally
that for $\delta \leq \gamma \leq \beta \leq \alpha$, the map
$$ \eta_{\delta, \gamma, \beta}: C_{\beta} \coprod_{ C_{\delta} } C'_{\delta} \rightarrow
C_{\beta} \coprod_{ C_{\gamma} } C'_{\gamma}$$
is a cofibration (trivial if either $f$ or $g$ is trivial).
The proof proceeds by induction on $\gamma$. If $\gamma$ is a limit ordinal,
then $\eta_{\delta, \gamma, \beta}$ can be obtained as a transfinite composition of the maps
$\{ \eta_{\epsilon, \epsilon+1, \beta} \}_{ \delta \leq \epsilon < \gamma}$, and the result follows from the inductive hypothesis. We may therefore assume that $\gamma = \gamma_0 +1$ is a successor ordinal. Since $\eta_{\delta, \gamma, \beta} = \eta_{ \gamma_0, \gamma, \beta} \circ \eta_{ \delta, \gamma_0, \beta}$, we can use the inductive hypothesis to reduce to the case where $\delta = \gamma_0$. Since $\eta_{\delta, \gamma, \beta}$ is a pushout of
$\eta_{\delta, \gamma, \gamma}$, we can assume also that $\beta = \gamma$. In other words, we are reduced to proving that the map
$$ h: C_{\gamma_0+1} \coprod_{ C_{\gamma_0} } C'_{\gamma_0} \rightarrow C'_{\gamma_0}$$
is a cofibration, which is trivial if either $f$ or $g$ is trivial. Let $J$ be the object of
$\calJ_{\gamma_0}$ which does not belong to $\calJ_{< \gamma_0}$. We now observe
that $h$ is a pushout of the evident map from
$$ ((F(J) \coprod_{L_J(F)} L_J(F') ) \otimes G'(J))
\coprod_{ (F(J) \coprod_{L_J(F)} L_J(F')) \otimes (G(J) \coprod_{L_J(G)} L_J(G'))}
(F'(J) \otimes (G(J) \coprod_{ L_J(G)} L_J(G'))$$
to $F'(J) \otimes G'(J),$
which is a cofibration (trivial if either $f$ or $g$ is trivial) by virtue of our assumptions on $f$, $g$, and the fact that $\otimes$ is a left Quillen bifunctor.
\end{proof}

\begin{remark}\label{cabler}
Proposition \ref{intreed} has an analogue for the model structures introduced in Proposition \ref{smurff}. That is, suppose that $\bfA$ and $\bfB$ are {\em combinatorial} model categories, and let $\calJ$ be an arbitrary small category. Then any left Quillen bifunctor
$\otimes: \bfA \times \bfB \rightarrow \bfC$ induces a left Quillen bifunctor
$$ \int_{ \calJ}: \Fun( \calJ, \bfA) \times \Fun( \calJ^{op}, \bfB) \rightarrow \bfC,$$
where we regard $\Fun(\calJ, \bfA)$ as endowed with the projective model structure
and $\Fun(\calJ^{op}, \bfB)$ with the injective model structure. To see this, we must show that 
for any projective cofibration $f: F \rightarrow F'$ in $\Fun(\calJ, \bfA)$ and any injective cofibration
$g: G \rightarrow G'$ in $\Fun(\calJ^{op}, \bfB)$, the induced map 
$$ h: \int_{\calJ}(F,G') \coprod_{ \int_{\calJ}(F,G) } \int_{\calJ}(F',G) \rightarrow \int_{\calJ}(F',G')$$
is a cofibration in $\bfC$, which is trivial if either $f$ or $g$ is trivial. Without loss of generality,
we may suppose that $f$ is a generating projective cofibration of the form
$\calF^{J}_{A} \rightarrow \calF^{J}_{A'}$ associated to an object $J \in \calJ$ and a
cofibration $i: A \rightarrow A'$ in $\bfA$, which is trivial if $f$ is trivial (see
the proof of Proposition \ref{smurff} for an explanation of this notation). Unwinding the definitions, we can identify $h$ with the map
$$ (A \otimes G'(J)) \coprod_{ A \otimes G(J) } (A' \otimes G(J)) \rightarrow
A' \otimes G'(J).$$
Since $i$ is a cofibration in $\bfA$ and the map $G(J) \rightarrow G'(J)$
is a cofibration in $\bfB$, we deduce that $h$ is a cofibration
in $\bfC$ (since $\otimes$ is a left Quillen bifunctor) which is trivial if either
$i$ or $h$ is trivial.
\end{remark}

\begin{example}\label{cabletome}
Let $\bfA$ be a simplicial model category, so that we have a left Quillen bifunctor
$$ \otimes: \bfA \times \sSet \rightarrow \bfA.$$
The coend construction determines a left Quillen bifunctor
$$ \int_{\cDelta}: \Fun( \cDelta, \bfA) \times \Fun( \cDelta^{op}, \sSet) \rightarrow \bfA.$$
where $\Fun(\cDelta, \bfA)$ and $\Fun( \cDelta^{op}, \sSet)$ are both endowed with the Reedy model structure. In particular, if we fix a cosimplicial object
$X^{\bigdot} \in \Fun(\cDelta, \bfA)$ which is Reedy cofibrant, then forming the coend
against $X^{\bigdot}$ determines a left Quillen functor from the category
of bisimplicial sets (with the Reedy model structure, which coincides with the
injective model structure by Example \ref{tetsu}) to $\bfA$.
\end{example}

\begin{example}\label{seventwo}
Let $\bfA$ be a simplicial model category, so that we have a left Quillen bifunctor
$$ \otimes: \bfA \times \sSet \rightarrow \bfA,$$
and consider the coend functor
$$ \int_{ \cDelta^{op} } \Fun( \cDelta^{op}, \bfA) \times \Fun( \cDelta, \sSet) \rightarrow \bfA.$$
Let $\Delta^{\bigdot} \in \Fun( \cDelta, \sSet)$ denote the standard simplex
(that is, the functor $[n] \mapsto \Delta^n$), and let ${\bf 1}$ denote the final object
of $\Fun(\cDelta, \sSet)$ (that is, the constant functor given by $[n] \mapsto \Delta^0$).
The unique map $\Delta^{\bigdot} \rightarrow {\bf 1}$ is a weak equivalence, 
and $\Delta^{\bigdot}$ is Reedy cofibrant: we may therefore regard $\Delta^{\bigdot}$ as a cofibrant replacement for the constant functor ${\bf 1}$.

The functor $X_{\bigdot} \mapsto \int_{ \cDelta^{op}}(X_{\bigdot}, {\bf 1})$ can be identified with the colimit functor $\Fun( \cDelta^{op}, \bfA) \rightarrow \bfA$. This is a left Quillen functor if
$\Fun( \cDelta^{op}, \bfA)$ is endowed with the projective model structure, but not the Reedy model structure. However, the {\it geometric realization} functor $X_{\bigdot} \mapsto |X_{\bigdot} | = \int_{ \cDelta^{op}}(X_{\bigdot}, \Delta^{\bigdot})$ is a left Quillen functor with respect to the Reedy model structure. 
\end{example}

\begin{corollary}\label{twinner}
Let $\bfA$ be a combinatorial simplicial model category, and let $X_{\bigdot}$ be a simplicial object
of $\bfA$. There is a canonical map
$$ \gamma: \hocolim X_{\bigdot} \rightarrow | X_{\bigdot} |$$
in the homotopy category of $\bfA$. This map is an equivalence if $X_{\bigdot}$ is Reedy cofibrant.
\end{corollary}

\begin{proof}
Let $\Delta^{\bigdot}$ and $\ast$ be the cosimplicial objects of $\sSet$ described in Example
\ref{seventwo}. Choose a weak equivalence of simplicial objects $X'_{\bigdot} \rightarrow X_{\bigdot}$, where $X'_{\bigdot}$ is projectively cofibrant. We then have a diagram
$$ \hocolim X_{\bigdot} \simeq \colim X'_{\bigdot}
\simeq \int_{ \cDelta^{op}}(X'_{\bigdot}, \ast)
\stackrel{\alpha}{\leftarrow} \int_{ \cDelta^{op}}( X'_{\bigdot}, \Delta^{\bigdot})
\stackrel{\beta}{\rightarrow} \int_{ \cDelta^{op}}( X_{\bigdot}, \Delta^{\bigdot}).$$
Since $X'_{\bigdot}$ is projectively cofibrant, Remark \ref{cabler} implies that the 
coend functor $\int_{\cDelta^{op}}( X'_{\bigdot}, \bigdot)$ preserves weak equivalences between
injectively cofibrant cosimplicial objects of $\sSet$; in particular, $\alpha$ is a weak equivalence in $\bfA$. This gives the desired map $\gamma$. Proposition \ref{intreed} implies that $\int_{ \cDelta^{op}}( \bigdot, \Delta^{\bigdot})$ preserves weak equivalences between Reedy cofibrant simplicial objects of $\bfA$, which proves that $\gamma$ is an isomorphism if $X_{\bigdot}$ is Reedy cofibrant.
\end{proof}

\begin{example}\label{swupt}
If $\bfA$ is the category of simplicial sets, then the map $\gamma$ of Corollary \ref{twinner}
is always an isomorphism; this follows from Example \ref{tetsu}. In other words, if
$X_{\bigdot, \bigdot}$ is a bisimplicial set, then we can identify the diagonal simplicial set
$[n] \mapsto X_{n,n}$ with the homotopy colimit of corresponding diagram
$\cDelta^{op} \rightarrow \sSet$.
\end{example}

\section{Simplicial Categories}\label{techapp}

Among the many different models for higher category theory, the theory of simplicial categories
is perhaps the most rigid. This can be either a curse or a blessing, depending on the situation. 
For the most part, we have chosen to use the less rigid theory of $\infty$-categories (see \S \ref{qqqc}) throughout this book. However, some arguments are substantially easier to carry out in the setting of simplicial categories. For this reason, we have devoted the final section of this appendix to giving a review of the theory of simplicial categories.

There exists a model structure on the category $\sCat$ of (small) simplicial categories, which was constructed by Bergner (\cite{bergner}). In \S \ref{compp4}, we will describe an analogous model structure on the category $\SCat$ of $\bfS$-enriched categories, where $\bfS$ is a suitable model category.
To formulate this generalization, we will need to employ the language of monoidal model categories, which we review in \S \ref{maymy}. Under mild assumptions on $\bfS$, one can show that a $\bfS$-enriched
category $\calC$ is fibrant if and only if each of the mapping objects $\bHom_{\calC}(X,Y)$ is a fibrant object of $\bfS$.\index{not}{CATsubS@$\SCat$}

In \S \ref{quasilimit3}, we will study the category $\bfA^{\calC}$ of diagrams $\calC \rightarrow \bfA$, where
$\calC$ is a small category and $\bfA$ a model category, both enriched over some
fixed model category $\bfS$. In the enriched setting we can again endow $\bfA^{\calC}$
with projective and injective model structures, which can be used to define homotopy limits and
colimits. 

Putting aside set-theoretic technicalities, every $\bfS$-enriched model category $\bfA$
gives rise to a fibrant object of $\SCat$: namely, the full subcategory
$\bfA^{\degree} \subseteq \bfA$ spanned by the fibrant-cofibrant objects. In \S \ref{pathspace}, we
will introduce a path object for $\bfA^{\degree}$, which will enable us to perform some calculations in
the homotopy category of $\SCat$. 

In \S \ref{hoco}, we will consider the problem of constructing homotopy colimits in
the category $\SCat$ of $\bfS$-enriched categories. Our main result, Theorem \ref{tubba}, asserts
that the formation of homotopy colimits in $\SCat$ is compatible with the formation of
(tensor) products in $\SCat$. We will apply this result in \S \ref{camper} to study the homotopy
theory of internal mapping objects in $\SCat$. 

We conclude this section with \S \ref{turka}, where we discuss localizations of (simplicial) model categories.


\subsection{Enriched and Monoidal Model Categories}\label{maymy}

Many of the model categories which arise naturally are {\em enriched} over the category
of simplicial sets. Our goal in this section to study enrichments of one model category over another.

\begin{definition}\index{gen}{Quillen bifunctor}\index{gen}{left Quillen bifunctor}\label{biquill}
Let $\bfA$, $\bfB$, and $\bfC$ be model categories. We will say that a functor
$F: \bfA \times \bfB \rightarrow \bfC$ is a {\em left Quillen bifunctor} if the following conditions
are satisfied:
\begin{itemize}
\item[$(a)$] Let $i: A \rightarrow A'$ and $j: B \rightarrow B'$ be cofibrations in $\bfA$ and
$\bfB$, respectively. Then the induced map
$$i \wedge j: F(A',B) \coprod_{F(A,B)} F(A,B') \rightarrow F(A',B')$$
is a cofibration in $\bfC$. Moreover, if either $i$ or $j$ is a trivial cofibration, then
$i \wedge j$ is also a trivial cofibration.
\item[$(b)$] The functor $F$ preserves small colimits separately in each variable.
\end{itemize}
\end{definition}

\begin{definition}\index{gen}{monoidal model category}\index{gen}{model category!monoidal}
A {\it monoidal model category} is a monoidal category $\bfS$ equipped with a model structure, which satisfies the following conditions:
\begin{itemize}
\item[$(i)$] The tensor product functor $\otimes: \bfS \times \bfS \rightarrow \bfS$ is a left
Quillen bifunctor.
\item[$(ii)$] The unit object ${\bf 1} \in \bfS$ is cofibrant.
\item[$(iii)$] The monoidal structure on $\bfS$ is closed.
\end{itemize}
\end{definition}

\begin{remark}
Some authors demand only a weakened form of axiom $(ii)$ in the preceding definition.
\end{remark}

\begin{example}\label{shootset}
The category of simplicial sets $\sSet$ is a monoidal model category, with
respect to the Cartesian product and the Kan model structure defined in \S \ref{simpset}.
\end{example}

\begin{definition}\index{gen}{model category!simplicial}\index{gen}{simplicial model category}
\index{gen}{model category!$\bfS$-enriched}
Let $\bfS$ be a monoidal model category. A {\it $\bfS$-enriched model category}
is a $\bfS$-enriched category $\bfA$ equipped with a model structure satisfying the following conditions:
\begin{itemize}
\item[$(1)$] The category $\bfA$ is tensored and cotensored over $\bfS$.
\item[$(2)$] The tensor product $\otimes: \bfA \times \bfS \rightarrow \bfA$ is a left Quillen bifunctor.
\end{itemize}
In the special case where $\bfS$ is the category of simplicial sets (regarded
as a monoidal model category as in Example \ref{shootset}), we will simply refer to
$\bfA$ as a {\it simplicial model category}.
\end{definition}

\begin{remark}\label{cyclor}
An easy formal argument shows that condition $(2)$ is equivalent to either of the following
statements:
\begin{itemize}
\item[$(2')$] Given any cofibration $i: D \rightarrow D'$ in $\bfA$ and any fibration
$j: X \rightarrow Y$ in $\bfA$, the induced map
$$k: \bHom_{\bfA}(D',X) \rightarrow \bHom_{\bfA}(D,X) \times_{ \bHom_{\bfA}(D,Y)} \bHom_{\bfA}(D',Y)$$ is fibration in $\bfS$, which is trivial if either $i$ or $j$ is a weak equivalence.
\item[$(2'')$] Given any cofibration $i: C \rightarrow C'$ in $\bfS$ and any fibration $j: X \rightarrow Y$ in $\bfA$, the induced map $$k: X^{C'} \rightarrow X^C \times_{ Y^C} Y^{C'}$$ is a fibration in $\bfA$, which is trivial if either $i$ or $j$ is trivial.
\end{itemize}
\end{remark}

The following provides a criterion for detecting simplicial model structures:

\begin{proposition}\label{testsimpmodel}
Let $\calC$ be a simplicial category that is equipped with a model structure $($not assumed to be compatible with the simplicial structure on $\calC$$)$. Suppose that every object of $\calC$ is cofibrant  and that the collection of weak equivalences in $\calC$ is stable under filtered colimits. Then $\calC$ is a simplicial model category if and only if the following conditions are satisfied:

\begin{itemize}
\item[$(1)$] As a simplicial category, $\calC$ is both tensored and cotensored over $\sSet$. 
\item[$(2)$] Given a cofibration $i: A \rightarrow B$ and a fibration $p: X \rightarrow Y$ in $\calC$, the induced map of simplicial sets
$$ q: \bHom_{\calC}(B, X) \rightarrow \bHom_{\calC}(A,X) \times_{ \bHom_{\calC}(A,Y)} \bHom_{\calC}(B,Y) $$
is a Kan fibration. 
\item[$(3)$] For every $n \geq 0$ and every object $C$ in $\calC$, the natural map
$$ C \otimes \Delta^n \rightarrow C \otimes \Delta^0 \simeq C$$
is a weak equivalence in $\calC$.
\end{itemize}
\end{proposition}

\begin{proof}
Suppose first that $\calC$ is a simplicial model category. It is clear that $(1)$ and $(2)$ are satisfied. To prove $(3)$, we note that the projection $\Delta^n \rightarrow \Delta^0$ admits a section $s: \Delta^0 \rightarrow \Delta^n$ which is a trivial cofibration. If $\calC$ is a simplicial model category, then since $C$ is cofibrant it follows
that $C \otimes \Delta^0 \rightarrow C \otimes \Delta^n$ is a trivial cofibration, and in particular a weak equivalence. Thus the projection
$C \otimes \Delta^n \rightarrow C \otimes \Delta^0$ is a weak equivalence by the two-out-of-three property.

Now suppose that $(1)$, $(2)$, and $(3)$ are satisfied. We wish to show that $\calC$ is a simplicial model category. 
We first show that the bifunctor
$$ (C,K) \mapsto C \otimes K$$
preserves weak equivalences separately in each variable separately.

Fix the object $C \in \calC$, and suppose that $f: K \rightarrow K'$ is a weak homotopy equivalence of simplicial sets. Choose a cofibration $K \rightarrow K''$, where $K''$ is a contractible Kan complex. Then we may factor $f$ as a composition
$$ K \stackrel{f'}{\rightarrow} K \times K'' \stackrel{f''}{\rightarrow} K.$$
To prove that $\id_{C} \otimes f$ is a weak equivalence, it suffices to prove that
$\id_{C} \otimes f'$ and $\id_{C} \otimes f''$ are weak equivalences. Note that the map
$f''$ has a section $s$, which is a trivial cofibration. Thus, to prove that $\id_{C} \otimes f''$ is a weak equivalence, it suffices to show that $\id_{C} \otimes s$ is a weak equivalence. In other words, we may reduce to the case where $f$ is itself a trivial cofibration of simplicial sets.

Consider the collection $A$ of all monomorphisms $f: K \rightarrow K'$ of simplicial sets having the property that $\id_{C} \otimes f$ is a weak equivalence in $\calC$. It is easy to see that this collection of morphisms is weakly saturated. Thus, to prove that it contains all trivial cofibrations of simplicial sets, it suffices to show that every horn inclusion $\Lambda^n_i \rightarrow \Delta^n$
belongs to $A$. We prove this by induction on $n > 0$. Choose a vertex $v$ belonging to
$\Lambda^n_i$. We note that the inclusion $\{v\} \rightarrow \Lambda^n_i$ is a pushout of horn inclusions in dimensions $< n$; by the inductive hypothesis, this inclusion belongs to $A$.
Thus, it suffices to show that $\{v\} \rightarrow \Delta^n$ belongs to $A$, which is equivalent to assumption $(3)$.

Now let us show that for each simplicial set $K$, the functor
$$ C \mapsto C \otimes K$$ preserves weak equivalences. We will prove this by induction on the (possibly infinite) dimension of $K$. Choose a weak equivalence
$g: C \rightarrow C'$ in $\calC$. Let $S$ denote the collection of all simplicial subsets $L \subseteq K$ such that $g \otimes \id_{L}$ is a weak equivalence. We regard $S$ as a partially ordered set with respect to inclusions of simplicial subsets. Clearly $\emptyset \in S$. Since weak equivalences in $\calC$ are stable under filtered colimits, the supremum of every chain in $S$ belongs to $S$.
By Zorn's lemma, $S$ has a maximal element $L$. We wish to show that $L = K$. If not, 
we may choose some nondegenerate simplex $\sigma$ of $K$ which does not belong to $L$. Choose $\sigma$ of the smallest possible dimension, so that all of the faces of $\sigma$ belong to $L$. Thus, there is an inclusion $L' = L \coprod_{ \bd \sigma} \sigma \subseteq K$. Since
$\calC$ is left proper, assumption $(2)$ implies that the diagram
$$ \xymatrix{ D \otimes \bd \sigma \ar[r] \ar[d] & D \otimes \sigma \ar[d] \\
D \otimes L \ar[r] & D \otimes L' }$$
is a homotopy pushout, for every object $D \in \calC$. We observe that
$g \otimes \id_{L}$ is a weak equivalence by assumption, $g \otimes \id_{\bd \sigma}$ is
a weak equivalence by the inductive hypothesis (since $\bd \sigma$ has dimension smaller than the dimension of $K$), and $g \otimes \id_{\sigma}$ is a weak equivalence in virtue of assumption $(3)$ and the fact that $g$ is a weak equivalence. It follows that $g \otimes \id_{L'}$ is a weak equivalence, which contradicts the maximality of $L$. This completes the proof that the bifunctor
$\otimes: \calC \times \sSet \rightarrow \calC$ preserves weak equivalences separately in each variable.

Now suppose given a cofibration
$i: C \rightarrow C'$ in $\calC$ and another cofibration $j: S \rightarrow S'$ in $\sSet$.
We wish to prove that the induced map 
$$i \wedge j: (C \otimes S') \coprod_{C \otimes S} (C' \otimes S) \rightarrow C' \otimes S'$$
is a cofibration in $\calC$, which is trivial if either $i$ or $j$ is trivial. The first point follows immediately from $(2)$. For the triviality, we will assume that $i$ is a weak equivalence (the case where $j$ is a weak equivalence follows using the same argument).
Consider the diagram
$$ \xymatrix{ C \otimes S \ar[rr]^{i \otimes \id_{S}} \ar[d] & & C' \otimes S \ar[d] & \\
C \otimes S' \ar[rr]^-{f} & & (C' \otimes S) \coprod_{ C \otimes S} (C \otimes S') \ar[r] & C' \otimes S'.}$$
The arguments above show that $i \otimes \id_{S}$ and $i \otimes \id_{S'}$ are weak equivalences.
The square in the diagram is a homotopy pushout, so Proposition \ref{propob} implies that $f$ is a weak equivalence as well. Thus $i \wedge j$ is a weak equivalence, by the two-out-of-three property.
\end{proof}

If $\calC$ is a simplicial model category, then there is automatically a strong relationship between the homotopy theory of the underlying model category and the homotopy theory of the simplicial sets $\bHom_{\calC}(\bigdot, \bigdot)$. For example, we have the following:

\begin{remark}
Let $\calC$ be a simplicial model category, let $X$ be a cofibrant object of $\calC$, and let $Y$ be a fibrant object of $\calC$. The simplicial set $K=\bHom_{\calC}(X,Y)$ is a Kan complex; moreover, 
there is a canonical bijection
$$ \pi_0 K \simeq \Hom_{ \h{ \calC}}(X,Y).$$
\end{remark}

We conclude this section by studying a situation which will arise in \S \ref{chap4}. 
Let $\calC$ and $\calD$ be model categories enriched over another model category $\bfS$, and suppose given a Quillen adjunction $$\Adjoint{F}{\calC}{\calD}{G}$$
between the underlying model categories. We wish to study the situation where $G$ (but not $F$) has the structure of a $\bfS$-enriched functor. Thus, for every triple of objects
$X \in \calC$, $Y \in \calD$, $S \in \bfS$, we have a canonical map
\begin{eqnarray*}
\Hom_{\calC}( S \otimes X, GY) & \simeq & \Hom_{\bfS}( S, \bHom_{\calC}( X, GY) ) \\
& \rightarrow &  \Hom_{\bfS}( S, \bHom_{\calD}(FX, FGY) ) \\
& \simeq & \Hom_{\calD}( S \otimes FX, FGY) \\
& \rightarrow & \Hom_{\calD}( S \otimes FX, Y). \end{eqnarray*}
Taking $Y = F(S \otimes X)$ and applying this map to the unit of the adjunction between
$F$ and $G$, we obtain a map
$S \otimes FX \rightarrow F(S \otimes X)$, which we will denote by $\beta_{X,S}$.
The collection of maps $\beta_{X,S}$ is simply another way of encoding the data of $G$
as a $\bfS$-enriched functor. If the maps $\beta_{X,S}$ are isomorphisms, then
$F$ is again a $\bfS$-enriched functor, and $(F,G)$ is an adjunction between $\bfS$-enriched categories. We wish to study an analogous situation, where the maps $\beta_{X,S}$ are only assumed to be weak equivalences.

\begin{remark}\label{tuccan}
Suppose that $\bfS$ is the category $\sSet$ of simplicial sets, with its usual model structure.
Then the map $\beta_{X,S}$ is automatically a weak equivalence for every cofibrant object
$X \in \calC$. To prove this, we consider the collection $\calK$ of all simplicial sets
$S$ such that $\beta_{S,X}$ is an equivalence. It is not difficult to show that $\calK$ is
closed under weak equivalences, homotopy pushout squares, and coproducts.
Since $\Delta^0 \in \calK$, we conclude that $\calK = \sSet$.
\end{remark}

\begin{proposition}\label{weakcompatequiv}
Let $\calC$ and $\calD$ be $\bfS$-enriched model categories. Let
$\Adjoint{F}{\calC}{\calD}{G}$ be a Quillen adjunction between the underlying model categories.
Assume that every object of $\calC$ is cofibrant,
and that the map $\beta_{X,S}: S \otimes F(X) \rightarrow F(S \otimes X)$ is a weak equivalence for every pair of cofibrant objects $X \in \calC$, $S \in \bfS$. The following are equivalent:
\begin{itemize}
\item[$(1)$] The adjunction $(F,G)$ is a Quillen equivalence.
\item[$(2)$] The restriction of $G$ determines a weak equivalence of $\bfS$-enriched categories
$\calD^{\degree} \rightarrow \calC^{\degree}$ (see \S \ref{compp4}).
\end{itemize}
\end{proposition}

\begin{remark}
Strictly speaking, in \S \ref{compp4} we only define weak equivalences between {\em small} $\bfS$-enriched categories; however, the definition extends to large categories in an obvious way.
\end{remark}

\begin{proof}
Since $G$ preserves fibrant objects, and every object of $\calC$ is cofibrant, it is clear that
$G$ carries $\calD^{\degree}$ into $\calC^{\degree}$. Condition $(1)$
is equivalent to the assertion that for every pair of fibrant-cofibrant objects
$C \in \calC$, $D \in \calD$, a map
$g: C \rightarrow GD$ is a weak equivalence in $\calC$ if and only if the adjoint map
$f: FC \rightarrow D$ is a weak equivalence in $\calD$. Choose a factorization of $f'$ as a composition
$FC \stackrel{f'}{\rightarrow} D' \stackrel{f''}{\rightarrow} D,$
where $f'$ is a trivial cofibration and $f''$ is a fibration. By the two-out-of-three property,
$f$ is a weak equivalence if and only if $f''$ is a weak equivalence. We note that
$g$ admits an analogous factorization as
$$ C \stackrel{g'}{\rightarrow} GD' \stackrel{g''}{\rightarrow} GD.$$
Using $(2)$, we deduce that $f''$ is an equivalence in $\calD^{\degree}$ if and only if
$g''$ is an equivalence in $\calC^{\degree}$. It will therefore suffice to show that
$g'$ is an equivalence in $\calC^{\degree}$. For this, it will suffice to show that
$C$ and $GD'$ corepresent the same functor on the homotopy category $\h{\calC}$.
Invoking $(2)$ again, it will suffice to show that for every fibrant-cofibrant object
$D'' \in \calD$, the induced map
$$\Hom_{ \h{\calC} }( GD', GD'') \rightarrow \Hom_{\h{\calC}}( C, GD'')
\simeq \Hom_{\h{\calD}}( FC, D'')$$ is bijective. Using $(2)$, 
we deduce that map $\Hom_{\h{\calD}}(D',D'') \rightarrow \Hom_{\h{\calD}}( GD', GD'')$
is bijective. The desired result now follows from the fact that $f'$ is a weak equivalence in $\calD$.

We now show that $(1) \Rightarrow (2)$.
The $\bfS$-enriched functor $G^{\degree}: \calD^{\degree} \rightarrow \calC^{\degree}$ is essentially surjective, since the right derived functor $RG$ is essentially surjective on homotopy categories. It suffices to show that $G^{\degree}$ is fully faithful: in other words, that for every pair of fibrant-cofibrant objects $X,Y \in \calD$, the induced map
$$ i: \bHom_{\calD}(X,Y) \rightarrow \bHom_{\calC}(G(X),G(Y))$$ 
is a weak equivalence in $\bfS$.

Since the left derived functor $LF$ is essentially surjective, there exists an object $X' \in \calC$
and a weak equivalence $FX' \rightarrow X$. We may regard $X$ as a fibrant replacement for $FX'$ in $\calD$; it follows that the adjoint map $X' \rightarrow GX$ may be identified with the
adjunction $X' \rightarrow (RG \circ LF) X'$, and is therefore a weak equivalence by $(1)$. Thus we have a diagram
$$ \xymatrix{ \bHom_{\calD}(X,Y) \ar[d] \ar[r]^{i} & \bHom_{\calC}(G(X),G(Y)) \ar[d] \\
\bHom_{\calD}(F(X'), Y) \ar[r]^{i'} & \bHom_{\calC}(X',G(Y)) }$$
in which the vertical arrows are homotopy equivalences; thus, to show that $i$ is a weak equivalence, it suffices to show that $i'$ is a weak equivalence. For this, it suffices to show that $i'$ induces a bijection from $[S, \bHom_{\calD}(F(X'), Y) ]$ to $[S, \bHom_{\calC}(X',G(Y))]$, for every cofibrant object $S \in \bfS$; here $[S,K]$ denotes the set of homotopy classes of maps from $S$ into $K$ in
the homotopy category $\h{\bfS}$. But we may rewrite this map of sets as
$$i'_S:  \bHom_{\h{\calD}}( F(X') \otimes S, Y) \rightarrow \bHom_{\h{\calC}}( X' \otimes S, G(Y) )
= \bHom_{\h{\calD}}( F( X' \otimes S), Y),$$
and it is given by composition with $\beta_{X',S}$. (Here $\h{\calC}$ and $\h{\calD}$ denote the
homotopy categories of $\calC$ and $\calD$ as {\em model categories}; these are equivalent
to the corresponding homotopy categories of $\calC^{\degree}$ and $\calD^{\degree}$ as $\bfS$-enriched categories). Since $\beta_{X',S}$ is an isomorphism in the homotopy category $\h{\calD}$, the map $i'_S$ is bijective and $(2)$ holds, as desired.
\end{proof}

\begin{corollary}\label{urchug}
Let $\Adjoint{F}{\calC}{\calD}{G}$
be a Quillen equivalence between simplicial model categories, where every object of $\calC$ is cofibrant. Suppose that $G$ is a simplicial functor. Then $G$ induces an equivalence of $\infty$-categories $\Nerve( \calD^{\degree}) \rightarrow \Nerve( \calC^{\degree})$.
\end{corollary}

\subsection{The Model Structure on $\bfS$-Enriched Categories}\label{compp4}

Throughout this section, we will fix a symmetric monoidal model category $\bfS$, and
study the category of $\bfS$-enriched categories. The main case of interest to us is that in which $\bfS$ is the category of simplicial sets (with its usual model structure and the Cartesian monoidal structure). However, the treatment of the general case requires little additional effort, and there are a number of other examples which arise naturally in other contexts:

\begin{itemize}
\item[$(i)$] The category $\sSet$ of simplicial sets, equipped with the Cartesian monoidal structure and the {\em Joyal} model structure defined in \S \ref{compp3}.
\item[$(ii)$] The category of complexes
$$ \ldots \rightarrow M_{n} \rightarrow M_{n} \rightarrow M_{n-1} \rightarrow \ldots, $$
of vector spaces over a field $k$, with its usual model structure (in which weak equivalences are quasi-isomorphisms, fibrations are epimorphisms, and cofibrations are monomorphisms) and monoidal structure given by the formation of tensor products of complexes.
\end{itemize}

Let $\bfS$ be an monoidal model category, and let $\SCat$ denote the category of (small) $\bfS$-enriched categories, in which morphisms are
given by $\bfS$-enriched functors. The goal of this section is to describe a model structure on $\SCat$. 
We first note that the monoidal structure on $\bfS$ induces a monoidal structure on its
homotopy category $\h{\bfS}$, which is determined up to (unique) isomorphism by the requirement that there exist a monoidal structure on the functor
$$ \bfS \rightarrow \h{\bfS} $$
given by inverting all weak equivalences. Consequently, 
we note that any $\bfS$-enriched category $\calC$ gives rise to an $\h{\bfS}$-enriched category $\h{\calC}$, having the same objects as $\calC$ and where mapping spaces are given by
$$ \bHom_{\h{\calC}}(X,Y) = [ \bHom_{\calC}(X,Y) ].$$
Here we let $[K]$ denote the image in $\h{\bfS}$ of an object $K \in \bfS$. We will refer to
$\h{\calC}$ as the {\it homotopy category} of $\calC$; the passage from
$\calC$ to $\h{\calC}$ is a special case of Remark \ref{laxcon}.\index{not}{hcalC@$\h{\calC}$}
\index{gen}{homotopy category!of a $\bfS$-enriched category}

\begin{definition}\label{equivdefequiv}\index{gen}{equivalence!of $\bfS$-enriched categories}
Let $\bfS$ be an monoidal model category.
We say that a functor $F: \calC \rightarrow \calC'$ in $\SCat$ is a {\it weak equivalence}
if the induced functor $\h{\calC} \rightarrow \h{\calC'}$ is an equivalence of $\h{\bfS}$-enriched categories. In other words, $F$ is a weak equivalence if and only if:

\begin{itemize}
\item[$(1)$] For every pair of objects $X,Y \in \calC$, the induced map
$$ \bHom_{\calC}(X,Y) \rightarrow \bHom_{\calC'}(F(X),F(Y))$$ is a weak equivalence in $\bfS$.
\item[$(2)$] Every object $Y \in \calC'$ is equivalent to $F(X)$ in the homotopy category $\h{\calC'}$, for some $X \in \calC$.
\end{itemize}
\end{definition}

\begin{remark}
If $\bfS$ is the category $\sSet$ (endowed with the Kan model structure), then Definition \ref{equivdefequiv} reduces to the definition given in \S \ref{stronghcat}.
\end{remark}

\begin{remark}\label{swunk}
Suppose that the collection of  weak equivalences in $\bfS$ is stable under filtered colimits.
Then it is easy to see that the collection of weak equivalences in $\SCat$ is also stable under filtered colimits. If $\bfS$ is also a combinatorial model category, then a bit more effort shows that
the class of weak equivalences in $\SCat$ is perfect, in the sense of Definition \ref{perfequiv}.
\end{remark}

We now introduce a bit of notation for working with $\bfS$-enriched categories. If $A$ is an object of $\bfS$, we will let $[1]_{A}$ denote the $\bfS$-enriched category having two objects $X$ and $Y$, with
$$ \bHom_{[1]_{A}}(Z,Z') = \begin{cases} {\bf 1}_{\bfS} & \text{if } Z=Z'=X \\
{\bf 1}_{\bfS} & \text{if } Z=Z'=Y \\
A & \text{if } Z=X, Z'=Y \\
\emptyset & \text{if } Z=Y, Z'=X.\end{cases}$$
Here $\emptyset$ denotes the initial object of $\bfS$ and ${\bf 1}_{\bfS}$ denotes the
unit object with respect to the monoidal structure on $\bfS$. We will denote
$[1]_{ {\bf 1}_{\bfS} }$ simply by $[1]_{\bfS}$. 
We let $[0]_{\bfS}$ denote the $\bfS$-enriched category having only a single object $X$, with $\bHom_{\ast}(X,X)={\bf 1}_{\bfS}$.\index{not}{1subS@$[1]_{A}$}\index{not}{1bfs@$[1]_{\bfS}$}
\index{not}{0bfs@$[0]_{\bfS}$}

We let $C_0$ denote the collection of all morphisms in $\bfS$ of the following types:
\begin{itemize}
\item[$(i)$] The inclusion $\emptyset \hookrightarrow [0]_{\bfS}$.
\item[$(ii)$] The induced maps $[1]_{S} \rightarrow [1]_{S'}$, where
$S \rightarrow S'$ ranges over a set of generators for the weakly saturated class of
cofibrations in $\bfS$.
\end{itemize}

\begin{proposition}\label{enrichcatper}\index{gen}{model category!of simplicial categories}
Let $\bfS$ be a combinatorial monoidal model category. Assume that every object of
$\bfS$ is cofibrant, and that the collection of weak equivalences in $\bfS$ is stable
under filtered colimits. Then
there exists a left proper, combinatorial model structure on $\SCat$, characterized by the following
conditions:
\begin{itemize}
\item[$(C)$] The class of cofibrations in $\SCat$ is the smallest weakly saturated class of morphisms
containing the set of morphisms $C_0$ appearing above.
\item[$(W)$] The weak equivalences in $\SCat$ are defined as in \S \ref{equivdefequiv}.
\end{itemize}
\end{proposition}

\begin{proof}
It suffices to verify the hypotheses of Proposition \ref{goot}. Condition $(1)$ follows from
Remark \ref{swunk}. For condition $(3)$, we must show that any functor $F: \calC \rightarrow \calC'$ having the right lifting property with respect to all morphisms in $C_0$ is a weak equivalence. Since $F$ has the right lifting property with respect to $i: \emptyset \rightarrow [0]_{\bfS}$, it is surjective on objects and therefore essentially surjective. The assumption that $F$ has the right lifting property with respect to the remaining morphisms of $C_0$ guarantees that for every $X, Y \in \calC$, the induced map
$$\bHom_{\calC}(X,Y) \rightarrow \bHom_{\calC'}(F(X),F(Y))$$ is a trivial fibration in
$\bfS$, and therefore a weak equivalence.

It remains to verify condition $(2)$: namely, that the class of weak equivalences is stable under pushout by the elements of $C_0$. We must show that given any pair of functors
$F: \calC \rightarrow \calD$, $G: \calC \rightarrow \calC'$ with $F$ a weak equivalence and $G$ a pushout of some morphism in $C_0$, the induced map $F': \calC' \rightarrow \calD' = \calD \coprod_{\calC} \calC'$ is a weak equivalence. There are two cases to consider.

First, suppose that $G$ is a pushout of the generating cofibration $i: \emptyset \rightarrow \ast$. In other words, the category $\calC'$ is obtained from $\calC$ by adjoining a new object $X$, which admits no morphisms to or from the objects of $\calC$ (and no endomorphisms other than the identity). The category $\calD'$ is obtained from $\calD$ by adjoining $X$ in the same fashion. It is easy to see that if $F$ is a weak equivalence, then $F'$ is also a weak equivalence.

The other basic case to consider is one in which $G$ is a pushout of one of the generating cofibrations $[1]_{S} \rightarrow [1]_{T}$, where $S \rightarrow T$ is a cofibration in $\bfS$. Let $H: [1]_{S} \rightarrow \calC$
denote the ``attaching map'', so that $H$ is determined by a pair of objects $x= H(X)$ and $y=H(Y)$
and a map of $h: S \rightarrow \bHom_{\calC}(x,y)$. 
By definition, $\calC'$ is universal
with respect to the property that it receives a map from $\calC$, and the map $h$ extends to a map
$\widetilde{h}: T \rightarrow \bHom_{\calC'}(x,y)$. To carry out the proof, we will give an explicit construction of a $\bfS$-enriched category $\calC'$ which has this universal property. 

For the remainder of the proof, we will assume that $\bfS$ is the category of simplicial sets.
This is purely for notational convenience; the same arguments can be employed without
change in the general case.

We begin by declaring that the objects of $\calC'$ are the objects of $\calC$. 
The definition of the morphisms in $\calC'$ is a bit more complicated. Let $w$ and $z$ be objects of $\calC$. We define a sequence of simplicial sets $M_{\calC'}^{k}$ as follows:
$$M_{\calC}^0 = \bHom_{\calC}(w,z)$$
$$M_{\calC}^1 = \bHom_{\calC}(y,z) \times T \times \bHom_{\calC}(w,x)$$
$$M_{\calC}^2 = \bHom_{\calC}(y,z) \times T \times \bHom_{\calC}(y,x) \times
T \times \bHom_{\calC}(w,x)$$ 
and so forth. More specifically, for $k \geq 1$, the $m$-simplices of $M_{\calC}^k$ are finite sequences
$$( \sigma_0, \tau_1, \sigma_1, \tau_2, \ldots, \tau_k, \sigma_k)$$ where 
$\sigma_0 \in \bHom_{\calC}(y,z)_m$, $\sigma_k \in \bHom_{\calC}(w,x)_m$, $\sigma_i \in \bHom_{\calC}(y,x)_m$ for $0 < i < k$, and $\tau_i \in T_m$ for $1 \leq i \leq k$.

We define $\bHom_{\calC'}(w,z)$ to be the quotient of the disjoint union
$ \coprod_{k} M_{\calC}^k$ by the equivalence relation which is generated by making the identification
$$( \sigma_0, \tau_1, \ldots, \sigma_k) \simeq (\sigma_0, \tau_1, \ldots, \tau_{j-1}, \sigma_{j-1} \circ h(\tau_j) \circ \sigma_{j}, \tau_{j+1}, \ldots, \sigma_k )$$ whenever the simplex $\tau_j$ belongs to $S_m \subseteq T_m$.

We equip $\calC'$ with an associative composition law, which is given on the level of simplices
by $$( \sigma_0, \tau_1, \ldots, \sigma_k ) \circ (\sigma'_0, \tau'_1, \ldots, \sigma'_l) =
(\sigma_0, \tau_1, \ldots, \tau_k, \sigma_k \circ \sigma'_0, \tau'_1, \ldots, \sigma'_l).$$
It is easy to verify that this composition law is well-defined (that is, compatible with the equivalence relation introduced above), associative, and that the identification $M_{\calC}^0 = \bHom_{\calC}(w,z)$ gives rise to an inclusion of categories $\calC \subseteq \calC'$. Moreover, the map 
$h: S \rightarrow \bHom_{\calC}(x,y)$ extends to $\widetilde{h}: T \rightarrow \bHom_{\calC'}(x,y)$, given by the composition
$$ T \simeq \{ \id_y \} \times T \times \{ \id_x\} \subseteq \bHom_{\calC}(y,y) \times T \times \bHom_{\calC}(x,x) = M_{\calC}^1 \rightarrow \bHom_{\calC'}(x,y).$$
Moreover, it is not difficult to see that $\calC'$ has the desired universal property.

We observe that, by construction, the simplicial sets $\bHom_{\calC'}(w,z)$ come equipped with a natural filtration. Namely, define $\bHom_{\calC'}(w,z)^k$ to be the image of
$$ \coprod_{0 \leq i \leq k} M_{\calC}^i $$ in $\bHom_{\calC'}(w,z)$. Then we have
$$ \bHom_{\calC}(w,z) = \bHom_{\calC'}(w,z)^0 \subseteq \bHom_{\calC'}(w,z)^1 \subseteq \ldots$$
and $\bigcup_{k} \bHom_{\calC'}(w,z)^k = \bHom_{\calC'}(w,z)$. Moreover, the inclusion
$ \bHom_{\calC'}(w,z)^k \subseteq \bHom_{\calC'}(w,z)^{k+1}$ is a pushout of the inclusion
$N_{\calC}^{k+1} \subseteq M_{\calC}^{k+1}$, where $N^{k+1}$ is the simplicial subset of $M_{\calC}^{k+1}$ whose
$m$-simplices consist of those $(2m+1)$-tuples $(\sigma_0, \tau_1, \ldots, \sigma_m)$ such that
$\tau_i \in S_m$ for at least one value of $i$.

Let us now return to the problem at hand: namely, we wish to prove that $F': \calC' \rightarrow \calD'$ is an equivalence. We note that the construction outlined above may also be employed to produce a model for $\calD'$, and an analogous filtration on its morphism spaces.

Since $G': \calD \rightarrow \calD'$ and $F: \calC \rightarrow \calD$ are essentially surjective, we deduce that $F'$ is essentially surjective. Hence it will suffice to show that, for any objects $w,z \in \calC'$, the induced map $$\phi: \bHom_{\calC'}(w,z) \rightarrow \bHom_{\calD'}(w,z)$$ is a weak homotopy equivalence.
For this, it will suffice to show that for each $i \geq 0$, the induced map $\phi_i: \bHom_{\calC'}(w,z)^i \rightarrow \bHom_{\calD'}(w,z)^i$ is a weak homotopy equivalence; then $\phi$, being a filtered colimit of weak homotopy equivalences $\phi_i$, will itself be a weak homotopy equivalence.

The proof now proceeds by induction on $i$. When $i=0$, $\phi_i$ is a weak homotopy equivalence by assumption (since $F$ is an equivalence of simplicial categories). For the inductive step, we note that
$\phi_{i+1}$ is obtained as a pushout
$$ \bHom_{\calC'}(w,z)^i \coprod_{ N^{i+1}_{\calC}} M^{i+1}_{\calC} \rightarrow \bHom_{\calD'}(w,z)^i \coprod_{ N^{i+1}_{\calD} } M^{i+1}_{\calD}.$$
Since $\bfS$ is left-proper, both of these pushouts are homotopy pushouts. Consequently, to show that $\phi_{i+1}$ is a weak equivalence, it suffices to show that $\phi_i$ is a weak equivalence and the each of the maps
$$ N^{i+1}_{\calC} \rightarrow N^{i+1}_{\calD}$$
$$ M^{i+1}_{\calC} \rightarrow M^{i+1}_{\calD}$$
are weak equivalences. These statements follow easily from the compatibility of the monoidal
structure of $\bfS$ with the model structure, and the assumption that every object of $\bfS$ is cofibrant.
\end{proof}

\begin{remark}
It follows from the proof of Proposition \ref{enrichcatper} that if
$f: \calC \rightarrow \calC'$ is a cofibration of $\bfS$-enriched categories, then
the induced map $\bHom_{\calC}(X,Y) \rightarrow \bHom_{\calC'}(fX,fY)$ is a cofibration
for every pair of objects $X,Y \in \calC$.
\end{remark}

\begin{remark}\label{cuttup}
The model structure of Proposition \ref{enrichcatper} enjoys the following functoriality:
suppose that $f: \bfS \rightarrow \bfS'$ is a monoidal left Quillen functor between
model categories satisfying the hypotheses of Proposition \ref{enrichcatper}, with right adjoint $g: \bfS' \rightarrow \bfS$. Then $f$ and $g$ induce a Quillen adjunction
$$ \Adjoint{ F}{\SCat}{\Cat_{\bfS'},}{G}$$
where $F$ and $G$ are as in Remark \ref{laxcon}. Moreover, if $(f,g)$ is a Quillen equivalence,
then $(F,G)$ is likewise a Quillen equivalence.
\end{remark}

In order for Proposition \ref{enrichcatper} to be useful in practice, we need to understand the fibrations in $\SCat$. For this, we first introduce a few definitions.

\begin{definition}\label{qfibb}
Let $F: \calC \rightarrow \calD$ be a functor between ordinary categories. We will say that
$F$ is a {\it quasi-fibration} if, for every object $C \in \calC$ and every isomorphism
$f: F(X) \rightarrow Y$ in $\calD$, there exists an isomorphism $\overline{f}: X \rightarrow \overline{Y}$ in $\calC$ such that $F( \overline{f} ) = f$.
\end{definition}

\begin{remark}
The relevance of Definition \ref{qfibb} is as follows: the category $\Cat$ admits a model structure
in which the weak equivalences are the equivalences of categories, and the fibrations are the quasi-fibrations. This is a special case of Theorem \ref{staycode}, which we will prove below
(namely, the special case where we take $\bfS = \Set$, endowed with the trivial model
structure of Example \ref{trivmodel}). 
\end{remark}

\begin{definition}\label{twurp}
Let $\bfS$ be a monoidal model category, and let $\calC$ be a
$\bfS$-enriched category. We will say that a morphism $f$
in $\calC$ is an {\it equivalence} if the homotopy class $[f]$
of $f$ is an isomorphism in $\h{\calC}$.\index{gen}{equivalence!in a $\bfS$-enriched category}

We will say that $\calC$ is {\it locally fibrant} if, for every pair of objects
$X,Y \in \calC$, the mapping space $\bHom_{\calC}(X,Y)$ is a fibrant object
of $\bfS$.\index{gen}{fibrant!locally}\index{gen}{locally fibrant}

We will say that a $\bfS$-enriched functor $F: \calC \rightarrow \calC'$ is a
{\it local fibration}\index{gen}{fibration!local} if the following conditions are satisfied:
\begin{itemize}
\item[$(i)$] For every pair of objects
$X,Y \in \calC$, the induced map $\bHom_{\calC}(X,Y) \rightarrow
\bHom_{\calC'}(FX,FY)$ is a fibration in $\bfS$. 
\item[$(ii)$] The induced map $\h{\calC} \rightarrow \h{\calC'}$ is a quasi-fibration
of categories.
\end{itemize}
\end{definition}

\begin{remark}\label{dinty}
Let $F: \calC \rightarrow \calC'$ be a functor between $\bfS$-enriched categories
which satisfies condition $(i)$ of Definition \ref{twurp}. Let $X \in \calC$
and $Y \in \calC'$ be objects. If $\calC'$
is locally fibrant, then every morphism $[f]: F(X) \rightarrow Y$ in
$\h{\calC'}$ can be represented by an equivalence $f: F(X) \rightarrow Y$ in $\calC'$.
Let $\overline{Y}$ be an object of $\calC$ such that $F( \overline{Y} = Y$. Since
${\bf 1}_{\bfS}$ is a cofibrant object of $\bfS$ and the map
$\bHom_{\calC}( X, \overline{Y} ) \rightarrow \bHom_{\calC}( F(X), Y)$
is a fibration, Proposition \ref{princex} implies that $[f]$ can be lifted to
an isomorphism $[\overline{f}]: X \rightarrow \overline{Y}$ in $\h{\calC}$ if
and only if $f$ can be lifted to an equivalence $\overline{f}: X \rightarrow \overline{Y}$
in $\calC$. Consequently, when $\h{\calC}$ is locally fibrant,
condition $(ii)$ is equivalent to the following analogous assertion:
\begin{itemize}
\item[$(ii')$] For every equivalence $f: F(X) \rightarrow Y$ in $\calC'$, there
exists an equivalence $\overline{f}: X \rightarrow \overline{Y}$ in $\calC$
such that $F( \overline{f} ) = f$.
\end{itemize}
\end{remark}

\begin{notation}\index{not}{[1]simS@$[1]^{\sim}_{\bfS}$}
We $[1]^{\sim}_{\bfS}$ denote the $\bfS$-enriched category
containing a pair of objects $X,Y$, with
$$ \bHom_{ [1]^{\sim}_{\bfS} }(Z, Z') = {\bf 1}_{\bfS}$$
for all $Z, Z' \in \{X, Y \}$. 
\end{notation}

\begin{definition}[Invertibility Hypothesis]\index{gen}{invertibility hypothesis}\label{inhyp}
Let $\bfS$ be a monoidal model category satisfying the hypotheses of Proposition
\ref{enrichcatper}. We will say that
$\bfS$ {\em satisfies the invertibility hypothesis} if the following condition is satisfied:
\begin{itemize}
\item[$(\ast)$] Let $i: [1]_{\bfS} \rightarrow \calC$ be a cofibration of $\bfS$-enriched
categories, classifying a morphism $f$ in $\calC$ which is invertible in the homotopy category $\h{\calC}$, and form a pushout diagram
$$ \xymatrix{ [1]_{\bfS} \ar[r]^{i} \ar[d] & \calC \ar[d]^{j} \\
[1]_{\bfS}^{\sim} \ar[r] & \calC\langle f^{-1} \rangle }$$ Then $j$ is an equivalence of $\bfS$-enriched categories.
\end{itemize}
\end{definition}

In other words, the invertibility hypothesis is the assertion that inverting a morphism
$f$ in a $\bfS$-enriched category $\calC$ does not change the homotopy type of $\calC$ when $f$ is already invertible up to homotopy. 

\begin{remark}\label{attaboy}
Let $\bfS$, $f$, and $\calC$ be as in Definition \ref{inhyp}, and choose a trivial cofibration
$F: \calC \rightarrow \calC'$, where $\calC'$ is a fibrant $\bfS$-enriched category.
Since $\SCat$ is left proper, the induced map
$\calC \langle f^{-1} \rangle \rightarrow \calC' \langle F(f)^{-1} \rangle$
is an equivalence of $\bfS$-enriched categories. Consequently,
assertion $(\ast)$ holds for $(\calC, f)$ if and only if it holds for
$( \calC', F(f) )$. In other words, to test whether $\bfS$ satisfies the invertibility
hypothesis, we only need to check $(\ast)$ in the case where $\calC$ is fibrant.
\end{remark}

\begin{remark}\label{uppa}
In Definition \ref{inhyp}, the condition that $i$ be a cofibration guarantees that
the construction $\calC \mapsto \calC \langle f^{-1} \rangle$ is homotopy invariant.
Alternatively, we can guarantee this by choosing a cofibrant replacement for the
map $j: [1]_{\bfS} \rightarrow [1]^{\sim}_{\bfS}$. Namely, choose a factorization for
$j$ as a composition
$$ [1]_{\bfS} \stackrel{j'}{\rightarrow} \calE \stackrel{j''}{\rightarrow} [1]^{\sim}_{\bfS},$$
where $j''$ is a weak equivalence and $j'$ is a cofibration. For every
$\bfS$-enriched category containing a morphism $f$, define
$\calC [f^{-1} ] = \calC \coprod_{ [1]_{\bfS} } \calE$. 
Then we have a canonical map $\calC [f^{-1}] \rightarrow \calC \langle f^{-1} \rangle$,
which is an equivalence whenever the map $[1]_{\bfS} \rightarrow \calC$ classifying $f$ is a cofibration. Moreover, the construction $\calC \mapsto \calC[f^{-1}]$ preserves weak equivalences in $\calC$. Consequently, we may reformulate the invertibility hypothesis as follows:
\begin{itemize}
\item[$(\ast')$] For every $\bfS$-enriched category $\calC$ containing an equivalence
$f$, the map $\calC \rightarrow \calC[f^{-1}]$ is a weak equivalence of $\bfS$-enriched categories.
\end{itemize}
\end{remark}

\begin{remark}\label{cuddle}
Let $\calC$ be a fibrant $\bfS$-enriched category containing an equivalence
$f: X \rightarrow Y$, and let $\calC[ f^{-1} ]$ be defined as in Remark \ref{uppa}. 
The canonical map $\calC \rightarrow \calC[f^{-1}]$ is a trivial cofibration, and
therefore admits a section. This section determines a map of $\bfS$-enriched
categories $h: \calE \rightarrow \calC$. We observe that $\calE$ is a mapping
cylinder for the object $[0]_{\SCat} \in \SCat$, so we can view
$h$ as a homotopy between the maps $[0]_{\SCat} \rightarrow \calC$
classifying the objects $X$ and $Y$.

More generally, the same argument shows that if $F: \calC \rightarrow \calD$
is a fibration of $\bfS$-enriched categories and $f: X \rightarrow Y$ is an equivalence
in $\calC$ such that $F(f) = \id_{D}$ for some object $D \in \calD$, then
the functors $[0]_{\bfS} \rightarrow \calC$ classifying the objects $X$ and $Y$
are homotopic in the model category $( \SCat)_{/ \calD}$. 
\end{remark}

\begin{definition}\index{gen}{excellent!model category}\index{gen}{model category!excellent}\label{modelexcellent}
We will say that a model category $\bfS$ is {\it excellent} if it is equipped with a symmetric
monoidal structure and satisfies the following conditions:
\begin{itemize}
\item[$(A1)$] The model category $\bfS$ is combinatorial.
\item[$(A2)$] Every monomorphism in $\bfS$ is a cofibration, and the collection of cofibrations
is stable under products.
\item[$(A3)$] The collection of weak equivalences in $\bfS$ is stable under filtered colimits.
\item[$(A4)$] The monoidal structure on $\bfS$ is compatible with the model structure.
In other words, the tensor product functor $\otimes: \bfS \times \bfS \rightarrow \bfS$ is a left Quillen bifunctor.
\item[$(A5)$] The monoidal model category $\bfS$ satisfies the invertibility hypothesis.
\end{itemize}
\end{definition}

\begin{remark}\label{dummu}
Axiom $(A2)$ of Definition \ref{modelexcellent} implies that every object of $\bfS$ is cofibrant.
In particular, $\bfS$ is left proper.
\end{remark}

\begin{example}[Dwyer, Kan]\label{cuppata}
The category of simplicial sets is an excellent model category, when endowed with the Kan model structure and the Cartesian product. The only nontrivial point is to show that $\sSet$ satisfies the invertibility hypothesis. This is one of the main theorems of \cite{dwyerkan}.
\end{example}

\begin{example}
Let $\bfS$ be a presentable category equipped with a closed symmetric monoidal structure.
Then $\bfS$ is an excellent model category with respect to the trivial model structure
of Example \ref{trivmodel}.
\end{example}

The following lemma guarantees a good supply of examples of excellent model categories:

\begin{lemma}\label{cuppat}
Suppose given a monoidal left Quillen functor $T: \bfS \rightarrow \bfS'$
between model categories $\bfS$ and $\bfS'$ satisfying axioms
$(A1)$ through $(A4)$ of Definition \ref{modelexcellent}.
If $\bfS$ satisfies axiom $(A5)$, then so does $\bfS'$.
\end{lemma}

\begin{proof}
As indicated in Remark \ref{cuttup}, the functor $T$ determines a Quillen adjunction
$$ \Adjoint{F}{\Cat_{\bfS}}{\Cat_{\bfS'}}{G}.$$
Let $\calC$ be a $\bfS'$-enriched category and $i: [1]_{\bfS'} \rightarrow \calC$
a cofibration classifying an equivalence $f$ in $\calC$. We wish to prove that the
map $\calC \rightarrow \calC \langle f^{-1} \rangle$ is an equivalence of
$\bfS'$-enriched categories. In view of Remark \ref{attaboy}, we may assume
that $\calC$ is fibrant.

Choose a factorization of the map $[1]_{\bfS} \rightarrow [1]_{\bfS}^{\sim}$
as a composition
$$ [1]_{\bfS} \stackrel{j}{\rightarrow} \calE \stackrel{j'}{\rightarrow} [1]_{\bfS}^{\sim}$$
as in Remark \ref{uppa}, so that we have an analogous factorization
$$ [1]_{\bfS'} \rightarrow F(\calE) \rightarrow [1]_{\bfS'}^{\sim}$$
in $\Cat_{\bfS'}$. Using the latter factorization, we can define
$\calC[ f^{-1}]$ as in Remark \ref{uppa}; we wish to show that the map
$h: \calC \rightarrow \calC[f^{-1}]$ is a trivial cofibration. 

Let $f_0$ be the morphism in $G(\calC)$ classified by $f$, and
let $G(\calC)[f_0^{-1}] \in \SCat$ be defined as in Remark \ref{uppa}. 
Using the fact that $\calC$ is locally fibrant (see Theorem \ref{staycode} below),
we conclude that $f_0$ is an equivalence in $G(\calC)$. Since $\bfS$
satisfies the invertibility hypothesis, the map
$h_0: G(\calC) \rightarrow G(\calC)[f_0^{-1}]$ is a trivial cofibration.
We now conclude by observing that $h$ is a pushout of $F(h_0)$.
\end{proof}

\begin{remark}\label{conrem}
Using a similar argument, we can prove a converse to Lemma \ref{cuppat} in the case where
$T$ is a Quillen equivalence.
\end{remark}

\begin{example}\label{supermat}
Let $\bfS$ be the category $\mSet$ of marked simplicial sets, endowed with
the Cartesian model structure defined in \S \ref{twuf}. Then the functor
$X \mapsto X^{\sharp}$ is a monoidal left Quillen functor 
$\sSet \rightarrow \bfS$. Combining Example \ref{cuppata} with Lemma \ref{cuppat}, we
conclude that $\bfS$ satisfies the invertibility hypothesis, so that $\bfS$ is an excellent model category (with respect to the Cartesian product).
\end{example}

\begin{example}
Let $\bfS$ denote the category of simplicial sets, endowed with the Joyal model structure.
The functor $X \mapsto X^{\flat}$ determines a monoidal left Quillen equivalence
$\bfS \rightarrow \mSet$. Using Remark \ref{conrem}, we deduce that $\bfS$ satisfies the invertibility hypothesis, so that $\bfS$ is an excellent model category (with respect to the Cartesian product).
\end{example}

We are now ready to state our main result:

\begin{theorem}\label{staycode}
Let $\bfS$ be an excellent model category.
Then:
\begin{itemize}
\item[$(1)$] An $\bfS$-enriched category $\calC$ is a fibrant object of $\SCat$ if and only if it is {\em locally fibrant}: that is, if and only if the mapping object $\bHom_{\calC}(X,Y) \in \bfS$ is fibrant for
every pair of objects $X,Y \in \calC$. 
\item[$(2)$] Let $F: \calC \rightarrow \calD$ be a $\bfS$-enriched functor, where $\calD$ is a
fibrant object of $\SCat$. Then $F$ is a fibration in $\SCat$ if and only if $F$ is a local fibration.
\end{itemize}
\end{theorem}

\begin{remark}
In the case where $\bfS$ is the category of simplicial sets (with its usual model structure),
Theorem \ref{staycode} is due to Bergner; see \cite{bergner}. Moreover, Bergner
proves a stronger result in this case: assertion $(2)$ holds without the assumption that
$\calD$ is fibrant.
\end{remark}

Before giving the proof of proof of Theorem \ref{staycode}, we need to establish some preliminaries.
Fix an excellent model category $\bfS$. We observe
that $\SCat$ is naturally {\em cotensored} over $\bfS$. That is, for
every $\bfS$-enriched category $\calC$ and every object $K \in \bfS$, we can define a new $\bfS$-enriched category $\calC^{K}$ as follows:

\begin{itemize}
\item[$(i)$] The objects of $\calC^{K}$ are the objects of $\calC$.
\item[$(ii)$] Given a pair of objects $X,Y \in \calC$, we have
$\bHom_{\calC^{K}}(X,Y) = \bHom_{\calC}(X,Y)^{K} \in \bfS$. 
\end{itemize} 

This construction does not endow $\SCat$ with the structure of a $\bfS$-enriched category, because
the construction $\calD \mapsto \calD^{K}$ is not compatible with colimits in $K$. However,
we can remedy the situation as follows. Let $\calC$ and $\calD$ be $\bfS$-enriched categories,
and let $\phi$ be a function from the set of objects of $\calC$ to the set of objects of $\calD$. Then
there exists an object $\bHom_{\SCat}^{\phi}( \calC, \calD) \in \bfS$, which is characterized by the following universal property: for every $K \in \bfS$, there is a natural bijection
$$ \Hom_{\bfS}(K, \bHom_{\SCat}^{u}(\calC, \calD) ) \simeq \Hom^{\phi}_{\SCat}(\calC, \calD^{K} ),$$
where $\Hom^{\phi}_{\SCat}(\calC, \calD^{K} )$ denotes the set of all functors from
$\calC$ to $\calD^{K}$ which is given on objects by the function $\phi$.

\begin{lemma}\label{stuttcat}
Let $\bfS$ be an excellent model category.
Fix a diagram of $\bfS$-enriched categories
$$ \xymatrix{ \calC \ar[d]^{F} \ar[r]^{u} & \calC' \ar[d]^{F'} \\
\calD \ar[r]^{u'} & \calD'. }$$ 
Assume that:
\begin{itemize}
\item[$(a)$] For every pair of objects $X, Y \in \calC$, the diagram
$$ \xymatrix{ \bHom_{\calC}(X,Y) \ar[r] \ar[d] & \bHom_{\calD}( FX, FY) \ar[d] \\
\bHom_{\calC'}(uX, uY) \ar[r] & \bHom_{\calD'}( u'FX, u'FY) }$$
is homotopy pullback square between fibrant objects of $\bfS$, and the vertical arrows are fibrations.
\end{itemize}

Let $G: \calA \rightarrow \calB$ be a functor between $\bfS$-enriched categories
which is a transfinite composition of pushouts of generating cofibrations of the form
$[1]_{S} \rightarrow [1]_{S'}$, where $S \rightarrow S'$ is a cofibration in
$\bfS$, and let $\phi$ be a function from the set of objects of $\calB$ (which
is isomorphic to the set of objects of $\calA$) to $\calC$. 
Then the diagram
$$ \xymatrix{  \bHom^{\phi}_{ \SCat}( \calB, \calC) \ar[r] \ar[d] &
\bHom^{F\phi}_{\SCat}( \calB, \calD) \times_{ \bHom^{F\phi}_{\SCat}(\calA, \calD)}
\bHom^{\phi}_{\SCat}(\calA, \calC) \ar[d] \\
 \bHom^{u\phi}_{\SCat}(\calB, \calC') \ar[r] & 
\bHom^{u'F\phi}_{\SCat}(\calB, \calD') \times_{ \bHom^{u'F\phi}_{\SCat}(\calA, \calD') }
\bHom^{u\phi}_{\SCat}( \calA, \calC') }$$
is a homotopy pullback square between fibrant objects of $\bfS$, and the vertical arrows
are fibrations.
\end{lemma}

\begin{proof}
It is easy to see that the collection of morphisms $G: \calA \rightarrow \calB$ which satisfy the conclusion of the lemma is weakly saturated. It will therefore suffice to show that $G$ contains every
morphism of the form $[1]_{S} \rightarrow [1]_{S'}$, where $S \rightarrow S'$
is a cofibration in $\bfS$. In this case, $\phi$ determines a pair of objects $X,Y \in \calC$, and we can rewrite the diagram of interest as
$$ \xymatrix{ \bHom_{\calC}(X,Y)^{S'} \ar[r] \ar[d] & 
\bHom_{\calC}(X,Y)^{S} \times_{ \bHom_{\calD}(FX,FY)^{S} }
\bHom_{\calD}(FX,FY)^{S'} \ar[d] \\
\bHom_{\calC'}( uX, uY)^{S'} \ar[r] & 
\bHom_{\calC'}(uX,uY)^{S} \times_{ \bHom_{\calD'}(u'FX, u'FY)^{S} }
\bHom_{\calD'}( u' FX, u'FY)^{S'}. }$$
The desired result now follows from $(a)$, since the map $S \rightarrow S'$ is a cofibration
between cofibrant objects of $\bfS$.
\end{proof}

\begin{proof}[Proof of Theorem \ref{staycode}]
Assertion $(1)$ is just a special case of $(2)$, where we take $\calD$ to be the final object
of $\SCat$. It will therefore suffice to prove $(2)$. 

We first prove the ``only if'' direction.
If $F$ is a fibration, then $F$ has the right lifting property with respect to every
trivial cofibration of the form $[1]_{S} \rightarrow [1]_{S'}$, where $S \rightarrow S'$ is
a trivial cofibration in $\bfS$. It follows that for every pair of objects $X, Y \in \calC$, 
the induced map $\bHom_{\calC}(X,Y) \rightarrow \bHom_{\calD}(FX, FY)$ is a fibration in $\bfS$.
In particular, $\calC$ is locally fibrant.

To complete the proof that $F$ is a local fibration, we will show that $F$ satisfies condition $(ii')$ of Remark \ref{dinty}. Suppose $X \in \calC$, and that $f: FX \rightarrow Y$ is an equivalence in
$\calD$. We wish to show that we can lift $f$ to an equivalence
$\overline{f}: X \rightarrow \overline{Y}$. Let $\calE$ and $\calD[f^{-1}]$ be defined as
in Remark \ref{uppa}. 
Since $\bfS$ satisfies the invertibility hypothesis, the map $h: \calD \rightarrow \calD[f^{-1}]$
is a trivial cofibration. Because we have assumed $\calD$ to be fibrant, the map $h$
admits a section. This section determines a map $s: \calE \rightarrow \calD$.
We now consider the lifting problem
$$ \xymatrix{ [0]_{\bfS} \ar[r]^{X} \ar[d] & \calC \ar[d]^{F} \\
\calE \ar[r]^{s} \ar@{-->}[ur] & \calD. }$$
Since $F$ is a fibration and the left vertical map is a trivial cofibration, there
exists a solution as indicated. This solution determines a morphism
$\overline{f}: X \rightarrow \overline{Y}$ in $\calC$ lifting $f$.
Moreover, $\overline{f}$ is the image of a morphism in $\calE$.
Since every morphism in $\calE$ is an equivalence, we deduce that
$\overline{f}$ is an equivalence in $\calC$.

Let us now suppose that $F$ is a local fibration. We wish to show that
$F$ is a fibration. Choose a factorization of $F$ as a composition
$$ \calC \stackrel{u}{\rightarrow} \calC' \stackrel{F'}{\rightarrow} \calD,$$
where $u$ is a weak equivalence and $F'$ is a fibration. We will prove the following:

\begin{itemize}
\item[$(\ast)$] Suppose given a commutative diagram of $\bfS$-enriched categories
$$ \xymatrix{ \calA \ar[r]^{v} \ar[d]^{G} & \calC \ar[d]^{F} \\
\calB \ar[r]^{v'} & \calD, }$$
where $G$ is a cofibration. If there exists a functor $\alpha: \calB \rightarrow \calC'$ such
that $\alpha G = u v$ and $F' \alpha = v'$, then there exists a functor $\beta: \calB \rightarrow \calC$
such that $\beta G = v$ and $F \beta = v'$. 
\end{itemize}

Since the map $F'$ has the right lifting property with respect to all trivial cofibrations, 
assertion $(\ast)$ implies that $F$ also has the right lifting property with respect to all trivial cofibrations, so that $F$ is a fibration as desired.

We now prove $(\ast)$. Using the small object argument, we deduce that 
the functor $G$ is a retract of some functor $G': \calA \rightarrow \calB'$, where
$G'$ is a transfinite composition of morphisms obtained as pushouts of 
generating cofibrations. It will therefore suffice to prove $(\ast)$ after replacing $G$ by $G'$.

Reordering the transfinite composition if necessary, we may assume that $G'$ factors as a
composition
$$ \calA \stackrel{G'_0}{\rightarrow} \calB'_0 \stackrel{G'_1}{\rightarrow} \calB',$$
where $\calB'_0$ is obtained from $\calA$ by adjoining a collection of new objects,
$\{B_i \}_{i \in I}$, and $\calB'$ is obtained from $\calB'_0$ by a transfinite sequence of pushouts
by generating cofibrations of the form $\calE_{S} \rightarrow \calE_{S'}$, where
$S \rightarrow S'$ is a cofibration in $\bfS$. Let $C'_i = \alpha(B_i)$ for each $i \in I$.
Since $u$ is an equivalence of $\bfS$-enriched categories, there exists a collection
of objects $\{ C_i \}_{i \in I}$ and equivalences$f_i: u C_i \rightarrow C'_i$. 
Let $g_i$ be the image of $f_i$ in $\calD$. Since $F$ is a local fibration, 
we can lift each $g_i$ to an equivalence
$f'_i: C_i \rightarrow C''_i$ in $\calC$. 
Since the maps $\bHom_{\calC'}( uC''_i, C'_i) \rightarrow \bHom_{\calD}( FC''_i, F'C'_i)$
are fibrations, we can choose morphisms $f''_i: uC''_i \rightarrow C'_i$
in $\calC'$ such that $F'( f'_i)$ is the identity for each $i$, and the diagrams
$$ \xymatrix{ & uC''_i \ar[dr]^{f''_i} & \\
uC_i \ar[rr]^{f_i} \ar[ur]^{f'_i} & & C'_i }$$
commute up to homotopy. Replacing $C_i$ by $C''_i$, we may assume
that each of the maps $f_i$ projects to the identity in $\calD$.

Let $\alpha_0 = \alpha | \calB'_0$, and let $\alpha'_0: \calB'_0 \rightarrow \calC'$ be defined by the formula $$ \alpha'_0(A) = \begin{cases} \alpha_0(A) & \text{if } A \in \calA \\
uC_i & \text{if } A = B_i, i \in I. \end{cases}$$
Remark \ref{cuddle} implies that the maps $\alpha_0$ and
$\alpha'_0$ are homotopic in the model category
$(\SCat)_{\calA/ \, / \calD}$. Applying Proposition \ref{princex}, we deduce the existence of a map
$\alpha': \calB' \rightarrow \calC$ which extends $\alpha_0$ and satisfies
$\alpha' G = u v$ and $F' \alpha' = v'$. We may therefore replace 
$\alpha$ by $\alpha'$, $v$ by $\alpha'_0$, and $\calA$ by $\calB'_0$, and thereby
reduce to the case where the functor $G: \calA \rightarrow \calB$ is a transfinite
composition of generating cofibrations of the form $\calE_{S} \rightarrow \calE_{S'}$, where
$S \rightarrow S'$ is a cofibration in $\bfS$.

Let $\phi$ be the map from the objects of $\calB$ to the objects of $\calC$ determined by $\alpha$.
Applying Lemma \ref{stuttcat}, we obtain a homotopy pullback diagram
$$ \xymatrix{  \bHom^{\phi}_{ \SCat}( \calB, \calC) \ar[r] \ar[d] & \bHom^{u\phi}{\SCat}(\calB, \calC') \ar[d] \\
\bHom^{F\phi}_{\SCat}( \calB, \calD) \times_{ \bHom^{F\phi}_{\SCat}(\calA, \calD)}
\bHom^{\phi}_{\SCat}(\calA, \calC) \ar[r] &
\bHom^{F\phi}_{\SCat}(\calB, \calD) \times_{ \bHom^{F\phi}_{\SCat}(\calA, \calD) }
\bHom^{u\phi}_{\SCat}( \calA, \calC'). }$$
in which the vertical arrows are fibrations. We therefore have a weak
equivalence $$\bHom^{\phi}_{\SCat}( \calB, \calC)
\rightarrow M = \bHom^{u \phi}_{\SCat}(\calB, \calC') \times_{ \bHom^{F \phi}_{\SCat}( \calB, \calD) }
\bHom^{\phi}_{\SCat}(\calA, \calC)$$
of fibrations over $N= \bHom^{F\phi}_{\SCat}(\calB, \calD) \times_{ \bHom^{F\phi}_{\SCat}(\calA, \calD) } \bHom^{u\phi}_{\SCat}( \calA, \calC')$. Moreover, the pair
$(\alpha,v)$ determines a map ${\bf 1}_{\bfS} \rightarrow M$ lifting the map
$(v', uv'): {\bf 1}_{\bfS} \rightarrow N$. Applying Proposition \ref{princex}, we deduce
that $(v,uv'): {\bf 1}_{\bfS} \rightarrow N$ can be lifted to a map
${\bf 1}_{\bfS} \rightarrow \bHom^{\phi}_{\SCat}(\calB, \calC)$, which is equivalent
to the existence of the desired map $\beta$.
\end{proof}

We conclude this section with a few easy results concerning homotopy limits in the
model category $\SCat$.

\begin{proposition}\label{scam}
Let $\bfS$ be an excellent model category, $\calJ$ a small category, and
$\{ \calC_{J} \}_{J \in \calJ}$ a diagram of $\bfS$-enriched categories.
Suppose given a compatible family of functors $\{ f_J: \calC \rightarrow \calC_{J} \}_{J \in \calJ}$ 
which exhibits
$\calC$ as a homotopy limit of the diagram $\{ \calC_{J} \}_{J \in \calJ}$ in $\SCat$. Then
for every pair of objects $X, Y \in \calC$, the maps
$\{ \bHom_{\calC}( X,Y) \rightarrow \bHom_{\calC_{J} }( f_{J} X, f_{J} Y) \}_{J \in \calJ}$
exhibit $\bHom_{\calC}(X,Y)$ as a homotopy limit of the diagram
$\{ \bHom_{ \calC_{J} }( f_J X, f_J Y) \}_{J \in \calJ}$ in $\bfS$.
\end{proposition}

\begin{proof}
Without loss of generality, we may assume that the diagram $\{ \calC_{J} \}_{J \in \calJ}$ is injectively fibrant, and that the maps $f_{J}$ exhibit $\calC$ as a limit of $\{ \calC_{J} \}_{J \in \calJ}$. 
It follows that $\bHom_{\calC}(X,Y)$ is a limit of the diagram $\{ \bHom_{\calC_J}( f_J X, f_J Y) \}_{J \in \calJ}$.
It will therefore suffice to show that the diagram $\{ \bHom_{ \calC_{J}}( f_J X, f_J Y) \}_{J \in \calJ}$
is injectively fibrant. For this, it will suffice to show that $\{ \bHom_{ \calC_{J}}( f_J X, f_J Y) \}_{J \in \calJ}$
has the right lifting property with respect to every weak trivial cofibration 
$\alpha: F \rightarrow F'$ of diagrams $F,F': \calJ \rightarrow \bfS$. Let
$G: \calJ \rightarrow \SCat$ be defined by the formula
$G(J) = [1]_{F(J)}$, and let $G': \calJ \rightarrow \SCat$ be defined likewise. The desired
result now follows from the observe that $\alpha$ induces a weak trivial cofibration
$G \rightarrow G'$ in $\Fun( \calJ, \SCat)$. 
\end{proof}

\begin{corollary}\label{wspin}
Let $\bfS$ be an excellent model category, $\calJ$ a small category, and
$\{ \calC_{J} \}_{J \in \calJ}$ a diagram of $\bfS$-enriched categories.
Suppose given $\bfS$-enriched functors
$$ \calD \stackrel{\beta}{\rightarrow} \calC \stackrel{\alpha}{\rightarrow} \lim \{ \calC_{J} \}_{J \in \calJ} $$
such that $\alpha \circ \beta$ exhibits $\calD$ as a homotopy limit of the diagram
$\{ \calC_{J} \}_{J \in \calJ}$. Then the following conditions are equivalent:
\begin{itemize}
\item[$(1)$] The functor $\alpha$ exhibits $\calC$ as a homotopy limit of the diagram
$\{ \calC_{J} \}_{J \in \calJ}$.
\item[$(2)$] For every pair of objects $X,Y \in \calC$, the functor $\alpha$ exhibits
$\bHom_{\calC}(X,Y)$ as a homotopy limit of the diagram $\{ \bHom_{\calC_{J}}( \alpha_{J} X,
\alpha_{J} Y ) \}_{J \in \calJ}$.
\end{itemize}
\end{corollary}

\begin{proof}
The implication $(1) \Rightarrow (2)$ follows from Proposition \ref{scam}. To prove the
converse, we may assume that the diagram $\{ \calC_{J} \}_{J \in \calJ}$ is injectively fibrant.
In view of $(2)$, Proposition \ref{scam} implies that $\alpha$ induces a fully faithful functor
between $\h{\bfS}$-enriched homotopy categories. It will therefore suffice to show that
$\alpha$ is essentially surjective on homotopy categories, which follows from our assumption that 
$\alpha \circ \beta$ is a weak equivalence. 
\end{proof}

\subsection{Model Structures on Diagram Categories}\label{quasilimit3}

In this section, we consider enriched analogues of the constructions presented in
\S \ref{qlim7}. Namely, suppose that $\bfS$ is an excellent model category,
$\bfA$ a combinatorial $\bfS$-enriched model category, and $\calC$ a small
$\bfS$-enriched category. Let $\bfA^{\calC}$ denote the category of
$\bfS$-enriched functors from $\calC$ to $\bfA$. In this section, we will study the associated
projective and injective model structures on $\bfA^{\calC}$. The ideas described here will be used in \S \ref{pathspace} to construct certain mapping objects in $\SCat$.

We begin with the analogue of Definition \ref{injproj}.

\begin{definition}\label{projinj}\index{gen}{cofibration!strong}\index{gen}{cofibration!weak}\index{gen}{fibration!strong}\index{gen}{fibration!weak}\index{gen}{strong!cofibration}\index{gen}{strong!fibration}\index{gen}{weak!cofibration}\index{gen}{weak!fibration}
Let $\calC$ be a small $\bfS$-category, and $\bfA$ a combinatorial $\bfS$-enriched model category.
A natural transformation $\alpha: F \rightarrow G$ in $\bfA^{\calC}$ is a:

\begin{itemize}
\item {\it injective cofibration} if the induced map $F(C) \rightarrow
G(C)$ is a cofibration in $\bfA$, for each $C \in \calC$.

\item {\it projective fibration} if the induced map $F(C) \rightarrow
G(C)$ is a fibration in $\bfA$, for each $C \in \calC$.

\item {\it weak equivalence} if the induced map $F(C) \rightarrow
G(C)$ is a weak equivalence in $\bfA$, for each $C \in \calC$.

\item {\it injective fibration} if it has the right lifting property
with respect to every morphism $\beta$ in $\bfA^{\calC}$ which is
simultaneously a weak equivalence and a injective cofibration.

\item {\it projective cofibration} if it has the left lifting property
with respect to every morphism $\beta$ in $\bfA^{\calC}$ which is
simultaneously a weak equivalence and a projective fibration.
\end{itemize}
\end{definition}

\begin{proposition}\label{smurf}\index{gen}{model category!projective}\index{gen}{model category!injective}
Let $\bfS$ be an excellent model category, $\bfA$ be a combinatorial $\bfS$-enriched model category, and let $\calC$ be a small $\bfS$-enriched category. Then there exist two combinatorial model structures on $\bfA^{\calC}$:

\begin{itemize}
\item The {\it projective model structure}, determined by the strong
cofibrations, weak equivalences, and projective fibrations.

\item The {\it injective model structure}, determined by the weak
cofibrations, weak equivalences, and injective fibrations.
\end{itemize}
\end{proposition}

The proof of Proposition \ref{smurf} is identical to that of Proposition \ref{smurff}, 
except that it requires the following more general form of Lemma \ref{mainerthyme}:

\begin{lemma}\label{mainertime}
Let $\bfA$ be a presentable category which is enriched, tensored, and cotensored over another presentable category $\bfS$. Let $S_0$ be a (small) set of morphisms of $\bfA$, and let $\overline{S}_0$ be the weakly saturated class of morphisms generated by $S_0$. Let $\calC$ be a small $\bfS$-enriched category. Let $\widetilde{S}$ be the collection of all morphisms $F \rightarrow G$ in $\bfA^{\calC}$ with the following property: for every $C \in \calC$, the map $F(C) \rightarrow G(C)$ belongs to $\overline{S}_0$. Then there exists a (small) set of morphisms $S$ of $\bfA^{\calC}$ which generates $\widetilde{S}$ (as a weakly saturated class of morphisms).
\end{lemma}

We will defer the proof until the end of this section.

\begin{remark}\label{suboteki}
In the situation of Proposition \ref{smurf}, the category $\bfA^{\calC}$ is again
enriched, tensored and cotensored over $\bfS$. Tensor product with an object $K \in \bfS$ is computed pointwise; in other words, if $\calF \in \bfA^{\calC}$, then we have the formula $$ (K \otimes \calF)(A) = K \otimes \calF(A).$$
Using criterion $(2')$ of Remark \ref{cyclor}, we deduce that $\bfA^{\calC}$ is a $\bfS$-enriched model category with respect to the injective model structure. A dual argument (using criterion $(2'')$ of Remark \ref{cyclor}) shows that $\bfA^{\calC}$ is also a $\bfS$-enriched model category with respect to the projective model structure.
\end{remark}

\begin{remark}\label{postsmurf}
For each object $C \in \calC$ and each $A \in \bfA$, let
$\calF^{C}_{A} \in \bfA^{\calC}$ be the functor given by
$$D \mapsto A \otimes \bHom_{\calC}(C,D).$$
As in the proof of Proposition \ref{smurff}, we learn that the
class of projective cofibrations in $\bfA^{\calC}$ is generated by
cofibrations of the form $j: \calF^{C}_{A} \rightarrow \calF^{C}_{A'}$, where
$A \rightarrow A'$ is a cofibration in $\bfA$. It follows that every projective cofibration is a injective cofibration; dually, every injective fibration is a projective fibration.
\end{remark}

As in \S \ref{qlim7}, the construction $(\calC, \bfA) \mapsto \bfA^{\calC}$ is functorial
in both $\calC$ and $\bfA$. We summarize the situation in the following
propositions, whose proofs are left to the reader:

\begin{proposition}
Let $\bfS$ be an excellent model category, $\calC$ a small $\bfS$-enriched model category, and
$\Adjoint{F}{\bfA}{\bfS}{G}$ a $\bfS$-enriched Quillen adjunction between
combinatorial $\bfS$-enriched model categories. The composition with $F$ and $G$
determines another $\bfS$-enriched Quillen adjunction
$$ \Adjoint{ F^{\calC} }{ \bfA^{\calC} }{ \bfB^{\calC} }{G^{\calC} },$$
with respect to either the projective or injective model structures.
Moreover, if $(F,G)$ is a Quillen equivalence, then $(F^{\calC}, G^{\calC})$ is also a Quillen equivalence.
\end{proposition}

Because the projective and injective model structures on
$\bfA^{\calC}$ have the same weak equivalences, the identity
functor $\id_{ \bfA^{\calC}}$ is a Quillen equivalence between them. However, it is
important to keep distinguish these two model structures, because
they have different variance properties as we now explain.

Let $f: \calC \rightarrow \calC'$ be a $\bfS$-enriched functor. Then
composition with $f$ yields a pullback functor $f^{\ast}:
\bfA^{\calC'} \rightarrow \bfA^{\calC}$. Since $\bfA$ has all
$\bfS$-enriched limits and colimits, $f^{\ast}$ has a
right adjoint which we shall denote by $f_{\ast}$ and a left
adjoint which we shall denote by $f_{!}$.

\begin{proposition}\label{colbin}
Let $\bfS$ be an excellent model category, $\bfA$ a combinatorial $\bfS$-enriched
model category, and $f: \calC \rightarrow \calC'$ a $\bfS$-enriched functor between
small $\bfS$-enriched categories. Let
$f^{\ast}: \bfA^{\calC'} \rightarrow \bfA^{\calC}$ be given
by composition with $f$. Then $f^{\ast}$ admits a right adjoint
$f_{\ast}$ and a left adjoint $f_{!}$. Moreover:

\begin{itemize}
\item[$(1)$] The pair $( f_{!}, f^{\ast} )$ determines a Quillen
adjunction between the {\em projective} model structures on
$\bfA^{\calC}$ and $\bfA^{\calC'}$.

\item[$(2)$] The pair $( f^{\ast}, f_{\ast} )$ determines a
Quillen adjunction between the {\em injective} model structures on
$\bfA^{\calC}$ and $\bfA^{\calC'}$.
\end{itemize}
\end{proposition}

We now study some aspects of the theory which are unique to the enriched context.

\begin{proposition}\label{lesstrick}
Let $\bfS$ be an excellent model category,
$\bfA$ a combinatorial $\bfS$-enriched model category, and 
$f: \calC \rightarrow \calC'$ an equivalence of small $\bfS$-enriched categories. 
Then:

\begin{itemize}
\item[$(1)$] The Quillen adjunction $(f_{!}, f^{\ast})$ determines a Quillen equivalence between
the projective model structures on $\bfA^{\calC}$ and $\bfA^{\calC'}$. 

\item[$(2)$] The Quillen adjunction $(f^{\ast}, f_{\ast})$ determines a Quillen equivalence between
the injective model structures on $\bfA^{\calC}$ and $\bfA^{\calC'}$.
\end{itemize}

\end{proposition}





\begin{proof}
We first note that $(1)$ and $(2)$ are equivalent: they are both equivalent to the assertion that
$f^{\ast}$ induces an equivalence on homotopy categories. It therefore suffices to prove $(1)$.
We first prove this under the following additional assumption:
\begin{itemize}
\item[$(\ast)$] For every pair of objects $C, D \in \calC'$, the map
$$ \bHom_{\calC'}( C,D) \rightarrow \bHom_{\calC}( f(C), f(D) )$$
is a cofibration in $\bfS$.
\end{itemize}
Let $Lf_{!}: \bfA^{\calC} \rightarrow \bfA^{\calC'}$ denote the left derived functor of $f_{!}$. We must show that the unit and counit maps
$$ h_F: F \mapsto f^{\ast} Lf_{!} F $$
$$ k_G: Lf_{!} f^{\ast} G \rightarrow G$$
are isomorphisms for all $F \in \h{\bfA^{\calC}}$, $G \in \h{\bfA^{\calC}}$. Since $f$ is essentially surjective on homotopy categories, a natural transformation $K \rightarrow K'$ of $\bfS$-enriched functors $K,K': \calC' \rightarrow \bfA$ is a weak equivalence if and only if $f^{\ast} K \rightarrow f^{\ast} K'$ is a weak equivalence. Consequently, to prove $k_G$ is an isomorphism, it suffices to show that
$h_{f^{\ast} G}$ is an isomorphism. 

Let us say that a map $F \rightarrow F'$ in $\bfA^{\calC}$ is {\it good} if the
induced map $f^{\ast} f_{!} F \coprod_{F} F' \rightarrow f^{\ast} f_{!} F'$
is a weak trivial cofibration. To complete the proof, it will suffice to show that
every projective cofibration is good. Since the collection of good transformations
is weakly saturated, it will suffice to show that each of the generating cofibrations
$\calF^{C}_{A} \rightarrow \calF^{C}_{A'}$ is good, where $C \in \calC'$ and
$j: A \rightarrow A'$ is a cofibration in $\bfA$. Unwinding the definitions, we must show that
for each $D \in \calC'$ the induced map
$$ \theta: (A' \otimes \bHom_{\calC'}(C, D)) \coprod_{ A \otimes \bHom_{\calC'}(C,D)}
(A \otimes \bHom_{\calC}(f(C), f(D))) \rightarrow A' \otimes \bHom_{\calC}( f(C), f(D) )$$
is a trivial cofibration. This follows from the fact that $j$ is a cofibration and our assumption $(\ast)$.

We now treat the general case. First, choose a trivial cofibration $g: \calC \rightarrow \calC''$, where
$\calC''$ is fibrant. Then $g$ satisfies $(\ast)$, so $g_!$ is a Quillen equivalence.
By a two-out-of-three argument, we see that $f_!$ is a Quillen equivalence if and only if
$(g \circ f)_!$ is a Quillen equivalence. Replacing $\calC$ by $\calC''$, we may
reduce to the case where $\calC$ is itself fibrant.

Choose a cofibration $j: \calC \coprod \calC' \rightarrow \calD$, where
$\calD$ is fibrant and equivalent to the final object of $\SCat$. Then $f$ factors as a composition
$$ \calC' \stackrel{f'}{\rightarrow} \calC \times \calD \stackrel{f''}{\rightarrow} \calC.$$
Since $\calC$ and $\calD$ are fibrant, the product $\calC \times \calD$ is equivalent
to $\calC$. Moreover, the map $f''$ admits a section $s: \calC \rightarrow \calC \times \calD$.
Using another two-out-of-three argument, it will suffice to show that $f'_{!}$ and $s_!$ are Quillen equivalences. For this, it will suffice to show that $f'$ and $s$ satisfy $(\ast)$.

We first show that $f'$ satisfies $(\ast)$. Fix a pair of objects $X,Y \in \calC'$.
Then $f'$ induces the composite map
$$ \bHom_{\calC'}(X,Y) \stackrel{u}{\rightarrow} \bHom_{\calC}(fX,fY) \times \bHom_{\calC'}(X,Y)
\stackrel{u'}{\rightarrow} \bHom_{\calC}(fX,fY) \times \bHom_{\calD}(jX,jY) \simeq \bHom_{\calC \times \calD}(f'X, f'Y).$$
The map $u$ is a monomorphism (since it admits a left inverse), and therefore a cofibration
in view of axiom $(A2)$ of Definition \ref{modelexcellent}. 
The map $u'$ is a product of cofibrations, and therefore a cofibration (again by axiom $(A2)$).

The proof that $s$ satisfies $(\ast)$ is similar: for every pair of objects $U,V \in \calC$,
the map $$ \bHom_{\calC}(U,V) \rightarrow \bHom_{ \calC \times \calD}( sU, sV)
\simeq \bHom_{ \calC}( U,V) \times \bHom_{\calD}(jU,jV)$$
is a monomorphism since it admits a left inverse, and is therefore a cofibration.
\end{proof}

In the special case where $f: \calC \rightarrow \calC'$ is a {\em cofibration} between
$\bfS$-enriched categories, we have some additional functoriality:

\begin{proposition}\label{sumner}
Let $\bfS$ be an excellent model category and let $f: \calC \rightarrow \calC'$ be a cofibration
of small $\bfS$-enriched categories. Then:
\begin{itemize}

\item[$(1)$] For every combinatorial $\bfS$-enriched model category $\bfA$, the pullback map
$f^{\ast}: \bfA^{\calC'} \rightarrow \bfA^{\calC}$ preserves projective cofibrations. 

\item[$(2)$] For every projectively cofibrant object $F \in \bfS^{\calC}$, the
unit map $F \rightarrow f^{\ast} f_{!} F$ is a projective cofibration.

\end{itemize}
\end{proposition}

\begin{lemma}\label{pseudopod}
Let $\bfS$ be an excellent model category, and suppose given a pushout diagram
$$ \xymatrix{ [1]_{S} \ar[r] \ar[d]^{i} & [1]_{S'} \ar[d] \\
\calC \ar[r]^{f} & \calC' }$$
of $\bfS$-enriched categories, where $j: S \rightarrow S'$ is a cofibration in $\bfS$.
Let $C$ be an object of $\calC$, and let $F \in \bfS^{\calC}$ be the functor
given by the formula $D \mapsto \bHom_{\calC}(C, D)$. Then
the unit map $F \rightarrow f^{\ast} f_{!} F$ is a projective cofibration in
$\bfS^{\calC}$.
\end{lemma}

\begin{proof}
The map $i$ determines a pair of objects $X,Y \in \calC$, and a map
$S \rightarrow \bHom_{\calC}(X,Y)$. The proof of Proposition \ref{enrichcatper}
shows that the functor $f^{\ast} f_{!} F$ is the colimit of a sequence
$$F = F(0) \stackrel{h_1}{\rightarrow} F(1) \stackrel{h_2}{\rightarrow} F(2) \rightarrow \ldots,$$
where each $h_k$ is a pushout of a map 
$\calF^{Y}_{A} \rightarrow \calF^{Y}_{A'}$ induced by a map
$t: A \rightarrow A'$ in $\bfS$. Moreover, the map
$t$ can be identified with the tensor product 
$$\id_{ \bHom_{\calC}(C, X)} \otimes \id_{ \bHom_{\calC}(Y,X) }^{\otimes k-1}
\otimes \wedge^{k}(j),$$ 
where $\wedge^{k}(j)$ denotes the $k$th pushout power of $j$. 
It follows that $t$ is a cofibration in $\bfS$, so that each $h_k$ is a projective cofibration
in $\bfS^{\calC}$. 
\end{proof}

\begin{proof}[Proof of Proposition \ref{sumner}]
The collection of $\bfS$-enriched functors $f$ which satisfy $(1)$ and $(2)$ is clearly
closed under the formation of retracts. We may there assume without loss of generality
that $f$ is a transfinite composition of pushouts of generating cofibrations
(see the discussion preceding Proposition \ref{enrichcatper}). Reordering
these pushouts if necessary, we can factor $f$ as a composition
$$ \calC \stackrel{f'}{\rightarrow} \overline{\calC} \stackrel{f''}{\rightarrow} \calC'$$
where $\overline{\calC}$ is obtained from $\calC$ by freely adjoining a collection of
new objects, and $f''$ is bijective on objects. Since $f'$ clearly satisfies
$(1)$ and $(2)$, it will suffice to prove that $f''$ satisfies $(1)$ and $(2)$.
Replacing $f$ by $f''$, we may assume that $f$ is bijective on objects.

We now show that $(2) \Rightarrow (1)$. Since the functor $f^{\ast}$ preserves colimits, the collection of morphisms $\alpha$ in $\bfA^{\calC'}$ such that $f^{\ast}$ is a projective cofibration in $\bfA^{\calC}$ is weakly saturated. It will therefore suffice to that for every object $X \in \calC'$ and every cofibration $A \rightarrow A'$ in $\bfA$, if $\alpha: \calF^{X}_{A} \rightarrow \calF^{X}_{A'}$ denotes the corresponding generating projective cofibration, then $f^{\ast}(\alpha)$ is a projective cofibration in $\bfS$.

There is a canonical left Quillen bifunctor
$$ \boxtimes: \bfS^{\calC} \times \bfA \rightarrow \bfA^{\calC}$$
described by the formula $(F \boxtimes A)(C) = F(C) \otimes A$.
(Here we regard $\bfS^{\calC}$ as endowed with the projective model structure.)
We observe that $f^{\ast}(\alpha)$ is the induced map
$(f^{\ast} F) \boxtimes A \rightarrow (f^{\ast} F) \boxtimes A'$, where
$F \in \bfS^{\calC'}$ is given by $F(C') = \bHom_{\calC'}( X, C')$. 
To prove $(1)$, it will suffice to show that $f^{\ast} F$ is projectively cofibrant.

Since $F$ is bijective on objects, we can choose an object $X_0 \in \calC$ such that
$fX_0 = X$. We now observe that $F \simeq f_{!} F_0$, where $F_0 \in \bfS^{\calC}$
is defined by the formula $F_0(C) = \bHom_{\calC}(X_0, C)$. If $(2)$ is satisfied, then
the unit map $F_0 \rightarrow f^{\ast} F$ is a projective cofibration in $\bfS^{\calC}$.
Since $F_0$ is projectively cofibrant, we conclude that $f^{\ast} F$ is projectively cofibrant as well.
This completes the proof that $(2) \Rightarrow (1)$.

To prove $(2)$, let us write $f$ as a transfinite composition of $\bfS$-enriched functors
$$ \calC = \calC_0 \rightarrow \calC_1 \rightarrow \ldots, $$
each of which is a pushout of a generating cofibration of the form $[1]_{S} \rightarrow [1]_{S'}$, where
$S \rightarrow S'$ is a cofibration in $\bfS$. 
For each $\alpha \leq \beta$, let $f^{\beta}_{\alpha}: \calC_{\alpha} \rightarrow \calC_{\beta}$
be the corresponding cofibration. We will prove that the following statement holds, for
every pair of ordinals $\alpha \leq \beta$:
\begin{itemize}
\item[$(2_{\alpha,\beta})$] For every projectively cofibrant object $F \in \bfS^{\calC_{\alpha} }$,
the unit map $u: F \rightarrow (f_{\alpha}^{\beta})^{\ast} (f_{\alpha}^{\beta})_{!} F$
is a projective cofibration.
\end{itemize}

The proof proceeds by induction on $\beta$. We observe that $u$ is a transfinite
composition of maps of the form
$$ u_{\gamma}: (f_{\alpha}^{\gamma})^{\ast} (f_{\alpha}^{\gamma})_{!} F \rightarrow
(f_{\alpha}^{\gamma})^{\ast} (f_{\gamma}^{\gamma+1})^{\ast}
(f_{\gamma}^{\gamma+1})_{!} (f_{\alpha}^{\gamma})_{!} F,$$
where $\gamma < \beta$. It will therefore suffice to show that each $u_{\gamma}$ is a projective cofibration. Our inductive hypothesis therefore guarantees that
$(2_{\alpha, \gamma})$ holds, so the first part of the proof shows that
$(f_{\alpha}^{\gamma})^{\ast}$ preserves trivial cofibrations. We are therefore reduced to proving
assertion $(2_{\gamma, \gamma+1})$. In other words, to prove $(2)$ in general, it will suffice to
treat the case in which $f$ is a pushout of a generating cofibration of the form $[1]_{S} \rightarrow [1]_{S'}$. 

We will in fact prove the following stronger version of $(2)$:
\begin{itemize}
\item[$(3)$] For every projective cofibration $\phi: F' \rightarrow F$ in $\bfS^{\calC}$, the induced map
$\phi': F \coprod_{F'} f^{\ast} f_{!} F' \rightarrow f^{\ast} f_{!} F$ is again a projective cofibration in 
$\bfS^{\calC}$. 
\end{itemize}
Consider the collection of {\em all} morphisms $\phi: F' \rightarrow F$ in $\bfS^{\calC}$ such that
the induced map $\phi': F \coprod_{F'} f^{\ast} f_{!} F' \rightarrow f^{\ast} f_{!} F$ is
a projective cofibration. It is easy to see that this collection is weakly saturated. Consequently, to prove
$(3)$ it suffices to treat the case where $\phi$ is a generating projective cofibration of the form
$\calF^{C}_{A} \rightarrow \calF^{C}_{A'}$, where $A \rightarrow A'$ is a cofibration in $\bfS$. 
In this case, we can identify $\phi'$ with the map 
$$ (F_C \boxtimes A') \coprod_{ F_C \boxtimes A } (f^{\ast} f_{!} F_C \boxtimes A) \rightarrow f^{\ast} f_{!} F_C \boxtimes A',$$
where $F_C \in \bfS^{\calC}$ is the functor $D \mapsto \bHom_{\calC}(C, D)$. Since
$\boxtimes$ is a left Quillen bifunctor, it will suffice to show that the unit map
$f_C \rightarrow f^{\ast} f_{!} F_C$ is a projective cofibration in $\bfS^{\calC}$. This is precisely the content of Lemma \ref{pseudopod}.
\end{proof}

In \S \ref{qlim7}, we introduced the definitions of homotopy limits and colimits in an arbitrary
combinatorial model category $\bfA$. We now discuss an analogous construction in the case where $\bfA$ is enriched over an excellent model category $\bfS$. To simplify the exposition, we will discuss only the case of homotopy limits; the case of homotopy colimits is entirely dual and left to the reader.

Fix an excellent model category $\bfS$ and a combinatorial $\bfS$-enriched model category $\bfA$. Let $f: \calC \rightarrow \calC'$ be a functor between small $\bfS$-enriched categories,
so that we have an induced Quillen adjunction
$$ \Adjoint{f^{\ast}}{\bfA^{\calC'}}{\bfA^{\calC}.}{f_{\ast}}.$$
We will refer to the right derived functor $Rf_{\ast}$ as the {\it homotopy right Kan extension} functor.\index{gen}{Kan extension!homotopy}\index{gen}{homotopy right Kan extension}\index{gen}{right Kan extension!homotopy} If we are given a pair of functors $F \in \bfA^{\calC}$, $G \in \bfA^{\calC'}$, and let $\eta: G \rightarrow f_{\ast} F$ be a map in $\bfA^{\calC'}$. We will say that $\eta$ {\it exhibits $G$ as the homotopy right Kan extension of $F$}
if, for some weak equivalence $F \rightarrow F'$ where $F'$ is injectively fibrant in $\bfA^{\calC}$, the
composite map $G \rightarrow f_{\ast} F \rightarrow f_{\ast} F'$ is a weak equivalence in
$\bfA^{\calC'}$. Since $f_{\ast}$ preserves weak equivalences between injectively fibrant objects, this condition is independent of the choice of $F'$.\index{gen}{Kan extension!homotopy}

\begin{remark}
In \S \ref{qlim7}, we defined homotopy right Kan extensions in the setting of
the diagram categories $\Fun( \calC, \bfA)$, where $\calC$ is an ordinary category.
In fact, this is a special case of the above construction. Namely, there is a unique colimit-preserving monoidal functor $F: \Set \rightarrow \bfS$, given by $F(S) = \coprod_{s \in S} {\bf 1}_{\bfS}$. We can therefore define a $\bfS$-enriched category $\overline{\calC}$ whose objects are the objects of $\calC$, with $\bHom_{ \overline{\calC} }(X,Y) = F \bHom_{\calC}(X,Y)$. We
now observe that we have an identification $\Fun( \calC, \bfA) \simeq \bfA^{\overline{\calC}}$, which is functorial in both $\calC$ and $\bfA$. This identification is compatible with the definition
of the injective model structures on both sides, so that either point of view gives rise to the same theory of homotopy right Kan extensions.
\end{remark}

We now discuss a special feature of the enriched theory of homotopy Kan extensions:
they can be reduced to the theory of homotopy Kan extensions in the model category $\bfS$:

\begin{proposition}\label{usecoinc}
Let $\bfS$ be an excellent model category, $\bfA$ a combinatorial model category
enriched over $\bfS$, and let $f: \calC \rightarrow \calC'$ be a functor between
small $\bfS$-enriched categories. Suppose
given objects $F \in \bfA^{\calC}$, $G \in \bfA^{\calC'}$, and a map $\eta: G \rightarrow f_{\ast} F$. 
Assume that $F$ and $G$ are projectively fibrant.
The following conditions are equivalent:

\begin{itemize}
\item[$(1)$]  The map $\eta$ exhibits $G$ as a homotopy right Kan extension of $F$.

\item[$(2)$] For each cofibrant object $A \in \bfA$, the induced map
$$ \eta_A: G_A \rightarrow f_{\ast} F_A$$
exhibits $G_A$ as a homotopy right Kan extension of $F_A$. Here $F_A \in \bfS^{\calC}$
and $G_A \in \bfS^{\calC'}$ are defined by $F_A(C) = \bHom_{\bfA}(A,F(C)), G_A(C) = \bHom_{\bfA}(A,G(C))$. 

\item[$(3)$] For every fibrant-cofibrant object $A \in \bfA$, the induced map
$$ \eta_A: G_A \rightarrow f_{\ast} F_A$$
exhibits $G_A$ as a homotopy right Kan extension of $F_A$.
\end{itemize}
\end{proposition}

\begin{proof}
Choose an equivalence $F \rightarrow F'$, where $F'$ is injectively fibrant. We note that
the induced maps $F_A \rightarrow F'_A$ are weak equivalences for any
cofibrant $A \in \bfA$, since $\bHom_{\bfA}(A, \bigdot)$ preserves weak equivalences between fibrant objects. Consequently, we may without loss of generality replace $F$ by $F'$ and thereby assume that $F$ is injectively fibrant.

Now suppose that $A$ is any cofibrant object of $\bfA$; we claim that
$F_A$ is injectively fibrant. To show that $F_A$ has
the right lifting property with respect to a trivial weak
cofibration $H \rightarrow H'$ of functors $\calC \rightarrow
\bfS$, one need only observe that $F$ has the right lifting
property with respect to trivial injective cofibration $A \otimes H
\rightarrow A \otimes H'$ in $\bfA^{\calC}$.

Now we note that $(1)$ is equivalent to the assertion that $\eta$ is a weak equivalence, $(2)$
is equivalent to the assertion that $\eta_{A}$ is a weak equivalence for any cofibrant object
$A$, and $(3)$ is equivalent to the assertion that $\eta_A$ is a weak equivalence whenever
$A$ is fibrant-cofibrant. Because $\bHom_{\bfA}(A, \bigdot)$ preserves weak equivalences between fibrant objects, we deduce that $(1) \Rightarrow (2)$. It is
clear that $(2) \Rightarrow (3)$. We will complete the proof by showing that $(3) \Rightarrow (1)$.
Assume that $(3)$ holds; we must show that
$\eta(C'): G(C') \rightarrow f_{\ast} F(C')$ is an isomorphism in the homotopy category $\h{\bfA}$,
for each $C' \in \calC'$.
For this, it suffices to show that $G(C')$ and $f_{\ast} F(C')$ represent the same $\calH$-valued functors on the homotopy category $\h{\bfA}$, which is precisely the content of $(3)$. 
\end{proof}

\begin{remark}\label{curble}
It follows from Proposition \ref{usecoinc} that we can make sense of homotopy right Kan extensions for diagrams
which do not take values in a model category. 
Let $f: \calC \rightarrow \calC'$ be a $\bfS$-enriched functor as in the discussion above, and let $\calA$ be an {\em arbitrary} locally fibrant $\bfS$-enriched category. Suppose given objects $F \in \calA^{\calC}$, $G \in \calA^{\calC'}$, and $\eta: f^{\ast} G \rightarrow F$
we say that $\eta$ {\it exhibits $G$ as a homotopy right
Kan extension of $F$} if, for each object $A \in \calA$, the induced map
$$\eta_A: G_A \rightarrow f_{\ast} F_A$$ exhibits $G_A \in \bfS^{\calC'}$ as a homotopy right Kan extension of $F_A \in \bfS^{\calC}$.\index{gen}{homotopy limit}\index{gen}{limit!homotopy}

Suppose that the monoidal structure on $\bfS$ is given by the Cartesian product, and take
$\calC'$ to be the final object of $\SCat$, so that we can identify $\calA^{\calC'}$ with $\calA$.
In this case, we can identify $G$ with a single object $B \in \calA$, and the map $\eta$ with
a collection of maps $ \{ B \rightarrow F(C) \}_{C \in \calC}$. We will say that $\eta$ 
{\it exhibits $B$ as a homotopy limit of $F$} if it identifies $G$ with a homotopy right Kan extension of $F$. In other words, $\eta$ exhibits $B$ as a homotopy limit of $F$ if, for every object
$A \in \calA$, the induced map
$$ \bHom_{\calA}(A,B) \rightarrow \lim \{ \bHom_{\calA}(A, F(C)) \}_{C \in \calC}$$
exhibits $\bHom_{\calA}(A,B)$ as a homotopy limit of the diagram
$\{ \bHom_{\calA}(A, F(C) \}_{C \in \calC}$ in the model category $\bfS$.

We also have a dual notion of {\em homotopy colimit} in an arbitrary fibrant
$\bfS$-enriched category $\calA$: a compatible family of maps $\{ F(C) \rightarrow B \}_{C \in \calC}$ {\it exhibits $B$ as a homotopy colimit of $F$} if, for every object $A \in \calA$, the
induced maps $\{ \bHom_{\calA}( B, A) \rightarrow \bHom_{\calA}( F(C), A) \}_{C \in \calC}$
exhibit $\bHom_{\calA}(B,A)$ as a homotopy limit of the diagram
$\{ \bHom_{\calA}( F(C), A) \}_{C \in \calC}$ in $\bfS$. \index{gen}{homotopy colimit}\index{gen}{colimit!homotopy}
\end{remark}

\begin{remark}
In view of Proposition \ref{usecoinc}, the terminology introduced in Remark \ref{curble} for general $\calA$ agrees with the terminology introduced for a combinatorial $\bfS$-enriched model category $\bfA$ if we set $\calA = \bfA^{\degree}$. We remark that, in general, the two notions do {\em not} agree if
we take $\calA = \bfA$, so that our terminology is potentially ambiguous; however, we feel that there is little danger of confusion.
\end{remark}

We conclude this section by giving the proof of Lemma \ref{mainertime}. 
Let $\bfA$ be a presentable category which is enriched, tensored, and cotensored over a presentable category $\bfS$. Let $\calC$ be a small $\bfS$-enriched category, and $\overline{S}_0$ a weakly saturated class of morphisms of $\bfA$ generated by a (small) set $S_0$. We regard this data as {\em fixed} for the remainder of this section.

Choose a regular cardinal $\kappa$ satisfying the following conditions:
\begin{itemize}
\item[$(i)$] The cardinal $\kappa$ is uncountable.

\item[$(ii)$] The category $\calC$ has fewer than $\kappa$-objects.

\item[$(iii)$] Let $X,Y \in \calC$, and let $K = \bHom_{\calC}(X,Y)$. Then the functor
from $\bfA$ to itself given by the formula $A \mapsto A^{K}$ preserves $\kappa$-filtered colimits.
This implies, in particular, that the collection of $\kappa$-compact objects of $\bfA$ is stable with respect to the functors $\bigdot \otimes K$.

\item[$(iv)$] The category $\bfA$ is $\kappa$-accessible. It follows also that $\bfA^{\calC}$ is
$\kappa$-accessible, and that an object $F \in \bfA^{\calC}$ is $\kappa$-compact if and only if
each $F(C) \in \bfA$ is $\kappa$-compact. We prove an $\infty$-category generalization of this
statement as Proposition \ref{horse1}. The same proof also works in the setting of ordinary categories.

\item[$(v)$] The source and target of every morphism in $S_0$ is a $\kappa$-compact object of $\bfA$.
\end{itemize}

Enlarging $S_0$ if necessary, we may assume that $S_0$ consists of {\em all} morphisms in
$f \in \overline{S}_0$ such that the source and target of $f$ are $\kappa$-compact.
Let $S$ be the collection of all injective cofibrations between $\kappa$-compact objects of $\bfA$ (in view of $(iv)$, we can equally well define $S$ to be the set of morphisms $F \rightarrow G$ in
$\bfA^{\calC}$ such that each of the induced morphisms $F(C) \rightarrow G(C)$ belongs to $S_0$). Let $\overline{S}$ be the weakly saturated class of morphisms in $\bfA^{\calC}$ generated by $S$, and choose a map $f: F \rightarrow G$ in $\bfA^{\calC}$ such that $f(C) \in \overline{S}_0$ for each $C \in \calC$. We wish to show that
$f \in \overline{S}$. Corollary \ref{unitape} implies that, for each $C \in \calC$, there exists a $\kappa$-good $S_0$-tree $\{ Y(C)_{\alpha} \}_{\alpha \in A(C)}$ with root $F(C)$ and colimit $G(C)$.

Let us define a {\it slice} to be the following data:
\begin{itemize}
\item[$(a)$] For each object $C \in \calC$, a downward-closed subset $B(C) \subseteq A(C)$.
\item[$(b)$] For every object $C \in \calC$, a morphism
$\eta_{C}: \coprod_{C' \in \calC} Y(C')_{B(C')} \otimes \bHom_{\bfA}(C',C) \rightarrow Y(C)_{B(C)}$, rendering the following diagrams commutative:
$$ \xymatrix{ \coprod_{C'',C' \in \calC} Y(C'')_{B(C'')} \otimes \bHom_{\bfA}(C'',C') \otimes \bHom_{\bfA}(C',C) \ar[r] \ar[d]^{\eta_{C'}} &
\coprod_{C'' \in \calC} Y(C'')_{B(C'')} \otimes \bHom_{\bfA}(C'',C) \ar[d]^{\eta_{C''}} \\
\coprod_{C' \in \calC} Y(C')_{B(C')} \otimes \bHom_{\bfA}(C',C) \ar[r]^{\eta_{C}} & Y(C)_{B(C)} }$$

$$ \xymatrix{ \coprod_{C' \in \calC} F(C') \otimes \bHom_{\bfA}(C',C) \ar[r] \ar[d] & F(C) \ar[d] \\
 \coprod_{C' \in \calC} Y(C')_{B(C')} \otimes \bHom_{\bfA}(C',C) \ar[d] \ar[r]^-{\eta_{C}} &
Y(C)_{B(C)} \ar[d] \\
\coprod_{C' \in \calC} G(C')  \otimes \bHom_{\bfA}(C',C) \ar[r] & G(C). }$$
\end{itemize}

We remark that $(b)$ is precisely the data needed to endow $C \mapsto Y(C)_{B(C)}$ with the structure of a $\bfS$-enriched functor $\calC \rightarrow \bfA$, lying between $F$ and $G$ in $\bfA^{\calC}$. 

\begin{lemma}\label{umpin}
Suppose given a collection of $\kappa$-small subsets $\{ B_0(C) \subseteq A(C) \}_{C \in \calC}$. Then there exists a slice $\{ (B(C), \eta_{C} \}_{C \in \calC}$
such that each $B(C)$ is a $\kappa$-small subset of $A(C)$ containing $B_0(C)$.
\end{lemma}

\begin{proof}
Enlarging each $B_0(C)$ if necessary, we may assume that each $B_0(C)$ is closed downwards. 
Note that because each $\{ Y(C)_{\alpha} \}_{\alpha \in A(C)}$ is a $\kappa$-good $S_0$-tree, if
$A' \subseteq A(C)$ is closed downward and $\kappa$-small, $Y(C)_{A'}$ is $\kappa$-compact when viewed as an object of $\bfA_{F(C)/}$. It follows from $(iii)$ that each $Y(C)_{B_0(C)} \otimes \bHom_{\bfA}(C,C')$ is a $\kappa$-compact object of $\bfA_{( F(C) \otimes \bHom_{\bfA}(C,C') )/}$. Consequently, each composition
$$ \coprod_{C' \in \calC} Y(C')_{B_0(C')} \otimes \bHom_{\bfA}(C',C)
{\rightarrow} \coprod_{C' \in \calC} G(C') \otimes \bHom_{\bfA}(C',C) 
\rightarrow G(C)$$
admits another factorization
$$ \coprod_{C' \in \calC} Y(C')_{B_0(C')} \otimes \bHom_{\bfA}(C',C)
\stackrel{\eta^1_{C}}{\rightarrow} Y(C)_{B_1(C)} \rightarrow G(C),$$
where $B_1(C)$ is downward closed and $\kappa$-small, and the diagram
$$ \xymatrix{ \coprod_{C' \in \calC} F(C') \otimes \bHom_{\bfA}(C',C) \ar[d] \ar[r] & 
\coprod_{C' \in \calC} Y(C')_{B_0(C')} \ar[d]^{\eta^1_C} \ar[d] \\
F(C) \ar[r] & Y(C)_{B_1(C)}}$$
commutes. Enlarging $B_1(C)$ if necessary, we may suppose that each $B_1(C)$ contains $B_0(C)$.

We now continue the preceding construction by defining, for each $C \in \calC$, a sequence of $\kappa$-small, downward closed subsets
$$ B_0(C) \subseteq B_1(C) \subseteq B_2(C) \subseteq \ldots $$
of $A(C)$, and maps 
$\eta^{i}_{C}: \coprod_{C' \in \calC} Y(C')_{B_{i-1}(C')} \otimes
\bHom_{\bfA}(C',C) \rightarrow Y(C)_{B_{i}(C)}$. 
Suppose that $i > 1$, and that the sets $B_{j}(C)$ and maps
$\eta^j_{C}$ have been constructed for $j < i$. Using a compactness argument, we conclude that the composition
$$ \coprod_{C' \in \calC} Y(C')_{B_{i-1}(C')} \otimes \bHom_{\bfA}(C',C)
\rightarrow \coprod_{C' \in \calC} G(C') \otimes \bHom_{\bfA}(C',C)
\rightarrow G(C)$$
coincides with
$$ \coprod_{ C' \in \calC} Y(C')_{B_{i-1}(C')} \otimes \bHom_{\bfA}(C',C)
\stackrel{\eta^i_{C}}{\rightarrow} Y(C)_{B_{i}(C)} \rightarrow G(C),$$
where $B_{i}(C)$ is $\kappa$-small and the diagram 
$$ \xymatrix{ \coprod_{C' \in \calC} F(C') \otimes \bHom_{\bfA}(C',C) \ar[d] \ar[r] & 
\coprod_{C' \in \calC} Y(C')_{B_{i-1}(C')} \otimes \bHom_{\bfA}(C',C) \ar[d]^{\eta^i_C} \ar[d] \\
F(C) \ar[r] & Y(C)_{B_i(C)}}$$
commutes. Enlarging $B_{i}(C)$ if necessary, we may suppose that $B_{i}(C)$ contains
$B_{i-1}(C)$ and that the following diagrams commute as well:
$$ \xymatrix{ \coprod_{C', C'' \in \calC} Y(C'')_{ B_{i-2}(C'')} \otimes \bHom_{\bfA}(C'',C')
\otimes \bHom_{\bfA}(C',C) \ar[r] \ar[d] & \coprod_{C'' \in \calC} Y(C'')_{B_{i-1}(C'') }
\otimes \bHom_{\bfA}(C'',C) \ar[d]^{\eta_{C}^{i}} \\
\coprod_{C' \in \calC} Y(C')_{B_{i-1}(C')} \otimes \bHom_{\bfA}(C',C) \ar[r]^{\eta^{i}_{C}} & 
Y(C)_{B_i(C)}}$$
$$ \xymatrix{ \coprod_{C' \in \calC} Y(C')_{B_{i-2}(C')} \otimes \bHom_{\bfA}(C',C) \ar[r] \ar[d]^{\eta^{i-1}_{C}} & \coprod_{C' \in \calC} Y(C')_{B_{i-1}(C')} \otimes \bHom_{\bfA}(C',C) \ar[d]^{\eta^{i}_{C}}
\\ Y(C)_{B_{i-1}(C)} \ar[r] & Y(C)_{B_{i}(C)}. }$$
We now define $B(C) = \bigcup B_i(C)$, and $\eta_{C}$ to be the amalgam of the compositions
$$ \coprod_{C' \in \calC} Y(C')_{B_{i-1}(C')} \otimes \bHom_{\bfA}(C',C)
\stackrel{\eta^i_{C}}{\rightarrow} Y(C)_{B_{i}}(C) \rightarrow Y(C)_{B(C)}.$$
\end{proof}

We now introduce a bit more terminology. Suppose given a pair of slices
$M = \{ ( B(C), \eta_C )\}_{C \in \calC}$, $M' = \{ ( B'(C), \eta'_{C} \}_{C \in \calC} \}$. We will say that $M$ is {\it $\kappa$-small} if each $B(C)$ is $\kappa$-small. We will say that $M'$ {\it extends} $M$
if $B(C) \subseteq B'(C)$ for each $C \in \calC$, and each diagram
$$ \xymatrix{ \coprod_{ C' \in \calC} Y(C')_{B(C')} \otimes \bHom_{\bfA}(C',C) \ar[r] \ar[d]^{\eta_C} &
\coprod_{C' \in \calC} Y(C')_{B'(C')} \otimes \bHom_{\bfA}(C',C) \ar[d]^{\eta'_{C}} \\
Y(C)_{B(C)} \ar[r] & Y(C)_{B'(C)} }$$
is commutative. Equivalently, $M'$ extends $M$ if $B(C) \subseteq B'(C)$ for each $C \in \calC$, and the induced maps $Y(C)_{B(C)} \rightarrow Y(C)_{B(C')}$ constitute a natural transformation of simplicial functors from $\calC$ to $\bfA$.

\begin{lemma}\label{goosebed}
Let $M' = \{ ( A'(C), \theta_{C} ) \}_{C \in \calC}$ be a slice, and let 
$\{ B_0(C) \subseteq A(C) \}_{C \in \calC}$ be a collection of $\kappa$-small subsets of
$A(C)$. Then there exists a pair of slices 
$N = \{ ( B(C), \eta_{C} ) \}_{C \in \calC}$, $N' = \{ (B(C) \cap A'(C), \eta'_{C}) \}$ where
$B(C)$ is $\kappa$-small, and $N'$ is compatible with both $N$ and $M'$.
\end{lemma}

\begin{proof}
Let $B'_0(C) = A'(C) \cap B_0(C)$. For every positive integer $i$, we will construct a pair of slices
$N_i = \{ ( B_i(C), \eta(i)_{C} ) \}$, $N'_{i} = \{ (B'_{i}(C), \eta'(i)_{C} ) \}$ so that the following conditions are satisfied:
\begin{itemize}
\item[$(a)$] Each $B_{i}(C)$ is $\kappa$-small and contains $B_{i-1}(C)$.  
\item[$(b)$] Each $B'_{i}(C)$ is $\kappa$-small, contains
$B'_{i-1}(C)$ and the intersection $B_{i}(C) \cap A'(C)$, and is contained in $A'(C)$.
\item[$(c)$] Each $N'_{i}$ is compatible with $M'$.
\item[$(d)$] If $i > 2$ and $C \in \calC$, then the diagram
$$ \xymatrix{ \coprod_{C' \in \calC} Y(C')_{B_{i-2}(C')} \otimes \bHom_{\bfA}(C',C)
\ar[r] \ar[d]^{ \eta(i-2)_{C} } & \coprod_{C' \in \calC} Y(C')_{B_{i-1}(C') } \otimes \bHom_{\bfA}(C',C) \ar[d]^{\eta(i-1)_{C}} \\
Y(C)_{B_{i-2}(C)} \ar[d] & Y(C)_{B_{i-1}}(C) \ar[d] \\
Y(C)_{B_{i}(C)} \ar@{=}[r] & Y(C)_{B_i(C)} }$$
commutes. 
\item[$(e)$] If $i > 2$ and $C \in \calC$, then the diagram
$$ \xymatrix{ \coprod_{C' \in \calC} Y(C')_{B'_{i-2}(C')} \otimes \bHom_{\bfA}(C',C)
\ar[r] \ar[d]^{ \eta'(i-2)_{C} } & \coprod_{C' \in \calC} Y(C')_{B'_{i-1}(C') } \otimes \bHom_{\bfA}(C',C) \ar[d]^{\eta('i-1)_{C}} \\
Y(C)_{B'_{i-2}(C)} \ar[d] & Y(C)_{B'_{i-1}}(C) \ar[d] \\
Y(C)_{B'_{i}(C)} \ar@{=}[r] & Y(C)_{B'_i(C)} }$$
commutes. 
\item[$(f)$] If $i > 1$ and $C \in \calC$, then the diagram
$$ \xymatrix{ \coprod_{C' \in \calC} Y(C')_{B'_{i-1}(C')} \otimes \bHom_{\bfA}(C',C) \ar[r] 
\ar[d]^{\eta'(i-1)_{C}} & \coprod_{C' \in \calC} Y(C')_{B_{i-1}(C')} \otimes \bHom_{\bfA}(C',C) \ar[d]^{\eta(i-1)_{C}} \\
Y(C)_{B'_{i-1}(C)} \ar[d] & Y(C)_{B_{i-1}(C)} \ar[d] \\
Y(C)_{B'_{i}(C) } \ar[r] & Y(C)_{B'_{i}(C)}  }$$
commutes.
\end{itemize}
The construction is by induction on $i$. 
The existence of $N_{i}$ satisfying $(a)$, $(d)$, and $(f)$ follows from Lemma \ref{umpin} (and a compactness argument). Similarly, the existence of $N'_{i}$ satisfying $(b)$, $(c)$, and $(e)$ follows by applying Lemma \ref{umpin} after replacing $G \in \bfA^{\calC}$ by
the functor $G'$ given by $G'(C) = Y(C)_{A'(C)}$, and the $S_0$-trees
$\{ Y(C)_{\alpha} \}_{\alpha \in A(C)}$ by the smaller trees $\{ Y(C)_{\alpha} \}_{ \alpha \in A'(C) }$.  

We now define $B(C) = \bigcup_{i} B_i(C)$. It follows from $(d)$ that the $\eta(i)_{C}$ assemble to
a map $$\eta_{C}: \coprod_{C' \in \calC} Y(C')_{B(C')} \otimes \bHom_{\bfA}(C',C) \rightarrow Y(C)_{B(C)}.$$ Taken together these maps determine a slice $N = \{ (B(C), \eta_{C} ) \}$. Similarly, $(e)$ implies that
the maps $\eta'(i)_{C}$ assemble to a slice $N' = \{ ( B(C) \cap A'(C), \eta'_{C} ) \}$. The compatibility of $N$ and $N'$ follows from $(f)$, while the compatibility of $M'$ and $N'$ follows from $(c)$.
\end{proof}

We now construct a transfinite sequence of compatible slices $\{ M(\gamma) = \{ (B(\gamma)(C), \eta(\gamma)_{C} ) \}_{C \in \calC} \}_{\gamma < \beta}$. The construction is by recursion. Assume that $M(\gamma')$ has been defined for $\gamma' < \gamma$, and let 
$M'(\gamma) = \{ ( B'(\gamma)(C), \eta'(\gamma)_{C} ) \}_{C \in \calC}$ be the slice
obtained by amalgamating the family $\{ M(\gamma') \}_{\gamma' < \gamma}$. If
$B'(\gamma)(C) = A(C)$ for all $C \in \calC$, we set $\beta = \gamma$ and conclude the construction. Otherwise, choose $C \in \calC$ and $a \in A(C) - B'(\gamma)(C)$. According to 
Lemma \ref{goosebed}, there exists a pair of slices
$N(\gamma) = \{ ( B''(C), \theta_C) \}_{C \in \calC}$, 
$N'(\gamma) = \{ ( B''(C) \cap B'(\gamma)(C), \theta'_{C} \}_{C \in \calC}$ such that
$N'(\gamma)$ is compatible with both $N(\gamma)$ and $M'(\gamma)$. We now define
$M(\gamma)$ to be the slice obtained by amalgamating $M'(\gamma)$ and $N(\gamma)$.

For $\gamma < \beta$, let $G(\gamma): \calC \rightarrow \bfA$ be the simplicial functor corresponding to the slice $M(\gamma)$. Then we have a transfinite sequence of composable morphisms 
$$ G(0) \rightarrow G(1) \rightarrow \ldots$$
in $(\bfA^{\calC})_{F/}$ having colimit $G \simeq \colim_{\gamma < \beta} G(\gamma)$. 
Since $\overline{S}$ is weakly saturated, to prove that the map $F \rightarrow G$ belongs to $\overline{S}$, it will suffice to show that for each $\gamma < \beta$, the map 
$$ f_{\gamma}: \colim_{\gamma' < \gamma} G(\gamma') \rightarrow G(\gamma)$$
belongs to $\overline{S}$. We observe that for each $C in \calC$, the map
$f_{\gamma}(C)$ can be identified with the map
$Y(C)_{ B'(\gamma)(C) } \rightarrow Y(C)_{ B(\gamma)(C) }$. Since $B(\gamma)(C) - B'(\gamma)(C)$ is $\kappa$-small, Remark \ref{relci}, Lemma \ref{tiura} and Lemma \ref{uper} imply that $f_{\gamma}$ is the pushout of a morphism belonging to $S_0$. We now conclude by applying the following result (replacing $G$ by $G(\gamma)$ and $F$ by $\colim_{\gamma' < \gamma} G(\gamma')$:

\begin{lemma}
Suppose that $f: F \rightarrow G$ has the property that, for each $C \in \calC$, there exists a pushout diagram
$$ \xymatrix{ X_{C} \ar[r]^{g_{C}} \ar[d] & Y_{C} \ar[d] \\
F(C) \ar[r]^{f(C)} & G(C) }$$
where $g_{C} \in S_0$. Then $f$ is the pushout of a morphism in $S$.
\end{lemma}

\begin{proof}
In view of $(iv)$, we can write $F$ as the colimit of a diagram $\{ F_{\lambda} \}_{ \lambda \in P}$
indexed by a $\kappa$-filtered partially ordered set $P$, where each $F_{\lambda}$ is a $\kappa$-compact object of $\bfA^{\calC}$, and is therefore a functor whose values are $\kappa$-compact objects of $\bfA$. Since each $X_{C} \in \bfA$ is $\kappa$-compact, the map $X_{C} \rightarrow F(C)$ factors through $F_{\lambda(C)}(C)$ for some sufficiently large $\lambda(C) \in P$. Since $\calC$ has fewer than $\kappa$ objects and $P$ is $\kappa$-filtered, we can choose a single
$\lambda \in P$ which works for every object $C \in \calC$. 

Consider, for each $C \in \calC$, the composite map
$$ \coprod_{C' \in \calC} Y_{C'} \otimes \bHom_{\bfA}(C',C)
\rightarrow \coprod_{C' \in \calC} G(C') \otimes \bHom_{\bfA}(C',C) \rightarrow G(C)
\simeq \colim_{\lambda' \in P} F_{\lambda'}(C) \coprod_{X_C} Y_{C}.$$ 
Using another compactness argument, we deduce that this map is equivalent to a composition
$$ \coprod_{C' \in \calC} Y_{C'} \otimes \bHom_{\bfA}(C',C)
\rightarrow F_{\lambda'(C)}(C) \coprod_{X_C} Y_{C} $$ 
for some sufficiently large $\lambda'(C) \in P$. Once again, because $P$ is $\kappa$-filtered we can choose a single $\lambda' \in P$ which works for all $C$. Enlarging $\lambda$ and $\lambda'$, we can assume $\lambda = \lambda'$. Using another compactness argument, we can (after enlarging $\lambda$ if necessary) assume that
each of the diagrams
$$ \xymatrix{ \coprod_{C' \in \calC} X_{C'} \otimes \bHom_{\bfA}(C',C) \ar[r] \ar[d]  &
F_{\lambda}(C) \ar[d] \\
\coprod_{C' \in \calC} Y_{C'} \otimes \bHom_{\bfA}(C',C) \ar[r] & F_{\lambda}(C) \coprod_{X_C} Y_{C} }$$
$$ \xymatrix{ \coprod_{C', C'' \in \calC} Y_{C''}
\otimes \bHom_{\bfA}(C'',C') \otimes \bHom_{\bfA}(C',C) \ar[r] \ar[d] & \coprod_{C'' \in \calC}
Y_{C''} \otimes \bHom_{\bfA}(C'',C) \ar[d] \\
\coprod_{C' \in \calC} ( F_{\lambda}(C') \coprod_{ X_{C'} } Y_{C'} ) \otimes \bHom_{\bfA}(C',C) \ar[r] & 
F_{\lambda}(C) \coprod_{ X_C} Y_C }$$
is commutative. Then the above maps allow us to define a $\bfS$-enriched functor
$G_{\lambda}: \calC \rightarrow \bfA$ by the formula $G_{\lambda}(C) = F_{\lambda}(C) \coprod_{ X_C} Y_{C}$. We now observe that there is a pushout diagram
$$ \xymatrix{ F_{\lambda} \ar[r]^{f_{\lambda}} \ar[d] & G_{\lambda} \ar[d] \\
F \ar[r]^{f} & G }$$
and that $f_{\lambda} \in S$.
\end{proof}

\subsection{Path Spaces in $\bfS$-Enriched Categories}\label{pathspace}

Let $\bfS$ be a excellent model category. We have seen that there exists a model structure on the category $\SCat$ of $\bfS$-enriched categories, whose fibrant objects are precisely
those categories which are enriched over the full subcategory $\bfS^{\degree}$ of fibrant
objects of $\bfS$. 

The theory of model categories provides a plethora of examples:
for every $\bfS$-enriched model category $\bfA$, the full subcategory $\bfA^{\degree} \subseteq \bfA$ of fibrant-cofibrant objects is a fibrant object of $\SCat$.
In other words, $\bfA^{\degree}$ is suitable to use for computing the homotopy set $[\calC, \bfA^{\degree}] = \Hom_{ \h{\SCat} }( \calC, \bfA^{\degree})$: if $\calC$ is cofibrant, then every map from $\calC$ to $\bfA^{\degree}$ in the homotopy category of $\SCat$ is represented by an actual $\bfS$-enriched functor from $\calC$ to $\bfA^{\degree}$. Moreover, two simplicial functors $F,F': \calC \rightarrow \bfA^{\degree}$ represent the same morphism
in $\h{\SCat}$ if and only if they are homotopic to one another. The relation of homotopy can be described either in terms of a cylinder object for $\calC$ or a path object for $\bfA^{\degree}$. Unfortunately, it is rather difficult to construct a cylinder object for $\calC$ explicitly, since the cofibrations in $\SCat$ are difficult to describe directly even when $\bfS = \sSet$ (the class of cofibrations of simplicial categories is {\em not} stable under products, so the usual procedure of
constructing mapping cylinders via ``product with an interval'' cannot be applied). 
On the other hand, Theorem \ref{staycode} gives a good understanding of the fibrations
in $\SCat$, which will allow us to give a very explicit construction of a path object for $\bfA^{\degree}$.

Let $\bfA$ be a $\bfS$-enriched model category. Our goal in this section is to give a direct construction of a path space object for $\bfA^{\degree}$ in $\SCat$. In other words, we wish to supply
a diagram of $\bfS$-enriched categories $$\bfA^{\degree} \rightarrow P( \bfA )
\rightarrow \bfA^{\degree} \times \bfA^{\degree}$$\index{gen}{path object!in simplicial categories}
where the composite map is the diagonal, the left map is a weak equivalence, and the right map is a fibration. For technical reasons, we will find it convenient to work not with the entire
category $\bfA$, but some (usually small) subcategory thereof. For this reason, we
introduce the following definition:

\begin{definition}\index{gen}{chunk}\label{defchunk}
Let $\bfS$ be an excellent model category, and let $\bfA$ be a
combinatorial $\bfS$-enriched model category. A {\it chunk of $\bfA$}
is a full subcategory $\calU \subseteq \bfA$ with the following properties:
\begin{itemize}
\item[$(a)$] Let $A$ be an object of $\calU$, and let $\{ \phi_i: A \rightarrow B_i \}_{i \in I}$ be a finite
collection of morphisms in $\calU$. Then there exists a factorization
$$ A \stackrel{p}{\rightarrow} \overline{A} \stackrel{q}{\rightarrow} \prod_{i \in I} B_i$$
of the product map $\prod_{i \in I} \phi_i$,
where $p$ is a trivial cofibration, $q$ a fibration, and $\overline{A} \in \calU$.
Moreover, this factorization can be chosen to depend functorially on the collection
$\{ \phi_i \}$, via a $\bfS$-enriched functor.

\item[$(b)$] Let $A$ be an object of $\calU$, and let $\{ \phi_i: B_i \rightarrow A \}_{i \in I}$ be a finite
collection of morphisms in $\calU$. Then there exists a factorization
$$ \coprod_{i \in I} B_i \stackrel{p}{\rightarrow} \overline{A} \stackrel{q}{\rightarrow} A$$
of the coproduct map $\coprod_{i \in I} \phi_i$,
where $p$ is a cofibration, $q$ a trivial fibration, and $\overline{A} \in \calU$.
Moreover, this factorization can be chosen to depend functorially on the collection
$\{ \phi_i \}$, via a $\bfS$-enriched functor.
\end{itemize}

If $\calU$ is a chunk of $\calA$, we let $\calU^{\degree}$\index{not}{Udeg@$\calU^{\degree}$}
denote the full subcategory $\bfA^{\degree} \cap \calU \subseteq \calU$ consisting
of fibrant-cofibrant objects of $\bfA$ which belong to $\calU$.

We will say that two chunks $\calU, \calU' \subseteq \bfA$ are {\it equivalent} if
they have the same essential image in the homotopy category $\h{\bfA}$.
\end{definition}

\begin{remark}
In particular, if $\calU$ is a chunk of $\bfA$, then each object $A \in \calU$ admits (functorial)
fibrant and cofibrant replacements which also belong to $\calU$ (take the set $I$ to be empty in
$(a)$ and $(b)$). 
\end{remark}

\begin{remark}
If $\calU \subseteq \calU' \subseteq \bfA$ are equivalent chunks of $\bfA$, then
the inclusion $\calU^{\degree} \subseteq {\calU'}^{\degree}$ is a weak equivalence
of $\bfS$-enriched categories.
\end{remark}

\begin{example}
Let $\bfS$ be an excellent model category, and $\bfA$ a combinatorial $\bfS$-enriched model category.
Then $\bfA$ is a chunk of itself; this follows from the small object argument.
\end{example}

\begin{example}
Let $\calU \subseteq \bfA$ be a chunk, and let $\{ X_{\alpha} \}$ be a collection of objects in $\bfA$.
Let $\calV \subseteq \calU$ be the full subcategory spanned by those objects $X \in \calU$
such that there exists an isomorphism $[X] \simeq [X_{\alpha}]$ in the homotopy category
$\h{\bfA}$. Then $\calV$ is also a chunk of $\bfA$.
\end{example}

We will prove a general existence theorem for chunks below (see Lemma \ref{exchunk}).

\begin{lemma}\label{tubble}
Let $\bfS$ be an excellent model category, and let $\calC$ be a small $\bfS$-enriched
category. Then there exists a weak equivalence of $\bfS$-enriched categories
$\calC \rightarrow \calU^{\degree}$, where $\calU$ is a chunk of a combinatorial
$\bfS$-enriched category $\bfA$.
\end{lemma}

\begin{proof}
Without loss of generality, we may suppose that $\calC$ is fibrant. Let
$\bfA = \bfS^{ \calC^{op}}$, endowed with the projective model structure.
We can identify $\calC$ with a full subcategory of $\bfA^{\degree}$ via the Yoneda embedding.
Using Lemma \ref{exchunk}, we can enlarge $\calC$ to a chunk in $\bfA$ having the same
image in the homotopy category $\h{\bfA}$.
\end{proof}

\begin{notation}\index{not}{Pbfa@$P(\bfA)$}
Let $\bfS$ be an excellent model category, let $\bfA$ be a combinatorial
$\bfS$-enriched model category, and let $\calU$ be a chunk of $\bfA$. We
define a new category
$P(\calU)$ as follows:
\begin{itemize}
\item[$(i)$] The objects of $P(\calU)$ are fibrations
$\phi: A \rightarrow B \times C$ in $\bfA$, where 
$A, B, C \in \calU^{\degree}$, and the composite maps
$A \rightarrow B$ and $A \rightarrow C$ are weak equivalences.
\item[$(ii)$] Morphisms in $P(\calU)$ are given by
maps of diagrams
$$ \xymatrix{ B \ar[d] & A \ar[l] \ar[r] \ar[d] & C \ar[d] \\
B' & A' \ar[l] \ar[r] & C'. }$$
\end{itemize}

We let $\pi, \pi': P(\calU) \rightarrow \calU^{\degree}$ be the functors described by the formulas
$$ \pi( \phi: A \rightarrow B \times C) = B \quad \quad \pi'( \phi: A \rightarrow B \times C) = C.$$
We observe that both $\pi$ and $\pi'$ have the structure of $\bfS$-enriched functors.
Invoking assumption $(a)$ of Proposition \ref{defchunk}, we deduce the
existence of another $\bfS$-enriched functor $\tau: \calU^{\degree} \rightarrow P(\calU)$, which
carries an object $A \in \bfA^{\degree}$ to the map $q$ appearing in a functorial
factorization
$$ A \stackrel{p}{\rightarrow} \overline{A} \stackrel{q}{\rightarrow} A \times A,$$
of the diagonal, where $p$ is a trivial cofibration and $q$ is a fibration.
\end{notation}

\begin{theorem}\label{catta}
Let $\bfS$ be an excellent model category, let $\bfA$ be a combinatorial $\bfS$-enriched model category, and let $\calU$ be a chunk of $\bfA$. 
Then the morphisms $\pi, \pi': P(\calU) \rightarrow \calU^{\degree}$ and
$\tau: \calU^{\degree} \rightarrow P(\calU)$
furnish $P(\calU)$ with the structure of a path object for
$\calU^{\degree}$ in $\SCat$.
\end{theorem}

\begin{proof}
We first show that $\pi \times \pi'$ is a fibration of $\bfS$-enriched
categories. In view of Theorem \ref{staycode}, it will suffice to show that
$\pi \times \pi'$ is a local fibration. Let $\phi: A \rightarrow B \times C$ and 
$\phi': A' \rightarrow B' \times C'$ be objects of $P(\calU)$.
We must show that the induced map
$$ \bHom_{P(\calU)}(\phi,\phi') \rightarrow \bHom_{\bfA}(B,B') \times
\bHom_{\bfA}(C,C')$$ is a fibration in $\bfS$. 
This map is a base change of 
$$ \bHom_{\bfA}(A,A') \rightarrow \bHom_{\bfA}(A, B' \times C'),$$
which is a fibration in virtue of the assumption that $\phi'$ is a fibration (since
$A$ is assumed to be cofibrant).

To complete the proof that $\pi \times \pi'$ is a quasi-fibration,
we must show that if $\phi: A \rightarrow B \times C$ is an object of
$P( \calU)$ and we are given weak equivalences
$f: B \rightarrow B'$, $g: C \rightarrow C'$, then we can lift
$f$ and $g$ to an equivalence in $P(\calU)$. To do so, we factor the composite map $A \rightarrow B' \times C'$ as a
trivial cofibration $A \rightarrow A'$ followed by a fibration
$\phi': A' \rightarrow B' \times C'$. Since $\calU$ is a chunk of $\bfA$, we may assume that
$A' \in \calU$ so that $\phi' \in P(\calU)$. We have an evident natural transformation
$\alpha: \phi \rightarrow \phi'$. We will show below that $\pi: P(\calU) \rightarrow \calU^{\degree}$
is an equivalence of $\bfS$-enriched categories; since $\pi(\alpha) = f$ is an isomorphism
in $\h{ \calU^{\degree}}$, we conclude that $\alpha$ is an isomorphism in $\h{P(\calU)}$ as
required.

To complete the proof, we must show that $\tau$ is a weak equivalence of $\bfS$-enriched categories.
By the two-out-of-three property, it will suffice to show that $\pi$ is a weak equivalence of
$\bfS$-enriched categories. Since $\tau$ is a section of $\pi$, it is clear that $\pi$ is essentially surjective. It remains only to prove that $\pi$ is fully faithful. Let
$\phi: A \rightarrow B \times C$ and $\phi': A' \rightarrow B' \times C'$ be objects
of $P(\calU)$; we wish to show that the induced map
$p: \bHom_{P(\calU)}(\phi, \phi') \rightarrow \bHom_{\bfA}(B,B')$ is a weak equivalence in
$\bfS$. We have a commutative diagram
$$ \xymatrix{ \bHom_{P(\calU)}(\phi,\phi') \ar[r] \ar[d] & \bHom_{\bfA}(A,A') \ar[d]^{u} \\
\bHom_{\bfA}(B,B') \times \bHom_{\bfA}(C,C') \ar[d] \ar[r] &  \bHom_{\bfA}(A, B' \times C') \ar[d] \\
\bHom_{\bfA}(B,B') \times \bHom_{\bfA}(A,C') \ar[d] \ar[r] & \bHom_{\bfA}(A,B') \times
\bHom_{\bfA}(A,C') \ar[d] \\
\bHom_{\bfA}(B,B') \ar[r] & \bHom_{\bfA}(A,B'). }$$
We note that, since the map $A \rightarrow B$ is a weak equivalence between cofibrant objects,
and $B'$ is fibrant, the bottom horizontal map is a weak equivalence in $\bfS$. Consequently, to
show that the top horizontal map is a weak equivalence in $\bfS$, it will suffice to show that each square
in the diagram is homotopy Cartesian. The bottom square is Cartesian and fibrant, so there is nothing to prove. The middle square is homotopy Cartesian because both of the middle vertical maps are weak equivalences. The upper square is a pullback square between fibrant objects of $\bfS$, and the
map $u$ is a fibration; we now complete the proof by invoking Proposition \ref{leftpropsquare}.
\end{proof}

Fix an excellent model category $\bfS$. The symmetric monoidal structure on $\bfS$
induces a symmetric monoidal structure on $\SCat$:
if $\calC$ and $\calD$ are $\bfS$-enriched categories, then we can define a new
$\bfS$-enriched category $\calC \otimes \calD$ as follows:
\begin{itemize}
\item[$(i)$] The objects of $\calC \otimes \calD$ are pairs $(C,D)$, where
$C \in \calC$ and $D \in \calD$.
\item[$(ii)$] Given a pair of objects $(C,D), (C',D') \in \calC \otimes \calD$, we have
$$ \bHom_{\calC \otimes \calD}( (C,D), (C',D')) = \bHom_{\calC}(C,C') \otimes \bHom_{\calD}(D,D') \in \bfS.$$
\end{itemize}
In the case where the tensor product on $\bfS$ is the Cartesian product, this simply
reduces to the usual product of $\bfS$-enriched categories.

Note that the operation $\otimes: \SCat \times \SCat \rightarrow \SCat$ is {\em not}
a left Quilen bifunctor, even when $\bfS = \sSet$: for example, a product of cofibrant
simplicial categories is generally not cofibrant. Nevertheless, $\otimes$ behaves
much like a left Quillen bifunctor at the level of homotopy categories. 
For example, the operation $\otimes$ respects weak equivalences in each argument, and therefore
induces a functor $\otimes: \h{ \SCat} \times \h{\SCat} \rightarrow \h{\SCat}$, which is
characterized by the existence of natural isomorphisms
$[ \calC \otimes \calD] \simeq [\calC] \otimes [\calD]$.

Our goal for the remainder of this section is to show that the monoidal structure
$\otimes$ on $\SCat$ is {\em closed}: that is, there exist internal mapping objects in
$\h{\SCat}$. This is not completely obvious. It is easy to see that the monoidal structure
$\otimes$ on $\SCat$ is closed: given a pair of $\bfS$-enriched categories $\calC$ and
$\calD$, the category of $\bfS$-enriched functors $\calD^{\calC}$ is itself $\bfS$-enriched,
and possesses the appropriate universal property. However, this is not necessarily
the ``correct'' mapping object, in the sense that the homotopy equivalence class
$[ \calD^{\calC} ]$ does not necessarily coincide with the internal mapping object
$[\calD]^{[ \calC] }$ in $\h{\SCat}$. Roughly speaking, the problem is that
$\calD^{\calC}$ consists of functors which are strictly compatible with composition; 
the correct mapping object should incorporate also functors which preserve composition only
up to (coherent) weak equivalence. However, when $\calD$ is the category of fibrant-cofibrant
objects of a $\bfS$-enriched {\em model} category $\bfA$, then we can proceed more directly.

\begin{definition}\label{cattusi}
Let $\bfS$ be an excellent model category, $\bfA$ a combinatorial $\bfS$-enriched model category,
and $\calC$ a cofibrant $\bfS$-enriched category. We will say that a full subcategory
$\calU \subseteq \bfA$ is a {\it $\calC$-chunk of $\bfA$} if it is a chunk of
$\bfA$ and the subcategory $\calU^{\calC}$ is a chunk of $\bfA^{\calC}$. 
Here we regard $\bfA^{\calC}$ as endowed with the {\em projective}
model structure.
\end{definition}


\begin{lemma}\label{tuff}
Let $\bfS$ be an excellent model category, $\bfA$ a combinatorial $\bfS$-enriched model category,
$\calC$ a $($small$)$ cofibrant $\bfS$-enriched category, and $\calU \subseteq \bfA$ a $\calC$-chunk.
Let $f,f': \calC \rightarrow \calU^{\degree}$ be a pair of maps. The following conditions are equivalent:
\begin{itemize}
\item[$(1)$] The homotopy classes $[f]$ and $[f']$ coincide in 
$\Hom_{ \h{\SCat}}( \calC, \calU^{\degree})$.
\item[$(2)$] The maps $f$ and $f'$ are weakly equivalent when regarded as objects
of $\bfA^{\calC}$.
\end{itemize}
\end{lemma}

\begin{proof}

Suppose first that $(1)$ is satisfied. Using Theorem \ref{catta}, we deduce
the existence of a homotopy $h: \calC \rightarrow P(\calU)$ from
$f = \pi \circ h$ to $f' = \pi' \circ h$. The map $h$
determines another simplicial functor $f'': \calC \rightarrow
\calU$, equipped with weak equivalences $f'' \rightarrow f$, $f''
\rightarrow f'$. This proves that $f$ and $f'$ are isomorphic in the homotopy category
of $\bfA^{\calC}$, so that $(2)$ is satisfied.

Now suppose that $(2)$ is satisfied. Since $\calU$ is a $\calC$-chunk, we can find
a projectively cofibrant $f'': \calU \rightarrow \calC$ equipped with a weak equivalence
$\alpha: f'' \rightarrow f$. Using $(2)$, we deduce that there is also a weak equivalence
$\beta: f'' \rightarrow f'$. Using again the assumption that $\calU^{\calC}$ is a chunk of
$\bfA^{\calC}$, we can choose a factorization of $\alpha \times \beta$ as a composition
$$ f'' \stackrel{u}{\rightarrow} f''' \stackrel{v}{\rightarrow} f \times f'$$ 
where $u$ is a trivial projective cofibration, $v$ is a projective fibration, and
$f''' \in \calU^{\calC}$. The map $v$ can be viewed as an object of $\calP(\calU)$, which
determines a right homotopy from $f$ to $f'$. 
\end{proof}

\begin{corollary}\label{suff}
Let $\bfS$ be an excellent model category, and let $f: \calC \rightarrow \calC'$ be a $\bfS$-enriched functor. Suppose that $f$ is fully faithful in the sense that for every pair of objects
$X,Y \in \calC$, the induced map $\bHom_{\calC}( X,Y) \rightarrow \bHom_{\calC'}(fX, fY)$ is a
weak equivalence in $\bfS$. Let $\calD$ be an arbitrary $\bfS$-enriched category. Then:
\begin{itemize}
\item[$(1)$] Composition with $f$ induces an injective map
$\phi: \Hom_{ \h{\SCat}} ( \calD, \calC) \rightarrow \Hom_{ \h{\SCat}}( \calD, \calC' )$.
\item[$(2)$] The image of $\phi$ consists of those maps
$g: \calD \rightarrow \calC'$ in $\h{\SCat}$ such that the essential image of
$[g]$ in $\h{\calC'}$ is contained in the essential image of $[f]$ in $\h{\calC'}$. 
\end{itemize}
\end{corollary}

\begin{proof}
Using Lemma \ref{tubble}, we may assume without loss of generality that
$\calC' = \calU^{\degree}$, where $\calU$ is a chunk of a $\bfA$-enriched model category.
Let $\calV \subseteq \calU$ be the full subcategory spanned by those objects which
are weakly equivalent to an object lying in the image of $f$. Since $f$ is fully faithful,
the induced map $\calC \rightarrow \calV^{\degree}$ is a weak equivalence. We may therefore
assume that $\calC = \calV^{\degree}$. 

Without loss of generality, we may suppose that $\calD$ is cofibrant. 
Enlarging $\calU$ and $\calV$ if necessary (using Lemma \ref{exchunk}), we may assume
that $\calU$ and $\calV$ are $\calD$-chunks. 
The desired results now follow immediately from Lemma \ref{tuff}.
\end{proof} 



Let $\pi_0 \bfA^{\calC}$ denote the collection of
weak equivalence classes of objects in $\bfA^{\calC}$. Every
equivalence class contains a fibrant-cofibrant representative, which
determines a $\bfS$-enriched functor $\calC \rightarrow \bfA^{\degree}$. 

\begin{proposition}\label{gumbaa}
Let $\bfS$ be an excellent model category, $\bfA$ a combinatorial
$\bfS$-enriched category, and $\calC$ a $($small$)$ cofibrant $\bfS$-enriched category.
Then the map
$$ \phi: \pi_0 \bfA^{\calC} \rightarrow \Hom_{\h{\SCat}}(\calC, \bfA^{\degree})$$
described above is bijective.
\end{proposition}

\begin{proof}
In view of Proposition \ref{lesstrick}, we may assume that $\calC$ is cofibrant.
Lemma \ref{tuff} shows that $\phi$ is well-defined and injective.
We show that $\phi$ is surjective. Let $[f] \in \Hom_{\h{\SCat}}(
\calC,\calU^{\degree})$. Since $\calC$ is cofibrant and $\bfA^{\degree}$
is fibrant in $\SCat$, we can find a $\bfS$-enriched functor $f: \calC
\rightarrow \bfA^{\degree}$ representing $[f]$. The simplicial
functor $f$ takes values in fibrant-cofibrant objects of $\bfA$,
but is not necessarily fibrant-cofibrant {\em as an object of}
$\bfA^{\calC}$. However, we can choose a 
weak trivial fibration $f' \rightarrow f$, where $f'$ is projectively cofibrant.
Consequently, it will suffice to show that a weak equivalence $u: f' \rightarrow f$
of $\bfS$-enriched functors $\calC \rightarrow \bfA^{\degree}$
guarantees that $[f] = [f'] \in \Hom_{\h{\SCat}}(\calC, \bfA^{\degree})$, which follows
from Lemma \ref{tuff}.
\end{proof}

\begin{proposition}\label{gumbarr}
Let $\bfS$ be an excellent model category, $\bfA$ a combinatorial $\bfS$-enriched model category,
and $\calC$ a small $\bfS$-enriched category.
Then the evaluation map
$e: ( \bfA^{\calC})^{\degree} \otimes \calC \rightarrow \bfA^{\degree}$ has
the following property: for every small $\bfS$-enriched category $\calD$, 
composition with $e$ induces a bijection
$$ \Hom_{ \h{\SCat} }( \calD, (\bfA^{\calC})^{\degree} )
\rightarrow \Hom_{ \h{\SCat} }( \calC \otimes \calD, \bfA^{\degree} ).$$
\end{proposition}

\begin{proof}
Using Proposition \ref{gumbaa}, we can identify both sides with
$\pi_0 \bfA^{ \calD \otimes \calC}$. 
\end{proof}

It is not clear that the conclusion of Proposition \ref{gumbarr} characterizes
$(\bfA^{\calC})^{\degree}$ up to equivalence, since
$(\bfA^{\calC})^{\degree}$ is a {\em large} $\bfS$-enriched category, and the
proof of the Proposition only applies in the case where $\calD$ is small.
To remedy this defect, we establish a more refined version:

\begin{corollary}\label{sniffle}
Let $\bfS$ be an excellent model category, $\bfA$ a combinatorial $\bfS$-enriched model category,
and $\calC$ a small cofibrant $\bfS$-enriched category. Let $\calU$ be a
$\calC$-chunk of $\bfA$. 
Then the evaluation map
$e: ( \calU^{\calC})^{\degree} \otimes \calC \rightarrow \calU^{\degree}$ has
the following property: for every small $\bfS$-enriched category $\calD$, 
composition with $e$ induces a bijection
$$ \Hom_{ \h{\SCat} }( \calD, (\calU^{\calC})^{\degree} )
\rightarrow \Hom_{ \h{\SCat} }( \calC \otimes \calD, \calU^{\degree} ).$$
\end{corollary}

\begin{proof}
Combine Proposition \ref{gumbarr} with Corollary \ref{suff}.
\end{proof}

We conclude this section with a technical result, which ensures the existence of a good supply of chunks of combinatorial model categories. 

\begin{lemma}\label{exchunk}
Let $\bfS$ be an excellent model category, $\bfA$ a combinatorial $\bfS$-enriched model category,
and $\{ \calC_\alpha \}_{ \alpha \in A}$ a $($small$)$ collection of $($small$)$ cofibrant
$\bfS$-enriched categories. Let $\calU$ be a small full subcategory of $\bfA$.
Then there exists a small subcategory $\calV \subseteq \bfA$ containing $\calU$, such
that $\calV$ is a $\calC_{\alpha}$-chunk for each $\alpha \in A$. Moreover, we may arrange
that $\calU$ and $\calV$ have the same essential image in the homotopy category
$\h{\bfA}$. 
\end{lemma}

\begin{proof}
Enlarging $A$ if necessary, we may suppose that the collection $\{ \calC_{\alpha} \}_{\alpha in A}$
includes the unit category $[0]_{\bfS}$. 
For each $\alpha \in A$, we can choose
$\bfS$-enriched functors $$F_{\alpha}: \bfA^{ \calC_{\alpha} \otimes [1]_{\bfS} } \rightarrow
\bfA^{ \calC_{\alpha} \otimes [2]_{\bfS} } \quad \quad
G_{\alpha}: \bfA^{ \calC_{\alpha} \otimes [1]_{\bfS} } \rightarrow
\bfA^{ \calC_{\alpha} \otimes [2]_{\bfS} },$$
such that $F$ carries each morphism $u: f \rightarrow g$ in $\bfA^{\calC}$ to a factorization
$$ f \stackrel{u'}{\rightarrow} f' \stackrel{u''}{\rightarrow} g$$
where $u'$ is a strong trivial cofibration and $u''$ is a projective fibration, and
$G$ carries $u$ to a factorization
$$ f \stackrel{v'}{\rightarrow} g' \stackrel{v''}{\rightarrow} g$$
where $v'$ is a projective cofibration and $v''$ is a weak trivial cofibration.
For $C \in \calC_{\alpha}$, let $F^{C}_{\alpha}$ be the
functor $u \mapsto f'(C)$, and define $G^{C}_{\alpha}$ likewise.

Choose a regular cardinal $\kappa$ such that each $\calC_{\alpha}$ is 
$\kappa$-small. We define a sequence of full subcategories
$\{ \calU_{\alpha} \subseteq \bfA \}_{ \alpha < \kappa }$ as follows:
\begin{itemize}
\item[$(i)$] If $\alpha = 0$, then $\calU_{\alpha} = \calU$.
\item[$(ii)$] If $\alpha$ is a nonzero limit ordinal, then 
$\calU_{\alpha} = \bigcup_{\beta < \alpha} \calU_{\beta}$.
\item[$(iii)$] If $\alpha = \beta + 1$, then
$\calU_{\alpha}$ is the full subcategory of $\bfA$ spanned by the following:
\begin{itemize}
\item[$(a)$] The objects which belong to $\calU_{\beta}$.
\item[$(b)$] The objects $F^{C}_{\alpha}(u) \in \bfA$, where
$\alpha \in A$, $C \in \calC_{\alpha}$, and $u: f \rightarrow g$
is a morphism from an object of $\calU_{\beta}^{\calC_{\alpha}}$ to
a finite product of object in $\calU_{\beta}^{\calC_{\alpha}}$. 
\item[$(c)$] The objects $G^{C}_{\alpha}(u) \in \bfA$, where
$\alpha \in A$, $C \in \calC_{\alpha}$, and $u: f \rightarrow g$
is a morphism from a finite coproduct of objects of $\calU_{\beta}^{\calC_{\alpha}}$ to
an object in $\calU_{\beta}^{\calC_{\alpha}}$. 
\end{itemize}
\end{itemize}

It is readily verified that the subcategory $\calV = \bigcup_{\alpha < \kappa} \calU_{\alpha}$
has the desired properties.
\end{proof}

\subsection{Homotopy Colimits of $\bfS$-Enriched Categories}\label{hoco}

Our goal in this section is to give an explicit construction of (certain) homotopy colimits
in the model category $\SCat$, where $\bfS$ is an excellent model category. We begin
with some general remarks concerning localization.

\begin{notation}\index{not}{calCWinv@$\calC[W^{-1}]$}\index{not}{Inv@$\Inv$}\label{localdef}
Consider the canonical map $\overline{i}: [1]_{\bfS} \rightarrow [1]^{\sim}_{\bfS}$. We
fix once and for all a factorization of $\overline{i}$ as a composition
$$ [1]_{\bfS} \stackrel{i}{\rightarrow} \calE \stackrel{i'}{\rightarrow} [1]_{\bfS}^{\sim},$$
where $i$ is a cofibration and $i'$ is a weak equivalence of $\bfS$-enriched categories.
For every $\bfS$-enriched category $\calC$ and every map of sets
$W \rightarrow \Hom_{\SCat}( [1]_{\bfS}, \calC$, we define a new $\bfS$-enriched
category $\calC[W^{-1}]$ by a pushout diagram
$$ \xymatrix{ \coprod_{w \in W} [1]_{\bfS} \ar[r] \ar[d] & \calC \ar[d] \\
\coprod_{w \in W} \calE \ar[r] & \calC[W^{-1} ]. }$$
\end{notation}

\begin{remark}\label{summat}
Since the model category $\SCat$ is left proper, the construction
$\calC \mapsto \calC[W^{-1}]$ preserves weak equivalences in $\calC$.
\end{remark}

We now characterize $\calC[W^{-1}]$ by a universal property in $\h{\SCat}$.

\begin{lemma}\label{pufft}
Let $\calC$ be a fibrant $\bfS$-enriched category, and 
$f$ be a morphism in $\calC$ classified by a map $j_0: [1]_{\bfS} \rightarrow \calC$.
The following conditions are equivalent:
\begin{itemize}
\item[$(1)$] The map $f$ is an equivalence in $\calC$.
\item[$(2)$] The extension problem depicted in the diagram
$$ \xymatrix{ [1]_{\bfS} \ar[d]^{i} \ar[r]^{j_0} & \calC \\
\calE \ar@{-->}[ur]^{j} & }$$
admits a solution.
\end{itemize}
\end{lemma}

\begin{proof}
The implication $(2) \Rightarrow (1)$ is clear, since every morphism in
$\calE$ is an equivalence. For the converse, we observe that
the desired lifting problem admits a solution if and only if the induced map
$i': \calC \rightarrow \calC \coprod_{[1]_{\bfS} } \calE$ admits a left inverse.
Since $\calC$ is fibrant, it suffices to show that $i'$ is a trivial cofibration.
The map $i'$ is a cofibration since it is a pushout of $i$, and a weak equivalence
because of the invertibility hypothesis. 
\end{proof}

Lemma \ref{pufft} immediately implies the following apparently stronger claim:

\begin{lemma}\label{canner}
Let $f_0: \calC \rightarrow \calD$ be a $\bfS$-enriched functor, where
$\calD$ is a fibrant $\bfS$-enriched category. Let $\psi: W \rightarrow \Hom_{\SCat}([1]_{\bfS},\calC)$ be a map of sets. The following conditions are equivalent:
\begin{itemize}
\item[$(1)$] The map $f_0$ extends to a map $f: \calC[W^{-1}] \rightarrow \calD$.
\item[$(2)$] For each $w \in W$, $f_0$ carries the morphism $\phi(w)$ to an equivalence in $\calD$.
\end{itemize}
\end{lemma}

\begin{proposition}\label{postcan}
Let $\calC$ and $\calD$ be $\bfS$-enriched categories, and let
$\psi: W \rightarrow \Hom_{\SCat}( [1]_{\bfS}, \calC)$ be a map of sets.
Then the induced map
$$ \phi: \Hom_{\h{\SCat}}( \calC[W^{-1}], \calD) \rightarrow \Hom_{\h{\SCat}}(\calC, \calD)$$
is injective, and its image is the subset $\Hom^{W}_{\h{\SCat}}(\calC, \calD)
\subseteq \Hom_{\h{\SCat}}(\calC, \calD)$ consisting of those homotopy classes of maps
which induce functors $\h{\calC} \rightarrow \h{\calD}$ carrying each element of $W$ to an isomorphism in $\h{\calD}$. 
\end{proposition}

\begin{proof}
Without loss of generality, we may suppose that $\calC$ is cofibrant and $\calD$ is fibrant.
The description of the image of $\phi$ follows immediately from Lemma \ref{canner}.
It will therefore suffice to show that $\phi$ is injective. Suppose we are given a pair of maps
$[f], [g] \in \Hom_{ \h{\SCat}}( \calC[W^{-1}], \calD)$ such that $\phi( [f] ) = \phi( [g] )$. Since
$\calC[ W^{-1}]$ is cofibrant, we may assume that $[f]$ and $[g]$ are represented by
actual $\bfS$-enriched functors $f,g: \calC[W^{-1}] \rightarrow \calD$. Moreover, the
condition that $\phi( [f] ) = \phi( [g] )$ guarantees that the restrictions
$f|\calC$ and $g|\calC$ are homotopic. We wish to show that $f$ and $g$ are homotopic.

Invoking Proposition \ref{princex}, we deduce that $g$ is homotopic to a map
$g': \calC[W^{-1}] \rightarrow \calD$ such that $g' | \calC = f | \calC$. Replacing $g$ by
$g'$ if necessary, we may assume that $g| \calC = f| \calC$. It will now suffice to show that
$f$ and $g$ are homotopic in the model category $( \SCat)_{\calC/}$. We
observe that $f$ and $g$ determine a map
$$ h: \calC[ ( W \coprod W)^{-1} ] \simeq \calC[W^{-1}] \coprod_{\calC} \calC[W^{-1}] \rightarrow \calD.$$
Using the invertibility hypothesis, we conclude that
$\calC[ (W \coprod W)^{-1}]$ is a cylinder object for
$\calC[W^{-1}]$ in the model category $( \SCat)_{\calC/}$, so that $h$ is the desired
homotopy from $f$ to $g$. 
\end{proof}

\begin{lemma}\label{swimcase}
Let $f: \calC \rightarrow \calD$ be a $\bfS$-enriched functor, and let
$\calM$ be the {\it categorical mapping cylinder} of $f$, defined as follows:
\begin{itemize}
\item[$(1)$] An object of $\calM$ is either an object of $\calC$ or an object of $\calD$.
\item[$(2)$] Given a pair of objects $X,Y \in \calM$, we have
$$ \bHom_{\calM}(X,Y) = \begin{cases} \bHom_{\calC}(X,Y) & \text{if } X,Y \in \calC \\
\bHom_{\calD}( X,Y) & \text{if } X,Y \in \calD \\
\bHom_{\calD}( fX, Y) & \text{if } X \in \calC, Y \in \calD \\
\emptyset & \text{if } X \in \calD, Y \in \calC. \end{cases} $$
Here $\emptyset$ denotes the initial object of $\bfS$.
\end{itemize}
We observe that there is a canonical retraction $j$ of $\calM$ onto $\calD$, described by the formula $$j(X) = \begin{cases} fX & \text{if } X \in \calC \\
X & \text{if } X \in \calD. \end{cases}$$
Let $W$ be a collection of morphisms in $\calM$ with the following properties:
\begin{itemize}
\item[$(i)$] For each $w \in W$, $j(w)$ is an identity morphism in $\calD$.
\item[$(ii)$] For every object $C \in \calC$, the morphism $C \rightarrow fC$
in $\calM$ classifying the identity map from $fC$ to itself belongs to $W$.
\end{itemize}
Assumption $(i)$ implies that the map $j$ canonically extends to a map
$\overline{j}: \calM[W^{-1}] \rightarrow \calD$. The map $\overline{j}$
is a weak equivalence of $\bfS$-enriched categories.
\end{lemma}

\begin{proof}
It will suffice to show that composition with $\overline{j}$ induces a bijection
$$ \Hom_{ \h{\SCat} }( \calD, \calA) \rightarrow \Hom_{ \h{\SCat}}( \calM[W^{-1}], \calA)$$
for every $\bfS$-enriched category $\calA$. Equivalently, we must show that the map $$t: \Hom_{\h{\SCat}}( \calD, \calA) \rightarrow \Hom^{W}_{\h{\SCat}}(\calM, \calA)$$
is bijective, where $\Hom_{\h{\SCat}}^{W}( \calM, \calA)$ is defined as in Proposition \ref{postcan}.
The map $t$ has a section $t'$, given by composition with the inclusion $\calD \rightarrow \calM$.
It will therefore suffice to show that $t \circ t'$ is the identity on $\Hom_{\h{\SCat}}^{W}(\calM, \calA)$. 

Using Lemma \ref{tubble} and Corollary \ref{suff}, we can reduce to the case where $\calA = \bfA^{\degree}$, where $\bfA$ is a combinatorial $\bfS$-enriched model category. Invoking Proposition \ref{gumbaa}, we deduce that every element $[g] \in \Hom_{\h{\SCat}}( \calM, \calA)$ can be represented
by a diagram $g: \calM \rightarrow \bfA^{\degree}$. We wish to prove that
$g$ and $g \circ i \circ j$ are homotopic. We observe that there is a canonical natural transformation $\alpha: g \rightarrow g \circ i \circ j$. Moreover, if $g$ carries
each element of $W$ to an equivalence in $\bfA^{\degree}$, then assumption
$(ii)$ guarantees that $\alpha$ is a weak equivalence in the model category $\bfA^{\calM}$. We now invoke Proposition \ref{gumbaa} to deduce that $g$ and $g \circ i \circ j$ are homotopic as desired.
\end{proof}

\begin{definition}
Let $A$ be a partially ordered set. An {\it $A$-filtered $\bfS$-enriched category}
is a $\bfS$-enriched category $\calC$ together with a map $r: \Ob(\calC) \rightarrow A$
with the following property: if $C, D \in \calC$ and $r(C) \nleq r(D)$, then
$\bHom_{\calC}(C,D) \simeq \emptyset$, where $\emptyset$ denotes an initial object of $\bfS$.

If $\calC$ is an $A$-filtered $\bfS$-enriched category and $a \in A$, then we let
$\calC_{\leq a}$ denote the full subcategory of $\calC$ spanned by those objects
$C \in \calC$ such that $r(C) \leq a$.
\end{definition}

\begin{remark}\label{sabreton}
Let $\calC$ be an $A$-filtered $\bfS$-enriched category, and let
$\psi: W \rightarrow \Hom_{\SCat}( [1]_{\bfS}, \calC)$ be a map of sets.
For each $a \in A$, we let $W_{a} \subseteq W$ be the subset
consisting of those elements $w \in W$ such that the morphism
$\psi(w)$ belongs to $\calC_{a}$. This data determines a diagram
$\chi^{W}: A \rightarrow \SCat$, described by the formula
$a \mapsto \calC_{\leq a}[ W_{a}^{-1} ]$. Moreover, we have a canonical isomorphism
of $\bfS$-enriched categories $\calC[W^{-1}] \simeq \colim(\chi)$.
\end{remark}

Using the small object argument, we easily deduce the following result:

\begin{lemma}\label{tooter}
Let $A$ be a partially ordered set, and let $\calC$ be an $A$-filtered
$\bfS$-enriched category. Then there exists a $\bfS$-enriched functor
$f: \calC' \rightarrow \calC$ with the following properties:
\begin{itemize}
\item[$(1)$] The functor $f$ is bijective on objects, and for every
pair of objects $C,D \in \calC'$, the map $\bHom_{\calC}(C,D)
\rightarrow \bHom_{\calC}(fC, fD)$ is a trivial fibration in $\bfS$. In particular,
$f$ is a weak equivalence of $\bfS$-enriched categories.

\item[$(2)$] The $A$-filtration on $\calC$ induces an $A$-filtration on $\calC'$. In other
words, if $C, D \in \calC'$ and $r(fC) \nleq r(fD)$, then $\bHom_{\calC'}(C,D)$ is
an initial object of $\bfS$.

\item[$(3)$] The diagram $A \rightarrow \SCat$ described by the formula
$a \mapsto \calC'_{\leq a}$ is projectively cofibrant.
\end{itemize}
\end{lemma}

\begin{proposition}\label{scun}
Let $A$ be a partially ordered set, $\calC$ be an $A$-filtered
$\bfS$-enriched category, and $\psi: W \rightarrow \Hom_{ \SCat}( [1]_{\bfS}, \calC)$
a map of sets. Let $\chi: A \rightarrow \SCat$ be defined as in Remark \ref{sabreton}. Then the
isomorphism $\colim \chi \simeq \calC[W^{-1}]$ exhibits $\calC$ as the homotopy colimit
of the diagram $\chi$.
\end{proposition}

\begin{proof}
Choose a map $\calC' \rightarrow \calC$ as in Lemma \ref{tooter}, and
choose a map $\psi': W \rightarrow \Hom_{\SCat}( [1]_{\bfS}, \calC' )$ lifting
$\psi$, and let $\chi': A \rightarrow \SCat$ be defined as in Remark \ref{sabreton}.
Then we have a canonical map $\chi' \rightarrow \chi$, which is a cofibrant
replacement for $\chi$ in the model category $\Fun( A, \SCat)$. It will therefore suffice to show that the induced map
$\calC'[W^{-1}] \simeq \colim \chi' \rightarrow \colim \chi \simeq \calC[W^{-1}]$ is a weak equivalence of
$\bfS$-enriched categories, which follows immediately from Remark \ref{summat}.
\end{proof}

\begin{definition}\index{gen}{Grothendieck construction}\index{not}{Groth@$\Groth(p)$}\label{swype}
\index{not}{@$W(p)$}
Let $A$ be a partially ordered set, and let $p: A \rightarrow \SCat$ be an $A$-indexed diagram of $\bfS$-enriched categories. Let us denote the image of $a \in A$ under $p$ by $\calC_{a}$.

The {\it Grothendieck construction on $p$} is a category
$\Groth(p)$ defined as follows:
\begin{itemize}
\item[$(1)$] The objects of $\Groth(p)$ are pairs $(a, C)$, where $a \in A$ and
$C \in \calC_{a}$.
\item[$(2)$] Given a pair of objects $(a,C), (a', C')$ in $\Groth(p)$, we set
$$ \bHom_{ \Groth(p) }( (a,C), (a',C') ) = \begin{cases} \bHom_{\calC_{a'}}( p^{a'}_{a} C,
C') & \text{if } a \leq a' \\
\emptyset & \text{otherwise.} \end{cases}$$
Here $p^{a'}_{a}$ denotes the functor $\calC_{a} \rightarrow \calC_{a'}$ determined by
$p$, and $\emptyset$ denotes an initial object of $\bfS$.
\item[$(3)$] Composition in $\Groth(p)$ is defined in the obvious way.
\end{itemize}

We observe that $\Groth(p)$ is $A$-filtered via the map
$r: \Ob( \Groth(p) ) \rightarrow A$ given by the formula
$r(a,C) = a$. We let $W(p)$ denote the collection of all morphisms in
$\Groth(p)$ of the form $\alpha: (a, C) \rightarrow (a', p^{a'}_{a} C)$, where
$a \leq a'$ and $\alpha$ corresponds to the identity in $\calC_{a'}$. 

For each $a \in A$, there is a canonical functor
$\pi_{a}: \Groth(p)_{\leq a} \rightarrow \calC_{a}$, given by the formula
$(C, a') \mapsto p^{a}_{a'}(C)$. We note that $\pi$ carries
each element of $W(p)_{a}$ to an identity map in $\calC_{a}$, so that
$\pi_{a}$ canonically extends to a map $\overline{\pi}_{a}: \Groth(p)_{\leq a}[ W(p)_{a}^{-1}] \rightarrow \calC_{a}$. The maps $\overline{\pi}_{a}$ are functorial in
$a$, and therefore determine a map of diagrams $\chi(p) \rightarrow p$, where
$\chi(p)$ is defined as in Remark \ref{sabreton}.
\end{definition}

\begin{lemma}\label{cutta}
Let $p: A \rightarrow \SCat$ be as in Definition \ref{swype}. Then for each
$a \in A$, the map $\overline{\pi}_{a}: \Groth(p)_{\leq a}[ W(p)_{a}^{-1}] \rightarrow \calC_{a}$
is a weak equivalence of $\bfS$-enriched categories.
\end{lemma}

\begin{proof}
This is a special case of Lemma \ref{swimcase}. 
\end{proof}

\begin{lemma}\label{toopa}
Let $p: A \rightarrow \SCat$ be as in Definition \ref{swype}. Then
there is a canonical isomorphism $\Groth(p)[ W(p)^{-1} ] \simeq \hocolim(p)$
in the homotopy category $\h{\SCat}$.
\end{lemma}

\begin{proof}
Combine Lemma \ref{cutta} with Proposition \ref{scun}.
\end{proof}

\begin{lemma}\label{twoface}
Let $\calC$ and $\calD$ be small $\bfS$-enriched categories.
Let $W$ be a collection of morphisms in $\calC$, and let
$W'$ be the collection of all morphisms in $\calC \otimes \calD$ of the form
$w \otimes \id_{D}$, where $w \in W$ and $D \in \calD$. Then the canonical map
$$ ( \calC \otimes \calD )[ {W'}^{-1} ] \rightarrow \calC[W^{-1}] \otimes \calD$$
is a weak equivalence of $\bfS$-enriched categories.
\end{lemma}

\begin{proof}
It will suffice to show that for every $\bfS$-enriched category $\calA$, 
the induced map
$$ \phi: \Hom_{ \h{\SCat} }( \calC[W^{-1}] \otimes \calD, \calA) \rightarrow
\Hom_{ \h{\SCat} }( ( \calC \otimes \calD)[ {W'}^{-1} ], \calA) $$
is bijective. Using Lemma \ref{tubble} and Corollary \ref{suff}, we can reduce
to the case where $\calA = \bfA^{\degree}$, where $\bfA$ is a combinatorial $\bfS$-enriched
model category. We now invoke Propositions \ref{gumbarr} and \ref{postcan} to get a chain of bijections
\begin{eqnarray*}
\Hom_{ \h{\SCat} }( \calC[W^{-1}] \otimes \calD, \bfA^{\degree}) & \simeq &
\Hom_{\h{\SCat} }( \calC[W^{-1}], ( \bfA^{\calD} )^{\degree} ) \\
& \simeq & \Hom_{\h{\SCat}}^{W}( \calC, ( \bfA^{\calD} )^{\degree} ) \\
& \simeq & \Hom_{\h{\SCat}}^{W'}( \calC \otimes \calD, \bfA^{\degree} ) \\
& \simeq & \Hom_{\h{\SCat}}( (\calC \otimes \calD)[W'], \bfA^{\degree} )
\end{eqnarray*}
whose composition is the map $\phi$.
\end{proof}

\begin{theorem}\label{tubba}
Let $A$ be a partially ordered set, and let $\calD$ be a $\bfS$-enriched category.
Then the functor $\calC \mapsto \calC \otimes \calD$ commutes with $A$-indexed homotopy colimits. In other words, if $p: A \rightarrow \SCat$ is a projectively cofibrant diagram,
and $p': A \rightarrow \SCat$ is defined by $p'(a) = p(a) \otimes \calD$, then the canonical isomorphism $\colim(p') \simeq \colim(p) \otimes \calD$ exhibits $\colim(p) \otimes \calD$
as a homotopy colimit of the diagram $p'$.
\end{theorem}

\begin{proof}
In view of Lemma \ref{toopa}, it will suffice to show that the canonical map
$h: \Groth(p')[ W(p')^{-1} ] \rightarrow \Groth(p)[W(p)^{-1}] \otimes \calD$
is a weak equivalence of $\bfS$-enriched categories. This is a special case of
Lemma \ref{twoface}.
\end{proof}

\subsection{Exponentiation in Model Categories}\label{camper}

Let $\calC$ be a category which admits finite products, containing a pair of objects $X$ and $Y$.
An {\it exponential of $X$ by $Y$} is an object $X^Y \in \calC$ together with a map
$e: X^Y \times Y \rightarrow X$, with the following universal property: for every
object $W \in \calC$, the composition
$$ \Hom_{\calC}( W, X^Y) \rightarrow \Hom_{\calC}( W \times Y, X^Y \times Y)
\stackrel{\circ e}{\rightarrow} \Hom_{\calC}(W \times Y, Z)$$
is bijective.\index{gen}{exponential}

Our goal in this section is to study the existence of exponentials in the homotopy category
of a model category $\bfA$. Suppose given a pair of objects $X, Y \in \bfA$, such that there
exists an exponential of $[X]$ by $[Y]$ in the homotopy category $\h{\bfA}$. 
We can then represent this exponential as $[Z]$, for some object $Z \in \bfA$. 
Without loss of generality, we may assume that $X$, $Y$, and $Z$ are fibrant and cofibrant,
so that we have a canonical identification $[Z] \times [Y] \simeq [Z \times Y]$. However,
we encounter a technical difficulty: the evaluation map $[Z] \times [Y] \rightarrow [X]$ need
not be representable by any morphism from $Z \times Y$ to $X$ in the category $\bfA$, because
$Z \times Y$ need not be cofibrant. We wish to work in certain contexts where this difficulty
does arise (for example, where $\bfA$ is the category of simplicial categories). For this reason we are forced to work with the following somewhat cumbersome definition:

\begin{definition}\index{gen}{weak exponential}\index{gen}{exponential!weak}\label{ki}
Let $\bfA$ be a model category. We will say that a diagram
$$ \xymatrix{ &  P \ar[dl]^{p} \ar[dr] & \\
Z \times Y & & X}$$
{\it exhibits $Z$ as a weak exponential of $X$ by $Y$} if the following conditions are satisfied:
\begin{itemize}
\item[$(1)$] The map $p$ exhibits $P$ as a homotopy product of $Z$ and $Y$; in other words,
the induced map $[p]: [P] \rightarrow [Z] \times [Y]$ is an isomorphism in the homotopy category
$\h{\bfA}$. 
\item[$(2)$] The composition
$[Z] \times [Y] \stackrel{ [p]^{-1} }{\rightarrow} [P] \rightarrow [X]$
exhibits $[Z]$ as an exponential of $[X]$ by $[Y]$ in the homotopy category $\h{\bfA}$.
\end{itemize}

We will say that a map $Z \times Y \rightarrow X$ {\it exhibits $Z$ as an exponential of
$X$ by $Y$} if the diagram
$$ \xymatrix{ & Z \times Y \ar[dl]^{\id} \ar[dr] & \\
Z \times Y & & X }$$
satisfies $(1)$ and $(2)$.
\end{definition}

\begin{remark}\index{gen}{standard diagram}
Suppose given a diagram
$$ \xymatrix{ & P \ar[dl]^{p} \ar[dr] &  \\
Z \times Y & & X}$$
as in Definition \ref{ki}. We will say that this diagram is {\em standard} if
$X, Y, Z \in \bfA$ are fibrant, and the map $p$ is a trivial fibration.

Suppose $X$ and $Y$ are fibrant objects of $\bfA$ such that there exists an exponential
of $[X]$ by $[Y]$ in the homotopy category $\h{\bfA}$. Without loss of generality, this
exponential has the form $[Z]$ where $Z$ is a fibrant object of $\bfA$. We can then
choose a trivial fibration $P \rightarrow Z \times Y$, where $P$ is cofibrant. The evaluation
map $[Z \times Y] \simeq [Z] \times [Y] \rightarrow [X]$ is then representable by
a map $P \rightarrow X$ in $\bfA$, so that we obtain a {\em standard} diagram which exhibits
$Z$ as a weak exponential of $X$ by $Y$.
\end{remark}

\begin{remark}
Suppose given a diagram
$$ \xymatrix{ & P \ar[dl]^{p} \ar[dr] & \\
Z \times Y & & X}$$
in a model category $\bfA$. The condition that this diagram exhibits
$Z$ as a weak exponential of $X$ by $Y$ depends only on the image of this diagram
in the homotopy category $\h{\bfA}$. We are therefore replace the above diagram by a weakly equivalent diagram when testing whether or not the conditions of Definition \ref{ki} are satisfied.
\end{remark}

\begin{remark}\label{toofus}
Let $\Adjoint{F}{\bfA}{\bfB}{G}$ be a Quillen equivalence of model categories.
Suppose given a standard diagram
$$ \xymatrix{&  P \ar[dl]^{p} \ar[dr] &  \\
Z \times Y & & X}$$
in $\calB$. Then this diagram exhibits $Z$ as a weak exponential of $X$ by $Y$ in
$\bfB$ if and only the associated diagram
$$ \xymatrix{ & GP \ar[dl] \ar[dr] & \\
GZ \times GY & & GX}$$
exhibits $GZ$ as a weak exponential of $GX$ by $GY$ in $\bfA$.
\end{remark}

To work effectively with weak exponentials, we need to introduce an additional assumption.

\begin{definition}\label{tumba}\index{gen}{multiplication!and homotopy colimits}
Let $\bfA$ be a combinatorial model category containing a fibrant object
$Y$. We will say that {\it multiplication by $Y$ preserves homotopy colimits}
if the following condition is satisfied:
\begin{itemize}
\item[$(\ast)$] Let $A$ be a (small) partially ordered set, let $F: A \rightarrow \bfA$
be a projectively cofibrant diagram, and let $F': A \rightarrow \bfA$ be another strongly
cofibrant diagram equipped with a natural transformation $F'(a) \rightarrow F(a) \times Y$
which is weak equivalence for each $a \in A$. Then the induced map
$\colim F' \rightarrow (\colim F) \times Y$ exhibits $\colim F'$ as a homotopy product
of $Y$ with $\colim F$ in $\bfA$.
\end{itemize}
We will say that {\it multiplication in $\bfA$ preserves homotopy colimits} if
condition $(\ast)$ is satisfied for every fibrant object $Y \in \bfA$.
\end{definition}

\begin{remark}
Definition \ref{tumba} refers only to homotopy colimits indexed by partially ordered sets.
However, every diagram indexed by an arbitrary category can be replaced by a diagram
indexed by a partially ordered set having the same homotopy colimit. We formulate and prove a precise statement to this effect (in the language of $\infty$-categories) in \S \ref{c4s2}. However, we will not need (or use) any such result in this appendix.
\end{remark}

\begin{remark}\label{canus}
Let $\Adjoint{F}{\bfA}{\bfB}{G}$ be a Quillen equivalence between combinatorial
model categories, and let $Y \in \bfB$ be a fibrant object. Then multiplication by
$Y$ preserves homotopy colimits in $\bfB$ if and only if multiplication by
$G(Y)$ preserves homotopy colimits in $\bfA$. Since the right derived functor $RG$
is essentially surjective on homotopy categories, we see that
multiplication in $\bfB$ preserves homotopy colimits if and only if multiplication in $\bfA$ preserves homotopy colimits.
\end{remark}

\begin{example}\label{canuss}
Let $\bfS$ be an excellent model category with respect to the symmetric monoidal structure
given by Cartesian product in $\bfS$. Then multiplication in $\SCat$ preserves homotopy colimits.
This is precisely the content of Theorem \ref{tubba}.
\end{example}

\begin{lemma}\label{alem}
Let $S$ be a collection of simplicial sets satisfying the following conditions:
\begin{itemize}
\item[$(i)$] The simplicial set $\Delta^0$ belongs to $S$.
\item[$(ii)$] If $f: X \rightarrow Y$ is a weak homotopy equivalence, then
$X \in S$ if and only if $Y \in S$.
\item[$(iii)$] For every small partially ordered set $A$, if
$F: A \rightarrow \sSet$ is a projectively cofibrant diagram such that
each $F(a) \in S$, then $\colim F \in S$.
\end{itemize}
\end{lemma}

\begin{proof}
Using $(ii)$ and $(iii)$, we deduce that if $F: A \rightarrow \sSet$ is
{\em any} diagram such that each $F(a)$ belongs to $S$, then
the homotopy colimit of $F$ belongs to $S$. In particular, $S$ is closed under the formation of coproducts and homotopy pushouts.

We now prove by induction on $n$ that every $n$-dimensional simplicial set $X$ belongs to $S$.
For this, we observe that there is a homotopy pushout diagram
$$ \xymatrix{ B \times \bd \Delta^n \ar[d] \ar[r] & B \times \Delta^n \ar[d] \\
\sk^{n-1} X \ar[r] & X, }$$
where $B$ denotes the set of $n$-simplices of $X$. The simplicial sets
$B \times \bd \Delta^n$ and $\sk^{n-1} X$ belong to $S$ by the inductive hypothesis.
The simplicial set $B \times \Delta^n$ is weakly equivalent to the constant simplicial
set $B$, which belongs to $S$ in view of $(i)$ and the fact that $S$ is stable under coproducts.
Since $S$ is stable under homotopy pushouts, we conclude that $X \in S$ as desired.

An arbitrary simplicial set $X$ can be written as the colimit of a projectively cofibrant diagram
$$ \sk^{0} X \subseteq \sk^{1} X \subseteq \sk^{2} X \subseteq \ldots$$
and therefore belongs to $S$ by assumption $(iii)$.
\end{proof}



\begin{proposition}\label{scat}
Let $\bfA$ be a combinatorial simplicial model category
containing a standard diagram
$$ \xymatrix{ & P \ar[dl]^{p} \ar[dr] & \\
Z \times Y & & X.}$$
Assume further that multiplication by $Y$ preserves homotopy colimits in $\bfA$.
The following conditions are equivalent:
\begin{itemize}
\item[$(i)$] The above diagram exhibits $Z$ as a weak exponential of 
$X$ by $Y$.

\item[$(ii)$] Let $W$ and $W'$ be cofibrant objects of $\bfA$, and
$W' \rightarrow W \times Y$ a map which exhibits $W'$ as a homotopy product
of $W$ and $Y$. Then the induced map
$$\bHom_{\bfA}(W,Z) \times_{ \bHom_{\bfA}( W', Z \times Y) } \bHom_{\bfA}(W', P)
\rightarrow \bHom_{\bfA}(W', X)$$
is a homotopy equivalence of Kan complexes.
\end{itemize}
\end{proposition}

\begin{remark}\label{kilpot}
In the situation of part $(ii)$ of Proposition \ref{scat}, the projection map
$\bHom_{\bfA}( W', P) \rightarrow \bHom_{\bfA}( W', Z \times Y)$ is a trivial Kan fibration, so
the fiber product
$\bHom_{\bfA}(W,Z) \times_{ \bHom_{\bfA}( W', Z \times Y) } \bHom_{\bfA}(W', P)$ is automatically
a Kan complex which is homotopy equivalent to $\bHom_{\bfA}(W,Z)$.
\end{remark}

\begin{proof}[Proof of Proposition \ref{scat}]
First suppose that $(ii)$ is satisfied. We wish to show that for every object
$[W] \in \h{\bfA}$, the map
$$ \Hom_{\h{\bfA}}( [W], [Z] ) \rightarrow
\Hom_{ \h{\bfA} }( [W] \times [Y], [Z] \times [Y])
\simeq \Hom_{ \h{\bfA}}( [W] \times [Y], [P] ) \rightarrow \Hom_{\h{\bfA}}( [W] \times [Y], [P] )$$
is bijective. Without loss of generality, we may assume that $[W]$ is the homotopy equivalence class
of a fibrant-cofibrant object $W \in \bfA$. Choose a cofibrant replacement $W' \rightarrow W \times Y$.
We observe that the map in question can be identified with the composition
$$ \pi_0 \bHom_{ \bfA}( W, Z) \rightarrow \pi_0 \bHom_{ \bfA}( W', Z \times Y)
\simeq \pi_0 \bHom_{\bfA}( W', P) \rightarrow \pi_0 \bHom_{\bfA}( W',X),$$
which is bijective in view of $(ii)$ and Remark \ref{kilpot}.

We now assume $(i)$ and prove $(ii)$. It will suffice to show that for every simplicial
set $K$, the induced map
$$\Hom_{\h{\sSet}}( K, \bHom_{\bfA}(W,Z) \times_{ \bHom_{\bfA}( W', Z \times Y) } \bHom_{\bfA}(W', P)
) \rightarrow \Hom_{\h{\sSet}}(K, \bHom_{\bfA}(W',X) )$$
is a bijection. Using Remark \ref{kilpot}, we can identify the left side with the
set $\Hom_{\h{\sSet}}(K, \bHom_{\bfA}(W,Z)) \simeq \Hom_{\h{\bfA}}( W \otimes K, Z)$.
Similarly, the right side can be identified with $\Hom_{\h{\bfA}}( W' \otimes K, X)$.
In view of assumption $(i)$, it will suffice to show that the map
$\beta_{K}: W' \otimes K \rightarrow Y \times (W \otimes K)$
exhibits $W' \otimes K$ as a homotopy product of $Y$ with $W \otimes K$.
The collection of simplicial sets $K$ with this property clearly contains
$\Delta^0$ and is stable under weak homotopy equivalence. The assumption
that multiplication by $Y$ preserves homotopy colimits guarantees that the hypotheses of
Lemma \ref{alem} are satisfied, so that the desired conclusion holds for {\em every} simplicial set $K$.
\end{proof}

\begin{lemma}\label{tukka}
Let $\bfA$ be a combinatorial model category and $i: B_0 \rightarrow B$
an inclusion of partially ordered sets. Suppose that there exists a retraction
$r: B \rightarrow B_0$, such that $r(b) \leq b$ for each $b \in B$.
Let $F: B \rightarrow \bfA$ be a diagram. Then a map
$\alpha: X \rightarrow \lim(F)$ in $\bfA$ exhibits $X$ as a homotopy limit of
$F$ if and only if $\alpha$ exhibits $X$ as a homotopy limit of $i^{\ast} F$.
\end{lemma}

\begin{proof}
Without loss of generality, we may assume that $F$ is injectively fibrant.
We have a canonical isomorphism $\lim(F) \simeq \lim( i^{\ast} F )$.
It will therefore suffice to show that the functor $i^{\ast}$ preserves injective fibrations.
It now suffices to observe that $i^{\ast}$ is right adjoint to $r^{\ast}$, and that
the functor $r^{\ast}$ preserves weak trivial cofibrations.
\end{proof}

\begin{lemma}\label{constancy}
Let $\bfA$ be a combinatorial model category, $\calC$ a small category, $F: \calC \rightarrow \bfA$ a diagram, and $\alpha: X \rightarrow \lim(F)$ a morphism in the category $\bfA$.
Suppose that:
\begin{itemize}
\item[$(i)$] For each $C \in \calC$, the induced map $X \rightarrow F(C)$ is a weak equivalence in $\bfA$.
\item[$(ii)$] The category $\calC$ has a final object $C_0$. 
\end{itemize}
Then $\alpha$ exhibits $X$ as a homotopy limit of the diagram $F$.
\end{lemma}

\begin{proof}
Without loss of generality, we may assume that the diagram $F$ is projectively fibrant.
Let $F': \calC \rightarrow \bfA$ be defined by the formula $F'(C) = F(C_0)$.
We observe that, for every $G \in \Fun(\calC, \bfA)$, we have
$\Hom_{ \Fun(\calC, \bfA) }( G, F') = \Hom_{\bfA}( G(C_0), F(C_0) )$. In particular,
we have a canonical map $\beta: F \rightarrow F'$. Condition $(i)$ guarantees that
$\beta$ is a weak equivalence. Since
$F(C_0) \in \bfA$, is fibrant, the diagram $F'$ is injectively fibrant. It therefore suffices to show that
the induced map $X \rightarrow \lim(F') \simeq F(C_0)$ is a weak equivalence, which follows from
$(i)$.
\end{proof}

\begin{lemma}\label{stride}
Let $\bfA$ be a combinatorial model category,
$A$ a partially ordered set, and set $B = \{ (a,b) \in A^{op} \times A: a \geq b \}$,
regarded as a partially ordered subset of $A^{op} \times A$. Let
$\pi: B \rightarrow A^{op}$ denote the projection onto the first factor.

Suppose given diagrams $F: B \rightarrow \bfA$, $G: A \rightarrow \bfA$,
and a natural transformation $\alpha: \pi^{\ast}(G) \rightarrow F$, which induces
weak equivalences $G(a) \rightarrow F(a,b)$ for each $(a,b) \in B$.
Then $\alpha$ exhibits $G$ as a homotopy right Kan extension of $F$.
\end{lemma}

\begin{proof}
In view of Proposition \ref{sabke}, it will suffice to show that
for each $a_0 \in A$, the transformation $\alpha$ exhibits
$G(a_0)$ as a homotopy limit of the diagram
$F | \{ (a,b) \in B: a \leq a_0 \}$. 
Let $B_0 = \{ (a,b) \in B: a = a_0 \}$. In view of Lemma \ref{tukka}, it will
suffice to show that $\alpha$ exhibits $G(a_0)$ as a limit of the diagram
$F_0 = F | B_0$. This follows immediately from Lemma \ref{constancy}.
\end{proof}

\begin{proposition}\label{psygood}
Let $\bfA$ be a combinatorial model category and $A = A_0 \cup \{ \infty \}$
a partially ordered set with a largest element $\infty$. Let 
$B = \{ (a,b) \in A^{op} \times A: a \geq b \}$, regarded as a partially ordered subset
of $A^{op} \times A$. 

Suppose given an object $X \in \bfA$ together with functors
$Y: A \rightarrow \bfA$, $Z: A^{op} \rightarrow \bfA$,
$P: B \rightarrow \bfA$ and diagrams $\sigma_{a,b}$: 
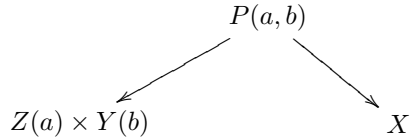
$$ \xymatrix{ & P(a,b) \ar[dl] \ar[dr] & \\
Z(a) \times Y(b) & & X}$$
which depend functorially on $(a,b) \in B$.
Suppose further that:
\begin{itemize}
\item[$(i)$] Each diagram $\sigma_{a,b}$ exhibits $P(a,b)$ as a homotopy product
of $Z(a)$ and $Y(b)$ in $\bfA$.

\item[$(ii)$] The diagrams $\sigma_{a,a}$ exhibit $Z(a)$ as a weak exponential
of $X$ by $Y(a)$.

\item[$(iii)$] Multiplication in $\bfA$ preserves homotopy colimits.

\item[$(iv)$] The diagram $Y$ exhibits $Y(\infty)$ as a homotopy colimit of
$Y_0 = Y|A_0$.
\end{itemize}
Then the diagram $Z$ exhibits $Z(\infty)$ as the homotopy limit of
the diagram $Z_0 = Z | A_0^{op}$. 
\end{proposition}

\begin{proof}
Making fibrant replacements if necessary, we may assume that each diagram
$\sigma_{a,b}$ is standard. 
According to the main result of \cite{combmodel}, there exists a Quillen equivalence
$\Adjoint{F}{\bfA'}{\bfA}{G}$ where $\bfA'$ is a combinatorial {\em simplicial} model category.
In view of Remark \ref{toofus}, we may replace $\bfA$ by $\bfA'$
and thereby reduce to the case where $\bfA$ is a simplicial model category. 

In view of Proposition \ref{usecoinc}, it will suffice to prove the following: for every
fibrant-cofibrant object $C \in \bfA$, if we define $G: A^{op} \rightarrow \sSet$
by the formula $G(a) = \bHom_{\bfA}( C, Z(a) )$, then
$G$ exhibits $G(\infty)$ as a homotopy limit of the diagram $G|A_0^{op}$.

Let $W: A \rightarrow \bfA$ be a cofibrant replacement for the functor
$a \mapsto C \times Y(a)$. Let $G': A^{op} \rightarrow \sSet$ be defined
by the formula $G'(a) = \bHom_{\bfA}( W(a), X)$. 

Define $G'': B \rightarrow \sSet$ by the formula
$$ G''(a,b) = \bHom_{\bfA}( C, Z(a) ) \times_{ \bHom_{\bfA}( W(a), Z(a) \times Y(b) ) }
\bHom_{\bfA}( W(a), P(a,b) ).$$
Let $\pi: B \rightarrow A^{op}$ denote projection onto the first factor, so that
we have natural transformations of diagrams
$$ \pi^{\ast} G \stackrel{\alpha}{\leftarrow} G'' \stackrel{\beta}{\rightarrow} \pi^{\ast} G'.$$
We observe that $\beta$ induces a trivial Kan fibration $G''(a,b) \rightarrow G'(a)$
for all $(a,b) \in B$. In particular, for $a \leq b$ the induced map
$G''(a,a) \rightarrow G''(a,b)$ is a homotopy equivalence.
Condition $(ii)$ guarantees that $\alpha$ induces a
homotopy equivalence $G''(a,b) \rightarrow G(a)$ if $a = b$, and therefore for all
$(a,b) \in B$.

Using Lemma \ref{stride}, we conclude that $\alpha$ and $\beta$ exhibit
$G$ and $G'$ as homotopy right Kan extensions of $G''$ along $\pi$. 
In particular, $G$ and $G'$ are equivalent in the homotopy category
$\h{ \Fun( A^{op}, \bfA)}$. Assumptions $(iii)$ and
$(iv)$ guarantee that $W$ exhibits $W(\infty)$ as the homotopy colimit
of $W | A_0$. Using Proposition \ref{usecoinc}, we deduce that
$G'$ exhibits $G'(\infty)$ as the homotopy limit of $G' | A_0^{op}$. 
It follows that $G$ exhibits $G(\infty)$ as the homotopy limit of
$G| A_0^{op}$, as desired.
\end{proof}

We conclude this section with an application of Proposition \ref{psygood}. 

\begin{proposition}\label{sturb}
Let $\bfS$ be an excellent model category in which the monoidal structure
is given by the Cartesian product. Let $\bfA$ be a combinatorial $\bfS$-enriched model category,
$A = A_0 \cup \{ \infty \}$ a partially ordered set with a largest element $\infty$, and $\{ \calC_{a} \}_{a \in A}$ a diagram of small $\bfS$-enriched categories indexed by $A$. Let
$\calU \subseteq \bfA$ be a chunk. For each $a \in A$, let $\calU^{\calC_{a}}_{f}$ denote the full
subcategory of $\calU^{\calC_a} \subseteq \bfA^{\calC_a}$ spanned by the projectively fibrant diagrams,
and let $W_{a}$ denote the collectinon of weak equivalences in $\calV_{a}$.
Assume that:
\begin{itemize}
\item[$(a)$] For each $a \in A$, the $\bfS$-enriched category $\calC_{a}$ is
cofibrant, and $\calU$ is a $\calC_{a}$-chunk of $\bfA$.

\item[$(b)$] The diagram $\{ \calC_{a} \}_{ a \in A}$ exhibits $\calC_{\infty}$ as a
homotopy colimit of the diagram $\{ \calC_{a} \}_{a \in A_0}$. 

\item[$(c)$] The chunk $\calU$ is small.
\end{itemize}

Then the induced diagram $\{ \calU^{\calC_{a}}_{f}[W_{a}^{-1}] \}_{a \in A}$ exhibits
$ \calU^{\calC_{\infty}}_{f}[ W_{\infty}^{-1} ]$ as a homotopy limit of the diagram
$ \{ \calU^{\calC_{a}}_{f}[ W_a^{-1} ] \}_{a \in A_0}$. 
\end{proposition}

Before proving Proposition \ref{sturb}, we need a simple lemma.

\begin{lemma}\label{kur}
Let $\bfS$ be an excellent model category, $\bfA$ a combinatorial $\bfS$-enriched
model category, and $\calU \subseteq \bfA$ a chunk. Let $\calU_{f}$ denote the
full subcategory of $\calU$ spanned by those objects which are fibrant in $\bfA$, and let
$W$ denote the collection of weak equivalences in $\calU_{f}$. Then the induced map
$\calU^{\degree} \rightarrow \calU_{f}[W^{-1}]$ is a weak equivalence of $\bfS$-enriched categories.
\end{lemma}

\begin{proof}
Let $W_0 = W \cap \calU^{\degree}$. Since every weak equivalence in $\calU^{\degree}$
is actually an equivalence, we conclude that the induced map
$\calU^{\degree} \rightarrow \calU^{\degree}[W_0^{-1}]$ is a weak equivalence.
It will therefore suffice to prove that the map
$i: \calU^{\degree}[W_0^{-1}] \rightarrow \calU_{f}[W^{-1}]$ is a weak equivalence.
Let $F$ be a $\bfS$-enriched fibrant replacement functor which carries $\calU$ to itself, so
that $F$ induces a map $j: \calU_{f}[W^{-1}]$ to $\calU^{\degree}[W_0^{-1}]$. We claim that $j$ is a homotopy inverse to $i$. To prove this, we observe that there is a natural transformation
$\alpha: \id \rightarrow F$, which we can identify with a map
$$\overline{h}: \calU_{f} \otimes [1]_{\bfS} \rightarrow \calU^{\degree}.$$
Let $W'_0$ be the collection of all morphisms in $\calU_{f} \otimes [1]_{\bfS}$
of the form $e \otimes \id$, where $e$ is an equivalence in $\calU_{f}$, and let
$W'_1$ be the collection of all morphisms of $\calU_f \otimes [1]_{\bfS}$ of the form
$\id \otimes g$, where $g: 0 \rightarrow 1$ is the tautological morphism in
$[1]_{\bfS}$. Let $W' = W'_0 \cup W'_1$, so that $\overline{h}$ determines a map
$$ h: ( \calU_{f} \otimes [1]_{\bfS} )[{W'}^{-1} ] \rightarrow \calU^{\degree}[W_0^{-1}].$$
We will prove that $h$ determines a homotopy from the identity to $j \circ i$, so that
$j$ is a left homotopy inverse to $i$. Applying the same argument to the
restriction $\overline{h} | \calU^{\degree} \otimes [1]_{\bfS}$ will show that
$j$ is a right homotopy inverse to $i$. 

To prove that $h$ gives the desired homotopy, it will suffice to show that the inclusions
$\{0\}, \{1\} \hookrightarrow [1]_{\bfS}$ induce weak equivalences
$$ \calU_{f}[W^{-1}] \rightarrow (\calU_{f} \otimes [1]_{\bfS})[{W'}^{-1}].$$
This follows immediately from Corollary \ref{suff} and Lemma \ref{twoface}. 
\end{proof}

\begin{proof}[Proof of Proposition \ref{sturb}]
Let $B = \{ (a,b) \in A^{op} \times A: a \geq b \}$. For
each $(a,b) \in B$, we define $P(a,b) = (\calU^{\calC_{a}}_{f} \times \calC_{b})[V_{a,b}^{-1}]$,
where $V_{a,b}$ is the collection of all morphisms of $\calU^{\calC_{a}}_{f} \times \calC_{b}$
of the form $e \otimes \id_{C}$, where $e \in W_a$ and $C \in \calC_{b}$. 
We have an evident family of diagrams $\sigma(a,b):$
$$ \xymatrix{ & P(a,b) \ar[dl] \ar[dr] & \\
\calU^{\calC_{a}}_{f}[W_{a}^{-1}] \times \calC_{b} & & \calU_{f}[W], }$$
where $\calU_{f}$ denotes the full subcategory of $\calU$ spanned by the fibrant objects,
and $W$ is the collection of weak equivalences in $\calU_{f} \subseteq \bfA$. 
To complete the proof, it will suffice to show that the hypotheses of Proposition
\ref{psygood} are satisfied. Condition $(i)$ follows from Lemma \ref{twoface},
condition $(iii)$ from Theorem \ref{tubba}, and condition $(iv)$ follows from $(b)$.
To prove $(ii)$, we observe that for each $a \in A$, the diagram
$\sigma(a,a)$ is weakly equivalent to the diagram
$$ \xymatrix{ & (\calU^{\calC_{a}})^{\degree} \times \calC_{a} \ar[dl]^{\sim} \ar[dr] & \\
( \calU^{\calC_{a}})^{\degree} \times \calC_{a} & & \calU^{\degree}. }$$
This diagram exhibits $(\calU^{\calC_{a}})^{\degree}$ as a weak exponential of
$\calU^{\degree}$ by $\calC_{a}$ by Corollary \ref{sniffle}. 
\end{proof}

\begin{corollary}\label{uspin}
Let $\bfS$ be an excellent model category in which the monoidal structure
is given by the Cartesian product. Let $\bfA$ be a combinatorial $\bfS$-enriched model category,
$A = A_0 \cup \{ \infty \}$ a partially ordered set with a largest element $\infty$, and $\{ \calC_{a} \}_{a \in A}$ a diagram of small $\bfS$-enriched categories indexed by $A$.

For each $a \in A$, let $\bfA^{\calC_{a}}_{f}$ denote the collection of projectively fibrant objects of
$\bfA^{\calC_a}$, and let $W_{a}$ denote the collection of weak equivalences in
$\bfA^{\calC_a}_{f}$. Assume that the diagram $\{ \calC_{a} \}_{a \in A}$ exhibits
$\calC_{\infty}$ as a homotopy colimit of the diagram $\{ \calC_{a} \}_{a \in A_0}$. Then
the induced diagram $\{ \bfA^{\calC_a}_{f}[W_a^{-1}] \}_{a \in A}$ exhibits
$\bfA^{\calC_{\infty}}[W_{\infty}^{-1}]$ as a homotopy limit of the diagram
$\{ \bfA^{\calC_a}_{f}[W_a^{-1}] \}_{a \in A_0}$.
\end{corollary}

\begin{proof}
Without loss of generality we may suppose that each $\calC_{a}$ is cofibrant.
The proof of Proposition \ref{smurff} shows that there exists a (small) regular cardinal $\kappa$
such that the collection of homotopy limit diagrams in $\Fun( A, \SCat)$ is stable under
$\kappa$-filtered colimits. This cardinal depends only on $A$ and $\bfS$, and remains invariant if we enlarge the universe. Using Lemma \ref{exchunk}, we can write
$\bfA$ as a $\kappa$-filtered union of full subcategories $\calU \subseteq \bfA$,
where $\calU$ is a $\calC_{a}$-chunk for each $a \in A$. We now conclude by
applying Proposition \ref{sturb}.
\end{proof}

\subsection{Localizations of Simplicial Model Categories}\label{turka}

Let $\bfA$ and $\bfA'$ be two model categories with the same underlying category. We say that $\bfA'$ is a {\it $($Bousfield$)$ localization} of $\bfA$ if:\index{gen}{localization!of a model category}\index{gen}{localization!Bousfield}
\begin{itemize}
\item[$(C)$] A morphism $f$ of $\calC$ is a cofibration in $\bfA$ if and only if $f$ is a cofibration in $\bfA'$.
\item[$(W)$] If a morphism $f$ of $\calC$ is a weak equivalence in $\bfA$, then $f$ is a weak equivalence in $\bfA'$.
\end{itemize}
It then follows also that:
\begin{itemize}
\item[$(F)$] If a morphism $f$ of $\calC$ is a fibration in $\bfA'$, then $f$ is a fibration in $\bfA$.
\end{itemize}

Our goal in this section is to study the localizations of a fixed model category $\bfA$, and to relate this to our study of localizations of presentable $\infty$-categories (\S \ref{invloc}). 

Let $\bfA$ be a simplicial model category. Let $\h{\bfA}$ be the homotopy category of $\bfA$, obtained from $\bfA$ by inverting all weak equivalences. Alternatively, we can obtain $\bfA$ by first passing to the full subcategory $\bfA^{\degree} \subseteq \bfA$ spanned by the fibrant-cofibrant objects, and then passing to the homotopy category of the simplicial category $\bfA^{\degree}$. From the second point of view, we see that $\h{\bfA}$ has a natural enrichment over the homotopy category $\calH$: if $X, Y \in \h{\bfA}$ are represented by fibrant-cofibrant objects
$\overline{X}, \overline{Y} \in \bfA$, then we let
$$ \bHom_{ \h{\bfA}}(X,Y) = [ \bHom_{\bfA}(\overline{X}, \overline{Y}) ].$$
Here $[K] \in \calH$ denotes the object of $\calH$ represented by a Kan complex $K$. In fact, this description is accurate if we assume only that $\overline{X}$ is cofibrant and $\overline{Y}$ fibrant.

Let $S$ be a collection of morphisms in $\h{\bfA}$. Then:
\begin{itemize}
\item[$(i)$] We will say that an object $Z \in \h{\bfA}$ is 
{\it $S$-local} if, for every morphism $f: X \rightarrow Y$ in $S$, the induced map
$$ \bHom_{\h{\bfA}}( Y, Z) \rightarrow \bHom_{\h{\bfA}}( X, Z)$$ is a homotopy equivalence.\index{gen}{$S$-local!object} We say that an object $\overline{Z} \in \bfA$ is $S$-local if its image in
$\h{\bfA}$ is $S$-local.
\item[$(ii)$] We will say that a morphism $f: X \rightarrow Y$ of $\h{\bfA}$ is an {\it $S$-equivalence} if, for every $S$-local object $Z \in \h{\bfA}$, the induced map
$$ \bHom_{\h{\bfA}}( Y, Z) \rightarrow \bHom_{\h{\bfA}}( X, Z)$$ is a homotopy equivalence. We say that a morphism $\overline{f}$ in $\bfA$ is an {\it $S$-equivalence} if its image in $\h{\bfA}$ is an $S$-equivalence.\index{gen}{$S$-equivalence}
\end{itemize}
If $\overline{S}$ is a collection of morphisms in $\h{\bfA}$, with image $S$ in $\h{\bfA}$, we will apply the same terminology: an object of $\bfA$ (or $\h{\bfA}$) is said to be {\it $\overline{S}$-local} if it is $S$-local, and morphism of $\bfA$ (or $\h{\bfA}$) will be said to be an {\it $\overline{S}$-equivalence}, if it is an $S$-equivalence.

\begin{lemma}\label{mtomp}
Let $\bfA$ be a left proper simplicial model category, let $S$ be a collection of morphisms
in $\h{\bfA}$, and let $i: A \rightarrow B$ be a cofibration in $\bfA$. The following conditions are equivalent:
\begin{itemize}
\item[$(1)$] The map $i$ is an $S$-equivalence.
\item[$(2)$] For every fibrant object $X \in \bfA$ which is $S$-local, the map
$\bHom_{\bfA}(B, X) \rightarrow \bHom_{\bfA}(A,X)$ is a trivial Kan fibration.
\end{itemize}
\end{lemma}

\begin{proof}
Choose a trivial fibration
$f: A' \rightarrow A$, where $A'$ is cofibrant, and choose a factorization
$$ A' \stackrel{i'}{\rightarrow} B' \stackrel{f'}{\rightarrow} B$$ of
$i \circ f$, where $i'$ is a cofibration and $f'$ is a trivial fibration. We have a commutative diagram
$$ \xymatrix{ A' \ar[r]^{i'} \ar[d]^{f} & B' \ar[d]^{g} \ar[dr]^{f'} & \\
A \ar[r]^{i} & A \coprod_{A'} B' \ar[r]^{j} & B. }$$
Since $f$ is a weak equivalence and $i'$ is a cofibration, the left properness of $\bfA$ guarantees that $g$ is a weak equivalence. It follows from the two-out-of-three property that $j$ is also a weak equivalence.

Suppose first that $(1)$ is satisfied. Let $X$ be an $S$-local fibrant object of $\bfA$. 
The map $p: \bHom_{\bfA}(B,X) \rightarrow \bHom_{\bfA}(A,X)$ is a Kan fibration. We wish to show that $p$ is a trivial Kan fibration. Our assumption that $X$ is $S$-local guarantees that the map
$q': \bHom_{\bfA}(B',X) \rightarrow \bHom_{\bfA}(A', X)$ is a homotopy equivalence, and therefore a trivial fibration (since $i'$ is a cofibration). The map
$$q: \bHom_{\bfA}( A \coprod_{A'} B', X) \rightarrow \bHom_{\bfA}(A)$$ is a pullback of $q'$, and therefore also a trivial fibration. To show that $p$ is a trivial Kan fibration, it will suffice to show that
for every $t: A \rightarrow X$, the fiber $p^{-1} \{ t\}$ is a contractible Kan complex. Since the corresponding fiber $q^{-1} \{t\}$ is contractible, it will suffice to show that composition with $j$ induces a homotopy equivalence $$ p^{-1} \{t\} \rightarrow q^{-1} \{t\}. $$
This is clear, since $j$ is a weak equivalence between cofibrant objects of the simplicial model category $\bfA_{A/}$.

Now suppose that $(2)$ is satisfied. We wish to show that $i$ is an $S$-equivalence. For this, it suffices to show that for every fibrant $S$-local object $X \in \bfA$, the map
$$ q': \bHom_{\bfA}( B', X) \rightarrow \bHom_{\bfA}(A',X)$$ is a trivial Kan fibration.
The preceding argument shows that the fiber of $q'$ over a morphism $t': A' \rightarrow X$ is contractible, provided that $t'$ factors as a composition
$$ A' \stackrel{f}{\rightarrow} A \stackrel{t}{\rightarrow} X.$$
To complete the proof, it suffices to show that the same result holds for an arbitrary vertex
$t'$ of $\bHom_{\bfA}(A', X)$. The map $t'$ factors as a composition
$$ A' \stackrel{u}{\rightarrow} Y \stackrel{v}{\rightarrow} X,$$
where $u$ is a cofibration and $v$ is a trivial fibration. We have a commutative diagram
$$ \xymatrix{ \bHom_{\bfA}(B', Y) \ar[r] \ar[d] & \bHom_{\bfA}(A', Y) \ar[d] \\
\bHom_{\bfA}( B', X) \ar[r] & \bHom_{\bfA}(A', X) }$$
in which the vertical arrows are trivial Kan fibrations. It will therefore suffice to show that the fiber 
$\bHom_{\bfA}(B', Y) \times_{ \bHom_{\bfA}(A',Y) } \{ u \}$ is a contractible.
Choose a trivial cofibration $A' \coprod_{A} Y \rightarrow Z$, where $Z$ is fibrant. We observe that
The map $Y \rightarrow A' \coprod_{A} Y$ is the pushout of a weak equivalence by a cofibration, and therefore a weak equivalence (since $\bfA$ is left proper). It follows that the map
$Y \rightarrow Z$ is a weak equivalence between fibrant objects of $\bfA$. We have a commutative diagram 
$$ \xymatrix{ \bHom_{\bfA}(B', Y) \ar[r] \ar[d] & \bHom_{\bfA}(A', Y) \ar[d] \\
\bHom_{\bfA}( B', Z) \ar[r]^{q''} & \bHom_{\bfA}(A', Z) }$$
in which the vertical maps are homotopy equivalences, and the horizontal maps are Kan fibrations. It will therefore suffice to show that the fiber of $q''$ is contractible, when taken over the composite map $t'': A' \stackrel{u}{\rightarrow} Y \rightarrow Z$. We now observe that $t''$ factors through
$A$, so that the desired result follows from the first part of the proof.
\end{proof}

\begin{corollary}\label{swask}
Let $\bfA$ and $\bfB$ be simplicial model categories, and suppose given a simplicial adjunction
$$ \Adjoint{F}{\bfA}{\bfB.}{G}$$
Assume that $\bfB$ is left proper. The following conditions are equivalent:
\begin{itemize}
\item[$(1)$] The adjunction between $F$ and $G$ is a Quillen adjunction.
\item[$(2)$] The functor $F$ preserves cofibrations and the functor $G$ preserves fibrant objects.
\end{itemize}
\end{corollary}

\begin{proof}
The implication $(1) \Rightarrow (2)$ is obvious. Conversely, suppose that $(2)$ is satisfied.
We wish to prove that $F$ is a left Quillen functor. Since $F$ preserves cofibrations, it will suffice
to show that for every trivial cofibration $u: A \rightarrow A'$ in $\bfA$, the image
$Fu$ is a weak equivalence in $\bfB$. Applying Lemma \ref{mtomp} in the case
$S = \emptyset$, it will suffice to prove the following: for every fibrant object $B \in \bfB$, the
induced map
$$ \bHom_{ \bfB}( FA', B) \rightarrow \bHom_{\bfB}( FA, B)$$
is a trivial Kan fibration. Since $F$ and $G$ are adjoint simplicial functors, this is equivalent
to the requirement that the map $\bHom_{\bfA}( A', GB) \rightarrow \bHom_{\bfA}(A, GB)$ be
a trivial Kan fibration, which follows from our assumption that $u$ is a trivial cofibration in
$\bfA$ and that $GB \in \bfA$ is fibrant.
\end{proof}

\begin{proposition}\label{suritu}
Let $\bfA$ be a left proper combinatorial simplicial model category, and let $S$ be a $($small$)$ set of cofibrations in $\bfA$. 
Let $S^{-1} \bfA$ denote the same category, with the following distinguished classes of morphisms:
\begin{itemize}
\item[$(C)$] A morphism $g$ in $S^{-1} \bfA$ is a {\it cofibration} if it is a cofibration when regarded as a morphism in $\bfA$.
\item[$(W)$] A morphism $g$ in $S^{-1} \bfA$ is a {\it weak equivalence} if it is an $S$-equivalence.
\end{itemize}
Then:
\begin{itemize}
\item[$(1)$] The above definitions endow $S^{-1} \bfA$ with the structure of a combinatorial simplicial model category.
\item[$(2)$] The model category $S^{-1} \bfA$ is left proper.
\item[$(3)$] An object $X \in \bfA$ is fibrant in $S^{-1} \bfA$ if and only if $X$ is $S$-local and fibrant in $\bfA$.
\end{itemize}
\end{proposition}

\begin{proof}
Enlarging $S$ if necessary, we may assume:
\begin{itemize}
\item[$(a)$] For every morphism $f: A \rightarrow B$ in $S$ and every $n \geq 0$, the induced map
$$(A \times \Delta^n) \coprod_{ A \times \bd \Delta^n } (B \times \bd \Delta^n)
\rightarrow B \times \Delta^n$$ belongs to $S$.  
\item[$(b)$] The set $S$ contains a collection of generating trivial cofibrations for $\bfA$.
\end{itemize}
It follows that an object $X \in \bfA$ is fibrant and $S$-local if and only if it has the extension property with respect to every morphism in $S$. Since $S \subseteq C \cap W$, we deduce that every fibrant object of $S^{-1} \bfA$ is $S$-local and fibrant in $\bfA$. The converse follows from Lemma \ref{mtomp}; this proves $(3)$. 

To prove $(1)$, it will suffice to show that the classes $C$ and $W$ satisfy the hypotheses of Proposition \ref{bigmaker} (the compatibility of the simplicial structure on $S^{-1} \bfA$ with its model structure
follows immediately from Proposition \ref{testsimpmodel}). We observe that Lemma \ref{mtomp} implies that $C \cap W$ is a weakly saturated class of morphisms in $\bfA$. The only other nontrivial point is to show that $W$ is an accessible subcategory of $\bfA^{[1]}$.

Proposition \ref{quillobj} implies the existence of a functor $T: \bfA \rightarrow \bfA$, equipped with a natural transformation $\id_{\bfA} \rightarrow T$, with the following properties:
\begin{itemize}
\item[$(i)$] For every $X \in \bfA$, the object $TX \in \bfA$ is fibrant and $S$-local.
\item[$(ii)$] For every $X \in \bfA$, the map $X \rightarrow TX$ belongs to the smallest weakly saturated class of morphisms containing $S$; in particular, it belongs to $W \cap C$ and is therefore an $S$-equivalence.
\item[$(iii)$] There exists a regular cardinal $\kappa$ such that $T$ commutes with $\kappa$-filtered colimits.
\end{itemize}

It follows that a morphism $f: X \rightarrow Y$ in $\bfA$ is an $S$-equivalence if and only if the induced map $Tf: TX \rightarrow TY$ is an $S$-equivalence. Since $TX$ and $TY$ are $S$-local, Yoneda's lemma (in the category $\h{\bfA})$) implies that $Tf$ is an $S$-equivalence if and only if
$Tf$ is a weak equivalence in $\bfA$. It follows that $W$ is the inverse image under $T$ of the collection of weak equivalenes in $\bfA$. Corollaries \ref{sundert} and \ref{smitty} imply that $W$ is an accessible subcategory of $\bfA^{[1]}$, as desired. This completes the proof of $(1)$.

We now prove $(2)$. We need to show that the collection of $S$-equivalences in $\bfA$ is stable under pushouts by cofibrations. We observe that every morphism $f: X \rightarrow Z$ admits a factorization 
$$ X \stackrel{f'}{\rightarrow} Y \stackrel{f''}{\rightarrow} Z$$
where $f'$ is a cofibration and $f''$ is a weak equivalence in $\bfA$ (in fact, we can choose $f''$ to be a trivial fibration in $\bfA$). If $f$ is an $S$-equivalence, then $f'$ is an $S$-equivalence, so that $f' \in C \cap W$. It will therefore suffice to show that $C \cap W$ and the class of weak equivalences in $\bfA$ are stable under pushouts by cofibrations. The first follows from the assertion that $C \cap W$ is weakly saturated, and the second from the assumption that $\bfA$ is left proper.
\end{proof}

\begin{proposition}\label{surito}
Let $\bfA$ be a left proper combinatorial simplicial model category. Then:
\begin{itemize}
\item[$(1)$] Every combinatorial localization of $\bfA$ has the form $S^{-1} \bfA$, where
$S$ is some $($small$)$ set of cofibrations in $\bfA$.
\item[$(2)$] Given two $($small$)$ sets of cofibrations $S$ and $T$, the localizations
$S^{-1} \bfA$ and $T^{-1} \bfA$ coincide if and only if the class of $S$-local objects of
$\h{\bfA}$ coincides with the class of $T$-local objects of $\h{\bfA}$.
\end{itemize}
\end{proposition}

\begin{proof}
The ``if'' direction of $(2)$ is obvious, and the converse follows from the characterization of the fibrant objects of $S^{-1} \bfA$ given in Proposition \ref{suritu}. We now prove $(1)$. Let $\bfB$ be a combinatorial model category which is a localization of $\bfA$, and let $S$ be a set of generating trivial cofibrations for $\bfB$. We claim that $\bfB = S^{-1} \bfA$. The cofibrations of $S^{-1} \bfA$ and $\bfB$ coincide. Moreover, the collection of trivial cofibrations in $S^{-1} \bfA$ is a weakly saturated class of morphisms which contains $S$, and therefore contains every trivial cofibration in $\bfB$. To complete the proof, it will suffice to show that every trivial cofibration $f: X \rightarrow Y$ in $S^{-1} \bfA$ is a trivial cofibration in $\bfB$.

Choose a diagram
$$ \xymatrix{ X' \ar[r]^{f'} \ar[d] & Y' \ar[d] \\
X \ar[r]^{f} & Y }$$
where $X'$ is cofibrant, $f'$ is a cofibration, and the vertical maps are weak equivalences in $\bfA$. Then $f'$ is a trivial cofibration in $S^{-1} \bfA$, and it will suffice to show that $f'$ is a trivial cofibration in $\bfB$. For this, it will suffice to show that for every fibrant object $Z \in \bfB$, the map
$$ \bHom_{\bfB}(Y', Z) \rightarrow \bHom_{\bfB}(X',Z)$$ is a trivial fibration. In view of Lemma \ref{mtomp}, it will suffice to show that $Z$ is $S$-local and fibrant as an object of $\bfA$. The second claim is obvious, and the first follows from the fact that $S$ consists of trivial cofibrations in $\bfB$. 
\end{proof}

\begin{remark}\label{surito2}
In the situation of Proposition \ref{surito}, we may assume that for every cofibration
$f: A \rightarrow B$ in $S$, the objects $A$ and $B$ are themselves cofibrant.
To see this, choose for each cofibration $f: A \rightarrow B$ in $S$ a diagram
$$ \xymatrix{ A' \ar[r]^{g_{f}} \ar[d]^{u} & B' \ar[d]^{v} \ar[dr]^{w} & \\
A \ar[r]^{f'} & A \coprod_{A'} B' \ar[r]^{f''} & B }$$
as in the proof of Lemma \ref{mtomp}, so that $u$ and $w$ are trivial cofibrations,
$f = f'' \circ f'$, and $g_f$ is a cofibration between cofibrant objects. Then $g_{f}$ is a trivial cofibration in $S^{-1} \bfA$. We claim that the localizations $S^{-1} \bfA$ and
$T^{-1} \bfA$ coincide, where $T = \{ g_{f} \}_{f \in S}$. To prove this, it will suffice to show that
for each $f \in S$, every $g_{f}$-local fibrant object $X \in \bfA$ is also $f$-local.

Suppose that $X$ is $g_{f}$-local. We wish to prove that the map
$p: \bHom_{\bfA}(B, X) \rightarrow \bHom_{\bfA}(A,X)$ is a trivial Kan fibration.
Since $p$ is automatically a Kan fibration, it will suffice to show that the fiber
$p^{-1} \{t\}$ is contractible, for every morphism $t: A \rightarrow X$. 
Since $X$ is $g_{f}$-local, we deduce that the fiber $q^{-1} \{t\}$ is contractible, where
$q$ is the projection map
$\bHom_{\bfA}( A \coprod_{A'} B', X) \rightarrow \bHom_{\bfA}(A, X)$. It will therefore
suffice to show that $f''$ induces a homotopy equivalence of fibers
$$ \bHom_{\bfA_{A/}}(B,X) \rightarrow \bHom_{ \bfA_{A/}}( A \coprod_{A'} B', X).$$
This is clear, since $f''$ is a weak equivalence between cofibrant objects
of the simplicial model category $\bfA_{A/}$.
\end{remark}

\begin{proposition}\label{notthereyet}
Let $\calC$ be an $\infty$-category. The following conditions are equivalent:
\begin{itemize}
\item[$(1)$] The $\infty$-category $\calC$ is presentable.
\item[$(2)$] There exists a combinatorial simplicial model category $\bfA$ and an equivalence
$\calC \simeq \Nerve( \bfA^{\degree})$. 
\end{itemize}
\end{proposition}

\begin{proof}
According to Theorem \ref{pretop} and Proposition \ref{local}, $\calC$ is presentable if and only if there exists a small simplicial set $K$, a small set $S$ of morphisms in $\calP(K)$, and an equivalence $\calC \simeq S^{-1} \calP(K)$. Let $\calD$ be the simplicial category $\sCoNerve[K]^{op}$, and let $\bfB$ be the category $\Set_{\Delta}^{\calD}$ of 
simplicial functors $\calD \rightarrow \sSet$, endowed with the injective model structure. Proposition \ref{gumby444} implies that there is an equivalence $\calP(K) \simeq \Nerve( \bfB^{\degree})$. Moreover, Propositions \ref{suritu} and \ref{surito} implies that there is a bijective correspondence between accessible localizations of
$\calP(K)$ (as a presentable $\infty$-category) and combinatorial localizations of $\bfB$ (as a model category). This proves the implication $(1) \Rightarrow (2)$. Moreover, it also shows that $(2) 
\Rightarrow (1)$ in the special case where $\bfA$ is a localization of a category of simplicial presheaves. 

We now complete the proof by invoking the following result, proven in \cite{combmodel}: for every combinatorial model category $\bfA$, there exists a small category $\calD$, a set $S$ of morphisms of $\Set^{\calD^{op}}_{\Delta}$, and a Quillen equivalence of $\bfA$ with $S^{-1} \Set_{\Delta}^{\calD}$.
Moreover, the proof given in \cite{combmodel} shows that when $\bfA$ is a {\em simplicial} model category, then $F$ and $G$ can be chosen to be simplicial functors.
\end{proof}

\begin{remark}
Let $\bfA$ and $\bfB$ be combinatorial simplicial model categories. Then the underlying $\infty$-categories $\Nerve( \bfA^{\degree})$ and $\Nerve( \bfB^{\degree})$ are equivalent if and only if $\bfA$ and $\bfB$ can be joined by a chain of simplicial Quillen equivalences. The ``only if'' assertion follows from Corollary \ref{urchug}, and the ``if'' direction can be proven using the methods described in \cite{combmodel}.
\end{remark}

\begin{proposition}\label{cabber}
Let $\bfA$ be a left proper combinatorial simplicial model category, and let
$\calC = \Nerve( \bfA^{\degree} )$ denote its underlying $\infty$-category.
Suppose that $\calC^{0} \subseteq \calC$ is an accessible localization of $\calC$, and let
$L: \calC \rightarrow \calC^{0}$ denote a left adjoint to the inclusion.

Then there exists a localization $\bfA'$ of $\bfA$ satisfying the following conditions:
\begin{itemize}
\item[$(1)$] An object $X \in \bfA'$ is fibrant if and only if it is fibrant in $\bfA$, and
the associated object of the homotopy category $\h{ \bfA} \simeq \h{ \calC}$ belongs to the full subcategory $\h{ \calC^{0}}$. 
\item[$(2)$] A morphism $f: X \rightarrow Y$ in $\bfA'$ is a weak equivalence if and only if
the functor $L: \h{ \calC} \rightarrow \h{ \calC^{0} }$ carries $f$ to an isomorphism in the homotopy category $\h{ \calC^{0} }$.
\end{itemize}
\end{proposition}

\begin{proof}
According to Proposition \ref{notthereyet}, the $\infty$-category $\calC$ is presentable. 
The results of \S \ref{invloc} imply that we can write $\calC^{0} = S^{-1} \calC$, for some
small collection of morphisms $S$ in $\calC$. We then take $\widetilde{S}$ be a collection of representatives for the elements of $S$ as cofibrations between cofibrant objects of $\bfA$, and let
$\bfA'$ denote the localization $\widetilde{S}^{-1} \bfA$. 
\end{proof}

We conclude this section by establishing a universal property enjoyed by the localization of a combinatorial simplicial model category.

\begin{proposition}\label{stake}
Suppose given a simplicial Quillen adjunction
$$ \Adjoint{F}{\bfA}{\bfB}{G}$$
between left proper combinatorial simplicial model categories, and let $\bfA'$ be a Bousfield localization of $\bfA$. The following conditions are equivalent:
\begin{itemize}
\item[$(1)$] The adjoint functors $F$ and $G$ determine a Quillen adjunction between
$\bfA'$ and $\bfB$.
\item[$(2)$] Let $\alpha$ be a morphism in $\bfA$ which is a weak equivalence in
$\bfA'$. Then the left derived functor $LF: \h{\bfA} \rightarrow \h{\bfB}$ carries
$\alpha$ to an isomorphism in the homotopy category $\h{\bfB}$.
\item[$(3)$] For every fibrant object $X \in \bfB$, the image $GX$ is a fibrant object of $\bfA'$.
\end{itemize}
\end{proposition}

\begin{proof}
The implication $(1) \Rightarrow (2)$ is obvious, and the implication $(3) \Rightarrow (1)$ follows
from Corollary \ref{swask}. We will complete the proof by showing that $(2) \Rightarrow (3)$.
According to Proposition \ref{surito} and Remark \ref{surito2}, we may suppose that $\bfA' = S^{-1} \bfA$, where $S$ is a small collection of cofibrations between cofibrant objects of $\bfA$. Let $X$ be a fibrant object of $\bfB$; we wish to show that $GX$ is a fibrant object of $\bfA'$. Since $GX$ is fibrant in
$\bfA$, it will suffice to show that $GX$ is $S$-local (Proposition \ref{suritu}). In other words, we
must show that if $\alpha: A \rightarrow B$ belongs to $S$, then the induced map
$p: \bHom_{\bfA}( B, GX) \rightarrow \bHom_{\bfA}( A, GX)$ is a weak homotopy equivalence.
Since $F$ and $G$ are simplicial functors, we can identify $p$ with the map
$\bHom_{\bfB}( FB, X) \rightarrow \bHom_{\bfB}( FA, X)$. To prove that $p$ is a weak homotopy equivalence, it will suffice to show that $F(\alpha)$ is a weak equivalence between cofibrant objects
of $\bfB$. This follows immediately from assumption $(2)$ (because $\alpha$ is a cofibration between cofibrant objects of $\bfA$, we can identify $F(\alpha)$ with the left derived functor $LF(\alpha)$.
\end{proof}

\begin{corollary}\label{swinker}
Let $\bfA$ and $\bfB$ be left proper combinatorial simplicial model categories, and suppose given
a simplicial Quillen adjunction
$$ \Adjoint{F}{\bfA}{\bfB.}{G}$$
Then:
\begin{itemize}
\item[$(1)$] There exists a new left proper combinatorial simplicial model structure $\bfA'$ on the category $\bfA$ with the following properties:
\begin{itemize}
\item[$(C)$] A morphism $\alpha$ in $\bfA'$ is a cofibration if and only if it is a cofibration in $\bfA$.
\item[$(W)$] A morphism $\alpha$ in $\bfA'$ is a weak equivalence if and only if the left
derived functor $LF$ carries $\alpha$ to an isomorphism in the homotopy category $\h{\bfB}$.
\item[$(F)$] A morphism $\alpha$ in $\bfA'$ is a fibration if and only if it has the right lifting
property with respect to every morphism in $\bfA'$ satisfying $(C)$ and $(W)$.
\end{itemize}
\item[$(2)$] The functors $F$ and $G$ determine a new simplicial Quillen adjunction
$$ \Adjoint{F'}{\bfA'}{\bfB.}{G'}$$
\item[$(3)$] Suppose that the right derived functor $RG$ is fully faithful. 
Then the adjoint pair $(F',G')$ is a Quillen equivalence.
\end{itemize}
\end{corollary}

\begin{proof}
The functors $F$ and $G$ determine a pair of adjoint functors 
$$ \Adjoint{f}{\Nerve \bfA^{\degree}}{\Nerve \bfB^{\degree}}{g}$$
between the underlying $\infty$-categories (see Proposition \ref{quiladj}), which
are themselves presentable (Proposition \ref{notthereyet}). Let $\overline{S}$
be the collection of all morphisms $u$ in $\Nerve \bfA^{\degree}$ such that
$f(u)$ is an equivalence in $\Nerve \bfB^{\degree}$. Proposition \ref{postbluse}
implies that $\overline{S}$ is generated (as a strongly saturated class of morphisms)
by a small subset $S \subseteq \overline{S}$. Without loss of generality, we may suppose
that the morphisms of $S$ are represented by some (small) collection $T$ of cofibrations between cofibrant objects of $\bfA$. Let $\bfA' = T^{-1} \bfA$. We claim that $\bfA'$ satisfies
the description given in $(1)$. In other words, we claim that a morphism
$\alpha$ in $\bfA$ is a $T$-equivalence if and only if the left derived functor
$LF$ carries $\alpha$ to an isomorphism in $\h{\bfB}$. Without loss of generality, we may suppose
that $\alpha$ is a morphism between fibrant-cofibrant objects of $\bfA$, so that we can
view $\alpha$ as a morphism in the $\infty$-category $\Nerve \bfA^{\degree}$. In this case, both conditions on $\alpha$ are equivalent to the requirement that $\alpha$ belongs to $\overline{S}$.
This completes the proof of $(1)$. Assertion $(2)$ follows immediately from Proposition \ref{stake}.

We now prove $(3)$. Note that the homotopy category $\h{\bfA'}$ can be identified with a full subcategory of the homotopy category $\h{\bfA}$, and that under this identification the left derived functor $LF$ restricts to the left derived functor $LF'$. It follows that for every fibrant object
$X \in \bfB$, the counit map
$$ (LF')(RG')X \simeq (LF')(GX)
\simeq (LF)(GX) \simeq (LF)(RG)X \simeq X$$
is an isomorphism in $\h{\bfB}$ (where the last equivalence follows from our assumption that $RG$
is fully faithful). It follows that the functor $RG'$ is fully faithful. To complete the proof, it will suffice to show that the left derived functor $LF'$ is conservative. In other words, we must show that
if $\alpha: X \rightarrow Y$ is a morphism in $\bfA'$, then $\alpha$ is a weak equivalence if and only if
$LF(\alpha)$ is an isomorphism in $\bfB$; this follows immediately from $(1)$.
\end{proof}

\printindex{gen}{Index}
\printindex{not}{Index of Notations}

\end{document}